\documentclass[11pt,reqno]{amsart}

\usepackage[margin=1.15in]{geometry}
\usepackage{amsmath,amssymb,amsfonts,mathtools,array}
\usepackage{mathrsfs}
\usepackage[expansion=false]{microtype}
\usepackage{xcolor}
\usepackage{tikz}
\usetikzlibrary{arrows.meta,calc,positioning}
\usepackage[
 colorlinks=true,
 linkcolor=blue!55!black,
 citecolor=blue!55!black,
 urlcolor=blue!60!black,
 linktoc=all,
 hypertexnames=false,
 bookmarksopen=true,
 bookmarksnumbered=true,
 pdftitle={Open full-metric formation and local marked first-order asymptotic moduli of FIK blowdown singularities},
 pdfauthor={Henry Shin}
]{hyperref}

\newcommand{\doi}[1]{%
  \href{https://doi.org/#1}{\nolinkurl{doi:#1}}}
\newcommand{\arxiv}[1]{%
  \href{https://arxiv.org/abs/#1}{\nolinkurl{arXiv:#1}}}

\theoremstyle{plain}
\newtheorem{theorem}{Theorem}[section]
\newtheorem{proposition}[theorem]{Proposition}
\newtheorem{lemma}[theorem]{Lemma}
\newtheorem{corollary}[theorem]{Corollary}
\newtheorem{maintheorem}{Theorem}

\theoremstyle{definition}
\newtheorem{definition}[theorem]{Definition}
\theoremstyle{remark}
\newtheorem{remark}[theorem]{Remark}

\newcommand{\Ric}{\operatorname{Ric}}
\newcommand{\Rm}{\operatorname{Rm}}

\newcommand{\supp}{\operatorname{supp}}

\newcommand{\Lie}{\mathcal L}
\newcommand{\A}{\mathcal A}
\newcommand{\Z}{\mathcal Z}
\newcommand{\Q}{\mathcal Q}
\newcommand{\E}{\mathcal E}
\newcommand{\B}{\mathcal B}
\newcommand{\R}{\mathbb R}

\newcommand{\ip}[2]{\left\langle #1,#2\right\rangle_f}
\newcommand{\norm}[1]{\left\lVert #1\right\rVert}
\newcommand{\abs}[1]{\left|#1\right|}

\title[Open full-metric FIK blowdown formation]
{Open full-metric formation and local marked first-order asymptotic moduli
of FIK blowdown singularities}

\author{Henry Shin}
\email{\href{mailto:hkshin@gmail.com}{hkshin@gmail.com}}
\date{August 2026}

\begin{document}
\raggedbottom

\begin{abstract}
We prove that the Feldman--Ilmanen--Knopf blowdown singularity has a
nonempty open nonlinear formation basin in the full space of smooth
Riemannian metrics and construct, on a smaller physical neighborhood, its
local marked first-order asymptotic moduli.  More precisely, given a
closed connected oriented Riemannian four-manifold, an implantation
point, and a positive time bound, its oriented blow-up admits a
relatively \(C^{2,\alpha}\)-open set of smooth metrics whose Ricci flows
form, before that time, a localized FIK singularity with global Type-I
curvature control.
No symmetry, K\"ahler condition, or finite-dimensional tuning is
imposed.  Curvature remains uniformly bounded outside the implantation
region, while fixed-convention marked parabolic rescalings converge
along the full singular-time sequence to the ancient FIK flow.  In the
same marking, the transported exceptional sphere collapses with the
sharp FIK asymptotics for area, intrinsic diameter, and curvature.  To our knowledge,
this is the first open full-metric formation theorem for a finite-time
singularity modeled on a noncompact, noncylindrical shrinker.

A self-contained spectral certification shows that the nonnegative
space of the weighted FIK operator is nine-dimensional and consists
entirely of scaling and diffeomorphism directions.  Exact modulation
removes this geometric block.  On a smaller little-H\"older
\(h^{2,\alpha}\) neighborhood whose
smooth locus lies inside the formation basin, one fixed positive-time
restart yields the marked first-profile coordinate
\[
 \mathfrak A_1=\lambda_\infty^{-\gamma_1}V_\infty\in E_1.
\]
Here \(E_1\) is the first stable eigenspace of the weighted operator.
The amplitude is a split \(C^1\) submersion and, in local \(C^1\) product
coordinates, is the projection onto \(E_1\).  On the smooth locus,
equality of amplitudes is precisely marked first-order asymptotic
agreement; a compactly
supported transverse disk realizes every sufficiently small amplitude
uniquely within that disk; and the normalized difference of two
rescaled flows converges to the Jacobi field of their amplitude
difference.  The same amplitude determines the first quadratic scale
and phase response.  All profile and comparison statements use one
fixed transported preparation, marking, and gauge convention.  Thus
FIK blowdown is both an open full-metric singularity mechanism and,
locally in that convention, a singularity type with a complete marked
first-order profile coordinate.
\end{abstract}

\subjclass[2020]{Primary 53E20; Secondary 53C44, 35K55}
\keywords{Ricci flow, FIK shrinker, Type-I singularity, open formation
basin, marked asymptotic moduli, first-profile scattering,
first-profile foliation, geometric modulation, complex blowdown}

\maketitle

\enlargethispage{4pt}
\tableofcontents

\section{Introduction}

\subsection{Open formation and local marked first-order asymptotic moduli}

Theorem~A establishes nonlinear stability of the FIK blowdown
singularity in the full space of smooth Riemannian metrics.  On the
oriented blow-up of an arbitrary closed connected oriented four-manifold at an
arbitrary prescribed point, it constructs a relatively
\(C^{2,\alpha}\)-open neighborhood of smooth metrics every one of whose
Ricci flows develops, within an arbitrarily prescribed positive time,
the same localized FIK singularity with global Type-I curvature
control.  With the marking
convention fixed, the construction supplies frozen base-time markings
in which the parabolic rescalings converge to the FIK model along the
full singular-time sequence, without passing to a subsequence.
Curvature remains uniformly bounded outside the chosen implantation
region.  No symmetry, K\"ahler condition, or
finite-dimensional tuning remains in the hypotheses.  To our
knowledge, this is the first such open basin on a closed manifold for a
finite-time singularity modeled by a noncompact, noncylindrical
shrinker.

This advances beyond both M\'aximo's open formation family, which
remains within \(U(2)\)-invariant K\"ahler metrics, and the selected
realization theorems of Stolarski and Hughes, which produce individual
trajectories rather than a full-metric basin~\cite{Maximo,Stolarski,Hughes}.
The detailed comparison with earlier formation and convergence results
appears below.

Theorem~B is the refined terminal synthesis.  It retains every
conclusion of Theorem~A on a formation basin and, after restriction near
a distinguished center \(G_{\rm ss}\) and one fixed positive-time restart, equips a
smaller open neighborhood of physical initial metrics with a
transported-restart-invariant marked first-profile coordinate.  This
coordinate is the coefficient of the sharp marked Jacobi expansion,
and its quadratic self-interaction determines the first scale and phase
response.  The associated map is a split \(C^1\) submersion.  Locally it
is the projection in a \(C^1\) product chart, so its fibers form a
finite-codimensional \(C^1\) first-order asymptotic foliation; the zero
fiber is the marked strong-stable leaf.  On the smooth locus, the
corresponding marked first-order leaf space is modeled on the first
stable eigenspace \(E_1\).  Every sufficiently small profile is realized uniquely on a
chosen compactly supported transverse disk, and the quantitative
two-state limit identifies the Jacobi field of the amplitude
difference.  Thus Theorems~A and~B identify FIK blowdown both as an open
full-metric singularity mechanism and, near a distinguished physical
trajectory, as a singularity type carrying a complete marked
first-order asymptotic coordinate.

The proof has one common analytic trunk and two branches.
Spectral certification, exact modulation, receding-domain trapping, and
global one-state continuation prove Theorem~A through the formation
branch alone.  Two-state comparison and first variation of those same
one-state trajectories produce the prepared scattering map in
Theorem~C(I).  Retained exact-core construction data and the
scattering and first-variation analysis yield the separately quantified
transversality statement in Theorem~C(II).  After choosing
\(G_{\rm ss}\), shrinking the physical source, and applying one fixed
positive-time restart into a Part~I ball, these inputs refine
Theorem~A to Theorem~B.  The terminal position of Theorem~B records this
local refinement; it does not enlarge the global formation domain of
Theorem~A.

In Figure~\ref{fig:intro-dependency-spine}, \(\mathscr U_0\) denotes the
large formation basin retained from Theorem~A, while
\(\widehat{\mathscr U}_1\) denotes the smaller physical neighborhood on
which the fixed-restart scattering and foliation conclusions hold.

\newcommand{\FIKdependencyfigure}{%
\begin{center}
\begin{minipage}{0.98\linewidth}
\centering
\begingroup
\hyphenpenalty=10000
\exhyphenpenalty=10000
\begin{tikzpicture}[
  >={Latex[length=2.1mm,width=1.4mm]},
  result/.style={
    draw,
    line width=0.35pt,
    align=center,
    inner xsep=4pt,
    inner ysep=3.5pt,
    font=\small
  },
  wide/.style={result,text width=0.86\textwidth},
  half/.style={
    result,
    text width=0.385\textwidth,
    minimum height=5.7em
  },
  cwide/.style={result,text width=0.78\textwidth},
  primary/.style={
    half,
    line width=0.70pt,
    fill=black!4
  },
  preparedhalf/.style={
    half,
    fill=black!1
  },
  preparedwide/.style={
    cwide,
    fill=black!1
  },
  synthesisbox/.style={
    wide,
    line width=0.45pt
  },
  terminal/.style={
    wide,
    line width=0.60pt,
    fill=black!3,
    inner ysep=3.5pt
  },
  edge/.style={->,line width=0.4pt},
  dataedge/.style={
    ->,
    dashed,
    dash pattern=on 2.2pt off 1.5pt,
    line width=0.35pt
  },
  datalabel/.style={
    font=\scriptsize,
    align=center,
    fill=white,
    inner sep=1pt
  }
]
\node[wide] (spectral) {
  \textbf{Spectral certification.}
  Exact FIK reduction, low-mode certificates, and uniform
  high-frequency comparison
  (Appendix~\ref{app:self-contained-FIK-spectrum})
};
\node[wide,below=4mm of spectral] (geometric) {
  \textbf{Common geometric reduction.}
  The certified nine-mode geometric spectrum yields exact modulation
  and coercive stable dynamics
  (Sections~\ref{sec:setup}--\ref{sec:phase})
};
\node[wide,below=4mm of geometric] (engine) {
  \textbf{Common one-state continuation theory.}
  Receding-domain trapping, adaptive grafts, prepared evolution, and
  endpoint-independent continuation
};

\coordinate (formation-center) at
  ($(engine.south)+(-0.215\textwidth,-7mm)$);
\coordinate (scattering-center) at
  ($(engine.south)+(0.215\textwidth,-7mm)$);
\node[half,anchor=north] (formation) at (formation-center) {
  \textbf{Formation branch (one state).}\\
  Exact-core implantation, strict entrance, and positive-time
  low-topology promotion.
};
\node[half,anchor=north] (scattering) at (scattering-center) {
  \textbf{Scattering and first-variation branch.}\\
  Stable profile and response for one state; two-state comparison and
  first variation of the resulting trajectories.
};

\node[primary,below=5mm of formation] (theoremA) {
  \textbf{Theorem A --- primary geometric theorem.}\\
  Universal full-metric formation basin, relatively
  \(C^{2,\alpha}\)-open among smooth metrics; localized Type-I FIK
  blowdown and full-sequence marked convergence.
};
\node[preparedhalf,below=5mm of scattering] (theoremCI) {
  \textbf{Theorem C(I) --- prepared scattering theorem.}\\
  Sharp \(C^1\) scattering on an arbitrary fixed common-margin
  strict-entrance ball, with buffered ambient extension.
};

\coordinate (output-midpoint) at
  ($(theoremA.south)!0.5!(theoremCI.south)$);
\node[preparedwide,below=16mm of output-midpoint] (theoremCII) {
  \textbf{Theorem C(II) --- exact-core transversality.}\\
  Only on its separately constructed domain: split submersion,
  unique realization on a fixed transverse disk, and a nonempty
  strong-stable zero leaf.
};
\node[synthesisbox,below=6mm of theoremCII] (synthesis) {
  \textbf{Fixed-restart physical synthesis near \(G_{\rm ss}\).}
  Retain Theorem~A on \(\mathscr U_0\), shrink to a smaller physical
  neighborhood, restart once into a Part~I ball, and combine
  Theorem~C(I) scattering with the separately retained
  Theorem~C(II) transverse data.
};
\node[terminal,below=4.5mm of synthesis] (theoremB) {
  \textbf{Theorem B --- physical first-profile theorem.}\\
  Theorem~A holds on \(\mathscr U_0\); on the smaller neighborhood
  near \(G_{\rm ss}\): prescribed profiles, sharp marked Jacobi normal
  form, quantitative two-state limit, and local marked first-order moduli.
};

\draw[edge] (spectral) -- (geometric);
\draw[edge] (geometric) -- (engine);
\coordinate (branch-split) at ($(engine.south)+(0,-2.5mm)$);
\draw[line width=0.4pt] (engine.south) -- (branch-split);
\draw[edge] (branch-split) -| (formation.north);
\draw[edge] (branch-split) -| (scattering.north);
\draw[edge] (formation) -- (theoremA);
\draw[edge] (scattering) -- (theoremCI);

\coordinate (data-left-target) at
  ($(theoremCII.north)+(-0.19\textwidth,0)$);
\coordinate (data-right-target) at
  ($(theoremCII.north)+(0.19\textwidth,0)$);
\coordinate (data-left-above) at
  ($(data-left-target)+(0,4mm)$);
\coordinate (data-right-above) at
  ($(data-right-target)+(0,4mm)$);
\draw[dataedge]
  (formation.south east) -- ++(1mm,-2mm) |-
  (data-left-above) -- (data-left-target);
\draw[dataedge]
  (scattering.south west) -- ++(-1mm,-2mm) |-
  (data-right-above) -- (data-right-target);
\node[datalabel,anchor=south] at
  ($(theoremCII.north)+(-0.20\textwidth,1mm)$)
  {retained exact-core center\\and strict-entrance data};
\node[datalabel,anchor=south] at
  ($(theoremCII.north)+(0.20\textwidth,1mm)$)
  {compact first-stable lift and\\scattering/first-variation analysis};

\draw[edge]
  (theoremA.south west) -- ++(-4mm,-2mm) |- (synthesis.west);
\draw[edge]
  (theoremCI.south east) -- ++(4mm,-2mm) |- (synthesis.east);
\draw[edge] (theoremCII) -- (synthesis);
\draw[edge] (synthesis) -- (theoremB);
\end{tikzpicture}
\endgroup
\refstepcounter{figure}
\label{fig:intro-dependency-spine}
\par\medskip
\parbox[t]{0.90\textwidth}{\raggedright\small
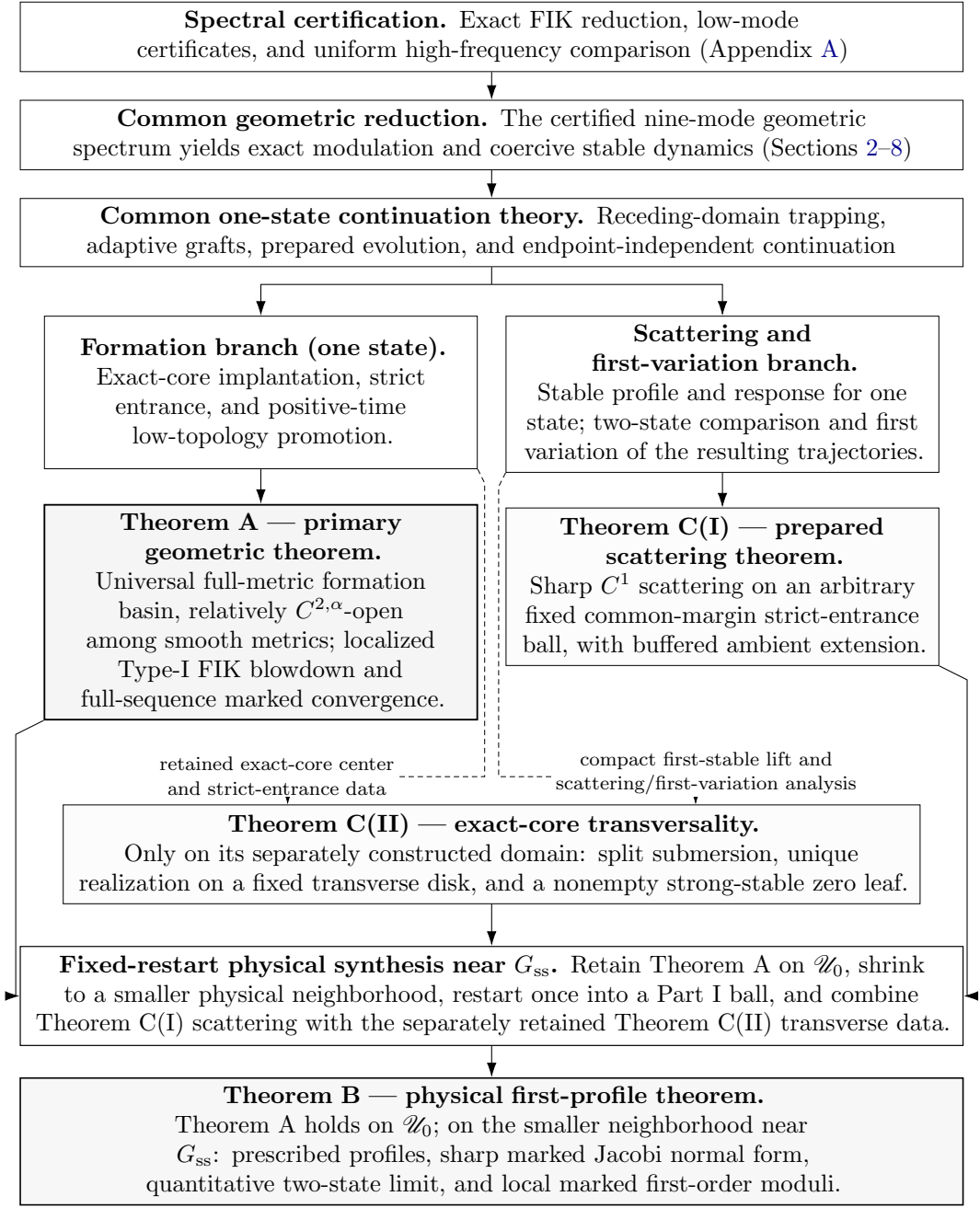
\textsc{Figure~\thefigure.} Proof architecture and theorem hierarchy.
Theorem~A is the primary geometric result and closes using only the
formation branch.  Theorem~C is the prepared analytic theorem: Part~I
gives scattering on an arbitrary fixed common-margin strict-entrance
ball satisfying the hypotheses of Part~I,
whereas Part~II gives transversality only on its independently
constructed exact-core domain.  After restriction near
\(G_{\rm ss}\), one fixed positive-time restart places the physical data
in a Part~I ball; these inputs yield Theorem~B, the physical
first-profile theorem.  Solid arrows trace the principal proof and
synthesis
dependencies.  Dashed arrows denote retained construction data or
analysis reused in Part~II, rather than applications of the preceding
theorem statements.  The restart image need
not lie in the Part~II domain,
and the terminal position of Theorem~B does not assert its local
scattering and foliation conclusions on all of \(\mathscr U_0\).}
\end{minipage}
\end{center}
\smallskip
}

\FIKdependencyfigure

The model has a direct topological meaning.  Let
\[
 M=\mathcal O_{\mathbb P^1}(-1)
\]
and let $(\bar g,\bar f)$ denote the complete asymptotically conical
gradient K\"ahler--Ricci shrinker constructed in~\cite{FIK}.  Its
zero section is an exceptional $(-1)$-sphere; the self-similar flow
collapses this sphere and is the local model associated with complex
blowdown.  The complete classification of K\"ahler--Ricci shrinker
surfaces identifies FIK as the nonflat asymptotically conical model
on the blow-up of $\mathbb C^2$~\cite{LiWang}, while the noncollapsed
K\"ahler-surface singularity theory singles it out as the Type-I
blow-up model~\cite{CifarelliConlonDeruelle,ConlonHallgrenMa}.

Type-I parabolic blow-ups admit nontrivial shrinking-soliton limits
subsequentially~\cite{EndersMullerTopping}; the problem resolved here is
the nonlinear formation, localization, and full-sequence marked
selection of one prescribed noncompact model.  Theorem~B further shows
that, on a smaller physical neighborhood inside this formation basin,
the first stable correction is encoded by
\[
 \mathfrak A_1=\lambda_\infty^{-\gamma_1}V_\infty\in E_1.
\]
Here \(E_1\) is the first stable eigenspace, \(\gamma_1>0\) its decay
exponent, \(\lambda_\infty\) the limiting scale, and \(V_\infty\) the
leading normalized stable coefficient.  The formation basin is open,
but no density assertion is made in the space of all metrics.

\subsection{Main results}

The three main theorems occupy distinct but connected
domains.  Theorem~A produces a relatively \(C^{2,\alpha}\)-open
formation basin \(\mathscr U\) in the space of smooth physical metrics,
using only the one-state formation branch.  Theorem~B retains every
conclusion of Theorem~A on a basin \(\mathscr U_0\) of this kind and, around a
distinguished \(G_{\rm ss}\in\mathscr U_0\), selects a smaller open
little-H\"older neighborhood \(\widehat{\mathscr U}_1\), with smooth
locus \(\mathscr U_1\subset\mathscr U_0\), on which one fixed
positive-time restart produces the marked first-order coordinate and
local asymptotic foliation.  Theorem~C is the prepared analytic theorem:
Part~I gives sharp \(C^1\) scattering on an arbitrary fixed common-margin
sliced ball of strict prepared entrances satisfying its hypotheses,
whereas Part~II gives transversality and profile
realization only on its separately constructed exact-core domain.  The
restart image used for Theorem~B lies in a Part~I ball, not necessarily
the Part~II domain; the physical foliation imports only the
transversality retained from Part~II.  Thus the local refinement in
Theorem~B does not enlarge the global formation domain of Theorem~A.

\paragraph{Metric-space notation.}
Throughout, whenever \(j,k,m,q\), or \(N\) denotes a
differentiability order or a discrete spectral index, it belongs to
\(\mathbb N_0=\{0,1,2,\ldots\}\), unless a different domain is stated
explicitly.  This convention does not apply when one of these letters
is separately declared to be a real-valued function or parameter.
For a closed smooth manifold \(\mathcal N\),
\(\operatorname{Met}^{\infty}(\mathcal N)\) denotes the cone of smooth
Riemannian metrics on \(\mathcal N\).  For an integer \(r\geq0\) and
\(0<\alpha<1\), let
\(h^{r,\alpha}(S^2T^*\mathcal N)\) be the closure of smooth symmetric
two-tensors in the \(C^{r,\alpha}\) norm, and let
\(\operatorname{Met}^{r,\alpha}(\mathcal N)\) denote its open
positive-definite cone.  Thus finite-order metric spaces without an
additional subscript are little-H\"older spaces.  When a big H\"older
metric space is needed, we write
\(\operatorname{Met}_{C}^{r,\alpha}(\mathcal N)\) for the
positive-definite cone in
\(C^{r,\alpha}(S^2T^*\mathcal N)\).  Every relative
\(C^{r,\alpha}\) topology on a smooth locus below is the subspace
topology induced by the corresponding little-H\"older cone.

\begin{maintheorem}[Open full-metric formation basin for the FIK blowdown]
\label{thm:intro-open-basin}
Let \((X^4,g_X)\) be a smooth closed connected oriented Riemannian
manifold, let \(x_*\in X\), let \(B_{\rm imp}\) be an oriented
coordinate neighborhood of \(x_*\), and let \(\varepsilon_T>0\).

\emph{Universal implantation.}
There are \(0<\alpha<1\), an oriented blow-up
\[
 \varpi:\widehat X\longrightarrow X,\qquad
 \widehat X\cong X\#\overline{\mathbb {CP}}^{\,2},
\]
a smooth center metric \(G_*\) containing an exact FIK core, a
  relatively \(C^{2,\alpha}\)-open neighborhood
\[
 \mathscr U\subset\operatorname{Met}^{\infty}(\widehat X)
\]
of \(G_*\), where ``relatively open'' means open for the
\(C^{2,\alpha}\) topology induced on the space of smooth metrics, and
constants
\[
 0<c\leq C<\infty,\qquad C_{\rm out}<\infty,
\]
such that
\begin{equation}\label{eq:intro-prescribed-exterior}
 G_*=\varpi^*g_X
 \quad\text{on}\quad
 \widehat X\setminus\varpi^{-1}(B_{\rm imp}).
\end{equation}

\emph{Open Type-I formation and localization.}
For every \(G_0\in\mathscr U\), let \(G(t;G_0)\) denote its maximal
Ricci flow.  Its singular time satisfies
\[
 0<T(G_0)<\varepsilon_T.
\]
There is a time \(0\leq t_{\rm I}(G_0)<T(G_0)\) such that, for every
\(t_{\rm I}(G_0)\leq t<T(G_0)\),
\begin{equation}\label{eq:intro-global-Type-I}
 c\leq
 (T(G_0)-t)
 \|\Rm_{G(t;G_0)}\|_{L^\infty(\widehat X,G(t;G_0))}
 \leq C.
\end{equation}
Moreover,
\[
 \sup_{G_0\in\mathscr U}
 \sup_{0\leq t<T(G_0)}
 \sup_{\widehat X\setminus\varpi^{-1}(B_{\rm imp})}
 |\Rm_{G(t;G_0)}|
 \leq C_{\rm out}.
\]
The buffered local derivative estimates consequently give smooth
convergence as \(t\uparrow T(G_0)\) on every compact subset of
\(\widehat X\setminus\varpi^{-1}(B_{\rm imp})\).  Thus every
curvature-singular point is confined to the prescribed implantation
region.

\emph{Full-sequence fixed-convention marked model.}
For each \(G_0\in\mathscr U\) and every \(t<T(G_0)\) sufficiently close
to \(T(G_0)\), there are relatively compact domains
\(\mathcal U_t(G_0)\Subset M\), exhausting \(M\) as
\(t\uparrow T(G_0)\), and marked embeddings
\[
 \Xi_t^{G_0}:\mathcal U_t(G_0)\longrightarrow\widehat X
\]
with the following property.  For every sequence
\(t_i\uparrow T(G_0)\), set
\(\delta_i=T(G_0)-t_i\) and freeze the base-time marking
\(\Xi_{t_i}^{G_0}\) along the \(i\)-th parabolically rescaled flow.
Then
\[
 \delta_i^{-1}(\Xi_{t_i}^{G_0})^*
 G\bigl(T(G_0)+s\delta_i;G_0\bigr)
 \longrightarrow g_{\mathrm{FIK}}(s)
 \quad\text{in}\quad
 C^\infty_{\mathrm{loc}}\bigl(M\times(-\infty,0)\bigr)
\]
without passing to a subsequence, where \(g_{\mathrm{FIK}}(s)\),
\(s<0\), is the canonical ancient self-similar flow generated by
\((\bar g,\bar f)\).

\emph{Marked exceptional-sphere collapse.}
Let \(E\subset M\) be the zero section.  For every
\(G_0\in\mathscr U\) there is a single time
\(t_*(G_0)<T(G_0)\), independent of the derivative order below, such
that, whenever \(t_*(G_0)\leq t<T(G_0)\), the zero section lies in
\(\mathcal U_t(G_0)\), and
\[
 \Sigma_t(G_0):=\Xi_t^{G_0}(E)\subset\widehat X
\]
is an embedded sphere with normal bundle of degree \(-1\).  Here and
below, if \(\Sigma\) is an embedded two-sphere in an oriented
four-manifold \(Y\), the degree of its normal bundle means its oriented
self-intersection number: choose either orientation on \(\Sigma\), orient
\(N_{\Sigma/Y}\) so that \(T\Sigma\oplus N_{\Sigma/Y}\) carries the
ambient orientation, and set
\[
 \deg N_{\Sigma/Y}
 :=\left\langle e(N_{\Sigma/Y}),[\Sigma]\right\rangle .
\]
Reversing the chosen orientation of \(\Sigma\) also reverses the induced
normal orientation, so this integer is independent of that choice.  In the same
marking, writing \(g_E=\bar g|_E\), there is \(\vartheta>0\), and for
every integer \(m\geq0\) there is a finite constant \(C_m\), such that
for every \(G_0\in\mathscr U\) and every
\(t_*(G_0)\leq t<T(G_0)\),
\begin{equation}\label{eq:intro-marked-bolt-collapse}
 \left\|
 (T(G_0)-t)^{-1}
 \left(\Xi_t^{G_0}|_E\right)^*
 G(t;G_0)
 -g_E
 \right\|_{C^m(E,g_E)}
\leq C_m(T(G_0)-t)^\vartheta .
\end{equation}
Consequently,
\[
 \operatorname{Area}_{G(t)}\Sigma_t(G_0)
 =(T(G_0)-t)
 \left(\operatorname{Area}_{g_E}E
       +O((T(G_0)-t)^\vartheta)\right),
\]
and the intrinsic diameter of \(\Sigma_t(G_0)\) equals
\[
 (T(G_0)-t)^{1/2}
 \left(\operatorname{diam}_{g_E}E
       +O((T(G_0)-t)^\vartheta)\right).
\]
Moreover, uniformly for \(z\in E\),
\[
 (T(G_0)-t)
 |\Rm_{G(t;G_0)}|\bigl(\Xi_t^{G_0}(z)\bigr)
 =
 |\Rm_{\bar g}|(z)+O((T(G_0)-t)^\vartheta).
\]
In particular, the ambient diameter of \(\Sigma_t(G_0)\) is
\(O((T(G_0)-t)^{1/2})\).
This is a marked, dynamically transported exceptional sphere; the
theorem does not identify it with one fixed submanifold for every
metric in the open basin.

\emph{Uniformity and absence of constraints.}
The late-time threshold \(t_*(G_0)\) may depend on the trajectory.
All numerical constants above, including every \(C_m\), and the size of
\(\mathscr U\) may depend on
\((X,g_X,x_*,B_{\rm imp},\varepsilon_T)\), the chosen normalization of
FIK, and the fixed implantation and preparation margins, but not on
\(G_0\in\mathscr U\).  Metrics in \(\mathscr U\) need not preserve the
exact core or the prescribed exterior and are subject to neither
\(U(2)\)-symmetry nor a K\"ahler condition.
\end{maintheorem}

\begin{remark}[Why the collapsing sphere is marked rather than fixed]
A common fixed exceptional sphere for every metric in the open basin
would contradict diffeomorphism covariance.  Indeed, a sufficiently
small diffeomorphism supported in the implantation region can move the
center sphere while keeping the pulled-back metric in the same open
neighborhood, and Ricci flow transports the singular geometry by the
inverse diffeomorphism.  The family \(\Sigma_t(G_0)\) above is therefore
the invariantly meaningful conclusion in the fixed marking convention;
it transforms equivariantly when the initial metric is pulled back.
\end{remark}

The exact-core metric itself is the center of the basin.  In particular,
Theorem~\ref{thm:intro-open-basin} localizes the blowdown
mechanism on the oriented blow-up of an arbitrary host and permits the
singular time to be as small as prescribed.
Taking \(X=S^4\) recovers the
\(\overline{\mathbb {CP}}^{\,2}\) ambient manifold of the classical
blowdown picture.

For the statements of Theorems~B and~C, let
\(\A=\bar\Delta_{\bar f}+2\overline{\Rm}\) denote the self-adjoint
weighted FIK operator fixed in
Section~\ref{sec:setup}, let \(\mathcal Z\) be its nine-dimensional
nonnegative geometric spectral space, and define
\[
 -\gamma_1
 :=\max\operatorname{spec}(\A|_{\mathcal Z^\perp})<0,
 \qquad E_1:=\ker(\A+\gamma_1 I).
\]
For \(A\in E_1\), \(\mathcal J_A\) denotes the corresponding
base-time-marked physical Jacobi field, normalized by
\(\mathcal J_A(-1)=A\); see
\eqref{eq:marked-Jacobi-field}.
For the introductory statements, set
\[
 d\nu:=(4\pi)^{-2}e^{-\bar f}\,dV_{\bar g},
 \qquad
 L^2_\nu:=L^2(M;S^2T^*M,d\nu).
\]
These conventions are recalled in Section~\ref{sec:setup} as part of
the global analytic framework.
In the fixed spectral normalization used throughout, choose an
\(L^2_\nu\)-orthonormal real basis \(Z_0,\ldots,Z_8\) of
\(\mathcal Z\), with \(Z_0\) the positively oriented Ricci mode.  For
the eight Hessian modes fix real potentials \(\phi_j\) and nonzero
constants \(\mathfrak c_j\), and define
\[
 Z_j=\mathfrak c_j\bar\nabla^2\phi_j,
 \qquad
 W_j:=\frac{\mathfrak c_j}{2}\bar\nabla\phi_j,
 \qquad
 Z_j=\Lie_{W_j}\bar g,\quad1\leq j\leq8.
\]
Theorem~\ref{thm:FIK-spectrum} certifies these modes, and
Section~\ref{sec:setup} records the same choices as a paper-wide
convention.
\begin{maintheorem}[Universal FIK formation with prescribed
profiles and local marked first-order asymptotic moduli]
\label{thm:intro-physical-amplitude-master}
Fix any host, implantation datum, and prescribed singular-time bound
covered by Theorem~\ref{thm:intro-open-basin}.  The construction has
the following additional structure.

\emph{I. Formation and the physical restart bridge.}
The formation neighborhood in Theorem~A may be chosen as a relative
\(C^{2,\alpha}\)-open basin \(\mathscr U_0\) on which every conclusion
of that theorem holds.  There are
\(G_{\rm ss}\in\mathscr U_0\), a physical time \(d>0\), and an open
neighborhood \(\widehat{\mathscr U}_1\) of \(G_{\rm ss}\) in the
positive cone of \(h^{2,\alpha}(S^2T^*\widehat X)\), with smooth locus
\[
 \mathscr U_1
 :=\widehat{\mathscr U}_1\cap\operatorname{Met}^{\infty}(\widehat X)
 \subset\mathscr U_0,
\]
such that a single DeTurck gauge, marking, graft, scale, and initial-map
convention, fixed independently of the initial metric and applied to
each time-\(d\) metric, gives a \(C^1\) restart map
\[
 \mathcal P_d:\widehat{\mathscr U}_1
 \longrightarrow
 \{\text{one common-margin sliced ball of strict prepared entrances}\}.
\]
For smooth data its tail is the original Ricci flow, up to the terminal
pullback produced by the fixed gauge convention and absorbed into this
common physical marking.

\emph{II. One transported-restart-invariant coordinate controls the
sharp normal form.}
For \(G\in\widehat{\mathscr U}_1\), evaluate all four components of
the Part~I scattering map in
Theorem~\ref{thm:intro-sharp-scattering} at \(\mathcal P_d(G)\), with the
tail physical clock normalized to vanish at the restart time \(d\), and
denote the resulting components by
\[
 T_{\rm tail}(G),\qquad \log\lambda_\infty(G),\qquad
 \Psi_\infty(G),\qquad V_\infty(G)\in E_1.
\]
Set \(T(G)=d+T_{\rm tail}(G)\).
Write \(\tau\) for its normalized clock,
\[
 \frac{dt}{d\tau}=\lambda(\tau)>0,
\]
write \(\tau(t)\) for the inverse clock and
\(\lambda(t):=\lambda(\tau(t))\), and let \(\Psi_\tau\) be the
transported phase in the fixed preparation convention.  On every
compact set, write \(\log_{\operatorname{Id}}\) for the inverse of the
fixed exponential chart at the identity used for phase differences.
Then
\[
 G\longmapsto
 \left(T(G),\log\lambda_\infty(G),\Psi_\infty(G),
       \mathfrak A_1(G)\right),
 \qquad
 \mathfrak A_1(G):=
 \lambda_\infty(G)^{-\gamma_1}V_\infty(G),
\]
is \(C^1\) componentwise: the scalar and \(E_1\)-valued components are
Banach-valued \(C^1\), while for every \(K\Subset M\) and \(m\geq0\)
the phase is \(C^1\) into the fixed \(C^m(K)\) exponential chart.
The amplitude \(\mathfrak A_1\) is unchanged by every transported
forward restart.  If normalized time is additively relabeled, every
time-typed cutoff, chart, column, and slice is relabeled with it; this
does not identify an independently prepared state.

There is \(\eta_1>0\) such that, for \(G\in\mathscr U_1\),
\(K\Subset M\), \(I\Subset(-\infty,0)\), and \(m\geq0\), the following
holds as \(\delta\downarrow0\).  For all sufficiently small
\(\delta>0\), put \(t=T(G)-\delta\); then \(K\) lies in the marking
domain and \(T(G)+s\delta\) lies in the flow interval for every
\(s\in I\), and the frozen base-time markings satisfy
\begin{equation}\label{eq:intro-physical-Jacobi-normal-form}
 \delta^{-1}\left(\Xi^G_{T(G)-\delta}\right)^*
 G\bigl(T(G)+s\delta;G\bigr)
 =
 g_{\rm FIK}(s)
 +\delta^{\gamma_1}\mathcal J_{\mathfrak A_1(G)}(s)
 +O_{C^m(K\times I)}
  \!\left(\delta^{\gamma_1+\eta_1}\right).
\end{equation}
Let
\[
 \mathbf c^{(2)}:E_1\longrightarrow\mathbb R^9,
 \qquad
 \mathbf c^{(2)}(V)
 =\bigl(a^{(2)}(V),b_1^{(2)}(V),\ldots,b_8^{(2)}(V)\bigr),
\]
be the continuous quadratic response map defined by
\eqref{eq:Q2-definition}--\eqref{eq:quadratic-feedback-coefficient}, and
define its phase field by
\[
 U^{(2)}(V):=\sum_{j=1}^8b_j^{(2)}(V)W_j.
\]
Thus \(a^{(2)}:E_1\to\mathbb R\) and
\(U^{(2)}:E_1\to\operatorname{span}\{W_1,\ldots,W_8\}\) are fixed
quadratic maps.  Then the same amplitude
determines the first nonlinear geometric response:
\begin{align}
 \frac{\lambda(t)}{T(G)-t}
 &=
 1+
 \frac{a^{(2)}(\mathfrak A_1(G))}{1+2\gamma_1}
 (T(G)-t)^{2\gamma_1}
 +o((T(G)-t)^{2\gamma_1}),
 \label{eq:intro-master-quadratic-scale}\\
 \log_{\operatorname{Id}}\!
 \left(\Psi_\infty\circ\Psi_{\tau(t)}^{-1}\right)
 &=
 \frac{(T(G)-t)^{2\gamma_1}}{2\gamma_1}
 U^{(2)}(\mathfrak A_1(G))
 +o_{C^m(K)}((T(G)-t)^{2\gamma_1}).
 \label{eq:intro-master-quadratic-phase}
\end{align}
The response law allows a component of the quadratic coefficient to
vanish.

\emph{III. Physical first-profile moduli and first-order completeness.}
The map
\[
 \mathfrak A_1:\widehat{\mathscr U}_1\longrightarrow E_1
\]
is a split \(C^1\) submersion, and
\(\mathfrak A_1(G_{\rm ss})=0\).  After shrinking
\(\widehat{\mathscr U}_1\) once more about \(G_{\rm ss}\), set
\[
 K_{\rm ss}:=\ker D\mathfrak A_1(G_{\rm ss}).
\]
There are neighborhoods \(\mathcal O\subset\widehat{\mathscr U}_1\)
of \(G_{\rm ss}\), \(\mathcal O_K\subset K_{\rm ss}\) of the origin,
and \(B_{\epsilon}^{E_1}(0)\subset E_1\), and a \(C^1\)
diffeomorphism
\[
 \Phi:\mathcal O
 \longrightarrow
 \mathcal O_K\times B_{\epsilon}^{E_1}(0)
\]
such that
\[
 \mathfrak A_1\!\left(\Phi^{-1}(k,A)\right)=A.
\]
Thus the amplitude is locally the projection onto \(E_1\).  Its
sufficiently small fibers are
codimension-\(\dim E_1\) split Banach submanifolds whose connected
components form a local \(C^1\) foliation.  The zero fiber is nonempty,
has the improved marked strong-stable rate, and has open dense
complement.  A compactly supported physical transverse disk realizes
every sufficiently small \(A\in E_1\) uniquely and agrees with the
implanted center outside the implantation region.

For \(G_1,G_2\in\mathscr U_1\), define in the common convention
\[
 \mathcal G_{i,\delta}(s)
 :=\delta^{-1}
 \left(\Xi^{G_i}_{T(G_i)-\delta}\right)^*
 G_i(T(G_i)+s\delta).
\]
Then, for every \(K\Subset M\), \(I\Subset(-\infty,0)\), and
\(m\geq0\),
\begin{equation}\label{eq:intro-master-quantitative-two-state}
 \delta^{-\gamma_1}
 \bigl(\mathcal G_{1,\delta}-\mathcal G_{2,\delta}\bigr)
 =
 \mathcal J_{\mathfrak A_1(G_1)-\mathfrak A_1(G_2)}
 +O_{C^m(K\times I)}(\delta^{\eta_1}).
\end{equation}
Consequently,
\begin{equation}\label{eq:intro-master-first-order-completeness}
 \mathfrak A_1(G_1)=\mathfrak A_1(G_2)
 \quad\Longleftrightarrow\quad
 \delta^{-\gamma_1}
 \bigl(\mathcal G_{1,\delta}-\mathcal G_{2,\delta}\bigr)
 \longrightarrow0
 \quad\text{in }C^\infty_{\rm loc}
\end{equation}
as \(\delta\downarrow0\); for the reverse implication one compact time
window containing \(s=-1\) suffices.  Hence, on the smooth relative
locus \(\mathcal O\cap\operatorname{Met}^{\infty}(\widehat X)\), the
local leaf space for marked first-order asymptotic agreement is
modeled, through \(\mathfrak A_1\), on
\(B_\epsilon^{E_1}(0)\).

\emph{IV. Scope.}
These are statements about actual initial metrics re-prepared at
positive time by this single fixed convention, and about the marking
and gauge transported along each resulting tail.  No quotient by
arbitrary diffeomorphisms, arbitrary time-dependent re-marking, or
global foliation of all of \(\mathscr U_0\) is asserted.  Neither open
physical neighborhood imposes symmetry, a K\"ahler condition,
exact-core identity, or finite-dimensional tuning.
The present spectral input identifies a six-dimensional pure-gauge
eigenspace at eigenvalue \(1-\sqrt2\), but neither proves that this is
the stable spectral edge nor identifies all of \(E_1\); no such
numerical identification is used.
\end{maintheorem}

\medskip
\noindent\textbf{Prepared theorem and physical bridge.}
The formal statement of Theorem~C appears in
Section~\ref{sec:global-two-state}.  Its two parts retain the distinct
domains summarized above and in
Figure~\ref{fig:intro-dependency-spine}: Part~I is the general prepared
scattering theorem, while the Part~II foliation and realization
conclusions belong only to the specially constructed exact-core domain.

The passage back to physical data is equally specific.  After choosing
\(G_{\rm ss}\) and shrinking the source, one fixed positive-time map
gives
\[
 \mathcal P_d:
 \widehat{\mathscr U}_1
 \longrightarrow
 \{\text{one common-margin prepared ball}\},
 \qquad
 \mathscr U_1
 =\widehat{\mathscr U}_1\cap
   \operatorname{Met}^{\infty}(\widehat X)
 \subset\mathscr U_0 .
\]
Positive-time smoothing makes this bridge \(C^1\) from the
\(h^{2,\alpha}\) physical chart, but the bridge does not identify the
physical and prepared topologies.  Theorem~A uses only the one-state
formation branch; the comparison and first-variation analysis refines
the resulting trajectories and enters Theorem~B only through this fixed
restart.

\subsection{From realization to robustness}

Several lines of work lead to these theorems.  Song--Weinkove proved
contraction of the exceptional divisor on Hirzebruch surfaces under an
invariance hypothesis~\cite{SongWeinkove}, and subsequently removed the
symmetry assumption for canonical contractions under the relevant
K\"ahler cohomological condition~\cite{SongWeinkoveContraction}.
M\'aximo constructed an open FIK-forming family of
\(U(2)\)-invariant K\"ahler metrics in real dimension
four~\cite{Maximo}.  Song proved Type-I formation for rotationally
symmetric K\"ahler flows on \(\mathbb{CP}^n\) blown up at one point and,
in the noncollapsing case, obtained complete nonflat shrinking K\"ahler
blow-up limits~\cite{SongTypeI};
Guo--Song subsequently
identified the noncollapsed blow-up limit along the exceptional divisor
as FIK in the \(U(n)\)-invariant setting~\cite{GuoSong}.  More recently,
Cifarelli--Conlon--Deruelle classified the noncollapsed Type-I
K\"ahler-surface model, and Conlon--Hallgren--Ma proved that every
noncollapsed finite-time singularity on a compact K\"ahler surface is
Type I~\cite{CifarelliConlonDeruelle,ConlonHallgrenMa}.  Taken together,
these results identify FIK, without a symmetry hypothesis, throughout
the noncollapsed finite-time K\"ahler-surface regime.

Outside the K\"ahler category, Isenberg--Knopf--\v{S}e\v{s}um developed a
\(U(2)\)-invariant picture whose final model identification depends on
two technical conjectures~\cite{IKS}, while the numerical evolutions of
Garfinkle--Isenberg--Knopf--Wu provide further evidence for FIK
formation under non-K\"ahler perturbations~\cite{GarfinkleIsenbergKnopfWu}.
In a different Type-II regime, Appleton proved Eguchi--Hanson formation
for a large class of \(U(2)\)-invariant non-K\"ahler flows on noncompact
four-manifolds~\cite{Appleton}.  None of these results gives a
neighborhood open in the full space of Riemannian metrics with a proved
FIK model.

Beyond symmetry, Stolarski proved that every asymptotically conical
shrinker can occur as a local Type-I singularity of a closed Ricci
flow~\cite{Stolarski}.  His construction uses a Wa\.{z}ewski exit
argument to select finitely many weighted-Lichnerowicz coefficients.
Hughes used the same strategy to find an arbitrarily
\(L^2\)-small compact perturbation of the noncompact Taub--Bolt metric
that develops FIK~\cite{Hughes}.  These are realization theorems:
they produce selected trajectories, not a neighborhood all of whose
flows have the prescribed model.  Classical open formation results do
exist for cylindrical neckpinches, but the relevant constructive
theorems are open inside preserved ansatz classes: rotational symmetry
for Angenent--Knopf~\cite{AngenentKnopf}, and the warped-Berger class
for Isenberg--Knopf--\v{S}e\v{s}um~\cite{IKSNeckpinch}.

The preceding work separates three levels.  Symmetry- or
K\"ahler-restricted formation results give open families within
preserved classes, while the classification results identify FIK once
a noncollapsed K\"ahler-surface singularity occurs.  The realization
results of Stolarski and Hughes select individual trajectories.
Theorem~A gives full-Riemannian openness: once the exact FIK core is
implanted, every sufficiently small Riemannian perturbation of the
resulting closed metric develops the same localized FIK model.  This
passage from realization to an unrestricted local basin is the central
advance of the paper.

The decisive input making that passage possible is spectral.
Naff--Ozuch discovered the nine-mode geometric picture for FIK and
organized its analysis by a tensor Wigner
decomposition~\cite{NaffOzuch}, in the stability framework originating
with Cao--Hamilton--Ilmanen~\cite{CaoHamiltonIlmanen}.

For logical independence, this paper proves the exact spectral
statement used by the nonlinear argument.  Appendix
\ref{app:self-contained-FIK-spectrum} reconstructs the FIK tensor
frame, connection and curvature action, establishes the Peter--Weyl
 reduction and the exact \(J\leq4\) block inventory, where
 \(J\in\mathbb N_0\) is the \(SU(2)\) Peter--Weyl highest-weight index,
 certifies the finite
low-mode algebra over \(\mathbb Q(\sqrt2)\), and proves the radial and
exceptional mode counts.  After the Wigner definitions and
cross-sectional reduction, the same appendix proves the remaining
uniform high-frequency comparison in
Lemma~\ref{lem:FIK-high-frequency}.  Together these arguments show that the
nonnegative spectral space of \(\A\) is nine-dimensional: one mode is
scaling and the other eight are Lie derivatives, while all remaining
spectrum is strictly negative.  Thus, modulo scale and diffeomorphisms,
there is no nonnegative mode to tune.

This spectral statement supplies the coercive input but does not by
itself imply nonlinear stability.  Nonpositivity of Perelman's entropy
Hessian~\cite{Perelman} is likewise insufficient: Fubini--Study is
neutrally linearly
stable but dynamically unstable~\cite{Kroencke}.  For FIK one must
control an expanding part of a noncompact model inside a closed flow, a
moving harmonic-map gauge, and a grafting annulus receding into the
conical end; one must then compare two such evolutions strongly enough
to differentiate their asymptotic data.  The paper converts the nine
geometric modes into exact feedback, closes this receding-domain
evolution, and identifies the resulting physical scattering geometry.

The basin argument has three essential interfaces: every strict
prepared entrance continues to the marked FIK limit; an exact implanted
core furnishes a strict prepared entrance; and strict entrance persists
under the prepared perturbation and positive-time low-topology
promotion used to construct the relative \(C^{2,\alpha}\)-open basin.
The scattering theory refines this same entrance class rather than
introducing a second construction.

Two nearby convergence theories have a different logical starting
point.  Colding--Minicozzi establish gauge-theoretic rigidity and
tangent-flow uniqueness for cylindrical Ricci-flow
singularities~\cite{ColdingMinicozzi}, while Choi--Lai obtain optimal
convergence after a flow is already known to approach a compact
integrable shrinker~\cite{ChoiLai}.  Neither result is a formation
theorem for a noncompact, noncylindrical model.  Here the conclusion is
full-sequence but marked: it makes no assertion about unrelated
basepoints or markings.  The relatively open basin is not claimed to
be open dense in the full metric space, and no continuation through
the singular time is asserted.

\subsection{Proof architecture}

Figure~\ref{fig:intro-dependency-spine} gives the global dependency
graph.  We now describe the analytic mechanisms along its common trunk
and its two branches.

The first step turns the nine-dimensional nonnegative spectral space
into exact nonlinear modulation, rather than nine parameters that must
be tuned.  Use the orthonormal eigenbasis
\(Z_0,\ldots,Z_8\) and the normalized generators fixed in
Section~\ref{sec:setup}.  The scale mode is diagonalized by
coupling scaling to the soliton radial flow, while the other eight modes
are generated by vector fields: for the fixed constants
\(\mathfrak c_j\) and
soliton potentials $\phi_j$,
\[
 \bar g-\Lie_{\bar\nabla\bar f}\bar g=2\Ric_{\bar g},
 \qquad
 Z_j=\mathfrak c_j\bar\nabla^2\phi_j=\Lie_{W_j}\bar g,
 \quad W_j=\frac{\mathfrak c_j}{2}\bar\nabla\phi_j
 \quad(1\leq j\leq8).
\]
For the normalized perturbation $h$, set
$\rho_\tau=\rho(e^{-\tau}\bar f)$ and $H=\rho_\tau h$, where $\rho$ is
the fixed drift-adapted cutoff.  Writing the selected scale and gauge
velocities as $\mathbf v=(a,b_1,\ldots,b_8)$, we impose the exact slice
\begin{equation}\label{eq:intro-exact-slice}
 \langle H,Z_\mu\rangle_{L^2(d\nu)}=0,
 \qquad 0\leq\mu\leq8.
\end{equation}
Differentiation gives a uniformly invertible $9\times9$ Gram system.
Self-adjointness of $\A$ removes the entire linear nonnegative spectrum
from its right-hand side, and the Gaussian first-moment estimate
controls the linearly growing gauge fields and the radial scale
transport.  Consequently $\mathbf v$ is quadratic in the stable
perturbation, up to the explicitly estimated modal part of the external
forcing; on the slice, the negative spectral gap yields the coercive
weighted energy inequality used throughout the continuation argument.

The second step closes that energy estimate on a domain receding into
the conical end.  The drift identity for
$\rho(e^{-\tau}\bar f)$ cancels the potentially order-one moving-cutoff
term, and the remaining annular errors are superexponentially small in
the Gaussian norms.  On a finite bootstrap interval
$[\tau_0,\tau_1]$, put
$b=(b_1,\ldots,b_8)$,
$q(\tau)=|a(\tau)|+|b(\tau)|$, and
 $P_{\tau_1}(\tau)=\int_\tau^{\tau_1}q(s)\,ds$.  The energy and Gram
 laws give, for bootstrap rates
 \(0<\sigma<\theta<\gamma_1\) and constants
\(c,C>0\),
$P_{\tau_1}(\tau)\leq
C e^{-2\theta\tau}+Ce^{-ce^\tau}$ uniformly in $\tau_1$.  If
$\mathcal B_0$ is an unmodulated scalar barrier in the intermediate or
outer region, choose fixed sufficiently large constants $K,K_0>0$ and
define the corrected barrier
\begin{equation}\label{eq:intro-corrected-barrier}
 \mathcal B(\tau,x)=
 \exp\!\left(K\int_{\tau_0}^{\tau}q(s)\,ds\right)
 \mathcal B_0(\tau,x)-K_0P_{\tau_1}(\tau),
\end{equation}
which absorbs both effects of the modulation: the exponential handles
transport and zeroth-order action,
whereas
$\partial_\tau[-K_0P_{\tau_1}]=K_0q$ handles the direct modulation
columns.  This closes the pointwise barriers without assuming a
pointwise estimate for $q$.  Conjugation by the accumulated modulation
flow then makes derivative recovery depend only on $\int q$; the
resulting three-region $C^2$ estimate upgrades the integrated
 dissipation to an instantaneous $H^1_\nu$ bound and returns, through
 the exact Gram law, the bootstrap pointwise velocity bound
\[
 q(\tau)\leq C e^{-2\sigma\tau}+Ce^{-ce^\tau}.
\]
After global continuation is closed, the stable spectral estimate
sharpens \(H^1_\nu\) from the three-region rate to
\(O(e^{-\theta\tau})\); only at that later stage does the same exact
Gram law yield
\(q(\tau)=O(e^{-2\theta\tau})+O(e^{-ce^\tau})\), as recorded in
\eqref{eq:master-bootstrap-H1}--\eqref{eq:master-bootstrap-velocity}.

The third step makes the noncompact argument compatible with a closed
ambient flow.  Let $\Theta$ be the adaptive soliton chart following the
selected scale and gauge, let $S=\lambda\Theta^*\bar g$ be its reference
metric, and let $\Phi$ be the controlled harmonic-map gauge.  If
$\acute G$ denotes the evolving extended metric on the model
manifold, the relative map $F=\Theta^{-1}\circ\Phi$ satisfies the
exact identities
\begin{equation}\label{eq:intro-relative-HMHF}
 \partial_tF=\Delta_{\acute G,S}F,\qquad
 (F^{-1})^*\acute G-S=\lambda\Theta^*h,
\end{equation}
where $\Delta_{\acute G,S}$ is the harmonic-map Laplacian from
$(M,\acute G)$ to $(M,S)$.  Thus the selected geometric phase is
factored out rather than counted as a large forcing.  The target has
finite total variation in scale-invariant bounded-geometry norms, the
relative map displaces the fixed graft annulus by
$O(e^{-\tau_0/2})$, and all pure graft and adaptive-column errors occur
where $\bar f\simeq e^\tau$ and are superexponentially small in fixed
Gaussian norms.  A uniform finite-dimensional implicit-function
theorem places every nearby prepared graft in
\eqref{eq:intro-exact-slice}.  For the center metric, however, no
correction is needed: implanting an exact FIK core and taking, for the
soliton radial flow $\varphi_\tau$ and an implantation parameter
$A>0$,
\begin{equation}\label{eq:intro-exact-core-zero}
 \lambda_0=Ae^{-\tau_0},\qquad
 \Theta_0=\Phi_0=\varphi_{\tau_0},\qquad
 F_0=\operatorname{Id},\qquad h_0=0,
\end{equation}
gives a complete prepared extension equal to
$\lambda_0\varphi_{\tau_0}^*\bar g$.  Here $A$ fixes the physical
implantation radius, while $\tau_0$ independently makes the singular
scale and remaining time small.  Since the three-region inequalities
and graft compatibility are strict, this exact zero-perturbation
center generates the full neighborhood in
Theorem~\ref{thm:intro-open-basin}.

The final step upgrades stability to a physical scattering geometry.
If \(E_1\) is the first stable eigenspace and \(-\gamma_1\) its
eigenvalue, the stable normal form gives, for some
\(\delta_{\rm prof}>0\),
\[
 H(\tau)=e^{-\gamma_1\tau}V_\infty
 +O_{H^1_\nu}\!\left(
   e^{-(\gamma_1+\delta_{\rm prof})\tau}\right),
 \qquad V_\infty\in E_1,
\]
but neither \(V_\infty\) nor the normalized clock alone is physical.
Under a transported forward restart in the fixed preparation, marking,
and gauge convention, the scale limit and stable coefficient transform
with opposite weights.  Consequently
\[
 \mathfrak A_1=\lambda_\infty^{-\gamma_1}V_\infty
\]
is unchanged by every such restart.  It is exactly the coefficient of
the frozen physical Jacobi field, and its quadratic self-interaction
determines the first scale and phase corrections.

The global two-state estimate,
Theorem~\ref{thm:global-two-state-estimate}, compares nearby prepared
solutions in a common relative gauge after using the two buffered
derivatives built into the prepared input.  The differentiated normal
form then makes the prepared scattering map \(C^1\).  Compactly
supported \(E_1\) seed tensors give an invertible derivative on a
physical transverse disk at its unique zero-profile point.  The
additional bridge is a single fixed physical time \(d>0\): positive-time
smoothing and the centered phase construction define one
\(C^1\) map from an open \(h^{2,\alpha}\) neighborhood of actual
metrics to one common prepared ball.  Along the smooth transverse disk,
the resulting prepared amplitude maps converge in \(C^1\) to the original
physical-amplitude map
\(v\mapsto\lambda_\infty(G_v)^{-\gamma_1}V_\infty(G_v)\) as
\(d\downarrow0\).  A one-sided parameterized implicit-function
argument moves its zero slightly and preserves the invertible disk
derivative.  The Banach submersion theorem
therefore produces the physical \(\mathfrak A_1\)-foliation of
Theorem~\ref{thm:intro-physical-amplitude-master}, and the sharp Jacobi
expansion proves that its leaves are exactly the fixed-marking
first-order asymptotic classes.

\subsection{Organization}

Appendix~\ref{app:self-contained-FIK-spectrum} supplies the exact FIK
reduction, finite low-frequency certification, and, in
 Lemma~\ref{lem:FIK-high-frequency}, the uniform \(J\geq5\) comparison
 in that same highest-weight index.
Sections
\ref{sec:setup}--\ref{sec:phase} turn the resulting spectral theorem
into a coercive geometric slice and exact modulation.

Sections~\ref{sec:tails}--\ref{sec:prepared-Banach-evolution} construct
the receding-domain evolution.  The one-state continuation and
exact-core promotion in Sections~\ref{sec:uniform-entrance} and
\ref{sec:prepared-open-basin} yield Theorem~A and provide the base
trajectories used later by the comparison theory.  Sections
\ref{sec:finite-horizon-two-state} and
\ref{sec:global-two-state} compare those trajectories and establish
Part~I of Theorem~C.  The same scattering and first-variation analysis,
combined with the specially constructed exact-core transverse seed,
gives profile realization and Part~II.  Section
\ref{sec:first-stable-profile} then combines Theorems~A and~C with the
fixed positive-time restart to prove Theorem~B.

Appendix~\ref{app:diagonal-scale-action} isolates the diagonal
scale-action calculation used in physical reconstruction, and
Appendix~\ref{app:FIK-Gaussian-tail} records the explicit Gaussian-tail
calculation used by the localized moving-frame estimates.

\subsection*{Acknowledgments}

The author is grateful to Bennett Chow, Michael Freedman, and Yongjia
Zhang for earlier collaboration on four-dimensional Ricci-flow
singularity models.

\part{The FIK operator and exact geometric modulation}

\section{Quantitative framework and the FIK operator}
\label{sec:setup}

\subsection{Global conventions and the prepared state space}
\label{subsec:global-prepared-conventions}

For the dynamical refinement, normalize FIK by
\begin{equation}\label{eq:shrinker}
 \Ric_{\bar g}+\bar\nabla^2\bar f=\frac12\bar g,
 \qquad
 \bar R+\abs{\bar\nabla\bar f}^2=\bar f,
\end{equation}
let \(\varphi_\tau\) be the complete flow
\begin{equation}\label{eq:radial-flow-convention}
 \partial_\tau\varphi_\tau
 =\bar\nabla\bar f\circ\varphi_\tau,
 \qquad
 \varphi_0=\operatorname{Id},
\end{equation}
and recall the Gaussian conventions fixed before Theorem~B:
\[
 d\nu=(4\pi)^{-2}e^{-\bar f}\,dV_{\bar g},
 \qquad
 L^2_\nu=L^2(M;S^2T^*M,d\nu).
\]
All Hilbert and Sobolev spaces in this paragraph consist of real
symmetric two-tensors on \(M\).  Set
\[
 \begin{split}
 H^1_\nu&:=
 \overline{C_c^\infty(M;S^2T^*M)}^{
  \left(\|u\|_{L^2_\nu}^2+
        \|\bar\nabla u\|_{L^2_\nu}^2\right)^{1/2}},
 \qquad H^{-1}_\nu:=(H^1_\nu)^* .
 \end{split}
\]
For \(h\in C_c^\infty(M;S^2T^*M)\), define
\[
 \A_0h=\bar\Delta_{\bar f}h+2\overline{\Rm}(h),\qquad
 \bar\Delta_{\bar f}h
 =\bar\Delta h-\bar\nabla_{\bar\nabla\bar f}h,
 \qquad
 (\overline{\Rm}(h))_{ij}
 =\bar R_{ikj\ell}h^{k\ell},
\]
and write
\(\overline{\Rm}(u,v)
 :=\langle\overline{\Rm}(u),v\rangle_{\bar g}\).
Our component convention is
\(R_{ijkl}=\langle R(e_j,e_i)e_k,e_l\rangle\); the order
\(ikj\ell\) in the displayed action is therefore essential and gives
\(\overline{\Rm}(\bar g)=\overline{\Ric}\).
The boundedness of the FIK curvature makes
\[
 \mathfrak a[u,v]
 :=-\int_M\langle\bar\nabla u,\bar\nabla v\rangle_{\bar g}\,d\nu
   +2\int_M\overline{\Rm}(u,v)\,d\nu
\]
a continuous symmetric form on \(H^1_\nu\).  Choose
\(C_{\rm Fr}\) so that
\[
 \mathfrak l(u,v)
 :=-\mathfrak a[u,v]
   +C_{\rm Fr}\langle u,v\rangle_{L^2_\nu}
\]
is closed and coercive, let \(L_0\) be its Friedrichs operator, and
define
\[
 \A:=C_{\rm Fr}I-L_0,\qquad D(\A):=D(L_0).
\]
Every spectral projection below refers to this self-adjoint
realization.  For an \(H^1_\nu\) argument, \(\A\) denotes the
associated \(H^1_\nu\)-to-\(H^{-1}_\nu\) form operator.  The
FIK-operator section records its domain and sign conventions in full.
The spectral theorem below gives
\[
 \mathcal Z
 :=\mathbf 1_{[0,\infty)}(\A)L^2_\nu
 =\operatorname{span}\{Z_0,\ldots,Z_8\}.
\]
Retain the \(L^2_\nu\)-orthonormal basis fixed before Theorem~B, in
which \(Z_0\) is the scale direction and \(Z_1,\ldots,Z_8\) are the
geometric Lie-derivative directions.  Its chosen orientation satisfies
\[
 Z_0=c_{\mathrm{Ric}}\Ric_{\bar g},\qquad
 c_{\mathrm{Ric}}>0,
\]
where \(c_{\mathrm{Ric}}\) is the \(L^2_\nu\)-normalizing constant.
Retain the real spectral potentials
\(\phi_1,\ldots,\phi_8\) and nonzero normalizations
\(\mathfrak c_1,\ldots,\mathfrak c_8\) fixed above, so that
\[
 Z_j=\mathfrak c_j\bar\nabla^2\phi_j,\qquad
 W_j:=\frac{\mathfrak c_j}{2}\bar\nabla\phi_j,\qquad
 Z_j=\Lie_{W_j}\bar g,\quad1\leq j\leq8.
\]
These potentials, constants, and vector fields are fixed before any
primitive Gram datum is selected.  Thus the nine moments and the
geometric generators appearing in the prepared slice below are defined
before they enter that construction.  Their completeness follows from the
linear-growth estimate \eqref{eq:W-all-order} with \(m=0\) on the
complete manifold \((M,\bar g)\); hence the generator flows used below
are global.
The spectral and form constructions in
Theorem~\ref{thm:FIK-spectrum} prove that the restriction of \(\A\) to
\(\mathcal Z^\perp\) has a strictly negative gap.  Fix throughout a
number \(\beta\) satisfying
\begin{equation}\label{eq:beta}
 0<\beta<\min\left\{\frac12,\gamma_1\right\},\qquad
 \mathfrak a[u,u]\leq-\beta\|u\|_{L^2_\nu}^2
 \quad\text{for }u\in H^1_\nu,\ u\perp\mathcal Z .
\end{equation}
Thus \(\beta\) is fixed before the entrance rates in
Theorem~\ref{thm:intro-sharp-scattering} are
chosen; it is not an input varying over the prepared ball.
Fix once and for all a smooth cutoff
\(\rho:[0,\infty)\to[0,1]\) with
\[
 \rho=1\ \hbox{on }[0,1],\qquad
 \rho=0\ \hbox{on }[2,\infty),
 \qquad
 \rho_\tau(x)=\rho(e^{-\tau}\bar f(x)).
\]
For every \(\sigma>0\), set
\begin{equation}\label{eq:omega-sigma}
 \omega_\sigma(\tau,x)
 =\min\left\{e^{-\sigma\tau}(1+\bar f(x))^\sigma,1\right\}.
\end{equation}
In every prepared H\"older statement below the exponent satisfies
\(0<\alpha<1\), and an unqualified prepared order satisfies \(r\geq3\).
The entrance time \(\tau_0\) is always a fixed finite normalized time.

Fix once and for all a number \(Q_{\rm har}>2\).  If \(g\) is a complete
\(C^{2,\alpha}\) metric on an \(n\)-manifold, a \(g\)-harmonic coordinate
chart of radius \(r>0\) at \(x\) is a diffeomorphism
\[
 u:B_g(x,r)\longrightarrow\Omega\subset\mathbb R^n,\qquad u(x)=0.
\]
No roundness of the coordinate image \(\Omega\) is required.  Put
\(\Omega_r=r^{-1}\Omega\), let
\(\psi_r(\xi)=u^{-1}(r\xi)\), and put
\(\widetilde g=r^{-2}\psi_r^*g\) on \(\Omega_r\).  We call the chart
\emph{admissible} when
\begin{equation}\label{eq:fixed-harmonic-radius-convention}
 Q_{\rm har}^{-1}\delta<\widetilde g<Q_{\rm har}\delta,
 \qquad
 \|\widetilde g\|_{C^{2,\alpha}(\Omega_r)}<Q_{\rm har},
\end{equation}
where the matrix inequalities and the norm inequality are uniform.
Define \(r_{\rm har}(g,x)\) to be the supremum of the radii \(r>0\)
for which there exists an admissible \(g\)-harmonic coordinate chart of
radius \(r\) at \(x\).  Thus the fixed ellipticity and
\(C^{2,\alpha}\) bounds in
\eqref{eq:fixed-harmonic-radius-convention} are part of the defining
existence condition, not an auxiliary property imposed afterward.
This metric-ball-domain convention is downward closed: restricting
\(u\) to \(B_g(x,r')\), \(0<r'<r\), is again an admissible harmonic
coordinate chart after the evident rescaling, although its Euclidean
image need not be a ball.

For quantitative perturbation arguments we retain the witness, not
merely the numerical value of the supremum.  An admissible chart has
\emph{coefficient reserve \(q\in(0,Q_{\rm har}-1)\)} when
\begin{equation}\label{eq:buffered-harmonic-chart-reserve}
 (Q_{\rm har}-q)^{-1}\delta
 <\widetilde g<(Q_{\rm har}-q)\delta,\qquad
 \|\widetilde g\|_{C^{2,\alpha}(\Omega_r)}
 <Q_{\rm har}-q.
\end{equation}
It has a \emph{domain buffer \(\zeta\in(0,1)\) at radius \(r\)} when it
is the restriction of such a reserved chart defined on a larger
metric ball.  If \(\Omega_+\) is the outer coordinate image after
normalization by \(r\), the buffer includes a bounded
\(C^{3,\alpha}\) Euclidean domain \(D\) satisfying
\[
 r^{-1}u\!\left(\overline{B_g(x,r)}\right)
 \Subset D\Subset\Omega_+,
\]
\[
 \operatorname{dist}_{\rm Euc}
 \!\left(r^{-1}u(\overline{B_g(x,r)}),\partial D\right)\geq\zeta,
 \qquad
 \operatorname{dist}_{\rm Euc}(D,\partial\Omega_+)\geq\zeta.
\]
We also require
\(\operatorname{diam}_{\rm Euc}D\leq\zeta^{-1}\) and a boundary atlas
of at most \(\lceil\zeta^{-1}\rceil\) charts, each of radius at least
\(\zeta\) and \(C^{3,\alpha}\)-character at most \(\zeta^{-1}\).
This is the quantitative domain on which the parameterized harmonic
Dirichlet problem is solved.  A
family has a buffered harmonic witness when the same \(q,\zeta>0\)
work at every center under consideration.
The fixed convention \eqref{eq:fixed-harmonic-radius-convention},
together with these explicitly recorded reserves whenever openness is
used, governs every occurrence of \(r_{\rm har}\) in the paper.

We also fix the singular-value convention used in the quantitative
faces.  If \(A\) is a real matrix, then
\[
 s_{\min}(A):=\inf_{|v|=1}|Av|
\]
in the displayed Euclidean bases; in particular, for a Gram matrix this
means the least singular value, not an unspecified eigenvalue.  If
\(F:(M,\acute G)\to(M,S)\) is a prepared map, then
\begin{equation}\label{eq:bundle-map-smallest-singular-value}
 s_{\min}(dF):=
 \inf_{x\in M}\inf_{v\ne0}
 \frac{|dF_xv|_{S,F(x)}}{|v|_{\acute G,x}},
\end{equation}
and \(s_{\min}(dF^{-1})\) is defined with the source and target metrics
reversed.  More explicitly, if \(A=[dF_x]\) and
\(\mathbf G_x,\mathbf S_{F(x)}\) are the coefficient matrices of
\(\acute G\) and \(S\) in a source--target coordinate pair, then the
coordinate representative whose Euclidean singular values equal the
intrinsic ones is the metric-weighted matrix
\[
 \mathbf S_{F(x)}^{1/2}A\mathbf G_x^{-1/2}.
\]
The unweighted matrix \(A\) is not asserted to have the same singular
values.  On every fixed common-margin scale-one atlas, ellipticity gives
uniform constants \(0<c_{\rm sv}\leq C_{\rm sv}<\infty\) such that
\[
 c_{\rm sv}s_{\min}(A)
 \leq
 s_{\min}\!\left(
   \mathbf S_{F(x)}^{1/2}A\mathbf G_x^{-1/2}\right)
 \leq C_{\rm sv}s_{\min}(A).
\]
Thus simultaneous parabolic rescaling of the source and target metrics
does not change the intrinsic quantity, while passage to an unweighted
coordinate Jacobian is only a uniformly controlled comparison.
For later first-exit margins, set
\begin{equation}\label{eq:typed-map-lower-margins}
 \mathfrak m_F(F;\acute G,S)
 :=\min\{s_{\min}(dF),s_{\min}(dF^{-1})\}.
\end{equation}
The coordinate-free scale-adapted functional for the relative marking
\(R\) is fixed below, immediately after its source and range metric
\(\widehat g_{\rm la}\) has been defined.

Fix once and for all a purely geometric number
\(\Gamma_{\rm atl}\geq1\) beyond which the quantitative AC charts of
Lemma~\ref{lem:FIK-AC-symbol}, their fixed buffers, and the compact
core cover are available.  Enlarge it once, still using only the fixed
background and cutoff profiles, past every exact-center cutoff-nesting
and support-separation threshold.  This number is fixed once and for all
as the base of the pre-radius atlas below and is never enlarged after
that atlas or any of its norms has been defined.  A distinct,
package-dependent lower
bound \(\Gamma_{\rm pre}\geq\Gamma_{\rm atl}\) will be produced by
Lemma~\ref{lem:pre-radius-low-order-closure}.

To type the primitive map ceiling before it is selected, fix now the
pre-radius core-plus-dyadic atlas.  Put
\begin{equation}\label{eq:pre-radius-dyadic-family}
 \mathscr L_{\rm pre}
 :=\{2^q\Gamma_{\rm atl}:q=0,1,2,\ldots\},
 \qquad
 A_L^{\rm pre}:=\{L/2<\bar f<4L\}
 \quad(L\in\mathscr L_{\rm pre}),
\end{equation}
and use the fixed compact core cover together with the uniformly
buffered scale-one AC covers of these annuli.  If
\((U,\kappa_U)\) and \((V,\kappa_V)\) are two members whose domains
overlap, rescale their source and target metrics by their recorded
radial levels.  For \(q\in\{6,14\}\), define the finite background
identity ceiling
\begin{equation}\label{eq:primitive-identity-map-ceilings}
 K_{{\rm Id},q}^{\rm pre}
 :=\sup_{U\cap V\ne\varnothing}
 \left(
  \|\kappa_V\circ\kappa_U^{-1}\|_{C^{q,\alpha}}
  +\|\kappa_U\circ\kappa_V^{-1}\|_{C^{q,\alpha}}
 \right)<\infty .
\end{equation}
The supremum uses the normalized overlap domains.  Its finiteness is
part of the fixed bounded-geometry and buffered-transition certificate;
it is independent of the eventual graft radius \(\Gamma\).

For later support faces, first define the scalar functional
\begin{equation}\label{eq:generic-support-separation-functional}
 \mathscr S_{\rm sep}(\tau;\Psi,A)
 :=\inf_{x\in A}\bigl(e^{-\tau}\bar f(\Psi(x))-2\bigr)
\end{equation}
for every map \(\Psi\) and set \(A\) for which the displayed infimum is
defined.  Thus support separation is a quantitative scalar condition,
not merely a qualitative disjointness assertion.

Using the complete geometric generators \(W_1,\ldots,W_8\) fixed in
the opening spectral convention above, set
\[
 Y_0:=2\Ric_{\bar g}
      =\bar g-\Lie_{\bar\nabla\bar f}\bar g,
 \qquad
 Y_j:=\Lie_{W_j}\bar g,\quad1\leq j\leq8.
\]
Fix the resulting reference background Gram matrix and its least
singular value:
\begin{equation}\label{eq:background-Gram-reserve}
 \begin{aligned}
  \mathbf G^{\rm bg}
  &:=\bigl(\langle Y_j,Z_\mu\rangle_{L^2_\nu}\bigr)_{\mu,j=0}^{8}\\
  &=\operatorname{diag}
    \left(\frac{2}{c_{\mathrm{Ric}}},1,\ldots,1\right),\\
  \gamma_{\rm Gram}^{\rm bg}
  &:=s_{\min}(\mathbf G^{\rm bg})
    =\min\left\{\frac{2}{c_{\mathrm{Ric}}},1\right\}>0 .
 \end{aligned}
\end{equation}
Thus the reserve is established explicitly from the already fixed
normalizations, before the primitive package is formed.  This is the
fixed matrix, rather than a moving cutoff matrix, relative to which
every primitive Gram margin is selected.

Define the primitive pre-radius datum by
\begin{equation}\label{eq:pre-radius-primitive-tuple}
 \begin{aligned}
 \mathfrak P_{\rm pre}^{\rm prim}:=
 (&\Lambda_{\rm ell},\Lambda_{\rm coef},\Lambda_{\rm map},
   \Lambda_{R,14}^{\rm pre},\Lambda_{F,6}^{\rm pre},
   \mu_R^{\rm pre},\mu_F^{\rm pre},K_{\rm gr},C_{\rm sc},\\
  &c_{\rm rad},C_{\rm rad},
   \kappa_{\rm sep},\kappa_{\rm Gram},
   \kappa_{\rm map},\kappa_{\rm har}).
 \end{aligned}
\end{equation}
The primitive tuple is admissible only when all of its entries have the
following finite ranges and compatibility:
\begin{equation}\label{eq:primitive-upper-slack-admissibility}
 \begin{gathered}
  1<\Lambda_{\rm ell}<\infty,\qquad
  0<\Lambda_{\rm coef},K_{\rm gr},C_{\rm sc}<\infty,\\
  0<2\mu_R^{\rm pre}<\Lambda_{R,14}^{\rm pre}
       <\Lambda_{\rm map}<\infty,\qquad
  0<2\mu_F^{\rm pre}<\Lambda_{F,6}^{\rm pre}
       <\Lambda_{\rm map}<\infty,\\
  0<c_{\rm rad}<1<C_{\rm rad}<\infty,\qquad
  0<\kappa_{\rm sep},\kappa_{\rm har}<\infty .
 \end{gathered}
\end{equation}
At the two map orders used in the primitive package, require explicitly
\begin{equation}\label{eq:primitive-identity-map-admissibility}
 \max\{K_{{\rm Id},6}^{\rm pre},K_{{\rm Id},14}^{\rm pre}\}
 <\Lambda_{\rm map}.
\end{equation}
Thus the atlas, domains, scales, orders, and finite identity bounds are
fixed components of the admissibility data.  The Gram margin is selected
against the fixed background reserve:
\begin{equation}\label{eq:primitive-Gram-admissibility}
 0<\kappa_{\rm Gram}<\gamma_{\rm Gram}^{\rm bg}.
\end{equation}
The map margin is normalized against the exact metric-isometry value:
\begin{equation}\label{eq:primitive-map-admissibility}
 0<\kappa_{\rm map}<1.
\end{equation}
Thus
\(\mu_{\rm Gram}^{\rm bg}:=
\gamma_{\rm Gram}^{\rm bg}-\kappa_{\rm Gram}>0\)
is a deterministic background reserve, not an additional package
choice.
Here \(C_{\rm ph}^{\rm pre}:=2C_{\rm sc}\).  For later use define the
following deterministic abbreviations from the
primitive radial constants:
\begin{equation}\label{eq:pre-radius-derived-radial-constants}
 \begin{gathered}
  c_R^\sharp:=\frac34c_{\rm rad},
  \qquad C_R^\sharp:=\frac43C_{\rm rad},\\
 c_{\widehat\Phi}^\sharp:=
   \frac34\min\{c_{\rm rad}^2,C_{\rm rad}^{-2}\},
  \qquad
 C_{\widehat\Phi}^\sharp:=
   \frac43\max\{C_{\rm rad}^2,c_{\rm rad}^{-2}\},\\
 c_F^\sharp:=c_R^\sharp c_{\widehat\Phi}^\sharp,
 \qquad C_F^\sharp:=C_R^\sharp C_{\widehat\Phi}^\sharp .
 \end{gathered}
\end{equation}
These constants are determined by the primitive radial constants.  Once
Lemma~\ref{lem:pre-radius-low-order-closure} has been applied to this
tuple after fixing \(0<\alpha<1\) and
\(0<\sigma<\theta<\beta\), put
\begin{equation}\label{eq:pre-radius-input-record}
 \mathfrak I_{\rm pre}:=
 (\alpha,\sigma,\theta,\mathfrak P_{\rm pre}^{\rm prim}).
\end{equation}
Its complete selected output is, by definition,
\begin{equation}\label{eq:pre-radius-derived-tuple}
 \begin{aligned}
 \mathfrak P_{\rm pre}^{\rm der}
 =\mathfrak S_{\rm pre}(\mathfrak I_{\rm pre}):=
 (&\Gamma_{\rm pre},\delta_{\rm atl}^{\rm pre},\delta_{\rm pre},
   \varepsilon_{\rm pre},\tau_{\rm pre},
   (C_{\rm ad}^{\rm pre}(m))_{m=0}^{5},
   c_{\rm ann}^{\rm pre},C_{\rm ann}^{\rm pre},C_{\rm gr}^{\rm pre},\\
  &C_{\mathcal E}^{\rm pre},c_{\mathcal E}^{\rm pre},
   \Lambda_{\rm C6}^{\rm pre},C_{\rm map,6}^{\rm pre},
   (K_{\mathcal Y,m}^{\rm pre})_{m=0}^{4},
   K_{\rm fb}^{\rm pre}).
 \end{aligned}
\end{equation}
Equations~\eqref{eq:pre-radius-input-record}--%
\eqref{eq:pre-radius-derived-tuple} specify the shorthand and complete
dependency list for the output furnished by
Lemma~\ref{lem:pre-radius-low-order-closure}.  No uniformity is asserted
as \(\sigma\downarrow0\), \(\theta-\sigma\downarrow0\), or
\(\theta\uparrow\beta\).

Define also the two radius functionals.  For
\(\Lambda,K>0\), let
\[
 \mathscr A_{\rm K}(\Lambda,K):=
 \left\{C>0:\ \exists\,\delta,c>0\ \hbox{for which
 Lemma~\ref{lem:uniform-two-state-Kato-ledger} holds with }
 (\Lambda,K,\delta,C,c)\right\}.
\]
With
\(\mathbf K_{\mathcal Y}:=(K_{\mathcal Y,m})_{m=0}^{4}\), and with the
FIK background \((\bar g,\bar f,\beta)\) fixed, define
\begin{align*}
 \mathscr D_{\rm 3reg}
 &:=
 \bigl\{(\sigma,\theta,C_{\mathcal E},c_{\mathcal E},C_{\rm gr},
          \mathbf K_{\mathcal Y},K_{\rm fb}):
        0<\sigma<\theta<\beta,\
        C_{\mathcal E},c_{\mathcal E},C_{\rm gr},K_{\rm fb}>0,\
        \mathbf K_{\mathcal Y}\in(0,\infty)^5\bigr\},\\
 \mathscr D_{\rm 2st}^{\rm adm}
 &:=
 \bigl\{(\theta_*,\Lambda_{\rm ell}^{(2)},
          K_{\mathcal Y,0}^{(2)},C_{\rm K}^{(2)}):
        0<\theta_*<\beta,\
         \Lambda_{\rm ell}^{(2)},K_{\mathcal Y,0}^{(2)}>0,\
         C_{\rm K}^{(2)}
         \in\mathscr A_{\rm K}
           (\Lambda_{\rm ell}^{(2)},K_{\mathcal Y,0}^{(2)})\bigr\}.
\end{align*}
Put
\(\mathscr D_{\rm 2st}:=\mathscr D_{\rm 2st}^{\rm adm}\), and define
\begin{equation}\label{eq:typed-radius-threshold-functionals}
 \mathfrak G_{\rm 3reg}:\mathscr D_{\rm 3reg}\longrightarrow[1,\infty),
 \qquad
 \mathfrak G_{\rm 2st}:\mathscr D_{\rm 2st}\longrightarrow[1,\infty).
\end{equation}
Lemma~\ref{lem:uniform-two-state-Kato-ledger} supplies an admissible
value of \(C_{\rm K}^{(2)}\) from
\((\Lambda_{\rm ell}^{(2)},K_{\mathcal Y,0}^{(2)})\), and
\(\mathfrak G_{\rm 2st}\) is evaluated only at such a value.
Theorem~\ref{thm:robust-modulated-three-region} and
Lemma~\ref{lem:uniform-two-state-package-radius} establish the
finiteness of \(\mathfrak G_{\rm 3reg}\) and
\(\mathfrak G_{\rm 2st}\), respectively, with precisely the displayed
arguments.  Once these arguments are fixed, set
\begin{align}
 \overline\Gamma_{\rm 3reg}
 &:=
 \mathfrak G_{\rm 3reg}
   (\sigma,\theta,C_{\mathcal E},c_{\mathcal E},C_{\rm gr},\notag\\
 &\hspace{9em}\mathbf K_{\mathcal Y},K_{\rm fb}),
 \label{eq:typed-three-region-threshold-value}\\
 \overline\Gamma_{\rm 2st}
 &:=
 \mathfrak G_{\rm 2st}
   (\theta_*,\Lambda_{\rm ell}^{(2)},
    K_{\mathcal Y,0}^{(2)},C_{\rm K}^{(2)}).
 \label{eq:typed-two-state-threshold-value}
\end{align}

The static clauses of
Lemmas~\ref{lem:uniform-two-state-Kato-ledger} and
\ref{lem:uniform-two-state-package-radius} are proved later, after the
exact two-state coefficient algebra is introduced.  Their proofs are
independent of the exact-core construction, the formation theorem, and
the global comparison conclusion.

The following dependency convention is fixed throughout.  In every
application of Theorem~\ref{thm:intro-sharp-scattering}, first fix
\(0<\alpha<1\) and
\begin{equation}\label{eq:fixed-entrance-rate-pair}
 0<\sigma<\theta<\beta .
\end{equation}
Next fix the primitive pre-radius ellipticity, coefficient, map,
graft, upper-scale, and positive-margin data with strict room,
including the quantitative background Gram reserve
\eqref{eq:primitive-Gram-admissibility}.

Lemma~\ref{lem:pre-radius-low-order-closure} then produces the package threshold
\(\Gamma_{\rm pre}\geq\Gamma_{\rm atl}\) and proves that the derived
low-order column, feedback, and graft-forcing ceilings needed in the
radius choice are uniform for every
\(\Gamma\geq\Gamma_{\rm pre}\); hence neither the fixed atlas nor
those ceilings are reselected when the eventual radius is enlarged.
Now fix the common two-state ellipticity ceiling
\(\Lambda_{\rm ell}^{(2)}\) and the direct-column ceiling
\(K_{\mathcal Y,0}^{(2)}\) strictly above the corresponding derived
 common faces.  Lemma~\ref{lem:uniform-two-state-Kato-ledger}
 furnishes one preliminary admissible
 \((\delta_{\rm K}^{(2)},C_{\rm K}^{(2)},c_{\rm K}^{(2)})\).
 Lemma~\ref{lem:uniform-two-state-package-radius} subsequently reduces
 \(\delta_{\rm K}^{(2)}\) by the uniform radial-absorption thresholds
 to obtain the final package entry \(\delta_{\rm box}^{(2)}\), still
 before the two-state radius value is evaluated.
Reserve the notation
\[
 C_P=C_P(\theta,K_{\rm fb},C_{\E},c_{\E},\bar g)\geq1,
 \qquad
 c_P=c_P(\theta,K_{\rm fb},C_{\E},c_{\E},\bar g)>0
\]
for the endpoint-independent constants produced by
Lemma~\ref{lem:future-phase-tail}.  Their construction uses only the
displayed data and is independent of every barrier radius, package
radius, outer support parameter, and finite endpoint.  Select these
named future-tail constants and the radius-affecting
barrier shape data as in
Theorem~\ref{thm:robust-modulated-three-region}, and put
\(\theta_*:=\theta\).  Evaluate
\[
 (C_P^*,c_P^*):=(C_P,c_P)\big|_{\theta=\theta_*},
\]
using those already frozen one-state feedback and forcing ceilings.
Then evaluate
\(\overline\Gamma_{\rm 3reg}\) by
\eqref{eq:typed-three-region-threshold-value}, and evaluate the
uniform subordinate two-state value
\(\overline\Gamma_{\rm 2st}\) by
\eqref{eq:typed-two-state-threshold-value}, with
\(\theta_*=\theta\).  Choose a graft radius
\begin{equation}\label{eq:global-compatible-package-radius}
 \Gamma_{\rm A}:=
 \max\{\Gamma_{\rm pre},\overline\Gamma_{\rm 3reg}\},
 \qquad
 \Gamma_{\rm B}:=
 \max\{\Gamma_{\rm A},\overline\Gamma_{\rm 2st}\},
 \qquad
 \Gamma\geq\Gamma_{\rm B}.
\end{equation}
Theorem~A alone requires only \(\Gamma\geq\Gamma_{\rm A}\).  In the
unified construction we impose the stronger
\(\Gamma\geq\Gamma_{\rm B}\) from the outset so that the same
exact-core center and basin also support the two-state refinement in
Theorem~B.  This reserves a static numerical ceiling; it does not use
a two-state evolution or difference estimate in the proof of
Theorem~A.
Only then construct the cutoff \(\eta_\Gamma\) and the associated
buffered atlases.  After the implantation scale in the exact-core
construction, one may still choose the lower scale margin
\(c_{\rm sc}\), which does not enter any radius threshold; all
remaining smallness and entrance-time data are then frozen in the order
of Remark~\ref{conv:authoritative-adaptive-order}.  Thus every
admissible rate pair remains available, but no uniformity as
\((\sigma,\theta)\) ranges over the open rate cone is asserted.
Rate-free static statements may use any one such fixed compatible pair
and package.  Dependence on the pair is usually suppressed from the
notation.  An admissible numerical prepared package is the following
ordered collection.  Its entries are chosen in the order specified in
Remark~\ref{conv:authoritative-adaptive-order}, using the named
construction results.  Every uniform prepared statement quantifying
over \(\mathfrak P_{\rm prep}\) refers to one completed choice of these
data:
\begin{equation}\label{eq:numerical-prepared-package}
 \begin{aligned}
 \mathfrak P_{\rm prep}
 :=(&\alpha,\sigma,\theta,
     \mathfrak P_{\rm pre}^{\rm prim},
     \mathfrak P_{\rm pre}^{\rm der},\\
   &
     \Lambda_{\rm ell}^{(2)},K_{\mathcal Y,0}^{(2)},
     C_{\rm K}^{(2)},c_{\rm K}^{(2)},
     K_{\rm J}^{(2)},K_0^{(2)},
     K_{\rm gr},c_{\rm supp}^{(2)},C_{\rm supp}^{(2)},\\
   & \Gamma_{\rm atl},\Gamma,
     c_{\rm rad},C_{\rm rad},
     \kappa_{\rm sep},\kappa_{\rm Gram},
     \kappa_{\rm map},\kappa_{\rm har},\\
   & c_{\rm sc},C_{\rm sc},c_{\rm scl},C_{\rm scl},
     m_{\rm ad},\delta_{\rm c2},C_{\rm raw},K_{h,5}^{\rm pre},
     \delta_{\rm box}^{(2)},
     \varepsilon_{\rm map}^{\rm HM},
     \varepsilon_{\rm hm}^{\rm HM},
     \varepsilon_{\rm ph}^{\rm HM},\\
   & \eta_{\rm core}^{(2)},\eta_{\rm ch}^{(2)},
     \varepsilon_{\rm ph,*}^{(2)},\varepsilon_{\rm ph},
     \varepsilon_{\rm ent},
     C_P^*,c_P^*,
     \tau_{\rm base}^{(2)},\tau_{\rm ad},C_\lambda,\\
   &
     N_{\rm cov},C_{\rm cov},R_{\rm in},
     c_{\rm atl},\ell_{\rm Leb},\Lambda_{\rm atl}).
 \end{aligned}
\end{equation}
The displayed occurrences of
\[
 K_{\rm gr},\ C_{\rm sc},\ c_{\rm rad},\ C_{\rm rad},\
 \kappa_{\rm sep},\ \kappa_{\rm Gram},\
 \kappa_{\rm map},\ \kappa_{\rm har}
\]
denote the corresponding projections of
\(\mathfrak P_{\rm pre}^{\rm prim}\); they are not additional
coordinates.
The first three entries are the fixed H\"older exponent and rate pair
\eqref{eq:fixed-entrance-rate-pair}.  The next two entries are the
complete primitive and derived pre-radius tuples
\eqref{eq:pre-radius-primitive-tuple}--%
\eqref{eq:pre-radius-derived-tuple}; in particular every existential
pre-radius ceiling used later is a frozen package entry.  The following
entries are the common two-state ellipticity, direct-column, Kato, and
absorption constants; the two support
constants are the fixed post-radius graft-annulus bounds.  In
particular,
\[
 1\leq K_{\rm J}^{(2)},K_0^{(2)}<\infty,
 \qquad
 0<c_{\rm supp}^{(2)}\leq C_{\rm supp}^{(2)}<\infty .
\]
The entries
\(K_{\rm gr}\), \(\Gamma_{\rm atl}\), and
\(\Gamma\) record, respectively, the high graft bound, the fixed
geometric atlas base, and the eventual graft radius; the derived tuple
already contains \(\Gamma_{\rm pre}\).  The next two are the common radial-comparison
constants.  They are
chosen admissibly so that
\[
 0<c_{\rm rad}<1<C_{\rm rad},
\]
and \(\Lambda_{\rm map}\) is chosen with positive room for the fixed
background maps at every displayed order in the pre-radius atlas
\eqref{eq:pre-radius-dyadic-family}; uniform atlas equivalence then
transfers it to the eventual \(\Gamma\)-atlas.  The next four are
positive lower margins for support separation, Gram invertibility, map
invertibility, and the dimensionless harmonic-radius functional
defined, on the fully typed prepared state, in
\eqref{eq:dimensionless-harmonic-radius-functional}.  The selected core
and annular chart radii are uniformly comparable with the intrinsic
scale in that definition by radial tracking.  Thus
\(\kappa_{\rm har}\) is scale-free; it is not a positive lower bound
for the unrescaled physical harmonic radius.  The entries
\[
 0<c_{\rm sc}<C_{\rm sc}<\infty,\qquad
 c_{\rm scl}:=\frac12c_{\rm sc},\qquad
 C_{\rm scl}:=2C_{\rm sc}
\]
are respectively the strict entrance-scale constants and the enlarged
dynamic scale bracket.  The upper constant \(C_{\rm sc}\) and
\(K_{\rm gr}\) are fixed before the radius.  The lower constant
\(c_{\rm sc}\) is fixed after \(\Gamma\), and in the exact-core
construction after the implantation scale \(A\), as explicitly allowed
in Remark~\ref{conv:authoritative-adaptive-order}; it enters no radius
threshold.  These constants are all frozen before the subsequent
adaptive smallness and entrance-time thresholds are chosen.  The
subsequent entries
\[
 \begin{aligned}
 &m_{\rm ad}=13,\quad
 \delta_{\rm c2},\quad C_{\rm raw}\geq1,\quad
 K_{h,5}^{\rm pre}\geq1,\\
 &\delta_{\rm box}^{(2)},\quad
 \varepsilon_{\rm map}^{\rm HM},\quad
 \varepsilon_{\rm hm}^{\rm HM},\quad
 \varepsilon_{\rm ph}^{\rm HM},\quad
 \eta_{\rm core}^{(2)},\quad
 \eta_{\rm ch}^{(2)},\quad
 \varepsilon_{\rm ph,*}^{(2)},\quad
 \varepsilon_{\rm ph},\\
 &\varepsilon_{\rm ent},\quad
 C_P^*,\quad c_P^*,\quad
 \tau_{\rm base}^{(2)},\quad
 \tau_{\rm ad},\quad
 C_\lambda=e^{\varepsilon_{\rm ph}}
 \end{aligned}
\]
are the derived harmonic-map, two-state, and specified adaptive
choices.
Lemma~\ref{lem:pre-radius-low-order-closure} constructs the pre-radius
entries.  Lemma~\ref{lem:uniform-two-state-package-radius} constructs
the two-state entries.  Remark~\ref{conv:authoritative-adaptive-order}
fixes their common order.  The final six entries
record the exterior overlap, scale comparison, and interface
separation constants.  They also record the auxiliary buffers, Lebesgue
number, and rescaled atlas bounds.  The remaining types are fixed as
follows: all coefficient, map, covering, comparison,
absorption, and time ceilings denoted by an upper-case
\(C,K,\Lambda,N,\tau\) are finite (with
\(N_{\rm cov}\in\mathbb N\)); all quantities denoted by
\(c,\delta,\varepsilon,\eta,\kappa,\ell_{\rm Leb}\) are strictly
positive in the ranges imposed by their named selection results;
\(\Gamma_{\rm atl},\Gamma\geq1\); and \(R_{\rm in}>0\).  The exact equalities
\[
 m_{\rm ad}=13,\qquad
 c_{\rm scl}=\tfrac12c_{\rm sc},\qquad
 C_{\rm scl}=2C_{\rm sc},\qquad
 C_\lambda=e^{\varepsilon_{\rm ph}}
\]
are part of the record.  Thus the compatible package,
including the adaptive entrance threshold, is fixed after the rate
pair and before the entrance ball and the statement of
Theorem~\ref{thm:intro-sharp-scattering} are
used.

The perturbative harmonic-coordinate modulus uses only the earlier
geometric projection
\begin{equation}\label{eq:reduced-prepared-harmonic-package}
 \mathfrak P_{\rm har}^{\rm geom}
 :=
 \bigl(
  \alpha,Q_{\rm har},\bar g,\bar f,
  \mathfrak P_{\rm pre}^{\rm prim},
  \mathfrak P_{\rm pre}^{\rm der},
  \Lambda_{\rm ell}^{(2)},
  \Gamma_{\rm atl},\Gamma,
  c_{\rm rad},C_{\rm rad},\kappa_{\rm har}
 \bigr).
\end{equation}
This record includes the fixed scale-one core and dyadic atlases and
their coefficient and domain constants through the two pre-radius
tuples.  It is completely fixed before
\(\delta_{\rm c2}\) and the later adaptive smallness data are chosen.
In particular,
\(\mathfrak P_{\rm har}^{\rm geom}\) contains neither
\(\delta_{\rm c2}\) nor any entrance, phase, or time-width threshold.
This is a subtuple of \(\mathfrak P_{\rm prep}\), not an additional
geometric choice.

We next fix the static phase-chart data needed for
Theorem~\ref{thm:intro-sharp-scattering}.  Fix the marked host and
implantation map \(\iota\), the graft cutoff \(\eta\), the polynomial
weight, the scale convention, and the same-output soliton-conjugated
map charts.
Precisely, fix a closed manifold \(\mathcal X\), an open marked set
\(\mathcal X''\subset\mathcal X\), and a diffeomorphism
\[
 \iota:\mathcal X''\longrightarrow M .
 \]
Let \(\mathfrak U_{\iota}^{\rm fix}=\{(U_a,R_a)\}_{a=1}^{A_\iota}\)
be the finite family of named compact marked buffers used below,
together with their fixed scales and normalized source and target
atlases.  Before using the notation, define
\(C_{R_a}^{q,\alpha}(U_a;M)\) to mean the ordinary coordinate
\(C^{q,\alpha}\) map norm after rescaling both the source and target
metrics by \(R_a^{-2}\); for the inverse-map norm the roles of those
atlases are reversed.  For every finite order \(q\) used below, record
once and for all the two-sided, scale-normalized marking-jet bound
\begin{equation}\label{eq:fixed-static-marking-jets}
 \max_{1\leq a\leq A_\iota}
 \left\{
  \|\iota\|_{C_{R_a}^{q,\alpha}(U_a;M)}
  +\|\iota^{-1}\|_{
       C_{R_a}^{q,\alpha}(\iota(U_a);\mathcal X'')}
 \right\}
 \leq K_{\iota,q}^{\rm fix}<\infty
\end{equation}
These are static data of the marked host, not evolving
prepared variables.
For this fixed \(\Gamma\), choose
\(\eta=\eta_\Gamma\in C^\infty(M;[0,1])\) with
\[
 \eta=1\ \hbox{on }\{\bar f\leq2\Gamma/3\},\qquad
 \eta=0\ \hbox{on }\{\bar f\geq5\Gamma/6\},\qquad
 \Omega_\eta:=\{0<\eta<1\}\Subset\{\bar f<\Gamma\}.
\]
For a metric \(G\) on \(\mathcal X\), the tensor \(\iota_*G\) is
globally defined on \(M\).  In particular the product
\(\eta\,\iota_*G\) is unambiguous, and the support compatibility needed
for the graft is automatic.

Fix once and for all a polynomial-loss exponent \(N\geq0\), suppressed
from the notation for the prepared manifold.  For every real \(L>0\)
put
\[
 A_L=\{L/2<\bar f<4L\}.
\]
The dyadic label sets merely select covers; they do not restrict this
continuous-level annulus notation.  Put
\[
 \mathscr L_\Gamma=\{2^q\Gamma:q=0,1,2,\ldots\}.
\]
Recall the full-end dyadic family and core cover fixed before the
primitive selection in \eqref{eq:pre-radius-dyadic-family}; together
they are called the \emph{pre-radius atlas}.  This atlas is used only for the low-order
coefficient, map, column, and parabolic estimates whose numerical
ceilings must be chosen before \(\Gamma\).  The
\(\Gamma\)-dependent family \(\mathscr L_\Gamma\) continues to define
the high-order prepared Banach topology.  On their common domain the
two atlases have uniformly bounded scale-one transition maps: every
\(\Gamma\)-collar meets only a bounded number of members of
\(\mathscr L_{\rm pre}\) having labels comparable with \(\Gamma\).
For every integer \(r\geq0\), define the scale-normalized tensor,
vector-field, and tame graph norms
\begin{align}
 \|u\|_{\mathfrak C_{{\rm sc},N}^{r,\alpha}}
 &:=
 \|u\|_{C^{r,\alpha}(\{\bar f<4\Gamma\},\bar g)}
 +\sup_{L\in\mathscr L_\Gamma}
   L^{-N}\|L^{-1}u\|_{C^{r,\alpha}(A_L,L^{-1}\bar g)},
 \label{eq:scaled-tensor-holder}\\
 \|V\|_{\mathfrak X_{\rm sc}^{r+1,\alpha}}
 &:=
 \|V\|_{C^{r+1,\alpha}(\{\bar f<4\Gamma\},\bar g)}
 +\sup_{L\in\mathscr L_\Gamma}
   \|V\|_{C^{r+1,\alpha}(A_L,L^{-1}\bar g)},
 \label{eq:scaled-vector-holder}\\
 \|u\|_{\mathfrak T_{{\rm sc},N}^{r,\alpha}}
 &:=
 \|u\|_{\mathfrak C_{{\rm sc},0}^{2,\alpha}}
 +\|u\|_{\mathfrak C_{{\rm sc},N}^{r,\alpha}} .
 \label{eq:scaled-tame-holder}
\end{align}
All three spaces mean the corresponding little-H\"older completions
in the fixed compact atlas and the complete dyadic atlas.
For every eventual radius \(\Gamma\), the low-order activation is
measured first in the radius-independent pre-atlas norm
\begin{equation}\label{eq:pre-radius-tensor-holder}
 \begin{split}
 \|u\|_{\mathfrak C_{{\rm pre},N}^{r,\alpha}}
 :=\;&
 \|u\|_{C^{r,\alpha}
   (\{\bar f<4\Gamma_{\rm atl}\},\bar g)}\\
 &+\sup_{L\in\mathscr L_{\rm pre}}
 L^{-N}\|L^{-1}u\|_
 {C^{r,\alpha}(A_L,L^{-1}\bar g)} .
 \end{split}
\end{equation}
The bounded-overlap and scale-one transition estimates between
\(\mathscr L_{\rm pre}\) and \(\mathscr L_\Gamma\) give the
\(\Gamma\)-uniform comparison used in the low-order activation,
\[
 \|u\|_{\mathfrak C_{{\rm sc},0}^{r,\alpha}}
 \leq C\|u\|_{\mathfrak C_{{\rm pre},0}^{r,\alpha}}.
\]
For each fixed \(\Gamma\), the reverse comparison, and the corresponding
comparisons with \(N>0\), also hold, but their constants may depend on
\(\Gamma\).  No uniform two-sided comparison with the expanding compact
piece is asserted.  Crucially, the norm in
\eqref{eq:pre-radius-tensor-holder} is fixed before \(\Gamma\).

Put
\[
 r_{\rm la}=(1+\bar f)^{1/2},\qquad
 \widehat g_{\rm la}=r_{\rm la}^{-2}\bar g .
\]
For a degree-one proper diffeomorphism \(R:M\to M\), define
\begin{align}
 s_{\min}^{\rm la}(dR)
 &:=
 \inf_{x\in M}\inf_{v\ne0}
 \frac{|dR_xv|_{\widehat g_{\rm la},R(x)}}
      {|v|_{\widehat g_{\rm la},x}},
 \label{eq:relative-marking-intrinsic-singular-value}\\
 \mathfrak m_R(R)
 &:=
 \min\{s_{\min}^{\rm la}(dR),
        s_{\min}^{\rm la}(dR^{-1})\}.
 \label{eq:relative-marking-lower-margin}
\end{align}
This is a coordinate-free dimensionless scalar and
\(\mathfrak m_R(\operatorname{Id})=1\).  In the fixed buffered core
charts and the dyadic source--range chart pairs of
\eqref{eq:pre-radius-dyadic-family}, the coefficient metrics are
uniformly elliptic relative to \(\widehat g_{\rm la}\).  Hence
\(\mathfrak m_R\) is uniformly comparable, with constants determined
only by the fixed atlas and radial-comparison package, to the infimum
of the unweighted coordinate singular values of \(dR\) and
\(dR^{-1}\).  It is continuous in the scale-adapted \(C^1\)
right-translated map charts and controls both local-invertibility faces
simultaneously.
The all-order AC symbol bounds imply that
\((M,\widehat g_{\rm la})\) is complete, has bounded geometry of every
fixed order, and has a positive injectivity-radius lower bound.  Indeed,
on the end the logarithmic radial coordinate has uniformly bounded
geometry, while the remaining compact set is harmless.  Fix
\(0<\varepsilon_{\rm la}\) below one quarter of that injectivity-radius
bound and use the \(\widehat g_{\rm la}\)-exponential to define the
scale-adapted smooth local addition
\begin{equation}\label{eq:scale-adapted-local-addition}
 \operatorname{Exp}:\mathscr U_{\rm la}\longrightarrow M,\qquad
 \mathscr U_{\rm la}
 =\{(x,v)\in TM:|v|_{\widehat g_{\rm la}}<
                    \varepsilon_{\rm la}\}.
\end{equation}
Equivalently, its fiber over \(x\) contains precisely the fixed
scale-relative ball
\[
|v|_{\bar g}<\varepsilon_{\rm la}r_{\rm la}(x).
\]
One has \(\operatorname{Exp}(0_x)=x\), and
\((\pi,\operatorname{Exp})\) is a diffeomorphism onto a
\(\widehat g_{\rm la}\)-uniform neighborhood of the diagonal.  On the
fixed core and on every scale-one annulus, the coordinate
representatives of this addition, its fiber inverse, and all their
derivatives through any prescribed finite order have one uniform
bound.  Thus \(\operatorname{Exp}(X)\) below denotes the globally
defined map \(x\mapsto\operatorname{Exp}_x(X(x))\) whenever the scaled
vector norm is sufficiently small; no unscaled uniform bound on
\(|X|_{\bar g}\) is required.

For a \(C^{r+1,\alpha}\) degree-one proper diffeomorphism
\(\psi:M\to M\), say first that \(\psi\) is
\emph{scale-radially admissible} if there are constants
\(0<c_\psi\leq C_\psi<\infty\) such that
\[
 c_\psi(1+\bar f)
 \leq 1+\bar f\circ\vartheta
 \leq C_\psi(1+\bar f),
 \qquad \vartheta\in\{\psi,\psi^{-1}\}.
\]
For such a map, call a buffered source--range chart pair
\(\psi\)-active if the image of the source chart meets the range
chart; use the analogous active pairs for \(\psi^{-1}\).  Define
\(\|\psi\|_{\operatorname{Map}_{\rm sc}^{r+1,\alpha}}\) to be the
supremum, over all active pairs in the fixed core cover and the
pre-radius family \eqref{eq:pre-radius-dyadic-family}, of the ordinary
\(C^{r+1,\alpha}\) coordinate norms of the restrictions
\[
 \psi:U\cap\psi^{-1}(V)\longrightarrow V,
 \qquad
 \psi^{-1}:V\cap\psi(U)\longrightarrow U,
\]
with the domain and range metrics rescaled by their respective
\(\bar f\)-levels.  The two-sided radial comparison implies that, for
each source level, active pairs involve only a uniformly bounded
number of adjacent dyadic range levels, depending on
\(c_\psi,C_\psi\) and the fixed atlas.  Thus the chart-pair
prescription is well-defined and has uniformly bounded multiplicity on
every fixed radial-comparison class.  If \(\psi\) is not
scale-radially admissible, set
\[
 \|\psi\|_{\operatorname{Map}_{\rm sc}^{r+1,\alpha}}:=+\infty.
\]
For \(q\in\{6,14\}\), this definition and
\eqref{eq:primitive-identity-map-ceilings} give
\[
 \|\operatorname{Id}_M\|_
      {\operatorname{Map}_{\rm sc}^{q,\alpha}}
 \leq K_{{\rm Id},q}^{\rm pre}<\Lambda_{\rm map}
\]
by \eqref{eq:primitive-identity-map-admissibility}.  Thus the identity
belongs to the fixed-package common-margin locus at every map order
used in the primitive closure.
At order \(r\), let
\[
 \operatorname{CtrlDiff}_{\rm prop}^{\,1}
   (M;c_{\rm rad},C_{\rm rad},\Lambda_{\rm map})
\]
denote the fixed-package \emph{controlled class} of degree-one proper
diffeomorphisms satisfying
\begin{equation}\label{eq:prepared-map-control-class}
 \begin{gathered}
  c_{\rm rad}(1+\bar f)
  \leq 1+\bar f\circ\vartheta
  \leq C_{\rm rad}(1+\bar f),
  \qquad \vartheta\in\{\psi,\psi^{-1}\},\\
  \|\psi\|_{\operatorname{Map}_{\rm sc}^{r+1,\alpha}}
  \leq\Lambda_{\rm map}.
 \end{gathered}
\end{equation}
This fixed non-strict sublevel is used for uniform estimates and is not
asserted to be open.  The ambient map manifold is
\[
 \operatorname{Diff}_{\rm prop,sc}^{\,1}(M)
 :=
 \bigcup_{\substack{0<c<C<\infty\\ \Lambda<\infty}}
 \operatorname{CtrlDiff}_{\rm prop}^{\,1}(M;c,C,\Lambda),
\]
with its right-translated little-H\"older charts.  For fixed package
constants, the notation
\[
 \operatorname{Diff}_{\rm prop}^{\,1}
   (M;c_{\rm rad},C_{\rm rad},\Lambda_{\rm map})
\]
means the interior of the corresponding controlled class in this
ambient manifold.  Equivalently, every member has a uniform
common-margin number \(\delta_\psi>0\) for which
\begin{equation}\label{eq:prepared-map-class}
 \begin{gathered}
  (c_{\rm rad}+\delta_\psi)(1+\bar f)
  \leq 1+\bar f\circ\vartheta
  \leq (C_{\rm rad}-\delta_\psi)(1+\bar f),
  \qquad \vartheta\in\{\psi,\psi^{-1}\},\\
  \|\psi\|_{\operatorname{Map}_{\rm sc}^{r+1,\alpha}}
  \leq\Lambda_{\rm map}-\delta_\psi,
 \end{gathered}
\end{equation}
and the corresponding scale-one local-invertibility margin is
positive.  The latter is equivalently a positive lower margin for the
Jacobians of \(\psi\) and \(\psi^{-1}\) in the buffered source and
target charts; it is already quantitatively controlled by the
map-and-inverse norm but is recorded to specify the open locus.
When the constants are those in
\(\mathfrak P_{\rm prep}\), they are suppressed from the notation.
If \(\psi_*\) is a prepared background and \(\psi\) is in its local
chart, the uniquely defined right-translated variables
\[
 \psi\circ\psi_*^{-1}=\operatorname{Exp}(X_\psi),
 \qquad
 \psi_*\circ\psi^{-1}=\operatorname{Exp}(X_{\psi^{-1}})
\]
belong to \(\mathfrak X_{\rm sc}^{r+1,\alpha}\).  On every
common-margin chart, the map-and-inverse norm above is uniformly
equivalent, up to the fixed background bound, to
 \(\|X_\psi\|_{\mathfrak X_{\rm sc}^{r+1,\alpha}}
 +\|X_{\psi^{-1}}\|_{\mathfrak X_{\rm sc}^{r+1,\alpha}}\).
For two maps in such a chart define the typed right-translated local
quasi-distance
\begin{equation}\label{eq:right-translated-map-distance}
 \begin{split}
 d_{\rm rt,sc}^{r+1,\alpha}(\psi,\psi_*):={}&
 \|X^+_{\psi,\psi_*}\|_{\mathfrak X_{\rm sc}^{r+1,\alpha}}
 +\|X^-_{\psi,\psi_*}\|_{\mathfrak X_{\rm sc}^{r+1,\alpha}},\\
 \psi\circ\psi_*^{-1}={}&\operatorname{Exp}
   (X^+_{\psi,\psi_*}),\qquad
 \psi_*\circ\psi^{-1}=\operatorname{Exp}
   (X^-_{\psi,\psi_*}).
 \end{split}
\end{equation}
After shrinking any fixed common-margin chart once, there is a
constant \(C_{\rm rt,\Delta}\geq1\), determined only by the common
scale-normalized map-and-inverse bounds and the fixed local-addition
radius,
such that any three maps \(\psi_1,\psi_2,\psi_3\) in it satisfy
\begin{equation}\label{eq:right-translated-quasi-triangle}
 d_{\rm rt,sc}^{r+1,\alpha}(\psi_1,\psi_3)
 \leq C_{\rm rt,\Delta}\left(
 d_{\rm rt,sc}^{r+1,\alpha}(\psi_1,\psi_2)
 +d_{\rm rt,sc}^{r+1,\alpha}(\psi_2,\psi_3)\right).
\end{equation}
Indeed,
\(\psi_1\circ\psi_3^{-1}
 =(\psi_1\circ\psi_2^{-1})
   \circ(\psi_2\circ\psi_3^{-1})\).
In an ordinary coordinate local addition, the displacement of this
composition is the inner displacement plus the outer displacement
composed with the inner map, with fixed smooth coefficient factors
coming from the exponential chart.  The same-order H\"older
composition and product estimates on the common-margin set therefore
give its \(C_{\rm sc}^{r+1,\alpha}\) norm by the sum of the two input
norms.  This is a tame size estimate; no same-order differentiability
of the composition operator is asserted.  Applying it to the inverse
factorization gives the inverse part of
\eqref{eq:right-translated-quasi-triangle}.  Hence
\(C_{\rm rt,\Delta}\) is uniform over every numerical prepared
package, in particular independent of the entrance time.  Every later
use of a triangle estimate for
\(d_{\rm rt,sc}\) refers to this local quasi-triangle inequality, not
to an assertion that a chart-log norm is a global metric.
In particular, a prepared map \(R\) or \(F\) has a typed map class,
two-sided radial comparison, and scale-normalized control of both the
map and its inverse before the prepared graph is formed.
Retain the complete geometric generators and the ordered basis fixed
in \eqref{eq:background-Gram-reserve}; thus
\[
 Y_0=2\Ric_{\bar g}
      =\bar g-\Lie_{\bar\nabla\bar f}\bar g,
 \qquad
 Y_j=\Lie_{W_j}\bar g,\quad 1\leq j\leq8,
\]
is the ordered geometric basis of \(\mathcal Z\) used below.  Fix also
a cutoff \(\chi:[0,\infty)\to[0,1]\), equal to one on
\([0,5/2]\) and zero on \([3,\infty)\), and put
\(\chi_\tau=\chi(e^{-\tau}\bar f)\).

Fix a prepared background
\(\mathbf z_*=(G_*,\lambda_*,R_*,F_*)\).  An order-\(r\) prepared tuple
has the fully typed components
\[
 \mathbf z=(G,\lambda,R,F),\qquad
 G\in\operatorname{Met}^{r,\alpha}(\mathcal X),\quad
 \lambda>0,\quad R,F\in\operatorname{Diff}_{\rm prop}^{\,1}(M),
\]
where the superscript \(1\) denotes the degree-one proper homotopy
class, not differentiability order.  For every normalized chart time
\(s\geq0\), put
\[
 \Theta(s,\mathbf z):=\varphi_s\circ R,\qquad
 \Phi(s,\mathbf z):=\Theta(s,\mathbf z)\circ F,
\]
and define the time-typed prepared graph by
\begin{equation}\label{eq:intro-prepared-state-graph}
 \acute G(s,\mathbf z)
 =\eta\,\iota_*G+(1-\eta)
       \lambda\Theta(s,\mathbf z)^*\bar g,
 \qquad
 h(s,\mathbf z)
 =\lambda^{-1}(\Phi(s,\mathbf z)^{-1})^*
       \acute G(s,\mathbf z)-\bar g.
\end{equation}
On the fixed chart \(\mathscr P_{\tau_0}^{r,\alpha}\) we suppress the
time argument:
\[
 \Theta:=\Theta(\tau_0,\mathbf z),\quad
 \Phi:=\Phi(\tau_0,\mathbf z),\quad
 \acute G(\mathbf z):=\acute G(\tau_0,\mathbf z),\quad
 h(\mathbf z):=h(\tau_0,\mathbf z).
\]
For a fully typed pair \((s,\mathbf z)\), put
\[
 S(s,\mathbf z):=\lambda\Theta(s,\mathbf z)^*\bar g,\qquad
 r_{\rm sol}(s,\mathbf z;x):=
 \bigl(\lambda[1+\bar f(\Phi(s,\mathbf z)(x))]\bigr)^{1/2},
\]
and define the dimensionless harmonic-radius functional
\begin{equation}\label{eq:dimensionless-harmonic-radius-functional}
 \mathfrak h_{\rm har}(s,\mathbf z):=
 \inf_{x\in M}\min\left\{
 \frac{r_{\rm har}(\acute G(s,\mathbf z),x)}
      {r_{\rm sol}(s,\mathbf z;x)},
 \frac{r_{\rm har}(S(s,\mathbf z),F(x))}
      {r_{\rm sol}(s,\mathbf z;x)}
 \right\}.
\end{equation}
For an evolution carrying
\((G(t),t(\tau),\lambda(\tau),\Theta_\tau,F_\tau)\), set first
\[
 R_\tau:=\varphi_{-\tau}\circ\Theta_\tau,\qquad
 \mathbf z(\tau):=
 \bigl(G(t(\tau)),\lambda(\tau),R_\tau,F_\tau\bigr).
\]
Define the prepared physical metric along this parametrized curve by
\begin{equation}\label{eq:dynamic-physical-prepared-metric}
 \acute G_{{\rm phys},\tau}
 :=
 \eta\,\iota_*G(t(\tau))
 +(1-\eta)\lambda(\tau)\Theta_\tau^*\bar g .
\end{equation}
Along this parametrized curve we use the unambiguous abbreviations
\begin{equation}\label{eq:dynamic-harmonic-and-separation-faces}
 \begin{gathered}
   \Theta_\tau:=\Theta(\tau,\mathbf z(\tau)),\quad
   \Phi_\tau:=\Phi(\tau,\mathbf z(\tau)),\quad
   \acute G_\tau:=\acute G(\tau,\mathbf z(\tau))
        =\acute G_{{\rm phys},\tau},\quad
  S_\tau:=S(\tau,\mathbf z(\tau)),\\
  r_{{\rm sol},\tau}(x):=
       r_{\rm sol}(\tau,\mathbf z(\tau);x),\qquad
  \mathfrak h_{\rm har}(\tau):=
       \mathfrak h_{\rm har}(\tau,\mathbf z(\tau)),\\
  \mathfrak s_{\rm sep}(\tau):=
  \mathscr S_{\rm sep}
   \bigl(\tau;\Phi_\tau,\operatorname{supp}(1-\eta)\bigr).
 \end{gathered}
\end{equation}
Thus the harmonic face \(\mathfrak h_{\rm har}(\tau)\geq\kappa_{\rm har}\)
is exactly the pair of pointwise inequalities
\begin{equation}\label{eq:dimensionless-harmonic-radius-face}
 \frac{r_{\rm har}(\acute G_\tau,x)}{r_{{\rm sol},\tau}(x)}
 \geq\kappa_{\rm har},
 \qquad
 \frac{r_{\rm har}(S_\tau,F_\tau(x))}
      {r_{{\rm sol},\tau}(x)}
 \geq\kappa_{\rm har}.
\end{equation}
These definitions use only the fixed convention
\eqref{eq:fixed-harmonic-radius-convention} and the prepared state; no
eventual-radius atlas enters them.
Put
\[
 \Theta_*:=\varphi_{\tau_0}\circ R_*,
 \qquad \Phi_*:=\Theta_*\circ F_*,
\]
and use the same-output soliton-conjugated relative maps
\begin{equation}\label{eq:intro-uniform-conjugated-map-coordinates}
 \begin{aligned}
  Q_\Theta&:=\Theta\circ\Theta_*^{-1}
   =\varphi_{\tau_0}\circ R\circ R_*^{-1}
      \circ\varphi_{-\tau_0}
   =\operatorname{Exp}(X_\Theta),\\
  Q_\Phi&:=\Phi\circ\Phi_*^{-1}
   =\varphi_{\tau_0}\circ R\circ F\circ F_*^{-1}
      \circ R_*^{-1}\circ\varphi_{-\tau_0}
   =\operatorname{Exp}(X_\Phi).
 \end{aligned}
\end{equation}
Recover
\[
 \Theta=Q_\Theta\circ\Theta_*,\qquad
 \Phi=Q_\Phi\circ\Phi_*,\qquad
 R=\varphi_{-\tau_0}\circ\Theta,\qquad
 F=\Theta^{-1}\circ\Phi .
\]
For each fixed \(\tau_0\), these coordinates are homeomorphically
equivalent, at every fixed finite prepared order, to the raw
right-translated \(R,F\) coordinates and hence define the same prepared
topology and open sets.  The conversions between raw and conjugated
coordinates are asserted to be \(C^1\) only with the derivative buffer
proved in Lemma~\ref{lem:prepared-chart-calculus}.  No
entrance-time-uniform equivalence with the raw product map norm is
asserted.
The Banach model is the product of the independent same-output
coordinates:
\[
 \mathscr E_{\rm prep}^{r,\alpha}
 :=h^{r,\alpha}(S^2T^*\mathcal X)\times\mathbb R\times
   \bigl(\mathfrak X_{\rm sc}^{r+1,\alpha}\bigr)^2.
\]
Writing
\[
 \widehat G:=\lambda_*^{-1}(G-G_*),\qquad
 \ell:=\log(\lambda/\lambda_*),
\]
give a coordinate increment the genuine vector-space norm
\begin{equation}\label{eq:prepared-model-Banach-norm}
 \|(\widehat G,\ell,X_\Theta,X_\Phi)\|_
      {\mathscr E_{\rm prep}^{r,\alpha}}
 :=
 \|\widehat G\|_{C^{r,\alpha}(\mathcal X)}+|\ell|
 +\|X_\Theta\|_{\mathfrak X_{\rm sc}^{r+1,\alpha}}
 +\|X_\Phi\|_{\mathfrak X_{\rm sc}^{r+1,\alpha}} .
\end{equation}
Here and below the norm on the metric factor is the usual
\(C^{r,\alpha}\) norm restricted to the closed little-H\"older space
\(h^{r,\alpha}\).  Every prepared ball, tangent space, and Fr\'echet
derivative is taken in this little-H\"older metric factor.
For finite-difference estimates we retain, rather than discard, the
dependent graph coordinate.  If \(\mathbf z_1,\mathbf z_2\) lie in one
fixed same-output chart with center \(\mathbf z_*\), set
\begin{equation}\label{eq:prepared-Banach-norm}
 \begin{split}
 d_{\rm prep}^{r,\alpha}(\mathbf z_1,\mathbf z_2)
 :=\;&
 \lambda_*^{-1}\|G_1-G_2\|_{C^{r,\alpha}(\mathcal X)}
 +\left|\log\frac{\lambda_1}{\lambda_2}\right|\\
 &+\|X_{\Theta,1}-X_{\Theta,2}\|_
       {\mathfrak X_{\rm sc}^{r+1,\alpha}}
  +\|X_{\Phi,1}-X_{\Phi,2}\|_
       {\mathfrak X_{\rm sc}^{r+1,\alpha}}\\
 &+\|h(\mathbf z_1)-h(\mathbf z_2)\|_
       {\mathfrak T_{{\rm sc},N}^{r,\alpha}} .
 \end{split}
\end{equation}
We continue to abbreviate this graph-augmented distance by
\(
 \|\mathbf z_1-\mathbf z_2\|_{\mathscr X_{\rm prep}^{r,\alpha}}
\) in finite-difference estimates.  It is not a vector-space norm: it
is the product distance after the graph embedding
\(
 \mathbf z\mapsto(\widehat G,\ell,X_\Theta,X_\Phi,h(\mathbf z))
\), and hence it does satisfy the triangle inequality.  At fixed finite
order, continuity of the graph map and the immediate domination of the
independent-coordinate distance show that it induces the same local
topology as \eqref{eq:prepared-model-Banach-norm}; no same-order
Lipschitz equivalence is asserted.  On a bounded coordinate-convex
common-margin order-\((k+2,\alpha)\) subball,
Lemma~\ref{lem:prepared-chart-calculus}
instead supplies the precise buffered comparison
\begin{equation}\label{eq:prepared-buffered-distance-comparison}
 d_{\rm prep}^{k,\alpha}(\mathbf z_1,\mathbf z_2)
 \leq C_{\mathscr B}
 \|\mathbf z_1-\mathbf z_2\|_{\mathscr E_{\rm prep}^{k+2,\alpha}} .
\end{equation}
The state-difference notation on the right means the norm of the
difference of the four independent coordinate vectors in
\eqref{eq:prepared-model-Banach-norm}.
Here and below Banach balls, Banach radii, tangent norms, and every
inverse- or submersion-theorem argument are taken in the independent
model norm \eqref{eq:prepared-model-Banach-norm}; the notation
\(\mathscr X_{\rm prep}\) for two states records the stronger
graph-augmented distance \eqref{eq:prepared-Banach-norm}.
The factor \(\lambda_*^{-1}\) is essential at the compact fixed set of
the soliton flow: a host-metric variation enters the normalized graph
as \(\lambda^{-1}(\Phi^{-1})^*(\eta\,\iota_*\dot G)\).
The conjugated map variables are equally essential: a raw
\(F\)-variation would enter through
\((\varphi_{\tau_0})_*X_F\), whose compact normal component need not be
uniformly bounded as \(\tau_0\) varies.  For each fixed entrance time
the independent model norm is equivalent to the unweighted host-metric
and raw-map product norm, so it changes neither the prepared topology
nor its open sets; its normalization makes the chart constants uniform
when entrance times vary.
The little-H\"older prepared Banach manifold
\(\mathscr P_{\tau_0}^{r,\alpha}\) is the open common-margin locus in
the independent conjugated variables
\((G,\log\lambda,X_\Theta,X_\Phi)\) on which \(G\) is a metric,
\(R,F\) lie in the interiors just defined, the graph identities hold,
and the common support, ellipticity, and bounded-geometry margins are
positive.  Its defining \(C^1\) Banach atlas consists of these
same-output soliton-conjugated local-addition charts, not of the raw
\(R,F\) charts.  Indeed, for two fixed chart centers denoted by
\(a,b\), their relative maps satisfy
\[
 Q_\Theta^{(b)}=Q_\Theta^{(a)}\circ
       (\Theta_a\circ\Theta_b^{-1}),\qquad
 Q_\Phi^{(b)}=Q_\Phi^{(a)}\circ
       (\Phi_a\circ\Phi_b^{-1}).
\]
The right factors are fixed controlled diffeomorphisms; together with
the fixed pointwise local-addition changes, they give same-order
\(C^1\) chart overlaps.  No variable--variable composition or
inversion is differentiated in an atlas transition.  At fixed order
the graph tensor \(h\) is a dependent continuous coordinate; it is
differentiated only with the buffer of
Lemma~\ref{lem:prepared-chart-calculus}, or after the separate scalar
integration-by-parts argument used for the moment map.  The fixed
controlled sublevels \eqref{eq:prepared-map-control-class} are used
only as uniform estimate sets and are not themselves called open.
Inside any compatible nested chart set
\begin{equation}\label{eq:intro-smooth-prepared-points}
 \mathscr P_{\tau_0}^{\infty}
 :=\bigcap_{m\geq3}\mathscr P_{\tau_0}^{m,\alpha}.
\end{equation}
This denotes smooth prepared points with finite scaled weighted norm at
every order for the fixed exponent \(N\); it carries no additional
Fr\'echet-differentiability assertion.

We next define the finite phase action which will be differentiated in
the proposition.  Let \(\mathfrak r_s\) be the complete flow of
\(-\bar\nabla\bar f\), and let
\(\psi_{j,s}^{(\tau_0)}\) be the complete flow of
\(\chi_{\tau_0}W_j\).  Encode a state equivalently by
\((G,\lambda,\Theta,\Phi)\) and set
\begin{align}
 \mathscr L_0(s)(G,\lambda,\Theta,\Phi)
 &=
 \bigl(G,e^{-s}\lambda,
       \mathfrak r_{-s}\circ\Theta,
       \mathfrak r_{-s}\circ\Phi\bigr),
 \label{eq:intro-prepared-scale-leg}\\
 \mathscr L_j(s)(G,\lambda,\Theta,\Phi)
 &=
 \bigl(G,\lambda,
       \psi_{j,-s}^{(\tau_0)}\circ\Theta,
       \psi_{j,-s}^{(\tau_0)}\circ\Phi\bigr),
 \quad1\leq j\leq8.
 \label{eq:intro-prepared-diffeomorphism-leg}
\end{align}
For \(p=(p_0,\ldots,p_8)\), in this fixed order, let
\begin{equation}\label{eq:intro-prepared-tuple-action}
 \mathbf A^{\rm prep}_{p,\tau_0}(\mathbf z)
 :=
 \mathscr L_8(p_8)\circ\cdots\circ
 \mathscr L_0(p_0)(\mathbf z)=:\mathbf z_p .
\end{equation}
Writing the resulting geometric components as
\((G,\lambda_p,\Theta_p,\Phi_p)\), set
\[
 R_p=\varphi_{-\tau_0}\circ\Theta_p,\qquad
 F_p=\Theta_p^{-1}\circ\Phi_p,
\]
and re-form the graph:
\begin{equation}\label{eq:intro-prepared-phase-action}
 \begin{split}
  \acute G_p
  &=\eta\,\iota_*G
    +(1-\eta)\lambda_p\Theta_p^*\bar g,\\
  \mathfrak A^{\rm prep}_{p,\tau_0}(\mathbf z)
  &:=
   \bar g+h\!\left(
      \mathbf A^{\rm prep}_{p,\tau_0}(\mathbf z)\right)
   =\lambda_p^{-1}(\Phi_p^{-1})^*\acute G_p .
 \end{split}
\end{equation}
Thus the phase action preserves the prepared graft graph rather than
acting only on the normalized inner tensor.

For comparison with the dynamical columns, define
\[
 K_{\tau_0,\mathbf z}(T)
 :=(\Phi^{-1})^*((1-\eta)\Theta^*T)
\]
and
\[
 \begin{aligned}
  \mathcal Y_{0,\tau_0}(\mathbf z)
  &:=Y_0-K_{\tau_0,\mathbf z}(Y_0),\\
  \mathcal Y_{j,\tau_0}(\mathbf z)
  &:=\Lie_{\chi_{\tau_0}W_j}\bar g
    -K_{\tau_0,\mathbf z}
       (\Lie_{\chi_{\tau_0}W_j}\bar g),
       \qquad1\leq j\leq8.
 \end{aligned}
\]
The \emph{full prepared phase columns} are
\begin{align}
 \mathscr C_{0,\tau_0}(\mathbf z)
 &:=
 D_{p_0}\!\left[
  \mathfrak A^{\rm prep}_{p,\tau_0}(\mathbf z)-\bar g
 \right]_{p=0}
 =
 \mathcal Y_{0,\tau_0}(\mathbf z)
 +h(\mathbf z)-\Lie_{\bar\nabla\bar f}h(\mathbf z),
 \label{eq:intro-full-column-zero}\\
 \mathscr C_{j,\tau_0}(\mathbf z)
 &:=
 D_{p_j}\!\left[
  \mathfrak A^{\rm prep}_{p,\tau_0}(\mathbf z)-\bar g
 \right]_{p=0}
 =
 \mathcal Y_{j,\tau_0}(\mathbf z)
 +\Lie_{\chi_{\tau_0}W_j}h(\mathbf z),
 \quad1\leq j\leq8.
 \label{eq:intro-full-column-j}
\end{align}
Accordingly, the matrix referred to below is the concrete matrix
\begin{equation}\label{eq:intro-full-phase-column-matrix}
 M^{\rm full}_{\mu j}(\mathbf z)
 :=
 \left\langle
  \rho_{\tau_0}\mathscr C_{j,\tau_0}(\mathbf z),Z_\mu
 \right\rangle_{L^2_\nu},
 \qquad0\leq\mu,j\leq8.
\end{equation}

On a common-margin neighborhood and for \(p\) in a fixed small ball,
the scale-one product, inverse, composition, pullback, and re-grafting
estimates give
\[
 (p,\mathbf z)\longmapsto
 \mathbf A^{\rm prep}_{p,\tau_0}(\mathbf z)
 \quad\text{of class }C^1:
 \quad
 \mathbb R^9\times\mathscr P_{\tau_0}^{r+2,\alpha}
 \longrightarrow\mathscr P_{\tau_0}^{r,\alpha},
\]
with locally uniform first-derivative bounds.  The tensor-valued action
in \eqref{eq:intro-prepared-phase-action} has the same buffered mapping
property.  After pairing against the compactly supported Gaussian test
tensors \(\rho_{\tau_0}Z_\mu\), one integration by parts removes the
derivative loss in the map-direction Lie derivative; hence the
scalar-valued moment map is \(C^1\) already on the unbuffered
\(\mathscr P_{\tau_0}^{r,\alpha}\) chart.  This scalar statement does
not turn the formal columns at a merely \(r\)-regular state into
order-\(r\) tangent vectors.  Split surjectivity below is instead
tested on nine genuine order-\(r\) tangents obtained at the smooth
center and continued as a fixed finite-dimensional complement in a
prepared Banach chart.  These are precisely the mapping facts proved
in the proposition and expanded quantitatively in
Lemma~\ref{lem:prepared-chart-calculus} and
Lemma~\ref{lem:unbuffered-Gaussian-moment-map}.

A \emph{common-margin} subset is one on which ellipticity, radial
comparison, support separation, Gram invertibility, graft compatibility,
bounded geometry, and the displayed coefficient and map norms have one
common quantitative package
\(\mathfrak P_{\rm prep}\), listed in
\eqref{eq:numerical-prepared-package}.
A \emph{uniformly interior common-margin subball}
\(\mathscr B'\Subset_{\rm u}\mathscr B\) means that, in the fixed
 prepared Banach chart, \(\mathscr B'\) has positive distance, measured
 in the independent model norm \eqref{eq:prepared-model-Banach-norm},
 from the boundary of \(\mathscr B\) and carries one uniform
numerical package.  The notation \(\Subset_{\rm u}\) records uniform
interiority only; no topological relative compactness in the
infinite-dimensional prepared space is asserted.

\begin{proposition}[Static sliced chart used in Theorem C]
\label{prop:intro-static-sliced-chart}
Fix \(r\geq3\), \(0<\alpha<1\), \(\tau_0\), and a smooth common-margin
prepared center
\(\mathbf z_c\in\mathscr P_{\tau_0}^{\infty}\) satisfying
\[
 \left\langle\rho_{\tau_0}h(\mathbf z_c),Z_\mu\right\rangle_{L^2_\nu}=0,
 \qquad0\leq\mu\leq8.
\]
Suppose that the \(9\times9\) matrix
\(M^{\rm full}(\mathbf z_c)\) in
\eqref{eq:intro-full-phase-column-matrix} is invertible.  There is a
sufficiently small open common-margin neighborhood
\[
 \mathscr N_{\mathrm{sl},\tau_0}^{r,\alpha}
 \subset\mathscr P_{\tau_0}^{r,\alpha}
\]
of \(\mathbf z_c\) such that the moment map
\[
 \mathfrak m_{\tau_0}:
 \mathscr N_{\mathrm{sl},\tau_0}^{r,\alpha}\longrightarrow\mathbb R^9,
 \qquad
 \mathfrak m_{\tau_0}(\mathbf z)
 =
 \bigl(
  \langle\rho_{\tau_0}h(\mathbf z),Z_\mu\rangle_{L^2_\nu}
 \bigr)_{\mu=0}^8,
\]
is \(C^1\), and \(D\mathfrak m_{\tau_0}(\mathbf z)\) is split
surjective at every
\(\mathbf z\in\mathscr N_{\mathrm{sl},\tau_0}^{r,\alpha}\).  Define
\[
 \Sigma_{\tau_0}^{r,\alpha}
 :=
 \left\{
  \mathbf z\in\mathscr N_{\mathrm{sl},\tau_0}^{r,\alpha}:
  \mathfrak m_{\tau_0}(\mathbf z)=0
 \right\}.
\]
Then \(\Sigma_{\tau_0}^{r,\alpha}\) is a split \(C^1\) submanifold of
codimension nine in
\(\mathscr N_{\mathrm{sl},\tau_0}^{r,\alpha}\), and hence a local split
submanifold of \(\mathscr P_{\tau_0}^{r,\alpha}\).  At the buffered
input order use the compatible restrictions
\[
 \mathscr N_{\mathrm{sl},\tau_0}^{r+2,\alpha}
 :=
 \mathscr N_{\mathrm{sl},\tau_0}^{r,\alpha}
 \cap\mathscr P_{\tau_0}^{r+2,\alpha},
 \qquad
 \Sigma_{\tau_0}^{r+2,\alpha}
 :=
 \Sigma_{\tau_0}^{r,\alpha}
 \cap\mathscr P_{\tau_0}^{r+2,\alpha}.
\]
There is an open neighborhood
\[
 \mathscr O_{\tau_0}^{r+2,\alpha}
 \subset
 \mathscr N_{\mathrm{sl},\tau_0}^{r+2,\alpha}
\]
and a centered \(C^1\) phase retraction
\[
 \Pi_{\rm sl}^{\,r+2\to r}:
 \mathscr O_{\tau_0}^{r+2,\alpha}
 \longrightarrow\Sigma_{\tau_0}^{r,\alpha}.
\]
Its restriction to
\(\Sigma_{\tau_0}^{r+2,\alpha}\cap
\mathscr O_{\tau_0}^{r+2,\alpha}\)
is the canonical inclusion into
\(\Sigma_{\tau_0}^{r,\alpha}\).
After shrinking \(\mathscr O_{\tau_0}^{r+2,\alpha}\) once more to a
convex model-ball inside the same common-margin prepared chart, the
quantitative implicit-function construction gives a derived constant
\(K_{\Pi,r}<\infty\) such that, for all
\(\mathbf z_1,\mathbf z_2\in
  \mathscr O_{\tau_0}^{r+2,\alpha}\),
\begin{align}
 \|\Pi_{\rm sl}^{\,r+2\to r}(\mathbf z_1)
      -\Pi_{\rm sl}^{\,r+2\to r}(\mathbf z_2)\|_{
      \mathscr X_{\rm prep}^{r,\alpha}}
 &\leq K_{\Pi,r}
 \|\mathbf z_1-\mathbf z_2\|_{
      \mathscr X_{\rm prep}^{r+2,\alpha}},
 \label{eq:intro-quantitative-phase-retraction-Lip}\\
 \sup_{\mathbf z\in\mathscr O_{\tau_0}^{r+2,\alpha}}
 \|D\Pi_{\rm sl}^{\,r+2\to r}(\mathbf z)\|_{
  \mathcal L(\mathscr E_{\rm prep}^{r+2,\alpha},
             \mathscr E_{\rm prep}^{r,\alpha})}
 &\leq K_{\Pi,r}.
 \label{eq:intro-quantitative-phase-retraction-D}
\end{align}
The constant depends only on the inverse phase-column bound and the
common prepared chart and composition package; it is not obtained from
compactness of a bounded set.
The chosen center and local slice neighborhood are suppressed from the
notation \(\Sigma_{\tau_0}^{r,\alpha}\).  This notation always denotes
the local sliced level set just defined; no assertion is made about
other zeros of \(\mathfrak m_{\tau_0}\) in the full prepared manifold.
If \(\mathscr B\subset\Sigma_{\tau_0}^{r,\alpha}\) is a common-margin
sliced ball and
\(\mathcal E_{\rm sl}:\mathscr B\to\mathbb B_{\rm tar}\) is \(C^1\) into one
of the Banach target charts used below, define
\[
 \mathcal E_{\rm amb}
 :=\mathcal E_{\rm sl}\circ\Pi_{\rm sl}^{\,r+2\to r}
\]
on any open
\(\mathscr O\subset\mathscr O_{\tau_0}^{r+2,\alpha}\) satisfying
\(\Pi_{\rm sl}^{\,r+2\to r}(\mathscr O)\subset\mathscr B\).
This composition is \(C^1\) and satisfies the chain rule.  If
\(\overline{\mathscr O}\subset
  \mathscr O_{\tau_0}^{r+2,\alpha}\),
\(\mathscr B'\Subset_{\rm u}\mathscr B\),
\(\Pi_{\rm sl}^{\,r+2\to r}(\overline{\mathscr O})\subset\mathscr B'\),
and there is \(L_{\mathcal E}<\infty\) such that, in the chosen target
chart, for all \(u,v\in\mathscr B'\),
\[
 \|\mathcal E_{\rm sl}(u)-\mathcal E_{\rm sl}(v)\|_{\mathbb B_{\rm tar}}
 \leq L_{\mathcal E}
       \|u-v\|_{\mathscr X_{\rm prep}^{r,\alpha}},
 \qquad
 \sup_{u\in\mathscr B'}\|D\mathcal E_{\rm sl}(u)\|
 \leq L_{\mathcal E},
\]
then
\[
 \operatorname{Lip}(\mathcal E_{\rm amb}|_{\mathscr O})
 +\sup_{z\in\mathscr O}\|D\mathcal E_{\rm amb}(z)\|
 \leq 2L_{\mathcal E}K_{\Pi,r}.
\]
The closure and uniformly-interior conditions retain the fixed domains
and prepared margins; the two displayed quantitative bounds, not
infinite-dimensional compactness, give uniformity.
\end{proposition}

\begin{proof}
On the compact core, work in a fixed finite harmonic atlas.  On each
dyadic end annulus rescale \(\bar g\) to unit size.  The common radial
comparison makes every background map send that annulus into one fixed
enlargement, with constants independent of its scale.  The ordinary
little-H\"older product and chain rules in these unit charts show that
composition, inversion, and pullback are \(C^1\) after the displayed
two-derivative input buffer; bounded overlap gives the global scaled
estimate.  Differentiating
\eqref{eq:intro-prepared-scale-leg}--%
\eqref{eq:intro-prepared-phase-action}, the simultaneous left
composition of \(\Theta\) and \(\Phi\) cancels on the pure outer target,
while on the physical summand it gives the corresponding pullback.
The radial scale leg similarly combines the factor \(e^{-s}\) with the
soliton pullback.  This proves the exact column formulas
\eqref{eq:intro-full-column-zero}--%
\eqref{eq:intro-full-column-j}.  The compact support of
\(\rho_{\tau_0}\) permits integration
by parts in every map-direction Lie derivative, so the scalar-valued
moment map itself is \(C^1\) at the unbuffered order.  Equations
\eqref{eq:intro-full-column-zero}--%
\eqref{eq:intro-full-phase-column-matrix} show that, at the smooth
center, its derivative on the nine genuine tangents
\[
 e_j:=
 \left.\partial_{p_j}
   \mathbf A^{\rm prep}_{p,\tau_0}(\mathbf z_c)\right|_{p=0}
 \in T_{\mathbf z_c}\mathscr P_{\tau_0}^{r,\alpha}
\]
is exactly \(M^{\rm full}(\mathbf z_c)\).  Choose a prepared Banach
chart at \(\mathbf z_c\) and continue the \(e_j\) as constant chart
tangent fields.  Continuity of \(D\mathfrak m_{\tau_0}\) makes its
restriction to their span remain invertible after the neighborhood is
shrunk.  Thus \(D\mathfrak m_{\tau_0}\) is genuinely split surjective
at every point of the unbuffered neighborhood; no derivative-losing
formal phase column at a nonsmooth point is used as a tangent.

The Banach implicit-function theorem now gives the asserted split
level set inside
\(\mathscr N_{\mathrm{sl},\tau_0}^{r,\alpha}\); it gives no conclusion
about zeros outside this neighborhood.
Applying the same theorem to the phase parameter, with the two
derivatives reserved for the tuple-valued pullback action, gives the
displayed retraction.  More explicitly, put
\[
 \mathcal F(\mathbf z,p)
 :=\mathfrak m_{\tau_0}\!\left(
   \mathbf A^{\rm prep}_{p,\tau_0}(\mathbf z)
 \right).
\]
On the shrunken neighborhood let \(p_{\tau_0}(\mathbf z)\) denote the
unique small \(C^1\) solution, centered by
\(p_{\tau_0}(\mathbf z_c)=0\), of
\[
 \mathcal F\bigl(\mathbf z,p_{\tau_0}(\mathbf z)\bigr)=0,
\]
and set
\[
 \Pi_{\rm sl}^{\,r+2\to r}(\mathbf z)
 :=\mathbf A^{\rm prep}_{p_{\tau_0}(\mathbf z),\tau_0}(\mathbf z).
\]
After the quantitative shrinking, \(D_p\mathcal F\) has one uniform
inverse bound and \(D_{\mathbf z}\mathcal F\) has one uniform prepared
chart bound.  Hence
\[
 Dp_{\tau_0}
 =-(D_p\mathcal F)^{-1}D_{\mathbf z}\mathcal F
\]
is uniformly bounded.  Differentiating the phase action gives
\[
 D\Pi_{\rm sl}
 =D_{\mathbf z}\mathbf A^{\rm prep}
  +D_p\mathbf A^{\rm prep}\,Dp_{\tau_0},
\]
so the prepared composition estimates give
\eqref{eq:intro-quantitative-phase-retraction-D}; the Lipschitz bound
follows on the chosen convex model-ball.  Centered uniqueness gives the
inclusion identity.  The final ambient assertions are now the chain
rule and the displayed product bound.  The later prepared-chart lemmas
expand these same estimates, but no evolution statement enters this
proof.
\end{proof}

\begin{remark}[Geometric content of the Theorem~C entrance conditions]
\label{rem:intro-Theorem-C-entrance-preview}
For orientation, the hypotheses in
Definition~\ref{def:strict-prepared-entrance} fall into the following
six groups.
Fix \(k_0\geq12\), \(0<\alpha<1\), and first fix
\[
 0<\sigma<\theta<\beta .
\]
Next fix one rate-compatible numerical prepared package
\(\mathfrak P_{\rm prep}\) from
\eqref{eq:numerical-prepared-package}, including its already fixed
strict scale constants \(0<c_{\rm sc}<C_{\rm sc}<\infty\).
In particular, the parameters
\((m_{\rm ad},\varepsilon_{\rm ph},\tau_{\rm ad},C_\lambda)\) are fixed.
Require
\[
 \tau_0\geq\tau_{\rm ad},
 \qquad t(\tau_0)=0,
\]
where the second equality is the physical-clock convention and does not
add a coordinate to the prepared tuple.
Use the single \(\varepsilon_{\rm ent}>0\) already selected in
Remark~\ref{conv:authoritative-adaptive-order}, no larger than every
smallness threshold in the feedback, three-region, graft,
harmonic-map, continuation, and uniform two-state estimates.  Fix
 \[
  0<\varepsilon\leq\varepsilon_{\rm ent}.
 \]
The exact definition concerns tuples
\(\mathbf z_0=(G_0,\lambda_0,R_0,F_0)
 \in\Sigma_{\tau_0}^{k_0+2,\alpha}\)
in the global weighted little topology.  Its six groups of conditions,
all measured in that single package with strict positive margins, are
previewed as follows:
\begin{enumerate}
\item it is obtained either from a raw
      \(\mathscr P_{\tau_0}^{k_0+4,\alpha}\) tuple by the buffered
      retraction
      \(\Pi_{\rm sl}^{\,k_0+4\to k_0+2}\), or is already sliced in
      \(\mathscr P_{\tau_0}^{k_0+2,\alpha}\) with invertible full
      phase matrix;
\item \(c_{\rm sc}e^{-\tau_0}<\lambda_0<
      C_{\rm sc}e^{-\tau_0}\), the nine moments vanish exactly, and
      \(h_0\) satisfies
      \[
       \|\rho_{\tau_0}h_0\|_{L^2_\nu}
        <\varepsilon e^{-\theta\tau_0},\qquad
       \sum_{\ell=0}^3|\bar\nabla^\ell h_0|
        <\varepsilon\omega_\sigma(\tau_0,\cdot),
      \]
      together with strict ellipticity, global \(C^2\), and
      historical-tail margins;
\item outside the fixed physical graft,
      \(\acute G_0=\lambda_0\Theta_0^*\bar g\), where
      \(\Theta_0=\varphi_{\tau_0}\circ R_0\);
      \(R_0^{\pm1}\), the adaptive columns, and the dyadic radial
      comparison satisfy their scale-normalized order-five bounds;
\item the fixed graft parameter and nested graft collars have the
      strict order-six interpolation, support, positivity, and
      pure-graft \(C^2\) margins;
\item \(F_0\) is a proper degree-one diffeomorphism with the strict
      scale-normalized order-five map margin, the adaptive Gram matrix
      is uniformly invertible, and the complete source and target
      have the stated scale-\(\lambda_0^{1/2}\) harmonic-radius,
      curvature, and Ricci-defect bounds through order ten;
\item the fixed core, interface, graft-input, exterior, and auxiliary
      collar chains satisfy the positive separation, bounded-overlap,
      Lebesgue-number, and coefficient-transfer bounds, including
      the one-time order-fourteen closed-metric input which yields the
      order-twelve retained physical coefficient package.
\end{enumerate}
Definition~\ref{def:strict-prepared-entrance} gives the exact
inequalities, collar inclusions, and constants.  These hypotheses apply
only to the marked high-regularity refinement and impose no restriction
on the full-metric neighborhood in Theorem~A.
\end{remark}

\subsection{Weighted spaces}

All norms and covariant derivatives in the paper are taken with
respect to $\bar g$ unless indicated otherwise.  Let \(E\to M\) be
any fixed natural real metric tensor bundle with its induced connection,
or its complexification with the induced Hermitian metric and
compatible connection.  Define
\[
 L^2_\nu(E):=L^2(M;E,d\nu),
\]
and, for every integer \(m\geq0\),
\[
 H^m_\nu(E):=
 \overline{C_c^\infty(M;E)}^{\|\cdot\|_{H^m_\nu(E)}},
 \qquad
 \|q\|_{H^m_\nu(E)}^2
 :=\sum_{\ell=0}^m
       \int_M|\bar\nabla^\ell q|_{\bar g}^2\,d\nu.
\]
Set
\[
 H^{-1}_\nu(E):=H^1_\nu(E)^*,
\]
with duality induced by the \(L^2_\nu(E)\) pairing.  We suppress the
bundle \(E\) whenever the tensor type is clear; unqualified spectral
and operator spaces default to \(E=S^2T^*M\).  For sections of the
same bundle, every unqualified pairing below means
\[
 \ip{q}{r}
 :=\int_M\langle q,r\rangle_{\bar g}\,d\nu.
\]
Smooth compactly supported sections are dense in the spaces used
below.  All integration-by-parts identities are first proved for such
sections and then extended by density.

We now record in detail the self-adjoint realization fixed before
Theorems~B and~C.  On
$C^\infty_c(S^2T^*M)$ consider the symmetric expression
\[
 \A h=\bar\Delta_{\bar f}h+2\overline{\Rm}(h),
 \qquad
 \bar\Delta_{\bar f}h
 =\bar\Delta h-\bar\nabla_{\bar\nabla\bar f}h.
\]
Our curvature sign convention is
\[
 \bar R(X,Y)Z
 =\bar\nabla_X\bar\nabla_YZ-\bar\nabla_Y\bar\nabla_XZ
  -\bar\nabla_{[X,Y]}Z,
\]
and our curvature action is
\[
 \bigl(\overline{\Rm}(h)\bigr)_{ij}
 =\bar R_{ikj\ell}h^{k\ell},
\]
where indices are raised with \(\bar g\).
Lemma~\ref{lem:self-contained-FIK-ledger} computes the curvature
components directly and shows that \(\overline{\Rm}\) is bounded.
Hence, for a fixed sufficiently large
Friedrichs-shift constant $C_{\rm Fr}$ the form
\[
 \mathfrak a[h,k]
 =-\ip{\bar\nabla h}{\bar\nabla k}
  +2\int_M\overline{\Rm}(h,k)\,d\nu ,
 \qquad h,k\in H^1_\nu,
\]
is a continuous symmetric form on \(H^1_\nu\), and
\[
 \mathfrak l(h,k)
 =-\mathfrak a[h,k]+C_{\rm Fr}\ip{h}{k}
\]
is closed and coercive on $H^1_\nu$.  We denote by $L_0$ its
Friedrichs operator and define the self-adjoint operator
$\A=C_{\rm Fr}-L_0$.  Thus \(L_0=-\A+C_{\rm Fr}\) has form
domain \(H^1_\nu\) and closed form \(\mathfrak l\).  Equivalently,
\(\A\) is represented on \(H^1_\nu\) by the continuous symmetric form
\begin{equation}\label{eq:A-form}
 \mathfrak a[h,h]
 =-\norm{\bar\nabla h}_{L^2_\nu}^2
   +2\int_M\overline{\Rm}(h,h)\,d\nu.
\end{equation}
For \(h\in D(\A)\), one has
\(\mathfrak a[h,k]=\ip{\A h}{k}\) for every \(k\in H^1_\nu\).
Whenever the operator notation is used with an \(H^1_\nu\) argument
below, it means this \(H^1_\nu\)--\(H^{-1}_\nu\) form pairing; no
claim that \(\A h\in L^2_\nu\) is then being made.
Equivalently, the operator domain is
\begin{equation}\label{eq:A-operator-domain}
 \begin{split}
 D(\A)=\bigl\{h\in H^1_\nu:\;&\text{there is }u\in L^2_\nu
 \text{ such that}\\
 &-\ip{\bar\nabla h}{\bar\nabla k}
   +2\int_M\overline{\Rm}(h,k)\,d\nu
   =\ip{u}{k}
 \quad\text{for every }k\in H^1_\nu\bigr\},
 \end{split}
\end{equation}
and then \(\A h=u\).  Thus the operator domain and its action have been
specified independently of the later spectral calculation.

\begin{lemma}[Gaussian first moment]\label{lem:first-moment}
For every fixed metric or Hermitian tensor bundle \(E\to M\) as above
and every \(q\in H^1_\nu(E)\),
\begin{equation}\label{eq:first-moment}
 \norm{\sqrt{\bar f}\,q}_{L^2_\nu}^2
 \leq4\left(\norm q_{L^2_\nu}^2+
             \norm{\bar\nabla q}_{L^2_\nu}^2\right).
\end{equation}
\end{lemma}

\begin{proof}
The soliton identities \eqref{eq:shrinker}, together with the
coordinate-independent inequality
\(\bar R\geq0\) proved in
Lemma~\ref{lem:self-contained-FIK-ledger}, give
\[
 \bar\Delta_{\bar f}\bar f=2-\bar f,
 \qquad
 \abs{\bar\nabla\bar f}^2\leq\bar f.
\]
For compactly supported $q$, weighted integration by parts yields
\[
 \int_M\bar f\abs q^2\,d\nu
 =2\int_M\abs q^2\,d\nu+
   \int_M\langle\bar\nabla\abs q^2,
                    \bar\nabla\bar f\rangle\,d\nu.
\]
The Kato inequality and $2ab\leq\frac12a^2+2b^2$ imply
\[
 \int_M\bar f\abs q^2\,d\nu
 \leq2\norm q_{L^2_\nu}^2+
 \frac12\norm{\sqrt{\bar f}\,q}_{L^2_\nu}^2+
 2\norm{\bar\nabla q}_{L^2_\nu}^2.
\]
Absorption proves \eqref{eq:first-moment}.  The general case follows
by density.
\end{proof}

\subsection{The spectral input}

We first complete the functional-analytic bridge for the spectral
calculation.  The exact coordinate formulas
in Lemma~\ref{lem:self-contained-FIK-ledger} use the completed FIK
radial coordinate
\[
 r\in[1,\infty),
\]
with \(r=1\) at the bolt, and give
\(\bar f=\sqrt2r^2+2-2\sqrt2\); hence \(\bar f\) is proper and satisfies
\(\bar f\simeq r^2\) on the asymptotically conical end; in particular,
its sublevel sets are compact.  The embedding of the form domain into
\(L^2_\nu\) is compact in the present Gaussian AC setting.  Indeed,
Lemma~\ref{lem:first-moment} gives, uniformly on bounded subsets of
\(H^1_\nu\),
\[
 \int_{\{\bar f>R\}}|h|^2\,d\nu
 \leq \frac{C}{R}\|h\|_{H^1_\nu}^2.
\]
On \(\{\bar f\leq R\}\) the weighted and unweighted norms are
equivalent, so local Rellich compactness, a diagonal subsequence, and
the displayed tail estimate prove compactness globally.  The
self-adjoint realization of \(\A\) therefore has compact resolvent.

We now state the precise normalization dictionary.  The mass-normalized
potential used in the Wigner calculation is
\[
 f_{\rm NO}=\sqrt2(r^2-1)-\log(2c_0),
 \qquad c_0=\sqrt2-1.
\]
It satisfies
\[
 \bar R
 +|\bar\nabla f_{\rm NO}|^2-f_{\rm NO}=C_{\rm NO},
 \qquad
 C_{\rm NO}=\log(2c_0)+\sqrt2\,c_0,
\]
Thus
\[
 \bar f=f_{\rm NO}+C_{\rm NO}.
\]
This additive shift leaves the drift Laplacian and Hessians unchanged,
and it multiplies their weighted measure by a positive constant:
\[
 d\nu=e^{-C_{\rm NO}}(4\pi)^{-2}e^{-f_{\rm NO}}\,dV_{\bar g}.
\]
For clarity, set
\[
 d\mu_f:=e^{-f_{\rm NO}}\,dV_{\bar g}.
\]
For every metric or Hermitian tensor bundle \(E\to M\) covered by the
bundlewise convention above, define the unnormalized Wigner spaces
\[
 L^2_f(E):=L^2(M;E,d\mu_f),
 \qquad
 H^1_f(E):=
 \overline{C_c^\infty(M;E)}^{
  \left(\int_M
   (|u|_{\bar g}^2+|\bar\nabla u|_{\bar g}^2)\,d\mu_f
  \right)^{1/2}}.
\]
We again suppress \(E\) when the section type is clear; in particular,
the scalar, one-form, and symmetric-two-tensor channels use
\(E=\mathbb C\), \((T^*M)_{\mathbb C}\), and
\((S^2T^*M)_{\mathbb C}\), respectively.  On each such bundle put
\[
 \Delta_f:=\bar\Delta-\bar\nabla_{\bar\nabla f_{\rm NO}}.
\]
Consequently \(L^2_f(E),H^1_f(E)\) and
\(L^2_\nu(E),H^1_\nu(E)\) have the same elements, with norms differing
only by fixed positive factors.  With the curvature convention fixed
above, the symmetric-two-tensor block operator
\(L_f=\Delta_f+2\overline{\Rm}\) is exactly \(\A\).
Appendix~\ref{app:self-contained-FIK-spectrum} proves the tensor
Peter--Weyl decomposition, the exact \(J\leq4\) block inventory, and
the exceptional and radial mode counts internally.  After the Wigner
blocks and cross-sectional reduction have been established there,
Lemma~\ref{lem:FIK-high-frequency} supplies the remaining uniform
coefficient comparison for \(J\geq5\).  Its proof uses no nonlinear
input.

For the exhaustion argument, we use the following terminology from
Appendix~\ref{app:self-contained-wigner}.  The integer
\(J\in\mathbb N_0\) is the \(SU(2)\) Peter--Weyl highest-weight index
in the normalization of that appendix; the radial sector is \(J=0\).
Among the connected zeroth-order Wigner blocks with
\(0\leq J\leq4\), the \emph{exceptional} blocks are the radial
two-component block, the two conjugate \(J=1\) anti-invariant blocks
carrying the positive Hessian family, and the \(J=2\) invariant block
carrying the zero Hessian family.  Every other connected block in this
range is called \emph{safe}.  The orthogonal reducing splitting,
entrywise matrices, and exhaustive block list are established in
\eqref{eq:self-mode-splitting},
\eqref{eq:self-Q-connected-block}--%
\eqref{eq:self-P-connected-block}, and
Appendix~\ref{app:exact-low-mode-certificates}.

\begin{theorem}[Certified FIK spectral decomposition]\label{thm:FIK-spectrum}
The nonnegative spectral subspace of $\A$ is nine-dimensional.  It has
an $L^2_\nu$-orthonormal real basis $Z_0,\ldots,Z_8$ with eigenvalues
\[
 1,\quad
 \underbrace{1-\frac1{\sqrt2},\ldots,
             1-\frac1{\sqrt2}}_{4\text{ times}},\quad
 \underbrace{0,\ldots,0}_{4\text{ times}}.
\]
The eigenvalue-$1$ space is spanned by $\Ric_{\bar g}$.  The remaining
eight tensors are Hessians of explicit functions.  Every other
eigenvalue is strictly negative.
\end{theorem}

\begin{proof}
The shrinker identities give
\(\A\Ric_{\bar g}=\Ric_{\bar g}\).
Lemma~\ref{lem:self-contained-scalar-modes} and
Corollary~\ref{cor:self-contained-Hessian-families} construct,
directly and with the present normalization, four real Hessian modes
of eigenvalue \(1-1/\sqrt2\), three real \(J=2\) Hessian modes of
eigenvalue \(0\), and the radial zero mode
\(\bar\nabla^2\bar f\).  Lemma~\ref{lem:self-contained-bolt-reality}
proves their smooth bolt extension and their exact real dimensions.
The growth estimates in
Lemma~\ref{lem:self-contained-geometric-growth} show that these tensors
and all their derivatives have at most polynomial growth.  Since the
weighted measure has Gaussian decay, all nine tensors belong to every
\(H^k_\nu\).  Their displayed distributional eigen-equations, tested
first against compactly supported tensors and then extended by the
 form-core cutoff, place them in \(D(\A)\) by the characterization in
 \eqref{eq:A-operator-domain}.

It remains to prove exhaustion.  The Peter--Weyl splitting
\eqref{eq:self-mode-splitting} is an orthogonal reducing decomposition
for both the form and the operator.  In the radial \(J=0\) sector, the
closed Friedrichs-domain reduction of
Proposition~\ref{prop:FIK-radial-domain-package}, together with
Proposition~\ref{prop:FIK-J0-inventory}, leaves exactly
\(\Ric_{\bar g}\) and \(\bar\nabla^2\bar f\) nonnegative.
For \(1\leq J\leq4\), Proposition
\ref{prop:exact-safe-low-mode-blocks} makes every safe zeroth-order
block pointwise negative definite; the exact mode identity
\eqref{eq:self-exact-mode-identity} has an additional nonpositive
radial-derivative term, so every corresponding operator eigenvalue is
strictly negative.  The only omitted blocks are the two conjugate
exceptional \(J=1\) families and the exceptional \(J=2\) family.
Corollary~\ref{cor:FIK-exceptional-count} allows at most one
nonnegative eigenvalue in each such complex block.  The explicit
Hessian modes already constructed attain that bound: after the
reality pairing they give precisely four real positive modes for
\(J=1\) and three real zero modes for \(J=2\).
Finally, Lemma~\ref{lem:FIK-high-frequency} excludes every block with
\(J\geq5\).

Thus no other nonnegative eigenvalue exists.  Compact resolvent
excludes non-point spectrum, and real orthonormalization within the
three finite-dimensional eigenspaces gives the stated basis.
\end{proof}

Henceforth
\begin{equation}\label{eq:Z-definition}
 \Z:=\operatorname{span}_{\mathbb R}\{Z_0,\ldots,Z_8\},
 \qquad
 \A Z_a=\lambda_aZ_a,\quad 0\leq a\leq8.
\end{equation}

The strict negativity in Theorem~\ref{thm:FIK-spectrum} is therefore separated
from zero by a spectral gap, which verifies the paper-wide choice
\eqref{eq:beta} made before Theorems~B and~C.

\begin{remark}
Theorem~\ref{thm:FIK-spectrum} supplies more than $\nu$-stability on the
weighted-divergence-free slice: we use the full nonnegative spectrum
of $\A$, including the pure-gauge eigenvectors outside that slice.
Our normalization can be checked directly:
\[
 \A(\Ric_{\bar g})=\Ric_{\bar g},\qquad
 \A(\bar\nabla^2\phi)=(1-\mu)\bar\nabla^2\phi
 \quad\text{when }\bar\Delta_{\bar f}\phi=-\mu\phi.
\]
This gives the eigenvalues $1$, $1-1/\sqrt2$, and $0$ listed above
with the present sign convention.
\end{remark}

\begin{lemma}[Quantitative asymptotically conical symbol estimate]
\label{lem:FIK-AC-symbol}
In the radial coframe of
Lemma~\ref{lem:self-contained-FIK-ledger}, let
\(g_{\rm C}\) be the \(U(2)\)-invariant cone obtained from the exact
FIK profile
\[
 F(r)=\frac1{\sqrt2}-\frac{\sqrt2-1}{r^2}
      -\frac{\sqrt2-1}{\sqrt2\,r^4}
\]
by replacing \(F\) with \(F_\infty=1/\sqrt2\).  For the cone
dilations \(\delta_R(r,\omega)=(Rr,\omega)\) and every \(j\geq0\),
\begin{equation}\label{eq:FIK-AC-quantitative}
 \left\|
  \delta_R^*(R^{-2}\bar g)-g_{\rm C}
 \right\|_{C^j_{g_{\rm C}}(\{1\leq r\leq2\})}
 \leq C_jR^{-2},\qquad R\geq2.
\end{equation}
\end{lemma}

\begin{proof}
Every radial coefficient of \(\bar g\) is a smooth algebraic
expression in \(F\) and its logarithmic derivatives.  Differentiating
the displayed exact profile with \(r\partial_r\) therefore gives the
same \(O(R^{-2})\) bound at every fixed order after dilation.  Thus
\eqref{eq:FIK-AC-quantitative} is an all-order symbol estimate, not
merely pointed \(C^0\) convergence.
\end{proof}

\subsection{Normalized Ricci--DeTurck perturbations}

Using DeTurck's strictly parabolic gauge \cite{DeTurck}, normalized
Ricci--DeTurck flow relative to $\bar g$, written in soliton
coordinates, has $\bar g$ as a stationary point and takes the form
\begin{equation}\label{eq:perturbation}
 \partial_\tau h=\A h+\Q(h)+\E.
\end{equation}
For later integration by parts we record the nonlinear remainder
exactly.  Write
\[
 g^{ab}=\bar g^{ab}-\widehat h^{ab},\qquad
 \widehat h^{ab}=\bar g^{ak}\bar g^{b\ell}h_{k\ell}
                         +\widetilde h^{ab},
 \qquad \widetilde h=O(h^2).
\]
In background coordinates,
\begin{equation}\label{eq:Q-exact}
\begin{split}
 \Q(h)_{ij}={}&
 -\widehat h^{ab}\bar\nabla_a\bar\nabla_bh_{ij}\\
 &+\bar R_{ja}{}^p{}_b
   \bigl(\bar g_{ip}\widetilde h^{ab}
         +\widehat h^{ab}h_{ip}\bigr)
  +\bar R_{ia}{}^p{}_b
   \bigl(\bar g_{jp}\widetilde h^{ab}
         +\widehat h^{ab}h_{jp}\bigr)\\
 &+g^{ab}g^{pq}\,
   \mathcal N(\bar\nabla h,\bar\nabla h)_{abpqij},
\end{split}
\end{equation}
where
\begin{equation}\label{eq:N-exact}
\begin{split}
 2\mathcal N_{abpqij}={}&
 (\bar\nabla_i h_{pa})(\bar\nabla_jh_{qb})
 +2(\bar\nabla_a h_{jp})(\bar\nabla_qh_{ib})
 -2(\bar\nabla_a h_{jp})(\bar\nabla_bh_{iq})\\
 &-2(\bar\nabla_jh_{pa})(\bar\nabla_bh_{iq})
 -2(\bar\nabla_i h_{pa})(\bar\nabla_bh_{jq}).
\end{split}
\end{equation}
This is the standard Ricci--DeTurck expansion used in
\cite[Section~4]{Stolarski}.  In particular,
\begin{equation}\label{eq:Q-schematic}
 \Q(h)
 =-\widehat h*\bar\nabla^2h
  +(\bar g+h)^{-2}*\bar\nabla h*\bar\nabla h
  +\overline{\Rm}*h*h.
\end{equation}
The term $\E$ records extension, cutoff, and chart errors.  The
advantage of retaining $\E$ is that the same perturbation equation applies to a
genuine complete perturbation, for which $\E=0$, and to a closed flow
extended across a moving outer annulus.

We will augment \eqref{eq:perturbation} by a finite-dimensional
geometric moving frame.  The precise abstract equation is
\eqref{eq:modulated}.

\section{Gaussian moments and the geometric modes}
\label{sec:moments}

\subsection{Growth of the geometric generators}

We use the Wigner conventions of
Appendix~\ref{app:self-contained-wigner}, summarized here.  For
\(J\in\mathbb N_0\), put
\[
 \mathcal J_J=\{-J,-J+2,\ldots,J-2,J\}.
\]
The complex harmonics \(D^J_{M,M'}\), \(M,M'\in\mathcal J_J\), are
zero when an index lies outside \(\mathcal J_J\), are normalized by
\[
 \int_{S^3}D^J_{M,M'}\overline{D^K_{N,N'}}\,dV_{S^3}
 =\frac{2\pi^2}{J+1}
   \delta_{JK}\delta_{MN}\delta_{M'N'},
\]
where \(dV_{S^3}\) denotes the unit-round Haar volume of total mass
\(2\pi^2\), in the \(SU(2)\) normalization fixed explicitly in
Appendix~\ref{app:self-contained-FIK-ledger},
and obey the angular derivative convention
\eqref{eq:self-wigner-derivatives}.  The real harmonics used below are
fixed normalized real or imaginary combinations of these complex
functions.  Appendix~\ref{app:self-contained-wigner} proves their
completeness, phase convention, and full tensor block reduction.

We use the radial coordinate \(r\in[1,\infty)\) and the real Wigner
harmonics \(D^J_{M,M'}\) fixed in
Appendix~\ref{app:self-contained-FIK-spectrum}.  The eigenspace
identifications are those proved in Theorem~\ref{thm:FIK-spectrum}.  In the
abbreviated displays below,
\(D^1\) denotes one of the four real \(J=1\) harmonics occurring in the
positive eigenspace, while \(D^2\) denotes one of the three real
\(J=2,\ M=0\) harmonics occurring in the zero eigenspace.  The fourth
zero-mode potential is \(\bar f\).  The eight real potentials and
their nonzero normalizations fixed in the opening spectral convention
are denoted by \(\phi_1,\ldots,\phi_8\) and
\(\mathfrak c_1,\ldots,\mathfrak c_8\); thus
\[
 Z_a=\mathfrak c_a\bar\nabla^2\phi_a,
 \qquad
 W_a:=\frac{\mathfrak c_a}{2}\bar\nabla\phi_a,
 \qquad 1\leq a\leq8,
\]
and keep these choices fixed.

The explicit functions in
Lemma~\ref{lem:self-contained-scalar-modes} have two asymptotic types.
The four positive modes have potentials
\[
 \phi=\widehat u\,D^1,\qquad
 \widehat u(r)\sim r^{\sqrt2},
\]
whereas the zero modes have quadratic potentials
\[
 \phi=\bar f
 \quad\text{or}\quad
 \phi=\widehat v\,D^2,\qquad
 \widehat v(r)=r^2
\]
in the normalization fixed here.  Since $\bar f\simeq r^2$ on the
asymptotically conical end, this gives the following bounds.

\begin{lemma}[Mode and generator growth]\label{lem:mode-growth}
With the preceding normalizations, for every integer \(m\geq0\),
\begin{align}
 |\bar\nabla^m W_a|
 &\leq C_m(1+\bar f)^{(1-m)/2},
 &&1\leq a\leq8,\label{eq:W-all-order}\\
 |\bar\nabla^m Z_a|
 &\leq C_m(1+\bar f)^{-m/2},
 &&0\leq a\leq8.\label{eq:Z-all-order}
\end{align}
In particular,
\begin{align}
 \abs{\bar\nabla\phi_a}
 &\leq C(1+\sqrt{\bar f}),&
 \abs{\bar\nabla^2\phi_a}
 &\leq C,\label{eq:V-growth}\\
 \abs{Z_a}+\sqrt{\bar f}\abs{\bar\nabla Z_a}
 &\leq C,&
 1&\leq a\leq8.\label{eq:Z-growth}
\end{align}
For \(Z_0=c_{\mathrm{Ric}}\Ric_{\bar g}\), where
\(c_{\mathrm{Ric}}>0\) is the fixed \(L^2_\nu\)-normalizing constant
chosen in the opening spectral convention and verified by
Theorem~\ref{thm:FIK-spectrum},
\eqref{eq:Z-all-order} follows from the differentiated quadratic
curvature decay of the FIK metric.  In particular, all
$Z_a$ belong to $H^k_\nu$ for every $k$.
\end{lemma}

\begin{proof}
We use the exact radial formulae in
Lemma~\ref{lem:self-contained-scalar-modes}, rather than
differentiating an undifferentiated asymptotic equivalence:
\[
 \widehat u
 =(r^2-1)^{1/2}(r^2+c_0)^{c_0/2}
 =r^{1+c_0}
   (1-r^{-2})^{1/2}(1+c_0r^{-2})^{c_0/2},
 \qquad c_0=\sqrt2-1,
\]
whereas \(\widehat v=r^2\).  Thus
\[
 \widehat u=r^{\sqrt2}\bigl(1+O_{\rm sym}(r^{-2})\bigr),
 \qquad
 |(r\partial_r)^jO_{\rm sym}(r^{-2})|
 \leq C_jr^{-2}
\]
for every \(j\).  The exact FIK profile in
Lemma~\ref{lem:self-contained-FIK-ledger},
\[
 F(r)=\frac1{\sqrt2}-\frac{c_0}{r^2}
      -\frac{c_0}{\sqrt2\,r^4},
\]
and the corresponding exact radial metric coefficients are
classical symbols under \(r\partial_r\).  The Wigner harmonics are
smooth on the compact principal orbit.  Consequently each covariant
derivative lowers radial order by one.  The positive potentials give
\[
 \bar\nabla\phi=O(r^{\sqrt2-1}),\quad
 \bar\nabla^2\phi=O(r^{\sqrt2-2}),\quad
 \bar\nabla^3\phi=O(r^{\sqrt2-3}),
\]
and the quadratic potentials give $O(r)$, $O(1)$, and $O(r^{-1})$,
respectively.  Repeated differentiation gives
\[
 |\bar\nabla^mW_a|\leq C_m(1+r)^{1-m},
 \qquad
 |\bar\nabla^mZ_a|\leq C_m(1+r)^{-m}.
\]
Since \(1+\bar f\simeq1+r^2\), these are
\eqref{eq:W-all-order}--\eqref{eq:Z-all-order}.  Smooth extension
across the exceptional divisor and the required core bounds follow
from Lemma~\ref{lem:self-contained-bolt-reality}.  For \(Z_0\), the curvature
components are algebraic combinations of the displayed exact profile
and its radial derivatives.  The same symbol calculation gives
\(|\bar\nabla^m\Ric_{\bar g}|=O(r^{-2-m})\), which is stronger than
\eqref{eq:Z-all-order}.
\end{proof}

\section{A localized geometric moving frame}
\label{sec:frame}

\subsection{Cutoff generators}

Retain the cutoffs \(\rho\) and \(\chi\) fixed before
Proposition~\ref{prop:intro-static-sliced-chart}; thus
\begin{equation}\label{eq:cutoffs}
 \begin{array}{lll}
 \rho=1\text{ on }[0,1],&
 \rho=0\text{ on }[2,\infty),\\
 \chi=1\text{ on }[0,5/2],&
 \chi=0\text{ on }[3,\infty).
 \end{array}
\end{equation}
For $\tau\geq1$, set
\[
 \rho_\tau(x)=\rho(e^{-\tau}\bar f(x)),\qquad
 \chi_\tau(x)=\chi(e^{-\tau}\bar f(x)).
\]
For $1\leq j\leq8$, use the fields \(W_j\) fixed above and define the
global geometric generators
\begin{equation}\label{eq:global-Y}
 Y_0=2\Ric_{\bar g}
     =\bar g-\Lie_{\bar\nabla\bar f}\bar g,\qquad
 Y_j=\Lie_{W_j}\bar g=\mathfrak c_j\bar\nabla^2\phi_j.
\end{equation}
By the opening normalization, the ordered list
$Y_0,\ldots,Y_8$ is a basis of $\Z$.  This choice
diagonalizes the scale direction: $Y_0$ spans the eigenvalue-$1$
space, while the $Y_j$, $j\geq1$, are the eight pure-gauge
eigentensors.

The corresponding direct generators and linear actions are
\begin{align}
 \widetilde Y_{0,\tau}&=Y_0,
 &\widetilde\B_{0,\tau}h
   &=h-\Lie_{\bar\nabla\bar f}h,
 \label{eq:localized-Y0}\\
 \widetilde Y_{j,\tau}
 &=\Lie_{\chi_\tau W_j}\bar g,
 &\widetilde\B_{j,\tau}h
 &=\Lie_{\chi_\tau W_j}h,
 \qquad1\leq j\leq8.
 \label{eq:localized-Yj}
\end{align}
On $\supp\rho_\tau$, the diffeomorphism fields agree with the global
ones.  The scale action need not be cut off analytically because its
direct tensor is $2\Ric_{\bar g}=O(\bar f^{-1})$.  In a closed-flow
graft it is localized geometrically by cancellation against the
adaptive outer model; see Proposition~\ref{prop:adaptive-graft}.

\begin{remark}[Scale--radial normalization]\label{rem:scale}
Let \(G(t)\) be an unnormalized Ricci flow, let \(\lambda(t)>0\), put
\(d\tau/dt=\lambda^{-1}\), and let
\(\Upsilon_t:M\to M\) be the diffeomorphism used for pullback.  Set
\[
 g=\lambda^{-1}\Upsilon_t^*G(t).
\]
If
\[
 Y_t^\Upsilon
 :=(\partial_t\Upsilon_t)\circ\Upsilon_t^{-1}
\]
is its Eulerian generator, then the normalized pullback generator is
\[
 \mathscr V
 :=\lambda\,\Upsilon_t^*Y_t^\Upsilon
 =\lambda\,(d\Upsilon_t)^{-1}\partial_t\Upsilon_t .
\]
Then
\[
 \partial_\tau g
 =-2\Ric_g-\lambda_tg+\Lie_{\mathscr V}g.
\]
The symbol \(\Upsilon_t\) denotes the pullback diffeomorphism.  In
Proposition~\ref{prop:controlled-HMH}, by contrast, \(\Phi_t\) denotes
the forward harmonic-map variable and
\[
 \Upsilon_t=\Phi_t^{-1}.
\]
With this identification,
\[
 \mathscr V
 =\lambda\,(\Phi_t)_*
   \bigl(\partial_t\Phi_t^{-1}\circ\Phi_t\bigr),
\]
which is denoted \(\lambda W\) in the proof of that proposition.  Thus
no inversion or Lie-derivative sign change is implicit between the two
conventions.

The soliton is stationary when \(\lambda_t=-1\) and
\(\mathscr V=-\bar\nabla\bar f\).  If
\(\lambda_t=-(1+a)\) and the radial part of \(\mathscr V\) is changed
from \(-\bar\nabla\bar f\) to
\(-(1+a)\bar\nabla\bar f\), the additional term is
\[
 a\bigl(g-\Lie_{\bar\nabla\bar f}g\bigr).
\]
At the background this is \(2a\Ric_{\bar g}\).  Thus the scale column is
the actual Ricci eigenvector rather than a mixture of that vector with
the radial Hessian mode.  The exact controlled construction is given
in Section~\ref{sec:controlled}.
\end{remark}

\subsection{Gaussian approximation of the global modes}

\begin{lemma}[Cutoff tails]\label{lem:cutoff-tails}
For each integer $m\geq0$ there are $C_m,c_m>0$ such that
\begin{equation}\label{eq:cutoff-tail}
 \norm{\widetilde Y_{j,\tau}-Y_j}_{H^m_\nu}
 \leq C_m e^{-c_m e^\tau},
 \qquad 0\leq j\leq8,\quad \tau\geq1.
\end{equation}
Moreover,
\begin{equation}\label{eq:gram-tail}
 \left|\ip{\widetilde Y_{j,\tau}-Y_j}{Z_a}\right|
 \leq C e^{-c e^\tau}.
\end{equation}
\end{lemma}

\begin{proof}
For $j=0$ the difference vanishes.  Every other difference in
\eqref{eq:cutoff-tail} is supported where $\bar f\geq2e^\tau$.
Lemma~\ref{lem:mode-growth}, the estimates
\[
 \abs{\bar\nabla^j\chi_\tau}
 \leq C_j e^{-j\tau/2}
 \quad\text{on }\{2e^\tau\leq\bar f\leq3e^\tau\},
\]
and the polynomial derivative bounds for the asymptotically conical
geometry show that the pointwise norm of each derivative grows at
most polynomially in $\bar f$.  The Gaussian tail estimate, verified
directly in Appendix~\ref{app:FIK-Gaussian-tail},
\begin{equation}\label{eq:general-Gaussian-tail}
 \int_{\{\bar f\geq R\}}(1+\bar f)^N\,d\nu
 \leq C_N(1+R)^{N+2}e^{-R/2}
\end{equation}
for $R$ sufficiently large then proves \eqref{eq:cutoff-tail}, after
changing the constants.  Equation~\eqref{eq:gram-tail} follows by
Cauchy--Schwarz.
\end{proof}

\subsection{The static finite-dimensional slice}

At a fixed time, let $\Theta_s$ be the flow of
$-\bar\nabla\bar f$ and define the scale--radial action
\begin{equation}\label{eq:finite-scale-action}
 \mathfrak S_sg=e^s\Theta_s^*g.
\end{equation}
It satisfies
\[
 \left.\frac d{ds}\right|_{s=0}\mathfrak S_s\bar g
 =\bar g-\Lie_{\bar\nabla\bar f}\bar g
 =2\Ric_{\bar g}.
\]
On the asymptotic cone, the homothety and radial pullback cancel.
Together with the complete flows of $W_1,\ldots,W_8$, this gives the
following static statement.

Indeed, \eqref{eq:W-all-order} and completeness of \((M,\bar g)\)
imply that every \(W_j\) is complete.  Let \(\Psi_j(s)\) denote its
flow and, in the fixed order displayed below, define
\begin{equation}\label{eq:static-action-definition}
 \mathfrak A_p g
 =\Psi_8(p_8)^*\cdots\Psi_1(p_1)^*
   \mathfrak S_{p_0}g.
\end{equation}

\begin{lemma}[Static FIK geometric slice]\label{lem:slice}
Fix $k\geq3$, $0<\alpha<1$, and a compact set $K\Subset M$.  In a
sufficiently small $C^{k,\alpha}$ neighborhood of $\bar g$ consisting
of metrics $\bar g+u$ with $\supp u\subset K$, there is a unique small
ordered scale--diffeomorphism parameter $p\in\R^9$ such that
\[
 h(p,u)=\mathfrak A_p(\bar g+u)-\bar g
\quad\text{satisfies}\quad
 \ip{h(p,u)}{Z_a}=0,\qquad0\leq a\leq8.
\]
\end{lemma}

\begin{proof}
For small $p$, the scalar moment map
\[
 F_a(p,u)=\ip{\mathfrak A_p(\bar g+u)-\bar g}{Z_a}
\]
is $C^1$.  Indeed, the generating fields have at most linear growth,
their small flows and the scale--radial action have uniformly
polynomially bounded derivatives, and the Gaussian weight supplies an
integrable majorant.  Its $p$-derivative at $(0,0)$ is the matrix
\[
 G_{aj}=\ip{Y_j}{Z_a}.
\]
It is invertible because $Y_0,\ldots,Y_8$ and
$Z_0,\ldots,Z_8$ are two bases of $\Z$.  The
finite-dimensional implicit function theorem proves the claim.
\end{proof}

\begin{remark}
Lemma~\ref{lem:slice} only prepares the initial slice.  The feedback
law in Theorem~\ref{thm:receding} propagates its receding-domain
analogue.  What remains outside that theorem is continuation of the
coupled geometric chart and the pointwise three-region box.
\end{remark}

\section{Controlled gauges and receding cutoffs}
\label{sec:controlled}

\subsection{The exact variable-scale equation}

We first record the geometric calculation which produces the diagonal
scale action.  It is an identity for a general metric family and does
not use linearization.

\begin{proposition}[Controlled harmonic-map-gauge identity]
\label{prop:controlled-HMH}
Let $\acute G(t)$ be a smooth complete metric family on $M$, let
$\lambda(t)>0$, and define
\[
 \frac{d\tau}{dt}=\lambda^{-1},\qquad
 \lambda_t=-(1+a).
\]
Let $U(\tau)$ be a time-dependent vector field and suppose that
$\Phi_t:M\to M$ is a family of diffeomorphisms satisfying
\begin{equation}\label{eq:controlled-HMH}
 \partial_t\Phi
 =\Delta_{\acute G(t),\bar g}\Phi
  +\lambda^{-1}\bigl((1+a)\bar\nabla\bar f-U\bigr)\circ\Phi.
\end{equation}
Set
\begin{equation}\label{eq:controlled-definitions}
 g=\lambda^{-1}(\Phi^{-1})^*\acute G,\qquad
 h=g-\bar g,\qquad
 \E=(\Phi^{-1})^*
       \bigl(\partial_t\acute G+2\Ric_{\acute G}\bigr).
\end{equation}
Our sign convention for the DeTurck vector field is
\begin{equation}\label{eq:DeTurck-vector-definition}
 B_{\bar g}(g)^k
 =g^{ij}\bigl(\Gamma(g)^k_{ij}-\Gamma(\bar g)^k_{ij}\bigr).
\end{equation}
Then
\begin{equation}\label{eq:controlled-g}
 \begin{split}
 \partial_\tau g={}&-2\Ric_g+\Lie_{B_{\bar g}(g)}g
 -\Lie_{\bar\nabla\bar f}g+g+\E\\
 &+a\bigl(g-\Lie_{\bar\nabla\bar f}g\bigr)+\Lie_Ug,
 \end{split}
\end{equation}
Equivalently,
\begin{equation}\label{eq:controlled-h}
 \partial_\tau h
 =\A h+\Q(h)+\E
 +a\bigl(Y_0+\widetilde\B_{0,\tau}h\bigr)
 +\Lie_U\bar g+\Lie_Uh.
\end{equation}
There is no factor of $\lambda$ multiplying $\E$ in the normalized
equation.
\end{proposition}

\begin{proof}
Let
\[
 W=\Phi_*\bigl(\partial_t\Phi^{-1}\circ\Phi\bigr)
   =-(\partial_t\Phi)\circ\Phi^{-1}.
\]
Covariance and scaling of the tension field give
\[
 (\Delta_{\acute G,\bar g}\Phi)\circ\Phi^{-1}
 =\Delta_{\lambda g,\bar g}\operatorname{Id}
 =-\lambda^{-1}B_{\bar g}(g).
\]
Consequently
\[
 \lambda W
 =B_{\bar g}(g)-(1+a)\bar\nabla\bar f+U.
\]
Differentiating \eqref{eq:controlled-definitions} with
$dt/d\tau=\lambda$ gives
\[
 \partial_\tau g
 =-\lambda_tg+\Lie_{\lambda W}g-2\Ric_g+\E.
\]
Substitution proves \eqref{eq:controlled-g}.  The shrinker equation
implies
\[
 -2\Ric_{\bar g}
 -\Lie_{\bar\nabla\bar f}\bar g+\bar g=0.
\]
The linearization of the first line of
\eqref{eq:controlled-g} at $\bar g$ is $\A$, and
\[
 \bar g-\Lie_{\bar\nabla\bar f}\bar g=2\Ric_{\bar g}=Y_0,
\]
which proves \eqref{eq:controlled-h}.
\end{proof}

Take
\begin{equation}\label{eq:velocity-cutoff}
 U(\tau)=\sum_{j=1}^8b_j(\tau)\chi_\tau W_j.
\end{equation}
The cutoff is inserted in the instantaneous velocity.  Therefore
\[
 \frac d{d\tau}\Psi_\tau^*g
 =\Psi_\tau^*\bigl(\partial_\tau g+\Lie_{U(\tau)}g\bigr)
\]
for the flow of $U$, and no term containing
$\partial_\tau\chi_\tau$ occurs.  Such a term would arise from
differentiating a finite-position parametrization, not from the
velocity formulation \eqref{eq:velocity-cutoff}.

\subsection{Cancellation at the moving boundary}

\begin{lemma}[Drift-adapted cutoff identity]
\label{lem:cutoff-cancellation}
For $\rho_\tau=\rho(e^{-\tau}\bar f)$,
\begin{equation}\label{eq:cutoff-cancellation}
 (\partial_\tau-\bar\Delta_{\bar f})\rho_\tau
 =-2e^{-\tau}\rho'(e^{-\tau}\bar f)
  -e^{-2\tau}\rho''(e^{-\tau}\bar f)
       |\bar\nabla\bar f|^2.
\end{equation}
In particular,
\[
 |(\partial_\tau-\bar\Delta_{\bar f})\rho_\tau|
 \leq Ce^{-\tau},\qquad
 |\bar\nabla\rho_\tau|\leq Ce^{-\tau/2},
\]
and both terms are supported in
$\{e^\tau\leq\bar f\leq2e^\tau\}$.
\end{lemma}

\begin{proof}
Writing $s=e^{-\tau}\bar f$ and using
$\bar\Delta_{\bar f}\bar f=2-\bar f$, one has
\[
 \partial_\tau\rho_\tau=-s\rho'(s)
\]
and
\[
 \bar\Delta_{\bar f}\rho_\tau
 =e^{-\tau}\rho'(s)(2-\bar f)
  +e^{-2\tau}\rho''(s)|\bar\nabla\bar f|^2.
\]
The two terms containing $s\rho'(s)$ cancel.  The stated bounds follow
from $|\bar\nabla\bar f|^2\leq\bar f$.
\end{proof}

Let $H=\rho_\tau h$.  On $\supp\rho_\tau$, the nesting of the
cutoffs gives $\chi_\tau=1$.  Equations~
\eqref{eq:controlled-h} and \eqref{eq:cutoff-cancellation} imply
\begin{equation}\label{eq:H-equation}
 \begin{split}
 \partial_\tau H={}&\A H+\rho_\tau(\Q(h)+\E)\\
 &+a\rho_\tau\bigl(Y_0+\widetilde\B_{0,\tau}h\bigr)
 +\sum_{j=1}^8b_j\rho_\tau
   \bigl(Y_j+\Lie_{W_j}h\bigr)
 +\mathcal C_\rho[h],
 \end{split}
\end{equation}
where
\begin{equation}\label{eq:C-rho}
 \mathcal C_\rho[h]
 =\bigl((\partial_\tau-\bar\Delta_{\bar f})\rho_\tau\bigr)h
  -2\bar\nabla_{\bar\nabla\rho_\tau}h.
\end{equation}

\begin{lemma}[Gaussian size of the boundary commutator]
\label{lem:C-rho}
Suppose $q$ and $\bar\nabla q$ have polynomial growth.  There are
$C,c>0$ such that
\begin{equation}\label{eq:C-rho-pairing}
 |\ip{\mathcal C_\rho[h]}q|
 \leq Ce^{-ce^\tau}
 \norm h_{C^1(\{e^\tau<\bar f<2e^\tau\})}.
\end{equation}
Moreover, with \(H=\rho_\tau h\),
\begin{equation}\label{eq:C-rho-energy-pairing}
 \left|\ip{\mathcal C_\rho[h]}H\right|
 \leq Ce^{-ce^\tau}
 \norm h_{C^1(\{e^\tau<\bar f<2e^\tau\})}^2.
\end{equation}
\end{lemma}

\begin{proof}
Both terms in \eqref{eq:C-rho} are supported where
$\bar f\simeq e^\tau$.  Lemma~\ref{lem:cutoff-cancellation} gives a
polynomial pointwise bound, while the measure contains
$e^{-\bar f}$.  The general Gaussian tail
\eqref{eq:general-Gaussian-tail} proves
\eqref{eq:C-rho-pairing}.  For the second assertion, put
\[
 A_\tau:=\{e^\tau<\bar f<2e^\tau\}.
\]
Since \(H=\rho_\tau h\), equations \eqref{eq:C-rho} and
\eqref{eq:cutoff-cancellation} give on \(A_\tau\)
\[
 |\mathcal C_\rho[h]|\,|H|
 \leq
 C\bigl(e^{-\tau}|h|+e^{-\tau/2}|\bar\nabla h|\bigr)|h|\,
 \mathbf 1_{A_\tau}.
\]
Consequently,
\[
 \left|\ip{\mathcal C_\rho[h]}H\right|
 \leq
 C\bigl(e^{-\tau}+e^{-\tau/2}\bigr)
 \|h\|_{C^1(A_\tau)}^2\,\nu(A_\tau)
 \leq
 Ce^{-ce^\tau}\|h\|_{C^1(A_\tau)}^2,
\]
where the last inequality follows from
\eqref{eq:general-Gaussian-tail}, after changing \(C,c\).  This proves
\eqref{eq:C-rho-energy-pairing}.
\end{proof}

\section{Nonlinear weighted estimates}
\label{sec:quadratic}

\begin{lemma}[Quadratic Ricci--DeTurck estimates]
\label{lem:quadratic}
There are $\varepsilon_0,C>0$ such that, if
\[
 \norm h_{C^1(\bar g)}\leq\varepsilon_0,\qquad
 h\in C^\infty(S^2T^*M)\cap H^1_\nu,
\]
then \(\Q(h)\) is understood as an \(H^{-1}_\nu\) distribution: its
second-order term is defined by one weighted integration by parts
against compactly supported tests and passage to the \(H^1_\nu\)
closure.  With every displayed pairing interpreted as this
\(H^{-1}_\nu\)--\(H^1_\nu\) duality, one has
\begin{align}
 \left|\ip{\Q(h)}h\right|
 &\leq C\|h\|_{C^1(\bar g)}\norm h_{H^1_\nu}^2,
 \label{eq:Q-energy}\\
 \left|\ip{\Q(h)}{Z_a}\right|
 &\leq C\norm h_{H^1_\nu}^2,
 \qquad0\leq a\leq8.
 \label{eq:Q-projection}
\end{align}
\end{lemma}

\begin{proof}
Let \(\zeta_R\) be a standard exhaustion cutoff with
\(|\bar\nabla\zeta_R|\leq C R^{-1}\), and first test against
\(\zeta_Rq\).  The \(C^1\) bound, the first-moment estimate, and
\(h\in H^1_\nu\) show that all terms in the resulting
integration-by-parts form are Cauchy as \(R\to\infty\); the terms
containing \(\bar\nabla\zeta_R\) tend to zero.  This defines the
claimed \(H^{-1}_\nu\) distribution and justifies the following
calculation without an unstated boundary term.

The lower-order terms in \eqref{eq:Q-schematic} follow from bounded
curvature, Cauchy--Schwarz, and $C^1$ smallness.  For the leading term,
integrate one derivative in
\[
 \int_M\left\langle
 -\widehat h^{ij}\bar\nabla_i\bar\nabla_jh,q\right\rangle d\nu
\]
by parts, with $q=h$ or $q=Z_a$.  When the derivative falls on
$\widehat h$, $h$, or $Z_a$, the claimed bound follows from
Lemma~\ref{lem:mode-growth}.  When it falls on the Gaussian density,
the new integral is bounded by
\[
 C\int_M\sqrt{\bar f}\,\abs h\,
                \abs{\bar\nabla h}\,\abs q\,d\nu.
\]
For $q=Z_a$, use \eqref{eq:Z-growth} and
Lemma~\ref{lem:first-moment}.  For $q=h$, take one factor of $h$ in
$L^\infty$ and use the same moment estimate.  This proves
\eqref{eq:Q-energy}--\eqref{eq:Q-projection}.
\end{proof}

\begin{lemma}[Energy of the geometric action]\label{lem:B-energy}
There is $C<\infty$ such that, for all $\tau\geq1$ and all smooth
compactly supported $h$,
\begin{align}
 \left|\ip{\widetilde\B_{a,\tau}h}{h}\right|
 &\leq C\norm h_{H^1_\nu}^2,\label{eq:B-energy}\\
 \left|\ip{\widetilde\B_{a,\tau}h}{Z_b}\right|
 &\leq C\norm h_{H^1_\nu},
 \label{eq:B-projection}
\end{align}
uniformly for $0\leq a,b\leq8$.
\end{lemma}

\begin{proof}
For $a\geq1$, expand the Lie derivative into
\[
 \Lie_{\chi_\tau W_a}h
 =\bar\nabla_{\chi_\tau W_a}h+
   (\bar\nabla(\chi_\tau W_a))*h.
\]
The second term is bounded because
$\abs{\bar\nabla(\chi_\tau W_a)}\leq C$.
For the transport term, direct Cauchy--Schwarz gives
\[
 \left|\int_M
 \langle\bar\nabla_{\chi_\tau W_a}h,h\rangle\,d\nu\right|
 \leq C\norm{\bar\nabla h}_{L^2_\nu}
       \norm{(1+\sqrt{\bar f})h}_{L^2_\nu}
 \leq C\norm h_{H^1_\nu}^2
\]
by Lemma~\ref{lem:first-moment}.  Similarly,
\[
 \left|\int_M
 \langle\bar\nabla_{\chi_\tau W_a}h,Z_b\rangle\,d\nu\right|
 \leq C\norm{\bar\nabla h}_{L^2_\nu}
       \norm{(1+\sqrt{\bar f})Z_b}_{L^2_\nu},
\]
and the last factor is finite by
Lemma~\ref{lem:mode-growth}.  The remaining zeroth-order terms follow
from \eqref{eq:V-growth}--\eqref{eq:Z-growth}.

For $a=0$, use the weighted divergence
\(
 \operatorname{div}_{\bar f}X
 :=e^{\bar f}\operatorname{div}(e^{-\bar f}X)
\), so that
\[
 \operatorname{div}_{\bar f}(\bar\nabla\bar f)
 =\bar\Delta_{\bar f}\bar f=2-\bar f.
\]
Weighted integration by parts therefore gives
\[
 \left|\int_M
 \langle\bar\nabla_{\bar\nabla\bar f}h,h\rangle\,d\nu\right|
 \leq C\int_M(1+\bar f)|h|^2\,d\nu
 \leq C\norm h_{H^1_\nu}^2.
\]
The terms in the Lie derivative involving
$\bar\nabla^2\bar f$ are zeroth order and bounded.  Pairing the
transport with $Z_b$ and integrating it off $h$ produces at worst
\[
 \int_M\bar f|h||Z_b|\,d\nu
 +\int_M\sqrt{\bar f}|h||\bar\nabla Z_b|\,d\nu
 +C\int_M|h||Z_b|\,d\nu.
\]
The first term is bounded by
$\norm{\sqrt{\bar f}h}_{L^2_\nu}
 \norm{\sqrt{\bar f}Z_b}_{L^2_\nu}$; all remaining mode factors have
finite Gaussian norm.  Lemma~\ref{lem:first-moment} proves the scale
estimates.  Thus both the linearly growing diffeomorphism transports
and the radial scale transport are controlled as bilinear forms on
$H^1_\nu$; no bounded-operator assertion
$H^1_\nu\to L^2_\nu$ is needed.
\end{proof}

\subsection{Localized tame estimates and modal forcing}
\label{subsec:localized-tame}

The cutoff argument requires localized versions of
Lemmas~\ref{lem:quadratic} and~\ref{lem:B-energy}.  We record them
explicitly.  Enlarge the transition annulus slightly by setting
\[
 \mathcal A_\tau^*
 =
 \left\{\frac12e^\tau<\bar f<3e^\tau\right\},
 \qquad
 \mathfrak a_m(h;\tau)
 =
 \sum_{\ell=0}^m
 \|\bar\nabla^\ell h\|_{L^\infty(\mathcal A_\tau^*)}.
\]
In this subsection, \(q_\tau\) denotes a quantity satisfying
\[
 0\leq q_\tau\leq Ce^{-ce^\tau};
\]
each occurrence below may denote a different nonnegative remainder
dominated by an envelope of this form, and the constants
\(C,c>0\) may change from line to line.  Define the uncut geometric
action operators on all of \(M\) by
\[
 \mathscr B_0h=h-\Lie_{\bar\nabla\bar f}h,
 \qquad
 \mathscr B_jh=\Lie_{W_j}h,\quad 1\leq j\leq8.
\]
Their restrictions to \(\operatorname{supp}\rho_\tau\) are the actions
entering the localized estimates below.

\begin{lemma}[Localized tame estimates]
\label{lem:localized-tame}
There are \(\delta_0,C,c>0\) with the following property.  Let
\(\tau\geq1\), let \(h\) be smooth on
\(\{\bar f<4e^\tau\}\), and suppose
\[
 \sum_{\ell=0}^2
 \|\bar\nabla^\ell h\|_{L^\infty(\{\bar f<4e^\tau\})}
 \leq\delta_0.
\]
Set \(H=\rho_\tau h\).  Then, for \(0\leq\mu,j\leq8\),
\begin{align}
 \left|\ip{\rho_\tau\Q(h)}{Z_\mu}\right|
 &\leq
 C\|H\|_{H^1_\nu}^2
+q_\tau\mathfrak a_2(h;\tau)^2,
 \label{eq:localized-Q-proj}\\
 \left|\ip{\rho_\tau\Q(h)}H\right|
 &\leq
 C\|h\|_{C^1(\{\bar f<4e^\tau\})}
 \|H\|_{H^1_\nu}^2
+q_\tau\mathfrak a_2(h;\tau)^2,
 \label{eq:localized-Q-energy}\\
 \left|\ip{\rho_\tau\mathscr B_jh}{Z_\mu}\right|
 &\leq
 C\|H\|_{H^1_\nu}
+q_\tau\mathfrak a_1(h;\tau),
 \label{eq:localized-B-proj}\\
 \left|\ip{\rho_\tau\mathscr B_jh}H\right|
 &\leq
 C\|H\|_{H^1_\nu}^2
+q_\tau\mathfrak a_1(h;\tau)^2,
 \label{eq:localized-B-energy}\\
 \left|\ip{\mathcal C_\rho[h]}{Z_\mu}\right|
 &\leq q_\tau\mathfrak a_1(h;\tau),
 \label{eq:localized-C-proj}\\
 \left|\ip{\mathcal C_\rho[h]}H\right|
 &\leq q_\tau\mathfrak a_1(h;\tau)^2.
 \label{eq:localized-C-energy}
\end{align}
The direct columns satisfy
\begin{equation}
 \left|
 \ip{\rho_\tau Y_j}{Z_\mu}
 -
 \ip{Y_j}{Z_\mu}
 \right|
 \leq q_\tau.
 \label{eq:localized-column-proj}
\end{equation}
If, in addition, \(H\perp\mathcal Z\), then
\begin{equation}
 \left|\ip{\rho_\tau Y_j}H\right|
 \leq q_\tau\mathfrak a_0(h;\tau).
 \label{eq:localized-column-energy}
\end{equation}
\end{lemma}

\begin{proof}
Choose \(\kappa\in C^\infty_c([0,1))\) such that
\(\kappa=1\) on \([0,\frac12]\), and set
\(\kappa_\tau=\kappa(e^{-\tau}\bar f)\).  On
\(\operatorname{supp}\kappa_\tau\) one has \(h=H\).  Split every
integral into its \(\kappa_\tau\)-part and its complement.

On the \(\kappa_\tau\)-part, integrate one derivative off the
quasilinear term
\[
 -\widehat h^{ij}\bar\nabla_i\bar\nabla_jh.
\]
The terms for which the derivative falls on \(h,H\), or \(Z_\mu\)
are estimated exactly as in Lemma~\ref{lem:quadratic}.  When it falls
on the Gaussian density, the resulting first moment is controlled by
Lemma~\ref{lem:first-moment}.  A derivative of \(\kappa_\tau\) is
supported in \(\mathcal A_\tau^*\), and hence contributes to the
Gaussian tail.  This proves the core contributions in
\eqref{eq:localized-Q-proj}--\eqref{eq:localized-Q-energy}.

The complementary part is supported where
\(\bar f\geq\frac12e^\tau\).  The schematic expansion
\eqref{eq:Q-schematic} and the \(C^2\) bound give
\[
 |\Q(h)|
 \leq
 C\bigl(|h|\,|\bar\nabla^2h|
       +|\bar\nabla h|^2+|h|^2\bigr).
\]
The mode and generator bounds in Lemma~\ref{lem:mode-growth}, together
with the Gaussian tail estimate
\eqref{eq:general-Gaussian-tail}, absorb every polynomial factor in
\(\bar f\) and give the second terms in
\eqref{eq:localized-Q-proj}--\eqref{eq:localized-Q-energy}.

The same decomposition applies to the geometric actions.  On the
core, use the integrations by parts in Lemma~\ref{lem:B-energy}.
Equivalently, one may use the exact identities
\begin{align*}
 \rho_\tau\Lie_{W_j}h
 &=
 \Lie_{W_j}H-(W_j\rho_\tau)h,
 \qquad 1\leq j\leq8,\\
 \rho_\tau
 \bigl(h-\Lie_{\bar\nabla\bar f}h\bigr)
 &=
 H-\Lie_{\bar\nabla\bar f}H
+(\bar\nabla_{\bar\nabla\bar f}\rho_\tau)h.
\end{align*}
The last terms are supported in the transition annulus.  There
\[
 |\bar\nabla^m\rho_\tau|
 \leq C_me^{-m\tau/2},\qquad m=1,2,
\]
whereas \(W_j\), \(\bar\nabla\bar f\), and their required derivatives
grow at most polynomially.  The Gaussian tail proves
\eqref{eq:localized-B-proj}--\eqref{eq:localized-B-energy}.
Equations~\eqref{eq:localized-C-proj}--\eqref{eq:localized-C-energy}
follow directly from Lemma~\ref{lem:cutoff-cancellation}, since
\(\mathcal C_\rho[h]\) is supported where
\(\bar f\simeq e^\tau\).

Finally, \((\rho_\tau-1)Y_j\) is relevant only where
\(\bar f\geq e^\tau\).  This proves
\eqref{eq:localized-column-proj}.  If \(H\perp\mathcal Z\), then
\[
 \ip{\rho_\tau Y_j}H
 =
 \ip{(\rho_\tau-1)Y_j}H,
\]
and the same argument proves \eqref{eq:localized-column-energy}.
\end{proof}

\begin{remark}[Lipschitz forms]
\label{rem:localized-Lipschitz}
The preceding proof also gives the difference estimates needed below.
If \(h_1,h_2\) satisfy the same \(C^2\) box, \(v=h_1-h_2\), and
\(V=\rho_\tau v\), then
\begin{align}
 \left|
 \ip{\rho_\tau(\Q(h_1)-\Q(h_2))}{Z_\mu}
 \right|
 &\leq
 C\delta_0\|V\|_{H^1_\nu}
+q_\tau\mathfrak a_2(v;\tau),
 \label{eq:localized-Q-Lip-proj}\\
 \left|
 \ip{\rho_\tau(\Q(h_1)-\Q(h_2))}V
 \right|
 &\leq
 C\delta_0\|V\|_{H^1_\nu}^2
+q_\tau\mathfrak a_2(v;\tau)^2.
 \label{eq:localized-Q-Lip-energy}
\end{align}
Indeed, expand the difference of the principal terms as
 \[
 \widehat h_1*\bar\nabla^2v
+(\widehat h_1-\widehat h_2)*\bar\nabla^2h_2
+\text{lower-order terms}
\]
and integrate one derivative in each summand.  The geometric actions
also satisfy the bilinear estimate
\begin{equation}
 \left|\ip{\mathscr B_ju}v\right|
 \leq C\|u\|_{H^1_\nu}\|v\|_{H^1_\nu},
 \label{eq:B-bilinear}
\end{equation}
first for compactly supported tensors and then by density.  The
localized version has the corresponding \(q_\tau\)-weighted annular
remainder.  Thus all difference remainders vanish when \(h_1=h_2\);
there is no solution-independent error in a two-solution estimate.
\end{remark}

\begin{lemma}[Coercivity on the geometric complement]
\label{lem:coercivity}
There is $c_{\A}>0$ such that
\begin{equation}\label{eq:H1-coercivity}
 -\mathfrak a[h,h]\geq c_{\A}\norm h_{H^1_\nu}^2
\end{equation}
whenever $h\in H^1_\nu$ and $h\perp\Z$.
\end{lemma}

\begin{proof}
Let $K=\norm{\overline{\Rm}}_{L^\infty}$.  Equations~
\eqref{eq:A-form} and \eqref{eq:beta} give
\[
 \norm{\bar\nabla h}_{L^2_\nu}^2
 \leq-\mathfrak a[h,h]+2K\norm h_{L^2_\nu}^2
 \leq\left(1+\frac{2K}{\beta}\right)
      \left[-\mathfrak a[h,h]\right].
\]
The $L^2_\nu$ term is bounded by
$\beta^{-1}[-\mathfrak a[h,h]]$.  Adding the two estimates proves the
claim.
\end{proof}

\section{The modulation system}
\label{sec:modulation}

\begin{theorem}[Projected energy and exact feedback]
\label{thm:abstract-energy}
Let $I$ be an interval and let $Y_j^0$, $0\leq j\leq8$, be any basis
of $\Z$.  Suppose $h$ is a smooth tensor satisfying
\begin{equation}\label{eq:slice}
 \ip{h(\tau)}{Z_a}=0,\qquad0\leq a\leq8,
\end{equation}
and, weakly in $H^{-1}_\nu$,
\begin{equation}\label{eq:modulated}
 \partial_\tau h
 =\A h+\Q(h)+\E+
 \sum_{j=0}^8b_j(\tau)
 \bigl(Y_j(\tau)+\B_j(\tau)h\bigr).
\end{equation}
Assume that
\begin{align}
 \norm{Y_j(\tau)-Y_j^0}_{L^2_\nu}
 &\leq\delta_*,\label{eq:Y-defect}\\
 \left|\ip{\B_j(\tau)h}{Z_a}\right|
 &\leq C\norm h_{H^1_\nu},\label{eq:abstract-B-proj}\\
 \left|\ip{\B_j(\tau)h}{h}\right|
 &\leq C\norm h_{H^1_\nu}^2,\label{eq:abstract-B-energy}
\end{align}
and that \eqref{eq:Q-energy}--\eqref{eq:Q-projection} hold.  Then, for
$\delta_*$ and
\[
 \sup_I\bigl(\norm h_{C^1}+\norm h_{H^1_\nu}\bigr)
\]
sufficiently small, differentiating the slice determines the full
coefficient vector \(\mathbf b:=(b_0,\ldots,b_8)\) uniquely and
\begin{equation}\label{eq:b-bound}
 \abs{\mathbf b}\leq C\left(\norm h_{H^1_\nu}^2+
                             \norm\E_{H^{-1}_\nu}\right).
\end{equation}
For every $0<\sigma<\beta$, after decreasing the smallness constants,
\begin{equation}\label{eq:energy}
 \frac{d}{d\tau}\norm h_{L^2_\nu}^2
 +2\sigma\norm h_{L^2_\nu}^2
 +c\norm h_{H^1_\nu}^2
 \leq C\norm{\E}_{H^{-1}_\nu}^2.
\end{equation}
\end{theorem}

\begin{proof}
Differentiate \eqref{eq:slice}.  Self-adjointness and
$\A Z_a=\lambda_aZ_a$ give the exact cancellation
\[
 \ip{\A h}{Z_a}
 =\lambda_a\ip h{Z_a}=0.
\]
Thus \(\mathbf b\) solves
\begin{equation}\label{eq:Gram-system}
 \sum_{j=0}^8M_{aj}(\tau,h)b_j
 =-\ip{\Q(h)+\E}{Z_a},
\end{equation}
where
\[
 M_{aj}(\tau,h)
 =\ip{Y_j(\tau)+\B_j(\tau)h}{Z_a}.
\]
At $\delta_*=0$ and $h=0$, this is the change-of-basis matrix from
$Y_j^0$ to $Z_a$ and is invertible.  Hypotheses
\eqref{eq:Y-defect} and \eqref{eq:abstract-B-proj} make $M$ a small
perturbation of that matrix.  Equations~\eqref{eq:Q-projection} and
duality prove \eqref{eq:b-bound}.

Pair \eqref{eq:modulated} with $h$.  Since $Y_j^0\in\Z$,
\[
 \ip{Y_j^0}{h}=0.
\]
Writing $x=\norm h_{H^1_\nu}$ and
$e=\norm\E_{H^{-1}_\nu}$, the direct-generator defect and action terms
are bounded by
\[
 C\abs{\mathbf b}\,(\delta_*x+x^2)
 \leq C(x^2+e)(\delta_*x+x^2).
\]
For $\delta_*$, $x$, and $\norm h_{C^1}$ small, Young's inequality
absorbs this expression into the coercive term, leaving $Ce^2$.
Lemma~\ref{lem:quadratic} treats $\Q(h)$ and duality treats $\E$.

Finally, \eqref{eq:beta} gives
$2\ip{\A h}{h}\leq-2\beta\norm h_{L^2_\nu}^2$, while
Lemma~\ref{lem:coercivity} gives $H^1_\nu$ coercivity.  A convex
combination leaves the two positive terms in \eqref{eq:energy}.
\end{proof}

\begin{corollary}[Complete-space feedback]
\label{cor:localized-energy}
Suppose the controlled equation \eqref{eq:controlled-h} is available
on the complete FIK manifold and $h(\tau_0)$ satisfies
\eqref{eq:slice}.  Assume all regularity, weak
\(H^{-1}_\nu\)-equation, nonlinear, and forcing hypotheses of
Theorem~\ref{thm:abstract-energy}; in particular,
\(\sup_I(\|h\|_{C^1}+\|h\|_{H^1_\nu})\) and the column defect have the
smallness required there.  For $\tau_0$ sufficiently large, there is a
unique choice of the eight velocities
\(b_1,\ldots,b_8\) in \eqref{eq:velocity-cutoff}, together with the
scale velocity $a$, which propagates the slice.  Writing
\(b_0=a\) and \(b=(b_1,\ldots,b_8)\), the estimate
\eqref{eq:b-bound} is equivalently
\(
 |a|+|b|\leq
 C(\norm h_{H^1_\nu}^2+\norm\E_{H^{-1}_\nu})
\)
after changing its fixed constant, and \eqref{eq:energy} also holds.
\end{corollary}

\begin{proof}
Use $Y_j^0=Y_j$ and the diagonal scale data
\eqref{eq:localized-Y0}.  Lemma~\ref{lem:cutoff-tails} verifies
\eqref{eq:Y-defect}, and Lemma~\ref{lem:B-energy} verifies
 \eqref{eq:abstract-B-proj}--\eqref{eq:abstract-B-energy}.  The
 feedback coefficients are the unique solution of
\eqref{eq:Gram-system}.  More explicitly, if
\(m_\mu=\ip h{Z_\mu}\), the same computation off the slice gives
\[
 m_\mu'=\lambda_\mu m_\mu
\]
after the feedback has cancelled the nonlinear, forcing, and geometric
column contributions.  Since \(m_\mu(\tau_0)=0\), uniqueness for
this scalar ODE propagates every slice condition.
\end{proof}

\begin{theorem}[Receding-domain modulation and phase control]
\label{thm:receding}
Fix $0<\sigma<\beta$.  Let
\[
 \mathcal D_{\tau_0,\tau_1}
 =\{(x,\tau):\tau_0\leq\tau<\tau_1,\ 
                  \bar f(x)<4e^\tau\}.
\]
Suppose $h$ is smooth on $\mathcal D_{\tau_0,\tau_1}$, satisfies
\eqref{eq:controlled-h} wherever $\bar f<3e^\tau$, with $U$ given by
\eqref{eq:velocity-cutoff}, and obeys
\begin{equation}\label{eq:receding-box}
 \sup_{\mathcal D_{\tau_0,\tau_1}}
 \bigl(|h|+|\bar\nabla h|+|\bar\nabla^2h|\bigr)
 \leq\delta.
\end{equation}
Set $H=\rho_\tau h$.  If $\delta$ is sufficiently small and
$\tau_0$ sufficiently large, the nine conditions
\begin{equation}\label{eq:receding-slice}
 \ip H{Z_\mu}=0,\qquad0\leq\mu\leq8,
\end{equation}
determine $(a,b_1,\ldots,b_8)$ uniquely at every time.  More
precisely, they solve
\begin{equation}\label{eq:receding-Gram}
 \sum_{j=0}^8M_{\mu j}c_j=-d_\mu,
 \qquad c=(a,b_1,\ldots,b_8),
\end{equation}
where
\begin{align}
 M_{\mu0}
 &=\ip{\rho_\tau
       (Y_0+\widetilde\B_{0,\tau}h)}{Z_\mu},
 \label{eq:receding-M0}\\
 M_{\mu j}
 &=\ip{\rho_\tau
       (Y_j+\Lie_{W_j}h)}{Z_\mu},
 \qquad1\leq j\leq8,\label{eq:receding-Mj}\\
 d_\mu
 &=\ip{\rho_\tau(\Q(h)+\E)+\mathcal C_\rho[h]}{Z_\mu}.
 \label{eq:receding-d}
\end{align}
The matrix $M$ is uniformly invertible, and
\begin{equation}\label{eq:receding-velocity}
 |a|+|b|
 \leq C\left(
   \norm H_{H^1_\nu}^2
  +\norm{\rho_\tau\E}_{H^{-1}_\nu}
  +e^{-ce^\tau}
       \mathfrak a_2(h;\tau)
       \right).
\end{equation}
If \eqref{eq:receding-slice} holds initially and the velocities are
chosen by \eqref{eq:receding-Gram}, it is propagated for as long as
the smooth controlled solution and \eqref{eq:receding-box} exist.
On that interval,
\begin{equation}\label{eq:receding-energy}
 \frac d{d\tau}\norm H_{L^2_\nu}^2
 +2\sigma\norm H_{L^2_\nu}^2
 +c\norm H_{H^1_\nu}^2
 \leq C\norm{\rho_\tau\E}_{H^{-1}_\nu}^2
      +Ce^{-ce^\tau}.
\end{equation}
\end{theorem}

\begin{proof}
Differentiate \eqref{eq:receding-slice} and insert the exact equation
\eqref{eq:H-equation}.  Since $\A$ is self-adjoint and each
\(Z_\mu\) is an eigentensor,
\[
 \ip{\A H}{Z_\mu}
 =\ip H{\A Z_\mu}=0.
\]
The remaining terms give
\eqref{eq:receding-Gram}--\eqref{eq:receding-d}.

At $h=0$, the matrix $M$ converges as $\tau\to\infty$ to
\[
 G_{\mu j}=\ip{Y_j}{Z_\mu},
\]
which is invertible because both ordered families are bases of $\Z$.
The cutoff tail is $O(e^{-ce^\tau})$.  The $h$-dependent matrix
terms are $O(\delta)$: for the Lie and radial transports this follows
by the weighted integrations by parts in
Lemma~\ref{lem:B-energy}, with the transition annulus estimated using
\eqref{eq:receding-box}.  Hence $M^{-1}$ is uniformly bounded.

Lemma~\ref{lem:localized-tame} gives
\[
 |\ip{\rho_\tau\Q(h)}{Z_\mu}|
 \leq C\norm H_{H^1_\nu}^2
      +Ce^{-ce^\tau}.
\]
The forcing is bounded by duality, uniformly because
\(\rho_\tau Z_\mu\) is bounded in \(H^1_\nu\).
Lemma~\ref{lem:C-rho} treats the last term in
\eqref{eq:receding-d}.  Uniform inversion proves
\eqref{eq:receding-velocity}.

Now pair \eqref{eq:H-equation} with $H$.  Spectral coercivity applies
because $H\perp\Z$.  The quadratic, direct-generator, and action terms
obey Lemma~\ref{lem:localized-tame}; the pointwise box absorbs the
quadratic energy term, and every transition error is
superexponentially small.  Insert
 \eqref{eq:receding-velocity}, use Young's inequality, and then combine
 the spectral gap with Lemma~\ref{lem:coercivity}.  After decreasing
 $\delta$ and increasing $\tau_0$, this yields
\eqref{eq:receding-energy}.  Finally, writing
\(m_\mu(\tau)=\ip{H(\tau)}{Z_\mu}\), the feedback cancels every
term in \(m_\mu'\) except its spectral part, and hence
\[
 m_\mu'=\lambda_\mu m_\mu .
\]
An initially sliced solution therefore remains sliced by uniqueness
for this finite-dimensional diagonal ODE.
\end{proof}

\begin{remark}[Why the positive modes disappear]
The eigenvalues $1$ and $1-1/\sqrt2$ do not occur in either
\eqref{eq:Gram-system} or \eqref{eq:receding-Gram}.  This is an exact
consequence of the slice, not a small-error estimate.
Theorem~\ref{thm:receding} constructs the algebraic feedback on the
pointwise existence interval.  By itself it does not prove existence
or continuation of the coupled harmonic-map PDE--feedback ODE; that
geometric continuation is established in
Section~\ref{sec:adaptive-continuation}.
\end{remark}

\begin{corollary}[Forced exponential contraction]
\label{cor:forced-decay}
Fix \(0<\sigma<\beta\), and let
\(I=[\tau_0,\tau_1)\subset[0,\infty)\), where
\(0\leq\tau_0<\tau_1\leq\infty\).  Assume the hypotheses of
Theorem~\ref{thm:abstract-energy} on \(I\), with its smallness
constants chosen for this \(\sigma\).  Suppose $\eta,\delta>0$ and
\begin{equation}\label{eq:forcing-decay}
 \norm{h(\tau_0)}_{L^2_\nu}
 \leq\delta e^{-\sigma\tau_0},
 \qquad
 \norm{\E(\tau)}_{H^{-1}_\nu}
 \leq\delta e^{-(\sigma+\eta)\tau}
 \quad(\tau\in I).
\end{equation}
Then, with \(C\) allowed to depend on the fixed
\(\sigma\) and \(\eta\),
\begin{equation}\label{eq:L2-decay}
 \norm{h(\tau)}_{L^2_\nu}
 \leq C\delta e^{-\sigma\tau}
\end{equation}
for $\tau\in[\tau_0,\tau_1)$.  Moreover,
\begin{equation}\label{eq:dissipation}
 \int_{\tau_0}^{\tau_1}
 \norm{h(\tau)}_{H^1_\nu}^2\,d\tau
 \leq C\delta^2e^{-2\sigma\tau_0}.
\end{equation}
\end{corollary}

\begin{proof}
For \(\tau\in I\), multiply \eqref{eq:energy} by
$e^{2\sigma\tau}$ and integrate from \(\tau_0\) to \(\tau\).
This gives \eqref{eq:L2-decay}.  Integrating the unweighted inequality
over \([\tau_0,s]\) and letting \(s\uparrow\tau_1\) gives
\eqref{eq:dissipation}.
\end{proof}

\begin{corollary}[Integrability of the moving frame]
\label{cor:b-integrable}
Under the hypotheses and conclusions of
Corollary~\ref{cor:forced-decay}, if in addition
$\E\in L^1(I;H^{-1}_\nu)$, then
\[
 \int_{\tau_0}^{\tau_1}\abs{\mathbf b(\tau)}\,d\tau<\infty.
\]
In the complete-space specialization this is equivalent to
\(\int_{\tau_0}^{\tau_1}(|a|+|b|)\,d\tau<\infty\).
Consequently the parameters obtained by integrating the renormalized
modulation velocity have finite limits, provided such parameters have
already been geometrically constructed and their true velocity is a
component of \(\mathbf b\).  No convergence of raw fixed-chart
diffeomorphism parameters is asserted.
\end{corollary}

\begin{proof}
Combine \eqref{eq:b-bound} with \eqref{eq:dissipation}.
\end{proof}

\begin{lemma}[Closure of the extinction scale]\label{lem:scale-closure}
Let $\lambda>0$ and $t=t(\tau)$ satisfy
\[
 \frac{dt}{d\tau}=\lambda,\qquad
 \frac{d\lambda}{dt}=-(1+a(\tau))
\]
on $[\tau_0,\infty)$, where $a\in L^1([\tau_0,\infty))$.
Then $t(\tau)$ converges to a finite time $T$, and
\[
 \lambda(\tau)e^\tau\longrightarrow\lambda_\infty\in(0,\infty),
 \qquad
 \frac{\lambda(\tau)}{T-t(\tau)}\longrightarrow1.
\]
\end{lemma}

\begin{proof}
The chain rule gives
\[
 \frac{d}{d\tau}\log\lambda=-(1+a).
\]
Consequently
\[
 \lambda(\tau)
 =\lambda(\tau_0)e^{-(\tau-\tau_0)}
  \exp\left(-\int_{\tau_0}^{\tau}a(s)\,ds\right),
\]
which proves the first limit and shows that
$\int_{\tau_0}^\infty\lambda\,d\tau<\infty$.  Thus
$T=\lim_{\tau\to\infty}t(\tau)<\infty$.  Finally,
\[
 \frac{T-t(\tau)}{\lambda(\tau)}
 =\int_0^\infty
   \exp\left(-u-\int_\tau^{\tau+u}a(s)\,ds\right)\,du.
\]
The $L^1$ tail of $a$ tends to zero.  The integrand is bounded by a
fixed constant times $e^{-u}$, so dominated convergence gives the
second limit.
\end{proof}

\section{Quantitative asymptotic phase}
\label{sec:phase}

The energy inequality controls more than convergence of the sliced
tensor.  It also controls the total future motion of the geometric
frame.  Only the nine modal components of the external forcing enter
that motion linearly.

For a forcing term in Theorem~\ref{thm:receding}, set
\begin{equation}\label{eq:zeta-Z}
 \zeta_{\Z}(\tau)
 =
 \left(
 \sum_{\mu=0}^8
 \left|\ip{\rho_\tau\E}{Z_\mu}\right|^2
 \right)^{1/2},
 \qquad
 F(\tau)=\norm{\rho_\tau\E}_{H^{-1}_\nu}.
\end{equation}
When $\E$ is only an $H^{-1}_\nu$ tensor, the brackets in
\eqref{eq:zeta-Z} denote the $H^{-1}_\nu$--$H^1_\nu$ duality pairing.

\begin{theorem}[Quantitative phase theorem]\label{thm:phase}
Assume the hypotheses of Theorem~\ref{thm:receding} on
$[\tau_0,\tau_1)$, assume that the slice
\eqref{eq:receding-slice} holds at \(\tau_0\), and choose
\(\mathbf c=(a,b_1,\ldots,b_8)\) by the exact Gram system
\eqref{eq:receding-Gram}.  Thus the slice is propagated and
\eqref{eq:receding-energy} holds.  Let \(\lambda>0\) and \(t\) solve
\begin{equation}\label{eq:phase-scale-odes}
 \frac{dt}{d\tau}=\lambda,\qquad
 \frac{d\lambda}{dt}=-(1+a),
\end{equation}
with prescribed \(\lambda(\tau_0)>0\) and \(t(\tau_0)\).  Let
\(\Psi_{\tau_0}=\operatorname{Id}\) and let \(\Psi_\tau\) be the
time-ordered flow generated by
\[
 U(\tau)=\sum_{j=1}^8b_j(\tau)\chi_\tau W_j.
\]
Write
\[
 y(\tau)=\norm H_{L^2_\nu}^2,\qquad
 D(\tau)=\norm H_{H^1_\nu}^2,\qquad
 \mathbf c(\tau)=(a,b_1,\ldots,b_8).
\]
Then
\begin{equation}\label{eq:modal-velocity}
 |\mathbf c(\tau)|
 \leq C\left(
 D(\tau)+\zeta_{\Z}(\tau)+e^{-ce^\tau}
 \right).
\end{equation}
For $\tau<s<\tau_1$, the phase variation satisfies
\begin{equation}\label{eq:finite-phase-tail}
 \int_\tau^s|\mathbf c(q)|\,dq
 \leq C\left[
 y(\tau)+
 \int_\tau^s
 \left(F(q)^2+\zeta_{\Z}(q)+e^{-ce^q}\right)dq
 \right].
\end{equation}
Moreover,
\begin{equation}\label{eq:energy-convolution}
 y(\tau)
 \leq
 e^{-2\sigma(\tau-\tau_0)}y(\tau_0)
 +C\int_{\tau_0}^{\tau}
 e^{-2\sigma(\tau-q)}
 \left(F(q)^2+e^{-ce^q}\right)dq.
\end{equation}

Suppose now that $\tau_1=\infty$,
\[
 F\in L^2([\tau_0,\infty)),\qquad
 \zeta_{\Z}\in L^1([\tau_0,\infty)),
\]
and put
\[
 P_\infty(\tau)=\int_\tau^\infty|\mathbf c(q)|\,dq.
\]
Then
\begin{equation}\label{eq:infinite-phase-tail}
 P_\infty(\tau)
 \leq C\left[
 y(\tau)+
 \int_\tau^\infty
 \left(F(q)^2+\zeta_{\Z}(q)+e^{-ce^q}\right)dq
 \right].
\end{equation}
Then Lemma~\ref{lem:scale-closure} defines
\[
 T=\lim_{\tau\to\infty}t(\tau)<\infty,\qquad
 \lambda_\infty=\lim_{\tau\to\infty}\lambda(\tau)e^\tau
 \in(0,\infty).
\]
The scale and the time-ordered modulation diffeomorphism have limits.
More precisely,
\begin{align}
 \left|
 \log\frac{\lambda(\tau)e^\tau}{\lambda_\infty}
 \right|
 &\leq P_\infty(\tau),\label{eq:scale-phase-rate}\\
 \left|
 \frac{\lambda(\tau)}{T-t(\tau)}-1
 \right|
 &\leq CP_\infty(\tau),\label{eq:physical-scale-rate}
\end{align}
where in the second estimate \(C\) depends only on the fixed tail bound
\(P_\infty(\tau_0)\).  In particular, \(C\) is uniform on any family
for which \(P_\infty(\tau_0)\) has a common upper bound.
There is also a global diffeomorphism $\Psi_\infty$ such that, for every
$K\Subset M$ and $k\geq0$,
\begin{equation}\label{eq:diffeomorphism-phase-rate}
 \norm{\Psi_\infty-\Psi_\tau}_{C^k(K)}
 \leq C_{K,k}\int_\tau^\infty|b(q)|\,dq.
\end{equation}
Here the difference of maps is measured in fixed bounded-geometry
exponential charts on a slightly larger compact set; equivalently,
\(\operatorname{Exp}_{\Psi_\tau}^{-1}\Psi_\infty\) is used wherever the exponential
chart is defined.  The estimates also give a constant \(C_\Psi\),
independent of \(\tau\), for which
\begin{equation}\label{eq:phase-radial-properness}
 C_\Psi^{-1}(1+\bar f(x))
 \leq 1+\bar f(\Psi_\tau(x))
 \leq C_\Psi(1+\bar f(x)),
\end{equation}
and the same inequalities hold for \(\Psi_\tau^{-1}\) and
\(\Psi_\infty^{\pm1}\).
\end{theorem}

\begin{proof}
In the right-hand side of the exact Gram system
\eqref{eq:receding-Gram}, the forcing occurs only through the nine
numbers
\[
 \ip{\rho_\tau\E}{Z_\mu}.
\]
Lemma~\ref{lem:localized-tame} treats the quadratic term, and
Lemma~\ref{lem:C-rho} treats the moving-boundary commutator.  Uniform
invertibility of $M$ therefore gives \eqref{eq:modal-velocity}.

Integrating \eqref{eq:receding-energy} from $\tau$ to $s$ and dropping
the nonnegative terminal energy and the $L^2_\nu$ damping gives
\[
 \int_\tau^sD(q)\,dq
 \leq C\left[
 y(\tau)+\int_\tau^s(F(q)^2+e^{-ce^q})\,dq
 \right].
\]
Integration of \eqref{eq:modal-velocity} proves
\eqref{eq:finite-phase-tail}.  The integrating-factor form of
\eqref{eq:receding-energy} proves
 \eqref{eq:energy-convolution}.  Letting $s\to\infty$ gives
\eqref{eq:infinite-phase-tail}.  In particular
\(\mathbf c\in L^1([\tau_0,\infty))\), so
Lemma~\ref{lem:scale-closure} applies to
\eqref{eq:phase-scale-odes}.

Since
\[
 \frac d{d\tau}\log(\lambda e^\tau)=-a,
\]
the first scale estimate follows.  For the second, use
\[
 \frac{T-t(\tau)}{\lambda(\tau)}
 =
 \int_0^\infty
 \exp\left(
 -u-\int_\tau^{\tau+u}a(q)\,dq
 \right)du.
\]
Since
\[
 \left|\int_\tau^{\tau+u}a(q)\,dq\right|\leq P_\infty(\tau)
 \leq P_\infty(\tau_0),
\]
the elementary inequality
\(\lvert e^{-x}-1\rvert\leq e^{|x|}|x|\) bounds the difference
between the integrand and \(e^{-u}\) by
\[
 e^{P_\infty(\tau_0)}e^{-u}P_\infty(\tau).
\]
Integration proves
\eqref{eq:physical-scale-rate}.

Finally, Lemma~\ref{lem:mode-growth} and the cutoff derivative bounds
give
\[
 |U|\leq C|b|(1+\sqrt{\bar f}),\qquad
 |\bar\nabla^kU|\leq C_k|b|,\quad k\geq1.
\]
The at-most-linear growth gives a complete time-dependent flow.
Along the flow and its inverse,
\[
 \left|\frac d{d\tau}\log(1+\bar f)\right|\leq C|b(\tau)|;
\]
integration proves \eqref{eq:phase-radial-properness}.  The $L^1$ bound
on $b$, Gronwall's inequality, and the differentiated flow equations
show that $\Psi_\tau$ is Cauchy in
$C^k$ on every compact set and prove
\eqref{eq:diffeomorphism-phase-rate}.  The same estimates for the
inverse flows give a local inverse for the limit, while
\eqref{eq:phase-radial-properness} makes both limiting maps proper.
They are therefore mutually inverse global diffeomorphisms.
\end{proof}

\begin{remark}[The sharp forcing spaces]\label{rem:phase-forcing}
The distinction in \eqref{eq:zeta-Z} is structural.  The full forcing
needs only to belong to $L^2_\tau H^{-1}_\nu$ in order to give
contraction and finite dissipation.  Only its nine geometric moments
need to be integrable in time in order for the phase to converge.  An
$L^2_\tau\setminus L^1_\tau$ forcing in a geometric direction can
have finite energy cost while producing an infinite phase drift.
\end{remark}

\begin{corollary}[Power-law physical reconstruction]
\label{cor:power-reconstruction}
Under the hypotheses of Theorem~\ref{thm:phase}, suppose that for some
$\mu>0$,
\[
 y(\tau)+
 \int_\tau^\infty
 \left(F(q)^2+\zeta_{\Z}(q)+e^{-ce^q}\right)dq
 \leq Ce^{-\mu\tau}.
\]
Then
\[
 \lambda(\tau)e^\tau
 =\lambda_\infty\bigl(1+O(e^{-\mu\tau})\bigr),
 \qquad
\lambda(t)
=(T-t)\bigl(1+O((T-t)^\mu)\bigr).
\]
\end{corollary}

\begin{proof}
The hypothesis and \eqref{eq:infinite-phase-tail} give
\(P_\infty(\tau)=O(e^{-\mu\tau})\); exponentiating
\eqref{eq:scale-phase-rate} proves the first assertion.  The first
assertion and \eqref{eq:physical-scale-rate} imply
\(T-t(\tau)\asymp e^{-\tau}\), so
\(e^{-\mu\tau}=O((T-t(\tau))^\mu)\); substituting this in
\eqref{eq:physical-scale-rate} proves the second assertion.
\end{proof}

\part{Receding-domain continuation and full-metric formation}

\section{Receding boundaries and Gaussian forcing}
\label{sec:tails}

\subsection{A general tail estimate}

\begin{lemma}[Superexponential receding tails]\label{lem:superexp}
Fix \(c_1>0\), an integer \(m\geq0\), and
\(N,A\geq0\).  Let $\E(\tau)$ be supported in
\[
 \{\bar f\geq c_1e^\tau\}
\]
and suppose
\[
 \sum_{j=0}^m\abs{\bar\nabla^j\E(\tau)}
 \leq A e^{A\tau}(1+\bar f)^N.
\]
Then there are $C,c>0$ such that
\begin{equation}\label{eq:superexp}
 \norm{\E(\tau)}_{H^m_\nu}
 \leq C e^{-c e^\tau}.
\end{equation}
The same conclusion holds in $H^{-1}_\nu$.
\end{lemma}

\begin{proof}
Square the pointwise bound and integrate over
$\{\bar f\geq c_1e^\tau\}$.  Equation~
\eqref{eq:general-Gaussian-tail} absorbs every polynomial factor and
the factor $e^{A\tau}$.  The $H^{-1}_\nu$ conclusion follows from the
continuous inclusion $L^2_\nu\hookrightarrow H^{-1}_\nu$.
\end{proof}

\subsection{Grafting errors}

Fix a closed manifold $\mathcal X$, a marked open set
$\mathcal X''\subset\mathcal X$, and a
diffeomorphism
\begin{equation}\label{eq:marked-identification}
 \iota:\mathcal X''\longrightarrow M.
\end{equation}
For a closed metric $G$ we write $\iota_*G$ for its pushforward to the
corresponding physical part of the soliton chart.  This convention is
kept fixed throughout the continuation and convergence arguments.

In Stolarski's normalization, for a fixed annular-width constant
\(\Gamma_0\geq1\), the extension error $\E_{\mathrm{gr}}$ satisfies,
schematically,
\begin{equation}\label{eq:Stolarski-error}
 \supp\E_{\mathrm{gr}}
 \subset\{e^\tau\leq\bar f\leq\Gamma_0 e^\tau\},
\qquad
 \abs{\E_{\mathrm{gr}}}\leq C e^{-\tau}.
\end{equation}
For the adaptive graft, the fixed-background derivative estimate used
below is \eqref{eq:adaptive-outer-forcing};
Remark~\ref{rem:fixed-background-graft-derivative} proves it with the
factor \(\Gamma_0^{-1}\).
Lemma~\ref{lem:superexp} gives
\begin{equation}\label{eq:graft-superexp}
 \norm{\E_{\mathrm{gr}}}_{H^{-1}_\nu}
 \leq C e^{-c e^\tau}.
\end{equation}
Thus the grafting region is large in the unweighted geometry but
invisible at the scale of the Gaussian spectral argument.

\begin{proposition}[Adaptive outer graft]
\label{prop:adaptive-graft}
Let \(\lambda(t)>0\), set \(d\tau/dt=\lambda^{-1}\), and assume
\begin{equation}\label{eq:adaptive-scale-ode}
 \lambda_t=-(1+a).
\end{equation}
Let $\Theta_t$ solve
\begin{equation}\label{eq:Theta-adaptive}
 \partial_t\Theta
 =\lambda^{-1}\bigl((1+a)\bar\nabla\bar f-U\bigr)\circ\Theta
\end{equation}
and set
\[
 S(t)=\lambda(t)\Theta_t^*\bar g.
\]
Then the outer model has the exact defect
\begin{equation}\label{eq:S-defect}
 \partial_tS+2\Ric_S
 =\Theta_t^*\bigl(-2a\Ric_{\bar g}-\Lie_U\bar g\bigr).
\end{equation}
Suppose $\eta$ is a fixed graft cutoff and
\[
 \acute G=\eta\,\iota_*G+(1-\eta)S,
\]
where $\iota_*G$ is a Ricci flow on the physical part of the chart.
Assume that \(\acute G\) is a smooth metric and let \(\Phi_t\) be a
family of diffeomorphisms solving \eqref{eq:controlled-HMH} with this
\(\acute G,\lambda,a,U\).  Define \(g,h\) by
\eqref{eq:controlled-definitions}.  Define
\begin{align}
 \mathcal G_{\mathrm{gr}}
 &=2\Ric_{\acute G}
   -2\eta\Ric_{\iota_*G}-2(1-\eta)\Ric_S,
 \label{eq:pure-graft-defect}\\
 K_\tau(T)
 &=(\Phi^{-1})^*((1-\eta)\Theta^*T).
 \label{eq:K-operator}
\end{align}
Then \eqref{eq:controlled-h} becomes
\begin{equation}\label{eq:adaptive-normalized}
 \begin{split}
 \partial_\tau h={}&\A h+\Q(h)
 +(\Phi^{-1})^*\mathcal G_{\mathrm{gr}}\\
 &+a\bigl[
   Y_0+\widetilde\B_{0,\tau}h-K_\tau(Y_0)\bigr]\\
 &+\Lie_U(\bar g+h)-K_\tau(\Lie_U\bar g).
 \end{split}
\end{equation}
Thus the non-Ricci defect of the adaptive exterior does not create an
independent forcing: it converts the global controls into spatially
localized generators.
\end{proposition}

\begin{proof}
Differentiate $S$ and use $\lambda_t=-(1+a)$:
\[
 \partial_tS+2\Ric_S
 =\Theta^*\left(
 -(1+a)\bar g
 +(1+a)\Lie_{\bar\nabla\bar f}\bar g
 -\Lie_U\bar g+2\Ric_{\bar g}\right).
\]
The shrinker identity reduces the bracket to
$-2a\Ric_{\bar g}-\Lie_U\bar g$, proving
\eqref{eq:S-defect}.  Differentiating the interpolated metric and using
the Ricci-flow equation for $\iota_*G$ gives
\[
 \partial_t\acute G+2\Ric_{\acute G}
 =\mathcal G_{\mathrm{gr}}
 +(1-\eta)\Theta^*
   \bigl(-2a\Ric_{\bar g}-\Lie_U\bar g\bigr).
\]
Insert this identity into Proposition~\ref{prop:controlled-HMH} and
use $Y_0=2\Ric_{\bar g}$.
\end{proof}

\begin{remark}
At the exact background, with $\Phi=\Theta$,
\[
 K_\tau(T)=((1-\eta)\circ\Theta^{-1})T.
\]
Hence each direct column in \eqref{eq:adaptive-normalized} is the
corresponding global geometric generator multiplied by
$\eta\circ\Theta^{-1}$.  Its transition recedes through
$\bar f\simeq e^\tau$.  The diagonal scale choice is particularly
favorable: its outer defect is $2a\Ric_{\bar g}=O(|a|\bar f^{-1})$,
gaining a factor $e^{-\tau}$ on the graft annulus and preserving the
asymptotic cone to leading order.
\end{remark}

\subsection{Scale conventions and cutoff defects}

Fix \(m\in\mathbb N_0\) and \(0<\alpha<1\).  We fix the scale
conventions used from this point onward.  On an
annulus
\[
 A_L^*=\{cL<1+\bar f<CL\}
\]
with fixed \(0<c<C<\infty\), put
\(\bar g_L=L^{-1}\bar g\).  A scale-normalized norm of a covariant
two-tensor \(u\) is the \(C^{m,\alpha}(\bar g_L)\) norm of
\(L^{-1}u\); equivalently its derivative terms are
\(L^{\ell/2}|\bar\nabla^\ell u|_{\bar g}\).
For a vector field we use its \(C^{m,\alpha}(\bar g_L)\) norm.
A scale-one ball means
\[
 B_{\bar g_L}(x,r)=B_{\bar g}(x,rL^{1/2})
\]
for fixed \(r\).  Map and inverse-map norms are computed in common
scale-one harmonic atlases for the rescaled domain and range metrics
on fixed enlarged annuli; differences are measured in the
right-translated exponential charts specified in
Subsection~\ref{subsec:prepared-Banach-chart}.  On a fixed compact core
we take \(L=1\).  These conventions also govern every earlier
informal occurrence of ``scale-normalized'' in the overview.

\begin{proposition}[Spatial cutoff defect]
\label{prop:cutoff-compatible}
For $1\leq j\leq8$, define the full direct-column defect
\begin{equation}\label{eq:direct-column-defect}
 \mathcal D_{j,\tau}
 :=\widetilde Y_{j,\tau}-Y_j
 =\Lie_{(\chi_\tau-1)W_j}\bar g.
\end{equation}
It is supported in \(\{\bar f\geq2e^\tau\}\), and for every integer
\(m\geq0\) there are \(C_m,c_m>0\) such that
\[
 \|\mathcal D_{j,\tau}\|_{H^m_\nu}
 +\|\mathcal D_{j,\tau}\|_{H^{-1}_\nu}
 \leq C_me^{-c_me^\tau}.
\]
The terms in \(\mathcal D_{j,\tau}\) containing a derivative of
\(\chi_\tau\) are supported in
\(\{2e^\tau\leq\bar f\leq3e^\tau\}\) and are uniformly bounded there
with all scale-normalized derivatives.  The scale column has no
analytic cutoff.
\end{proposition}

\begin{proof}
Equation~\eqref{eq:W-all-order} gives
\(|\bar\nabla^mW_j|\leq C_m(1+\bar f)^{(1-m)/2}\).
On $\{\bar f\simeq e^\tau\}$,
\[
 \abs{\bar\nabla\chi_\tau}=O(e^{-\tau/2}),
 \qquad
 \abs{W_j}=O(e^{\tau/2}),
\]
so every term containing a cutoff derivative is uniformly
scale-normalized on the transition annulus.  Outside that annulus,
\(\mathcal D_{j,\tau}=-Y_j\) on
\(\{\bar f\geq3e^\tau\}\).  Thus every derivative of
\(\mathcal D_{j,\tau}\) has polynomial growth on its full support
\(\{\bar f\geq2e^\tau\}\).  The Gaussian tail estimate
\eqref{eq:general-Gaussian-tail} proves the \(H^m_\nu\) bound, and
\(L^2_\nu\hookrightarrow H^{-1}_\nu\) proves the last assertion.
\end{proof}

\section{Modulated three-region trapping}
\label{sec:three-region}

The argument is organized in four stages.  First, the projected energy
estimate is converted into a future phase tail and then into scalar
barriers that absorb the direct modulation columns.  Second, derivatives
are recovered using only the accumulated phase.  Third, these inputs are
assembled into the three-region improvement theorem.  Finally, a
weighted graph estimate upgrades the post-bootstrap decay to the
instantaneous spectral rate.  This ordering prevents the barrier
argument from using a pointwise velocity estimate whose proof depends
on that argument.

\subsection{Future-tail absorption and scalar barriers}
\label{subsec:phase-tail-three-region}

We now prove that the pointwise part of Stolarski's trapping argument is
stable under the feedback controls.  The essential point is that the
weighted estimate controls the \emph{future total motion} of the frame.
That quantity, rather than the instantaneous velocity, is the correct
coefficient in the pointwise argument.

Throughout this subsection, set
\begin{equation}\label{eq:q-and-P}
 q(\tau)=|a(\tau)|+|b(\tau)|,\qquad
 P_{\tau_1}(\tau)=\int_\tau^{\tau_1}q(s)\,ds .
\end{equation}
We use the scale-adapted weight
\(\omega_\sigma\) fixed in \eqref{eq:omega-sigma}.

The model perturbation equation used throughout this subsection is
\begin{equation}\label{eq:target-equation}
 \begin{split}
 \partial_\tau h={}&\A h+\Q(h)+\E\\
 &+a\bigl(\mathcal Y_{0,\tau}
          +h-\Lie_{\bar\nabla\bar f}h\bigr)\\
 &+\sum_{j=1}^8b_j
 \bigl(\mathcal Y_{j,\tau}
       +\Lie_{\chi_\tau W_j}h\bigr).
 \end{split}
\end{equation}
Here
\(\mathcal Y_{j,\tau}\), \(0\leq j\leq8\), are smooth
time-dependent symmetric two-tensor fields.  In the abstract lemmas
of this section they are arbitrary families subject to the displayed
scale and support bounds.  In the adaptive construction they are the
effective columns
\[
 \mathcal Y_{0,\tau}=Y_0-K_\tau(Y_0),\qquad
 \mathcal Y_{j,\tau}
 =\Lie_{\chi_\tau W_j}\bar g
   -K_\tau(\Lie_{\chi_\tau W_j}\bar g),
 \quad1\leq j\leq8,
\]
where \(K_\tau(T)=(\Phi^{-1})^*((1-\eta)\Theta^*T)\).
Thus the correction to the global column is supported in the union of
the cutoff exterior and the transported graft region; on every fixed
compact set the columns equal the global \(Y_j\) after both supports
have exited.  The formulas,
their precise support identity, and their scale-normalized and Gaussian
tail bounds are recorded in
\eqref{eq:effective-column-zero}--\eqref{eq:effective-column-j} and
Proposition~\ref{prop:effective-column-tails}.
Throughout this section the feedback controls \(a,b\) are continuous
in normalized time, and the column families are continuous in time
with values in the indicated spatial classes.  Profiles obtained by
dividing by \(q=|a|+|b|\), and set equal to zero where \(q=0\), are
understood as strongly Bochner-measurable families with the displayed
uniform spatial bounds.  Differential inequalities involving these
profiles hold for almost every time; equivalently, they follow by time
mollification and passage to the limit.  No time derivative of such a
profile is used.
We call its direct columns \emph{low-order controlled} when, for every
integer \(0\leq m\leq4\), they satisfy the endpoint-independent scale
bounds
\begin{equation}\label{eq:controlled-column-bounds}
 \sup_{\tau_0\leq\tau\leq\tau_1}\sup_M
 \sum_{\ell=0}^m(1+\bar f)^{\ell/2}
 |\bar\nabla^\ell\mathcal Y_{j,\tau}|_{\bar g}
 \leq C_{\mathcal Y,m},
 \qquad0\leq j\leq8.
\end{equation}
No all-order hypothesis on abstract column families is imposed here.
The compact all-order smoothing used for the adaptive solution later
comes instead from the exact support identity
\(\mathcal Y_{j,\tau}=Y_j\) on a fixed compact cylinder; see
Corollary~\ref{cor:quantitative-core-smoothing}.

\begin{lemma}[Future phase tail]\label{lem:future-phase-tail}
Fix \(0<\theta<\beta\), an amplitude \(\varepsilon>0\), forcing
constants \(C_{\E},c_{\E}>0\), and input constants
\(C_{\rm fb}^{\rm in},c_{\rm fb}^{\rm in}>0\).  Suppose on
\([\tau_0,\tau_1]\) that \(H\perp\mathcal Z\) and that
\begin{align*}
 \frac d{d\tau}\norm H_{L^2_\nu}^2
 +2\theta\norm H_{L^2_\nu}^2
 +c_{\rm fb}^{\rm in}\norm H_{H^1_\nu}^2
 &\leq
 C_{\rm fb}^{\rm in}\norm{\rho_\tau\E}_{H^{-1}_\nu}^2
 +C_{\rm fb}^{\rm in}e^{-c_{\rm fb}^{\rm in}e^\tau},\\
 q(\tau)
 &\leq
 C_{\rm fb}^{\rm in}\norm H_{H^1_\nu}^2
 +C_{\rm fb}^{\rm in}e^{-c_{\rm fb}^{\rm in}e^\tau}.
\end{align*}
Assume
\begin{equation}\label{eq:tail-energy-data}
 \|H(\tau_0)\|_{L^2_\nu}
 \leq\varepsilon e^{-\theta\tau_0},
 \qquad
 \|\rho_\tau\E(\tau)\|_{H^{-1}_\nu}
 \leq C_{\E}e^{-c_{\E}e^\tau}.
\end{equation}
Put
\[
 C_{\rm hist}:=C_{\rm fb}^{\rm in}(1+C_{\E}^2),
 \qquad
 c_{\rm hist}:=\frac12\min\{c_{\rm fb}^{\rm in},2c_{\E}\}>0.
\]
Increase \(\tau_0\), depending on \(\varepsilon\), so that
\begin{equation}\label{eq:historical-tail-absorption}
 C_{\rm hist}e^{-c_{\rm hist}e^{\tau_0}}
 \leq\varepsilon^2e^{-2\theta\tau_0}.
\end{equation}
Assume also the coarse \(C^1\) bound in
\eqref{eq:receding-box}.  There are
\(C_{\rm ft}\geq1\) and \(c_{\rm ft}>0\), independent of
\(\varepsilon,\tau_0,\tau_1\), such that
\begin{align}
 \|H(\tau)\|_{L^2_\nu}^2
 &\leq
 C_{\rm ft}\varepsilon^2e^{-2\theta\tau}
 +C_{\rm ft}e^{-c_{\rm ft}e^\tau},
 \label{eq:tail-L2}\\
 \int_\tau^{\tau_1}\|H(s)\|_{H^1_\nu}^2\,ds
 &\leq
 C_{\rm ft}\varepsilon^2e^{-2\theta\tau}
 +C_{\rm ft}e^{-c_{\rm ft}e^\tau},
 \label{eq:tail-diss}\\
 P_{\tau_1}(\tau)
 &\leq
 C_{\rm ft}\varepsilon^2e^{-2\theta\tau}
 +C_{\rm ft}e^{-c_{\rm ft}e^\tau}.
 \label{eq:tail-P}
\end{align}
Here \(C_{\rm ft},c_{\rm ft}\) depend only on
\[
 \theta,\ C_{\rm fb}^{\rm in},c_{\rm fb}^{\rm in},\
 C_{\E},c_{\E},\ \bar g .
\]
Quantitatively, if the assumed package has every upper constant at
most \(K_{\rm fb}\) and every favorable coercivity or Gaussian-decay
constant at least \(K_{\rm fb}^{-1}\), one may take
\[
 C_{\rm fb}^{\rm in}=K_{\rm fb},
 \qquad
 c_{\rm fb}^{\rm in}=K_{\rm fb}^{-1}.
\]
After enlarging \(C_{\rm ft}\) and decreasing \(c_{\rm ft}\) once,
denote the resulting common constants by
\begin{equation}\label{eq:named-future-tail-constants}
 C_P=C_P(\theta,K_{\rm fb},C_{\E},c_{\E},\bar g)\geq1,
 \qquad
 c_P=c_P(\theta,K_{\rm fb},C_{\E},c_{\E},\bar g)>0.
\end{equation}
The same \(C_P,c_P\) apply in
\eqref{eq:tail-L2}--\eqref{eq:tail-P}; in particular, they are fixed
before any barrier or package-radius threshold is evaluated.
\end{lemma}

\begin{proof}
The two Gaussian terms on the right-hand side of the energy inequality
are bounded by
\[
 C_{\rm hist}e^{-2c_{\rm hist}e^\tau}.
\]
The integrating-factor formula therefore contains a convolution
bounded by
\[
 C_{\rm hist}\int_{\tau_0}^{\tau}
 e^{-2\theta(\tau-s)}e^{-2c_{\rm hist}e^s}\,ds .
\]
Its early-time part is bounded by
\[
 C e^{-2\theta(\tau-\tau_0)}e^{-c_{\rm hist}e^{\tau_0}}
\]
and is absorbed by \eqref{eq:historical-tail-absorption}; its late-time
part is \(O(e^{-c'e^\tau})\).  This proves \eqref{eq:tail-L2}.
Integrating the same energy inequality from \(\tau\) to \(\tau_1\),
retaining the energy at \(\tau\), proves \eqref{eq:tail-diss}.
Integrating the assumed feedback estimate and using
\eqref{eq:tail-diss} proves \eqref{eq:tail-P}.  The coarse \(C^1\) box
supplies the stated feedback bound when the concrete estimate
\eqref{eq:receding-velocity} is used.
\end{proof}

\begin{lemma}[Squared-norm inequality with the geometric drift]
\label{lem:modulated-kato}
Fix a numerical bound \(K_{\mathcal Y,0}<\infty\).  There are
\(\varepsilon_{\rm K},c_{\rm K},C_{\rm K}>0\), depending only on this
bound and the fixed background, with the following property.  Let
$h$ be a smooth solution of \eqref{eq:target-equation} on a spacetime
domain, let $g=\bar g+h$, and suppose
\[
 |h|\leq\varepsilon_{\rm K}.
\]
Assume that $\E$ is continuous and that the direct columns satisfy the
zeroth-order case of \eqref{eq:controlled-column-bounds}, with
\(C_{\mathcal Y,0}\leq K_{\mathcal Y,0}\).  Put
\begin{equation}\label{eq:modulated-drift}
 V_{a,b}
 =(1+a)\bar\nabla\bar f
   -\sum_{j=1}^8b_j\chi_\tau W_j .
\end{equation}
Then \(v=|h|_{\bar g}^2\) satisfies the classical pointwise inequality
\begin{equation}\label{eq:modulated-kato}
 \begin{split}
 \partial_\tau v
 +\bar\nabla_{V_{a,b}}v
 \leq
 {}&g^{ij}\bar\nabla_i\bar\nabla_jv
 +\left(\frac{C_{\rm K}}{1+\bar f}
              +C_{\rm K}q(\tau)\right)v\\
 &+C_{\rm K}q(\tau)v^{1/2}
  +2v^{1/2}|\E|_{\bar g}
  -c_{\rm K}g^{ij}
       \langle\bar\nabla_i h,\bar\nabla_jh\rangle_{\bar g}.
 \end{split}
\end{equation}
Equivalently, with the linear scalar operator
\begin{equation}\label{eq:Kato-operator}
 \mathscr P_{a,b}w
 =\partial_\tau w
  -g^{ij}\bar\nabla_i\bar\nabla_jw
  +\bar\nabla_{V_{a,b}}w
  -\left(\frac{C_{\rm K}}{1+\bar f}
              +C_{\rm K}q\right)w,
\end{equation}
one has
\[
 \mathscr P_{a,b}v
 \leq C_{\rm K}qv^{1/2}
      +2v^{1/2}|\E|_{\bar g}
      -c_{\rm K}g^{ij}
       \langle\bar\nabla_i h,\bar\nabla_jh\rangle_{\bar g}.
\]
\end{lemma}

\begin{proof}
The leading part of $\A h+\Q(h)$ is
\[
 g^{ij}\bar\nabla_i\bar\nabla_jh
 -\bar\nabla_{\bar\nabla\bar f}h .
\]
Because the norm is taken with respect to the fixed background
\(\bar g\), we first define the exact groupings used below.  With
\(\mathcal N\) as in \eqref{eq:N-exact}, put
\[
 \begin{split}
 \mathcal R(h,\bar\nabla h)_{ij}:={}&
 2\bar R_{ikj\ell}h^{k\ell}
 +\bar R_{ja}{}^p{}_b
   \bigl(\bar g_{ip}\widetilde h^{ab}
         +\widehat h^{ab}h_{ip}\bigr)\\
 &+\bar R_{ia}{}^p{}_b
   \bigl(\bar g_{jp}\widetilde h^{ab}
         +\widehat h^{ab}h_{jp}\bigr)
 +g^{ab}g^{pq}\mathcal N(\bar\nabla h,\bar\nabla h)_{abpqij}.
 \end{split}
\]
Writing \(X_{j'}=\chi_\tau W_{j'}\), put
\[
 \begin{split}
 \mathcal M_{a,b}(h)_{ij}:={}&
 a(\mathcal Y_{0,\tau})_{ij}+a h_{ij}
 -a\bigl((\bar\nabla_i\bar\nabla^k\bar f)h_{kj}
          +(\bar\nabla_j\bar\nabla^k\bar f)h_{ik}\bigr)\\
 &+\sum_{j'=1}^8 b_{j'}\bigl(
    (\mathcal Y_{j',\tau})_{ij}
   +(\bar\nabla_iX_{j'}^k)h_{kj}
   +(\bar\nabla_jX_{j'}^k)h_{ik}\bigr).
 \end{split}
\]
Thus all second derivatives of \(h\) are in the principal coefficient
\(g^{ij}\), and all first-order control transports are in
\(\bar\nabla_{V_{a,b}}\).  Direct differentiation of
\(v=|h|_{\bar g}^2\) now gives
\[
 \begin{split}
 &\left(\partial_\tau-g^{ij}\bar\nabla_i\bar\nabla_j
       +\bar\nabla_{V_{a,b}}\right)v
 ={}-2g^{ij}
       \langle\bar\nabla_i h,\bar\nabla_jh\rangle_{\bar g}\\
 &+2\langle h,\mathcal R(h,\bar\nabla h)
                    +\E+\mathcal M_{a,b}(h)\rangle_{\bar g}.
 \end{split}
\]
The exact expansion
\eqref{eq:Q-exact} implies
\[
 2|\langle h,\mathcal R(h,\bar\nabla h)\rangle|
 \leq C|h|\,g^{ij}
       \langle\bar\nabla_i h,\bar\nabla_jh\rangle_{\bar g}
      +\frac{C}{1+\bar f}|h|^2.
\]
The quadratic curvature decay of the FIK metric was used in the last
term.  After decreasing \(\varepsilon_{\rm K}\), the first term is
absorbed by the full negative gradient term.  This coercive step uses
the evolution of \(|h|^2\) and does not require a strict Kato
inequality for \(|h|\).

The scale action decomposes as
\[
 a\bigl(h-\Lie_{\bar\nabla\bar f}h\bigr)
 =-a\bar\nabla_{\bar\nabla\bar f}h
   +a\bigl(h-2(\bar\nabla^2\bar f)*h\bigr),
\]
and, for $j\geq1$,
\[
 b_j\Lie_{\chi_\tau W_j}h
 =b_j\bar\nabla_{\chi_\tau W_j}h
   +b_j\bigl(\bar\nabla(\chi_\tau W_j)\bigr)*h .
\]
Lemma~\ref{lem:mode-growth} and the cutoff estimates give
\[
 |\bar\nabla^2\bar f|
 +|\bar\nabla(\chi_\tau W_j)|\leq C.
\]
Thus the first terms in these two displays produce precisely the drift
in \eqref{eq:modulated-drift}, and the remaining terms are bounded by
$Cqv$.  The direct columns are uniformly bounded, so their pairing
with \(2h\) is at most \(Cqv^{1/2}\), while the forcing contributes
\(2v^{1/2}|\E|\).  Uniform ellipticity of \(g^{ij}\) and one final
decrease of \(\varepsilon_{\rm K}\) prove
\eqref{eq:modulated-kato}.
\end{proof}

\begin{lemma}[Phase-tail corrected barriers]
\label{lem:phase-tail-corrected-barriers}
Fix \(0<\sigma<\theta\), \(C_P,c_P,C_{\rm gr}>0\), and a numerical
zeroth-order column ceiling \(K_{\mathcal Y,0}<\infty\).  Let
\(\varepsilon_{\rm K},c_{\rm K},C_{\rm K}\) be the constants furnished
by Lemma~\ref{lem:modulated-kato} for this ceiling.  There are constants
\[
 0<16\gamma_-<\gamma_+<\frac14,\qquad
 K,K_0,A_{\mathrm I},D_{\mathrm I},
 A_{\mathrm O},D_{\mathrm O}>0
\]
and a finite package-radius threshold
\begin{equation}\label{eq:barrier-package-radius-threshold}
 \overline\Gamma_{\rm bar}
 =\overline\Gamma_{\rm bar}
   (\sigma,\theta,C_P,c_P,C_{\rm gr},
    K_{\mathcal Y,0},\text{background})<\infty .
\end{equation}
For every already fixed package radius
\(\Gamma\geq\overline\Gamma_{\rm bar}\), there are named thresholds
\[
 \varepsilon_{\rm bar}>0,
 \qquad
 \eta_{\rm ph}^{\rm bar}>0,
\]
depending only on the displayed data, the fixed \(\Gamma\), the
background, and \(K_{\mathcal Y,0}\).  Neither threshold, nor any
barrier constant, depends on the distinct outer support parameter
\(\Gamma_0\geq1\).

For any such \(\Gamma_0\), fix
\[
 0<\varepsilon\leq\varepsilon_{\rm bar},
 \qquad
 0\leq\eta_{\rm ph}\leq\eta_{\rm ph}^{\rm bar},
\]
let \(h\) be a smooth solution of
\eqref{eq:target-equation} on \(M\times[\tau_0,\tau_1]\), put
\(g=\bar g+h\), and assume
\[
 a,b\in C^0([\tau_0,\tau_1]),
 \qquad
 |h|+|\bar\nabla h|\leq\varepsilon_{\rm K},
\]
with \(\varepsilon_{\rm K}\) as in Lemma~\ref{lem:modulated-kato}.
Assume that the direct columns satisfy the zeroth-order case of
\eqref{eq:controlled-column-bounds}, with
\(C_{\mathcal Y,0}\leq K_{\mathcal Y,0}\), that \(\E\) is continuous,
and that
\begin{equation}\label{eq:barrier-phase-hyp}
 \int_{\tau_0}^{\tau_1}q\,d\tau\leq\eta_{\rm ph},\qquad
 P_{\tau_1}(\tau)
 \leq C_P\varepsilon^2e^{-2\theta\tau}
       +C_Pe^{-c_Pe^\tau}.
\end{equation}
If \(\tau_0\) is sufficiently large, the following assertions hold.
The lower bound for \(\tau_0\) may depend on \(\varepsilon\), but is
uniform in \(\Gamma_0\).  Define
\begin{align}
 J(\tau)
 &=\exp\left(K\int_{\tau_0}^{\tau}q(s)\,ds\right),
 \label{eq:barrier-J}\\
 \mathcal B_{\mathrm I}^0(\tau,x)
 &=\varepsilon e^{-\sigma\tau}
   \left(A_{\mathrm I}\bar f^\sigma
         -D_{\mathrm I}\bar f^{\sigma-1}\right),
 \label{eq:barrier-I0}\\
 \mathcal B_{\mathrm O}^0(x)
 &=A_{\mathrm O}\varepsilon-\frac{D_{\mathrm O}}{\bar f},
 \label{eq:barrier-O0}\\
 \mathcal B_{\mathrm I}(\tau,x)
 &=J(\tau)\mathcal B_{\mathrm I}^0(\tau,x)
   -K_0P_{\tau_1}(\tau),
 \label{eq:barrier-I}\\
 \mathcal B_{\mathrm O}(\tau,x)
 &=J(\tau)\mathcal B_{\mathrm O}^0(x)
   -K_0P_{\tau_1}(\tau).
 \label{eq:barrier-O}
\end{align}
Then $\mathcal B_{\mathrm I}>0$ on
\[
 \{\Gamma<\bar f<\gamma_+e^\tau\},
\]
and
\begin{equation}\label{eq:inner-squared-supersolution}
 \mathscr P_{a,b}(\mathcal B_{\mathrm I}^2)
 \geq C_{\rm K}q\,\mathcal B_{\mathrm I}
\end{equation}
there.  In particular, \(\mathcal B_{\mathrm I}^2\) is a supersolution
of \eqref{eq:modulated-kato} wherever \(\E=0\).  If, in addition,
\begin{equation}\label{eq:outer-forcing-hyp}
 \supp\E(\tau)\subset
 \{e^\tau\leq\bar f\leq\Gamma_0e^\tau\},\qquad
 \sum_{\ell=0}^2|\bar\nabla^\ell\E|
 \leq\frac{C_{\rm gr}}{\Gamma_0}e^{-\tau},
\end{equation}
then $\mathcal B_{\mathrm O}>0$ on
\[
 \{\gamma_-e^\tau<\bar f<\infty\}
\]
and
\begin{equation}\label{eq:outer-squared-supersolution}
 \mathscr P_{a,b}(\mathcal B_{\mathrm O}^2)
 \geq C_{\rm K}q\,\mathcal B_{\mathrm O}
      +2\mathcal B_{\mathrm O}|\E|_{\bar g}
\end{equation}
there.  Moreover,
the constants can be chosen so that
\begin{equation}\label{eq:barrier-crossing}
 \mathcal B_{\mathrm O}<\mathcal B_{\mathrm I}
 \quad\hbox{on }\{\bar f=\gamma_+e^\tau\},
 \qquad
 \mathcal B_{\mathrm I}<\mathcal B_{\mathrm O}
 \quad\hbox{on }\{\bar f=\gamma_-e^\tau\}.
\end{equation}
Consequently, on the complete exterior region
\(\{\Gamma<\bar f<\infty\}\), the locally Lipschitz function
\begin{equation}\label{eq:glued-barrier}
 \mathcal B=\begin{cases}
 \mathcal B_{\mathrm I},&\bar f\leq\gamma_-e^\tau,\\
 \min\{\mathcal B_{\mathrm I},\mathcal B_{\mathrm O}\},
   &\gamma_-e^\tau\leq\bar f\leq\gamma_+e^\tau,\\
 \mathcal B_{\mathrm O},&\bar f\geq\gamma_+e^\tau
 \end{cases}
\end{equation}
is positive there, and its square is a viscosity supersolution in the precise
sense that
\begin{equation}\label{eq:glued-squared-supersolution}
 \mathscr P_{a,b}(\mathcal B^2)
 \geq C_{\rm K}q\,\mathcal B+2\mathcal B|\E|_{\bar g}
\end{equation}
on the complete region \(\{\Gamma<\bar f<\infty\}\).
\end{lemma}

\begin{proof}
We give the calculation because it is the point at which an
instantaneous estimate for $q$ would otherwise appear to be needed.
For every smooth one-variable function $F$,
\begin{equation}\label{eq:radial-drift-calculation}
 \bar\Delta_{\bar f}F(\bar f)
 =F'(\bar f)(2-\bar f)
  +F''(\bar f)|\bar\nabla\bar f|^2,
 \qquad
 |\bar\nabla\bar f|^2=\bar f-\bar R .
\end{equation}
Consequently, before replacing $\bar g^{-1}$ by $g^{-1}$,
\begin{align*}
 &\left(\partial_\tau-\bar\Delta_{\bar f}\right)
 \left[\varepsilon e^{-\sigma\tau}
       \left(A_{\mathrm I}\bar f^\sigma
             -D_{\mathrm I}\bar f^{\sigma-1}\right)\right]\\
 &\quad=\varepsilon e^{-\sigma\tau}
 \left[
   \left(D_{\mathrm I}
         -A_{\mathrm I}\sigma(1+\sigma)\right)
      \bar f^{\sigma-1}
   +O_{\sigma,D_{\mathrm I}}(\bar f^{\sigma-2})
   +O_{\sigma,A_{\mathrm I}}
       (|\bar R|\bar f^{\sigma-2})
 \right].
\end{align*}
The FIK shrinker has \(\bar R\geq0\), so
\(|\bar\nabla\bar f|^2\leq\bar f\).  Write
\[
 |\bar\nabla w|_g^2
 :=g^{ij}(\bar\nabla_iw)(\bar\nabla_jw).
\]
There is a fixed gradient constant \(C_{\rm grad}\), depending only on
\(\sigma\) and the background ellipticity bounds, such that, once
\(\Gamma\geq4D_{\mathrm I}/A_{\mathrm I}\),
\begin{equation}\label{eq:inner-barrier-gradient}
 \mathcal B_{\mathrm I}^0
 \geq\frac12A_{\mathrm I}\varepsilon
       e^{-\sigma\tau}\bar f^\sigma,
 \qquad
 \frac{|\bar\nabla\mathcal B_{\mathrm I}^0|_g^2}
      {\mathcal B_{\mathrm I}^0}
 \leq C_{\rm grad}A_{\mathrm I}\varepsilon
       e^{-\sigma\tau}\bar f^{\sigma-1}.
\end{equation}
Choose, explicitly,
\begin{equation}\label{eq:D-I-choice}
 D_{\mathrm I}>
 A_{\mathrm I}\left[
 \sigma(1+\sigma)+8(C_{\rm K}+C_{\rm grad}+1)\right].
\end{equation}
Next choose the finite threshold
\(\overline\Gamma_{\rm bar}\) in
\eqref{eq:barrier-package-radius-threshold} at least
\(4D_{\mathrm I}/A_{\mathrm I}\), beyond the compact critical region of
\(\bar f\), and large enough for every lower-order absorption in this
proof.  The threshold may be increased to a nearby regular value; this
is possible because
\(|\bar\nabla\bar f|^2=\bar f-\bar R>0\) for all sufficiently large
\(\bar f\).  Fix henceforth the package value
\(\Gamma\geq\overline\Gamma_{\rm bar}\).  The lower-order terms in the
radial calculation are then absorbed by the threshold choice, and the
\((g^{-1}-\bar g^{-1})*\bar\nabla^2\mathcal B_{\mathrm I}^0\)
term is absorbed by decreasing \(\varepsilon_{\rm K}\).  The resulting
margin is
\begin{equation}\label{eq:intermediate-margin}
 \begin{split}
 &\left(\partial_\tau
  -g^{ij}\bar\nabla_i\bar\nabla_j
  +\bar\nabla_{\bar\nabla\bar f}
  -\frac{C_{\rm K}}{1+\bar f}\right)
       \mathcal B_{\mathrm I}^0\\
 &\hspace{18mm}\geq
 4\frac{|\bar\nabla\mathcal B_{\mathrm I}^0|_g^2}
        {\mathcal B_{\mathrm I}^0}
  +
  c_{\mathrm I}\varepsilon e^{-\sigma\tau}\bar f^{\sigma-1}
 \end{split}
\end{equation}
on $\{\Gamma<\bar f<\gamma_+e^\tau\}$.  The two order
$\bar f^\sigma$ terms cancel, and the negative
$D_{\mathrm I}\bar f^{\sigma-1}$ correction gives the displayed
positive leading coefficient.  Thus \eqref{eq:D-I-choice} absorbs both
the full zeroth-order potential and the gradient cost arising when the
barrier is squared; it is stronger than the coefficient needed merely
for positivity of the unsquared radial profile.

The growth bounds for the geometric fields imply
\begin{equation}\label{eq:phase-on-radial-barrier}
 \left|
 \bar\nabla_{V_{a,b}-\bar\nabla\bar f}
 \mathcal B_{\mathrm I}^0
 \right|
 \leq Cq\mathcal B_{\mathrm I}^0 .
\end{equation}
Here one uses
$|\bar\nabla\bar f|^2\leq\bar f$ and
$|W_j(\bar f)|\leq C\bar f$.  Since $J'=KqJ$, choosing $K$ large
absorbs \eqref{eq:phase-on-radial-barrier} and the
\(C_{\rm K}q\)-potential in \eqref{eq:Kato-operator}; choose it at this
stage large enough also to absorb the \(q(1+\bar f)\) terms in the
exhaustion corrector used below.

The direct columns in \eqref{eq:modulated-kato} give the term
\(C_{\rm K}qv^{1/2}\).  At a contact with
\(v=\mathcal B_{\mathrm I}^2\), this is
\(C_{\rm K}q\mathcal B_{\mathrm I}\); at the unsquared barrier level
it is supplied by an additive \(C_{\rm K}q\) margin.  That margin is
absorbed by the last term in \eqref{eq:barrier-I}, because
\begin{equation}\label{eq:negative-tail-derivative}
 \partial_\tau\bigl[-K_0P_{\tau_1}(\tau)\bigr]
 =K_0q(\tau)
\end{equation}
for every $\tau$.  Notice that the zeroth-order terms applied
to the negative spatial constant have the favorable sign.  Choose
$K_0$ after $K$, large relative to \(C_{\rm K}\).  Finally,
\eqref{eq:barrier-phase-hyp} and
$\theta>\sigma$ imply
\[
 K_0P_{\tau_1}(\tau)
 =o\left(\varepsilon e^{-\sigma\tau}\bar f^\sigma\right)
\]
uniformly for $\bar f\geq\Gamma$.  Thus
$\mathcal B_{\mathrm I}\geq\frac12J\mathcal B_{\mathrm I}^0>0$ after
decreasing \(\varepsilon\) and increasing \(\tau_0\).  Since its
spatial gradient is \(J\bar\nabla\mathcal B_{\mathrm I}^0\),
\eqref{eq:intermediate-margin} and
\eqref{eq:inner-barrier-gradient} now give
\begin{equation}\label{eq:inner-root-margin}
 \mathscr P_{a,b}\mathcal B_{\mathrm I}
 \geq
 2\frac{|\bar\nabla\mathcal B_{\mathrm I}|_g^2}
         {\mathcal B_{\mathrm I}}
 +C_{\rm K}q.
\end{equation}
For every positive smooth \(w\), if
\[
 \mathcal V
 =\frac{C_{\rm K}}{1+\bar f}+C_{\rm K}q,
\]
then the exact product identity is
\begin{equation}\label{eq:squared-barrier-identity}
 \mathscr P_{a,b}(w^2)
 =2w\,\mathscr P_{a,b}w
   -2|\bar\nabla w|_g^2+\mathcal Vw^2.
\end{equation}
Equations~\eqref{eq:inner-root-margin} and
\eqref{eq:squared-barrier-identity} prove
\eqref{eq:inner-squared-supersolution}.

For completeness, the outer calculation is
\[
 \bar\Delta_{\bar f}(\bar f^{-1})
 =\bar f^{-1}-2\bar f^{-2}
   +2\bar f^{-3}|\bar\nabla\bar f|^2
\]
in dimension four.  Hence, for $\bar f$ large,
\begin{equation}\label{eq:outer-margin}
 \left(
  -g^{ij}\bar\nabla_i\bar\nabla_j
  +\bar\nabla_{\bar\nabla\bar f}
  -\frac{C_{\rm K}}{1+\bar f}\right)\mathcal B_{\mathrm O}^0
 \geq
 \frac{c_{\mathrm O}D_{\mathrm O}
       -C_{\rm K}'A_{\mathrm O}\varepsilon}{\bar f}.
\end{equation}
On the support of the outer forcing,
\[
 \frac{C_{\rm gr}}{\Gamma_0}e^{-\tau}
 \leq\frac{C_{\rm gr}}{\bar f}.
\]
This estimate removes \(\Gamma_0\) from the constants: only the
location of the outer edge of the support changes with \(\Gamma_0\).
Choose $D_{\mathrm O}$, depending on the fixed $C_{\rm gr}$, so that
the right-hand side of
\eqref{eq:outer-margin} absorbs \(2C_{\rm gr}/\bar f\).  Beyond
$\Gamma_0e^\tau$ the forcing vanishes.  Increase $\tau_0$ until
$D_{\mathrm O}/\bar f\leq A_{\mathrm O}\varepsilon/2$ throughout the
outer region.  On that region,
\[
 \frac{|\bar\nabla\mathcal B_{\mathrm O}^0|_g^2}
      {\mathcal B_{\mathrm O}^0}
 \leq
 \frac{CD_{\mathrm O}^2}{A_{\mathrm O}\varepsilon}\bar f^{-3}.
\]
After one further increase of \(\tau_0\), this is absorbed by one
quarter of the \(D_{\mathrm O}\bar f^{-1}\) margin.  The integrating
factor and the negative phase-tail correction treat the modulation
terms exactly as above and ensure
\(\mathcal B_{\mathrm O}\geq\frac12J\mathcal B_{\mathrm O}^0>0\).
Consequently,
\begin{equation}\label{eq:outer-root-margin}
 \mathscr P_{a,b}\mathcal B_{\mathrm O}
 \geq
 2\frac{|\bar\nabla\mathcal B_{\mathrm O}|_g^2}
         {\mathcal B_{\mathrm O}}
 +C_{\rm K}q+|\E|_{\bar g}.
\end{equation}
The squared identity \eqref{eq:squared-barrier-identity} proves
\eqref{eq:outer-squared-supersolution}.

At $\bar f=\gamma e^\tau$, the leading values of the two uncorrected
barriers are
\[
 \mathcal B_{\mathrm I}^0
 =\varepsilon A_{\mathrm I}\gamma^\sigma+O(e^{-\tau}),
 \qquad
 \mathcal B_{\mathrm O}^0
 =A_{\mathrm O}\varepsilon+O(e^{-\tau}).
\]
Choose the ratio $A_{\mathrm O}/A_{\mathrm I}$ and then
$0<16\gamma_-<\gamma_+<1/4$ so that the unique leading crossover lies
strictly between $\gamma_-$ and $\gamma_+$.  Because both barriers are
multiplied by the same $J$ and contain the same correction
$-K_0P_{\tau_1}$,
\eqref{eq:barrier-crossing} follows for $\tau_0$ large.  The minimum
of the positive barriers has positive square and
\[
 \mathcal B^2
 =\min\{\mathcal B_{\mathrm I}^2,\mathcal B_{\mathrm O}^2\}
\]
on the overlap.  The downward spatial corner of this minimum
contributes a nonnegative measure to
\(-g^{ij}\bar\nabla_i\bar\nabla_j\mathcal B^2\).  Equivalently,
regularize the minimum and pass to the limit.  At the crossing the
right-hand side of \eqref{eq:glued-squared-supersolution} has the same
value on both branches.  The weak and viscosity notions agree for this
continuous uniformly parabolic scalar operator, proving
\eqref{eq:glued-squared-supersolution}.
\end{proof}

\subsection{Derivative recovery under accumulated modulation}

\begin{lemma}[Derivative recovery from a quantitative phase tail]
\label{lem:L1-phase-derivatives}
Fix constants
\[
 0<\sigma<\theta,\qquad
 C_0,C_P,c_P,C_{\rm gr}>0.
\]
There are thresholds
\[
 \varepsilon_{\rm der}>0,
 \qquad
 \eta_{\rm ph}^{\rm der}>0,
\]
depending only on the background, the displayed fixed constants, and
the numerical column ceilings through order four, with the following
property.  For any outer support parameter \(\Gamma_0\geq1\), fix
\[
 0<\varepsilon\leq\varepsilon_{\rm der},
 \qquad
 0\leq\eta_{\rm ph}\leq\eta_{\rm ph}^{\rm der},
\]
and let \(h\) be a smooth solution of
\eqref{eq:target-equation} on all of
\(M\times[\tau_0,\tau_1]\).  Suppose
\begin{align}
 |h(\tau,x)|
 &\leq C_0\varepsilon\omega_\sigma(\tau,x),
 \label{eq:derivative-C0-hyp}\\
 \int_{\tau_0}^{\tau_1}q(s)\,ds
 &\leq\eta_{\rm ph},\qquad
 P_{\tau_1}(s)
 \leq C_P\varepsilon^2e^{-2\theta s}
      +C_Pe^{-c_Pe^s}
 \quad(\tau_0\leq s\leq\tau_1).
 \label{eq:derivative-phase-tail-hyp}
\end{align}
Assume also the coarse derivative box
\begin{equation}\label{eq:derivative-coarse-C2-box}
 \sup_{M\times[\tau_0,\tau_1]}
 \sum_{\ell=0}^2|\bar\nabla^\ell h|
 \leq\delta_{\rm B},
\end{equation}
where the structural constant \(\delta_{\rm B}>0\) is chosen, after
\(\eta_{\rm ph}^{\rm der}\) is fixed, below the
ellipticity-dependent absorption threshold.  This is the coarse
\(C^2\) bootstrap hypothesis; the
conclusion below improves its size and spatial weight.
Assume the entrance bound
\begin{equation}\label{eq:derivative-entrance-C3}
 \sum_{\ell=0}^3|\bar\nabla^\ell h(\tau_0,x)|
 \leq C_0\varepsilon\omega_\sigma(\tau_0,x),
\end{equation}
the forcing bound \eqref{eq:outer-forcing-hyp}, and the column bounds
\eqref{eq:controlled-column-bounds}.  If \(\tau_0\) is sufficiently
large, with its threshold independent of \(\Gamma_0\), so that
\begin{equation}\label{eq:derivative-late-start}
 C_{\rm gr}e^{-\tau_0}\leq\varepsilon,\qquad
 C_P^{1/2}e^{-\frac12c_Pe^{\tau_0}}
 \leq\varepsilon e^{-\sigma\tau_0},
\end{equation}
then, throughout $M\times[\tau_0,\tau_1]$,
\begin{equation}\label{eq:L1-phase-derivative-conclusion}
 |\bar\nabla h|+|\bar\nabla^2h|
 \leq C\varepsilon\omega_\sigma .
\end{equation}
The constant depends only on the background and the displayed fixed
constants (including the column constants through order four), but
not on \(\Gamma_0\), $\sup q$, or $\tau_1$.
\end{lemma}

\begin{proof}
All derivatives in this proof are taken in raw fixed-background norms.
The hypotheses \eqref{eq:derivative-coarse-C2-box}--%
\eqref{eq:derivative-entrance-C3} therefore supply no factors
\((1+\bar f)^{\ell/2}\).

Extracting the covariant transport from \eqref{eq:target-equation}
writes that equation in the form
\begin{equation}\label{eq:raw-transport-equation}
 \partial_\tau h+\bar\nabla_{V_{a,b}}h
 =g^{ij}\bar\nabla_i\bar\nabla_jh
   +\mathcal L_0(h,\bar\nabla h)
   +\E+qS+qC*h ,
\end{equation}
where \(V_{a,b}\) is defined in \eqref{eq:modulated-drift}.
When \(q>0\), \(S\) is the full direct linear combination divided by
\(q\), and \(C*h\) is the remaining zeroth-order control combination
divided by \(q\); both are set to zero when \(q=0\).  The
finite-dimensional norm equivalence and
\eqref{eq:controlled-column-bounds} give, in raw norms,
\begin{equation}\label{eq:raw-control-profile-bounds}
 \sum_{\ell=0}^4
 \bigl(|\bar\nabla^\ell S|+|\bar\nabla^\ell C|\bigr)\leq C .
\end{equation}
The coefficients in \(\mathcal L_0\) are the fixed curvature terms and
the explicit quasilinear terms from \eqref{eq:Q-exact}.  On the coarse
\(C^2\) box their raw derivatives needed below are bounded, with the
top quasilinear factors retained for absorption rather than placed in
a Schauder constant.

The geometric generators satisfy
\begin{equation}\label{eq:raw-generator-bounds}
 |V_{a,b}|\leq C(1+q)(1+\bar f)^{1/2},\qquad
 |\bar\nabla^rV_{a,b}|\leq C_r(1+q),\quad r\geq1 .
\end{equation}
Fix a terminal point \((x_*,\tau_*)\), put
\[
 I_*=[\tau_-,\tau_*],\qquad
 \tau_-=\max\{\tau_0,\tau_*-1\},
\]
and let \(F_{\tau,\tau_*}\) be the terminal-value flow
\begin{equation}\label{eq:raw-drift-following-flow}
 \partial_\tau F_{\tau,\tau_*}
 =V_{a,b}(\tau)\circ F_{\tau,\tau_*},\qquad
 F_{\tau_*,\tau_*}=\operatorname{Id}.
\end{equation}
The flow is complete on \(I_*\).  Differentiating
\eqref{eq:raw-drift-following-flow} and using
\eqref{eq:raw-generator-bounds} gives, through every fixed order,
\begin{equation}\label{eq:raw-flow-distortion}
 \|F_{\tau,\tau_*}^{\pm1}\|_{C^r_{\rm loc}}
 \leq C_r\exp\left(C_r\int_{\tau}^{\tau_*}q(s)\,ds\right),
\end{equation}
where the norms are taken in fixed-radius \(\bar g\)-balls.  Although
\(|V_{a,b}|\) grows like \((1+\bar f)^{1/2}\), it does not enter this
derivative estimate.  Along the same flow,
\begin{equation}\label{eq:raw-flow-radial-comparison}
 C^{-1}(1+\bar f(x))
 \leq1+\bar f(F_{\tau,\tau_*}(x))
 \leq C(1+\bar f(x)),
\end{equation}
because the logarithmic derivative is bounded by \(C(1+q)\) and
\(|I_*|\leq1\).  Thus both the raw geometry and
\(\omega_\sigma\) are uniformly comparable on a recent
drift-following cylinder.  Accordingly, these estimates depend on the
phase through \(\int q\), not through \(\sup q\).

The FIK background has globally bounded curvature derivatives and a
uniform positive harmonic radius.  Choose once a small raw
\(\bar g\)-radius \(r_{\rm B}>0\) and concentric balls
\[
 B_0=B_{\bar g}(x_*,r_{\rm B}),\qquad
 B_1=B_{\bar g}(x_*,2r_{\rm B}),\qquad
 B_2=B_{\bar g}(x_*,3r_{\rm B}).
\]
Pull all tensors on the drift-following cylinder back to \(B_2\) by
\(F_{\tau,\tau_*}\), and write
\[
 \widehat h=F_{\tau,\tau_*}^*h,\quad
 \widehat g_0=F_{\tau,\tau_*}^*\bar g,\quad
 \widehat\E=F_{\tau,\tau_*}^*\E,\quad
 \widehat q=q(\tau).
\]
Let \(\widehat\nabla\) be the connection of \(\widehat g_0\).
Because
\(\Lie_{V_{a,b}}h-\bar\nabla_{V_{a,b}}h
  =(\bar\nabla V_{a,b})*h\),
\eqref{eq:raw-transport-equation} becomes
\begin{equation}\label{eq:raw-drift-conjugated-equation}
 \partial_\tau\widehat h
 -\widehat a^{ij}\widehat\nabla_i\widehat\nabla_j\widehat h
 =
 \widehat A*\widehat\nabla\widehat h
 +\widehat B_0*\widehat h
 +\widehat q\,\widehat C*\widehat h
 +\widehat q\,\widehat S+\widehat\E .
\end{equation}
The raw flow bounds, the coarse \(C^2\) box, and
\eqref{eq:raw-control-profile-bounds} make this equation uniformly
parabolic on
\[
 B_2\times I_*.
\]
For the Bernstein calculation use the translated time
\[
 s=\tau-\tau_*,\qquad s_0=\tau_- -\tau_*,
\]
write \(x_0=x_*\), set \(Q_2=B_2\times[s_0,0]\), and put
\[
 \mathcal H(s)=\widehat h(\tau_*+s),\qquad
 \widehat q(s)=q(\tau_*+s),
\]
and regard \(\widehat\E\) as reparametrized by the same translated
time.
All fixed-background and flow coefficient jets needed below are
uniformly bounded after Gronwall, with \(q\) retained as the
integrable coefficient displayed in
\eqref{eq:raw-drift-conjugated-equation}.  The top quasilinear
derivatives are expanded and absorbed in the Bernstein inequalities
below; no time derivative of \(q\) is required.

Put
\[
 u(s):=\mathcal H(s),\qquad
 D:=\widehat\nabla(s),\qquad
 \gamma_s:=F_{\tau_*+s,\tau_*},\qquad
 g_0(s):=\gamma_s^*\bar g,
\]
and let \(a^{ij}\) denote the inverse of \(g_0+u\), with its indices
written relative to \(g_0\).  All norms in the following calculation
are taken with \(g_0(s)\).  Naturality of the pullback connection gives
\[
 D^j\partial_sg_0
 =\gamma_s^*\bar\nabla^j(\Lie_{V_{a,b}}\bar g).
\]
The variation formula for the Levi--Civita connection,
\[
 2g_0\bigl((\partial_sD)_XY,Z\bigr)
 =(D_X\partial_sg_0)(Y,Z)
 +(D_Y\partial_sg_0)(X,Z)
 -(D_Z\partial_sg_0)(X,Y),
\]
and \eqref{eq:raw-generator-bounds} imply, at every finite order used
below,
\begin{equation}\label{eq:Bernstein-pulled-background-ledger}
 \sum_{j=0}^3\left(
  |D^j\partial_sg_0|+|D^j\partial_sD|
 \right)
 \leq C(1+\widehat q(s)).
\end{equation}
The constant includes the harmless factor
\(\exp(C\int_{s_0}^0\widehat q)\) from the spatial jets of
\(\gamma_s^{\pm1}\), and hence contains no \(\sup\widehat q\).

Since \(Dg_0=0\), differentiation of the inverse principal coefficient
gives the exact identities
\begin{equation}\label{eq:Bernstein-inverse-coefficient-ledger}
 Da=-a*(Du)*a,
 \qquad
 D^2a=-a*(D^2u)*a+a*a*(Du)*(Du)*a.
\end{equation}
Set
\[
 \widehat{\mathscr P}
 :=\partial_s-a^{ij}D_iD_j.
\]
Using \eqref{eq:Q-exact}, collect the fixed curvature and pulled-back
background terms into smooth bounded tensors.  Commuting the equation
once and twice with \(D\), including the connection variation in
\eqref{eq:Bernstein-pulled-background-ledger}, gives
\begin{align}
 \widehat{\mathscr P}(Du)
 &=\mathcal R_1
   +\widehat q\,D(\widehat C*u+\widehat S)+D\widehat\E,
 \label{eq:Bernstein-first-commuted-ledger}\\
 \widehat{\mathscr P}(D^2u)
 &=\mathcal R_2
   +\widehat q\,D^2(\widehat C*u+\widehat S)
   +D^2\widehat\E,
 \label{eq:Bernstein-second-commuted-ledger}
\end{align}
where, pointwise on \(Q_2\),
\begin{align}
 |\mathcal R_1|
 &\leq C(1+\widehat q)(|u|+|Du|)
   +C|D^2u|+C|Du|\,|D^2u|,
 \label{eq:Bernstein-R1-ledger}\\
 |\mathcal R_2|
 &\leq C(1+\widehat q)(|u|+|Du|+|D^2u|)\notag\\
 &\quad+C|D^3u|+C\left(
    |Du|\,|D^3u|+|D^2u|^2+|Du|^2|D^2u|
   \right).
 \label{eq:Bernstein-R2-ledger}
\end{align}
These displays account for every top-order term.  The factors
\((Da)*D^2u\), \((Da)*D^3u\), and \((D^2a)*D^2u\) are exactly the
principal commutators exposed by
\eqref{eq:Bernstein-inverse-coefficient-ledger}; the differentiated
quadratic term in \eqref{eq:Q-exact} gives the remaining products.
Curvature commutators and \(\partial_sD\) are covered by
\eqref{eq:Bernstein-pulled-background-ledger}.  Thus all
principal-coefficient derivatives are displayed explicitly, without
invoking an a priori Schauder bound for them.

The phase-flow pullback changes the raw coarse box by at most
\(C\exp(C\int\widehat q)\).  Decrease the already structural
\(\delta_{\rm B}\), after the upper bound
\(\eta_{\rm ph}^{\rm der}\) is fixed, so
that on \(Q_2\)
\begin{equation}\label{eq:Bernstein-pulled-small-box}
 |u|+|Du|+|D^2u|\leq\delta_1,
\end{equation}
where \(\delta_1\) lies below the ellipticity-dependent absorption
threshold.  The Bochner identity for the evolving connection,
\eqref{eq:Bernstein-first-commuted-ledger}--%
\eqref{eq:Bernstein-R2-ledger}, and Young's inequality then give
\begin{align}
 \widehat{\mathscr P}|Du|^2
 &\leq-c|D^2u|^2
   +C(1+\widehat q)(|u|^2+|Du|^2)+C\widehat q
   +C\sum_{j=0}^1|D^j\widehat\E|^2,
 \label{eq:Bernstein-unlocalized-first}\\
 \widehat{\mathscr P}|D^2u|^2
 &\leq-c|D^3u|^2
   +C(1+\widehat q)\sum_{j=0}^2|D^ju|^2+C\widehat q
   +C\sum_{j=0}^2|D^j\widehat\E|^2.
 \label{eq:Bernstein-unlocalized-second}
\end{align}
For example, the nonlinear products paired with \(D^2u\) are bounded
by
\[
 C\left(
 |D^2u|^3+|Du|\,|D^2u|\,|D^3u|
 +|Du|^2|D^2u|^2
 \right),
\]
and \eqref{eq:Bernstein-pulled-small-box} absorbs them into the
negative derivative and lower-order terms.  For the direct source, at
each \(p=0,1,2\), use
\begin{equation}\label{eq:Bernstein-unsquared-source}
 2\widehat q\,|D^p\widehat S|\,|D^pu|
 \leq C\widehat q(1+|D^pu|^2),
\end{equation}
whereas
\(\widehat qD^p(\widehat C*u)\) is bounded by
\(C\widehat q\sum_{j\leq p}|D^ju|^2\).  Hence neither a time
derivative nor a square of \(\widehat q\) occurs.

The cylinder is truncated precisely when \(\tau_-=\tau_0\).
In either case, \eqref{eq:derivative-phase-tail-hyp} gives
\begin{equation}\label{eq:recent-phase-bound}
 \left(\int_{\tau_-}^{\tau_*}q(\tau)\,d\tau\right)^{1/2}
 \leq C\varepsilon e^{-\theta\tau_-}
      +Ce^{-c e^{\tau_-}} .
\end{equation}
When the cylinder is not truncated, \(\tau_*-\tau_-=1\), so the
right side is comparable with the same expression at \(\tau_*\);
when it is truncated, \(\tau_*-\tau_0<1\), and the entrance-time
weight is comparable throughout the cylinder.
The local Bernstein calculation for the first two spatial derivatives
therefore yields
\begin{equation}\label{eq:Bernstein-L1}
 \begin{split}
 |\widehat\nabla\mathcal H|
 +|\widehat\nabla^2\mathcal H|(x_0,0)
 \leq C\bigg(&
 \|\mathcal H\|_{C^0(Q_2)}
 +\mathbf 1_{\{\tau_*-\tau_0<1\}}
   \|\mathcal H(s_0)\|_{C^3(B_2)}\\
 &+\left(\int_{s_0}^0\widehat q(s)\,ds\right)^{1/2}
 +\|\widehat\E\|_{L^\infty([s_0,0];C_x^2(B_2))}\bigg).
 \end{split}
\end{equation}
Here \(L^\infty([s_0,0];C_x^2(B_2))\) uses only the fixed-time
spatial derivatives through order two supplied by
\eqref{eq:outer-forcing-hyp}; no time derivative of
\(\widehat\E\) is asserted or used.
To prove \eqref{eq:Bernstein-L1}, choose concentric balls
\[
 B_0\Subset B_1\Subset B_2
\]
with \(x_0\in B_0\), and spatial cutoffs
\(\zeta_0,\zeta_1\) such that \(\zeta_i\) is supported in \(B_{i+1}\)
and equals one on \(B_i\).  In a cylinder which does not meet the
initial face, choose nested time cutoffs \(\xi_1,\xi_0\) which vanish
at the back face, with \(\xi_1=1\) wherever \(\xi_0\) or
\(\xi_0'\) is nonzero.  If the initial face is present, take both
time cutoffs equal to one and retain the initial derivative norm.
 With \(a=\widehat a\), \(D=\widehat\nabla\), and
 \(\widehat{\mathscr P}\) as defined above, introduce the augmented
 quantities
\[
 \mathcal B_1
 =\xi_1\zeta_1^2|\widehat\nabla\mathcal H|^2
   +A_1|\mathcal H|^2,\qquad
 \mathcal B_2
 =\xi_0\zeta_0^2|\widehat\nabla^2\mathcal H|^2
   +A_2\mathcal B_1 .
\]
Here \(A_1\) is chosen first and \(A_2\) second, both sufficiently
large in terms of the fixed cutoff and coefficient bounds.  The
commutator identities for the uniformly controlled background, the
differentiated conjugated equation, and
\eqref{eq:derivative-coarse-C2-box} then give
\begin{equation}\label{eq:Bernstein-differential}
\begin{aligned}
 \widehat{\mathscr P}\mathcal B_1
 &\leq
 -c\xi_1\zeta_1^2|\widehat\nabla^2\mathcal H|^2
 +C(1+\widehat q)\mathcal B_1
 +C\|\mathcal H\|_{C^0(Q_2)}^2+C\widehat q
 +C\sum_{\ell=0}^1|\widehat\nabla^\ell\widehat\E|^2,
 \\
 \widehat{\mathscr P}\mathcal B_2
 &\leq
 -c\xi_0\zeta_0^2|\widehat\nabla^3\mathcal H|^2
 +C(1+\widehat q)\mathcal B_2
 +C\|\mathcal H\|_{C^0(Q_2)}^2+C\widehat q
 +C\sum_{\ell=0}^2|\widehat\nabla^\ell\widehat\E|^2 .
\end{aligned}
\end{equation}
The full unbounded transport has already been removed by
\eqref{eq:raw-drift-following-flow}.  The bounded
\(\widehat A*\widehat\nabla\mathcal H\) term and its differentiated
commutators are handled by Young's inequality and the lower-order
augmentation.  Derivatives of the drift flow contribute only
\(C(1+\widehat q)\mathcal B_i\), by
\eqref{eq:raw-generator-bounds}.  The term
\(\widehat q\,\widehat C*\mathcal H\) is retained as an integrable
zeroth-order coefficient and contributes the displayed
\(C\widehat q\,\mathcal B_i\) terms.  Thus only
\(\int\widehat q\), not \(\sup\widehat q\), enters Gronwall.
The augmentation and nesting ensure that no transition term is
discarded.  The negative gradient term in
\(\widehat{\mathscr P}(A_1|\mathcal H|^2)\) absorbs the spatial and
temporal cutoff errors in the first pass.  On the support of
\(\bar\nabla\zeta_0\) or \(\xi_0'\), one has
\(\zeta_1=\xi_1=1\); after \(A_2\) is chosen, the negative
\(|\widehat\nabla^2\mathcal H|^2\) term from
\(A_2\widehat{\mathscr P}\mathcal B_1\) absorbs every cutoff error in
the second pass.  In particular, transition-region derivatives are
never estimated by a quantity carrying the vanishing inner cutoff.

The genuinely quasilinear commutators are bounded schematically by
\[
 C\left(
 |\widehat\nabla^2\mathcal H|^3
 +|\widehat\nabla\mathcal H|\,
   |\widehat\nabla^2\mathcal H|\,
   |\widehat\nabla^3\mathcal H|
 +|\widehat\nabla\mathcal H|^2
   |\widehat\nabla^2\mathcal H|^2
 \right).
\]
The drift-flow pullback changes the raw coarse box by at most
\(C e^{C\int_{I_*}q}\).
Choose \(\delta_{\rm B}\) sufficiently small.  Young's inequality
absorbs these terms into the negative derivative terms in the
appropriate line of \eqref{eq:Bernstein-differential}, with the
remainder bounded by \(C\mathcal B_i\).  This is the required nonlinear
absorption; \(C^0\)-smallness alone would not justify
\eqref{eq:Bernstein-differential}.
 The direct-source contribution is exactly
 \eqref{eq:Bernstein-unsquared-source}; using Young's inequality with a
 constant coefficient instead would create a spurious
 \(\widehat q^2\) term.  Apply the scalar maximum principle to the first
and then the second line of \eqref{eq:Bernstein-differential}, followed
each time by Gronwall.  On the lateral boundary only the augmented
\(C^0\) term remains.  At the back face either the time cutoffs vanish
or the initial \(C^3\) norm controls both quantities.  This gives
\eqref{eq:Bernstein-L1}, with dependence on \(\widehat q\) only through
\((\int_{s_0}^0\widehat q\,ds)^{1/2}\).  The square root is unavoidable
because the Bernstein quantities are squared derivative quantities,
but it is harmless below.  In particular, no bound for \(\sup q\) has
been used.

On a backward time interval of length at most \(1\), the weight
$\omega_\sigma$ changes by at most a fixed factor on $Q_2$; the same
is true spatially by \eqref{eq:raw-flow-radial-comparison} and the
fixed raw radius of \(B_2\).  Thus
\eqref{eq:derivative-C0-hyp} controls the first term on the right of
\eqref{eq:Bernstein-L1} by
$C\varepsilon\omega_\sigma(\tau_*,x_0)$.  The initial term is bounded
in the same way by \eqref{eq:derivative-entrance-C3}.  By
\eqref{eq:recent-phase-bound}, $\theta>\sigma$, and
\eqref{eq:derivative-late-start}, the phase term is also bounded by
$C\varepsilon\omega_\sigma(\tau_*,x_0)$.

The graft forcing has the scale asserted in
\eqref{eq:outer-forcing-hyp}.  If a doubled cylinder meets its support,
then \(1+\bar f\geq e^\tau\) somewhere on that cylinder.  The raw
radial-flow comparison therefore gives
\(1+\bar f\geq c e^\tau\) throughout a fixed smaller cylinder, and
hence \(\omega_\sigma\simeq1\) there; no upper comparison involving
\(\Gamma_0e^\tau\) is used.
Fixed-background derivatives through order two of size
\(O(\Gamma_0^{-1}e^{-\tau})\) therefore contribute at most
\(C C_{\rm gr}\Gamma_0^{-1}e^{-\tau_*}
\leq C C_{\rm gr}e^{-\tau_*}\), with a constant independent of
\(\Gamma_0\), which is absorbed using
\eqref{eq:derivative-late-start}.  If the cylinder does not meet the
support, this term is zero.  Finally,
\(F_{\tau_*,\tau_*}=\operatorname{Id}\), so the left side of
\eqref{eq:Bernstein-L1} is exactly the raw fixed-\(\bar g\) derivative
norm of \(h\) at \((x_*,\tau_*)\), with no rescaling to undo.
The terminal point was arbitrary, including on the compact core.
This proves \eqref{eq:L1-phase-derivative-conclusion} on all of \(M\).
\end{proof}

\begin{remark}[Scale-rescaled variant used at higher order]
\label{rem:scale-rescaled-phase-variant}
The preceding proof uses raw fixed-background cylinders.  The
higher-order weighted estimates use the scale-normalized variant below,
for which the hypotheses already supply scale-normalized jets.  We
record it here to keep the two norm conventions separate.

Put
\[
 V_{\rm ph}=-a\bar\nabla\bar f
       +\sum_{j=1}^8b_j\chi_\tau W_j
\]
and, for any starting time \(t_-\), let
\(\mathcal P_{\tau;t_-}\) be the two-parameter flow generated
by \(-V_{\rm ph}\), with
\(\mathcal P_{t_-;t_-}=\operatorname{Id}\).  On every annular scale,
\begin{equation}\label{eq:scaled-generator-bounds}
 |\bar\nabla^rV_{\rm ph}|
 \leq C_rq(\tau)(1+\bar f)^{(1-r)/2},\qquad r\geq0.
\end{equation}
For \(\widetilde h=\mathcal P_{\tau;t_-}^*h\), exact cancellation of
the phase Lie derivative gives
\begin{equation}\label{eq:exact-phase-conjugated-equation}
 \partial_\tau\widetilde h
 =
 \mathcal P_{\tau;t_-}^*(\A h+\Q(h)+\E)
 +a\widetilde h
 +a\mathcal P_{\tau;t_-}^*\mathcal Y_{0,\tau}
 +\sum_{j=1}^8b_j
       \mathcal P_{\tau;t_-}^*\mathcal Y_{j,\tau}.
\end{equation}
No derivative of \(q\) occurs.  To record the remaining control terms
with the pullback explicit, define, where \(q(\tau)>0\),
\begin{equation}\label{eq:phase-conjugated-control-profiles}
 \begin{split}
  C_{\rm ph}(\tau)
  &:=\frac{a(\tau)}{q(\tau)}
       \operatorname{Id}_{S^2T^*M},\\
  S_{\rm ph}(\tau)
  &:=\frac1{q(\tau)}\,
     \mathcal P_{\tau;t_-}^*
       \left(a(\tau)\mathcal Y_{0,\tau}
             +\sum_{j=1}^8b_j(\tau)\mathcal Y_{j,\tau}\right),
 \end{split}
\end{equation}
and set both profiles equal to zero where \(q=0\).  Thus
\eqref{eq:exact-phase-conjugated-equation} is equivalently
\[
 \partial_\tau\widetilde h
 =
 \mathcal P_{\tau;t_-}^*(\A h+\Q(h)+\E)
 +q\,C_{\rm ph}[\widetilde h]+q\,S_{\rm ph}.
\]
In particular, the direct source is conjugated as a covariant
two-tensor by the phase flow; it is not the unpulled profile from the
raw transport equation.

For a terminal time \(\tau_*\), put
\[
 \tau(s)=\tau_*-\log(1-s),\qquad
 \phi_s=\varphi_{-\log(1-s)},\qquad
 s_-=\max\{-1,1-e^{\tau_*-t_-}\},
\]
and define, for \(s_-\leq s\leq0\),
\begin{equation}\label{eq:phase-rescaled-H}
 \mathcal H(s)
 =(1-s)\phi_s^*\widetilde h(\tau(s)).
\end{equation}
Relative to
\(\widehat g_0=(1-s)\phi_s^*
(\mathcal P_{\tau(s);t_-}^*\bar g)\), this tensor satisfies
\begin{equation}\label{eq:one-state-phase-rescaled-equation}
 \partial_s\mathcal H
 -\widehat a^{ij}\widehat\nabla_i\widehat\nabla_j\mathcal H
 =
 \Lie_{\widehat Z_s}\mathcal H
 +\widehat A*\widehat\nabla\mathcal H
 +\widehat B_0*\mathcal H
 +\widehat q\,\widehat C*\mathcal H
 +\widehat q\,\widehat S+\widehat\E ,
\end{equation}
where
\[
 \widehat q(s)=(1-s)^{-1}q(\tau(s)),\qquad
 \widehat Z_s=(1-s)^{-1}\phi_s^*
 \bigl(\bar\nabla\bar f
  -\mathcal P_{\tau(s);t_-}^*\bar\nabla\bar f\bigr).
\]
Here, with pullback of a bundle endomorphism understood in the precise
sense
\[
 (\phi_s^*C_{\rm ph})[\phi_s^*u]
 :=\phi_s^*(C_{\rm ph}[u]),
\]
\[
 \widehat C=\phi_s^*C_{\rm ph}(\tau(s)),\qquad
 \widehat S=(1-s)\phi_s^*S_{\rm ph}(\tau(s)),\qquad
 \widehat\E=\phi_s^*\mathcal P_{\tau(s);t_-}^*\E,
\]
and \(\widehat A,\widehat B_0\) collect the remaining pulled-back
quasilinear lower-order coefficients.
The scale-normalized coefficient bounds follow from
\eqref{eq:scaled-generator-bounds} and the corresponding target-jet
bounds.  The rescaled time coefficient has exactly the recent \(L^1\)
mass
\begin{equation}\label{eq:phase-rescaled-q-integral}
 \int_{s_-}^{0}\widehat q(s)\,ds
 =\int_{\tau(s_-)}^{\tau_*}q(\tau)\,d\tau .
\end{equation}
This variant is used only where scale-normalized entrance and
coefficient bounds are explicitly available; it is not used to
deduce raw derivatives from the raw \(C^2\) box.
\end{remark}

\begin{lemma}[Phase-conjugated nested-cylinder smoothing]
\label{lem:phase-conjugated-nested-smoothing}
Fix \(m\geq0\), put \(N_{\rm sm}=m+4\), and let
\[
 K_0\Subset K_1\Subset\cdots\Subset K_{N_{\rm sm}}\Subset M
\]
be fixed compact sets.  Let \(I=[\tau_0,\tau_1]\), suppose that
\(\E=0\) on \(K_{N_{\rm sm}}\times I\), and let \(h\) solve
\eqref{eq:target-equation} on a neighborhood of that cylinder.  Assume
the coarse \(C^2\) box \eqref{eq:derivative-coarse-C2-box}, and suppose
that, for some scalar \(\eta_{\rm ph}\geq0\),
\[
 \int_Iq\leq\eta_{\rm ph}.
\]
There is a threshold
\[
 \eta_{\rm ph}^{\rm sm}
 =\eta_{\rm ph}^{\rm sm}(K_0,\ldots,K_{N_{\rm sm}})>0
\]
such that the following holds whenever
\(\eta_{\rm ph}\leq\eta_{\rm ph}^{\rm sm}\).  Assume that the direct
columns have uniform
spatial \(C^{m+2}\) bounds on \(K_{N_{\rm sm}}\).  For
\[
 I_\tau=[\max\{\tau_0,\tau-\tfrac34\},\tau]
\]
one has
\begin{equation}\label{eq:nested-cylinder-smoothing}
\begin{aligned}
 \|h(\tau)\|_{C^m(K_0)}
 \leq C_m\bigg[
 &\sup_{s\in I_\tau}\|h(s)\|_{L^2(K_{N_{\rm sm}})}
 +\left(\int_{I_\tau}q(s)\,ds\right)^{1/2}\\
 &+\mathbf1_{\{\tau-\tau_0<3/4\}}
      \|h(\tau_0)\|_{C^{m+1}(K_{N_{\rm sm}})}
 \bigg].
\end{aligned}
\end{equation}
The constants depend only on \(m\), the separation distances of the
nested sets, the fixed background jets through order \(m+4\) on
\(K_{N_{\rm sm}}\), the ellipticity and coarse \(C^2\) box, the displayed
\(C^{m+2}\) column bounds, and the
upper bound \(\eta_{\rm ph}^{\rm sm}\).  In particular, no a priori derivative of \(h\)
above order two and no column jet above order \(m+2\) enters the
constant.  It is independent of \(\sup q\) and of the terminal endpoint.
\end{lemma}

\begin{proof}
Let \(t_-=\max\{\tau_0,\tau-\frac34\}\), and conjugate on
\([t_-,\tau]\) by the phase flow which is the identity at \(t_-\), as
in \eqref{eq:exact-phase-conjugated-equation}.  If
\[
 K_0\Subset K_0'\Subset K_0''\Subset K_1
\]
are fixed intermediate domains, put
\[
\begin{aligned}
 d_*=\min\bigl\{&
 \operatorname{dist}_{\bar g}(K_0,M\setminus K_0'),
 \operatorname{dist}_{\bar g}(K_0',M\setminus K_0''),\\
 &\operatorname{dist}_{\bar g}(K_0'',M\setminus K_1),
 \min_{0\leq j<N_{\rm sm}}
 \operatorname{dist}_{\bar g}(K_j,M\setminus K_{j+1})
 \bigr\},
\end{aligned}
\]
then \eqref{eq:scaled-generator-bounds} on the fixed compact
\(K_{N_{\rm sm}}\)
gives
\[
 \sup_{\substack{t_-\leq s\leq\tau\\x\in K_{N_{\rm sm}-1}}}
 \bigl(d_{\bar g}(\mathcal P_{s;t_-}(x),x)
       +d_{\bar g}(\mathcal P_{s;t_-}^{-1}(x),x)\bigr)
 \leq C_{K_{N_{\rm sm}}}\int_{t_-}^{\tau}q .
\]
Choose \(\eta_{\rm ph}^{\rm sm}\) so that the right-hand side is less
than \(d_*/4\).
The conjugated image of each \(K_j\) then lies in \(K_{j+1}\), for the
flow and its inverse.  In particular, the pulled-back graft forcing
vanishes on the cylinders used below; this is the precise point at
which the hypothesis \(\E=0\) on \(K_{N_{\rm sm}}\times I\) enters.

Write \(\widetilde h\) for the conjugated tensor.  Its principal part
has the form
\[
 \mathsf A^{ij}(\widetilde h,s,x)
 \nabla_i\nabla_j\widetilde h
\]
with a uniformly elliptic matrix \(\mathsf A\).  The remaining
modulation terms are \(qC*\widetilde h+qS\), with the required spatial
derivatives of \(C,S\) uniformly bounded and with no time derivative
of \(q\).  In the untruncated case set
\[
 t_r=\tau-\frac34+\frac{r}{8N_{\rm sm}},
 \qquad J_r=[t_r,\tau],\qquad0\leq r\leq N_{\rm sm}.
\]
Choose spatial cutoffs supported in \(K_{N_{\rm sm}-r-1}\) and equal to
one on \(K_{N_{\rm sm}-r-2}\), together with temporal cutoffs supported
in \(J_r\)
and equal to one on \(J_{r+1}\).  Thus the first cutoff is supported
inside \(K_{N_{\rm sm}-1}\), whose phase image remains in
\(K_{N_{\rm sm}}\) by the
 preceding displacement estimate.  We now give the tame estimates used
 in the iteration.  Take all derivatives with one fixed
 bounded-geometry connection \(D\) on the largest conjugated cylinder
 and write
 \[
  u:=\widetilde h,\qquad
 \mathscr P_{\rm loc}
  :=\partial_s-\mathsf A^{ij}(u,s,x)D_iD_j
       +\mathsf B^i(u,Du,s,x)D_i .
\]
On the cylinder under consideration, where the pulled-back graft
forcing vanishes, the equation is
\[
 \mathscr P_{\rm loc}u
 =
 \mathsf N(u,Du,s,x)+q\,C(s,x)*u+q\,S(s,x),
 \qquad
 \mathsf N(0,0,s,x)=0,
\]
where \(\mathsf N\) contains the quadratic-gradient and zeroth-order
terms.
 Changing between \(D\) and the pulled-back background connection
 contributes only fixed coefficient terms.  By \eqref{eq:Q-exact},
 \(\mathsf A\) is a smooth function of \(u\); the remaining
 solution-dependent nonlinearity is a smooth coefficient times
 \(Du*Du\), plus terms of order at most zero.  The phase pullback has
 spatial jets through every fixed order bounded by
 \(C_r\exp(C_r\int q)\), while its time derivative and the time
 derivative of its pulled-back connection have spatial jets bounded by
 \(C_rq(s)\).  They therefore contribute only
 \(C_r(1+q)\) times the level-\(r\) energy.

 If \(\zeta\) is supported in one member of the nested family and is
 one on the next smaller member, then, for every integer \(r\geq1\)
 and every \(\epsilon_0>0\),
 \begin{align}
  &\left|\int\zeta^2
   \left\langle D^ru,
    [D^r,\mathsf A^{ij}D_iD_j]u\right\rangle\right|
  \notag\\
  &\qquad\leq
  \epsilon_0\|\zeta D^{r+1}u\|_{L^2}^2
  +C_{r,\epsilon_0}\|u\|_{H^r(\operatorname{supp}\zeta)}^2
  +C_{r,\epsilon_0,d_*}
    \|u\|_{H^r(\operatorname{supp}D\zeta)}^2,
  \label{eq:nested-tame-principal-commutator}\\
  &\left|\int\zeta^2
   \left\langle D^ru,
    D^r\mathsf N(u,Du)\right\rangle\right|
  \notag\\
  &\qquad\leq
  \epsilon_0\|\zeta D^{r+1}u\|_{L^2}^2
  +C_{r,\epsilon_0}\|u\|_{H^r(\operatorname{supp}\zeta)}^2
  +C_{r,\epsilon_0,d_*}
    \|u\|_{H^r(\operatorname{supp}D\zeta)}^2,
  \label{eq:nested-tame-lower-commutator}
 \end{align}
 The first-order coefficient in \(\mathscr P_{\rm loc}\) produces the
 additional commutator
 \[
 \begin{aligned}
  &\left|\int\zeta^2
  \left\langle D^ru,
  [D^r,\mathsf B^i(u,Du,s,x)D_i]u\right\rangle\right|\\
  &\qquad\leq
  \epsilon_0\|\zeta D^{r+1}u\|_{L^2}^2
  +C_{r,\epsilon_0}\|u\|_{H^r(\operatorname{supp}\zeta)}^2\\
  &\hspace{38mm}
  +C_{r,\epsilon_0,d_*}
   \|u\|_{H^r(\operatorname{supp}D\zeta)}^2 .
 \end{aligned}
 \]
 Indeed, write
 \[
  \mathsf B=\mathsf B_0+\mathsf B_1,\qquad
  \mathsf B_0(s,x)=\mathsf B(0,0,s,x),\qquad
  \mathsf B_1(0,0,s,x)=0.
 \]
 The fixed part is a linear lower-order commutator, while the
 solution-dependent part satisfies schematically
 \[
  \|[D^r,\mathsf B_1^iD_i]u\|_{L^2}
  \leq
  C_r(1+\|D^2u\|_{L^\infty})\|u\|_{H^r}
  +C_r\|Du\|_{L^\infty}\|u\|_{H^{r+1}}.
 \]
 Pairing with \(D^ru\), localizing, and applying Young's inequality
 gives the displayed estimate.  These are the standard integer-order
 Moser estimates under the stated hypotheses.  In a coordinate
 extension to \(\mathbb R^4\), use
 \[
  \|[D^r,F]G\|_{L^2}
  \leq C_r\left(
   \|DF\|_{L^\infty}\|G\|_{H^{r-1}}
   +\|F\|_{H^r}\|G\|_{L^\infty}
  \right)
 \]
 with \(F=\mathsf A(u,s,x)\), \(G=D^2u\).  To preserve the homogeneous
 right-hand side, first split
 \[
  \mathsf A=\mathsf A_0+\mathsf A_1,\qquad
  \mathsf A_0(s,x)=\mathsf A(0,s,x),\qquad
  \mathsf A_1(0,s,x)=0.
 \]
 The commutator with \(\mathsf A_0\) is a fixed-coefficient linear
 commutator and is bounded, after one integration by parts, by
 \[
  \epsilon_0\|\zeta D^{r+1}u\|_{L^2}^2
  +C_{r,\epsilon_0}\|u\|_{H^r(\operatorname{supp}\zeta)}^2
  +C_{r,\epsilon_0,d_*}
    \|u\|_{H^r(\operatorname{supp}D\zeta)}^2.
 \]
 For the remaining part, smooth composition and
 \(\mathsf A_1(0,s,x)=0\) give
 \[
  \|\mathsf A_1\|_{H^r}\leq C_r\|u\|_{H^r},
  \qquad
  \|D\mathsf A_1\|_{L^\infty}
  \leq C(|u|+|Du|).
 \]
 The standard commutator inequality gives the same bound without an
 additive constant.  Pairing, integrating the principal term once, and
 applying Young's inequality proves
 \eqref{eq:nested-tame-principal-commutator}.  Likewise,
 \[
  \|D^r(F(u,s,x)Du*Du)\|_{L^2}
  \leq C_r\left(
   \|Du\|_{L^\infty}\|u\|_{H^{r+1}}
   +(1+\|D^2u\|_{L^\infty})\|u\|_{H^r}
  \right)
 \]
 proves \eqref{eq:nested-tame-lower-commutator}.  This accounts for
 every intermediate derivative split without placing an intermediate
 derivative in \(L^\infty\).

 Let \(\zeta_r,\xi_r\) be the spatial and temporal cutoffs specified
 above, with
 \[
  |D^j\zeta_r|\leq C_{m,j}d_*^{-j},
  \qquad |\xi_r'|\leq C_mN_{\rm sm},
 \]
 and each cutoff identically one on the support of the next.  Apply the
 three commutator estimates with \(\epsilon_0\) below one quarter of the
 ellipticity constant.  Keep the differentiated direct source in the
 form
 \begin{equation}\label{eq:nested-unsquared-direct-source}
  2q\,|D^rS|\,|D^ru|
  \leq C_mq(1+|D^ru|^2),
 \end{equation}
 and bound the Leibniz expansion of \(qC*u\) by
 \(C_mq\sum_{j\leq r}|D^ju|^2\).  Thus no \(q^2\) and no time
 derivative of \(q\) occurs.  With
 \[
  E_r=\sum_{j=0}^r
       \|\xi_r^{1/2}\zeta_rD^ju\|_{L^2}^2,
  \qquad
  D_{r+1}=\|\xi_r^{1/2}\zeta_rD^{r+1}u\|_{L^2}^2,
 \]
 for \(1\leq r\leq m+1\), let \(E_{r-1}^{+}\) denote the preceding
 energy and let \(D_r^{+}\) denote its order-\(r\) dissipation density,
 both evaluated with the next larger pair of cutoffs.  Then the
 localized energy recursion is
 \begin{equation}\label{eq:nested-local-energy-recursion}
  \frac d{ds}E_r+cD_{r+1}
  \leq C_m(1+q)E_r
      +C_mE_{r-1}^{+}+C_mD_r^{+}+C_mq.
 \end{equation}
 The \(r=0\) Caccioppoli inequality is separate.  Starting with that
 base estimate, integrate the recursion, use
 the nesting, and apply Gronwall.  Since
\(\int(1+q)\leq3/4+\eta_{\rm ph}^{\rm sm}\), induction for
 \(0\leq r\leq m+1\) gives
 \begin{equation}\label{eq:nested-energy-induction}
 \begin{split}
  &\sup_{s\in J_{r+1}}
    \|u(s)\|_{H^r(K_{N_{\rm sm}-r-2})}^2
  +\int_{J_{r+1}}
    \|u(s)\|_{H^{r+1}(K_{N_{\rm sm}-r-2})}^2\,ds\\
  &\qquad\leq C_m\left[
    \sup_{s\in J_0}\|u(s)\|_{L^2(K_{N_{\rm sm}-1})}^2
    +\int_{J_0}q(s)\,ds\right].
 \end{split}
 \end{equation}
 The largest direct-source derivative used is \(D^{m+1}S\), below the
 assumed \(C^{m+2}\) column ceiling; the spare derivative covers the
 tensor pullback.  The curvature, connection, and cutoff commutators at
 level \(m+1\) use background and phase-flow jets of order at most
 \(m+4\).  Thus the stated \(C^{m+2}\) column bounds and
 order-\((m+4)\) background and phase-flow jet bounds suffice.

 Because the phase image of \(K_{N_{\rm sm}-1}\) lies in
 \(K_{N_{\rm sm}}\), its uniformly controlled Jacobian gives
 \[
  \sup_{s\in J_0}\|u(s)\|_{L^2(K_{N_{\rm sm}-1})}
  \leq C\sup_{s\in J_0}\|h(s)\|_{L^2(K_{N_{\rm sm}})}.
 \]

We spell out the finite-jet endpoint bootstrap, since an abstract
quasilinear \(L^2\)-to-\(C^m\) invocation would conceal its constant
dependence.  In the untruncated case set
\[
 \mathscr E_0
 :=
 \sup_{s\in J_0}
   \|\widetilde h(s)\|_{L^2(K_{N_{\rm sm}-1})}^2
 +\int_{J_0}q(s)\,ds .
\]
The coarse \(C^2\) box, the finite volume of the fixed cylinder, and
the selected upper bound for \(\eta_{\rm ph}^{\rm sm}\) give
\begin{equation}\label{eq:nested-E0-small-ledger}
 \mathscr E_0\leq M_0
=M_0(K_{N_{\rm sm}},\delta_{\rm B},\eta_{\rm ph}^{\rm sm})<\infty .
\end{equation}
Consequently every fixed power \(p\geq1\) obeys
\(\mathscr E_0^p\leq M_0^{p-1}\mathscr E_0\).  This is not an
additional smallness assumption: \(M_0\) is absorbed into the
dependencies already declared in the statement.
After including the uniformly bounded first-order term in the local
parabolic operator, the tensor equation and the coarse \(C^2\) box give
\begin{equation}\label{eq:nested-Moser-scalar-inequality}
 \mathscr P_{\rm loc}|\widetilde h|^2
 \leq
 C(1+q)|\widetilde h|^2+Cq
\end{equation}
on the relevant nested cylinders.  To retain the vanishing estimate
when both the tensor and \(q\) vanish, introduce the future-tail
majorant
\[
 w_+(s,x)
 :=
 |\widetilde h(s,x)|^2
 +C_0\int_s^\tau q(r)\,dr .
\]
For \(C_0\) fixed sufficiently large,
\eqref{eq:nested-Moser-scalar-inequality} gives
\[
 \mathscr P_{\rm loc}w_+
 \leq C(1+q)w_+ .
\]
Multiply this inequality by a spatial--temporal cutoff and by
\(w_+^{p-1}\).  The local Sobolev inequality in real dimension four
and the usual iteration
\(p\mapsto\frac32p\) give, after Gronwall in the spatially constant
coefficient \(q(s)\),
starting from the spacetime \(L^1\) norm of \(w_+\), which is bounded
by \(C\mathscr E_0\),
\begin{equation}\label{eq:nested-Moser-C0}
 \sup_{K_0''\times J_{m+2}}
 |\widetilde h|^2
 \leq C\mathscr E_0 .
\end{equation}
At no point is \(q\) squared: the direct term is absorbed by the
future-tail majorant and the zeroth-order term is retained as the
integrable coefficient \(q(s)w_+\).  Hence the iteration constant
contains only \(\exp(C\int q)\).

Starting from \eqref{eq:nested-Moser-C0}, recover the derivatives by a
finite augmented Bernstein induction.  For \(1\leq r\leq m\), choose
successively smaller space--time cutoffs \(\zeta_r,\xi_r\), beginning
inside \(K_0''\) and ending identically one on \(K_0'\).  Define
\[
 \mathcal B_r
 =
  \xi_r\zeta_r^2|D^ru|^2
  +A_r\mathcal B_{r-1},
  \qquad
  \mathcal B_0=|u|^2,
\]
where \(A_1,\ldots,A_m\) are chosen recursively.  We record every
principal derivative split before applying the induction.  Commuting
\(D^r\) with the equation gives
\[
 (D\mathsf A)*D^{r+1}u,\qquad
 (D^r\mathsf A)*D^2u,\qquad
 (D^p\mathsf A)*D^{r-p+2}u,\quad 2\leq p\leq r-1.
\]
The first term is absorbed into the negative
\(|D^{r+1}u|^2\) term.  Smooth composition gives
\[
 D^r\mathsf A
 =\mathsf A_u*D^ru
  +\mathfrak P_r(Du,\ldots,D^{r-1}u)
  +D_x^r\mathsf A .
\]
The term containing \(D^ru\) is multiplied by the coarse-box factor
\(D^2u\).  For \(2\leq p\leq r-1\), both solution derivatives have
order at most \(r-1\), except at \(p=2\), where the \(D^2u\) factor is
in the coarse box and the other factor has order \(r\).  Thus every
strict intermediate split is controlled by the preceding Bernstein
level.  The identical decomposition applies to
 \(D^r\mathsf N(u,Du)\): a possible \(D^{r+1}u\) is multiplied by
 \(Du\), an order-\(r\) factor is multiplied by \(D^2u\), and all other
 factors have order at most \(r-1\).

The omitted first-order commutator has the same tame structure:
\[
 [D^r,\mathsf B^iD_i]u
 =
 \mathsf B_{Du}*Du*D^{r+1}u
 +(1+D^2u)*D^ru
 +\mathfrak Q_r(Du,\ldots,D^{r-1}u)
\]
schematically, after separating the fixed coefficient part.  The first
term is absorbed by the negative \(|D^{r+1}u|^2\) term, the second uses
the coarse \(C^2\) box, and every strict intermediate factor is
controlled by the preceding Bernstein level.

Assume inductively that
\(\sup\mathcal B_{r-1}\leq C_{r-1}\mathscr E_0\).  The preceding
splitting, \eqref{eq:nested-E0-small-ledger}, and
\eqref{eq:nested-unsquared-direct-source} give
\begin{equation}\label{eq:nested-Bernstein-induction}
 \mathscr P_{\rm loc}\mathcal B_r
 \leq-c\xi_r\zeta_r^2|D^{r+1}u|^2
   +C_m(1+q)\mathcal B_r
   +C_m(q+\mathscr E_0).
\end{equation}
Indeed, every product containing already controlled intermediate
derivatives is a fixed polynomial in \(\mathscr E_0^{1/2}\).  After
pairing with \(D^ru\), \eqref{eq:nested-E0-small-ledger} bounds it by
\(C_m\mathscr E_0+C_m\mathcal B_r\); no constant uses an a priori
higher derivative.  Differentiating \(qC*u+qS\) uses spatial
derivatives of \(C,S\) only through order \(r\leq m\), leaving two
derivatives in the assumed column ceiling.  On the support of a
level-\(r\) cutoff derivative the level-\((r-1)\) cutoffs are one.
Choose \(A_r\) so that the negative \(|D^ru|^2\) term in
\(A_r\mathscr P_{\rm loc}\mathcal B_{r-1}\) absorbs every spatial and
temporal cutoff error at level \(r\).

The scalar maximum principle, followed by Gronwall and induction in
\(r\), now yields
\begin{equation}\label{eq:nested-finite-jet-endpoint}
 \|\widetilde h(\tau)\|_{C^m(K_0')}^2
 \leq C_m\mathscr E_0 .
\end{equation}
This proves the needed endpoint regularization using only the displayed
finite jets.  Combining \eqref{eq:nested-finite-jet-endpoint} with the
phase-image \(L^2\) comparison above and taking square roots gives the
first two terms on the right of
\eqref{eq:nested-cylinder-smoothing}.

If \(\tau-\tau_0<3/4\), take the temporal cutoffs to be identically
one.
In this truncated case set
\[
 \mathscr E_{\rm tr}
 :=
 \sup_{s\in I_\tau}
   \|\widetilde h(s)\|_{L^2(K_{N_{\rm sm}-1})}^2
 +\int_{I_\tau}q(s)\,ds
 +\|h(\tau_0)\|_{C^{m+1}(K_{N_{\rm sm}})}^2 .
\]
Repeating the scalar cutoff--Moser calculation with the temporal
cutoffs identically one and retaining the back-face term gives
\[
 \sup_{K_0''\times I_\tau}|\widetilde h|^2
 \leq C_m\mathscr E_{\rm tr}.
\]
This is the \(r=0\) base for the augmented Bernstein induction; at
orders \(1\leq r\leq m\), retain the corresponding initial derivative
terms.  At each order the negative
derivative term of the preceding augmented quantity absorbs the
spatial cutoff error, exactly as in
\eqref{eq:Bernstein-differential}; hence no inverse power of
\(\tau-\tau_0\) occurs.  This gives the additional initial-trace term
displayed in \eqref{eq:nested-cylinder-smoothing}.  Finally, undoing
the phase pullback is legitimate because, after one further decrease
of \(\eta_{\rm ph}^{\rm sm}\),
\(\mathcal P_{s;t_-}^{-1}(K_0)\subset K_0'\) throughout the recent
cylinder.  It changes the constants only by
\(e^{C_m\int_{I_\tau}q}\) and gives
\eqref{eq:nested-cylinder-smoothing} on the stated domain \(K_0\).
\end{proof}

\subsection{The three-region improvement}

\begin{theorem}[Robust modulated three-region improvement]
\label{thm:robust-modulated-three-region}
Fix
\[
 0<\sigma<\theta<\beta .
\]
Fix also forcing constants
$C_{\E},c_{\E},C_{\rm gr}>0$.
Fix numerical column bounds
\[
 K_{\mathcal Y,m}<\infty,\qquad0\leq m\leq4,
\]
and a numerical feedback constant \(K_{\rm fb}<\infty\).  In the
statement below every upper constant in the energy and feedback
estimates assumed in Lemma~\ref{lem:future-phase-tail} is required to
be at most \(K_{\rm fb}\), while every favorable coercivity or decay
constant is required to be at least \(K_{\rm fb}^{-1}\).  The constants
in
\eqref{eq:controlled-column-bounds} are required to satisfy
\(C_{\mathcal Y,m}\leq K_{\mathcal Y,m}\).  These bounds are fixed
before any threshold is chosen.
In the adaptive application they are precisely the pre-radius values
in \eqref{eq:pre-radius-three-region-inputs}; in particular, none is
evaluated from an atlas or evolution constant depending on the
eventual package radius.
Let \(C_P,c_P\) be the named constants in
\eqref{eq:named-future-tail-constants}.  Apply
Lemma~\ref{lem:modulated-kato} with the ceiling
\(K_{\mathcal Y,0}\), and then apply
Lemma~\ref{lem:phase-tail-corrected-barriers} with
\[
 (\sigma,\theta,C_P,c_P,C_{\rm gr},K_{\mathcal Y,0}).
\]
At this stage fix every shape constant furnished by the barrier lemma,
including \(\gamma_\pm\) and all ratios among the four provisional
barrier amplitudes.  Leave their common positive normalization free,
since it enters no radius inequality.  Let
\(\overline\Gamma_{\rm bar}\) denote the resulting threshold.
Let \(\delta_{\rm rec}>0\) be an admissible receding-box threshold in
Theorem~\ref{thm:receding}.  Retain for the moment the constant
\(\varepsilon_{\rm K}\) furnished by
Lemma~\ref{lem:modulated-kato}.  Its final decrease will be made only
after the actual compatible radius and the common normalization of the
barriers have fixed the numerical constant \(C_0\) entering
Lemma~\ref{lem:L1-phase-derivatives}.  That later decrease enters no
radius inequality and preserves every conclusion obtained with the
present value of \(\varepsilon_{\rm K}\).
There is a finite threshold
\begin{equation}\label{eq:three-region-package-radius-threshold}
 \overline\Gamma_{\rm 3reg}
 =\mathfrak G_{\rm 3reg}
  (\sigma,\theta,C_{\E},c_{\E},C_{\rm gr},
   (K_{\mathcal Y,m})_{m=0}^{4},K_{\rm fb})
 \geq\overline\Gamma_{\rm bar}.
\end{equation}
The admissibility condition is
\begin{equation}\label{eq:three-region-compatible-package-radius}
 \Gamma\geq\overline\Gamma_{\rm 3reg}.
\end{equation}
For every already fixed package radius satisfying
\eqref{eq:three-region-compatible-package-radius}, first take the named
thresholds \(\varepsilon_{\rm bar}\) and
\(\eta_{\rm ph}^{\rm bar}\) furnished by
Lemma~\ref{lem:phase-tail-corrected-barriers} for this value of
\(\Gamma\).  Fix also the \(m=2\) nested core collars used below and
let \(\eta_{\rm ph}^{\rm sm}>0\) be the corresponding phase threshold
in Lemma~\ref{lem:phase-conjugated-nested-smoothing}.  Choose a
provisional structural number \(\delta_{\rm sm}>0\) below the
quasilinear absorption threshold in that fixed nested-core calculation,
and decrease, if necessary, the presently retained
\(\varepsilon_{\rm K}\) so that
\(\varepsilon_{\rm K}\leq\delta_{\rm sm}\).  This choice depends on no
barrier normalization and enters no radius inequality.  Using this
value as the coarse \(C^2\)-box ceiling, fix the resulting \(m=2\)
core-smoothing constant, and hence its zeroth-order component; it
remains valid after every later decrease of
\(\varepsilon_{\rm K}\).  Increase, if necessary,
one common normalization of the four barrier amplitudes, retaining all
previously frozen ratios, so that the entrance and inner-boundary
comparisons below hold.  Since
\(J\leq e^{K\eta_{\rm ph}^{\rm bar}}\), choose a numerical
\(C_{0,{\rm bar}}\geq1\), depending only on the data already fixed and
on this normalized barrier package, which dominates both the fixed
zeroth-order core-smoothing constant and the normalized glued-barrier
comparison constant.  Thus the zeroth-order argument below has the
form
\[
 |h|\leq C_{0,{\rm bar}}\varepsilon\omega_\sigma .
\]
Apply Lemma~\ref{lem:L1-phase-derivatives} with
\(C_0=C_{0,{\rm bar}}\), the named \(C_P,c_P\), the fixed
\(C_{\rm gr}\), and the prescribed column ceilings.  Let
\(\varepsilon_{\rm der}>0\) and
\(\eta_{\rm ph}^{\rm der}>0\) be its thresholds, and choose its
structural \(\delta_{\rm B}>0\), after
\(\eta_{\rm ph}^{\rm der}\), below the absorption thresholds in both
that lemma and the fixed nested-core smoothing calculation.  Now
decrease \(\varepsilon_{\rm K}\), if necessary, so that
\begin{equation}\label{eq:ordered-three-region-thresholds}
 0<\varepsilon_{\rm K}
 \leq\min\{\delta_{\rm sm},\delta_{\rm B},\delta_{\rm rec}\}.
\end{equation}
This decrease changes no radius or barrier constant and preserves the
Kato and barrier conclusions.

There are constants
\[
 C_*<\infty,\qquad c_*>0,
\]
depending only on the fixed background, rate pair, column and feedback
ceilings, forcing package, and the fixed compatible package radius
\(\Gamma\), and independent of the outer support parameter
\(\Gamma_0\geq1\) and of the finite endpoint, with the following
property.  Choose \(C_*\) larger than every barrier, core-smoothing,
and derivative-recovery output constant for the normalized package.
After these constants and the package radius are fixed, there is
\(\varepsilon_*>0\), still independent of \(\Gamma_0\), chosen below
\(\varepsilon_{\rm bar}\) and \(\varepsilon_{\rm der}\) and so that
\begin{equation}\label{eq:three-region-strict-improvement-ledger}
 2C_*\varepsilon_*\leq\varepsilon_{\rm K}.
\end{equation}
For any \(\Gamma_0\geq1\), let
$0<\varepsilon\leq\varepsilon_*$, let $\tau_0$ be sufficiently large
with a lower bound independent of \(\Gamma_0\),
  and let $h$ be a smooth solution on
  $[\tau_0,\tau_1]$ of \eqref{eq:target-equation},
  where $(a,b)$ is selected by the exact receding feedback system
  associated with the direct columns
  $\mathcal Y_{0,\tau},\ldots,\mathcal Y_{8,\tau}$.  Assume that these
  columns are low-order controlled in the sense of
  \eqref{eq:controlled-column-bounds}
and that the exact feedback package supplies the two estimates in
Lemma~\ref{lem:future-phase-tail}, with upper constants bounded by
\(K_{\rm fb}\), favorable constants bounded below by
\(K_{\rm fb}^{-1}\), and all constants independent of \(\tau_1\).
Assume:
\begin{enumerate}
\item $H=\rho_\tau h$ satisfies the exact slice
      \eqref{eq:receding-slice};
\item the coarse box
      \[
       \sup_{\tau_0\leq\tau\leq\tau_1}
       \sup_M
       \sum_{\ell=0}^2|\bar\nabla^\ell h|
       \leq\varepsilon_{\rm K}
      \]
      holds;
\item the weighted entrance bound is
      \[
       \|H(\tau_0)\|_{L^2_\nu}
       \leq\varepsilon e^{-\theta\tau_0};
      \]
\item the entrance tensor is arbitrary, subject only to the slice and
      the scale-adapted inequalities
      \begin{equation}\label{eq:arbitrary-stable-entrance}
       \sum_{\ell=0}^3|\bar\nabla^\ell h(\tau_0)|
       \leq\varepsilon\omega_\sigma(\tau_0,\cdot);
      \end{equation}
\item the forcing satisfies the second inequality in
      \eqref{eq:tail-energy-data}, as well as
      \eqref{eq:outer-forcing-hyp}, and
      vanishes on $\{\bar f<e^\tau\}$, and $\tau_0$ is chosen so that
      \eqref{eq:historical-tail-absorption} and
      \eqref{eq:derivative-late-start} hold.  In the former take
      \[
       C_{\rm hist}=K_{\rm fb}(1+C_{\E}^2),
       \qquad
       c_{\rm hist}=\frac12\min\{K_{\rm fb}^{-1},2c_{\E}\};
      \]
      in the latter use the named constants \(C_P,c_P\) furnished by
      Lemma~\ref{lem:future-phase-tail}.
\end{enumerate}
Then, with the already fixed \(C_*\) and independently of \(\tau_1\),
\begin{align}
 \|H(\tau)\|_{L^2_\nu}
 &\leq C_*\varepsilon e^{-\theta\tau}
       +C_*e^{-c_*e^\tau},\label{eq:three-region-L2}\\
 \sum_{\ell=0}^2|\bar\nabla^\ell h(x,\tau)|
 &\leq C_*\varepsilon
 e^{-\sigma\tau}(1+\bar f(x))^\sigma,
 &&\bar f(x)\leq e^\tau,\label{eq:three-region-middle}\\
 \sum_{\ell=0}^2|\bar\nabla^\ell h(x,\tau)|
 &\leq C_*\varepsilon,
 &&x\in M,\label{eq:three-region-global}\\
 |a(\tau)|+|b(\tau)|
 &\leq C_*\varepsilon^2e^{-2\sigma\tau}
       +C_*e^{-c_*e^\tau}.\label{eq:three-region-velocity}
\end{align}
In particular, the entrance tensor in
\eqref{eq:arbitrary-stable-entrance} need not be a cutoff finite sum
of eigentensors.  It may contain an arbitrary small component in
$\Z^\perp$.
\end{theorem}

\begin{proof}
We record the order of choices, since it is part of the strict
bootstrap improvement.  First fix the background, exponents, forcing
package, column bounds, and feedback package in the statement.  Next
fix the named future-tail constants, the Kato and receding-energy
thresholds, and all barrier shape data specified before
\eqref{eq:three-region-package-radius-threshold}.  In particular, every
barrier shape ratio which enters a radial lower bound is already fixed.
The numerical thresholds in Lemma~\ref{lem:L1-phase-derivatives} are
fixed only after the compatible radius, because their displayed input
\(C_0\) is determined by the common barrier normalization.
Then choose the finite number
\(\overline\Gamma_{\rm 3reg}\) in
 \eqref{eq:three-region-package-radius-threshold}.  Fix an arbitrary
 package value satisfying
 \eqref{eq:three-region-compatible-package-radius}; the resulting radius
 is never reassigned or enlarged in the proof.  Retain the already fixed
 Gaussian exponents in the theorem-level estimates as follows.  Let
 \(c_{\rm cut}>0\) be a common exponent in
 \eqref{eq:cutoff-tail} for \(0\leq m\leq2\) and in
 \eqref{eq:gram-tail}, decreased once to cover the finitely many cutoff
 commutators used below.  By the quantitative convention in the
 statement, every Gaussian decay exponent in the assumed receding-energy
 and feedback package is at least \(K_{\rm fb}^{-1}\).  Define
 \begin{equation}\label{eq:three-region-common-Gaussian-exponent}
  c_*:=
  \frac{1}{4e}
  \min\left\{
   c_P,\ c_{\E},\ c_{\rm cut},\ K_{\rm fb}^{-1}
  \right\}>0 .
 \end{equation}
 The factor \(1/(4e)\) absorbs the square-root loss in passing from the
 squared energy estimate to its norm form and every fixed normalized-time
 translation of length at most one below.  Consequently every Gaussian
 remainder used in this proof may, after increasing its multiplicative
 constant, be written as \(Ce^{-c_*e^\tau}\).  The dependencies of
 \(c_*\) are exactly those declared in the theorem; in particular it is
 independent of \(\Gamma_0\) and of the finite endpoint.  Neither \(c_*\)
 nor the compatible package radius is reassigned later.  With the fixed
core collars, select \(\eta_{\rm ph}^{\rm sm}\) and
\(\delta_{\rm sm}\), make the preliminary structural decrease of
\(\varepsilon_{\rm K}\), and fix the corresponding core-smoothing
constant.  Retain the already fixed \(\gamma_\pm\) and amplitude ratios
and choose the common normalization prescribed in the statement.  Then
fix \(C_{0,{\rm bar}}\), \(\varepsilon_{\rm der}\),
\(\eta_{\rm ph}^{\rm der}\), and \(\delta_{\rm B}\) in exactly the
order stated above, and make the final decrease of
\(\varepsilon_{\rm K}\).  Only then choose \(C_*\) larger than every
output constant in the estimates below.  After \(C_*\) is fixed choose
\(\varepsilon_*\) below \(\varepsilon_{\rm bar}\) and
\(\varepsilon_{\rm der}\), and so that
\eqref{eq:three-region-strict-improvement-ledger} holds.  Finally
choose \(\tau_0\), depending on \(\varepsilon\) when required by the
displayed late-start inequalities.  No later step enlarges \(C_*\).

Lemma~\ref{lem:future-phase-tail}, applied with
\[
 C_{\rm fb}^{\rm in}=K_{\rm fb},
 \qquad
 c_{\rm fb}^{\rm in}=K_{\rm fb}^{-1},
\]
and with its Gaussian remainders weakened to the already fixed exponent
\(c_*\), proves
\eqref{eq:three-region-L2}, the dissipation estimate
\eqref{eq:tail-diss}, and the phase-tail bound \eqref{eq:tail-P}.
For the fixed core collars used below, let
\[
 \eta_{\rm ph}^*
 :=\min\left\{
  \eta_{\rm ph}^{\rm bar},
  \eta_{\rm ph}^{\rm der},
  \eta_{\rm ph}^{\rm sm}
 \right\}>0.
\]
By decreasing \(\varepsilon_*\) and increasing \(\tau_0\), while
retaining the already fixed package and radius, arrange that
\[
 \int_{\tau_0}^{\tau_1}q\,d\tau\leq\eta_{\rm ph}^*.
\]
Thus the phase hypotheses of
Lemmas~\ref{lem:phase-tail-corrected-barriers},
\ref{lem:L1-phase-derivatives}, and
\ref{lem:phase-conjugated-nested-smoothing} hold simultaneously.

On a fixed core $\{\bar f<2\Gamma\}$, the forcing vanishes.  After the
first unit of normalized time,
Lemma~\ref{lem:phase-conjugated-nested-smoothing}, with \(m=2\) and
fixed nested sublevel sets between
\(\{\bar f\leq\Gamma\}\) and \(\{\bar f<2\Gamma\}\), gives
\begin{equation}\label{eq:core-improvement}
 \sum_{\ell=0}^2|\bar\nabla^\ell h|
 \leq C\varepsilon e^{-\theta\tau}
 \qquad\hbox{on }\{\bar f\leq\Gamma\},\quad
 \tau\geq\tau_0+1.
\end{equation}
On the initial unit interval, use the initial-face clause of
Lemma~\ref{lem:phase-conjugated-nested-smoothing} with \(m=2\).
The entrance \(C^3\) term in
\eqref{eq:nested-cylinder-smoothing}, the displayed bound on the
accumulated phase, and comparability of the time weights then give
\[
 \sum_{\ell=0}^2|\bar\nabla^\ell h|
 \leq C\varepsilon e^{-\sigma\tau}
 \qquad\hbox{on }\{\bar f\leq\Gamma\}.
\]
Indeed, the accumulated phase on that interval obeys
\[
 \int_{\tau_0}^{\tau}q(s)\,ds
 \leq P_{\tau_1}(\tau_0)
 \leq C\varepsilon^2e^{-2\theta\tau_0}
      +Ce^{-ce^{\tau_0}}.
\]
For $\tau-\tau_0\geq1$, apply the interior estimate on
$[\tau-1,\tau]$; then
\[
 \int_{\tau-1}^{\tau}q(s)\,ds
 \leq P_{\tau_1}(\tau-1)
 \leq C_\theta\varepsilon^2e^{-2\theta\tau}
      +Ce^{-ce^\tau}.
\]
Together with
$\sup_{s\in[\tau-1,\tau]}\|H(s)\|_{L^2_\nu}
\leq C_\theta\varepsilon e^{-\theta\tau}+Ce^{-ce^\tau}$,
Lemma~\ref{lem:phase-conjugated-nested-smoothing} gives
\eqref{eq:core-improvement}.
Thus the direct columns contribute only the recent $L^1$ phase tail,
not $\sup q$.

Apply Lemma~\ref{lem:modulated-kato} to \(v=|h|_{\bar g}^2\) and
compare \(v\) with the squared glued barrier
\(\mathcal B^2\).  The two-regime core estimate
supplies the inner spatial boundary inequality at $\bar f=\Gamma$;
after the initial unit it improves that boundary because
$\theta>\sigma$.  Assumption
\eqref{eq:arbitrary-stable-entrance} supplies the initial inequality,
because the common normalization of
\(A_{\mathrm I},D_{\mathrm I},A_{\mathrm O},D_{\mathrm O}\)
was chosen after the radius and before \(C_*\).  That normalization
changes none of the ratios used in the radius condition or the
crossover geometry and only enlarges the strict supersolution margins.
No amplitude is changed after \(C_*\) is frozen.
There is no outer spatial boundary.  To justify comparison on the
complete end, use
\[
 \mathcal B_\delta
 =\mathcal B+z_\delta,\qquad
 z_\delta
 =\delta J(\tau)e^{L(\tau-\tau_0)}(1+\bar f).
\]
The factor \(J\) is essential here.  The identities
\(\bar\Delta_{\bar f}\bar f=2-\bar f\) and \(J'=KqJ\), the coarse
ellipticity box, and
\[
 |\bar\nabla_{V_{a,b}-\bar\nabla\bar f}\bar f|
 \leq Cq(1+\bar f)
\]
show, with that fixed choice of \(K\) and then a fixed sufficiently
large \(L\), that
\begin{equation}\label{eq:exhaustion-corrector}
 \mathscr P_{a,b}z_\delta
 \geq
 2\frac{|\bar\nabla z_\delta|_g^2}{z_\delta}+z_\delta
 \qquad\text{on }\{\bar f>\Gamma\}.
\end{equation}
Thus the exhaustion argument uses no bound on \(\sup q\).
The weighted Cauchy inequality
\[
 \frac{|\bar\nabla(w+z)|_g^2}{w+z}
 \leq\frac{|\bar\nabla w|_g^2}{w}
     +\frac{|\bar\nabla z|_g^2}{z},
 \qquad w,z>0,
\]
together with \eqref{eq:inner-root-margin},
\eqref{eq:outer-root-margin}, and
\eqref{eq:squared-barrier-identity}, shows that every smooth branch
of \(\mathcal B_\delta^2\) is a strict supersolution of the scalar
inequality \eqref{eq:modulated-kato}.  The minimum construction
preserves this in the viscosity sense.  Since \(h\) is bounded while
\(\mathcal B_\delta\to\infty\) as \(\bar f\to\infty\), the inequality
\(v<\mathcal B_\delta^2\) holds on the outer boundary of a sufficiently
large exhaustion.  Apply the scalar maximum principle, let the
exhaustion radius tend to infinity, and then let \(\delta\downarrow0\).
This yields \(v\leq\mathcal B^2\), and hence
\begin{equation}\label{eq:C0-three-region}
 |h(x,\tau)|
 \leq C_{0,{\rm bar}}\varepsilon\omega_\sigma(\tau,x).
\end{equation}
In Stolarski's proof the specially prepared eigensum supplies both
this initial comparison and the initial derivative bounds used in the
next step.  Here both are supplied directly by
\eqref{eq:arbitrary-stable-entrance}; no spectral representation of
$h(\tau_0)$ is used.

Lemma~\ref{lem:L1-phase-derivatives}, with the already fixed input
\(C_0=C_{0,{\rm bar}}\), upgrades
\eqref{eq:C0-three-region} to the first two derivative bounds.  On
$\bar f\leq e^\tau$ this is
\eqref{eq:three-region-middle}; on the complete exterior it is
\eqref{eq:three-region-global}.

It remains to close the instantaneous velocity.  Since
$H=\rho_\tau h$ and the cutoff annulus is Gaussian
superexponentially small, the pointwise estimate just proved gives
\[
 \|H(\tau)\|_{H^1_\nu}^2
 \leq
 C\varepsilon^2e^{-2\sigma\tau}
 \int_M(1+\bar f)^{2\sigma}\,d\nu
 +Ce^{-ce^\tau}
 \leq C\varepsilon^2e^{-2\sigma\tau}+Ce^{-ce^\tau}.
\]
Insert this into \eqref{eq:receding-velocity} to obtain
\eqref{eq:three-region-velocity}.  The barrier and derivative arguments
use only the integral phase tail, while the pointwise velocity estimate
enters only at this final step.
Finally, \eqref{eq:three-region-strict-improvement-ledger} gives
\[
 \sup_M\sum_{\ell=0}^2|\bar\nabla^\ell h|
 \leq C_*\varepsilon\leq\frac12\varepsilon_{\rm K}.
\]
Thus a bootstrap opened with the doubled face
\(\sup\sum_{\ell\leq2}|\bar\nabla^\ell h|
 \leq2C_*\varepsilon\) is returned to \(C_*\varepsilon\), strictly
inside both that face and the structural coarse box.
\end{proof}

\subsection{Post-bootstrap spectral recovery}

\begin{lemma}[Weighted graph estimate and operator core]
\label{lem:weighted-graph-estimate}
There is a background constant $C$ such that every
$V\in C^\infty_c(S^2T^*M)$ satisfies
\begin{equation}\label{eq:A-graph-estimate}
 \|\bar\nabla^2V\|_{L^2_\nu}
 +\|\sqrt{\bar f}\,\bar\nabla V\|_{L^2_\nu}
 \leq C\bigl(\|LV\|_{L^2_\nu}
              +\|V\|_{H^1_\nu}\bigr),
 \qquad L=-\A .
\end{equation}
The operator \(L\) is the closure of its restriction to
\(C^\infty_c(S^2T^*M)\); in particular this space is an operator core
for the Friedrichs realization.  Estimate
\eqref{eq:A-graph-estimate} therefore extends to every \(V\in D(L)\).
Every spectral projection of \(L\) preserves \(D(L)\), so the same
estimate applies to the projected tensor.
\end{lemma}

\begin{proof}
The weighted Bochner formula for the covariant tensor
$\bar\nabla V$, integrated against $d\nu$, gives
\begin{equation}\label{eq:tensor-weighted-Bochner}
 \begin{split}
 \|\bar\nabla^2V\|_{L^2_\nu}^2
 \leq{}&
 \|\bar\Delta_{\bar f}V\|_{L^2_\nu}^2
 +C\|\bar\nabla V\|_{L^2_\nu}^2\\
 &+C\int_M
 |\bar\nabla\overline{\Rm}|\,|V|\,|\bar\nabla V|\,d\nu .
 \end{split}
\end{equation}
Indeed, commute $\bar\nabla$ with $\bar\Delta_{\bar f}$ in
$\frac12\bar\Delta_{\bar f}|\bar\nabla V|^2$ and use
$\overline{\Ric}+\bar\nabla^2\bar f=\frac12\bar g$; the commutators
are
$\overline{\Rm}*\bar\nabla V$ and
$\bar\nabla\overline{\Rm}*V$.  The FIK shrinker has bounded
$\overline{\Rm}$ and $\bar\nabla\overline{\Rm}$, so Cauchy--Schwarz in
\eqref{eq:tensor-weighted-Bochner} yields
\[
 \|\bar\nabla^2V\|_{L^2_\nu}
 \leq C\bigl(
 \|\bar\Delta_{\bar f}V\|_{L^2_\nu}
 +\|V\|_{H^1_\nu}\bigr).
\]
Since
\[
 \bar\Delta_{\bar f}V=-LV-2\overline{\Rm}(V),
\]
the first term on the left of \eqref{eq:A-graph-estimate} has the
claimed bound.  Lemma~\ref{lem:first-moment}, applied to the compactly
supported tensor $\bar\nabla V$, then gives
\[
 \|\sqrt{\bar f}\,\bar\nabla V\|_{L^2_\nu}
 \leq2\bigl(
 \|\bar\nabla V\|_{L^2_\nu}
 +\|\bar\nabla^2V\|_{L^2_\nu}\bigr),
\]
which proves \eqref{eq:A-graph-estimate}.

For completeness, since
\(d\nu=(4\pi)^{-2}e^{-\bar f}\,dV_{\bar g}\), conjugation by the
unitary multiplication
\[
 \mathcal U_\nu V:=(4\pi)^{-1}e^{-\bar f/2}V
\]
from \(L^2_\nu(S^2T^*M)\) to
\(L^2(dV_{\bar g};S^2T^*M)\) changes \(L\) into a generalized
Schr\"odinger operator on the complete bundle \(S^2T^*M\): its
principal part is the connection Laplacian and its smooth endomorphism
potential is bounded below (its scalar part is quadratic in
\(|\bar\nabla\bar f|\), up to bounded shrinker terms).
The vector-bundle essential-self-adjointness result
\cite[Corollary~2.9]{BravermanMilatovicShubin} for semibounded
generalized Schr\"odinger operators on complete manifolds shows that
\(C^\infty_c(S^2T^*M)\) is an operator core.

Let \(V_n\to V\) in the graph norm of \(L\).  Semiboundedness and
\[
 \|V_n-V_m\|_{H^1_\nu}^2
 \leq
 C\bigl\langle
  (L+C_{\rm Fr})(V_n-V_m),V_n-V_m
 \bigr\rangle_f
\]
show that \(V_n\) is also Cauchy in \(H^1_\nu\).  Applying
\eqref{eq:A-graph-estimate} to \(V_n-V_m\) and passing to the limit
extends the estimate to \(D(L)\).  Finally, the spectral calculus
gives \(LP_\Lambda=P_\Lambda L\) on \(D(L)\) for every spectral
projection \(P_\Lambda\), proving the last assertion.
\end{proof}

\begin{proposition}[Post-bootstrap \(H^1\) recovery at every strict bootstrap rate]
\label{prop:instantaneous-H1-spectral-rate}
Fix $0<\theta<\beta$.  In addition to the hypotheses of
Theorem~\ref{thm:robust-modulated-three-region}, assume the concrete
adaptive-graft forcing bounds
\begin{equation}\label{eq:H1-recovery-tails}
 \sum_{j=0}^8
 \|\rho_\tau\mathcal Y_{j,\tau}-Y_j\|_{L^2_\nu}
 +\|\rho_\tau\E(\tau)\|_{L^2_\nu}
 +\|\mathcal C_\rho[h](\tau)\|_{L^2_\nu}
 +\|\mathcal C_\rho[h](\tau)\|_{H^1_\nu}
 \leq Ce^{-ce^\tau}.
\end{equation}
Then, after increasing $\tau_0$ and decreasing the bootstrap
smallness once, for $\tau_0+1\leq\tau\leq\tau_1$,
\begin{align}
 \|H(\tau)\|_{H^1_\nu}
 &\leq C\varepsilon e^{-\theta\tau}
       +Ce^{-ce^\tau},
 \label{eq:instantaneous-H1-spectral-rate}\\
 |a(\tau)|+|b(\tau)|
 &\leq C\varepsilon^2e^{-2\theta\tau}
       +Ce^{-ce^\tau}.
 \label{eq:instantaneous-velocity-spectral-rate}
\end{align}
The constants are independent of the finite endpoint $\tau_1$.
\end{proposition}

\begin{proof}
Put $L=-\A$ and
\[
 X(\tau)=\langle LH,H\rangle_{L^2_\nu}.
\]
Because $H\perp\mathcal Z$, Lemma~\ref{lem:coercivity} and
\eqref{eq:A-form} give
\begin{equation}\label{eq:X-H1-equivalence}
 c\|H\|_{H^1_\nu}^2
 \leq X
 \leq C\|H\|_{H^1_\nu}^2.
\end{equation}
The tensor $H=\rho_\tau h$ is smooth and compactly supported, so
Lemma~\ref{lem:weighted-graph-estimate} applies to it without a
graph-core density assertion.

On $\supp\rho_\tau$, the nesting of the two cutoffs gives
$\chi_\tau=1$.  The exact equation for $H$ may therefore be written
\begin{equation}\label{eq:H-graph-decomposition}
 \partial_\tau H
 =-LH+\sum_{j=0}^8c_jY_j+\mathcal R,
 \qquad c=(a,b_1,\ldots,b_8),
\end{equation}
where $\mathcal R$ contains the localized quasilinear term, the
action terms on $h$, the defects of the effective columns, the
moving-cutoff commutator, and the graft forcing.

The schematic expansion \eqref{eq:Q-schematic}, the identities
\[
 \rho_\tau\bar\nabla h
 =\bar\nabla H-h\bar\nabla\rho_\tau,
 \qquad
 \rho_\tau\bar\nabla^2h
 =\bar\nabla^2H
  -2\bar\nabla\rho_\tau*\bar\nabla h
  -h\bar\nabla^2\rho_\tau,
\]
and \eqref{eq:A-graph-estimate} imply
\begin{equation}\label{eq:Q-graph-bound}
 \|\rho_\tau\Q(h)\|_{L^2_\nu}
 \leq C\delta
 \bigl(\|LH\|_{L^2_\nu}+X^{1/2}\bigr)
 +Ce^{-ce^\tau},
\end{equation}
where $\delta$ is the size of the global $C^2$ box.  All terms
containing a derivative of $\rho_\tau$ are supported where
$\bar f\simeq e^\tau$ and are Gaussian-superexponentially small.
The first-moment estimate and \eqref{eq:A-graph-estimate} likewise
give, for each geometric action,
\begin{equation}\label{eq:B-graph-bound}
 \|\rho_\tau\mathscr B_jh\|_{L^2_\nu}
 \leq C\bigl(\|LH\|_{L^2_\nu}+X^{1/2}\bigr)
      +Ce^{-ce^\tau}.
\end{equation}
Using the already obtained small bound for $|c|$ and
\eqref{eq:H1-recovery-tails}, and then decreasing $\delta$, we obtain
\begin{equation}\label{eq:R-graph-pairing}
 2|\langle LH,\mathcal R\rangle|
 \leq \|LH\|_{L^2_\nu}^2+CX+Ce^{-ce^\tau}.
\end{equation}

The global geometric columns make no contribution:
\[
 \langle LH,Y_j\rangle
 =\langle H,LY_j\rangle=0,
\]
because $L\mathcal Z\subset\mathcal Z$ and
$H\perp\mathcal Z$.  Since $L$ is time independent,
\eqref{eq:H-graph-decomposition} and
\eqref{eq:R-graph-pairing} yield
\[
 X'\leq CX+Ce^{-ce^\tau}.
\]
For any $s\in[\tau-1,\tau]$, Gronwall gives
$X(\tau)\leq CX(s)+Ce^{-ce^\tau}$.  Average in $s$ and use
\eqref{eq:X-H1-equivalence} together with the endpoint-independent
dissipation estimate \eqref{eq:tail-diss}:
\[
 X(\tau)
 \leq C\int_{\tau-1}^{\tau}
       \|H(s)\|_{H^1_\nu}^2\,ds+Ce^{-ce^\tau}
 \leq C\varepsilon^2e^{-2\theta\tau}+Ce^{-ce^\tau}.
\]
This proves \eqref{eq:instantaneous-H1-spectral-rate}.  Substitution
in the exact receding Gram estimate
\eqref{eq:receding-velocity} proves
\eqref{eq:instantaneous-velocity-spectral-rate}.
\end{proof}

\begin{remark}
The concrete $L^2_\nu$ forcing hypothesis in
\eqref{eq:H1-recovery-tails} is essential here.  Under only an
abstract $H^{-1}_\nu$ forcing bound, the pairing
$\langle LH,\E\rangle$ is unavailable and integrated dissipation
does not exclude instantaneous $H^1_\nu$ spikes.  Thus
Proposition~\ref{prop:instantaneous-H1-spectral-rate} is a
post-bootstrap upgrade for the adaptive graft, not an extra
  conclusion of Theorem~\ref{thm:abstract-energy}.
\end{remark}

\begin{corollary}[Quantitative core smoothing]
\label{cor:quantitative-core-smoothing}
Under the hypotheses of
Theorem~\ref{thm:robust-modulated-three-region}, suppose that the
direct columns are the adaptive effective columns
\eqref{eq:effective-column-zero}--\eqref{eq:effective-column-j} and that
the adaptive interval carries the exact support-separation identities
\eqref{eq:coarse-Gram-support-separation}.
Then, for every
$K\Subset K^+\Subset M$ and every integer $m\geq0$ there is
$C_{K,K^+,m}$, independent of $\tau_1$, such that
\begin{equation}\label{eq:quantitative-core-smoothing}
 \|h(\tau)\|_{C^m(K)}
 \leq C_{K,K^+,m}\varepsilon e^{-\theta\tau}
      +C_{K,K^+,m}e^{-c_*e^\tau}
\end{equation}
whenever $\tau\geq\tau_0+1$ and
$K^+\subset\{\bar f<e^{\tau-3/4}\}$.
\end{corollary}

\begin{proof}
For \(s\in[\tau-\tfrac34,\tau]\), the condition on $K^+$ gives
$K^+\subset\{\bar f<e^s\}\subset\supp\rho_s$, so
$\rho_s=\chi_s=1$ and $\E(s)=0$ on $K^+$.  Moreover,
\eqref{eq:coarse-Gram-support-separation} and the exact formulas
\eqref{eq:effective-column-zero}--\eqref{eq:effective-column-j} give
\[
 \mathcal Y_{j,s}=Y_j\quad\hbox{on }K^+,
 \qquad0\leq j\leq8.
\]
Hence $h=H$ there, and
\eqref{eq:tail-L2} gives
\begin{equation}\label{eq:core-cylinder-L2}
 \sup_{s\in[\tau-\frac34,\tau]}
 \|h(s)\|_{L^2(K^+)}
 \leq C_{K^+}\varepsilon e^{-\theta\tau}
      +C_{K^+}e^{-ce^\tau}.
\end{equation}
Put \(N=m+4\), and choose a fixed nested chain
\[
 K=K_0\Subset K_1\Subset\cdots\Subset K_N\Subset K^+.
\]
Thus every spatial derivative of every direct column required below is
a fixed-background derivative of \(Y_j\).  No smoothing of the
ODE-carried finite-order entrance data is being asserted.

Theorem~\ref{thm:robust-modulated-three-region} has already been proved,
so its instantaneous estimate \eqref{eq:three-region-velocity} applies
on the recent cylinder.  In particular, there is an endpoint-independent
constant \(Q_*<\infty\), depending only on the fixed theorem package,
such that
\begin{equation}\label{eq:core-smoothing-q-upper-bound}
 0\leq q(s)\leq Q_*,
 \qquad
 q(s)^2\leq Q_*q(s).
\end{equation}
We may therefore perform the nested-cylinder argument directly in the
fixed background coordinates, without conjugating by the phase flow.
This fixed-coordinate formulation introduces no compact-set-dependent
phase threshold into the present all-order corollary.

For completeness, record the additional estimate needed in these fixed
coordinates.  On the fixed compact set \(K^+\), put
\[
 X_0=-\bar\nabla\bar f,
 \qquad X_j=W_j\quad(1\leq j\leq8).
\]
Although these fields have scale-normalized growth on the complete end,
every spatial derivative of every \(X_j\) is bounded on the fixed set
\(K^+\).  After commuting \(r\) fixed-background derivatives with the
equation, the top part of the transport is
\(c_j(s)\bar\nabla_{X_j}\bar\nabla^rh\), where
\(|c(s)|\leq q(s)\).  If \(\zeta\) is one of the nested spatial cutoffs,
integration by parts gives
\[
 \begin{split}
 \left|
  \int_{K^+}\zeta^2
   \left\langle
    c_j\bar\nabla_{X_j}\bar\nabla^rh,
    \bar\nabla^rh
   \right\rangle d\mu_{\bar g}
 \right|
 \leq{}& C_{m,K^+}q(s)
       \|\zeta\bar\nabla^rh\|_{L^2}^2\\
 &+C_{m,K^+}q(s)
       \|h\|_{H^r(\operatorname{supp}\bar\nabla\zeta)}^2.
 \end{split}
\]
All commutators in which a derivative falls on \(X_j\) have order at
most \(r\) in \(h\) and satisfy the same lower-order bound.  The
zeroth-order pieces of the Lie derivatives are handled identically.
For the direct source, the fixed all-order bounds for \(Y_j\), Young's
inequality, and \eqref{eq:core-smoothing-q-upper-bound} give, at every
order used in the induction,
\[
 \left|
  \int_{K^+}\zeta^2
   \left\langle c_j\bar\nabla^rY_j,
   \bar\nabla^rh\right\rangle d\mu_{\bar g}
 \right|
 \leq C_{m,K^+}q(s)
       \bigl(1+\|\zeta\bar\nabla^rh\|_{L^2}^2\bigr).
\]
Equivalently, if a standard Young estimate first produces \(q^2\), the
second inequality in \eqref{eq:core-smoothing-q-upper-bound} converts it
to \(Q_*q\).  No derivative of \(q\) and no endpoint-dependent bound
for it is used.

The forcing vanishes on the cylinder.  The principal quasilinear
commutators are exactly those estimated in
\eqref{eq:nested-tame-principal-commutator}--%
\eqref{eq:nested-tame-lower-commutator}; the coarse \(C^2\) box absorbs
the top-order coefficient term, and the remaining products are treated
by the same tame induction.  Consequently the fixed-coordinate local
energies satisfy the recursion
\eqref{eq:nested-local-energy-recursion}, with constants enlarged only
by \(C_{m,K^+}(1+Q_*)\).  Gronwall uses
\[
 \int_{\tau-\frac34}^{\tau}(1+q(s))\,ds
 \leq\frac34+P_{\tau_1}\!\left(\tau-\frac34\right),
\]
and hence is uniform in the finite endpoint.  The cutoff--Moser base
step and the augmented Bernstein induction
\eqref{eq:nested-Moser-scalar-inequality}--%
\eqref{eq:nested-Bernstein-induction} now give
\begin{equation}\label{eq:core-all-order-interior}
 \|h(\tau)\|_{C^m(K)}
 \leq C_{K,K^+,m}\left[
  \sup_{s\in[\tau-\frac34,\tau]}
       \|h(s)\|_{L^2(K^+)}
  +\left(\int_{\tau-\frac34}^{\tau}q(s)\,ds\right)^{1/2}
 \right].
\end{equation}
Here no initial-face term occurs because
\(\tau\geq\tau_0+1\).  The constant depends only on the fixed compact
chain, the background jets on \(K^+\), the coarse ellipticity and
\(C^2\) box, the fixed all-order jets of the tensors \(Y_j\), and the
endpoint-independent number \(Q_*\).

Finally,
\[
 \int_{\tau-\frac34}^{\tau}q(s)\,ds
 \leq P_{\tau_1}\!\left(\tau-\frac34\right)
 \leq C\varepsilon^2e^{-2\theta\tau}+Ce^{-ce^\tau}.
\]
Substitution of this estimate and
\eqref{eq:core-cylinder-L2} into
\eqref{eq:core-all-order-interior} proves
\eqref{eq:quantitative-core-smoothing}.  The factor \(1/(4e)\) in
\eqref{eq:three-region-common-Gaussian-exponent} absorbs both the square
root and the fixed translation
\(\tau\mapsto\tau-\tfrac34\), so the Gaussian exponent is the already
fixed \(c_*\).  All constants are determined on one fixed-length
cylinder from endpoint-independent theorem constants; hence they are
independent of \(\tau_1\).  This argument is independent of the proof of
Theorem~\ref{thm:robust-modulated-three-region} and alters none of its
thresholds or conclusions.
\end{proof}

\begin{remark}[A priori theorem versus geometric continuation]
Theorem~\ref{thm:robust-modulated-three-region} is endpoint-independent
and therefore improves the pointwise faces of every finite bootstrap
interval.  The construction and simultaneous control of the adaptive
chart, graft, Gram matrix, and underlying closed Ricci flow are proved
together later in
Theorem~\ref{thm:prepared-entrance-continuation}; no continuation
alternative is being assumed here.
\end{remark}

\begin{remark}[Relation to Stolarski's estimates]
\label{rem:honest-Stolarski-import}
The unmodulated strict barriers and the scale-invariant local
rescaling are taken from
\cite[Theorem~6.1, Lemmas~6.5--6.9]{Stolarski}.  The argument above
modifies these inputs in three respects:
\begin{enumerate}
\item the specially prepared initial eigensum is replaced by the open
      scale-adapted entrance condition
      \eqref{eq:arbitrary-stable-entrance};
\item the direct feedback columns are absorbed by the negative future
      phase tail $-K_0P_{\tau_1}$, while their transports are absorbed
      by $J$;
\item derivative recovery is performed in raw fixed-background norms
      after conjugating the full transport, so
      its constants depend on $\int q$ and not on an unavailable
      prior bound for $\sup q$.
\end{enumerate}
Thus Theorem~\ref{thm:robust-modulated-three-region} is a complete
pointwise a priori theorem for a controlled solution satisfying that
theorem's stated hypotheses.  It
does not by itself construct the adaptive harmonic-map chart or prove
its graft estimates; those are coupled to this theorem by the
finite-endpoint coupled bootstrap.
\end{remark}

\section{Adaptive grafts and geometric continuation}
\label{sec:adaptive-continuation}

The purpose of this section is to close the geometric part of the
continuation argument.  There are three points.  First, the prescribed
radial and gauge motion can be factored out of the harmonic-map heat
flow exactly.  Second, on the fixed physical grafting annulus the
resulting target has finite total variation in every scale-invariant
$C^m$ norm.  Third, the defect of the interpolated metric is supported
on an annulus which recedes through $\bar f\simeq e^\tau$ and hence is
Gaussian-superexponentially small.  None of these statements uses the
spectral gap.

\begin{lemma}[Buffered local Ricci control]
\label{lem:buffered-local-Ricci-control}
Fix \(n\in\mathbb N\), \(0<\alpha<1\), \(J\geq0\),
\(\Lambda<\infty\), and
$\upsilon_{\rm har}>0$.  There are
$\delta_{\rm RF},C_{\rm RF}>0$ with the following property.  Let
$U\Subset U^+$ be open subsets of a closed manifold $\mathcal X$ and
let $R>0$
satisfy
\[
 d_{G_0}(U,\mathcal X\setminus U^+)\geq10R.
\]
Suppose that on $U^+$
\[
 \inf_{x\in U^+}r_{\rm har}(G_0,x)
 \geq\upsilon_{\rm har}R,\qquad
 R^{2+\ell}|(\nabla^{G_0})^\ell\Rm_{G_0}|
 \leq\Lambda,\qquad0\leq\ell\leq J.
\]
Here \(r_{\rm har}\) is precisely the fixed
\(C^{2,\alpha},Q_{\rm har}\) harmonic radius of
\eqref{eq:fixed-harmonic-radius-convention}; no tolerance or
normalization is reselected in this lemma.
If the closed Ricci flow with initial value $G_0$ exists on
$[0,t_*]$, then on
$U\times[0,\min\{t_*,\delta_{\rm RF}R^2\}]$,
\begin{equation}\label{eq:buffered-local-Ricci-control}
 R^{2+\ell}|(\nabla^{G(t)})^\ell\Rm_{G(t)}|
 \leq C_{\rm RF},\qquad0\leq\ell\leq J.
\end{equation}
The same conclusion follows if the harmonic-radius hypothesis is
replaced by the quantitative bound
\[
 \operatorname{inj}_{G_0}\geq\upsilon_{\rm inj}R
 \quad\text{on }U^+
\]
together with the displayed curvature bounds.  In that formulation
\(\delta_{\rm RF}\) and \(C_{\rm RF}\) may also depend on the fixed
constant \(\upsilon_{\rm inj}>0\).
\end{lemma}

\begin{proof}
If \(n=1\), the Riemann curvature tensor vanishes identically and the
Ricci flow is stationary.  Hence
\eqref{eq:buffered-local-Ricci-control} holds in every derivative order
with, for example, \(C_{\rm RF}=1\) and any fixed
\(\delta_{\rm RF}>0\).  We may therefore assume \(n\geq2\).
Rescale \(R\) to one.  In the harmonic-radius alternative choose
\[
 \rho_{\rm har}
 :=c_{\rm har}\min\left\{
       1,\upsilon_{\rm har},(1+\Lambda)^{-1/2}\right\},
 \qquad 0<c_{\rm har}<\frac18,
\]
where \(c_{\rm har}=c_{\rm har}(n,Q_{\rm har})\) is fixed below.
For every \(x\in U\), one has
\(4\rho_{\rm har}<\upsilon_{\rm har}\leq r_{\rm har}(G_0,x)\).
Since \(r_{\rm har}\) is a supremum over admissible radii and the
convention is downward closed, there is an admissible chart on
\(B_{G_0}(x,4\rho_{\rm har})\); no attainment of the supremum is being
used.  Restrict once more to the admissible chart on
\(B_{G_0}(x,\rho_{\rm har})\).  In its normalized coordinates one has
\[
 B_{\mathbb R^n}(0,Q_{\rm har}^{-1/2})
 \subset \Omega_{\rho_{\rm har}},
 \qquad
 \sqrt{\det\widetilde G_0}\geq Q_{\rm har}^{-n/2}.
\]
Indeed, if a Euclidean radial segment first left
\(\Omega_{\rho_{\rm har}}\) at radius less than
\(Q_{\rm har}^{-1/2}\), its \(\widetilde G_0\)-length would be less
than one, contradicting that the chart domain is the unit
\(\widetilde G_0\)-metric ball.  Consequently
\[
 \operatorname{Vol}_{G_0}B_{G_0}(x,\rho_{\rm har})
 \geq \omega_nQ_{\rm har}^{-n}\rho_{\rm har}^{\,n}
 =:v_{\rm har}\rho_{\rm har}^{\,n}.
\]
Moreover \(\Lambda\leq\rho_{\rm har}^{-2}\), and the ten-unit buffer
places every ball just used inside \(U^+\).

In the injectivity-radius alternative set
\[
 \rho_{\rm inj}
 :=c_{\rm inj}\min\left\{
       1,\upsilon_{\rm inj},(1+\Lambda)^{-1/2}\right\},
\]
with \(0<c_{\rm inj}=c_{\rm inj}(n)<1/8\) sufficiently small.  The
injectivity lower bound and the curvature bound, applied in exponential
coordinates and combined with the Rauch comparison estimates, give
\[
 \operatorname{Vol}_{G_0}B_{G_0}(x,\rho_{\rm inj})
 \geq v_{\rm inj}\rho_{\rm inj}^{\,n},
 \qquad
 v_{\rm inj}=v_{\rm inj}(n,\Lambda,\upsilon_{\rm inj})>0,
\]
and again \(\Lambda\leq\rho_{\rm inj}^{-2}\).  Thus in either
alternative there are fixed numbers
\(\rho>0\) and \(v_0>0\), with precisely the dependencies stated in
the lemma, for which the initial curvature and volume hypotheses of
Peng Lu's local doubling-time theorem
\cite[Theorem~1.2]{PengLu} hold at scale \(\rho\).  That theorem
therefore gives, with
\(\varepsilon_0=\varepsilon_0(n,v_0)>0\),
\[
 \sup_{B_{G(t)}(x,\varepsilon_0\rho)}
 |\Rm_{G(t)}|\leq(\varepsilon_0\rho)^{-2},
 \qquad
 0\leq t\leq
 \min\{t_*,(\varepsilon_0\rho)^2\}.
\]
The completeness and bounded-curvature hypotheses of that theorem hold
because the present flow is closed.  Since \(x\in U\) was arbitrary,
this is the required \(\ell=0\) estimate, after absorbing the fixed
\(\rho\) into \(\delta_{\rm RF}\) and \(C_{\rm RF}\); unlike the raw
pseudolocality estimate, it has no \(t^{-1}\) term.

The base curvature estimate gives uniform metric comparison on the
same time interval.  Choose a finite-overlap cover by nested metric
balls with a fixed positive separation.  A Calabi-smoothed evolving
distance cutoff \(\zeta\), equal to one on the smaller ball and
supported in the next larger ball, can be chosen so that
\[
 |\nabla^{G(t)}\zeta|^2
 +|(\partial_t-\Delta_{G(t)})\zeta|\leq C
\]
in the barrier sense; the curvature bound, metric comparison, and
fixed buffer make \(C\) uniform.  These are the only cutoff bounds in
the following Bernstein argument.  Apply the local
Shi--Bernstein induction~\cite{Shi} to
\[
 \partial_t\nabla^\ell\Rm
 =\Delta\nabla^\ell\Rm
  +\sum_{p+q=\ell}\nabla^p\Rm*\nabla^q\Rm .
\]
At induction level \(\ell\), the initial term is bounded by the
hypothesis for that \(\ell\), while all lower-order terms have already
been controlled.  The cutoff errors are uniform by the base curvature
bound, metric comparison, and the fixed buffer.  Thus the Bernstein
maximum principle gives a bound depending on the initial
\(\ell\)-jet, rather than the positive-time factor \(t^{-\ell/2}\).
This proves the displayed estimate for every \(0\leq\ell\leq J\),
including \(t=0\).
A finite-overlap nested covering absorbs the loss from $U^+$ to $U$.
\end{proof}

\begin{lemma}[Buffered Ricci--DeTurck coefficients and transported markings]
\label{lem:buffered-Ricci-DeTurck-coefficients}
Fix the data in Lemma~\ref{lem:buffered-local-Ricci-control}, an integer
\(r\geq4\) with \(J\geq r-2\), and a smooth reference metric
\(\widehat G\) on \(U^+\)
whose scale-\(R\) coefficient bounds through order \(r+1\) are fixed.
Here and below \(C_R^{j,\alpha}\) denotes the dimensionless
coefficient norm in the scale-\(R\) atlas.
Choose smooth nested buffers
\begin{equation}\label{eq:Ricci-DeTurck-nested-buffers}
 U=U^0\Subset U^1\Subset U^2\Subset U^3\Subset U^4
 \Subset U^5=U^+,
 \qquad
 d_{G_0}(U^\ell,\mathcal X\setminus U^{\ell+1})
 \geq c_{\rm buf}R,\quad 0\leq\ell\leq4 .
\end{equation}
where \(c_{\rm buf}>0\) depends only on the original \(10R\) buffer.
Assume, in addition, that \(G_0\) has a scale-normalized
\(C^{r,\alpha}\) coefficient bound \(K_r\) in one finite atlas
subordinate to these buffers.  After decreasing \(\delta_{\rm RF}\),
put
\[
 T=\min\{t_*,\delta_{\rm RF}R^2\}.
\]
There are diffeomorphisms \(\chi(t)\), defined on \(U^3\) and onto
their images, with the uniform buffered containment
\[
 U^1\Subset\chi(t)(U^2)\Subset U^3,
\]
such that on \(U^2\times[0,T]\)
\begin{equation}\label{eq:buffered-DeTurck-representation}
 G(t)=\chi(t)^*\widetilde G(t),\qquad
 \chi(0)=\operatorname{Id}.
\end{equation}
Here \(\widetilde G\) satisfies, on the corresponding image,
\begin{equation}\label{eq:buffered-Ricci-DeTurck-equation}
 \partial_t\widetilde G
 =-2\Ric_{\widetilde G}
   +\Lie_{B_{\widehat G}(\widetilde G)}\widetilde G,
 \qquad
 B_{\widehat G}(\widetilde G)^k
 =\widetilde G^{ij}
   \bigl(\Gamma(\widetilde G)^k_{ij}
        -\Gamma(\widehat G)^k_{ij}\bigr),
\end{equation}
and
\begin{equation}\label{eq:buffered-DeTurck-coefficients}
 \sup_{0\leq t\leq T}
 \left(
  \|\widetilde G(t)\|_{C_R^{r,\alpha}(U^1)}
  +\|\chi(t)\|_{C_R^{r-1,\alpha}(U^1)}
  +\|\chi(t)^{-1}\|_{C_R^{r-1,\alpha}(U^1)}
 \right)
 \leq C_{\rm RD}.
\end{equation}
Here and below a map norm on a moving image is read in the fixed
buffered atlas after restriction to the common image of \(U^1\).
The constant \(C_{\rm RD}\) depends only on
\(n,r,\alpha,\Lambda,\upsilon_{\rm har},K_r,c_{\rm buf}\) and the
fixed coefficient bounds through order \(r+1\) for \(\widehat G\).
The construction is used only through the asserted interior
restrictions and estimates.  Their constants depend on
\(\widehat G\) and the auxiliary localization only through the recorded
coefficient and buffer bounds.  The particular maps
\(\chi,\chi^{-1}\), and hence
\(\widetilde G\), may depend on the auxiliary extension of
\(\widehat G\) and on the chosen remote Dirichlet data; no
extension-independence or canonical-gauge assertion is made.

The construction is uniform for two flows in a common coefficient
ball.  If their existence endpoints are \(t_{1,*}\) and \(t_{2,*}\),
set
\[
 T_{12}=\min\{t_{1,*},t_{2,*},\delta_{\rm RF}R^2\}.
\]
For \(4\leq j\leq r\), let \(\mathfrak D_j(V,t)\) be the sum on \(V\)
of the \(C_R^{j,\alpha}\) norm of
\(\widetilde G_1-\widetilde G_2\) and the
\(C_R^{j-1,\alpha}\) norms of
\(\chi_1-\chi_2\) and
\(\chi_1^{-1}-\chi_2^{-1}\).  Then the genuinely local estimate is
\begin{equation}\label{eq:buffered-DeTurck-difference-local}
 \mathfrak D_j(U^0,t)
 \leq
 C_j\|G_{1,0}-G_{2,0}\|_{C_R^{j,\alpha}(U^2)}
 +C_jR^{-2}\int_0^t
       \mathfrak D_j(U^1,q)\,dq
\end{equation}
for \(0\leq t\leq T_{12}\).  In particular, suppose that a fixed
finite collection of triples
\[
 U_a^0\Subset U_a^1\Subset U_a^2,\qquad
 1\leq a\leq N_{\rm loc},
\]
with comparable scales \(R_a\), covers the closed manifold
\(\mathcal X\), every \(U_a^1\) is covered by members of
\(\{U_b^0\}_{b=1}^{N_{\rm loc}}\) with uniformly bounded overlap, and the
preceding hypotheses hold in each \(U_a^2\).  Assume here that the
local reference metrics are restrictions of one fixed smooth metric
on \(\mathcal X\), and take the local constructions to be restrictions
of the resulting common global DeTurck gauge.  Summing
\eqref{eq:buffered-DeTurck-difference-local} and applying Gronwall gives
the closed estimate
\begin{equation}\label{eq:buffered-DeTurck-difference}
 \sup_{0\leq q\leq t}
 \sum_{a=1}^{N_{\rm loc}}\mathfrak D_{j,a}(U_a^0,q)
 \leq
 C_j\sum_{a=1}^{N_{\rm loc}}
 \|G_{1,0}-G_{2,0}\|_{C_{R_a}^{j,\alpha}(U_a^2)},
 \qquad 0\leq t\leq T_{12}.
\end{equation}
Here \(\mathfrak D_{j,a}\) is the preceding quantity measured at scale
\(R_a\).
Thus no norm on an unestimated larger buffer remains on the
right-hand side.
The endpoint \(j=r\) is intentional.  For the metric component,
\(J\geq r-2\) controls the differentiated Ricci terms, while an
order-\(j\) principal-coefficient difference is needed only in
\(C^{j-2,\alpha}\) and multiplies the second derivative of the bounded
order-\(j\) metric.  For the map component at \(j=r\), use exactly one
integrable derivative of Ricci--DeTurck smoothing on the unused
buffer:
\[
 \left(\frac{t}{R^2}\right)^{1/2}
 \|\widetilde G_i(t)\|_{C_R^{r+1,\alpha}(U^1)}
 \leq C,\qquad 0<t\leq T .
\]
Thus \(B_{\widehat G}(\widetilde G_i)\) is \(C_R^{r,\alpha}\) with
coefficient \(O(t^{-1/2})\).  Subtraction of the two flow equations
then estimates \(\chi_1-\chi_2\) in \(C_R^{r-1,\alpha}\) by a Volterra
inequality with the integrable kernel \(t^{-1/2}\); weakly singular
Gronwall gives the asserted endpoint bound.  The same smoothing gives
the individual \(C_R^{r,\alpha}\) map bound needed when the inverse
identity is differentiated.  Thus no order-\((r+1)\), and a fortiori
no order-\((r+2)\), initial metric jet is concealed at the displayed
endpoint.
This last finite-cover conclusion is used only when the smallest sets
really cover the entire closed manifold.  It is not invoked for the
proper noncollapsing exterior of a singular flow; that case is
terminated at a separated normalized interface in
Lemma~\ref{lem:inner-terminated-exterior-DeTurck}.

If a marking \(\iota:\mathcal X''\to M\) is fixed and
\(\overline{U^+}\Subset\mathcal X''\), assume, whenever the marking
conclusion is invoked, the fixed scale-normalized jet bound
\begin{equation}\label{eq:buffered-fixed-marking-jet}
 \|\iota\|_{C_R^{r,\alpha}(U^+;M)}
 +\|\iota^{-1}\|_{C_R^{r,\alpha}(\iota(U^+);\mathcal X'')}
 \leq K_{\iota,r}<\infty
\end{equation}
in the recorded source and target atlases.  Then the transported maps
are typed as
\begin{equation}\label{eq:transported-marking-definition}
 \begin{aligned}
  \widetilde\iota(t)&:
  \chi(t)(U^2)\longrightarrow\iota(U^2),
  &\widetilde\iota(t)&=\iota\circ\chi(t)^{-1},\\
  \widetilde\iota(t)^{-1}&:
  \iota(U^2)\longrightarrow\chi(t)(U^2),
  &\widetilde\iota(t)^{-1}&=\chi(t)\circ\iota^{-1}.
 \end{aligned}
\end{equation}
obey the corresponding two-sided \(C_R^{j-1,\alpha}\) bounds and, for
two flows using this same fixed marking,
\begin{equation}\label{eq:buffered-transported-marking-difference}
 \begin{split}
 &\|\widetilde\iota_1(t)-\widetilde\iota_2(t)\|_
       {C_R^{j-1,\alpha}(U^0;M)}
 +\|\widetilde\iota_1(t)^{-1}-\widetilde\iota_2(t)^{-1}\|_
       {C_R^{j-1,\alpha}(\iota(U^0);\mathcal X'')}\\
 &\qquad\leq C_{j,\iota}\left(
 \|\chi_1(t)^{-1}-\chi_2(t)^{-1}\|_
       {C_R^{j-1,\alpha}(U^1)}
 +\|\chi_1(t)-\chi_2(t)\|_
       {C_R^{j-1,\alpha}(U^1)}\right),
 \qquad4\leq j\leq r .
 \end{split}
\end{equation}
Here \(C_{j,\iota}\) depends, in addition to the displayed geometric
package, on \(K_{\iota,r}\); the constants in the metric and gauge
conclusions \eqref{eq:buffered-DeTurck-coefficients}--%
\eqref{eq:buffered-DeTurck-difference} do not.  Moreover,
\begin{equation}\label{eq:buffered-marking-covariance}
 (\widetilde\iota(t))_*\widetilde G(t)=\iota_*G(t)
\end{equation}
on the common smaller marked set.
\end{lemma}

\begin{proof}
Rescale \(R=1\).  The harmonic-coordinate normalization at \(t=0\),
the curvature bound in Lemma~\ref{lem:buffered-local-Ricci-control},
the bounds through \(J\geq r-2\), and
\(\partial_tG=-2\Ric_G\) first give uniform bounded-geometry
coefficients, coordinate ellipticity, and metric equivalence on
\(U^4\times[0,T]\).  Fix a smooth extension of \(\widehat G\) across
the closed manifold and fixed remote Dirichlet data in order to pose
the auxiliary localized problem.  These choices may change the
resulting gauge, but only the interior restrictions and estimates
stated in the lemma will be used.  On a smooth domain between \(U^4\)
and \(U^+\), solve the localized harmonic-map heat equation
\begin{equation}\label{eq:localized-DeTurck-HMHF}
 \partial_t\psi=\Delta_{G(t),\widehat G}\psi,\qquad
 \psi(\,\cdot\,,0)=\operatorname{Id},
\end{equation}
with fixed Dirichlet data on the remote lateral boundary.  Standard
initial-boundary quasilinear theory constructs \(\psi\) for a uniform
short time.  Possible corner incompatibility at the remote boundary is
irrelevant to the initial-face interior estimates on \(U^4\).

After a further uniform decrease of \(\delta_{\rm RF}\), the
interior \(C^1\) estimate keeps \(\psi\) and its inverse in the
successive buffers in \eqref{eq:Ricci-DeTurck-nested-buffers}.  Set
\[
 \chi=\psi,\qquad
 \widetilde G=(\psi^{-1})^*G .
\]
This proves \eqref{eq:buffered-DeTurck-representation} by definition.
The Eulerian velocity of \(\psi\) is
\(-B_{\widehat G}(\widetilde G)\); differentiating the pushforward
therefore gives \eqref{eq:buffered-Ricci-DeTurck-equation}.  This also
identifies the locally constructed system with the prescribed Ricci
flow, rather than with an unrelated solution of a Dirichlet
Ricci--DeTurck problem.

In the fixed reference atlas, the principal part of
\eqref{eq:buffered-Ricci-DeTurck-equation} is
\(\widetilde G^{ab}\partial_a\partial_b\).  The ellipticity just proved,
the initial \(C^{r,\alpha}\) bound, and the initial-face interior
quasilinear Schauder estimate propagate the \(C^{r,\alpha}\) bound on
successively smaller buffers.  Since the DeTurck vector contains one
derivative of \(\widetilde G\), the equations
\[
 \partial_t\chi
 =-B_{\widehat G}(\widetilde G)\circ\chi
 \quad\text{and}\quad
  \partial_t\chi^{-1}
  =d\chi^{-1}\,B_{\widehat G}(\widetilde G)
  \]
  give the asserted \(C^{r-1,\alpha}\) estimates.  If \(T'<T\) were the
maximal time on which the construction, the buffer inclusion, and
these estimates held, the same interior estimates at \(T'\) would
retain strict ellipticity and \(C^1\)-invertibility margins and would
  restart \eqref{eq:localized-DeTurck-HMHF}.  Hence \(T'=T\).  This is
  the required construction, identification, and continuation
  argument.  Only the single integrable \(t^{-1/2}\) derivative gain
  recorded above is used at the top two-state map order; no
  nonintegrable full-smoothing estimate or extra initial jet is used.

For the finite-cover assertion, solve
\eqref{eq:localized-DeTurck-HMHF} on the closed manifold, with the
single global reference metric and no lateral boundary.  The finite
cover supplies the same uniform construction time and shows that its
restrictions satisfy all of the preceding local estimates.  Thus the
gauges agree on overlaps.  For two flows, subtract
  \eqref{eq:buffered-Ricci-DeTurck-equation} in the common reference
  gauge.  At order \(j\), the principal-coefficient difference is
  estimated in \(C_R^{j-2,\alpha}\) and multiplies the second derivatives
  of one metric in the same space.  Thus the quasilinear source is
  bounded by the \(C_R^{j,\alpha}\) solution difference and requires no
  derivatives above order \(j\).  The localized initial-face Schauder
  estimate gives the metric part on
  \(U^0\Subset U^1\Subset U^2\).  For \(j\leq r-1\), subtraction of the
  triangular flow equations gives the map and inverse-map parts
  directly.  At \(j=r\), the one-derivative interior smoothing estimate
  on the unused buffer gives
  \[
   \|B_{\widehat G}(\widetilde G_i(t))\|_{C_R^{r,\alpha}}
   \leq Ct^{-1/2}.
  \]
  The difference \(B_{\widehat G}(\widetilde G_1)
  -B_{\widehat G}(\widetilde G_2)\) is already controlled in
  \(C_R^{r-1,\alpha}\) by the metric part.  Hence the subtracted flow
  equations give a weakly singular Volterra inequality for the
  \(C_R^{r-1,\alpha}\) map difference.  Its \(t^{-1/2}\) kernel is
  integrable and is absorbed, after the fixed short-time reduction, by
  the standard weakly singular Gronwall lemma.  Differentiating the
  inverse identity and using the corresponding individual
  \(C_R^{r,\alpha}\) flow bounds gives the inverse-map estimate.
  Together these estimates give
  \eqref{eq:buffered-DeTurck-difference-local}.  Restoring scale inserts
the factor \(R^{-2}\) in the time integral.  On the stated finite
cover, a partition of unity and bounded overlap bound every larger-set
norm by the sum of the smaller-set norms.  Gronwall then proves
\eqref{eq:buffered-DeTurck-difference}.  Finally, the scale-one
Fa\`a di Bruno formula for
\(\iota\circ\chi^{-1}\) uses derivatives of \(\iota\) through order
\(j\) to estimate a \(C_R^{j-1,\alpha}\) difference.  Thus
\eqref{eq:buffered-fixed-marking-jet}, together with the individual and
difference bounds for \(\chi^{-1}\), gives
\eqref{eq:buffered-transported-marking-difference} and the asserted
individual marking bounds.  The identity
\(\widetilde\iota^{-1}=\chi\circ\iota^{-1}\) and the same estimate give
the inverse-marking terms with the identical derivative count.
Naturality of pushforward gives
\eqref{eq:buffered-marking-covariance}.
\end{proof}

\subsection{The adaptive target}

Recall the soliton radial flow \(\varphi_\tau\) fixed in
\eqref{eq:radial-flow-convention}.
Use the graft radius \(\Gamma\) and cutoff \(\eta=\eta_\Gamma\) already
fixed in \(\mathfrak P_{\rm prep}\), with
\[
 \eta=1\quad\text{on }\{\bar f\leq2\Gamma/3\},\qquad
 \eta=0\quad\text{on }\{\bar f\geq5\Gamma/6\},
\]
and whose transition set
$\Omega_\eta=\{0<\eta<1\}$ is compactly contained in
\[
 \left\{\frac23\Gamma<\bar f<\frac56\Gamma\right\}.
\]
We use fixed enlargements
\[
 \Omega_\eta\Subset\Omega_\eta^+
 \Subset\Omega_\eta^{++}
 \Subset\left\{\frac12\Gamma<\bar f<\Gamma\right\}.
\]
Fix with these enlargements numerical constants
\begin{equation}\label{eq:fixed-enlarged-collar-radial-range}
 0<c_\eta\leq C_\eta<\infty,
 \qquad
 c_\eta\Gamma\leq1+\bar f\leq C_\eta\Gamma
 \quad\hbox{on }\Omega_\eta^{++},
\end{equation}
uniformly for the allowed graft radii.  After increasing the fixed
lower graft-radius threshold, one may, for example, take any fixed
\(c_\eta<1/2\) and any fixed \(C_\eta>1\).
The constants below are uniform in time; they may depend on $\Gamma$
and on the indicated number of derivatives.

All differences of the maps \(R_\tau\) on the compact collar
\(\Omega_\eta^+\) are measured in one fixed finite background atlas,
or equivalently by the target exponential field
\[
 \mathfrak r_{s,\tau}(x)
 =\exp_{R_\tau(x)}^{-1}(R_s(x)).
\]
The phase budget is chosen so that this field, and the analogous
field for the inverse maps on their common enlarged domain, remain in
one normal neighborhood.  Thus
\(\|R_s-R_\tau\|_{C^m(\Omega_\eta^+)}\) below denotes the norm of
\(\mathfrak r_{s,\tau}\), after the fixed parallel-transport
identification.  On the AC end all map norms use the common rescaled
domain and range atlases fixed above.  No subtraction of
manifold-valued maps is intended.

Fix a numerical adaptive-position package
\begin{equation}\label{eq:adaptive-position-package}
 \Lambda_{\rm ad}\geq1,\qquad
 0<c_{\rm rad}\leq C_{\rm rad}<\infty,\qquad
 0<c_{\rm scl}<C_{\rm scl}<\infty .
\end{equation}

\begin{proposition}[Adaptive target and annulus tracking]
\label{prop:adaptive-target-tracking}
For every integer $m\geq0$ and every fixed package
\eqref{eq:adaptive-position-package} there are
\(\varepsilon_{\rm ph}\in(0,1]\), \(\tau_\Gamma<\infty\), and
\(C_m<\infty\) with the following property.  Choose these constants
before \(\tau_0,\tau_1\) and the entrance data, and set
\begin{equation}\label{eq:fixed-C-lambda}
 C_\lambda=e^{\varepsilon_{\rm ph}}.
\end{equation}

Let
\[
 U_\tau=\sum_{j=1}^8b_j(\tau)\chi_\tau W_j,\qquad
 \chi_\tau=\chi(e^{-\tau}\bar f),
\]
and fix a finite bootstrap endpoint $\tau_1>\tau_0$.  Suppose on
\([\tau_0,\tau_1)\) that \(\tau_0\geq\tau_\Gamma\) and
\begin{equation}\label{eq:adaptive-phase-budget}
 \lambda_\tau=-(1+a)\lambda,\qquad
 c_{\rm scl}e^{-\tau}\leq\lambda(\tau)\leq C_{\rm scl}e^{-\tau},\qquad
 P_{\tau_1}(\tau_0):=
 \int_{\tau_0}^{\tau_1}(|a|+|b|)\,d\tau
 \leq\varepsilon_{\rm ph}.
\end{equation}
The scale identity then gives, on every such finite interval,
\begin{equation}\label{eq:finite-interval-physical-width}
 \int_{\tau_0}^{\tau_1}\lambda(\tau)\,d\tau
 \leq C_\lambda\lambda(\tau_0).
\end{equation}
Define $\Theta_\tau$ by
\begin{equation}\label{eq:adaptive-Theta-tau}
 \partial_\tau\Theta_\tau
 =\bigl((1+a)\bar\nabla\bar f-U_\tau\bigr)\circ\Theta_\tau
\end{equation}
and put
\begin{equation}\label{eq:adaptive-target-S}
 S_\tau=\lambda(\tau)\Theta_\tau^*\bar g.
\end{equation}
At the already fixed derivative order \(m\), normalize the initial
position by requiring that
\begin{equation}\label{eq:adaptive-initial-position}
 R_{\tau_0}:=\varphi_{-\tau_0}\circ\Theta_{\tau_0}
\end{equation}
is a global proper diffeomorphism (equivalently,
\(\Theta_{\tau_0}\) is), and that it and its inverse have
\(C^{m+1}\) norm at most
\(\Lambda_{\rm ad}\) on \(\Omega_\eta^+\), and
that on the explicit outer region \(\{\bar f\geq\Gamma/2\}\)
\begin{equation}\label{eq:adaptive-initial-radial-comparison}
 c_{\rm rad}(1+\bar f(x))
 \leq1+\bar f(R_{\tau_0}(x))
 \leq C_{\rm rad}(1+\bar f(x)).
\end{equation}
We also require uniform scale-normalized asymptotic-identity bounds:
on every dyadic tracking annulus
$A_L^\circ=\{L<\bar f<4L\}$, the maps $R_{\tau_0}$ and
$R_{\tau_0}^{-1}$ have \(C^{m+1}\) norm at most
\(\Lambda_{\rm ad}\) after the domain and range metrics are rescaled by
$L^{-1}$.  This condition is independent of $L\geq\Gamma$.
For the exact-core prepared center used in the formation construction
one has $R_{\tau_0}=\operatorname{Id}$.

If \(\tau_0\geq\tau_\Gamma\) and
\eqref{eq:adaptive-initial-position}--%
\eqref{eq:adaptive-initial-radial-comparison} hold, together with the
uniform scale-normalized dyadic bounds stated immediately after
\eqref{eq:adaptive-initial-radial-comparison}, then on
$[\tau_0,\tau_1)$:
\begin{enumerate}
\item for every $x\in\Omega_\eta^+$ and
      $\tau_0\leq\tau<\tau_1$,
      \begin{equation}\label{eq:Theta-f-tracking}
       C_m^{-1}e^{\tau-\tau_0}\bar f(\Theta_{\tau_0}(x))
       \leq \bar f(\Theta_\tau(x))
       \leq
       C_me^{\tau-\tau_0}\bar f(\Theta_{\tau_0}(x));
      \end{equation}
\item the relative maps
      \[
       R_\tau=\varphi_{-\tau}\circ\Theta_\tau
      \]
      and their inverses have uniformly bounded $C^{m+1}$ norms on
      $\Omega_\eta^+$.  For $\tau_0\leq\tau\leq s<\tau_1$,
      \begin{equation}\label{eq:relative-map-Cauchy}
       \norm{R_s-R_\tau}_{C^m(\Omega_\eta^+)}
       \leq C_m\left(
        \int_\tau^s(|a|+|b|)\,dr+e^{-\tau}\right);
      \end{equation}
      the same maps and their inverses retain uniform
      scale-normalized $C^{m+1}$ bounds on the whole AC outer region.
      In particular, for $\bar f(x)\geq\Gamma/2$,
      \begin{equation}\label{eq:R-global-radial-comparison}
       C_m^{-1}(1+\bar f(x))
       \leq1+\bar f(R_\tau(x))
       \leq C_m(1+\bar f(x));
      \end{equation}
\item on $\Omega_\eta^+$ the metrics $S_\tau$ are uniformly
      equivalent and have uniformly bounded geometry.  More
      precisely,
      \begin{equation}\label{eq:target-relative-derivatives}
       \norm{\partial_\tau S_\tau}_{C^m(S_\tau;\Omega_\eta^+)}
       \leq C_m\bigl(e^{-\tau}
                    +e^{-\tau}|a(\tau)|+|b(\tau)|\bigr).
      \end{equation}
\end{enumerate}
There is also a finite-order upgrade which uses no additional
geometric smallness.  Fix an integer \(m_*\geq1\), and choose
\(\varepsilon_{\rm ph}\) and \(\tau_\Gamma\) from the preceding
statement at order \(m_*\).  Under the same hypotheses at that order,
let \(q\geq m_*\), and suppose in addition that
\(R_{\tau_0}^{\pm1}\) have ordinary \(C^{q+1}\) bounds on
\(\Omega_\eta^+\) and uniform scale-normalized \(C^{q+1}\) bounds on
every dyadic tracking annulus, with some finite ceiling
\(\Lambda_q\).  Then assertions \textup{(2)} and \textup{(3)} hold
through order \(q\), with a constant
\[
 C_q=C_q\bigl(q,\Gamma,\Lambda_q,
          \eqref{eq:adaptive-position-package},\text{background}\bigr)
 <\infty,
\]
but with the same \(\varepsilon_{\rm ph}\) and \(\tau_\Gamma\).
In particular, no further reduction of the phase budget and no
increase of the geometric entrance time is required to propagate a
fixed finite higher order.  If the estimates hold globally, the
corresponding \(C^q\) convergence conclusions hold as well.  No
uniformity as \(q\to\infty\) is asserted.
If the estimates hold for every finite endpoint with constants
independent of $\tau_1$, then $R_\tau$ and $S_\tau$ converge in
$C^m(\Omega_\eta^+)$ as $\tau\to\infty$.
If $\Gamma$ is chosen beyond the support of the conjugated cutoff,
then the term $|b|$ in \eqref{eq:target-relative-derivatives} vanishes
on $\Omega_\eta^+$.
\end{proposition}

\begin{proof}
Along an orbit of \eqref{eq:adaptive-Theta-tau}, write
$y(\tau)=\bar f(\Theta_\tau(x))$.  The shrinker identities and the
growth estimates for the fields $W_j$ give, whenever $y$ is large,
\[
 y'=(1+a)|\bar\nabla\bar f|^2-U_\tau(\bar f),\qquad
 \left|\frac{y'}y-1\right|
 \leq C\bigl(y^{-1}+|a|+|b|\bigr).
\]
A preliminary differential inequality gives
$y(\tau)\geq c e^{\tau-\tau_0}y(\tau_0)$.  Hence
\[
 \int_{\tau_0}^{\tau_1} y(s)^{-1}\,ds
 \leq C y(\tau_0)^{-1},
\]
with the same bound if $\tau_1$ is replaced by any smaller endpoint.
Integration of the logarithmic inequality proves
\eqref{eq:Theta-f-tracking}.

Differentiating $R_\tau=\varphi_{-\tau}\circ\Theta_\tau$ gives the
exact equation
\begin{equation}\label{eq:relative-target-flow}
 \partial_\tau R_\tau
 =
 \left[(\varphi_{-\tau})_*
       \bigl(a\bar\nabla\bar f-U_\tau\bigr)\right]\circ R_\tau.
\end{equation}
On a fixed enlargement of $\Omega_\eta^+$, the fields in brackets
obey
\[
 \norm{(\varphi_{-\tau})_*
       (a\bar\nabla\bar f-U_\tau)}_{C^{m+1}}
  \leq C_m(|a|+|b|).
\]
Indeed, $\bar\nabla\bar f$ is invariant under its own flow, while
$e^{-\tau}\bar f\circ\varphi_\tau$, the conjugated cutoffs, and the
conjugated fields $W_j$ have bounded derivatives on every fixed
annulus.  The asymptotically conical error is $O(e^{-\tau})$ and is
harmless.  The standard variational equations for a flow, followed
by Gronwall, prove the uniform derivative bounds and
\eqref{eq:relative-map-Cauchy}.  If the estimates hold for arbitrary
finite endpoints, the resulting $L^1$ tail makes the maps Cauchy in
$C^m$.
The inequality
$|(\varphi_{-\tau})_*U_\tau(\bar f)|\leq C|b|(1+\bar f)$
and its scale analogue preserve
\eqref{eq:adaptive-initial-radial-comparison} on the whole outer
region.  Applying the differentiated flow equations on each dyadic
annulus after rescaling its metric by the annular $\bar f$-level
preserves the uniform scale-normalized derivative bounds, with a
constant independent of the annulus.  The same at-most-linear growth
prevents finite-time escape, so $\Theta_\tau$ remains a complete
diffeomorphism.  This is the outer-region form of the same argument.

Finally,
\begin{equation}\label{eq:S-tau-exact}
 \partial_\tau S_\tau
 =\lambda\Theta_\tau^*
   \bigl(-2(1+a)\Ric_{\bar g}-\Lie_{U_\tau}\bar g\bigr).
\end{equation}
On the tracked annulus,
\[
 |\bar\nabla^\ell\Ric_{\bar g}|
 \leq C_\ell\bar f^{-1-\ell/2},\qquad
 |\bar\nabla^\ell\Lie_{W_j}\bar g|
 \leq C_\ell\bar f^{-\ell/2}.
\]
Tensor scaling in \eqref{eq:S-tau-exact}, together with
\eqref{eq:Theta-f-tracking}, yields
\eqref{eq:target-relative-derivatives}.  Its right-hand side is
integrable on each bootstrap interval, proving uniform equivalence.
After global continuation its tail is integrable, which proves
convergence of the target metrics.  If the physical graft annulus
lies beyond the conjugated
support of $\chi_\tau$, then $U_\tau$ vanishes there.

We prove the finite-order upgrade separately, since this is where the
distinction between geometric smallness and higher-order boundedness is
used.  Put
\[
 \mathcal V_\tau
 :=
 (\varphi_{-\tau})_*
 \bigl(a\bar\nabla\bar f-U_\tau\bigr),
 \qquad
 \partial_\tau R_\tau=\mathcal V_\tau\circ R_\tau .
\]
The order-\(m_*\) conclusion already keeps \(R_\tau^{\pm1}\) in the
common compact-collar and scale-normalized dyadic coordinate domains
and preserves the two-sided radial comparison.  On each of those
domains, the all-order AC symbol estimates and the fixed cutoff profile
give, for every fixed \(q\),
\[
 \|\mathcal V_\tau\|_{C^{q+1}_{\rm sc}}
 \leq C_q\bigl(|a(\tau)|+|b(\tau)|+e^{-\tau}\bigr),
\]
where the subscript ``\({\rm sc}\)'' denotes either the fixed collar
norm or the corresponding scale-normalized dyadic norm.  Hence
\[
 \int_{\tau_0}^{\tau_1}
 \|\mathcal V_\tau\|_{C^{q+1}_{\rm sc}}\,d\tau
 \leq C_q\bigl(\varepsilon_{\rm ph}+e^{-\tau_0}\bigr),
\]
uniformly in the finite endpoint.

Differentiating
\(\partial_\tau R_\tau=\mathcal V_\tau\circ R_\tau\) spatially
\(j\leq q+1\) times gives a triangular variational system: its
top-order term is
\(D\mathcal V_\tau(R_\tau)D^jR_\tau\), while every remaining term
contains only lower derivatives of \(R_\tau\) and derivatives of
\(\mathcal V_\tau\) of order at most \(j\).  Induction in \(j\) and
Gronwall therefore propagates the stated \(C^{q+1}\) bounds.  Applying
the same argument to the inverse-flow equation gives the corresponding
bounds for \(R_\tau^{-1}\), and integration between \(\tau\) and \(s\)
gives the order-\(q\) version of
\eqref{eq:relative-map-Cauchy}.  Finally,
\[
 S_\tau
 =\lambda(\tau)R_\tau^*\varphi_\tau^*\bar g
\]
and the fixed scale bracket give the spatial \(C^q\) coefficient bounds
for \(S_\tau\); the exact identity \eqref{eq:S-tau-exact} gives the
order-\(q\) version of
\eqref{eq:target-relative-derivatives}.  Only the already established
low-order range and radial controls use the small phase threshold.
Higher derivatives require bounded initial jets, but no additional
smallness.

The constants in the order-\(m\) statement depend only on the fixed
background, \(m,\Gamma\), and the numerical package
\eqref{eq:adaptive-position-package}.
In the finite-order upgrade, the additional dependence on the
quantified initial ceiling \(\Lambda_q\) is exactly the one displayed
in the statement.  For the low orders used before the radius choice,
the same calculation
is instead performed in the fixed core-plus-pre-radius atlas
\eqref{eq:pre-radius-dyadic-family}; Lemma~\ref{lem:pre-radius-low-order-closure}
records the resulting \(\Gamma\)-independent majorants.  The displayed
\((m,\Gamma)\)-dependence belongs only to the later fixed-package
high-regularity output.
\end{proof}

\begin{remark}[Fixed adaptive parameter order and phase budget]
\label{conv:authoritative-adaptive-order}
Fix \(0<\alpha<1\).  For every continuation, restart, and two-state
argument below, fix once and for all
\[
 m_{\rm ad}=13.
\]
The dependency order is the following.  First fix
\(0<\sigma<\theta<\beta\).  Before choosing the radius, fix the
primitive pre-radius ellipticity, coefficient, \(R^{\pm1}\),
\(F^{\pm1}\), graft, upper-scale, and positive-margin data in
\(\mathfrak P_{\rm pre}^{\rm prim}\), including
\[
 K_{\rm gr},\quad C_{\rm sc}
\]
with positive room.  Apply
Lemma~\ref{lem:pre-radius-low-order-closure} once.  It produces
\(\Gamma_{\rm pre}\), the annulus-tracking and robust forcing
constants, and all effective-column and feedback ceilings in
\(\mathfrak P_{\rm pre}^{\rm der}\).  In particular, the product in
\eqref{eq:pre-radius-robust-graft-ceiling} supplies the robust forcing
constant in \eqref{eq:adaptive-robust-forcing-constants}.  The fixed
cutoff profile, the AC estimates, and the direct scale-normalized proof
of that lemma make precisely these derived low-order ceilings uniform
for \(\Gamma\geq\Gamma_{\rm pre}\); they do not use the later
\(\Gamma\)-dependent high-regularity atlas.  Fix
\(\Lambda_{\rm ell}^{(2)}\) and
\(K_{\mathcal Y,0}^{(2)}\) strictly above the corresponding common
 faces in this derived package, and obtain the named
 preliminary admissible triple
 \((\delta_{\rm K}^{(2)},C_{\rm K}^{(2)},c_{\rm K}^{(2)})\)
 from Lemma~\ref{lem:uniform-two-state-Kato-ledger}.  Apply
 Lemma~\ref{lem:uniform-two-state-package-radius} to reduce the
 preliminary witness by the two radial-absorption thresholds and record
 the resulting \(\delta_{\rm box}^{(2)}\).  Now evaluate
\(\overline\Gamma_{\rm 3reg}\) from the already fixed ceilings and
\(\overline\Gamma_{\rm 2st}\) from
\eqref{eq:uniform-two-state-package-radius}, with
\(\theta_*=\theta\), and choose \(\Gamma\) satisfying
\eqref{eq:global-compatible-package-radius}.
Only at this point construct \(\eta_\Gamma\), the graft collars, and
the \(\Gamma\)-dependent atlases.  Record the fixed positive lower and
upper constants in the resulting two-state graft-support annulus as
\[
 0<c_{\rm supp}^{(2)}\leq C_{\rm supp}^{(2)}<\infty.
\]
These are post-radius geometric constants and do not enter either
radius functional.  Then choose the positive lower scale
margin \(0<c_{\rm sc}<C_{\rm sc}\) and set
\[
 c_{\rm scl}=\frac12c_{\rm sc},\qquad
 C_{\rm scl}=2C_{\rm sc}.
\]
The lower margin \(c_{\rm sc}\) does not enter the upper forcing,
column, or feedback ceilings and hence does not enter
\(\overline\Gamma_{\rm 3reg}\) or
\(\overline\Gamma_{\rm 2st}\).  In the exact-core construction it may
therefore be chosen after the implantation scale \(A\), with
\(c_{\rm sc}<A\).  Freeze the static package at this point.  Thus the
enlarged-bracket constants are package entries, while the forcing
ceiling used to choose \(\Gamma\) is not chosen from a
\(\Gamma\)-dependent estimate.
At this same post-radius stage freeze the three named harmonic-map
thresholds
\(\varepsilon_{\rm map}^{\rm HM}\),
\(\varepsilon_{\rm hm}^{\rm HM}\), and
\(\varepsilon_{\rm ph}^{\rm HM}\) from
\eqref{eq:named-HMHF-smallness-thresholds}.  Their dependencies are
the already fixed geometric package and collar separations; none
 enters either radius functional.
At this post-radius stage also choose the domains \(K_i^\Gamma\), the
core and chart thresholds
\(\eta_{\rm core}^{(2)},\eta_{\rm ch}^{(2)}\), and the resulting
\(\varepsilon_{\rm ph,*}^{(2)}\) in
\eqref{eq:uniform-two-state-phase-threshold}.  These are frozen before
any subordinate pair.
With \(\theta_*=\theta\), take the rate-independent two-state
smallness constants from
\eqref{eq:uniform-two-state-smallness-package}.  Choose one common
 coarse threshold \(\delta_{\rm c2}>0\).  Let
\(C_{\rm raw}\geq1\) be the fixed comparison from the pre-radius
activation norm to the raw harmonic-map discrepancy norm, characterized
by
\begin{equation}\label{eq:pre-atlas-to-raw-C2-comparison}
 \sup_M\sum_{j=0}^{2}(1+\bar f)^{j/2}
 |\bar\nabla^jh|_{\bar g}
 \leq C_{\rm raw}
 \|h\|_{\mathfrak C_{{\rm pre},0}^{2,\alpha}} .
\end{equation}
It depends only on the fixed background, \(\alpha\), and the fixed
pre-radius atlas and norm conventions; it is independent of
\(\Gamma\), the rate pair, the evolution, and the finite endpoint.
Require
\begin{equation}\label{eq:authoritative-coarse-C2-threshold}
 2\delta_{\rm c2}
 \leq
 \min\{\delta_{\rm atl}^{\rm pre},\delta_{\rm pre},
       C_{\rm raw}^{-1}\varepsilon_{\rm hm}^{\rm HM},
       \varepsilon_{\rm K},\delta_{\rm rec},
       \delta_{\rm box}^{(2)}\}
\end{equation}
and below every remaining fixed structural \(C^2\) threshold in the
physical, target, and harmonic-map packages.  Thus the doubled
bootstrap face is already an activation face for every pre-radius and
two-state estimate; no later argument decreases it.
For orders \(0\leq m\leq5\), use the
\(\Gamma\)-uniform comparison constant
\(C_{\rm ad}^{\rm pre}(m)\) constructed in
Lemma~\ref{lem:pre-radius-low-order-closure}; these are the only map
constants used in the pre-radius column and feedback ceilings.  For the
later high-regularity continuation, let
\(C_{\rm ad}=C_{\rm ad}(m_{\rm ad},\Gamma)\) be the fixed comparison
constant from the prepared scale-normalized map atlas to the
compact-collar and dyadic coordinate norms used in
\eqref{eq:adaptive-initial-position} and
\eqref{eq:adaptive-initial-radial-comparison}, and set
\begin{equation}\label{eq:authoritative-adaptive-map-bound}
 \Lambda_{\rm ad}
 :=\max\{1,C_{\rm ad}\Lambda_{\rm map}\}.
\end{equation}
Thus the adaptive-position package
\eqref{eq:adaptive-position-package}, through every order
\(0\leq m\leq m_{\rm ad}\), is determined by
\(\mathfrak P_{\rm prep}\) and the fixed background atlases; no
additional map bound is chosen later.
The doubled order-twelve source--target coefficient faces and the
order-six map-and-inverse faces, inserted in the prepared graph formula,
also determine one endpoint-independent post-radius ceiling: for every
evolution carrying those displayed doubled faces,
\begin{equation}\label{eq:pre-atlas-C5-interpolation-ceiling}
 \sup_\tau
 \|h(\tau)\|_{\mathfrak C_{{\rm pre},0}^{5}}
 \leq K_{h,5}^{\rm pre}<\infty .
\end{equation}
This ceiling is a completed package entry but does not enter either
radius functional.  Let
\[
 \vartheta_{\rm int}:=\frac{2+\alpha}{5}\in(0,1)
\]
and let \(C_{\rm int}\) be the uniform scale-one interpolation constant
in the fixed pre-radius atlas.  Let
\(C_{\rm 3reg}^{(0)}\) and \(C_{\rm 3reg}^{(2)}\) denote,
respectively, the fixed global zeroth-order and through-order-two
output constants in the three-region theorem for this package.  Every
later entrance threshold is also selected against these constants.
Next let
\(C_{6,K_{\rm gr}}\) denote the constant
\(C_{m,K_{\rm gr}}\) in
\eqref{eq:graft-compatibility-propagated} at \(m=6\), for these fixed
geometric data.  Only after those choices define
\(\varepsilon_{\rm ph}\) to be the minimum of the thresholds in
Proposition~\ref{prop:adaptive-target-tracking} for
\(0\leq m\leq m_{\rm ad}\), every later strict-improvement threshold
which uses the phase budget for this fixed rate pair, and the number
required to ensure
\[
 C_{6,K_{\rm gr}}\varepsilon_{\rm ph}
 <\frac18K_{\rm gr}.
\]
Include also the requirement
\begin{equation}\label{eq:authoritative-scale-budget}
 \varepsilon_{\rm ph}
 \leq
 \min\{\log2,\varepsilon_{\rm pre},
               \varepsilon_{\rm ph}^{\rm HM},
               \varepsilon_{\rm ph,*}^{(2)}\}.
\end{equation}
In particular,
\[
 \varepsilon_{\rm ph}
 \leq\varepsilon_{\rm ph,*}^{(2)}
 \leq
 \min\left\{\frac{\log2}{K_{\rm J}^{(2)}},
             \eta_{\rm core}^{(2)},\eta_{\rm ch}^{(2)}\right\}.
\]
Indeed, the exact scale equation and the strict entrance bracket give
\begin{equation}\label{eq:authoritative-scale-bracket-closure}
 e^\tau\lambda(\tau)
 =
 e^{\tau_0}\lambda(\tau_0)
 \exp\!\left(-\int_{\tau_0}^{\tau}a(q)\,dq\right).
\end{equation}
Consequently, on every interval on which the exact scale equation
holds and the accumulated phase budget is at most
\(\varepsilon_{\rm ph}\), the entrance
inequality
\[
 c_{\rm sc}e^{-\tau_0}
 <\lambda(\tau_0)
 <C_{\rm sc}e^{-\tau_0}
\]
and \eqref{eq:authoritative-scale-budget} imply the strict propagated
bracket
\[
 c_{\rm scl}e^{-\tau}
 <\lambda(\tau)
 <C_{\rm scl}e^{-\tau}.
\]
Thus the scale comparison in
\eqref{eq:adaptive-phase-budget} is a closed bootstrap face determined
by the prepared entrance, not an additional unsupplied hypothesis.
Choose \(\varepsilon_{\rm ent}>0\) once, no larger than every fixed
one-state feedback, three-region, graft, harmonic-map, continuation,
two-state common-box entrance threshold, and the polarized two-state
energy-margin threshold determined by \(\beta-\theta\) for this
package; in particular require
\[
 \varepsilon_{\rm ent}
 \leq\min\{\varepsilon_{\rm map}^{\rm HM},
            \varepsilon_{\rm hm}^{\rm HM}\}.
\]
In addition impose the strict interpolation and raw-face choices
\begin{equation}\label{eq:authoritative-pre-Holder-interpolation-choice}
 \begin{gathered}
  C_{\rm 3reg}^{(2)}\varepsilon_{\rm ent}
  <\delta_{\rm c2},\\
  C_{\rm int}
  (C_{\rm 3reg}^{(0)}\varepsilon_{\rm ent})^{1-\vartheta_{\rm int}}
  (K_{h,5}^{\rm pre})^{\vartheta_{\rm int}}
  <\delta_{\rm c2}.
 \end{gathered}
\end{equation}
The
 named future-tail constants evaluated at the frozen top rate,
\((C_P^*,c_P^*)=(C_P,c_P)|_{\theta=\theta_*}\), the fixed positive
lower graft-support constant
\(c_{\rm supp}^{(2)}\), the now fixed
\(\varepsilon_{\rm ph}\), and this
\(\varepsilon_{\rm ent}\) determine the rate-independent two-state
base time
\[
 \tau_{\rm base}^{(2)}
 =\tau_{\rm base}^{(2)}
  (\theta,\Gamma,\Lambda_{\rm ell}^{(2)},
   K_{\mathcal Y,0}^{(2)},\varepsilon_{\rm ph},
    \varepsilon_{\rm ent},C_P^*,c_P^*,
   c_{\rm supp}^{(2)},\text{background})
\]
from \eqref{eq:uniform-two-state-base-time}.  Thus this time is fixed
before the common entrance ball and before any subordinate pair is
quantified.
Then put
\[
 \tau_{\rm gr}
 :=1+\max\!\left\{
  0,\log\!\left(
   \frac{8C_{6,K_{\rm gr}}C_{\rm sc}}
        {K_{\rm gr}\Gamma}\right)\right\}.
\]
Define
\begin{equation}\label{eq:common-adaptive-entrance-time}
 \tau_{\rm ad}
 :=\max\left\{
   \tau_{\rm pre},
   \tau_{\rm base}^{(2)},
   \max_{0\leq m\leq m_{\rm ad}}\tau_\Gamma(m),
   \tau_{\rm gr}\right\}.
\end{equation}
Finally define the single constant
\(C_\lambda=e^{\varepsilon_{\rm ph}}\) by
\eqref{eq:fixed-C-lambda}.  Every later invocation of that proposition
uses these common choices and requires
\(\tau_0\geq\tau_{\rm ad}\).  Its uniform scale-one \(C^{14}\) bounds for
\(R^{\pm1}\) imply all lower \(C^{j,\alpha}\) bounds required below,
including the order-\(12,\alpha\) target-coefficient package, by the
fixed scale-one interpolation inequalities.  For every later
continuation output order \(k_0\geq12\), the prepared input order
\(k_0+2\geq14\) has map components with one additional derivative and
supplies the corresponding finite higher-order initial data.  The
finite-order upgrade in
Proposition~\ref{prop:adaptive-target-tracking} propagates precisely
that fixed order without changing \(m_{\rm ad}=13\),
\(\varepsilon_{\rm ph}\), or \(\tau_{\rm ad}\).  Thus order thirteen
remains the fixed order for every geometric, graft, cutoff,
target-defect, and restart threshold; higher prepared orders enter only
as boundedness data for the requested continuation output order.
\end{remark}

\subsection{Factoring the controlled harmonic-map heat flow}

The following identity is the main geometric simplification.

\begin{proposition}[Relative harmonic-map heat flow]
\label{prop:relative-HMHF}
Let $\acute G(t)$, $\lambda$, $\Theta$, and $\Phi$ be as in
Propositions~\ref{prop:controlled-HMH} and
\ref{prop:adaptive-graft}, with
\[
 \frac{d\tau}{dt}=\lambda^{-1}.
\]
Set
\begin{equation}\label{eq:relative-map-F}
 F_t=\Theta_t^{-1}\circ\Phi_t.
\end{equation}
Then $F$ solves the ordinary time-dependent harmonic-map heat flow
\begin{equation}\label{eq:relative-HMHF}
 \partial_tF=\Delta_{\acute G(t),S(t)}F,
 \qquad S(t)=\lambda(t)\Theta_t^*\bar g.
\end{equation}
 Moreover, if
\[
 q=(F^{-1})^*\acute G,\qquad
 g=\lambda^{-1}(\Phi^{-1})^*\acute G=\bar g+h,
\]
then
\begin{equation}\label{eq:relative-metric-identity}
 q-S=\lambda\Theta^*h.
\end{equation}
Wherever \(q\) and \(S\) are \(\Lambda_0\)-uniformly equivalent,
the Eulerian velocity $V_F=(\partial_tF)\circ F^{-1}$ satisfies
\begin{equation}\label{eq:F-speed}
 |V_F|_S
 \leq C(\Lambda_0)\lambda^{-1/2}
       (|\bar\nabla h|_{\bar g}\circ\Theta),
 \qquad
 |\lambda V_F|_S
 \leq C(\Lambda_0)\lambda^{1/2}
       (|\bar\nabla h|_{\bar g}\circ\Theta).
\end{equation}
\end{proposition}

\begin{proof}
Write
\[
 V=\lambda^{-1}\bigl((1+a)\bar\nabla\bar f-U\bigr).
\]
Thus $\partial_t\Theta=V\circ\Theta$.  Since
$\Theta:(M,\Theta^*\bar g)\to(M,\bar g)$ is an isometry,
naturality of the tension field gives
\[
 \Delta_{\acute G,\bar g}(\Theta\circ F)
 =d\Theta\bigl(\Delta_{\acute G,\Theta^*\bar g}F\bigr).
\]
The target connections of $\Theta^*\bar g$ and
$S=\lambda\Theta^*\bar g$ agree.  Substituting
$\Phi=\Theta\circ F$ into \eqref{eq:controlled-HMH}, the two terms
$V\circ\Theta\circ F$ cancel, proving \eqref{eq:relative-HMHF}.

The relation $\Phi=\Theta\circ F$ also gives
\[
 (F^{-1})^*\acute G
 =\lambda\Theta^*(\bar g+h),
\]
which is \eqref{eq:relative-metric-identity}.  In target coordinates,
the harmonic-map velocity is the negative DeTurck vector:
\[
 V_F=-B_S(q).
\]
For $q$ uniformly equivalent to $S$,
$|B_S(q)|_S\leq C|\nabla^S(q-S)|_S$.  The scaling in
\eqref{eq:relative-metric-identity} yields
\[
 |\nabla^S(q-S)|_S
 =\lambda^{-1/2}(|\bar\nabla h|_{\bar g}\circ\Theta),
\]
and \eqref{eq:F-speed} follows.
\end{proof}

\paragraph{Local solution class.}
In a prepared chart frozen at the left endpoint \(s\), write
\[
 \widehat t_s(\tau):=\frac{t(\tau)-t(s)}{\lambda(s)}.
\]
The notation \(\mathbb S^{r,\alpha}(J)\) means the class in which the
zero-trace affine increment of the scale-normalized compact
Ricci--DeTurck representative has little anisotropic regularity
\(h^{(r+\alpha)/2,r+\alpha}\) in the frozen scale-adapted compact
atlas, and the scalar pair
\((\widehat t_s,\ell)\) is \(C^{1,\alpha/2}\).  The same-output
\(\Theta\)-coordinate is \(C^{1,\alpha/2}\) with values in the weighted
spatial class of order \(r+1\), and the same-output \(\Phi\)-coordinate
has, uniformly in the frozen source-adapted charts,
\[
 L^\infty C_x^{r+1,\alpha}\cap
 W^{1,\infty}C_x^{r-1,\alpha},
\]
together with its prescribed order-\((r+1)\) left trace.  The dependent
graph tensor is
\[
 C^{\alpha/2}_\tau\mathfrak T_{{\rm sc},N}^{r,\alpha}
 \cap C^{1,\alpha/2}_\tau
        \mathfrak T_{{\rm sc},N}^{r-2,\alpha}.
\]
All coordinate differences have the stated left trace, every endpoint
has the canonical two-order-lower trace, and the compact and noncompact
equations hold in their Bochner mild senses.
Subsection~\ref{subsec:time-dependent-prepared-spaces} gives the
corresponding frozen clocks, atlas-supremum norms, traces, and mild
formulations.

\begin{definition}[Admissible adaptive and first-exit intervals]
\label{def:admissible-first-exit-interval}
For \(r\geq12\) and \(0<\alpha<1\), an
\emph{admissible prepared evolution of order \((r,\alpha)\)} on a
normalized-time interval \(I\) is a tuple
\[
 (G(t),t(\tau),\lambda(\tau),\Theta_\tau,F_\tau,
   c(\tau)=(a(\tau),b(\tau)))
\]
which has the following local-in-time meaning.  Every compact
subinterval \(K\subset I\), with compactness understood in the relative
topology of \(I\), has a finite cover by intervals
\[
 J=[s_J,s_J+\delta_J]\subset I,\qquad \delta_J>0,
\]
 such that, in the trial-independent prepared chart frozen at the left
 endpoint \(s_J\), the restricted tuple is represented by an element of
 \(\mathbb S^{r,\alpha}(J)\) in the convention just fixed, with its
 precise norm as constructed in
Subsection~\ref{subsec:time-dependent-prepared-spaces}.  Once the
interval has been fixed, we suppress it from the notation.  Thus the
compact Ricci--DeTurck representative, scalar variables, prepared maps,
and graph tensor have exactly the spatial and temporal regularity
specified in that solution space; no simultaneous all-order bound is
part of admissibility.  In particular, the nine moments are locally
absolutely continuous, and
every ordinary bootstrap face defined by the prepared Banach variables
is continuous.  A normalized or physical harmonic-radius condition is
not an ordinary continuous face: whenever one is declared, it is
interpreted through its recorded buffered operative and reserve
witnesses and the one-sided lower-stability and reserve-to-operative
clauses of Lemmas~\ref{lem:prepared-harmonic-radius-lower-stability}
and~\ref{lem:finite-physical-harmonic-openness}.  No continuity of a
bare supremal-radius functional is part of admissibility.  The ODE and
parabolic equations below hold in the strong, mild, or
absolutely-continuous sense encoded by \(\mathbb S^{r,\alpha}(J)\),
classically for smooth data.  On a finite interval on which this
solution norm stays bounded, all components have endpoint traces in the
two-order-lower prepared topology.

The tuple satisfies the closed Ricci flow, \(t_\tau=\lambda\), the exact
scale equation \(\lambda_\tau=-(1+a)\lambda\), the adaptive
\(\Theta\)-equation \eqref{eq:adaptive-Theta-tau}, and the relative
harmonic-map equation \eqref{eq:relative-HMHF}.  It also satisfies,
identically,
\[
 S=\lambda\Theta^*\bar g,\qquad
 \acute G=\eta\,\iota_*G+(1-\eta)S,\qquad
 \Phi=\Theta\circ F,\qquad
 h=\lambda^{-1}(\Phi^{-1})^*\acute G-\bar g,
\]
the adaptive normalized equation \eqref{eq:adaptive-normalized}, the
exact nine receding-slice identities, and the exact adaptive Gram
feedback: at each time, \(c\) is the unique coefficient vector obtained
by differentiating those nine identities in
\eqref{eq:adaptive-normalized} on the stated invertible-Gram locus.
All prepared maps and inverses are
degree-one proper diffeomorphisms, the structural graph identities
hold, and the tuple remains in the stated common-margin prepared
regularity class.

Given a finite, explicitly declared list of ordinary continuous
bootstrap faces, together with any explicitly declared witnessed
harmonic-radius conditions, an \emph{admissible first-exit interval} is
the maximal half-open interval \([\tau_0,\tau_1)\) on which an
admissible prepared evolution exists, the ordinary faces retain their
displayed weak bounds, and the recorded harmonic witnesses retain their
operative bounds and positive reserve-to-operative gaps.  Admissibility
comprises only the evolution, slice, feedback, graph, and regularity
identities listed above.  Quantitative coefficient, margin, smallness,
buffer, and witnessed-harmonic hypotheses are imposed separately.

Every strict quantitative requirement among the ordinary faces is
understood as a continuous positive margin: an upper face is its
ceiling minus the controlled quantity, a lower face is the controlled
quantity minus its floor, and an invertibility face is expressed by
the corresponding least-singular-value or inverse-bound margin.  An
ordinary face \(\mathfrak m_j>0\) occurs at a finite right endpoint if
\[
 \liminf_{\tau\uparrow\tau_1}\mathfrak m_j(\tau)=0.
\]
A witnessed harmonic-radius condition occurs at a finite right endpoint
if no terminal subinterval carries one common buffered witness package
with positive componentwise reserve-to-operative gaps.  Since the
ordinary face list and the witnessed harmonic list are finite, if
neither kind of condition occurs, then there are
\(\sigma<\tau_1\) and \(\mu_*>0\) such that every ordinary declared
margin is at least \(\mu_*\) on \([\sigma,\tau_1)\), and every declared
harmonic condition carries one common operative-and-reserve witness
package there with all componentwise gaps at least \(\mu_*\), after
decreasing \(\mu_*\) if necessary.

An unqualified admissible prepared evolution means one of a fixed
finite order sufficient for the separately displayed hypotheses of
the result in which it occurs.  For continuation order \(k_0\), the
constructed evolution has order \((k_0+2,\alpha)\).
\end{definition}

\begin{remark}
The conjugation in Proposition~\ref{prop:relative-HMHF} is the
receding, noncompact analogue of the modified harmonic-map gauge used
by Choi--Lai~\cite{ChoiLai}.  Their long-time theorem assumes a compact
integrable shrinker and an already asymptotic Ricci flow.  It therefore
does not supply Theorem~\ref{thm:adaptive-HMHF-continuation}, but its
modified-HMHF calculation is consistent with
\eqref{eq:relative-HMHF}.
\end{remark}

\begin{remark}[The precise Bamler--Kleiner input]
\label{rem:BK-A9-specialization}
The result used below is the published general-form
\cite[Proposition~A.9]{BamlerKleiner}.  After translating the initial
time to zero and parabolically rescaling, it applies to smooth metric
families on arbitrary, possibly noncompact, manifolds whose initial
metrics are complete.  Its hypotheses are: uniform curvature bounds;
the stated derivative bounds through order ten for curvature and the
metric time derivatives (with the permitted \(t^{-j/2}\) initial-face
behavior); two-sided tensor bounds for
\(\partial_tg+2\Ric_g\) and \(\partial_ts+2\Ric_s\); and a
diffeomorphic initial map with small \(C^0\) metric discrepancy.  Its
conclusion is a short-time harmonic-map heat flow consisting of
diffeomorphisms, with a prescribed larger strict \(C^0\) discrepancy
bound and the corresponding restart alternative.

The hypotheses below are stronger in exactly the needed directions.
The common harmonic-radius and coefficient package gives the curvature
and metric-time derivative bounds after scale normalization; the
displayed curvature and defect derivative bounds are uniform and hence
imply the allowed initial-face bounds; properness strengthens the
initial diffeomorphism hypothesis; and
\eqref{eq:moving-target-HMHF-initial-closeness} implies the required
\(C^0\) discrepancy.  Thus neither compactness nor an exact Ricci-flow
identity is imported from the citation.  Displacement, uniqueness in
the controlled class, and the higher local derivative estimates used
here are established separately in the proof below.
Most importantly for iteration, the admissible discrepancy cap
\(\bar\eta_n\) in the cited proposition is dimensional.  The finite
curvature, coefficient, metric-time-derivative, and Ricci-defect
ceilings determine its short lifespan \(\tau\), but do not decrease
\(\bar\eta_n\).  We retain this separation explicitly below.
\end{remark}

\begin{proposition}[Uniform moving-target HMHF restart]
\label{prop:moving-target-HMHF-restart}
Fix \(n\in\mathbb N\), \(0<\alpha<1\), and
\(m\in\{4,5\}\).  There is a dimensional discrepancy cap
\[
 \epsilon_{\rm BK}=\epsilon_{\rm BK}(n)>0 .
\]
For every finite bounded-geometry, coefficient, and defect ceiling
\(\Lambda\), and every
\(0<\epsilon<\epsilon'<\epsilon_{\rm BK}\), there is
\(\delta=\delta(n,m,\alpha,\Lambda,\epsilon,\epsilon')>0\) with the following
property.  Let \((X^n,g(t))\) and \((Y^n,s(t))\) be smooth complete
metric families on \([t_*,t_*+\delta r^2]\), smooth jointly in space
and time.  Suppose, at scale $r$, that both families
have a uniform \(C^{2,\alpha}\) harmonic-radius lower bound and
scale-normalized metric-coefficient bounds through order
\(m+7\leq12\).  Thus, in particular, the coefficient jets through
order \(m+6\leq11\) used in the map estimates below are uniformly
controlled.  Suppose also that the curvature derivatives through
order ten and the derivatives through order ten of
\[
 \mathfrak D_g=\partial_tg+2\Ric_g,\qquad
 \mathfrak D_s=\partial_ts+2\Ric_s
\]
are bounded by $\Lambda$ in scale-normalized norms.  These assumptions
imply the curvature and metric-time-derivative hypotheses of
\cite[Proposition~A.9]{BamlerKleiner}.  If $F_*:X\to Y$ is a proper
diffeomorphism and
\begin{equation}\label{eq:moving-target-HMHF-initial-closeness}
 \|(F_*^{-1})^*g(t_*)-s(t_*)\|_{C^1_r(s(t_*))}\leq\epsilon,
\end{equation}
where the first spatial derivative in $C^1_r$ is multiplied by $r$,
then
\[
 \partial_tF=\Delta_{g(t),s(t)}F,\qquad F(t_*)=F_*,
\]
has a solution on this interval.  It is unique in the following
precise class.  A competitor \(\widetilde F\) is a
\(C^{2,1}_{\rm loc}\) solution such that every
\(\widetilde F(t):X\to Y\) is a proper diffeomorphism and, for some
finite constant \(L_{\widetilde F}\) fixed for that solution,
\[
 \sup_{t_*\leq t\leq t_*+\delta r^2}
 \left(
  \|d\widetilde F(t)\|_{L^\infty(g(t),s(t))}
  +\|d\widetilde F(t)^{-1}\|_{L^\infty(s(t),g(t))}
 \right)
 \leq L_{\widetilde F}.
\]
The constant \(L_{\widetilde F}\) is not an input to the existence
time \(\delta\); it enters only the pair-dependent local argument used
to prove uniqueness.  Relative to fixed basepoints, the displayed
bound, short-time metric equivalence, and the bounded trajectory of
one basepoint imply at-most-linear growth of the map and inverse, so no
additional unquantified polynomial-growth constants are part of the
class.  Every
$F(t)$ is a proper diffeomorphism and
\begin{equation}\label{eq:moving-target-HMHF-closeness}
 \|(F(t)^{-1})^*g(t)-s(t)\|_{C^0(s(t))}<\epsilon'.
\end{equation}
The scale-normalized local derivatives of \(F\) through order
\(m+1\), including their \(C^\alpha\) seminorms, obey uniform estimates
on the second half of every restart interval.  The same is true for
\(F^{-1}\).  If the initial \(C^{m+1,\alpha}_r\) map and inverse-map
norms are bounded, these estimates hold up to the initial face as
well.
After decreasing \(\delta\), the solution also obeys the uniform
displacement estimate
\begin{equation}\label{eq:moving-target-uniform-displacement}
 \sup_{x\in X}
 d_{s(t_*)}\bigl(F(x,t),F_*(x)\bigr)
 \leq C(n,\Lambda,\epsilon')\sqrt{t-t_*},
 \qquad t_*\leq t\leq t_*+\delta r^2 .
\end{equation}
The existence time and estimates are uniform.  Here persistence of a
strict version of \eqref{eq:moving-target-HMHF-closeness} at a
candidate finite endpoint \(T^\dagger\) means that, for some
\(\epsilon''<\epsilon'\),
\[
 \limsup_{t\uparrow T^\dagger}
 \|(F(t)^{-1})^*g(t)-s(t)\|_{C^0(s(t))}
 \leq\epsilon''.
\]
Consequently, as long as the displayed geometry bounds extend past
\(T^\dagger\) and this strict terminal margin persists, loss of the
harmonic-map chart cannot be the first finite-time breakdown.
\end{proposition}

\begin{proof}
Scale \(r=1\) and push the source family forward by \(F_*\), so that
the initial map becomes the identity on \(Y\).
The curvature and defect assumptions imply the curvature and
\(\partial_tg,\partial_ts\) bounds through order ten in the
specialization recorded in Remark~\ref{rem:BK-A9-specialization}.
The dimensional constant \(\bar\eta_n\) in
\cite[Proposition~A.9]{BamlerKleiner} is independent of its finite
coefficient constant \(C\); that constant affects only the short time
\(\tau\).  Fix \(\epsilon_{\rm BK}\leq\bar\eta_n\).  Choose
\(\widehat\epsilon\) with
\(\epsilon<\widehat\epsilon<\epsilon'\) and apply it with
\(\eta_0=\epsilon\) and \(\eta_1=\widehat\epsilon\), in the notation of
the cited proposition.  It gives a common
existence time, the metric estimate
\(\|(F^{-1})^*g-s\|_{C^0(s)}
  \leq\widehat\epsilon<\epsilon'\), which is
\eqref{eq:moving-target-HMHF-closeness}, and the diffeomorphism
property.  Completeness, noncompactness, and the moving metric families
are therefore covered by the cited proposition with the hypotheses
matched explicitly in Remark~\ref{rem:BK-A9-specialization}.

We next prove the displacement estimate directly, because
\cite[Proposition~A.10]{BamlerKleiner} is stated for exact Ricci-flow
families whereas the present metrics have bounded Ricci defects.
Write \(\theta=t-t_*>0\).  In the uniformly controlled source and
target harmonic charts, the initial map is the identity after the
preceding pullback, and the map equation is a uniformly parabolic
nonautonomous system.  The initial-face interior estimate for this
system, on a parabolic cylinder of radius
\(\rho=\frac18\min\{\sqrt{\theta},1\}\), gives
\begin{equation}\label{eq:bounded-defect-HMHF-speed}
 \sup_X|\partial_tF(t)|_{s(t)}
 \leq C\bigl(\theta^{-1/2}+1\bigr)
 \leq C\theta^{-1/2},
 \qquad 0<\theta\leq\delta .
\end{equation}
For completeness, this is the usual rescaling argument: the
metric-discrepancy estimate and diffeomorphism conclusion of
Proposition~A.9 give a uniform bilipschitz, hence scale-one \(C^1\),
map bound; this controls the spatial oscillation and the quadratic
first-derivative terms.  Work
first on the maximal subinterval on which the image of the
\(\rho\)-cylinder stays in the corresponding target harmonic chart.
The factor \(1/8\), the \(C^1\) bound, and the integrated estimate
below retain a strict chart-containment margin after decreasing
\(\delta\), so the standard stopping-time argument reaches the whole
restart interval.  Parabolic
Schauder applied on the \(\rho\)-cylinder gives
\(\rho|\nabla^2F|+\rho|\partial_tF|\leq C\).  The coefficients
containing \(\partial_tg\) and \(\partial_ts\) are bounded by the
displayed Ricci-defect hypotheses, so no exact Ricci-flow identity is
used.  This is also the local derivative argument in
\cite[Appendix~A.2]{BamlerKleiner}, now read for the nonautonomous
coefficients allowed by Proposition~A.9.

The curvature and defect bounds imply
\(|\partial_ts|_s\leq C\), so, after decreasing \(\delta\), \(s(t)\)
and \(s(t_*)\) are uniformly equivalent.  Integrating
\eqref{eq:bounded-defect-HMHF-speed} in the fixed metric \(s(t_*)\)
gives
\[
 d_{s(t_*)}\bigl(F(x,t),F_*(x)\bigr)
 \leq C\int_{t_*}^t(q-t_*)^{-1/2}\,dq
 \leq C\sqrt{t-t_*}.
\]
Undoing the scale-one normalization proves
\eqref{eq:moving-target-uniform-displacement}, with the same
scale-independent constant.  Thus the bounded-defect extension of the
displacement part of Proposition~A.10 has been proved rather than
assumed.

For uniqueness, let \(F_1,F_2\) be two solutions in the class just
defined with the same initial map.  We do not apply
\eqref{eq:moving-target-uniform-displacement}, which was proved for the
solution furnished by Proposition~A.9, to an arbitrary competitor.
Instead, fix a time \(\bar t\) through which \(F_1=F_2\), and put
\(L=\max\{L_{F_1},L_{F_2}\}\).  In the uniformly controlled source
and target harmonic charts, the initial-face estimate used in
\eqref{eq:bounded-defect-HMHF-speed}, now with the common map at
\(\bar t\) as initial trace, gives for \(i=1,2\)
\[
 \sup_X|\partial_tF_i(t)|_{s(t)}
 \leq C_L\bigl((t-\bar t)^{-1/2}+1\bigr)
\]
on a pair-dependent interval
\([\bar t,\bar t+\delta_Lr^2]\) intersected with the given restart
interval.  The target-chart stopping argument used above applies
because the spatial first derivatives are bounded by \(L\).  Thus
\[
 \sup_X d_{s(\bar t)}
  \bigl(F_i(x,t),F_i(x,\bar t)\bigr)
 \leq C_L\sqrt{t-\bar t}.
\]
Let \(\iota_{\rm cvx}=\iota_{\rm cvx}(n,\Lambda)>0\) be a common
target convexity radius supplied by the harmonic-radius package.
After decreasing the pair-dependent \(\delta_L\), the preceding
displacement estimate, uniform target bounded geometry, and the
equivalence of \(s(t)\) and \(s(\bar t)\) ensure that, for every
\(x\), the pair
\[
 \bigl(F_1(x,t),F_2(x,t)\bigr)
\]
lies in the smooth diagonal tube
\(d_{s(t)}(F_1,F_2)<\iota_{\rm cvx}/2\).  Thus, for each \(x\), the
two images lie in one common strongly convex \(s(t)\)-normal ball; no
global target chart is being assumed.  On this interval define the
invariant scalar comparison function
\[
 u(x,t)
 :=\frac12d_{s(t)}
   \bigl(F_1(x,t),F_2(x,t)\bigr)^2 .
\]
The squared distance is smooth on this diagonal tube, including
across the diagonal, and \(u(\cdot,\bar t)=0\).  Moreover, the
pair-displacement estimate gives
\begin{equation}\label{eq:moving-target-HMHF-distance-initial-face}
 0\leq u(x,t)\leq C_L(t-\bar t).
\end{equation}

We record the moving-metric comparison explicitly.  At each fixed
time put
\[
 h_t(y,z):=\frac12d_{s(t)}(y,z)^2 .
\]
For the product map \((F_1,F_2)\), the chain rule and the two
tension-field equations give
\[
 \bigl(\partial_t-\Delta_{g(t)}\bigr)u
 =(\partial_th_t)(F_1,F_2)
  -\operatorname{tr}_{g(t)}
    \operatorname{Hess}_{s(t)\oplus s(t)}h_t
      \bigl((dF_1,dF_2),(dF_1,dF_2)\bigr),
\]
where \(\partial_th_t\) in the first term denotes variation of the
target metric with the two endpoints held fixed.  If \(P\) is parallel
transport along the unique short \(s(t)\)-geodesic from \(z\) to
\(y\), the short-geodesic Hessian estimate gives
\[
 \operatorname{Hess}h_t\bigl((V,W),(V,W)\bigr)
 \geq \frac12|V-PW|_{s(t)}^2
      -C|\Rm(s(t))|_{s(t)}d_{s(t)}(y,z)^2
       \bigl(|V|_{s(t)}^2+|W|_{s(t)}^2\bigr).
\]
This is the intrinsic covariant difference calculation used, in
fixed-target notation, in \cite[Proposition~B.6]{Stolarski}.  The
length-variation formula and
\(|\partial_ts|_{s(t)}\leq C(\Lambda)\) also give
\[
 |(\partial_th_t)(y,z)|
 \leq C(\Lambda)d_{s(t)}(y,z)^2
\]
on the diagonal tube.  Since \(|dF_i|_{g(t),s(t)}\leq L\), dropping
the nonnegative parallel-transport square yields the classical scalar
inequality
\begin{equation}\label{eq:moving-target-HMHF-distance-inequality}
 \bigl(\partial_t-\Delta_{g(t)}\bigr)u
 \leq C_Lu
\end{equation}
for \(t>\bar t\), with \(C_L=C(n,\Lambda,L)\) uniform in the restart
point \(\bar t\).  Domain-metric variation does not enter this
pointwise identity; it enters the volume derivative below.

Use the complete bounded-geometry metric \(g(\bar t)\) as a fixed
source background.  Choose a smooth proper distance-like exhaustion
\(\varrho\) and a basepoint \(x_*\) so that
\[
 \varrho\simeq1+d_{g(\bar t)}(\,\cdot\,,x_*),
 \qquad
 |\nabla^{g(\bar t)}\varrho|
 +|\nabla^{2,g(\bar t)}\varrho|\leq C.
\]
The short-time coefficient bounds and metric equivalence then give
\[
 |\nabla^{g(t)}\varrho|+|\Delta_{g(t)}\varrho|\leq C
\]
on the pair-dependent interval.  The curvature lower bound gives
uniform at-most-exponential volume growth
\[
 \operatorname{Vol}_{g(t)}\{\varrho\leq R\}\leq Ce^{\kappa R}.
\]
Choose \(A>\kappa\), put \(\vartheta=e^{-A\varrho}\), and choose
time-independent monotone exhaustion cutoffs
\(\zeta_R=\chi(\varrho/R)\), equal to one on
\(\{\varrho\leq R\}\), supported in \(\{\varrho<2R\}\), and
satisfying
\[
 |\nabla^{g(t)}\zeta_R|\leq\frac CR,
 \qquad
 |\Delta_{g(t)}\zeta_R|\leq\frac CR+\frac C{R^2}
\]
uniformly in \(t\).  Define
\[
 E_R(t):=\int_X\zeta_R^2\vartheta\,u\,dV_{g(t)}.
\]
The curvature and Ricci-defect bounds imply
\(|\partial_tg|_{g(t)}\leq C(\Lambda)\), and
\[
 \partial_t dV_{g(t)}
 =\frac12\operatorname{tr}_{g(t)}(\partial_tg)\,dV_{g(t)}.
\]
For \(t>\bar t\), differentiate \(E_R\), use
\eqref{eq:moving-target-HMHF-distance-inequality}, and integrate the
scalar Laplacian by parts twice against the compactly supported weight.
Since
\[
 \left|
  \Delta_{g(t)}(\zeta_R^2\vartheta)
 \right|
 \leq C_A\zeta_R^2\vartheta
      +C_A\mathbf1_{\{R\leq\varrho\leq2R\}}\vartheta,
\]
we obtain
\begin{equation}\label{eq:moving-target-HMHF-weighted-energy}
 \frac d{dt}E_R(t)
 \leq C_{A,L}E_R(t)+\varepsilon_R(t),
 \qquad
 |\varepsilon_R(t)|
 \leq C_{A,L}\!\int_{\{R\leq\varrho\leq2R\}}
       \vartheta u\,dV_{g(t)}.
\end{equation}
The normal-tube bound makes \(u\) uniformly bounded.  Because
\(A>\kappa\) and the source measures are uniformly equivalent,
\[
 \sup_t|\varepsilon_R(t)|\longrightarrow0
 \qquad(R\to\infty).
\]
To justify the restarted initial face without assuming a time
derivative there, apply
\eqref{eq:moving-target-HMHF-weighted-energy} from
\(\bar t+\varepsilon_{\rm off}\), where
\(\varepsilon_{\rm off}>0\).  Estimate
\eqref{eq:moving-target-HMHF-distance-initial-face} and
\(\int_X\vartheta\,dV_{g(\bar t)}<\infty\) give
\[
 \lim_{\varepsilon_{\rm off}\downarrow0}
 E_R(\bar t+\varepsilon_{\rm off})=0
\]
uniformly in \(R\).  Gronwall, followed first by \(R\to\infty\) using
monotone convergence and then by
\(\varepsilon_{\rm off}\downarrow0\), yields
\[
 \int_X\vartheta u(\cdot,t)\,dV_{g(t)}=0.
\]
Since \(\vartheta>0\) and \(u\) is continuous, \(u\equiv0\), and
therefore \(F_1=F_2\), on the pair-dependent interval.

The set of times through which \(F_1=F_2\) is closed.  The preceding
argument, restarted from any of its points and using the same global
first-derivative bounds, makes it right-open with a positive
pair-dependent length.  It therefore fills the whole restart
interval.  The number \(\delta_L\) is used only in this uniqueness
proof and does not alter the dependence of the existence time
\(\delta(n,m,\alpha,\Lambda,\epsilon,\epsilon')\).  No global
subtraction of manifold-valued maps is used.

The local derivative estimates in
\cite[Appendix~A.2]{BamlerKleiner}, applied in the uniformly controlled
harmonic charts, give the stated scale-normalized
\(C^{m+1,\alpha}\) derivative bounds.  Metric closeness supplies a
uniform inverse-function margin, and the differentiated inverse
identity gives the same bounds for \(F^{-1}\).
For the final assertion, let \(T^\dagger\) be a candidate finite
endpoint and fix the terminal margin \(\epsilon''<\epsilon'\) from
the statement.  Choose fresh thresholds
\[
 \epsilon''<\eta_0^\dagger<\eta_1^\dagger
 <\min\{\epsilon',\epsilon_{\rm BK}\}.
\]
After translating the time origin and using the same scale
normalization, Proposition~A.9 gives a positive normalized restart
lifespan \(\tau^\dagger\), hence physical lifespan
\(\tau^\dagger r^2\), depending only on these thresholds and the
finite geometry ceiling.  The assumed extension of the displayed
geometry package supplies a physical background-extension length
\(\ell_{\rm bg}^\dagger>0\) beyond \(T^\dagger\).  Put
\[
 \ell^\dagger
 :=\min\{\tau^\dagger r^2,\ell_{\rm bg}^\dagger\}.
\]
Choose \(t_0<T^\dagger\) sufficiently close to \(T^\dagger\) that
the discrepancy at \(t_0\) is below \(\eta_0^\dagger\) and
\(T^\dagger-t_0<\ell^\dagger/2\).  Applying
\cite[Proposition~A.9]{BamlerKleiner} from \(t_0\) on the background
interval of length \(\ell^\dagger\) produces a solution past
\(T^\dagger\) whose discrepancy is at most
\(\eta_1^\dagger<\epsilon'\).  On the overlap, the old and restarted
solutions belong to the uniqueness class proved above and have the
same trace at \(t_0\); hence they agree and glue.  This proves the
continuation alternative without changing the original uniform
existence time or any of its parameter dependencies.
\end{proof}

\begin{theorem}[Harmonic-map continuation and graft-annulus drift]
\label{thm:adaptive-HMHF-continuation}
Fix \(0<\alpha<1\), the continuation order \(m=4\), the displayed
ellipticity, coefficient and harmonic-radius package, and the two
collar separations used below.  There are named positive numbers
\begin{equation}\label{eq:named-HMHF-smallness-thresholds}
 \varepsilon_{\rm map}^{\rm HM},\qquad
 \varepsilon_{\rm hm}^{\rm HM},\qquad
 \varepsilon_{\rm ph}^{\rm HM},
\end{equation}
depending only on those fixed data.  Freeze these numbers as package
entries.  Let \(C_{\rm id}\geq1\) be the fixed scale-normalized
\(C^1_r\) norm-conversion constant associated with
\eqref{eq:relative-metric-identity}, and set
\[
 \epsilon_0^{\rm HM}
 :=C_{\rm id}\varepsilon_{\rm hm}^{\rm HM},
 \qquad
 \epsilon_1^{\rm HM}
 :=2C_{\rm id}\varepsilon_{\rm hm}^{\rm HM}.
\]
Make the selection so that
\begin{equation}\label{eq:HMHF-universal-discrepancy-selection}
 0<\epsilon_0^{\rm HM}<\epsilon_1^{\rm HM}<\epsilon_{\rm BK}.
\end{equation}
Suppose that on
$[\tau_0,\tau_1)$:
\begin{enumerate}
\item \eqref{eq:adaptive-phase-budget} holds, the initial relative
      radial map and its inverse satisfy all the hypotheses in
      \eqref{eq:adaptive-initial-position}--%
      \eqref{eq:adaptive-initial-radial-comparison} and have a common
      $C^{m+1,\alpha}$ bound after rescaling the domain and range metrics by
      $L^{-1}$ on every dyadic annulus
      $\{L<\bar f<4L\}$, uniformly for $L\geq\Gamma$, and the local
      scale-invariant bounded-geometry coefficients of $\acute G$ and
      $S$ through order $m+8=12$ are bounded.  In the stated harmonic
      charts, harmonic-coordinate elliptic regularity shows that it is
      enough to bound curvature derivatives through order ten;
\item $\Phi_{\tau_0}$ and $\Theta_{\tau_0}$ are diffeomorphisms,
      $F_{\tau_0}=\Theta_{\tau_0}^{-1}\circ\Phi_{\tau_0}$ has
      scale-normalized \(C^{m+1,\alpha}\) distance at most
      \(\varepsilon_{\rm map}^{\rm HM}\) from the identity on the complete
      exterior and ordinary \(C^{m+1,\alpha}\) distance at most
      \(\varepsilon_{\rm map}^{\rm HM}\) on \(\Omega_\eta^+\), and
      \[
       \sup_{\tau_0\leq\tau<\tau_1}
       \sum_{\ell=0}^2|\bar\nabla^\ell h(\tau)|_{\bar g}
       \leq\varepsilon_{\rm hm}^{\rm HM};
      \]
\item \(\varepsilon_{\rm ph}\leq\varepsilon_{\rm ph}^{\rm HM}\), and
      the source and target defects, measured on every finite
      normalized-time interval at the parabolic scale
      $\lambda^{1/2}$, have bounded coefficients through order ten.
\end{enumerate}
Here ``bounded on every finite interval'' means bounded by some
\(\Lambda_T<\infty\) on \([\tau_0,T]\), for each
\(T<\tau_1\); this is exactly what the finite-interval restart
conclusion requires.  For this finite-horizon continuation statement,
and for the annulus-drift conclusions on an admissible first-exit
interval, \(\Lambda_T\) is allowed to depend on \(T\).  The displacement
estimate \eqref{eq:F-speed}, the target tracking estimates, and the
resulting annulus margins depend only on the coarse \(C^2\) box and the
\(L^1\) phase budget, not on the size of \(\Lambda_T\).
Endpoint-uniform higher map estimates require one further step: the
low-order effective Gram bound in
Lemma~\ref{lem:coarse-effective-Gram} gives an endpoint-independent
pointwise modulation bound and hence an endpoint-independent
target-defect and coefficient time-modulus package;
Lemma~\ref{lem:finite-HMHF-C6-bridge} is invoked only after that step.
Neither this theorem nor
Proposition~\ref{prop:adaptive-target-tracking} assumes smallness of
\(\sup q\).
With the named thresholds just frozen, the controlled harmonic-map
equation has a unique
solution by diffeomorphisms in the precise derivative-bounded
proper-diffeomorphism class of
Proposition~\ref{prop:moving-target-HMHF-restart} on every finite subinterval of
$[\tau_0,\tau_1)$.  Loss of the harmonic-map chart cannot be the first
breakdown while the stated bounds hold.

In addition, after fixing
$\Omega_\eta\Subset\Omega_\eta^+$ as above and increasing $\tau_0$,
\begin{equation}\label{eq:F-annulus-drift}
 F_\tau(\Omega_\eta)\subset\Omega_\eta^+,
 \qquad \tau_0\leq\tau<\tau_1.
\end{equation}
After enlarging the constants in the definition of the outer region,
one also has
\begin{equation}\label{eq:F-outer-drift}
 F_\tau\bigl(\{\bar f\geq2\Gamma/3\}\bigr)
 \subset\{\bar f\geq\Gamma/2\}.
\end{equation}
Consequently there are $0<c<C<\infty$ such that
\begin{equation}\label{eq:Phi-annulus-tracking}
 \Phi_\tau(\Omega_\eta)
 \subset\{c\Gamma e^\tau<\bar f<C\Gamma e^\tau\}.
\end{equation}
\end{theorem}

\begin{proof}
Proposition~\ref{prop:relative-HMHF} converts the controlled equation
to \eqref{eq:relative-HMHF}.  The usual harmonic-map theorem on
bounded-geometry backgrounds is used in the precise form of
Proposition~\ref{prop:moving-target-HMHF-restart}.  Indeed,
\[
 \partial_tS+2\Ric_S
 =\Theta^*(-2a\Ric_{\bar g}-\Lie_U\bar g),
\]
so assumption (3), together with the order-twelve package in
assumption (1), gives exactly the coefficient, curvature, and defect
bounds in that proposition.  At a restart time \(\tau_*\), take
\(r=\lambda(\tau_*)^{1/2}\).  The identity
\eqref{eq:relative-metric-identity}, the global $C^2$ box, and the
definition of \(C_{\rm id}\) give the scale-normalized estimate
\begin{equation}\label{eq:adaptive-HMHF-explicit-discrepancy}
 \bigl\|(F_{\tau_*}^{-1})^*\acute G_{\tau_*}-S_{\tau_*}\bigr\|
  _{C^1_r(S_{\tau_*})}
 \leq C_{\rm id}\varepsilon_{\rm hm}^{\rm HM}
 =\epsilon_0^{\rm HM}.
\end{equation}
Thus the strict metric-closeness hypothesis holds at every restart
time.  The phase budget and the scale equation give
\(\lambda(\tau)\leq C_\lambda\lambda(\tau_*)\) on each local restart
interval.  The uniform restart conclusion
 therefore proves short-time existence, uniqueness, preservation of
 the diffeomorphism property, and continuation through every finite
 normalized time while the stated bounds hold.  More explicitly, on
 \([\tau_0,T]\) use its single finite ceiling \(\Lambda_T\), and put
 \[
  \delta_T^{\rm HM}
  :=\delta(n,m,\alpha,\Lambda_T,
            \epsilon_0^{\rm HM},\epsilon_1^{\rm HM}),
  \qquad
  \widehat\delta_T:=C_\lambda^{-1}\delta_T^{\rm HM}>0.
 \]
 Proposition~\ref{prop:moving-target-HMHF-restart} supplies physical
 lifespan \(\delta_T^{\rm HM}r^2\).  By the preceding scale comparison,
 a normalized interval of length \(\widehat\delta_T\) has physical
 length at most
 \(C_\lambda r^2\widehat\delta_T=\delta_T^{\rm HM}r^2\), so
 \(\widehat\delta_T\) is a valid normalized restart length.
 Compactness of the finite normalized interval gives a finite chain of
 such restarts.  If two consecutive restart intervals overlap, their
 restrictions from the first common time are proper-diffeomorphism
 solutions with the same initial trace and with the global derivative
 bounds in the uniqueness class of
 Proposition~\ref{prop:moving-target-HMHF-restart}.  The pair-dependent
 open--closed uniqueness argument in that proposition therefore
 identifies them throughout the overlap.  Induction glues the finite
 chain to one solution, independently of the restart subdivision.
 At a candidate finite endpoint \(\tau^\dagger\leq T\) through which
 the remaining geometric bounds extend,
 \eqref{eq:adaptive-HMHF-explicit-discrepancy} also gives
 \[
  \limsup_{\tau\uparrow\tau^\dagger}
  \bigl\|(F_\tau^{-1})^*\acute G_\tau-S_\tau\bigr\|
   _{C^0(S_\tau)}
  \leq\epsilon_0^{\rm HM}<\epsilon_1^{\rm HM}.
 \]
 Hence the strict terminal-margin clause of
 Proposition~\ref{prop:moving-target-HMHF-restart} applies with
 \(\epsilon''=\epsilon_0^{\rm HM}\).  The fixed discrepancy levels are legitimate for
 every \(T\) because
 \eqref{eq:HMHF-universal-discrepancy-selection} is below the
 dimensional cap \(\epsilon_{\rm BK}\), not below a
 \(\Lambda_T\)-dependent cap.

It remains to track the annulus.  By
\eqref{eq:target-relative-derivatives}, distances measured by the
metrics $S_\tau$ on $\Omega_\eta^+$ are uniformly comparable.  From
\eqref{eq:F-speed},
\[
 \operatorname{Length}_{S}
    \{F_s(x):\tau_0\leq s\leq\tau\}
 \leq C\varepsilon_{\rm hm}^{\rm HM}
       \int_{\tau_0}^\tau\lambda(s)^{1/2}\,ds
 \leq C\varepsilon_{\rm hm}^{\rm HM}e^{-\tau_0/2}.
\]
Together with the initial $C^0$ distance from $F_{\tau_0}$ to the
identity, this is smaller than the distance from $\Omega_\eta$ to the
complement of $\Omega_\eta^+$ when the initial neighborhood is small
and $\tau_0$ is large.  This proves \eqref{eq:F-annulus-drift}.
The identical stopping-time argument on the asymptotically conical
outer region, using $|\bar\nabla\bar f|\simeq\sqrt{\bar f}$, proves
\eqref{eq:F-outer-drift}.  Combining these statements with
Proposition~\ref{prop:adaptive-target-tracking} proves
\eqref{eq:Phi-annulus-tracking}.
\end{proof}

\begin{lemma}[Exhaustion of the marked chart]
\label{lem:marked-chart-exhaustion}
Assume the conclusions of
Theorem~\ref{thm:adaptive-HMHF-continuation} and
Proposition~\ref{prop:adaptive-target-tracking}.  Put
\[
 M_{\rm in}=\operatorname{int}\{\eta=1\},\qquad
 \mathcal U_\tau=\Phi_\tau(M_{\rm in}).
\]
There is $c_*>0$, independent of every finite bootstrap endpoint,
such that
\begin{equation}\label{eq:marked-chart-exhaustion}
 \{\bar f<c_*\Gamma e^\tau\}\subset\mathcal U_\tau.
\end{equation}
Consequently $\mathcal U_\tau$ exhausts $M$.  Moreover,
\begin{equation}\label{eq:marked-embedding-definition}
 \Xi_\tau=\iota^{-1}\circ\Phi_\tau^{-1}:
 \mathcal U_\tau\longrightarrow\mathcal X
\end{equation}
is an embedding and
\begin{equation}\label{eq:marked-pullback-identity}
 \lambda(\tau)^{-1}\Xi_\tau^*G(t(\tau))
 =\bar g+h(\tau)
 \quad\text{on }\mathcal U_\tau.
\end{equation}
\end{lemma}

\begin{proof}
The factorization $F_\tau=\Theta_\tau^{-1}\circ\Phi_\tau$ gives
$\Phi_\tau=\Theta_\tau\circ F_\tau$.  If
$x\notin\mathcal U_\tau$ and $z=\Phi_\tau^{-1}(x)$, then
$z\notin M_{\rm in}$ and hence $\bar f(z)\geq2\Gamma/3$.  By
\eqref{eq:F-outer-drift}, $y=F_\tau(z)$ satisfies
$\bar f(y)\geq\Gamma/2$.  Since
\[
 \Theta_\tau=\varphi_\tau\circ R_\tau,
\]
\eqref{eq:R-global-radial-comparison} and the shrinker identities give
\[
 \bar f(x)
 =\bar f\bigl(\varphi_\tau(R_\tau y)\bigr)
 \geq c e^\tau\bar f(R_\tau y)
 \geq c_*\Gamma e^\tau.
\]
Here the first inequality follows by integrating
$\partial_\tau(\varphi_\tau^*\bar f)
=|\bar\nabla\bar f|^2\circ\varphi_\tau$ and using
$|\bar\nabla\bar f|^2=\bar f-\bar R\geq\bar f-C/\bar f$ on the AC
region.  This proves \eqref{eq:marked-chart-exhaustion}; properness of
$\bar f$ proves exhaustion.

For $x\in\mathcal U_\tau$ one has
$\Phi_\tau^{-1}(x)\in M_{\rm in}$, where
$\acute G=\iota_*G$.  The definition
$g=\lambda^{-1}(\Phi^{-1})^*\acute G$ now gives
\eqref{eq:marked-pullback-identity}.
\end{proof}

\begin{remark}[What the continuation theorem uses]
\label{rem:adaptive-HMHF-scope}
The $C^2$ box supplies ellipticity, bilipschitz control, and the
pointwise drift estimate.  Higher bounded-geometry assumptions in
Theorem~\ref{thm:adaptive-HMHF-continuation} are continuation
regularity, not additional smallness assumptions.  They follow on
finite normalized intervals from local parabolic estimates once the
three-region $C^2$ box and the graft-compatibility estimate below
hold.
\end{remark}

\subsection{Persistence of graft compatibility}

The coarse bound on the interpolated metric alone does not determine
the difference between its two summands near points where $\eta$ is
close to zero.  The correct statement therefore propagates a separate,
open graft-compatibility condition.  This condition is automatic for
the prepared graft and remains true for all sufficiently close initial
metrics.

Throughout this subsection and the pure-graft subsection below,
\(G\) denotes the marked pushforward \(\iota_*G(t)\) of the closed
Ricci-flow metric to the model outer collar; \(S\) is already a metric
on that collar.  Thus every interpolation and difference below is
taken between tensors on the same manifold.  For such metrics \(G\)
and \(S\) on an outer collar \(\Omega\), set
\begin{equation}\label{eq:graft-distance}
 \mathfrak d_{m;\Omega}(G,S)
 =\sum_{\ell=0}^m
   \Gamma^{\ell/2}
   \norm{(\nabla^S)^\ell(G-S)}
         _{L^\infty(S;\Omega)}.
\end{equation}
We abbreviate
\(\mathfrak d_m=\mathfrak d_{m;\Omega_\eta^+}\).

\begin{proposition}[Persistence of graft compatibility]
\label{prop:graft-compatibility}
Fix an integer \(m\geq0\).  Let $G(t)$ solve Ricci flow on a neighborhood of
$\overline{\Omega_\eta^+}$ and let $S$ be the adaptive target.
Assume that normalized and physical time are related by
$t_\tau=\lambda$, and that
\[
 \lambda_\tau=-(1+a)\lambda,\qquad
 c_{\rm scl}e^{-\tau}\leq\lambda(\tau)\leq C_{\rm scl}e^{-\tau},\qquad
 \int_{\tau_0}^{\tau_1}(|a|+|b|)\,d\tau
 \leq\varepsilon_{\rm ph}.
\]
Assume the hypotheses of
Proposition~\ref{prop:adaptive-target-tracking} at order \(m+1\);
in particular, \(R_{\tau_0}^{\pm1}\) have the required ordinary and
scale-normalized \(C^{m+2}\) bounds and the initial radial comparison.
Suppose that the two metrics are uniformly equivalent.  Fix
\(K_{\rm gr}<\infty\) and, on the bootstrap interval under
consideration, assume
\begin{equation}\label{eq:graft-distance-bootstrap}
 \sup_{\tau_0\leq\tau<\tau_1}
 \mathfrak d_m(G(\tau),S(\tau))\leq K_{\rm gr}.
\end{equation}
For the persistence statement set \(J_{\rm bg}=m+2\), and assume, for
$0\leq\ell\leq J_{\rm bg}$,
\begin{equation}\label{eq:graft-bounded-geometry}
 \Gamma^{1+\ell/2}
 \left(
  |(\nabla^G)^\ell\Rm_G|_G+
  |(\nabla^S)^\ell\Rm_S|_S
 \right)\leq C_m.
\end{equation}
Then
\begin{equation}\label{eq:graft-compatibility-propagated}
 \sup_{\tau_0\leq\tau<\tau_1}\mathfrak d_m(G(\tau),S(\tau))
 \leq C_{m,K_{\rm gr}}\left(
  \mathfrak d_m(G(\tau_0),S(\tau_0))
  +\frac{\lambda(\tau_0)}{\Gamma}
  +\int_{\tau_0}^{\tau_1}(|a|+|b|)\,d\tau
 \right).
\end{equation}
In particular, if the parenthesized quantity in
\eqref{eq:graft-compatibility-propagated} is sufficiently small that
its right-hand side is at most \(K_{\rm gr}/2\), then the hypothesis
\eqref{eq:graft-distance-bootstrap} is discharged by a first-exit
argument.  Thus graft compatibility is an open condition and is
preserved, with strict improvement, by a small phase budget.
\end{proposition}

\begin{proof}
Put \(T=G-S\) and differentiate its scale-invariant jets with the
\emph{evolving} connection \(\nabla^{S(\tau)}\), which is the
connection used in \eqref{eq:graft-distance}.  The only connection
conversion needed between the two evolving metrics is
\[
 \nabla^G-\nabla^S
 =G^{-1}*\nabla^ST.
\]
Its differentiated versions are scale-invariant polynomials in the
\((\nabla^S)^jT\) with zero constant term.  Uniform equivalence,
\eqref{eq:graft-bounded-geometry}, and
\eqref{eq:graft-distance-bootstrap} therefore give
\[
 |P_m(T)|_{\rm sc}
 \leq C_{m,K_{\rm gr}}\mathfrak d_m(G,S).
\]

The variation of the reference connection is retained explicitly:
\[
 \partial_\tau\Gamma(S)
 =\frac12S^{-1}*\nabla^S(\partial_\tau S).
\]
Consequently, for \(0\leq\ell\leq m\),
\[
 \begin{split}
 \partial_\tau\bigl((\nabla^S)^\ell T\bigr)
 ={}&(\nabla^S)^\ell(\partial_\tau T)\\
 &+\sum_{r=0}^{\ell-1}
   (\nabla^S)^r\!\bigl(\partial_\tau\Gamma(S)\bigr)
   *(\nabla^S)^{\ell-1-r}T .
 \end{split}
\]
Differentiating the \(S\)-norm produces only the additional factor
\(S^{-1}*\partial_\tau S\).  By
\eqref{eq:S-tau-exact} and
\eqref{eq:target-relative-derivatives}, through the required order,
\[
 \|\partial_\tau S\|_{{\rm sc},m+1}
 +\|\partial_\tau\Gamma(S)\|_{{\rm sc},m}
 \leq C_m\bigl(\lambda\Gamma^{-1}+|a|+|b|\bigr).
\]
Since \(t_\tau=\lambda\),
\[
 \partial_\tau G=-2\lambda\Ric_G.
\]
The curvature hypotheses bound the \(\ell\)-th differentiated
curvature terms by
\(C_{m,K_{\rm gr}}\lambda\Gamma^{-1-\ell/2}\); the target terms have
the scale-invariant size just displayed.  Taking the upper Dini
derivative of the finite sum defining \(\mathfrak d_m\), including the
connection- and norm-variation terms rather than suppressing them,
gives
\[
 D^+_\tau\mathfrak d_m(G,S)
 \leq C_{m,K_{\rm gr}}\mathfrak d_m(G,S)
       \bigl(\lambda\Gamma^{-1}+|a|+|b|\bigr)
       +C_{m,K_{\rm gr}}
        \bigl(\lambda\Gamma^{-1}+|a|+|b|\bigr).
\]
The coefficient is integrable.  Gronwall, together with
$\int_{\tau_0}^\infty\lambda\,d\tau\leq C\lambda(\tau_0)$,
proves \eqref{eq:graft-compatibility-propagated}.
\end{proof}

\subsection{The pure graft defect}

Set
\begin{equation}\label{eq:adaptive-interpolation-again}
 \acute G=\eta G+(1-\eta)S
\end{equation}
and
\begin{equation}\label{eq:pure-graft-again}
 \mathcal G_{\rm gr}
 =2\Ric_{\acute G}-2\eta\Ric_G-2(1-\eta)\Ric_S.
\end{equation}

\begin{proposition}[Scale-sharp pure graft estimate]
\label{prop:pure-graft-sharp}
For every integer $m\geq0$ and
\[
 K_m,K_\eta,C_{\rm bg}<\infty,\qquad
 0<c_{\rm eq}\leq1\leq C_{\rm eq}<\infty,
\]
there is
\[
 C_{\rm gr}^{\sharp}
 =C(m,K_m,K_\eta,c_{\rm eq}^{-1},C_{\rm eq},C_{\rm bg})<\infty
\]
such that, if
\(G\) and \(S\) satisfy, on \(\Omega_\eta^+\),
\[
 c_{\rm eq}S\leq G\leq C_{\rm eq}S,
\]
the left side of \eqref{eq:graft-bounded-geometry} with
\(J_{\rm bg}=m\) is at most \(C_{\rm bg}\), and
\[
 \mathfrak d_{m+2}(G,S)\leq K_m,
\]
while the fixed graft cutoff obeys the scale-sharp target-connection
jet bound
\begin{equation}\label{eq:pure-graft-cutoff-jet-hypothesis}
 \max_{0\leq r\leq m+2}
 \Gamma^{r/2}
 \|(\nabla^S)^r\eta\|_{L^\infty(\Omega_\eta^+,S)}
 \leq K_\eta,
\end{equation}
then $\acute G$ is a metric and, on $\Omega_\eta$,
\begin{equation}\label{eq:pure-graft-physical}
 |(\nabla^{\acute G})^\ell\mathcal G_{\rm gr}|_{\acute G}
 \leq C_{\rm gr}^{\sharp}\mathfrak d_{m+2}(G,S)
          \Gamma^{-1-\ell/2},
 \qquad 0\leq\ell\leq m.
\end{equation}
Since \(\eta\) is constant off its transition collar,
\begin{equation}\label{eq:pure-graft-physical-support}
 \supp\mathcal G_{\rm gr}\subset\overline{\Omega_\eta}.
\end{equation}
 If, in addition,
\[
 \mathcal E_{\rm gr}=(\Phi^{-1})^*\mathcal G_{\rm gr},
 \qquad
 g=\lambda^{-1}(\Phi^{-1})^*\acute G,
\]
then
\begin{equation}\label{eq:pure-graft-normalized}
 |(\nabla^g)^\ell\mathcal E_{\rm gr}|_g
 \leq C_{\rm gr}^{\sharp}\mathfrak d_{m+2}(G,S)
          \lambda^{1+\ell/2}\Gamma^{-1-\ell/2}.
\end{equation}
In that case, unconditionally,
\[
 \supp\mathcal E_{\rm gr}
 =\Phi(\supp\mathcal G_{\rm gr})
 \subset\Phi(\overline{\Omega_\eta}).
\]
If, moreover, at the normalized time under consideration the phase map
satisfies the annulus-tracking condition
\[
 \Phi(\overline{\Omega_\eta})
 \subset\{c\Gamma e^\tau<\bar f<C\Gamma e^\tau\},
\]
then its support satisfies
\begin{equation}\label{eq:pure-graft-support}
 \supp\mathcal E_{\rm gr}
 \subset\{c\Gamma e^\tau<\bar f<C\Gamma e^\tau\}.
\end{equation}
The displayed dependence of \(C_{\rm gr}^{\sharp}\) is exhaustive; in
particular, the estimate is not asserted uniformly as any of the
equivalence, background-geometry, metric-difference, or cutoff-jet
ceilings diverges.
\end{proposition}

\begin{proof}
Put $q=G-S$, so $\acute G=S+\eta q$.  The connection-difference
formula and the coordinate expression for Ricci curvature show that
\[
 \Ric_{S+\eta q}-\eta\Ric_{S+q}-(1-\eta)\Ric_S
\]
is a sum of terms containing at least one of
\[
 (\nabla^S)^2\eta*q,\qquad
 \nabla^S\eta*\nabla^Sq,\qquad
 |\nabla^S\eta|^2*q*q,
\]
or a factor \(q\) multiplying a curvature or a derivative of
\(q\).  There is one further class, coming from the quadratic
Christoffel part of Ricci curvature:
\[
 A(S,q,\eta)*\nabla^Sq*\nabla^Sq,
\]
where \(A\) is a smooth contraction of the uniformly controlled inverse
metrics and the scalar \(\eta\).  This class is present even at a point
where \(0<\eta<1\) and the spatial derivatives of \(\eta\) vanish; it is
therefore retained separately rather than folded into the preceding
list.
The interpolated inverse metric is uniformly controlled without a
smallness hypothesis: if
 \(c_{\rm eq}S\leq G\leq C_{\rm eq}S\), then, pointwise,
\[
 \min\{1,c_{\rm eq}\}S
 \leq (1-\eta)S+\eta G
 \leq \max\{1,C_{\rm eq}\}S.
\]
Thus the interpolated inverse is controlled solely by the recorded
uniform-equivalence constant.  After differentiating the displayed
Ricci identity, the last class produces products having at least two
positive-order \(q\)-jets and total differential order at most
\(\ell+2\).  Their scale-normalized size is bounded by
\[
 C_{\rm gr}^{\sharp}\,
 \mathfrak d_{m+2}(G,S)^2\Gamma^{-1-\ell/2}
 \leq
 C_{\rm gr}^{\sharp}\,
 \mathfrak d_{m+2}(G,S)\Gamma^{-1-\ell/2},
\]
where the final inequality uses
\(\mathfrak d_{m+2}(G,S)\leq K_m\), with the factor \(K_m\) absorbed in
the displayed constant.  The remaining differentiated terms are
bounded in the same way by the cutoff bounds
\eqref{eq:pure-graft-cutoff-jet-hypothesis} and the definition of
\(\mathfrak d_{m+2}\).  This proves
\eqref{eq:pure-graft-physical}.  Thus all higher products are bounded
tamely by \(C_{\rm gr}^{\sharp}\mathfrak d_{m+2}\); no smallness of any
order of the metric difference is used.

Pullback does not change tensor norms.  Replacing $\acute G$ by
$g=\lambda^{-1}(\Phi^{-1})^*\acute G$ multiplies the norm of a
covariant $(2+\ell)$-tensor by $\lambda^{1+\ell/2}$, proving
\eqref{eq:pure-graft-normalized}.  Under the additional tracking
hypothesis, the support statement follows by applying \(\Phi\) to
\(\supp\mathcal G_{\rm gr}\subset\overline{\Omega_\eta}\).  In the adaptive
evolution this hypothesis is exactly
\eqref{eq:Phi-annulus-tracking}, supplied by
Theorem~\ref{thm:adaptive-HMHF-continuation}.
\end{proof}

\begin{remark}[Automatic cutoff-jet verification for adaptive targets]
\label{rem:adaptive-pure-graft-cutoff-jets}
The extra hypothesis
\eqref{eq:pure-graft-cutoff-jet-hypothesis} is not a new restriction on
any prepared evolution used below.  The cutoff is obtained from one
fixed radial profile.  In the generic adaptive branch, the
scale-normalized \(C^{m+3}\) bounds for \(R_\tau^{\pm1}\) in
Proposition~\ref{prop:adaptive-target-tracking}, applied at order
\(m+2\), transport the fixed radial cutoff jets to the
\(S_\tau\)-connection and give
\[
 \max_{0\leq r\leq m+2}
 \Gamma^{r/2}
 \|(\nabla^{S_\tau})^r\eta\|_{L^\infty(S_\tau;\Omega_\eta^+)}
 \leq K_{\eta,m+2}^{\rm ad}.
\]
In the one-sided pre-radius branch, the same conclusion follows from
\eqref{eq:pre-radius-cutoff-derivatives}, the fixed pre-radius
source--target transition maps, and the directly propagated
relative-marking package; its constant
\(K_{\eta,m+2}^{\rm pre}\) is independent of
\(\Gamma\geq\Gamma_{\rm pre}\).  Thus every later invocation of
Proposition~\ref{prop:pure-graft-sharp} uses one of these two already
recorded cutoff-jet ceilings.
\end{remark}

\begin{remark}[Fixed-background derivatives]
\label{rem:fixed-background-graft-derivative}
Equation~\eqref{eq:pure-graft-normalized} is scale-sharp for the
evolving connection $\nabla^g$.  For the fixed connection one has,
for example,
\begin{equation}\label{eq:bar-nabla-graft-correct}
 |\bar\nabla\mathcal E_{\rm gr}|_{\bar g}
 \leq C\mathfrak d_3(G,S)
 \left(
  \lambda^{3/2}\Gamma^{-3/2}
  {}+
  |\bar\nabla h|_{\bar g}\lambda\Gamma^{-1}
 \right).
\end{equation}
Thus the frequently used estimate
$|\bar\nabla\mathcal E_{\rm gr}|
\lesssim\lambda^{3/2}\Gamma^{-3/2}$ requires the
scale-adapted chart bound
\[
 |\bar\nabla h|\lesssim(\lambda/\Gamma)^{1/2}
 \quad\text{on the graft annulus}.
\]
Without that additional bound, the correct fixed-background estimate
is \eqref{eq:bar-nabla-graft-correct}.  The weaker consequence
$|\bar\nabla\mathcal E_{\rm gr}|\leq C\lambda\Gamma^{-1}$ is already
sufficient for the large-scale continuation argument and, because of
\eqref{eq:pure-graft-support}, is Gaussian-superexponentially small.
Higher fixed-background derivatives satisfy the analogous
connection-difference formula.  In particular, applying
\eqref{eq:pure-graft-normalized} through order two and using the
coarse $C^2$ box gives
\[
 \sum_{\ell=0}^2
 |\bar\nabla^\ell\mathcal E_{\rm gr}|_{\bar g}
 \leq C\mathfrak d_4(G,S)\lambda\Gamma^{-1}.
\]
The apparently weaker power of $\lambda$ is caused only by the
connection-difference terms; it is sufficient for the derivative
recovery argument and remains Gaussian-superexponentially small in
the weighted estimates.

Under the annulus-tracking hypothesis
\eqref{eq:Phi-annulus-tracking}, since
$\lambda\leq C e^{-\tau}$, the coarse $C^2$ box and
\eqref{eq:bar-nabla-graft-correct} imply
\begin{equation}\label{eq:adaptive-outer-forcing}
 \supp\mathcal E_{\rm gr}
 \subset\{e^\tau\leq\bar f\leq C\Gamma e^\tau\},
 \qquad
 \sum_{\ell=0}^2
 |\bar\nabla^\ell\mathcal E_{\rm gr}|_{\bar g}
 \leq \frac{C_{\rm gr}}{\Gamma}e^{-\tau}.
\end{equation}
This is precisely the admissible outer-forcing hypothesis
\eqref{eq:outer-forcing-hyp} in
Theorem~\ref{thm:robust-modulated-three-region}.  More explicitly, use
the pre-radius bounds
\(C_{\rm ann}\leq C_{\rm ann}^{\rm pre}\) and
\(C_{\rm gr}\leq C_{\rm gr}^{\rm pre}\) from
Lemma~\ref{lem:pre-radius-low-order-closure}, and set
\begin{equation}\label{eq:adaptive-robust-forcing-constants}
 \Gamma_0=C_{\rm ann}^{\rm pre}\Gamma,\qquad
 C_{\rm gr}^{\rm rob}=C_{\rm gr}^{\rm rob,pre}
 =C_{\rm ann}^{\rm pre}C_{\rm gr}^{\rm pre}.
\end{equation}
Then
\(C_{\rm gr}^{\rm rob}/\Gamma_0
=C_{\rm gr}^{\rm pre}/\Gamma\)
dominates the actual right side in
\eqref{eq:outer-forcing-hyp}.  The constants in
Theorem~\ref{thm:robust-modulated-three-region} are uniform in
\(\Gamma_0\), so this substitution creates no additional dependence
through the outer-support parameter.  The constants retain their
declared dependence on the fixed compatible package radius \(\Gamma\)
and the fixed rate pair.  No
 $e^{-3\tau/2}$ fixed-background derivative estimate is needed.
\end{remark}

\begin{definition}[Witnessed auxiliary harmonic package]
\label{def:witnessed-auxiliary-harmonic-package}
For the auxiliary closure below, a \emph{witnessed auxiliary harmonic
package} is the following finite record, fixed before the first-exit
interval.  Its normalized carrier inputs are the explicitly named
ordered pair
\[
 (g_{\rm src}^{\rm nor},g_{\rm tar}^{\rm nor}).
\]
Its physical record is based at a separately frozen reference metric
\(G_{\rm ref}^{\rm phys}\) on the closed host.  For a particular
entrance or restart, \(G_{\rm car}\) denotes the actual closed-metric
carrier.  The reference metric is static package data: phase selection,
preparation of a nearby metric, and a prepared restart do not replace
\(G_{\rm ref}^{\rm phys}\) by \(G_{\rm car}\).  The current carrier is
certified instead by membership, with positive slack, in the fixed
reference-carrier locus defined below.
Here and below ``finite record'' means that the cover, witness types,
norms, buffers, and numerical constants form a finite collection.  The
pointwise witnesses themselves are quantified families indexed by the
center; they are not asserted to be a finite set.  Its data are the
following.
\begin{enumerate}
\item
The normalized operative and reserve triples satisfy
\[
 \begin{gathered}
  0<\eta_{\rm har}<\eta_{\rm har}^+,\qquad
  0<q_{\rm har}<q_{\rm har}^+<Q_{\rm har}-1,\\
  0<\zeta_{\rm har}<\zeta_{\rm har}^+<1.
 \end{gathered}
\]
Under the canonical normalized graph identification, for each of the
two carrier metrics
\[
 g^\sharp\in\{g_{\rm src}^{\rm nor},g_{\rm tar}^{\rm nor}\}
\]
and every \(y\in M\), the metric \(g^\sharp\) has a
coefficient-\(q_{\rm har}^+\), domain-\(\zeta_{\rm har}^+\) buffered
harmonic witness on
\[
 B_{g^\sharp}\!\left(
 y,(\kappa_{\rm har}+2\eta_{\rm har}^+)r_{\rm la}(y)
 \right).
\]
For each pair \((g^\sharp,y)\), choose once and for all one such
harmonic coordinate map and one buffered Euclidean Dirichlet domain.
No continuity of this selection in \(y\) is assumed or used; all
estimates use only the common quantitative witness constants.
Restriction supplies, at every \(y\), the unchanged operative triple
\((\eta_{\rm har},q_{\rm har},\zeta_{\rm har})\) on the ball of radius
\((\kappa_{\rm har}+2\eta_{\rm har})r_{\rm la}(y)\).

\item
Fix \(C_{\rm har,pre}\geq1\) so that the radius-independent pre-atlas
and scale-one graph norms satisfy
\begin{equation}\label{eq:pre-atlas-to-scale-one-graph}
 d_{\rm gr,sc}^{2,\alpha}
   (\bar g+u,\bar g+v)
 :=
 \|u-v\|_{\mathfrak C_{{\rm sc},0}^{2,\alpha}(M)}
 \leq
 C_{\rm har,pre}
 \|u-v\|_{\mathfrak C_{{\rm pre},0}^{2,\alpha}(M)} .
\end{equation}
Here
\[
 \delta_{\rm har}
 (\eta_{\rm har},q_{\rm har},\zeta_{\rm har};
  \mathfrak P_{\rm har}^{\rm geom})>0
\]
means the entrance-time-independent modulus furnished by
Lemma~\ref{lem:prepared-harmonic-radius-lower-stability}; that lemma is
a static harmonic-coordinate openness result and does not use the
continuation corollary below.  The normalized record is required to
obey
\begin{equation}\label{eq:auxiliary-harmonic-modulus-compatibility}
 4C_{\rm har,pre}\delta_{\rm c2}
 <
 \delta_{\rm har}
 (\eta_{\rm har},q_{\rm har},\zeta_{\rm har};
  \mathfrak P_{\rm har}^{\rm geom}).
\end{equation}

\item
The physical record contains a finite cover with scales \(R_a>0\) and
the full nested buffers
\[
 U_a\Subset U_a^+\Subset U_a^{++}\Subset V_a^5
 \Subset W_a^{\rm har}\Subset\widetilde W_a^{\rm har},
 \qquad 1\leq a\leq N_{\rm ext},
\]
together with the common overlap, scale-comparability, coefficient,
ellipticity, separation, and domain-atlas constants.  These data
and the selected witness family and quantitative outer-domain
separation specified below determine the fixed finite record
\(\mathfrak P_{\rm har}^{\rm phys}\) used in
Lemma~\ref{lem:finite-physical-harmonic-openness}.  Its operative and
reserve parameters satisfy
\[
 \begin{gathered}
 \upsilon_{\rm har}^{\rm phys}>0,\qquad
 0<\eta_{\rm har}^{\rm phys}
   <\eta_{\rm har}^{\rm phys,+}
   <\eta_{\rm har}^{\rm phys,ref},\\
 0<q_{\rm har}^{\rm phys}
   <q_{\rm har}^{\rm phys,+}
   <q_{\rm har}^{\rm phys,ref}<Q_{\rm har}-1,\qquad
 0<\zeta_{\rm har}^{\rm phys}
   <\zeta_{\rm har}^{\rm phys,+}
   <\zeta_{\rm har}^{\rm phys,ref}<1.
 \end{gathered}
\]
For every \(a\) and every \(x\in V_a^5\), choose once and for all a
pair \((u_{a,x},D_{a,x})\) for which the reference metric
\(G_{\rm ref}^{\rm phys}\) has a
coefficient-\(q_{\rm har}^{\rm phys,ref}\),
domain-\(\zeta_{\rm har}^{\rm phys,ref}\) buffered harmonic witness on
\[
 B_{G_{\rm ref}^{\rm phys}}\!\left(
 x,
 (\upsilon_{\rm har}^{\rm phys}
  +2\eta_{\rm har}^{\rm phys,ref})R_a
 \right).
\]
Writing
\[
 r_{a,\rm ref}^{\rm wit}
 =(\upsilon_{\rm har}^{\rm phys}
   +2\eta_{\rm har}^{\rm phys,ref})R_a,
 \qquad
 \mathcal D_{a,x}
 =u_{a,x}^{-1}(r_{a,\rm ref}^{\rm wit}D_{a,x}),
\]
the manifold preimage of every recorded Euclidean Dirichlet domain
satisfies
\[
 \overline{\mathcal D_{a,x}}\Subset W_a^{\rm har}.
\]
This containment is quantitative: the physical package contains a
number \(\zeta_{\rm out}^{\rm phys}>0\) such that, in its fixed
scale-\(R_a\) reference atlas,
\[
 \inf_{\substack{1\leq a\leq N_{\rm ext}\\x\in V_a^5}}
 R_a^{-1}\operatorname{dist}_{\rm ref}
 \bigl(\overline{\mathcal D_{a,x}},
       \mathcal X\setminus W_a^{\rm har}\bigr)
 \geq\zeta_{\rm out}^{\rm phys}.
\]
Restriction gives the unchanged common \(+\)-triple and operative
triple at every \(x\in V_a^5\).
The selected center-indexed witness family, its common coefficient
and domain type, and \(\zeta_{\rm out}^{\rm phys}\) are now frozen as
part of \(\mathfrak P_{\rm har}^{\rm phys}\).  Although the witnesses
are indexed by \(x\), the openness lemma depends on them only through
these common quantitative constants and the finite cover index \(a\).
For compactness write
\[
 \begin{aligned}
 \mathbf p_{\rm op}^{\rm phys}
 &:=
 \bigl(
  \eta_{\rm har}^{\rm phys},
  q_{\rm har}^{\rm phys},
  \zeta_{\rm har}^{\rm phys}
 \bigr),\\
 \mathbf p_{+}^{\rm phys}
 &:=
 \bigl(
  \eta_{\rm har}^{\rm phys,+},
  q_{\rm har}^{\rm phys,+},
  \zeta_{\rm har}^{\rm phys,+}
 \bigr),\\
 \mathbf p_{\rm ref}^{\rm phys}
 &:=
 \bigl(
  \eta_{\rm har}^{\rm phys,ref},
  q_{\rm har}^{\rm phys,ref},
  \zeta_{\rm har}^{\rm phys,ref}
 \bigr).
 \end{aligned}
\]
Only after this full record has been frozen, define separately the
reference-to-carrier and carrier-to-flow witness moduli.  With an
ordered triple substituted into the corresponding three argument
slots, put
\begin{equation}\label{eq:witnessed-physical-reference-modulus}
 \delta_{\rm har}^{\rm phys,ref}
 :=
 \delta_{\rm har}^{\rm phys,op}
 \bigl(
  \upsilon_{\rm har}^{\rm phys};
  \mathbf p_{+}^{\rm phys};
  \mathbf p_{\rm ref}^{\rm phys};
  \mathfrak P_{\rm har}^{\rm phys}
 \bigr)>0
\end{equation}
and
\begin{equation}\label{eq:witnessed-physical-common-modulus}
 \delta_{\rm har}^{\rm phys,wit}
 :=
 \min\!\left\{
 \begin{gathered}
 \delta_{\rm har}^{\rm phys}
 \bigl(
  \upsilon_{\rm har}^{\rm phys},
  \mathbf p_{\rm op}^{\rm phys};
  \mathfrak P_{\rm har}^{\rm phys}
 \bigr),\\
 \delta_{\rm har}^{\rm phys,op}
 \bigl(
  \upsilon_{\rm har}^{\rm phys};
  \mathbf p_{\rm op}^{\rm phys};
  \mathbf p_{+}^{\rm phys};
  \mathfrak P_{\rm har}^{\rm phys}
 \bigr)
 \end{gathered}
 \right\}>0.
\end{equation}
The reference modulus
\(\delta_{\rm har}^{\rm phys,ref}\) is used once in passing from the
frozen reference metric to every carrier in the coefficient ball.  The
separate common modulus
\(\delta_{\rm har}^{\rm phys,wit}\) is reserved for the subsequent
short Ricci-flow motion from that carrier's common \(+\)-tier to the
operative tier.

\item
Let \(\mathfrak A_{\rm phys}\) be the fixed finite physical reference
atlas on \(\bigcup_a\widetilde W_a^{\rm har}\), with the recorded
scales \(R_a\), and put
\[
 \|K\|_{\mathcal C_{\rm phys}^{14,\alpha}}
 :=
 \max_{1\leq a\leq N_{\rm ext}}
 \|K\|_{C_{R_a}^{14,\alpha}
       (\widetilde W_a^{\rm har};\mathfrak A_{\rm phys})}.
\]
Choose \(0<\varepsilon_{\rm coeff}^{\rm phys}<1\) and
\(\Lambda_{\rm coeff}^{\rm phys}>1\), with
\(\|G_{\rm ref}^{\rm phys}\|_{\mathcal C_{\rm phys}^{14,\alpha}}
 <\Lambda_{\rm coeff}^{\rm phys}\), and define the named physical
coefficient ball
\begin{equation}\label{eq:physical-coefficient-ball}
 \begin{split}
 \mathscr B_{\rm coeff}^{\rm phys}
 (G_{\rm ref}^{\rm phys};
       \varepsilon_{\rm coeff}^{\rm phys},
       \Lambda_{\rm coeff}^{\rm phys})
 :=\bigl\{G_\circ:\;&G_\circ\text{ is a closed }
 C^{14,\alpha}\text{ metric on }\mathcal X,\\
 &\|G_\circ-G_{\rm ref}^{\rm phys}\|_
       {\mathcal C_{\rm phys}^{14,\alpha}}
       <\varepsilon_{\rm coeff}^{\rm phys},\\
 &(\Lambda_{\rm coeff}^{\rm phys})^{-1}
       G_{\rm ref}^{\rm phys}
       <G_\circ<
       \Lambda_{\rm coeff}^{\rm phys}G_{\rm ref}^{\rm phys}
       \quad\text{on }\bigcup_a\widetilde W_a^{\rm har},\\
 &\|G_\circ\|_{\mathcal C_{\rm phys}^{14,\alpha}}
       <\Lambda_{\rm coeff}^{\rm phys}\bigr\}.
 \end{split}
\end{equation}
Let \(C_{\rm emb}^{\rm phys}\) be the fixed finite-atlas
\(C^{14,\alpha}\)-to-\(C^{2,\alpha}\) restriction constant.  The radius
is chosen so that
\[
 C_{\rm emb}^{\rm phys}\varepsilon_{\rm coeff}^{\rm phys}
 <\frac14\delta_{\rm har}^{\rm phys,ref}
\]
holds.  Hence every member of the outer ball carries the common
\(+\)-tier, and therefore the operative tier, at every
\(x\in V_a^5\).

For \(G_\circ\) in this ball, define its package-face distance by
\begin{equation}\label{eq:physical-coefficient-package-distance}
\begin{split}
\operatorname{dist}_{\rm pkg}
\!\left(G_\circ,\partial\mathscr B_{\rm coeff}^{\rm phys}\right)
:=\min\Biggl\{&
1-\frac{\|G_\circ-G_{\rm ref}^{\rm phys}\|_
 {\mathcal C_{\rm phys}^{14,\alpha}}}
 {\varepsilon_{\rm coeff}^{\rm phys}},\\
&\inf_{\bigcup_a\widetilde W_a^{\rm har}}
 \lambda_{\min}\!\left(
 (G_{\rm ref}^{\rm phys})^{-1}G_\circ\right)
 -(\Lambda_{\rm coeff}^{\rm phys})^{-1},\\
&\Lambda_{\rm coeff}^{\rm phys}
 -\sup_{\bigcup_a\widetilde W_a^{\rm har}}
 \lambda_{\max}\!\left(
 (G_{\rm ref}^{\rm phys})^{-1}G_\circ\right),\\
&1-\frac{\|G_\circ\|_{\mathcal C_{\rm phys}^{14,\alpha}}}
 {\Lambda_{\rm coeff}^{\rm phys}}
\Biggr\}.
\end{split}
\end{equation}
Choose \(\mu_{\rm coeff}^{\rm phys}>0\) below one quarter of the
corresponding distance for \(G_{\rm ref}^{\rm phys}\), and define
\begin{equation}\label{eq:physical-coefficient-inner-locus}
\begin{aligned}
&\mathscr B_{\rm coeff,in}^{\rm phys}
\bigl(G_{\rm ref}^{\rm phys};
 \varepsilon_{\rm coeff}^{\rm phys},
 \Lambda_{\rm coeff}^{\rm phys},
 \mu_{\rm coeff}^{\rm phys}\bigr)\\
&\quad:=
\Bigl\{G_\circ\in
\mathscr B_{\rm coeff}^{\rm phys}
(G_{\rm ref}^{\rm phys};
 \varepsilon_{\rm coeff}^{\rm phys},
 \Lambda_{\rm coeff}^{\rm phys}):\\
&\hspace{8em}
 \operatorname{dist}_{\rm pkg}
 (G_\circ,\partial\mathscr B_{\rm coeff}^{\rm phys})
 >\mu_{\rm coeff}^{\rm phys}\Bigr\}.
\end{aligned}
\end{equation}
Every actual entrance carrier and every carrier used in one common
prepared restart package is required to belong to this inner locus.
The reference metric and the outer ball are not recentered when the
actual carrier changes.

After the preceding physical witnesses, scaled separations, and this
coefficient ball have been frozen, choose
\(0<\mu_{\rm RF}<1\) and decrease \(\delta_{\rm RF}>0\), if necessary,
so that every carrier \(G_\circ\) in that coefficient ball and its
actual fixed-marking closed Ricci flow satisfy the following.  Write
\[
 G(\,\cdot\,;G_\circ):
 [0,t_*(G_\circ))\longrightarrow
 \operatorname{Met}_{C}^{14,\alpha}(\mathcal X),
 \qquad t_*(G_\circ)\in(0,\infty],
\]
for its maximal forward lifetime.  Then
\begin{equation}\label{eq:witnessed-physical-quarter-modulus}
 \begin{aligned}
 &\sup_{\substack{0\leq t<t_*(G_\circ)\\
                  t\leq\delta_{\rm RF}R_a^2}}
  \|G(t;G_\circ)-G_\circ\|_
   {C_{R_a}^{2,\alpha}(W_a^{\rm har})}\\
 &\hspace{8em}\leq
 \frac14\delta_{\rm har}^{\rm phys,wit}
 \end{aligned}
\end{equation}
for every \(a\).  Here the two positive moduli entering
\(\delta_{\rm har}^{\rm phys,wit}\) are the basic and
reserve-to-operative moduli furnished by
Lemma~\ref{lem:finite-physical-harmonic-openness}, and
\(\delta_{\rm RF}\) is selected after both.  The number
\(\mu_{\rm RF}\) is the residual time-width fraction used in
\eqref{eq:auxiliary-buffered-time-width}.
\end{enumerate}
The normalized operative and reserve triples, the physical operative,
\(+\)-, and reference triples, all their strict gaps, both families of
quantified witnesses, every recorded Dirichlet preimage, the fixed
reference metric, outer ball, inner-locus margin, and residual
time-width fraction are part of the static package.  The two
lower-stability lemmas cited above are proved from these fixed data and
do not invoke the continuation corollary.
\end{definition}

\begin{corollary}[Closure of the auxiliary geometric bootstrap]
\label{cor:adaptive-auxiliary-closure}
Let $k_0\geq12$ and \(0<\alpha<1\), and consider an admissible
first-exit interval in the sense of
Definition~\ref{def:admissible-first-exit-interval}, on which the
  controlled chart is smooth, the robust three-region \(C^2\) box and
  the radius-independent pre-atlas face
  \(\|h(\tau)\|_{\mathfrak C_{{\rm pre},0}^{2,\alpha}}
    \leq2\delta_{\rm c2}\) hold, with
  \(\|h(\tau_0)\|_{\mathfrak C_{{\rm pre},0}^{2,\alpha}}
    <\delta_{\rm c2}\), and the initial closed metric, graft,
  and relative map belong to one bounded prepared
  \(C^{k_0+2,\alpha}\) package.  Assume in addition the full quantitative
  hypotheses of Proposition~\ref{prop:adaptive-target-tracking} at the
  fixed order \(m_{\rm ad}=13\).  Set
  \[
   m_{\rm aux}:=\max\{m_{\rm ad},k_0\}.
  \]
  The order-\(m_{\rm ad}\) instance is the operative geometric
  instance.  The finite-order upgrade in
  Proposition~\ref{prop:adaptive-target-tracking} will be used at
  \(m_{\rm aux}\) only to obtain the fixed higher-order target
  coefficient bound asserted below.  Thus use the
  adaptive-position package \eqref{eq:adaptive-position-package} and
  the choices fixed in
  Remark~\ref{conv:authoritative-adaptive-order}, require
  \[
   \tau_0\geq\tau_{\rm ad}\geq\tau_\Gamma(m_{\rm ad}),
  \]
  and assume on the whole first-exit interval that
  \begin{equation}\label{eq:auxiliary-adaptive-position-hypotheses}
   \lambda_\tau=-(1+a)\lambda,\qquad
   c_{\rm scl}e^{-\tau}\leq\lambda(\tau)
   \leq C_{\rm scl}e^{-\tau},\qquad
   \int_{\tau_0}^{\tau_1}(|a|+|b|)\,d\tau
   \leq\varepsilon_{\rm ph}.
  \end{equation}
  Here \(C_\lambda=e^{\varepsilon_{\rm ph}}\), with the same
  \(\varepsilon_{\rm ph}\) used below in
  \eqref{eq:auxiliary-buffered-time-width}.  The admissibility
  identities already supply the exact \(\Theta\)-equation
  \eqref{eq:adaptive-Theta-tau} and
  \(S_\tau=\lambda(\tau)\Theta_\tau^*\bar g\).

  Require also that
  \[
   R_{\tau_0}=\varphi_{-\tau_0}\circ\Theta_{\tau_0}
  \]
  is a global proper diffeomorphism and that
  \(R_{\tau_0}\) and \(R_{\tau_0}^{-1}\) satisfy, with the common ceiling
  \(\Lambda_{\rm ad}\), the ordinary \(C^{14}\) bounds on
  \(\Omega_\eta^+\), the scale-normalized \(C^{14}\) bounds on every
  dyadic annulus \(\{L<\bar f<4L\}\), uniformly for \(L\geq\Gamma\),
  and the two-sided radial comparison
  \eqref{eq:adaptive-initial-radial-comparison}.  These are explicit
  hypotheses, not conclusions supplied by the word ``admissible'';
  they imply every lower adaptive-target-tracking order used below.
  In the global first-exit argument they impose no new restriction:
  the strict entrance conditions supply the initial data, while the
  propagated bracket and phase budget are the faces recorded in
  \eqref{eq:master-pre-radius-activation-ledger}.
  The bounded prepared \(C^{k_0+2,\alpha}\) package also determines a
  finite ceiling
  \(\Lambda_{\rm ad}^{\langle k_0\rangle}\) for the ordinary and
  scale-normalized \(C^{m_{\rm aux}+1}\) norms of
  \(R_{\tau_0}^{\pm1}\).  Indeed, the map components of an order
  \(k_0+2\) prepared tuple have \(k_0+3\) derivatives and
  \[
   m_{\rm aux}+1\leq k_0+3.
  \]
  The fixed comparison from the prepared map atlas to the
  compact-collar and dyadic atlases supplies this ceiling.  It is a
  bounded high-regularity input, not a smallness condition, and enters
  neither radius functional nor the fixed phase threshold.

  Assume that the initial normalized
  source and target branches carry the normalized part of one
  witnessed auxiliary harmonic package in the sense of
  Definition~\ref{def:witnessed-auxiliary-harmonic-package}, with
  normalized carriers
  \(g_{\rm src}^{\rm nor}=\bar g+h(\tau_0)\) and
  \(g_{\rm tar}^{\rm nor}=\bar g\), and with
  operative triple
  \[
   (\eta_{\rm har},q_{\rm har},\zeta_{\rm har})
  \]
  and compatibility
  \eqref{eq:auxiliary-harmonic-modulus-compatibility}.  At
  \(r_0=\lambda(\tau_0)^{1/2}\), assume also curvature derivatives
  through order ten and Ricci-defect derivatives through order ten
  bounded in scale-normalized norms.  Assume moreover that the
  physical graft collar and the complementary outer region admit a
  finite buffered cover carrying the physical part of that same
  witnessed package, with the same frozen physical reference metric
  \(G_{\rm ref}^{\rm phys}\) and actual entrance carrier
  \[
   G_{\rm car}:=G(t(\tau_0))
   \in\mathscr B_{\rm coeff,in}^{\rm phys}
   \bigl(G_{\rm ref}^{\rm phys};
    \varepsilon_{\rm coeff}^{\rm phys},
    \Lambda_{\rm coeff}^{\rm phys},
    \mu_{\rm coeff}^{\rm phys}\bigr).
  \]
  In particular,
  \[
   U_a\Subset U_a^+\Subset U_a^{++}\Subset V_a^5
   \Subset W_a^{\rm har}\Subset\widetilde W_a^{\rm har},
   \qquad 1\leq a\leq N_{\rm ext},
  \]
  with scales \(R_a>0\), common parameters
  \((\upsilon_{\rm har}^{\rm phys},
    \eta_{\rm har}^{\rm phys},
    q_{\rm har}^{\rm phys},
    \zeta_{\rm har}^{\rm phys})\), and recorded manifold preimages of
  the outer Dirichlet domains in \(W_a^{\rm har}\).  The hypotheses of
  Lemma~\ref{lem:buffered-local-Ricci-control} hold on every auxiliary
  outer member through order \(k_0\), the reduced
  \(\delta_{\rm RF}\) has the witness-preserving meaning fixed in
  Definition~\ref{def:witnessed-auxiliary-harmonic-package}, and
  \begin{equation}\label{eq:auxiliary-buffered-time-width}
   2C_\lambda\lambda(\tau_0)
   \leq(1-\mu_{\rm RF})\delta_{\rm RF}
   \min_{1\leq a\leq N_{\rm ext}}R_a^2.
  \end{equation}
  Use the already fixed pre-radius low-order graft threshold
  \(K_{\rm gr}\), let
  \(C_{6,K_{\rm gr}}\) be the constant in
  \eqref{eq:graft-compatibility-propagated}, and assume the complete
  strict-improvement margin
  \begin{equation}\label{eq:auxiliary-complete-graft-margin}
   C_{6,K_{\rm gr}}\left(
    \mathfrak d_6(G(\tau_0),S(\tau_0))
    +\frac{\lambda(\tau_0)}{\Gamma}
    +\int_{\tau_0}^{\tau_1}(|a|+|b|)\,d\tau
   \right)
   \leq\frac12K_{\rm gr}.
  \end{equation}
  In the strict-entrance application this is enforced by the common
  positive entrance margin: \(K_{\rm gr}\) and its induced forcing and
  annulus ceilings are fixed first, then
  \eqref{eq:global-compatible-package-radius} fixes \(\Gamma\).
  Afterward choose only the initial graft and phase thresholds, and
  enlarge the
  admissible lower bound for \(\tau_0\) so that
  \(\lambda(\tau_0)/\Gamma\) uses the remaining margin.  Assume that the prepared
  \(C^{k_0+2,\alpha}\) bound supplies finite scale-normalized
  order-\((k_0+2)\) coefficient bounds \(K_{k_0+2}\) for the initial closed
  metric and finite
  \(\mathfrak d_{k_0}(G(\tau_0),S(\tau_0))\leq K_{k_0}\).
  Fix one smooth reference metric on the closed host whose restriction
  to every \(U_a^+\) has fixed scale-\(R_a\) coefficient bounds through
  order \(k_0+3\).  Thus the hypotheses of
  Lemma~\ref{lem:buffered-Ricci-DeTurck-coefficients} hold with
  \(r=k_0+2\).
  Then the closed-flow, graft, curvature, and spatial source and target
coefficient bounds required below hold with constants independent of
the finite normalized endpoint.  In addition, the harmonic faces obey
\begin{align}
 \inf_{\tau_0\leq\tau<\tau_1}
 \mathfrak h_{\rm har}(\tau)
 &\geq\kappa_{\rm har}+\eta_{\rm har},
 \label{eq:auxiliary-normalized-harmonic-improvement}\\
 \inf_{\substack{\tau_0\leq\tau<\tau_1\\
                  1\leq a\leq N_{\rm ext}\\
                  x\in U_a^{++}}}
 \frac{r_{\rm har}(G(t(\tau)),x)}{R_a}
 &\geq
 \upsilon_{\rm har}^{\rm phys}
 +\eta_{\rm har}^{\rm phys}.
 \label{eq:auxiliary-physical-harmonic-improvement}
\end{align}
Both conclusions retain positive residual coefficient and domain
reserves, with radius slack
\(\eta_{\rm har}/2\) and
\(\eta_{\rm har}^{\rm phys}/2\), respectively.

On every finite normalized-time subinterval, both Ricci-defect
hypotheses in Theorem~\ref{thm:adaptive-HMHF-continuation} hold, with
the finite-horizon constant permitted there.  At this stage no
endpoint-independent target-defect constant is
asserted: by \eqref{eq:S-defect} that defect contains the instantaneous
coefficients \(a,b\).  The endpoint-independent defect and coefficient
time-modulus package is obtained from the coarse Gram estimate in
Lemma~\ref{lem:coarse-effective-Gram}.  The low-order hypotheses of
Proposition~\ref{prop:graft-compatibility} hold, and the
graft-compatibility norm remains
 strictly below its bootstrap threshold.  In the fixed finite family
 of exterior Ricci--DeTurck gauges, the closed metric coefficients
 through order \(k_0+2\) and the gauge and inverse-gauge coefficients
 through order \(k_0+1\) remain bounded on the retained cover.  On the
 marked graft chart the transported-marking coefficients through order
 \(k_0+1\) remain bounded as well, and the metric in the fixed marking
 remains bounded through order \(k_0\).  Thus these higher hypotheses cannot
 be an independent first exit face.  These independently localized
 gauges certify one-state invariant and coefficient bounds only; they
 are not identified across overlaps and are not used as the common
 two-state gauge.  That gauge is constructed once on the prepared
 exterior in Lemma~\ref{lem:anchored-exterior-interface}.
\end{corollary}

\begin{proof}
Use the fixed physical enlargement
 $\Omega_\eta^+\Subset\Omega_\eta^{++}$ of the graft annulus.  The
 finite buffered cover may be refined once, without changing its
uniform scale or overlap constants, so that each member contains a
nested chain as in \eqref{eq:Ricci-DeTurck-nested-buffers}.  The
smallest members still cover the physical graft collar and the
complementary noncollapsing outer region.  No member is required to
cover the collapsing deep core, and no largest member is required to
be covered by smaller members.  The finite-interval width estimate
\eqref{eq:finite-interval-physical-width} and
\eqref{eq:auxiliary-buffered-time-width} allow
Lemma~\ref{lem:buffered-local-Ricci-control} to be applied on every
member of the buffered cover for the entire remaining physical
  interval.  It gives, in particular, the curvature derivatives through
  order \(k_0\) required for the coefficient transfer below, with
  constants independent of the finite normalized endpoint.
  Independently of these curvature estimates, first apply the
  fixed order-\(m_{\rm ad}\) instance of
  Proposition~\ref{prop:adaptive-target-tracking}; it supplies the
  endpoint-independent range, radial-comparison, and low-order
  relative-map bounds.  Then apply the finite-order upgrade in
  Proposition~\ref{prop:adaptive-target-tracking} with
  \(q=m_{\rm aux}\) and initial ceiling
  \(\Lambda_{\rm ad}^{\langle k_0\rangle}\).  It gives
  endpoint-independent ordinary and scale-normalized
  \(C^{m_{\rm aux}+1}\) bounds for \(R_\tau^{\pm1}\) and
  \(C^{m_{\rm aux}}\) coefficient bounds for \(S_\tau\) on the graft
  collar and its fixed enlargements.  In particular, the target has the
  required endpoint-independent \(C^{k_0}\) coefficient package.  This
  step uses only the exact relative-map equation and the preceding
  order-\(m_{\rm ad}\) geometric tracking; it does not use the
  Ricci--DeTurck coefficient conclusion below.

  The initial
  \(C^{k_0+2,\alpha}\) package, \(k_0\geq12\), controls the initial
  order-\((k_0+2)\) metric coefficients.  Apply
  Lemma~\ref{lem:buffered-Ricci-DeTurck-coefficients} with
  \(r=k_0+2\) independently on this fixed finite family, using the
  static marking jet \(K_{\iota,k_0+2}^{\rm fix}\) from
  \eqref{eq:fixed-static-marking-jets}.  It propagates the
  \(C^{k_0+2,\alpha}\) coefficient bound up to the initial face in each
  buffered Ricci--DeTurck chart and controls the gauge and inverse
  gauge in \(C^{k_0+1,\alpha}\).  Separately, choose one compact marked
  domain containing \(\overline{\Omega_\eta^{++}}\), cover it by a
  fixed finite atlas, and solve one localized harmonic-map problem on
  that entire domain.  Thus no single coordinate chart is assumed to
  contain the whole graft collar.  In this marked gauge,
  \eqref{eq:transported-marking-definition} controls the transported
  marking in \(C^{k_0+1,\alpha}\).  Local
  Shi--Bernstein estimates~\cite{Shi} give the corresponding
  positive-time curvature formulation.  The covariance identity
  \eqref{eq:buffered-marking-covariance}, together with the
  \(C^{k_0+1,\alpha}\) composition and inverse estimates, therefore
  transfers the metric back to the fixed physical marking in
  \(C^{k_0,\alpha}\).  This is the one-derivative marking buffer that
  would be lost if the DeTurck lemma were used only with \(r=k_0\).
  Combine this fixed-marking source bound with the \(C^{k_0}\) target
  bound obtained above.  Uniform equivalence and the coordinate formula
  for \(\nabla^{S_\tau}\) then bound every term in
  \eqref{eq:graft-distance}, and hence
  \[
   \sup_{\tau_0\leq\tau<\tau_1}
   \mathfrak d_{k_0}(G(\tau),S(\tau))
   \leq
   C_{k_0}\bigl(
    K_{k_0+2},K_{k_0},
    \Lambda_{\rm ad}^{\langle k_0\rangle},
    K_{\iota,k_0+2}^{\rm fix},
    \mathfrak P_{\rm prep}\bigr),
  \]
  with a constant independent of the finite normalized endpoint and
  without asserting that this high norm is small.  These estimates
  concern the actual closed Ricci flow on fixed buffered physical sets
  and do not use the harmonic-map chart.

We now close the two harmonic faces rather than treating them as
consequences of the preceding coefficient bounds.  On the normalized
region, the first-exit and entrance faces give
\begin{equation}\label{eq:auxiliary-scale-one-graph-distance}
\begin{split}
 d_{\rm gr,sc}^{2,\alpha}
  \bigl(\bar g+h(\tau),\bar g+h(\tau_0)\bigr)
 &=
 \|h(\tau)-h(\tau_0)\|_
  {\mathfrak C_{{\rm sc},0}^{2,\alpha}(M)}\\
 &\leq
 C_{\rm har,pre}
 \|h(\tau)-h(\tau_0)\|_
  {\mathfrak C_{{\rm pre},0}^{2,\alpha}(M)}\\
 \leq
 3C_{\rm har,pre}\delta_{\rm c2}
 <
 \delta_{\rm har}
 (\eta_{\rm har},q_{\rm har},\zeta_{\rm har};
  \mathfrak P_{\rm har}^{\rm geom}).
\end{split}
\end{equation}
The equality is the canonical model identification, the first
inequality is \eqref{eq:pre-atlas-to-scale-one-graph}, and the second
uses the \(2\delta_{\rm c2}\) first-exit face and the strict
\(\delta_{\rm c2}\) entrance face.  Apply the cross-time clause of
Lemma~\ref{lem:prepared-harmonic-radius-lower-stability} to the
recorded witnesses at \(\tau_0\).  The exact graph identity
\eqref{eq:harmonic-radius-normalized-graph-identity} removes the maps,
scale, and normalized-time labels; the target branch reduces to the
fixed background branch.  This proves
\eqref{eq:auxiliary-normalized-harmonic-improvement} and leaves the
residual normalized witness package with radius slack
\(\eta_{\rm har}/2\).

On the physical cover,
\eqref{eq:finite-interval-physical-width} and
\eqref{eq:auxiliary-buffered-time-width} keep the entire remaining
flow inside the reduced local Ricci--DeTurck time window.  By the
fixed-marking quarter-modulus condition
\eqref{eq:witnessed-physical-quarter-modulus} in item~\textup{(4)} of
Definition~\ref{def:witnessed-auxiliary-harmonic-package}, the nested-buffer
Ricci--DeTurck estimate has already been transferred back to the actual
metric in the fixed physical marking on \(W_a^{\rm har}\), with norm at
most one quarter of both the basic and reserve-to-operative moduli of
Lemma~\ref{lem:finite-physical-harmonic-openness}.  That lemma is a
deferred static input proved in
Section~\ref{sec:uniform-entrance}; its proof uses only the fixed
finite-cover Dirichlet data and does not invoke this continuation
argument.  Apply it
with
\[
 U_a^{++}\Subset V_a^5\Subset W_a^{\rm har}
 \Subset\widetilde W_a^{\rm har}.
\]
The centers are therefore the same physical points \(x\), rather than
untracked points moved by a local gauge.  The result is precisely
\eqref{eq:auxiliary-physical-harmonic-improvement} for the actual
closed metric \(G(t(\tau))\), together with radius slack
\(\eta_{\rm har}^{\rm phys}/2\) and common residual coefficient and
domain reserves.  Both harmonic conclusions are uniform in the finite
normalized endpoint.

The explicit order-\(m_{\rm ad}\) adaptive-position hypotheses in the
statement permit the specified instance of
Proposition~\ref{prop:adaptive-target-tracking} to be applied without
importing an unnamed first-exit face.  That instance gives every
lower-order spatial bound used in the graft, cutoff, target-defect, and
harmonic-map restart arguments below.  Its finite-order upgrade at
\(m_{\rm aux}\), applied above, is used only for the displayed
order-\(k_0\) source--target coefficient bound; it requires no new
phase or geometric bootstrap face.
  Apply Proposition~\ref{prop:graft-compatibility} only with \(m=6\),
  on the first-exit bootstrap
  \(\mathfrak d_6\leq K_{\rm gr}\), where \(K_{\rm gr}\) is the fixed
  graft threshold.  Its curvature requirement stops at order eight, so
  the preceding package applies; the complete margin
  \eqref{eq:auxiliary-complete-graft-margin}, including the
  \(\lambda(\tau_0)/\Gamma\) term, makes its conclusion at most
  \(K_{\rm gr}/2\), discharging the bootstrap and keeping the strict
  low-order graft margin small.
The target-chart bounds just obtained, together with the fixed radial
cutoff profile, verify
\eqref{eq:pure-graft-cutoff-jet-hypothesis} through order twelve, with
the adaptive ceiling recorded in
Remark~\ref{rem:adaptive-pure-graft-cutoff-jets}.
Apply Proposition~\ref{prop:pure-graft-sharp} with \(m=2\) to obtain
the small \(C^2\) outer forcing used in the barriers, and with
  \(m=10\), \(K_m=C(K_{k_0+2},K_{k_0})\), to obtain a bounded order-ten
continuation defect.  This separates low-order smallness from
high-order boundedness and does not require curvature derivatives
above order \(k_0\).  The identity
\[
 \partial_tS+2\Ric_S
 =\Theta^*(-2a\Ric_{\bar g}-\Lie_U\bar g)
\]
does the same for the target defect on each fixed finite interval.
The spatial target coefficients are endpoint-uniform by
Proposition~\ref{prop:adaptive-target-tracking}.  On each fixed finite
interval the continuity of \(a,b\) gives the finite defect bound required
for harmonic-map restart, but we do not promote that bound to an
endpoint-independent one here.  Hence the source
$\acute G=\eta G+(1-\eta)S$, the target $S$, and both of their defects
have bounded geometry after rescaling by the local parabolic scale
$\lambda^{1/2}$ on that finite interval.

There is no demand for bounds of all spatial orders simultaneously at
the initial face.  For each fixed continuation order \(k_0\), the
prepared \(C^{k_0+2,\alpha}\) theory supplies the finite source and
target input jets used above.  The buffered theory supplies the
order-\((k_0+2)\) DeTurck coefficient and order-\((k_0+1)\) marking
bounds used to recover the fixed-marking order-\(k_0\) coefficient
package, while the finite-order target upgrade propagates the matching
target bound.  The continuation theorem uses only its fixed
order-twelve subpackage, and the source and target Ricci defects are
needed only through order ten.  After any positive physical time,
parabolic smoothing of the gauge-fixed metric variables supplies every
higher order needed for \(C^\infty_{\rm loc}\) convergence of the
normalized metric; no simultaneous all-order bound or derivative gain
for the ODE-carried coordinate maps is asserted here.

To match the preceding finite bounds to the restart theorem explicitly,
\eqref{eq:relative-metric-identity} and the three-region \(C^2\) box
give uniform ellipticity and the scale-normalized \(C^1\) metric
discrepancy at every restart, while admissibility supplies the proper
prepared diffeomorphism class.  The harmonic-radius hypotheses are
\eqref{eq:auxiliary-normalized-harmonic-improvement} and
\eqref{eq:auxiliary-physical-harmonic-improvement}.  The order-twelve
source and target coefficient hypotheses follow from
Lemma~\ref{lem:buffered-Ricci-DeTurck-coefficients}, used with
\(r=k_0+2\), and
Proposition~\ref{prop:adaptive-target-tracking}; the required
order-ten curvature and Ricci-defect bounds follow from
Lemma~\ref{lem:buffered-local-Ricci-control},
Proposition~\ref{prop:pure-graft-sharp} with \(m=10\), and the target
identity \eqref{eq:S-defect}.  These are precisely the finite-order
hypotheses of Proposition~\ref{prop:moving-target-HMHF-restart} with
\(m=4\).  That proposition is posed on complete manifolds, so no
lateral boundary trace or compatibility condition is being invoked.
Ordinary interior parabolic regularity is used only to recover higher
smoothness on shorter positive-time subintervals, not as an additional
restart hypothesis.  The strict improvement of the \(C^2\) box,
\eqref{eq:graft-compatibility-propagated},
\eqref{eq:auxiliary-normalized-harmonic-improvement}, and
\eqref{eq:auxiliary-physical-harmonic-improvement} excludes
simultaneous equality at any analytic, graft, or harmonic auxiliary
threshold.  The usual first-exit argument therefore closes the
auxiliary bootstrap.
\end{proof}

\subsection{The effective geometric columns}

\paragraph{Dependency structure.}
The estimates in this subsection are arranged in the following
dependency order.  We first fix
the direct pre-\(C^6\) input package and estimate the localized
geometric columns without using any later bootstrap conclusion.  We
then obtain the coarse Gram and phase controls, pass to the
scale-adapted column and forcing bounds, and finally close the
pre-radius low-order estimates uniformly in the eventual graft radius.
Only the outputs of that closure are fed into the one-state and
two-state radius thresholds.  Each later constant therefore has an
explicit earlier input, and all constants needed to choose the radius
are fixed independently of that radius.

Recall
\[
 K_\tau(T)=(\Phi^{-1})^*((1-\eta)\Theta^*T).
\]
The direct columns in the adaptive equation are
\begin{align}
 \mathcal Y_{0,\tau}
 &=Y_0-K_\tau(Y_0),\label{eq:effective-column-zero}\\
\mathcal Y_{j,\tau}
 &=\Lie_{\chi_\tau W_j}\bar g
   -K_\tau(\Lie_{\chi_\tau W_j}\bar g),
 \qquad1\leq j\leq8.\label{eq:effective-column-j}
\end{align}

\paragraph{Direct pre-\(C^6\) input package.}
For the one-sided pre-radius versions of the next results, let
\(\mathscr P_{\rm pre}^{(6)}\) denote the following hypotheses on the
given admissible first-exit interval in the sense of
Definition~\ref{def:admissible-first-exit-interval}.  These hypotheses
are independent of all subsequent bootstrap-closure conclusions.
\begin{enumerate}
\item[(P1)] In the fixed core-plus-dyadic pre-radius atlas
      \eqref{eq:pre-radius-dyadic-family}, the source metric
      \(\acute G\) and target metric \(S\) have one scale-normalized
      ellipticity, harmonic-radius, doubling, and spatial
      \(C^{12,\alpha}\) coefficient package.  The physical
      Ricci--DeTurck defect is bounded through order ten.  The fixed
      prepared transition maps and the normalized derivatives of the
      radial and graft cutoffs through order fifteen have one bound.
      Every upper bound and favorable lower bound in this item is
      independent of \(\Gamma\geq\Gamma_{\rm pre}\).  Target-defect
      time regularity is not included here; when needed it is supplied
      separately by Lemma~\ref{lem:coarse-effective-Gram}.
\item[(P2)] The entrance relative marking
      \(R_{\tau_0}=\varphi_{-\tau_0}\circ\Theta_{\tau_0}\) and its
      inverse have the ordinary and scale-normalized
      \(C^{14,\alpha}\) bounds in the fixed pre-radius atlas, and obey
      the two-sided radial comparison
      \eqref{eq:adaptive-initial-radial-comparison}.  The propagated
      relative markings \(R_\tau^{\pm1}\), radial comparison, and
      annulus-tracking bounds retain the same fixed
      \(C^{14,\alpha}\) package.
\item[(P3)] For some fixed
      \(0<\varepsilon_{\rm ph}^{\rm pre}\leq\log2\),
      with the already frozen primitive upper-scale constant
      \[
       C_{\rm ph}^{\rm pre}:=2C_{\rm sc}<\infty
      \]
      independent of \(\Gamma\),
      \begin{equation}\label{eq:direct-pre-C6-phase}
       \lambda_\tau=-(1+a)\lambda,\qquad
       \lambda>0,\qquad
       \lambda(\tau)\leq C_{\rm ph}^{\rm pre}e^{-\tau},\qquad
       \int_{\tau_0}^{\tau_1}(|a|+|b|)\,d\tau
       \leq\varepsilon_{\rm ph}^{\rm pre}.
      \end{equation}
      Thus, directly from the first and last relations,
      \begin{equation}\label{eq:direct-pre-C6-relative-scale}
       e^{-\varepsilon_{\rm ph}^{\rm pre}}e^{-(\tau-s)}
       \leq\frac{\lambda(\tau)}{\lambda(s)}
       \leq
       e^{\varepsilon_{\rm ph}^{\rm pre}}e^{-(\tau-s)}
       \quad(\tau_0\leq s\leq\tau<\tau_1).
      \end{equation}
\item[(P4)] The entrance map and inverse have a common
      \(\operatorname{Map}_{\rm sc}^{6,\alpha}\) bound in the same
      pre-radius atlas.  On the interval,
      \begin{equation}\label{eq:direct-pre-C6-atlas-smallness}
       \sup_{\tau_0\leq\tau<\tau_1}
       \|h(\tau)\|_{\mathfrak C_{{\rm pre},0}^{2,\alpha}}
       \leq\delta_{\rm atl}^{\rm pre},
      \end{equation}
      where \(\delta_{\rm atl}^{\rm pre}>0\) is the output selected in
      the pre-radius branch of
      Lemma~\ref{lem:C6-persistent-source-atlas}.  It depends only on
      the primitive geometric package and the fixed buffer margin, not
      on any effective-column or Kato ceiling.  Properness and degree
       one hold, and the local-invertibility, metric-discrepancy,
       radial-comparison, and raw \(C^2\) faces have fixed positive
       margins.
\end{enumerate}

\begin{lemma}[Persistent paired source--target entry atlas]
\label{lem:C6-persistent-source-atlas}
Let \([\tau_0,\tau_1)\) be an admissible first-exit interval.  Assume
the endpoint-independent physical coefficient, curvature,
harmonic-radius, and ellipticity package in the physical part of
Corollary~\ref{cor:adaptive-auxiliary-closure}; the order-twelve target
spatial package of
Proposition~\ref{prop:adaptive-target-tracking}; the strict
harmonic-map chart, radial-comparison, and metric-discrepancy faces; and
the exact scale equation and accumulated phase budget
\eqref{eq:adaptive-phase-budget}.  There are constants
\[
 \begin{gathered}
 \delta_{\rm src}>0,\qquad
 \vartheta_{\rm src}>0,\qquad
 c_{\rm src}>0,\qquad
 C_{\rm src}<\infty,\\
 C_{\rm clock}<\infty,\qquad
 \Lambda_{\rm src}<\infty,\qquad
 0<\epsilon_{\rm atl}<1,\qquad
 0<c_{\rho,-}\leq L_\rho<1<C_{\rho,+},\\
 0<c_{\rho,{\rm ret}}<1<C_{\rho,{\rm ret}}<\infty,\qquad
 N_{\rm tr},C_{\rm tr}<\infty,
 \end{gathered}
\]
depending only on that package, with the following property.
For the pre-radius invocation, replace the preceding generic spatial
and phase hypotheses by \(\mathscr P_{\rm pre}^{(6)}\).  In this branch
there are, in addition, primitive-dependent constants
\[
 C_{\rm disp}^{\rm pre}<\infty,\qquad
 C_{\rm clock}^{\rm pre}<\infty,\qquad
 \delta_{\rm atl}^{\rm pre}>0.
\]
The last is an output of the lemma, selected before
\(\mathscr P_{\rm pre}^{(6)}\) is imposed with that value.  It may be
decreased to preserve all finitely many primitive metric, map, and
buffer margins, without changing either preceding constant.  The proof
then
uses \eqref{eq:direct-pre-C6-relative-scale}, the upper physical-width
bound obtained from \eqref{eq:direct-pre-C6-phase}, and ratios of chart
radii to the actual positive value of \(\lambda\); no numerical lower
bound for \(e^\tau\lambda(\tau)\) enters.  For the pre-radius
construction, all the displayed upper constants and reciprocals of
positive lower constants
have one bound independent of
\(\Gamma\geq\Gamma_{\rm pre}\).  A fixed-\(\Gamma\)
high-regularity continuation may retain a larger package constant, but
that larger value is not used in the pre-radius construction.

The atlas is selected from a single admissible-radius field.  Put
\begin{equation}\label{eq:C6-prepared-physical-scale}
 r_{{\rm sol},s}(x):=
 \left\{\lambda(s)
 \bigl(1+\bar f(\Theta_s(F_s(x)))\bigr)\right\}^{1/2}.
\end{equation}
There is a canonically regularized function
\(\varrho_s:M\to(0,\infty)\), Lipschitz with respect to
\(\acute G(s)\), such that
\begin{equation}\label{eq:C6-admissible-radius}
 c_{\rho,-}r_{{\rm sol},s}(x)\leq\varrho_s(x)
 \leq C_{\rho,+}r_{{\rm sol},s}(x),
 \qquad
 |\varrho_s(x)-\varrho_s(y)|
 \leq L_\rho\,d_{\acute G(s)}(x,y).
\end{equation}
If \(x_{\mathcal U}\) is the center of a member entered at \(s\), its
radius is
\begin{equation}\label{eq:C6-atlas-radius-selection}
 r_{\mathcal U}=\epsilon_{\rm atl}\varrho_s(x_{\mathcal U}).
\end{equation}
For any active member \(\mathfrak a\), entered at
\(s_{\mathfrak a}\), let \(I_{\mathfrak a}\) denote its active
interval: one prescribed short slab for a near member, and its entire
remaining lifetime for a retained remote member.  The current
admissible radius on its outer buffer obeys
\begin{equation}\label{eq:C6-admissible-radius-persistence}
 c_{\rho,{\rm ret}}\epsilon_{\rm atl}^{-1}r_{\mathfrak a}
 \leq\varrho_\tau(x)\leq
 C_{\rho,{\rm ret}}\epsilon_{\rm atl}^{-1}r_{\mathfrak a},
 \qquad
 x\in\mathcal U_{\mathfrak a}^4,\quad
 \tau\in I_{\mathfrak a}.
\end{equation}

At every refresh time \(s\in[\tau_0,\tau_1)\) there is a uniformly
locally finite six-level buffered paired source--target atlas
\[
 \mathfrak A_s^{(6)}
 =
 \left\{
 \mathcal U^0\Subset\mathcal U^1\Subset
  \mathcal U^2\Subset\mathcal U^3\Subset\mathcal U^4
  \Subset\mathcal U^5;\quad
  \mathcal V^0\Subset\mathcal V^1\Subset
  \mathcal V^2\Subset\mathcal V^3\Subset\mathcal V^4
  \Subset\mathcal V^5;\quad
  r_{\mathcal U},L_{\mathcal U},L_{\mathcal V}
 \right\}
\]
whose \(\mathcal U^0\)-members cover the source and whose
\(\mathcal V^0\)-members cover the target, with uniform overlap; the
corresponding \(\mathcal U^1\)- and \(\mathcal V^1\)-enlargements
have a uniform scale-relative Lebesgue number.  Moreover,
\begin{equation}\label{eq:C6-source-atlas-lower-scale}
 r_{\mathcal U}^{\,2}\geq c_{\rm src}\lambda(s).
\end{equation}
Each member has an explicit prepared scale label:
\[
 L_{\mathcal U}=1
 \quad\hbox{only when \(\mathcal U^5\) lies in the fixed core},\qquad
 L_{\mathcal U}\in\mathscr L_{\rm pre}
 \quad\hbox{otherwise}.
\]
For a dyadic member, \(\mathcal U^5\) lies in a fixed enlargement of
\(A_{L_{\mathcal U}}^{\rm pre}\).  Members whose outer buffers overlap have
comparable radii, and their nonunit dyadic labels are comparable, with
factor at most \(C_{\rm src}\).  The finitely many core--dyadic
interfaces obey the same assertion after increasing \(C_{\rm src}\),
with \(C_{\rm src}\) independent of
 \(\Gamma\geq\Gamma_{\rm pre}\) when the pre-radius normalized package
is fixed.
The \(\mathcal U^5\)-coordinates are harmonic for \(\acute G(s)\), and
the \(\mathcal V^5\)-coordinates are harmonic for \(S(s)\), both at the
common scale \(r_{\mathcal U}\).  In these fixed coordinates the
metrics \(r_{\mathcal U}^{-2}\acute G(s)\) and
\(r_{\mathcal U}^{-2}S(s)\), and their spatial coefficients through
order twelve, have one \(\Lambda_{\rm src}\)-bounded geometry package.
For a Riemannian metric \(g\), a subset \(A\subset M\), and
\(\rho>0\), write
\begin{equation}\label{eq:metric-neighborhood-convention}
 \mathcal N_\rho^g(A)
 :=\{y\in M:\operatorname{dist}_g(y,A)<\rho\}.
\end{equation}
Thus \(\mathcal N_\rho^g(A)\Subset B\) means that
\(\overline{\mathcal N_\rho^g(A)}\) is a compact subset of \(B\).
There is a package number \(\kappa_{\rm buf}>0\) such that, for
\(0\leq j\leq4\),
\begin{equation}\label{eq:C6-paired-entry-containment}
 \mathcal N_{\kappa_{\rm buf}r_{\mathcal U}}^{S(s)}
   \bigl(F_s(\overline{\mathcal U^j})\bigr)
 \Subset\mathcal V^j,\qquad
 \mathcal N_{\kappa_{\rm buf}r_{\mathcal U}}^{S(s)}
   \bigl(\overline{\mathcal V^j}\bigr)
 \Subset F_s(\mathcal U^{j+1}).
\end{equation}
Thus both the forward and inverse coordinate representatives have a
quantitative one-buffer margin.  Moreover, at the entry time,
\(\Theta_s(\mathcal V^5)\) lies in a fixed enlargement of
\(A_{L_{\mathcal V}}^{\rm pre}\), with the core convention
\(L_{\mathcal V}=1\) and otherwise
\(L_{\mathcal V}\in\mathscr L_{\rm pre}\), and
\begin{equation}\label{eq:C6-target-output-scale}
 c_{\rm src}r_{\mathcal U}^{\,2}
 \leq\lambda(s)L_{\mathcal V}
 \leq C_{\rm src}r_{\mathcal U}^{\,2}.
\end{equation}
The chart norms on all six levels are uniformly equivalent to the
core and dyadic norms defining
\(\mathfrak X_{\rm sc}^{q,\alpha}\), \(2\leq q\leq6\).

For a chart entered at time \(s\), use the exact effective clock
\begin{equation}\label{eq:C6-source-effective-clock}
 \vartheta_{\mathcal U}^{\,s}(\tau)
 :=
 \frac1{r_{\mathcal U}^{\,2}}
 \int_s^\tau\lambda(q)\,dq .
\end{equation}
If
\begin{equation}\label{eq:C6-atlas-persistence-clock}
 \vartheta_{\mathcal U}^{\,s}(v)
 \leq\vartheta_{\rm src},
\end{equation}
then the same fixed source and target coordinate buffers remain valid
on \([s,v]\), with the same uniform ellipticity, spatial coefficient,
overlap, and norm-equivalence constants.  Measured with the frozen
entry metric \(S(s)\), they retain the half-margin containments
\begin{equation}\label{eq:C6-paired-persistent-containment}
 \mathcal N_{\kappa_{\rm buf}r_{\mathcal U}/2}^{S(s)}
   \bigl(F_\tau(\overline{\mathcal U^j})\bigr)
 \Subset\mathcal V^j,\qquad
 \mathcal N_{\kappa_{\rm buf}r_{\mathcal U}/2}^{S(s)}
   \bigl(\overline{\mathcal V^j}\bigr)
 \Subset F_\tau(\mathcal U^{j+1})
\end{equation}
for \(0\leq j\leq4\) and \(s\leq\tau\leq v\).  In particular, after
decreasing
\(\delta_{\rm src}\),
\[
 \vartheta_{\mathcal U}^{\,s}
 \bigl(\min\{s+\delta_{\rm src},\tau_1\}\bigr)
 \leq\vartheta_{\rm src}
\]
for every newly refreshed chart.
The scale equation also gives the named terminal-clock estimate
\begin{equation}\label{eq:C6-source-terminal-clock}
 \vartheta_{\mathcal U}^{\,s}(\tau)
 \leq
 C_{\rm clock}\frac{\lambda(s)}{r_{\mathcal U}^{\,2}},
 \qquad s\leq\tau<\tau_1.
\end{equation}

Now fix any \(K_{\rm rem}\geq1\) satisfying
\begin{equation}\label{eq:C6-atlas-regime-admissibility}
 \frac{C_{\rm clock}}{K_{\rm rem}}
 \leq\vartheta_{\rm src}.
\end{equation}
A member freshly entered at time \(s\) is called \(K_{\rm rem}\)-near
when
\[
 r_{\mathcal U}^{\,2}\leq K_{\rm rem}\lambda(s),
\]
and \(K_{\rm rem}\)-remote otherwise.  A near member is used for one
normalized-time slab of length at most \(\delta_{\rm src}\) and is then
refreshed.  At every positive refresh time, retained remote pairs are
carried unchanged, and new centers are added only on
\[
 \mathscr R_s
 :=M\setminus
   \bigcup_{\mathfrak r\ {\rm retained\ remote}}
       \mathcal U_{\mathfrak r}^{0}.
\]
Each newly entered near pair, and each pair entering the remote regime
at that handoff, has a \emph{trace-level star refinement} subordinate
to the preceding near cover: there is a preceding near parent
\(\mathfrak p(\mathfrak a_{\rm new})\) such that
\begin{equation}\label{eq:C6-refresh-trace-star}
  \overline{\mathcal U_{\rm new}^{2}}
   \Subset\mathcal U_{\mathfrak p(\mathfrak a_{\rm new})}^{3},
  \qquad
  \overline{\mathcal V_{\rm new}^{2}}
   \Subset\mathcal V_{\mathfrak p(\mathfrak a_{\rm new})}^{3}.
\end{equation}
The child and parent have comparable radii, source labels, and invariant
target scales.
The source and target transition maps, and their inverses, have one
uniform \(C^{12,\alpha}\) scale-one bound on the trace pair.
Consequently direct restriction and coordinate change transfer the
parent's \(C^{q,\alpha}\)-reset trace to the child's level-two
entrance trace, without a factor depending on the refresh number.  No
subordination of
\(\mathcal U_{\rm new}^{4}\) or \(\mathcal V_{\rm new}^{4}\) is
asserted or used.

A region which is remote at \(\tau_0\), or first becomes remote at the
end of a completed near slab, is covered by a \emph{freshly entered}
subordinate remote pair at that time.  Its clock therefore starts at
that handoff time; the pair is thereafter retained with its source and
target buffers and is not reclassified.  By
\eqref{eq:C6-source-terminal-clock} and
\eqref{eq:C6-atlas-regime-admissibility}, every retained remote pair
satisfies \eqref{eq:C6-atlas-persistence-clock} for its entire
remaining life.  At every time the active near pairs together with
 the retained remote pairs form uniformly locally finite source and
 target level-zero covers with the preceding norm equivalences and
 containment margins.  Across retained and newly entered generations,
 the corresponding level-one enlargements have one uniform source and
 target scale-relative Lebesgue number.  No such Lebesgue-number claim
 is made for the mixed-generation level-zero cover itself.

 More precisely, write
 \(\xi_{\mathfrak a}^{\mathcal U}\) and
 \(\xi_{\mathfrak a}^{\mathcal V}\) for the scale-one source and target
 coordinate maps of an active pair, and set
 \begin{equation}\label{eq:C6-active-neighbor-family}
 \mathscr B_{\mathfrak a}(\tau)
  :=
  \left\{\mathfrak b:\ \mathfrak b\ {\rm is\ active\ at}\ \tau,\ 
  \mathcal U_{\mathfrak b}^0
  \cap\mathcal U_{\mathfrak a}^3\ne\varnothing
  \ \hbox{or}\
  \mathcal V_{\mathfrak b}^0
  \cap\mathcal V_{\mathfrak a}^3\ne\varnothing\right\}.
 \end{equation}
 Then, including when the two members entered in different
 generations,
 \begin{equation}\label{eq:C6-active-cross-generation-packing}
  \#\mathscr B_{\mathfrak a}(\tau)\leq N_{\rm tr},\qquad
  C_{\rm tr}^{-1}r_{\mathfrak a}
  \leq r_{\mathfrak b}\leq
  C_{\rm tr}r_{\mathfrak a}
  \quad(\mathfrak b\in\mathscr B_{\mathfrak a}(\tau)).
 \end{equation}
 Their source labels are comparable by \(C_{\rm tr}\), with the fixed
 core convention at a core--dyadic interface.  The invariant target
 scales, rather than the unweighted moving labels, satisfy
 \[
  C_{\rm tr}^{-1}
  \leq
  \frac{\lambda(s_{\mathfrak b})L_{\mathcal V_{\mathfrak b}}}
       {\lambda(s_{\mathfrak a})L_{\mathcal V_{\mathfrak a}}}
  \leq C_{\rm tr}.
 \]
 On every nonempty overlap, in scale-one coordinates,
 \begin{align}
  &\|\xi_{\mathfrak b}^{\mathcal U}
       \circ(\xi_{\mathfrak a}^{\mathcal U})^{-1}\|_{C^{12,\alpha}}
   +\|\xi_{\mathfrak a}^{\mathcal U}
       \circ(\xi_{\mathfrak b}^{\mathcal U})^{-1}\|_{C^{12,\alpha}}
   \nonumber\\
  &\quad+
   \|\xi_{\mathfrak b}^{\mathcal V}
       \circ(\xi_{\mathfrak a}^{\mathcal V})^{-1}\|_{C^{12,\alpha}}
   +\|\xi_{\mathfrak a}^{\mathcal V}
       \circ(\xi_{\mathfrak b}^{\mathcal V})^{-1}\|_{C^{12,\alpha}}
  \leq C_{\rm tr}.
  \label{eq:C6-active-direct-transition}
 \end{align}
 Each norm in \eqref{eq:C6-active-direct-transition} is restricted to
 the relevant coordinate overlap.  These are direct transitions
 through a common member of the fixed prepared atlas (the pre-radius
 atlas in the one-sided version), not products along a chain of
 refresh generations.
\end{lemma}

\begin{proof}
We first make the variable scale and the packing construction
quantitative.  Fix a small structural number \(\epsilon_0>0\) and a
dimensionless coarse paired-containment margin
\(\mu_{\rm pair}>0\), both determined only by the input package and its
coarse source--target quasi-isometry.  For \(x\in M\), let
\(\widehat\varrho_s(x)\) be the supremum of
\(0<r\leq\epsilon_0 r_{{\rm sol},s}(x)\) for which the following
admissibility
property holds on the \(40r\) source ball: there are source and paired
target harmonic coordinates with the input-package ellipticity and
order-twelve coefficient bounds and one coarse forward-and-inverse
paired containment on the full \(40r\) envelope with normalized margin
at least \(\mu_{\rm pair}\).  No member
of the eventual six-level profile occurs in this definition.  The
scale-normalized harmonic-radius bound, the strict discrepancy face,
and the radial and target tracking bounds give
\begin{equation}\label{eq:C6-raw-admissible-radius}
 c_{\rm adm} r_{{\rm sol},s}(x)\leq\widehat\varrho_s(x)
 \leq\epsilon_0 r_{{\rm sol},s}(x)
\end{equation}
after decreasing \(\epsilon_0\), with \(c_{\rm adm}>0\) uniform.
Notice that
this definition uses only fixed numerical thresholds, so
\(\widehat\varrho_s\) is determined by the geometric data at time
\(s\).

The function \(r_{{\rm sol},s}\) is uniformly Lipschitz in the source
metric.
Indeed, this is the scale-one gradient estimate for
\((1+\bar f)^{1/2}\), transported by \(\Theta_s\circ F_s\) using the
radial, ellipticity, and discrepancy faces.  In particular,
\begin{equation}\label{eq:C6-prepared-scale-temperate}
 r_{{\rm sol},s}(x)
 \leq C_{\rm tem}
 \bigl(r_{{\rm sol},s}(y)+d_{\acute G(s)}(x,y)\bigr),
\end{equation}
where \(C_{\rm tem}\geq1\) is fixed by the input package.  Choose
\begin{equation}\label{eq:C6-radius-regularization-constants}
 0<L_\rho:=\epsilon_\rho
 \leq\frac12\min\{c_{\rm adm},1\},
 \qquad
 c_{\rho,-}:=\frac{\min\{c_{\rm adm},L_\rho\}}{C_{\rm tem}}
 =\frac{L_\rho}{C_{\rm tem}},
 \qquad
 C_{\rho,+}:=\max\{2,\epsilon_0\},
\end{equation}
and take the metric
infimal convolution
\begin{equation}\label{eq:C6-admissible-radius-regularization}
 \varrho_s(x):=
 \inf_{y\in M}
 \left\{\widehat\varrho_s(y)
       +L_\rho d_{\acute G(s)}(x,y)\right\}.
\end{equation}
It is \(L_\rho\)-Lipschitz and is no larger than
\(\widehat\varrho_s\).  Conversely, for every \(y\in M\),
\eqref{eq:C6-raw-admissible-radius} and
\eqref{eq:C6-prepared-scale-temperate} imply
\begin{align*}
 \widehat\varrho_s(y)
 +L_\rho d_{\acute G(s)}(x,y)
 &\geq
 \min\{c_{\rm adm},L_\rho\}
 \bigl(r_{{\rm sol},s}(y)+d_{\acute G(s)}(x,y)\bigr)\\
 &\geq c_{\rho,-}r_{{\rm sol},s}(x).
\end{align*}
Taking the infimum gives the lower bound in
\eqref{eq:C6-admissible-radius}; taking \(y=x\) gives its upper bound.
Thus \(c_{\rho,-}\) is the lower-comparability factor and \(L_\rho\)
is the Lipschitz constant; the two constants are kept distinct.

Before choosing any six-level profile, make the outer envelope
explicit.  For every \((s,x_0)\), choose a half-maximal raw radius
from the admissible set in the definition of
\(\widehat\varrho_s(x_0)\), with
\[
 \frac12\widehat\varrho_s(x_0)
 <r_{\rm env}(s,x_0)\leq\widehat\varrho_s(x_0);
\]
such a choice exists even when the defining supremum is not attained.
The raw \(40r_{\rm env}\) paired chart contains fixed source and target
subenvelopes
\(\mathcal E_{\mathcal U}(s,x_0)\) and
\(\mathcal E_{\mathcal V}(s,x_0)\) of radius
\(\sigma_{\rm env}\varrho_s(x_0)\).  Let
\(\sigma_{\rm pair}>0\) be the fixed fraction supplied by the coarse
paired containment and its quasi-isometry constant, and choose
\begin{equation}\label{eq:C6-envelope-fraction-choice}
 0<\sigma_{\rm env}\leq
 \min\left\{\sigma_{\rm pair},
 \frac{1}{2C_{\rm tem}C_{\rho,+}}\right\}.
\end{equation}
In particular,
\(\sigma_{\rm env}\) is independent of all buffer factors and atlas
radius factors selected below.

The radial- and target-tracking faces, the exact relative-scale law,
\eqref{eq:C6-prepared-scale-temperate}, and the coefficient modulus on
these raw envelopes give package constants
\[
 0<\vartheta_{\rm env},\qquad
 1\leq C_{\rm met}^{\rm env}<\infty,\qquad
 0<\kappa_-\leq1\leq\kappa_+<\infty
\]
with the following property.  If a future member centered at
\((s,x_0)\) has radius \(r\), lies in these subenvelopes, and its
effective clock is at most \(\vartheta_{\rm env}\), then the current
source and target metrics are
\(C_{\rm met}^{\rm env}\)-equivalent to their entry metrics throughout
the subenvelopes and
\begin{equation}\label{eq:C6-preselected-envelope-scale-comparison}
 \kappa_-r_{{\rm sol},s}(x_0)
 \leq r_{{\rm sol},\tau}(x)
 \leq\kappa_+r_{{\rm sol},s}(x_0),
 \qquad x\in\mathcal E_{\mathcal U}(s,x_0).
\end{equation}
Indeed, spatial comparison at \(s\) follows from
\eqref{eq:C6-prepared-scale-temperate} and the fixed envelope fraction:
for \(x\in\mathcal E_{\mathcal U}(s,x_0)\),
\[
 d_{\acute G(s)}(x,x_0)
 \leq\sigma_{\rm env}C_{\rho,+}r_{{\rm sol},s}(x_0)
 \leq\frac{r_{{\rm sol},s}(x_0)}{2C_{\rm tem}}.
\]
Applying \eqref{eq:C6-prepared-scale-temperate} in both directions
therefore gives
\((2C_{\rm tem})^{-1}r_{{\rm sol},s}(x_0)
\leq r_{{\rm sol},s}(x)
\leq C r_{{\rm sol},s}(x_0)\).
The exact relative-scale identity, target radial tracking, and the
source--target map-speed estimate control the time variation; choosing
\(\vartheta_{\rm env}\) below their finitely many strict thresholds
gives the displayed constants.  This argument uses the raw package and
the preselected envelopes only.

Once the level-five source and target members are fitted inside these
envelopes by \eqref{eq:C6-radius-factor-smallness} below, the preceding
already selected constants give, for every active member while its
persistence clock holds,
\begin{equation}\label{eq:C6-retained-soliton-scale-comparison}
 \kappa_-r_{{\rm sol},s_{\mathfrak a}}(x_{\mathfrak a})
 \leq r_{{\rm sol},\tau}(x)
 \leq\kappa_+r_{{\rm sol},s_{\mathfrak a}}(x_{\mathfrak a}),
 \qquad
 x\in\mathcal U_{\mathfrak a}^4,\quad
 \tau\in I_{\mathfrak a}.
\end{equation}
Fix
\begin{equation}\label{eq:C6-retained-radius-constants}
 c_{\rho,{\rm ret}}
 :=\frac{c_{\rho,-}\kappa_-}{C_{\rho,+}},
 \qquad
 C_{\rho,{\rm ret}}
 :=\frac{C_{\rho,+}\kappa_+}{c_{\rho,-}}.
\end{equation}
These constants automatically obey the inequalities stated in the
lemma.

Before fixing the six levels, put
\begin{equation}\label{eq:C6-handoff-radius-ratio-bound}
 C_{\rm hand}:=\frac97
  \max\{C_{\rho,{\rm ret}},c_{\rho,{\rm ret}}^{-1}\}.
\end{equation}
The factor \(9/7\) is the worst ratio allowed by
\(q_{\rm atl}\leq1/8\).  Enlarge \(C_{\rm hand}\), by one fixed
source--target quasi-isometry factor and by
\(C_{\rm met}^{\rm env}\), to \(C_\star\).
Choose concentric source ball factors \(b_j^-<b_j^+\) and target
factors \(\widetilde b_j^-<\widetilde b_j^+\),
\(0\leq j\leq5\), with the usual nesting and, in addition,
\begin{equation}\label{eq:C6-handoff-buffer-profile}
 b_0^+ + C_\star b_2^+ < b_3^-,
 \qquad
 \widetilde b_0^+
       +C_\star\widetilde b_2^+
 <\widetilde b_3^- .
\end{equation}
Thus a member of radius \(r\) contains the corresponding
\(b_j^-r\)-ball and lies in the \(b_j^+r\)-ball, and similarly on the
target.  The finite inequalities
\eqref{eq:C6-handoff-buffer-profile} are achieved by choosing the
dimensionless buffer gaps after \(C_\star\) is fixed.  Fix this
six-level profile once and for all; its now finite outer factors will
be fitted into the preselected envelopes by decreasing
\(\epsilon_{\rm atl}\) below.
Uniform ellipticity, the buffer gaps, and the already selected
envelope constant \(C_{\rm met}^{\rm env}\) determine constants
\(\beta_0>0\) and \(A_*\geq1\), independent of the
two small radius factors, such that \(A_*\) dominates every fixed
enlargement used in the source and target covering arguments and
\(B_{\acute G(s)}(x_{\mathfrak a},\beta_0r_{\mathfrak a})
\subset\mathcal U_{\mathfrak a}^0\) at entry.  The corresponding
current-time inclusions are recorded below in
\eqref{eq:C6-active-radius-geometry}.  No persistence constant is
selected at this stage.

Choose \(\epsilon_{\rm atl}\) first and then
\(0<\epsilon_{\rm wit}<\epsilon_{\rm atl}\), once and for all, so that
\begin{equation}\label{eq:C6-radius-factor-smallness}
 \begin{gathered}
  \epsilon_{\rm atl}
  \max\{b_5^+,\widetilde b_5^+\}
  \leq\sigma_{\rm env},\qquad
  q_{\rm atl}:=
  \frac{A_*L_\rho\epsilon_{\rm atl}}{c_{\rho,{\rm ret}}}
  \leq\frac18,\\
  q_{\rm wit}:=
  \frac{A_*L_\rho\epsilon_{\rm wit}}{c_{\rho,{\rm ret}}}
  \leq\frac18,\qquad
  \frac{2A_*\epsilon_{\rm wit}}
  {c_{\rho,{\rm ret}}(1-q_{\rm wit})}
  \leq\frac{\beta_0}{2}\epsilon_{\rm atl}.
 \end{gathered}
\end{equation}
By \eqref{eq:C6-radius-regularization-constants} and
\eqref{eq:C6-retained-radius-constants},
\(L_\rho/c_{\rho,{\rm ret}}
=C_{\rm tem}C_{\rho,+}/\kappa_-\) is already a fixed package
constant.  Choose the initial parameters small enough also for the
finite buffer-containment inequalities below; neither radius factor is
reselected after the atlas is constructed.

At the first refresh, choose a maximal family of pairwise disjoint
witness balls
\[
 D_{\mathfrak a}
 :=B_{\acute G(s)}
       (x_{\mathfrak a},\epsilon_{\rm wit}
        \varrho_s(x_{\mathfrak a})),
\]
and set \(r_{\mathfrak a}\) by
\eqref{eq:C6-atlas-radius-selection}.  If a candidate witness ball at
\(x\) meets a selected ball centered at \(x_{\mathfrak a}\), then,
with \(q_0:=L_\rho\epsilon_{\rm wit}\leq q_{\rm wit}\),
\[
 d_{\acute G(s)}(x,x_{\mathfrak a})
 \leq\epsilon_{\rm wit}
 \bigl(\varrho_s(x)+\varrho_s(x_{\mathfrak a})\bigr),
 \qquad
 |\varrho_s(x)-\varrho_s(x_{\mathfrak a})|
 \leq q_0
 \bigl(\varrho_s(x)+\varrho_s(x_{\mathfrak a})\bigr).
\]
Consequently
\begin{equation}\label{eq:C6-variable-radius-ratio}
 \frac{1-q_0}{1+q_0}
 \leq\frac{\varrho_s(x)}{\varrho_s(x_{\mathfrak a})}
 \leq\frac{1+q_0}{1-q_0},
 \qquad
 d_{\acute G(s)}(x,x_{\mathfrak a})
 \leq\frac{2\epsilon_{\rm wit}}{1-q_0}
       \varrho_s(x_{\mathfrak a}).
\end{equation}
The last inequality, \eqref{eq:C6-radius-factor-smallness}, and the
definition of \(\beta_0\) put \(x\) in
\(\mathcal U_{\mathfrak a}^0\).  Maximality therefore proves that the
inner members cover \(M\).  More generally, if two \(A\)-enlarged
variable-radius witness balls meet and
\(q_A:=AL_\rho\epsilon_{\rm wit}<1\), the same calculation gives
\begin{equation}\label{eq:C6-enlarged-witness-radius-ratio}
 \frac{1-q_A}{1+q_A}
 \leq\frac{\varrho_s(x)}{\varrho_s(y)}
 \leq\frac{1+q_A}{1-q_A}.
\end{equation}
Every witness-scale enlargement used in the maximal-selection and
disjointness argument is dominated by
\(A_*/c_{\rho,{\rm ret}}\), so \(q_A\leq q_{\rm wit}\).  Outer atlas
members are at the distinct scale
\(\epsilon_{\rm atl}\varrho_s\), not at the witness scale.  If two
entry outer members centered at \(x_{\mathfrak a}\) and
\(x_{\mathfrak b}\) meet, the fixed ball profile instead gives
\[
 d_{\acute G(s)}(x_{\mathfrak a},x_{\mathfrak b})
 \leq
 \frac{A_*\epsilon_{\rm atl}}{c_{\rho,{\rm ret}}}
 \bigl(\varrho_s(x_{\mathfrak a})
       +\varrho_s(x_{\mathfrak b})\bigr).
\]
The \(L_\rho\)-Lipschitz estimate and the first inequality in
\eqref{eq:C6-radius-factor-smallness} therefore give
\begin{equation}\label{eq:C6-entry-outer-radius-ratio}
 \frac{1-q_{\rm atl}}{1+q_{\rm atl}}
 \leq
 \frac{\varrho_s(x_{\mathfrak b})}
      {\varrho_s(x_{\mathfrak a})}
 \leq
 \frac{1+q_{\rm atl}}{1-q_{\rm atl}}.
\end{equation}
Thus witness-scale intersections use \(q_{\rm wit}\), whereas outer
atlas overlaps use \(q_{\rm atl}\).  This proves the variable-radius
Vitali covering and same-generation radius comparability without
identifying \(L_\rho\) with \(c_{\rho,-}\).  Apply the
admissibility property at each center.  The strict
metric-discrepancy and bilipschitz faces make \(F_s\) a uniform
quasi-isometry from \((M,\acute G(s))\) to \((M,S(s))\); hence one may
choose an \(S(s)\)-harmonic target ball about \(F_s(x_{\mathfrak a})\)
at the same scale.  Choosing six fixed levels between the witness and
admissibility radii according to
\eqref{eq:C6-handoff-buffer-profile}, and then choosing
\(\kappa_{\rm buf}>0\) below both \(\mu_{\rm pair}/2\) and the finitely
many normalized gaps between consecutive source and target levels,
gives the displayed buffers,
\eqref{eq:C6-paired-entry-containment}, and uniform source and target
overlap and scale-relative Lebesgue-number constants.  In particular,
the images of the source inner cover lie in the target inner members,
so the latter cover the target.
The atlas- and witness-to-admissibility ratios are both fixed before
any center is selected.
At a selected center, \eqref{eq:C6-admissible-radius} and radial
tracking give
\begin{equation}\label{eq:C6-entry-radius-label-comparison}
 r_{\mathfrak a}^{\,2}
 \asymp
 \lambda(s)
 \bigl(1+\bar f(\Phi_s(x_{\mathfrak a}))\bigr).
\end{equation}
In particular, independently of the labeling convention, the lower
comparison is quantitative:
\begin{equation}\label{eq:C6-source-lower-scale-constant-choice}
 r_{\mathfrak a}^{\,2}
 =\epsilon_{\rm atl}^{\,2}
   \varrho_s(x_{\mathfrak a})^2
 \geq
 \epsilon_{\rm atl}^{\,2}c_{\rho,-}^{\,2}
 r_{{\rm sol},s}(x_{\mathfrak a})^2
 \geq
 \epsilon_{\rm atl}^{\,2}c_{\rho,-}^{\,2}\lambda(s).
\end{equation}
Thus \(c_{\rm src}\) may, and henceforth does, satisfy
\(c_{\rm src}\leq
\epsilon_{\rm atl}^{\,2}c_{\rho,-}^{\,2}\).
The one-sided pre-radius scale package does not permit replacement of
the last factor by a two-sided expression involving the source label.
For a dyadically tagged source member, radial tracking and
\(\lambda(s)e^s\leq C\) give only
\(r_{\mathfrak a}^{\,2}\leq C L_{\mathcal U_{\mathfrak a}}\);
source-label comparability is instead read directly from overlapping
fixed pre-radius annuli.  Under a later two-sided dynamic scale
bracket this upper bound improves to
\(r_{\mathfrak a}^{\,2}\asymp L_{\mathcal U_{\mathfrak a}}\), but that
stronger comparison is not used here.  Under the core convention
\(L_{\mathcal U_{\mathfrak a}}=1\), the same comparison gives
\[
 c\lambda(s)\leq r_{\mathfrak a}^{\,2}\leq C.
\]
Members for which \(\Phi_s(x_{\mathfrak a})\) stays in the compact
soliton core satisfy \(r_{\mathfrak a}^{\,2}\asymp\lambda(s)\); other
fixed-source-core members may have larger radii.  Thus
\(L_{\mathcal U}=1\) serves only as a fixed-atlas label and does not
assert \(r_{\mathfrak a}^{\,2}\asymp1\).
In all cases this proves
\eqref{eq:C6-source-atlas-lower-scale}.  The order-twelve physical
coefficient package, target tracking, and the two-sided radial
comparison identify the resulting source and target chart norms with
the fixed core and dyadic prepared norms.  Tag only the members lying
in the fixed core by \(L_{\mathcal U}=1\).  Tag every other member by an
\(L_{\mathcal U}\in\mathscr L_{\rm pre}\), unique up to one adjacent
dyadic choice, whose fixed enlargement contains its source outer
buffer.  Neighboring members then have comparable labels uniformly in
\(\Gamma\).  Apply the target radial-tracking comparison
to \(\Theta_s(\mathcal V^5)\) and call the resulting core or dyadic
label \(L_{\mathcal V}\).  Since \(S=\lambda\Theta^*\bar g\), the
common target harmonic scale is exactly comparable to
\((\lambda(s)L_{\mathcal V})^{1/2}\); this proves
\eqref{eq:C6-target-output-scale}.  Maximal separation,
quasi-isometry, and neighboring-annulus overlap give the asserted
radius and label comparability and the constant \(C_{\rm src}\).
For the asserted independence from \(\Gamma\), do not regard
\(\{\bar f<4\Gamma\}\) as one expanding compact chart.  Use the fixed
 core \(\{\bar f<4\Gamma_{\rm atl}\}\) from
\eqref{eq:pre-radius-dyadic-family} and cover its complement by the
scale-normalized dyadic AC annuli.  The finitely many
core--annulus transitions are fixed, while
all annulus--annulus transitions are dilates of one bounded AC
package.  The cutoff derivatives satisfy the same normalized bounds
on the single annular scale \(L\simeq\Gamma\).  Hence maximal
separation, overlap, Lebesgue numbers, radius comparability, and
transition norms have one bound for every
 \(\Gamma\geq\Gamma_{\rm pre}\).  Because the pre-radius coefficient
ceiling already carries the source and target spatial coefficients
through order twelve, the displayed
\(\Lambda_{\rm src}\) has the same uniform bound.  Only optional
continuation estimates above the displayed order twelve may use larger
constants depending on the already frozen high-regularity package;
those constants are absent from the direct pre-\(C^6\) package.

In physical time, the width of \([s,v]\) relative to a chart is exactly
\[
 \frac{t(v)-t(s)}{r_{\mathcal U}^{\,2}}
 =
 \vartheta_{\mathcal U}^{\,s}(v).
\]
The buffered Ricci--DeTurck coefficient estimate, the target spatial
tracking estimate, and the fixed harmonic-radius margin therefore
preserve all six fixed source and target coordinate levels whenever
\eqref{eq:C6-atlas-persistence-clock} holds.  It remains to verify that
their map-containment margins persist.  Let
\(C_{\rm disp}^{\rm pre}\) be the uniform primitive-package constant in
the normalized map-speed estimate, and let
\(C_{\rm clock}^{\rm pre}\) be the corresponding constant in the
integrated scale estimate below.  The normalized-time form of
\eqref{eq:F-speed} and the raw \(C^1\) small box give
\[
 \left|(\partial_\tau F_\tau)\circ F_\tau^{-1}\right|_{S(\tau)}
 \leq C_{\rm disp}^{\rm pre}\lambda(\tau)^{1/2}
       |\bar\nabla h(\tau)|_{\bar g}\circ\Theta_\tau
 \leq C_{\rm disp}^{\rm pre}
       \delta_{\rm atl}^{\rm pre}\lambda(\tau)^{1/2}.
\]
The scale equation and phase budget imply
\[
 \int_s^\tau\lambda(q)^{1/2}\,dq
 \leq C_{\rm clock}^{\rm pre}\lambda(s)^{1/2}
 \leq C_{\rm clock}^{\rm pre} r_{\mathcal U}.
\]
The number \(\delta_{\rm atl}^{\rm pre}\) is selected from the
primitive metric, map, and buffer package before any effective-column
 or Kato ceiling.  Choose it so that
\begin{equation}\label{eq:pre-radius-atlas-displacement-choice}
 C_{\rm disp}^{\rm pre}C_{\rm clock}^{\rm pre}
 \delta_{\rm atl}^{\rm pre}
 \leq\frac14\kappa_{\rm buf}.
\end{equation}
Then the displacement is less than
\(\kappa_{\rm buf}r_{\mathcal U}/4\).  Uniform equivalence of \(S(\tau)\)
and \(S(s)\) converts the estimate to the frozen entry distance.  The
same conclusion for \(F_\tau^{-1}\) follows from the quantitative
lower singular-value face.  Hence
\eqref{eq:C6-paired-persistent-containment} holds.  By
\eqref{eq:C6-source-atlas-lower-scale}, scale comparison on a fixed
normalized-time interval gives
\[
 \vartheta_{\mathcal U}^{\,s}(s+\delta)
 \leq C\delta .
\]
This proves the near-slab assertion after fixing
\(\delta_{\rm src}\) and \(\vartheta_{\rm src}\) sufficiently small,
with
\(\vartheta_{\rm src}\leq\vartheta_{\rm env}\).

The exact identity
\[
 \lambda(q)
 =\lambda(s)e^{-(q-s)}
   \exp\!\left(-\int_s^q a(u)\,du\right)
\]
and the phase budget give
\[
 \vartheta_{\mathcal U}^{\,s}(\tau)
 \leq
 C_{\rm clock}\frac{\lambda(s)}{r_{\mathcal U}^{\,2}},
 \qquad \tau\geq s.
\]
This proves \eqref{eq:C6-source-terminal-clock}; for a member classified
as remote at its entry time, the right side is at most
\(C_{\rm clock}/K_{\rm rem}\).
We finish with the point that involves different refresh generations.
Let \(\mathfrak a\) be active.  The fit condition in
\eqref{eq:C6-radius-factor-smallness} places its level-five source and
target pair inside the envelopes selected before the buffer profile.
For a near member, the prescribed short-slab clock is at most
\(\vartheta_{\rm src}\leq\vartheta_{\rm env}\).  For a retained remote
member, the terminal-clock and regime choices give the same inequality
for its entire remaining life.  The preselected envelope bounds
\eqref{eq:C6-preselected-envelope-scale-comparison} therefore give
the entry-to-current metric equivalence and
\eqref{eq:C6-retained-soliton-scale-comparison}; the target radial
tracking estimate on
\(F_\tau(\mathcal U_{\mathfrak a}^4)\Subset
 \mathcal V_{\mathfrak a}^4\) verifies the application of the already
fixed constants.  At entry,
\eqref{eq:C6-admissible-radius} and
\(r_{\mathfrak a}
=\epsilon_{\rm atl}
\varrho_{s_{\mathfrak a}}(x_{\mathfrak a})\) give
\[
 \frac{\epsilon_{\rm atl}^{-1}r_{\mathfrak a}}{C_{\rho,+}}
 \leq r_{{\rm sol},s_{\mathfrak a}}(x_{\mathfrak a})
 \leq
 \frac{\epsilon_{\rm atl}^{-1}r_{\mathfrak a}}{c_{\rho,-}}.
\]
Combining this display, \eqref{eq:C6-retained-soliton-scale-comparison},
and \eqref{eq:C6-admissible-radius} at time \(\tau\) yields
\[
 \frac{c_{\rho,-}\kappa_-}{C_{\rho,+}}\,
 \epsilon_{\rm atl}^{-1}r_{\mathfrak a}
 \leq\varrho_\tau(x)\leq
 \frac{C_{\rho,+}\kappa_+}{c_{\rho,-}}\,
 \epsilon_{\rm atl}^{-1}r_{\mathfrak a}.
\]
This is \eqref{eq:C6-admissible-radius-persistence}, with the explicit
factor \(\epsilon_{\rm atl}^{-1}\) and the constants in
\eqref{eq:C6-retained-radius-constants}.

Keep, as fixed subsets, the witness balls of all retained remote
pairs.  Inductively they are pairwise disjoint.  For every active pair
put \(\varrho_{\mathfrak a}:=\varrho_\tau(x_{\mathfrak a})\).
The fixed ball profile, metric persistence, and
\eqref{eq:C6-admissible-radius-persistence} give, with the constants
fixed before \eqref{eq:C6-radius-factor-smallness},
\begin{align}
 B_{\acute G(\tau)}
 \bigl(x_{\mathfrak a},
       \beta_0\epsilon_{\rm atl}\varrho_{\mathfrak a}\bigr)
 &\subset\mathcal U_{\mathfrak a}^0,
 &
 \mathcal U_{\mathfrak a}^4
 &\subset
 B_{\acute G(\tau)}
 \left(x_{\mathfrak a},
       \frac{A_*\epsilon_{\rm atl}}{c_{\rho,{\rm ret}}}
       \varrho_{\mathfrak a}\right),
 \nonumber\\
 B_{\acute G(\tau)}
 \left(x_{\mathfrak a},
       \frac{\epsilon_{\rm wit}}
            {A_*C_{\rho,{\rm ret}}}
       \varrho_{\mathfrak a}\right)
 &\subset D_{\mathfrak a},
 &
 D_{\mathfrak a}
 &\subset
 B_{\acute G(\tau)}
 \left(x_{\mathfrak a},
       \frac{A_*\epsilon_{\rm wit}}{c_{\rho,{\rm ret}}}
       \varrho_{\mathfrak a}\right).
 \label{eq:C6-active-radius-geometry}
\end{align}
The same inclusions, with room to spare, hold for a newly entered
current witness ball.

These estimates first justify the maximal extension across
generations.  If a current candidate witness ball centered at \(x\)
meets a retained \(D_{\mathfrak a}\), then
\[
 d_{\acute G(\tau)}(x,x_{\mathfrak a})
 \leq
 \frac{A_*\epsilon_{\rm wit}}{c_{\rho,{\rm ret}}}
 \bigl(\varrho_\tau(x)+\varrho_{\mathfrak a}\bigr).
\]
The \(L_\rho\)-Lipschitz estimate and
\eqref{eq:C6-radius-factor-smallness} therefore give
\[
 \frac{1-q_{\rm wit}}{1+q_{\rm wit}}
 \leq\frac{\varrho_\tau(x)}{\varrho_{\mathfrak a}}
 \leq\frac{1+q_{\rm wit}}{1-q_{\rm wit}},
 \qquad
 d_{\acute G(\tau)}(x,x_{\mathfrak a})
 \leq
 \frac{2A_*\epsilon_{\rm wit}}
      {c_{\rho,{\rm ret}}(1-q_{\rm wit})}
 \varrho_{\mathfrak a}.
\]
Hence \(x\in\mathcal U_{\mathfrak a}^0\).  A candidate centered in the
complement of the retained inner buffers is consequently disjoint
from every retained witness set.  Extending the retained sets by a
maximal disjoint family of current witness balls is therefore
legitimate, and the calculation
\eqref{eq:C6-variable-radius-ratio} covers the complementary region.
Thus disjointness and coverage hold across all generations.

It remains to count active neighbors.  First suppose that the source
overlap in \eqref{eq:C6-active-neighbor-family} holds.  The first line
of \eqref{eq:C6-active-radius-geometry} gives
\[
 d_{\acute G(\tau)}(x_{\mathfrak a},x_{\mathfrak b})
 \leq
 \frac{A_*\epsilon_{\rm atl}}{c_{\rho,{\rm ret}}}
 \bigl(\varrho_{\mathfrak a}+\varrho_{\mathfrak b}\bigr).
\]
Thus the \(L_\rho\)-Lipschitz bound and
\eqref{eq:C6-radius-factor-smallness} give the explicit
cross-generation comparison
\begin{equation}\label{eq:C6-cross-generation-radius-ratio}
 \frac{1-q_{\rm atl}}{1+q_{\rm atl}}
 \leq\frac{\varrho_{\mathfrak b}}{\varrho_{\mathfrak a}}
 \leq\frac{1+q_{\rm atl}}{1-q_{\rm atl}}.
\end{equation}
Together with \eqref{eq:C6-admissible-radius-persistence}, this makes
\(r_{\mathfrak a}\) and \(r_{\mathfrak b}\) uniformly comparable.
If instead only the target overlap in
\eqref{eq:C6-active-neighbor-family} holds, the identical calculation
in \(S(\tau)\), followed by the current uniform quasi-isometry
\(F_\tau^{\pm1}\), gives the same comparison after enlarging its fixed
constant.
By the second line of \eqref{eq:C6-active-radius-geometry}, the
pairwise disjoint witness sets contain current balls of radii
\[
 \frac{\epsilon_{\rm wit}}{A_*C_{\rho,{\rm ret}}}
 \varrho_{\mathfrak b}.
\]
For \(\mathfrak b\in\mathscr B_{\mathfrak a}(\tau)\), these radii are
uniformly comparable by
\eqref{eq:C6-cross-generation-radius-ratio}, and all their centers
lie in one fixed multiple of the
\(\epsilon_{\rm atl}\varrho_{\mathfrak a}\)-ball.  The input
scale-one doubling bound (item {\rm(P1)} in the pre-radius version)
therefore gives a uniform upper bound for their number.  Applying the
current quasi-isometry \(F_\tau\) gives the target packing bound, so
the union of the source- and target-overlap subfamilies still has
uniform cardinality after increasing \(N_{\rm tr}\).
Maximality gives the active level-zero source cover, and paired
containment gives the active level-zero target cover.  Since every
closed level-zero member lies inside its level-one member with a fixed
scale-relative gap, the level-one enlargements of these covers have a
common source and target Lebesgue number, even at a retained/new
interface.  This proves the level-one analytic Lebesgue-number
statement and
\eqref{eq:C6-active-cross-generation-packing}, including the source
label and invariant target-scale comparability by radial tracking.

At a positive refresh, apply this maximal-extension construction with
new centers in \(\mathscr R_s\).  At the instant before refresh the
active level-zero members cover \(M\).  Since a point of
\(\mathscr R_s\) lies in no retained remote level-zero member, it lies
in a preceding near level-zero member.  This is a statement about the
center only; we do not infer that the whole new outer buffer lies in
the near-covered residual.
The outgoing near pairs are regarded as active through the handoff
instant \(s\).  Their persistence estimates hold on the closed
terminal face, so the cross-generation radius and transition
comparisons apply simultaneously to an outgoing parent and an incoming
child at \(s\).

We instead use the fixed gaps
\[
 \overline{\mathcal U_{\rm old}^{0}}
 \Subset\mathcal U_{\rm old}^{3},
 \qquad
 \overline{\mathcal V_{\rm old}^{0}}
 \Subset\mathcal V_{\rm old}^{3}.
\]
The cross-generation radius comparison just proved shows that the
radius of a new member centered in
\(\mathcal U_{\rm old}^{0}\) is uniformly comparable with the old
radius, with ratio bounded by
\eqref{eq:C6-handoff-radius-ratio-bound}.  The star-refinement
inequality \eqref{eq:C6-handoff-buffer-profile} therefore places the
child's closed level-two source--target
pair in the level-three pair of any chosen preceding near parent
whose level-zero member contains the child center.  Allowing all such
parents would give a uniformly bounded star by the packing estimate;
choose one and call it
\(\mathfrak p(\mathfrak a_{\rm new})\).  Direct transition through the
fixed prepared atlas gives the uniform \(C^{12,\alpha}\) transition
bounds.  This proves
\eqref{eq:C6-refresh-trace-star}.

Notice that \(\mathcal U_{\rm new}^{4}\) may cross a retained remote
inner buffer.  That harmless possibility is precisely why the trace
claim is made only on the new level-two pair.  Retained remote members
are carried unchanged and never serve as refresh-trace parents.
Restriction and one direct coordinate change now transfer the old
level-three reset trace to the new level-two
entrance trace with one fixed constant, rather than a product indexed
by refresh generations.

For completeness, the transition estimate is direct.  If active
pairs \(\mathfrak a,\mathfrak b\) overlap, the packing conclusion puts
both source outer members in a uniformly bounded collection of fixed
prepared charts having comparable labels.  The entry harmonic
coordinates of each member have a uniform \(C^{12,\alpha}\)
transition to one common fixed prepared chart.  Composing these two
maps gives the first line of
\eqref{eq:C6-active-direct-transition}; it is a two-map composition,
not a composition through the intervening refreshes.  For the target,
use the two entrance relative markings
\(R_{s_{\mathfrak a}}\) and \(R_{s_{\mathfrak b}}\) to pass directly
to the corresponding fixed pre-radius charts.  Their scale-normalized
\(C^{14,\alpha}\) bounds, the exact soliton dilation, and the invariant
target-scale comparison above give a common prepared chart with
uniform constants.  Composing back gives the second line and the
inverse estimates.  Thus no product of target transitions over
intervening generations occurs.  A uniformly bounded partition of
unity in the common fixed charts has the asserted scale-one
derivative bounds.

A new remote pair at a near-to-remote handoff is selected as part of
this subordinate maximal extension, so its entry clock begins at the
handoff and does not include the preceding near slab.  A retained pair
is never refreshed, while every active near pair belongs only to the
current slab.  The disjoint-witness packing and maximality prove the
final active source--target cover assertion.  The resulting constants
\(c_{\rm src},C_{\rm clock},N_{\rm tr},C_{\rm tr}\) are chosen only
after
\[
 c_{\rho,-},\quad C_{\rho,+},\quad L_\rho,\quad
 c_{\rho,{\rm ret}},\quad C_{\rho,{\rm ret}},\quad
 \epsilon_{\rm atl},\quad\epsilon_{\rm wit}
\]
have been fixed, and depend on them only through the displayed package
bounds.  This completes the radius and packing constant choices.
\end{proof}

\begin{lemma}[Positive-clock buffered interior reset]
\label{lem:C6-positive-clock-reset}
Fix \(2\leq q\leq6\), \(0<\alpha<1\), and four scale-one buffers
\[
 \mathcal U^0\Subset\mathcal U^1\Subset
 \mathcal U^2\Subset\mathcal U^3
\]
with one bounded-geometry and separation package.  Consider on
\(\mathcal U^3\times[0,\Theta]\) a scalar-principal-symbol bundle
equation
\[
 \partial_\vartheta X
 -A^{ij}\nabla_i\nabla_jX-B^i\nabla_iX-CX=\mathscr G .
\]
Assume uniform ellipticity, the spatial order-\((q-2,\alpha)\)
coefficient bounds, and the common parabolic time-oscillation modulus
of Lemma~\ref{lem:weighted-prepared-Schauder}.  If
\(\Theta\geq2\vartheta_*>0\), then
\begin{equation}\label{eq:C6-positive-clock-reset}
 \sup_{\vartheta_*\leq\vartheta\leq\Theta}
 \|X(\vartheta)\|_{C^{q,\alpha}(\mathcal U^2)}
 \leq C_{\vartheta_*}\left(
  \|X\|_{L^\infty C^{q-1,\alpha}
                (\mathcal U^3\times[0,\Theta])}
  +\|\mathscr G\|_{L^\infty C^{q-2,\alpha}
                (\mathcal U^3\times[0,\Theta])}
 \right).
\end{equation}
For \(q=2\), the first norm on the right is the \(C^{1,\alpha}\) norm.
The estimate is an interior reset: no \(C^{q,\alpha}\) norm of the
initial trace occurs.
\end{lemma}

\begin{proof}
Let \(\delta_0\) be the short-interval length supplied by
Lemma~\ref{lem:weighted-prepared-Schauder} for this coefficient
package, and put
\[
 d_*:=\min\left\{\frac{\vartheta_*}{4},
                  \frac{\delta_0}{4}\right\}.
\]
Fix an arbitrary output time
\(\bar\vartheta\in[\vartheta_*,\Theta]\) and use the backward window
\[
 I_{\bar\vartheta}
 :=[\bar\vartheta-2d_*,\bar\vartheta]\subset[0,\Theta].
\]
Choose a cutoff \(\kappa_{\bar\vartheta}\) which vanishes at the left
face of this window and equals one on
\([\bar\vartheta-d_*,\bar\vartheta]\).  Choose also a spatial cutoff
\(\zeta\), equal to one on \(\mathcal U^2\) and supported in
\(\mathcal U^3\).  The section
\(Y=\kappa_{\bar\vartheta}\zeta X\) has zero trace at the left face.
Writing
\[
 D:=A^{ij}\nabla_i\nabla_j+B^i\nabla_i,\qquad
 [D,\zeta]X:=D(\zeta X)-\zeta DX,
\]
its exact source is
\[
 \kappa_{\bar\vartheta}\zeta\mathscr G
 +\dot\kappa_{\bar\vartheta}\zeta X
 -\kappa_{\bar\vartheta}[D,\zeta]X .
\]
Thus the forcing carries the required time cutoff, and the sign agrees
with the displayed commutator convention.  The commutator contains at
most one spatial derivative of \(X\), so its
\(C^{q-2,\alpha}\) norm is controlled by the
\(C^{q-1,\alpha}\) norm in
\eqref{eq:C6-positive-clock-reset}.

The window length is \(2d_*\leq\delta_0\).  Apply the zero-trace local
form of Lemma~\ref{lem:weighted-prepared-Schauder} on this single
buffered window and evaluate at \(\bar\vartheta\).  The cutoff
derivatives cost only powers of the fixed number \(d_*>0\); in
particular the resulting constant is independent of
\(\bar\vartheta\) and of \(\Theta\).  Taking the supremum over all
\(\bar\vartheta\in[\vartheta_*,\Theta]\) proves the stated estimate
without iterating a short-time constant.
\end{proof}

Let \(\Lambda_{\rm C6}\) dominate the ellipticity constants and all
spatial operator-coefficient norms through order six furnished by the
preceding physical and target packages and by the coarse pointwise
modulation bound in Lemma~\ref{lem:coarse-effective-Gram} below.  The
pre-radius construction gives
\(\Lambda_{\rm C6}\leq\Lambda_{\rm C6}^{\rm pre}\) on the low-order
first-exit box, uniformly in \(\Gamma\).  The
proof of that lemma shows that this spatial bound is fixed independently
of \(K_{\rm rem}\); only its optional near-chart time modulus depends on
the chosen regime.  For \(2\leq q\leq6\), let
\(\Theta_{\rm rem}(q)\) denote the short-clock threshold furnished by
Lemma~\ref{lem:short-effective-clock-propagation} when that lemma is
invoked with the explicit triple
\((q,\alpha,\Lambda_{\rm C6})\), and set
\[
 \Theta_{\rm rem}
 :=
 \min\left\{
  \vartheta_{\rm src},
  \min_{2\leq q\leq6}\Theta_{\rm rem}(q)
 \right\}.
\]
Fix once a regime constant \(K_{\rm rem}\gg1\), depending only on the
prepared coefficient package, so large that
\begin{equation}\label{eq:C6-remote-regime-threshold}
 \frac{C_{\rm clock}}{K_{\rm rem}}\leq\Theta_{\rm rem}.
\end{equation}
This single choice supplies both
\eqref{eq:C6-atlas-regime-admissibility} and every short-clock
threshold used below.

\begin{lemma}[Endpoint-uniform coarse Gram and target-time package]
\label{lem:coarse-effective-Gram}
Let \([\tau_0,\tau_1)\) be an admissible first-exit interval on which
the coarse global \(C^2\) and ellipticity box and the strict
harmonic-map metric-discrepancy, diffeomorphism, support-separation,
and annulus-tracking faces hold.  In the generic version assume the
scale and phase comparisons \eqref{eq:adaptive-phase-budget} and the
endpoint-independent closed-flow, graft, curvature, and spatial
coefficient package supplied by the physical part of
Corollary~\ref{cor:adaptive-auxiliary-closure}.  In the pre-radius
version assume directly items {\rm(P1)}--{\rm(P3)} of
\(\mathscr P_{\rm pre}^{(6)}\) instead.  In particular, in either
version the fixed nested cutoffs
and support margin give
\begin{equation}\label{eq:coarse-Gram-support-separation}
 \Phi_\tau(\supp(1-\eta))\cap\supp\rho_\tau=\varnothing,
 \qquad
 \chi_\tau\equiv1\quad\hbox{on a neighborhood of }\supp\rho_\tau .
\end{equation}
In the pre-radius version, the Gram estimate uses the exact support
separation and the actual positive scale, while the defect-time
package uses only
\eqref{eq:direct-pre-C6-relative-scale}; no lower bound for
\(e^\tau\lambda(\tau)\) is used.  No spatial bound for \(F^{\pm1}\)
above scale-normalized order one is assumed.

Put
\[
 \mathscr C_{0,\tau}
 =\mathcal Y_{0,\tau}+\widetilde\B_{0,\tau}h,
 \qquad
 \mathscr C_{j,\tau}
 =\mathcal Y_{j,\tau}+\Lie_{\chi_\tau W_j}h,
 \quad1\leq j\leq8,
\]
and
\[
 M^{\rm low}_{\mu j}(\tau)
 :=\ip{\rho_\tau\mathscr C_{j,\tau}}{Z_\mu}.
\]
If the coarse \(C^2\) constant is sufficiently small and \(\tau_0\)
is sufficiently large, then
\begin{equation}\label{eq:coarse-effective-Gram}
 \left\|M^{\rm low}(\tau)
       -\bigl(\ip{Y_j}{Z_\mu}\bigr)_{\mu,j=0}^8\right\|
 \leq C\|h(\tau)\|_{C^1}+Ce^{-ce^\tau},
 \qquad
 \|(M^{\rm low})^{-1}\|\leq C,
\end{equation}
with constants independent of \(\tau_1\).  Consequently the exact
slice Gram system gives the coarse pointwise estimate
\begin{equation}\label{eq:coarse-pointwise-modulation}
 q(\tau):=|a(\tau)|+|b(\tau)|
 \leq C\left(
   \|\rho_\tau h(\tau)\|_{H^1_\nu}^{2}
   +\|\rho_\tau\mathcal E_{\rm gr}(\tau)\|_{H^{-1}_\nu}
   +e^{-ce^\tau}\right)
 \leq C_{\rm coarse},
\end{equation}
where \(C_{\rm coarse}\) depends only on the fixed bootstrap package,
not on \(\tau_1\).  In particular, the source and target Ricci defects
 through order ten have one common endpoint-independent bound.  On the
 near-parabolic charts
 \(r_{\mathcal U}^2\leq K_{\rm rem}\lambda\), the rescaled source and target
coefficients through order twelve also have a common
source-adapted-time oscillation modulus.  On the remote charts the same
spatial package is endpoint-independent; no uniform annular-time
modulus is asserted or needed there.  This is the target-time and
defect package used in
Lemma~\ref{lem:finite-HMHF-C6-bridge}.
\end{lemma}

\begin{proof}
The relative metric identity \eqref{eq:relative-metric-identity}, the
strict metric-discrepancy face, and uniform source and target
ellipticity give scale-normalized order-one bounds for \(F^{\pm1}\).
In the generic branch,
Proposition~\ref{prop:adaptive-target-tracking} supplies the
corresponding bounds for \(R_\tau^{\pm1}\); in the pre-radius branch
they are item {\rm(P2)} of \(\mathscr P_{\rm pre}^{(6)}\), proved
directly from the relative-marking ODE in the application below.
Together with
\[
 J_\tau=\Theta_\tau\circ\Phi_\tau^{-1}
 =\varphi_\tau\circ R_\tau\circ F_\tau^{-1}
       \circ R_\tau^{-1}\circ\varphi_{-\tau},
\]
these facts would give the required scale-normalized order-one bounds
for \(J_\tau\) in either branch.  The Gram estimate itself does not use
even this consequence: the first identity in
\eqref{eq:coarse-Gram-support-separation} gives
\(\rho_\tau K_\tau(T)=0\) exactly.  The second gives
\(\rho_\tau\Lie_{\chi_\tau W_j}\bar g=\rho_\tau Y_j\).
The omitted tail \((1-\rho_\tau)Y_j\) is
Gaussian-superexponentially small.  Hence, using no map derivative
above order one,
\[
 \|\rho_\tau\mathcal Y_{j,\tau}-Y_j\|_{L^2_\nu}
 \leq Ce^{-ce^\tau},\qquad0\leq j\leq8.
\]
The terms \(\widetilde\B_{0,\tau}h\) and
\(\Lie_{\chi_\tau W_j}h\) perturb the fixed background Gram matrix by
\(O(\|h\|_{C^1})\), by the same weighted integration by parts as in
Lemma~\ref{lem:B-energy}.  This proves
\eqref{eq:coarse-effective-Gram}.

Differentiate the exact slice and use this inverse matrix.  The linear
spectral term vanishes, the nonlinear term is quadratic, and the
moving-cutoff term is Gaussian-superexponentially small.  The modal
pure-graft forcing is bounded by duality.  This gives the first
inequality in \eqref{eq:coarse-pointwise-modulation}, exactly as in
\eqref{eq:receding-velocity}.  The global coarse \(C^1\) box and the
finite Gaussian mass control
\(\|\rho_\tau h\|_{H^1_\nu}\).  The physical graft estimate and the
uniform \(C^1\) diffeomorphism and ellipticity faces control the
displayed forcing term by the unconditional normalized estimate
\eqref{eq:pure-graft-normalized}, without using transported support.
The second inequality follows with an endpoint-independent constant.
This bound is only a coarse coefficient-regularity estimate; it is not
the sharp quadratic velocity decay proved after the three-region
argument.

Finally,
\[
 \partial_tS+2\Ric_S
 =\Theta^*(-2a\Ric_{\bar g}-\Lie_U\bar g),
\]
and the differentiated target-flow equations bound every required
spatial derivative by \(C_m(1+q)\).  In the generic branch this is the
spatial part of Proposition~\ref{prop:adaptive-target-tracking}.  In
the pre-radius branch, item {\rm(P2)} supplies the scale-normalized
\(R_\tau^{\pm1}\)-jets, item {\rm(P1)} supplies the target coefficient
and cutoff jets, and item {\rm(P3)} supplies the exact scale and phase
relations.  Thus this differentiation uses no numerical lower bound
for \(e^\tau\lambda(\tau)\).  Consequently
\eqref{eq:coarse-pointwise-modulation} gives an endpoint-independent
target-defect bound.  Equations
\eqref{eq:relative-target-flow} and \eqref{eq:S-tau-exact} give total
coefficient oscillation at most \(C\delta\) on every normalized-time
 interval of length \(\delta\).  On charts with
 \(r_{\mathcal U}^2\leq K_{\rm rem}\lambda\), the source-adapted clock
 satisfies
 \[
  \frac{d}{d\tau}\vartheta_{\mathcal U}^{\,s}(\tau)
  =
  \frac{\lambda(\tau)}{r_{\mathcal U}^{\,2}}
  \geq c(K_{\rm rem})>0;
 \]
 hence the same estimate gives
a common modulus in \(\vartheta_{\mathcal U}\).  The closed
Ricci--DeTurck equation, the physical coefficient package, and the
 pure-graft identity give the corresponding source modulus.  On remote
 charts \(r_{\mathcal U}^2>K_{\rm rem}\lambda\), only the endpoint-independent
spatial coefficient bound is retained; the sixth-order bridge uses
same-order propagation over their finite total local clock and does
not freeze coefficients there.  Thus no finite-horizon
\(\Lambda_T\) is promoted in either regime.
\end{proof}

For a persistent paired member, an integer \(q\geq1\), and
\(0\leq j\leq4\), write
\begin{equation}\label{eq:typed-local-map-norm}
 \|F\|_{C_{\rm map}^{q,\alpha}
              (\mathcal U^j,\mathcal V^j)}
 :=
 \|F|_{\mathcal U^j}\|_
      {C^{q,\alpha}(\mathcal U^j\to\mathcal V^j)}
 +\|F^{-1}\|_{C^{q,\alpha}
       (\mathcal V^j\to\mathcal U^{j+1})} .
\end{equation}
The two terms mean the representatives in the fixed entry source and
target coordinates, restricted to the indicated buffers; they are
well-defined for the whole active lifespan by
\eqref{eq:C6-paired-persistent-containment}.  On a right-translated
prepared chart they
are equivalently the norms of the corresponding local displacement
vectors for \(F\) and \(F^{-1}\).  Thus, for the radially admissible
maps under consideration,
\(\operatorname{Map}_{\rm sc}^{q,\alpha}\) is precisely the global
map-and-inverse norm defined above, obtained here by taking the
corresponding fixed-pair buffered core and dyadic suprema.  No moving
domain \(F(\mathcal U)\) occurs in this local representation.  This
notation is used below so that a vector-field norm is never applied
directly to a manifold-valued map.

\begin{lemma}[Finite-horizon uniform sixth-order harmonic-map bridge]
\label{lem:finite-HMHF-C6-bridge}
Let \(k_0\geq12\), \(0<\alpha<1\), and let
\([\tau_0,\tau_1)\) be an admissible first-exit interval with finite
endpoint.  In
the generic version, assume the endpoint-independent closed-flow and
spatial bounds in the physical part of
Corollary~\ref{cor:adaptive-auxiliary-closure}, the target spatial
bounds of Proposition~\ref{prop:adaptive-target-tracking}, and
\eqref{eq:adaptive-phase-budget}.  In the one-sided pre-radius version,
replace precisely those hypotheses by the direct package
\(\mathscr P_{\rm pre}^{(6)}\) above.  In either version assume the
two-regime target-defect and coefficient package of
Lemma~\ref{lem:coarse-effective-Gram}, and the strict
metric-discrepancy and diffeomorphism faces of
Theorem~\ref{thm:adaptive-HMHF-continuation}.  If the entrance values
of \(F\) have a common
\(\operatorname{Map}_{\rm sc}^{6,\alpha}\) bound, then
\begin{equation}\label{eq:finite-HMHF-C6-bridge}
 \sup_{\tau_0\leq\tau<\tau_1}
 \|F(\tau)\|_{\operatorname{Map}_{\rm sc}^{6,\alpha}}
 \leq C,
\end{equation}
where \(C\) depends only on the entrance and coefficient packages and
is independent of \(\tau_1\).  More precisely, \(C\) is a monotone
function of the displayed ellipticity, coefficient, buffer,
effective-clock, and incoming-trace bounds.  For the one-sided
pre-radius version, the proof uses
\eqref{eq:direct-pre-C6-relative-scale} and never a positive lower bound
for \(e^\tau\lambda(\tau)\).
\end{lemma}

\begin{proof}
Use the entry atlas and refresh schedule of
Lemma~\ref{lem:C6-persistent-source-atlas}.  On each active near chart,
that is, on a chart satisfying
\[
 r_{\mathcal U}^{\,2}
 \leq K_{\rm rem}\lambda(s)
\]
at its current entry time \(s\), take a near slab of length
\(\delta_{\rm src}\), or its truncated terminal version.  Decrease
\(\delta_{\rm src}\) once so that the upper local clock of every such
slab is at most \(\Theta_{\rm rem}\).  Scale comparison over the slab
also gives
\[
 r_{\mathcal U}^{\,2}\leq C K_{\rm rem}\lambda(\tau)
 \quad\text{and}\quad
 \vartheta_{\mathcal U}^{\,s}(s+\delta_{\rm src})
 \geq8\vartheta_*
\]
for one \(\vartheta_*>0\) depending only on the fixed package.  Thus a
full near slab has both a uniformly bounded total effective clock and a
uniformly positive terminal effective-clock reserve.  The factor eight
allows the five reset orders below to be staggered rather than
incorrectly using a newly smoothed lower norm before it is available.

We first justify the \(C^{1,\alpha}\) base used by the reset; it does
not follow from bilipschitz control alone.  Put
\[
 \mathfrak q:=(F^{-1})^*\acute G.
\]
Then
\[
 \mathfrak q-S=\lambda\Theta^*h,\qquad
 F^*\mathfrak q=\acute G.
\]
On a near slab,
\[
 c_{\rm src}\lambda(s)
 \leq r_{\mathcal U}^{\,2}
 \leq C K_{\rm rem}\lambda(\tau),
\]
with all three scales comparable throughout the slab.  The global raw
\(C^2\) box for \(h\), the order-twelve target package, and the
physical source package therefore give endpoint-uniform scale-one
\(C^{1,\alpha}\) coefficient bounds for \(\mathfrak q\) and
\(\acute G\) on the auxiliary outer buffers; the \(C^\alpha\) bound for
the first derivatives follows from the uniform scale-one second
derivative bound.  Since
\[
 F:(M,\acute G)\longrightarrow(M,\mathfrak q)
\]
is an isometry, choose the target coordinates harmonic for
\(\mathfrak q\).  Such buffered charts have one uniform radius and
\(C^{1,\alpha}\) coefficient package by the preceding bounds.  In the
source \(\acute G\)-harmonic coordinates the components of \(F\)
therefore satisfy the uniformly elliptic harmonic-coordinate system.
The transition from these \(\mathfrak q\)-harmonic coordinates to the
fixed \(S\)-harmonic target atlas has a uniform \(C^{2,\alpha}\) bound,
by the same buffered elliptic coordinate estimate applied to
\(\mathfrak q-S=\lambda\Theta^*h\).  The order-one bound from the strict
metric-discrepancy/diffeomorphism face and buffered interior elliptic
estimates consequently give
\begin{equation}\label{eq:C6-near-C1alpha-base}
 \sup_{\text{near slabs}}
 \|F\|_{C_{\rm map}^{1,\alpha}
              (\mathcal U^4,\mathcal V^4)}
 \leq C.
\end{equation}
The same argument applied to the inverse gives the second term in the
typed norm.  Here the level-five pair is the unused elliptic buffer
for the displayed level-four estimate.  This derivation uses no map
jet of order two or higher and is independent of the finite endpoint.

Use the paired buffers in the order
\[
 \mathcal U^0\Subset\mathcal U^1\Subset
 \mathcal U^2\Subset\mathcal U^3\Subset\mathcal U^4,
 \qquad
 \mathcal V^0\Subset\mathcal V^1\Subset
 \mathcal V^2\Subset\mathcal V^3\Subset\mathcal V^4.
\]
At \(\tau_0\), the incoming fixed-pair
\(C^{q,\alpha}(\mathcal U^2,\mathcal V^2)\) traces are the prepared
traces.  At every later refresh they are supplied by the paired
trace-level star refinement
\eqref{eq:C6-refresh-trace-star} from the preceding level-three reset
traces.  Apply Lemma~\ref{lem:short-effective-clock-propagation}, with
an auxiliary intermediate buffer between \(\mathcal U^1\) and
\(\mathcal U^2\), successively for \(q=2,\ldots,6\).  It controls the
whole slab on \(\mathcal U^1\), including a terminally truncated slab,
from that single level-two incoming trace.

On every full near slab, apply
Lemma~\ref{lem:C6-positive-clock-reset} on the final positive-clock
portion as follows.  The half-margin containment
\eqref{eq:C6-paired-persistent-containment}, the lower
singular-value face, and the uniform \(C^1\) bound allow one fixed
auxiliary source buffer \(\mathcal Z^3\) such that
\[
 \overline{\mathcal U^3}
 \cup
 \overline{\bigcup_{\tau\ {\rm in\ the\ terminal\ reset\ tail}}
                  F_\tau^{-1}(\mathcal V^3)}
 \Subset\mathcal Z^3\Subset\mathcal U^4 .
\]
Choose once a nested reset chain
\[
 \mathcal Z^3\Subset\mathcal W_5\Subset\mathcal W_4
 \Subset\mathcal W_3\Subset\mathcal W_2
 \Subset\mathcal U^4
\]
and its paired target chain, with the additional auxiliary buffers
required by that lemma.  At \(q=2\), use
\eqref{eq:C6-near-C1alpha-base} on \(\mathcal U^4\) and reset onto
\(\mathcal W_2\).  This bound is available after one
\(\vartheta_*\)-clock unit.  Apply the order-three reset only on the
remaining tail, use that \(C^{2,\alpha}(\mathcal W_2)\) bound there,
and reset onto \(\mathcal W_3\).  Inductively, start the order-\(q\)
reset on the tail beginning at clock
\((q-2)\vartheta_*\) and reset onto \(\mathcal W_q\), with the
convention \(\mathcal W_6:=\mathcal Z^3\).  Every such tail has length
at least \(2\vartheta_*\) by the eight-unit reserve.  At every step the
order-\((q-2)\) source contains only already controlled lower map jets
and the fixed order-twelve geometric coefficients.  Hence at the
terminal face of every full near slab,
\begin{equation}\label{eq:C6-near-reset-trace}
 \|F\|_{C_{\rm map}^{q,\alpha}
              (\mathcal U^3,\mathcal V^3)}
 \leq C_q,\qquad2\leq q\leq6,
\end{equation}
where \(C_q\) is independent of the incoming
\(C^{q,\alpha}\)-trace and of the refresh number.  This reset, rather
than the same-order Schauder estimate, prevents multiplication of a
 fixed constant over successive near slabs.  The quantitative lower
singular-value face and the differentiated inverse identities give
the inverse bound on \(\mathcal V^3\), since its entire terminal
preimage lies in \(\mathcal Z^3\), while the forward estimate is
restricted from \(\mathcal Z^3\) to \(\mathcal U^3\).  This is why
\eqref{eq:C6-near-reset-trace} is a \(C_{\rm map}\)-estimate rather
than only a coordinate estimate for \(F\).  A terminally truncated
 slab ending before its reset time is already controlled on
 \(\mathcal U^0\subset\mathcal U^1\) by the same-order propagation from
 its incoming trace, which is either the prepared trace at \(\tau_0\)
 or the preceding reset trace.

If a region is declared remote at the end of a full near slab, the
atlas lemma enters a fresh subordinate remote pair at that handoff
time.  Thus
\eqref{eq:C6-near-reset-trace} supplies exactly the
\(C_{\rm map}^{q,\alpha}(\mathcal U^2,\mathcal V^2)\) entrance trace
needed to apply
Lemma~\ref{lem:short-effective-clock-propagation} with output on
\(\mathcal U^1\), after inserting one auxiliary intermediate buffer
between levels one and two.  A
pair which is already remote at \(\tau_0\) uses the prepared initial
trace instead.

Once a pair is declared remote at an entry time \(s\), it is retained
and is not reclassified, even if the scalar comparison
\(r_{\mathcal U}^{\,2}\lessgtr K_{\rm rem}\lambda(\tau)\) later
recrosses.  The exact scale equation and the phase budget give
\begin{equation}\label{eq:finite-C6-retained-remote-clock}
 \vartheta_{\mathcal U}^{\,s}(\tau_1)
 =
 \int_s^{\tau_1}
   \frac{\lambda(q)}{r_{\mathcal U}^{\,2}}\,dq
 \leq
 C_{\rm clock}\frac{\lambda(s)}{r_{\mathcal U}^{\,2}}
 \leq
 \frac{C_{\rm clock}}{K_{\rm rem}}
 \leq\Theta_{\rm rem}.
\end{equation}
The initially remote \(q=2\) step requires one additional argument:
bilipschitz control alone does not make \(DF\) H\"older.  For an active
paired member \(\mathfrak a\), write \(s_{\mathfrak a}\) for its entry
time, \(I_{\mathfrak a}\) for its active lifespan, and
\(r_{\mathfrak a}\) for its source harmonic radius.  Let
\(\mathfrak A_{\rm rem}\) be the retained pairs declared remote at
their entry times, and let \(\mathfrak A_{\rm near}\) be the pairs
belonging to the successive near slabs.  In the frozen entry source
coordinates and the effective clock
\eqref{eq:C6-source-effective-clock}, the principal matrix is
\begin{equation}\label{eq:C6-remote-principal-matrix}
 \mathsf A_{\mathfrak a}^{ij}(\tau,x)
 :=r_{\mathfrak a}^{\,2}\acute G(\tau)^{ij}(x).
\end{equation}
Here \(\acute G^{ij}\) denotes the inverse coordinate matrix in the
frozen entry source chart.  The source-connection first-order
coefficient in these coordinates is
\begin{equation}\label{eq:C6-remote-source-first-order}
 \mathsf B_{\mathfrak a}^{k}
 :=-\mathsf A_{\mathfrak a}^{ij}
       \Gamma(\acute G)^{k}_{ij}.
\end{equation}
The physical coefficient package gives, uniformly in every active pair
and every \(\tau\in I_{\mathfrak a}\),
\begin{equation}\label{eq:C6-remote-principal-bound}
 \Lambda_{\rm C6}^{-1}|\xi|^2
 \leq \mathsf A_{\mathfrak a}^{ij}\xi_i\xi_j
 \leq \Lambda_{\rm C6}|\xi|^2,
 \qquad
 \|\mathsf A_{\mathfrak a}\|_{C^{0,\alpha}}
 \leq\Lambda_{\rm C6},
 \qquad
 \|\mathsf B_{\mathfrak a}\|_{C^{0,\alpha}}
 \leq\Lambda_{\rm C6}.
\end{equation}
Move this source-connection first-order block to the left-hand side.
The remaining order-two target block is
\begin{equation}\label{eq:C6-remote-q2-block}
 \mathcal N_{2,\mathfrak a}
 =\mathsf A_{\mathfrak a}^{ij}
   \Gamma(S)\bigl(F\bigr)(\partial_iF,\partial_jF)
   +\mathcal R_{2,\mathfrak a}(F,DF),
\end{equation}
where the coefficients in \(\mathcal R_{2,\mathfrak a}\) have the same
uniform \(C^{0,\alpha}\) bound and every term contains only the
displayed lower jets.  The global \(C^1\) face, the target
\(C^{2,\alpha}\) coefficient package, and buffered interpolation
therefore give, for every \(0<\delta_{\rm int}\leq1\),
\begin{equation}\label{eq:C6-remote-q2-interpolation}
 \|\mathcal N_{2,\mathfrak a}(\tau)\|_
      {C^{0,\alpha}(\mathcal U_{\mathfrak a}^2)}
 \leq C_{\delta_{\rm int}}
   +\delta_{\rm int}
    \|F_\tau\|_{C^{2,\alpha}
      (\mathcal U_{\mathfrak a}^3\to
       \mathcal V_{\mathfrak a}^3)} .
\end{equation}
Here \(C_{\delta_{\rm int}}\) is independent of
\(\mathfrak a\), \(\tau\), and the finite endpoint.

We record explicitly the buffered outer-to-inner transition used
below.  The
family \(\mathscr B_{\mathfrak a}(\tau)\) in
\eqref{eq:C6-active-neighbor-family} covers
\(\overline{\mathcal U_{\mathfrak a}^3}\) by its level-zero source
members and covers
\(\overline{\mathcal V_{\mathfrak a}^3}\) by its level-zero target
members.  It has cardinality at most \(N_{\rm tr}\), and every one of
its entry radii lies between
\(C_{\rm tr}^{-1}r_{\mathfrak a}\) and
\(C_{\rm tr}r_{\mathfrak a}\), by
\eqref{eq:C6-active-cross-generation-packing}.  The direct source and
target coordinate changes
\eqref{eq:C6-active-direct-transition}, together with the fixed gaps
\(\overline{\mathcal U_{\mathfrak b}^0}\Subset
\mathcal U_{\mathfrak b}^1\) and
\(\overline{\mathcal V_{\mathfrak b}^0}\Subset
\mathcal V_{\mathfrak b}^1\), therefore give
\begin{equation}\label{eq:C6-remote-outer-to-inner-transition}
 \|F_\tau\|_{C^{2,\alpha}
      (\mathcal U_{\mathfrak a}^3\to
       \mathcal V_{\mathfrak a}^3)}
 \leq C_{\rm tr}
 \max_{\substack{\mathfrak b\ {\rm active\ at}\ \tau\\
                  \mathfrak b\in
                   \mathscr B_{\mathfrak a}(\tau)}}
 \|F_\tau\|_{C^{2,\alpha}
      (\mathcal U_{\mathfrak b}^1\to
       \mathcal V_{\mathfrak b}^1)} .
\end{equation}
The local-to-global H\"older estimate implicit in
\eqref{eq:C6-remote-outer-to-inner-transition} uses the buffered
level-one Lebesgue margin, not a level-zero Lebesgue number for the
mixed-generation cover.  Hence its constant is independent of the
particular finite subcover and of the number of refresh generations.

Define
\begin{align}
 X^F_{2,\rm rem}
 &:={\sup}_{\mathfrak a\in\mathfrak A_{\rm rem}}
    {\sup}_{\tau\in I_{\mathfrak a}}
    \|F_\tau\|_{C^{2,\alpha}
      (\mathcal U_{\mathfrak a}^1\to
       \mathcal V_{\mathfrak a}^1)},
 \label{eq:C6-X2-remote-definition}\\
 X_{2,\rm near}
 &:={\sup}_{\mathfrak a\in\mathfrak A_{\rm near}}
    {\sup}_{\tau\in I_{\mathfrak a}}
    \|F_\tau\|_{C^{2,\alpha}
      (\mathcal U_{\mathfrak a}^1\to
       \mathcal V_{\mathfrak a}^1)},
 \label{eq:C6-X2-near-definition}\\
 E_2
 &:={\sup}_{\mathfrak a\in\mathfrak A_{\rm rem}}
    \|F_{s_{\mathfrak a}}\|_{C^{2,\alpha}
      (\mathcal U_{\mathfrak a}^2\to
       \mathcal V_{\mathfrak a}^2)} .
 \label{eq:C6-E2-definition}
\end{align}
Each supremum over an empty class is understood to be zero.
The prepared trace bounds the terms in \(E_2\) for members remote at
\(\tau_0\).  At a near-to-remote handoff, the subordinate
trace-transfer clause and \eqref{eq:C6-near-reset-trace} bound the
corresponding term; thus \(E_2<\infty\) with a bound independent of the
refresh number and the finite endpoint.  The near reset argument
already gives \(X_{2,\rm near}<\infty\).

On each retained remote pair apply the local short-effective-clock
estimate of Lemma~\ref{lem:short-effective-clock-propagation} with
level-two input, an auxiliary intermediate buffer, and level-one output,
and then take the supremum over retained pairs; each clock is bounded by
\eqref{eq:finite-C6-retained-remote-clock}.  Equations
\eqref{eq:C6-remote-q2-interpolation} and
\eqref{eq:C6-remote-outer-to-inner-transition} give
\begin{equation}\label{eq:C6-remote-q2-absorption}
 \begin{split}
 X^F_{2,\rm rem}
 &\leq C_{\delta_{\rm int}}
       (E_2+1+X_{2,\rm near})
   +C\delta_{\rm int}C_{\rm tr}X^F_{2,\rm rem}.
 \end{split}
\end{equation}
Choose \(\delta_{\rm int}>0\) after the uniform constants have been
fixed, so that
\(C\delta_{\rm int}C_{\rm tr}\leq\tfrac12\).  This closes the remote
\(C^{2,\alpha}\), hence \(C^{1,\alpha}\), base without assuming it.
The strict lower singular-value face and the differentiated inverse
identity give the corresponding \(C^{2,\alpha}\) inverse bound in the
fixed target buffers.

Freeze this order-two bound.  In the covariantly commuted
right-translated tension equation at each order \(q=3,\ldots,6\), the
current top map jet occurs only linearly, with coefficient controlled
by the uniform \(C^1\) map and geometric packages.  Every other term
contains strictly lower map jets and is controlled by the preceding
induction level.  At each such order the prepared trace at
\(\tau_0\), or \eqref{eq:C6-refresh-trace-star} together with
\eqref{eq:C6-near-reset-trace} at a near-to-remote handoff, supplies
the uniform level-two incoming norm.  Because the direct transition
maps have order twelve, the buffered outer-to-inner estimate
\eqref{eq:C6-remote-outer-to-inner-transition} has the identical
\(C^{q,\alpha}\) form for every \(q\leq6\); it is used at each
induction step, and on \(\mathcal U^2\) before applying the
differentiated inverse identity on \(\mathcal V^1\).
Lemma~\ref{lem:short-effective-clock-propagation}
therefore propagates each incoming norm over the entire remaining life
of the remote pair with the single factor
\[
 C\exp\!\left(
   C\vartheta_{\mathcal U}^{\,s}(\tau_1)
 \right).
\]
The order-twelve package controls all spatial coefficients required
through \(q=6\); no sixth-order map coefficient is assumed.  The global
\(C^1\) face supplies the outer \(C^0\) term in
\eqref{eq:short-effective-clock-propagation}.  At every higher order,
the differentiated inverse identity gives the same estimate for
\(F^{-1}\).  There is at most one near-to-remote handoff for each
retained pair, and its clock begins at that handoff; the number of
preceding near refreshes is irrelevant because each near slab ends
with the same reset estimate.

At every time, the active near and retained remote level-zero pairs
cover the source and target with uniformly bounded overlap, by
Lemma~\ref{lem:C6-persistent-source-atlas}.  Taking their supremum and
using the uniform fixed-pair typed-norm equivalence for the
right-translated map and inverse-map fields,
then restricting the level-one estimates to level zero, equivalently
the global
\(\operatorname{Map}_{\rm sc}^{6,\alpha}\) norm, proves
\eqref{eq:finite-HMHF-C6-bridge}, with no dependence on \(\tau_1\).
\end{proof}

\begin{proposition}[Effective-column tails]
\label{prop:effective-column-tails}
Assume the conclusions of
Theorem~\ref{thm:adaptive-HMHF-continuation}.  Assume in addition one
endpoint-independent bootstrap coefficient package at some
continuation order \(k_0\geq12\).  Concretely, require one constant
\(\Lambda_{\rm col}\), independent of the finite endpoint \(\tau_1\),
which bounds, in every scale-one dyadic source and target chart, the
metric and defect coefficients through order six and the
right-translated maps
\[
 R^{\pm1},\quad F^{\pm1},\quad \Theta,\quad\Phi
\]
through order six.  On near-parabolic charts require the corresponding
common coefficient time-oscillation moduli; on remote charts require
the endpoint-independent spatial package used in the same-order
propagation argument.  In the global first-exit argument these are
explicit bootstrap
faces.  The physical coefficient part is discharged by
Corollary~\ref{cor:adaptive-auxiliary-closure}, the target spatial part
by Proposition~\ref{prop:adaptive-target-tracking}, the target-time and
defect part by Lemma~\ref{lem:coarse-effective-Gram}, and the order-six
map and inverse-map part by
Lemma~\ref{lem:finite-HMHF-C6-bridge}.  They are stated here as
hypotheses so that the finite-horizon constant \(\Lambda_T\) from
Theorem~\ref{thm:adaptive-HMHF-continuation} is not treated as
endpoint-independent.  For every integer
$0\leq m\leq4$, the effective columns obey, including at the initial
face, the global scale-normalized bounds
\begin{equation}\label{eq:effective-column-scale-bounds}
 \sup_{\tau_0\leq\tau<\tau_1}\sup_M
 \sum_{\ell=0}^m
 (1+\bar f)^{\ell/2}
 |\bar\nabla^\ell\mathcal Y_{j,\tau}|_{\bar g}
 \leq C_m,
 \qquad0\leq j\leq8.
\end{equation}
All constants are independent of the finite endpoint \(\tau_1\).
In the pre-radius application, the direct package
\(\mathscr P_{\rm pre}^{(6)}\), together with the conclusions of
Lemmas~\ref{lem:coarse-effective-Gram} and
\ref{lem:finite-HMHF-C6-bridge}, supplies constants
\(K_{\mathcal Y,m}^{\rm dir}<\infty\),
\(0\leq m\leq4\), depending only on that package and uniformly in
\(\Gamma\geq\Gamma_{\rm pre}\), such that
\(C_m\leq K_{\mathcal Y,m}^{\rm dir}\).  These direct constants may be
used as the numerical pre-radius column ceilings in the later closure.
Set
\[
 \mathscr S_{\tau}^{\rm col}
 :=\supp(1-\chi_\tau)
   \cup\Phi_\tau(\supp(1-\eta)).
\]
Then
\[
 \supp(\mathcal Y_{j,\tau}-Y_j)
 \subset\mathscr S_{\tau}^{\rm col},
 \qquad0\leq j\leq8.
\]
Consequently, if \(K\Subset M\) satisfies
\(K\cap\mathscr S_{\tau}^{\rm col}=\varnothing\), then
\(\mathcal Y_{j,\tau}=Y_j\) on \(K\) for \(0\leq j\leq8\).  For an arbitrary
\(K\Subset K^+\Subset M\), there is a finite support-exit time
\(\tau_{\rm exit}(K^+)\), independent of the terminal endpoint, such
that
\begin{equation}\label{eq:effective-column-support-exit}
 \mathcal Y_{j,\tau}=Y_j\quad\hbox{on }K^+,\qquad
 \tau\geq\tau_{\rm exit}(K^+),\quad0\leq j\leq8.
\end{equation}
Hence arbitrary spatial orders are available after support exit.  On
the finite prefix before \(\tau_{\rm exit}(K^+)\), column bounds hold
only through the order carried by the entrance metric, map, and ODE
package.  Smooth entrances carry every fixed order, whereas finite-order
entrance data yield no simultaneous all-order global scale-normalized
estimate.
In particular, the direct columns and the scale derivatives needed in
the Kato and Bernstein arguments are uniformly bounded on every
receding transition annulus.  For \(0\leq m\leq4\), globally there are
$A_m,N_m<\infty$ such that the chart and target derivatives which
enter \eqref{eq:effective-column-zero}--%
\eqref{eq:effective-column-j} are bounded by
\[
 C_me^{A_m\tau}(1+\bar f)^{N_m}.
\]
Consequently there are $c_m,C_m>0$ such that
\begin{equation}\label{eq:effective-column-tail}
 \norm{\mathcal Y_{j,\tau}-Y_j}_{H^m_\nu}
 \leq C_me^{-c_me^\tau},
 \qquad0\leq j\leq8.
\end{equation}
Here and in \eqref{eq:effective-column-tail},
\(0\leq m\leq4\).  These orders contain every derivative and pairing
used in the receding Gram matrix.  In particular,
\begin{equation}\label{eq:effective-Gram}
 \ip{\rho_\tau\mathcal Y_{j,\tau}}{Z_\mu}
 =\ip{Y_j}{Z_\mu}+O(e^{-ce^\tau}).
\end{equation}
\end{proposition}

\begin{proof}
The cutoff part
$\Lie_{\chi_\tau W_j}\bar g-Y_j$ is covered by
Lemma~\ref{lem:cutoff-tails}.  Every $K_\tau(T)$ is supported in
\[
 \Phi_\tau(\supp(1-\eta))
 \subset\{\bar f\geq c\Gamma e^\tau\},
\]
by \eqref{eq:F-outer-drift} and
Proposition~\ref{prop:adaptive-target-tracking}.  On that set, the
tensors $Y_j$ and
their derivatives have polynomial growth.

We first prove the stronger annular statement.  Write
\[
 J_\tau=\Theta_\tau\circ\Phi_\tau^{-1}.
\]
On the support of $(1-\eta)\circ\Phi_\tau^{-1}$,
\[
 K_\tau(T)
 =\bigl((1-\eta)\circ\Phi_\tau^{-1}\bigr)J_\tau^*T.
\]
The identity \eqref{eq:relative-metric-identity}, the global
$C^2$ bootstrap, and the endpoint-independent order-six package in
the hypotheses, obtained order by order from
Lemmas~\ref{lem:weighted-prepared-Schauder} and
\ref{lem:short-effective-clock-propagation}, give uniform
scale-normalized derivative bounds through order \(m+1\) for $J_\tau$ and
$J_\tau^{-1}$ on
$\{\bar f\geq c\Gamma e^\tau\}$.  Here one uses the global rescaled-annulus
bounds in the initial-position hypothesis and
Proposition~\ref{prop:adaptive-target-tracking}; the parabolic
estimates have the same constant on every rescaled dyadic annulus.
Annulus tracking gives the same bounds for the transported cutoff.
Finally, the AC
estimates give
\[
 (1+\bar f)^{\ell/2}
 |\bar\nabla^\ell Y_j|_{\bar g}\leq C_\ell,
 \qquad
 (1+\bar f)^{\ell/2}
 |\bar\nabla^\ell\Lie_{\chi_\tau W_j}\bar g|_{\bar g}
 \leq C_\ell
\]
on the complete exterior.  The product and pullback formulae prove
the asserted bound on the support of $K_\tau$.  Below its receding
support $K_\tau=0$.  The global AC estimates for $Y_0$ and the
scale-normalized cutoff estimates for
$\Lie_{\chi_\tau W_j}\bar g$ give the same bound there.  This proves
\eqref{eq:effective-column-scale-bounds} on all of $M$.  Moreover, the
definitions give the exact support inclusion in the statement: for
$j\geq1$ the difference from $Y_j$ is
$\Lie_{(\chi_\tau-1)W_j}\bar g-K_\tau(\Lie_{\chi_\tau W_j}\bar g)$,
while for $j=0$ it is $-K_\tau(Y_0)$.  Let $K^+\Subset M$ and put
$R=\sup_{K^+}\bar f$.  Since
\[
 \supp(1-\chi_\tau)
 \subset\{\bar f\geq(5/2)e^\tau\},
 \qquad
 \Phi_\tau(\supp(1-\eta))
 \subset\{\bar f\geq c\Gamma e^\tau\},
\]
choose $\tau_{\rm exit}(K^+)$ so that
\[
 R<\min\{5/2,c\Gamma\}e^\tau
 \quad\text{for every }\tau\geq\tau_{\rm exit}(K^+).
\]
Then $K_\tau(T)=0$ for every tensor $T$ occurring in the columns and
$\chi_\tau=1$ on $K^+$, which proves
\eqref{eq:effective-column-support-exit}.  If \(K^+\) meets the
 transported support on the preceding compact time prefix, the
 initial-value estimates for the relative harmonic-map equation, the
 buffered Ricci--DeTurck estimates, and the \(R,\Theta\) variational
 equations propagate exactly the finite spatial order supplied at the
 entrance.  They do not create missing jets of the ODE-carried target
 or marking.  For a smooth entrance the same argument may be run at
 any prescribed finite order.

It remains, for \(0\leq m\leq4\), to obtain a sufficiently coarse
global bound for the Gaussian tail.  The harmonic-map equation and
the $C^2$ box, now with one common coefficient, defect, map, and
inverse-map package, the near-parabolic time modulus, and the remote
short-clock estimate \eqref{eq:short-effective-clock-propagation} give
the required estimates on each
unit normalized-time interval with a constant independent of the
terminal endpoint.  Iterating the local estimates
and the variational equations for $\Theta$ and $\Phi$ gives at worst
$C_me^{A_m\tau}(1+\bar f)^{N_m}$ growth.  Uniform boundedness is not
needed: Lemma~\ref{lem:superexp} absorbs every exponential in $\tau$
and every polynomial in $\bar f$ on
$\{\bar f\geq c\Gamma e^\tau\}$.  This proves
\eqref{eq:effective-column-tail}; Cauchy--Schwarz gives
\eqref{eq:effective-Gram}.
\end{proof}

\begin{remark}
With the cutoffs nested so that
$\Phi_\tau(\supp(1-\eta))\cap\supp\rho_\tau=\varnothing$ and
$\chi_\tau=1$ on $\supp\rho_\tau$, the adaptive columns agree exactly
with $Y_0,\ldots,Y_8$ in the receding Gram system.
Proposition~\ref{prop:effective-column-tails} is retained because it is stable
under small changes of the graft and cutoff positions.
\end{remark}

\begin{corollary}[Adaptive feedback package]
\label{cor:adaptive-feedback-package}
Let \([\tau_0,\tau_1)\) be an admissible first-exit interval with
finite endpoint satisfying the hypotheses of
Proposition~\ref{prop:effective-column-tails} and the adaptive
counterparts of all hypotheses of
Theorems~\ref{thm:receding} and \ref{thm:phase}.  In particular, the
normalized equation, robust \(C^2\) box, forcing bounds, initial exact
slice, and scale and phase equations are part of the hypotheses, and
the velocities are chosen by the exact modified Gram system.  Then all
conclusions of
Theorem~\ref{thm:receding} and
Theorem~\ref{thm:phase} remain valid for the adaptive normalized
equation when $Y_j$ is replaced by
$\mathcal Y_{j,\tau}$ and the velocities are defined by the exact
modified Gram system.  In particular, the modified matrix is
uniformly invertible, the exact receding slice is propagated, and
\eqref{eq:receding-velocity}, \eqref{eq:receding-energy}, and
\eqref{eq:finite-phase-tail} hold with constants independent of the
finite endpoint.  The same replacement is valid in the
future-phase-tail Lemma~\ref{lem:future-phase-tail}.
On the pre-radius first-exit box, all upper constants in this package
are at most \(K_{\rm fb}^{\rm pre}\) and every favorable constant is at
least \((K_{\rm fb}^{\rm pre})^{-1}\); this is the final assertion of
Lemma~\ref{lem:pre-radius-low-order-closure}.
\end{corollary}

\begin{proof}
The pointwise and scale-derivative hypotheses used in the Kato and
local parabolic estimates are exactly
\eqref{eq:effective-column-scale-bounds}.  In the modified Gram
matrix and in every Gaussian pairing, the difference from the global
columns is bounded by
\eqref{eq:effective-column-tail}--\eqref{eq:effective-Gram}.
These errors are Gaussian-superexponential and therefore are absorbed
by the same smallness choices as the moving-cutoff commutator.  The
action terms involving $h$ obey the original first-moment estimates
because the effective fields have the same scale-normalized growth.
Thus the proofs of the two cited theorems apply term by term.  All
estimates are made on an arbitrary finite endpoint; no infinite-time
existence is used.
\end{proof}

\begin{lemma}[Uniform pre-radius low-order closure]
\label{lem:pre-radius-low-order-closure}
Fix the FIK background, the radial and moving-cutoff profiles,
\(0<\alpha<1\), \(m_{\rm ad}=13\), and one rate pair
\[
 0<\sigma<\theta<\beta .
\]
 The geometric atlas base \(\Gamma_{\rm atl}\) and the pre-radius
 family \eqref{eq:pre-radius-dyadic-family} have already been frozen
 and will not be changed.  Before choosing a graft radius, fix the
 finite primitive tuple \(\mathfrak P_{\rm pre}^{\rm prim}\) defined,
 with its complete dependency list, in
 \eqref{eq:pre-radius-primitive-tuple}.
Here
\(\Lambda_{\rm coef}\) bounds the scale-one spatial metric
coefficients through order twelve and the physical Ricci-defect
coefficients through order ten.  The two distinct map ceilings mean
\[
 \|R_{\tau_0}^{\pm1}\|_{\operatorname{Map}_{\rm sc}^{14,\alpha}}
 \leq\Lambda_{R,14}^{\rm pre}-2\mu_R^{\rm pre},\qquad
 \|F_0^{\pm1}\|_{\operatorname{Map}_{\rm sc}^{6,\alpha}}
 \leq\Lambda_{F,6}^{\rm pre}-2\mu_F^{\rm pre},
\]
and both are chosen below the umbrella ceiling
\(\Lambda_{\rm map}\), with positive room.  The positive margins
\(\mu_R^{\rm pre}\) and \(\mu_F^{\rm pre}\) reserve explicit upper-norm
slack in these two entrance faces.  The margins
 \(\kappa_{\rm map}\) and \(\kappa_{\rm har}\) are, respectively, the
 scale-one lower singular-value/inverse-map margin and the
 dimensionless harmonic-radius margin
 \(\mathfrak h_{\rm har}\), equivalently
 \(r_{\rm har}/r_{{\rm sol},\tau}\); see
 \eqref{eq:dimensionless-harmonic-radius-face}.  This primitive tuple
 contains no forcing, effective-column, feedback, or eventual-radius
 atlas constant.

Using only the complete input record
\(\mathfrak I_{\rm pre}
=(\alpha,\sigma,\theta,\mathfrak P_{\rm pre}^{\rm prim})\), the
 fixed pre-radius atlas, and the fixed background, there exist a
package-dependent threshold
\[
 \Gamma_{\rm pre}
 =\Gamma_{\rm pre}(\mathfrak I_{\rm pre})
 \geq\Gamma_{\rm atl}
\]
and numbers
 \[
  \begin{gathered}
   \delta_{\rm atl}^{\rm pre}>0,\qquad
   0<\delta_{\rm pre}\leq\delta_{\rm atl}^{\rm pre},\qquad
   \varepsilon_{\rm pre}\in(0,\log2],\qquad
  \tau_{\rm pre}<\infty,\\
  C_{\rm ad}^{\rm pre}(m)<\infty\quad(0\leq m\leq5),\qquad
   0<c_{\rm ann}^{\rm pre}\leq C_{\rm ann}^{\rm pre}<\infty,\qquad
   C_{\rm gr}^{\rm pre},C_{\mathcal E}^{\rm pre}<\infty,\qquad
  c_{\mathcal E}^{\rm pre}>0,\\
  \Lambda_{\rm C6}^{\rm pre}<\infty,\qquad
  C_{\rm map,6}^{\rm pre}<\infty,\qquad
  K_{\mathcal Y,m}^{\rm pre}<\infty\quad(0\leq m\leq4),\qquad
  K_{\rm fb}^{\rm pre}\geq1 ,
 \end{gathered}
\]
all independent of \(\Gamma\).  Set
\begin{equation}\label{eq:pre-radius-robust-graft-ceiling}
 C_{\rm gr}^{\rm rob,pre}
 :=C_{\rm ann}^{\rm pre}C_{\rm gr}^{\rm pre}.
\end{equation}
 The displayed witness defines the derived tuple
 \(\mathfrak P_{\rm pre}^{\rm der}\) in
 \eqref{eq:pre-radius-derived-tuple}.  It is fixed before any
 first-exit interval or radius functional is evaluated; its primitive
 inputs are exactly those specified above.
 The construction is uniform in
 \(\Gamma\geq\Gamma_{\rm pre}(\mathfrak I_{\rm pre})\) for each fixed
 input record, but no uniformity is asserted at the boundary of the
 open rate cone.

For every \(\Gamma\geq\Gamma_{\rm pre}\), construct
\(\eta_\Gamma\) from the fixed radial profile.  Then
\begin{equation}\label{eq:pre-radius-cutoff-derivatives}
 \sup_{\Gamma\geq\Gamma_{\rm pre}}
 \sup_{\{\Gamma/2<\bar f<\Gamma\}}
 \sum_{\ell=0}^{m_{\rm ad}+2}
 \Gamma^{\ell/2}
 |\bar\nabla^\ell\eta_\Gamma|_{\bar g}
 <\infty .
\end{equation}
Consider an admissible first-exit interval
\([\tau_0,\tau_1)\), in the sense of
Definition~\ref{def:admissible-first-exit-interval}, with finite
endpoint and \(\tau_0\geq\tau_{\rm pre}\), on which:
\begin{enumerate}
\item the activation face satisfies
      \[
       \sup_{[\tau_0,\tau_1)}
       \|h(\tau)\|_{\mathfrak C_{{\rm pre},0}^{2,\alpha}}
       \leq\delta_{\rm pre}.
      \]
      The complete source and target metrics remain in the
      \(2\Lambda_{\rm ell}\) ellipticity box.  Their physical
       harmonic face satisfies
       \(\mathfrak h_{\rm har}(\tau)\geq\kappa_{\rm har}/2\),
       their order-twelve spatial coefficient
       and order-ten Ricci-defect faces are at most
       \(2\Lambda_{\rm coef}\), and the radial-comparison and primitive
       graft-compatibility faces lie at their doubled upper thresholds
       and half lower margins.  The weak annulus face is the explicit
       inclusion
       \begin{equation}\label{eq:pre-radius-weak-annulus-face}
        \Phi_\tau(\Omega_\eta^{++})\subset
        \left\{\tfrac12c_{\rm ann}^{\rm pre}\Gamma e^\tau
        <1+\bar f<2C_{\rm ann}^{\rm pre}\Gamma e^\tau\right\}.
       \end{equation}
       The
       quantitative support face is explicitly
       \begin{equation}\label{eq:pre-radius-support-separation-activation}
        \inf_{\tau_0\leq\tau<\tau_1}
        \mathfrak s_{\rm sep}(\tau)
        \geq\tfrac12\kappa_{\rm sep}.
       \end{equation}
       The
      coarse cutoff Gram matrix formed from the exact background
       columns on the receding support has least singular value at least
      \(\kappa_{\rm Gram}/2\);
\item
      the one-sided scale-and-phase package
      \begin{equation}\label{eq:pre-radius-one-sided-phase-package}
       \lambda_\tau=-(1+a)\lambda,\qquad
       0<\lambda(\tau)\leq2C_{\rm sc}e^{-\tau},\qquad
       \int_{\tau_0}^{\tau_1}(|a|+|b|)\,d\tau
       \leq\varepsilon_{\rm pre}
      \end{equation}
      holds.  In particular, for
      \(\tau_0\leq s\leq\tau<\tau_1\),
      \begin{equation}\label{eq:pre-radius-relative-scale}
       e^{-\varepsilon_{\rm pre}}e^{-(\tau-s)}
       \leq\frac{\lambda(\tau)}{\lambda(s)}
       \leq
       e^{\varepsilon_{\rm pre}}e^{-(\tau-s)}.
      \end{equation}
\item the entrance representatives obey the distinct static faces
      \[
       \|R_{\tau_0}^{\pm1}\|_
         {\operatorname{Map}_{\rm sc}^{14,\alpha}}
       \leq\Lambda_{R,14}^{\rm pre}-2\mu_R^{\rm pre},\qquad
       \|F_0^{\pm1}\|_
         {\operatorname{Map}_{\rm sc}^{6,\alpha}}
       \leq\Lambda_{F,6}^{\rm pre}-2\mu_F^{\rm pre},
      \]
       together with the common entrance radial comparison
       \begin{equation}\label{eq:pre-radius-entrance-radial-comparison}
        c_{\rm rad}(1+\bar f(x))
        \leq1+\bar f(\vartheta(x))
        \leq C_{\rm rad}(1+\bar f(x)),
        \qquad
        \vartheta\in
        \{R_{\tau_0},R_{\tau_0}^{-1},F_0,F_0^{-1}\},
       \end{equation}
       and the strict local-invertibility margins
       \[
        \mathfrak m_R(R_{\tau_0})>\kappa_{\rm map},
        \qquad
        \mathfrak m_F(F_0;\acute G_{\tau_0},S_{\tau_0})
        >\kappa_{\rm map},
       \]
       and
       \[
        \mathfrak s_{\rm sep}(\tau_0)\geq\kappa_{\rm sep};
       \]
\item throughout the interval, \(R_\tau^{\pm1}\) and
      \(F(\tau)^{\pm1}\) remain degree-one proper diffeomorphisms in
      the fixed buffered harmonic-map charts.  Their weak radial
      bootstrap faces are, explicitly,
      \begin{align}
       \frac12c_{\rm rad}(1+\bar f(x))
       &\leq1+\bar f(\vartheta(x))
       \leq2C_{\rm rad}(1+\bar f(x)),
       &&\vartheta\in\{R_\tau,R_\tau^{-1}\},
       \label{eq:pre-radius-weak-R-radial-face}\\
       \frac12c_F^\sharp(1+\bar f(x))
       &\leq1+\bar f(\vartheta(x))
       \leq2C_F^\sharp(1+\bar f(x)),
       &&\vartheta\in\{F_\tau,F_\tau^{-1}\}.
       \label{eq:pre-radius-weak-F-radial-face}
      \end{align}
      Their bootstrap
      lower-margin face is
      \[
       \min\{\mathfrak m_R(R_\tau),
       \mathfrak m_F(F_\tau;\acute G_\tau,S_\tau)\}
       \geq\frac12\kappa_{\rm map};
      \]
      the source--target metric discrepancy
      and the coefficients of the harmonic-map system lie in the
      doubled \((\Lambda_{\rm ell},\Lambda_{\rm coef})\) box; the
       target radial faces lie in the doubled
       \(C_{\rm ad}^{\rm pre}(m)\), \(0\leq m\leq5\), faces, and the
       annulus face is \eqref{eq:pre-radius-weak-annulus-face}; the order-six map and
      inverse-map face is placed at
      \(2C_{\rm map,6}^{\rm pre}\).  In addition, the propagated
      relative markings \(R_\tau^{\pm1}\) retain the
      \(2\Lambda_{R,14}^{\rm pre}\)
      \(\operatorname{Map}_{\rm sc}^{14,\alpha}\) face, which the
      conclusion improves to
      \(\Lambda_{R,14}^{\rm pre}-\mu_R^{\rm pre}\).
\end{enumerate}
No effective-column bound and no sharp receding-feedback estimate is
assumed here.  The coarse \(C^2\) and accumulated-phase inequalities
are independent activation hypotheses rather than conclusions of the
closure.
 The physical ellipticity, physical harmonic-radius, physical
 coefficient, physical Ricci-defect, and primitive
 \(K_{\rm gr}\)-compatibility face in item~(1)
 are activation faces only; their independent improvement is supplied
 by the low-order physical and graft-compatibility closures invoked in
 the coupled bootstrap.  The present lemma improves the remaining
 derived atlas, map, annulus, Gram, effective-column, and feedback
 faces, together with the displayed derived graft-forcing estimates,
 to their undoubled ceilings.  In particular, it yields the strict
 map-margin improvement
 \begin{equation}\label{eq:pre-radius-map-margin-improvement}
  \inf_{\tau_0\leq\tau<\tau_1}
  \min\{\mathfrak m_R(R_\tau),
  \mathfrak m_F(F_\tau;\acute G_\tau,S_\tau)\}
  \geq\frac34\kappa_{\rm map}.
 \end{equation}
 It improves every other named lower margin to at least three quarters
 of the corresponding primitive margin, but
 does not claim a new improvement of the primitive graft-compatibility
 face.  In particular,
independently of both \(\Gamma\) and \(\tau_1\),
\begin{align}
 \supp\mathcal E_{\rm gr}(\tau)
 &\subset
 \{e^\tau\leq\bar f
      \leq C_{\rm ann}^{\rm pre}\Gamma e^\tau\},
 \label{eq:pre-radius-uniform-graft-support}\\
 \sum_{\ell=0}^{2}
 |\bar\nabla^\ell\mathcal E_{\rm gr}|_{\bar g}
 &\leq
 \frac{C_{\rm gr}^{\rm pre}}{\Gamma}e^{-\tau},
 \label{eq:pre-radius-uniform-graft-pointwise}\\
 \|\rho_\tau\mathcal E_{\rm gr}(\tau)\|_{H^{-1}_\nu}
 &\leq
 C_{\mathcal E}^{\rm pre}
 e^{-c_{\mathcal E}^{\rm pre}e^\tau},
 \label{eq:pre-radius-uniform-graft-gaussian}\\
 \sup_{\tau_0\leq\tau<\tau_1}
 \|F(\tau)\|_{\operatorname{Map}_{\rm sc}^{6,\alpha}}
 &\leq C_{\rm map,6}^{\rm pre},
 \label{eq:pre-radius-uniform-C6-map}\\
 \sup_{\tau_0\leq\tau<\tau_1}\sup_M
 \sum_{\ell=0}^{m}
 (1+\bar f)^{\ell/2}
 |\bar\nabla^\ell\mathcal Y_{j,\tau}|_{\bar g}
 &\leq K_{\mathcal Y,m}^{\rm pre},
 \quad
 0\leq m\leq4,\quad0\leq j\leq8.
 \label{eq:pre-radius-uniform-columns}
\end{align}
The paired source--target atlas and its refresh maps may simultaneously
be chosen with coefficient, overlap, label-comparison, containment,
and effective-clock constants at most
\(\Lambda_{\rm C6}^{\rm pre}\).  Moreover, every upper constant in the
receding energy and feedback estimates used in
Lemma~\ref{lem:future-phase-tail} is at most
\(K_{\rm fb}^{\rm pre}\), and every favorable coercivity,
inverse-Gram, or decay constant is at least
\((K_{\rm fb}^{\rm pre})^{-1}\).

Consequently the choices
\begin{equation}\label{eq:pre-radius-three-region-inputs}
 \begin{gathered}
 C_{\mathcal E}=C_{\mathcal E}^{\rm pre},\qquad
 c_{\mathcal E}=c_{\mathcal E}^{\rm pre},\qquad
 C_{\rm gr}=C_{\rm gr}^{\rm rob,pre},\\
 K_{\mathcal Y,m}=K_{\mathcal Y,m}^{\rm pre}
 \quad(0\leq m\leq4),\qquad
 K_{\rm fb}=K_{\rm fb}^{\rm pre}
 \end{gathered}
\end{equation}
are legitimate numerical inputs to
\(\overline\Gamma_{\rm 3reg}\), fixed before the actual radius
\(\Gamma\).  The lower scale margin \(c_{\rm sc}\) and every
\(\Gamma\)-dependent high-regularity atlas constant are absent from
these dependencies.
\end{lemma}

\begin{proof}
Use the fixed core together with the pre-radius dyadic family
\eqref{eq:pre-radius-dyadic-family}.  Choose
\(\Gamma_{\rm pre}\geq\Gamma_{\rm atl}\) once so that every fixed
enlargement of a possible graft-transition collar lies in the
quantitative AC region and the fixed radial, separation, and Gaussian
tail margins used below have positive room.  This changes neither
\(\mathscr L_{\rm pre}\) nor any norm or ceiling already defined on
 it.  Choose \(R_{\rm AC}<\infty\) and fixed constants
 \(0<c_{\rm sol}\leq C_{\rm sol}<\infty\) so that the FIK radial-flow
 comparison
 \[
  c_{\rm sol}e^\tau(1+\bar f(y))
  \leq1+\bar f(\varphi_\tau y)
  \leq C_{\rm sol}e^\tau(1+\bar f(y))
 \]
 holds whenever \(\bar f(y)\geq R_{\rm AC}\).  This is an AC-end
 statement, not a global assertion at the fixed zero section.
 Increase \(\Gamma_{\rm pre}\), still using only the primitive input
 record, so that the weak lower radial faces place every intermediate
 point used for \(\Omega_\eta^{++}\) and
 \(\operatorname{supp}(1-\eta)\) in
 \(\{\bar f\geq R_{\rm AC}\}\).  Choose, with room for the fixed
 core--end overlap,
  \[
   0<c_{\rm ann}^{\rm pre}
     <c_{\rm sol}c_\eta c_{\widehat\Phi}^\sharp,
   \qquad
   C_{\rm ann}^{\rm pre}
     >C_{\rm sol}C_\eta C_{\widehat\Phi}^\sharp,
  \]
 and use this \(C_{\rm ann}^{\rm pre}\) as the derived annulus
 ceiling.  Put \(c_\Phi^{\rm pre}:=c_{\rm ann}^{\rm pre}\).  These
 choices use only the primitive input record.  Increase
 \(\Gamma_{\rm pre}\), using
only the primitive input record, so that
\begin{equation}\label{eq:pre-radius-direct-separation-choice}
 c_\Phi^{\rm pre}
 \left(1+\frac23\Gamma_{\rm pre}\right)-3
 \geq\frac34\kappa_{\rm sep}.
\end{equation}
If
\(2^q\Gamma_{\rm atl}\leq\Gamma<2^{q+1}\Gamma_{\rm atl}\), then the
transition collar of \(\eta_\Gamma\), including every fixed
enlargement used below, meets only a bounded number of pre-radius
annuli whose labels are comparable with \(\Gamma\).  The AC estimates
give one scale-one \(C^{m_{\rm ad}+2,\alpha}\) metric package,
ellipticity bound, harmonic-radius lower bound, and transition-map
bound on all these annuli.  The chain rule applied to the fixed radial
profile gives \eqref{eq:pre-radius-cutoff-derivatives}.  The two-sided
 radial comparison and the pre-radius map bound transport the same
 estimates to the source and target charts.  This proves the asserted
 uniform \(C_{\rm ad}^{\rm pre}(m)\) and atlas constants.  Use the same
 \(\delta_{\rm atl}^{\rm pre}>0\) output by the pre-radius branch of
 Lemma~\ref{lem:C6-persistent-source-atlas}.  Its permitted
 primitive-margin decrease has already been made so that the six-level
 buffered containments and the scale-one coefficient and
 lower-singular-value faces persist
 whenever
 \(\|h\|_{\mathfrak C_{{\rm pre},0}^{2,\alpha}}
 \leq\delta_{\rm atl}^{\rm pre}\).  Thereafter choose
 \(\delta_{\rm pre}\leq\delta_{\rm atl}^{\rm pre}\).
 The number of charts never enters, since all global estimates take
 suprema over a uniformly locally finite family.

The exact scale equation and
\(\varepsilon_{\rm pre}\leq\log2\) give the upper scale comparison
needed in the physical and target estimates.  Repeating the low-order
proof of Proposition~\ref{prop:graft-compatibility} in the pre-radius
atlas gives the indicated graft estimate.  We do not invoke
Proposition~\ref{prop:adaptive-target-tracking} here, because its
generic statement contains a two-sided absolute scale bracket.
Instead we now derive the relative-marking part of item {\rm(P2)}
directly from its ODE.  No growing compact atlas and no positive lower
scale margin is used.
Let \(r_{\rm rt}>0\) be a common right-translated chart radius and
\(C_{\rm rt}\geq1\) the map-norm equivalence constant on the primitive
common-margin set; these are deterministic consequences of its
map-and-inverse and local-invertibility margins.  Shrink this chart
once so that singular values are quantitatively stable there: for a
primitive constant \(C_{\rm sing,R}<\infty\),
\begin{equation}\label{eq:pre-radius-R-singular-stability}
 \mathfrak m_R(R)
 \geq\mathfrak m_R(R_{\tau_0})
 -C_{\rm sing,R}
  d_{\rm rt,sc}^{14,\alpha}(R,R_{\tau_0}).
\end{equation}
This is the scale-one \(C^1\) perturbation inequality for \(dR\) and
the differentiated inverse identity, taken uniformly over the fixed
core and dyadic charts.  Let
\(\tau^\sharp\) be the first time at which \(R_\tau\) or its inverse
leaves the common chart about \(R_{\tau_0}\).  On
\([\tau_0,\tau^\sharp)\), the direct relative-marking ODE in
\eqref{eq:relative-target-flow} gives, through order fourteen,
\[
 \sup_{\tau_0\leq\tau<\tau^\sharp}
 d_{\rm rt,sc}^{14,\alpha}(R_\tau,R_{\tau_0})
 \leq C(\varepsilon_{\rm pre}+e^{-\tau_0}).
\]
Choose \(\varepsilon_{\rm pre}\) smaller and
\(\tau_{\rm pre}\) larger, if necessary, so that
\begin{equation}\label{eq:pre-radius-right-chart-and-slack-choice}
 C(\varepsilon_{\rm pre}+e^{-\tau_0})
 \leq
 \min\left\{\tfrac12r_{\rm rt},
             \frac{\mu_R^{\rm pre}}{C_{\rm rt}},
             \frac{\kappa_{\rm map}}
                  {4C_{\rm sing,R}}\right\}.
\end{equation}
The first inequality contradicts a finite chart exit, so
\(\tau^\sharp=\tau_1\).  The second and the local norm equivalence give
\[
 \|R_\tau^{\pm1}\|_{\operatorname{Map}_{\rm sc}^{14,\alpha}}
 \leq
 \|R_{\tau_0}^{\pm1}\|_{\operatorname{Map}_{\rm sc}^{14,\alpha}}
 +
 C_{\rm rt}d_{\rm rt,sc}^{14,\alpha}
       (R_\tau,R_{\tau_0})
 \leq\Lambda_{R,14}^{\rm pre}-\mu_R^{\rm pre}.
\]
Moreover,
\eqref{eq:pre-radius-R-singular-stability}, the strict entrance
margin, and the third choice in
\eqref{eq:pre-radius-right-chart-and-slack-choice} give
\begin{equation}\label{eq:pre-radius-R-margin-closure}
 \inf_{\tau_0\leq\tau<\tau_1}\mathfrak m_R(R_\tau)
 >\frac34\kappa_{\rm map}.
\end{equation}
On every pre-radius dyadic annulus, the same differentiated
relative-marking ODE, after rescaling the domain and range metrics by
their actual annular levels, gives the full
\(\operatorname{Map}_{\rm sc}^{14,\alpha}\) bounds for
\(R_\tau^{\pm1}\).  Its radial component obeys
\[
 \left|\frac{d}{d\tau}
  \log\frac{1+\bar f(R_\tau(x))}{1+\bar f(x)}\right|
 \leq C(|a(\tau)|+|b(\tau)|)+Ce^{-\tau},
\]
uniformly in the dyadic label; the inverse variational equation gives
the identical estimate for \(R_\tau^{-1}\).  Decrease
\(\varepsilon_{\rm pre}\) and increase \(\tau_{\rm pre}\), using only
the primitive input record, so that the accumulated right side is at
most \(\log(4/3)\).  The entrance comparison
\eqref{eq:pre-radius-entrance-radial-comparison} then gives
\begin{equation}\label{eq:pre-radius-R-radial-improvement}
 \frac34c_{\rm rad}(1+\bar f(x))
 \leq1+\bar f(\vartheta(x))
 \leq\frac43C_{\rm rad}(1+\bar f(x)),
 \qquad
 \vartheta\in\{R_\tau,R_\tau^{-1}\}.
\end{equation}
This proves the \(R\)-part of item {\rm(P2)} directly from
\eqref{eq:relative-target-flow}; neither this argument nor its inverse
variational equation uses a lower bound for \(e^\tau\lambda\).  The
annulus assertion, which also contains the harmonic-map factor, has
not been used and will be proved below.

Put
\[
 q_\tau:=(F_\tau^{-1})^*\acute G_\tau .
\]
The exact identity
\(q_\tau-S_\tau=\lambda\Theta_\tau^*h(\tau)\), the primitive
ellipticity box, and the activated
\(\mathfrak C_{{\rm pre},0}^{2,\alpha}\)-bound give a primitive
constant \(C_{\rm sing,F}<\infty\) such that
\begin{equation}\label{eq:pre-radius-F-metric-discrepancy}
 (1-\epsilon_F)S_\tau\leq q_\tau\leq(1+\epsilon_F)S_\tau,
 \qquad
 \epsilon_F:=C_{\rm sing,F}\delta_{\rm pre}<1.
\end{equation}
This estimate is obtained before the persistent atlas and before the
order-six bridge.

We next close the radial faces without assuming the annulus conclusion.
Put
\[
 \widehat\Phi_\tau
 :=\varphi_{-\tau}\circ\Phi_\tau
 =R_\tau\circ F_\tau,
 \qquad
 \ell:=\log(1+\bar f).
\]
We use a radial covector aligned with the target metric, rather than a
tensor lower bound on the collapsing core.  On \(\{r>1\}\), in the exact
FIK radial orthonormal frame, the calculation in
Lemma~\ref{lem:self-contained-FIK-ledger} gives
\[
 (\bar\nabla^2\bar f)_{00}
 =\frac12+\frac{c_0}{2r^4}\geq\frac12.
\]
At the bolt \(d\ell=0\), and the estimates below extend there by
continuity.  The flow \(\varphi_\tau\) preserves the radial splitting, so its
differential expands the radial line by at least \(e^{\tau/2}\).
Since \(d\ell\) is radial and
\(\sup_M|d\ell|_{\bar g^{-1}}<\infty\),
\begin{equation}\label{eq:pre-radius-aligned-FIK-covector}
 \sup_M|d\ell|_{(\varphi_\tau^*\bar g)^{-1}}
 \leq Ce^{-\tau/2}.
\end{equation}
Using \(S_\tau=\lambda R_\tau^*\varphi_\tau^*\bar g\) gives the exact
pullback identity
\begin{equation}\label{eq:pre-radius-aligned-target-covector}
 |d(\ell\circ R_\tau)|_{S_\tau^{-1}}
 =\lambda(\tau)^{-1/2}
   \bigl(|d\ell|_{(\varphi_\tau^*\bar g)^{-1}}
          \circ R_\tau\bigr)
 \leq C(\lambda(\tau)e^\tau)^{-1/2}.
\end{equation}
No fixed-core ellipticity or annulus conclusion is used here.

Let
\[
 Y_\tau:=(\partial_\tau F_\tau)\circ F_\tau^{-1}.
\]
The normalized-time form of \eqref{eq:F-speed}, the metric comparison
\eqref{eq:pre-radius-F-metric-discrepancy}, the isometry
\(F_\tau:(M,\acute G_\tau)\to(M,q_\tau)\), and
\eqref{eq:pre-atlas-to-raw-C2-comparison} give
\begin{equation}\label{eq:pre-radius-forward-speed}
 \sup_M|Y_\tau|_{S_\tau}
 \leq C_{\rm drift}\delta_{\rm pre}\lambda(\tau)^{1/2}.
\end{equation}
Here and below \(C_{\rm drift}\) depends only on the primitive
pre-radius package.

For fixed \(x\), the chain rule gives
\begin{equation}\label{eq:pre-radius-composite-radial-speed}
 \frac d{d\tau}\ell(\widehat\Phi_\tau(x))
 =d\ell_{R_\tau(F_\tau(x))}
      \bigl[(\partial_\tau R_\tau)(F_\tau(x))\bigr]
  +d(\ell\circ R_\tau)_{F_\tau(x)}
       \bigl[Y_\tau(F_\tau(x))\bigr].
\end{equation}
The first term is bounded by
\(C(|a|+|b|)+Ce^{-\tau}\) by the relative-marking calculation leading
to \eqref{eq:pre-radius-R-radial-improvement}; on the compact part this
is the same calculation in the fixed core charts.  Equations
\eqref{eq:pre-radius-aligned-target-covector} and
\eqref{eq:pre-radius-forward-speed} bound the second term by
\(C\delta_{\rm pre}e^{-\tau/2}\).  Hence
\begin{equation}\label{eq:pre-radius-composite-radial-drift}
 \sup_{\tau_0\leq\tau<\tau_1}\sup_{x\in M}
 |\ell(\widehat\Phi_\tau(x))-\ell(\widehat\Phi_{\tau_0}(x))|
 \leq C\bigl(\varepsilon_{\rm pre}+e^{-\tau_0}
              +\delta_{\rm pre}e^{-\tau_0/2}\bigr).
\end{equation}
Decrease \(\varepsilon_{\rm pre},\delta_{\rm pre}\) and increase
\(\tau_{\rm pre}\), using only the primitive package, so that
\begin{equation}\label{eq:pre-radius-composite-radial-choice}
 C\bigl(\varepsilon_{\rm pre}+e^{-\tau_0}
              +\delta_{\rm pre}e^{-\tau_0/2}\bigr)
 \leq\log\frac43.
\end{equation}

The entrance faces for \(R_{\tau_0}^{\pm1}\) and \(F_0^{\pm1}\)
give
\[
 c_{\rm rad}^2(1+\bar f(x))
 \leq1+\bar f(\widehat\Phi_{\tau_0}(x))
 \leq C_{\rm rad}^2(1+\bar f(x)).
\]
It follows from \eqref{eq:pre-radius-composite-radial-drift}--%
\eqref{eq:pre-radius-composite-radial-choice} that
\[
 \frac34c_{\rm rad}^2(1+\bar f(x))
 \leq1+\bar f(\widehat\Phi_\tau(x))
 \leq\frac43C_{\rm rad}^2(1+\bar f(x)).
\]
Applying this inequality at \(\widehat\Phi_\tau^{-1}(x)\), taking reciprocal
bounds, and using \eqref{eq:pre-radius-derived-radial-constants} gives
the common forward--inverse estimate
\begin{equation}\label{eq:pre-radius-hat-Phi-radial-improvement}
 c_{\widehat\Phi}^\sharp(1+\bar f(x))
 \leq1+\bar f(\vartheta(x))
 \leq C_{\widehat\Phi}^\sharp(1+\bar f(x)),
 \qquad
 \vartheta\in\{\widehat\Phi_\tau,\widehat\Phi_\tau^{-1}\}.
\end{equation}
Finally,
\[
 F_\tau=R_\tau^{-1}\circ\widehat\Phi_\tau,
 \qquad
 F_\tau^{-1}=\widehat\Phi_\tau^{-1}\circ R_\tau.
\]
Combining \eqref{eq:pre-radius-R-radial-improvement} and
\eqref{eq:pre-radius-hat-Phi-radial-improvement} therefore yields
\begin{equation}\label{eq:pre-radius-F-radial-improvement}
 c_F^\sharp(1+\bar f(x))
 \leq1+\bar f(\vartheta(x))
 \leq C_F^\sharp(1+\bar f(x)),
 \qquad
 \vartheta\in\{F_\tau,F_\tau^{-1}\}.
\end{equation}
This strictly improves
\eqref{eq:pre-radius-weak-F-radial-face} and excludes a first exit
through either \(F\)-radial face.  The entire argument precedes and is
independent of annulus tracking.

Finally, \eqref{eq:fixed-enlarged-collar-radial-range} and
\eqref{eq:pre-radius-hat-Phi-radial-improvement} give
\[
 c_{\widehat\Phi}^\sharp c_\eta\Gamma
 \leq1+\bar f(\widehat\Phi_\tau(x))
 \leq C_{\widehat\Phi}^\sharp C_\eta\Gamma,
 \qquad x\in\Omega_\eta^{++}.
\]
The choice of \(\Gamma_{\rm pre}\) places these intermediate points in
the quantified AC region.  Applying the AC-end FIK estimate directly
to
\(\Phi_\tau=\varphi_\tau\circ\widehat\Phi_\tau\), and using the already
selected inequalities for
\(c_{\rm ann}^{\rm pre},C_{\rm ann}^{\rm pre}\), gives
\begin{equation}\label{eq:pre-radius-annulus-improvement}
 \Phi_\tau(\Omega_\eta^{++})
 \subset
 \{c_{\rm ann}^{\rm pre}\Gamma e^\tau
      <1+\bar f<
    C_{\rm ann}^{\rm pre}\Gamma e^\tau\}.
\end{equation}
Thus the annulus face is proved, with constants independent of
\(\Gamma\) and \(\tau_1\), before it is used in any transported-support
or Gaussian estimate.

With item {\rm(P2)} now established, the fixed radial profile,
\eqref{eq:pre-radius-cutoff-derivatives}, and the scale-one transition
maps transport the cutoff jets to the target connection.  The explicit
formula \(S_\tau=\lambda\Theta_\tau^*\bar g\), differentiated in these
charts, gives the required target spatial package through order
twelve.  The corresponding physical estimates belong to the separate
physical low-order closure and enter only through the activation faces
in item~(1).  Thus items {\rm(P1)}--{\rm(P3)} needed by the pre-radius
coarse-Gram branch have now been verified without a two-sided absolute
scale estimate.
The support face is controlled directly, including the harmonic-map
factor.  Recall that
\[
 \widehat\Phi_\tau=\varphi_{-\tau}\circ\Phi_\tau
 =R_\tau\circ F_\tau .
\]
The improved composite radial face
\eqref{eq:pre-radius-hat-Phi-radial-improvement}, followed by the fixed
lower radial estimate for \(\varphi_\tau\), gives
\[
 e^{-\tau}\bigl(1+\bar f(\Phi_\tau(x))\bigr)
 \geq c_\Phi^{\rm pre}(1+\bar f(x)).
\]
Because \(\eta_\Gamma=1\) on
\(\{\bar f\leq2\Gamma/3\}\), every
\(x\in\operatorname{supp}(1-\eta_\Gamma)\) satisfies
\(\bar f(x)\geq2\Gamma/3\).  Hence, for \(\tau\geq\tau_0\geq0\),
\begin{equation}\label{eq:pre-radius-support-separation-direct}
 \begin{split}
  e^{-\tau}\bar f(\Phi_\tau(x))-2
  &=
  e^{-\tau}\bigl(1+\bar f(\Phi_\tau(x))\bigr)
       -e^{-\tau}-2\\
  &\geq
  c_\Phi^{\rm pre}\left(1+\frac23\Gamma\right)-3
  \geq\frac34\kappa_{\rm sep}.
 \end{split}
\end{equation}
Taking the infimum proves the strict improvement of
\eqref{eq:pre-radius-support-separation-activation}.  This argument
uses the full factor \(R_\tau\circ F_\tau\); the relative-marking
equation alone does not provide a uniform absolute-displacement
estimate on the unbounded support.

The transported cutoff-jet estimate just proved is precisely the
pre-radius verification of
\eqref{eq:pure-graft-cutoff-jet-hypothesis} recorded in
Remark~\ref{rem:adaptive-pure-graft-cutoff-jets}.
The normalized pure-graft identity
\eqref{eq:pure-graft-normalized}, on the fixed \(K_{\rm gr}\)-face,
then gives \eqref{eq:pre-radius-uniform-graft-pointwise}, while
\eqref{eq:pre-radius-annulus-improvement} gives
\eqref{eq:pre-radius-uniform-graft-support}.  Gaussian
integration on
\(\{\bar f\geq c\Gamma e^\tau\}\), uniformly for
\(\Gamma\geq\Gamma_{\rm pre}\), gives
\[
 C(1+\Gamma e^\tau)^N e^{-c\Gamma e^\tau}
 \leq C_{\mathcal E}^{\rm pre}
       e^{-c_{\mathcal E}^{\rm pre}e^\tau}
\]
after increasing \(\tau_{\rm pre}\), proving
\eqref{eq:pre-radius-uniform-graft-gaussian}.

Next apply the pre-radius branch of
Lemma~\ref{lem:coarse-effective-Gram}, whose items
{\rm(P1)}--{\rm(P3)} have just been verified, and use the exact support
separation
\eqref{eq:coarse-Gram-support-separation}.  It gives
\(\rho_\tau K_\tau(T)=0\) and
\(\rho_\tau\mathcal Y_{j,\tau}=\rho_\tau Y_j\) without any derivative
bound for \(F\) above order one.  That branch yields a
\(\Gamma\)-uniform inverse-Gram, coarse pointwise-modulation,
target-defect, and coefficient-time-modulus package before any sharp
future-tail estimate is invoked.  Decrease
\(\delta_{\rm pre},\varepsilon_{\rm pre}\) and increase
\(\tau_{\rm pre}\), if necessary, so that the perturbation from the
fixed reference matrix in \eqref{eq:background-Gram-reserve} satisfies
\[
 \|M^{\rm low}(\tau)-\mathbf G^{\rm bg}\|
 \leq\frac14\kappa_{\rm Gram}.
\]
The standard perturbation inequality for \(s_{\min}\), together with
\eqref{eq:primitive-Gram-admissibility}, then gives
\[
 s_{\min}M^{\rm low}(\tau)
 \geq\gamma_{\rm Gram}^{\rm bg}
      -\frac14\kappa_{\rm Gram}
 >\frac34\kappa_{\rm Gram}.
\]
Thus the half-margin face is strictly improved without assuming a full
effective-column bound.

It remains to close the differential \(F\)-margin before invoking any
atlas or harmonic-map bridge that uses it.  In the same finite initial
choice of \(\delta_{\rm pre}\), require
\begin{equation}\label{eq:pre-radius-F-margin-choice}
 \min\{(1+\epsilon_F)^{-1/2},
           (1-\epsilon_F)^{1/2}\}
 >\frac34\kappa_{\rm map}.
\end{equation}
This is possible precisely because
\eqref{eq:primitive-map-admissibility} fixes
\(0<\kappa_{\rm map}<1\).  Since
\[
 F_\tau:(M,\acute G_\tau)\longrightarrow(M,q_\tau)
\]
is an isometry, \eqref{eq:pre-radius-F-metric-discrepancy} gives
\[
 s_{\min}(dF_\tau)\geq(1+\epsilon_F)^{-1/2},
 \qquad
 s_{\min}(dF_\tau^{-1})\geq(1-\epsilon_F)^{1/2}.
\]
Thus
\begin{equation}\label{eq:pre-radius-F-margin-closure}
 \inf_{\tau_0\leq\tau<\tau_1}
 \mathfrak m_F(F_\tau;\acute G_\tau,S_\tau)
 >\frac34\kappa_{\rm map}.
\end{equation}
Together with \eqref{eq:pre-radius-R-margin-closure}, this proves
\eqref{eq:pre-radius-map-margin-improvement} without using the bridge.

Construct the persistent paired atlas using
\(\mathscr L_{\rm pre}\).  The near-chart estimates are scale-one
parabolic estimates with one positive reset clock.  On every retained
remote pair, the entire remaining effective clock is bounded by
\eqref{eq:finite-C6-retained-remote-clock}.  The explicitly closed
order-two absorption
\eqref{eq:C6-remote-q2-absorption} and the successive
order-three through order-six propagation therefore involve only the
uniform package just obtained.  The level-two incoming trace is
either the prepared trace or the restriction of a level-three near
reset trace through \eqref{eq:C6-refresh-trace-star}.  Thus the proof of
 Lemma~\ref{lem:finite-HMHF-C6-bridge} applies directly.  The physical
 components of \({\mathscr P}_{\rm pre}^{(6)}\), items (P1)--(P4), are
 precisely the independent activation faces in item~(1), whereas the
 target, entrance-map, phase, discrepancy, properness, and
 lower-singular-value components have been closed from the primitive
 package by \eqref{eq:pre-radius-R-margin-closure} and
 \eqref{eq:pre-radius-F-margin-closure}.  No conclusion of the bridge
 is used to verify its
 hypotheses.  The bridge therefore yields
\eqref{eq:pre-radius-uniform-C6-map} with a right-hand side independent
of the assumed doubled order-six face.

With \eqref{eq:pre-radius-uniform-C6-map} fixed, the pullback formulas
\eqref{eq:effective-column-zero}--%
\eqref{eq:effective-column-j}, the uniform cutoff estimates, and the
scale-one product and composition estimates give
\eqref{eq:pre-radius-uniform-columns}.  On the receding Gram support
the correction vanishes exactly; off it the same Gaussian estimate
makes every column and Gram tail uniformly superexponential.  Finally,
the receding energy and feedback proofs contain only a finite list of
constants depending on the fixed spectral gap, the background Gram
inverse, \eqref{eq:pre-radius-uniform-columns}, and
\eqref{eq:pre-radius-uniform-graft-gaussian}.  Choose
\(K_{\rm fb}^{\rm pre}\) larger than every upper constant and the
reciprocal of every favorable lower constant in this list.

Every step was performed on the coarse first-exit box before invoking
the sharp three-region improvement.  Thus this is a bootstrap
improvement, not an estimate obtained from the conclusion whose radius
it is used to select.  All constants are uniform for
\(\Gamma\geq\Gamma_{\rm pre}\), while the fixed atlas base
\(\Gamma_{\rm atl}\) has not changed.  This proves the lemma.
\end{proof}

\section{Prepared Banach geometry and local evolution}
\label{sec:prepared-Banach-evolution}

Using the adaptive geometric package of
Section~\ref{sec:adaptive-continuation}, this section constructs the
Banach manifold and coupled local flow and obtains the differentiable
finite-time prepared solution map.

\subsection{The prepared Banach chart and static columns}
\label{subsec:prepared-Banach-chart}

The openness and parameter-dependence statements below use the fixed
Banach chart introduced before
Theorem~\ref{thm:intro-sharp-scattering}.  In particular,
\(\mathscr L_\Gamma,A_L,N\), the tensor and vector norms, and the tame
graph norm are those in
\eqref{eq:scaled-tensor-holder}--\eqref{eq:scaled-tame-holder}, with
their little-H\"older completions.  We retain the detailed graph
formula here because it is differentiated repeatedly below.

Fix a marked chart, the cutoff \(\eta\), a normalized time \(\tau_0\),
and a background prepared state
\(\mathbf z_*=(G_*,\lambda_*,R_*,F_*)\).  A nearby state is written
\[
 \mathbf z=(G,\lambda,R,F),\qquad
 \Theta=\varphi_{\tau_0}\circ R,\qquad
 \Phi=\Theta\circ F,
\]
and determines
\begin{equation}\label{eq:prepared-state-graph}
 \acute G(\mathbf z)
 =\eta\,\iota_*G+(1-\eta)\lambda\Theta^*\bar g,
 \qquad
 h(\mathbf z)
 =\lambda^{-1}(\Phi^{-1})^*\acute G(\mathbf z)-\bar g.
\end{equation}
Use the soliton-conjugated charts
\eqref{eq:intro-uniform-conjugated-map-coordinates}.  Thus
\(X_\Theta,X_\Phi\) are genuine vector fields on the common normalized
output; equivalently the untranslated tangent variables are sections
of the appropriate pullback bundles.  All scaled norms below are
pulled back through the fixed backgrounds \(\Theta_*,\Phi_*\).  The
graph finite differences are equipped with the already defined
augmented distance \eqref{eq:prepared-Banach-norm}; the underlying
Banach model is \eqref{eq:prepared-model-Banach-norm}.  Raw \(R,F\)
coordinates remain
available at a single fixed \(\tau_0\), but their change of coordinates
is not used for any entrance-time-uniform operator estimate.
The structural identities in \eqref{eq:prepared-state-graph}, the
fixed marked chart, and the degree-one properness conditions are part
of the chart, rather than additional equations imposed later.  The
following lemma records the functional-analytic calculus used below.

\begin{lemma}[Shifted-annulus soliton calculus]
\label{lem:shifted-annulus-soliton-calculus}
Fix a finite order \(m\), the radial-comparison and map-and-inverse
ceilings in \(\mathfrak P_{\rm prep}\), and numbers
\[
 0<c_\lambda<C_\lambda<\infty .
\]
For \(\tau\geq0\), \(L\in\mathscr L_{\rm pre}\), and every relative
map \(\psi\) in this controlled class, the soliton flow satisfies
\begin{equation}\label{eq:soliton-shifted-annuli}
 \varphi_\tau\circ\psi(A_L)
 \subset
 \{c e^\tau L<\bar f<C e^\tau L\},
 \qquad
 \psi\circ\varphi_{-\tau}(A_{e^\tau L})
 \subset
 \{cL<\bar f<CL\},
\end{equation}
after fixed enlargements of the annuli.  In harmonic coordinates for
\((A_L,L^{-1}\bar g)\) and the corresponding output annulus with
metric \((e^\tau L)^{-1}\bar g\), the maps in
\eqref{eq:soliton-shifted-annuli}, their inverses, and their derivatives
through order \(m\) have one uniform bound.

If
\begin{equation}\label{eq:shifted-annulus-scale-bracket}
 c_\lambda\leq\lambda e^\tau\leq C_\lambda,
\end{equation}
then \(\lambda(\varphi_\tau\circ\psi)^*\bar g\) has a uniform
scale-one coefficient package on \(A_L\).  If \(R,F\) and their
inverses lie in the same controlled class, put
\begin{equation}\label{eq:prepared-conjugated-map-factors}
 \begin{aligned}
  \psi_{\rm rel}&:=R\circ F^{-1}\circ R^{-1},\\
  J_\tau(\mathbf z)&:=
    \Theta(\tau,\mathbf z)\circ\Phi(\tau,\mathbf z)^{-1}
    =\varphi_\tau\circ\psi_{\rm rel}\circ\varphi_{-\tau},\\
  \zeta_\tau(\mathbf z)&:=
    (1-\eta)\circ\Phi(\tau,\mathbf z)^{-1}.
 \end{aligned}
\end{equation}
On \(\supp\zeta_\tau(\mathbf z)\), the map \(J_\tau(\mathbf z)\) and
its inverse preserve the current annular level up to fixed factors, and
\(J_\tau(\mathbf z)^{\pm1}\) and \(\zeta_\tau(\mathbf z)\) have
uniform scale-one \(C^{m,\alpha}\) bounds.  More precisely, an output
annulus \(A_{L_{\rm out}}\) meeting
\(\supp\zeta_\tau(\mathbf z)\) has
\[
 L_{\rm out}\geq c\Gamma e^\tau,\qquad
 L_{\rm in}:=e^{-\tau}L_{\rm out}\geq c\Gamma,
\]
and all estimates are obtained directly by passing through the
pre-radius charts with label comparable to \(L_{\rm in}\), not by a
chain of annular transitions.
\end{lemma}

\begin{proof}
Along the radial flow,
\[
 \frac d{d\tau}(\bar f\circ\varphi_\tau)
 =|\bar\nabla\bar f|^2\circ\varphi_\tau
 =(\bar f-\bar R)\circ\varphi_\tau .
\]
The FIK curvature is bounded and its AC derivatives have the symbol
decay fixed in Lemma~\ref{lem:FIK-AC-symbol}.  Integration gives
\[
 c e^\tau(1+\bar f(x))
 \leq1+\bar f(\varphi_\tau x)
 \leq C e^\tau(1+\bar f(x))
\]
on the end.  Differentiating the flow equation in the AC
scale-one charts, or equivalently integrating its variational
equations after the domain and range rescalings, gives the asserted
all-order symbol bounds.  Composition with a controlled relative map
preserves them by the two-sided radial comparison.  The tensor
scaling in \eqref{eq:shifted-annulus-scale-bracket} cancels the
factor \(e^\tau\) in the output metric scale.

The factorization \eqref{eq:prepared-conjugated-map-factors} is exact.
The two occurrences of the soliton flow shift from
\(L_{\rm out}\) to \(e^{-\tau}L_{\rm out}\) and back, whereas
\(\psi_{\rm rel}\) is a same-level controlled map.  Thus
\(J_\tau(\mathbf z)\) is a same-output-level map even though neither
\(\Theta(\tau,\mathbf z)\) nor \(\Phi(\tau,\mathbf z)\) is.
Finally, \(\zeta_\tau(\mathbf z)\ne0\) implies that
\(\Phi(\tau,\mathbf z)^{-1}(A_{L_{\rm out}})\) meets
\(\{\bar f\geq c\Gamma\}\); the preceding shifted comparison gives the
lower bounds for \(L_{\rm out}\) and \(L_{\rm in}\).  The normalized
cutoff bounds and the direct pre-radius coordinate changes give the
claimed estimates for \(\zeta_\tau(\mathbf z)\) and complete the proof.
\end{proof}

\begin{lemma}[Prepared-chart calculus]
\label{lem:prepared-chart-calculus}
Fix \(k\geq3\), \(0<\alpha<1\), and the polynomial-loss index \(N\).
 Let
 \(\mathbf z_*=(G_*,\lambda_*,R_*,F_*)
 \in\mathscr P_{\tau_0}^{k+2,\alpha}\) be a prepared background
 state.  Thus \(G_*\in C^{k+2,\alpha}\), while
 \(R_*^{\pm1},F_*^{\pm1}\) have the scaled
 \(C^{k+3,\alpha}\) bounds required by that prepared class, and
 \(R_*,F_*\) are degree-one proper diffeomorphisms.
Assume that, for some constants \(0<c_{\rm rad}\leq C_{\rm rad}<\infty\),
\[
 c_{\rm rad}(1+\bar f)
 \leq 1+\bar f\circ\psi_*
 \leq C_{\rm rad}(1+\bar f),
 \qquad
 \psi_*\in\{R_*,F_*,R_*^{-1},F_*^{-1}\},
\]
and that all ellipticity, support-separation, graft, and
bounded-geometry conditions hold with positive margin.  For a single
fixed \(\tau_0\) no uniformity in that parameter is needed.  Whenever
the conclusion is used uniformly over entrance times, assume in
addition the package scale bracket
\begin{equation}\label{eq:prepared-chart-physical-scale-bracket}
 0<c_{\rm scl}\leq\lambda_*e^{\tau_0}
 \leq C_{\rm scl}<\infty ;
\end{equation}
after shrinking the \(\ell\)-chart the same inequalities hold with
relaxed fixed constants.  Set
\[
 \begin{aligned}
 \mathscr X_{\rm in}^{k+2,\alpha}
 &:=
 h^{k+2,\alpha}(S^2T^*\mathcal X)
 \times\mathbb R
 \times
 \bigl(\mathfrak X_{\rm sc}^{k+3,\alpha}\bigr)^2,\\
 \mathscr X_{\rm out}^{k,\alpha}
 &:=
 h^{k,\alpha}(S^2T^*\mathcal X)
 \times\mathbb R
 \times
 \bigl(\mathfrak X_{\rm sc}^{k+1,\alpha}\bigr)^2
 \times\mathfrak T_{{\rm sc},N}^{k,\alpha}.
 \end{aligned}
\]
 Write
 \[
  \widehat G:=\lambda_*^{-1}(G-G_*),\qquad
  G=G_*+\lambda_*\widehat G .
 \]
 After shrinking a neighborhood of
 \((0,0,0,0)\) in \(\mathscr X_{\rm in}^{k+2,\alpha}\), the assignment
 \[
  (\widehat G,\ell,X_\Theta,X_\Phi)
  \longmapsto
  \bigl(\widehat G,\ell,X_\Theta,X_\Phi,h(\mathbf z)\bigr),
 \qquad
 \lambda=\lambda_*e^\ell,
\]
with \(R,F,\Theta,\Phi\) and \(h(\mathbf z)\) defined by
\eqref{eq:prepared-state-graph}, is a \(C^1\) map into
\(\mathscr X_{\rm out}^{k,\alpha}\).  In the same charts:
\begin{enumerate}
\item inversion and composition of the charted diffeomorphisms are
      \(C^1\) from vector-field input order \(k+3\) to output order
      \(k+1\);
\item pullback is \(C^1\) from a charted diffeomorphism of order
      \(k+3\) and a scaled tensor of order \(k+2\) to a scaled tensor
      of order \(k\);
\item multiplication maps
      \[
       \mathfrak C_{{\rm sc},N}^{k,\alpha}
       \times\mathfrak C_{{\rm sc},0}^{k,\alpha}
       \longrightarrow
       \mathfrak C_{{\rm sc},N}^{k,\alpha}
      \]
      continuously, and the analogous tame estimates hold with
      \(\mathfrak T_{{\rm sc},N}^{k,\alpha}\) in the first factor.
 \end{enumerate}
For a fixed \(\tau_0\), the same assertion in the unweighted raw
\(G\)-coordinate is \(C^1\), but its operator constant and chart
radius may depend on \((\tau_0,\lambda_*)\).  Uniform entrance-time
families are measured in the normalized coordinate above, equivalently
in the conjugated product tangent norm
\begin{equation}\label{eq:uniform-prepared-graph-tangent-norm}
 \|\dot{\mathbf z}\|_{\rm con,*}:=
 \lambda_*^{-1}\|\dot G\|_{C^{k+2,\alpha}}
 +|\dot\ell|
 +\|\dot X_\Theta\|_{\mathfrak X_{\rm sc}^{k+3,\alpha}}
 +\|\dot X_\Phi\|_{\mathfrak X_{\rm sc}^{k+3,\alpha}} .
\end{equation}
More precisely,
\begin{equation}\label{eq:uniform-prepared-graph-derivative-bound}
 \|Dh(\mathbf z_*)[\dot{\mathbf z}]\|
 _{\mathfrak T_{{\rm sc},N}^{k,\alpha}}
 \leq C\|\dot{\mathbf z}\|_{\rm con,*},
\end{equation}
with \(C\) independent of the entrance time on the displayed package
family.  Consequently the full graph tangent norm obtained by adding
the left side of
\eqref{eq:uniform-prepared-graph-derivative-bound} is uniformly
equivalent to \(\|\cdot\|_{\rm con,*}\); this equivalence is a
conclusion, not its definition.
No entrance-time-uniform operator bound from the unweighted raw
\(G,R,F\) product norm is asserted.
For every \(L\in\mathscr L_\Gamma\), the operator norms of these maps
and of their first derivatives, after restriction to an enlarged
annulus and rescaling by \(L^{-1}\bar g\), are bounded independently of
\(L\).  The bounds are uniform on any family lying in one fixed
prepared chart for which the positive margins and the constants in the
two-sided radial comparisons are common, the scaled
\(C^{k+3,\alpha}\) norms of every background map and inverse entering
\eqref{eq:prepared-state-graph} are bounded by one common constant, and
the scaled \(C^{k+2,\alpha}\) norms of the background tensor
coefficients are bounded by one common constant.  The coordinate
operations in \textup{(1)}--\textup{(3)} are \(C^1\) with the displayed
input and output orders.  The two-order input buffer also controls the
first two derivatives with respect to the finite-dimensional phase
parameter used below.  No \(C^\infty\) assertion is made on a fixed
finite-regularity Banach space.

On every bounded common-margin ball these same finite coordinate
operations satisfy the mixed tame Taylor estimate
\begin{equation}\label{eq:prepared-mixed-tame-remainder}
 \|\mathcal N(\mathbf z+\mathbf u)-\mathcal N(\mathbf z)
      -D\mathcal N(\mathbf z)[\mathbf u]\|_{{\rm out},k}
 \leq C_K
   \|\mathbf u\|_{{\rm in},k+2}
   \|\mathbf u\|_{{\rm in},k}.
\end{equation}
Here \(\mathcal N\) is any graph, local-addition, inversion,
composition, pullback, or product operation occurring in the prepared
state; the high input norm is the one displayed above, and the low
input norm is obtained by lowering the tensor order from \(k+2\) to
\(k\) and each vector-field order from \(k+3\) to \(k+1\).
For a multi-input operation the right side is the sum of the analogous
one-high--one-low products.  Estimate
\eqref{eq:prepared-mixed-tame-remainder} is not a same-order
\(C^2\) assertion: it is precisely the one-high--one-low remainder
bound allowed by the two-order buffer.

 For every order \(r\geq3\) for which
 \(\mathbf z_*\in\mathscr P_{\tau_0}^{r,\alpha}\), the
 same-output soliton-conjugated coordinate variables define that Banach manifold
 near \(\mathbf z_*\).  With the displayed two-order buffer, its coordinate
graph from input order \(k+2\) into output order \(k\) is \(C^1\).
At every fixed order \(r\), the same graph operations are continuous
from \(C^{r,\alpha}\) prepared inputs to
\(\mathscr P_{\tau_0}^{r,\alpha}\); the two-order loss is needed for
the asserted Fr\'echet derivative, not for continuity.  Consequently
a bounded \(C^{r,\alpha}\) input ball with common margins maps to a
bounded common-margin ball in \(\mathscr P_{\tau_0}^{r,\alpha}\).
Within each \(\mathscr P_{\tau_0}^{r,\alpha}\), the loci on which the
two-sided radial comparisons, map-and-inverse bounds, ellipticity,
support separation, Gram invertibility, graft compatibility, and
bounded-geometry inequalities retain positive common margins are open.
The associated fixed non-strict controlled sublevels need not be open
and are not used as manifold domains.  Every map in a sufficiently
small common-margin chart is properly homotopic to its background map
and therefore remains in the degree-one proper homotopy class.  In
particular, the strict prepared states form an open subset.
\end{lemma}

\begin{proof}
Choose constants \(0<c_*<1/2<4<C_*<\infty\), depending only on the
common constants \(c_{\rm rad},C_{\rm rad}\) and the radius of the
prepared chart, so that the relative maps
\[
 R^{\pm1},\quad F^{\pm1},\quad
 R\circ F,\quad F^{-1}\circ R^{-1}
\]
occurring in \eqref{eq:prepared-state-graph}, and every corresponding
map in the shrunk chart, send \(A_L\) into
\[
 A_{L,*}^+=\{c_*L<\bar f<C_*L\}.
\]
This same-annulus assertion is deliberately not made for
\(\Theta=\varphi_{\tau_0}\circ R\) or
\(\Phi=\varphi_{\tau_0}\circ R\circ F\).  Those maps shift the level by
the factor \(e^{\tau_0}\) and are handled by
Lemma~\ref{lem:shifted-annulus-soliton-calculus}.
The asymptotically conical bounds give uniformly bounded geometry to
\((A_{L,*}^+,L^{-1}\bar g)\).  These annuli admit scale-one harmonic
coordinate covers of bounded multiplicity, with uniform
\(C^{k+2,\alpha}\) bounds for the metric coefficients.  The compact region
\(\{\bar f<4\Gamma\}\) is handled by one fixed finite atlas.  A small
ball in \(\mathfrak X_{\rm sc}^{k+3,\alpha}\) moves \(A_L\) only inside
\(A_{L,*}^+\).  Thus every calculation may be made in a fixed Euclidean
ball after rescaling, with constants independent of \(L\).

In one such chart, the product estimate follows directly from the
Leibniz rule and the H\"older product inequality.  Composition follows
by differentiating \(u\circ\psi\) and applying the chain rule at each
order.  Inversion follows inductively from
\[
 D(\psi^{-1})
 =\bigl((D\psi)^{-1}\bigr)\circ\psi^{-1},
\]
and pullback from the coordinate formula expressing \(\psi^*u\) as
\(u\circ\psi\) multiplied by copies of \(D\psi\).  Their first
variations are obtained by the same formulas; in particular, along a
curve \(\psi_s\),
\[
 \frac d{ds}\bigg|_{s=0}\psi_s^*u_s
 =\psi_0^*\bigl(\dot u_0+\Lie_Wu_0\bigr),
 \qquad
 W=\dot\psi_0\circ\psi_0^{-1}.
\]
 Consequently an output of order \(k\) uses at most input order \(k+2\)
 in every twice-differentiated phase expression occurring below.  The
 same estimates hold in the little-H\"older completions by approximation
 and closure.

For \eqref{eq:prepared-mixed-tame-remainder}, subtract the displayed
first-variation formula from the exact coordinate formula.  Every
remaining summand contains at least two increment factors.  After
differentiating to output order \(k\), put the factor carrying the
largest number of derivatives in the high norm and the other increment
factor in the two-order-lower norm; all remaining factors are bounded
in the fixed \(C^2\) box.  The H\"older product inequality gives
\[
 C_K\|\mathbf u\|_{{\rm in},k+2}
       \|\mathbf u\|_{{\rm in},k}.
\]
For inversion, use the exact identity for the difference of two matrix
inverses before differentiating; for composition and pullback, use the
integral form of the chain rule along the local-addition segment.
These identities have the same one-high--one-low allocation, uniformly
in every rescaled chart.  This proves
\eqref{eq:prepared-mixed-tame-remainder}; density passes it to the
little-H\"older completions.

On \(A_L\), the factor \(L^{-N}\) in
\eqref{eq:scaled-tensor-holder} is constant.  The relative maps just
listed preserve its annular label up to a fixed factor.  For the full
prepared graph use instead the exact decomposition
\begin{equation}\label{eq:prepared-graph-shifted-decomposition}
 h
 =\lambda^{-1}(\Phi^{-1})^*(\eta\,\iota_*G)
   +\zeta J^*\bar g-\bar g,
 \qquad
 \zeta=(1-\eta)\circ\Phi^{-1},
 \quad J=\Theta\circ\Phi^{-1}.
\end{equation}
In the conjugated coordinates, put
\[
 J_*:=\Theta_*\circ\Phi_*^{-1},\qquad
 \zeta_*:=(1-\eta)\circ\Phi_*^{-1}.
\]
Then the variable map dependence factors exactly as
\begin{equation}\label{eq:prepared-graph-conjugated-variable-factors}
 \begin{gathered}
  J=Q_\Theta\circ J_*\circ Q_\Phi^{-1},\qquad
  \zeta=\zeta_*\circ Q_\Phi^{-1},\\
  \lambda^{-1}(\Phi^{-1})^*(\eta\,\iota_*G)
  =
  \lambda^{-1}(Q_\Phi^{-1})^*
   (\Phi_*^{-1})^*(\eta\,\iota_*G).
 \end{gathered}
\end{equation}
Thus every variable composition, inverse, and pullback in the graph is
same-output-level.  The soliton shifts occur only in the fixed
background coefficients \(J_*,\zeta_*,\Phi_*\).
The compact and end estimates play different roles.  On the compact
 support of \(\eta\), direct differentiation gives
 \[
  D_Gh(\mathbf z_*)[\dot G]
  =\lambda_*^{-1}(\Phi_*^{-1})^*
      (\eta\,\iota_*\dot G).
 \]
 Hence the raw compact-core operator norm may grow like
 \(\lambda_*^{-1}\); it is harmless at fixed \(\tau_0\), and is
 uniform in entrance-time families precisely after replacing
 \(\dot G\) by the normalized variable
 \(\dot{\widehat G}=\lambda_*^{-1}\dot G\), as in
 \eqref{eq:uniform-prepared-graph-tangent-norm}.  No shifted-annulus
 argument is used to conceal this compact contribution.
For the map directions,
\eqref{eq:prepared-graph-conjugated-variable-factors} and the ordinary
same-level inverse, composition, and pullback estimates give
\[
 \|D_{(X_\Theta,X_\Phi)}h[\dot X_\Theta,\dot X_\Phi]\|
 _{\mathfrak T_{{\rm sc},N}^{k,\alpha}}
 \leq
 C\bigl(
  \|\dot X_\Theta\|_{\mathfrak X_{\rm sc}^{k+3,\alpha}}
 +\|\dot X_\Phi\|_{\mathfrak X_{\rm sc}^{k+3,\alpha}}\bigr).
\]
At an exact center, the preceding identity contains no term
\(\Lie_{(\varphi_{\tau_0})_*X_F}\bar g\).

 On the end, the first term is estimated between the shifted physical
 source and normalized output charts.  There the uniform estimate is
 exactly the shifted tensor-pullback estimate supplied by
 \eqref{eq:prepared-chart-physical-scale-bracket} and
 Lemma~\ref{lem:shifted-annulus-soliton-calculus}.  The second term uses
 the same-output-level conjugated map \(J_{\tau_0}(\mathbf z)\) from
\eqref{eq:prepared-conjugated-map-factors}.  Thus the weight is
preserved in each term of \eqref{eq:prepared-graph-shifted-decomposition}.
In every product the second factor is controlled in the unweighted
\(\mathfrak C_{{\rm sc},0}^{k,\alpha}\) norm, so no additional
polynomial loss is created.  Taking the dyadic supremum and adding the
fixed-core estimates proves the asserted global bounds.  In
particular, no estimate in this proof treats \(\Theta\) or \(\Phi\) as
a same-annulus variable map before the fixed soliton factors have been
removed.  Together with the scale derivative, the preceding estimates
prove \eqref{eq:uniform-prepared-graph-derivative-bound}.

Finally, retain \(r_{\rm la}=(1+\bar f)^{1/2}\).  The soliton identity
\(|\bar\nabla\bar f|^2\leq\bar f\) gives
\[
 |d\log r_{\rm la}|_{\widehat g_{\rm la}}
 =|dr_{\rm la}|_{\bar g}\leq\frac12 .
\]
If \(\psi_*\) denotes either background map and
\[
 \psi_s=\operatorname{Exp}(sX)\circ\psi_*,
 \qquad 0\leq s\leq1,
\]
then smallness in \(\mathfrak X_{\rm sc}^{1,\alpha}\) gives
\[
 |X(y)|_{\widehat g_{\rm la}}
 =r_{\rm la}(y)^{-1}|X(y)|_{\bar g}\leq C\varepsilon .
\]
The defining curve
\(s\mapsto\operatorname{Exp}_y(sX(y))\) is a
\(\widehat g_{\rm la}\)-geodesic lying in
\(\mathscr U_{\rm la}\).  Integration of the preceding logarithmic
gradient bound gives
\[
 e^{-C\varepsilon}r_{\rm la}(y)
 \leq r_{\rm la}(\operatorname{Exp}_y(sX(y)))
 \leq e^{C\varepsilon}r_{\rm la}(y).
\]
After decreasing \(\varepsilon\), the background radial comparison
therefore yields constants independent of \(s\) such that
\[
 c(1+\bar f(x))
 \leq1+\bar f(\psi_s(x))
 \leq C(1+\bar f(x)).
\]
Thus \((s,x)\mapsto\psi_s(x)\) is a proper homotopy.  In particular,
\(\deg\psi_1=\deg\psi_*=1\).  Smallness in the scaled \(C^1\) norm
also preserves local invertibility, so \(\psi_1\) is a proper local
diffeomorphism.  It is therefore a finite covering of the connected
model \(M\); degree one makes the covering one-sheeted, and hence
\(\psi_1\) is a global diffeomorphism.

Given a center in the ambient map class, first relax its radial and
map-norm constants and retain half of the resulting radial,
map-and-inverse, and local-invertibility margins.  The estimates above
show that a sufficiently small right-translated ball preserves those
relaxed bounds, properness, and degree.  All remaining common-margin
conditions are strict inequalities in the displayed Banach norms and
are therefore open.  Thus these common-margin variables, rather than
the boundary of a fixed controlled sublevel, form an open subset of the
corresponding Banach product and give the asserted manifold charts.
The graph formula
\eqref{eq:prepared-state-graph} is a finite composition of the maps
just estimated, which proves the \(C^1\) graph assertion with the
two-order buffer.
\end{proof}

In the differentiability statements below, \(k_0\geq12\) denotes the
continuation output order and the prepared input order is \(k_0+2\);
when no differentiation is taken we write \(k=k_0\).  A
\emph{common-margin prepared ball} at displayed order \(r\) means a
ball
\[
 \mathscr B\subset\mathscr P_{\tau_0}^{r,\alpha}
\]
on which the ellipticity, two-sided radial comparison, support
separation, Gram invertibility, graft, and bounded-geometry
inequalities retain common positive margins, and for which one number
\[
 K_{\rm init}^{r,\alpha}=K_{\rm init}^{r,\alpha}(\mathscr B)\geq1
\]
has been fixed that bounds all scaled metric, coefficient, graph,
map, inverse-map, and local-addition coordinate norms occurring at
the displayed order.  The choice of
\(K_{\rm init}^{r,\alpha}\) is part of the data whenever uniformity on
\(\mathscr B\) is asserted.  All maps in the ball lie in the
degree-one proper homotopy class of its center.
More quantitatively, every uniform statement uses the numerical
prepared package \(\mathfrak P_{\rm prep}\) fixed, after its rate pair,
before the entrance ball and
Theorem~\ref{thm:intro-sharp-scattering} are invoked, as prescribed in
\eqref{eq:numerical-prepared-package}.  Its entries have the meanings
listed there; in particular the adaptive phase budget, entrance
threshold, and scale-comparison constant are part of the package, not
choices made during a later restart.
Constants described as uniform on a common-margin ball may depend on
\[
 k,\alpha,N,\quad\text{the fixed background},\quad
 \mathfrak P_{\rm prep},\quad K_{\iota,r+2}^{\rm fix},\quad\text{and}\quad
 K_{\rm init}^{r,\alpha},
\]
but not on any further feature of the ball or on \(\tau_0\).

\begin{lemma}[Static prepared columns]
\label{lem:static-prepared-columns}
Fix \(k\geq3\) and one numerical prepared package
\(\mathfrak P_{\rm prep}\).  There is \(\tau_*>0\), depending only on
that package and the fixed background data, with the following
package-local property.  For every \(\tau_0\geq\tau_*\), every
prepared center satisfying the structural graph identities and
\(\mathfrak P_{\rm prep}\), and every
sufficiently small common-margin ball
\(\mathscr B^k\subset\mathscr P_{\tau_0}^{k,\alpha}\) about that center
which retains the same package, form the effective columns at the
 single time \(\tau_0\) by
 \eqref{eq:effective-column-zero}--%
 \eqref{eq:effective-column-j}.  Then, for every \(0\leq j\leq8\),
 the \(j\)-th column map is continuous on
 \(\mathscr B^k\), and, for every \(\mathbf z\in\mathscr B^k\) and
 \(0\leq m\leq3\), with \(0\leq\mu\leq8\) in the last estimate,
\begin{align}
 \sup_M\sum_{\ell=0}^m(1+\bar f)^{\ell/2}
 |\bar\nabla^\ell\mathcal Y_{j,\tau_0}(\mathbf z)|
 &\leq C_m,
 \label{eq:static-column-scale-bound}\\
 \|\mathcal Y_{j,\tau_0}(\mathbf z)-Y_j\|_{H^m_\nu}
 &\leq C_m e^{-c_me^{\tau_0}},
 \label{eq:static-column-Gaussian-tail}\\
 \left|
 \ip{\rho_{\tau_0}\mathcal Y_{j,\tau_0}(\mathbf z)}{Z_\mu}
 -\ip{Y_j}{Z_\mu}
 \right|
 &\leq Ce^{-ce^{\tau_0}}.
 \label{eq:static-column-Gram}
\end{align}
with constants independent of the center and of
\(\tau_0\geq\tau_*\).

For the differentiability assertion, let
\(\mathscr B^{k+2}\subset\mathscr P_{\tau_0}^{k+2,\alpha}\) be any
common-margin ball satisfying the same package whose canonical
 inclusion is contained in \(\mathscr B^k\).  For every
 \(0\leq j\leq8\), the \(j\)-th column map is \(C^1\)
 from \(\mathscr B^{k+2}\) to \(H^m_\nu\) for each
 \(0\leq m\leq3\), and, for
\(\mathbf z_1,\mathbf z_2\in\mathscr B^{k+2}\),
\begin{equation}\label{eq:static-column-Lipschitz}
 \|\mathcal Y_{j,\tau_0}(\mathbf z_1)
       -\mathcal Y_{j,\tau_0}(\mathbf z_2)\|_{H^m_\nu}
 \leq
 C_me^{-c_me^{\tau_0}}
 \|\mathbf z_1-\mathbf z_2\|_{\mathscr X_{\rm prep}^{k+2,\alpha}} .
\end{equation}
\end{lemma}

\begin{proof}
At the background state the transported outer term
\[
 K_{\tau_0}(T)
 =(\Phi^{-1})^*((1-\eta)\Theta^*T)
\]
is supported in
\(\{\bar f\geq c\Gamma e^{\tau_0}\}\).  This support separation,
with a slightly smaller \(c\), is one of the open conditions defining
  \(\mathscr B^k\).  Use the exact factorization
\[
 J=\Theta\circ\Phi^{-1},\qquad
 \zeta=(1-\eta)\circ\Phi^{-1},\qquad
 K_{\tau_0}(T)=\zeta J^*T .
\]
Lemma~\ref{lem:shifted-annulus-soliton-calculus} gives uniform bounds
for \(J^{\pm1}\), \(\zeta\), and all their derivatives through order
\(k\) on every rescaled output annulus meeting this support.  The two
soliton-flow shifts cancel inside \(J\), and the estimate passes
directly through the pre-radius input annulus with
\(L_{\rm in}\simeq e^{-\tau_0}L_{\rm out}\geq c\Gamma\); it does not
use a same-annulus estimate for either \(\Theta\) or \(\Phi\).
More precisely, the map factors are used through order
\(m+1\leq4\), as required by the \(m\)-fold pullback formula; these
derivatives are available from the \(k+1\) map component of
\(\mathscr P_{\tau_0}^{k,\alpha}\).  The polynomial FIK bounds for
\(Y_j\), together with
\eqref{eq:scaled-tensor-holder}--%
\eqref{eq:scaled-vector-holder}, prove
\eqref{eq:static-column-scale-bound}.  Gaussian integration on the
receding support controls the transported term.  The complete
difference, however, has two terms:
\begin{equation}\label{eq:static-column-complete-tail-decomposition}
 \begin{aligned}
  \mathcal Y_{0,\tau_0}-Y_0
  &=-K_{\tau_0}(Y_0),\\
  \mathcal Y_{j,\tau_0}-Y_j
  &=
  \bigl(\Lie_{\chi_{\tau_0}W_j}\bar g-Y_j\bigr)
  -K_{\tau_0}
    \bigl(\Lie_{\chi_{\tau_0}W_j}\bar g\bigr),
  \qquad1\leq j\leq8 .
 \end{aligned}
\end{equation}
Lemma~\ref{lem:cutoff-tails}, including its normalized
cutoff-derivative estimates, controls the first summand for \(j\geq1\);
the preceding shifted-annulus calculation controls the second.  Hence
\eqref{eq:static-column-Gaussian-tail} follows for every \(j\).  For the
Gram estimate use the exact identity
\[
 \left\langle\rho_{\tau_0}\mathcal Y_{j,\tau_0},Z_\mu\right\rangle
 -\langle Y_j,Z_\mu\rangle
 =
 \left\langle
  \rho_{\tau_0}(\mathcal Y_{j,\tau_0}-Y_j),Z_\mu
 \right\rangle
 -
 \left\langle(1-\rho_{\tau_0})Y_j,Z_\mu\right\rangle .
\]
Both terms are Gaussian tails, proving
\eqref{eq:static-column-Gram}.

The preceding value estimates use only the continuous fixed-order
prepared calculus.  On \(\mathscr B^{k+2}\),
Lemma~\ref{lem:prepared-chart-calculus} gives the inverse, pullback,
product, and composition operations as \(C^1\) maps from input order
\(k+2\) to output order \(k\).
Their first derivatives have the same
scale-invariant bounds and are supported in the same receding region.
The direct moving-cutoff summand in
\eqref{eq:static-column-complete-tail-decomposition} is independent of
\(\mathbf z\) and cancels in differences.
The Gaussian tail therefore multiplies the ordinary Banach
Lipschitz estimate and gives
\eqref{eq:static-column-Lipschitz}.  No finite phase action,
harmonic-map evolution, or bootstrap estimate is used.
\end{proof}

\subsection{Local existence of the coupled feedback system}
\label{subsec:coupled-local-existence}

The algebraic feedback theorem, Theorem~\ref{thm:receding}, does not by
itself construct the geometric evolution to which it is applied.  There
are two functional-analytic points which are slightly nonstandard here.
First, the harmonic-map equation is posed on a complete asymptotically
conical manifold and its coefficients depend on the moving graft.
Second, the nine feedback coefficients are nonlocal functionals of the
current state.  We record the linear estimate and the regularity of the
feedback functional before carrying out the fixed-point argument.

\subsubsection{Time-dependent prepared spaces}
\label{subsec:time-dependent-prepared-spaces}

Fix hereafter \(k\geq12\), \(0<\alpha<1\), and the polynomial-loss
exponent \(N\geq0\).  Fix one base entrance in the common-margin ball.
On the first local interval take \(s=\tau_0\) and use that entrance as
the reference state.  At a later restart time \(s>\tau_0\), the
construction is inductive: a distinguished reference solution
\(\mathbf z^\circ(\tau)\) has already been constructed through \(s\),
and we freeze
\[
 \mathbf z_s^\circ:=\mathbf z^\circ(s),\qquad
 \lambda^\circ(s):=\lambda(\mathbf z_s^\circ),\qquad
 \acute G_\circ(s):=\acute G(\mathbf z_s^\circ).
\]
The source metric, harmonic atlas, bundle trivializations, and
same-output soliton-conjugated map charts at this single frozen state are used for
every nearby trial state on \(I=[s,s+\delta]\).  Thus none of the clocks
or identifications below depends on the trial state.  This convention
also specifies the ``center path'' at every later occurrence.
Fix a common-margin prepared ball
\[
 \mathscr B_s\subset\mathscr P_s^{k+2,\alpha}
\]
about \(\mathbf z_s^\circ\) (at the first interval
\(\mathscr B_s=\mathscr B_{\tau_0}\)); write simply \(\mathscr B\) when
the restart time is clear.  Let
\(I=[s,s+\delta]\), \(0<\delta\leq1\).  In addition to the state distance
in \eqref{eq:prepared-Banach-norm}, for every integer
\(0\leq r\leq k+2\) we use
\begin{equation}\label{eq:extended-prepared-distance}
 \begin{split}
  \|\mathbf z_1-\mathbf z_2\|_{\mathscr Y_{\rm prep}^{r,\alpha}}
  :={}&
  \frac{|t_1-t_2|}{\lambda^\circ(s)}
  +\|G_1-G_2\|_{C^{r,\alpha}(\mathcal X)}
  +\left|\log\frac{\lambda_1}{\lambda_2}\right|\\
 &+\|X_{\Theta,1}-X_{\Theta,2}\|_{\mathfrak X_{\rm sc}^{r+1,\alpha}}
  +\|X_{\Phi,1}-X_{\Phi,2}\|_{\mathfrak X_{\rm sc}^{r+1,\alpha}}\\
  &+\|h(\mathbf z_1)-h(\mathbf z_2)\|
       _{\mathfrak T_{{\rm sc},N}^{r,\alpha}} .
 \end{split}
\end{equation}
Here and below all map differences are taken in the one frozen
same-output exponential chart, as in \eqref{eq:prepared-Banach-norm}.  On a
common-margin ball the resulting distances are uniformly equivalent
when the center of the exponential chart is changed.
The factor \((\lambda^\circ(s))^{-1}\) in the compact clock records
dimensionless elapsed compact time.  The raw compact metric retains its
unscaled physical-host norm, which controls the order-one exterior; the
scale-covariant collapsing-core control is instead carried by the affine
normalized compact increment in the frozen scale-adapted atlas defined
below.
The normalized graph term \(h\) records the remaining
\(\lambda^{-1}\)-sensitive graft contribution.  This dual tier avoids
placing the entire order-one exterior metric in a
\(\lambda^{-1}\)-weighted fixed-host norm.

We next define the parabolic norms used for the
noncompact map component.  The clock must follow the physical source
scale; the raw \(F\)-equation is not uniformly parabolic in the
annular clock \((\tau-s)/L\).

At the left endpoint \(s\), use the trial-independent reference scale
\[
 \lambda_s^\circ(\tau)
 :=\lambda^\circ(s)e^{-(\tau-s)},\qquad
 t_s^\circ(\tau):=\int_s^\tau\lambda_s^\circ(q)\,dq .
\]
Here \(\lambda^\circ(s)\) is the scale of the frozen reference state
\(\mathbf z_s^\circ\) specified above; at the first interval it is the
scale of the base entrance.  The superscript \({}^\circ\) denotes only
this restart reference and never a trial-dependent quantity.
After the fixed-point ball and \(\delta\) are decreased, every trial
scale satisfies
\begin{equation}\label{eq:reference-clock-comparison}
 C^{-1}\leq\frac{\lambda(\tau)}{\lambda_s^\circ(\tau)}
 \leq C,\qquad \tau\in I,
\end{equation}
with the same statement for first variations after the ratio is
differentiated.

The fixed base source metric \(\acute G_\circ(s)\) and its
scale-normalized bounded geometry supply a
uniformly locally finite source-harmonic cover
\[
 \mathfrak A_s
 =\{(\mathcal U,r_{\mathcal U},L_{\mathcal U})\}.
\]
The radius \(r_{\mathcal U}\) is the local physical harmonic scale,
\(L_{\mathcal U}=1\) on the finite collapsing-core subcover, and
\(L_{\mathcal U}\in\mathscr L_\Gamma\) labels the corresponding
rescaled model annulus elsewhere.  The cover can be chosen so that
\begin{equation}\label{eq:source-adapted-atlas}
 r_{\mathcal U}^2\geq c\lambda_s^\circ(s),\qquad
 r_{\mathcal U}^{-2}\acute G_\circ(s)
 \ \text{has one common harmonic-chart package},
\end{equation}
neighboring radii and labels are comparable, and the chart norms are
uniformly equivalent to the core and annular norms in
\eqref{eq:scaled-vector-holder}.  These facts follow from the complete
source harmonic-radius package, the two-sided radial comparison, and
the prepared metric identity; they are retained as part of the
numerical prepared package.  In particular, on every trial state,
\begin{equation}\label{eq:source-adapted-ellipticity}
 r_{\mathcal U}^2\frac{\lambda}{\lambda_s^\circ}\acute G^{-1}
\end{equation}
is uniformly elliptic in the
\(r_{\mathcal U}^{-2}\acute G_\circ(s)\)-harmonic coordinates.  This
single statement includes normalized-scale charts near the collapsing
core and physical-scale charts on the noncollapsing exterior.
After the ball is shrunk, all trial source metrics, target bundles, and
same-output soliton-conjugated map charts are uniformly equivalent to
these single base choices, including for first variations.

Use the frozen, trial-independent clock
\begin{equation}\label{eq:source-adapted-clock}
 \vartheta_{\mathcal U}
 :=\frac{t_s^\circ(\tau)-t_s^\circ(s)}{r_{\mathcal U}^2},
 \qquad
 J_{\mathcal U}
 :=\left[
  0,\frac{t_s^\circ(s+\delta)-t_s^\circ(s)}
          {r_{\mathcal U}^2}
 \right].
\end{equation}
Since \(\lambda_s^\circ>0\), this frozen clock is strictly increasing.
We denote its inverse by
\begin{equation}\label{eq:source-adapted-clock-inverse}
 \tau_{\mathcal U}:J_{\mathcal U}\longrightarrow[s,s+\delta],
 \qquad
 \vartheta_{\mathcal U}
 =\frac{t_s^\circ(\tau_{\mathcal U}(\vartheta_{\mathcal U}))
        -t_s^\circ(s)}{r_{\mathcal U}^2}.
\end{equation}
By \eqref{eq:source-adapted-atlas},
\(|J_{\mathcal U}|\leq C\delta\).  Let
\(\mathscr S_{\mathcal U}\) denote the coordinate representative in
the corresponding scale-one harmonic chart.  For a normalized-time
forcing vector field \(F\), set
\begin{equation}\label{eq:source-adapted-forcing-scaling}
 \mathscr S_{\mathcal U}^{(2)}F
 :=\frac{r_{\mathcal U}^2}{\lambda_s^\circ(\tau)}
      \mathscr S_{\mathcal U}F .
\end{equation}
This is the exact forcing scaling obtained after changing from
\(\tau\) to \(\vartheta_{\mathcal U}\).

For an integer \(q\geq2\), put
\[
 \|u\|_{C_{{\rm Sch},\infty}^{q,\alpha}(\Omega\times J)}
 :=
 \|u\|_{L^\infty(J;C^{q,\alpha}(\Omega))}
 +\|\partial_\vartheta u\|_
       {L^\infty(J;C^{q-2,\alpha}(\Omega))}.
\]
Here \(\partial_\vartheta u\) is the Bochner weak derivative; thus the
displayed space is
\[
 L^\infty(J;C^{q,\alpha})
 \cap W^{1,\infty}(J;C^{q-2,\alpha}).
\]
These are the big spatial H\"older spaces; \(h^{r,\alpha}\) denotes the
corresponding little-H\"older completion, with a uniform vanishing
high-frequency modulus over the normalized atlas.  This definition
does \emph{not} impose spatial tightness.  For later comparison, call a
bounded family \(\mathcal F\) \emph{annularly tight} if it also satisfies
\begin{equation}\label{eq:weighted-little-annular-c0}
 \lim_{J\to\infty}\ 
 \sup_{X\in\mathcal F}
 \sup_{\mathcal U\subset\{\bar f\geq2^J\}}
  L_{\mathcal U}^{-N}
  \|\mathscr S_{\mathcal U}X\|_{C^{q,\alpha}}=0,
\end{equation}
with the tensor normalization of
\eqref{eq:scaled-tensor-holder} inserted when appropriate.
The optional condition \eqref{eq:weighted-little-annular-c0} is useful
for direct compactness, but it is not part of the entrance class in
Theorem~\ref{thm:intro-sharp-scattering}.  The global passage there is
instead made by a direct
two-state Cauchy estimate, which excludes annular escape without
restricting the weighted Banach domain.  Prepared initial traces lie in
the high-frequency little space, and the common high-frequency modulus
imposed in the fixed-point set below recovers that solution subspace.
This distinction is necessary:
arbitrary strongly measurable \(L^\infty(J;h^{q-2,\alpha})\) forcing
need not have a uniform high-frequency tail in time, and a bounded
family on the dyadic atlas need not satisfy
\eqref{eq:weighted-little-annular-c0}.

We shall also use the following elementary sequential property.  If
\(u_n\to u\) in any one of the global little-H\"older
atlas-supremum spaces used here, then
\(\{u,u_1,u_2,\ldots\}\) has one common vanishing top-order spatial
modulus.  Indeed, if \(\omega_v(\varrho)\) denotes the corresponding
truncated top-order seminorm, then
\[
 \omega_{u_n}(\varrho)
 \leq \omega_u(\varrho)
       +C\|u_n-u\|_{h^{q,\alpha}} .
\]
First make the last term small for all sufficiently large \(n\), and
then decrease \(\varrho\) for the remaining finite set.  The same
argument applies to the finite product coefficient, forcing, and trace
spaces below.  This is a property of norm-convergent sequences, not of
arbitrary bounded subsets of a little-H\"older space.

The displayed norm is the
initial-time-inclusive Bochner--Schauder norm.  We do not put a time
H\"older seminorm on either the forcing or the coefficients: after the
annular change of time no common such seminorm exists.  Spatial
H\"older regularity and a uniform modulus of total coefficient
oscillation are sufficient for the Cauchy estimate below.  Define
\begin{align}
 \|X\|_{\mathbb E_{\rm sc}^{q,\alpha}(I)}
 &:=
 \sup_{\mathcal U\in\mathfrak A_s}
   \|\mathscr S_{\mathcal U}X\|_
   {C_{{\rm Sch},\infty}^{q,\alpha}
    (\mathcal U\times J_{\mathcal U})},
 \label{eq:weighted-parabolic-spaces}\\
 \|F\|_{\mathbb F_{\rm sc}^{q-2,\alpha}(I)}
 &:=
 \sup_{\mathcal U\in\mathfrak A_s}
   \operatorname*{ess\,sup}_{\vartheta_{\mathcal U}\in J_{\mathcal U}}
   \|\mathscr S_{\mathcal U}^{(2)}
          F(\vartheta_{\mathcal U})\|_
    {C^{q-2,\alpha}(\mathcal U)} .
\end{align}
\(\mathbb F_{\rm sc}^{q-2,\alpha}(I)\) is the full strongly measurable
Bochner \(L^\infty\) forcing space defined by this norm; no common
modulus in the displayed source-adapted times is imposed.  This distinction
is essential: a normalized-time scalar \(a(\tau)\) is represented by
\(a(\tau_{\mathcal U}(\vartheta_{\mathcal U}))\), whose displayed-time
moduli need not be uniform as the physical scale varies.  Classical
time regularity will be recovered
scale-by-scale for the actual geometric solution after the mild
construction.

Replacing each unweighted chart term by
\(L_{\mathcal U}^{-N}\) times that term defines
\(\mathbb E_{{\rm sc},N}^{q,\alpha}\) and
\(\mathbb F_{{\rm sc},N}^{q-2,\alpha}\).  For symmetric two-tensors we
also insert the tensor normalization
\(u\mapsto L_{\mathcal U}^{-1}u\), exactly as
in \eqref{eq:scaled-tensor-holder}.  These conventions give
\[
 \sup_{\tau\in I}\|X(\tau)\|_{\mathfrak X_{\rm sc}^{q,\alpha}}
 \leq C\|X\|_{\mathbb E_{\rm sc}^{q,\alpha}(I)} .
\]

We now fix the compact affine factor and all trace conventions in the
coupled solution class.  Denote by \(G_s^{\circ}\) the compact-metric
component of the center state \(\mathbf z_s^{\circ}\), and freeze
\[
 g_{s,0}^{\circ}:=
   (\lambda^{\circ}(s))^{-1}G_s^{\circ}
\]
and a unit harmonic atlas \(\mathfrak A_{\mathcal X,s}\) for
\(g_{s,0}^{\circ}\), with doubled buffers, uniformly bounded overlap,
and the uniform transition bounds supplied by the normalized carrier
package.  All compact tensor norms below are atlas-supremum norms in
this one frozen atlas, with covariant two-tensor components normalized
by \(g_{s,0}^{\circ}\).  Thus put
\[
 \mathcal H_{\mathcal X,s}^{j,\alpha}
 :=h^{j,\alpha}_{\mathfrak A_{\mathcal X,s}}
       (\mathcal X;S^2T^*\mathcal X).
\]
On the frozen local construction, let
\[
 \widehat I_{\rm ph}=[0,\widehat\delta_{\rm ph}]
\]
be the common dimensionless compact-time interval supplied by the
scale-normalized compact quasilinear theorem.  At a restart \(s\), this
corresponds to the physical interval
\([t(s),t(s)+\lambda(s)\widehat\delta_{\rm ph}]\).  The
normalized-time interval will be chosen so that the clock of the actual
fixed point takes values in \(\widehat I_{\rm ph}\).  Trial clocks in the
open fixed-point neighborhood are evaluated against the auxiliary
two-sided coefficient extension constructed in Step~1 below; no
backward Ricci flow is asserted.  For \(r\geq2\), define
\[
 \begin{aligned}
 \mathbb E_{\mathcal X,s}^{r,\alpha}(\widehat I_{\rm ph})
 &:=h^{(r+\alpha)/2,\,r+\alpha}
   _{\mathfrak A_{\mathcal X,s}}
   (\widehat I_{\rm ph}\times\mathcal X;S^2T^*\mathcal X),
 \\
 \mathbb F_{\mathcal X,s}^{r-2,\alpha}(\widehat I_{\rm ph})
 &:=h^{(r-2+\alpha)/2,\,r-2+\alpha}
   _{\mathfrak A_{\mathcal X,s}}
   (\widehat I_{\rm ph}\times\mathcal X;S^2T^*\mathcal X).
 \end{aligned}
\]
where the first exponent is temporal, the second is spatial, and both
spaces are the closures of smooth tensor increments in the corresponding
anisotropic parabolic atlas-supremum norms, including the initial face.
The subscript \(0\) denotes the closed zero-initial-trace subspace.  These norms
are, by definition, the norms denoted below by
\(\|\cdot\|_{\mathbb E_{\mathcal X,s}^{r,\alpha}}\) and
\(\|\cdot\|_{\mathbb F_{\mathcal X,s}^{r-2,\alpha}}\).  The endpoint
trace maps from \(\mathbb E_{\mathcal X,s}^{r,\alpha}\) are bounded into
\(\mathcal H_{\mathcal X,s}^{r,\alpha}\).  Since \(\mathcal X\) is
closed, there are no lateral-boundary compatibility conditions; for
solutions, the equation determines the time jets at the initial face.

For a nearby initial state set
\[
 \lambda_s:=\lambda(s),\qquad
 g_{s,0}:=\lambda_s^{-1}G(t(s)).
\]
The compact component is always measured affinely from \(g_{s,0}\),
never by putting the absolute tensor \(\lambda_s^{-1}G\) in a fixed
host-atlas norm.  Indeed, in the common frozen prepared chart,
\[
 g_{s,0}-g_{s,0}^{\circ}
 =\lambda_s^{-1}
     \bigl(G(t(s))-G_s^{\circ}\bigr)
  +\bigl(\lambda_s^{-1}-(\lambda^{\circ}(s))^{-1}\bigr)
     G_s^{\circ}.
\]
Hence the prepared metric displacement together with
\(|\log(\lambda_s/\lambda^{\circ}(s))|\) controls this difference in
\(\mathfrak A_{\mathcal X,s}\).  After shrinking the common-margin
ball, every compact coefficient considered below satisfies, for fixed
constants \(0<c<C<\infty\),
\begin{equation}\label{eq:compact-affine-ellipticity}
 c\,g_{s,0}^{\circ}
 \leq g_{s,0}+H
 \leq C\,g_{s,0}^{\circ}.
\end{equation}
The constants and all atlas transition bounds are independent of the
collapsing scale.  On the exterior, constant rescaling only decreases
curvature and increases the available harmonic radius; on the implanted
core and transition region this is precisely part of the normalized
carrier package.

For the noncompact factor, the Bochner norm alone gives a canonical
endpoint trace two spatial orders lower.  To record the prescribed
high-order left trace without asserting top-order time continuity, set
\[
 \begin{split}
 \mathbb E_{{\rm sc},{\rm tr}}^{q,\alpha}(I)
 :=\bigl\{(X,X_s):\;&X\in\mathbb E_{\rm sc}^{q,\alpha}(I),\quad
 X_s\in\mathfrak X_{\rm sc}^{q,\alpha},\\
 &\operatorname{tr}_s^{q-2}X=X_s
   \text{ in }\mathfrak X_{\rm sc}^{q-2,\alpha}\bigr\},
 \end{split}
\]
with norm
\[
 \|(X,X_s)\|_{\mathbb E_{{\rm sc},{\rm tr}}^{q,\alpha}(I)}
 :=\|X\|_{\mathbb E_{\rm sc}^{q,\alpha}(I)}
   +\|X_s\|_{\mathfrak X_{\rm sc}^{q,\alpha}} .
\]
Here \(\operatorname{tr}_s^{q-2}\) is the strong trace supplied by the
\(W^{1,\infty}(C^{q-2,\alpha})\) part of
\(\mathbb E_{\rm sc}^{q,\alpha}\).  The same argument gives a trace at
every other endpoint in \(\mathfrak X_{\rm sc}^{q-2,\alpha}\), but no
top-order endpoint continuity is claimed.

The solution space \(\mathbb S^{r,\alpha}(I)\) consists of tuples
\[
 (G,\widehat t_s,\ell,X_\Theta,X_\Phi),\qquad
 \ell=\log(\lambda/\lambda(s)),\qquad
 \widehat t_s=\frac{t-t(s)}{\lambda(s)},
\]
recorded in the chart frozen at \(s\), with
\(t=t(s)+\lambda(s)\widehat t_s\) reconstructed, and with the following
properties.
Choose as part of that frozen chart one fixed smooth closed reference
metric \(\widehat G_s\), uniformly equivalent to
\(g_{s,0}^{\circ}\), whose coefficients and inverse have the required
uniform buffered order-\((k+3,\alpha)\) bounds in
\(\mathfrak A_{\mathcal X,s}\).  It is obtained by fixed-scale smoothing
of the center carrier package.  Put
\(B_s(g):=B_{\widehat G_s}(g)\).  For a physical compact-time argument
\(u\), set
\(\widehat\sigma=(u-t(s))/\lambda_s\), and let
\[
 \widetilde G_s(\widehat\sigma)
 :=g_{s,0}+H_s(\widehat\sigma).
\]
The zero-trace affine increment \(H_s\) and the reconstructing
diffeomorphism are the unique solutions
\[
 \begin{aligned}
  \partial_{\widehat\sigma}H_s
   &=-2\Ric_{g_{s,0}+H_s}
     +\Lie_{B_s(g_{s,0}+H_s)}(g_{s,0}+H_s),
   &H_s(0)&=0,\\
  \partial_{\widehat\sigma}\chi_s
   &=-B_s(g_{s,0}+H_s)\circ\chi_s,
  &\chi_s(0)&=\operatorname{Id}.
 \end{aligned}
\]
Consequently,
\[
 \frac d{d\widehat\sigma}
    \bigl(\chi_s^*(g_{s,0}+H_s)\bigr)
 =\chi_s^*\!\left(
   \partial_{\widehat\sigma}H_s
   +\Lie_{-B_s(g_{s,0}+H_s)}(g_{s,0}+H_s)\right)
 =-2\Ric_{\chi_s^*(g_{s,0}+H_s)},
\]
and constant-scale invariance of the Ricci tensor gives the exact
physical reconstruction
\begin{equation}\label{eq:compact-affine-reconstruction}
 G(t(s)+\lambda_s\widehat\sigma)
 =\lambda_s\chi_s(\widehat\sigma)^*
       (g_{s,0}+H_s(\widehat\sigma)).
\end{equation}
Compact
quasilinear uniqueness and ODE uniqueness make the reconstruction
unique.  On a common short interval it depends \(C^1\) on buffered
initial data from order \(r+2\) to reconstructed output order \(r\);
no same-order differentiability of the pullback reconstruction is
asserted.  The differential of the affine center is
\[
 \delta g_{s,0}
 =\lambda_s^{-1}\delta G(t(s))
   -g_{s,0}\,\delta\log\lambda_s,
\]
which is uniformly controlled by the prepared tangent norm in the
frozen atlas.  Thus the preceding parameter statement is uniform on
the common-margin ball.  In the compact factor we require
\[
 H_s\in
 \mathbb E_{\mathcal X,s,0}^{r,\alpha}(\widehat I_{\rm ph}),
 \qquad \widehat t_s,\ell\in C^{1,\alpha/2}(I),
\]
\[
 X_\Theta\in C^{1,\alpha/2}
   (I;\mathfrak X_{\rm sc}^{r+1,\alpha}),
 \qquad
 (X_\Phi,X_{\Phi,s})\in
   \mathbb E_{{\rm sc},{\rm tr}}^{r+1,\alpha}(I).
\]
The dependent graph tensor determined by the metric and map variables
belongs to
\[
 \mathbb H_{{\rm sc},N}^{r,\alpha}(I)
 :=C^{\alpha/2}\bigl(I;
   \mathfrak T_{{\rm sc},N}^{r,\alpha}\bigr)
 \cap
 C^{1,\alpha/2}\bigl(I;
   \mathfrak T_{{\rm sc},N}^{r-2,\alpha}\bigr),
\]
equipped with the sum of the two displayed norms.  We deliberately use
normalized time, without annular time rescaling, for this last norm:
the drift in the normalized \(h\)-equation is a scale-one transport
term on the AC end.  In the one frozen chart, set
\[
 \begin{aligned}
 \|\mathbf z\|_{\mathbb S^{r,\alpha}(I)}:={}&
 \|H_s\|_{
    \mathbb E_{\mathcal X,s}^{r,\alpha}(\widehat I_{\rm ph})}
 +\|\widehat t_s\|_{C^{1,\alpha/2}(I)}
 +\|\ell\|_{C^{1,\alpha/2}(I)}\\
 &+\|X_\Theta\|_{
    C^{1,\alpha/2}(I;\mathfrak X_{\rm sc}^{r+1,\alpha})}
 +\|(X_\Phi,X_{\Phi,s})\|_{
    \mathbb E_{{\rm sc},{\rm tr}}^{r+1,\alpha}(I)}\\
 &+\|h\|_{\mathbb H_{{\rm sc},N}^{r,\alpha}(I)}
 +\sup_{\tau\in I}
   \|\mathbf z(\tau)-\mathbf z_s^\circ\|_
      {\mathscr Y_{\rm prep}^{r,\alpha}} .
 \end{aligned}
\]
The compact metric factor in this formula is already the affine
zero-trace increment; the remaining affine scalar and map factors are
understood after subtracting the frozen chart center.  The subscript \(0\) denotes
the closed subspace whose coordinate differences have zero left trace;
in particular, \(X_{\Phi,s}=0\) in the centered local construction.
The compact and prepared-coordinate factors use their little-H\"older
completions.  The linear \(X_\Phi\)-factor is the big Bochner space just
defined, augmented by its prescribed high-order left trace.  The
nonlinear construction restricts it to a common-spatial-tail subset
and then takes the lower-contraction-metric completion defined in
Step~2 below.
The common-modulus requirement is not part of the ambient affine class
\(\mathbb S^{r,\alpha}(I)\); it is an additional restriction imposed
only at the fixed-point indices.  Step~2 introduces the notation for
that subset and its completion and proves that the completion retains
the full \(\mathbb S^{k+2,\alpha}(I)\) regularity and traces.

Membership in \(\mathbb S^{r,\alpha}(I)\) fixes the ambient regularity
and traces.  The term \emph{mild solution} additionally means the
following integral identities.  With the fixed reference metric above,
\[
 \mathcal R_{{\rm DT},s}(g)
 :=-2\Ric_g+\Lie_{B_{\widehat G_s}(g)}g,
\]
then the compact affine increment satisfies
\[
 H_s(\widehat\sigma)
 =\int_0^{\widehat\sigma}
    \mathcal R_{{\rm DT},s}(g_{s,0}+H_s(\zeta))\,d\zeta
 \quad\text{in }\mathcal H_{\mathcal X,s}^{r-2,\alpha},
 \qquad0\leq\widehat\sigma\leq\widehat\delta_{\rm ph}.
\]
For a noncompact component of spatial order \(q\geq2\) satisfying
\(\mathscr L_\tau X=F\), put
\[
 u_{\mathcal U}(\vartheta)
 :=\mathscr S_{\mathcal U}X
       (\tau_{\mathcal U}(\vartheta)),
 \qquad
 f_{\mathcal U}(\vartheta)
 :=\mathscr S_{\mathcal U}^{(2)}F
       (\tau_{\mathcal U}(\vartheta)).
\]
Write its fully rescaled local equation as
\[
 \partial_\vartheta u_{\mathcal U}
 -\widehat A_{\mathcal U}^{ij}D_iD_j u_{\mathcal U}
 -\widehat B_{\mathcal U}^{i}D_i u_{\mathcal U}
 -\widehat C_{\mathcal U}u_{\mathcal U}
 =f_{\mathcal U}.
\]
Here the hatted coefficients are precisely the local coefficients after
multiplication of \(\mathscr L_\tau\) by
\(r_{\mathcal U}^2/\lambda_s^\circ(\tau)\), including the fixed
coordinate-connection lower-order terms.
Then mildness means the Bochner identity
\[
 \begin{aligned}
 u_{\mathcal U}(\vartheta)
 ={}&u_{\mathcal U}(0)\\
 &+\int_0^\vartheta
 \bigl(\widehat A_{\mathcal U}^{ij}D_iD_j u_{\mathcal U}
       +\widehat B_{\mathcal U}^{i}D_i u_{\mathcal U}
       +\widehat C_{\mathcal U}u_{\mathcal U}
       +f_{\mathcal U}\bigr)(\zeta)\,d\zeta
 \end{aligned}
\]
in \(C^{q-2,\alpha}\) on every scale-one source chart.  These
identities agree on overlaps.  They require only strong Bochner
measurability of the integrands and impose no common time-H\"older
modulus in the source-adapted clocks.

The contraction itself is taken in the two-order-lower metric
\begin{equation}\label{eq:coupled-contraction-metric}
 d_{\mathbb D^{r,\alpha}(I)}(\mathbf w_1,\mathbf w_2)
 :=
 \sup_{\tau\in I}
 \|\mathbf z(\mathbf w_1)(\tau)
      -\mathbf z(\mathbf w_2)(\tau)\|
      _{\mathscr Y_{\rm prep}^{r,\alpha}} .
\end{equation}
At the output index \(r=k\) we use the mixed refinement
\begin{equation}\label{eq:coupled-mixed-contraction-metric}
 \begin{split}
 d_{\mathbb D_{\rm mix}^{k,\alpha}(I)}
   (\mathbf w_1,\mathbf w_2)
 :={}&d_{\mathbb D^{k,\alpha}(I)}(\mathbf w_1,\mathbf w_2)\\
 &+\sup_{\tau\in I}
 \|X_{\Theta,1}(\tau)-X_{\Theta,2}(\tau)\|
       _{\mathfrak X_{\rm sc}^{k+2,\alpha}} .
 \end{split}
\end{equation}
The single extra \(\Theta\)-derivative controls the \(C^{k,\alpha}\)
difference of the target Christoffel symbols
\(\Gamma(\lambda\Theta^*\bar g)\).  Since \(\Theta\) obeys an ODE, its
output in this additional norm acquires an explicit factor
\(|I|\), so the refinement does not obstruct contraction.
The full parabolic norm is used for the invariant high-norm bound, not
as the contraction metric.  This distinction is necessary: a
zero-trace lower-order forcing is small on a short interval after two
spatial orders are lowered, but not in the top parabolic norm.

\begin{lemma}[Uniform weighted parabolic estimate]
\label{lem:weighted-prepared-Schauder}
Let \(q\geq2\).  On a common-margin prepared ball, consider a linear
system on a vector bundle \(E\to M\),
\begin{equation}\label{eq:prepared-linear-parabolic}
 \mathscr L_\tau X
 :=\partial_\tau X-
   A^{ij}(\tau,x)\bar\nabla_i\bar\nabla_jX
   -B^i(\tau,x)\bar\nabla_iX-C(\tau,x)X=F .
\end{equation}
Assume that the principal coefficients act scalarly on the bundle
fibers:
\[
 A^{ij}=a^{ij}\operatorname{Id}_{E},\qquad
 a^{ij}=a^{ji}.
\]
The uniform ellipticity hypothesis below is imposed on the rescaled
coefficient
\[
 \frac{r_{\mathcal U}^{\,2}}{\lambda_s^\circ}\,a^{ij},
\]
not on the unrescaled coefficient \(a^{ij}\).  On a finite product of
bundles, assume that the second-order part is block diagonal with
scalar-principal diagonal blocks; triangular coupling is confined to
terms of order at most one.
Assume that, after passage to the charts and scalings used in
\eqref{eq:weighted-parabolic-spaces}, the operator
\[
 \frac{r_{\mathcal U}^2}{\lambda_s^\circ(\tau)}
 \mathscr L_\tau
\]
in the variable \(\vartheta_{\mathcal U}\) has ellipticity constants in
\([\Lambda^{-1},\Lambda]\), and spatial coefficient
\(C^{q-2,\alpha}\) norms bounded by \(\Lambda\), independently of
\(\mathcal U\).  Assume also that the oscillation of these rescaled
coefficients over \(I\), in the same spatial norms, is bounded by a
common modulus \(\omega(\delta)\) with
\(\omega(\delta)\to0\) as \(\delta\downarrow0\).  Suppose
 \(X_s\in\mathfrak X_{\rm sc}^{q,\alpha}\) and
 \(F\in\mathbb F_{\rm sc}^{q-2,\alpha}(I)\), with no time-continuity
 hypothesis on \(F\) beyond strong Bochner measurability.  Then there is a
number \(\delta_0>0\), depending only on \(q,\alpha,\Lambda\), the
common modulus \(\omega\), and the fixed prepared margins, such that
 the initial-value problem
\[
 \mathscr L_\tau X=F,\qquad X(s)=X_s,
\]
has
 a unique mild solution in \(\mathbb E_{\rm sc}^{q,\alpha}(I)\), and
\begin{equation}\label{eq:weighted-prepared-Schauder}
 \|X\|_{\mathbb E_{\rm sc}^{q,\alpha}(I)}
 \leq C\left(
       \|X_s\|_{\mathfrak X_{\rm sc}^{q,\alpha}}
       +\|F\|_{\mathbb F_{\rm sc}^{q-2,\alpha}(I)}
       \right),
 \qquad0<\delta\leq\delta_0.
\end{equation}
If the localized coefficient fields have a common vanishing top-order
high-frequency H\"older tail, and the localized initial traces and
forcing have the corresponding common high-frequency tails at orders
\(q\) and \(q-2\), respectively, then the solution has the same
property and lies in the global
little-H\"older subspace used in
Theorem~\ref{thm:intro-sharp-scattering}.  If the family is
additionally such that the initial traces and forcing are annularly
tight in the sense of \eqref{eq:weighted-little-annular-c0}, at their
respective orders and with the source-adapted forcing normalization,
that optional second modulus is also preserved.  No annular vanishing
is required of the nonzero principal coefficients.  The assertions are
uniform for the supplied moduli.

We shall use the following quantitative version of the spatial-tail
assertion.  Fix concentric members
\(\mathcal U^0\Subset\mathcal U^1\) of one normalized buffered chart
and use one of the fixed uniformly bounded extension operators from
\(\mathcal U^1\) to a Euclidean ball.  For a chart representative \(u\)
and an integer \(m\geq0\), put
\[
 \omega_{m,\alpha}^{\mathcal U}(u;\varrho)
 :=
 \operatorname*{ess\,sup}_{\vartheta}
 \sup_{\substack{x,y\in\mathcal U^0\\
                  0<|x-y|\leq\varrho}}
 \frac{|D^m u(\vartheta,x)-D^m u(\vartheta,y)|}
      {|x-y|^\alpha}.
\]
For an initial trace the essential supremum in \(\vartheta\) is
omitted.  The quantity
\(\operatorname{Tail}_{m,\alpha}(u;\varrho)\) is the supremum of these
seminorms over the normalized atlas, after inserting the fixed tensor,
bundle, and polynomial-weight normalizations.  Different choices of
the uniformly buffered extensions give equivalent quantities, with
one fixed dilation of \(\varrho\).

The coefficient-tail hypothesis below is understood in the strong
Bochner essential supremum over the displayed time variable.  It
includes the order-\((q-2)\) tails of \(A,B,C\), the tails at the
highest order in which the fixed cutoff coefficients occur, and the
tails of the induced transition operators on the chosen bundle
representatives.  Thus it concerns the transition matrices acting on
components and their derivatives, not a fictitious finite-band action
of a nonlinear coordinate map.  Let \(\eta_{\mathscr L}\) be any
bounded nondecreasing function, vanishing at zero, which dominates all
these tails uniformly over the complete normalized atlas.  The fixed
atlas package supplies the corresponding uniform big-H\"older bounds,
buffer widths, Lipschitz constants, and overlap multiplicity.

There are constants \(C_{\rm tail},C_{\rm dil}\), depending only on
the numerical data in this lemma and not on the rate at which
\(\eta_{\mathscr L}\) vanishes, such that
\begin{equation}\label{eq:quantitative-tail-Schauder}
\begin{split}
 \operatorname{Tail}_{q,\alpha}(X;\varrho)
 \leq C_{\rm tail}\bigl\{&
 \operatorname{Tail}_{q,\alpha}
      (X_s;C_{\rm dil}\varrho)
 +\operatorname{Tail}_{q-2,\alpha}
      (F;C_{\rm dil}\varrho)\\
 &+\bigl(\eta_{\mathscr L}(C_{\rm dil}\varrho)
          +\varrho^{1-\alpha}\bigr)
   \bigl(\|X_s\|_{\mathfrak X_{\rm sc}^{q,\alpha}}
         +\|F\|_{\mathbb F_{\rm sc}^{q-2,\alpha}(I)}\bigr)
 \bigr\}
\end{split}
\end{equation}
for \(0<\varrho\leq C_{\rm dil}^{-1}\).  The same statement holds with
 the fixed weights, tensor normalizations, trace extensions, finite
 product systems, and localized transition operators above.  This is an
 estimate for each fixed operator; a family with one coefficient-tail
 modulus has the same right-hand side uniformly.

The same estimate holds with the polynomial weight \(N\), with tensor
scaling, and for finite product systems having the block-diagonal
scalar-principal and triangular lower-order structure just specified.
The initial datum is then measured in the correspondingly weighted
trace norm.  The constant depends only on
\(q,\alpha,\Lambda\), the fixed background bounded-geometry and atlas
data, and the common prepared margins; for the weighted version it may
also depend on the fixed exponent \(N\).  It is independent of the
outer annulus and of the exhaustion of \(M\).

If \(X(s)=0\), then two-order lowering gives the uniform small-time
estimate
\begin{equation}\label{eq:weighted-zero-trace}
 \sup_{\tau\in I}
 \|X(\tau)\|_{\mathfrak X_{\rm sc}^{q-2,\alpha}}
 \leq C\delta^{\alpha/2}
       \|X\|_{\mathbb E_{\rm sc}^{q,\alpha}(I)} .
\end{equation}
Interpolation with the top norm gives the one-order lowering
\begin{equation}\label{eq:weighted-zero-trace-one-order}
 \sup_{\tau\in I}
 \|X(\tau)\|_{\mathfrak X_{\rm sc}^{q-1,\alpha}}
 \leq C\delta^{\alpha/4}
       \|X\|_{\mathbb E_{\rm sc}^{q,\alpha}(I)} .
\end{equation}
For two solutions, both zero-trace estimates hold for their difference
provided their initial traces agree:
\[
 X_1(s)=X_2(s).
\]
For a first variation they hold provided the differentiated initial
trace vanishes, \(D_pX_s[\dot p]=0\).  With a nonzero difference or
variation of the initial trace, one first subtracts its homogeneous
trace evolution and applies the zero-trace estimates to the remainder;
the trace evolution itself is controlled by
\eqref{eq:weighted-prepared-Schauder}.

The solution operator is \(C^1\) with the same bounds for
parameter-dependent families under the following topology.  Let the
parameter range be an open subset of a Banach space, or a Banach
manifold in one fixed local chart.  Require the rescaled coefficient
map
\[
 p\longmapsto(A_p,B_p,C_p)
\]
to be \(C^1\) into the strong-Bochner
\(L^\infty(I;C_x^{q-2,\alpha})\) product space defined by the same
global atlas-supremum coefficient norm, and require bounded parameter
neighborhoods to carry the same coefficient time-oscillation modulus
\(\omega\).  Suppose also that the forcing is \(C^1\) into
\(\mathbb F_{\rm sc}^{q-2,\alpha}(I)\) and that the initial datum is
\(C^1\) into \(\mathfrak X_{\rm sc}^{q,\alpha}\).  Then the solution
operator is \(C^1\) with the same bounds.  For a parameter \(p\) and
tangent \(\dot p\), its derivative is the unique mild solution of the
differentiated problem
\[
 \mathscr L_\tau\dot X
 =
 \dot F+\dot A^{ij}\bar\nabla_i\bar\nabla_jX
       +\dot B^i\bar\nabla_iX+\dot C X,
 \qquad
 \dot X(s)=D_pX_s[\dot p].
\]
In particular, a parameter-independent initial trace gives
\(\dot X(s)=0\).

The solution operator also inherits mixed tame remainders.  Precisely,
suppose the coefficient, forcing, and trace maps on a bounded high-order
ball have one-high--one-low Taylor remainders, and let \(p_n\to p\) in
their high topology.  If their increments and the corresponding
solution increments are \(O(\epsilon_n)\) in the two-order-lower
topology, then, after subtracting the differentiated problem, the
residual forcing \(\mathfrak r_n\) satisfies
\begin{equation}\label{eq:Schauder-mixed-tame-residual}
\begin{split}
 \|\mathfrak r_n\|_{\mathbb F_{\rm sc}^{q-2,\alpha}}
 \leq C_K\bigl(&
   \|\delta p_n\|_{\rm high}\|\delta X_n\|_{\rm low}
  +\|\delta p_n\|_{\rm low}\|\delta X_n\|_{\rm high}\\
 &+\|\delta p_n\|_{\rm high}\|\delta p_n\|_{\rm low}\bigr)
 =o(\epsilon_n).
\end{split}
\end{equation}
The residual initial trace obeys the analogous estimate.  Applying
\eqref{eq:weighted-prepared-Schauder} to the residual equation gives
the same \(o(\epsilon_n)\) conclusion in the corresponding
two-order-lower solution norm.  This clause uses convergence in the
high topology but requires a rate only in the lower topology.

The numerical Schauder and differentiated estimates are independent
of any spatial high-frequency modulus.  If the coefficient, forcing,
and trace maps take values \(C^1\)-smoothly in the corresponding
closed uniform-spatial-tail little-H\"older subspaces of the displayed
Bochner spaces, then the solution operator restricts to a \(C^1\) map
between those open little-H\"older parameter and solution spaces.  Here
the forcing hypothesis is stronger than mere membership in
\(L^\infty(I;h^{q-2,\alpha})\): one spatial tail is required uniformly
over the time variable.  No single modulus is imposed on the whole
open parameter neighborhood.  A common modulus is needed only for a
simultaneous tail assertion about an entire family; each individual
value in this closed uniform-tail subspace supplies its own modulus,
and every norm-convergent sequence in that subspace supplies a common
one by the sequential observation above.
If, in each fixed core or annular chart, the coefficients and forcing
are continuous in the original normalized time, the mild solution is
classical there.  Higher normalized-time and spatial regularity, when
present in the data, is recovered by the usual local parabolic
bootstrap.  Neither recovery statement requires a modulus uniform in
the displayed annular time.
\end{lemma}

\begin{proof}
On the compact core use the initial-value spatial H\"older estimate
with strongly measurable \(L^\infty\)-in-time forcing in a fixed
finite harmonic atlas.  Equivalently, the frozen Euclidean heat
 operator maps
\[
 L^\infty_\vartheta C^{q-2,\alpha}_x
 \longrightarrow
 L^\infty_\vartheta C^{q,\alpha}_x
 \cap W^{1,\infty}_\vartheta C^{q-2,\alpha}_x
\]
 with its initial trace included.  In a bundle trivialization the
 frozen principal part is
\(a_0^{ij}\partial_i\partial_j\otimes\operatorname{Id}_{E}\);
connection terms and every triangular coupling are of order at most
one.  Hence the scalar Fourier heat multiplier acts componentwise.
Indeed, if \(\Delta_j\) is a spatial dyadic block and \(A_0\) is this
frozen uniformly elliptic scalar-principal operator, then
 \[
 \begin{split}
 &2^{(q+\alpha)j}\left\|
  \Delta_j\int_0^\vartheta
       e^{(\vartheta-r)A_0}F(r)\,dr
 \right\|_{L^\infty_x}\\
 &\qquad\leq
 C\,2^{(q-2+\alpha)j}
 \|\Delta_jF\|_{L^\infty_\vartheta L^\infty_x}
  \int_0^\vartheta
  2^{2j}e^{-c2^{2j}(\vartheta-r)}\,dr\\
 &\qquad\leq
 C\,2^{(q-2+\alpha)j}
 \|\Delta_jF\|_{L^\infty_\vartheta L^\infty_x}.
 \end{split}
 \]
 The low-frequency block is bounded directly, and
  \(\partial_\vartheta X=A_0X+F\) gives the
  \(W^{1,\infty}_\vartheta C^{q-2,\alpha}_x\) term.  Taking the dyadic
  supremum proves the spatial-H\"older Bochner--Schauder estimate; the
  full \(L^\infty\)-in-time norm alone does not assert a uniform
  vanishing high-frequency tail.  If the data, coefficients, and
 forcing are supplied with one common high-frequency modulus, the
 identical block estimate, followed by a finite low-frequency
 truncation, preserves that modulus and proves the global
 little-H\"older clause.  The optional annular-tightness assertion is
  proved below directly on the global parametrix; no off-diagonal
  estimate from a different realization is being invoked.  No
 time modulus of the right-hand side is used.  On each
\(\mathcal U\in\mathfrak A_s\), pass to the scale-one
\(r_{\mathcal U}^{-2}\acute G_\circ(s)\)-harmonic chart, use
\(\vartheta_{\mathcal U}\), and multiply
\eqref{eq:prepared-linear-parabolic} by
\(r_{\mathcal U}^2/\lambda_s^\circ\).  By hypothesis, one
obtains a uniformly parabolic system on a Euclidean ball with
ellipticity and coefficient bounds independent of \(\mathcal U\).  The
 localized Bochner--Schauder estimate on a concentric smaller ball,
 with the initial trace included, therefore has one common constant.

We next prove the quantitative tail assertion.  No approximation in
time is made.  Apply fixed spatial Littlewood--Paley truncations (or
fixed spatial mollifiers in the extended charts) at each time to the
coefficients and forcing, and the same spatial truncations to the
initial trace.  These operators preserve strong Bochner measurability.
Under the common uniform spatial little-H\"older tail assumed in the
statement, the truncated fields converge in the relevant
\(L^\infty_\vartheta C_x^{m,\alpha}\) norms, while uniform ellipticity
persists for all sufficiently large truncation levels.  We may
therefore carry out the following calculation for representatives
smooth in the spatial variables, with their original merely measurable
time dependence, and pass to the stated little-H\"older closures.  Let
\(\Delta_\ell\), \(\ell\geq-1\), be a
fixed inhomogeneous spatial Littlewood--Paley decomposition on the
extended Euclidean chart.  For \(m\geq0\), set
\[
 b_\ell^m(u)
 :=
 2^{(m+\alpha)\ell}
 \|\Delta_\ell u\|_{L^\infty_{\vartheta,x}},
\]
with the evident omission of the time supremum for a trace, and define
the two-sided H\"older frequency envelope
\[
 \mathcal W_j^m(u)
 :=
 \max\left\{
   \sup_{\ell\geq j} b_\ell^m(u),\
   \sup_{-1\leq\ell<j}
      2^{-(1-\alpha)(j-\ell)}b_\ell^m(u)
 \right\}.
\]
The second term is essential: blocks below \(j\) are not absent from a
H\"older quotient at pair-distance \(2^{-j}\); their contribution is
suppressed by the displayed geometric factor.

The difference characterization of \(C^{m,\alpha}\), applied to the
kernels of \(\Delta_\ell\), gives fixed integers \(J_0\) and constants
\(c_0,C_0\), depending only on the chosen Euclidean extensions, such
that
\begin{equation}\label{eq:tail-frequency-envelope-equivalence}
 \begin{split}
 \operatorname{Tail}_{m,\alpha}(u;2^{-j})
 &\leq
 C_0\mathcal W_{j-J_0}^m(u)
 +C_0\,2^{-(1-\alpha)j}\|u\|_{C^{m,\alpha}},\\
 \mathcal W_j^m(u)
 &\leq
 C_0\operatorname{Tail}_{m,\alpha}
      (u;c_0\,2^{-j})
 +C_0\,2^{-(1-\alpha)j}\|u\|_{C^{m,\alpha}}.
 \end{split}
\end{equation}
Indeed, a block of frequency at least \(2^j\) is estimated using the
cancellation of its kernel and a difference at its own scale.  For a
block of frequency below \(2^j\), Bernstein's inequality gives the
factor \(2^{-(1-\alpha)(j-\ell)}\).  Summing the two geometric series
proves the first inequality; applying the block kernels to differences
of \(u\) proves the second.  The inhomogeneous block contributes only
the displayed \(2^{-(1-\alpha)j}\)-term.  This proves the required
quantitative equivalence with the pair-distance tail, rather than only
the qualitative equivalence of the associated little-H\"older
closures.

Work now on one normalized chart.  After subdividing its fixed outer
buffer into a uniformly bounded number of smaller buffered balls,
choose a spatial center \(x_0\) and freeze only in space:
\[
 A_0^{ij}(\vartheta):=A^{ij}(\vartheta,x_0).
\]
The scalar-principal hypothesis implies that the evolution family
\(U_0(\vartheta,r)\) of
\[
 \partial_\vartheta
   -A_0^{ij}(\vartheta)\partial_i\partial_j
\]
acts componentwise on the bundle fibers.  Its Fourier symbol has
covariance \(\int_r^\vartheta A_0(\zeta)\,d\zeta\); uniform ellipticity
therefore gives
\[
 \|\Delta_\ell U_0(\vartheta,r)\|_{L^\infty\to L^\infty}
 \leq C e^{-c2^{2\ell}(\vartheta-r)}.
\]
Only strong measurability in \(\vartheta\) is used here.  Consequently
the homogeneous term and the Duhamel term satisfy
\[
 b_\ell^q(X)
 \leq C b_{\ell+O(1)}^q(X_s)
      +C b_{\ell+O(1)}^{q-2}(F)
      +C b_{\ell+O(1)}^{q-2}(\mathcal Q),
\]
where \(\mathcal Q\) contains the coefficient-freezing error,
connection terms, lower-order coefficients, and localization
commutators.  Here and below \(\ell+O(1)\) denotes a uniformly bounded
set of neighboring blocks and therefore produces only one fixed shift
in the envelope index.

We record all top-order possibilities in \(\mathcal Q\).  Write
\(Y=\zeta X\) for a fixed inner cutoff.  Denote the inhomogeneous
paraproduct and resonant remainder in Bony's decomposition by
\(\mathsf T_f g\) and \(\mathsf R(f,g)\), respectively.
For the second-order freezing error,
\[
 (A-A_0)D^2Y
 =
 \mathsf T_{A-A_0}D^2Y
 +\mathsf T_{D^2Y}(A-A_0)
 +\mathsf R(A-A_0,D^2Y).
\]
The low-coefficient--high-solution part of the first term has norm at
most \(\varepsilon_{\rm sp}\mathcal W_j^q(Y)\), where
\(\varepsilon_{\rm sp}\) is made uniformly small by the preceding
choice of spatial radius.  Its commutator with \(\Delta_\ell\) has
either a top coefficient block or one unused spatial derivative.  The
former is bounded by
\[
 C\eta_{\mathscr L}(C2^{-j})
   \|Y\|_{L^\infty C^{q,\alpha}},
\]
and the latter by
\[
 C2^{-(1-\alpha)j}
   \|Y\|_{L^\infty C^{q,\alpha}}.
\]
In \(\mathsf T_{D^2Y}(A-A_0)\) and in the resonant term, a high block lies on
the coefficient unless the solution factor has one unused derivative.
These give the same two alternatives.  This exhausts the high--low,
low--high, and high--high interactions of the principal coefficient.

For \(B\,DY\) and \(CY\), a top block on \(B\) or \(C\) is controlled
by \(\eta_{\mathscr L}\), while a top block on \(Y\) gains at least one
power of \(2^{-\ell}\) because these operators have order at most one.
The same statement applies to the fixed connection coefficients.
Finally,
\[
 [\mathscr L,\zeta]X
 =
 -2A^{ij}(D_i\zeta)D_jX
 -A^{ij}(D_iD_j\zeta)X
 -B^i(D_i\zeta)X
\]
is of order at most one.  A top block on a cutoff or on a coefficient
is included in \(\eta_{\mathscr L}\); otherwise there is an unused
solution derivative.  Thus no localization term has been omitted.
Taking the two-sided envelope gives, for one fixed integer \(J_1\),
\begin{equation}\label{eq:quantitative-tail-block}
 \begin{split}
 \mathcal W_j^q(Y)
 \leq {}&
 C\mathcal W_{j-J_1}^q(X_s)
 +C\mathcal W_{j-J_1}^{q-2}(F)\\
 &+C\bigl(
     \eta_{\mathscr L}(C2^{-j})
     +2^{-(1-\alpha)j}\bigr)
   \bigl(
     \|X\|_{\mathbb E_{\rm sc}^{q,\alpha}}
     +\|X_s\|_{\mathfrak X_{\rm sc}^{q,\alpha}}
     +\|F\|_{\mathbb F_{\rm sc}^{q-2,\alpha}}\bigr)
 +\frac18\mathcal W_j^q(Y).
 \end{split}
\end{equation}
The fraction \(1/8\) includes the absorbed principal freezing error;
its precise value is immaterial.

We next pass between overlapping normalized charts.  If \(\Psi\) is a
fixed coordinate change and \(P_\Psi\) is the induced bundle
transition matrix, the transition operator is
\[
 Tu=P_\Psi\,(u\circ\Psi).
\]
It is not frequency banded.  The uniform transition bounds and the
buffer imply instead the almost-diagonal estimate
\[
 2^{(m+\alpha)\ell}
 \|\Delta_\ell T\Delta_r u\|_{L^\infty}
 \leq
 C_M\,2^{-M(|\ell-r|-J_2)_+}
 2^{(m+\alpha)r}\|\Delta_r u\|_{L^\infty}
\]
for one fixed \(J_2\) and some
\(M>\max\{\alpha,1-\alpha\}\), with the top-order remainders bounded
by the prescribed transition tail.  Equivalently, differentiating
\(P_\Psi(u\circ\Psi)\) and using the pair-distance chain rule gives
the following estimate.  Let \(L_{\rm tr}\geq1\) be the common
Lipschitz bound for the transition maps and their inverses, and let
\(\operatorname{Tail}_{\rm tr}\) denote the maximum of the prescribed
top-order transition-map and transition-matrix tails:
\[
 \begin{split}
 \operatorname{Tail}_{m,\alpha}(Tu;\varrho)
 \leq C\bigl\{&
 \operatorname{Tail}_{m,\alpha}
       (u;L_{\rm tr}\varrho)\\
 &+\bigl(
    \operatorname{Tail}_{\rm tr}(\Psi,P_\Psi;
                                  L_{\rm tr}\varrho)
    +\varrho^{1-\alpha}\bigr)
   \|u\|_{C^{m,\alpha}}\bigr\}.
 \end{split}
\]
The lower Fa\`a di Bruno terms have one unused
derivative.  The same estimate, with \(\Psi=\operatorname{Id}\),
treats multiplication by a cutoff.  Hence every transition and cutoff
preserves \(\mathcal W_j^m\) with one fixed index shift and a remainder
dominated by
\(\eta_{\mathscr L}(C2^{-j})+2^{-(1-\alpha)j}\); no assertion of a
literal finite frequency shift is used.

The cover has bounded multiplicity, each cutoff is supported in one
uniformly bounded finite star, and overlapping members have comparable
radii and dyadic labels.  Taking the normalized atlas supremum
therefore changes only the constants and replaces \(1/8\) by, say,
\(1/4\).  The factor \(L_{\mathcal U}^{-N}\) and the tensor
normalization are constant on each normalized chart and change by a
bounded factor on overlaps, so the same estimate holds in the weighted
spaces.  For a finite product system, apply the argument to each
scalar-principal diagonal block.  Every triangular coupling has order
at most one and is already covered by the lower-order case above.

Absorb the final \(\frac14\mathcal W_j^q(X)\), and use the ordinary
Schauder bound \eqref{eq:weighted-prepared-Schauder}, established
independently by the parametrix construction below, to replace
\(\|X\|_{\mathbb E_{\rm sc}^{q,\alpha}}\) by the trace and forcing
norms.  The two envelope--tail inequalities, with
\(\varrho\simeq2^{-j}\), now give
\eqref{eq:quantitative-tail-Schauder}, with one fixed
\(C_{\rm dil}\).

The parametrix and Neumann series constructed below are used only to
produce the exact solution and its ordinary Schauder bound.  Once that
solution has been constructed, the preceding argument is applied
directly to the exact equation that it satisfies.  Thus the fixed
chart dilation occurs once in the a priori estimate; no tail estimate
is iterated through the powers of the parametrix error operator, and
no accumulated dilation is introduced.

 A partition of unity subordinate to the source-adapted cover has
uniformly bounded derivatives in its scale-one metrics;
the commutators with the partition are lower order and are absorbed,
after decreasing \(\delta_0\), by interpolation and the zeroth-order
term in the estimate.

Neighboring charts have comparable radii and labels, and the cover has
bounded multiplicity.  Hence taking the chart supremum proves
\eqref{eq:weighted-prepared-Schauder}.  The factor
\(L_{\mathcal U}^{-N}\) changes by at most a fixed factor on an
overlap, which proves the weighted assertion.  The tensor
normalization is likewise fixed on each chart.  This is the reason for
the source-adapted forcing scaling
\eqref{eq:source-adapted-forcing-scaling}.

We now specify the global operator and trace structure.  Set
\[
 \mathbb E_{{\rm sc},0}^{q,\alpha}(I)
 :=\bigl\{X\in\mathbb E_{\rm sc}^{q,\alpha}(I):
       \operatorname{tr}_{s}^{q-2}X=0\bigr\}.
\]
This is a closed subspace of \(\mathbb E_{\rm sc}^{q,\alpha}(I)\).
The coefficient bounds make
\[
 \mathscr L_0:\mathbb E_{{\rm sc},0}^{q,\alpha}(I)
 \longrightarrow \mathbb F_{\rm sc}^{q-2,\alpha}(I),
 \qquad \mathscr L_0X:=\mathscr L X,
\]
a bounded operator.  Notice that these are spaces of global bundle
sections; the atlas is used to define their norms, not to replace a
global section by an arbitrary family of chart representatives.

Choose a uniformly locally finite partition of unity
\(\{\zeta_{\mathcal U}\}\) subordinate to the source-adapted cover,
together with outer cutoffs equal to one on
\(\operatorname{supp}\zeta_{\mathcal U}\).  In the normalized outer
chart choose a spatial center \(x_{\mathcal U}\) and set
 \[
 A_{\mathcal U,0}^{ij}(\vartheta)
 :=\widehat A_{\mathcal U}^{ij}(\vartheta,x_{\mathcal U}).
\]
This is the same spatially frozen, strongly measurable,
scalar-principal family used above.  Its evolution family
\(U_{\mathcal U,0}(\vartheta,r)\) is defined by the Fourier multiplier
with covariance
\(\int_r^\vartheta A_{\mathcal U,0}(\zeta)\,d\zeta\), and the preceding
ellipticity estimate gives the uniform zero-trace Cauchy solver
 \[
 (\mathcal S_{\mathcal U}F)(\vartheta)
 :=\int_0^\vartheta
   U_{\mathcal U,0}(\vartheta,r)F(r)\,dr .
\]
Thus no value of a merely measurable coefficient is selected at a
time endpoint, and no time-freezing error is introduced.  Restrict a
global forcing to an outer chart, use the fixed bounded Euclidean
extension, apply \(\mathcal S_{\mathcal U}\), pull the result back as a
bundle section, and multiply by \(\zeta_{\mathcal U}\).  The locally
finite sum defines a bounded operator
\[
 \mathcal P_0:\mathbb F_{\rm sc}^{q-2,\alpha}(I)
 \longrightarrow \mathbb E_{{\rm sc},0}^{q,\alpha}(I).
\]
Thus \(\mathcal P_0F\) is an actual global section and has exactly zero
initial trace.  Bounded overlap, the transition estimates proved above,
and the fixed buffers show that the sum and its Bochner weak time
derivative belong to the displayed global spaces, and that
\(\mathscr L_0\mathcal P_0F\) is the distributional application of
\(\mathscr L\) to that section.

Put
\[
 \mathcal E_0:=I-\mathscr L_0\mathcal P_0
 \quad\hbox{on }\mathbb F_{\rm sc}^{q-2,\alpha}(I).
\]
Its terms are precisely the spatial principal-coefficient freezing
errors, the lower-order and connection terms, and the extension,
localization, and partition commutators.  After first fixing one
sufficiently small harmonic-chart radius, the undifferentiated spatial oscillation of the
principal coefficient has norm \(\varepsilon_{\rm sp}\), uniformly in
the measurable time variable, whereas every term in which a derivative
hits a coefficient and every partition commutator is lower order.
Since \(\mathcal P_0F\) has zero trace, spatial interpolation followed
by time integration gives
\[
 \|\mathcal E_0\|_{\mathcal L(\mathbb F_{\rm sc}^{q-2,\alpha}(I))}
 \leq \varepsilon_{\rm sp}
       +\varepsilon_{\rm com}(\delta),
 \qquad \varepsilon_{\rm com}(\delta)\longrightarrow0.
\]
Choose \(\varepsilon_{\rm sp}<1/4\), and then choose \(\delta_0\) so
that
\[
 \sup_{0<r\leq\delta_0}\varepsilon_{\rm com}(r)<\frac14.
\]
Then \(\|\mathcal E_0\|<1/2\), and
\[
 \mathcal Q_0:=\mathcal P_0(I-\mathcal E_0)^{-1}
 :\mathbb F_{\rm sc}^{q-2,\alpha}(I)
 \longrightarrow\mathbb E_{{\rm sc},0}^{q,\alpha}(I)
\]
is a bounded right inverse of \(\mathscr L_0\).  Because
\(\mathbb E_{{\rm sc},0}^{q,\alpha}(I)\) is complete and
\(\mathscr L_0\) is bounded, the operator-norm Neumann series converges
in the stated operator domain, not merely chartwise.

It remains to insert the prescribed trace.  The same buffered local
construction, now using the homogeneous nonautonomous spatially frozen
Cauchy evolution \(U_{\mathcal U,0}(\vartheta,0)\) of the localized
initial data, gives a bounded trace extension
\[
 \mathcal T:\mathfrak X_{\rm sc}^{q,\alpha}
 \longrightarrow\mathbb E_{\rm sc}^{q,\alpha}(I),
 \qquad
 \operatorname{tr}_{s}^{q-2}\mathcal T X_s=X_s,
 \qquad
 \mathscr L\mathcal T X_s\in
 \mathbb F_{\rm sc}^{q-2,\alpha}(I).
\]
The trace identity is exact because each local extension has the
localized trace and \(\sum_{\mathcal U}\zeta_{\mathcal U}=1\).
For \((X_s,F)\) define
\[
 X:=\mathcal T X_s
   +\mathcal Q_0\bigl(F-\mathscr L\mathcal T X_s\bigr).
\]
Then \(\operatorname{tr}_{s}^{q-2}X=X_s\), \(\mathscr L X=F\), and
the preceding bounds give \eqref{eq:weighted-prepared-Schauder}.  In
equivalent augmented-operator notation,
\[
 \mathbf L:\mathbb E_{{\rm sc},{\rm tr}}^{q,\alpha}(I)
 \longrightarrow
 \mathfrak X_{\rm sc}^{q,\alpha}\times
 \mathbb F_{\rm sc}^{q-2,\alpha}(I),
 \qquad
 \mathbf L(X,X_s)=(X_s,\mathscr L X),
\]
has the bounded right inverse just constructed.  Applying the a priori
estimate to the difference of two zero-trace solutions proves
uniqueness.

We finally verify the optional annular \(c_0\) clause.  Identify the
weighted normalized chart data with the corresponding closed subspace
of the atlas-supremum sequence space.  The trace extension \(\mathcal T\)
and every localized zero-trace Cauchy solver entering \(\mathcal P_0\)
use only one uniformly buffered finite star.  Bounded overlap and
comparability of labels therefore imply
\[
 \mathcal T:c_0(\mathfrak A_s)_{\rm tr}
       \longrightarrow c_0(\mathfrak A_s)_{\rm sol},
 \qquad
 \mathcal P_0:c_0(\mathfrak A_s)_{\rm for}
       \longrightarrow c_0(\mathfrak A_s)_{\rm sol}.
\]
Restriction, transition, multiplication by uniformly bounded
coefficients, and the partition commutators preserve these closed
subspaces.  Hence \(\mathcal E_0\) preserves
\(c_0(\mathfrak A_s)_{\rm for}\), its Neumann series restricts to that
subspace, and \(\mathcal Q_0\) preserves annular tightness.  The formula
for \(X\) above then proves the asserted annular \(c_0\) statement.
The equation then gives the corresponding order-\((q-2)\) tightness
for \(\partial_\vartheta X\).  Notice that only uniform boundedness,
not annular decay, was used for the coefficients.

If the initial trace is zero, the fundamental theorem of calculus in
each source-adapted clock gives
\[
 \sup_{\vartheta_{\mathcal U}\in J_{\mathcal U}}
 \|X(\vartheta_{\mathcal U})\|_{C^{q-2,\alpha}}
 \leq |J_{\mathcal U}|
 \sup_{\vartheta_{\mathcal U}\in J_{\mathcal U}}
 \|\partial_{\vartheta_{\mathcal U}}X\|_{C^{q-2,\alpha}}.
\]
Since \(|J_{\mathcal U}|\leq C\delta\) uniformly, this proves the weaker
but convenient estimate
\eqref{eq:weighted-zero-trace}.  Interpolation between this
\(C^{q-2,\alpha}\) estimate and the top
\(C^{q,\alpha}\) bound gives
\eqref{eq:weighted-zero-trace-one-order}; the source-adapted interval
has one common upper bound on every chart.
Subtraction proves the full difference estimate.  If the two initial
traces agree, the difference has zero trace and the preceding
fundamental-theorem-of-calculus argument proves
\eqref{eq:weighted-zero-trace}--%
\eqref{eq:weighted-zero-trace-one-order}; otherwise subtract the
homogeneous evolution of their trace difference first.

For the parameter assertion, first extend the varying initial trace
with the same bounded trace extension used in the parametrix.  Work in
one fixed Banach chart and let \(p_n=p+\eta_n\), where
\(\eta_n\to0\) is an arbitrary sequence of increments for which
\(p_n\) remains in the parameter domain.  Let \(\dot X_n\) be the
solution of the differentiated equation with direction \(\eta_n\),
and put
\[
 R_n:=X(p_n)-X(p)-\dot X_n .
\]
The ordinary difference equation and
\eqref{eq:weighted-prepared-Schauder} first give
\[
 \|X(p_n)-X(p)\|_{\mathbb E_{\rm sc}^{q,\alpha}}
 \leq C\|\eta_n\|.
\]
Thus every product of a first-order coefficient increment with the
solution increment is \(O(\|\eta_n\|^2)\).
After subtracting the three equations, the residual initial trace and
residual forcing are \(o(\|\eta_n\|)\) in their stated trace and
Bochner norms, because the coefficient, forcing, and initial-trace
maps are Fr\'echet \(C^1\) in those norms.  The
modulus-independent estimate
\eqref{eq:weighted-prepared-Schauder} therefore gives
\[
 \|R_n\|_{\mathbb E_{\rm sc}^{q,\alpha}}
 =o(\|\eta_n\|).
\]
This proves Fr\'echet differentiability for arbitrary increments, not
merely differentiability along a chosen curve.  Subtracting the
linearized equations at two base points, taking the supremum over unit
tangent directions, and applying the same estimate proves continuity
of the derivative in operator norm.  Each individual solution and
 variation lies in the little-H\"older subspace by the tail-preservation
 clause; the constants in the preceding remainder estimate do not
 depend on its individual modulus.

For the mixed clause, perform the same subtraction but estimate every
product by the H\"older tame inequality
\[
 \|uv\|_{\rm out}
 \leq C\bigl(\|u\|_{\rm high}\|v\|_{\rm low}
             +\|u\|_{\rm low}\|v\|_{\rm high}\bigr).
\]
The assumed mixed Taylor bounds handle the pure coefficient, forcing,
and trace remainders.  High convergence makes the high factors
\(o(1)\), while the lower difference estimate makes the low factors
\(O(\epsilon_n)\).  This proves
\eqref{eq:Schauder-mixed-tame-residual}, and the lower-order version of
\eqref{eq:weighted-prepared-Schauder} gives the asserted residual
solution estimate.

Finally, on a fixed chart, time-continuous coefficients and forcing
may be approximated in their own chartwise modulus.  The mild identity
and the local estimate then give the classical equation, and
differentiation proves the stated bootstrap.
\end{proof}

\begin{lemma}[Short effective-clock same-order propagation]
\label{lem:short-effective-clock-propagation}
Fix \(q\geq2\), \(0<\alpha<1\), and \(\Lambda\geq1\).
There is \(\Theta_{\rm rem}>0\), depending only on these data and the
fixed buffer and bounded-geometry package, with the following
 property.  Let \(J=[0,\Theta]\), with
 \(0<\Theta\leq\Theta_{\rm rem}\), and let
 \(\mathcal U\Subset\mathcal U^+\Subset\mathcal U^{++}\) be a
 scale-one triply buffered bounded-geometry chart.  Consider on
 \(\mathcal U^{++}\)
\[
 \partial_\vartheta X
 -A^{ij}(\vartheta,x)\nabla_i\nabla_jX
 -B^i(\vartheta,x)\nabla_iX-C(\vartheta,x)X=\mathscr G .
\]
Assume that the principal matrix acts scalarly on the bundle fibers,
that \(A,B,C\) are strongly Bochner measurable in \(\vartheta\), that
\(A\) has ellipticity constants in \([\Lambda^{-1},\Lambda]\), and that
\[
 \operatorname*{ess\,sup}_{\vartheta\in J}
 \left(
  \|A(\vartheta)\|_{C^{q-2,\alpha}}
  +\|B(\vartheta)\|_{C^{q-2,\alpha}}
  +\|C(\vartheta)\|_{C^{q-2,\alpha}}
 \right)\leq\Lambda .
\]
Let \(X\) be the restriction to
\(\mathcal U^{++}\times J\) of a mild solution on an ambient cylinder
containing
\(\overline{\mathcal U^{++}}\times J\).  If
\(X(0)\in C^{q,\alpha}(\mathcal U^{++})\) and
\(\mathscr G\in
L^\infty(J;C^{q-2,\alpha}(\mathcal U^{++}))\), then \(X\) satisfies
\begin{equation}\label{eq:short-effective-clock-propagation}
 \begin{split}
 &\|X\|_{L^\infty(J;C^{q,\alpha}(\mathcal U))}
 +\|\partial_\vartheta X\|_
       {L^\infty(J;C^{q-2,\alpha}(\mathcal U))}\\
 &\qquad\leq
 C e^{C\Theta}\left(
  \|X(0)\|_{C^{q,\alpha}(\mathcal U^{++})}
  +
  \|\mathscr G\|_
       {L^\infty(J;C^{q-2,\alpha}(\mathcal U^{++}))}
  +\|X\|_{L^\infty(J;C^0(\mathcal U^{++}))}
 \right).
 \end{split}
\end{equation}
Here \(C\) depends only on \(q,\alpha,\Lambda\), the dimension, and the
fixed buffer and bounded-geometry data.  No time-continuity or common
time modulus of the coefficients is required.  This is an a priori
interior estimate for the restriction of an ambient solution, not a
locally determined initial-value problem on
\(\mathcal U^{++}\) without lateral data; the outer \(C^0\) term
records precisely that possible lateral influx.

The same estimate holds for one global section on a uniformly locally
finite triply buffered atlas, after taking the supremum of the
localized norms on the indicated inner and outer members, and with the
polynomial weights and tensor normalization used in
\eqref{eq:weighted-parabolic-spaces}.  In this atlas form the outer
\(C^0\) term is the global weighted \(C^0\) norm.  It also holds for
finite triangular systems, differences, and first variations, provided
their displayed spatial coefficient and forcing bounds are uniform.
More explicitly, if \(X(\tau)\) is one global normalized-time
solution, then on each atlas member \(\mathcal U\) entered at time
\(s_{\mathcal U}\) one applies the estimate to
\[
 X_{\mathcal U}(\vartheta)
 :=
 \mathscr S_{\mathcal U}
 X\bigl(\tau_{\mathcal U}(\vartheta)\bigr),
 \qquad
 \vartheta\in J_{\mathcal U},
 \qquad |J_{\mathcal U}|\leq\Theta_{\rm rem}.
\]
No common numerical effective clock is asserted across different
members; compatibility comes from restricting the same global
normalized-time solution.
\end{lemma}

\begin{proof}
After scaling, work first on a Euclidean unit ball.  At a fixed spatial
center freeze only in space:
\[
 A_0^{ij}(\vartheta):=A^{ij}(\vartheta,x_0).
\]
The matrices \(A_0(\vartheta)\) remain strongly measurable and
uniformly elliptic.  Their evolution family has the standard dyadic
bound
\[
 \|\Delta_jU_0(\vartheta,s)\|_{L^\infty\to L^\infty}
 \leq C e^{-c2^{2j}(\vartheta-s)},
 \qquad0\leq s\leq\vartheta\leq\Theta .
\]
This follows directly from the Fourier multiplier with covariance
\(\int_s^\vartheta A_0(r)\,dr\); only ellipticity and measurability in
time enter.  The same spatial-block calculation as in the proof of
Lemma~\ref{lem:weighted-prepared-Schauder} therefore gives the
Bochner \(C^{q,\alpha}\) estimate for the spatially frozen operator,
with no time freezing.

Choose the spatial localization radius uniformly small, take a fixed
buffered cutoff \(\zeta\), and write \(Y=\zeta X\).  We first record
explicitly the endpoint \(q=2\).  No spatial derivative is then applied
to any coefficient.  Uniformly for almost every \(\vartheta\), on a
ball \(B_r(x_0)\) the ordinary H\"older product estimate gives
\[
 \begin{split}
 \bigl\|(A-A_0)D^2Y\bigr\|_{C^\alpha}
 \leq{}&
 C\|A-A_0\|_{C^0(B_r)}
       \|Y\|_{C^{2,\alpha}(B_r)}\\
 &+C[A]_{C^\alpha(B_r)}
       \|D^2Y\|_{C^0(B_r)} .
 \end{split}
\]
Since
\(\|A-A_0\|_{C^0(B_r)}\leq\Lambda r^\alpha\), the first term is
absorbed after fixing \(r=r(\alpha,\Lambda)\) sufficiently small.
Spatial interpolation on the fixed buffered pair gives, for every
\(\epsilon>0\),
\[
 \|D^2Y\|_{C^0}
 +\|DY\|_{C^\alpha}
 +\|Y\|_{C^\alpha}
 \leq
 \epsilon\|Y\|_{C^{2,\alpha}}
 +C_\epsilon\|Y\|_{C^0}.
\]
Consequently,
\[
 \bigl\|(A-A_0)D^2Y\bigr\|_{C^\alpha}
 +\|BDY\|_{C^\alpha}
 +\|CY\|_{C^\alpha}
 \leq
 \epsilon\|Y\|_{C^{2,\alpha}}
 +C_{\epsilon,\Lambda}\|Y\|_{C^0}.
\]
Thus the endpoint uses only
\(A,B,C\in C^\alpha\), and no classical derivative of a coefficient is
taken.

For \(q\geq3\), commute only the derivatives \(D^\beta\) with
\(|\beta|\leq q-2\).  The term in which no derivative lands on \(A\) is
treated by the same freezing and absorption.  Every remaining
principal term has the form
\[
 (D^\gamma A)D^{\beta-\gamma+2}Y,
 \qquad
 1\leq|\gamma|\leq|\beta|\leq q-2,
\]
so the coefficient derivative is available under the stated
\(C^{q-2,\alpha}\) hypothesis and the derivative of \(Y\) is strictly
below order \(q\).  Spatial interpolation treats it.  The analogous
terms from \(B\,DY\), \(CY\), and the fixed partition commutators have
still lower differential order.  Hence no derivative of \(A,B,C\)
beyond order \(q-2\) is used.  After
\(\Theta_{\rm rem}\) is fixed sufficiently small, a finite buffered
localization followed by Gronwall gives the factor \(Ce^{C\Theta}\) in
\eqref{eq:short-effective-clock-propagation}.  The outer \(C^0\) term
accounts for lateral influx into an individual inner chart.  For a
global section, take the supremum over the triply buffered atlas;
bounded overlap, spatial interpolation, and the outer \(C^0\) term
close all cutoff commutators.  Neighboring polynomial weights are
comparable, and tensor normalization is constant on each scaled chart,
proving the weighted statement.  Subtraction and differentiation of
the equation give the final two assertions.
\end{proof}

\begin{remark}[Why the source-adapted physical clock is required]
\label{rem:annular-time-scaling}
In the right-translated chart for the relative harmonic-map flow,
\[
 \partial_\tau F=\lambda\Delta_{\acute G,S}F,
 \qquad S=\lambda\Theta^*\bar g.
\]
In local coordinates,
\[
 (\partial_\tau F)^A
 =\lambda\acute G^{ij}
 \left(
  \partial_i\partial_jF^A
  -\Gamma(\acute G)_{ij}^{p}\partial_pF^A
  +\Gamma(S)_{BC}^{A}(F)\partial_iF^B\partial_jF^C
 \right).
\]
At an exact outer-model state,
\(\acute G=\lambda\Theta^*\bar g\) and
\(\Theta\simeq\varphi_\tau\).  On a fixed physical annulus,
\(\lambda\acute G^{-1}\) therefore retains the small physical-time
factor; multiplying by \(L\) in the clock \((\tau-s)/L\) does not
remove it.  Near the collapsing core, on the other hand, the natural
source harmonic scale is \(r\simeq\lambda^{1/2}\).  A single raw
annular clock cannot cover both regimes.

In the frozen source-adapted chart
\((\mathcal U,r_{\mathcal U})\), multiplication by
\(r_{\mathcal U}^2/\lambda_s^\circ\) and the change
\eqref{eq:source-adapted-clock} turn the principal matrix into
\[
 r_{\mathcal U}^2\frac{\lambda}{\lambda_s^\circ}
 \acute G^{-1},
\]
which is uniformly elliptic by
\eqref{eq:reference-clock-comparison} and
\eqref{eq:source-adapted-ellipticity}.  Constant metric rescaling does
not alter Christoffel symbols, and
Lemma~\ref{lem:prepared-chart-calculus} controls the remaining
coefficients.  The clock is fixed from the center path, rather than
from a trial state, so differences of
\(\lambda/\lambda_s^\circ\) appear explicitly as ordinary
principal-coefficient differences in the contraction and variational
equations.  On a bounded \(\mathbb S^{k+2,\alpha}\) ball the ODE,
Ricci--DeTurck, and map bounds give a common coefficient modulus
\(\omega(\delta)\leq C_K\delta^{\alpha/2}\).
\end{remark}

\subsubsection{The feedback functional}

\begin{lemma}[No derivative loss in the adaptive Gram feedback]
\label{lem:feedback-functional-regularity}
On a common-margin prepared ball on which the adaptive Gram matrix is
invertible, the exact feedback law
\begin{equation}\label{eq:algebraic-feedback-map}
 c=\mathscr C(\tau,\mathbf z)
   =-M(\tau,\mathbf z)^{-1}d(\tau,\mathbf z)
\end{equation}
is a \(C^1\) map of \((\tau,\mathbf z)\) into \(\mathbb R^9\).
For every \(3\leq r\leq k\) and two states in that ball,
\begin{equation}\label{eq:feedback-local-Lipschitz}
 |\mathscr C(\tau,\mathbf z_1)
       -\mathscr C(\tau,\mathbf z_2)|
 \leq C
 \|\mathbf z_1-\mathbf z_2\|_{\mathscr Y_{\rm prep}^{r,\alpha}},
\end{equation}
and, for every tangent state \(\dot{\mathbf z}\),
\begin{equation}\label{eq:feedback-C1-bound}
 |D_{\mathbf z}\mathscr C(\tau,\mathbf z)
       [\dot{\mathbf z}]|
 \leq C\|\dot{\mathbf z}\|_{\mathscr Y_{\rm prep}^{r,\alpha}} .
\end{equation}
On a bounded order-\((k+2,\alpha)\) prepared ball it also satisfies
the mixed tame remainder estimate
\begin{equation}\label{eq:feedback-mixed-tame-remainder}
\begin{split}
 &|\mathscr C(\tau,\mathbf z+\mathbf u)
       -\mathscr C(\tau,\mathbf z)
       -D\mathscr C(\tau,\mathbf z)[\mathbf u]|\\
 &\hspace{18mm}\leq
 C_K\|\mathbf u\|_{\mathscr Y_{\rm prep}^{k+2,\alpha}}
      \|\mathbf u\|_{\mathscr Y_{\rm prep}^{k,\alpha}} .
\end{split}
\end{equation}
The constants are uniform on finite normalized-time intervals and on
common-margin prepared balls.  In particular, after substitution in
the coupled equations, the feedback is a finite-rank lower-order term
and costs no derivative beyond the second-order parabolic principal
part.  The harmless threshold \(r\geq3\) is the threshold for \(C^1\)
composition in Lemma~\ref{lem:prepared-chart-calculus}; the
differential expression itself uses only the two-jet described in the
proof.
\end{lemma}

\begin{proof}
The entries of \(M\) are Gaussian pairings of the effective columns
and their terms linear in \(h\).  Lemma~\ref{lem:static-prepared-columns}
and Lemma~\ref{lem:prepared-chart-calculus} show that these entries are
\(C^1\) functions of the prepared state at the displayed orders.
Uniform Gram invertibility makes \(M\mapsto M^{-1}\) a \(C^1\) map
with uniformly bounded derivative.

It remains to check the terms in \(d\) which appear to contain two
derivatives.  Work first with a smooth state.  Every such term is
paired with \(\rho_\tau Z_\mu\) against \(d\nu\).  A typical
quadratic term has the form
\[
 \int_M
 \langle h*\bar\nabla^2h,\rho_\tau Z_\mu\rangle\,d\nu .
\]
Because \(\rho_\tau\) is compactly supported, integration by parts
has no boundary contribution and writes this as a sum of integrals of
\[
 \bar\nabla h*\bar\nabla h*\rho_\tau Z_\mu,\qquad
 h*\bar\nabla h*
 \bar\nabla(\rho_\tau Z_\mu e^{-\bar f}).
\]
The derivatives of \(Z_\mu\), \(\rho_\tau\), and \(\bar f\) have
fixed polynomial growth on their supports, and the Gaussian measure
absorbs every such polynomial.  The same calculation applied to the
principal Ricci part of the pure graft defect moves one derivative
from the second jet of the metric difference onto the fixed cutoff,
the Gaussian test tensor, or the coefficient of that difference.
Thus the resulting functional uses at most first derivatives of the
varying metric tensors and at most the second derivatives of the
charted maps which occur in pullback and inversion.  The
moving-cutoff commutator is already first order.

The scaled product and composition estimates in
Lemma~\ref{lem:prepared-chart-calculus}, followed by the Gaussian tail
estimate, now give
\[
 |M_1-M_2|+|d_1-d_2|
 \leq C
 \|\mathbf z_1-\mathbf z_2\|_{\mathscr Y_{\rm prep}^{r,\alpha}}.
\]
The same integration-by-parts formula may be differentiated in a
tangent direction.  It gives the identical bound for \(DM\) and
\(Dd\).  Density of smooth states in the little-H\"older completion
extends the identities to the entire prepared ball.  Finally,
\[
 D\mathscr C[\dot{\mathbf z}]
 =M^{-1}(DM[\dot{\mathbf z}])M^{-1}d
  -M^{-1}Dd[\dot{\mathbf z}],
\]
which proves \eqref{eq:feedback-local-Lipschitz} and
\eqref{eq:feedback-C1-bound}.  Subtracting the displayed linearization
from the exact integrated formulas leaves terms with at least two
state increments.  The same one-high--one-low allocation as in
\eqref{eq:prepared-mixed-tame-remainder}, followed by the ordinary
finite-dimensional quadratic remainder for \(M\mapsto M^{-1}\), proves
\eqref{eq:feedback-mixed-tame-remainder}.
\end{proof}

\subsubsection{Construction and continuation}

Recall Definition~\ref{def:admissible-first-exit-interval}.  For later
reference, its exact system on an interval \(I\) is the tuple
\[
 (G(t),t(\tau),\lambda(\tau),\Theta_\tau,F_\tau,
   c(\tau)=(a(\tau),b(\tau)))
\]
which satisfies the exact system
\[
 \begin{aligned}
  \partial_tG&=-2\Ric_G,&
  t_\tau&=\lambda,&
  \lambda_\tau&=-(1+a)\lambda,\\
  \partial_\tau\Theta
  &=\bigl((1+a)\bar\nabla\bar f-U_b\bigr)\circ\Theta,&
  \partial_\tau F&=\lambda\Delta_{\acute G,S}F,
 \end{aligned}
\]
where
\[
 U_b=\sum_{j=1}^8b_j\chi_\tau W_j,\qquad
 S=\lambda\Theta^*\bar g,\qquad
 \acute G=\eta\,\iota_*G+(1-\eta)S,
\]
\[
 \Phi=\Theta\circ F,\qquad
 h=\lambda^{-1}(\Phi^{-1})^*\acute G-\bar g,\qquad
 c=\mathscr C(\tau,\mathbf z)
\]
is the exact adaptive Gram feedback
\eqref{eq:algebraic-feedback-map}.  In addition it must satisfy the
structural graph identities, have degree-one proper prepared maps, and
remain in one stated common-margin prepared ball.  The same system is
repeated with equation labels in
\eqref{eq:coupled-scale-system}--\eqref{eq:coupled-map-system}.
When two evolutions are called admissible, all numerical margins,
buffered physical covers, reference gauges, and scaled coefficient
bounds are common.

For the fixed slice chart used in the evolution, choose once a
sufficiently small neighborhood
\(\mathscr N_{\mathrm{sl},\tau_0}^{3,\alpha}\) furnished by
Proposition~\ref{prop:intro-static-sliced-chart}.  For every \(r\geq3\)
use the compatible higher-regularity restrictions
\[
 \mathscr N_{\mathrm{sl},\tau_0}^{r,\alpha}
 :=
 \mathscr N_{\mathrm{sl},\tau_0}^{3,\alpha}
 \cap\mathscr P_{\tau_0}^{r,\alpha},
\]
and define
\begin{equation}\label{eq:preliminary-sliced-space}
 \Sigma_{\tau_0}^{r,\alpha}
 :=
 \left\{
  \mathbf z\in\mathscr N_{\mathrm{sl},\tau_0}^{r,\alpha}:
  \left\langle\rho_{\tau_0}h(\mathbf z),Z_\mu\right\rangle=0
  \ \text{for }0\leq\mu\leq8
 \right\}.
\end{equation}
Thus the sliced spaces at different regularity orders refer to the
same local component and satisfy the canonical nesting
\[
 \Sigma_{\tau_0}^{s,\alpha}
 =
 \Sigma_{\tau_0}^{r,\alpha}
 \cap\mathscr P_{\tau_0}^{s,\alpha},
 \qquad s\geq r\geq3.
\]
On a two-derivative-buffered common-margin neighborhood,
\(\Pi_{\rm sl}^{\,r+2\to r}\) denotes the unique centered prepared
phase retraction imposing these nine moments.  Its split \(C^1\)
submanifold structure, existence, uniqueness, and mapping properties
were established, independently of the evolution, in
Proposition~\ref{prop:intro-static-sliced-chart}.  The detailed
tuple-valued phase formula is
\eqref{eq:intro-prepared-scale-leg}--%
\eqref{eq:intro-prepared-phase-action}; its scale-normalized estimates
and the corresponding intrinsic manifold statement are expanded later
in Proposition~\ref{prop:sliced-prepared-manifold}.  Neither the local
evolution construction nor its proof is used to obtain the static
retraction.

\begin{proposition}[Local coupled feedback evolution]
\label{prop:coupled-local-feedback}
Fix \(k\geq12\) and \(0<\alpha<1\).  Let
\[
 \mathbf z_0=(G_0,\lambda_0,R_0,F_0)
 \in\mathscr P_{\tau_0}^{k+2,\alpha},
 \qquad
 \Theta_0=\varphi_{\tau_0}\circ R_0,
\]
belong to a common-margin prepared ball with recorded ceiling
\(K_{\rm init}^{k+2,\alpha}\).  Equivalently, \(G_0\) is a
\(C^{k+2,\alpha}\) closed metric, \(R_0,F_0\) and their inverses have
the scaled \(C^{k+3,\alpha}\) bounds of the prepared chart, and the
following strict conditions hold:
\begin{enumerate}
\item the adaptive extension \(\acute G_0\) and
      \(S_0=\lambda_0\Theta_0^*\bar g\) are complete, uniformly
      equivalent metrics with the required bounded geometry;
\item
      \[
       h_0=\lambda_0^{-1}
          ((\Theta_0\circ F_0)^{-1})^*\acute G_0-\bar g
      \]
      belongs to the small \(C^2\) box, satisfies the exact receding
      slice, and has invertible adaptive Gram matrix;
\item the graft compatibility, support separation, radial comparison,
      and outer bounded geometry inequalities all hold with positive
      margin.
\end{enumerate}
Then there is \(\delta_\tau>0\) and a unique tuple
\[
 \bigl(G(t),t(\tau),\lambda(\tau),\Theta_\tau,F_\tau,
       c(\tau)=(a(\tau),b(\tau))\bigr),
 \qquad \tau_0\leq\tau\leq\tau_0+\delta_\tau,
\]
in the time-dependent prepared space above, with
\[
 G(0)=G_0,\qquad t(\tau_0)=0,\qquad
 \lambda(\tau_0)=\lambda_0,\qquad
 \Theta_{\tau_0}=\Theta_0,\qquad F_{\tau_0}=F_0,
\]
and
\begin{align}
 \partial_tG&=-2\Ric_G,&
 t_\tau&=\lambda,&
 \lambda_\tau&=-(1+a)\lambda,
 \label{eq:coupled-scale-system}\\
 \partial_\tau\Theta
 &=\bigl((1+a)\bar\nabla\bar f-U_b\bigr)\circ\Theta,
 &
 \partial_\tau F
 &=\lambda\Delta_{\acute G,S}F .
 \label{eq:coupled-map-system}
\end{align}
Here \(U_b=\sum_{j=1}^8b_j\chi_\tau W_j\),
\[
 S=\lambda\Theta^*\bar g,\qquad
 \acute G=\eta\,\iota_*G+(1-\eta)S,
\]
and \(c=\mathscr C(\tau,\mathbf z)\) is the unique solution of the
exact adaptive Gram system obtained from \eqref{eq:receding-Gram} by
using the effective columns
\eqref{eq:effective-column-zero}--\eqref{eq:effective-column-j}.
Equivalently, regarding \(F\) as a function of physical time,
\(\partial_tF=\Delta_{\acute G,S}F\).

The maps \(\Theta_\tau\) and \(F_\tau\) remain global
diffeomorphisms.  Moreover,
\[
 \Phi_\tau=\Theta_\tau\circ F_\tau,\qquad
 h(\tau)=\lambda(\tau)^{-1}
          (\Phi_\tau^{-1})^*\acute G_\tau-\bar g
\]
solve the adaptive normalized equation, and the exact receding slice
is propagated.  Smooth initial data give a smooth solution.  On every
common compact subinterval of existence the solution depends
 continuously on the initial state.  In the prepared Banach charts the
 local solution map on the common existence interval is
\[
 C^1:\Sigma_{\tau_0}^{k+2,\alpha}
        \longrightarrow\mathscr P_{\tau}^{k,\alpha}.
\]
More precisely, on every fixed uniformly interior common-margin
subball with one recorded value of \(K_{\rm init}^{k+2,\alpha}\),
after choosing one common short existence interval, there is a
constant \(C\) such that
\begin{equation}\label{eq:local-global-prepared-Lipschitz}
 \sup_{\tau_0\leq\tau\leq\tau_0+\delta_\tau}
 d_{\rm prep}^{k,\alpha}
 \bigl(\mathcal S_{\tau,\tau_0}(\mathbf z_{1,0}),
       \mathcal S_{\tau,\tau_0}(\mathbf z_{2,0})\bigr)
 \leq
 C\|\mathbf z_{1,0}-\mathbf z_{2,0}\|
       _{\mathscr E_{\rm prep}^{k+2,\alpha}} .
\end{equation}
The left side retains the global graph-augmented prepared distance,
including the supremum over the complete source atlas; the right side
is the independent model increment norm.  The preceding \(C^1\)
assertion, and hence its domain tangent, uses the independent model
norm \eqref{eq:prepared-model-Banach-norm}.
Here the domain is the local sliced manifold of
Proposition~\ref{prop:sliced-prepared-manifold}; its tangent vectors
satisfy the linearized nine moment constraints.  The corresponding
ambient map is obtained only after phase projection:
\begin{equation}\label{eq:local-ambient-solution-map}
 \operatorname{dom}\Pi_{\rm sl}^{\,k+4\to k+2}
 \xrightarrow{\ \Pi_{\rm sl}^{\,k+4\to k+2}\ }
 \Sigma_{\tau_0}^{k+2,\alpha}
 \xrightarrow{\ \mathcal S_{\tau,\tau_0}\ }
 \mathscr P_{\tau}^{k,\alpha}.
\end{equation}
By \eqref{eq:quantitative-phase-retraction-Lipschitz}--%
\eqref{eq:quantitative-phase-retraction-derivative}, its local ambient
Lipschitz and differential constants are bounded by the corresponding
sliced constants times \(K_{\Pi,k+2}\).

Let \([\tau_0,\tau^\dagger)\) be the maximal half-open continuation of
this local coupled solution before the desired physical singular time.
If \(\tau^\dagger<\infty\), then at least one of the following occurs:
\begin{enumerate}
\item the closed Ricci flow has unbounded curvature at
      \(t(\tau^\dagger)\);
\item the coarse global \(C^2\) extension box or uniform ellipticity
      is lost;
\item the finite-endpoint accumulated phase budget reaches
      \(\varepsilon_{\rm ph}\);
\item the adaptive Gram matrix loses invertibility;
\item graft compatibility or the required source/target
      bounded-geometry package is lost;
\item \(F\) ceases to be a controlled diffeomorphism, meaning that
      the proper radial comparison, local-invertibility bound, or
      the corresponding inverse-map bound degenerates.
\end{enumerate}
Here item \textup{(1)} occurs when
\[
 \limsup_{\tau\uparrow\tau^\dagger}
 \|\Rm_{G(t(\tau))}\|_{L^\infty(\mathcal X,G(t(\tau)))}=\infty,
\]
 and items \textup{(2)}--\textup{(6)} occur in the terminal ordinary-margin
 or witnessed-package sense fixed in
 Definition~\ref{def:admissible-first-exit-interval}.  Thus failure of
 all six alternatives supplies, on one terminal subinterval, common
 positive margins for the ordinary faces and common positive
 reserve-to-operative gaps for every witnessed harmonic condition
 contained in the bounded-geometry package.
On a sufficiently small admissible box, item \textup{(4)} is excluded
by Corollary~\ref{cor:adaptive-feedback-package}.  Under
Corollary~\ref{cor:adaptive-auxiliary-closure}, the
higher-regularity portion of item \textup{(5)} is not an independent
exit face.
\end{proposition}

\begin{proof}
For notational economy the construction is written on the first
interval.  At a later restart, every occurrence below of
\(\tau_0,\lambda_0,\Theta_0,\Phi_0\) and every integral beginning at
\(\tau_0\) means, respectively,
\(s,\lambda(s),\Theta_s^\circ,\Phi_s^\circ\) and an integral beginning
at \(s\), in the frozen charts of
\eqref{eq:extended-prepared-distance}.  The common-margin hypotheses
and the recorded ceiling \(K_{\rm init}^{k+2,\alpha}\) make all
constants uniform under this substitution.  In a chart frozen at
\(s\), put \(\lambda_s:=\lambda(s)\).  The scalar used in the local
fixed-point system is the dimensionless elapsed compact clock
\[
 \theta(\tau)
 :=\frac{t(\tau)-t(s)}{\lambda_s},
\]
and hence has left value zero; the corresponding global physical time
is \(t(s)+\lambda_s\theta\).

\emph{Step 1: reduction to a strictly parabolic--ODE system.}
Use the frozen compact atlas and smooth reference metric specified in
the definition of \(\mathbb S^{r,\alpha}(I)\), set
\[
 g_{s,0}:=\lambda_s^{-1}G(t(s)),
\]
and solve for the zero-trace affine increment
\(H_s=\widetilde G_s-g_{s,0}\).  In the dimensionless compact time
\(\widehat\sigma=(u-t(s))/\lambda_s\), its equation is
\[
 \partial_{\widehat\sigma}H_s
 =\mathcal R_{{\rm DT},s}(g_{s,0}+H_s),
 \qquad H_s(0)=0.
\]
The uniform ellipticity
\eqref{eq:compact-affine-ellipticity} and the normalized carrier bounds
allow the compact quasilinear Schauder theorem to be applied in the
single frozen atlas.  It gives a unique \(H_s\) on
\(\widehat I_{\rm ph}=[0,\widehat\delta_{\rm ph}]\), with
\(\widehat\delta_{\rm ph}>0\) depending only on the recorded normalized
prepared bound and margins.  Equivalently, the forward physical
lifespan is
\[
 \delta_{\rm ph}(s)=\lambda_s\widehat\delta_{\rm ph}.
\]
Its solution map is \(C^1\) from input order \(k+2\) to output order
\(k\).  The reconstruction ODE and its sign are those fixed in the
definition of \(\mathbb S^{r,\alpha}(I)\); the two-order buffer controls
that ODE and its first variation.

Let \(\chi_s\) be the reconstructing diffeomorphism and define the
zero-trace reconstructed normalized increment
\[
 K_s(\widehat\sigma)
 :=\chi_s(\widehat\sigma)^*
      (g_{s,0}+H_s(\widehat\sigma))-g_{s,0},
 \qquad0\leq\widehat\sigma\leq\widehat\delta_{\rm ph}.
\]
Then \(K_s(0)=0\), and
\[
 \lambda_sK_s(\widehat\sigma)
 =G(t(s)+\lambda_s\widehat\sigma)-G(t(s)).
\]
We now establish the same-order reconstruction estimate needed for the
uniform little-H\"older bounds.  Put \(r:=k+2\).  The transition maps and inverse
transition maps of \(\mathfrak A_{\mathcal X,s}\), the subordinate
cutoffs on the doubled buffers, and the local-addition maps used to
represent compact diffeomorphisms are taken with the uniform
order-\((r+1,\alpha)\) bounds supplied by the normalized carrier
package.  The same bounds hold for the coefficients and inverse
coefficients of \(\widehat G_s\).  These are bounds on fixed
background objects; no order-\((r+1,\alpha)\) bound is imposed on the
varying tensors \(g_{s,0}\) or \(H_s\).
Write \(M_{\rm rec}\geq1\) for one common numerical bound for the
compact ellipticity, these fixed order-\((r+1,\alpha)\) background
norms, and the recorded order-\((r,\alpha)\) norms of
\(g_{s,0}\) and \(H_s\).  The subscript distinguishes this constant
from the manifold \(M\).  All absorbable paradifferential remainders
below are written with a parameter \(\varepsilon_{\rm rec}>0\).  The
frequency separations and spatial freezing radii are first chosen so
that those remainders have coefficient \(\varepsilon_{\rm rec}\);
after all fixed paracomposition and block-reduction constants have
been determined, \(\varepsilon_{\rm rec}\) is chosen small enough to
close the final absorption.

\medskip
\noindent\emph{Same-order reconstruction sublemma.}
On the common normalized compact-time interval
\(\widehat I_{\rm ph}\), the increment \(K_s\) defined above belongs to
\[
 \mathbb E_{\mathcal X,s,0}^{r,\alpha}(\widehat I_{\rm ph}).
\]
Its norm is bounded by a constant depending only on the recorded
normalized compact package, the ellipticity constants in
\eqref{eq:compact-affine-ellipticity}, and the common compact
quasilinear bound.

More generally, consider a family in the one frozen chart for which
these bounds are common.  Suppose that the bounded nondecreasing
functions \(\eta_0,\eta_H\), both vanishing at zero and extended
constantly for arguments at least \(1\), satisfy
\[
 \operatorname{Tail}_{r,\alpha}(g_{s,0};\varrho)
 \leq\eta_0(\varrho),
 \qquad
 \sup_{\widehat\sigma\in\widehat I_{\rm ph}}
 \operatorname{Tail}_{r,\alpha}
       (H_s(\widehat\sigma);\varrho)
 \leq\eta_H(\varrho).
\]
There are constants \(C_{\rm rec}\geq1\) and \(C_K<\infty\),
depending only on the preceding common bounds and not on the rates at
which \(\eta_0,\eta_H\) vanish, such that
\begin{equation}\label{eq:compact-reconstruction-tail}
 \operatorname{Tail}_{k+2,\alpha}(K_s;\varrho)
 \leq C_K\bigl(
     \eta_0(C_{\rm rec}\varrho)
    +\eta_H(C_{\rm rec}\varrho)
    +\min\{\varrho,1\}^{1-\alpha}\bigr).
\end{equation}
The tail on the left includes the supremum over
\(\widehat\sigma\in\widehat I_{\rm ph}\).

\smallskip
\noindent\emph{Proof of the sublemma.}
Put \(m:=r-2=k\), write

\[
 t:=\widehat\sigma,\qquad
 T:=\widehat\delta_{\rm ph},\qquad
 \widetilde g:=g_{s,0}+H_s,\qquad
 g:=\chi_s^*\widetilde g=g_{s,0}+K_s,
\]

and let \(D:=\nabla^{\widehat G_s}\) be the Levi--Civita connection of the
fixed smooth reference metric \(\widehat G_s\).  We first treat smooth
data, so every top-order envelope introduced below is finite.  All
constants are uniform under fixed-atlas smoothing and depend
only on the package recorded before the statement of the sublemma.

We use the spatial Littlewood--Paley operators from the proof of
Lemma~\ref{lem:weighted-prepared-Schauder}.  For a time-dependent
tensor \(u\), set

\[
 b_\ell^q(u;t)
 :=2^{(q+\alpha)\ell}
   \|\Delta_\ell u\|_{L^\infty([0,t]\times\mathcal X)},
\]

and define

\[
 \mathcal W_j^q(u;t):=
 \max\left\{
  \sup_{\ell\geq j}b_\ell^q(u;t),
  \sup_{-1\leq\ell<j}
   2^{-(1-\alpha)(j-\ell)}b_\ell^q(u;t)
 \right\}.
\]

For a time-independent tensor the time supremum is omitted.  Notice
that for every fixed integer \(J\geq0\),

\begin{equation}\label{eq:reconstruction-envelope-shift}
 \mathcal W_{j-J}^q(u;t)
 \leq 2^{(1-\alpha)J}\mathcal W_j^q(u;t).
\end{equation}

Thus fixed frequency shifts of an unknown never cause an iterated
spatial dilation.  The two inequalities in
\eqref{eq:tail-frequency-envelope-equivalence} will be used only at
the end.

We begin with the genuinely lower-order reconstruction, for which
ordinary factorwise calculus has a two-derivative buffer.  The
DeTurck field \(B_s(\widetilde g)\) is of order \(r-1\).  Differentiating
its flow equation through order \(r-1\), and using Gronwall, controls
\(\chi_s\) at that order.  In the top differentiated term either the
top block lies on \(B_s(\widetilde g)\), and hence on an order-\(r\)
block of \(g_{s,0}\) or \(H_s\), or it lies on \(\chi_s\) and is
multiplied by a bounded first derivative of the vector field.  Every
other Fa\`a di Bruno term has one unused spatial derivative.  Applying
the same alternatives to

\[
 g=\chi_s^*\widetilde g,\qquad
 \mathcal R:=\Rm(g)=\chi_s^*\Rm(\widetilde g),\qquad
 \mathcal A:=\nabla^g-D,
\]

gives, for one fixed \(J_0\),

\begin{equation}\label{eq:reconstruction-buffered-lower}
\begin{split}
 &\|g\|_{L^\infty_tC_x^{m,\alpha}}
 +\|\mathcal R\|_{L^\infty_tC_x^{m,\alpha}}
 +\|\mathcal A\|_{L^\infty_tC_x^{m-1,\alpha}}
 \leq C M_{\rm rec},\\
 &\mathcal W_j^m(g;T)
  +\mathcal W_j^m(\mathcal R;T)
  +\mathcal W_j^{m-1}(\mathcal A;T)\\
 &\qquad\leq C\left(
   \mathcal W_{j-J_0}^r(g_{s,0})
  +\mathcal W_{j-J_0}^r(H_s;T)
  +M_{\rm rec}2^{-(1-\alpha)j}
 \right).
\end{split}
\end{equation}

This step uses only the order-\((r-1,\alpha)\) diffeomorphism and the
order-\((r-2,\alpha)\) pullback; it is not the same-order conclusion
being proved.  The cutoff and chart-transition terms in this
calculation are covered by the almost-diagonal transition estimate
proved in Lemma~\ref{lem:weighted-prepared-Schauder}.  Since the fixed
background objects have one further derivative, their top blocks give
the final term in \eqref{eq:reconstruction-buffered-lower}.  The same
calculation at \(t=0\), now without loss because
\(g(0)=g_{s,0}\), gives

\begin{equation}\label{eq:reconstruction-initial-curvature-connection}
 \mathcal W_j^m(\mathcal R(0))
 +\mathcal W_j^{m+1}(\mathcal A(0))
 \leq C\left(
   \mathcal W_{j-J_0}^r(g_{s,0})
  +M_{\rm rec}2^{-(1-\alpha)j}
 \right).
\end{equation}
Here, in particular,
\[
 \mathcal A(0)=\nabla^{g_{s,0}}-D;
\]
it is controlled by the initial metric in
\eqref{eq:reconstruction-initial-curvature-connection} and is not
being set equal to zero.

We next recover the two missing derivatives without requiring
higher-order control of the reconstructing diffeomorphism.  For
\(0\leq t\leq T\), define

\[
 \mathbf A_j(t):=\mathcal W_j^{m+1}(\mathcal A;t),
 \qquad
 \mathbf G_j(t):=\mathcal W_j^{m+2}(g;t).
\]

By the Moser expansion for inversion on the uniformly elliptic set,
the order-\((m+2,\alpha)\) envelope of \(g^{-1}\) is bounded by a fixed
multiple of \(\mathbf G_j(t)+M_{\rm rec}2^{-(1-\alpha)j}\); thus no
separate inverse-metric unknown is needed below.

For a localized curvature block put

\[
 z_\ell(t):=
 2^{(m+1+\alpha)\ell}
 \int_0^t2^\ell
   \|\Delta_\ell\mathcal R(\theta)\|_{L^\infty_x}\,d\theta
\]

and let \(\mathbf Z_j(t)\) be its two-sided envelope, after taking the
supremum over the finite frozen atlas and its fixed inner cutoffs,

\[
 \mathbf Z_j(t):=
 \max\left\{
  \sup_{\ell\geq j}z_\ell(t),
  \sup_{-1\leq\ell<j}
   2^{-(1-\alpha)(j-\ell)}z_\ell(t)
 \right\}.
\]

The integral is taken after the \(L^\infty_x\)-norm.  This stronger
quantity is what permits the principal freezing error below to be
absorbed without using cancellation in time.

The three exact geometric identities used in the recovery are

\begin{align}
 \partial_t\mathcal R
 &=g^{-1}*D^2\mathcal R
   +\mathcal P_1(g^{-1},\mathcal A)*D\mathcal R
   +\mathcal P_0(g^{-1},\mathcal A,D\mathcal A,
                  \mathcal R,\Rm(D))*\mathcal R,
 \label{eq:reconstruction-curvature-system}\\
 \partial_t\mathcal A^k{}_{ij}
 &=-g^{k\ell}\bigl(
   \nabla_i\Ric_{j\ell}
  +\nabla_j\Ric_{i\ell}
  -\nabla_\ell\Ric_{ij}\bigr)
  =\mathcal Q_1(g^{-1})*D\mathcal R
   +\mathcal Q_0(g^{-1},\mathcal A)*\mathcal R,
 \label{eq:reconstruction-connection-system}\\
 D_i g_{jk}
 &=\mathcal A^p{}_{ij}g_{pk}
   +\mathcal A^p{}_{ik}g_{jp}.
 \label{eq:reconstruction-metric-compatibility}
\end{align}

Here \(\Rm(D):=\Rm(\widehat G_s)\), and \(\mathcal P_a\) and
\(\mathcal Q_a\) are universal smooth contractions on the uniformly
elliptic set.  In
\eqref{eq:reconstruction-curvature-system} the only second-order term
is \(g^{ab}D_aD_b\), acting scalarly on the tensor components.  The
display follows from
\((\partial_t-\Delta_g)\Rm=\Rm*\Rm\) after replacing \(\nabla^g\) by
\(D+\mathcal A\).  The second display is the standard variation
formula for the Levi--Civita connection under
\(\partial_tg=-2\Ric_g\).  Thus these identities contain all
top-order terms, not merely their principal symbols.

First use metric compatibility.  Localizing
\eqref{eq:reconstruction-metric-compatibility}, applying
\(\Delta_\ell\), and using Bony's decomposition gives, once a fixed
low-frequency threshold has been increased,

\begin{equation}\label{eq:reconstruction-metric-from-connection}
 \mathbf G_j(t)
 \leq C\mathbf A_{j-J_1}(t)
   +CM_{\rm rec}2^{-(1-\alpha)j}
   +\frac18\mathbf G_j(t).
\end{equation}

Indeed, the high block of
\(\mathcal A*g\) either lies on \(\mathcal A\), which is the first
term, or lies on \(g\).  In the latter case the derivative on the
left supplies \(2^\ell\), so the low-
\(\mathcal A\)--high-\(g\) paraproduct has coefficient
\(C2^{-\ell}\|\mathcal A\|_{C^0}\), which is at most \(1/16\) above
the fixed threshold.  The resonant term is split by the same two
alternatives.  The finitely many remaining blocks are controlled by
\eqref{eq:reconstruction-buffered-lower}; in the two-sided envelope
they give the displayed \(2^{-(1-\alpha)j}\) term.  A top block of the
fixed connection \(D\), a cutoff, or a transition matrix has one
unused derivative and gives the same remainder.

We now establish the only parabolic estimate needed here.  Work on
one doubled buffered chart, insert an inner cutoff, and freeze in
space, but not in time,

\[
 a_0^{ab}(\theta):=g^{ab}(\theta,x_0).
\]

Let \(U_0(t,\theta)\) be the evolution family of
\(\partial_t-a_0^{ab}(t)\partial_a\partial_b\).  Since the principal
matrix is scalar on the tensor fiber and independent of \(x\), its
Fourier symbol is explicit.  Hence the precise estimate used below is
the sandwiched estimate

\begin{equation}\label{eq:reconstruction-sandwiched-frozen-heat}
 \|\Delta_\ell U_0(t,\theta)\Delta_p\|_{L^\infty\to L^\infty}
 \leq C\mathbf 1_{\{|p-\ell|\leq2\}}
 e^{-c2^{2\ell}(t-\theta)},
 \qquad0\leq \theta\leq t\leq T.
\end{equation}

No estimate for an unsandwiched variable-coefficient evolution is
being asserted.  Duhamel's formula, Fubini, and

\[
 \int_\theta^t2^\ell e^{-c2^{2\ell}(v-\theta)}\,dv
 \leq C2^{-\ell}
\]

give, for the localized curvature and its complete remainder
\(\mathcal E\),

\begin{equation}\label{eq:reconstruction-integrated-block}
 \begin{split}
 2^\ell\int_0^t
  \|\Delta_\ell\mathcal R(v)\|_{L^\infty}\,dv
 \leq{}&C2^{-\ell}
  \|\Delta_{\ell+O(1)}\mathcal R(0)\|_{L^\infty}\\
 &+C2^{-\ell}\int_0^t
  \|\Delta_{\ell+O(1)}\mathcal E(\theta)\|_{L^\infty}\,d\theta.
 \end{split}
\end{equation}

We list all terms in \(\mathcal E\), because the derivative count is
the point of the argument.

\begin{enumerate}
\item The only order-two term is
\((g^{-1}-g^{-1}(\theta,x_0))*D^2\mathcal R\).  In its
low-coefficient--high-curvature paraproduct, the coefficient is made
at most \(\varepsilon_{\rm sp}\) by the same fixed spatial
subdivision used in the proof of
Lemma~\ref{lem:weighted-prepared-Schauder}.  After
\eqref{eq:reconstruction-integrated-block}, its contribution is
exactly
\[
 C\varepsilon_{\rm sp}
 \int_0^t2^\ell
  \|\Delta_{\ell+O(1)}\mathcal R(\theta)\|_{L^\infty}\,d\theta,
\]
and hence is an absorbable multiple of \(\mathbf Z_j(t)\).  If the
top block lies on \(g^{-1}\), it is bounded by
\(\int_0^t\mathbf G_{j-J_2}(\theta)\,d\theta\); if a derivative falls on the
coefficient in the block commutator, the same bound has one additional
frequency gain.  The resonant term has precisely these two
alternatives.

\item In
\(\mathcal P_1(g^{-1},\mathcal A)*D\mathcal R\), a high curvature
block has only one derivative.  The factor \(2^{-\ell}\) in
\eqref{eq:reconstruction-integrated-block} therefore leaves an
additional \(2^{-\ell}\) relative to \(\mathbf Z_j\).  A high block
of \(g^{-1}\) is routed to \(\mathbf G_j\), and a high block of
\(\mathcal A\) is routed to \(\mathbf A_j\); both occur under the
time integral.

\item In the zeroth-order term, a high curvature factor is again
controlled by \(\mathbf Z_j\), now with two unused powers of
\(2^{-\ell}\).  A high \(D\mathcal A\) loses one of those powers and
is bounded by \(\int_0^t\mathbf A_{j-J_2}(\theta)\,d\theta\).  High blocks of
\(g^{-1}\) or \(\mathcal A\) are treated as in item~2.  The fixed
curvature \(\Rm(D)\) has the recorded extra derivative.

For clarity, above the fixed low-frequency threshold the resonant
alternatives just invoked are absolutely
summable, rather than being hidden in a formal product rule.  If
\(\widetilde\Delta_p\) denotes a fixed enlargement of \(\Delta_p\),
then the forcing weight left after
\eqref{eq:reconstruction-integrated-block} satisfies, for one fixed
\(C_0\),
\[
 \begin{split}
 &2^{(m+\alpha)\ell}
  \sum_{p\geq\ell-C_0}
  \Bigl\|\Delta_\ell\bigl(
    (\Delta_pg^{-1})D^2\widetilde\Delta_p\mathcal R
    +(D\Delta_p\mathcal A)\widetilde\Delta_p\mathcal R
  \bigr)\Bigr\|_{L^\infty}\\
 &\qquad\leq
  CM_{\rm rec}2^{-(m+\alpha)\ell}
  \bigl(
    \mathbf G_{\ell-C_0}(t)
   +\mathbf A_{\ell-C_0}(t)+M_{\rm rec}
  \bigr).
 \end{split}
\]
Here \(m=k\geq12\), so the sum is geometric.  Every other resonant
product in items~1--3 has at least the same unused frequency power.

\item The localization commutator is
\[
 -2g^{ab}(D_a\zeta)D_b\mathcal R
 -g^{ab}(D_aD_b\zeta)\mathcal R
 +\text{fixed-connection terms},
\]
so it is covered by items~2 and 3.  A top block of a cutoff,
transition matrix, or fixed connection coefficient has one unused
derivative and contributes
\(CM_{\rm rec}2^{-(1-\alpha)j}\).  Passage between charts uses the
almost-diagonal transition estimate already proved in
Lemma~\ref{lem:weighted-prepared-Schauder}; no finite-band property of
a coordinate change is used.
\end{enumerate}

After multiplying \eqref{eq:reconstruction-integrated-block} by
\(2^{(m+1+\alpha)\ell}\), taking both halves of the envelope, using
\eqref{eq:reconstruction-envelope-shift}, and choosing
\(\varepsilon_{\rm sp}\) so that the principal error is at most
\(\frac18\mathbf Z_j(t)\), these four alternatives give

\begin{equation}\label{eq:reconstruction-Z-estimate}
 \mathbf Z_j(t)
 \leq C\mathbf I_j
 +C\int_0^t
   \bigl(\mathbf A_j(\theta)+\mathbf G_j(\theta)\bigr)\,d\theta
 +\frac18\mathbf Z_j(t),
\end{equation}

where, after increasing one fixed shift \(J_*\),

\begin{equation}\label{eq:reconstruction-input-envelope}
 \mathbf I_j:=
  \mathcal W_{j-J_*}^r(g_{s,0})
 +\mathcal W_{j-J_*}^r(H_s;T)
 +M_{\rm rec}2^{-(1-\alpha)j}.
\end{equation}

Finally integrate
\eqref{eq:reconstruction-connection-system} and apply the same Bony
decomposition.  A low \(g^{-1}\)-coefficient multiplying a high
\(D\mathcal R\) is exactly controlled by \(\mathbf Z_j\).  A high
block of \(g^{-1}\) multiplying a lower curvature derivative is
controlled by \(\int_0^t\mathbf G_j(\theta)\,d\theta\), with one unused
frequency power.  In
\(\mathcal Q_0(g^{-1},\mathcal A)*\mathcal R\), a high
\(\mathcal A\) block and a low curvature block give
\(\int_0^t\mathbf A_j(\theta)\,d\theta\), whereas a high curvature block has
one unused power relative to \(\mathbf Z_j\).  Resonant terms split
between the same alternatives.  Together with
\eqref{eq:reconstruction-initial-curvature-connection}, this proves

\begin{equation}\label{eq:reconstruction-A-estimate}
 \mathbf A_j(t)
 \leq C\mathbf I_j+C\mathbf Z_j(t)
 +C\int_0^t
   \bigl(\mathbf A_j(\theta)+\mathbf G_j(\theta)\bigr)\,d\theta.
\end{equation}

Absorb the last term in
\eqref{eq:reconstruction-metric-from-connection} and substitute the
result into \eqref{eq:reconstruction-Z-estimate}.  Absorb the
\(\frac18\mathbf Z_j\) term there and then substitute the resulting
bound into \eqref{eq:reconstruction-A-estimate}.  Gronwall's
inequality on the fixed interval \([0,T]\) yields

\begin{equation}\label{eq:reconstruction-coupled-envelope}
 \mathbf A_j(T)+\mathbf Z_j(T)+\mathbf G_j(T)
 \leq C\mathbf I_j.
\end{equation}

The same argument without the high-frequency tail weights, together
with the finite low block from
\eqref{eq:reconstruction-buffered-lower}, gives the uniform
\(C_x^{m+2,\alpha}=C_x^{r,\alpha}\) bound for \(g\), and hence for
\(K_s=g-g_{s,0}\).  Since \(K_s(0)=0\) and

\[
 \partial_tK_s=-2\Ric_{g_{s,0}+K_s},
\]

the spatial estimate and the differentiated Ricci identity give, for
every \(q\) with \(2q\leq r\),

\[
 \partial_t^{\,q}K_s
 \in L^\infty_tC_x^{r-2q,\alpha}.
\]
For \(q\geq1\), its initial value is
\[
 \partial_t^{\,q}K_s(0)
 =\mathcal P_q\bigl(
  g_{s,0}^{-1},Dg_{s,0},\ldots,D^{2q}g_{s,0};
  \Rm(D),\ldots,D^{2q-2}\Rm(D)\bigr).
\]

Here \(\mathcal P_q\) denotes a universal smooth tensorial
contraction.  The case \(q=0\) is the already recorded identity
\(K_s(0)=0\).  One further distributional differentiation of the
Ricci identity, followed by the same Bony product estimates, gives
\[
 \partial_t^{\,q+1}K_s
 \in L^\infty_tC_x^{r-2q-2,\alpha},
\]
where the last space is interpreted in the Littlewood--Paley sense
when its exponent is negative.  Consequently, for every
spatial block and every \(2q+|\beta|=r\),
\[
 \bigl\|\Delta_\ell D^\beta\partial_t^{\,q}K_s(t)
       -\Delta_\ell D^\beta\partial_t^{\,q}K_s(t')\bigr\|_{L^\infty}
 \leq C\min\bigl\{2^{-\alpha\ell},
       |t-t'|2^{(2-\alpha)\ell}\bigr\}.
\]
Splitting at \(2^{2\ell}|t-t'|=1\) gives the temporal
\(\alpha/2\)-seminorm at every top parabolic derivative; the identical
argument below the top order gives all mixed seminorms.  This is the
standard dyadic parabolic interpolation characterization, and proves
that the smooth estimates are uniform in the full
\(C^{(r+\alpha)/2,r+\alpha}\) norm, including the initial face.  The
displayed initial identities are precisely the required compatibility
relations there; no lateral compatibility condition occurs because
\(\mathcal X\) is closed.

Approximate \(g_{s,0}\) in the frozen atlas and solve the corresponding
smooth Ricci--DeTurck systems.  The compact quasilinear estimate gives
convergence of \(H_s\) in its stated space, and the buffered
factorwise reconstruction gives convergence of \(K_s\) at order
\(r-2\).  A fixed datum, and likewise a norm-convergent family together
with its limit, has one common input envelope in
\eqref{eq:reconstruction-input-envelope}.  Estimate
\eqref{eq:reconstruction-coupled-envelope} makes the order-\(r\)
envelopes of all reconstructed approximants uniformly small as
\(j\to\infty\).  After \(j\) is fixed, their finitely many remaining
blocks converge by the order-\((r-2,\alpha)\) reconstruction estimate.
They are therefore Cauchy at order \(r\).  Applying the physical Ricci
identity, the time-jet estimates, and the preceding dyadic parabolic
interpolation estimate to pairwise differences gives convergence in
the full anisotropic space.  The limit is the already defined tensor
\(\chi_s^*(g_{s,0}+H_s)-g_{s,0}\).

It remains to translate the envelope estimate.  From
\eqref{eq:reconstruction-coupled-envelope}, metric compatibility, and
\(K_s=g-g_{s,0}\),

\[
 \mathcal W_j^r(K_s;T)
 \leq C\left(
  \mathcal W_{j-J_*}^r(g_{s,0})
 +\mathcal W_{j-J_*}^r(H_s;T)
 +M_{\rm rec}2^{-(1-\alpha)j}
 \right).
\]

Apply the second inequality in
\eqref{eq:tail-frequency-envelope-equivalence} to the two inputs and
the first inequality to \(K_s\).  The fixed atlas passage and the
single index shift are absorbed into one dilation
\(C_{\rm rec}\).  Enlarging \(C_K\) to absorb the common bound
\(M_{\rm rec}\) gives exactly
\eqref{eq:compact-reconstruction-tail}.
\par\hfill\(\square\)
\medskip

For a singleton family, the sublemma gives the full fixed-value
order-\((k+2,\alpha)\) anisotropic regularity and one uniform top
spatial little-H\"older tail for \(K_s\) on the whole compact time
slab.  For a norm-convergent family, the compact quasilinear estimate
first gives a common tail for \(H_s\), and
\eqref{eq:compact-reconstruction-tail} transfers it through the
physical reconstruction.  This is a fixed-value continuity statement:
it does not assert same-order Fr\'echet differentiability.  Parameter
differentiability remains the buffered \(k+2\to k\) assertion above.

For use only in the trial operator, extend this already reconstructed
normalized increment across its left endpoint.  Set
\[
 m=1+\left\lfloor\frac{k+2+\alpha}{2}\right\rfloor,
\]
set \(b_j:=j\) for \(1\leq j\leq m\), solve the Vandermonde system
\[
 \sum_{j=1}^m a_j(-b_j)^\ell=1,
 \qquad 0\leq\ell\leq m-1,
\]
and fix \(0<c_{\rm ext}\leq m^{-1}\).  Now set
\[
 \varepsilon_{\rm ext}
 :=c_{\rm ext}\widehat\delta_{\rm ph},
\]
and consider
a bounded linear zero-trace extension operator
\[
 \operatorname{Ext}_0:
 \left\{
  u\in\mathbb E_{\mathcal X,s}^{k+2,\alpha}
       (\widehat I_{\rm ph}):u(0)=0
 \right\}
 \longrightarrow
 \mathbb E_{\mathcal X,s}^{k+2,\alpha}
       ([-\varepsilon_{\rm ext},\widehat\delta_{\rm ph}])
\]
which agrees with \(u\) on \(\widehat I_{\rm ph}\).  Set
\[
 (\operatorname{Ext}_0u)(-\widehat\sigma)
 :=\sum_{j=1}^m a_j u(b_j\widehat\sigma)
 \quad(0\leq\widehat\sigma\leq\varepsilon_{\rm ext}),
 \qquad
 (\operatorname{Ext}_0u)(\widehat\sigma):=u(\widehat\sigma)
 \quad(\widehat\sigma\geq0).
\]
The moment identities match all temporal derivatives required by the
anisotropic exponent.  The formula is bounded in the anisotropic norm
with norm depending only on \(k,\alpha,c_{\rm ext}\), acts only in
time, and passage from smooth tensors to their closure shows that it
preserves both the little-H\"older class and its uniform top spatial
tail.  Put
\[
 K_s^{\rm tr}(\widehat\sigma)
 :=(\operatorname{Ext}_0K_s)(\widehat\sigma),
 \qquad
 -\varepsilon_{\rm ext}\leq\widehat\sigma
 \leq\widehat\delta_{\rm ph}.
\]
Choose \(\widehat\delta_{\rm ph}\) sufficiently small at the outset,
uniformly in the normalized prepared data.  Then
\(g_{s,0}+K_s^{\rm tr}\) remains positive definite.  On
\(\widehat I_{\rm ph}\) one has \(K_s^{\rm tr}=K_s\), so that
\[
 G(t(s))+\lambda_sK_s^{\rm tr}(\widehat\sigma)
 =G(t(s)+\lambda_s\widehat\sigma).
\]
The time-only extension is fixed and linear, while the forward
reconstruction is \(C^1\) with the stated two-order buffer.  Thus
\(K_s^{\rm tr}\) is a \(C^1\) coefficient of the initial
prepared state from input order \(k+2\) to output order \(k\), and at
the same time has the full fixed-value order-\((k+2,\alpha)\) spatial
class used in the invariant high tier.  Its
negative-\(\widehat\sigma\) part is auxiliary and is not asserted to
solve Ricci flow.

Put
\[
 R_\tau=\varphi_{-\tau}\circ\Theta_\tau,\qquad
 \ell(\tau)=\log\frac{\lambda(\tau)}{\lambda_0}.
\]
Proposition~\ref{prop:adaptive-target-tracking}, or a direct
differentiation, gives the exact relative system
\begin{equation}\label{eq:relative-coupled-system}
 \begin{aligned}
  \theta_\tau&=e^\ell,&
  \ell_\tau&=-(1+a),\\
  \partial_\tau R_\tau
  &=\left[(\varphi_{-\tau})_*
       (a\bar\nabla\bar f-U_b)\right]\circ R_\tau,&
  \partial_\tau F_\tau
  &=\lambda\Delta_{\acute G,S}F_\tau ,
 \end{aligned}
\end{equation}
with
\[
 t(\tau)=t(s)+\lambda_s\theta(\tau).
\]
Here \(c=(a,b)=\mathscr C(\tau,\mathbf z)\).
Equivalently, with
\(V_\Theta:=(1+a)\bar\nabla\bar f-U_b\), the two map equations used in
the same-output chart are
\begin{equation}\label{eq:same-output-coupled-map-system}
 \partial_\tau\Theta=V_\Theta\circ\Theta,\qquad
 \partial_\tau\Phi
 =V_\Theta\circ\Phi
  +\lambda\,d\Theta|_F
       \bigl(\Delta_{\acute G,S}F\bigr),
 \quad F=\Theta^{-1}\circ\Phi .
\end{equation}
The second identity is the chain rule for \(\Phi=\Theta\circ F\);
because \(\Theta:(M,S)\to(M,\lambda\bar g)\) is a target isometry up
to a constant factor, its second term is precisely the transformed
harmonic-map principal part.  Lemma~\ref{lem:feedback-functional-regularity}
shows that the ODE right-hand sides are \(C^1\), locally Lipschitz maps
in the prepared spaces.  In the same-output exponential chart at
\(\Phi_0\), the second equation in
\eqref{eq:same-output-coupled-map-system} is a quasilinear system whose
linearization satisfies
Lemma~\ref{lem:weighted-prepared-Schauder} in its Bochner mild
formulation; the scaling assertion is precisely
Remark~\ref{rem:annular-time-scaling}.  Normalized-time continuity of
the trial state makes every resulting forcing strongly measurable,
but no common modulus in all displayed source-adapted clocks is
asserted.  The clocks and atlas are frozen from the center state, so
the fixed-point map is compared in one Banach space; the trial factor
\(\lambda/\lambda_s^\circ\) is an ordinary \(C^1\) coefficient.

\emph{Step 2: the contraction.}
Write
\[
 \Theta_\tau=\operatorname{Exp}(X_\Theta(\tau))\circ\Theta_0,\qquad
 \Phi_\tau=\operatorname{Exp}(X_\Phi(\tau))\circ\Phi_0,
\]
with \(X_\Theta(\tau_0)=X_\Phi(\tau_0)=0\), and recover
\(R_\tau=\varphi_{-\tau}\circ\Theta_\tau\) and
\(F_\tau=\Theta_\tau^{-1}\circ\Phi_\tau\).  We use a
lower-triangular Picard map.  This is essential at the top spatial
order: the transformed harmonic-map operator contains a second-order
cross term from \(\Theta\) into \(\Phi\), and that term must act on the
already computed \(\Theta\)-output rather than on an arbitrary
\(\Theta\)-trial.

For a trial tuple
\[
 \mathbf w=(\theta,\ell,X_\Theta,X_\Phi),
\]
put
\[
 G_{\mathbf w}^{\rm tr}(\tau)
 :=G(t(s))+\lambda_sK_s^{\rm tr}(\theta(\tau)),
\]
form the preliminary trial graph state
\(\mathbf z^{\rm tr}(\mathbf w)\), and set
\[
 c^{\rm tr}(\mathbf w)
 =(a^{\rm tr}(\mathbf w),b^{\rm tr}(\mathbf w))
 :=\mathscr C\bigl(\tau,\mathbf z^{\rm tr}(\mathbf w)\bigr).
\]
Let \(\mathscr V_\Theta^{\rm tr}(\tau,\mathbf w)\) denote the
same-output coordinate form of \(V_\Theta\circ\Theta\), evaluated with
this preliminary trial state and feedback.  First define the three
ODE outputs by
\begin{equation}\label{eq:coupled-fixed-point-ode-stage}
 \begin{aligned}
  \widehat\theta(\tau)
   &=\int_{\tau_0}^{\tau}e^{\ell(q)}\,dq,\\
  \widehat\ell(\tau)
   &=-\int_{\tau_0}^{\tau}
       \bigl(1+a^{\rm tr}(\mathbf w)(q)\bigr)\,dq,\\
  \widehat X_\Theta(\tau)
   &=\int_{\tau_0}^{\tau}
       \mathscr V_\Theta^{\rm tr}(q,\mathbf w(q))\,dq .
 \end{aligned}
\end{equation}
The compact affine increment \(H_s\) remains the fixed Step~1
solution.

Retain the trial \(X_\Phi\), but replace the other three entries by
the outputs just obtained:
\begin{equation}\label{eq:coupled-sharp-trial}
 \mathbf w^\sharp(\mathbf w)
 :=
 \bigl(\widehat\theta,\widehat\ell,
       \widehat X_\Theta,X_\Phi\bigr).
\end{equation}
Equivalently, the geometric quantities in the sharp trial state are
\[
 \begin{gathered}
  G^\sharp
   =G(t(s))+\lambda_sK_s^{\rm tr}(\widehat\theta),
  \qquad
  \lambda^\sharp=\lambda_s e^{\widehat\ell},\\
  \Theta^\sharp
   =\operatorname{Exp}(\widehat X_\Theta)\circ\Theta_0,
  \qquad
  \Phi^{\rm tr}
   =\operatorname{Exp}(X_\Phi)\circ\Phi_0,\\
  R^\sharp=\varphi_{-\tau}\circ\Theta^\sharp,
  \qquad
  F^\sharp=(\Theta^\sharp)^{-1}\circ\Phi^{\rm tr},\\
  S^\sharp=\lambda^\sharp(\Theta^\sharp)^*\bar g,
  \qquad
  \acute G^\sharp
   =\eta\,\iota_*G^\sharp+(1-\eta)S^\sharp .
 \end{gathered}
\]
Let
\[
 \mathbf z^\sharp(\mathbf w)
 :=(G^\sharp,\lambda^\sharp,R^\sharp,F^\sharp),
 \qquad
 c^\sharp(\mathbf w)
 :=\mathscr C\bigl(\tau,\mathbf z^\sharp(\mathbf w)\bigr).
\]
Thus every occurrence of \(\Theta\) in the coefficients or forcing of
the \(\Phi\)-equation is evaluated coherently at
\(\Theta^\sharp\); in particular, \(S,\acute G,F\), the target
connection, and every second-order cross term use the same
\(\Theta^\sharp\).  Only \(X_\Phi\) remains a trial variable.

Freeze the \(X_\Phi\)-principal coefficients at
\(\mathbf z^\sharp(\mathbf w)\), keep \(\widehat X_\Phi\) as the
unknown, and denote the resulting linear operator by
\(\mathscr L_{\mathbf w^\sharp}\).  Let
\(\mathscr R_\Phi(\tau,\mathbf w^\sharp)\) be the remaining coordinate
form of the second equation in
\eqref{eq:same-output-coupled-map-system}, evaluated with
\(c^\sharp(\mathbf w)\).  Complete the definition by
\begin{equation}\label{eq:coupled-fixed-point}
 \mathscr L_{\mathbf w^\sharp}\widehat X_\Phi
   =\mathscr R_\Phi(\tau,\mathbf w^\sharp),
 \qquad
 \widehat X_\Phi(\tau_0)=0,
\end{equation}
and set
\[
 \mathfrak K(\mathbf w)
 :=
 \bigl(\widehat\theta,\widehat\ell,
       \widehat X_\Theta,\widehat X_\Phi\bigr).
\]
Throughout Step~2, an \(\mathbb S\)-norm of a four-tuple
\(\mathbf v\) abbreviates the corresponding norm of
\((H_s,\mathbf v)\).

The construction has exactly the desired fixed points.  Indeed, if
\(\mathbf w=\mathfrak K(\mathbf w)\), then
\[
 \mathbf w^\sharp(\mathbf w)=\mathbf w,\qquad
 \mathbf z^\sharp(\mathbf w)=\mathbf z^{\rm tr}(\mathbf w),\qquad
 c^\sharp(\mathbf w)=c^{\rm tr}(\mathbf w).
\]
Consequently the fixed-point equations are precisely
\eqref{eq:relative-coupled-system} and
\eqref{eq:same-output-coupled-map-system}.  Conversely, every solution
of those equations is a fixed point of this triangular map.

Put \(\beta:=\alpha/4\), let \(\pi_{\rm ode}\) denote projection onto
the first three tuple components, and let \(d_{\rm ode}\) be the sum
of the corresponding component seminorms occurring in
\(d_{\mathbb D_{\rm mix}^{k,\alpha}(I)}\).  Direct integration in
\eqref{eq:coupled-fixed-point-ode-stage}, together with
Lemma~\ref{lem:feedback-functional-regularity}, gives
\begin{equation}\label{eq:coupled-sharp-trial-difference}
 \begin{split}
 d_{\rm ode}\bigl(
   \pi_{\rm ode}\mathbf w_1^\sharp,
   \pi_{\rm ode}\mathbf w_2^\sharp\bigr)
 &\leq C_K\delta^\beta
   d_{\mathbb D_{\rm mix}^{k,\alpha}(I)}
      (\mathbf w_1,\mathbf w_2),\\
 d_{\mathbb D_{\rm mix}^{k,\alpha}(I)}
   \bigl(\mathbf w_1^\sharp,\mathbf w_2^\sharp\bigr)
 &\leq C_K
   d_{\mathbb D_{\rm mix}^{k,\alpha}(I)}
      (\mathbf w_1,\mathbf w_2).
 \end{split}
\end{equation}
If \(Y_i=\widehat X_{\Phi,i}\) are the corresponding parabolic
outputs, subtraction gives
\[
 \mathscr L_{\mathbf w_1^\sharp}(Y_1-Y_2)
 =
 \bigl(\mathscr L_{\mathbf w_2^\sharp}
       -\mathscr L_{\mathbf w_1^\sharp}\bigr)Y_2
 +\mathscr R_\Phi(\mathbf w_1^\sharp)
  -\mathscr R_\Phi(\mathbf w_2^\sharp).
\]
The prepared-chart calculus,
Lemma~\ref{lem:feedback-functional-regularity}, and
\eqref{eq:coupled-sharp-trial-difference} imply
\begin{equation}\label{eq:coupled-map-difference-forcing}
 \left\|
  \bigl(\mathscr L_{\mathbf w_2^\sharp}
       -\mathscr L_{\mathbf w_1^\sharp}\bigr)Y_2
  +\mathscr R_\Phi(\mathbf w_1^\sharp)
       -\mathscr R_\Phi(\mathbf w_2^\sharp)
 \right\|_{\mathbb F_{\rm sc}^{k,\alpha}}
 \leq C_K
 d_{\mathbb D_{\rm mix}^{k,\alpha}(I)}
    (\mathbf w_1,\mathbf w_2).
\end{equation}
The difference of principal matrices multiplies the uniformly bounded
high derivatives of \(Y_2\).  It is therefore estimated in
\(\mathbb F_{\rm sc}^{k,\alpha}\), exactly as before.  The
\(\Theta\)-part of that coefficient difference is smaller by the ODE
factor in \eqref{eq:coupled-sharp-trial-difference}; the retained
\(X_\Phi\)-part is handled by the one-order zero-trace estimate after
the Schauder solve.  All feedback terms satisfy the same estimate by
\eqref{eq:feedback-local-Lipschitz}.

The map
\[
 \mathbf w\longmapsto\mathbf w^\sharp(\mathbf w)
\]
is \(C^1\) in the same buffered prepared topologies: it is the product
of the three \(C^1\) ODE solution maps and the identity on \(X_\Phi\).
The parameter-dependent Schauder theorem therefore makes
\(\mathfrak K\) \(C^1\) with the original two-order buffer.  Along a
curve of fixed points, differentiation of
\(\mathbf w^\sharp(\mathbf w)=\mathbf w\) gives
\(\delta\mathbf w^\sharp=\delta\mathbf w\); substituting this identity
in the differentiated triangular equations gives exactly the original
coupled linearized system.

Let
\[
 \mathbf w_0(\tau)\equiv(0,0,0,0)
\]
be the constant extension of the initial coordinate state
\((\theta,\ell,X_\Theta,X_\Phi)=(0,0,0,0)\).  Fix provisionally
\[
 C_{\rm clk}:=
  \sup_{\mathbf z_s\in\mathscr B_s}
       \frac{\lambda^\circ(s)}{\lambda_s}<\infty,
 \qquad
 0<\rho_{\rm fp}
 <\frac1{2C_{\rm clk}}
    \min\{\varepsilon_{\rm ext},
             \widehat\delta_{\rm ph}\}.
\]
Consider the two-tier ball on which the
\(\mathbb S^{k+2,\alpha}\) norm is at most \(K\) and the
\(\mathbb D_{\rm mix}^{k,\alpha}\) displacement from the constant initial
extension is at most \(\rho_{\rm fp}\).  Indeed,
\(|\theta|\leq C_{\rm clk}\rho_{\rm fp}\).  Thus every trial clock lies in
the domain of \(K_s^{\rm tr}\), so the trial operator is
defined on an open Banach neighborhood of the two-tier ball.  The tame
product estimates (one factor is
always controlled in the low \(C^2\) box), the high-order Schauder
estimate, and the ODE estimates give
\begin{equation}\label{eq:coupled-high-self-map}
 \|\mathfrak K(\mathbf w)\|_{\mathbb S^{k+2,\alpha}}
 \leq
 C_0+C\bigl(\rho_{\rm fp}+\delta^\beta\bigr)K
    +C_K\delta^{\alpha/2}.
\end{equation}
where
\[
 C_0=C_0(k,\alpha,N,\mathfrak P_{\rm prep},
              K_{\rm init}^{k+2,\alpha}).
\]
Thus \(C_0\), and hence the later choice \(K=2C_0+1\), depend only on
the recorded uniformity data.  The invariant high self-map uses
\(q=k+3\).  For differences,
apply Lemma~\ref{lem:weighted-prepared-Schauder} with \(q=k+2\) to
\eqref{eq:coupled-map-difference-forcing}, and then use the
one-order zero-trace estimate
\eqref{eq:weighted-zero-trace-one-order}.  Direct integration of the
ODE variables gives the same small factor, including in the extra
\(\Theta\)-norm.  Therefore
\begin{align}
 d_{\mathbb D_{\rm mix}^{k,\alpha}(I)}
       (\mathfrak K(\mathbf w),\mathbf w_0)
 &\leq C_K\delta^{\alpha/4},
 \label{eq:coupled-fixed-point-self-map}\\
 d_{\mathbb D_{\rm mix}^{k,\alpha}(I)}
       (\mathfrak K(\mathbf w_1),\mathfrak K(\mathbf w_2))
 &\leq
 C_K\delta^{\alpha/4}
 d_{\mathbb D_{\rm mix}^{k,\alpha}(I)}
   (\mathbf w_1,\mathbf w_2).
\label{eq:coupled-fixed-point-contraction}
\end{align}
Choose \(K=2C_0+1\).  First decrease the provisional fixed-point
displacement radius so that \(C\rho_{\rm fp}<1/8\), and then decrease
\(\delta_\tau\) so that
\[
 C\delta_\tau^\beta<\frac18,\qquad
 C_K\delta_\tau^\beta<\frac12,
\]
\eqref{eq:coupled-fixed-point-self-map} is at most
\(\rho_{\rm fp}\), and the last term in
\eqref{eq:coupled-high-self-map} is at most \(K/4\).
After this choice every sharp trial
\(\mathbf w^\sharp(\mathbf w)\) remains in the same prepared chart and
coefficient package.  Equations
\eqref{eq:coupled-fixed-point-self-map} and
\eqref{eq:coupled-fixed-point-contraction} follow with the same
contraction metric: the retained \(X_\Phi\)-difference is treated by
the one-order zero-trace Schauder estimate, while every difference in
the first three sharp entries already carries the ODE factor from
\eqref{eq:coupled-sharp-trial-difference}.
Decrease \(\delta_\tau\) once more, using only the same recorded
uniformity data, so that
\[
 0\leq\widehat\theta(\tau)
 =\int_s^\tau e^{\ell(q)}\,dq
 \leq\frac12\widehat\delta_{\rm ph}
 \qquad(s\leq\tau\leq s+\delta_\tau).
\]
The two-tier ball is invariant and
\eqref{eq:coupled-fixed-point-contraction} is a contraction.  Its
iterates are taken in the following explicitly specified low-topology
completion in the fixed common-margin chart.  For a trial
\(\mathbf w\), let
\(\Omega_{\mathrm{hf}}^{k+2,\alpha}(\mathbf w;\varrho)\) be the supremum,
after every scale, tensor, and polynomial-weight normalization
prescribed in the component norms, of the truncated
\(\alpha\)-H\"older seminorms at spatial pair-distance at most
\(\varrho\) of the top spatial derivatives in every compact-core and
rescaled dyadic-chart representative occurring in
\(\mathbb S^{k+2,\alpha}(I)\), followed by the supremum over the
relevant time variable and over all charts.  For \(\varrho\geq1\),
extend this quantity constantly by its value at \(1\).

We now separate the numerical lifespan from the individual
little-H\"older tail.  After \(K\) has been fixed, decrease
\(\rho_{\rm fp}\) and then choose \(\delta_{\rm tail}>0\), depending
only on
\[
 k,\alpha,N,\quad \mathfrak P_{\rm prep},\quad
 K_{\rm init}^{k+2,\alpha},\quad\text{and the fixed common margins},
\]
so that there are constants
\[
 C_{\rm hf}\geq1,\qquad 0<\vartheta_{\rm hf}<1,
\]
with the following property.  For each fixed initial state
\(\mathbf z_0\), there is a bounded nondecreasing function
\[
 \nu_{\mathbf z_0}:[0,\infty)\longrightarrow[0,\infty),
 \qquad
 \lim_{\varrho\downarrow0}\nu_{\mathbf z_0}(\varrho)=0,
\]
extended constantly for arguments at least \(1\), such that every
trial tuple in the two-tier ball and every
\(0<\delta\leq\delta_{\rm tail}\) satisfy
\begin{equation}\label{eq:coupled-tail-self-map}
 \Omega_{\mathrm{hf}}^{k+2,\alpha}
   \bigl(\mathfrak K_{\mathbf z_0}(\mathbf w);\varrho\bigr)
 \leq
 \nu_{\mathbf z_0}(C_{\rm hf}\varrho)
 +\vartheta_{\rm hf}
  \Omega_{\mathrm{hf}}^{k+2,\alpha}
    (\mathbf w;C_{\rm hf}\varrho).
\end{equation}
The constants \(C_{\rm hf}\), \(\vartheta_{\rm hf}\), and
\(\delta_{\rm tail}\) are independent of the particular function
\(\nu_{\mathbf z_0}\).

Here is the required top-frequency estimate.  Let
\(\omega_{m,\alpha}(u;\varrho)\) denote the truncated
\(\alpha\)-H\"older seminorm of \(D^m u\) at pair-distance at most
\(\varrho\).  In every normalized chart, for
\(0<\varrho\leq C_{\rm hf}^{-1}\), Leibniz and Fa\`a di Bruno give
\[
 \begin{aligned}
  \omega_{m,\alpha}(uv;\varrho)
  &\leq C_K\Bigl(
    \omega_{m,\alpha}(u;\varrho)
    +\omega_{m,\alpha}(v;\varrho)
    +\varrho^{1-\alpha}\Bigr),\\
  \omega_{m,\alpha}(u\circ\Psi;\varrho)
  &\leq C_K\Bigl(
    \omega_{m,\alpha}(u;C_{\rm hf}\varrho)
    +\omega_{m,\alpha}(\Psi;\varrho)
    +\varrho^{1-\alpha}\Bigr).
 \end{aligned}
\]
The last terms contain precisely the lower-derivative products: those
factors have one unused derivative and hence are uniformly Lipschitz
under the fixed \(K\)-bound.  Enlarge
\(\nu_{\mathbf z_0}\) by
\(C_K\min\{\varrho,1\}^{1-\alpha}\) to absorb these terms.

For \(C_{\rm hf}^{-1}\leq\varrho\leq1\), enlarge
\(\nu_{\mathbf z_0}\) on that compact scale range by the uniform
high-norm bound for
\(\mathfrak K_{\mathbf z_0}(\mathbf w)\), and then extend it constantly
for arguments at least \(1\).  This changes neither boundedness nor
vanishing at zero and makes
\eqref{eq:coupled-tail-self-map} valid for every
\(0<\varrho\leq1\).

Apply these estimates and
\eqref{eq:quantitative-tail-Schauder} to the triangular construction.
We record the top-order routing explicitly.  No small factor will be
attributed to a same-order zero-trace parabolic inverse.

First consider the ODE stage.  Once the scalar feedback is fixed, the
coordinate equation for \(X_\Theta\) is pointwise in space: it contains
no spatial derivative of \(X_\Theta\).  Differentiating its coordinate
right-hand side to the top spatial order can therefore place the top
derivative only on a trial \(X_\Theta\), with a uniformly bounded
coefficient, or on a fixed background coefficient.  The former term
is integrated over normalized time, while the latter is included in
the datum modulus.  The feedback coefficients are scalar functions of
time; multiplying the fixed profiles
\(\bar\nabla\bar f,\chi_\tau W_1,\ldots,\chi_\tau W_8\) creates only a
fixed spatial tail.  Hence, after enlarging
\(\nu_{\mathbf z_0}\),
\begin{equation}\label{eq:coupled-triangular-ode-tail}
 \Omega_{\rm hf}^{k+2,\alpha}
 \bigl((\widehat\theta,\widehat\ell,\widehat X_\Theta);
       \varrho\bigr)
 \leq
 \nu_{\mathbf z_0}(C_{\rm hf}\varrho)
 +C_K\delta\,
  \Omega_{\rm hf}^{k+2,\alpha}
     (\mathbf w;C_{\rm hf}\varrho).
\end{equation}
The scalar components have no spatial tail; they are included in the
left side only to emphasize that the entire first stage has already
been computed.

The geometric reason for the triangular ordering is the exact
naturality identity
\begin{equation}\label{eq:triangular-transformed-tension}
 d\Theta^\sharp\!
 \left[
  \Delta_{\acute G^\sharp,S^\sharp}
   \bigl((\Theta^\sharp)^{-1}\circ\Psi\bigr)
 \right]
 =
 \Delta_{\acute G^\sharp,\bar g}\Psi .
\end{equation}
Here the constant rescaling from \(\bar g\) to
\(\lambda^\sharp\bar g\) does not change the target connection.
Thus the principal \(\Phi\)-operator depends on \(\Theta^\sharp\)
through \(\acute G^\sharp\).  On the pure exterior,
\(\acute G^\sharp=\lambda^\sharp(\Theta^\sharp)^*\bar g\);
linearization at the identity contains the cross-principal term
\(-\Delta X_\Theta^\sharp\).  Formula
\eqref{eq:coupled-triangular-ode-tail} is precisely what gives this
term a small tail.  It would have no small same-order factor if
\(X_\Theta^\sharp\) were replaced by an arbitrary trial.

Put \(q:=k+3\), the spatial order of each map component in the high
tier, and write \(Z:=X_\Phi\) for the sole retained map trial.  In one
normalized chart, after the three sharp entries have been fixed, let
\(A_\Phi[Z]\) denote the principal matrix of the coordinate equation
for \(Y=\widehat X_\Phi\), and put
\(A_\Phi^{\rm ret}[Z]:=A_\Phi[Z]-A_\Phi[0]\).  The coordinate
linearization used to define
\(\mathscr L_{\mathbf w^\sharp}\) shows that this principal matrix
depends smoothly and at differential order zero on \(Z\); occurrences
of \(DZ\) belong to the lower-order coefficient tensors.  This is the
usual quasilinear, rather than fully nonlinear, structure of the
harmonic-map system in a fixed exponential chart.

The complete top-order routing for the parabolic stage is as follows.
\begin{enumerate}
\item The initial state, compact coefficient path, frozen atlas,
      cutoffs, transitions, and the value of the triangular map at the
      constant trial are fixed.  Their tails are included in
      \(\nu_{\mathbf z_0}\).

\item Every second-order dependence on \(\Theta\), including that
      carried through \(S,\acute G,F\), and the transformed target
      connection, is evaluated at \(\Theta^\sharp\).  A same-order
      parabolic estimate may transmit this tail with a bounded
      coefficient, but
      \eqref{eq:coupled-triangular-ode-tail} has already supplied the
      factor \(C_K\delta\).

\item Every term containing two spatial derivatives of the unknown
      \(\widehat X_\Phi\) is contained in
      \(\mathscr L_{\mathbf w^\sharp}\); none is treated as a
      prescribed source.  More precisely, Bony decomposition of the
      part depending on the retained trial gives
      \[
       A_\Phi^{\rm ret}[Z]D^2Y
       =\mathsf T_{A_\Phi^{\rm ret}[Z]}D^2Y
        +\mathsf T_{D^2Y}A_\Phi^{\rm ret}[Z]
        +\mathsf R(A_\Phi^{\rm ret}[Z],D^2Y).
      \]
      The low-coefficient--high-output term is retained in the
      principal operator.  The tame composition estimate and the
      two-sided frequency-envelope product estimate give, with one
      fixed integer \(J_\Phi\),
      \begin{equation}\label{eq:triangular-retained-Phi-block}
       \begin{split}
       &\mathcal W_j^{q-2}\!\left(
          \mathsf T_{D^2Y}A_\Phi^{\rm ret}[Z]
          +\mathsf R(A_\Phi^{\rm ret}[Z],D^2Y)
        \right)\\
       &\qquad\leq
        C_K\|Y\|_{L^\infty
             \mathfrak X_{\rm sc}^{q-1,\alpha}}
             \mathcal W_{j-J_\Phi}^{q}(Z)
        +C_K2^{-(1-\alpha)j}.
       \end{split}
      \end{equation}
      Indeed, in the first displayed paraproduct the high block lies
      on the order-zero composition \(A_\Phi^{\rm ret}[Z]\), while the
      low factor contains only \(D^2Y\); in the resonant term the two
      comparable high blocks leave at least one unused derivative.
      The proof of
      \eqref{eq:weighted-zero-trace-one-order} is an interpolation
      estimate for every zero-left-trace element of
      \(\mathbb E_{\rm sc}^{q,\alpha}(I)\), not an additional
      equation for that element.  Since \(Y(\tau_0)=0\) and the
      invariant high-tier bound is \(K\), it therefore yields
      \[
       \|Y\|_{L^\infty
          \mathfrak X_{\rm sc}^{q-1,\alpha}}
       \leq C_K\delta^\beta.
      \]
      Thus \eqref{eq:triangular-retained-Phi-block} transmits a top
      \(Z\)-block only with the factor \(C_K\delta^\beta\).  Spatial
      freezing of the first paraproduct contributes instead
      \(\varepsilon_\Phi\mathcal W_j^q(Y)\), where
      \(\varepsilon_\Phi\) is made uniformly small by choosing the
      normalized freezing radius after the coefficient bounds have
      been fixed.

\item Let \(\mathscr R_{\Phi,{\rm ret}}^{\rm lo}[Z]\) be the sum of
      the remaining lower-order and coordinate terms after their value
      at \(Z=0\) has been subtracted.  In a top block, either the high
      derivative lies on \(Z\) and its companion is an
      undifferentiated centered factor, or at most \(q-1\) derivatives
      land on \(Z\).  The first case gives \(C_K\rho_{\rm fp}\); the
      second has one unused derivative.  Terms in which a lower
      derivative of the zero-trace output is the companion are covered
      by \eqref{eq:weighted-zero-trace-one-order}.  Consequently,
      after enlarging the fixed datum modulus,
      \begin{equation}\label{eq:triangular-retained-Phi-source}
       \mathcal W_j^{q-2}
          (\mathscr R_{\Phi,{\rm ret}}^{\rm lo}[Z])
       \leq
       C_K(\rho_{\rm fp}+\delta^\beta)
          \mathcal W_{j-J_\Phi}^{q}(Z)
       +\nu_{\mathbf z_0}(C2^{-j}).
      \end{equation}
      The finite-dimensional Gram inverse introduces no spatial tail.
      Its scalar output multiplies fixed spatial profiles and
      subsequently passes through the ODE stage.

\item After the four evolution outputs have been obtained, the maps
      \(F\), the inverse maps, the prepared graph, the graft, the
      moving support, and the effective columns are formed only by the
      finite same-output inverse, composition, pullback, product, and
      localization operations of
      Lemma~\ref{lem:prepared-chart-calculus}.  Their fixed
      linearizations transmit the tails of outputs already obtained;
      their centered remainders carry
      \(C_K\rho_{\rm fp}\), and their lower
      Fa\`a di Bruno terms are absorbed into the datum modulus.
\end{enumerate}

Let
\[
 Y_{\mathbf w}:=\widehat X_\Phi(\mathbf w).
\]
For a tuple \(\mathbf v\), write
\(\mathcal W_j^{\rm mix}(\mathbf v)\) for the maximum of the
normalized two-sided envelopes of its components at their respective
top spatial orders; in particular, both map components are measured
at order \(q\).  Apply the dyadic proof of the quantitative Schauder
tail estimate, rather than only its coarser final display, to
\[
 \mathscr L_{\mathbf w^\sharp}Y_{\mathbf w}
 =\mathscr R_\Phi(\mathbf w^\sharp),
 \qquad Y_{\mathbf w}(\tau_0)=0.
\]
The low-principal--high-output interaction remains in the frozen
operator.  Item~\textup{(2)} and
\eqref{eq:coupled-triangular-ode-tail} give
\(C_K\delta\mathcal W_{j-J_\Phi}^{\rm mix}(\mathbf w)\);
\eqref{eq:triangular-retained-Phi-block} and
\eqref{eq:triangular-retained-Phi-source} give the retained-\(Z\)
contribution; centered nonlinearities give \(C_K\rho_{\rm fp}\); and
all fixed or one-unused-derivative blocks belong to
\(\nu_{\mathbf z_0}\).  Hence, after increasing \(J_\Phi\) once,
\begin{equation}\label{eq:triangular-Phi-dyadic}
 \begin{split}
 \mathcal W_j^q(Y_{\mathbf w})
 \leq{}&\nu_{\mathbf z_0}(C2^{-j})
 +C_K(\rho_{\rm fp}+\delta^\beta)
       \mathcal W_{j-J_\Phi}^{\rm mix}(\mathbf w)\\
 &+\varepsilon_\Phi\mathcal W_j^q(Y_{\mathbf w}),
 \end{split}
\end{equation}
where \(0<\delta\leq1\) was used to replace \(\delta\) by
\(\delta^\beta\).  Choose the common freezing radius so that
\(\varepsilon_\Phi\leq\frac12\), and absorb the final term.  Taking
the normalized-atlas supremum and using
\eqref{eq:tail-frequency-envelope-equivalence} in both directions,
with its one fixed index/scale dilation, gives
\begin{equation}\label{eq:triangular-Phi-tail}
 \Omega_{\rm hf}^{k+2,\alpha}
   (Y_{\mathbf w};\varrho)
 \leq
 \nu_{\mathbf z_0}(C_{\rm hf}\varrho)
 +C_K\bigl(\rho_{\rm fp}+\delta^\beta\bigr)
  \Omega_{\rm hf}^{k+2,\alpha}
    (\mathbf w;C_{\rm hf}\varrho).
\end{equation}
In particular, the factor \(\delta^\beta\) is used only in
\eqref{eq:triangular-retained-Phi-block}, where a top coefficient
block is paired with a lower derivative of a zero-trace output, and in
the already triangular ODE stage.  It is not a same-order zero-trace
smoothing assertion.

Combining
\eqref{eq:coupled-triangular-ode-tail},
\eqref{eq:triangular-Phi-tail}, and item~\textup{(5)}, and using
\(\delta\leq\delta^\beta\) for \(0<\delta\leq1\), gives
\[
 \Omega_{\mathrm{hf}}^{k+2,\alpha}
   \bigl(\mathfrak K_{\mathbf z_0}(\mathbf w);\varrho\bigr)
 \leq
 \nu_{\mathbf z_0}(C_{\rm hf}\varrho)
 +C_*\bigl(\rho_{\rm fp}+\delta^\beta\bigr)
  \Omega_{\mathrm{hf}}^{k+2,\alpha}
    (\mathbf w;C_{\rm hf}\varrho),
\]
where \(C_*\) is the sum of the finitely many transmission constants
in the preceding routing.  Choose
\[
 C_*\rho_{\rm fp}\leq\frac14,
 \qquad
 C_*\delta_{\rm tail}^{\beta}\leq\frac14,
 \qquad
 \vartheta_{\rm hf}:=
 C_*\bigl(\rho_{\rm fp}+\delta_{\rm tail}^{\beta}\bigr)<1.
\]
This proves \eqref{eq:coupled-tail-self-map}.  Both choices depend only
on the recorded numerical package and not on the rate at which an
individual little-H\"older tail vanishes.

The construction is quantitative in the fixed inputs:
\(\nu_{\mathbf z_0}\) is a fixed \(C_K\)-multiple of the sum of the
top-order tails of the frozen initial coefficients, the fixed compact
coefficient path, \(\mathfrak K_{\mathbf z_0}(\mathbf w_0)\), and
\(\varrho^{1-\alpha}\).  The compact quasilinear parametrix and
Lemma~\ref{lem:weighted-prepared-Schauder} preserve a common input
tail with constants depending only on the recorded uniformity data;
\eqref{eq:compact-reconstruction-tail} transfers the compact tail from
the affine DeTurck increment to the reconstructed coefficient used in
the graph.
Consequently, if a family of initial prepared states carries one
common datum modulus, the functions \(\nu_{\mathbf z_0}\) are
dominated by one common bounded nondecreasing function \(\nu\)
vanishing at zero.

Enlarge \(\nu_{\mathbf z_0}\), if necessary, so that it also dominates
the tail of the fixed element
\((H_s,\mathbf w_0)\).  For a norm-convergent family, use
the common function \(\nu\) supplied by the preceding paragraph.
Define
\begin{equation}\label{eq:invariant-tail-envelope}
 \mu_{\mathbf z_0}(\varrho)
 :=
 \sum_{j=0}^{\infty}
  \vartheta_{\rm hf}^{\,j}
  \nu_{\mathbf z_0}
   \bigl(C_{\rm hf}^{\,j+1}\varrho\bigr),
 \qquad \varrho\geq0 .
\end{equation}
This function is finite, nondecreasing, and vanishes at zero.  Indeed,
bounded convergence applies to the summable geometric series.
Moreover,
\[
 \nu_{\mathbf z_0}(C_{\rm hf}\varrho)
 +\vartheta_{\rm hf}
  \mu_{\mathbf z_0}(C_{\rm hf}\varrho)
 =\mu_{\mathbf z_0}(\varrho).
\]
Consequently \eqref{eq:coupled-tail-self-map} makes the class with
\(\mu=\mu_{\mathbf z_0}\) invariant.  Put
\[
 \begin{split}
 \mathbb A_{K,\rho_{\rm fp},\mu}:=\bigl\{
   \mathbf w=(\theta,\ell,X_\Theta,X_\Phi):\;&
   (H_s,\mathbf w)
      \in\mathbb S_0^{k+2,\alpha}(I),\\
  &\|(H_s,\mathbf w)\|_
       {\mathbb S^{k+2,\alpha}(I)}\leq K,\quad
   d_{\mathbb D_{\rm mix}^{k,\alpha}(I)}
     (\mathbf w,\mathbf w_0)\leq\rho_{\rm fp},\\
  &\Omega_{\mathrm{hf}}^{k+2,\alpha}
     (\mathbf w;\varrho)\leq\mu(\varrho)
     \quad(0<\varrho\leq1)\bigr\},\\
 \mathbb X_{K,\rho_{\rm fp},\mu}:={}&
  \operatorname{Comp}_{d_{\mathbb D_{\rm mix}^{k,\alpha}(I)}}
   \bigl(\mathbb A_{K,\rho_{\rm fp},\mu}\bigr) .
 \end{split}
\]
Thus the common-tail condition is purely a common small-scale spatial
modulus.  It contains no time derivative, no common time-H\"older
modulus in the source-adapted clocks, and no annular-tightness condition
\eqref{eq:weighted-little-annular-c0}.  The choices of
\(K,\rho_{\rm fp}\), and the final lifespan
\[
 \delta_\tau\leq\min\{\delta_{\rm tail},\delta_0\}
\]
are made before, and independently of, the individual function
\(\mu_{\mathbf z_0}\).  The latter is used only to identify the
low-topology completion with the intended little-H\"older solution
class.  If a family of initial states has one common datum modulus,
\eqref{eq:invariant-tail-envelope} supplies one common invariant
solution modulus.  The space
\(\mathbb X_{K,\rho_{\rm fp},\mu}\) is complete by construction.  To
identify its points with genuine prepared tuples,
take a Cauchy sequence in \(\mathbb A_{K,\rho_{\rm fp},\mu}\).
It converges in the global lower-order atlas-supremum norm, which
already excludes annular escape.  Split each normalized chart function
into low and high spatial frequencies.  Global lower-order convergence
controls the finitely many low-frequency blocks uniformly over the
entire atlas, while the common modulus \(\mu\) controls the
high-frequency remainder.  Hence the spatial top derivatives converge
in their \(C^\alpha\) norms in the same global atlas-supremum topology,
and the limit retains \(\mu\).  The common-margin graph, initial-trace,
and map constraints pass to the limit by the prepared chart calculus.

Estimate \eqref{eq:coupled-fixed-point-contraction} extends
\(\mathfrak K\) uniquely and contractively to
\(\mathbb X_{K,\rho_{\rm fp},\mu}\), and the self-map estimate keeps
its image in
that completion.  Banach's theorem therefore gives a fixed point
\(\mathbf w\) there.  This use of the completion does not discard the
temporal regularity.  Choose Picard iterates
\(\mathbf w_{n+1}=\mathfrak K(\mathbf w_n)\).  Every
\(\mathbf w_{n+1}\) lies in the full
\(\mathbb S^{k+2,\alpha}\) ball and obeys its uniform
time-H\"older and weak-time-derivative bounds.  The strong spatial
convergence just proved is uniform in normalized time and lets the
time-H\"older inequalities pass to the limit.

We recover the time derivatives strongly from the equations, so no
weak compactness in a H\"older space is used.  For the scalar and
\(X_\Theta\) components, local Lipschitz continuity of the feedback and
the prepared ODE right-hand sides gives strong uniform convergence of
\[
 \partial_\tau\theta_{n+1},\qquad
 \partial_\tau\ell_{n+1},\qquad
 \partial_\tau X_{\Theta,n+1}
\]
in their defining lower-order spaces.  For each \(n\), put
\(\mathbf w_n^\sharp:=\mathbf w^\sharp(\mathbf w_n)\).  On every
source-adapted chart,
the scaled parabolic identity for the next output has the form
\[
 \begin{aligned}
 \partial_{\vartheta_{\mathcal U}}
  \mathscr S_{\mathcal U}X_{\Phi,n+1}
 ={}&
 A_{\mathcal U}(\mathbf w_n^\sharp)*
  \nabla^2\mathscr S_{\mathcal U}X_{\Phi,n+1}\\
 &+B_{\mathcal U}(\mathbf w_n^\sharp)*
  \nabla\mathscr S_{\mathcal U}X_{\Phi,n+1}\\
 &+C_{\mathcal U}(\mathbf w_n^\sharp)*
  \mathscr S_{\mathcal U}X_{\Phi,n+1}
 +\mathscr R_{\mathcal U}(\mathbf w_n^\sharp).
 \end{aligned}
\]
 The common spatial modulus gives strong convergence of
 \(X_{\Phi,n+1}\) through its invariant high-tier order
 \(C^{k+3,\alpha}\), uniformly in the global atlas-supremum Bochner
 norm.  At this completion step we use the high self-map bounds
 \(q=k+3\), not the two-order-lower contraction estimate:
 \[
  (A_{\mathcal U}(\mathbf w_n^\sharp),
    B_{\mathcal U}(\mathbf w_n^\sharp),
    C_{\mathcal U}(\mathbf w_n^\sharp),
    \mathscr R_{\mathcal U}(\mathbf w_n^\sharp))
 \]
 converges in the corresponding strong-Bochner
 \(L^\infty C^{k+1,\alpha}\) coefficient product norm.  Indeed, the
 common top-order modulus controls the high-frequency tails, the lower
 contraction topology controls the finitely many low frequencies, and
 the high-tier prepared coefficient calculus supplies exactly these
 one-order-higher bounds.  Consequently every term in the displayed
 parabolic right-hand side converges strongly in
 \(L^\infty C^{k+1,\alpha}\).  Passing to the limit in the Bochner
 integral identities shows that the limiting \(X_\Phi\) belongs to
 \(W^{1,\infty}C^{k+1,\alpha}\), exactly the time-derivative component
 of \(\mathbb E_{\rm sc}^{k+3,\alpha}\), and that its weak derivative
 is the limiting right-hand side.  The normalized graph identity and
 the prepared composition calculus separately give the required
 \(C^{k,\alpha}\)-valued time derivative for \(h\).
Consequently all required temporal traces and derivatives belong to
the asserted solution spaces, and passing to the limit in the ODE and
parabolic identities gives
\[
 \mathbf w=\mathfrak K(\mathbf w)
 \quad\hbox{in the full }\mathbb S^{k+2,\alpha}\hbox{ solution class}.
\]
Thus the equations themselves restore the required temporal traces
and derivatives, without placing them in the contraction metric.
At the fixed point,
\(\mathbf w^\sharp(\mathbf w)=\mathbf w\); hence the limiting
coefficients and forcing are those of the original exact geometric
system, rather than of an auxiliary modified equation.
The fixed point satisfies \eqref{eq:relative-coupled-system}.  In
particular, its local elapsed clock is
\[
 \theta(\tau)=\int_s^\tau e^{\ell(q)}\,dq\geq0,
 \qquad
 t(\tau)-t(s)=\lambda_s\theta(\tau).
\]
It therefore samples only
\(K_s^{\rm tr}|_{[0,\widehat\delta_{\rm ph}]}=K_s\), and hence, by
\eqref{eq:compact-affine-reconstruction}, the genuine forward physical
Ricci-flow coefficient.  The
auxiliary negative-time
extension disappears from every geometric conclusion.  The same
difference estimate proves uniqueness in the prepared class and
continuous dependence.

\emph{Step 3: preservation of the diffeomorphisms.}
The vector field generating \(\Theta_\tau\) has at most linear growth
and uniformly bounded scale derivatives on the short interval.
Forward and backward integral curves therefore cannot escape to
infinity in finite time, so its flow is a global diffeomorphism.
For the same-output variable \(Q_\Phi\), the fixed-point estimate gives
\[
 \|X_\Phi(\tau)\|_{\mathfrak X_{\rm sc}^{1,\alpha}}\longrightarrow0
 \quad\text{as }\tau\downarrow\tau_0.
\]
The radial argument in the proof of
Lemma~\ref{lem:prepared-chart-calculus} shows that
\[
 (\sigma,x)\longmapsto
 \operatorname{Exp}(\sigma X_\Phi(\tau))\circ\Phi_0(x),
 \qquad0\leq\sigma\leq1,
\]
is a proper homotopy.  Smallness of the scaled \(C^1\) displacement
preserves local invertibility.  Hence \(\Phi_\tau\) is a proper local
diffeomorphism with the same degree, of absolute value one, as
\(\Phi_0\).  It is consequently a one-sheeted covering and thus a
global diffeomorphism.  Since \(\Theta_\tau\) is already a global
diffeomorphism, \(F_\tau=\Theta_\tau^{-1}\circ\Phi_\tau\) is one as
well.

\emph{Step 4: the normalized equation and the slice.}
Propositions~\ref{prop:relative-HMHF} and
\ref{prop:adaptive-graft} apply to the solution just constructed and
give the displayed adaptive normalized equation for \(h\).  Initially
the mild regularity makes the nine moments
\[
 m_\mu(\tau):=\ip{\rho_\tau h(\tau)}{Z_\mu},
 \qquad 0\leq\mu\leq8,
\]
absolutely continuous.  Their almost-everywhere derivatives satisfy
\[
 m_\mu'
 =\lambda_\mu m_\mu+
   \bigl(M(\tau,\mathbf z)c+d(\tau,\mathbf z)\bigr)_\mu
 =\lambda_\mu m_\mu,
\]
where the last equality is
\eqref{eq:algebraic-feedback-map}.  Since
\(m_\mu(\tau_0)=0\), uniqueness for these scalar linear ODEs gives
\(m_\mu\equiv0\).  Hence the exact receding slice is propagated.
For the fixed point just obtained, the coefficients and forcing are
continuous in normalized time in each fixed chart.  The classical
recovery clause of Lemma~\ref{lem:weighted-prepared-Schauder} therefore
upgrades the mild equation chartwise.  If the initial state is smooth,
the compact Ricci--DeTurck bootstrap,
Lemma~\ref{lem:weighted-prepared-Schauder} at successively higher
orders, and the \(R,\Theta,\lambda\) ODEs give all spatial and time
derivatives on every shorter interval.  Pulling back by the smooth
DeTurck flow then proves the asserted smoothness of the geometric
solution.

\emph{Step 5: \(C^1\) dependence.}
Work in one fixed chart of the sliced prepared manifold.  By
Lemmas~\ref{lem:prepared-chart-calculus} and
\ref{lem:feedback-functional-regularity}, and by the
parameter-dependent clause of
Lemma~\ref{lem:weighted-prepared-Schauder}, the coefficient, forcing,
trace, and ODE maps in the fixed-point system are Fr\'echet \(C^1\)
from prepared input order \(k+2\) to output order \(k\).  Their
operator estimates and the common contraction factor depend only on
the recorded high bound and prepared margins, not on an individual
spatial tail modulus.

Let \(\xi_n\to0\) be arbitrary chart increments with
\(\mathbf z_0+\xi_n\) in the sliced domain.  Write
\[
 \mathbf w_n=\mathcal S(\mathbf z_0+\xi_n),\qquad
 \mathbf w=\mathcal S(\mathbf z_0),
\]
and let \(\mathbf v_n\) be the solution of the coupled linearized
fixed-point system at \((\mathbf z_0,\mathbf w)\) with initial
increment \(\xi_n\).  The uniform local Lipschitz estimate gives
\[
 \|\mathbf w_n-\mathbf w\|_{\mathbb D_{\rm mix}^{k,\alpha}(I)}
 \leq C\|\xi_n\|_{\mathscr E_{\rm prep}^{k+2,\alpha}} .
\]
Since the initial increments converge in the order-\((k+2,\alpha)\)
little-H\"older topology, the sequential observation preceding
\eqref{eq:weighted-parabolic-spaces} gives the initial states one common
spatial tail.  The quantitative construction
\eqref{eq:coupled-tail-self-map}--%
\eqref{eq:invariant-tail-envelope} therefore gives
\(\{\mathbf w,\mathbf w_1,\mathbf w_2,\ldots\}\) one common invariant
high-order tail.  Splitting into finitely many low frequencies and a
common high-frequency remainder, as in the proof that
\(\mathbb X_{K,\rho_{\rm fp},\mu}\) retains its high regularity, now
shows
\[
 \mathbf w_n-\mathbf w\longrightarrow0
 \quad\text{in the high spatial atlas-supremum coefficient topology}.
\]
This is exactly the high factor occurring in the mixed product
remainders; no rate in that high spatial topology, and no convergence
claim in a stronger top time-H\"older norm, is used.

Subtract the two nonlinear systems and the linearized system.  The
mixed tame estimates
\eqref{eq:prepared-mixed-tame-remainder} and
\eqref{eq:feedback-mixed-tame-remainder}, followed by
\eqref{eq:Schauder-mixed-tame-residual}, make the residual trace and
forcing
\[
 o\!\left(
   \|\xi_n\|_{\mathscr E_{\rm prep}^{k+2,\alpha}}
 \right)
\]
in the output-order trace and Bochner norms: the high factor tends to
zero, while every low factor is bounded by the preceding Lipschitz
estimate.  If
\(\mathbf r_n:=\mathbf w_n-\mathbf w-\mathbf v_n\), the same
modulus-independent contraction and Schauder estimates give
\[
 \|\mathbf r_n\|_{\mathbb D_{\rm mix}^{k,\alpha}(I)}
 \leq\frac12
       \|\mathbf r_n\|_{\mathbb D_{\rm mix}^{k,\alpha}(I)}
 +o\!\left(
   \|\xi_n\|_{\mathscr E_{\rm prep}^{k+2,\alpha}}
 \right).
\]
Absorption first gives
\[
 \|\mathbf r_n\|_{\mathbb D_{\rm mix}^{k,\alpha}(I)}
 =o\!\left(
   \|\xi_n\|_{\mathscr E_{\rm prep}^{k+2,\alpha}}
 \right).
\]
For clarity, the compact metric component of this upgrade does not
come from the unscaled \(\mathscr Y_{\rm prep}\)-term alone.  Differentiate
the exact affine reconstruction
\[
 G(t(\tau))=G(t(s))+\lambda_sK_s(\theta(\tau)).
\]
For an initial tangent increment this gives
\[
 \delta G(t(\tau))
 =\delta G(t(s))+\delta\lambda_s K_s
  +\lambda_sDK_s[\delta\mathbf z_s]
  +\lambda_s(\partial_{\widehat\sigma}K_s)\,\delta\theta,
\]
with the last three terms evaluated at \(\theta(\tau)\).  The buffered
affine compact estimates control \(K_s\), \(DK_s\), and
\(\partial_{\widehat\sigma}K_s\) in the output physical-host
\(C^{k,\alpha}\) tier.  Since
\(\lambda_s/\lambda_s^\circ(\tau)\) is uniformly comparable on the local
interval and the input prepared norm controls
\((\lambda^\circ(s))^{-1}\delta G(t(s))\) and
\(\delta\log\lambda_s\), one obtains
\[
 \frac{\|\delta G(t(\tau))\|_{C^{k,\alpha}(\mathcal X)}}
      {\lambda_s^\circ(\tau)}
 \leq C
 \|\delta\mathbf z_s\|_{\mathscr E_{\rm prep}^{k+2,\alpha}}.
\]
The mixed affine compact remainder gives the analogous little-
\(o(\|\xi_n\|)\) estimate.  Substitution into the zero-trace Schauder
and ODE estimates then
upgrades this lower-metric remainder to the full output-order prepared
norm and proves
\[
 \|\mathbf w_n-\mathbf w-\mathbf v_n\|
      _{\mathscr E_{\rm prep}^{k,\alpha}}
 =
 o\!\left(
   \|\xi_n\|_{\mathscr E_{\rm prep}^{k+2,\alpha}}
 \right).
\]
Thus the local solution map is Fr\'echet differentiable for arbitrary
increments.  Subtracting the coupled linearized systems at two base
states and taking the supremum over unit tangent directions proves
continuity of its derivative in operator norm.  The individual
invariant moduli supplied by
\eqref{eq:invariant-tail-envelope} ensure little-H\"older membership;
their values do not enter the remainder estimate.  Consequently,
\[
 \mathcal S_{\tau,\tau_0}:
 \Sigma_{\tau_0}^{k+2,\alpha}
 \longrightarrow
 \mathscr P_\tau^{k,\alpha}
\]
is \(C^1\) on the asserted open sliced manifold, although the principal
coefficient is not \(C^{k+1,\alpha}\)-Lipschitz in the lower contraction
metric.  Finite-horizon differentiability is therefore not obtained by
composing these derivative-losing local maps.  Instead,
Proposition~\ref{prop:two-state-prepared-evolution} gives same-order
two-state and first-variation propagation, and hence the finite-horizon
statement without loss at restart times.  On the
present local interval, composing with the buffered phase retraction gives
\eqref{eq:local-ambient-solution-map}; the additional two derivatives
belong to the phase projection, not to the sliced evolution.

\emph{Step 6: continuation.}
Suppose that \(\tau^\dagger<\infty\) and none of
\textup{(1)}--\textup{(6)} occurs.  Let \(\sigma<\tau^\dagger\) and
\(\mu_*>0\) be the common terminal margin supplied by
Definition~\ref{def:admissible-first-exit-interval}.  The phase budget
bounds
\(\int_{\tau_0}^{\tau^\dagger}|a|\,d\tau\), and hence
\[
 0<c\leq\lambda(\tau)e^\tau\leq C<\infty
 \quad\text{on }[\tau_0,\tau^\dagger).
\]
Thus \(t_\tau=\lambda\) is controlled above and below on every
terminal compact subinterval.  Hamilton's extension theorem
\cite{Hamilton} extends the closed Ricci flow to
\(t(\tau^\dagger)\), and the compact Shi--Schauder
estimates give uniform restart bounds at order \(k+2\) away from the
initial physical time.  On \([\sigma,\tau^\dagger)\), the ellipticity,
radial, graft, bounded-geometry, Gram-inverse, and
controlled-diffeomorphism margins are bounded below by \(\mu_*\).
The corresponding metric, Gram, geometric-package, and map bounds
therefore give the uniform coefficient estimates required in
Lemma~\ref{lem:weighted-prepared-Schauder}.  Applying the lemma on a
fixed terminal normalized-time interval gives uniform
\(C^{k+3,\alpha}\) restart bounds for \(F\); the ODEs give the
corresponding bounds for \(R,\lambda\), and their inverses.
When Corollary~\ref{cor:adaptive-auxiliary-closure} applies, it
supplies precisely the higher-order portion of this
bounded-geometry package, so that portion is not a separate exit face.

The equations and these uniform bounds give the two-order-lower
endpoint trace asserted in
Definition~\ref{def:admissible-first-exit-interval}; no top-order
endpoint trace and no compactness of the terminal family in the
top-order little-H\"older norm are being inferred.  The terminal
restart bounds above determine one finite ceiling
\[
 K_{\rm init,term}^{k+2,\alpha}
\]
for all states \(\mathbf z(s)\), \(\sigma\leq s<\tau^\dagger\),
together with one fixed collection of positive geometric margins.  For
each such state, the tail estimate
\eqref{eq:coupled-tail-self-map} supplies an individual invariant
modulus \(\mu_{\mathbf z(s)}\), while the admissible lifespan furnished
by the local construction depends only on
\(K_{\rm init,term}^{k+2,\alpha}\) and those margins, not on
\(\mu_{\mathbf z(s)}\).  Hence there is one
\(\delta_*>0\), independent of \(s\).  Choose \(s\) so close to
\(\tau^\dagger\) that \(s+\delta_*>\tau^\dagger\).  Uniqueness on the
overlap then extends the original solution past
\(\tau^\dagger\), a contradiction.  This proves the continuation
alternative.  Finally, at every time the adaptive Gram matrix differs
from its model matrix by an
\(O(\|h\|_{C^1})\) term and the tail in
\eqref{eq:effective-column-tail}; hence item \textup{(4)} is absent in
the sufficiently small box.  Corollary~\ref{cor:adaptive-auxiliary-closure}
gives the last assertion.
\end{proof}

\begin{lemma}[Finite-horizon approximation and diagonal passage]
\label{lem:finite-regularity-prepared-approximation}
Fix \(k\geq12\) and \(0<\alpha<1\).  Let
\(\mathbf z_0\in\Sigma_{\tau_0}^{k+2,\alpha}\) be a strict prepared
state in the little-H\"older chart, with a fixed positive
common-margin package.  Write first
\[
 \Sigma_{\tau_0}^{\infty}:=\bigcap_{r\geq3}
 \Sigma_{\tau_0}^{r,\alpha}
\]
for the smooth sliced class.  There are smooth sliced prepared states
\[
 \mathbf z_0^{(n)}\in\Sigma_{\tau_0}^{\infty},
 \qquad
 \mathbf z_0^{(n)}\longrightarrow\mathbf z_0
 \quad\hbox{in }\mathscr P_{\tau_0}^{k+2,\alpha},
\]
For all sufficiently large \(n\), these states retain one common
package with at least half of every ordinary numerical margin, together
with the same operative normalized triple, a new strictly stronger
normalized reserve, the unchanged common physical \(+\)- and operative
tiers from the fixed reference-carrier certificate, and the common
future-window certificate.

Let \(S>\tau_0\) be finite.  Suppose that, on the maximal admissible
intervals of these smooth solutions, one has endpoint-independent
package bounds \emph{with a common positive exit-face margin} for the
quantities used in the continuation argument:
the exact slice and structural identities; the scale and phase budget;
uniform ellipticity and the global \(C^2\) box; the order-\((k+2)\)
prepared coefficient, map, inverse-map, graft, and bounded-geometry
package; and the \(L^2_\nu\), \(H^1_\nu\)-dissipation, three-region, and
pure-graft inequalities.  Concretely, after all scale-normalized
thresholds have been fixed, ellipticity stays uniformly above its lower
face, the \(C^2\), phase, graft, and analytic bootstrap quantities stay
uniformly below their upper faces, the Gram inverse remains uniformly
bounded, and the radial-comparison, proper local-invertibility,
inverse-map, support-separation, and source/target bounded-geometry
margins retain one fixed positive slack.  Thus ``package bounds'' here
means strict improvement of every continuation inequality, not mere
boundedness of its defining functional.  Then every sufficiently large
smooth approximant exists on \([\tau_0,S]\), the finite-regularity solution
exists there.  Define
\(\mathscr P^{k,\alpha}_{\rm loc}\) to mean convergence of every
prepared component in the fixed compact charts at the displayed
order.  Then
\begin{equation}\label{eq:finite-regularity-finite-horizon-convergence}
 \mathbf z^{(n)}\longrightarrow\mathbf z
 \quad\text{in }C^0\bigl([\tau_0,S];
       \mathscr P^{k,\alpha}_{\rm loc}\bigr).
\end{equation}
The convergence also holds in the global weighted prepared output
topology at every strictly
lower H\"older exponent.  In the little-H\"older output topology the
convergence is also strong.  Indeed, the entrance sequence converges
in the little-H\"older norm, so the sequential observation preceding
\eqref{eq:weighted-parabolic-spaces} supplies one common datum modulus;
then \eqref{eq:invariant-tail-envelope} supplies one common invariant
solution modulus.  The output modulus need not equal the original
datum modulus.  No
annular \(c_0\) condition is imposed.  Instead, the proof below uses
the global prepared Banach distance in the local solution theorem to
make the approximants globally Cauchy, so annular escape is excluded
without shrinking the entrance class.

The exact identities, scale and phase estimates, \(L^2_\nu\) estimates,
and pointwise estimates through order \(k\) pass by strong convergence.
The spacetime \(H^1_\nu\)-dissipation inequality passes by weak lower
semicontinuity.  No estimate in a topology stronger than the displayed
uniform package is inferred merely from
\eqref{eq:finite-regularity-finite-horizon-convergence}.  Instead, for
every \(\delta>0\) and \(m\geq0\), nested interior parabolic estimates
give for the gauge-fixed closed Ricci--DeTurck component
\begin{equation}\label{eq:finite-regularity-positive-time-diagonal}
 \sup_n\sup_{\tau_0+\delta\leq\tau\leq S}
 \|\widetilde G^{(n)}(\tau)\|_{C^m(\mathcal X)}
 <\infty .
\end{equation}
For \(h^{(n)}\), the identical all-order estimate holds on every recent
compact parabolic cylinder on which both the cutoff exterior and the
transported graft support have exited and hence
\(\mathcal Y_{j,\tau}=Y_j\); before support exit it
is asserted only through the finite order carried by the prepared
package.  In particular, for each \(K\Subset M\) all orders are
available for sufficiently late \(\tau\), and on the moving core they
are available whenever the backward cylinder lies in
\(\{\bar f<e^{\tau-1}\}\), exactly the range used below.  After a
diagonal extraction these parabolic components converge smoothly on
the stated positive-time cylinders.  Uniqueness of the
finite-regularity mild solution removes the extraction.  The
ODE-carried maps are not
claimed to gain derivatives beyond their entrance order.  If the
package bounds hold with constants independent of every finite \(S\),
a diagonal argument over \(S\to\infty\) gives the global
finite-regularity evolution with the same low-order estimates.  The
same statement holds at a prepared restart time after translating
\(\tau_0\) to that time.
\end{lemma}

\begin{proof}
Smooth the tensor variables and the same-output soliton-conjugated map
variables in the fixed core
atlas and in the uniformly locally finite rescaled dyadic atlas,
using the common tail modulus from the little-H\"older definition.
This is approximation of the one fixed entrance state, not a
compactness assertion for an arbitrary bounded family: its
little-H\"older frequency modulus makes the mollification error
uniform on every annulus.  The bounded overlap and fixed weights
therefore give convergence in the global prepared
\(C^{k+2,\alpha}\) norm.  Properness, degree one, ellipticity, radial
comparison, support separation, graft compatibility, and bounded
geometry persist by the openness statement in
Lemma~\ref{lem:prepared-chart-calculus}.

The smoothed states need not satisfy the nine moments exactly.  Their
moment errors tend to zero.  The derivative of the static moment map
in the nine centered phase directions is the invertible prepared Gram
matrix.  The finite-dimensional implicit-function correction
therefore changes the phase parameters by \(o(1)\), preserves
smoothness, and produces smooth states in
\(\Sigma_{\tau_0}^{\infty}\).  Continuity of the prepared graph at the
same finite order, also recorded in
Lemma~\ref{lem:prepared-chart-calculus}, gives the asserted
\(C^{k+2,\alpha}\) convergence; the two-order buffer is needed for the
Fr\'echet derivative of the phase retraction, not for this continuity
statement.

Apply the reserve-to-operative clause of
Lemma~\ref{lem:prepared-harmonic-radius-lower-stability} to the
corrected normalized graph tensors.  Because the corrected states
converge at the required scale-one \(C^{2,\alpha}\) order, this
preserves the normalized operative triple and leaves a new common
reserve strictly above it.  The actual carriers of the corrected
states remain, with positive slack, in the same inner locus
\eqref{eq:physical-coefficient-inner-locus} based at the original
frozen \(G_{\rm ref}^{\rm phys}\).  Hence
Lemma~\ref{lem:reference-carrier-physical-certificate} gives every one
the fixed common physical \(+\)-tier and operative tier, and they lie in
the outer ball.
 The uniform half-open estimate
 \eqref{eq:witnessed-physical-quarter-modulus} is then the precomputed
 assertion for every member of that outer ball; it is not inferred by
 continuity of Ricci flow with respect to the smoothed input.
 Continuity is used only to retain inner-locus membership and the
 \(\mu_{\rm RF}\)-slack in the lifetime-width inequality.  Thus the
 approximation retains the complete dynamic witness certificate and
 does not replace it by bare harmonic-radius lower bounds.

All corrected states eventually lie in one smaller common-margin
ball with one recorded \(K_{\rm init}^{k+2,\alpha}\).
Proposition~\ref{prop:coupled-local-feedback}, whose lifespan is
independent of the individual invariant tail modulus, then gives a
common lifespan and continuous dependence in the output
\(C^{k,\alpha}\) topology.

Fix \(S>\tau_0\).  For a smooth approximant, let \(S_n\leq S\) be the
largest time through which it remains in the common package.  The
assumed endpoint-independent estimates keep every continuation
quantity a fixed positive distance from its exit face on
\([\tau_0,S_n]\).  The continuation alternative in
Proposition~\ref{prop:coupled-local-feedback}, applied at \(S_n\) with
the same package constants and recorded high-norm ceiling, supplies a
restart interval of uniform length independently of the individual
tail modulus.  Hence \(S_n=S\) for all sufficiently large \(n\).

The package bounds give uniform boundedness at the top prepared order
and strong compactness two orders lower on each fixed chart.  We do not
upgrade that chartwise statement by bounded overlap: such an argument
would permit annular escape.  Instead use the genuinely global
prepared distance in
the global Lipschitz estimate
\eqref{eq:local-global-prepared-Lipschitz}.  Since the smoothed
entrances are Cauchy in
\(\mathscr P_{\tau_0}^{k+2,\alpha}\), that proposition makes their
solutions Cauchy in the single fixed same-output restart chart used in
\eqref{eq:prepared-Banach-norm}, in the precise sense that
\[
 \sup_{\tau_0\leq\tau\leq\tau_0+\delta_0}
 d_{\rm prep}^{k,\alpha}
   \bigl(\mathbf z_n(\tau),\mathbf z_m(\tau)\bigr)
 \longrightarrow0 .
\]
Here the difference of two time-\(\tau\) states is measured by the
global graph-augmented distance \eqref{eq:prepared-Banach-norm} after the common
same-output identification of their \(\Theta\)- and \(\Phi\)-blocks; in
particular, its end component is the supremum over the entire source
atlas, not over a fixed finite collection of charts.  The uniform
\(k+2,\alpha\) prepared package and global interpolation
then give strong convergence with respect to
\(d_{\rm prep}^{k+2,\alpha'}\) for every
\(\alpha'<\alpha\).
Choose a finite common restart partition
\[
 \tau_0=s_0<s_1<\cdots<s_J=S
\]
whose mesh is at most one half of the uniform local lifespan supplied
by the preceding continuation argument.  We propagate both convergence
and a common top spatial modulus across this partition.  At
\(s_0=\tau_0\), norm convergence in the little-H\"older prepared space
gives the entrance sequence one common datum modulus.  The invariant
tail estimate \eqref{eq:invariant-tail-envelope}, whose constants are
independent of that modulus, gives all solutions one common modulus on
\([s_0,s_1]\).  The global Lipschitz estimate gives convergence two
orders lower there.  In each normalized chart split a representative
into finitely many low spatial frequencies and its high-frequency
remainder.  The lower-order convergence controls the former uniformly
over the complete atlas, while the common modulus controls the latter.
Thus the endpoint states converge in the full high spatial coefficient
topology.

The same-output restart coordinates and all fixed-order prepared
coordinate changes are continuous at this same spatial order by
Lemma~\ref{lem:prepared-chart-calculus}; the two-order loss there is
needed for Fr\'echet differentiation, not for continuity.  Consequently
the endpoint sequence, expressed in the chart frozen at \(s_1\), again
has one common datum modulus and remains in the same common-margin
ball.  Repeating the preceding argument on
\([s_1,s_2]\), and then finitely many times, propagates the common
modulus and the high spatial convergence through \(S\).  At the same
time, applying the global local theorem along a finite descending
exponent ladder gives global convergence in the weighted prepared
output topology at every prescribed \(\alpha'<\alpha\).  This direct
global Cauchy argument, rather than local compactness, rules out an
annularly translated bump.

The global limit passes every coefficient and the finite-dimensional
Gram feedback, so it solves the exact coupled mild system with
entrance \(\mathbf z_0\).  Local uniqueness identifies it with
\(\mathbf z\) and removes every subsequence.  At the original exponent,
the same finite-frequency/common-tail split on the finite restart cover
upgrades the convergence to the stated global little-H\"older output
topology, without any annular-tightness hypothesis.

Strong convergence passes the exact identities and the stated
finite-order inequalities.  Uniform \(L^2_\tau H^1_\nu\) control gives
weak compactness of the spatial gradients and hence the dissipation
inequality by lower semicontinuity.  On
\([\tau_0+\delta,S]\), uniform ellipticity and the compact
Ricci--DeTurck equation permit the usual nested parabolic Schauder
estimates at arbitrary order, proving
\eqref{eq:finite-regularity-positive-time-diagonal}.  The normalized
equation has the same property on a recent compact cylinder after both
the cutoff exterior and the transported graft support have exited that
cylinder; before that time its coefficient package permits iteration
only through its stated finite
order.  Arzel\`a--Ascoli on a countable exhaustion of the permitted
cylinders gives the asserted smooth convergence.  This is a new
interior estimate for the limit, not a passage of an unbounded entrance
norm.  Finally take
\(S=\tau_0+1,\tau_0+2,\ldots\) and diagonalize.  Translating the
construction proves the restart version.
\end{proof}

\section{Shared exterior control, finite-horizon comparison, and first
variations}
\label{sec:finite-horizon-two-state}

The one-state flow has now been constructed.  The terminated exterior
atlas and anchored-gauge clauses recorded here are shared one-state
tools used by the formation branch.  The difference and first-variation
estimates are proved alongside them but are not used until Part~III,
where they enter the proofs of the prepared scattering theorem and
Theorem~B.

\begin{lemma}[Dynamic effective-column differences]
\label{lem:dynamic-effective-column-differences}
Fix a finite endpoint \(S>\tau_0\) and a common-margin prepared ball at
order \(k+2\).  If \(\mathbf z_i(\tau)\), \(i=1,2\), are two admissible
states in that ball for \(\tau_0\leq\tau\leq S\), then, for
\(0\leq m\leq3\),
\begin{equation}\label{eq:dynamic-column-Lipschitz}
 \|\mathcal Y_{j,\tau}(\mathbf z_1(\tau))
       -\mathcal Y_{j,\tau}(\mathbf z_2(\tau))\|_{H^m_\nu}
 \leq C_S e^{-ce^\tau}
       \|\mathbf z_1(\tau)-\mathbf z_2(\tau)\|
          _{\mathscr X_{\rm prep}^{k,\alpha}},
 \qquad0\leq j\leq8.
\end{equation}
The same estimate holds after pairing with
\(\rho_\tau Z_\mu\), for \(0\leq\mu\leq8\).
In addition, for the fixed weight exponent \(N\) of the prepared ball,
\begin{equation}\label{eq:dynamic-column-high-order-Lipschitz}
 \|\mathcal Y_{j,\tau}(\mathbf z_1(\tau))
       -\mathcal Y_{j,\tau}(\mathbf z_2(\tau))\|
       _{\mathfrak T_{{\rm sc},N}^{k-2,\alpha}}
 \leq C_{S,N}
 \left(
  d_{\rm rt,sc}^{k+1,\alpha}(R_1,R_2)
  +d_{\rm rt,sc}^{k+1,\alpha}(F_1,F_2)
 \right),
 \qquad0\leq j\leq8 .
\end{equation}
The high-order bound, with the corresponding tangent norm, also holds
for sliced first variations; its Fr\'echet remainder is \(o(1)\) times
the increment norm on the common prepared ball.
\end{lemma}

\begin{proof}
The global part of each direct column is independent of the state.
Consequently the difference consists only of differences of the
pullback \(K_\tau(T)\), and is supported where
\(\bar f\geq ce^\tau\).  For the asserted same-order estimate, apply
Lemma~\ref{lem:prepared-chart-calculus} at output order \(3\) and input
order \(5\), rather than at output order \(k\).  The standing assumption
\(k\geq12\) implies that the
\(\mathscr X_{\rm prep}^{k,\alpha}\)-distance controls every required
order-five input component.  The differentiated column formula uses at
most four derivatives of its map factors.  Hence its derivative, as a
map into \(H^m_\nu\) for \(0\leq m\leq3\), is bounded by the displayed
same-order prepared distance.  The mean-value formula in the common
same-output soliton-conjugated prepared chart gives a polynomial-growth difference
on the receding support.  Gaussian integration absorbs that growth and
proves \eqref{eq:dynamic-column-Lipschitz}.  Multiplication by
\(\rho_\tau Z_\mu\) and Cauchy--Schwarz give the pairing statement.
For \eqref{eq:dynamic-column-high-order-Lipschitz}, apply the same
 prepared-chart calculus at output order \(k-2\) and input order \(k\)
 in every fixed core and scale-normalized end chart.  The global column
 cancels, while the remaining pullback depends only on the \(R\)- and
 \(F\)-blocks (including their induced \(\Theta\)- and \(\Phi\)-maps);
 because both states lie in one common-margin right-translated chart,
 the two typed distances in
 \eqref{eq:dynamic-column-high-order-Lipschitz} control every map,
 inverse-map, and composition input factor.  The prescribed
factor \(L^{-N}\) absorbs the polynomial column growth uniformly over
the dyadic atlas.  This gives the displayed
\(\mathfrak T_{{\rm sc},N}^{k-2,\alpha}\) bound without using Gaussian
smallness.  The prepared composition maps are \(C^1\) with the
two-derivative input buffer, so differentiation and Taylor's formula
give the final two assertions.
\end{proof}

\subsection{Terminated exterior interface and hybrid distance}
\label{subsubsec:terminated-interface-hybrid-distance}

Fix as part of the common prepared package four smooth physical core
domains
\[
 K_-\Subset K_0\Subset K_1\Subset K_2\Subset\mathcal X'',
 \qquad
 \iota\bigl(K_2\setminus\overline{K_-}\bigr)
 \Subset\operatorname{int}\{\eta=1\}.
\]
Put
\[
 E^{++}:=\mathcal X\setminus\overline{K_-},\qquad
 E^+:=\mathcal X\setminus\overline{K_0},\qquad
 E:=\mathcal X\setminus\overline{K_2}.
\]
Thus \(E\Subset E^+\Subset E^{++}\).  Fix a smooth compact
codimension-zero collar
\[
 \mathcal A_{\rm in}\Subset K_1\setminus\overline{K_0}
\]
and a cutoff \(\zeta\in C^\infty(E^{++};[0,1])\) which vanishes on a
neighborhood of \(\partial E^{++}\), equals one on a neighborhood of
\(\overline E\), and satisfies
\[
 \supp d\zeta\Subset\operatorname{int}\mathcal A_{\rm in}.
\]
This is the normalized inner interface at which the exterior
localization terminates.  Cover
\(\overline E\) by an ordinary finite buffered atlas
\[
 U_a\Subset U_a^+\Subset U_a^{++},\qquad
 1\leq a\leq N_{\rm ext},
\]
whose largest members are relatively compact in \(E^{++}\), with local
noncollapsing scales \(R_a\).  Fix one exterior reference metric
\(g_{\rm ext}\), uniformly equivalent to the initial physical metrics
in the common package.  For each retained largest member fix a
auxiliary six-level enlargement
\[
 U_a^{++}\Subset V_a^0\Subset V_a^1\Subset\cdots
 \Subset V_a^5\Subset E^{++},
\]
with every successive \(g_{\rm ext}\)-separation bounded below by
\(c_{\rm atl}R_a\).  These auxiliary sets carry the one-state local Ricci
and harmonic-map estimates; no difference norm is taken on them.
Require only
\begin{equation}\label{eq:prepared-exterior-termination}
\begin{aligned}
 \overline E&\subset\bigcup_{a=1}^{N_{\rm ext}}U_a,&
 \overline {E^+}&\subset\bigcup_{a=1}^{N_{\rm ext}}U_a^{++},\\
 \#\{b:U_a^{++}\cap U_b^{++}\ne\varnothing\}
 &\leq N_{\rm cov},&
 R_a&\asymp_{C_{\rm cov}}R_b
 \quad\text{on overlaps},\\
 \min\!\left\{
 d_{g_{\rm ext}}\!\left(\supp d\zeta,
              \bigcup_{a=1}^{N_{\rm ext}}\overline{U_a}\right),
 d_{g_{\rm ext}}(\supp d\zeta,\partial E^{++})
 \right\}
 &\geq R_{\rm in}>0 .
\end{aligned}
\end{equation}
In addition, the fixed atlas has a \(g_{\rm ext}\)-Lebesgue number
\(\ell_{\rm Leb}>0\) on \(\overline E\), and its rescaled chart,
partition, and auxiliary-buffer constants are bounded by
\(\Lambda_{\rm atl}\).  The retained largest members carry the
one-state coefficient bounds and all
members are expressed in the single anchored exterior Ricci--DeTurck
gauge constructed below.
Fix, once with this atlas, a cutoff
\(\vartheta_{\rm ext}\in C_c^\infty(\bigcup_aU_a^+;[0,1])\) which is
one on a neighborhood of \(\bigcup_a\overline{U_a}\), and decrease a
fixed \(R_{{\rm ext},0}>0\), if necessary, so that
\begin{equation}\label{eq:exterior-initial-corridor-separation}
 d_{g_{\rm ext}}\!\left(
   \bigcup_a\overline{U_a},
   \supp(1-\vartheta_{\rm ext})\right)
 \geq4R_{{\rm ext},0}.
\end{equation}
This cutoff is used only to split the homogeneous initial face into a
buffered local trace and a genuinely separated corridor trace.
There is deliberately no cyclic requirement that the largest buffers
be covered by smaller buffers.  The hierarchy terminates instead at
the separated normalized interface \(\mathcal A_{\rm in}\).
The cutoff and all atlas data are fixed once with the package.  Hence a
H\"older quotient whose two points
are closer than \(\ell_{\rm Leb}\) is computed in one enlarged chart; if the
points are farther apart, it is bounded by the corresponding
scale-normalized supremum norm.  Thus the local H\"older estimates glue
with constants depending only on the recorded exterior certificate.

We keep two geometrically distinct collar systems.  First fix an
\emph{interface collar}
\[
 \mathcal W_{\rm in}\Subset\mathcal W_{\rm in}^+
 \Subset\mathcal W_{\rm in}^{++}\Subset E^+\cap\mathcal X'',
 \qquad
 \mathcal A_{\rm in}\Subset\operatorname{int}\mathcal W_{\rm in}.
\]
It crosses the artificial inner boundary of the terminated exterior
problem and is used only to read its normalized trace.  Independently,
after choosing \(K_2\) strictly inside the region on which \(\eta=1\),
fix the \emph{graft-input collars}
\begin{equation}\label{eq:separated-graft-collars}
 \mathcal W_{\rm gr}\Subset\mathcal W_{\rm gr}^+
 \Subset\mathcal W_{\rm gr}^{++}\Subset E\cap\mathcal X'',
 \qquad
 \Omega_\eta^{++}\Subset\operatorname{int}\iota(\mathcal W_{\rm gr}).
\end{equation}
Thus the physical transition of \(\eta\) lies in the smallest graft
collar, whereas \(\mathcal A_{\rm in}\subset K_1\subset K_2\) lies
strictly outside it.  In particular, no collar contained in \(E\) is
ever required to contain \(\mathcal A_{\rm in}\).  Fix an auxiliary
graft-input chain
\[
 \mathcal W_{\rm gr}^+
 \Subset\mathcal W_{\rm gr}^{0}
 \Subset\mathcal W_{\rm gr}^{1}
 \Subset\cdots
 \Subset\mathcal W_{\rm gr}^{5}
 \Subset\mathcal W_{\rm gr}^{++},
\]
with fixed positive scale-normalized separations.  These auxiliary collars
are used by one-time initial-face and off-diagonal interior estimates;
they are not restarted with the unknown trace on the next collar.
Put the two closed flows in
one time-independent reference gauge,
\[
 G_i(t)=\chi_i(t)^*\widetilde G_i(t),\qquad
 \widetilde\iota_i(t)=\iota\circ\chi_i(t)^{-1}.
\]
Here \(\iota_1=\iota_2=\iota\) is the fixed marking of the common
prepared chart; the marking is not an additional two-state variable.
The common displacement margin is chosen so that the transported
middle interface and graft sets remain compactly inside
\(\mathcal W_{\rm in}^{++}\) and
\(\mathcal W_{\rm gr}^{++}\), respectively, throughout the interval;
hence both transported markings are defined with one unused spatial
buffer in each collar system.
The prepared package supplies the one-state coefficient bounds on the
largest physical and marked sets, the difference norms below are taken
on the smallest sets covering \(E\), interface comparisons are made on
\(\mathcal W_{\rm in}\), and graft comparisons are made only on
\(\mathcal W_{\rm gr}\).  On \(\mathcal A_{\rm in}\), define the
derived interface map
\[
 P_i:=\Phi_i\circ\widetilde\iota_i .
\]
For the rest of this subsection write
\[
 c_i=(a_i,b_i),\qquad
 \delta c=c_1-c_2=(a_1-a_2,b_1-b_2).
\]
The exact identities
\begin{equation}\label{eq:inner-terminal-graph-identity}
 (\widetilde\iota_i)_*\widetilde G_i
 =\iota_*G_i
 =\lambda_i\Phi_i^*(\bar g+h_i),
 \qquad
 \widetilde G_i=P_i^*\bigl(\lambda_i(\bar g+h_i)\bigr)
\end{equation}
recover the low-order physical interface from the normalized graph,
scale, and map blocks.  In the fixed scale-normalized collar charts set
the following convention before using the corresponding norms:
\(C_{\rm sc}^{r,\alpha}\) means the ordinary covariant
\(C^{r,\alpha}\) norm after rescaling the fixed collar atlas to its
recorded graft scale; map differences are read in the fixed
right-translated exponential charts and tensor norms include the
corresponding scale factors.  More precisely, throughout the
two-state argument the shorthand
\[
 \|\psi_1-\psi_2\|_{\mathfrak X_{\rm sc}^{r,\alpha}}
 :=d_{\rm rt,sc}^{r,\alpha}(\psi_1,\psi_2)
\]
means the two-sided forward-plus-inverse distance in
\eqref{eq:right-translated-map-distance}; on a restricted star it
means the corresponding localized two-sided distance.
\begin{equation}\label{eq:hybrid-interface-block}
 \mathfrak I_{\rm in}(\tau):=
 \left|\log\frac{\lambda_1}{\lambda_2}\right|
 +\|P_1-P_2\|_{C_{\rm sc}^{2,\alpha}(\mathcal A_{\rm in})}
 +\|h_1-h_2\|_{C_{\rm sc}^{1,\alpha}
      (P_1(\mathcal A_{\rm in})\cup P_2(\mathcal A_{\rm in}))}.
\end{equation}
Put
\begin{equation}\label{eq:hybrid-interface-algebraic-data}
\begin{split}
 \mathfrak J_{\rm in}:={}&
 \left|\log\frac{\lambda_1}{\lambda_2}\right|
 +\|R_1-R_2\|_{\mathfrak X_{\rm sc}^{3,\alpha}}
 +\|F_1-F_2\|_{\mathfrak X_{\rm sc}^{3,\alpha}}\\
 &+\|\chi_1^{-1}-\chi_2^{-1}\|
       _{C_{\rm sc}^{2,\alpha}(\mathcal W_{\rm in}^+)}
 +\|h_1-h_2\|_{C_{\rm sc}^{1,\alpha}
       (P_1(\mathcal A_{\rm in})\cup P_2(\mathcal A_{\rm in}))}.
\end{split}
\end{equation}
For an integer \(m\geq4\), define the hybrid distance before using it
in any estimate.  With the analogous conventions on the physical
atlas and model end, set
\[
\begin{aligned}
\mathfrak D_m^{\rm hyb,full}:={}&
 |t_1-t_2|+\left|\log\frac{\lambda_1}{\lambda_2}\right|
 +\sum_a\|\widetilde G_1-\widetilde G_2\|
                 _{C_{R_a}^{m,\alpha}(U_a)}\\
&+\sum_a\bigl(
 \|\chi_1-\chi_2\|_{C_{R_a}^{m-1,\alpha}(U_a^+)}
 +\|\chi_1^{-1}-\chi_2^{-1}\|
                 _{C_{R_a}^{m-1,\alpha}(U_a^+)}\bigr)\\
&+\|\chi_1-\chi_2\|_{C_{\rm sc}^{m-1,\alpha}(\mathcal W_{\rm in}^+)}
 +\|\chi_1^{-1}-\chi_2^{-1}\|
                 _{C_{\rm sc}^{m-1,\alpha}(\mathcal W_{\rm in}^+)}
 +\|\widetilde\iota_1-\widetilde\iota_2\|
                 _{C^{m-1,\alpha}(\mathcal W_{\rm in}^+)}\\
&+\|\widetilde G_1-\widetilde G_2\|
                 _{C_{\rm sc}^{m+2,\alpha}(\mathcal W_{\rm gr}^+)}
 +\|\chi_1-\chi_2\|_{C_{\rm sc}^{m+1,\alpha}(\mathcal W_{\rm gr})}
 +\|\chi_1^{-1}-\chi_2^{-1}\|
                 _{C_{\rm sc}^{m+1,\alpha}(\mathcal W_{\rm gr})}\\
&+\|\widetilde\iota_1-\widetilde\iota_2\|
                 _{C^{m+1,\alpha}(\mathcal W_{\rm gr})}
 +\|R_1-R_2\|_{\mathfrak X_{\rm sc}^{m+1,\alpha}}
 +\|F_1-F_2\|_{\mathfrak X_{\rm sc}^{m+1,\alpha}}\\
&+\|h_1-h_2\|_{\mathfrak T_{{\rm sc},N}^{m,\alpha}}
 +\mathfrak I_{\rm in}.
\end{aligned}
\]
This is the \emph{full finite-horizon distance}; in particular its
\(F\)-summand is global.  The localized uniform distance is defined
separately after the gauge construction.  The block decomposition and
its numbered formula are recorded after
the gauge construction in
\eqref{eq:hybrid-physical-block}--\eqref{eq:two-state-hybrid-distance}.
Right translation is understood in the map term.  Since
\(\Phi=\varphi_\tau\circ R\circ F\), the prepared composition calculus
controls the algebraic composition.  To close its time propagation, we
use the following exact inner-collar identity.

\begin{lemma}[Closed normalized inner-trace map]
\label{lem:closed-normalized-inner-trace}
On a fixed marked collar
\(\mathcal C_{\rm in}^+\Subset\mathcal X''\) containing
\(\overline{K_2\setminus K_-}\) and satisfying
 \(\iota(\mathcal C_{\rm in}^+)\Subset\{\eta=1\}\), put
 \[
 Q_i:=\Phi_i\circ\iota,\qquad
 g_i:=\bar g+h_i,\qquad
 U_i:=U_{b_i}:=\sum_{j=1}^8b_{i,j}\chi_\tau W_j,
 \qquad i=1,2.
 \]
Then
\begin{equation}\label{eq:closed-normalized-inner-trace}
 G_i=\lambda_iQ_i^*g_i,\qquad
 \partial_\tau Q_i=
 \bigl((1+a_i)\bar\nabla\bar f-U_i-B_{\bar g}(g_i)\bigr)\circ Q_i .
\end{equation}
Consequently, fix a finite normalized interval \([\tau_0,S]\) and a
collar
\[
 \mathcal A_{\rm in}\Subset
 \mathcal C_{\rm in}\Subset\mathcal C_{\rm in}^+.
\]
For \(\tau_0\leq s\leq\tau\leq S\), put
\[
 \mathcal O_{\rm in}(s,\tau)
 :=
 \bigcup_{i=1}^2\ \bigcup_{q\in[s,\tau]}
 Q_i(q)(\mathcal C_{\rm in})\Subset M .
\]
In the scale-normalized charts one has
\begin{equation}\label{eq:closed-normalized-inner-trace-difference}
 \begin{split}
 &\|Q_1(\tau)-Q_2(\tau)\|_
      {C_{\rm sc}^{3,\alpha}(\mathcal C_{\rm in})}\\
 &\quad\leq
 C_{S,\mathcal C_{\rm in}}
 \|Q_1(s)-Q_2(s)\|_
      {C_{\rm sc}^{3,\alpha}(\mathcal C_{\rm in})}\\
 &\qquad+
 C_{S,\mathcal C_{\rm in}}\int_s^\tau\left(
  |\delta c(q)|+
  \|h_1(q)-h_2(q)\|_
   {C_{\rm sc}^{4,\alpha}(\mathcal O_{\rm in}(s,\tau))}
 \right)dq .
 \end{split}
\end{equation}
The same quantified estimate, on the same domains, holds for sliced
first variations.
This finite-horizon estimate is not used to obtain a uniform
unweighted \(C^{4,\alpha}\) trace on the receding collar.
\end{lemma}

\begin{proof}
On \(\{\eta=1\}\), the controlled harmonic-map equation and naturality
give
\[
 \lambda_i\Delta_{G_i,\bar g}Q_i
 =-B_{\bar g}(g_i)\circ Q_i .
\]
Combining this with the prescribed drift in
\eqref{eq:controlled-HMH} proves
\eqref{eq:closed-normalized-inner-trace}.  The fixed-collar
composition estimate and Gronwall give
\eqref{eq:closed-normalized-inner-trace-difference}; differentiating
the same identity gives the first-variation statement.
\end{proof}

\subsection{Anchored exterior gauge and interface trace}
\label{subsubsec:anchored-exterior-gauge}

\begin{lemma}[Anchored exterior gauge and algebraic interface trace]
\label{lem:anchored-exterior-interface}
Fix an integer \(r\geq6\).  Suppose that the auxiliary exterior chains
\(V_a^0\Subset\cdots\Subset V_a^5\) carry the hypotheses of
Lemmas~\ref{lem:buffered-local-Ricci-control} and
\ref{lem:buffered-Ricci-DeTurck-coefficients} through output order
\(r\), with one common physical-time-width margin.  Suppose also that
the exact normalized graph identity on
\(K_2\setminus\overline{K_-}\) supplies the corresponding
scale-adapted order-\(r\) initial geometry, while the coarse normalized
graph and radial-comparison bounds persist there.  Use the reference
fixed in the exterior certificate,
\(\widehat G_{\rm ext}:=g_{\rm ext}\), on the compact manifold with
boundary \(\overline{E^{++}}\).
The numerical package is chosen so that the strict time-width ratio in
\eqref{eq:auxiliary-buffered-time-width} is below the threshold in the
Dirichlet harmonic-map estimate below, and the common entrance time is
above its corresponding collar threshold.  These choices depend only
on the fixed package, not on either state or on a terminal time.

Then, throughout the remaining physical interval, there is one
fixed-reference exterior Ricci--DeTurck gauge
\[
 G_i(t)=\chi_i(t)^*\widetilde G_i(t)
 \quad\text{on }E^{++},\qquad
 \chi_i(t):E^{++}\longrightarrow E^{++}
 \ \text{a diffeomorphism},\qquad
 \chi_i|_{\partial E^{++}}=\operatorname{Id},
 \]
and its restrictions to all retained sets satisfy
\begin{equation}\label{eq:anchored-exterior-one-state-bounds}
 \sup_t\sum_{a=1}^{N_{\rm ext}}\left(
  \|\widetilde G_i(t)\|_{C_{R_a}^{r,\alpha}(U_a^{++})}
  +\|\chi_i(t)^{\pm1}\|_{C_{R_a}^{r-1,\alpha}(U_a^{++})}
 \right)\leq C_r .
\end{equation}
In the effective physical-time variables on the retained charts, the
Ricci--DeTurck coefficients have one common
 \(C^{\alpha/2}\) modulus (equivalently, the common parabolic
 time-oscillation modulus used in
 Lemma~\ref{lem:anchored-buffered-Davies}).
The same reference metric and boundary convention may be used for
every member of a common-margin prepared ball.  Hence two-state
differences and sliced first variations are tensorial in this one
gauge, rather than comparisons of independently chosen local gauges.

On the common marked collar
\(P_i=Q_i\circ\chi_i^{-1}\), and in local coordinates its normalized
equation is
\[
 \begin{split}
 \partial_\tau P_i^A={}&
 (P_i^*g_i)^{ab}\bigl(
  \partial_{ab}P_i^A-\widehat\Gamma^c_{ab}\partial_cP_i^A
  +\bar\Gamma^A_{BC}(P_i)\partial_aP_i^B\partial_bP_i^C\bigr)\\
 &+\bigl((1+a_i)\bar\nabla\bar f-U_i\bigr)^A(P_i),
 \end{split}
 \]
 where \(\widehat\Gamma\) is the fixed common exterior reference
 connection.  In particular, for \(4\leq m\leq r-2\),
\begin{equation}\label{eq:interface-block-controlled}
 \mathfrak I_{\rm in}\leq C\mathfrak J_{\rm in}.
\end{equation}
This is the full finite-horizon algebraic estimate.  For uniform
theory the exact identity for \(P\), together with the current global
low \(F\)-block propagated on each short interval from its actual
global left-endpoint trace, gives
\eqref{eq:interface-block-algebraic-closure}.  The graft source star is
not enlarged across the collapsing interface, and no arbitrary-time
same-collar parabolic restart is asserted.  The same estimates hold
for sliced first variations.  This
low-order derived block is retained precisely so that no high-order
boundary trace of a pullback metric is asserted.
\end{lemma}

\begin{proof}
Fix an arbitrary finite normalized endpoint \(S\) in the common
remaining interval and construct the gauge on
\([\tau_0,S]\).  Every constant below depends only on the common
package and is independent of this choice of \(S\).
For \(i=1,2\), set \(t_{i,0}:=t_i(\tau_0)\).
On each finite physical endpoint solve the single Dirichlet
harmonic-map problem
\[
 \partial_t\psi_i=\Delta_{G_i(t),\widehat G_{\rm ext}}\psi_i,\qquad
 \psi_i(\cdot,t_{i,0})=\operatorname{Id},\qquad
 \psi_i|_{\partial E^{++}}=\operatorname{Id},
\]
and set
\[
 \chi_i=\psi_i,\qquad
 \widetilde G_i=(\psi_i^{-1})^*G_i .
 \]
For the construction only, extend the target and the reference metric
smoothly across a doubled boundary collar and solve the map equation
into that extension.  The identity boundary condition and the strict
\(C^1\)-closeness-to-identity margin preserve the inward side of the
collar; side preservation is included in the same first-exit
alternative as local invertibility.  Thus the resulting map has image
in \(E^{++}\), rather than this self-map property being assumed in the
construction.
Only zeroth-order compatibility of the exact initial and lateral data
is automatic at the initial-boundary corner.  We use compatible
approximation without imposing any additional condition on the
physical data.  Let \(d_\partial\) be a fixed smooth boundary defining
function on the doubled collar.  For \(\varepsilon\downarrow0\), the
standard noncharacteristic boundary-jet recursion for the uniformly
parabolic map system gives a smooth initial map
\(\Phi_{i,\varepsilon}\) into the doubled target such that
\[
 \Phi_{i,\varepsilon}|_{\partial E^{++}}
 =\operatorname{Id},\qquad
 \operatorname{supp}
   (\Phi_{i,\varepsilon}-\operatorname{Id})
 \subset\{d_\partial<2\varepsilon\},\qquad
 \|\Phi_{i,\varepsilon}-\operatorname{Id}\|_
       {C^1(\overline{E^{++}})}
 \longrightarrow0,
\]
and the compatibility conditions obtained by differentiating the
fixed Dirichlet problem at \(t=t_{i,0}\) hold through order \(r\).
After the boundary value and first normal jet are fixed, uniform
parabolicity makes the coefficient of the next even normal derivative
invertible; the compatibility identities determine the required even
normal jets successively.  A fixed Borel extension in the collar
realizes those jets, and a shrinking collar cutoff gives the displayed
support and \(C^1\)-convergence.  The recursion and extension depend
smoothly on the finite source-coefficient jets, so the construction
also applies to two-state differences and sliced first variations.

Let \(\psi_{i,\varepsilon}\) solve the compatible problem
\[
 \partial_t\psi_{i,\varepsilon}
 =\Delta_{G_i(t),\widehat G_{\rm ext}}\psi_{i,\varepsilon},\qquad
 \psi_{i,\varepsilon}(\cdot,t_{i,0})
 =\Phi_{i,\varepsilon},\qquad
 \psi_{i,\varepsilon}|_{\partial E^{++}}
 =\operatorname{Id}.
\]
The boundary \(C^1\) estimate depends only on the common coefficient
package and the uniformly bounded \(C^1\) norms of
\(\Phi_{i,\varepsilon}\), and is therefore uniform in
\(\varepsilon\).  Every retained set has a fixed positive distance
from \(\partial E^{++}\); for sufficiently small \(\varepsilon\), the
initial map equals \(\operatorname{Id}\) on a neighborhood of that set,
and the interior estimates through order \(r-1\) are uniform there.
Boundary compactness for the spatial first derivatives and interior
Schauder compactness for the higher derivatives give a limit
\(\psi_i\).  The limit solves the original problem and satisfies
\[
 \psi_i(\cdot,t_{i,0})=\operatorname{Id},\qquad
 \psi_i|_{\partial E^{++}}=\operatorname{Id}.
\]
It inherits the stated global \(C^1\) and retained interior
\(C^{r-1,\alpha}\) bounds, and uniqueness removes subsequence
dependence.  The strict \(C^1\) side-preservation and
local-invertibility margins also pass to the limit, so the existing
degree-one argument yields the required self-diffeomorphism.

The auxiliary chains supply the initial-face and positive-time
Ricci--DeTurck estimates on the noncollapsing part of \(E^{++}\).
On the remaining marked collar, use
\[
 G_i=\lambda_iQ_i^*(\bar g+h_i)
\]
and the recorded scale-adapted collar charts.  On this fixed physical
collar \(Q_i\) lands in a receding annulus
\(1+\bar f\simeq e^\tau\), so its natural physical scale satisfies
\[
 \begin{aligned}
 r_{\rm in}(\tau)^2
 &\asymp\lambda_i(\tau)(1+\bar f\circ Q_i)
  \asymp\lambda_i(\tau)e^\tau,\\
 0\leq t_i(S)-t_i(\tau)
 &=\int_\tau^S\lambda_i(q)\,dq
  \leq C\lambda_i(\tau),\qquad \tau\leq S,\\
 C\lambda_i(\tau)&\ll r_{\rm in}(\tau)^2 .
 \end{aligned}
\]
Here \(S<\infty\) is the endpoint fixed at the start of the proof.
The estimate follows directly from the common scale bracket and its
constant is independent of \(S\); no as-yet-unconstructed terminal
time is used.
In these charts the source metric, boundary data, and map equation
have the same uniformly parabolic bounds.  The normalized collar
restart and the fixed exterior auxiliary-buffer restart therefore have
one common continuation alternative.  Standard scale-one Dirichlet
estimates, first for compatible approximants and then by passage to
the limit, give on a restart chart of physical scale \(r_{\rm loc}\),
for some fixed \(\vartheta>0\),
\[
 r_{\rm loc}^{-1}\|\psi_i-\psi_i(t_*)\|_{C^0}
 +\|d\psi_i-d\psi_i(t_*)\|_{C^0}
 \leq C\left(\frac{t-t_*}{r_{\rm loc}^2}\right)^\vartheta
\]
on each restart chart.  On the auxiliary exterior chains the ratio on the
right is small by the strict time-width certificate; on the normalized
collar it is \(O(e^{-\tau})\) by the preceding display.  Decreasing the
fixed time-width constant and increasing the entrance threshold make
 this estimate strictly smaller than the displacement and
 \(C^1\)-invertibility margins.  Thus those margins genuinely exclude a
 first exit.  Summing the \(C^1\) increments over the finitely many
 noncollapsing exterior restarts and the geometrically summable
 normalized-collar restarts keeps
 \(\|d\psi_i-\operatorname{Id}\|_{C^0}\) below one fixed constant
 \(<1\) on all of \(\overline{E^{++}}\).  Hence \(\psi_i\) is a local
 diffeomorphism.  It fixes the boundary and is homotopic relative to
 the boundary to the identity; as a proper local diffeomorphism of the
 compact manifold with boundary it is a covering.  Componentwise, its
 homotopy to the identity makes the induced fundamental-group map
 surjective, whereas the image of the fundamental group of a connected
 covering has index equal to its number of sheets.  Every component
 therefore has one sheet, so \(\psi_i\) is a global diffeomorphism of
 \(E^{++}\).  No orientation convention is needed.  The identity
 \[
  d\psi_i^{-1}=(d\psi_i\circ\psi_i^{-1})^{-1}
 \]
 gives the uniform inverse \(C^1\) margin.  Differentiating this
 identity and using the retained-set estimates gives all inverse bounds
 in \eqref{eq:anchored-exterior-one-state-bounds}.
 The same local Ricci and harmonic-map estimates propagate
 the initial order-\(r\) collar geometry, so no future order-\(r\) bound
 for \(h\) is being assumed here.  Uniqueness on
overlapping finite endpoints produces one gauge on the whole remaining
interval.  Applying the interior initial-face estimate from
Lemma~\ref{lem:buffered-Ricci-DeTurck-coefficients} to this single
solution on each auxiliary chain proves
 \eqref{eq:anchored-exterior-one-state-bounds}.  This is where the
auxiliary levels beyond \(U_a^{++}\) are used.

We next record the two-state gauge system rather than suppressing it.
Write
\[
 \widehat\chi_i(\tau)=\chi_i(t_i(\tau)),\qquad
 \mathcal T_i(\psi)=
 \Delta_{G_i(t_i(\tau)),\widehat G_{\rm ext}}\psi .
\]
Then, exactly,
\[
 \partial_\tau\widehat\chi_i
  =\lambda_i\mathcal T_i(\widehat\chi_i),\qquad
 \widehat\chi_i|_{\partial E^{++}}=\operatorname{Id}.
\]
For \(\delta\chi=\widehat\chi_1-\widehat\chi_2\), the mean-value
linearization in the map variable gives
\begin{equation}\label{eq:anchored-gauge-difference-equation}
\begin{split}
 \partial_\tau\delta\chi
 -\lambda_1\mathcal L_{\chi,12}(\tau)\delta\chi
 ={}&(\lambda_1-\lambda_2)
       \mathcal T_2(\widehat\chi_2)\\
 &+\lambda_1(\mathcal T_1-\mathcal T_2)
       [\widehat\chi_2],
 \qquad
 \delta\chi|_{\partial E^{++}}=0,
\end{split}
\end{equation}
where
\[
 \mathcal L_{\chi,12}
 =\int_0^1D_\psi\mathcal T_1
   (\widehat\chi_2+\vartheta\delta\chi)\,d\vartheta .
\]
For every
\[
 \tau_0\leq s\leq\tau\leq S,
\]
the coefficient difference in the last source term is a universal
linear combination of the metric and first metric-derivative
differences in the common exterior gauge.  Boundary Schauder estimates
on the compact collar, the auxiliary-chain interior estimates outside it,
and the clock estimate therefore give, through the orders used below,
\[
 \begin{split}
 &\sum_{a=1}^{N_{\rm ext}}
   \|\delta\chi(\tau)\|_{C_{R_a}^{m-1,\alpha}(U_a^+)}
  +\|\delta\chi(\tau)\|_
     {C_{\rm sc}^{m-1,\alpha}(\mathcal W_{\rm in}^+)}
 \\
 &\quad\leq C\left[
  \sum_{a=1}^{N_{\rm ext}}
   \|\delta\chi(s)\|_{C_{R_a}^{m-1,\alpha}(U_a^+)}
  +\|\delta\chi(s)\|_
     {C_{\rm sc}^{m-1,\alpha}(\mathcal W_{\rm in}^+)}
 \right]
 +C\int_s^\tau\!
    \left(\mathfrak D_m^{\rm hyb,full}
          +|\delta c|\right)dq .
 \end{split}
\]
with \(C\) depending only on the common package and not on
\(s,\tau\), or \(S\).
Differentiating
\(\widehat\chi_i\circ\widehat\chi_i^{-1}=\operatorname{Id}\)
gives the identical estimates, on the same listed domains and at the
same scales, for the inverse difference and then for
\(\widetilde\iota_1-\widetilde\iota_2\) on
\(\mathcal W_{\rm in}^+\).  The corresponding higher-order
graft-collar estimates are not inferred from this interface estimate;
they are supplied by the one-time buffered memory estimate below.
The same equations in effective physical time give a coefficient time
modulus uniform on the common-margin ball.

The domain DeTurck contribution in the tension field cancels the
\(-B_{\bar g}(g_i)\)-term in
 \eqref{eq:closed-normalized-inner-trace}, which gives the displayed
 \(P_i\)-equation.  More importantly for the uniform estimate,
\[
 P_i=\varphi_\tau\circ R_i\circ F_i\circ\iota\circ\chi_i^{-1}
 \quad\text{on }\mathcal A_{\rm in}.
 \]
The scale-adapted composition estimate applied directly to
\[
 P_i=\varphi_\tau\circ R_i\circ F_i\circ\iota\circ\chi_i^{-1}
 \quad\text{on }\mathcal A_{\rm in}
\]
and \eqref{eq:hybrid-interface-algebraic-data} gives
\eqref{eq:interface-block-controlled}.  This identity, rather than a
localized initial-boundary problem for the displayed \(P\)-equation,
also supplies every lateral value used later.  Differentiating the
algebraic identity and the anchored gauge equation gives the
first-variation statement with the same derivative count.
\end{proof}

All differences below are taken in the fixed finite atlases and in
right-translated exponential charts.  At normalized time \(\tau\), every
occurrence of \(\chi_i\), \(\chi_i^{-1}\), and
\(\widetilde\iota_i\) in the following blocks means, respectively,
\(\chi_i(t_i(\tau))\), \(\chi_i(t_i(\tau))^{-1}\), and
\(\widetilde\iota_i(t_i(\tau))\).  Define
\begin{equation}\label{eq:hybrid-physical-block}
 \mathfrak G_m(\tau)
 := 
 \sum_{a=1}^{N_{\rm ext}}
 \|\widetilde G_1(t_1(\tau))-\widetilde G_2(t_2(\tau))\|
       _{C_{R_a}^{m,\alpha}(U_a)},
\end{equation}
and, for \(m\geq4\), the order-dependent gauge and marking block
\begin{equation}\label{eq:hybrid-marking-block}
 \begin{split}
 \mathfrak M_{{\rm gr},m}(\tau):={}&
 \sum_{a=1}^{N_{\rm ext}}\left(
  \|\chi_1-\chi_2\|_{C_{R_a}^{m-1,\alpha}(U_a^+)}
  +\|\chi_1^{-1}-\chi_2^{-1}\|_{C_{R_a}^{m-1,\alpha}(U_a^+)}
  \right)\\
  &+\|\chi_1-\chi_2\|
       _{C_{\rm sc}^{m-1,\alpha}(\mathcal W_{\rm in}^+)}
   +\|\chi_1^{-1}-\chi_2^{-1}\|
       _{C_{\rm sc}^{m-1,\alpha}(\mathcal W_{\rm in}^+)}\\
  &+\|\widetilde\iota_1-\widetilde\iota_2\|
       _{C^{m-1,\alpha}(\mathcal W_{\rm in}^+)}.
 \end{split}
\end{equation}
We retain the abbreviation
\(\mathfrak M_{\rm gr}:=\mathfrak M_{{\rm gr},4}\) only in arguments
whose output order is \(m_\#=4\).
The compact physical graft costs two input derivatives.  We therefore
record a metric input on the larger graft collar and map outputs on the
strictly smaller collar:
\begin{equation}\label{eq:hybrid-graft-buffer-block}
\begin{split}
 \mathfrak B_{{\rm gr},m}(\tau):={}&
 \|\widetilde G_1-\widetilde G_2\|
   _{C_{\rm sc}^{m+2,\alpha}(\mathcal W_{\rm gr}^+)}
 +\|\chi_1-\chi_2\|
   _{C_{\rm sc}^{m+1,\alpha}(\mathcal W_{\rm gr})}\\
 &+\|\chi_1^{-1}-\chi_2^{-1}\|
   _{C_{\rm sc}^{m+1,\alpha}(\mathcal W_{\rm gr})}
 +\|\widetilde\iota_1-\widetilde\iota_2\|
   _{C^{m+1,\alpha}(\mathcal W_{\rm gr})} .
\end{split}
\end{equation}
This is a graft-input block, not a global output claim on the
noncompact end.  Its tensor component is already in the single anchored
exterior gauge.  The gauge, inverse-gauge, and marking components are
then propagated one order below by their triangular equations.  The
larger tensor collar and the unused chain in
\eqref{eq:separated-graft-collars} supply the lateral memory which a
same-collar parabolic restart would otherwise miss.
For later uniform estimates fix nested finite source-adapted stars
\[
 \mathfrak S_{\rm gr}\Subset\mathfrak S_{\rm gr}^+
 \Subset\mathfrak S_{\rm gr}^{++},
\]
where \(\mathfrak S_{\rm gr}\) contains every chart whose physical
image meets \(\mathcal W_{\rm gr}^{++}\), and successive stars are
obtained by adding one layer of the fixed source atlas.  The star is
deliberately not enlarged across the collapsing interface collar:
that would destroy the \(\lambda\)-weighted noncollapsing estimate
used below.  Define the
localized source-adapted map block
\begin{equation}\label{eq:localized-graft-F-block}
 \mathfrak F_{{\rm gr},m}(\tau)
 :=
 \|F_1(\tau)-F_2(\tau)\|_
 {\mathfrak X_{\rm sc}^{m,\alpha}(\mathfrak S_{\rm gr})}.
\end{equation}
For comparison, the notation
\begin{equation}\label{eq:global-low-F-block}
 \mathfrak F_m^{\rm glob}(\tau)
 :=
 \|F_1(\tau)-F_2(\tau)\|_
 {\mathfrak X_{\rm sc}^{m,\alpha}}
\end{equation}
will be used only for a coarse global map estimate; it is not a
summand of the localized uniform hybrid distance below.  The notation
\(\mathfrak F_m^{\rm glob}[\xi](\tau)\) denotes its linearized
seminorm.
Its typed larger-star entrance trace is
\begin{equation}\label{eq:localized-graft-F-initial-buffer}
 d_{{\rm Fgr},m,0}^{++}
   (\mathbf z_{1,0},\mathbf z_{2,0})
 :=
 \|F_1(\tau_0)-F_2(\tau_0)\|_
 {\mathfrak X_{\rm sc}^{m,\alpha}(\mathfrak S_{\rm gr}^{++})}.
\end{equation}
The genuinely global low initial trace is
\begin{equation}\label{eq:global-F-initial-trace}
 d_{{\rm F},m,0}^{\rm glob}
   (\mathbf z_{1,0},\mathbf z_{2,0})
 :=
 \|F_1(\tau_0)-F_2(\tau_0)\|_
 {\mathfrak X_{\rm sc}^{m,\alpha}} .
\end{equation}
We suppress its two entrance arguments when the pair is fixed, and
use the bracket notation
\(d_{{\rm Fgr},m,0}^{++}[\xi]\) and
\(d_{{\rm F},m,0}^{\rm glob}[\xi]\) for the corresponding linearized
seminorms.
Thus \(\mathfrak F_{{\rm gr},m}\) records exactly the portion of \(F\)
which enters the graft pullback; it is dominated by the global
\(\mathfrak X_{\rm sc}^{m,\alpha}\) norm on finite horizons, but no
uniform future global \(F\)-estimate is inferred from a local graft
argument.  The global norm in
\eqref{eq:global-F-initial-trace} is retained at the original entrance
face, where it is controlled by the prepared entrance norm; subsequent
short intervals restart the global equation from its actual global
trace.
With these abbreviations, define the \emph{localized uniform hybrid
distance}
\begin{equation}\label{eq:two-state-hybrid-distance}
 \begin{split}
 \mathfrak D_m^{\rm hyb}(\tau):={}&
 |t_1-t_2|
 +\left|\log\frac{\lambda_1}{\lambda_2}\right|
 +\mathfrak G_m(\tau)+\mathfrak M_{{\rm gr},m}(\tau)
  +\mathfrak B_{{\rm gr},m}(\tau)\\
  &+\|R_1-R_2\|_{\mathfrak X_{\rm sc}^{m+1,\alpha}}
  +\mathfrak F_{{\rm gr},m+1}
  +\|h_1-h_2\|_{\mathfrak T_{{\rm sc},N}^{m,\alpha}}
  +\mathfrak I_{\rm in}.
\end{split}
\end{equation}
It differs from the full finite-horizon distance
\(\mathfrak D_m^{\rm hyb,full}\) defined above only by replacing the
global \(F\)-norm with \(\mathfrak F_{{\rm gr},m+1}\).  Therefore
\begin{equation}\label{eq:full-local-hybrid-comparison}
 \mathfrak D_m^{\rm hyb}
 \leq C\mathfrak D_m^{\rm hyb,full},\qquad
 \mathfrak D_m^{\rm hyb}(\tau_0)
 +\mathfrak D_m^{\rm hyb,full}(\tau_0)
 \leq C\|\mathbf z_{1,0}-\mathbf z_{2,0}\|
 _{\mathscr X_{\rm prep}^{m+2,\alpha}} .
\end{equation}
No reverse comparison is asserted at future times.  Let
\(\mathfrak D_m^{\rm full,base}\) denote the full finite-horizon sum
with the final \(\mathfrak I_{\rm in}\)-summand omitted, and put
\[
 \mathfrak D_m^{\rm red}
 :=\mathfrak D_m^{\rm hyb}-\mathfrak I_{\rm in}.
\]
The exact identity
\[
 P_i=\varphi_\tau\circ R_i\circ F_i\circ
       \iota\circ\chi_i^{-1}
 \quad\text{on }\mathcal A_{\rm in}
\]
and the scale-adapted composition calculus then give
\begin{equation}\label{eq:interface-block-algebraic-closure}
 \mathfrak J_{\rm in}\leq C\mathfrak D_m^{\rm full,base},
 \qquad
 \mathfrak I_{\rm in}
 \leq C\left(
   \mathfrak D_m^{\rm red}+\mathfrak F_{m+1}^{\rm glob}\right)
 \leq C\mathfrak D_m^{\rm full,base}.
\end{equation}
Thus the interface block records the terminating trace explicitly while
remaining controlled by \(\mathfrak D_m^{\rm full,base}\).  The global
map block in the middle expression is restarted only as a global
parabolic unknown from its global trace; it is neither restarted from
nor replaced by an interface-collar trace.
The full unscaled \(C^{m,\alpha}(\mathcal X)\) norm of the closed metric
is intentionally absent.  The collapsing marked core is measured by
the normalized graph tensor \(h\); the closed metric is measured only
on the fixed buffered noncollapsing physical sets.  The ordinary
marking block is deliberately measured at order \(m-1\), rather than
at one fixed low order.  This is the order supplied by the triangular
DeTurck ODE from an order-\(m\) metric difference and is sufficient
away from the graft.  On the compact graft collar,
\(\mathfrak B_{{\rm gr},m}\) uses the two derivatives already present
in the prepared input to evaluate the complete second-order graft
operator as an ordinary \(C^{m-2,\alpha}\) source.  Together the two
blocks keep the exact covariance
\((\widetilde\iota_i)_*\widetilde G_i=\iota_*G_i\) visible and ensure
that the moving DeTurck gauge remains separate from the fixed graft
cutoff.

\begin{lemma}[Anchored buffered off-diagonal estimate]
\label{lem:anchored-buffered-Davies}
Fix an integer \(m\geq0\), \(0<\alpha<1\), and an interval
\[
 I=[\tau_*,S),\qquad \tau_*<S\leq\infty .
\]
Work on a fixed Riemannian domain \((\mathcal D,g_*)\), either complete
without boundary or with smooth boundary and homogeneous Dirichlet
condition, and on a finite-rank metric bundle \(E\to\mathcal D\) with
compatible connection and fixed positive smooth density \(d\mu_*\).
For this lemma put
\[
 V_\partial(\mathcal D;E):=
 \begin{cases}
  H^1(\mathcal D;E),&\partial\mathcal D=\varnothing,\\
  H^1_0(\mathcal D;E),&\text{homogeneous Dirichlet boundary},
 \end{cases}
 \qquad
 H^{-1}_\partial(\mathcal D;E):=
 V_\partial(\mathcal D;E)^*,
\]
with duality induced by the \(L^2(d\mu_*)\) pairing, and abbreviate
\(H^{-1}:=H^{-1}_\partial(\mathcal D;E)\) within this lemma.

Let \(W\Subset W^+\Subset\mathcal D\) satisfy
\[
 \operatorname{dist}_{g_*}
 \bigl(\overline W,\mathcal D\setminus W^+\bigr)\geq4R
\]
for some \(R>0\).  The pair may be contained in the anchored physical
domain \(E^{++}\), or may be a pair of fixed source-adapted stars with
the same recorded bounded-geometry certificate.  Let
\(\lambda:I\to(0,\infty)\) be locally integrable, and let
\(\mathcal U(\tau,q)\) be the evolution family of
\[
 \partial_\tau-\lambda(\tau)\mathcal L_\tau,\qquad
 \mathcal L_\tau
 =A_\tau^{ab}\nabla_a\nabla_b+B_\tau*\nabla+C_\tau .
\]
Assume that this evolution and its adjoint are realized on the form
scale
\[
 V_\partial\subset L^2(d\mu_*)\subset H^{-1}_\partial
\]
and extend consistently to initial data in \(H^{-1}_\partial\).
Assume that \(A_\tau\) is a real symmetric scalar bundle symbol,
uniformly elliptic relative to \(g_*\), and that, in the scale-one
buffered source atlas,
\[
 \|A\|_{C^{m,\alpha}}
 +\|B\|_{C^{(m-1)_+,\alpha}}
 +\|C\|_{C^{(m-2)_+,\alpha}}
 \leq\Lambda_m ,
\]
with the same spatial jets carrying one common effective-time
\(C^{\alpha/2}\) modulus.  Here \(j_+=\max\{j,0\}\).
Assume uniform forward and adjoint Caccioppoli estimates on every
buffered chart meeting the admissible source region.  Using the
\(4R\) buffer, choose the smoothing chain so that
\[
 W\Subset W_1\Subset W_2
 \Subset\mathcal N_{R/2}^{g_*}(W)\Subset W^+,
\]
and assume uniform interior parabolic H\"older estimates on this
chain.  Thus a source separated from \(W\) by \(R\) is separated from
\(W_2\) by at least \(R/2\).
In the geometric applications the second derivatives of the
background metric occur only in the zeroth-order coefficient \(C\);
consequently a one-state metric package through the displayed output
order \(m\) supplies exactly these coefficient bounds.

Require the following forward and adjoint form inequalities.  There
exist \(\kappa>0\) and a locally integrable function
\(\beta:I\to\mathbb R\) such that, for every \(v\in V_\partial\),
\begin{equation}\label{eq:anchored-Davies-energy}
\begin{aligned}
 2\operatorname{Re}
 \langle\lambda(\tau)\mathcal L_\tau v,v\rangle_
        {H^{-1}_\partial,V_\partial}
 &\leq
 -\kappa\lambda(\tau)\|\nabla v\|_{L^2}^2
 +2\beta(\tau)\|v\|_{L^2}^2,\\
 2\operatorname{Re}
 \langle\lambda(\tau)\mathcal L_\tau^\dagger v,v\rangle_
        {H^{-1}_\partial,V_\partial}
 &\leq
 -\kappa\lambda(\tau)\|\nabla v\|_{L^2}^2
 +2\beta(\tau)\|v\|_{L^2}^2 .
\end{aligned}
\end{equation}
Put
\[
 \mathfrak s(q,\tau):=\int_q^\tau\lambda(r)\,dr
\]
and assume the uniform bounds
\begin{equation}\label{eq:anchored-Davies-clock-growth}
 \sup_{\tau_*\leq q\leq\tau<S}\mathfrak s(q,\tau)
 \leq\Theta_0<\infty,
 \qquad
 \sup_{\tau_*\leq q\leq\tau<S}
 \int_q^\tau\beta_+(r)\,dr
 \leq B_0<\infty .
\end{equation}
The constants below are allowed to depend on
\(\Theta_0\) and \(B_0\), but not on \(s,q,\tau\), or \(S\).

Fix
\[
 \tau_*\leq s\leq\tau<S.
\]
If
\(\mathcal F\in L^1_{\rm loc}(I;H^{-1}_\partial)\) is strongly
measurable and
\[
 \operatorname{dist}_{g_*}(\supp\mathcal F(q),W)\geq R
\]
for almost every \(q\in[s,\tau]\), then there are \(C_m,c_m>0\) and an integer
\(q_m\geq0\), depending only on the displayed package, such that
\begin{equation}\label{eq:anchored-Davies-Duhamel}
\begin{aligned}
 \left\|\int_s^\tau
  \mathcal U(\tau,q)\mathcal F(q)\,dq
 \right\|_{C^{m,\alpha}(W)}
 &\leq C_m\int_s^\tau
  e^{B_0+C_m\mathfrak s(q,\tau)}
  \bigl(1+\mathfrak s(q,\tau)^{-q_m}\bigr)
  e^{-c_mR^2/\mathfrak s(q,\tau)}
  \|\mathcal F(q)\|_{H^{-1}_\partial}\,dq\\
 &\leq C_{m,B_0,\Theta_0}\int_s^\tau
  \bigl(1+\mathfrak s(q,\tau)^{-q_m}\bigr)
  e^{-c_mR^2/\mathfrak s(q,\tau)}
  \|\mathcal F(q)\|_{H^{-1}_\partial}\,dq .
\end{aligned}
\end{equation}
Every scalar kernel in this statement is defined to be zero at zero
effective time.  If \(g\in H^{-1}_\partial\) is supported a distance
at least \(R\) from \(W\), then
\[
\begin{aligned}
 \|\mathcal U(\tau,s)g\|_{C^{m,\alpha}(W)}
 &\leq
 C_m e^{B_0+C_m\mathfrak s(s,\tau)}
 \bigl(1+\mathfrak s(s,\tau)^{-q_m}\bigr)
 e^{-c_mR^2/\mathfrak s(s,\tau)}
 \|g\|_{H^{-1}_\partial}\\
 &\leq
 C_{m,B_0,\Theta_0}
 \bigl(1+\mathfrak s(s,\tau)^{-q_m}\bigr)
 e^{-c_mR^2/\mathfrak s(s,\tau)}
 \|g\|_{H^{-1}_\partial}.
\end{aligned}
\]
The same assertion holds for a cutoff commutator
\[
 \mathcal F=\nabla_iP^i+Q,\qquad P,Q\in L^2,
\]
provided \(P,Q\) have the same support separation; in that case
\[
 \|\mathcal F\|_{H^{-1}_\partial}
 \leq C\bigl(\|P\|_{L^2}+\|Q\|_{L^2}\bigr).
\]
\end{lemma}

\begin{proof}
For a homogeneous datum supported a distance at least \(R\) from
\(W\), the refined chain gives distance at least \(R/2\) from
\(W_2\).  Choose a smooth real weight \(\phi\) which is zero on the
datum support, equals \(R/4\) on a neighborhood of \(W_2\), and satisfies
\[
 |\nabla\phi|+R|\nabla^2\phi|\leq C
\]
in the fixed bounded-geometry atlas.  The Dirichlet form domain is
preserved by multiplication by \(e^{\vartheta\phi}\).  Apply
\eqref{eq:anchored-Davies-energy} to
\(e^{\vartheta\phi}\mathcal U(\tau,q)g\).  The scalar principal symbol
absorbs the bundle-valued first-order cross term, and the remaining
weight commutators give
\[
 \frac d{d\tau}
 \|e^{\vartheta\phi}\mathcal U(\tau,q)g\|_{L^2}^2
 \leq
 \bigl(
  2\beta_+(\tau)+C(1+\vartheta^2)\lambda(\tau)
 \bigr)
 \|e^{\vartheta\phi}\mathcal U(\tau,q)g\|_{L^2}^2 .
\]
Gronwall and the choice
\(\vartheta=cR/\mathfrak s(q,\tau)\) give
\[
 \|\mathcal U(\tau,q)g\|_{L^2(W_2)}
 \leq
 C e^{B_0+C\mathfrak s(q,\tau)}
 e^{-cR^2/\mathfrak s(q,\tau)}
 \|g\|_{L^2}.
\]
For effective times comparable to or larger than \(R^2\), the same
bound follows by taking \(\vartheta=0\) and decreasing \(c\); hence the
constant is uniform over the full interval in
\eqref{eq:anchored-Davies-clock-growth}.

Apply the identical weighted form estimate to the adjoint evolution.
Duality on
\(V_\partial\subset L^2\subset H^{-1}_\partial\), followed by the
adjoint Caccioppoli estimate on the source side, gives the corresponding
\(H^{-1}_\partial\)-to-\(L^2(W_2)\) estimate.  Interior parabolic
H\"older estimates first on \(W_1\) and then on \(W\) give the
\(C^{m,\alpha}\) output and introduce only the finite factor
\(1+\mathfrak s^{-q_m}\).  The bounds in
\eqref{eq:anchored-Davies-clock-growth} absorb the displayed Gronwall
factor into \(C_{m,B_0,\Theta_0}\).  This proves the homogeneous
estimate, and integration at each source time proves
\eqref{eq:anchored-Davies-Duhamel}.  Finally,
\(\nabla_iP^i+Q\) acts continuously on \(V_\partial\), which proves
the commutator clause.
\end{proof}

\begin{lemma}[Endpoint effective-time maximal regularity]
\label{lem:effective-time-endpoint-maximal-regularity}
Let \(m\geq0\), and let \(\mathcal U(\tau,q)\) be an evolution family
with the scalar-principal-symbol, bounded-geometry, and effective-time
coefficient package of
Lemma~\ref{lem:anchored-buffered-Davies}, including its fixed form
realization and forward/adjoint energy package.  For every two-derivative
conclusion in this lemma, require the full \(q=m+2\) Schauder package:
in every retained scale-one chart the rescaled operator is uniformly
parabolic and
\begin{equation}\label{eq:endpoint-MR-full-coefficient-ledger}
 A,\ B,\ C\in C^{m,\alpha},
\end{equation}
with one uniform bound and one common effective-time
\(C^{\alpha/2}\) modulus in those spatial norms.  The anchored estimate
at output order \(m+2\) is a stronger sufficient package and is the
one supplied in the graft application below; invoking the anchored
estimate at index \(m\) alone is not sufficient for the
\(C^{m+2,\alpha}\) output.
The same full coefficient assumptions are used for both the global
one-order kernel \eqref{eq:effective-time-one-order-kernel} and its
short source-atlas form below; no staggered-coefficient variant is
invoked.

For the H\"older conclusions, an output region is understood either in
a boundaryless realization (in particular, on the fixed closed doubles
used below) or on an interior region with one fixed positive buffer from
the boundary of the form realization.  The Dirichlet realization in
Lemma~\ref{lem:anchored-buffered-Davies} is sufficient for the energy
and off-diagonal estimates but is not, by itself, being promoted to an
up-to-the-boundary endpoint Schauder assertion.  An up-to-boundary
version is valid under a uniform \(C^{m+2,\alpha}\) boundary atlas,
homogeneous Dirichlet data, and the standard initial--boundary
compatibility conditions through the asserted order.  No application
of this lemma below uses that additional version.

Put
 \(\mathfrak s(q,\tau)=\int_q^\tau\lambda(r)\,dr\).
 Fix an initial time
\[
 \tau_*\leq s<S
\]
and use the inherited half-open interval \([s,S)\).  Assume that the
total effective time
satisfies
\[
 \sup_{s\leq\tau<S}\mathfrak s(s,\tau)
 \leq\Theta_0<\infty .
\]
If \(u(s)=0\) and
\[
 \partial_\tau u-\lambda(\tau)\mathcal L_\tau u
   =\lambda(\tau)F,
 \qquad
 F\in L^\infty\!\left(
       [s,S);C^{m,\alpha}_{\rm sc}\right),
\]
then, on every retained compact collar or uniformly buffered
source-adapted chart,
\begin{equation}\label{eq:effective-time-two-order-endpoint}
 \sup_{s\leq\tau<S}
 \|u(\tau)\|_{C_{\rm sc}^{m+2,\alpha}}
 \leq C_{\Theta_0}
 \operatorname*{ess\,sup}_{s\leq q<S}
 \|F(q)\|_{C_{\rm sc}^{m,\alpha}} .
\end{equation}
Precisely, the closed-interval estimate is first applied on every
\([s,\widehat S]\) with \(s<\widehat S<S\); its constant is independent
of \(\widehat S\), and the displayed assertion follows by
\(\widehat S\uparrow S\).  If \(S<\infty\) and the coefficients,
forcing, evolution, and solution extend through \(S\), the same
formula holds with both half-open inequalities replaced by closed
ones.
The same zero-trace contribution obeys the one-order estimate
\begin{equation}\label{eq:effective-time-one-order-kernel}
 \|u(\tau)\|_{C_{\rm sc}^{m+1,\alpha}}
 \leq C\int_s^\tau
  \lambda(q)\bigl(1+\mathfrak s(q,\tau)^{-1/2}\bigr)
  \|F(q)\|_{C_{\rm sc}^{m,\alpha}}\,dq .
\end{equation}
The last display is an improper, equivalently Lebesgue, integral at
\(q=\tau\); the value assigned to the integrand at that single endpoint
is irrelevant, and the \(\mathfrak s^{-1/2}\) singularity is integrable.

In the full source-adapted atlas, the corresponding statement is the
following short-interval form.  There is
\(\delta_0>0\), depending only on the fixed source-atlas,
coefficient, and effective-clock package, with the following property.
Let
\[
 J=[\sigma,\sigma+\delta]\Subset[s,S),
 \qquad 0<\delta\leq\delta_0,
\]
and let \(J_{\mathcal U}\) denote the image of \(J\) in the displayed
source-adapted clock of a retained chart.  Then
\(|J_{\mathcal U}|\leq C\delta\), uniformly in the chart,
\(\sigma\), \(\delta\), and the terminal endpoint.

Under the source-atlas coefficient bounds
\[
 A\in C^{m+1,\alpha},\qquad
 B\in C^{m,\alpha},\qquad
 C\in C^{m,\alpha},
\]
the chartwise version of
\eqref{eq:effective-time-one-order-kernel}, with
\(\Theta_{\mathcal U}(q,\tau)\) equal to elapsed displayed time in the
source-adapted chart, gives every solution with \(u(\sigma)=0\)
\begin{equation}\label{eq:source-atlas-one-order-Abel-block}
 \sup_{\tau\in J}
 \|u(\tau)\|_{\mathfrak X_{\rm sc}^{m+1,\alpha}}
 \leq C\delta^{1/2}
       \|F\|_{\mathbb F_{\rm sc}^{m,\alpha}(J)}.
\end{equation}
Indeed,
\[
 \int_0^\Theta(\Theta-r)^{-1/2}\,dr=2\sqrt\Theta,
\]
uniformly over the atlas.

Under the same full \(q=m+2\) package
\eqref{eq:endpoint-MR-full-coefficient-ledger}, after
decreasing \(\delta_0\) if necessary, the same interval \(J\) and
zero-trace condition give the stronger two-order form
\begin{equation}\label{eq:source-atlas-zero-trace-maximal-block}
 \|u\|_{\mathbb E_{\rm sc}^{m+2,\alpha}(J)}
 \leq C\|F\|_{\mathbb F_{\rm sc}^{m,\alpha}(J)},
 \qquad
 \sup_{\tau\in J}
 \|u(\tau)\|_{\mathfrak X_{\rm sc}^{m+1,\alpha}}
 \leq C\delta^{\alpha/4}
       \|F\|_{\mathbb F_{\rm sc}^{m,\alpha}(J)}.
\end{equation}
The constants in both displays are uniform over all such \(J\), all
retained source-adapted charts, and all terminal subintervals.
All these estimates hold in the weighted and tensor-normalized spaces.
For a parameter \(p\) in a Banach space there are two admissible
parameter routes.  In the coefficient route, the rescaled
coefficients \((A_p,B_p,C_p)\) depend \(C^1\) on \(p\) into the
atlas-supremum strong-Bochner product
\[
 L^\infty([s,S);C_x^{m,\alpha})^3
\]
corresponding to
\eqref{eq:endpoint-MR-full-coefficient-ledger}, and \(F_p\) depends
\(C^1\) into the corresponding strong-Bochner forcing space
\(L^\infty([s,S);C_x^{m,\alpha})\).  The coefficient family and its
first parameter derivative have one common effective-time modulus
and, whenever the little-H\"older clause is used, one common spatial
high-frequency tail.  At each fixed base parameter the derivatives
are continuous and the Fr\'echet remainders are little-oh in these
exact Bochner norms, uniformly on terminal subintervals
\([s,\widehat S]\).  Uniformity with respect to the base parameter is
asserted only on a compact parameter subset, or on a set on which a
common modulus of continuity for the first derivatives and a common
little-oh remainder modulus are assumed explicitly; boundedness of a
subset of an infinite-dimensional Banach space alone is not used for
that conclusion.  If the clock varies, fix one reference parameter
\(p_0\) and require
\[
 p\longmapsto\log(\lambda_p/\lambda_{p_0})
\]
to be \(C^1\) in \(L^\infty([s,S))\).  The clocks are positive and
locally integrable and have a common effective-time bound.

The actual initial trace may vary: assume
\[
 u_p(s)=u_{s,p},\qquad
 p\longmapsto u_{s,p}\in\mathfrak X_{\rm sc}^{m+2,\alpha}
\]
is \(C^1\).  The zero trace is the special case \(u_{s,p}=0\).  Then,
in the coefficient route, \(v=D_pu[\dot p]\) has initial trace
\(v(s)=D_pu_s[\dot p]\) and satisfies
\[
 \begin{split}
 \partial_\tau v-\lambda\mathcal L_\tau v
 =\lambda\bigl(&D_pF[\dot p]
 +D_pA[\dot p]*\bar\nabla^2u
 +D_pB[\dot p]*\bar\nabla u
 +D_pC[\dot p]*u\\
 &+D_p(\log\lambda)[\dot p]\,
       (F+\mathcal L_\tau u)\bigr).
 \end{split}
\]
The displayed maximal-regularity and Abel bounds apply to the
zero-trace remainder obtained by subtracting the fixed bounded trace
extension constructed in the proof of
Lemma~\ref{lem:weighted-prepared-Schauder}; the trace extension and its
induced forcing are added separately.  A two-state difference has the
corresponding secant trace and source.

In the residual-source route no parameter differentiability of
\((A_p,B_p,C_p)\) is required.  Fix a base parameter \(p_0\), retain
only the fixed operator \(\mathcal L_{p_0,\tau}\) satisfying
\eqref{eq:endpoint-MR-full-coefficient-ledger}, and suppose the exact
family can be written on a neighborhood of \(p_0\) as
\begin{equation}\label{eq:endpoint-MR-fixed-operator-residual-route}
 \partial_\tau u_p-\lambda_p\mathcal L_{p_0,\tau}u_p
 =\lambda_p F_p^{(p_0)},
\end{equation}
where \(p\mapsto F_p^{(p_0)}\) is \(C^1\) in the same
strong-Bochner forcing norm, locally at \(p_0\), and has the same
terminal-subinterval uniformity.  Then
\begin{equation}\label{eq:endpoint-MR-fixed-operator-variation}
 \partial_\tau v-\lambda_{p_0}\mathcal L_{p_0,\tau}v
 =\lambda_{p_0}\left(
   D_pF^{(p_0)}_{p_0}[\dot p]
   +D_p\log\lambda_{p_0}[\dot p]\,
      \bigl(F^{(p_0)}_{p_0}+\mathcal L_{p_0,\tau}u_{p_0}\bigr)
  \right),
\end{equation}
with the same differentiated trace.  Subtracting this equation from
the secant equation leaves a little-oh source in the exact
strong-Bochner \(C^{m,\alpha}\) forcing norm.  The preceding
local/compact-modulus convention applies to both routes.  Thus the
endpoint estimate applies either to a differentiable coefficient
family or to a differentiable exact residual; in the latter route no
top-order derivative of the coefficient family is required.
An ordinary \(L^1\)-in-time same-order Duhamel estimate is reserved for
sources already smooth at the asserted output order; it cannot replace
\eqref{eq:effective-time-two-order-endpoint} at a two-derivative
endpoint.
\end{lemma}

\begin{proof}
Fix first \(s<\widehat S<S\) and change variables to effective time
\[
 \vartheta(\tau)=\int_s^\tau\lambda(r)\,dr,
 \qquad0\leq\vartheta\leq\widehat\Theta\leq\Theta_0.
\]
We verify the local variable-coefficient estimate, including its
realization.  On a uniformly buffered source-adapted chart, with the
fixed bundle trivialization and tensor normalization, the equation has
the scale-one form
\begin{equation}\label{eq:endpoint-MR-local-equation}
 \partial_\vartheta u
 -A^{ij}(\vartheta,x)\partial_i\partial_j u
 -B^i(\vartheta,x)\partial_i u-C(\vartheta,x)u=F.
\end{equation}
The principal symbol is scalar on the fibers and uniformly elliptic,
and \eqref{eq:endpoint-MR-full-coefficient-ledger} is exactly the
coefficient hypothesis of
Lemma~\ref{lem:weighted-prepared-Schauder} with \(q=m+2\).
That lemma freezes the principal coefficient only in space,
\[
 A_0^{ij}(\vartheta):=A^{ij}(\vartheta,x_0),
\]
so its componentwise Fourier evolution has covariance
\(\int_r^\vartheta A_0(\zeta)\,d\zeta\) and obeys the uniform dyadic
heat estimate even though \(A_0\) depends on time.  Its proof treats
explicitly
\[
 (A-A_0)D^2(\zeta u),\qquad
 [\Delta_j,A]D^2(\zeta u),\qquad
 [\mathcal L_\vartheta,\zeta]u,
\]
all connection and \(B,C\) terms, and every overlap transition by the
paraproduct, interpolation, and bounded-overlap estimates in the proof
of \eqref{eq:weighted-prepared-Schauder}.  After the spatial
localization radius is fixed, the resulting parametrix error has the
form
\[
 \|\mathcal E\|
 \leq\varepsilon_{\rm sp}
      +C_{\rm fr}\omega(\delta_\vartheta)
      +\varepsilon_{\rm com}(\delta_\vartheta),
 \qquad
 \varepsilon_{\rm com}(r)\longrightarrow0,
\]
where \(\omega\) is the common effective-time coefficient modulus.
Choose \(\varepsilon_{\rm sp}<1/6\) and then
\(\vartheta_*>0\) so that the other two terms sum to less than
\(1/3\).  Neumann inversion gives, on every effective interval
\(J_\vartheta\) of length at most \(\vartheta_*\),
\begin{equation}\label{eq:endpoint-MR-local-two-order}
 \|u\|_{L^\infty C^{m+2,\alpha}_{\rm sc}}
 +\|\partial_\vartheta u\|_{L^\infty C^{m,\alpha}_{\rm sc}}
 \leq C_0\left(
  \|u(\vartheta_0)\|_{C^{m+2,\alpha}_{\rm sc}}
  +\|F\|_{L^\infty C^{m,\alpha}_{\rm sc}}
 \right).
\end{equation}
The constants depend only on the displayed coefficient, buffer, and
bounded-geometry package.  Thus spatially variable principal
coefficients, localization commutators, bundle transitions, weights,
and tensor normalizations are all included.

On a fixed compact collar used below, the equation is first extended
to the recorded fixed closed double and the same finite-atlas argument
applies.  On an interior output region in a realization with boundary,
all cutoffs are supported inside the fixed boundary buffer, so the same
interior parametrix applies and its cutoff commutator is already among
the displayed terms.  No unknown boundary trace is introduced.  The
optional up-to-boundary version is the compatible Dirichlet Schauder
realization stated in the boundary convention of the lemma.

Partition \([0,\widehat\Theta]\) into at most
\[
 N\leq1+\left\lceil\Theta_0/\vartheta_*\right\rceil
\]
successive effective intervals.  Iterating
\eqref{eq:endpoint-MR-local-two-order}, beginning with \(u(s)=0\),
proves \eqref{eq:effective-time-two-order-endpoint} with a constant
independent of \(\widehat S\).  Letting \(\widehat S\uparrow S\)
proves the half-open assertion; the same last local estimate proves
the closed-endpoint version when the data extend through \(S\).

We derive the one-order kernel from the preceding two-order estimate
under these full coefficient assumptions, without introducing a
separate parametrix norm.  First
the same localized construction gives the same-order evolution bound
\begin{equation}\label{eq:endpoint-MR-same-order-evolution}
 \sup_{r\leq\zeta\leq\vartheta}
 \|\mathcal U(\zeta,r)v\|_{C^{m,\alpha}_{\rm sc}}
 \leq C_{\Theta_0}\|v\|_{C^{m,\alpha}_{\rm sc}}.
\end{equation}
For \(m\geq2\), this is
Lemma~\ref{lem:weighted-prepared-Schauder} at order \(q=m\), localized
and chained over the same finite effective-time subdivision.  At
\(m=0\), it is the base-space \(C^{0,\alpha}\)-estimate in the same
\(q=2\) localized parametrix: the scalar-principal Kato inequality
controls its supremum component and the localized parabolic H\"older
estimate controls its \(\alpha\)-seminorm.  At \(m=1\), interpolate this
order-zero bound with the order-two bound furnished by that parametrix.
Thus all integer values \(m\geq0\) in the statement are covered.

Let \(0<d:=\vartheta-r\leq\vartheta_*\), put
\(u(\zeta)=\mathcal U(\zeta,r)v\), and choose a scalar cutoff
\(\eta\) which vanishes on the first quarter of \([r,\vartheta]\), is
one on its last half, and satisfies \(|\eta'|\leq C d^{-1}\).  The
section \(\eta u\) has zero trace at \(r\) and solves the same equation
with forcing \(\eta'u\).  Hence
\eqref{eq:effective-time-two-order-endpoint} and
\eqref{eq:endpoint-MR-same-order-evolution} give
\[
 \|u(\vartheta)\|_{C^{m+2,\alpha}_{\rm sc}}
 \leq C d^{-1}\|v\|_{C^{m,\alpha}_{\rm sc}}.
\]
The scale-one H\"older interpolation inequality between orders \(m\)
and \(m+2\), uniform over the buffered atlas, now yields
\begin{equation}\label{eq:endpoint-MR-variable-one-order-smoothing}
 \|\mathcal U(\vartheta,r)\|_
 {C^{m,\alpha}_{\rm sc}\to C^{m+1,\alpha}_{\rm sc}}
 \leq C_{\Theta_0}\bigl(1+(\vartheta-r)^{-1/2}\bigr).
\end{equation}
For a larger gap, factor the evolution at the beginning of its last
\(\vartheta_*/2\)-piece and use the finite same-order chain
\eqref{eq:endpoint-MR-same-order-evolution} on the earlier pieces.
Duhamel's formula and
\eqref{eq:endpoint-MR-variable-one-order-smoothing} prove
\eqref{eq:effective-time-one-order-kernel}; its terminal singularity is
integrable.

On the full source atlas choose \(\delta_0>0\) so that every displayed
source clock has length at most \(\vartheta_*\) whenever
\(|J|\leq\delta_0\).  Integrating
\eqref{eq:endpoint-MR-variable-one-order-smoothing} gives
\eqref{eq:source-atlas-one-order-Abel-block}.  Applying
Lemma~\ref{lem:weighted-prepared-Schauder} with \(q=m+2\) gives the
first inequality in
\eqref{eq:source-atlas-zero-trace-maximal-block}; its zero-trace
one-order-lowering estimate
\eqref{eq:weighted-zero-trace-one-order} gives the factor
\(\delta^{\alpha/4}\) in the second.  The atlas proof of that lemma
already includes the partition, polynomial weights, tensor
normalizations, and uniform finite overlap, so none is omitted here.

It remains to justify the parameter clauses.  Fix a base parameter
and fix once and for all the bounded linear trace extension
\(\mathcal T\) constructed in the proof of
Lemma~\ref{lem:weighted-prepared-Schauder}.  Replace
\(u_p\) by
\[
 \widetilde u_p:=u_p-\mathcal T u_{s,p}.
\]
This has exactly zero initial trace.  Since \(\mathcal T\) is fixed and
bounded, the induced forcing is \(C^1\) in the same trace and
strong-Bochner forcing norms.  Its first variation contains
\(\mathcal T D_pu_s[\dot p]\), and its Taylor remainder is the bounded
image of the trace Taylor remainder.  Thus the zero-trace estimates
just proved apply without discarding a varying initial datum.

In the coefficient route all parameters use the same effective
subdivision because the coefficient family and its first derivative
have the common modulus stated in the lemma.  The finite-horizon
difference quotient satisfies the secant version of the displayed
differentiated equation in the statement.  Subtract from it the
solution \(v\) of that differentiated equation.  The strong-Bochner
Fr\'echet remainders of \(A,B,C,F\), together with
\[
 D_p\log\lambda[\dot p]\,(F+\mathcal L_\tau u),
\]
the exact clock-variation term, make the residual source
\(o(\|\dot p\|)\) in
\(L^\infty C^{m,\alpha}_{\rm sc}\), uniformly on
\([s,\widehat S]\) and uniformly as \(\widehat S\uparrow S\).

In the residual-source route, subtract
\eqref{eq:endpoint-MR-fixed-operator-variation} from the secant form of
\eqref{eq:endpoint-MR-fixed-operator-residual-route}.  The fixed
operator cancels exactly, the clock derivative is precisely the term
displayed in \eqref{eq:endpoint-MR-fixed-operator-variation}, and the
only remaining nonlinear term is the Fr\'echet remainder of
\(F_p^{(p_0)}\).  Hence this route gives the same little-oh forcing
without differentiating the coefficients of
\(\mathcal L_{p_0,\tau}\).

The chained two-order estimate sends either residual to
\(o(\|\dot p\|)\) in the asserted output norm.  Repeating the argument
for two nearby base parameters proves continuity of the derivative.
Uniformity in the base parameter follows on compact subsets, or under
the explicit common derivative and remainder moduli stipulated in the
statement; it is not inferred from boundedness alone.  Thus the
difference, first-variation, Fr\'echet-remainder, and endpoint claims
all follow in the exact strong-Bochner topology stated above, with the
trace, forcing, and clock variations accounted for in the route in
which they occur.
\end{proof}

\begin{lemma}[Fixed-bottleneck Davies--Gaffney estimate]
\label{lem:fixed-bottleneck-Davies}
Fix an integer \(m\geq0\) and \(0<\alpha<1\).  Work on a fixed
Riemannian domain \((\mathcal D,g_*)\), either complete without boundary
or with smooth boundary equipped with a homogeneous Dirichlet or
homogeneous conormal realization, and on a finite-rank bundle
\(E\to\mathcal D\) with fixed bundle metric, compatible connection
\(\nabla\), and positive smooth density \(d\mu_*\), all uniformly
controlled in the source atlas.  Define the form domain
\begin{equation}\label{eq:boundary-form-space}
 V_\partial(\mathcal D;E):=
 \begin{cases}
  H^1(\mathcal D;E),
   &\partial\mathcal D=\varnothing,\\
  H^1_0(\mathcal D;E),
   &\text{homogeneous Dirichlet realization},\\
  H^1(\mathcal D;E),
   &\text{homogeneous conormal realization},
 \end{cases}
\end{equation}
and define
\[
 H^{-1}_\partial(\mathcal D;E):=
 V_\partial(\mathcal D;E)^*,
 \qquad
 H^{-1}(\mathcal D;E):=
 H^{-1}_\partial(\mathcal D;E),
\]
using the \(L^2(d\mu_*)\) pairing.  Thus every occurrence of
\(H^{-1}\) within this lemma and its proof means the dual of the
boundary-realization form domain, not automatically
\((H^1_0)^*\).

In the conormal case, assume that the operator and its adjoint are
realized by uniformly quasi-coercive closed forms on
\(V_\partial=H^1\), that smooth sections are a common form core, and
that the integration-by-parts boundary contributions are nonpositive
for both forms.  In the Dirichlet case multiplication by every weight
used below preserves \(H^1_0\); in the conormal case it preserves
\(H^1\), and the preceding form assumption controls the boundary
contribution.  Assume throughout that the evolution family and its
adjoint act consistently on the Gelfand triple
\[
 V_\partial\subset L^2(d\mu_*)\subset H^{-1}_\partial .
\]

Fix \(s<S\).  Let \(\lambda:[s,S)\to(0,\infty)\) be locally
integrable, put
\[
 \mathfrak s(q,\tau)=\int_q^\tau\lambda(r)\,dr,
\]
and let \(\mathcal U(\tau,q)\) be the evolution family of the original,
unnormalized operator
\[
 \partial_\tau-\mathcal L_\tau,\qquad
 \mathcal L_\tau
 =a_\tau^{ij}\nabla_i\nabla_j+b_\tau*\nabla+c_\tau
\]
on a complete uniformly locally finite source atlas.  Let
\(\mathcal O\) be an output union with uniform local
\(H^{-1}_\partial\)-to-\(C^{m,\alpha}\) smoothing.  Assume that the
original principal coefficient \(a_\tau\) is real symmetric and
nonnegative as a scalar quadratic form globally.
Before changing time, and after only the stated spatial rescaling,
assume on the fixed bottleneck and every buffered output chart that
\[
 \kappa\lambda(\tau)|\xi|_{g_*}^2
 \leq a_\tau(\xi,\xi)
 \leq\kappa^{-1}\lambda(\tau)|\xi|_{g_*}^2
\]
in each fixed-scale chart, with the scale-adapted analogue below.

To distinguish the two coefficient families, define the
clock-normalized operator in a fixed-scale chart by
\begin{equation}\label{eq:fixed-bottleneck-normalized-operator}
 \widehat{\mathcal L}_\vartheta
 :=\lambda(\tau(\vartheta))^{-1}\mathcal L_{\tau(\vartheta)}
 =\widehat a_\vartheta^{ij}\nabla_i\nabla_j
  +\widehat b_\vartheta*\nabla+\widehat c_\vartheta,
 \qquad d\vartheta=\lambda\,d\tau .
\end{equation}
Thus, in a fixed-scale chart,
\[
 \widehat a_\vartheta=\lambda^{-1}a_\tau,\qquad
 \widehat b_\vartheta=\lambda^{-1}b_\tau,\qquad
 \widehat c_\vartheta=\lambda^{-1}c_\tau .
\]
For a recorded time-independent scale-\(r_{\mathcal V}\) chart
\(\Phi_{\mathcal V}\), put
\[
 d\vartheta_{\mathcal V}
 =\lambda r_{\mathcal V}^{-2}\,d\tau
\]
and define the scale-adapted normalized operator
\(\widehat{\mathcal L}^{(\mathcal V)}_{\vartheta_{\mathcal V}}\)
by the exact identity
\begin{equation}\label{eq:fixed-bottleneck-scale-normalized-operator}
 \left(
  \partial_{\vartheta_{\mathcal V}}
  -\widehat{\mathcal L}^{(\mathcal V)}_{\vartheta_{\mathcal V}}
 \right)(\Phi_{\mathcal V}^*u)
 =
 \frac{r_{\mathcal V}^{\,2}}{\lambda(\tau)}
 \Phi_{\mathcal V}^*
 \bigl((\partial_\tau-\mathcal L_\tau)u\bigr).
\end{equation}
This definition incorporates both the spatial dilation and the clock
change; in particular, no unrecorded convention
\(\widehat b=\lambda^{-1}b\),
\(\widehat c=\lambda^{-1}c\) is imposed after a nontrivial spatial
rescaling.

Assume that the principal coefficient of every normalized operator is
uniformly elliptic in its scale-one chart and that the normalized
coefficients satisfy the full output-order bounds
\[
 \widehat a\in C^{m,\alpha},\qquad
 \widehat b\in C^{(m-1)_+,\alpha},\qquad
 \widehat c\in C^{(m-2)_+,\alpha},
\]
with one uniform bound and one effective-time
\(C^{\alpha/2}\) modulus.
Record also the resulting uniform local smoothing package: for some
integer \(p_m\geq0\), every buffered output
member \((\mathcal V,r_{\mathcal V})\) has the uniform local smoothing
bound
\[
 \|\mathcal U(\tau,q)g\|_
  {C^{m,\alpha}_{r_{\mathcal V}}(\mathcal V)}
 \leq
 C\left[
  1+\left(
   \frac{\mathfrak s(q,\tau)}{r_{\mathcal V}^{\,2}}
  \right)^{-p_m}\right]\|g\|_{H^{-1}(\mathcal D;E)}
\]
for \(\mathfrak s(q,\tau)>0\), with the bounded-scale interpretation
when \(r_{\mathcal V}\asymp1\).
Suppose that a smooth
function \(\psi\) has the following fixed-bottleneck properties:
\[
\begin{aligned}
 \psi&=R_{\rm b}&&\text{near }\mathcal O,\\
 \psi&=0&&\text{on the source region}\\
 &&&\text{and on the collapsing core},\\
 \supp d\psi&\Subset\mathcal B_{\rm b}.&
\end{aligned}
\]
where \(R_{\rm b}>0\) and \(\mathcal B_{\rm b}\) is a fixed
noncollapsing buffer.  Suppose that, in the fixed scale-one charts on
\(\mathcal B_{\rm b}\),
\begin{equation}\label{eq:bottleneck-Davies-coefficient-condition}
\begin{split}
 a_\tau(d\psi,d\psi)&\leq C\lambda(\tau),\\
 |b_\tau*d\psi|
 +|(\nabla a_\tau)*d\psi|
 +|a_\tau*\nabla^2\psi|&\leq C\lambda(\tau),
\end{split}
\end{equation}
In addition, in those fixed-scale charts, require that, for
\(0\leq\ell\leq m\), the \(C^{0,\alpha}\) norms of
\(\nabla^\ell\) of the normalized weight-commutator coefficients
\[
 \widehat a_\vartheta(d\psi,d\psi),\qquad
 \widehat b_\vartheta*d\psi
 +(\nabla\widehat a_\vartheta)*d\psi
 +\widehat a_\vartheta*\nabla^2\psi
\]
are uniformly bounded in those charts and carry the same
effective-time modulus.

Assume also the precise forward and adjoint form inequalities
\[
\begin{aligned}
 2\operatorname{Re}
 \langle\mathcal L_\tau v,v\rangle_
        {H^{-1}_\partial,V_\partial}
 &\leq
 -\kappa\lambda(\tau)\|\nabla v\|_{L^2(d\mu_*)}^2
 +2\beta(\tau)\|v\|_{L^2(d\mu_*)}^2,\\
 2\operatorname{Re}
 \langle\mathcal L_\tau^\dagger v,v\rangle_
        {H^{-1}_\partial,V_\partial}
 &\leq
 -\kappa\lambda(\tau)\|\nabla v\|_{L^2(d\mu_*)}^2
 +2\beta(\tau)\|v\|_{L^2(d\mu_*)}^2
\end{aligned}
\]
for every \(v\in V_\partial(\mathcal D;E)\), where
\(\kappa>0\) is fixed and
\[
 \sup_{q\leq\tau}\int_q^\tau\beta_+(r)\,dr\leq B_0<\infty .
\]
Assume finally that the normalized coefficients on the bottleneck and
the output charts have the prepared chartwise time-oscillation moduli.
No global comparison
\(a_\tau\simeq\lambda(\tau)g_*^{-1}\) is assumed.
For the first conclusion below, assume that the retained output charts
have bounded scale, equivalently that their local parabolic clocks are
uniformly comparable to
\(\mathfrak s(q,\tau)=\int_q^\tau\lambda(r)\,dr\).  This includes every
fixed finite buffered output atlas.  The scale-adapted clause below
treats an unbounded output whose local clocks are
\(\mathfrak s/r_{\mathcal V}^{2}\).

If
\[
 \mathcal F\in
 L^1_{\rm loc}\bigl(
  [s,S);H^{-1}_\partial(\mathcal D;E)\bigr)
\]
is strongly measurable and is supported where \(\psi=0\) for almost
every \(q\), then
\begin{equation}\label{eq:fixed-bottleneck-Duhamel}
\left\|\int_s^\tau
 \mathcal U(\tau,q)\mathcal F(q)\,dq\right\|_
 {C^{m,\alpha}(\mathcal O)}
\leq C\int_s^\tau
 (1+\mathfrak s(q,\tau)^{-p_m})
 e^{-cR_{\rm b}^2/\mathfrak s(q,\tau)}
 \|\mathcal F(q)\|_{H^{-1}}\,dq .
\end{equation}
If \(g\in H^{-1}(\mathcal D;E)\) is supported where \(\psi=0\), the
separated homogeneous estimate is
\begin{equation}\label{eq:fixed-bottleneck-homogeneous}
 \|\mathcal U(\tau,s)g\|_{C^{m,\alpha}(\mathcal O)}
 \leq
 C\bigl(1+\mathfrak s(s,\tau)^{-p_m}\bigr)
 e^{-cR_{\rm b}^2/\mathfrak s(s,\tau)}
 \|g\|_{H^{-1}} .
\end{equation}
The same kernel applies to a cutoff commutator written
\(\mathcal F=\nabla_iP^i+Q\), provided \(P,Q\in L^2\), both are
supported on the zero side of the bottleneck, and
\begin{equation}\label{eq:bottleneck-commutator-Hminusone}
 \|\mathcal F\|_{H^{-1}}
 \leq C\bigl(\|P\|_{L^2}+\|Q\|_{L^2}\bigr).
\end{equation}
For \(v\in V_\partial\), the divergence-form term in this clause means
\[
 \langle\nabla_iP^i,v\rangle_{H^{-1}_\partial,V_\partial}
 :=
 -\int_{\mathcal D}\langle P^i,\nabla_i v\rangle\,d\mu_* .
\]
Thus no undeclared boundary-trace contribution is present in the
conormal realization.
Every displayed kernel is defined to be zero at zero effective time.

There is also a scale-adapted unbounded-output version.  Suppose that
\(\mathcal O\) has a uniformly locally finite buffered cover
\[
 \mathfrak A_{\mathcal O}
 =\{(\mathcal V,r_{\mathcal V})\}
\]
and that its scaled \(C^{m,\alpha}\)-norm is the supremum of the
ordinary norms in the
\(r_{\mathcal V}^{-2}g_*\)-scale-one charts.  Assume that, for every
retained output member, there is a smooth weight
\(\psi_{\mathcal V}\) which vanishes on the source region and on the
collapsing core, satisfies
\begin{equation}\label{eq:bottleneck-radial-weight}
 \psi_{\mathcal V}\geq
 c(R_{\rm b}+r_{\mathcal V})
 \quad\hbox{near }\mathcal V,
\end{equation}
and whose gradient is supported in the fixed bottleneck together with
a pure-outer radial corridor on which
\eqref{eq:bottleneck-Davies-coefficient-condition}, including its
required spatial jets, holds uniformly with
\(\psi\) replaced by \(\psi_{\mathcal V}\).  Suppose also
that on the buffered output chart the local parabolic clock is
uniformly comparable to
\(\mathfrak s(q,\tau)/r_{\mathcal V}^{2}\).  Assume in addition that
\begin{equation}\label{eq:fixed-bottleneck-total-effective-time}
 \Theta_{\rm b}
 :=
 \sup_{s\leq q\leq\tau<S}\mathfrak s(q,\tau)
 <\infty .
\end{equation}
The constants in the scale-adapted conclusions may depend on
\(\Theta_{\rm b}\).  Then
\begin{equation}\label{eq:fixed-bottleneck-radial-Duhamel}
\left\|\int_s^\tau
 \mathcal U(\tau,q)\mathcal F(q)\,dq\right\|_
 {C_{\rm sc}^{m,\alpha}(\mathcal O)}
\leq C\int_s^\tau
 (1+\mathfrak s(q,\tau)^{-p_m})
 e^{-cR_{\rm b}^2/\mathfrak s(q,\tau)}
 \|\mathcal F(q)\|_{H^{-1}}\,dq .
\end{equation}
For the separated datum \(g\) above, each retained output chart first
obeys the sharper estimate
\[
 \|\mathcal U(\tau,s)g\|_
 {C^{m,\alpha}_{r_{\mathcal V}}(\mathcal V)}
 \leq
 C\left[
  1+\left(
    \frac{\mathfrak s(s,\tau)}{r_{\mathcal V}^{\,2}}
   \right)^{-p_m}\right]
 \exp\!\left(
  -\frac{c(R_{\rm b}+r_{\mathcal V})^2}
         {\mathfrak s(s,\tau)}\right)
 \|g\|_{H^{-1}} .
\]
Taking the atlas supremum gives
\begin{equation}\label{eq:fixed-bottleneck-radial-homogeneous}
 \|\mathcal U(\tau,s)g\|_{C_{\rm sc}^{m,\alpha}(\mathcal O)}
 \leq
 C\bigl(1+\mathfrak s(s,\tau)^{-p_m}\bigr)
 e^{-cR_{\rm b}^2/\mathfrak s(s,\tau)}
 \|g\|_{H^{-1}} .
\end{equation}
The same scale-adapted conclusion holds for the commutator satisfying
\eqref{eq:bottleneck-commutator-Hminusone}.  Accordingly, no
bottleneck-time gain is inferred from the remote source-adapted clock
\(\mathfrak s/r_{\mathcal V}^{2}\).
\end{lemma}

\begin{proof}
Conjugate by \(e^{\theta\psi}\).  The unweighted part of the covariant
energy form contributes \(\beta(\tau)\).  Every new first- or
zeroth-order term contains \(d\psi\) or \(\nabla^2\psi\), hence is
supported in \(\mathcal B_{\rm b}\); completing the square and using
\eqref{eq:bottleneck-Davies-coefficient-condition} gives
\[
 \frac{d}{d\tau}\|e^{\theta\psi}u\|_{L^2}^2
 \leq
 \bigl(
  2\beta_+(\tau)
  +C(\theta+\theta^2)\lambda(\tau)
 \bigr)
 \|e^{\theta\psi}u\|_{L^2}^2 .
\]
If \(\mathfrak s(q,\tau)<R_{\rm b}^2\), choose
\(\theta=cR_{\rm b}/\mathfrak s(q,\tau)\); the linear contribution
\(C\theta\mathfrak s\) is then a fixed \(O(R_{\rm b})\) factor, while
the quadratic contribution yields the Gaussian
\(e^{-cR_{\rm b}^2/\mathfrak s}\).
If \(\mathfrak s(q,\tau)\geq R_{\rm b}^2\), take \(\theta=0\);
the unweighted energy estimate and
\(e^{-cR_{\rm b}^2/\mathfrak s}\geq e^{-c}\) give the same conclusion
after changing the constant.  The adjoint form inequality on
\(V_\partial\), duality with
\(H^{-1}_\partial=V_\partial^*\), and the boundary-compatible
Caccioppoli estimate give the \(H^{-1}_\partial\)-to-\(L^2\) version.
Uniform local
source-adapted smoothing on the output charts supplies only the finite
factor \(1+\mathfrak s^{-p_m}\), because the bounded-scale output
clocks are uniformly comparable to \(\mathfrak s\).  Bounded overlap
glues the chart estimates, and Duhamel proves
\eqref{eq:fixed-bottleneck-Duhamel}.  A cutoff commutator is put in
\(H^{-1}\) in divergence form, exactly as in
Lemma~\ref{lem:anchored-buffered-Davies}.

For the scale-adapted clause, put
\(L_{\mathcal V}:=R_{\rm b}+r_{\mathcal V}\) and repeat the conjugated
energy calculation with \(\psi_{\mathcal V}\).  Write
\(\varsigma=\mathfrak s(q,\tau)\).  If
\(\varsigma<L_{\mathcal V}^{2}\), choose
\(\theta=cL_{\mathcal V}/\varsigma\).  The quadratic weight term is absorbed
by the separation exponent
\(\theta L_{\mathcal V}=cL_{\mathcal V}^{2}/\varsigma\), while the linear term
obeys
\[
 C\theta\varsigma
 \leq
 \varepsilon\frac{L_{\mathcal V}^{2}}{\varsigma}
 +C_\varepsilon\varsigma .
\]
After fixing \(\varepsilon>0\) sufficiently small, the first term is
absorbed into the Gaussian and the second is bounded by
\eqref{eq:fixed-bottleneck-total-effective-time}.  If
\(\varsigma\geq L_{\mathcal V}^{2}\), take \(\theta=0\); then the unweighted
forward and adjoint estimates and
\(e^{-cL_{\mathcal V}^{2}/\varsigma}\geq e^{-c}\) give the same conclusion.
Consequently the \(H^{-1}_\partial\)-to-\(L^2\) estimate on the
twice-enlarged output member carries the uniform factor
\[
 C_{B_0,\Theta_{\rm b}}
 \exp\!\left(
  -\frac{c(R_{\rm b}+r_{\mathcal V})^2}
         {\mathfrak s(q,\tau)}
 \right).
\]
Interior smoothing in the
\(r_{\mathcal V}^{-2}g_*\)-scale-one chart costs at most
\[
 1+\left(
   \frac{\mathfrak s(q,\tau)}{r_{\mathcal V}^{2}}
  \right)^{-p_m}.
\]
After decreasing \(c>0\), the elementary bound
\[
 \left[
  1+\left(\frac{r_{\mathcal V}^{2}}{\varsigma}\right)^{p_m}
  \right]
 e^{-c(R_{\rm b}+r_{\mathcal V})^2/\varsigma}
 \leq
 C(1+\varsigma^{-p_m})e^{-c'R_{\rm b}^2/\varsigma},
 \qquad \varsigma>0,
\]
is uniform in \(\mathcal V\).  Taking the supremum over the bounded
overlap cover and then integrating in \(q\) proves
\eqref{eq:fixed-bottleneck-radial-Duhamel}.  The homogeneous and
commutator variants are identical.
\end{proof}

For the fixed exterior reference, write
\[
 \mathscr R_{\rm D}(G)
 :=-2\Ric_G+\Lie_{B_{\widehat G_{\rm ext}}(G)}G
\]
for the Ricci--DeTurck right-hand side.

\begin{lemma}[Inner-terminated exterior Ricci--DeTurck estimate]
\label{lem:inner-terminated-exterior-DeTurck}
Fix \(m\geq4\), \(0<\alpha<1\), and two prepared evolutions beginning at
\(\tau_0\).  Suppose they share the exterior certificate
\eqref{eq:prepared-exterior-termination}, the anchored common DeTurck
gauge of Lemma~\ref{lem:anchored-exterior-interface} at order
\(r=m+2\), and its common effective-time coefficient modulus.  Put
\[
 \widehat G_i(\tau)=\widetilde G_i(t_i(\tau)),\qquad
 u=\widehat G_1-\widehat G_2 .
\]
Let \(\mathfrak D_{\rm ext}\) be the fixed double used below, and
extend by zero every tensor carrying the factor \(\zeta\).  The two
typed pieces of the homogeneous initial face are
\begin{align}
 d_{{\rm ext},m,0}^{+}
 &:=
 \sum_{a=1}^{N_{\rm ext}}
 \|u(\tau_0)\|_{C_{R_a}^{m,\alpha}(U_a^+)},
 \label{eq:exterior-buffered-initial-trace}\\
 d_{{\rm ext},-1,0}^{\rm corr}
 &:=
 \|(1-\vartheta_{\rm ext})\zeta u(\tau_0)\|_
   {H^{-1}(\mathfrak D_{\rm ext};g_{\rm ext})}.
 \label{eq:exterior-initial-corridor-memory}
\end{align}
Thus the first datum is a same-order trace on one unused atlas buffer,
whereas the second is only a low norm and is separated from every
retained output chart by
\eqref{eq:exterior-initial-corridor-separation}.  Bracketed versions
of these symbols denote the same seminorms with \(u(\tau_0)\) replaced
by its sliced first variation.  Put
\[
 K_{{\rm ext},m}^{0}(s)
 :=(1+s^{-q_m})e^{-cR_{{\rm ext},0}^2/s},
 \qquad K_{{\rm ext},m}^{0}(0):=0 .
\]
The common dynamic scale bracket gives the automatic ratio bound
\[
 \sup_\tau
 \left|\log\frac{\lambda_1(\tau)}{\lambda_2(\tau)}\right|
 \leq
 \log\frac{C_{\rm scl}}{c_{\rm scl}}
 =:C_{\rm rat}.
\]
Let
\[
 \mathfrak s(q,\tau):=\int_q^\tau\lambda_1(r)\,dr,\qquad
 K_{m,R_{\rm in}}(s):=(1+s^{-q_m})e^{-cR_{\rm in}^2/s},
 \quad K_{m,R_{\rm in}}(0):=0 .
\]
Here \(K_{m,R_{\rm in}}:[0,\infty)\to[0,\infty)\) is a scalar kernel,
not an
operator or tensor norm, and \(c>0\) and the integer \(q_m\) are the
constants furnished by the anchored buffered off-diagonal estimate
Lemma~\ref{lem:anchored-buffered-Davies}.  After increasing the common entrance
threshold, the effective-time range lies in an interval on which
\(K_{m,R_{\rm in}}\) is increasing.  In arguments using only an upper bound, we
instead use its finite supremum on that range.
Then, for every finite endpoint \(S\) and
\(\tau_0\leq\tau\leq S\),
\begin{equation}\label{eq:inner-terminated-exterior-DeTurck}
 \begin{split}
 \mathfrak G_m(\tau)\leq{}&
 C_Sd_{{\rm ext},m,0}^{+}
 +C_SK_{{\rm ext},m}^{0}\bigl(
      \mathfrak s(\tau_0,\tau)\bigr)
      d_{{\rm ext},-1,0}^{\rm corr}\\
 &+C_S\int_{\tau_0}^{\tau}
       \lambda_1(q)\left|\log\frac{\lambda_1}{\lambda_2}(q)\right|\,dq\\
 &+C_S\int_{\tau_0}^{\tau}
       \lambda_1(q)K_{m,R_{\rm in}}\bigl(\mathfrak s(q,\tau)\bigr)
       \mathfrak I_{\rm in}(q)\,dq .
 \end{split}
\end{equation}
The same estimate holds for sliced first variations.  If the
one-state coefficient package is uniform for all future time, then
\(C_S\) may be replaced by one constant \(C\).  In particular, if
\(d_0\geq0\),
\[
 d_{{\rm ext},m,0}^{+}
 +d_{{\rm ext},-1,0}^{\rm corr}\leq d_0,
\]
\(\mathfrak I_{\rm in}(q)\leq C d_0e^{A(q-\tau_0)}\) and the scale and
phase tails give
\(\sup_q|\log(\lambda_1/\lambda_2)(q)|\leq C_Dd_0\), then
\begin{equation}\label{eq:inner-terminated-exterior-global}
 \sup_{\tau\geq\tau_0}\mathfrak G_m(\tau)\leq C_Dd_0 .
\end{equation}
At \(m=4\), the two typed quantities
\eqref{eq:exterior-buffered-initial-trace} and
\eqref{eq:exterior-initial-corridor-memory} are the complete
homogeneous exterior entrance memory.  They are retained as explicit
low hybrid entrance summands when that entrance distance is introduced;
 no single order-six graph-augmented prepared distance
 \(d_{\rm prep}^{6,\alpha}\) is substituted for them.
\end{lemma}

\begin{proof}
In the common gauge the exact asynchronous difference equation on
\(E^{++}\) is
\begin{equation}\label{eq:terminated-asynchronous-DeTurck}
 \partial_\tau u-\lambda_1\mathcal A^{ab}\nabla_a\nabla_bu
 =
 \lambda_1(\mathcal B*\nabla u+\mathcal C*u)
 +(\lambda_1-\lambda_2)\mathscr R_{\rm D}(\widehat G_2).
\end{equation}
Use the fixed cutoff \(\zeta\) in the exterior certificate.  Its support
is contained in the retained largest-set cover of \(\overline{E^+}\),
and its derivative is supported in
\(\operatorname{int}\mathcal A_{\rm in}\).  Extend
\(v=\zeta u\) by zero to \(\mathfrak D_{\rm ext}\), extending the
coefficients from their uniformly controlled neighborhood of
\(\supp\zeta\) by one fixed uniformly parabolic package.  The cutoff
commutator is used only in divergence form:
\[
 [\mathcal A^{ab}\nabla_a\nabla_b,\zeta]u
 =\operatorname{div}(\mathcal P_\zeta u)+\mathcal Q_\zeta u .
\]
The domain \(\mathfrak D_{\rm ext}\) is the fixed closed double, so its
form space is \(H^1\).  Integration by parts for the scalar-symbol
polarized Ricci--DeTurck operator and for its adjoint gives
\eqref{eq:anchored-Davies-energy} with
\[
 \beta_+(\tau)\leq C\lambda_1(\tau).
\]
The scale comparison \(\lambda_1(\tau)\leq Ce^{-\tau}\) therefore gives
both
\[
 B_0\leq C\int_{\tau_0}^{\infty}\lambda_1<C,
 \qquad
 \Theta_0\leq\int_{\tau_0}^{\infty}\lambda_1<C,
\]
uniformly in the finite endpoint.  Thus the absorbed off-diagonal
kernel below has exactly the stated endpoint-independent constant.
Consequently it is bounded in \(H^{-1}\) by
\(\|u\|_{L^2(\mathcal A_{\rm in})}\).  From
\eqref{eq:inner-terminal-graph-identity} and the prepared composition
calculus,
\[
 \|u(q)\|_{L^2(\mathcal A_{\rm in})}
 \leq C\mathfrak I_{\rm in}(q).
\]
This is the decisive point: no \(C^{m,\alpha}\) trace of a pullback
metric, and hence no unavailable \((m+1)\)-st map derivative, is used.

 The \(H^{-1}\) Duhamel estimate
\eqref{eq:anchored-Davies-Duhamel}, followed on the separated domain
\(E\) by interior parabolic Schauder estimates, gives the kernel
\(K_{m,R_{\rm in}}(\mathfrak s(q,\tau))\).  The clock term in
\eqref{eq:terminated-asynchronous-DeTurck} is bounded by
\[
 C\lambda_1\left|\log\frac{\lambda_1}{\lambda_2}\right|.
\]
For the homogeneous initial face use the exact split
\[
 v(\tau_0)=
 \vartheta_{\rm ext}\zeta u(\tau_0)
 +(1-\vartheta_{\rm ext})\zeta u(\tau_0).
\]
The first summand, extended by zero, has
\(C^{m,\alpha}\)-norm bounded by
\(Cd_{{\rm ext},m,0}^{+}\), because the derivative of
\(\vartheta_{\rm ext}\) is contained in
\(\bigcup_aU_a^+\).  The second summand is separated from
\(\bigcup_a\overline{U_a}\) by \(4R_{{\rm ext},0}\); the homogeneous
part of Lemma~\ref{lem:anchored-buffered-Davies} therefore bounds it by
\[
 CK_{{\rm ext},m}^{0}\bigl(\mathfrak s(\tau_0,\tau)\bigr)
 d_{{\rm ext},-1,0}^{\rm corr}.
\]
This is precisely the initial corridor contribution which cannot be
recovered from \(\mathfrak G_m(\tau_0)\) on the smaller output atlas.
Duhamel's formula therefore proves
\eqref{eq:inner-terminated-exterior-DeTurck}; differentiating the
identities and using the same linear evolution proves its
first-variation version.

For the final assertion,
\(\mathfrak s(q,\tau)\leq Ce^{-q}\).  Hence
\(K_{m,R_{\rm in}}(\mathfrak s(q,\tau))\leq
C e^{-c'e^q}\), up to a fixed polynomial in \(e^q\), and this absorbs
the assumed coarse exponential growth of \(\mathfrak I_{\rm in}\).
The total initial effective-time range is bounded, so
\(K_{{\rm ext},m}^{0}(\mathfrak s(\tau_0,\tau))\) has one finite
package supremum; the two typed initial terms are therefore
\(O(d_0)\).
The clock integral is \(O_D(d_0)\) because
\(\lambda_1(q)\leq Ce^{-q}\).  This proves
\eqref{eq:inner-terminated-exterior-global}.
\end{proof}

\begin{lemma}[One-time graft-input propagation with exterior memory]
\label{lem:compact-graft-buffer-propagation}
Fix \(m\geq4\).  Let two admissible prepared evolutions share the
separated graft-input chain
\[
 \mathcal W_{\rm gr}\Subset\mathcal W_{\rm gr}^+
 \Subset\mathcal W_{\rm gr}^{0}\Subset\cdots
 \Subset\mathcal W_{\rm gr}^{5}
 \Subset\mathcal W_{\rm gr}^{++}\Subset E
\]
and the common anchored exterior gauge.  Suppose that their one-state
closed-metric, gauge, and marking bounds are available through metric
order \(m+2\) and map order \(m+1\).  In the Ricci--DeTurck
polarization this controls \(A\) through order \(m+2\), \(B\) through
order \(m+1\), and \(C\) through order \(m\), exactly the coefficient
bounds in Lemma~\ref{lem:anchored-buffered-Davies} at output order
\(m+2\).  Define the typed entrance
metric trace on the collar actually supporting the localization cutoff
by
\begin{equation}\label{eq:typed-graft-larger-collar-metric-trace}
 d_{{\rm Ggr},m+2,0}^{0}
 :=
 \|\widetilde G_1(t_1(\tau_0))
       -\widetilde G_2(t_2(\tau_0))\|_
   {C_{\rm sc}^{m+2,\alpha}(\mathcal W_{\rm gr}^{0})}.
\end{equation}
Its bracketed version denotes the corresponding sliced tangent
seminorm.  Define the full typed entrance quantity
\begin{equation}\label{eq:typed-graft-initial-distance}
 d_{{\rm gr},m,0}:=
 d_{{\rm Ggr},m+2,0}^{0}
 +\mathfrak B_{{\rm gr},m}(\tau_0)
 +\mathfrak G_m(\tau_0)+\mathfrak I_{\rm in}(\tau_0)
 +\left|\log\frac{\lambda_1}{\lambda_2}(\tau_0)\right|
 +|t_1(\tau_0)-t_2(\tau_0)|.
\end{equation}
Thus \(d_{{\rm gr},m,0}\) is a scalar distance in the anchored hybrid
chart; it is not an undefined higher prepared norm.  Let
\(R_{\rm gr}>0\) be one quarter of the least separation between
successive graft collars and put
\[
 K_{{\rm gr},m}(\sigma)
 :=(1+\sigma^{-q_{m+2}})
   \exp(-cR_{\rm gr}^2/\sigma),\qquad
 K_{{\rm gr},m}(0):=0.
\]
Define the effective-time clock-source amplitude
\begin{equation}\label{eq:graft-clock-source-amplitude}
 \mathfrak L_{12}(\tau)
 :=\sup_{\tau_0\leq q\leq\tau}
   \left|\log\frac{\lambda_1}{\lambda_2}(q)\right|
 \leq
 \left|\log\frac{\lambda_1}{\lambda_2}(\tau_0)\right|
 +\int_{\tau_0}^{\tau}|a_1-a_2|(q)\,dq .
\end{equation}
The inequality is the exact consequence of
\(\partial_\tau\log(\lambda_1/\lambda_2)=-(a_1-a_2)\).
This \(K_{{\rm gr},m}\) is a nonnegative scalar kernel.  With
\(\mathfrak s(q,\tau)=\int_q^\tau\lambda_1(r)\,dr\), one has on every
finite normalized horizon
\begin{equation}\label{eq:compact-graft-buffer-finite}
\begin{split}
 \mathfrak B_{{\rm gr},m}(\tau)\leq{}&
 C_S d_{{\rm gr},m,0}
 +C_S\mathfrak L_{12}(\tau)\\
 &+C_S\int_{\tau_0}^{\tau}
   \lambda_1(q)K_{{\rm gr},m}\bigl(\mathfrak s(q,\tau)\bigr)
   \bigl(\mathfrak G_m(q)+\mathfrak I_{\rm in}(q)\bigr)\,dq .
\end{split}
\end{equation}
The initial face is always \(\tau_0\).  In particular,
\eqref{eq:compact-graft-buffer-finite} is not a same-collar unit
restart: the final integral is the off-diagonal lateral and exterior
memory omitted by such a restart.
Under the uniform future coefficient and effective-time package used
in Lemma~\ref{lem:uniform-weighted-Schauder-restart}, the constant
\(C_S\) in \eqref{eq:compact-graft-buffer-finite} can be chosen
independently of the terminal normalized time.

The same estimate holds for sliced first variations, with every
quantity replaced by its linearized block.  On a bounded
common-margin chart it also holds for every difference quotient, and
the remainder after subtraction of the linearized estimate is
\(o(1)\) times the corresponding hybrid increment.  If
\(\lambda_1(q)\leq Ce^{-q}\) and
\(\mathfrak G_m+\mathfrak I_{\rm in}\leq
C d\,e^{A(q-\tau_0)}\), while the clock ratio is bounded by \(Cd\),
then the Gaussian in \(K_{{\rm gr},m}\) absorbs the coarse exponential
growth and
\begin{equation}\label{eq:compact-graft-buffer-global-memory}
 \sup_{\tau\geq\tau_0}\mathfrak B_{{\rm gr},m}(\tau)
 \leq C d_{{\rm gr},m,0}+C_Dd .
\end{equation}
\end{lemma}

\begin{proof}
Put both closed flows first in the single anchored gauge of
Lemma~\ref{lem:anchored-exterior-interface}; no second compact gauge is
introduced.  Their metric difference \(u\) then satisfies the exact
polarized equation \eqref{eq:terminated-asynchronous-DeTurck}.  Choose
\(\xi_{\rm gr}\) equal to one on
\(\mathcal W_{\rm gr}^+\) and supported in
\(\mathcal W_{\rm gr}^{0}\).  The equation for
\(\xi_{\rm gr}u\) has the original initial trace at \(\tau_0\), the
clock source
\[
 (\lambda_1-\lambda_2)\mathscr R_{\rm D}(\widehat G_2),
\]
and a divergence-form commutator supported a distance at least
\(R_{\rm gr}\) from \(\mathcal W_{\rm gr}^+\).  The clock source is
bounded, in the scale-normalized order-\(m\) source norm, by
\[
 C\lambda_1\left|\log(\lambda_1/\lambda_2)\right|.
\]
Extend this compactly supported localization and its coefficients to
the same fixed closed double used for the anchored exterior equation.
The scalar-symbol principal form and its adjoint obey
\eqref{eq:anchored-Davies-energy} there with
\(\beta_+\leq C\lambda_1\).  Hence the exponential scale comparison
supplies common finite constants \(B_0,\Theta_0\), independently of the
terminal endpoint.  This verifies the realization and energy
hypotheses whenever the separated commutator is inserted into
Lemma~\ref{lem:anchored-buffered-Davies}.
After division by the effective clock \(\lambda_1\,d\tau\), this is an
\(L^\infty_{\mathfrak s}C_{\rm sc}^{m,\alpha}\) source with norm
\(C\mathfrak L_{12}(\tau)\).  Apply the two-order endpoint estimate
\eqref{eq:effective-time-two-order-endpoint} from
Lemma~\ref{lem:effective-time-endpoint-maximal-regularity} to this
nonseparated clock contribution.  This is the step which recovers the
output \(C_{\rm sc}^{m+2,\alpha}\) norm from an order-\(m\) clock
source; no ordinary \(L^1\) same-order endpoint estimate is used.
The homogeneous metric datum
\(\xi_{\rm gr}u(\tau_0)\) is controlled in
\(C_{\rm sc}^{m+2,\alpha}\) by
\(Cd_{{\rm Ggr},m+2,0}^{0}\).  This is why the entrance trace is taken
on \(\mathcal W_{\rm gr}^{0}\), rather than only on the output collar
\(\mathcal W_{\rm gr}^{+}\).
Apply Lemma~\ref{lem:anchored-buffered-Davies} at output order
\(m+2\) only to the separated commutator
and then the two-step interior Schauder estimate on the remaining
auxiliary collars.  Its \(H^{-1}\)-to-\(C^{m+2,\alpha}\) form gives
exactly the last integral in
\eqref{eq:compact-graft-buffer-finite}; the commutator norm is
controlled by the exterior block, and any portion meeting the
terminating cutoff is controlled by the exact low trace
\(\mathfrak I_{\rm in}\).  This proves the metric component without
unknown boundary data on \(\partial\mathcal W_{\rm gr}^{0}\).

The gauge equation \eqref{eq:anchored-gauge-difference-equation} and
its inverse equation are triangular.  Interior parabolic estimates on
\(\mathcal W_{\rm gr}\Subset\mathcal W_{\rm gr}^+\), followed by
\(\widetilde\iota_i=\iota\circ\chi_i^{-1}\), bound their order-\(m+1\)
differences by the just-obtained order-\(m+2\) metric term, the typed
initial gauge terms, and the same clock and off-diagonal memory.  This
proves \eqref{eq:compact-graft-buffer-finite}.

Finally \(\mathfrak s(q,\tau)\leq Ce^{-q}\).  Hence
\(K_{{\rm gr},m}(\mathfrak s(q,\tau))\) is bounded by a fixed
polynomial in \(e^q\) times \(e^{-c'e^q}\), proving the global clause.
Differentiation gives the variational estimate.  For a difference
quotient, the mean-value polarizations of the metric, gauge, inverse,
and composition maps are \(C^1\) at the displayed buffered orders;
Taylor's integral formula and dominated convergence give the stated
remainder estimate for all four blocks, not only for the final graft
operator.
\end{proof}

\begin{lemma}[Coefficient and parameter bounds for the polarized
relative harmonic-map equation]
\label{lem:polarized-HMHF-full-coefficient-ledger}
Fix an integer \(m\geq0\).  Work in one of the fixed buffered
source-adapted charts, with the tensor and map normalizations used in
\eqref{eq:source-adapted-atlas}--\eqref{eq:source-adapted-clock}.
Suppose the source and target metrics, their inverses, and the coordinate
representatives of the two maps and their inverses are uniformly bounded
through order \(m+2\), with a fixed ellipticity and buffer margin.  After
right translation to one map chart, the polarized relative harmonic-map
operator has the form
\[
 A^{ij}\partial_i\partial_j+B^i\partial_i+C
\]
and obeys
\begin{equation}\label{eq:polarized-HMHF-full-coefficient-ledger}
 \|A\|_{C^{m+1,\alpha}}
 +\|B\|_{C^{m,\alpha}}
 +\|C\|_{C^{m,\alpha}}
 \leq C_m .
\end{equation}
The same assertion holds in the global atlas-supremum strong-Bochner
norm, with the common effective-time modulus and spatial tail supplied
by the prepared package.

For top-order parameter dependence, define the unpolarized coordinate
tension map
\[
 \mathscr H(g,S,F):=\Delta_{g,S}F.
\]
Let \(P\) be an open subset of a Banach space, or of a Banach manifold
in one fixed chart, and let the preceding metrics, maps, charts, and
clock depend Fr\'echet \(C^1\) on \(p\in P\) in their displayed
order-\((m+2,\alpha)\) prepared norms.  Then \(\mathscr H\), the
difference of two such tension fields after the fixed right
translations, and the resulting exact nonseparated residual are
Fr\'echet \(C^1\) into
\[
 L^\infty\bigl(I;C_x^{m,\alpha}\bigr).
\]
This assertion is local at each base parameter.  It is uniform on a
compact parameter subset, or on a set carrying an explicitly assumed
common modulus for the first derivative and its Taylor remainder; it
is not promoted to a uniform assertion merely because a
common-margin subset is bounded.

No \(C^1\) assertion is made here for the top-order map
\(p\mapsto C_p\) in \(C^{m,\alpha}\).  Instead, at a base parameter
\(p_0\), freeze the fixed-state polarized operator
\(\mathcal L_{p_0}\) from
\eqref{eq:polarized-HMHF-full-coefficient-ledger} and write the exact
right-translated relative equation in the residual-source form
\[
 \partial_\tau w_p-\lambda_p\mathcal L_{p_0}w_p
   =\lambda_p\mathcal R_p^{(p_0)},
\]
with any separately recorded phase and separated cutoff sources left
outside this displayed nonseparated block.  The combined residual
\(p\mapsto\mathcal R_p^{(p_0)}\) is \(C^1\) in the preceding forcing
space by the nonlinear tension-map assertion; its definition is by
this exact nonlinear identity, not by separately differentiating the
polarized zeroth-order coefficient.  Together with the \(C^1\) clock,
chart, and actual initial-trace maps, this verifies the fixed-operator
residual-source and varying-trace clauses of
Lemma~\ref{lem:effective-time-endpoint-maximal-regularity}.  Hence all
first-variation and Fr\'echet-remainder conclusions hold at the stated
prepared order without assuming an unavailable extra target-metric
derivative.
\end{lemma}

\begin{proof}
In source coordinates \(x^i\) and target coordinates \(y^\gamma\),
the tension field is
\begin{equation}\label{eq:relative-HMHF-coordinate-ledger}
 (\Delta_{g,S}F)^\gamma
 =g^{ij}\left(
   \partial_i\partial_jF^\gamma
   -\Gamma(g)^k_{ij}\partial_kF^\gamma
   +\Gamma(S)^\gamma_{\beta\delta}(F)
      \partial_iF^\beta\partial_jF^\delta
 \right).
\end{equation}
Polarize this identity by the mean-value formula along the segment in
the fixed right-translated chart.  The principal coefficient is an
inverse source metric.  Every first-order coefficient is a finite sum
of products of an inverse source metric, one source or target
Christoffel symbol, and first derivatives of an interpolated map.
Every zeroth-order coefficient is a finite sum of the same factors and
terms of the form
\[
 g^{ij}\,
 D\Gamma(S)(F)[\,\cdot\,]\,
 \partial_iF\,\partial_jF .
\]
Thus the worst zeroth-order factor uses two derivatives of the target
metric and no more than first derivatives of either map.  Since
\(C^{m,\alpha}\) is a Banach algebra, buffered inversion and composition
give exactly
\eqref{eq:polarized-HMHF-full-coefficient-ledger}; in particular, the
zeroth-order coefficient has order \(m\), not merely order \(m-1\).
The fixed transition maps, tensor normalizations, and polynomial weights
have the recorded bounded jets, so taking the atlas supremum preserves
the estimate and its common moduli.

For the parameter assertion, do not differentiate the polarized
zeroth-order coefficient at top order.  Indeed, differentiating its
factor
\[
 D\Gamma(S)(F)[\,\cdot\,]
\]
with respect to the map would produce
\(D^2\Gamma(S)(F)[\dot F,\cdot]\), which generally uses three spatial
derivatives of \(S\).  That calculation is neither needed nor licensed
by the order-\((m+2,\alpha)\) hypothesis.

Differentiate instead the unpolarized identity
\eqref{eq:relative-HMHF-coordinate-ledger} before polarization.  Its
first derivative is a sum of the variations of \(g^{-1}\),
\(\partial^2F\), and \(\Gamma(g)\partial F\), together with
\[
 \begin{split}
 g^{-1}\bigl(&D_S\Gamma(S)[\dot S](F)
                 (\partial F)^2
 +D_y\Gamma(S)(F)[\dot F](\partial F)^2\\
 &+\Gamma(S)(F)\,\partial\dot F\,\partial F
 +\Gamma(S)(F)\,\partial F\,\partial\dot F\bigr).
 \end{split}
\]
Here \(D_S\Gamma(S)[\dot S]\) uses at most one spatial derivative of
the metric variation, while \(D_y\Gamma(S)\) uses at most two spatial
derivatives of \(S\).  The other terms use at most two derivatives of
\(F\) or \(\dot F\) and at most one derivative of the source metric.
The H\"older algebra, buffered inversion, and buffered composition
estimates therefore place this derivative in \(C^{m,\alpha}\) from
order-\((m+2,\alpha)\) inputs.  The same calculation, in difference
form, gives the Fr\'echet little-oh remainder in the exact
strong-Bochner forcing norm.  Fixed right translations and chart
changes preserve this conclusion by
Lemma~\ref{lem:prepared-chart-calculus}.

At a base parameter, move only the unknown derivative represented by
the fixed operator \(\mathcal L_{p_0}\) to the left side.  All remaining
terms form the combined residual \(\mathcal R_p^{(p_0)}\).  The
preceding nonlinear calculation proves that this residual is \(C^1\)
in the required forcing norm, so
\eqref{eq:endpoint-MR-fixed-operator-residual-route}--%
\eqref{eq:endpoint-MR-fixed-operator-variation} apply.  The actual
initial trace is handled by the fixed trace extension in that lemma,
and the scalar clock assertion follows by differentiating its positive
ratio.  Continuity and little-oh are local at the base parameter;
uniformity over base parameters is used only on compact subsets or
under the explicit common-modulus hypothesis stated above.  Thus the
parameter conclusion uses no derivative of \(C_p\) and no hidden third
target-metric derivative.
\end{proof}

Every subsequent fixed-state invocation of the source-atlas Abel block
for the relative harmonic-map equation is made under
\eqref{eq:polarized-HMHF-full-coefficient-ledger}; no later occurrence
of that phrase denotes an unproved staggered-coefficient estimate.
Every top-order parameter or remainder invocation uses instead the
nonlinear fixed-operator residual-source route just proved and does not
differentiate the polarized zeroth-order coefficient.

\begin{lemma}[Localized \(F\)-propagation: separated memory and Abel block]
\label{lem:localized-graft-F-coarse-memory}
Fix \(m\geq4\).  On a finite horizon, suppose that in the fixed
scale-one charts on \(\mathfrak S_{\rm gr}^{++}\) the polarized
source-adapted \(F\)-operator has
\[
 \|A\|_{C^{m+1,\alpha}}
 +\|B\|_{C^{m,\alpha}}
 +\|C\|_{C^{m,\alpha}}\leq\Lambda_{{\rm F},m},
\]
with one common effective-time \(C^{\alpha/2}\) modulus.  In the
harmonic-map application these bounds, including the full
zeroth-order bound, are supplied by
Lemma~\ref{lem:polarized-HMHF-full-coefficient-ledger} from the source
and target metric coefficients through \(C^{m+2,\alpha}\) and the
right-translated maps and inverses through \(C^{m+2,\alpha}\).  Thus a
prepared input at order \(m+2\), whose map components in fact have one
additional derivative, supplies the stated package.  The
nonseparated target-connection difference below is therefore measured
in the auxiliary
\(\mathfrak X_{\rm sc}^{m+2,\alpha}\) \(R\)-norm; this is one order
above the final \(R\)-output but remains within the prepared map input.
Let
\(R_{\rm F}>0\) be one quarter of
the recorded separation of \(\mathfrak S_{\rm gr}\) from the cutoff
corridor outside \(\mathfrak S_{\rm gr}^{+}\).  Then
\begin{equation}\label{eq:localized-graft-F-coarse-memory}
\begin{split}
 \mathfrak F_{{\rm gr},m+1}(\tau)\leq{}&
 C_S d_{{\rm Fgr},m+1,0}^{++}
 +C_S\int_{\tau_0}^{\tau}|\delta c(q)|\,dq\\
 &+C_S\int_{\tau_0}^{\tau}\lambda_1(q)
 K_{\rm MR}^{(1)}\bigl(\mathfrak s(q,\tau)\bigr)
 \left(
  \mathfrak G_m+\mathfrak B_{{\rm gr},m}
  +\left|\log\frac{\lambda_1}{\lambda_2}\right|
  +\|R_1-R_2\|_{\mathfrak X_{\rm sc}^{m+2,\alpha}}
 \right)(q)\,dq\\
 &+C_S\int_{\tau_0}^{\tau}
 \lambda_1(q)
 K_{{\rm Fgr},m}\bigl(\mathfrak s(q,\tau)\bigr)
 \left(
  \mathfrak D_m^{\rm hyb}
  +\mathfrak F_{m+1}^{\rm glob}
 \right)(q)\,dq ,
\end{split}
\end{equation}
where
\[
 K_{\rm MR}^{(1)}(\sigma):=1+\sigma^{-1/2},
 \qquad \sigma>0,
\]
is used only under the locally integrable effective-time integral, and
\[
 K_{{\rm Fgr},m}(\sigma)
 =(1+\sigma^{-q_{m+1}})e^{-cR_{\rm F}^2/\sigma},
 \qquad K_{{\rm Fgr},m}(0)=0,
\]
is the scalar off-diagonal kernel associated with the separation of
\(\mathfrak S_{\rm gr}\) and
\(\partial\mathfrak S_{\rm gr}^+\).  For the later restart estimate,
denote by
\begin{equation}\label{eq:localized-graft-F-hs-memory}
 \begin{split}
 \mathfrak F_{{\rm gr},m+1}^{\rm hs}(\tau)
 :={}&C_S d_{{\rm Fgr},m+1,0}^{++}\\
 &+C_S\int_{\tau_0}^{\tau}
  \lambda_1(q)
  K_{{\rm Fgr},m}\bigl(\mathfrak s(q,\tau)\bigr)
  \left(\mathfrak D_m^{\rm hyb}
        +\mathfrak F_{m+1}^{\rm glob}\right)(q)\,dq .
 \end{split}
\end{equation}
This is the \emph{homogeneous/separated localized \(F\)-memory}.
Equation~\eqref{eq:localized-graft-F-coarse-memory} bounds the full
localized block by
\(\mathfrak F_{{\rm gr},m+1}^{\rm hs}\), the smooth phase \(L^1\)
term, and the nonseparated Abel term.  Only
\(\mathfrak F_{{\rm gr},m+1}^{\rm hs}\) is a one-time memory; the
nonseparated zero-trace contribution is estimated on each short
interval by
\eqref{eq:source-atlas-one-order-Abel-block}.
The final global-\(F\) summand is
the low lateral value created by cutting the unknown; crucially, both
it and the other separated sources retain the physical factor
\(\lambda_1\).  The homogeneous/separated memory is always anchored at
the initial face \(\tau_0\).  The global nonseparated convolution in
\eqref{eq:localized-graft-F-coarse-memory} is justified by the
coefficient bounds above and
\eqref{eq:effective-time-one-order-kernel}.  In the later restart
argument the same contribution may instead be recentered with zero
trace on each short interval and estimated by
\eqref{eq:source-atlas-one-order-Abel-block}.  Under the uniform future
coefficient package the constant in this formula is independent of the
terminal time.

For the parameterized assertions, let \(P\) be an open subset of the
independent prepared model \(\mathscr E_{\rm prep}^{m+2,\alpha}\), and
let \(p\mapsto\mathbf z_p\) be a Fr\'echet \(C^1\) family in one bounded
common-margin chart, with the common time modulus and spatial tail of
the prepared package.  Assume that its actual initial trace is \(C^1\)
in the typed larger-star norm and use the positive clock
\(\lambda_p\).  At each base parameter \(p_0\), freeze the fixed-state
operator supplied by
Lemma~\ref{lem:polarized-HMHF-full-coefficient-ledger}.  The nonlinear
tension-map part of that lemma, together with
Lemma~\ref{lem:prepared-chart-calculus}, writes the exact
nonseparated equation with this fixed operator and a \(C^1\) residual
source in the strong-Bochner order-\(m\) forcing norm.  The actual
initial trace is removed by the fixed bounded trace extension, and the
clock is treated by
\eqref{eq:endpoint-MR-fixed-operator-variation}.  Thus the same formula
holds for sliced first variations without differentiating the
top-order polarized zeroth-order coefficient.  Subtracting that first
variation from a difference quotient leaves an \(o(1)\)-multiple of
the corresponding hybrid increment on the right in the exact forcing
norm.  This little-oh assertion is local at \(p_0\); it is uniform over
base parameters only on compact subsets or under an explicitly common
derivative and remainder modulus.
\end{lemma}

\begin{proof}
Subtract the exact \(F\)-equations in their separate source-adapted
clocks and polarize in one right-translated map chart.  Let \(w\)
denote the resulting coordinate difference.  Choose a fixed cutoff
\(\xi_{\rm F}\) which is one on
\(\mathfrak S_{\rm gr}^{+}\), is supported in
\(\mathfrak S_{\rm gr}^{++}\), and whose derivative is separated from
\(\mathfrak S_{\rm gr}\) by \(R_{\rm F}>0\).  We cut the unknown:
\[
 z:=\xi_{\rm F}w .
\]
Every chart meeting \(\supp\xi_{\rm F}\) lies in the genuinely
noncollapsing physical graft buffer.  In its fixed scale-one atlas the
polarized equation for \(z\) has the form
\begin{equation}\label{eq:localized-F-cutoff-equation}
 \partial_\tau z-\lambda_1(\tau)
 \bigl(A^{ab}\nabla_a\nabla_b+B*\nabla+C\bigr)z
 =\xi_{\rm F}\mathcal Q
  +\lambda_1[\mathcal L,\xi_{\rm F}]w .
\end{equation}
Here the principal symbol is scalar and uniformly elliptic.  The
one-state graft-buffer package and
Lemma~\ref{lem:polarized-HMHF-full-coefficient-ledger} give the
displayed \(C^{m+1,\alpha}/C^{m,\alpha}/C^{m,\alpha}\) spatial bounds
for \(A,B,C\), respectively; the common scale bracket controls the
polarization of the two clocks, and the source-adapted coefficient
package gives the asserted common H\"older time modulus in
\(\mathfrak s=\int\lambda_1\).  Thus all principal, lower-order, and
time-modulus hypotheses of
Lemma~\ref{lem:anchored-buffered-Davies} hold on this fixed
noncollapsing star.  Nothing is asserted about a scalar-\(\lambda\)
factor for the global operator on the collapsing core.
More precisely, extend the compactly supported star equation, its
coefficient bundle, and its coefficients to one fixed closed double
\(\mathfrak D_{\rm F}\) of \(\mathfrak S_{\rm gr}^{++}\), and denote
the extended bundle by \(E_{\rm F}\to\mathfrak D_{\rm F}\).  The
forward and adjoint scalar-symbol
forms satisfy \eqref{eq:anchored-Davies-energy} with
\(\beta_+\leq C\lambda_1\); since
\(\lambda_1(\tau)\leq Ce^{-\tau}\), the constants \(B_0\) and
\(\Theta_0\) in \eqref{eq:anchored-Davies-clock-growth} are uniform.
Thus the later separated-memory uses invoke the fully typed
off-diagonal lemma.

The homogeneous datum \(z(\tau_0)\) is supported in
\(\mathfrak S_{\rm gr}^{++}\) and is bounded at the required order by
\(Cd_{{\rm Fgr},m+1,0}^{++}\).  The commutator in
\eqref{eq:localized-F-cutoff-equation} is written in divergence form.
Because \(\supp d\xi_{\rm F}\) is separated from the output and remains
in the noncollapsing buffer,
\[
 \|\lambda_1[\mathcal L,\xi_{\rm F}]w\|_
 {H^{-1}(\mathfrak D_{\rm F};E_{\rm F})}
 \leq C\lambda_1\mathfrak F_{m+1}^{\rm glob}.
\]
This is the only global lateral map value used by the localized
argument, and it retains the factor \(\lambda_1\).

Decompose the exact raw right-hand side as
\[
 \xi_{\rm F}\mathcal Q
   =\mathcal Q_{\rm ph}+\mathcal Q_{\rm ns}
     +\mathcal Q_{\rm cor},
\]
where \(\mathcal Q_{\rm ph}\) is the smooth phase source,
\(\mathcal Q_{\rm ns}\) is supported on
\(\mathfrak S_{\rm gr}^{+}\), and \(\mathcal Q_{\rm cor}\) is
supported in the cutoff corridor.  Normalize the nonseparated source
by the exact definition
\begin{equation}\label{eq:localized-F-normalized-nonseparated-source}
 \mathcal F_{\rm ns}:=\lambda_1^{-1}\mathcal Q_{\rm ns};
 \qquad
 \mathcal Q_{\rm ns}=\lambda_1\mathcal F_{\rm ns}.
\end{equation}
Thus there is no implicit convention about whether \(\mathcal Q\)
already contains the physical clock.  The source calculation gives
\[
 \|\mathcal F_{\rm ns}\|_{\mathbb F_{\rm sc}^{m,\alpha}}
 \leq C\left(
   \mathfrak G_m+\mathfrak B_{{\rm gr},m}
   +\left|\log\frac{\lambda_1}{\lambda_2}\right|
   +\|R_1-R_2\|_{\mathfrak X_{\rm sc}^{m+2,\alpha}}
 \right).
\]
This source has order \(m\), while the asserted map output has order
\(m+1\); hence its contribution is estimated by the one-order Abel
kernel \eqref{eq:effective-time-one-order-kernel}, with the factor
\(\lambda_1\) supplied exactly by
\eqref{eq:localized-F-normalized-nonseparated-source}, not by an
ordinary same-order \(L^1\) bound.  The displayed full coefficient
bounds \(A\in C^{m+1,\alpha}\), \(B,C\in C^{m,\alpha}\) are exactly
the fixed-state hypotheses of that one-derivative parametrix estimate.
Phase profiles are already uniformly \(C^{m+1,\alpha}\), so
\(\mathcal Q_{\rm ph}\) remains in exact \(L^1\)-Duhamel form and gives
the first integral.  The corridor source retains its physical factor
and, since the corridor is separated from the output and remains
noncollapsing, obeys
\[
 \|\mathcal Q_{\rm cor}\|_{H^{-1}(\mathfrak D_{\rm F};E_{\rm F})}
 \leq C\lambda_1\mathfrak D_m^{\rm hyb}.
\]
Applying Lemma~\ref{lem:anchored-buffered-Davies} to that source and to
the commutator gives the final integral, without using either a global
scalar-clock evolution or a same-star restart.

For the parameter family specified in the statement, fix a base
parameter \(p_0\), write \(\mathcal L_p\) for the fixed-state polarized
operator in \eqref{eq:localized-F-cutoff-equation}, and set
\(\mathcal L_0:=\mathcal L_{p_0}\).  Rewrite the nonseparated block
exactly as
\begin{equation}\label{eq:localized-F-fixed-operator-residual}
 \partial_\tau z_p-\lambda_p\mathcal L_0z_p
 =\lambda_p\mathcal F_{{\rm ns},p}^{(p_0)}
   +\mathcal Q_{{\rm ph},p}+\mathcal Q_{{\rm cor},p}
   +\lambda_p[\mathcal L_p,\xi_{\rm F}]w_p,
\end{equation}
where, as an exact identity,
\begin{equation}\label{eq:localized-F-fixed-operator-source-definition}
 \mathcal F_{{\rm ns},p}^{(p_0)}
 :=\lambda_p^{-1}\mathcal Q_{{\rm ns},p}
   +(\mathcal L_p-\mathcal L_0)z_p.
\end{equation}
The right side of
\eqref{eq:localized-F-fixed-operator-source-definition} is regarded as
the single combined residual obtained from the exact nonlinear tension
identity.  Its \(C^1\) regularity in
\(\mathbb F_{\rm sc}^{m,\alpha}\) is the nonlinear tension-map
conclusion of
Lemma~\ref{lem:polarized-HMHF-full-coefficient-ledger}; we do not
differentiate its two displayed summands separately.  At \(p=p_0\),
it reduces to the normalized source
\(\mathcal F_{\rm ns}\) in
\eqref{eq:localized-F-normalized-nonseparated-source}.

If \(v=D_pz_{p_0}[\dot p]\), differentiation of
\eqref{eq:localized-F-fixed-operator-residual} gives the exact equation
\[
 \begin{split}
 \partial_\tau v-\lambda_{p_0}\mathcal L_0v
 ={}&\lambda_{p_0}\Bigl(
    D_p\mathcal F_{{\rm ns},p_0}^{(p_0)}[\dot p]
    +D_p\log\lambda_{p_0}[\dot p]\,
       (\mathcal F_{{\rm ns},p_0}^{(p_0)}+\mathcal L_0z_{p_0})
    \Bigr)\\
 &+D_p\mathcal Q_{{\rm ph},p_0}[\dot p]
  +D_p\mathcal Q_{{\rm cor},p_0}[\dot p]
  +D_p\!\left(
      \lambda_p[\mathcal L_p,\xi_{\rm F}]w_p
    \right)_{p=p_0}[\dot p].
 \end{split}
\]
Thus the \(D\log\lambda\) term multiplies exactly the normalized source
plus \(\mathcal L_0z\); it neither omits nor double-counts a clock
factor.  The phase term and both separated terms are differentiated in
their displayed raw form, so every occurrence of \(\lambda_p\) in
them is retained.  The commutator contains no zeroth-order coefficient,
and the prepared chart calculus places its variation and the corridor
variation in the same separated \(H^{-1}\) source class.

The fixed-operator residual-source clause of
Lemma~\ref{lem:effective-time-endpoint-maximal-regularity}, its fixed
trace extension for the varying actual initial trace, and the
parameter-uniform Davies estimate now give the first-variation formula.
Subtracting it from the secant equation leaves the nonlinear residual,
phase, separated-source, chart, trace, and clock Taylor remainders, all
little-oh in their exact forcing norms.  The same two estimates give
the asserted Fr\'echet remainder, locally at \(p_0\), and uniformly on
compact parameter subsets or under the explicit common-modulus
hypothesis.  No derivative of the top-order polarized zeroth-order
coefficient is used.
\end{proof}

\subsection{Same-order graft propagation}
\label{subsubsec:same-order-graft-propagation}

\begin{lemma}[Same-order pure-graft source from the compact input buffer]
\label{lem:same-order-pure-graft-difference}
Fix \(m\geq4\) and a bounded common prepared package whose low-order
graft margin keeps every interpolated metric uniformly positive.  On
the common marked collar put
\[
 G_i^{\rm mk}:=(\widetilde\iota_i)_*\widetilde G_i=\iota_*G_i,\qquad
 \acute G_i=\eta G_i^{\rm mk}+(1-\eta)S_i ,
\]
Set additionally
\[
 \lambda_{12}(\tau):=\max\{\lambda_1(\tau),\lambda_2(\tau)\}.
\]
Form \(\mathcal E_{{\rm gr},i}\) by
\eqref{eq:pure-graft-again} and
\eqref{eq:pure-graft-normalized}.  Then, after the fixed
right-translation identifications of the two normalized charts,
\begin{equation}\label{eq:pure-graft-buffered-difference}
\begin{split}
 \|\mathcal E_{{\rm gr},1}-\mathcal E_{{\rm gr},2}\|
      _{C_{\rm sc}^{m-2,\alpha}}
 \leq C\Bigg(&
  \mathfrak B_{{\rm gr},m}
  +\left|\log\frac{\lambda_1}{\lambda_2}\right|\\
 &+\|R_1-R_2\|_{\mathfrak X_{\rm sc}^{m+1,\alpha}}
  +\mathfrak F_{{\rm gr},m+1}
  \Bigg)
 \leq C\mathfrak D_m^{\rm hyb}.
\end{split}
\end{equation}
The norm is taken on the receding normalized image of the fixed graft
collar and is zero off that image.  The same estimate holds for sliced
first variations.  At each base state, the Fr\'echet remainder for
difference quotients is \(o(1)\) times the corresponding tangent norm;
the little-\(o\) is uniform when the base state ranges over a fixed
compact parameter trajectory.  No package-wide uniform remainder on an
arbitrary bounded infinite-dimensional set is asserted.
Keeping the physical cutoff scale explicit at arbitrary order gives,
on every dyadic annulus \(A_L=\{L/2<\bar f<4L\}\),
\begin{equation}\label{eq:pure-graft-scale-normalized-all-orders}
 L\,
 \|\mathcal E_{{\rm gr},1}-\mathcal E_{{\rm gr},2}\|
 _{C_{\rm sc}^{m-2,\alpha}(A_L)}
 \leq C\frac{L\lambda_{12}}{\Gamma}\left(
  \mathfrak B_{{\rm gr},m}
  +\left|\log\frac{\lambda_1}{\lambda_2}\right|
  +\|R_1-R_2\|_{\mathfrak X_{\rm sc}^{m+1,\alpha}}
  +\mathfrak F_{{\rm gr},m+1}\right).
\end{equation}
The left side vanishes unless \(L\simeq\Gamma e^\tau\); on that
 support \(L\lambda_{12}/\Gamma\asymp1\).  Thus
\eqref{eq:pure-graft-scale-normalized-all-orders} covers every
finite order \(m\geq4\), not only the order-four uniform closure.

At \(m=4\), retain the physical cutoff scale rather than absorbing it
into the common-package constant.  If
\begin{equation}\label{eq:phase-controlled-graft-input-block}
 \mathfrak B_{{\rm gr},4}
 +\left|\log\frac{\lambda_1}{\lambda_2}\right|
  +\|R_1-R_2\|_{\mathfrak X_{\rm sc}^{5,\alpha}}
  +\mathfrak F_{{\rm gr},5}
 \leq Dd_0
\end{equation}
and \(\lambda_i\simeq e^{-\tau}\), then
\begin{equation}\label{eq:pure-graft-buffered-scale-sharp}
 \supp(\mathcal E_{{\rm gr},1}-\mathcal E_{{\rm gr},2})
 \subset\{c\Gamma e^\tau\leq\bar f\leq C\Gamma e^\tau\},
 \qquad
 \sum_{\ell=0}^2
 |\bar\nabla^\ell
   (\mathcal E_{{\rm gr},1}-\mathcal E_{{\rm gr},2})|_{\bar g}
 \leq C_Dd_0e^{-\tau}.
\end{equation}
The same scale-sharp conclusion holds for sliced first variations.
\end{lemma}

\begin{proof}
On a fixed compact collar the map
\[
 \mathscr G:(G,S)\longmapsto
 2\Ric_{\eta G+(1-\eta)S}-2\eta\Ric_G-2(1-\eta)\Ric_S
\]
is a smooth second-order differential operator from the open set of
uniformly positive \(C^{m,\alpha}\) metric pairs to
\(C^{m-2,\alpha}\) tensors.  Its derivative is bounded on every
bounded subset with a fixed ellipticity margin.  The mean-value formula
therefore gives
\[
 \|\mathscr G(G_1^{\rm mk},S_1)
       -\mathscr G(G_2^{\rm mk},S_2)\|_{C^{m-2,\alpha}}
 \leq C\bigl(
   \|G_1^{\rm mk}-G_2^{\rm mk}\|_{C^{m,\alpha}}
  +\|S_1-S_2\|_{C^{m,\alpha}}\bigr).
\]
No high-order smallness is asserted or needed here; the low
\(C^2\)-smallness is used only for the uniform ellipticity margin.

We do not use a same-regularity smoothness assertion for composition.
With the derivatives used as displayed in
\(\mathfrak B_{{\rm gr},m}\), the standard tame composition estimate
is \(C^1\) from a bounded set of
\(C^{m+1,\alpha}\) diffeomorphisms and inverse diffeomorphisms together
with \(C^{m+2,\alpha}\) tensors into \(C^{m,\alpha}\).  Applied to the
combined marking map, it gives
\[
 \|G_1^{\rm mk}-G_2^{\rm mk}\|_{C^{m,\alpha}}
 \leq C\mathfrak B_{{\rm gr},m}.
\]
Likewise \(S_i=\lambda_i\Theta_i^*\bar g\) and
\(\Phi_i=\varphi_\tau\circ R_i\circ F_i\).  Every composition in the
graft operator is evaluated on \(\mathfrak S_{\rm gr}\).
The same tame estimate, with the fixed smooth soliton data, therefore
bounds the target difference in \(C^{m,\alpha}\) by the scale, \(R\),
and localized \(\mathfrak F_{{\rm gr},m+1}\) terms on the right of
\eqref{eq:pure-graft-buffered-difference}.  No value of \(F_1-F_2\)
on the collapsing core enters this calculation.

For the final normalization we again avoid a same-order composition
claim.  Pull the two physical metrics, the target, and the scalar
cutoff through \(\Phi_i^{-1}\) first.  Ricci naturality gives, exactly,
\[
 \mathcal E_{{\rm gr},i}
 =
 \mathscr G_{\eta_i^{\rm nor}}
   (G_i^{\rm nor},S_i^{\rm nor}),
 \qquad
 (G_i^{\rm nor},S_i^{\rm nor},\eta_i^{\rm nor})
 :=
 (\Phi_i^{-1})^*(G_i^{\rm mk},S_i,\eta).
\]
Here the pullback of a scalar has its usual compositional meaning.
The combined marked-normalized pullback is a tame \(C^1\) map from the
\(C^{m+2,\alpha}\) tensor and \(C^{m+1,\alpha}\) map buffer into
\(C^{m,\alpha}\).  Hence the difference of the three normalized inputs
in \(C^{m,\alpha}\) is bounded by the right-hand side of
\eqref{eq:pure-graft-buffered-difference}.  The operator
\(\mathscr G_\eta(G,S)\), now including its smooth dependence on the
cutoff, is second order and \(C^1\) from that bounded
\(C^{m,\alpha}\) set to \(C^{m-2,\alpha}\).  Its mean-value formula
proves \eqref{eq:pure-graft-buffered-difference} without a
derivative-deficient terminal pullback.  Differentiating these tame
 maps once proves the variation statement.  Continuity of the derivative
and Taylor's formula with integral remainder prove the stated base-local,
and compact-trajectory-uniform, difference-quotient remainder.

For \eqref{eq:pure-graft-buffered-scale-sharp}, repeat the same
mean-value calculation before suppressing the physical cutoff scale.
After the Ricci terms are expanded, every summand in
\(\mathscr G_\eta(G,S)\) either contains a derivative of \(\eta\), or
contains a coefficient which vanishes when \(\eta=0\) or \(\eta=1\).
Thus the physical difference is supported in the fixed transition
collar and its \(\ell\)-th covariant derivative, \(0\leq\ell\leq2\),
has the cutoff-scale factor
\(\Gamma^{-1-\ell/2}\), up to a fixed common-package constant.  The
marked-normalized pullback and tensor rescaling contribute
\(\lambda_i^{1+\ell/2}\).  The \(C^{6,\alpha}\) tensor and
\(C^{5,\alpha}\) map buffers in
\eqref{eq:phase-controlled-graft-input-block} make all of these
composition estimates tame at \(\ell=2\).  Since
\(\lambda_i\simeq e^{-\tau}\) and \(\Gamma\) is fixed, this gives the
pointwise bound in
\eqref{eq:pure-graft-buffered-scale-sharp}.  The two-sided annulus
tracking of the transition collar gives its support statement.
Differentiating this physical-scale calculation proves the identical
claim for sliced first variations.
The identical calculation with \(m\) derivatives retained, and with
the norms expressed in the scale-\(L\) atlas before the factor
\(L^{-1}\) is suppressed, gives
\eqref{eq:pure-graft-scale-normalized-all-orders}.  Each pair of
derivatives of the cutoff costs one factor \(\Gamma^{-1}\), while
tensor normalization contributes at most \(\lambda_{12}\); multiplying by the
two-derivative parabolic cost \(L\) gives precisely
\(L\lambda_{12}/\Gamma\).  The input derivatives recorded in
\(\mathfrak B_{{\rm gr},m}\) make this argument valid at every declared
finite order.
\end{proof}

For later two-state formulae, fix once and for all the modulation
operators
\begin{equation}\label{eq:two-state-modulation-operators}
 \mathscr T_0u=u-\Lie_{\bar\nabla\bar f}u,\qquad
\mathscr T_ju=\Lie_{\chi_\tau W_j}u
\quad(1\leq j\leq8).
\end{equation}

\subsection{Finite-horizon hybrid closure and first variations}
\label{subsubsec:finite-horizon-hybrid-closure}

\begin{lemma}[Finite-horizon sequential common-tail continuity]
\label{lem:finite-horizon-sequential-common-tail}
Fix \(k\geq12\), \(0<\alpha<1\), and a finite normalized endpoint
\(S>\tau_0\).  Let
\(\mathbf z_0\in\Sigma_{\tau_0}^{k+2,\alpha}\), put
\[
 \mathbf z(\tau):=\mathcal S_{\tau,\tau_0}(\mathbf z_0),
\]
and suppose that this coupled trajectory is defined on
\([\tau_0,S]\), remains in one fixed prepared package, and has
positive distance from every exit face on that interval.  Let
\[
 \mathbf z_{0,n}\in\Sigma_{\tau_0}^{k+2,\alpha},
 \qquad
 d_{\rm prep}^{k+2,\alpha}(\mathbf z_{0,n},\mathbf z_0)
 \longrightarrow0.
\]
Then, for all sufficiently large \(n\), the trajectories
\[
 \mathbf z_n(\tau)
 :=\mathcal S_{\tau,\tau_0}(\mathbf z_{0,n})
\]
exist on \([\tau_0,S]\) and remain in one common interior subpackage
of the same numerical prepared package.  In particular, their
exit-face slack has one positive lower bound independent of all
sufficiently large \(n\).

Moreover, there is a finite common restart cover of
\([\tau_0,S]\) such that, in every frozen same-output chart of that
cover, the coordinate tuples of
\[
 \mathbf z,\mathbf z_1,\mathbf z_2,\ldots
\]
carry one common top spatial little-H\"older modulus.  Uniformly in
\(\tau\in[\tau_0,S]\),
\[
 \mathbf z_n(\tau)-\mathbf z(\tau)\longrightarrow0
\]
in the high spatial atlas-supremum coefficient topology carried by
the prepared order-\((k+2,\alpha)\) package.  Thus the tensor
components converge at their order-\((k+2,\alpha)\) spatial tier and
the map and inverse-map components at their corresponding
order-\((k+3,\alpha)\) tier.  No convergence in a stronger top
time-H\"older norm is asserted.

The same conclusion holds, at the buffered spatial orders declared in
the package, for every finite coefficient family obtained from these
states by the prepared graph, inverse, composition, pullback, graft,
moving-support, Gram, gauge, inverse-gauge, marking, and interface
operations.  It also holds in the associated finite-time Bochner and
\(L^1_\tau\) coefficient norms needed for Taylor remainders.
\end{lemma}

\begin{proof}
Let \(m_*>0\) be one quarter of the least, along the base trajectory on
\([\tau_0,S]\), among the ordinary normalized exit-face slacks and the
recorded reserve-to-operative gaps for all witnessed normalized or
physical harmonic conditions.  Let \(K_*\) be a common ceiling for its
order-\((k+2,\alpha)\) prepared package.
The local lifespan in
Proposition~\ref{prop:coupled-local-feedback} depends only on
\(K_*\), the numerical package, and \(m_*\), and is independent of an
individual little-H\"older modulus.  After increasing \(K_*\) by a
fixed amount, choose a corresponding
\(\delta_*>0\), and fix a finite partition
\[
 \tau_0=s_0<s_1<\cdots<s_J=S,
 \qquad
 s_{j+1}-s_j\leq\frac12\delta_*.
\]

We argue inductively.  At \(s_0=\tau_0\), norm convergence in the
little-H\"older prepared space gives the entrance sequence one common
top spatial modulus.  All sufficiently large entrances lie in the
same prepared ball with margin at least \(2m_*\).  Apply
\eqref{eq:coupled-tail-self-map} on \([s_0,s_1]\).  Its constants are
independent of the common datum modulus, and
\eqref{eq:invariant-tail-envelope} therefore gives the entire family
of local solutions one common invariant modulus on that interval.

The local global prepared Lipschitz estimate
\eqref{eq:local-global-prepared-Lipschitz} gives convergence two
orders lower, uniformly in normalized time.  In every normalized
chart, split a representative into finitely many low spatial
frequencies and a high-frequency remainder.  The lower-order
convergence controls the finite-frequency portion uniformly over the
complete atlas, while the common modulus controls the high-frequency
remainder.  It follows that
\[
 \mathbf z_n-\mathbf z\longrightarrow0
\]
uniformly on \([s_0,s_1]\) in the full high spatial coefficient
topology.  Every ordinary normalized exit functional is continuous in
the lower prepared topology.  The witnessed normalized and physical
harmonic conditions persist instead by their recorded
reserve-to-operative gaps and the one-sided lower-stability clauses of
Lemmas~\ref{lem:prepared-harmonic-radius-lower-stability}
and~\ref{lem:finite-physical-harmonic-openness}.  Since the base
ordinary margins and witness gaps are at least \(4m_*\), all
sufficiently large trajectories retain margins and gaps at least
\(2m_*\) throughout this interval.

In particular the endpoint states converge at the spatial order
required for a prepared restart.  The fixed-order prepared coordinate
operations are continuous at that same spatial order by
Lemma~\ref{lem:prepared-chart-calculus}; the two-order loss there is
needed for Fr\'echet differentiation, not for continuity.
Consequently the endpoint sequence, expressed in the chart frozen at
\(s_1\), again has one common datum modulus and, for all sufficiently
large \(n\), retains margin at least \(2m_*\).  The same argument
applies on \([s_1,s_2]\).  Iterating finitely proves existence,
the common-tail assertion, and high spatial convergence through
\(S\).  A single index \(n_0\) suffices because the cover is finite.

The derived coefficient families are treated in their displayed
triangular order.  The compact reconstructed coefficient and the
same-output map variables have already been included in the local
tail argument.  The graph, inverse, pullback, graft, moving-support,
Gram, and interface blocks are finite prepared-chart operations and
therefore preserve the common spatial modulus and strong convergence.
The anchored metric block is parabolic; the gauge, inverse-gauge, and
transported marking equations are triangular; and their coefficients
are among those already shown to converge strongly at the required
buffered spatial orders.  Their integral equations consequently give
the same conclusion for those blocks.  Finite-time integration gives
the corresponding \(L^1_\tau\) history convergence.

No time derivative is inserted into the common modulus.  The required
Bochner derivatives are recovered from the exact equations after the
spatial coefficients and solutions have converged, exactly as in the
completion argument following
\eqref{eq:invariant-tail-envelope}.
\end{proof}

\begin{lemma}[Same-order weighted normalized-tensor difference]
\label{lem:same-order-normalized-h-difference}
Fix \(k\geq12\), \(0<\alpha<1\), \(N\geq0\), a finite normalized
horizon \([\tau_0,S]\), and an interval
\(I=[s,s+\ell]\subset[\tau_0,S]\), \(0<\ell\leq1\).  Let two admissible
prepared evolutions lie in one common-margin ball with one-state
coefficient bounds two orders above \(k\).  Put
\[
 w=h_1-h_2,\qquad H_i=\rho_\tau h_i,\qquad
 V=H_1-H_2 .
\]
Separate the one-time blocks
\begin{equation}\label{eq:finite-horizon-memory-block}
 \mathfrak M_{k,S}^{\rm mem}:=
 \mathfrak G_k+\mathfrak B_{{\rm gr},k}
 +\mathfrak F_{{\rm gr},k+1}^{\rm hs}
\end{equation}
from the locally restarted distance
\begin{equation}\label{eq:finite-horizon-local-distance}
 \begin{split}
 \mathfrak D_k^{\rm loc}:={}&
 |t_1-t_2|
  +\left|\log\frac{\lambda_1}{\lambda_2}\right|
  +\mathfrak M_{{\rm gr},k}
  +\|R_1-R_2\|_{\mathfrak X_{\rm sc}^{k+2,\alpha}}\\
 &+\mathfrak F_{k+1}^{\rm glob}
  +\mathfrak F_{{\rm gr},k+1}
  +\|h_1-h_2\|_{\mathfrak T_{{\rm sc},N}^{k,\alpha}} .
 \end{split}
\end{equation}
Here \(\mathfrak F_{{\rm gr},k+1}^{\rm hs}\) is retained in
\(\mathfrak M_{k,S}^{\rm mem}\) because its larger-star homogeneous
trace and separated corridor forcing are the typed map input of the
graft operator.  The full \(\mathfrak F_{{\rm gr},k+1}\) block,
including its nonseparated zero-trace part, is locally propagated.
Likewise, \(\mathfrak F_{k+1}^{\rm glob}\) is a
genuinely global parabolic unknown and is restarted from its actual
global trace at \(s\); no localized trace is substituted for it.  The
order-\((k+2)\) \(R\)-term is an auxiliary ODE block, supplied by the
map regularity in the prepared input and used only to type the
nonseparated \(F\)-source.  Define
\begin{align}
 d_{k,0}^{\rm mem}
 &:=
 d_{{\rm ext},k,0}^{+}
 +d_{{\rm ext},-1,0}^{\rm corr}
 +d_{{\rm Ggr},k+2,0}^{0}\notag\\
 &\quad
 +\mathfrak B_{{\rm gr},k}(\tau_0)
 +\mathfrak G_k(\tau_0)+\mathfrak I_{\rm in}(\tau_0)\notag\\
 &\quad
 +\left|\log\frac{\lambda_1}{\lambda_2}(\tau_0)\right|
  +|t_1(\tau_0)-t_2(\tau_0)|
  +d_{{\rm Fgr},k+1,0}^{++},
 \label{eq:finite-horizon-memory-entrance}\\
 \mathscr H_{k,S}(\tau)^2
 &:=
 \bigl(d_{k,0}^{\rm mem}\bigr)^2
 +\int_{\tau_0}^{\tau}\bigl(\notag\\[-2pt]
 &\hspace{35mm}
   \bigl(\mathfrak D_k^{\rm hyb,full}(q)\bigr)^2
   +\|V(q)\|_{H^1_\nu}^2\bigr)dq .
 \label{eq:finite-horizon-history-functional}
\end{align}
Thus \(\mathscr H_{k,S}\) is monotone and is anchored at the original
entrance face; it is not reconstructed from a smaller-collar trace at
time \(s\).
If
\[
 |\delta c|\leq
 C\bigl(\|V\|_{H^1_\nu}+\mathfrak D_k^{\rm hyb,full}\bigr),
\]
then
\begin{equation}\label{eq:same-order-normalized-h-difference}
 \|w(\tau)\|_{\mathfrak T_{{\rm sc},N}^{k,\alpha}}^2
 \leq
 C\bigl(\mathfrak D_k^{\rm loc}(s)\bigr)^2
 +C\mathscr H_{k,S}(\tau)^2
 +C\int_s^\tau\left(
   \bigl(\mathfrak D_k^{\rm loc}(q)\bigr)^2
   +\|V(q)\|_{H^1_\nu}^2\right)dq .
\end{equation}
On a fixed finite horizon the constant is independent of the dyadic
annulus.  The same estimate holds for sliced first variations.  For a
difference quotient minus its sliced linearization, the corresponding
forced estimate holds with the mixed Taylor trace, Bochner, and
history forcing norms added to the right side.  No little-o assertion
relative to the entrance-increment norm is made for those forcing
terms in this lemma; that size is established after the finite-horizon
two-state estimate and the sequential common-tail argument are both
available.
\end{lemma}

\begin{proof}
Subtract the two exact controlled equations before suppressing any
outer term.  In the notation of
\eqref{eq:controlled-h}, the result is
\begin{equation}\label{eq:finite-horizon-exact-h-difference}
 \begin{split}
 \partial_\tau w={}&
 \A w+\Q(h_1)-\Q(h_2)
 +\sum_{j=0}^8c_{1,j}\mathscr T_jw\\
 &+\sum_{j=0}^8\delta c_j
      \bigl(\mathcal Y_{j,\tau}(\mathbf z_1)
             +\mathscr T_jh_2\bigr)
 +\sum_{j=0}^8c_{2,j}
      \bigl(\mathcal Y_{j,\tau}(\mathbf z_1)
            -\mathcal Y_{j,\tau}(\mathbf z_2)\bigr)
 +(\E_1-\E_2).
 \end{split}
\end{equation}
Use the exact mean-value polarization, rather than treating a
principal-coefficient difference as a prescribed source.  In a fixed
core chart, and in a scale-one soliton chart on the end,
\eqref{eq:Q-exact} gives
\[
 \partial_\tau w-\mathsf a_1^{ij}\bar\nabla_i\bar\nabla_jw
 -\bar\nabla_{\mathsf V_1}w
 =
 \mathsf B_{12}*\bar\nabla w+\mathsf C_{12}*w
 +\mathsf F_{\rm fb}+\mathsf L_{\rm coup}[\delta\mathbf Z].
\]
The coarse \(C^2\) box makes \(\mathsf a_1\) uniformly elliptic.  The
apparently second-order term is polarized exactly as
\[
 (\mathsf a_1-\mathsf a_2)*\bar\nabla^2h_2
 =
 \left(\int_0^1
  D\mathsf a_{h_2+\theta w}[w]\,d\theta\right)
 *\bar\nabla^2h_2,
\]
and is therefore a zeroth-order coefficient applied to \(w\), not an
endpoint forcing.  This is the same identity recorded later in
\eqref{eq:two-state-principal-polarization}.
Here \(\delta\mathbf Z\) denotes the full physical, normalized, map,
gauge, marking, scale, and clock difference;
\(\mathsf L_{\rm coup}\) collects the remaining bounded linearized
couplings, while \(\mathsf F_{\rm fb}\) consists only of the smooth
same-order profiles multiplied by \(\delta c\).

On the end, first conjugate by the first-state phase flow and then use
the soliton rescaling
\eqref{eq:phase-rescaled-H}.  The resulting cylinders are uniformly
parabolic in a scale-one metric, exactly as in
\eqref{eq:one-state-phase-rescaled-equation}; this is a normalized
tensor clock and is not the raw harmonic-map clock
\eqref{eq:source-adapted-clock}.  Phase distortion is bounded by the
\(L^1\) phase budget.  The dyadic factor \(L^{-N}\) changes by at most
a fixed constant on enlarged overlapping cylinders, so the chartwise
same-order estimates sum in
\(\mathfrak T_{{\rm sc},N}^{k,\alpha}\).

We now estimate only the \(h\)-equation in its natural normalized
clock.  The exact polarization above defines a uniformly parabolic
nonautonomous operator \(\mathscr L_{h,12}(\tau)\) on each
phase-conjugated soliton cylinder.  Its principal coefficients are
those of the first state; all polarized quasilinear terms involving
\(w\) are lower-order coefficients of this same operator.  The
order-\(k+2\) one-state buffer controls those coefficients uniformly on
the finite interval.

The endpoint estimate must now be closed simultaneously with the local
non-\(h\) blocks; an integral of a merely
\(C^{k-2,\alpha}\) source is not, by itself, a
\(C^{k,\alpha}\) endpoint estimate.  Put
\[
 \mathfrak W_k(q):=
 \|w(q)\|_{\mathfrak T_{{\rm sc},N}^{k,\alpha}},
 \qquad
 \mathfrak Z_k(q):=\mathfrak D_k^{\rm loc}(q)-\mathfrak W_k(q).
\]
The exterior, graft-input, and homogeneous/separated localized \(F\)
blocks are not restarted at \(s\).  The one-time components of
\eqref{eq:inner-terminated-exterior-DeTurck},
\eqref{eq:compact-graft-buffer-finite}, and
\eqref{eq:localized-graft-F-coarse-memory} retain their original
entrance traces.  Their scalar off-diagonal kernels are bounded on the
present finite effective-time range.  The global \(F\)-block is
excluded from this memory and is propagated below from its actual
global trace on each short interval.  The assumed feedback estimate,
Cauchy--Schwarz, the effective-time clock bound
\eqref{eq:graft-clock-source-amplitude}, and the triangular dependence
of the localized \(F\)-block give the coupled scalar inequality
\begin{equation}\label{eq:finite-horizon-coupled-memory-inequality}
\begin{split}
 \mathfrak M_{k,S}^{\rm mem}(\tau)
 \leq{}& C_S d_{k,0}^{\rm mem}\\
 &+C_S\int_{\tau_0}^{\tau}
 \left(
  \mathfrak M_{k,S}^{\rm mem}(q)
  {}+\mathfrak D_k^{\rm hyb,full}(q)
  {}+\|V(q)\|_{H^1_\nu}
 \right)dq .
\end{split}
\end{equation}
Here the first term retains all original homogeneous entrance traces;
the current global \(F\)-value enters only through
\(\mathfrak D_k^{\rm hyb,full}\), while its original trace belongs to
\(\mathfrak D_k^{\rm loc}(\tau_0)\).  Thus every localized memory term
remains anchored at the original trace \(\tau_0\).  Scalar
Gronwall, followed by Cauchy--Schwarz on
\([\tau_0,\tau]\subset[\tau_0,S]\), yields
\[
 \mathfrak M_{k,S}^{\rm mem}(\tau)
 \leq C_S\left(
  d_{k,0}^{\rm mem}
  {}+\left[
   \int_{\tau_0}^{\tau}
   \left(
    \bigl(\mathfrak D_k^{\rm hyb,full}(q)\bigr)^2
    {}+\|V(q)\|_{H^1_\nu}^2
   \right)dq
  \right]^{1/2}\right).
\]
By the definition
\eqref{eq:finite-horizon-history-functional}, this proves
\begin{equation}\label{eq:finite-horizon-memory-by-history}
 \mathfrak M_{k,S}^{\rm mem}(q)
 \leq C_S\mathscr H_{k,S}(q),
 \qquad \tau_0\leq q\leq S .
\end{equation}
In particular, no value on
\(\mathcal W_{\rm gr}^{0}\setminus\mathcal W_{\rm gr}^{+}\), and no
pre-\(s\) lateral forcing, has been replaced by the smaller current
graft block.  The algebraic interface identity gives
\begin{equation}\label{eq:finite-horizon-interface-memory-splitting}
 \mathfrak I_{\rm in}
 \leq C\bigl(\mathfrak D_k^{\rm loc}
             +\mathfrak M_{k,S}^{\rm mem}\bigr),
 \qquad
 \mathfrak D_k^{\rm hyb,full}
 \leq C\bigl(\mathfrak D_k^{\rm loc}
             +\mathfrak M_{k,S}^{\rm mem}\bigr).
\end{equation}

Let \(\delta_{\rm Ab}>0\) denote the short-interval constant
\(\delta_0\) from
Lemma~\ref{lem:effective-time-endpoint-maximal-regularity}.  Recenter
the simultaneous contraction from the proof of
Proposition~\ref{prop:coupled-local-feedback} at time \(s\), but only
for the genuinely local variables.  Its mixed difference estimate
\eqref{eq:coupled-map-difference-forcing}, together with the
one-order Abel block
\eqref{eq:source-atlas-one-order-Abel-block}, is valid on every
\(J=[s,s+\delta]\subset I\) with
\(0<\delta\leq\delta_{\rm Ab}\), with constants uniform on the present
finite horizon.  Propagate the clock, scale, and the auxiliary
order-\((k+2)\) \(R\) ODE block, followed by the local gauge,
inverse-gauge, and marking blocks, in triangular order.
The interface trace is recovered from
\eqref{eq:finite-horizon-interface-memory-splitting}; the exterior and
graft-input blocks enter through one-time memory.

Split the global \(F\)-difference on \(J\) into its homogeneous
evolution from the actual global value
\(F_1(s)-F_2(s)\), the smooth finite-rank feedback part, and a
zero-trace nonseparated part.  By
\eqref{eq:coupled-map-difference-forcing}, the last has
\(\mathbb F_{\rm sc}^{k,\alpha}(J)\)-norm bounded by
\[
 C\sup_J\bigl(
   \mathfrak D_k^{\rm loc}+\mathfrak M_{k,S}^{\rm mem}\bigr).
\]
In particular, the coefficient difference multiplying
\(\bar\nabla^2F_2\) is measured at source order \(k\), with the
auxiliary \(R\)-derivative controlling the target connection.
Equation~\eqref{eq:source-atlas-one-order-Abel-block} makes its
contribution to \(\mathfrak F_{k+1}^{\rm glob}\) at most
\(C\delta^{1/2}\) times this source norm.  The same decomposition
treats the nonseparated portion of the localized \(F\)-block; its
homogeneous larger-star datum and its separated corridor source stay
in \(\mathfrak M_{k,S}^{\rm mem}\).
Thus
\begin{equation}\label{eq:finite-horizon-non-h-top-ledger}
 \begin{split}
 \sup_{s\leq q\leq\tau}\mathfrak Z_k(q)
 \leq{}&
 C\mathfrak D_k^{\rm loc}(s)
 +C\sup_{s\leq q\leq\tau}\mathfrak M_{k,S}^{\rm mem}(q)\\
  &+C\int_s^\tau
    \bigl(\mathfrak D_k^{\rm loc}(q)+|\delta c(q)|\bigr)\,dq\\
  &+C\bigl(\delta^{\alpha/4}+\delta^{1/2}\bigr)
    \sup_{s\leq q\leq\tau}
    \bigl(\mathfrak W_k(q)+\mathfrak Z_k(q)
          +\mathfrak M_{k,S}^{\rm mem}(q)\bigr).
 \end{split}
\end{equation}
Here every nonfeedback top-order source has first been placed in its
declared time-supremum norm.
\eqref{eq:dynamic-column-high-order-Lipschitz} supplies the moving
columns in \(C_{\rm sc}^{k-2,\alpha}\).  The last term in
\eqref{eq:finite-horizon-non-h-top-ledger} contains both the
two-order-lowered zero-trace coupling and the one-derivative
global/local \(F\) Abel block.  After \(\delta\) is fixed, its
\(\mathfrak Z_k\) summand is absorbed on the left and its memory
summand is merged with the preceding memory term; only
\(\mathfrak W_k\) is returned to the \(h\)-estimate.

We spell out the graft contribution to that block.  Define
\[
 \Delta\mathcal E_{\rm gr}
 :=\mathcal E_{{\rm gr},1}-\mathcal E_{{\rm gr},2},
 \qquad
 \lambda_*(q):=\max\{\lambda_1(q),\lambda_2(q)\}.
\]
For a phase-rescaled cylinder meeting the moving graft support, let
\(L(q)\) be the dyadic scale of its enlarged annulus at source time
\(q\); set the corresponding integrand to zero when the enlarged
annulus misses the support.  The common scale bracket and the support
statement in \eqref{eq:pure-graft-scale-normalized-all-orders} give
\[
 L(q)\simeq\Gamma e^q,\qquad
 L(q)\frac{\lambda_*(q)}{\Gamma}\asymp1 .
\]
Apply the weighted prepared Schauder estimate to the \(h\)-specific
evolution family
\(\mathcal U_{h,L}(\tau,q)\) of \(\mathscr L_{h,12}\) on that cylinder.
This family is not reused for the raw \(F\)- or physical-metric
equations.  The two-derivative Duhamel gain costs \(L(q)\), and the
all-order scale-sharp graft estimate gives the correctly typed bound
\begin{equation}\label{eq:h-specific-scaled-graft-Duhamel}
 \begin{split}
 &\left\|\int_s^\tau\mathcal U_{h,L}(\tau,q)
       \Delta\mathcal E_{\rm gr}(q)\,dq
 \right\|_{C_{\rm sc}^{k,\alpha}}\\
 &\quad\leq C\sup_{s\leq q\leq\tau}
   L(q)\|\Delta\mathcal E_{\rm gr}(q)\|
       _{C_{\rm sc}^{k-2,\alpha}}\\
 &\quad\leq C\sup_{s\leq q\leq\tau}\left(
   \mathfrak B_{{\rm gr},k}
   +\left|\log\frac{\lambda_1}{\lambda_2}\right|
   +\|R_1-R_2\|_{\mathfrak X_{\rm sc}^{k+1,\alpha}}
   +\mathfrak F_{{\rm gr},k+1}\right)(q)\\
 &\quad\leq C\sup_{s\leq q\leq\tau}
   \bigl(\mathfrak D_k^{\rm loc}(q)
         +\mathfrak M_{k,S}^{\rm mem}(q)\bigr).
 \end{split}
\end{equation}
This is the scale cancellation which prevents a hidden loss from the
\(L^{-1}\)-diffusion.  Cylinders disjoint from the graft have no such
source.  Every other nonfeedback source is bounded in
\(L^\infty_qC_{\rm sc}^{k-2,\alpha}\) by the blocks in
\eqref{eq:finite-horizon-non-h-top-ledger}.
For the finite-rank feedback source in
\eqref{eq:finite-horizon-exact-h-difference}, write
\(\mathsf F_{{\rm fb},L}(q)\) for its phase-conjugated,
scale-\(L\) representative.  Uniformity in \(L\) here uses the exact
prepared-graph structure; it is not obtained by discarding the factor
\(L^{-N}\) for an arbitrary element of
\(\mathfrak T_{{\rm sc},N}^{k+2,\alpha}\).  Indeed,
\eqref{eq:prepared-graph-shifted-decomposition} writes the graph as the
sum of the compact-host term
\[
 \lambda^{-1}(\Phi^{-1})^*(\eta\,\iota_*G)
\]
and the outer term \(\zeta J^*\bar g-\bar g\).  The latter has
unweighted scale-normalized bounds through every fixed order carried
by the prepared map buffer, by
Lemma~\ref{lem:shifted-annulus-soliton-calculus}.  The former is
supported in
\[
 \Phi_q(\operatorname{supp}\eta)
 \subset\{\bar f\leq C\Gamma e^q\}
 \subset\{\bar f\leq C\Gamma e^S\},
\]
so on the fixed horizon its weighted order-\((k+2)\) bound gives the
corresponding unweighted bound, at the allowed cost
\((C\Gamma e^S)^N\).  The explicit effective-column formulas
\eqref{eq:effective-column-zero}--\eqref{eq:effective-column-j} and the
finite-order shifted-annulus calculation in the proof of
Proposition~\ref{prop:effective-column-tails}, rerun with the present
order-\((k+2)\) buffer, give the same unweighted conclusion at every
fixed order \(m\leq k\) used here; no endpoint-independent all-order
column assertion is invoked.  Finally, the
all-order FIK scale-symbol bounds and
Lemma~\ref{lem:mode-growth}, together with the support of \(\chi_q\),
show that \(\mathscr T_0\) preserves the outer scale-symbol bounds and
that \(\mathscr T_j\), \(1\leq j\leq8\), is confined to
\(\{\bar f\leq3e^q\}\).  Thus the order-\((k+2)\) buffer, the uniform
phase distortion, and the bounded change of dyadic labels give, for
every integer \(0\leq m\leq k\),
\[
 \|\mathsf F_{{\rm fb},L}(q)\|_{C_{\rm sc}^{m,\alpha}}
 \leq C_{m,S}|\delta c(q)|
\]
uniformly in \(L\) and in the selected cylinder.  Apply the preceding
same-order evolution estimate
\eqref{eq:endpoint-MR-same-order-evolution} at the same order \(m\) to the
\(h\)-specific family \(\mathcal U_{h,L}\) in its phase-rescaled
scale-one clock.  The effective-time ceiling is fixed on the selected
cylinder.  Hence
\begin{equation}\label{eq:h-feedback-L1-Duhamel}
 \left\|
  \int_s^\tau
   \mathcal U_{h,L}(\tau,q)
   \mathsf F_{{\rm fb},L}(q)\,dq
  \right\|_{C_{\rm sc}^{m,\alpha}}
 \leq C_{m,S}\int_s^\tau|\delta c(q)|\,dq,
 \qquad 0\leq m\leq k .
\end{equation}
Here \(C_{m,S}\) may depend on the fixed finite horizon, in addition
to \(m,\alpha,N\) and the declared common coefficient and phase
package, but is independent of \(L\) and of the selected cylinder.
No endpoint-independent unweighted high-order estimate is asserted in
this finite-horizon step.  Since every dyadic label satisfies
\(L\geq\Gamma\geq1\), taking the compact-core maximum and the dyadic
supremum controls simultaneously the unweighted order-two tier and the
weighted order-\(m\) tier of
\(\mathfrak T_{{\rm sc},N}^{m,\alpha}\); comparable transported labels
change \(L^{-N}\) by only a fixed factor.  Taking \(m=k\) in this
estimate gives
\begin{equation}\label{eq:finite-horizon-h-top-ledger}
 \begin{split}
 \sup_{s\leq q\leq\tau}\mathfrak W_k(q)
 \leq{}&
 C\mathfrak D_k^{\rm loc}(s)
 +C\sup_{s\leq q\leq\tau}\mathfrak Z_k(q)
 +C\sup_{s\leq q\leq\tau}\mathfrak M_{k,S}^{\rm mem}(q)\\
 &+C\int_s^\tau
   \bigl(\mathfrak D_k^{\rm loc}(q)+|\delta c(q)|\bigr)\,dq .
 \end{split}
\end{equation}
This constructs no evolution family on a heterogeneous product of
physical, source-adapted, and normalized clocks.

Choose one
\(0<\delta_\#\leq\min\{1,\delta_{\rm Ab}\}\), depending only on the
finite-horizon coefficient package, so that both the
\(\delta^{\alpha/4}\) mixed zero-trace term and the
\(\delta^{1/2}\) global/local \(F\) Abel term in
\eqref{eq:finite-horizon-non-h-top-ledger} are absorbed after
\eqref{eq:finite-horizon-h-top-ledger} is substituted.  For
\(\delta\leq\delta_\#\) this gives the history-aware triangular estimate
\begin{equation}\label{eq:finite-horizon-triangular-top-ledger}
 \begin{split}
 \sup_{s\leq q\leq\tau}\mathfrak D_k^{\rm loc}(q)
 \leq C\Bigg(&
  \mathfrak D_k^{\rm loc}(s)
  +\sup_{s\leq q\leq\tau}\mathfrak M_{k,S}^{\rm mem}(q)\\
 &+\int_s^\tau
   \bigl(\mathfrak D_k^{\rm loc}(q)+|\delta c(q)|\bigr)\,dq\Bigg).
 \end{split}
\end{equation}
By \eqref{eq:finite-horizon-memory-by-history}, monotonicity of
\(\mathscr H_{k,S}\), the feedback hypothesis, and
\eqref{eq:finite-horizon-interface-memory-splitting},
\[
 \left(\int_s^\tau|\delta c(q)|\,dq\right)^2
 \leq C\delta\int_s^\tau\left(
  \bigl(\mathfrak D_k^{\rm loc}(q)\bigr)^2
  +\|V(q)\|_{H^1_\nu}^2\right)dq
 +C\delta^2\mathscr H_{k,S}(\tau)^2 .
\]
Squaring \eqref{eq:finite-horizon-triangular-top-ledger} proves
\eqref{eq:same-order-normalized-h-difference} on a short interval.
A finite subdivision of \(I\), followed by the discrete
Volterra--Gronwall inequality and monotonicity of
\(\mathscr H_{k,S}\), proves it for every \(0<\ell\leq1\).

For sliced first variations, differentiation of the simultaneous
system gives the same estimates.

For a difference quotient minus its sliced linearization, apply
Taylor's integral formula block by block.  The clock and scale blocks
are finite-dimensional ODEs; \(R\) and \(F\) are differentiated in
their right-translated source-adapted charts; the anchored metric,
gauge, inverse-gauge, and marking blocks use the buffered equations in
\eqref{eq:compact-graft-buffer-finite}; the interface block uses the
differentiated algebraic \(P\)-identity in
\eqref{eq:interface-block-controlled}; and the normalized tensor,
moving-column, and pure-graft blocks use, respectively, the prepared
chart calculus and
Lemma~\ref{lem:same-order-pure-graft-difference}.  These calculations
give the forced analogue of
\eqref{eq:same-order-normalized-h-difference}, with the corresponding
Taylor trace, Bochner, and history norms on its right side.

At this stage we do not declare those Taylor terms to be little-o
relative to the entrance-increment norm.  Their required size follows in
Proposition~\ref{prop:two-state-prepared-evolution} by combining the
two-state lower-order bound with
Lemma~\ref{lem:finite-horizon-sequential-common-tail}.  This ordering
avoids using finite-horizon high-spatial continuity before it has been
proved.
\end{proof}

\begin{theorem}[Gauge-covariant finite-horizon hybrid propagation]
\label{thm:finite-horizon-weighted-parabolic}
Fix \(k\geq12\), \(0<\alpha<1\), \(N\geq0\), and a common-margin
prepared ball.  Let two admissible prepared evolutions with initial
states \(\mathbf z_{1,0},\mathbf z_{2,0}\) be defined on
\([\tau_0,S]\), \(S<\infty\), and put
\[
 H_i=\rho_\tau h_i,\qquad V=H_1-H_2,\qquad
 \delta c=(a_1-a_2,b_1-b_2),\qquad
 d_0=\|\mathbf z_{1,0}-\mathbf z_{2,0}\|
       _{\mathscr X_{\rm prep}^{k+2,\alpha}}.
\]
Define the complete typed entrance distance
\begin{equation}\label{eq:finite-horizon-typed-entrance-distance}
 \begin{split}
 d_{k,0}^{\rm typ}:={}&
 \mathfrak D_k^{\rm hyb,full}(\tau_0)
 +d_{{\rm ext},k,0}^{+}
 +d_{{\rm ext},-1,0}^{\rm corr}
 +d_{{\rm Ggr},k+2,0}^{0}\\
 &+d_{{\rm Fgr},k+1,0}^{++}
 +\|R_1(\tau_0)-R_2(\tau_0)\|
      _{\mathfrak X_{\rm sc}^{k+2,\alpha}} .
 \end{split}
\end{equation}
Every term is a trace of the prepared input on its declared global or
larger-buffer domain, and
\[
 d_{k,0}^{\rm typ}
 \leq C\|\mathbf z_{1,0}-\mathbf z_{2,0}\|
          _{\mathscr E_{\rm prep}^{k+2,\alpha}}
 \leq C d_0.
\]
Assume the one-state coefficients are bounded two orders above the
output order on the largest buffered physical and marked sets, as they
are on the prepared ball, and suppose
\begin{equation}\label{eq:finite-horizon-feedback-hypothesis}
 |\delta c|
 \leq C\bigl(\|V\|_{H^1_\nu}
              +\mathfrak D_k^{\rm hyb,full}\bigr).
\end{equation}
Then
\begin{equation}\label{eq:finite-horizon-weighted-parabolic}
 \mathfrak D_k^{\rm hyb,full}(\tau)^2
 \leq C_S\bigl(d_{k,0}^{\rm typ}\bigr)^2
 +C_S\int_{\tau_0}^\tau
   \left(\mathfrak D_k^{\rm hyb,full}(q)^2
         +\|V(q)\|_{H^1_\nu}^2\right)\,dq,
 \qquad \tau_0\leq\tau\leq S.
\end{equation}
The constant \(C_S\) depends only on \(S\), the common prepared bounds,
and the fixed exterior certificate, and is independent of the outer
dyadic annulus.  The normalized and map subproblems may be estimated
on unit subintervals without derivative loss; the physical exterior
block is recovered once from \(\tau_0\) by
Lemma~\ref{lem:inner-terminated-exterior-DeTurck}, not by a cyclic
spatial restart.
\end{theorem}

\begin{proof}
Let \(\delta_{\rm Ab}>0\) denote the short-interval constant from
Lemma~\ref{lem:effective-time-endpoint-maximal-regularity}, and let
\(0<\delta_\#\leq\min\{1,\delta_{\rm Ab}\}\) be the finite-horizon
absorption length fixed in the proof of
Lemma~\ref{lem:same-order-normalized-h-difference}.  Fix first a
normalized subinterval
\[
 I=[s,s+\ell]\subset[\tau_0,S],\qquad0<\ell\leq\delta_\#.
\]
Cover the model core by finitely many fixed harmonic charts and the AC
end by uniformly enlarged scale-one charts for \(L^{-1}\bar g\).
Lemma~\ref{lem:prepared-chart-calculus} gives uniform coefficients for
inverse, pullback, composition, and moving support.  Estimate the
components in their natural clocks.

For \(F_1-F_2\), use
the source-adapted one-order Abel estimate
\eqref{eq:source-atlas-one-order-Abel-block} in the frozen clocks
\eqref{eq:source-adapted-clock}.  The global map is restarted from its
actual global trace \(F_1(s)-F_2(s)\); this is legitimate because the
equation and norm are global and create no lateral boundary datum.
Only the larger-star homogeneous datum and separated corridor source
of the localized graft map remain anchored at \(\tau_0\) through
\eqref{eq:localized-graft-F-coarse-memory}.  On \(I\), recenter the
nonseparated inhomogeneous contribution with zero trace.  Differences of
\(\lambda/\lambda_s^\circ\) are retained among the
principal-coefficient differences.  Each such difference
multiplies the bounded \(k+2\) one-state jet and belongs to
\(\mathbb F_{\rm sc}^{k,\alpha}\), exactly as recorded in
\eqref{eq:coupled-map-difference-forcing}; the auxiliary
order-\((k+2)\) \(R\)-difference types the target-connection term.
The Abel kernel supplies the missing one spatial derivative and the
factor \(C\ell^{1/2}\).  Separate the terms containing \(\delta c\).
They are finite sums of uniformly \(C^{k+1,\alpha}\) coefficient
profiles multiplied by the scalar functions \(\delta c_j(q)\).  If
\(\mathcal U_F(\tau,q)\) is the evolution family for the frozen
linearized map equation, its same-order bound gives
\begin{equation}\label{eq:feedback-L1-Duhamel}
 \left\|
  \int_s^\tau\mathcal U_F(\tau,q)
       \bigl[\delta c_j(q)\mathcal P_j(q)\bigr]\,dq
 \right\|_{\mathfrak X_{\rm sc}^{k+1,\alpha}}
 \leq C\int_s^\tau|\delta c(q)|\,dq .
\end{equation}
Thus the energy information is used in its legitimate form:
\[
 \int_s^\tau|\delta c(q)|\,dq
 \leq C\ell^{1/2}
 \left(\int_s^\tau\|V(q)\|_{H^1_\nu}^2\,dq\right)^{1/2}
 +C\int_s^\tau\mathfrak D_k^{\rm hyb,full}(q)\,dq .
\]
No \(L^\infty_\tau H^1_\nu\) bound is inserted into the Schauder
forcing norm.

For the closed metrics the two physical clocks must be retained.  On the
fixed exterior domain set
\[
 \widehat G_i(\tau)=\widetilde G_i(t_i(\tau)),\qquad
 u=\widehat G_1-\widehat G_2 .
\]
If \(\mathscr R_{\rm D}\) is the Ricci--DeTurck right-hand side in the
fixed common reference gauge, then
\[
 \partial_\tau\widehat G_i
 =\lambda_i\mathscr R_{\rm D}(\widehat G_i).
\]
Consequently the exact difference equation has the form
\begin{equation}\label{eq:finite-horizon-unequal-clock-DeTurck}
 \partial_\tau u-\lambda_1\mathcal A^{ab}\nabla_a\nabla_bu
 =
 \lambda_1(\mathcal B*\nabla u+\mathcal C*u)
 +\mathcal F_{\rm clk},
 \qquad
 \mathcal F_{\rm clk}
 =(\lambda_1-\lambda_2)\mathscr R_{\rm D}(\widehat G_2).
\end{equation}
On a fixed normalized interval \([\tau_0,S]\), scale comparability and
the one-state \(C^{k+2,\alpha}\) bounds give
\[
 \|\mathcal F_{\rm clk}\|_{C^{k,\alpha}(E^{++})}
 \leq C_S\left|\log\frac{\lambda_1}{\lambda_2}\right|
 \leq C_S\mathfrak D_k^{\rm hyb,full}.
\]
Lemma~\ref{lem:inner-terminated-exterior-DeTurck}, applied to
\eqref{eq:finite-horizon-unequal-clock-DeTurck}, therefore propagates
\(\mathfrak G_k\); its two homogeneous entrance terms satisfy
\[
 d_{{\rm ext},k,0}^{+}
 +d_{{\rm ext},-1,0}^{\rm corr}\leq C d_{k,0}^{\rm typ} .
\]
The required inner trace is not estimated by
iterating to a larger physical buffer: it is the exact normalized
graph trace \eqref{eq:inner-terminal-graph-identity}.  The triangular
DeTurck ODE and inverse ODE, written in the same normalized time,
contain the identical integrable clock mismatch and propagate
\(\mathfrak M_{{\rm gr},k}\) at order \(k-1\).

On the separated graft collars we use the larger input buffer and
retain the exterior trace.  Apply
Lemma~\ref{lem:compact-graft-buffer-propagation} at order \(k\).
Its one-time estimate includes the asynchronous clock term from
\eqref{eq:finite-horizon-unequal-clock-DeTurck}, the separated
interface/exterior memory, and the triangular gauge, inverse-gauge,
and transported-marking equations.  Since
 \(d_{{\rm gr},k,0}\leq Cd_{k,0}^{\rm typ}\), the endpoint clock
 contribution is
first bounded by
\eqref{eq:graft-clock-source-amplitude}; Cauchy--Schwarz is applied to
\(\int|\delta a|\) and to the separated memory integrals, and no
top-order clock Duhamel estimate is used.  This gives
\begin{equation}\label{eq:graft-buffer-finite-propagation}
 \mathfrak B_{{\rm gr},k}(\tau)^2
 \leq C_S\bigl(d_{k,0}^{\rm typ}\bigr)^2+
 C_S\int_{\tau_0}^{\tau}\left(
       \bigl(\mathfrak D_k^{\rm hyb,full}(q)\bigr)^2
       +\|V(q)\|_{H^1_\nu}^2\right)\,dq .
\end{equation}
This is exactly where the two extra prepared derivatives are used.
It is a larger-input/smaller-output graft-collar estimate for the
closed flow, is independent of the \(h\)-equation, and has neither an
uncontrolled next-collar norm nor an omitted lateral datum on its
right side.

The transported-marking identity
\eqref{eq:buffered-marking-covariance} then preserves the exact
physical graft on the smallest marked set.  The remaining variables
\(t,\log\lambda,R\) obey locally Lipschitz ODEs and are controlled by
direct integration.

Finally, apply
Lemma~\ref{lem:same-order-normalized-h-difference} with the history
functional \(\mathscr H_{k,S}\) from
\eqref{eq:finite-horizon-history-functional}.  The original entrance
traces in \eqref{eq:finite-horizon-memory-entrance} are bounded by
\(C_Sd_{k,0}^{\rm typ}\).  Hence, for every \(\tau\leq S\),
\begin{equation}\label{eq:finite-horizon-history-by-Volterra-data}
 \mathscr H_{k,S}(\tau)^2
 \leq C_S\bigl(d_{k,0}^{\rm typ}\bigr)^2+
 C_S\int_{\tau_0}^{\tau}\left(
  \bigl(\mathfrak D_k^{\rm hyb,full}(q)\bigr)^2
  +\|V(q)\|_{H^1_\nu}^2\right)dq .
\end{equation}
The lemma constructs only the characteristic normalized tensor
cylinders, uses the full four-term graft bracket in
\eqref{eq:h-specific-scaled-graft-Duhamel}, records all weight
commutators, and keeps the finite-rank feedback in its
\(L^1\)-in-time Duhamel form.  Its local estimate excludes the exterior,
graft-input, and homogeneous/separated localized \(F\) memories, which
are controlled from
\(\tau_0\) by \eqref{eq:finite-horizon-memory-by-history}.

Subdivide \([\tau_0,\tau]\) into finitely many intervals of length at
most \(\delta_\#\).  On each interval use
\eqref{eq:finite-horizon-triangular-top-ledger}; then use
\eqref{eq:finite-horizon-memory-by-history},
\eqref{eq:finite-horizon-interface-memory-splitting}, and
\eqref{eq:finite-horizon-history-by-Volterra-data}.  Discrete
Volterra--Gronwall gives
\[
 \bigl(\mathfrak D_k^{\rm hyb,full}(\tau)\bigr)^2
 \leq C_S\bigl(d_{k,0}^{\rm typ}\bigr)^2+
 C_S\int_{\tau_0}^{\tau}\left(
  \bigl(\mathfrak D_k^{\rm hyb,full}(q)\bigr)^2
  +\|V(q)\|_{H^1_\nu}^2\right)dq ,
\]
which is \eqref{eq:finite-horizon-weighted-parabolic}.  Every history
integral begins at \(\tau_0\) and has upper limit \(\tau\).  No
heterogeneous product propagator has been used: only zero-trace local
contributions are recentered.  The exterior, graft-input, and
homogeneous/separated localized \(F\) terms retain their original
entrance face and separated spatial memory; the global \(F\)-block and
the nonseparated localized \(F\)-block are treated separately.  The
former restarts from its actual global trace, while the latter is
recentered as a zero-trace Abel contribution on the same source star.
\end{proof}

\begin{proposition}[Two-state and first-variation prepared evolution]
\label{prop:two-state-prepared-evolution}
Fix \(k\geq12\), \(0<\alpha<1\), a common-margin prepared ball
 \(\mathscr B\subset\mathscr P_{\tau_0}^{k+2,\alpha}\), its sliced part
 \(\mathscr B_{\rm sl}=\mathscr B\cap
 \Sigma_{\tau_0}^{k+2,\alpha}\), and a finite
normalized endpoint \(S>\tau_0\).  For
\(\mathbf z_0\in\mathscr B_{\rm sl}\), let
\(\mathbf z(\,\cdot\,;\mathbf z_0)\) denote its unique maximal coupled
trajectory, equivalently
\[
 \mathbf z(\tau;\mathbf z_0)
 =\mathcal S_{\tau,\tau_0}(\mathbf z_0)
\]
wherever the solution operator is defined.  Set
\begin{equation}\label{eq:finite-horizon-survival-domain}
 \mathscr B_{\rm sl}(S)
 :=
 \left\{\mathbf z_0\in\mathscr B_{\rm sl}:
 \begin{array}{l}
 [\tau_0,S]\text{ lies in the maximal existence interval, and}\\
 \mathbf z(\tau;\mathbf z_0)\text{ remains in the fixed common}\\
 \text{prepared package with positive trajectory-specific}\\
 \text{distance from its exit faces for }\tau_0\leq\tau\leq S
 \end{array}
 \right\}.
\end{equation}
Here the survival condition has two distinct quantitative components.
Every trajectory under consideration remains in the same completed
numerical prepared package, with the common ellipticity, Gram-inverse,
coefficient, map, gauge, atlas, and buffer bounds fixed by
\(\mathscr B\) and the exterior certificate.  For each individual
entrance its compact trajectory also has positive distance from the
exit faces of that package.  This latter distance may depend on the
entrance and is used only to make \(\mathscr B_{\rm sl}(S)\) open by
continuation and continuous dependence; no estimate below divides by,
or otherwise depends on, that trajectory-specific distance.  Let
\(\mathbf z_i(\tau)\), \(i=1,2\), be two coupled feedback solutions on
\([\tau_0,S]\) whose entrances belong to
\(\mathscr B_{\rm sl}(S)\).
Use the smooth time-independent exterior reference metric
\(\widehat G_{\rm ext}\) determined by the center entrance of
\(\mathscr B\), and put the restrictions of both closed flows to
\(E^{++}\) in the single anchored Ricci--DeTurck gauge of
Lemma~\ref{lem:anchored-exterior-interface}:
\[
 G_i(t)=\chi_i(t)^*\widetilde G_i(t).
\]
On the fixed marked graft/interface collar contained in \(E^{++}\),
transport each physical marking by
\[
 \widetilde\iota_i=\iota\circ\chi_i(t)^{-1}.
\]
Then the pushforward convention in \eqref{eq:marked-identification}
gives the exact covariance identity
\begin{equation}\label{eq:two-state-marking-covariance}
 (\widetilde\iota_i)_*\widetilde G_i
 =\iota_*G_i .
\end{equation}
Thus the normalized tensors \(h_i\) are unchanged and the physical
graft identity is covariant.  The gauge diffeomorphism is retained
explicitly: the gauge-fixed metric is measured by
\eqref{eq:hybrid-physical-block}, and the diffeomorphism and
transported marking are measured by
\eqref{eq:hybrid-marking-block}.  Consequently the fixed cutoff
\(\eta\) creates no unmeasured commutator.
Write
\[
 H_i=\rho_\tau h_i,\qquad
 V=H_1-H_2,\qquad
 c_i=(a_i,b_i).
\]
All differences of maps are taken in the exponential charts used in
\eqref{eq:prepared-Banach-norm}.  We use the hybrid distance
\(\mathfrak D_k^{\rm hyb,full}\) defined above.
There is
\[
 C_S=C\!\left(
 S,k,\alpha,\mathfrak P_{\rm prep},
 \text{the fixed exterior certificate}
 \right)<\infty,
\]
uniform for every pair of entrances in
\(\mathscr B_{\rm sl}(S)\).  It depends on the fixed common package
bounds, but not on either trajectory's additional exit-face distance.
Moreover,
\begin{equation}\label{eq:two-state-finite-horizon}
 \begin{split}
 &\sup_{\tau_0\leq\tau\leq S}
   \left(\mathfrak D_k^{\rm hyb,full}(\tau)
         +\|V(\tau)\|_{L^2_\nu}\right)
 +\left(
   \int_{\tau_0}^{S}\|V(\tau)\|_{H^1_\nu}^2\,d\tau
  \right)^{1/2}\\
 &\hspace{35mm}\leq
 C_S\,d_{k,0}^{\rm typ}
 \leq C_S\,
 \|\mathbf z_1(\tau_0)-\mathbf z_2(\tau_0)\|
       _{\mathscr X_{\rm prep}^{k+2,\alpha}} .
 \end{split}
\end{equation}
The feedback difference satisfies
\begin{equation}\label{eq:two-state-feedback}
 |c_1-c_2|
 \leq C_S\left(
  \|V\|_{H^1_\nu}+\mathfrak D_k^{\rm hyb,full}\right).
\end{equation}
The same estimates hold for first variations.  Consequently the
finite-horizon coupled solution map is \(C^1\) from
\(\mathscr B_{\rm sl}(S)\) to the state space two derivatives lower,
and this
assertion includes the solution-dependent pullbacks and moving
physical graft support.
\end{proposition}

\begin{proof}
Subtract the two exact Gram systems.  With
\(\delta M=M_1-M_2\), \(\delta d=d_1-d_2\), and
\(\delta c=c_1-c_2\),
\[
 M_1\delta c=-\delta d-\delta M\,c_2.
\]
Uniform Gram inversion, the localized difference estimates in
Remark~\ref{rem:localized-Lipschitz}, and the dynamic composition
estimate \eqref{eq:dynamic-column-Lipschitz} give
\eqref{eq:two-state-feedback}.  In particular, every term vanishes
when the two states agree; there is no solution-independent annular
remainder.

Subtract the two cutoff equations for \(H_i\) and pair with
\(V\).  Since both \(H_i\) satisfy the same exact slice,
\(V\perp\mathcal Z\).  Coercivity, polarization of
\eqref{eq:localized-Q-Lip-energy}, and
\eqref{eq:B-bilinear} yield
\begin{equation}\label{eq:two-state-energy-proof}
 \frac d{d\tau}\|V\|_{L^2_\nu}^2
 +c\|V\|_{H^1_\nu}^2
 \leq
  C_S\left(\|V\|_{L^2_\nu}^2
          +(\mathfrak D_k^{\rm hyb,full})^2\right).
\end{equation}
Here the difference of the outer terms is controlled as follows.
On every dyadic annulus write
\[
 K_i(T)
 =
 \bigl((1-\eta)\circ\Phi_i^{-1}\bigr)
 (\Theta_i\circ\Phi_i^{-1})^*T.
\]
The scale-one inverse and composition estimates in
Lemma~\ref{lem:prepared-chart-calculus} bound \(K_1-K_2\), including the
change of its support, by \(C_S\mathfrak D_k^{\rm hyb,full}\).  On its receding
support the Gaussian factor gives the stronger
\(e^{-ce^\tau}\mathfrak D_k^{\rm hyb,full}\) bound in every pairing used in
\eqref{eq:two-state-energy-proof}.  The same calculation treats the
pure graft defect.

The equations for \(t,\log\lambda\), and \(R\) are ODEs with
locally Lipschitz right-hand sides in the scaled chart.  Apply
Theorem~\ref{thm:finite-horizon-weighted-parabolic} on
\([\tau_0,S]\), using
\eqref{eq:two-state-feedback}.  It gives
\begin{equation}\label{eq:two-state-auxiliary-proof}
 \bigl(\mathfrak D_k^{\rm hyb,full}(\tau)\bigr)^2
 \leq C_S\bigl(d_{k,0}^{\rm typ}\bigr)^2
 +C_S\int_{\tau_0}^{\tau}
   \left(\bigl(\mathfrak D_k^{\rm hyb,full}(s)\bigr)^2
         +\|V(s)\|_{H^1_\nu}^2\right)ds .
\end{equation}
Adding a small multiple of
\eqref{eq:two-state-auxiliary-proof} to the integral of
\eqref{eq:two-state-energy-proof}, and applying Gronwall, proves
\eqref{eq:two-state-finite-horizon}.

We now verify, rather than assume, the smallness of every Taylor
forcing.  Let
\[
 \xi_n\longrightarrow0
\]
be arbitrary sliced initial increments and put
\[
 \mathbf z_n(\tau)
 :=\mathcal S_{\tau,\tau_0}(\mathbf z_0+\xi_n),
 \qquad
 \mathbf z(\tau)
 :=\mathcal S_{\tau,\tau_0}(\mathbf z_0),
 \qquad
 \epsilon_n
 :=\|\xi_n\|_{\mathscr E_{\rm prep}^{k+2,\alpha}}.
\]
Lemma~\ref{lem:finite-horizon-sequential-common-tail} applies to this
entrance sequence.  Let \(\eta_n\) be the sum of the finitely many
high spatial atlas-supremum differences which occur as high factors
in the prepared, Schauder, marking, gauge, interface, graft,
moving-support, history, Gram, and coefficient Taylor formulas.
The lemma gives
\[
 \eta_n\longrightarrow0.
\]
No rate for \(\eta_n\) is required.

For this paragraph, let
\(\|\mathbf z_n-\mathbf z\|_{\rm low}\) denote the finite sum of the
lower-order supremum component and coefficient norms occurring as low
factors in the mixed tame formulas, and let
\(\|\mathbf z_n-\mathbf z\|_{\rm low,Bochner}\) denote the analogous
finite sum of their source-adapted parabolic forcing and
time-derivative norms.  The two-state estimate
\eqref{eq:two-state-finite-horizon} controls the hybrid and energy
members of these sums.  The global \(F\)-block and its target
coefficients are controlled by
\eqref{eq:source-atlas-one-order-Abel-block} and
\eqref{eq:coupled-map-difference-forcing}; the anchored metric block is
controlled from its exact unequal-clock equation
\eqref{eq:finite-horizon-unequal-clock-DeTurck}; and the gauge,
inverse-gauge, and transported-marking blocks are controlled by
Lemma~\ref{lem:compact-graft-buffer-propagation} and
\eqref{eq:graft-buffer-finite-propagation}.  Direct integration of the
triangular equations controls \(t,\log\lambda,R\), the algebraic
Lipschitz bounds in \eqref{eq:two-state-feedback} and
\eqref{eq:interface-block-controlled} control the Gram/feedback and
interface blocks, and
\eqref{eq:finite-horizon-history-by-Volterra-data} controls the history
terms.  Thus these componentwise estimates, together with
\eqref{eq:two-state-finite-horizon}, give
\[
 \sup_{\tau_0\leq\tau\leq S}
   \|\mathbf z_n(\tau)-\mathbf z(\tau)\|_{\rm low}
 +\|\mathbf z_n-\mathbf z\|_{\rm low,Bochner}
 \leq C_S\epsilon_n .
\]
The mixed tame estimates
\eqref{eq:prepared-mixed-tame-remainder},
\eqref{eq:feedback-mixed-tame-remainder}, and
\eqref{eq:Schauder-mixed-tame-residual}, together with Taylor's
integral formula for the buffered gauge, marking, interface, and
graft blocks, therefore place one factor in the high spatial
difference and the other in the preceding lower-order difference.
Consequently there is a sequence \(\gamma_n\downarrow0\) such that the
sum of all residual trace, Bochner, and \(L^1_\tau\) history norms is
bounded by
\begin{equation}\label{eq:finite-horizon-Taylor-forcing-small}
 C_S\gamma_n\epsilon_n .
\end{equation}
For terms involving an \(L^2_\tau\) lower factor, this follows from
the uniform-in-time high factor and Cauchy--Schwarz on the fixed
finite horizon.  The finite-dimensional Gram remainder and the
history integrals obey the same bound by uniform continuity of their
first derivatives on the compact base trajectory and dominated
convergence.

Subtract from
\(\mathbf z_n-\mathbf z\) the solution of the coupled linearized
system with initial value \(\xi_n\), and call the result
\(\mathcal R_n\).  Choose \(r_n\) to dominate the finitely many
residual norms above and so that
\[
 0\leq r_n\leq C_S\gamma_n\epsilon_n.
\]
Then
\begin{equation}\label{eq:finite-horizon-r-small}
 r_n=o(\epsilon_n).
\end{equation}
Thus \(r_n=o(\epsilon_n)\) is a consequence of the finite-horizon
common-tail passage and the two-state estimate, not an additional
hypothesis.
Write \(V_{R,n}\) for the cutoff-tensor component of
\(\mathcal R_n\), \(D_{R,n}\) for the analogue of
\(\mathfrak D_k^{\rm hyb,full}\) formed from its remaining components,
and \(c_{R,n}\) for its feedback residual.

The residual has zero initial value in the prepared chart.  It also
retains the slice orthogonality: the nonlinear difference and the
linearized cutoff tensor each have zero pairing with every \(Z_\mu\),
and the slice functional is linear in that cutoff tensor.  Repeating
the derivations of
\eqref{eq:two-state-feedback},
\eqref{eq:two-state-energy-proof}, and
\eqref{eq:two-state-auxiliary-proof} for the forced residual equations
therefore gives
\begin{align}
 |c_{R,n}|
 &\leq C_S\bigl(\|V_{R,n}\|_{H^1_\nu}+D_{R,n}\bigr)+r_n,
 \label{eq:forced-residual-feedback}\\
 \frac d{d\tau}\|V_{R,n}\|_{L^2_\nu}^2
 +c\|V_{R,n}\|_{H^1_\nu}^2
 &\leq C_S\bigl(\|V_{R,n}\|_{L^2_\nu}^2+D_{R,n}^2\bigr)
       +C_Sr_n^2,
 \label{eq:forced-residual-energy}\\
 D_{R,n}(\tau)^2
 &\leq C_Sr_n^2
   +C_S\int_{\tau_0}^{\tau}
       \bigl(D_{R,n}(q)^2
             +\|V_{R,n}(q)\|_{H^1_\nu}^2\bigr)\,dq .
 \label{eq:forced-residual-auxiliary}
\end{align}
The additional terms are exactly the residual-forcing pairings or
Duhamel integrals.  Cauchy--Schwarz and Young give \(r_n^2\) in the
last two displays.  The first display follows by subtracting and
linearizing the exact Gram system and applying
\eqref{eq:feedback-mixed-tame-remainder};
\eqref{eq:feedback-L1-Duhamel} controls the feedback contribution in
the auxiliary Duhamel estimate.
Thus these are forced analogues of the earlier estimates, not an
application of exact-solution-difference inequalities to a
non-solution.

Add a small multiple of
\eqref{eq:forced-residual-auxiliary} to the integral of
\eqref{eq:forced-residual-energy}, use
\eqref{eq:forced-residual-feedback}, and apply Gronwall.  Since
\(r_n=o(\epsilon_n)\), this gives
\[
 \sup_{\tau_0\leq\tau\leq S}
 \|\mathcal R_n(\tau)\|_{\mathscr E_{\rm prep}^{k,\alpha}}
 =
 o(\|\xi_n\|_{\mathscr E_{\rm prep}^{k+2,\alpha}}).
\]
Thus the finite-horizon derivative is Fr\'echet, not merely
curvewise.  The first derivatives of all chart operations and
coefficient maps depend continuously on the base state from prepared
input order \(k+2\) to output order \(k\).  Subtracting the two
corresponding linearized systems, taking the supremum over unit tangent
directions, and applying the same finite-horizon estimate proves
operator-norm continuity of the derivative.  Hence the finite-horizon
solution map is \(C^1\) with the two-derivative input--output buffer
asserted in the proposition.
\end{proof}

\section{The uniform entrance map}
\label{sec:uniform-entrance}

The initial phase must be selected with constants which do not
deteriorate as the support radius recedes.  This is a finite-dimensional
statement.  It is important, however, to use the prepared phase
action rather than the raw complete-space action
$\mathfrak A_p$.  The latter generally changes the asymptotic cone.

Let
\[
 \mathbf z=(G,\lambda,\Theta,\Phi)
\]
denote a prepared geometric state at time $\tau_0$.  We now recall, in
the detailed notation used below, the finite action already defined in
\eqref{eq:intro-prepared-scale-leg}--%
\eqref{eq:intro-prepared-phase-action}.  Let $\mathfrak r_s$ be the complete flow of
$-\bar\nabla\bar f$, and, for $1\leq j\leq8$, let
$\psi_{j,s}^{(\tau_0)}$ be the complete flow of
$\chi_{\tau_0}W_j$.  Define nine one-parameter maps on prepared states
by
\begin{align}
 \mathscr L_0(s)(G,\lambda,\Theta,\Phi)
 &=
 \bigl(G,e^{-s}\lambda,
       \mathfrak r_{-s}\circ\Theta,
       \mathfrak r_{-s}\circ\Phi\bigr),
 \label{eq:prepared-scale-leg}\\
 \mathscr L_j(s)(G,\lambda,\Theta,\Phi)
 &=
 \bigl(G,\lambda,
       \psi_{j,-s}^{(\tau_0)}\circ\Theta,
       \psi_{j,-s}^{(\tau_0)}\circ\Phi\bigr),
 \qquad1\leq j\leq8.
 \label{eq:prepared-diffeomorphism-leg}
\end{align}
In each diffeomorphism leg the relative map
$F=\Theta^{-1}\circ\Phi$ is unchanged.  In the scale leg the factor
$e^{-s}$ in $\lambda$ and the radial pullback combine to give the
diagonal action \eqref{eq:finite-scale-action} on the inner metric.

Fix the order $0,1,\ldots,8$ and put
\[
 \mathbf z_p
 =\mathscr L_8(p_8)\circ\cdots\circ
   \mathscr L_0(p_0)(\mathbf z).
\]
Write its geometric components as
\[
 (G,\lambda_p,\Theta_p,\Phi_p),\qquad
 R_p=\varphi_{-\tau_0}\circ\Theta_p,\qquad
 F_p=\Theta_p^{-1}\circ\Phi_p,
\]
and re-form
\[
 S_p=\lambda_p\Theta_p^*\bar g,\qquad
 \acute G_p=\eta\,\iota_*G+(1-\eta)S_p.
\]
Thus the nine legs define the tuple-valued prepared phase action
\begin{equation}\label{eq:prepared-tuple-action}
 \mathbf A^{\rm prep}_{p,\tau_0}(\mathbf z)
 :=\mathbf z_p .
\end{equation}
For every output order \(r\geq3\), this action is \(C^2\) from the
two-derivative-buffered input space
\(\mathscr P_{\tau_0}^{r+2,\alpha}\) into
\(\mathscr P_{\tau_0}^{r,\alpha}\).  Completeness and the
scale-normalized generator bounds make all constants uniform on a
common-margin prepared ball.  Equivalently, after the Gaussian
integrations by parts used below, its finite-dimensional scalar moment
map is \(C^1\) at the unbuffered order by
Lemma~\ref{lem:unbuffered-Gaussian-moment-map}.  The quantitative
second-derivative bounds used for the phase implicit-function theorem
are taken in the displayed two-derivative-buffered space.  We do not
assert that pullback by a varying diffeomorphism is \(C^2\) on one
unbuffered \(C^{r,\alpha}\) space.

The normalized-metric realization of
\eqref{eq:prepared-tuple-action} is
\begin{equation}\label{eq:prepared-phase-action}
 \mathfrak A^{\rm prep}_{p,\tau_0}(\mathbf z)
 :=
 \bar g+h\!\left(\mathbf A^{\rm prep}_{p,\tau_0}(\mathbf z)\right)
 =\lambda_p^{-1}(\Phi_p^{-1})^*\acute G_p.
\end{equation}
Thus \(\mathbf A^{\rm prep}_{p,\tau_0}\) is tuple-valued, whereas
\(\mathfrak A^{\rm prep}_{p,\tau_0}\) is tensor-valued.  The inner
metric action is the usual \(\mathfrak A_p\), while the outer metric
is re-grafted to the correspondingly changed adaptive target.
In particular, the prepared asymptotic identity is preserved.  To
check the linearization, differentiate
\eqref{eq:prepared-scale-leg}--%
\eqref{eq:prepared-diffeomorphism-leg}.  The simultaneous left
composition of $\Theta$ and $\Phi$ cancels on the pure outer target,
whereas on the physical part it produces the corresponding pullback.
Consequently, at every prepared state,
\begin{align}
 D_{p_j}\left[
  \mathfrak A^{\rm prep}_{p,\tau_0}(\mathbf z)-\bar g
 \right]_{p=0}
 &=\mathscr C_{j,\tau_0}(\mathbf z),
 \label{eq:prepared-action-linearization}\\
 \mathscr C_{0,\tau_0}(\mathbf z)
 &=\mathcal Y_{0,\tau_0}(\mathbf z)
   +\widetilde\B_{0,\tau_0}h(\mathbf z),
 \label{eq:prepared-full-column-zero}\\
 \mathscr C_{j,\tau_0}(\mathbf z)
 &=\mathcal Y_{j,\tau_0}(\mathbf z)
   +\Lie_{\chi_{\tau_0}W_j}h(\mathbf z),
 \qquad1\leq j\leq8.
 \label{eq:prepared-full-column-j}
\end{align}
Thus the derivative equals the direct effective column only at an
exact prepared background \(\mathbf z_0\) with \(h(\mathbf z_0)=0\).

\begin{lemma}[Buffered full prepared columns]
\label{lem:buffered-full-prepared-columns}
Let \(k\geq3\), let \(\mathbf z_0\) be an exact prepared state with
\(h(\mathbf z_0)=0\), and let
\(\mathscr B\subset\mathscr P_{\tau_0}^{k+2,\alpha}\) be a sufficiently
small common-margin ball about \(\mathbf z_0\).  For
\(\tau_0\geq\tau_*\), uniformly in \(\mathbf z\in\mathscr B\),
\begin{equation}\label{eq:static-full-phase-column}
 \left|
 \ip{\rho_{\tau_0}\mathscr C_{j,\tau_0}(\mathbf z)}{Z_\mu}
 -\ip{Y_j}{Z_\mu}
 \right|
 \leq
 C\|\mathbf z-\mathbf z_0\|_{\mathscr X_{\rm prep}^{k+2,\alpha}}
 +Ce^{-ce^{\tau_0}}.
\end{equation}
Moreover the scalar maps
\[
 \mathbf z\longmapsto
 \ip{\rho_{\tau_0}\mathscr C_{j,\tau_0}(\mathbf z)}{Z_\mu}
\]
have uniformly bounded first derivatives on \(\mathscr B\), and the
moment map
\[
 (p,\mathbf z)\longmapsto
 \ip{\rho_{\tau_0}
  \bigl(\mathfrak A^{\rm prep}_{p,\tau_0}(\mathbf z)-\bar g\bigr)}
   {Z_\mu}
\]
has uniformly bounded first and second derivatives from a fixed
neighborhood of \(\{0\}\times\mathscr B\) to \(\mathbb R\).
\end{lemma}

\begin{proof}
Equations \eqref{eq:prepared-full-column-zero}--%
\eqref{eq:prepared-full-column-j} split each full column into the
direct effective column and a term linear in \(h\).  The direct term is
controlled by Lemma~\ref{lem:static-prepared-columns}.  In the remaining
Gaussian pairing, one integration by parts moves the derivative in the
Lie derivative off \(h\); the derivatives of
\(\rho_{\tau_0}Z_\mu e^{-\bar f}\) have fixed polynomial growth
times the Gaussian weight.  Hence that pairing is bounded by the
buffered prepared norm of \(h\).  Since \(h(\mathbf z_0)=0\), this gives
\eqref{eq:static-full-phase-column}.  The scale-one inverse, pullback,
product, and composition estimates in
Lemma~\ref{lem:prepared-chart-calculus} give the first derivatives.  A
second differentiation uses precisely the two spatial derivatives in
the passage from input order \(k+2\) to output order \(k\); the same
integration by parts supplies a common Gaussian majorant.  This proves
the asserted uniform derivative bounds.
\end{proof}

\begin{lemma}[Witnessed scale-uniform lower stability of the
harmonic-radius face]
\label{lem:prepared-harmonic-radius-lower-stability}
Fix \(r\geq3\), \(0<\alpha<1\), and one reduced geometric package
\(\mathfrak P_{\rm har}^{\rm geom}\) from
\eqref{eq:reduced-prepared-harmonic-package}, at the stage before
\(\delta_{\rm c2}\) and the later entrance thresholds are selected.
Every completed numerical prepared package used below is required only
to extend this reduced record.  For every
\[
 \eta>0,\qquad q\in(0,Q_{\rm har}-1),\qquad \zeta\in(0,1)
\]
there is a number
\[
 \delta_{\rm har}
 (\eta,q,\zeta;\mathfrak P_{\rm har}^{\rm geom})>0,
\]
with the following property; no monotonicity of this modulus in any
parameter is asserted or needed.  The number is independent of the
entrance time.  Suppose
\(\mathbf z,\mathbf z'\in\mathscr P_{\tau_c}^{r,\alpha}\) lie in one
uniformly interior common-margin chart carrying the same reduced
geometric package.  Assume that, after the exact graph reduction
\eqref{eq:harmonic-radius-normalized-graph-identity} below, both
normalized metric branches of \(\mathbf z\), at every \(y\in M\),
possess coefficient-\(q\) harmonic witnesses on the ball of radius
\[
 (\kappa_{\rm har}+2\eta)r_{\rm la}(y)
\]
with domain buffer \(\zeta\), in the sense of
\eqref{eq:buffered-harmonic-chart-reserve}.  If
\[
 \|\mathbf z'-\mathbf z\|_{\mathscr X_{\rm prep}^{r,\alpha}}
 \leq\delta_{\rm har}
       (\eta,q,\zeta;\mathfrak P_{\rm har}^{\rm geom}),
\]
then
\begin{equation}\label{eq:prepared-harmonic-radius-lower-stability}
 \mathfrak h_{\rm har}(\tau_c,\mathbf z')
 \geq\kappa_{\rm har}+\eta.
\end{equation}
Moreover the primed charts at radius
\((\kappa_{\rm har}+\eta)r_{\rm la}(y)\) may be chosen with a common
positive coefficient reserve and a common positive domain buffer
depending only on
\((\eta,q,\zeta;\mathfrak P_{\rm har}^{\rm geom})\).  Thus the lemma propagates
the witness needed at the next perturbative step: with
\(\eta':=\eta/2\), the primed state carries an
\((\eta',q',\zeta')\)-witness for some common
\(q'>0\), \(\zeta'\in(0,1)\).  It asserts
one-sided lower stability under a recorded witness, not continuity or
a Lipschitz estimate for the harmonic-radius functional.

There is also a reserve-to-operative form used when membership in a
strict entrance class must itself persist.  Fix two triples satisfying
\begin{equation}\label{eq:prepared-harmonic-two-tier-order}
 \begin{gathered}
  0<\eta_{\rm op}<\eta_+,\qquad
  0<q_{\rm op}<q_+<Q_{\rm har}-1,\\
  0<\zeta_{\rm op}<\zeta_+<1.
 \end{gathered}
\end{equation}
If the unprimed branches carry the \(+\)-witness on radius
\((\kappa_{\rm har}+2\eta_+)r_{\rm la}\), there is
\begin{equation}\label{eq:prepared-harmonic-reserve-modulus}
 \delta_{\rm har}^{\rm op}
 \bigl(
  \eta_{\rm op},q_{\rm op},\zeta_{\rm op};
  \eta_+,q_+,\zeta_+;
  \mathfrak P_{\rm har}^{\rm geom}
 \bigr)>0
\end{equation}
such that the same closeness hypothesis with this modulus gives the
\emph{same operative witness}
\[
 (\eta_{\rm op},q_{\rm op},\zeta_{\rm op})
\]
for the primed branches.  Moreover the primed branches retain a new
recorded \(+\)-triple strictly larger componentwise than the operative
one.  Thus the operative triple, and every numerical condition imposed
on it, is unchanged under sufficiently small perturbation.

The same conclusion is uniform on a family of centers precisely when
the same \(\eta,q,\zeta\) and numerical package work on that family.
In particular, compact-family uniformity follows from a finite
subcover once actual buffered witnesses have been chosen at every
center.  A bare lower bound for \(r_{\rm har}\), without a common
coefficient reserve, is not used as a substitute for this hypothesis.

The same conclusion holds when the two states carry different
normalized-time labels.  In that case the hypothesis is the
scale-one \(C^{2,\alpha}\) closeness of their normalized graph tensors
\(\bar g+h\), under the canonical model identification, together with
the same fixed background branch and the displayed common witnesses.
No comparison of the time-typed cutoffs, maps, or slices is required:
\eqref{eq:harmonic-radius-normalized-graph-identity} removes all of
them before the openness argument is applied.
\end{lemma}

\begin{proof}
The graph identities remove every raw map and scale from this
functional.  Indeed, pullback and constant-scale covariance of the
fixed-convention harmonic radius give, with \(y=\Phi(x)\),
\[
 \frac{r_{\rm har}(\acute G,x)}{r_{\rm sol}(x)}
 =
 \frac{r_{\rm har}(\bar g+h,y)}{(1+\bar f(y))^{1/2}}.
\]
Since \(\Theta(F(x))=\Phi(x)\), the same covariance gives
\[
 \frac{r_{\rm har}(S,F(x))}{r_{\rm sol}(x)}
 =
 \frac{r_{\rm har}(\bar g,y)}{(1+\bar f(y))^{1/2}}.
\]
Both \(\Phi\) and \(F\) are diffeomorphisms.  Taking the global infima
therefore yields the exact normalized identity
\begin{equation}\label{eq:harmonic-radius-normalized-graph-identity}
 \mathfrak h_{\rm har}(\tau_c,\mathbf z)
 =
 \min\left\{
  \inf_{y\in M}
   \frac{r_{\rm har}(\bar g+h(\mathbf z),y)}{r_{\rm la}(y)},
  \inf_{y\in M}
   \frac{r_{\rm har}(\bar g,y)}{r_{\rm la}(y)}
 \right\}.
\end{equation}
Thus only the normalized metric \(\bar g+h\) varies.  Its scale-one
\(C^{2,\alpha}\) distance is controlled directly by the prepared norm,
uniformly in \(\tau_c\).

Fix one center \(y\), divide lengths by \(r_{\rm la}(y)\), and use
the reserved unprimed chart on the buffered outer ball.  In these
coordinates the prepared norm makes the primed and unprimed
coefficient metrics \(C^{2,\alpha}\)-close, with a bound independent
of \(y\) and \(\tau_c\).  Metric-ball comparison first places the
primed ball of radius \(\kappa_{\rm har}+\eta\) a definite distance
inside the buffered unprimed coordinate domain.  On an intermediate
domain solve the primed harmonic-coordinate Dirichlet problem with
the unprimed coordinate functions as boundary data.  Uniform
ellipticity, the coefficient reserve \(q\), the domain buffer
\(\zeta\), and the scale-one package give a common Dirichlet inverse
and a common Schauder constant.  The new coordinate functions are
therefore \(C^{3,\alpha}\)-close to the old ones.  For a sufficiently
small bound depending only on
\((\eta,q,\zeta;\mathfrak P_{\rm har}^{\rm geom})\), they remain a
diffeomorphism on the primed ball and their coefficient metric obeys
\eqref{eq:fixed-harmonic-radius-convention}, with positive residual
coefficient and domain reserves.  Here is the global injectivity step.
In the old coordinates the new harmonic map
\(w:\overline D\to\mathbb R^4\) has boundary value
\(\operatorname{Id}\) on \(\partial D\), and its \(C^1\)-closeness to
the identity gives \(\det Dw>0\) on \(D\).  For \(z\in D\), homotopy
relative to the boundary gives \(\deg(w,D,z)=1\).  Every preimage has
local degree \(+1\), so \(z\) has exactly one preimage.  For
\(z\notin\overline D\), the degree is zero and positivity of every
local degree excludes a preimage.  Finally, no interior point can map
to \(\partial D\): otherwise openness of the local diffeomorphism
would produce nearby image points outside \(\overline D\), a
contradiction.  Thus \(w(D)=D\), with exactly one preimage for every
point, and \(w\) is a diffeomorphism.  After the injectivity argument,
replace it by
\[
 \widetilde w:=w-w(y).
\]
Then \(\widetilde w(y)=0\), as required in the definition of
\(r_{\rm har}\) preceding
\eqref{eq:fixed-harmonic-radius-convention}.  Translation of the
Euclidean target preserves harmonicity, the Jacobian, global
injectivity, and every coefficient and domain estimate; it merely
translates the allowed Euclidean image.  Restricting to a slightly
larger primed metric ball before passing to the desired
\((\kappa_{\rm har}+\eta)\)-ball leaves the residual domain buffer.

This is the standard parameterized harmonic-coordinate openness
argument on a buffered domain.  The metric-ball-domain definition is
essential: after restriction the Euclidean image is allowed to be
nonround, so no reparametrization that could destroy harmonicity is
made.  The constants are common at every \(y\) because the hypotheses
record the same \(q,\zeta\) and the prepared package is scale-uniform.
Applying the argument to the varying branch in
\eqref{eq:harmonic-radius-normalized-graph-identity}, while observing
that the background branch is unchanged, proves
\eqref{eq:prepared-harmonic-radius-lower-stability}.  The residual
reserves give the last assertion.  For a compact family, choose
witnesses locally and take a finite subcover before taking the
minimum perturbation radius.  No monotone-minorant construction is
involved.  The cross-time clause follows from the same proof after
starting directly with
\eqref{eq:harmonic-radius-normalized-graph-identity}; only the two
normalized graph metrics are compared.

For the reserve-to-operative clause, insert the fixed intermediate
triple
\[
 \eta_{\rm mid}=\frac{\eta_{\rm op}+\eta_+}{2},\qquad
 q_{\rm mid}=\frac{q_{\rm op}+q_+}{2},\qquad
 \zeta_{\rm mid}=\frac{\zeta_{\rm op}+\zeta_+}{2}.
\]
Run the same Dirichlet argument using only the gaps between
the \(+\)- and intermediate triples.  For a sufficiently small
perturbation the primed charts satisfy the intermediate certificate.
It is both a reserve strictly above the operative triple and, after
restriction, an operative witness.  Taking the uniform minimum of the
three positive gaps proves
\eqref{eq:prepared-harmonic-reserve-modulus}; no monotone dependence of
either modulus is used.
\end{proof}

\begin{lemma}[Witnessed harmonic openness on a finite physical cover]
\label{lem:finite-physical-harmonic-openness}
Fix \(0<\alpha<1\), let \(\mathcal X\) be a closed four-manifold, and
let
\[
 U_a\Subset V_a\Subset W_a\Subset W_a^+,\qquad
 R_a>0,\qquad1\leq a\leq N,
\]
be one finite buffered physical cover with fixed overlap,
scale-comparability, coefficient, ellipticity, separation, and
normalized boundary-atlas constants.  The coefficient record is taken
on \(W_a^+\) and includes a fixed scaled separation from \(W_a\) to its
boundary.  Fix
\[
 \upsilon>0,\qquad \eta>0,\qquad
 q\in(0,Q_{\rm har}-1),\qquad
 \zeta,\zeta_{\rm out}\in(0,1).
\]
Suppose a \(C^{2,\alpha}\) metric \(G\) has, for every \(a\) and
\(x\in V_a\), a coefficient-\(q\), domain-\(\zeta\) harmonic witness
on the ball
\[
 B_G\!\left(x,r_a^{\rm wit}\right),
 \qquad
 r_a^{\rm wit}:=(\upsilon+2\eta)R_a,
\]
and let
\(\mathcal D_{a,x}
  =u_{a,x}^{-1}(r_a^{\rm wit}D_{a,x})\)
be the manifold preimage of its Euclidean Dirichlet domain.  Require
\(\overline{\mathcal D_{a,x}}\Subset W_a\), quantitatively in the
fixed physical reference atlas:
\[
 \inf_{\substack{1\leq a\leq N\\x\in V_a}}
 R_a^{-1}\operatorname{dist}_{\rm ref}
 \bigl(\overline{\mathcal D_{a,x}},
       \mathcal X\setminus W_a\bigr)\geq\zeta_{\rm out}.
\]
Freeze these common quantitative data as
\begin{equation}\label{eq:finite-physical-harmonic-package}
 \begin{aligned}
 \mathfrak P_{\rm har}^{\rm phys}
 :=\bigl(&4,\alpha,Q_{\rm har},N,\zeta_{\rm out},\\
  &\{U_a\Subset V_a\Subset W_a\Subset W_a^+\}_{a=1}^N,\\
  &\{R_a\}_{a=1}^N,\\
  &\text{the displayed cover, coefficient, and ellipticity constants},\\
  &\text{the separation and quantitative domain constants}
  \bigr).
 \end{aligned}
\end{equation}
The continuously center-indexed witness pairs are not asserted to be a
finite set; the finite record contains their common quantitative type
and constants.
This record is fixed before any Ricci--DeTurck time width is selected;
it contains neither \(\delta_{\rm RF}\), a lifetime width, nor a later
perturbation size.  The norm
\(C^{2,\alpha}_{R_a}(W_a)\) below is the dimensionless coefficient norm
in the fixed physical reference atlas included in this record.  There is
\[
 \delta_{\rm har}^{\rm phys}
 (\upsilon,\eta,q,\zeta;\mathfrak P_{\rm har}^{\rm phys})>0
\]
such that, if another metric \(G'\) obeys
\[
 \max_a
 \|G'-G\|_{C^{2,\alpha}_{R_a}(W_a)}
 \leq\delta_{\rm har}^{\rm phys},
\]
then, for every \(a\) and \(x\in V_a\),
\[
 r_{\rm har}(G',x)\geq(\upsilon+\eta)R_a.
\]
At every such center the primed charts retain common positive
coefficient and domain reserves; equivalently, with
\(\eta'=\eta/2\), they form a new witness package at radius
\((\upsilon+2\eta')R_a\).

The reserve-to-operative form has the same center set.  If
\[
 0<\eta_{\rm op}<\eta_+,\qquad
 0<q_{\rm op}<q_+<Q_{\rm har}-1,\qquad
 0<\zeta_{\rm op}<\zeta_+<1,
\]
and \(G\) carries, for every \(a\) and \(x\in V_a\), the
\(+\)-certificate at radius
\((\upsilon+2\eta_+)R_a\), with the same recorded outer-domain
separation \(\zeta_{\rm out}\), there is a positive modulus
\begin{equation}\label{eq:physical-harmonic-reserve-modulus}
 \delta_{\rm har}^{\rm phys,op}
 \bigl(
  \upsilon;
  \eta_{\rm op},q_{\rm op},\zeta_{\rm op};
  \eta_+,q_+,\zeta_+;
  \mathfrak P_{\rm har}^{\rm phys}
 \bigr)
\end{equation}
such that every metric satisfying the corresponding displayed
\(C^{2,\alpha}_{R_a}(W_a)\) closeness retains, at every
\(x\in V_a\), the same operative physical witness triple and a new
reserve triple strictly larger componentwise.  This modulus is
independent of
\(\delta_{\rm RF}\) and every later time-width choice.

The same assertion holds for a compact family when the displayed
witness and finite-cover constants are common.  It also holds
pointwise for a time-dependent metric \(G(t)\) whenever the actual
metrics, expressed in the fixed physical reference atlas, satisfy the
displayed scale-\(R_a\) \(C^{2,\alpha}\) closeness on the recorded outer
witness domains.  No conclusion is inferred here merely from a
DeTurck representative at a moved center.
\end{lemma}

\begin{proof}
Fix an arbitrary pair \((a,x)\) with \(x\in V_a\), rescale \(R_a\) to
one, and use the
selected witness \((u_{a,x},D_{a,x})\).  Solve the \(G'\)-harmonic
coordinate Dirichlet problem with the old coordinates as boundary
data.  The Dirichlet inverse, Schauder, metric-ball comparison,
Jacobian, degree, and restriction constants depend only on
\[
 (\upsilon,\eta,q,\zeta,\zeta_{\rm out};
   \mathfrak P_{\rm har}^{\rm phys}),
\]
not on \(x\) or on the particular selected witness: uniform
ellipticity and coefficient control give the common Dirichlet and
Schauder bounds, the two quantitative domain buffers give the common
boundary and restriction bounds, and the fixed radius gap is
\(\eta\).  Thus the new coordinates are uniformly
\(C^{3,\alpha}\)-close to the old ones.  The positive-Jacobian,
degree-one argument in the proof of
Lemma~\ref{lem:prepared-harmonic-radius-lower-stability} gives global
injectivity, and restriction from the buffered outer domain leaves
residual coefficient and domain reserves.  Since the estimate is
already independent of the center, one takes a minimum only over the
finitely many cover indices \(a\), not over the continuously indexed
witness family.  The compact-family assertion follows by first making
the displayed quantitative package common on a finite subcover.  The
time-dependent version is the same local argument in the fixed
physical coordinates.
For the reserve-to-operative version, take the componentwise arithmetic
midpoint of the operative and \(+\)-triples, as in the proof of
\eqref{eq:prepared-harmonic-reserve-modulus}.  The three positive gaps
give the new reserve, whose restriction is the unchanged operative
certificate.
\end{proof}

\begin{lemma}[Reference-carrier physical coefficient certificate]
\label{lem:reference-carrier-physical-certificate}
Fix \(0<\alpha<1\), a finite buffered physical cover, its scales and
reference atlas,
and three strictly nested physical witness tiers
\[
 \mathbf p_{\rm op}^{\rm phys}
 <\mathbf p_{+}^{\rm phys}
 <\mathbf p_{\rm ref}^{\rm phys}
\]
componentwise.  Let \(G_{\rm ref}^{\rm phys}\) be a smooth metric whose
scale-\(R_a\) coefficients through order fourteen are uniformly bounded
on every \(\widetilde W_a^{\rm har}\), and suppose that it carries the
reference-tier center-indexed witnesses and quantitative Dirichlet
containments specified in item~(3) of
Definition~\ref{def:witnessed-auxiliary-harmonic-package}, prior to the
coefficient-ball and Ricci-time choices in item~(4).
Then one may choose
\[
 \varepsilon_{\rm coeff}^{\rm phys}>0,\qquad
 1<\Lambda_{\rm coeff}^{\rm phys}<\infty,\qquad
 \mu_{\rm coeff}^{\rm phys}>0,\qquad
 \delta_{\rm RF}>0
\]
so that the following hold.
\begin{enumerate}
\item The reference metric has more than
 \(4\mu_{\rm coeff}^{\rm phys}\) package-face slack in the outer ball
 \eqref{eq:physical-coefficient-ball}.
\item Every metric in that outer ball carries the same common
 \(+\)-tier physical witness package at every \(x\in V_a^5\), and hence
 also the operative tier.
\item For every \(G_\circ\) in the outer ball, its actual fixed-marking
 closed Ricci flow satisfies, for every \(a\),
 \[
 \sup_{\substack{0\leq t<t_*(G_\circ)\\
                 t\leq\delta_{\rm RF}R_a^2}}
 \|G(t;G_\circ)-G_\circ\|_
 {C_{R_a}^{2,\alpha}(W_a^{\rm har})}
 \leq\frac14\delta_{\rm har}^{\rm phys,wit}.
 \]
\end{enumerate}
All constants depend only on the frozen finite reference package.  The
assertion is uniform for a compact reference family carrying the same
quantitative cover and witness data.
\end{lemma}

\begin{proof}
Take the maximum of the scale-\(R_a\)
\(C^{14,\alpha}\) coefficient norms of the reference metric and the
finite-atlas ellipticity constants.  Choose
\(\Lambda_{\rm coeff}^{\rm phys}\) strictly above these numbers and
strictly above \(1\).  The
reference-to-\(+\) modulus
\eqref{eq:witnessed-physical-reference-modulus} is positive by
Lemma~\ref{lem:finite-physical-harmonic-openness}.  Choose
\(\varepsilon_{\rm coeff}^{\rm phys}\) below the finite coefficient,
ellipticity, buffer, and domain-stability radii and so that
\[
 C_{\rm emb}^{\rm phys}\varepsilon_{\rm coeff}^{\rm phys}
 <\frac14\delta_{\rm har}^{\rm phys,ref}.
\]
The reserve-to-operative clause of that lemma then gives the common
\(+\)-tier, uniformly for every member of the outer ball.  All four
faces in \eqref{eq:physical-coefficient-package-distance} have positive
reference-center slack, so choose
\(\mu_{\rm coeff}^{\rm phys}\) below one quarter of their minimum.

The outer ball has uniform order-fourteen coefficients, ellipticity,
physical buffers, curvature bounds, and harmonic-radius lower bounds
on the fixed finite atlas.  Apply
Lemmas~\ref{lem:buffered-local-Ricci-control} and
\ref{lem:buffered-Ricci-DeTurck-coefficients} on the recorded nested
buffers, with the constants made common over the finite cover.  The
transported-marking covariance in the latter lemma transfers the
coefficient estimate back to the actual fixed physical marking.
Equivalently, after that transfer one integrates
\(\partial_tG=-2\Ric_G\) in the fixed atlas.  Consequently there are
uniform \(C<\infty\) and \(\delta_0>0\) such that
\[
 \|G(t;G_\circ)-G_\circ\|_
 {C_{R_a}^{2,\alpha}(W_a^{\rm har})}
 \leq C\,t/R_a^2
\]
whenever \(0\leq t<t_*(G_\circ)\) and
\(t\leq\delta_0R_a^2\).  Decrease \(\delta_0\), over the finite cover,
until the right side is at most
\(\delta_{\rm har}^{\rm phys,wit}/4\), and call the result
\(\delta_{\rm RF}\).  The restriction \(t<t_*(G_\circ)\) is retained:
no lower bound for the global lifetime is asserted here.  A compact
reference family is handled by making the displayed finite package
common on a finite subcover.
\end{proof}

For use by the phase theorem, we now define the static strict faces
without invoking either a phase map or the later entrance
construction.  For any prepared tuple with invertible algebraic
feedback matrix, form
\[
 \mathcal D_{\acute G,0}:=
  (\partial_t\acute G+2\Ric_{\acute G})|_{t=0},\qquad
 \mathcal D_{S,0}:=
  (\partial_tS+2\Ric_S)|_{t=0},
\]
where the time derivatives are the algebraic values prescribed by the
feedback and adaptive equations at that tuple, and put
\[
 \mathscr J_0:=
 (\acute G_0,S_0,\mathcal D_{\acute G,0},\mathcal D_{S,0}).
\]
For a fixed pre-radius source--target chart \(\mathcal U\), let
\(\mathcal N_{12,10}(\mathscr J_0;\mathcal U)\) be the sum of the
scale-one \(C^{12,\alpha}\) coefficient norms of
\(\acute G_0,S_0\), their scale-normalized curvature-derivative norms
through order ten, and the scale-normalized spatial derivative norms
through order ten of
\(\mathcal D_{\acute G,0},\mathcal D_{S,0}\).  With
\(\mathscr U_{\rm core}\) the fixed finite core cover and
\(\mathscr U_L^{\rm pre}\) the fixed cover of \(A_L^{\rm pre}\), set
\begin{equation}\label{eq:typed-entrance-coefficient-functional}
 \mathfrak C_{\rm ent}^{12,10}(\mathscr J_0):=
 \max\left\{
  \max_{\mathcal U\in\mathscr U_{\rm core}}
       \mathcal N_{12,10}(\mathscr J_0;\mathcal U),\
  \sup_{L\in\mathscr L_{\rm pre}}
  \sup_{\mathcal U\in\mathscr U_L^{\rm pre}}
       \mathcal N_{12,10}(\mathscr J_0;\mathcal U)
 \right\}.
\end{equation}
This full-state finite-jet functional is independent of the eventual
radius; no maximum over an infinite atlas is asserted.

\begin{definition}[Phase-independent strict-face certificate]
\label{def:phase-independent-strict-face-certificate}
Fix a prepared numerical package and an output order \(k\geq12\).
A prepared tuple carries the \emph{phase-independent strict-face
certificate} when, without reference to how its nine moments are or
will be imposed, it carries all of the following static data with
positive slack:
\begin{enumerate}
\item the scale bracket; weighted \(L^2_\nu\), pointwise,
pre-atlas, ellipticity, coefficient/defect
\(\mathfrak C_{\rm ent}^{12,10}\), Gram-inverse, and support-separation
faces;
\item the exact graft identity and the finite graft-compatibility,
radial-comparison, raw \(R^{\pm1}\), raw \(F^{\pm1}\), properness,
degree, right-translated-distance, and metric-dependent lower
singular-value faces, at the fixed orders fourteen for \(R\), six for
\(F\), and twelve/ten for the coefficient/curvature bounds;
\item the normalized source and target operative and reserve harmonic
witnesses, including their modulus-compatibility inequality;
\item the finite physical cover, the frozen reference metric
\(G_{\rm ref}^{\rm phys}\), its reference-tier center-indexed witness
pairs and quantitative Dirichlet containments, the common operative
and \(+\)-tiers, the named outer ball
\eqref{eq:physical-coefficient-ball}, membership of the tuple's actual
closed carrier in the inner locus
\eqref{eq:physical-coefficient-inner-locus}, the uniform half-open
quarter-modulus certificate
\eqref{eq:witnessed-physical-quarter-modulus} on the outer ball, and
the strict \(\mu_{\rm RF}\)-width inequality.
\end{enumerate}
This is a static predicate on a prepared tuple.  It contains neither
the nine slice equalities nor a phase-selection route, and therefore
does not depend on Proposition~\ref{prop:uniform-receding-phase}, its
centered corollary, or the later definition of a strict entrance.  The
later entrance definition inserts the exact
numerical thresholds and assembles this certificate with one of the
two independent ways of imposing the nine moments.  For \(3\leq k<12\)
we use only an explicitly chosen finite subfamily of the ordinary
conjugated prepared faces and no raw or physical certificate.
\end{definition}

\begin{proposition}[Uniform receding phase map]
\label{prop:uniform-receding-phase}
Fix $k\geq3$, $0<\alpha<1$, the indices \(N,\Gamma\), and one numerical
prepared package \(\mathfrak P_{\rm prep}\) as in
\eqref{eq:numerical-prepared-package}.  There are
$\tau_*>0$, $\delta>0$, and $C<\infty$ such that the following holds
for every $\tau_0\geq\tau_*$.  Let \(\mathbf z_0\) be an exact prepared
state with \(h(\mathbf z_0)=0\), and let
\(\mathscr B\subset\mathscr P_{\tau_0}^{k+2,\alpha}\) be a
common-margin prepared ball about \(\mathbf z_0\) satisfying this same
numerical package.  Let
\(\mathbf z_u\in\mathscr B\) satisfy
\[
 \norm{\mathbf z_u-\mathbf z_0}
       _{\mathscr X_{\rm prep}^{k+2,\alpha}}\leq\delta .
\]
There is a unique
$p=p_{\tau_0}(\mathbf z_u)\in\R^9$, $|p|<C\delta$,
such that
\begin{equation}\label{eq:uniform-phase-slice}
 \ip{\rho_{\tau_0}
   \bigl(\mathfrak A^{\rm prep}_{p,\tau_0}(\mathbf z_u)
         -\bar g\bigr)}
   {Z_\mu}=0,\qquad0\leq\mu\leq8.
\end{equation}
The maps $p_{\tau_0}$ are $C^1$ and
\begin{equation}\label{eq:uniform-phase-Lipschitz}
 |p_{\tau_0}(\mathbf z_u)-p_{\tau_0}(\mathbf z_v)|
 \leq C\norm{\mathbf z_u-\mathbf z_v}
              _{\mathscr X^{k+2,\alpha}_{\rm prep}}
\end{equation}
for every pair \(\mathbf z_u,\mathbf z_v\in\mathscr B\) lying in the
\(\delta\)-neighborhood of \(\mathbf z_0\), where the graph-augmented
prepared distance is \eqref{eq:prepared-Banach-norm}.  The phase-adjusted prepared state
\[
 \mathbf A^{\rm prep}_{p_{\tau_0}(\mathbf z_u),\tau_0}(\mathbf z_u)
\]
is controlled at output order \(k\), with the bound determined by
\(\mathfrak P_{\rm prep}\).
The same constant works for every $\tau_0\geq\tau_*$.  More
quantitatively, let $\mathbf z_*$ be a strict sliced state in this
ball, so that \(p_{\tau_0}(\mathbf z_*)=0\), and suppose it satisfies
a fixed finite family \(\mathscr L_{\rm prep}^{(k)}\) of strict
prepared-coordinate entrance and graft inequalities at time
\(\tau_0\).  Apart from the harmonic-radius face, require every member
of this family to have a normalized signed defining functional which
is uniformly locally Lipschitz, on the numerical package, in the
conjugated prepared output-\(k\) norm.  When the harmonic-radius face is
present, it is controlled instead by the lower-stability modulus of
Lemma~\ref{lem:prepared-harmonic-radius-lower-stability}; in that case
the center is required to carry recorded common witness parameters
\(q_{\rm har,*}>0\), \(\zeta_{\rm har,*}\in(0,1)\).  No Lipschitz regularity of
\(\mathfrak h_{\rm har}\) is assumed.
When the family is the phase-independent strict-face certificate of
Definition~\ref{def:phase-independent-strict-face-certificate}, the
center instead carries its recorded normalized and physical operative
and \(+\)-triples.  Preservation uses normalized
reserve-to-operative openness together with physical inner-locus
membership and the fixed reference-carrier certificate.  Thus both
operative certificates and the normalized modulus-compatibility
inequality are unchanged without using the carrier-to-flow modulus.
For \(3\leq k<12\), this family contains no raw map, inverse-map,
right-translated, radial-comparison, or metric-dependent
local-invertibility inequality.

When \(k\geq12\), the family may be enlarged by exactly the fixed raw
part of the phase-independent strict-face certificate: its \(R\)-faces have
order fourteen, its \(F\)-faces have order six, and its
metric-dependent lower singular-value faces have the displayed fixed
orders.  No arbitrary higher-order raw inequality is included in this
quantitative assertion.  Let
\(C_{\rm ph,raw}(\tau_0;\mathscr B)<\infty\) be a fixed-time operator
bound for the structured phase leg in the raw entrance-map charts,
chosen on the displayed common-margin ball so that
\begin{equation}\label{eq:phase-raw-R-leg-bound}
 d_{\rm rt,sc}^{14,\alpha}(R_p,R)
 \leq C_{\rm ph,raw}(\tau_0;\mathscr B)|p|.
\end{equation}
Such a finite bound follows from the ordinary right-translated chart
calculus at fixed \(\tau_0\); no entrance-time-uniform equivalence
with the conjugated prepared norm is asserted.  Define the typed
entrance-control distance by
\begin{equation}\label{eq:phase-entrance-control-distance}
 \mathfrak d_{{\rm ph},{\rm ent},\tau_0;\mathscr B}^{k+2,\alpha}
   (\mathbf z_u,\mathbf z_*):=
 \begin{cases}
  \|\mathbf z_u-\mathbf z_*\|_{\mathscr X_{\rm prep}^{k+2,\alpha}},
       &3\leq k<12,\\[2mm]
  (1+C_{\rm ph,raw}(\tau_0;\mathscr B))
  \|\mathbf z_u-\mathbf z_*\|_{\mathscr X_{\rm prep}^{k+2,\alpha}}
  +d_{\rm rt,sc}^{14,\alpha}(R_u,R_*)\\
  \qquad
  +d_{\rm rt,sc}^{6,\alpha}(F_u,F_*),
       &k\geq12 .
 \end{cases}
\end{equation}
For every ordinary face \(\mathcal F\), take a recorded nonnegative
local Lipschitz modulus \(L_{\mathcal F}\) in the applicable distance
\eqref{eq:phase-entrance-control-distance} and replace it by the
positive majorant
\[
 L_{\mathcal F}^+:=\max\{1,L_{\mathcal F}\}.
\]
If \(\operatorname{slack}_{\mathcal F}(\mathbf z_*)>0\) is its
normalized signed slack, define the always nonempty minimum
\[
 \mathfrak r_{\rm ent}^{\rm Lip}(\tau_0)
 :=\min\left(
 \{1\}\cup
 \left\{
  \frac{\operatorname{slack}_{\mathcal F}(\mathbf z_*)}
       {L_{\mathcal F}^+}:\mathcal F\text{ is an ordinary face}
 \right\}\right)>0.
\]
Thus an empty ordinary-face family and an invariant face with zero
optimal Lipschitz modulus cause no undefined minimum or division.  The
faces included are,
for \(3\leq k<12\), all faces of
\(\mathscr L_{\rm prep}^{(k)}\), and for \(k\geq12\), those faces
together with every face in the fixed raw family just added.
If the harmonic-radius face is among them but the family is not the
phase-independent strict-face certificate, put
\begin{equation}\label{eq:phase-harmonic-radius-slack}
 0<\eta_{\rm har,*}\leq
 \frac14\bigl(
  \mathfrak h_{\rm har}(\tau_0,\mathbf z_*)
  -\kappa_{\rm har}\bigr)
\end{equation}
so small that the recorded witness is valid on every normalized ball
of radius
\((\kappa_{\rm har}+2\eta_{\rm har,*})r_{\rm la}\), and define
\begin{equation}\label{eq:phase-entrance-radius-with-harmonic-modulus}
 \mathfrak r_{\rm ent}(\tau_0):=
 \min\left\{
  \mathfrak r_{\rm ent}^{\rm Lip}(\tau_0),
  \delta_{\rm har}(\eta_{\rm har,*},
                   q_{\rm har,*},\zeta_{\rm har,*};
                   \mathfrak P_{\rm har}^{\rm geom})
 \right\}.
\end{equation}
If the harmonic-radius face is absent, set
\(\mathfrak r_{\rm ent}:=\mathfrak r_{\rm ent}^{\rm Lip}\).
If the full phase-independent strict-face certificate is present,
define instead
\begin{equation}\label{eq:phase-entrance-radius-full-harmonic-modulus}
 \mathfrak r_{\rm ent}(\tau_0):=
 \min\left\{
  \mathfrak r_{\rm ent}^{\rm Lip}(\tau_0),
  \delta_{\rm har}^{\rm op}
  \bigl(
   \eta_{\rm har},q_{\rm har},\zeta_{\rm har};
   \eta_{\rm har}^+,q_{\rm har}^+,\zeta_{\rm har}^+;
   \mathfrak P_{\rm har}^{\rm geom}
   \bigr)
 \right\}.
\end{equation}
No face is omitted merely because the phase leg later leaves its value
invariant: the typed distance transfers every Lipschitz margin, while
the displayed modulus controls the harmonic-radius margin along the
complete phase adjustment.
There is a constant \(c_*>0\), independent of \(\tau_0\), such that
the family just specified remains strict after applying the phase map
whenever
\begin{equation}\label{eq:phase-margin-radius}
 \mathfrak d_{{\rm ph},{\rm ent},\tau_0;\mathscr B}^{k+2,\alpha}
   (\mathbf z_u,\mathbf z_*)
 \leq c_*\mathfrak r_{\rm ent}(\tau_0).
\end{equation}
The radius in \eqref{eq:phase-margin-radius} is allowed to depend on
the entrance time through the actual margin
\(\mathfrak r_{\rm ent}(\tau_0)\).  When \(3\leq k<12\), the uniform
quantitative clause applies only to
\(\mathscr L_{\rm prep}^{(k)}\).  It excludes every raw map,
inverse-map, right-translated, radial-comparison, and
metric-dependent local-invertibility face, whether or not the phase
leg would leave that face invariant.  If \(\mathbf z_u\) is separately
assumed to satisfy an invariant raw face with its own quantified
margin, that face persists, but it is not inferred from
\eqref{eq:phase-margin-radius}.  A merely continuous
finite family still persists after a qualitative fixed-time
shrinking, but no entrance-time-uniform linear radius is claimed for
it.  Every prepared phase leg leaves the closed metric \(G\) itself
unchanged and leaves \(F=\Theta^{-1}\circ\Phi\) unchanged.  Hence, when
\(k\geq12\), the pre-existing order-twelve closed-curvature face and
the fixed raw \(F\) map, inverse-map, and right-translated-distance
faces persist along the phase leg after their center margins have
first been transferred by the typed distance.  The weighted term involving
\(C_{\rm ph,raw}(\tau_0;\mathscr B)\) controls the fixed order-fourteen
raw \(R\) and \(R^{-1}\) faces.  The prepared metric calculus together
with these typed map bounds controls the metric-dependent lower
singular-value faces.  This is precisely the phase-independent
certificate needed
by the continuation theorem.  A generic harmonic-radius member is
preserved by \eqref{eq:prepared-harmonic-radius-lower-stability} and
\eqref{eq:phase-entrance-radius-with-harmonic-modulus}, rather than by
 an unproved linear estimate for \(r_{\rm har}\).  For the full
 phase-independent certificate, the normalized reserve-to-operative
 modulus in
\eqref{eq:phase-entrance-radius-full-harmonic-modulus} preserves the
normalized operative triple and a new reserve.  The typed ordinary-face
distance preserves membership in the fixed physical inner locus, so
Lemma~\ref{lem:reference-carrier-physical-certificate} supplies the
unchanged common physical \(+\)- and operative tiers.  Hence the
recorded normalized modulus-compatibility inequality is preserved
verbatim without using the carrier-to-flow modulus.
\end{proposition}

\begin{proof}
Consider
\[
 \mathcal F_\mu(p,\mathbf z_u)
 =\ip{\rho_{\tau_0}
   \bigl(\mathfrak A^{\rm prep}_{p,\tau_0}(\mathbf z_u)
         -\bar g\bigr)}
   {Z_\mu}.
\]
At the background prepared state $\mathbf z_0$,
\[
 D_{p_j}\mathcal F_\mu(0,\mathbf z_0)
 =\ip{\rho_{\tau_0}\mathscr C_{j,\tau_0}(\mathbf z_0)}{Z_\mu}
 =\ip{\rho_{\tau_0}\mathcal Y_{j,\tau_0}(\mathbf z_0)}{Z_\mu}
 =\ip{Y_j}{Z_\mu}+O(e^{-ce^{\tau_0}}).
\]
The limiting matrix is invertible.
Lemma~\ref{lem:buffered-full-prepared-columns}, applied before any
evolution is constructed, gives uniform first and second derivative
bounds for the scalar moment map, with constants determined by
\(\mathfrak P_{\rm prep}\).  The quantitative implicit-function theorem
therefore gives
a uniform ball, uniqueness, and
\eqref{eq:uniform-phase-Lipschitz}.  The entrance functionals controlled
by the conjugated prepared calculus are locally Lipschitz at output
order \(k\), uniformly on this ball.  This proves the quantitative
claim for \(\mathscr L_{\rm prep}^{(k)}\).  When \(k\geq12\), the
additional fixed raw product faces require the separate one-sided
control recorded in
\eqref{eq:phase-raw-R-leg-bound}: the structured phase action has
\(G_p=G\) and \(F_p=F\) exactly, while the local quasi-triangle
\eqref{eq:right-translated-quasi-triangle},
\eqref{eq:uniform-phase-Lipschitz}, and
\eqref{eq:phase-raw-R-leg-bound} control \(R_p\) and its inverse by
\eqref{eq:phase-entrance-control-distance}.  The same typed distance,
combined with the phase-adjusted source and target metric bounds,
controls the lower singular values in
\eqref{eq:relative-marking-lower-margin} and
\eqref{eq:typed-map-lower-margins}.  The prepared phase-leg
bound and \eqref{eq:uniform-phase-Lipschitz} also give
\[
 \left\|
  \mathbf A^{\rm prep}_{p_{\tau_0}(\mathbf z_u),\tau_0}
       (\mathbf z_u)-\mathbf z_*
 \right\|_{\mathscr X_{\rm prep}^{k,\alpha}}
 \leq C\,
 \mathfrak d_{{\rm ph},{\rm ent},\tau_0;\mathscr B}^{k+2,\alpha}
       (\mathbf z_u,\mathbf z_*).
\]
If a generic harmonic-radius face is present without the full
phase-independent certificate, its center value is at least
\(\kappa_{\rm har}+4\eta_{\rm har,*}\); hence
Lemma~\ref{lem:prepared-harmonic-radius-lower-stability} and
\eqref{eq:phase-entrance-radius-with-harmonic-modulus} preserve it with
strict room after decreasing the universal factor \(c_*\).  When the
full phase-independent certificate is present, use instead
\eqref{eq:phase-entrance-radius-full-harmonic-modulus} and
\eqref{eq:prepared-harmonic-reserve-modulus} for the normalized
certificate.  The typed entrance distance places the input closed
carrier in the inner locus of the center's fixed reference-carrier
ball.  Lemma~\ref{lem:reference-carrier-physical-certificate} then
supplies, at every recorded center \(x\in V_a^5\), the unchanged common
physical \(+\)-tier and operative tier.  Since the phase leg has
\(G_p=G\), it changes neither those witnesses, inner-locus membership,
nor the frozen physical reference.  The signed face
\[
 \operatorname{dist}_{\rm pkg}
 (G_{\rm car},\partial\mathscr B_{\rm coeff}^{\rm phys})
 -\mu_{\rm coeff}^{\rm phys}
\]
and the strict width inequality are ordinary faces controlled by
\(\mathfrak r_{\rm ent}^{\rm Lip}\).  Once inner-locus membership has
been preserved, the uniform half-open quarter-window conclusion is
the precomputed consequence
\eqref{eq:witnessed-physical-quarter-modulus}, not a new finite-jet
functional.  In either
case, because
\(p_{\tau_0}(\mathbf z_*)=0\), reducing this typed distance by a fixed
multiple of the actual strict margin proves
\eqref{eq:phase-margin-radius}.  No assertion is made about the
unadjusted margin of an arbitrary uncentered state.  This last radius need
not be uniform in the unweighted raw product norm as \(\tau_0\) varies:
in addition to entrance thresholds which decay as
\(e^{-\sigma\tau_0}\), in the high-order raw clause the finite factor
\(C_{\rm ph,raw}(\tau_0;\mathscr B)\) is retained explicitly.  The
fixed factor \(C_{\rm rt,\Delta}\) is absorbed into \(c_*\).
\end{proof}

\begin{corollary}[Phase map about a sliced prepared center]
\label{cor:sliced-center-phase-map}
Fix \(k\geq3\) and \(0<\alpha<1\).  Let
\(\mathbf z_c\in\mathscr P_{\tau_c}^{k+2,\alpha}\) be a common-margin
prepared state at time \(\tau_c\geq\tau_*\) satisfying
\[
 \ip{\rho_{\tau_c}h(\mathbf z_c)}{Z_\mu}=0,
 \qquad0\leq\mu\leq8,
\]
and suppose its full phase-column matrix
\[
 M^c_{\mu j}
 =
 \ip{\rho_{\tau_c}
       \mathscr C_{j,\tau_c}(\mathbf z_c)}{Z_\mu}
\]
is invertible with \(\|(M^c)^{-1}\|\leq C_M\).  First choose a
preliminary common-margin ball
\(\mathscr B_0\subset\mathscr P_{\tau_c}^{k+2,\alpha}\) about
\(\mathbf z_c\) on which the prepared column bounds hold.  Then on a
 sufficiently small ball
 \(\mathscr B\Subset_{\rm u}\mathscr B_0\) there is a unique
\(C^1\) prepared phase map \(p_{\tau_c}\), centered by
\(p_{\tau_c}(\mathbf z_c)=0\), which imposes all nine moments in
\eqref{eq:uniform-phase-slice}.  Its Lipschitz constant and admissible
implicit-function radius depend only on \(C_M\) and the common
prepared bounds.  Shrink \(\mathscr B\) once more so that the complete
phase segment satisfies
\begin{equation}\label{eq:centered-phase-segment-containment}
 \mathbf A^{\rm prep}_{s p_{\tau_c}(\mathbf z),\tau_c}(\mathbf z)
 \in\mathscr B_0,\qquad
 \mathbf z\in\mathscr B,\quad0\leq s\leq1 .
\end{equation}

For \(k\geq12\), select on this new preliminary ball its own finite raw
phase-leg constant
\(C_{\rm ph,raw}^{c}(\tau_c;\mathscr B_0)\) such that
\begin{equation}\label{eq:centered-phase-raw-leg-bound}
 d_{\rm rt,sc}^{14,\alpha}(R_p,R)
 \leq C_{\rm ph,raw}^{c}(\tau_c;\mathscr B_0)|p|
\end{equation}
whenever the phase leg remains in \(\mathscr B_0\).  Define the
center-local typed distance
\begin{equation}\label{eq:centered-phase-entrance-distance}
 \mathfrak d_{{\rm ph},c}^{k+2,\alpha}
   (\mathbf z_u,\mathbf z_c):=
 \begin{cases}
  \|\mathbf z_u-\mathbf z_c\|_
    {\mathscr X_{\rm prep}^{k+2,\alpha}},
      &3\leq k<12,\\[2mm]
  (1+C_{\rm ph,raw}^{c}(\tau_c;\mathscr B_0))
  \|\mathbf z_u-\mathbf z_c\|_
    {\mathscr X_{\rm prep}^{k+2,\alpha}}\\
  \qquad
  +d_{\rm rt,sc}^{14,\alpha}(R_u,R_c)
  +d_{\rm rt,sc}^{6,\alpha}(F_u,F_c),
      &k\geq12 .
 \end{cases}
\end{equation}
If \(\mathbf z_c\) satisfies the fixed quantitative family specified
in Proposition~\ref{prop:uniform-receding-phase}, define
\(\mathfrak r_{{\rm ent},c}>0\) by the same rule: use the baseline
\(\{1\}\), divide each ordinary normalized slack by
\(\max\{1,L_{\mathcal F}\}\), and take the resulting nonempty minimum.
When a generic
harmonic-radius face is present without the phase-independent
strict-face certificate, include
the lower-stability threshold formed from its recorded center witness
\((\eta_{{\rm har},c},q_{{\rm har},c},\zeta_{{\rm har},c})\),
\[
 \delta_{\rm har}\!\left(
  \eta_{{\rm har},c},q_{{\rm har},c},\zeta_{{\rm har},c};
 \mathfrak P_{\rm har}^{\rm geom}\right).
\]
When \(\mathbf z_c\) carries the full phase-independent strict-face
certificate, use instead its normalized reserve-to-operative modulus
\[
 \delta_{\rm har}^{\rm op}
 \bigl(
  \eta_{\rm har},q_{\rm har},\zeta_{\rm har};
  \eta_{\rm har}^+,q_{\rm har}^+,\zeta_{\rm har}^+;
 \mathfrak P_{\rm har}^{\rm geom}
 \bigr)
\]
in the minimum with the ordinary-face slack, which includes the signed
inner-locus face
\[
 \operatorname{dist}_{\rm pkg}-\mu_{\rm coeff}^{\rm phys}.
\]
Thus the same normalized operative triple and compatibility inequality,
and the same physical common \(+\)- and operative tiers supplied by the
fixed reference-carrier certificate, are part of the preserved family
without using the carrier-to-flow modulus.
A further shrinking depending on this radius preserves the family
whenever
\[
 \mathfrak d_{{\rm ph},c}^{k+2,\alpha}
   (\mathbf z_u,\mathbf z_c)
 \leq c_c\mathfrak r_{{\rm ent},c}.
\]
Here \(c_c>0\) depends on \(C_M\), the common prepared bounds, and the
fixed local quasi-triangle constant; the raw neighborhood size also
depends, through the weighted distance, on the explicitly selected
\(C_{\rm ph,raw}^{c}(\tau_c;\mathscr B_0)\).

In particular, the conclusion applies at every strict center obtained
by evolving a sliced feedback solution: there \(M^c\) is the exact
receding Gram matrix and is uniformly invertible.
When \(3\leq k<12\), the preservation-of-margins conclusion is
restricted to the uniformly Lipschitz conjugated prepared functionals
specified in the proposition; no raw-map, radial, or
local-invertibility face is inferred unless it is imposed separately
on the input with its own margin.  Preservation
of the complete fixed order-twelve/ten coefficient package and the
order-fourteen \(R\)/order-six \(F\) raw conditions is asserted only for
\(k\geq12\), exactly as in
Proposition~\ref{prop:uniform-receding-phase}.
\end{corollary}

\begin{proof}
Use the moment map \(\mathcal F\) from the proof of
Proposition~\ref{prop:uniform-receding-phase}.  The propagated slice
gives \(\mathcal F(0,\mathbf z_c)=0\), while
\[
 D_p\mathcal F(0,\mathbf z_c)=M^c.
\]
The prepared-chart calculus and the proof of
Lemma~\ref{lem:buffered-full-prepared-columns} give uniform first and
second derivatives of this scalar moment map on a common-margin ball;
that argument only uses the prepared bounds and does not require
\(h(\mathbf z_c)=0\).  The quantitative implicit-function theorem
 therefore gives the asserted centered phase map on a ball
 \(\mathscr B\Subset_{\rm u}\mathscr B_0\).  Continuity of the phase curve and
its uniform prepared bound allow this ball to be shrunk so that
\eqref{eq:centered-phase-segment-containment} holds.  The
prepared-chart calculus gives
the quantitative preservation statement for the uniformly Lipschitz
conjugated functionals.  When \(k\geq12\), ordinary fixed-time raw
chart calculus on \(\mathscr B_0\) gives
\eqref{eq:centered-phase-raw-leg-bound}; the local quasi-triangle
\eqref{eq:right-translated-quasi-triangle} and the phase-map Lipschitz
bound then control the raw \(R\)-faces by
\eqref{eq:centered-phase-entrance-distance}.  The closed \(G\)-face
and the raw \(F\) map, inverse-map, and right-translated-distance faces
are invariant, while the prepared source--target metric calculus
controls the metric-dependent
\(\mathfrak m_F(F;\acute G,S)\) face.  The segment containment
\eqref{eq:centered-phase-segment-containment}, the prepared phase-leg
bound, and
Lemma~\ref{lem:prepared-harmonic-radius-lower-stability} preserve the
 harmonic-radius face by its displayed center-local modulus.  At a
 center carrying the full phase-independent certificate, the same
 shrinking keeps the input carrier in the center's fixed inner
 coefficient locus.  Lemma~%
 \ref{lem:reference-carrier-physical-certificate} therefore supplies,
 at every \(x\in V_a^5\), the common physical \(+\)- and operative tiers
 before the phase leg.  The identity \(G_p=G\) preserves the carrier,
 those witnesses, and the frozen reference metric along the leg.
 Shrinking by
the actual normalized slack proves the stated center-local
preservation radius.  Thus no raw-leg
constant selected on the different exact-center ball is reused here.
\end{proof}

\begin{lemma}[Unbuffered Gaussian moment calculus]
\label{lem:unbuffered-Gaussian-moment-map}
Fix \(r\geq3\), \(0<\alpha<1\), a normalized time \(\tau_c\), and a
common-margin neighborhood in
\(\mathscr P_{\tau_c}^{r,\alpha}\) with numerical package
\eqref{eq:numerical-prepared-package}.  The scalar moment map
\begin{equation}\label{eq:unbuffered-moment-map}
\mathfrak m_\mu(\mathbf z)
 =\ip{\rho_{\tau_c}h(\mathbf z)}{Z_\mu},
 \qquad0\leq\mu\leq8,
\end{equation}
is \(C^1\) on this unbuffered \(r\)-th order prepared manifold.
Its differential is uniformly bounded in the
\(\mathscr E_{\rm prep}^{r,\alpha}\) model tangent norm and depends
continuously on the base state.  If
\(\mathbf z\in\mathscr P_{\tau_c}^{r+2,\alpha}\), so that the
tuple-valued phase curve is differentiable into the order-\(r\)
chart, then
\begin{equation}\label{eq:unbuffered-moment-phase-column}
 D\mathfrak m_\mu(\mathbf z)
  [\partial_{p_j}\mathbf A^{\rm prep}_{p,\tau_c}(\mathbf z)|_{p=0}]
 =
 \ip{\rho_{\tau_c}\mathscr C_{j,\tau_c}(\mathbf z)}{Z_\mu}.
\end{equation}
At a merely \(r\)-regular state the right side still denotes the
integrated-by-parts scalar column, but it is not asserted to arise
from an order-\(r\) tangent phase curve.  In particular, this lemma is
only scalar-valued; it does not assert unbuffered differentiability of
the tensor-valued pullback map.
\end{lemma}

\begin{proof}
For a smooth curve of prepared states, differentiate the exact graph
\[
 h=\lambda^{-1}(\Phi^{-1})^*\acute G-\bar g.
\]
The derivative is the pullback of the direct tensor variation plus a
Lie derivative in the Eulerian variation field of \(\Phi^{-1}\), and
the scalar scale variation.  Pair it with
\(\rho_{\tau_c}Z_\mu e^{-\bar f}\,dV_{\bar g}\).  Because
\(\rho_{\tau_c}\) is compactly supported, integration by parts moves
the derivative in the Lie term onto the fixed Gaussian test tensor,
the density, and the bounded prepared coefficients.  There is no
boundary term.  The resulting expression contains no derivative of a
varying tensor at its top order and is bounded by
\[
 C(\mathfrak P_{\rm prep})
 \|\dot{\mathbf z}\|_{\mathscr E_{\rm prep}^{r,\alpha}} .
\]
The same formula for the difference of two base states, together with
the scale-one product and composition bounds, proves continuity of
the differential.  Polynomial coefficient growth on the receding
part is harmless by the Gaussian tail estimate.  Approximation in the
little-H\"older chart extends the bounded differential formula, for
genuine order-\(r\) tangent curves, from smooth states to the entire
unbuffered neighborhood.  When the base state has the displayed
two-derivative buffer, the prepared phase leg is such a curve and
gives \eqref{eq:unbuffered-moment-phase-column}.  No phase leg at a
merely \(r\)-regular base state is used in this conclusion.
\end{proof}

\begin{proposition}[Sliced prepared manifold and buffered phase projection]
\label{prop:sliced-prepared-manifold}
Fix \(r\geq3\), and let the smooth center \(\mathbf z_c\) satisfy the
hypotheses of
Corollary~\ref{cor:sliced-center-phase-map}.  There is a sufficiently
small open common-margin neighborhood
\[
 \mathscr N_{\mathrm{sl},\tau_c}^{r,\alpha}
 \subset\mathscr P_{\tau_c}^{r,\alpha}
\]
of \(\mathbf z_c\) on which
\[
 \mathfrak m(\mathbf z)
 =
 \left(
  \ip{\rho_{\tau_c}h(\mathbf z)}{Z_\mu}
 \right)_{\mu=0}^8
\]
is \(C^1\) and has split-surjective differential.  Define
\begin{equation}\label{eq:sliced-prepared-manifold}
 \Sigma_{\tau_c}^{r,\alpha}
 :=
 \left\{
  \mathbf z\in\mathscr N_{\mathrm{sl},\tau_c}^{r,\alpha}:
  \mathfrak m(\mathbf z)=0
 \right\}.
\end{equation}
Then \(\Sigma_{\tau_c}^{r,\alpha}\) is a split \(C^1\) Banach
submanifold of codimension nine, and
\[
 T_{\mathbf z}\Sigma_{\tau_c}^{r,\alpha}
 =\ker D\mathfrak m(\mathbf z).
\]
Put
\[
 \mathscr N_{\mathrm{sl},\tau_c}^{r+2,\alpha}
 :=
 \mathscr N_{\mathrm{sl},\tau_c}^{r,\alpha}
 \cap\mathscr P_{\tau_c}^{r+2,\alpha},
 \qquad
 \Sigma_{\tau_c}^{r+2,\alpha}
 :=
 \Sigma_{\tau_c}^{r,\alpha}
 \cap\mathscr P_{\tau_c}^{r+2,\alpha}.
\]
On a sufficiently small two-derivative-buffered neighborhood
\[
 \mathscr O_{\tau_c}^{r+2,\alpha}
 \subset\mathscr N_{\mathrm{sl},\tau_c}^{r+2,\alpha},
\]
the centered phase map defines a \(C^1\) phase projection
\begin{equation}\label{eq:phase-retraction}
 \Pi_{\rm sl}^{\,r+2\to r}(\mathbf z)
 :=
 \mathbf A^{\rm prep}_{p_{\tau_c}(\mathbf z),\tau_c}(\mathbf z)
 \in\Sigma_{\tau_c}^{r,\alpha}.
\end{equation}
On the higher-regularity sliced part of its domain it is the canonical
inclusion:
\[
 \left.\Pi_{\rm sl}^{\,r+2\to r}\right|_{
   \Sigma_{\tau_c}^{r+2,\alpha}
   \cap\mathscr O_{\tau_c}^{r+2,\alpha}}
 =
 \iota_{r+2,r}:
 \Sigma_{\tau_c}^{r+2,\alpha}
 \cap\mathscr O_{\tau_c}^{r+2,\alpha}
 \hookrightarrow\Sigma_{\tau_c}^{r,\alpha}.
\]
After shrinking \(\mathscr O_{\tau_c}^{r+2,\alpha}\) to a convex
model-ball inside the same common-margin chart, there is a derived
constant \(K_{\Pi,r}<\infty\) such that, for all
\(\mathbf z_1,\mathbf z_2\in
  \mathscr O_{\tau_c}^{r+2,\alpha}\),
\begin{align}
 \|\Pi_{\rm sl}^{\,r+2\to r}(\mathbf z_1)
      -\Pi_{\rm sl}^{\,r+2\to r}(\mathbf z_2)\|_{
      \mathscr X_{\rm prep}^{r,\alpha}}
 &\leq K_{\Pi,r}
 \|\mathbf z_1-\mathbf z_2\|_{
      \mathscr X_{\rm prep}^{r+2,\alpha}},
 \label{eq:quantitative-phase-retraction-Lipschitz}\\
 \sup_{\mathbf z\in\mathscr O_{\tau_c}^{r+2,\alpha}}
 \|D\Pi_{\rm sl}^{\,r+2\to r}(\mathbf z)\|_{
  \mathcal L(\mathscr E_{\rm prep}^{r+2,\alpha},
             \mathscr E_{\rm prep}^{r,\alpha})}
 &\leq K_{\Pi,r}.
 \label{eq:quantitative-phase-retraction-derivative}
\end{align}
This constant is supplied by the quantitative phase-column inverse and
prepared-chart calculus, not by compactness of a bounded subset.
We suppress the superscript \(r+2\to r\) whenever the input and output
orders are fixed by context.

A \emph{common-margin sliced ball} means a relatively open bounded
ball in \(\Sigma_{\tau_c}^{r,\alpha}\) on which the prepared margins
are common.  If
\(\mathcal E_{\rm sl}:\mathscr B_{\rm sl}\to\mathbb B_{\rm tar}\)
is a \(C^1\) map into one of the Banach target charts used below, and
if an ambient neighborhood
\(\mathscr O\subset\mathscr P_{\tau_c}^{r+2,\alpha}\) satisfies
\[
 \mathscr O\subset\operatorname{dom}
     \Pi_{\rm sl}^{\,r+2\to r},
 \qquad
 \Pi_{\rm sl}^{\,r+2\to r}(\mathscr O)
     \subset\mathscr B_{\rm sl},
\]
its buffered ambient extension is
\[
 \mathcal E_{\rm amb}
 =
 \mathcal E_{\rm sl}\circ\Pi_{\rm sl}^{\,r+2\to r}
\]
on \(\mathscr O\), and
\[
 D\mathcal E_{\rm amb}
 =D\mathcal E_{\rm sl}\circ
 D\Pi_{\rm sl}^{\,r+2\to r}.
\]
If \(\overline{\mathscr O}\subset\operatorname{dom}
\Pi_{\rm sl}^{\,r+2\to r}\),
\(\mathscr B'\Subset_{\rm u}\mathscr B_{\rm sl}\),
\(\Pi_{\rm sl}^{\,r+2\to r}(\overline{\mathscr O})
\subset\mathscr B'\), and there is \(L_{\mathcal E}<\infty\) such
that, in the chosen target chart, for all \(u,v\in\mathscr B'\),
\[
 \|\mathcal E_{\rm sl}(u)-\mathcal E_{\rm sl}(v)\|_{\mathbb B_{\rm tar}}
 \leq L_{\mathcal E}
       \|u-v\|_{\mathscr X_{\rm prep}^{r,\alpha}},
 \qquad
 \sup_{u\in\mathscr B'}\|D\mathcal E_{\rm sl}(u)\|
 \leq L_{\mathcal E},
\]
then
\begin{equation}\label{eq:quantitative-ambient-extension-bound}
 \operatorname{Lip}(\mathcal E_{\rm amb}|_{\mathscr O})
 +\sup_{z\in\mathscr O}\|D\mathcal E_{\rm amb}(z)\|
 \leq2L_{\mathcal E}K_{\Pi,r}.
\end{equation}
Thus the closure and uniformly-interior hypotheses retain domains and
prepared margins; they are not used as a compactness substitute.
\end{proposition}

\begin{proof}
Lemma~\ref{lem:unbuffered-Gaussian-moment-map} proves that
\(\mathfrak m\) is \(C^1\) on the unbuffered
\(\mathscr P_{\tau_c}^{r,\alpha}\) neighborhood.  Since the center is
smooth, the nine vectors
\[
 e_j=\left.\partial_{p_j}
 \mathbf A^{\rm prep}_{p,\tau_c}(\mathbf z_c)\right|_{p=0}
\]
are genuine members of
\(T_{\mathbf z_c}\mathscr P_{\tau_c}^{r,\alpha}\), and the matrix of
\(D\mathfrak m(\mathbf z_c)\) on their span is the full phase-column
matrix.  Continue this span as a constant finite-dimensional
subbundle in a prepared Banach chart.  Continuity of
\(D\mathfrak m\), followed by shrinking, keeps its restriction
invertible.  Hence \(D\mathfrak m\) is split surjective on the
unbuffered neighborhood without invoking a derivative-losing phase
curve at a nonsmooth point.  The Banach implicit-function theorem
proves \eqref{eq:sliced-prepared-manifold}.  The centered
implicit phase parameter is the unique solution of
\[
 \mathfrak m\!\left(
   \mathbf A^{\rm prep}_{p,\tau_c}(\mathbf z)
 \right)=0.
\]
The buffered prepared-chart calculus makes the resulting tuple-valued
map \(C^1\) from input order \(r+2\) to output order \(r\).  With
\[
 \mathcal F(\mathbf z,p)
 :=\mathfrak m\!\left(
   \mathbf A^{\rm prep}_{p,\tau_c}(\mathbf z)
 \right),
\]
the quantitative inverse-column bound from
Corollary~\ref{cor:sliced-center-phase-map} and the uniform prepared
action bounds give
\[
 Dp_{\tau_c}=-(D_p\mathcal F)^{-1}D_{\mathbf z}\mathcal F,
 \qquad
 D\Pi_{\rm sl}
 =D_{\mathbf z}\mathbf A^{\rm prep}
  +D_p\mathbf A^{\rm prep}\,Dp_{\tau_c}.
\]
After the stated shrinking these formulas yield the finite constant
\(K_{\Pi,r}\) in
\eqref{eq:quantitative-phase-retraction-Lipschitz}--%
\eqref{eq:quantitative-phase-retraction-derivative}.  If
\[
 \mathbf z\in
 \Sigma_{\tau_c}^{r+2,\alpha}
 \cap\mathscr O_{\tau_c}^{r+2,\alpha},
\]
then \(p=0\) is a solution; uniqueness and
\(\mathbf A^{\rm prep}_{0,\tau_c}=\operatorname{Id}\) give the
canonical-inclusion identity.  The ambient derivative formula and
\eqref{eq:quantitative-ambient-extension-bound} follow from the chain
rule and the product of the sliced and retraction bounds.
\end{proof}

\begin{corollary}[Fixed-convention preparation map]
\label{cor:fixed-convention-preparation-map}
Fix \(r\geq3\), a normalized time \(\tau_c\), and one marked host,
graft, scale, and initial-map convention.  On a sufficiently small
\(C^{r+2,\alpha}\) neighborhood of a closed metric for which these
objects lie in a common prepared chart, the convention defines a
\(C^1\) raw lift
\[
 G\longmapsto\mathbf z_{\rm raw}(G)
 \in\mathscr P_{\tau_c}^{r+2,\alpha}.
\]
If the centered phase-column matrix is uniformly invertible, then
\begin{equation}\label{eq:fixed-convention-preparation-map}
 \mathcal P_{\tau_c}(G)
 :=
 \Pi_{\rm sl}^{\,r+2\to r}
 \bigl(\mathbf z_{\rm raw}(G)\bigr)
 \in\Sigma_{\tau_c}^{r,\alpha}
\end{equation}
is \(C^1\).  Every strict finite-jet or weighted prepared inequality
whose defining norms use at most output order \(r\) persists after
shrinking the metric neighborhood.  If the center carries recorded
operative and \(+\)-harmonic witness triples, then the dimensionless
harmonic-radius face persists by the reserve-to-operative clause of
Lemma~\ref{lem:prepared-harmonic-radius-lower-stability}; the same
operative triple and hence its modulus-compatibility inequality are
retained.
Persistence of the full phase-independent strict-face certificate in
Definition~\ref{def:phase-independent-strict-face-certificate} is
asserted only for \(r\geq12\); at those orders the finite
physical-cover certificate persists because the nearby closed carrier
remains in the same inner coefficient locus
\eqref{eq:physical-coefficient-inner-locus}, based at the center's
frozen \(G_{\rm ref}^{\rm phys}\), and
Lemma~\ref{lem:reference-carrier-physical-certificate} supplies its
common \(+\)- and operative tiers.  The phase projection does not alter
that carrier or use the carrier-to-flow reserve.  Hence the half-open
quarter-modulus future-window
estimate \eqref{eq:witnessed-physical-quarter-modulus}, already uniform
on the outer ball, and the \(\mu_{\rm RF}\)-width slack persist as well.
If the center
metrics range over a compact already-sliced family with a uniform
Gram inverse and common margins, finitely many local maps
\eqref{eq:fixed-convention-preparation-map} glue, after shrinking, to
one \(C^1\) preparation map on a neighborhood of that family.
\end{corollary}

\begin{proof}
With the host, cutoff, scale, and maps fixed, the raw lift is assembled
from affine dependence on \(G\) and the inverse, pullback,
composition, and graft operations controlled by
Lemma~\ref{lem:prepared-chart-calculus}; hence it is \(C^1\) with the
two-derivative buffer displayed above.  Proposition~%
\ref{prop:sliced-prepared-manifold} and the chain rule prove the first
claim.  Strict inequalities are open in the stated finite-jet and
weighted norms.  The raw lift and phase projection are continuous in
the scale-one \(C^{2,\alpha}\) metric and map variables controlled by
the prepared norm.  When the stated two-tier harmonic certificate is
present, the reserve-to-operative clause of
Lemma~\ref{lem:prepared-harmonic-radius-lower-stability} applies after
shrinking by its recorded operative and reserve gaps.  At output order
\(r\geq12\), for a tuple carrying the full phase-independent
strict-face certificate, continuity of the raw lift and phase
projection places the nearby closed carrier in the same inner
coefficient locus with positive residual slack.  The phase projection
does not alter that carrier or the frozen reference metric.
Lemma~\ref{lem:reference-carrier-physical-certificate} then supplies,
at every recorded center \(x\in V_a^5\), the same common physical
\(+\)-tier and operative tier.  Since the
quarter-window certificate was proved uniformly for every metric in
the outer ball, it follows from membership rather than from a new
continuity argument for Ricci flows.  Continuity of the scale component
preserves the
\(\mu_{\rm RF}\)-strict width inequality.  On a compact center family both the
prepared and physical witness reserves are made uniform by a finite
subcover before the source neighborhood is chosen.  On
overlaps, centered uniqueness identifies the
local phase parameters, so a finite compact subcover glues the maps.
\end{proof}

\section{Prepared entrances and localized basins}
\label{sec:prepared-open-basin}

This section proves Theorem~A in its natural dependency order.  It
first closes one-state continuation from a strict entrance, then
implants an exact FIK core, promotes the resulting high-regularity
neighborhood through positive-time smoothing to the relative
\(C^{2,\alpha}\) topology, and finally states the conclusion as an
invariant marked basin.  The evolutionary closure proving Theorem~A
uses only one-state estimates.  The unified construction nevertheless
reserves the static coefficient and radius constants required later for
Theorem~B, including \(\Gamma_{\rm B}\) in
\eqref{eq:global-compatible-package-radius}, from the outset.  These
entries are not used in the formation argument: no two-state
trajectory, difference estimate, or variational conclusion is used
there.  The finite-horizon comparison theory has already been
established; its
global scattering and asymptotic-differentiability use begins only in
Section~\ref{sec:global-two-state}, after the proof of Theorem~A is
complete.

\subsection[Strict entrances and global continuation]
{Strict prepared entrances and global continuation}

The analytic theorem is stated first for an arbitrary strict prepared
entrance.  Its nonemptiness will then follow from a direct exact-core
implant, rather than from a pre-existing tuned singular trajectory.

First fix a continuation output order \(k_0\geq12\) and
\(0<\alpha<1\).  Next fix
\[
 0<\sigma<\theta<\beta,
\]
then fix a rate-compatible numerical prepared package satisfying
\eqref{eq:global-compatible-package-radius}, and hence
\eqref{eq:three-region-compatible-package-radius}.  With those data
frozen, fix
\[
 0<\varepsilon\leq\varepsilon_{\rm ent},
\]
where \(\varepsilon_{\rm ent}\) is no larger than every smallness
threshold in the feedback, three-region, graft, and continuation
arguments for this fixed pair and package.  The number \(\varepsilon\)
is part of the entrance data;
all statements below are uniform for
\(0<\varepsilon\leq\varepsilon_{\rm ent}\) once the remaining common
margins are fixed.
Strict entrances are points of the finite
little-H\"older prepared manifold at input order \(k_0+2\), not merely
points of its smooth subclass.  Smooth strict entrances form the dense
subclass used for the closed geometric flows in Theorem~A.  The
finite-regularity class is the actual open Banach domain used for
continuous and \(C^1\) dependence in
Theorem~\ref{thm:intro-sharp-scattering}; its solutions satisfy the
compact and chartwise Bochner integral identities specified in
Subsection~\ref{subsec:time-dependent-prepared-spaces}.  The closed
Ricci--DeTurck and harmonic-map components acquire their usual
positive-time parabolic regularity, but the ODE-carried
\(R,\Theta\), and target data retain the spatial regularity of their
entrance values.  Hence the full prepared tuple is smooth at positive
time only for a smooth entrance.  Smoothness never substitutes for a
weighted high-order hypothesis on the noncompact prepared end.
The common numerical constants are fixed by
\eqref{eq:numerical-prepared-package}, and the single adaptive order,
phase threshold, and scale-comparison constant are those of
Remark~\ref{conv:authoritative-adaptive-order}.
The support scalar is already the dynamic functional in
\eqref{eq:dynamic-harmonic-and-separation-faces}.  Thus
\(\mathfrak s_{\rm sep}(\tau_0)>0\) quantitatively implies the
first disjointness in
\eqref{eq:coarse-Gram-support-separation}.

\begin{definition}[Quantitative conditions for a strict prepared entrance]
\label{def:strict-prepared-entrance}
This is the operative definition of the input class in
Theorem~\ref{thm:intro-sharp-scattering}.  A tuple is a strict prepared
entrance if and only if it satisfies the following full quantitative
conditions.  The preceding
Remark~\ref{rem:intro-Theorem-C-entrance-preview} gives their geometric
content; the list below is the complete numerical expansion of
the phase-independent strict-face certificate in
Definition~\ref{def:phase-independent-strict-face-certificate},
assembled with the exact nine slice equalities only after one of the
two independent routes below has been chosen.  A point of the finite
little-H\"older prepared chart, whose
closed-metric component is a
\(C^{k_0+2,\alpha}\) metric $G_0$ on a closed manifold \(\mathcal X\),
is a \emph{strict prepared FIK entrance of continuation order \(k_0\)}
if it admits a marked FIK chart,
an initial scale
$\lambda_0$, an adaptive exterior, nested graft collars
\[
 \Omega_\eta\Subset\Omega_\eta^+
 \Subset\Omega_\eta^{++},
\]
and initial maps $\Theta_0,F_0$, all at one normalized time
$\tau_0$.  The initial slice is prepared by either of the following
two routes.
\begin{enumerate}
\item[(a)] The raw tuple lies in the specified neighborhood
\[
 \operatorname{dom}
 \Pi_{\rm sl}^{\,k_0+4\to k_0+2}
 \subset\mathscr P_{\tau_0}^{k_0+4,\alpha}
\]
of a fixed sliced center and has the common margins required by that
phase chart, including the phase-independent strict-face certificate.
Apply the centered phase map of
Corollary~\ref{cor:sliced-center-phase-map}, equivalently the buffered
phase retraction of
Proposition~\ref{prop:sliced-prepared-manifold},
\[
 \Pi_{\rm sl}^{\,k_0+4\to k_0+2},
\]
once and rename the resulting phase-adjusted tuple
\((G_0,\lambda_0,\Theta_0,\Phi_0)\).  Proposition
\ref{prop:uniform-receding-phase} and
Corollary~\ref{cor:sliced-center-phase-map} preserve that static
certificate, including the normalized two-tier and physical three-tier
harmonic witness data.
\item[(b)] The raw tuple already satisfies the exact nine slice
equalities, has the uniformly invertible adaptive Gram matrix and the
numerical prepared package required below, belongs to a common-margin
ball in \(\mathscr P_{\tau_0}^{k_0+2,\alpha}\), and directly satisfies
 all of the strict inequalities below, equivalently the same
phase-independent strict-face certificate.  In this already-sliced case
 set \(p=0\).  No additional smallness relative to an exact \(h=0\)
 center is required beyond the listed strict faces.
\end{enumerate}
In both cases \(F_0=\Theta_0^{-1}\circ\Phi_0\).  We set
\(\Theta_{\tau_0}:=\Theta_0\), \(\Phi_{\tau_0}:=\Phi_0\), and
\(F_{\tau_0}:=F_0\).  Every condition below refers to this
 prepared sliced tuple, which is required to be a point of one fixed
 open common-margin ball in
\(\mathscr P_{\tau_0}^{k_0+2,\alpha}\).  The physical clock is
normalized by
 \[
  t(\tau_0):=0,\qquad G(t(\tau_0)):=G_0.
 \]
 This normalization is part of the entrance convention, not an
 additional prepared-state coordinate.  For the pre-radius energy,
 forcing, and feedback constants fixed in
 \eqref{eq:pre-radius-three-region-inputs}, put
 \begin{equation}\label{eq:strict-entrance-historical-constants}
  C_{\rm hist}^{\rm pre}
  :=K_{\rm fb}^{\rm pre}
     \bigl(1+(C_{\mathcal E}^{\rm pre})^2\bigr),
  \qquad
  c_{\rm hist}^{\rm pre}
  :=\frac12\min\bigl\{(K_{\rm fb}^{\rm pre})^{-1},
                         2c_{\mathcal E}^{\rm pre}\bigr\}>0.
 \end{equation}
 These are the worst-case pre-package instance of the constants in
 Lemma~\ref{lem:future-phase-tail}; they are fixed before the entrance
 time and before the actual package radius.  We also require
 \begin{equation}\label{eq:strict-entrance-adaptive-time}
  \tau_0\geq\tau_{\rm ad},
 \end{equation}
 with the single threshold fixed in
 \eqref{eq:common-adaptive-entrance-time}.  All inequalities in the
 following list hold with a common positive margin:
\begin{enumerate}
\item the exact receding slice holds, and, for fixed constants
      \(0<c_{\rm sc}<C_{\rm sc}<\infty\),
      \[
       c_{\rm sc}e^{-\tau_0}<\lambda_0<C_{\rm sc}e^{-\tau_0};
      \]
\item for fixed $0<\sigma<\theta<\beta$, the initial tensor satisfies,
      with strict inequalities,
      \[
       \norm{\rho_{\tau_0}h_0}_{L^2_\nu}
       <\varepsilon e^{-\theta\tau_0},\qquad
       \sum_{\ell=0}^3|\bar\nabla^\ell h_0|
       <\varepsilon\omega_\sigma(\tau_0,\cdot),
      \]
       and its radius-independent pre-atlas activation norm obeys
       \[
        \|h_0\|_{\mathfrak C_{{\rm pre},0}^{2,\alpha}}
        <\delta_{\rm c2},
       \]
       while its ellipticity lies strictly inside
      the fixed \(\Lambda_{\rm ell}\) box.  The entrance time satisfies,
      with strict room, the pre-package instance of
      \eqref{eq:historical-tail-absorption}:
      \begin{equation}\label{eq:strict-entrance-historical-tail-absorption}
       C_{\rm hist}^{\rm pre}
       e^{-c_{\rm hist}^{\rm pre}e^{\tau_0}}
       <\varepsilon^2e^{-2\theta\tau_0};
      \end{equation}
\item the extension is genuinely prepared: outside the fixed
      physical graft it equals
      $S_0=\lambda_0\Theta_0^*\bar g$, the initial adaptive radial map
      $R_{\tau_0}=\varphi_{-\tau_0}\circ\Theta_{\tau_0}$ has the
      radial-comparison bound
      \eqref{eq:adaptive-initial-radial-comparison};
      \(R_{\tau_0}\) and \(R_{\tau_0}^{-1}\) have bounded
      \(C^{14,\alpha}\) norm on \(\Omega_\eta^+\) and a common
      \(C^{14,\alpha}\) bound
      after rescaling domain and range by \(L^{-1}\) on every dyadic
      annulus \(\{L<\bar f<4L\}\), uniformly for \(L\geq\Gamma\).
      Equivalently, with strict room,
      \[
       \|R_{\tau_0}^{\pm1}\|_
        {\operatorname{Map}_{\rm sc}^{14,\alpha}}
       <\Lambda_{R,14}^{\rm pre}-2\mu_R^{\rm pre}.
       \]
       Its forward-and-inverse local-invertibility face is, explicitly,
       \begin{equation}\label{eq:strict-entrance-R-singular-face}
        \mathfrak m_R(R_{\tau_0})>\kappa_{\rm map}.
       \end{equation}
       This is the explicit quantitative \(R\)-component of the strict
       entrance package.  It strengthens the unspecified positive
       local-invertibility margin in the ambient common-margin locus to
       the fixed primitive threshold \(\kappa_{\rm map}\); it is not
       inferred from qualitative common-margin membership.
       In addition, for every \(0\leq m\leq4\),
       \begin{equation}\label{eq:strict-entrance-column-faces}
        \max_{0\leq j\leq8}\sup_M
        \sum_{\ell=0}^{m}(1+\bar f)^{\ell/2}
        |\bar\nabla^\ell\mathcal Y_{j,\tau_0}|_{\bar g}
        <K_{\mathcal Y,m}^{\rm pre}.
       \end{equation}
       The exact adaptive Gram and scaled support-separation faces are
       \begin{equation}\label{eq:strict-entrance-Gram-separation-faces}
        s_{\min}M^{\rm low}(\tau_0)>\kappa_{\rm Gram},
        \qquad
        \mathfrak s_{\rm sep}(\tau_0)>\kappa_{\rm sep}.
       \end{equation}
\item the graft parameter $\Gamma$ is the fixed rate-compatible package
      value satisfying
      \eqref{eq:global-compatible-package-radius}.  The graft
      discrepancy obeys the complete strict-improvement margin
      \begin{equation}\label{eq:strict-entrance-complete-graft-margin}
       C_{6,K_{\rm gr}}\left(
        \mathfrak d_{6;\Omega_\eta^{++}}(\iota_*G_0,S_0)
        +\frac{\lambda_0}{\Gamma}
        +\varepsilon_{\rm ph}\right)
       <\frac12K_{\rm gr},
      \end{equation}
      and the interpolation defect computed
      from $(\iota_*G_0,S_0,\eta)$ satisfies the support and
      fixed-background $C^2$ bounds in
      \eqref{eq:outer-forcing-hyp}, with strict margin;
\item the initial relative map
      $F_0$ is a proper diffeomorphism whose scale-normalized
       distance from the identity satisfies
       \[
        d_{\rm rt,sc}^{6,\alpha}(F_0,\operatorname{Id})
        <\varepsilon_{\rm map}^{\rm HM},
       \]
       and its ordinary \(C^{6,\alpha}\) right-translated distance on
       \(\Omega_\eta^+\) is also strictly below
       \(\varepsilon_{\rm map}^{\rm HM}\).  Moreover,
       \begin{equation}\label{eq:strict-entrance-harmonic-coefficient-faces}
       \mathfrak h_{\rm har}(\tau_0)>\kappa_{\rm har},
       \qquad
       \mathfrak C_{\rm ent}^{12,10}(\mathscr J_0)
       <\Lambda_{\rm coef},
      \end{equation}
      and the harmonic inequality is witnessed quantitatively: there
      are recorded operative and reserve triples
      \[
       \begin{gathered}
        0<\eta_{\rm har}<\eta_{\rm har}^+,\qquad
        0<q_{\rm har}<q_{\rm har}^+<Q_{\rm har}-1,\\
        0<\zeta_{\rm har}<\zeta_{\rm har}^+<1,
       \end{gathered}
      \]
      such that, after the normalized graph reduction
      \eqref{eq:harmonic-radius-normalized-graph-identity}, both metric
      branches at every center possess coefficient-\(q_{\rm har}^+\)
      harmonic charts on the ball of radius
      \((\kappa_{\rm har}+2\eta_{\rm har}^+)r_{\rm la}\), with domain
      buffer \(\zeta_{\rm har}^+\).  By restriction these charts also
      carry the unchanged operative triple
      \((\eta_{\rm har},q_{\rm har},\zeta_{\rm har})\).
      Both triples and the strict componentwise gaps are part of the
      common positive margin.  The reserve-to-operative clause of
      Lemma~\ref{lem:prepared-harmonic-radius-lower-stability} is the
      openness statement for this certificate.  The bare strict inequality
      for the supremal radius is not used later as an unrecorded
      perturbation reserve.  For the already fixed comparison constant
      \(C_{\rm har,pre}\geq1\) in
      \eqref{eq:pre-atlas-to-scale-one-graph} from the
      radius-independent pre-atlas \(C^{2,\alpha}\) norm to the
      scale-one graph distance, the witness package is required to
      satisfy
      \begin{equation}\label{eq:strict-entrance-harmonic-modulus-compatibility}
       4C_{\rm har,pre}\delta_{\rm c2}
       <
       \delta_{\rm har}
       (\eta_{\rm har},q_{\rm har},\zeta_{\rm har};
        \mathfrak P_{\rm har}^{\rm geom}).
      \end{equation}
      Thus item~(5) is exactly the normalized part of
      Definition~\ref{def:witnessed-auxiliary-harmonic-package};
      \eqref{eq:strict-entrance-harmonic-modulus-compatibility} repeats
      \eqref{eq:auxiliary-harmonic-modulus-compatibility} among the
      entrance conditions.
      This is a compatibility of already fixed numerical margins, not
      a continuity assertion for the radius functional;
      and, intrinsically---equivalently, by the metric-weighted
      coordinate convention in
       \eqref{eq:bundle-map-smallest-singular-value}---one has
       \begin{equation}\label{eq:strict-entrance-map-singular-face}
        \mathfrak m_F(F_0;\acute G_0,S_0)
        =\min\{s_{\min}(dF_0),s_{\min}(dF_0^{-1})\}
        >\kappa_{\rm map}.
       \end{equation}
       With the same common map and local-invertibility margins,
      \[
       \|F_0^{\pm1}\|_
        {\operatorname{Map}_{\rm sc}^{6,\alpha}}
       <\Lambda_{F,6}^{\rm pre}-2\mu_F^{\rm pre}.
      \]
\item there are fixed smooth core domains
      \[
       K_-\Subset K_0\Subset K_1\Subset K_2\Subset\mathcal X''
       \]
       with
       \(\iota(\overline{K_2})\Subset\{\bar f<\Gamma/2\}\),
       \(\overline{\Omega_\eta^{++}}\Subset\{\bar f>\Gamma/2\}\), and
       \(\iota(\overline{K_-})\cap
         \overline{\Omega_\eta^{++}}=\varnothing\), and
       \(\iota(K_2\setminus\overline{K_-})
        \Subset\operatorname{int}\{\eta=1\}\).  Set
      \[
       E^{++}:=\mathcal X\setminus\overline{K_-},\quad
       E^+:=\mathcal X\setminus\overline{K_0},\quad
       E:=\mathcal X\setminus\overline{K_2}.
       \]
       The package fixes a compact collar
       \(\mathcal A_{\rm in}\Subset
       K_1\setminus\overline{K_0}\) and the terminating cutoff
       \(\zeta\) described above
       \eqref{eq:prepared-exterior-termination}.
       It also fixes two collars with disjoint analytic roles.  The
       interface collar satisfies
       \[
        \mathcal W_{\rm in}\Subset\mathcal W_{\rm in}^+
        \Subset\mathcal W_{\rm in}^{++}
        \Subset E^+\cap\mathcal X'',
        \qquad
        \mathcal A_{\rm in}
        \Subset\operatorname{int}\mathcal W_{\rm in}.
       \]
       The physical graft-input collar lies strictly farther out:
       \[
        \mathcal W_{\rm gr}\Subset\mathcal W_{\rm gr}^+
        \Subset\mathcal W_{\rm gr}^{++}
        \Subset E\cap\mathcal X'',
        \qquad
        \Omega_\eta^{++}
        \Subset\operatorname{int}\iota(\mathcal W_{\rm gr}).
       \]
       Between the last two retained graft collars the package fixes
       the auxiliary chain
       \[
        \mathcal W_{\rm gr}^+
        \Subset\mathcal W_{\rm gr}^{0}\Subset\cdots
        \Subset\mathcal W_{\rm gr}^{5}
        \Subset\mathcal W_{\rm gr}^{++},
       \]
       whose successive scaled separations have one positive lower
       bound among the fixed package constants.  No collar contained
       in \(E\) is required to contain the interior interface
       \(\mathcal A_{\rm in}\).
       There are finitely many
       ordinary buffered physical subsets
      \[
       U_a\Subset U_a^+\Subset U_a^{++}\Subset E^{++},
       \qquad1\leq a\leq N_{\rm ext},
       \]
       with scales \(R_a>0\), and auxiliary chains
       \[
       U_a^{++}\Subset V_a^0\Subset\cdots\Subset V_a^5
       \Subset E^{++}.
      \]
       For each \(a\) the package also fixes a harmonic witness
       enlargement
       \[
        V_a^5\Subset W_a^{\rm har}
        \Subset\widetilde W_a^{\rm har}\Subset E^{++}
       \]
       and common operative, carrier-reserve, and reference physical
       witness parameters
       \[
        \begin{gathered}
        \upsilon_{\rm har}^{\rm phys}>0,\qquad
        0<\eta_{\rm har}^{\rm phys}
          <\eta_{\rm har}^{\rm phys,+}
          <\eta_{\rm har}^{\rm phys,ref},\\
        0<q_{\rm har}^{\rm phys}
          <q_{\rm har}^{\rm phys,+}
          <q_{\rm har}^{\rm phys,ref}<Q_{\rm har}-1,\qquad
        0<\zeta_{\rm har}^{\rm phys}
          <\zeta_{\rm har}^{\rm phys,+}
          <\zeta_{\rm har}^{\rm phys,ref}<1.
        \end{gathered}
       \]
       The physical part of the entrance certificate carries one frozen
       reference metric \(G_{\rm ref}^{\rm phys}\).  For every
       \(x\in V_a^5\), choose once and for all a reference-tier witness
       pair \((u_{a,x},D_{a,x})\) for which
       \(G_{\rm ref}^{\rm phys}\) has a
       coefficient-\(q_{\rm har}^{\rm phys,ref}\),
       domain-\(\zeta_{\rm har}^{\rm phys,ref}\) harmonic witness on
       \[
        B_{G_{\rm ref}^{\rm phys}}\!\left(
          x,
          (\upsilon_{\rm har}^{\rm phys}
           +2\eta_{\rm har}^{\rm phys,ref})R_a
        \right).
       \]
       Its recorded Euclidean Dirichlet domain has manifold preimage
       compactly contained in \(W_a^{\rm har}\), with the common
       scale-\(R_a\) lower separation
       \(\zeta_{\rm out}^{\rm phys}>0\).

       The named outer coefficient ball, its inner locus, and the
       common \(+\)-tier throughout the outer ball are furnished by
       Lemma~\ref{lem:reference-carrier-physical-certificate}.  The
       actual closed carrier is required to satisfy
       \[
        G_0\in
        \mathscr B_{\rm coeff,in}^{\rm phys}
        \bigl(
        G_{\rm ref}^{\rm phys};
        \varepsilon_{\rm coeff}^{\rm phys},
        \Lambda_{\rm coeff}^{\rm phys},
        \mu_{\rm coeff}^{\rm phys}
        \bigr).
       \]
       Thus, for every \(x\in V_a^5\), \(G_0\) carries the transferred
       common \(+\)-tier and, by restriction, the operative tier, with
       the corresponding transferred witness pairs and quantitative Dirichlet
       containments.  These triples and their strict componentwise gaps,
       like the two normalized prepared triples in item~(5), are part
       of the common positive margin.  No recentering at \(G_0\) is
       performed: every nearby or restarted carrier in this prepared
       package uses the same \(G_{\rm ref}^{\rm phys}\), outer ball, and
       inner locus.

       The same reference-carrier lemma selects \(\delta_{\rm RF}\)
       after the physical harmonic moduli and proves, for every metric
       \(G_\circ\) in the named outer ball, the following uniform
       half-open estimate.  In the fixed physical reference atlas on
       every \(W_a^{\rm har}\), write its actual fixed-marking Ricci
       flow as
       \[
        G(\,\cdot\,;G_\circ):
        [0,t_*(G_\circ))\longrightarrow
        \operatorname{Met}_{C}^{14,\alpha}(\mathcal X),
        \qquad t_*(G_\circ)\in(0,\infty].
       \]
       Then
       \[
        \begin{aligned}
        &\sup_{\substack{0\leq t<t_*(G_\circ)\\
                         t\leq\delta_{\rm RF}R_a^2}}
         \|G(t;G_\circ)-G_\circ\|_
          {C_{R_a}^{2,\alpha}(W_a^{\rm har})}\\
        &\hspace{8em}\leq
        \frac14\delta_{\rm har}^{\rm phys,wit}.
        \end{aligned}
       \]
       The reference-to-carrier modulus is used only to obtain the common
       \(+\)-tier on the outer ball; the displayed quarter-window uses
       the separate \(+\)-to-operative modulus.  Thus the two moduli
       enter in this order.
       Their smallest members cover both the physical graft collar and
       the complementary outer region, which lie in
       \(E=\mathcal X\setminus\overline{K_2}\); their retained largest
       members cover \(\overline{E^+}\).  They satisfy
       the overlap, scale-comparability, and interface-separation
       certificate \eqref{eq:prepared-exterior-termination}, with common
       \(N_{\rm cov},C_{\rm cov},R_{\rm in},c_{\rm atl},
       \ell_{\rm Leb},\Lambda_{\rm atl}\).  No cyclic cross-cover is
       imposed.  The exterior gauge is the single gauge anchored at its
       inner boundary in
       Lemma~\ref{lem:anchored-exterior-interface}.
       The hypotheses of
       Lemma~\ref{lem:buffered-local-Ricci-control}
       hold on every auxiliary outer member \(V_a^5\) through order twelve.
       Lemma~\ref{lem:anchored-exterior-interface} therefore propagates
       the coefficient bounds to the retained \(U_a^{++}\).  For one
       recorded \(0<\mu_{\rm RF}<1\),
      \[
       2C_\lambda\lambda_0
       \leq(1-\mu_{\rm RF})\delta_{\rm RF}
       \min_{1\leq a\leq N_{\rm ext}}R_a^2.
      \]
      The prepared input order \(k_0+2\geq14\) supplies the required
      order-fourteen coefficient bound and hence these
      order-twelve initial curvature bounds.  Here \(C_\lambda\) is the constant fixed in
       \eqref{eq:fixed-C-lambda}.  Thus the graft-collar and outer
      controls required on every admissible remaining physical-time
      interval are both certified by one-time data.  This is exactly
      the physical part of
      Definition~\ref{def:witnessed-auxiliary-harmonic-package}.
\end{enumerate}
The exact exterior-graft identity is imposed by the fixed preparation
convention.  The nine slice equalities are imposed either by the
centered phase projection in route~(a) or directly in route~(b).
Properness, degree, and radial comparison are
structural conditions of the prepared chart and persist by its scaled
\(C^1\) and radial margins.  Apart from the two normalized
harmonic-coordinate certificates and the finite physical-cover
certificate, the remaining analytic conditions are strict
inequalities in finite jets or weighted prepared norms.
Lemma~\ref{lem:prepared-harmonic-radius-lower-stability} makes the
normalized witnessed condition open and furnishes a new strictly
stronger normalized reserve.  The physical condition is open because
the actual carrier has positive inner-locus slack; throughout that
locus Lemma~\ref{lem:reference-carrier-physical-certificate} supplies
the same common \(+\)- and operative tiers.
The uniform quarter-modulus fixed-marking Ricci-flow estimate is
\eqref{eq:witnessed-physical-quarter-modulus} on the same named
coefficient ball, while the recorded
\(\mu_{\rm RF}\)-slack makes the lifetime-width inequality open.
Thus the entrance set is relatively open in the finite
little-H\"older sliced prepared manifold constructed above; the
equalities themselves are defining constraints, not open inequalities.
Its intersection with the smooth prepared points is dense.  No
conclusion on an ensuing \emph{global} bootstrap interval is part of
the definition; the quarter-modulus physical estimate above is the
precomputed local Ricci--DeTurck consequence of the recorded one-time
coefficient and buffer data.
When the continuation order has been fixed in the surrounding
argument, we suppress the words ``of continuation order \(k_0\).''
\end{definition}

Recall that $\varphi_u$ is the complete flow of $\bar\nabla\bar f$,
normalized by $\varphi_0=\operatorname{Id}$.  The canonical FIK shrinking flow
based at time $-1$ is
\begin{equation}\label{eq:canonical-FIK-flow}
 g_{\mathrm{FIK}}(s)
 =(-s)\varphi_{-\log(-s)}^*\bar g,
 \qquad -\infty<s<0.
\end{equation}
The shrinker equation gives
$\partial_sg_{\mathrm{FIK}}=-2\Ric_{g_{\mathrm{FIK}}}$ and
$g_{\mathrm{FIK}}(-1)=\bar g$.

\begin{theorem}[Prepared-entrance continuation theorem]
\label{thm:prepared-entrance-continuation}
Fix \(k_0\geq12\) and \(0<\alpha<1\), then
$0<\sigma<\theta<\beta$, then one rate-compatible numerical prepared
package at this order, and finally
\(0<\varepsilon\leq\varepsilon_{\rm ent}\).
Every strict prepared FIK entrance of continuation order \(k_0\) in
the sense of Definition~\ref{def:strict-prepared-entrance} at
$\tau_0$ generates a unique coupled mild feedback evolution for all
$\tau\geq\tau_0$.  The closed Ricci--DeTurck representative and the
parabolic prepared components acquire positive-time interior
regularity.  The ODE-carried components and the diffeomorphism used to
return to the ungauged Ricci flow retain their finite entrance order.
Consequently the ungauged closed flow \(G(t)\) is smooth for positive
physical time when the entrance, including its coordinate
diffeomorphism data, is smooth.  For a general finite little-H\"older
entrance, the conclusions below comprise gauge-fixed smoothing and the
coordinate-invariant curvature estimates, but not a gain of coordinate
regularity in the unsmoothed pullback.  The entire prepared
evolution is smooth up to, and away from, the initial face when the
entrance is smooth.  Put
\[
 H=\rho_\tau h,\qquad q=|a|+|b|,\qquad
 P_\infty(\tau)=\int_\tau^\infty q(s)\,ds.
\]
For every $\tau\geq\tau_0$,
\begin{align}
 \|H(\tau)\|_{L^2_\nu}
 &\leq C\varepsilon e^{-\theta\tau}
       +Ce^{-ce^\tau},
 \label{eq:master-bootstrap-L2}\\
 \sum_{\ell=0}^2|\bar\nabla^\ell h|
 &\leq C\varepsilon e^{-\sigma\tau}(1+\bar f)^\sigma
       +Ce^{-ce^\tau}
 &&\text{on }\{\bar f\leq e^\tau\},
 \label{eq:master-bootstrap-inner}\\
 \sup_M\sum_{\ell=0}^2|\bar\nabla^\ell h|
 &\leq C\varepsilon,
 \label{eq:master-bootstrap-global}\\
 P_\infty(\tau)
 &\leq C\varepsilon^2e^{-2\theta\tau}
       +Ce^{-ce^\tau}.
 \label{eq:master-bootstrap-phase-tail}
\end{align}
For $\tau\geq\tau_0+1$ one has the sharper instantaneous estimates
\begin{align}
 \|H(\tau)\|_{H^1_\nu}
 &\leq C\varepsilon e^{-\theta\tau}
       +Ce^{-ce^\tau},
 \label{eq:master-bootstrap-H1}\\
 q(\tau)
 &\leq C\varepsilon^2e^{-2\theta\tau}
       +Ce^{-ce^\tau}.
 \label{eq:master-bootstrap-velocity}
\end{align}
On every $K\Subset M$ and for every $m\geq0$,
\begin{equation}\label{eq:master-core-rate}
 \|h(\tau)\|_{C^m(K)}
 \leq C_{K,m}\varepsilon e^{-\theta\tau}
      +C_{K,m}e^{-ce^\tau}.
\end{equation}
This holds whenever $\tau\geq\tau_0+1$ and
$K\subset\{\bar f<e^{\tau-1}\}$.

There are limits $\lambda_\infty>0$ and
$\Psi_\infty\in\operatorname{Diff}(M)$ such that
\begin{align}
 \left|\log\frac{\lambda(\tau)e^\tau}{\lambda_\infty}\right|
 &\leq C\varepsilon^2e^{-2\theta\tau}+Ce^{-ce^\tau},
 \label{eq:master-scale-rate}\\
 \left|\frac{\lambda(\tau)}{T-t(\tau)}-1\right|
 &\leq C\varepsilon^2e^{-2\theta\tau}+Ce^{-ce^\tau},
 \label{eq:master-physical-scale-rate}\\
 \|\Psi_\infty-\Psi_\tau\|_{C^m(K)}
 &\leq C_{K,m}\varepsilon^2e^{-2\theta\tau}
       +C_{K,m}e^{-ce^\tau}.
 \label{eq:master-phase-rate}
\end{align}
For each \(m\) for which the entrance carries the
\(C^{m+1}\) map bounds and \(C^m\) target bounds in
Proposition~\ref{prop:adaptive-target-tracking}, the relative radial
maps \(R_\tau\), their inverses, and the adaptive targets \(S_\tau\)
converge in the corresponding norms on \(\Omega_\eta^+\), with tails
\(O(e^{-2\theta\tau})\).

For a finite-order entrance, the marked embeddings below are asserted
only at the finite map order carried by the prepared package.  The
\(C^\infty_{\mathrm{loc}}\) conclusions concern the pulled-back metric
tensors, not the regularity of the absolute marking maps.  Their
all-order regularity is established by the smooth-approximation and
positive-time diagonal argument in the proof below.  After the theorem,
Lemma~\ref{lem:finite-order-relative-marking-cancellation} records the
exact relative-marking identity which explains why this conclusion
does not assert a derivative gain for the absolute markings.

The underlying closed Ricci flow has maximal time
\begin{equation}\label{eq:master-singular-time}
 T=t(\tau_0)+\int_{\tau_0}^{\infty}\lambda(\tau)\,d\tau<\infty
\end{equation}
and develops a global Type-I singularity:
\begin{equation}\label{eq:master-global-Type-I}
0<c\leq
 (T-t)\|\Rm_{G(t)}\|_{L^\infty(\mathcal X,G(t))}
 \leq C
\end{equation}
for all sufficiently late $t<T$.  There are exhausting domains
$\mathcal U_t\subset M$ and marked embeddings
$\Xi_t:\mathcal U_t\to\mathcal X$ for which
\begin{equation}\label{eq:master-marked-convergence}
 \lambda(t)^{-1}\Xi_t^*G(t)\longrightarrow\bar g
 \quad\text{in }C^\infty_{\mathrm{loc}}(M),
\end{equation}
and
\begin{equation}\label{eq:master-physical-power}
 \lambda(t)=(T-t)\bigl(1+O((T-t)^{2\theta})\bigr).
\end{equation}
More strongly, let $t_i\uparrow T$ be arbitrary, put
\[
 \delta_i=T-t_i,\qquad \tau_i=\tau(t_i),\qquad
 \Xi_i=\Xi_{t_i},
\]
and freeze $\Xi_i$ for the entire $i$th rescaled flow.  For every
$K\Subset M$, $I\Subset(-\infty,0)$, and $m\geq0$,
\begin{equation}\label{eq:master-marked-parabolic-convergence}
 \left\|
  \delta_i^{-1}\Xi_i^*G(T+s\delta_i)
  -g_{\mathrm{FIK}}(s)
 \right\|_{C^m(K\times I)}
 \leq C_{K,I,m}\delta_i^\theta
\end{equation}
for all sufficiently large $i$.  Thus the full sequence, without
subsequence extraction, has the canonical FIK flow as its marked
parabolic limit.

All constants in the bootstrap and continuation estimates are
uniform on prepared-coordinate neighborhoods with
common strict margins and uniform entrance and scale bounds.
\end{theorem}

\begin{proof}
We first assume that the strict entrance is smooth.  Fix a finite test
endpoint $\tau_1>\tau_0$, without allowing any
constant to depend on it.  Let \([\tau_0,\tau_*)\) be the maximal
interval on which the following faces hold simultaneously with doubled
bootstrap constants:
the analytic box, with its global raw scale-normalized \(C^2\) face at
\(2\delta_{\rm c2}\), the separate radius-independent activation face
\begin{equation}\label{eq:master-pre-Holder-first-exit-face}
 \sup_{\tau_0\leq\tau<\tau_*}
 \|h(\tau)\|_{\mathfrak C_{{\rm pre},0}^{2,\alpha}}
 \leq2\delta_{\rm c2},
\end{equation}
the quantitative support-separation face
\begin{equation}\label{eq:master-support-separation-first-exit-face}
 \inf_{\tau_0\leq\tau<\tau_*}
 \mathfrak s_{\rm sep}(\tau)
 \geq\tfrac12\kappa_{\rm sep},
\end{equation}
the quantitative forward-and-inverse map-margin face
\begin{equation}\label{eq:master-map-margin-first-exit-face}
 \inf_{\tau_0\leq\tau<\tau_*}
 \min\{\mathfrak m_R(R_\tau),
 \mathfrak m_F(F_\tau;\acute G_\tau,S_\tau)\}
 \geq\tfrac12\kappa_{\rm map},
\end{equation}
and the ellipticity box; accumulated phase and scale
comparability; controlled HMHF diffeomorphism and annulus tracking;
the coarse pointwise modulation and endpoint-independent target-defect
and coefficient time-modulus package; the endpoint-independent
order-six map, inverse-map, and target spatial-coefficient package;
effective-column estimates and Gram invertibility; graft
  compatibility, the order-fourteen coefficient-transfer package, and
  the order-twelve closed-curvature package; and the support
and $C^2$ bounds for the pure graft forcing.  The one-time strict
entrance conditions and
Proposition~\ref{prop:coupled-local-feedback} make this interval
nonempty.
The fixed choices give, on this interval,
\begin{equation}\label{eq:master-pre-radius-activation-ledger}
 \begin{aligned}
 \tau_0\geq\tau_{\rm ad}
 &\geq\max\{\tau_{\rm pre},\tau_{\rm base}^{(2)}\},\\
 \int_{\tau_0}^{\tau_*}(|a|+|b|)\,d\tau
 &\leq\varepsilon_{\rm ph}
 \leq\min\{\varepsilon_{\rm pre},
            \varepsilon_{\rm ph}^{\rm HM},
            \varepsilon_{\rm ph,*}^{(2)}\},\\
 2\delta_{\rm c2}
 &\leq\min\{\delta_{\rm atl}^{\rm pre},\delta_{\rm pre},
             C_{\rm raw}^{-1}\varepsilon_{\rm hm}^{\rm HM},
             \delta_{\rm box}^{(2)}\}.
 \end{aligned}
\end{equation}
The entrance conditions supply the
\(\operatorname{Map}_{\rm sc}^{14,\alpha}\) bounds for
\(R_{\tau_0}^{\pm1}\), the
\(\operatorname{Map}_{\rm sc}^{6,\alpha}\) bounds for
\(F_0^{\pm1}\), and their strict radial and local-invertibility
margins; explicitly, these norms are below
\(\Lambda_{R,14}^{\rm pre}-2\mu_R^{\rm pre}\) and
\(\Lambda_{F,6}^{\rm pre}-2\mu_F^{\rm pre}\), respectively.
The two quantitative lower margins are precisely
\eqref{eq:strict-entrance-R-singular-face} and
\eqref{eq:strict-entrance-map-singular-face}; by continuity they also
make \eqref{eq:master-map-margin-first-exit-face} strict at the entrance.
 The witnessed normalized and physical harmonic-radius conditions and
 the ordinary order-twelve coefficient face are present.  So are the source--target discrepancy and
proper-diffeomorphism faces, and the order-six map-and-inverse face at
\(2C_{\rm map,6}^{\rm pre}\).  Thus items (P1)--(P4) of
\(\mathscr P_{\rm pre}^{(6)}\), and every activation hypothesis of
Lemma~\ref{lem:pre-radius-low-order-closure}, are present before that
 lemma is invoked.
The undoubled order-six map, effective-column, and feedback faces are
those furnished by Lemma~\ref{lem:pre-radius-low-order-closure}; their
doubled versions define the corresponding first-exit faces.  We close
these faces in the following order.  The physical auxiliary closure
uses only the stated low-order conclusions of
Lemma~\ref{lem:pre-radius-low-order-closure}.  The fixed-\(\Gamma\)
high-regularity constants are chosen afterward and do not affect the
radius.

Use first only the coarse phase, scale, \(C^2\), graft, and physical
geometry faces.  Proposition~\ref{prop:graft-compatibility} and the
physical part of Corollary~\ref{cor:adaptive-auxiliary-closure} improve
the closed-flow, graft, curvature, harmonic-radius, and spatial
coefficient package without using a transported support conclusion.
The unconditional physical support is
\eqref{eq:pure-graft-physical-support}.

Before constructing a persistent paired atlas or invoking the
order-six map bridge, close the two quantitative map margins.  The
relative-marking argument
\eqref{eq:pre-radius-R-singular-stability}--%
\eqref{eq:pre-radius-R-margin-closure} uses only the strict entrance
margin, the right-translated \(R\)-equation, and the accumulated phase
budget.  The metric argument
\eqref{eq:pre-radius-F-metric-discrepancy}--%
\eqref{eq:pre-radius-F-margin-closure} uses only the exact prepared
graph identity, the activated pre-radius \(C^{2,\alpha}\) bound, and
the primitive ellipticity box.  Neither argument uses the persistent
atlas or the bridge.  They therefore apply on the present first-exit
interval and give
\begin{equation}\label{eq:master-map-margin-improvement}
 \inf_{\tau_0\leq\tau<\tau_*}
 \min\{\mathfrak m_R(R_\tau),
 \mathfrak m_F(F_\tau;\acute G_\tau,S_\tau)\}
 \geq\tfrac34\kappa_{\rm map}.
\end{equation}
Thus \eqref{eq:master-map-margin-first-exit-face} is strictly improved
before any result depending on its lower singular-value bounds is used.

For each \(S<\tau_*\), the remaining harmonic-map defects have a finite
constant \(\Lambda_S\), so
Theorem~\ref{thm:adaptive-HMHF-continuation} applies on
\([\tau_0,S]\).  By the overlap clause in its proof, the gauges
obtained for different test endpoints \(S\) are restrictions of the
same uniquely determined finite-horizon solution.  Its displacement
and tracking proof uses only the
coarse \(C^2\) box and the \(L^1\) phase budget; letting
\(S\uparrow\tau_*\) therefore improves
\eqref{eq:F-outer-drift}--\eqref{eq:Phi-annulus-tracking} with constants
independent of the test endpoint.  The scale-normalized target package
also supplies \eqref{eq:pure-graft-cutoff-jet-hypothesis} by
Remark~\ref{rem:adaptive-pure-graft-cutoff-jets}.
Proposition~\ref{prop:pure-graft-sharp}
now supplies the transported support, but that support was not used to
construct the preceding physical coefficient package.

Lemma~\ref{lem:coarse-effective-Gram} next improves the low-order Gram
face and gives an endpoint-independent coarse pointwise bound for
\(q\).  Hence the target defect and the source and target coefficient
time moduli are endpoint-independent.
With the lower singular-value inputs already closed in
\eqref{eq:master-map-margin-improvement},
Lemma~\ref{lem:finite-HMHF-C6-bridge} then improves the order-six map
and inverse-map face, while
Proposition~\ref{prop:adaptive-target-tracking} supplies the target
spatial part.  Together with the doubled source coefficient face and
the prepared graph formula, these bounds give the
\(\Gamma\)-independent estimate
\eqref{eq:pre-atlas-C5-interpolation-ceiling}.  Every hypothesis of
Proposition~\ref{prop:effective-column-tails} is now available, so the
full effective-column and Gram faces improve.
Corollary~\ref{cor:adaptive-feedback-package} may now be applied with endpoint
$\tau_*-\delta$; letting $\delta\downarrow0$ gives
\[
 \|H(\tau)\|_{L^2_\nu}^2
 +\int_\tau^{\tau_*}\|H(s)\|_{H^1_\nu}^2\,ds
 +\int_\tau^{\tau_*}q(s)\,ds
 \leq C\varepsilon^2e^{-2\theta\tau}+Ce^{-ce^\tau}.
\]
Thus the future phase tail is available before any sharp decaying
pointwise velocity estimate is used.
Theorem~\ref{thm:robust-modulated-three-region} then strictly improves
the three-region and global raw \(C^2\) faces; explicitly, the global
through-order-two norm is at most
\[
 C_{\rm 3reg}^{(2)}\varepsilon
 \leq C_{\rm 3reg}^{(2)}\varepsilon_{\rm ent}
 <\delta_{\rm c2}
\]
by \eqref{eq:authoritative-pre-Holder-interpolation-choice}.
On each scale-one
pre-radius chart apply the interpolation inequality to
\(u_L=L^{-1}h\):
\[
 \|u_L\|_{C^{2,\alpha}}
 \leq
 C_{\rm int}\|u_L\|_{C^0}^{1-\vartheta_{\rm int}}
 \|u_L\|_{C^5}^{\vartheta_{\rm int}}.
\]
Taking the compact-core maximum and dyadic supremum, using the
zeroth-order three-region output and
\eqref{eq:pre-atlas-C5-interpolation-ceiling}, gives
\[
 \|h\|_{\mathfrak C_{{\rm pre},0}^{2,\alpha}}
 \leq
 C_{\rm int}
 (C_{\rm 3reg}^{(0)}\varepsilon)^{1-\vartheta_{\rm int}}
 (K_{h,5}^{\rm pre})^{\vartheta_{\rm int}}
 <\delta_{\rm c2}
\]
by \eqref{eq:authoritative-pre-Holder-interpolation-choice}.  Hence the
pre-atlas H\"older face
\eqref{eq:master-pre-Holder-first-exit-face} is also strictly improved;
it cannot be the first exit.  With these improvements established, the
exact Gram system is combined with the recovered \(H^1_\nu\) control to improve the
sharp instantaneous velocity, accumulated phase, and scale faces; the
earlier coarse Gram bound served only to close coefficient regularity.

With the endpoint-independent defect package now available,
Proposition~\ref{prop:moving-target-HMHF-restart}, applied with the
fixed levels \(\epsilon_0^{\rm HM}<\epsilon_1^{\rm HM}\), excludes loss
of the map chart with one uniform restart constant.  The annulus and
transported-support margins were already improved by
Theorem~\ref{thm:adaptive-HMHF-continuation} and
Proposition~\ref{prop:pure-graft-sharp} above.  Thus the dependency
order is
\[
 \begin{aligned}
 \text{activation}
 &\longrightarrow \text{physical auxiliary closure}\\
 &\longrightarrow \text{quantitative \(R/F\) map-margin closure}\\
 &\longrightarrow \text{finite-horizon HMHF tracking}\\
 &\longrightarrow \text{coarse Gram and pointwise }q\\
 &\longrightarrow \text{uniform target-time/defect package}\\
 &\longrightarrow \text{uniform \(C^6\) map bridge}\\
 &\longrightarrow \text{full columns and feedback},
 \end{aligned}
\]
which closes the continuation in the displayed order.

If $\tau_*<\tau_1$, the only remaining exit is termination of the
closed flow.  The normalized \(C^2\) estimate controls the marked inner
region.  On every member of the fixed cover of the physical graft
collar and the outer complement,
Lemma~\ref{lem:buffered-local-Ricci-control} applies because
\[
 t(\tau_*)-t(\tau_0)
 \leq C_\lambda\lambda_0
 <\tfrac12\delta_{\rm RF}
 \min_{1\leq a\leq N_{\rm ext}}R_a^2.
\]
Together with the marked inner region, these sets cover \(\mathcal X\),
 so curvature remains globally
 bounded and the closed Ricci-flow extension theorem continues the
 flow.  Every ordinary face therefore retains a strict terminal margin,
 while the witnessed harmonic-radius conditions retain the common
 operative and reserve packages established above.  Local existence
 continues the evolution, a contradiction.  Hence
$\tau_*\geq\tau_1$.  Letting $\tau_1\to\infty$ gives the global
normalized evolution and the infinite future tail.

The $L^2_\nu$ and pointwise estimates just obtained give
\eqref{eq:master-bootstrap-L2}--%
\eqref{eq:master-bootstrap-global}.  The pure-graft estimate, the
effective-column tails, and the $C^2$ cutoff calculation give
\eqref{eq:H1-recovery-tails}.
Proposition~\ref{prop:instantaneous-H1-spectral-rate} supplies that improvement
and the sharp quadratic velocity estimate;
Lemma~\ref{lem:future-phase-tail} supplies
\eqref{eq:master-bootstrap-phase-tail}.
The admissible first-exit interval carries
\eqref{eq:coarse-Gram-support-separation}.  For a fixed
\(K\Subset M\), put
\[
 R_K:=\max\left\{e^{\tau_0},\sup_K\bar f\right\}
\]
and choose one fixed smooth buffer \(K\Subset K^+\Subset M\) with
\(\sup_{K^+}\bar f<e^{1/4}R_K\).  Whenever
\(\tau\geq\tau_0+1\) and
\(K\subset\{\bar f<e^{\tau-1}\}\), one has
\[
 \sup_{K^+}\bar f<e^{\tau-3/4}.
\]
Thus the same fixed buffer works throughout the full range in the
theorem, and
Corollary~\ref{cor:quantitative-core-smoothing} gives
\eqref{eq:master-core-rate}.  The finite-regularity passage at the end of
the proof uses the support-exited moving-core cylinders in
Lemma~\ref{lem:finite-regularity-prepared-approximation}; no uncarried
jet of an ODE component is invoked.  Theorem~\ref{thm:phase} gives
\eqref{eq:master-scale-rate}--\eqref{eq:master-phase-rate}, while
Proposition~\ref{prop:adaptive-target-tracking} gives convergence of
$R_\tau$, $R_\tau^{-1}$, and $S_\tau$.  Its
$e^{-\tau}$ geometric error is faster than $e^{-2\theta\tau}$ because
$\theta<\beta<1/2$.

The scale bounds make \eqref{eq:master-singular-time} finite.  Cover
the closed manifold at every $t<T$ by the marked inner region, the
fixed physical graft collar, and the complementary outer region.
Their curvature bounds give $|\Rm_G|\leq C\lambda^{-1}$.  Conversely,
\eqref{eq:master-core-rate} and nonflatness of $\bar g$ give
$|\Rm_G|\geq c\lambda^{-1}$ at a marked point.  Hence the maximal
time is exactly $T$, and
\eqref{eq:master-physical-scale-rate} gives
\eqref{eq:master-global-Type-I} and
\eqref{eq:master-physical-power}.

 Lemma~\ref{lem:marked-chart-exhaustion} supplies the exhausting
 normalized-time domains and embeddings.  In physical time set
 \(\mathcal U_t:=\mathcal U_{\tau(t)}\) and
 \(\Xi_t:=\Xi_{\tau(t)}\); they obey the exact identity
 $\lambda^{-1}\Xi_t^*G-\bar g=h$.  The core estimates and local
 parabolic smoothing prove \eqref{eq:master-marked-convergence}.

 We finally freeze the marking and recover the spacetime limit.  For
 $s\in I\Subset(-\infty,0)$ define $\tau_i(s)$ by
 \[
  t(\tau_i(s))=T+s\delta_i.
 \]
 This is well defined for all large $i$, uniformly on $I$.  From
 \eqref{eq:master-physical-scale-rate},
 \[
  \frac{d}{d\tau}\log(T-t(\tau))
  =-\frac{\lambda(\tau)}{T-t(\tau)}
  =-1+O(e^{-2\theta\tau}).
 \]
 Since $t_i=T-\delta_i$, integration over bounded normalized-time
 windows gives, uniformly for $s\in I$,
 \begin{align}
  \tau_i(s)-\tau_i
  &=-\log(-s)+O_I(e^{-2\theta\tau_i}),
  \label{eq:master-normalized-time-shift}\\
  \frac{\lambda(\tau_i(s))}{\delta_i}
  &=(-s)\bigl(1+O_I(e^{-2\theta\tau_i})\bigr).
  \label{eq:master-parabolic-scale-ratio}
 \end{align}

 Put
 \begin{equation}\label{eq:master-marking-transition}
  D_{i,s}=\Phi_{\tau_i(s)}\circ\Phi_{\tau_i}^{-1}.
 \end{equation}
 The domains in Lemma~\ref{lem:marked-chart-exhaustion} have the form
 $\mathcal U_\tau=\Phi_\tau(M_{\rm in})$.  Hence $D_{i,s}$ maps
 $\mathcal U_{\tau_i}$ to $\mathcal U_{\tau_i(s)}$, and the composition
 order in
 \[
  \Xi_{\tau_i}
  =\Xi_{\tau_i(s)}\circ D_{i,s}
 \]
 is exact.  Therefore
 \begin{equation}\label{eq:master-frozen-pullback}
 \delta_i^{-1}\Xi_{\tau_i}^*G(T+s\delta_i)
 =
 \frac{\lambda(\tau_i(s))}{\delta_i}
 D_{i,s}^*\bigl(\bar g+h(\tau_i(s))\bigr).
\end{equation}
For every fixed \(K\Subset M\) and \(I\Subset(-\infty,0)\), enlarge
\(K\) once so that all trajectories \(D_{i,s}(K)\), \(s\in I\), lie
in that enlargement for all large \(i\).  Equation
\eqref{eq:master-frozen-pullback} then identifies the entire
frozen-window family with sections of the single tensor bundle
\(S^2T^*K\); all \(C^m(K\times I)\) norms and time derivatives below
are taken in this fixed bundle, using one fixed background connection
on the enlarged compact set.  Thus no subtraction of tensors based at
different markings or tangent spaces is implicit.

 In normalized time the Eulerian generator of $\Phi_\tau$ is
 \begin{equation}\label{eq:master-Phi-generator}
  (\partial_\tau\Phi_\tau)\circ\Phi_\tau^{-1}
  =
  -B_{\bar g}(g(\tau))
  +(1+a(\tau))\bar\nabla\bar f-U(\tau).
 \end{equation}
 On every fixed compact set, the core estimate, the instantaneous
 velocity estimate, and local smoothing imply in every $C^m$ norm
 \[
  (\partial_\tau\Phi_\tau)\circ\Phi_\tau^{-1}
  =\bar\nabla\bar f+O_{K,m}(e^{-\theta\tau}).
 \]
 The generator has the uniform at-most-linear growth bound established
 for the controlled chart.  Because
 \eqref{eq:master-normalized-time-shift} gives a bounded time window,
 trajectories issuing from $K$ remain in one fixed compact set.  ODE
 stability, together with
 \eqref{eq:master-normalized-time-shift}, now yields
 \begin{equation}\label{eq:master-transition-limit}
  D_{i,s}
  =\varphi_{-\log(-s)}
   +O_{K,I,m}(e^{-\theta\tau_i})
 \end{equation}
 in $C^m(K\times I)$.  Hence
 \eqref{eq:master-core-rate} gives the same rate for
 $h(\tau_i(s))$ there.  Substituting these estimates and
 \eqref{eq:master-parabolic-scale-ratio} into
 \eqref{eq:master-frozen-pullback}, and using
 $\delta_i\asymp e^{-\tau_i}$, proves
 \eqref{eq:master-marked-parabolic-convergence} for spatial
 derivatives.  The tensors on the left solve Ricci flow in the fixed
 marking $\Xi_i$:
 \[
  \partial_s\bigl(\delta_i^{-1}\Xi_i^*
       G(T+s\delta_i)\bigr)
  =-2\Ric_{\delta_i^{-1}\Xi_i^*G(T+s\delta_i)}.
 \]
 Uniform spatial curvature-derivative bounds on $K\times I$ in the
 fixed bundle just specified, together with this
 equation inductively give the asserted mixed spacetime derivatives.
 Every step above is uniform on a common-margin prepared family with
 uniform entrance and scale bounds.

 We now remove the auxiliary smoothness assumption.  Let the entrance
 be a finite little-H\"older point of
 \(\Sigma_{\tau_0}^{k_0+2,\alpha}\), and choose the smooth sliced
 approximants furnished by
 Lemma~\ref{lem:finite-regularity-prepared-approximation}.  The proof
 just completed gives, on every smooth approximant and every finite
 endpoint \(S\), exactly the package-local estimates listed in that
 lemma, with constants independent of the approximant and of \(S\).
 The finite-horizon induction in the lemma therefore continues all
 approximants through \([\tau_0,S]\), passes the slice, scale, phase,
 \(L^2_\nu\), dissipation, three-region, graft, and Type-I estimates to
 the unique finite-regularity mild solution.  The gauge-fixed estimate
 \eqref{eq:finite-regularity-positive-time-diagonal}, together with the
 normalized-equation support-exit clause of that lemma on the recent
 core cylinders, yields the local \(C^m\) estimates; they are not
 obtained by pretending that \(C^{k_0,\alpha}\) convergence controls
 arbitrary derivatives.  Diagonalizing over finite \(S\) gives the global mild
 evolution.  The same positive-time diagonal argument on the frozen
 compact spacetime windows used above proves
 \eqref{eq:master-marked-convergence} and
 \eqref{eq:master-marked-parabolic-convergence} in every fixed order.
 The qualification concerning the unsmoothed ungauged coordinate
 pullback is precisely the one stated in the theorem.

\end{proof}

\begin{lemma}[Finite-order markings and smooth relative pullbacks]
\label{lem:finite-order-relative-marking-cancellation}
Let a finite-order strict prepared entrance generate the evolution in
Theorem~\ref{thm:prepared-entrance-continuation}.  The absolute
markings \(\Xi_t\) and the ODE-carried prepared maps are asserted only
at their carried finite orders.  Nevertheless, on every fixed compact
set and for every sufficiently late time,
\[
 \lambda(t)^{-1}\Xi_t^*G(t)=\bar g+h(\tau(t))
\]
has a canonical smooth representative.  Moreover, if
\(t_i\uparrow T\), \(\delta_i=T-t_i\), and
\(t(\tau_i(s))=T+s\delta_i\), then on every
\(K\Subset M\) and \(I\Subset(-\infty,0)\),
\[
 \delta_i^{-1}\Xi_{\tau_i}^*G(T+s\delta_i)
 =
 \frac{\lambda(\tau_i(s))}{\delta_i}
 D_{i,s}^*\bigl(\bar g+h(\tau_i(s))\bigr),
 \qquad
 D_{i,s}=\Phi_{\tau_i(s)}\circ\Phi_{\tau_i}^{-1},
\]
and the right-hand side is smooth to every fixed order for all
sufficiently large \(i\).  Thus every
\(C^m(K\times I)\) estimate for the displayed pulled-back tensors is
meaningful for arbitrary fixed \(m\), without asserting any derivative
gain for the absolute markings or the ODE-carried data.
\end{lemma}

\begin{proof}
The first identity is
\eqref{eq:marked-pullback-identity}, and the frozen identity is
\eqref{eq:master-frozen-pullback}.  On every late compact parabolic
cylinder, the positive-time diagonal estimate
\eqref{eq:finite-regularity-positive-time-diagonal} and the
support-exit clause of
Lemma~\ref{lem:finite-regularity-prepared-approximation} give
all-order regularity of \(h\).  The Eulerian generator
\eqref{eq:master-Phi-generator} is therefore smooth to every fixed
order there, so ordinary ODE regularity gives the same conclusion for
the relative transition \(D_{i,s}\).  The fixed-marking Ricci-flow
identity at the end of the continuation proof supplies the mixed time
derivatives.
\end{proof}

\begin{remark}[Marked scope of the parabolic limit]
\label{rem:marked-parabolic-scope}
The embeddings in
\eqref{eq:master-marked-parabolic-convergence} are frozen from the
receding-slice marking at the base times $t_i$.  The conclusion is
quantitative uniqueness of the marked parabolic limit.  It does not
compare blow-ups based at unrelated points or formed in unrelated
charts, and it makes no unmarked tangent-flow uniqueness assertion.
\end{remark}

\begin{lemma}[One-state continuity of terminal geometric data]
\label{lem:one-state-terminal-continuity}
Fix \(k_{\rm ent}\geq12\) and a common-margin open set
\(\mathscr B_{\rm ent}\subset
\Sigma_{\tau_0}^{k_{\rm ent}+2,\alpha}\) of finite-regularity strict
prepared entrances at one entrance time.  The scalar maps
\[
 \mathbf z_0\longmapsto T(\mathbf z_0),
 \qquad
 \mathbf z_0\longmapsto\lambda_\infty(\mathbf z_0)
\]
are continuous.  For every \(K\Subset M\) and every fixed carried
order \(\ell\geq0\), the map
\[
 \mathbf z_0\longmapsto\Psi_\infty(\mathbf z_0)|_K
\]
is continuous into the fixed \(C^\ell(K;M)\) exponential chart.

Fix now an integer \(0\leq m\leq k_{\rm ent}+1\) for which the
initial \(R^{\pm1}\)-coordinate projections are continuous into
\(C^{m+1}\), and for which the prepared set carries uniform
\(C^{m+1}\) bounds for \(R_\tau\) and \(R_\tau^{-1}\), and \(C^m\)
bounds for \(S_\tau\), in the sense of
Proposition~\ref{prop:adaptive-target-tracking}.  Then
\[
 \mathbf z_0\longmapsto
 \left(
  R_\infty|_{\Omega_\eta^+},
  R_\infty^{-1}|_{\Omega_\eta^+},
  S_\infty|_{\Omega_\eta^+}
 \right)
\]
is continuous into
\[
 C^m(\Omega_\eta^+;M)
 \times C^m(\Omega_\eta^+;M)
 \times C^m(\Omega_\eta^+;S^2T^*M).
\]
The two map topologies use the fixed exponential and
parallel-transport identifications of
Proposition~\ref{prop:adaptive-target-tracking}; the tensor topology
uses the fixed background atlas on \(\Omega_\eta^+\).
These conclusions use only finite-time continuous dependence and the
uniform one-state tails in
Theorem~\ref{thm:prepared-entrance-continuation}.
\end{lemma}

\begin{proof}
Continuity is local, so fix an arbitrary entrance and shrink to a
bounded uniformly interior common-margin prepared-coordinate ball
\[
 \mathscr B_0\Subset_{\rm u}\mathscr B_{\rm ent}.
\]
Thus \(\mathscr B_0\) has positive independent-model Banach-radius
distance from the exit
faces and carries one uniform numerical package; no topological
relative compactness is used.  Fix a terminal accuracy
\(\varepsilon_*>0\).  Uniformly on \(\mathscr B_0\), choose \(S\) so
large that
\[
\begin{aligned}
 |T-t(S)|
 &+\left|\log\frac{\lambda(S)e^S}{\lambda_\infty}\right|
 +\|\Psi_\infty-\Psi_S\|_{C^\ell(K)}
 +\|R_\infty-R_S\|_{C^m(\Omega_\eta^+)}\\
 &+\|R_\infty^{-1}-R_S^{-1}\|_{C^m(\Omega_\eta^+)}
 +\|S_\infty-S(S)\|_{C^m(\Omega_\eta^+)}
 <\frac{\varepsilon_*}{3}.
\end{aligned}
\]
This is possible by
\eqref{eq:master-scale-rate}--\eqref{eq:master-phase-rate}, the bound
\(|T-t(S)|\leq Ce^{-S}\), and the adaptive tails following
\eqref{eq:master-phase-rate}; all constants are uniform by the last
paragraph of Theorem~\ref{thm:prepared-entrance-continuation}.
After shrinking \(\mathscr B_0\) once more, the global continuation
theorem keeps every resulting trajectory in one admissible ball on
\([\tau_0,S]\).  Apply
Proposition~\ref{prop:two-state-prepared-evolution} only at the fixed
base order \(12\).  The continuous inclusion from the
\(k_{\rm ent}+2\) entrance chart gives
\begin{equation}\label{eq:terminal-data-finite-horizon-continuity}
 \begin{split}
 &\sup_{\tau_0\leq\tau\leq S}
 \left(
  \mathfrak D_{12}^{\rm hyb,full}(\tau)
  +\|H_1(\tau)-H_2(\tau)\|_{L^2_\nu}
 \right)\\
 &\qquad+
 \left(\int_{\tau_0}^{S}
 \|H_1-H_2\|_{H^1_\nu}^2\,d\tau\right)^{1/2}
 \leq C_S
 \|\mathbf z_{1,0}-\mathbf z_{2,0}\|_
      {\mathscr X_{\rm prep}^{k_{\rm ent}+2,\alpha}} .
 \end{split}
\end{equation}
Combining \eqref{eq:two-state-feedback} with Cauchy--Schwarz on this
fixed interval gives
\begin{equation}\label{eq:terminal-data-feedback-L1}
 \int_{\tau_0}^{S}|c_1-c_2|\,d\tau
 \leq C_S
 \|\mathbf z_{1,0}-\mathbf z_{2,0}\|_
      {\mathscr X_{\rm prep}^{k_{\rm ent}+2,\alpha}},
 \qquad c_i=(a_i,b_i).
\end{equation}
In particular \(t(S)\), \(\lambda(S)\), and the feedback coefficients
are continuous from the stated entrance topology.  The exact phase ODE
\[
 \partial_\tau\Psi_\tau
 =\left(\sum_{j=1}^8b_j(\tau)\chi_\tau W_j\right)\circ\Psi_\tau
\]
then gives continuity of \(\Psi_S\) in \(C^\ell(K;M)\) for every fixed
\(\ell\): on the compact set swept out up to time \(S\), the
generators are smooth, and the ordinary differentiated
parameter-dependent ODE estimates use only the already controlled
scalar coefficients.

It remains to recover the top carried derivative of the adaptive
target.  This derivative is obtained from the exact relative-marking
ODE, not from a higher-order invocation of the hybrid parabolic
estimate.  Put
\[
 Z_{i,\tau}
 :=
 (\varphi_{-\tau})_*
 \bigl(a_i(\tau)\bar\nabla\bar f-U_{b_i(\tau)}\bigr),
 \qquad
 \partial_\tau R_i=Z_{i,\tau}\circ R_i ,
\]
as in \eqref{eq:relative-target-flow}.  On the fixed enlargement of
\(\Omega_\eta^+\) swept out by the two finite-horizon trajectories,
the explicit smooth generators and their conjugated cutoffs satisfy
\[
 \|Z_{1,\tau}-Z_{2,\tau}\|_{C^{m+2}}
 \leq C_{S,m}|c_1(\tau)-c_2(\tau)|,
 \qquad
 \|Z_{i,\tau}\|_{C^{m+2}}
 \leq C_{S,m}|c_i(\tau)|.
\]
These are estimates for fixed background fields; they require no
higher prepared metric norm.  The differentiated flow equations,
Gronwall, and \eqref{eq:terminal-data-feedback-L1} give
\begin{equation}\label{eq:terminal-top-relative-marking-bridge}
 \begin{split}
 d_{C^{m+1}}\bigl(R_1(S),R_2(S)\bigr)
 \leq C_{S,m}\bigg[
 &d_{C^{m+1}}\bigl(R_1(\tau_0),R_2(\tau_0)\bigr)\\
 &+\int_{\tau_0}^{S}|c_1-c_2|\,d\tau
 \bigg]\longrightarrow0 .
 \end{split}
\end{equation}
The restriction \(m\leq k_{\rm ent}+1\) makes the initial
\(C^{m+1}\) term continuous in the entrance topology induced by
\(d_{\rm prep}^{k_{\rm ent}+2,\alpha}\).
The common proper-diffeomorphism margin and
\[
 D(R_S^{-1})[\dot R_S]
 =-DR_S^{-1}\bigl(\dot R_S\circ R_S^{-1}\bigr)
\]
give the corresponding \(C^m\) continuity of \(R_S^{-1}\).
Finally the exact target formula
\[
 S(S)=\lambda(S)(\varphi_S\circ R_S)^*\bar g
\]
and the fixed composition charts give
\begin{equation}\label{eq:terminal-top-target-metric-bridge}
 \|S_1(S)-S_2(S)\|_{C^m(\Omega_\eta^+)}
 \leq C_{S,m}\left(
 \left|\log\frac{\lambda_1(S)}{\lambda_2(S)}\right|
 +d_{C^{m+1}}\bigl(R_1(S),R_2(S)\bigr)\right).
\end{equation}
This explicitly closes the borderline pair
\(k_{\rm ent}=12\), \(m=13\), without asking for an unavailable
\(\mathscr P_{\tau_0}^{15,\alpha}\) entrance bound.  The
finite-horizon parabolic derivative loss is paid once at base order;
the missing carried derivative is supplied by the triangular smooth
ODE.

Nearby entrances now make every corresponding finite-time value
differ by less than \(\varepsilon_*/3\).  The triangle inequality,
with the two uniform terminal tails, proves all asserted continuity
statements.  No global two-state decay estimate is used.
\end{proof}

\begin{lemma}[Static compatibility of the formation and scattering
packages]\label{lem:static-formation-scattering-compatibility}
Fix \(0<\theta_*<\beta\) and finite positive ceilings
\(\Lambda_{\rm ell}^{(2)}\) and \(K_{\mathcal Y,0}^{(2)}\) strictly
dominating the corresponding common faces of the already frozen
one-state coefficient record.
There are finite positive numbers
\[
 \delta_{\rm box}^{(2)},\quad
 C_{\rm K}^{(2)},\quad c_{\rm K}^{(2)},\quad
 K_{\rm J}^{(2)},\quad K_0^{(2)},\quad
 \overline\Gamma_{\rm 2st}
\]
which, when adjoined to that record, form a common A/B package
compatible with every subordinate pair
\(0<\sigma'<\theta'<\theta_*\).  This is a static coefficient and
radius selection; it assumes no second trajectory and invokes no
two-state continuation or decay estimate.
\end{lemma}

\begin{proof}
The exact polarization and local Kato calculation in
Lemma~\ref{lem:uniform-two-state-Kato-ledger} uses only the two finite
common ceilings and the fixed background, and therefore supplies
\(C_{\rm K}^{(2)},c_{\rm K}^{(2)}\) and a positive box size without
assuming a comparison solution.  The parameter-selection part of
Lemma~\ref{lem:uniform-two-state-package-radius} then chooses
\(K_{\rm J}^{(2)},K_0^{(2)}\) and
\(\overline\Gamma_{\rm 2st}\) before any subordinate rate pair.
Its later barrier conclusion is not used here.
\end{proof}

\subsection{Exact-core implantation and the high-regularity basin}

\begin{proposition}[Quantitative exact-core implantation]
\label{prop:exact-core-implant}
Let $(X^4,g_X)$ be a smooth closed connected oriented Riemannian
manifold, let $x_*\in X$, and let $B_{\rm imp}$ be an oriented
coordinate neighborhood of $x_*$.  Fix $k\geq12$, $0<\alpha<1$, and
choose the graft parameter $\Gamma$ sufficiently large for the
estimates of Sections~\ref{sec:three-region} and
\ref{sec:adaptive-continuation}.  In every continuation or
Theorem~\ref{thm:intro-sharp-scattering} application, this means the
already chosen package value
\(\Gamma\) satisfying
\eqref{eq:global-compatible-package-radius},
obtained after fixing the rate pair and the pre-radius ceilings in
Remark~\ref{conv:authoritative-adaptive-order}; the static implantation
itself uses only this fixed radius.  There is $A_*>0$ and, for every
$0<A<A_*$, a number $\tau_*(A)<\infty$ with the following property.
The constants may depend on
$(X,g_X,B_{\rm imp},\Gamma,k,\alpha)$; no uniformity over the class of
hosts is asserted.

For every $\tau_0\geq\tau_*(A)$, let
\begin{equation}\label{eq:exact-core-initial-parameters}
 \lambda_0=Ae^{-\tau_0},\qquad
 \Theta_0=\Phi_0=\varphi_{\tau_0},\qquad
 F_0=\operatorname{Id},\qquad
 S_0=\lambda_0\Theta_0^*\bar g.
\end{equation}
There exist an oriented blow-up
\[
 \varpi:\widehat X\longrightarrow X,\qquad
 \widehat X\cong X\#\overline{\mathbb {CP}}^{\,2},
\]
a marked open set $\widehat X''\subset\widehat X$, a diffeomorphism
$\iota:\widehat X''\to M$, and a smooth metric
$G_{A,\tau_0}$ on $\widehat X$ such that:
\begin{enumerate}
\item
      $\varpi$ is a diffeomorphism off the exceptional sphere and
      \begin{equation}\label{eq:implant-prescribed-exterior}
       G_{A,\tau_0}=\varpi^*g_X
       \quad\text{on }
       \widehat X\setminus\varpi^{-1}(B_{\rm imp});
      \end{equation}
\item
      on a fixed relatively compact open neighborhood
      $\mathcal K_\Gamma\Subset M$ of
      $\supp\eta\cup\overline{\Omega_\eta^{++}}$,
      \begin{equation}\label{eq:implant-exact-core}
       \iota_*G_{A,\tau_0}=S_0;
      \end{equation}
\item
      the prepared noncompact extension
      \begin{equation}\label{eq:implant-prepared-extension}
       \acute G_0
       =\eta\,\iota_*G_{A,\tau_0}+(1-\eta)S_0
      \end{equation}
      equals the complete metric $S_0$ on all of $M$, and hence
      \begin{equation}\label{eq:implant-zero-perturbation}
       h_0
       =\lambda_0^{-1}(\Phi_0^{-1})^*\acute G_0-\bar g
       =0;
      \end{equation}
\item
      there are fixed smooth core domains
       \[
       \begin{gathered}
        K_-\Subset K_0\Subset K_1\Subset K_2\Subset\widehat X'',\\
        \iota(\overline{K_2})\Subset\{\bar f<\Gamma/2\},
        \qquad
        \iota(\overline{K_-})\cap\overline{\Omega_\eta^{++}}
        =\varnothing,\\
        \iota(K_2\setminus\overline{K_-})
        \Subset\operatorname{int}\{\eta=1\},
       \end{gathered}
       \]
      and the associated physical exterior domains
      \[
       U_{\rm ext}^{++}:=\widehat X\setminus\overline{K_-},\qquad
       U_{\rm ext}^{+}:=\widehat X\setminus\overline{K_0},\qquad
       U_{\rm ext}:=\widehat X\setminus\overline{K_2}.
       \]
       There exist a compact collar
       \(\mathcal A_{\rm in}\Subset
       K_1\setminus\overline{K_0}\) and a fixed cutoff
       \(\zeta\) satisfying the support and separation conditions in
       \eqref{eq:prepared-exterior-termination}.
       In addition, there exist two separately buffered collar triples.  The
       interface collar is
       \[
        \mathcal W_{\rm in}\Subset\mathcal W_{\rm in}^+
        \Subset\mathcal W_{\rm in}^{++}
        \Subset U_{\rm ext}^+\cap\widehat X'',
        \qquad
        \mathcal A_{\rm in}
        \Subset\operatorname{int}\mathcal W_{\rm in}.
       \]
       The physical graft-input collar is
       \[
        \mathcal W_{\rm gr}\Subset\mathcal W_{\rm gr}^+
        \Subset\mathcal W_{\rm gr}^{++}
        \Subset U_{\rm ext}\cap\widehat X'',
        \qquad
        \Omega_\eta^{++}
        \Subset\operatorname{int}\iota(\mathcal W_{\rm gr}).
       \]
       Only the graft triple carries the auxiliary six-level chain
       \[
        \mathcal W_{\rm gr}^+
        \Subset\mathcal W_{\rm gr}^{0}\Subset\cdots
        \Subset\mathcal W_{\rm gr}^{5}
        \Subset\mathcal W_{\rm gr}^{++}
       \]
       with fixed positive scaled separations.
       An ordinary retained finite family
       \(U_a\Subset U_a^+\Subset U_a^{++}\Subset U_{\rm ext}^{++}\)
       with scales \(R_a\) covers \(\overline{U_{\rm ext}}\) by its
       smallest members and \(\overline{U_{\rm ext}^+}\) by its
       retained largest members.  Each \(U_a^{++}\) has an auxiliary chain
       \[
       U_a^{++}\Subset V_a^0\Subset\cdots\Subset V_a^5
       \Subset U_{\rm ext}^{++}.
       \]
       The family may be chosen with harmonic witness enlargements
       \[
        V_a^5\Subset W_a^{\rm har}
        \Subset\widetilde W_a^{\rm har}
        \Subset U_{\rm ext}^{++}
       \]
       and common constants
       \[
        \begin{gathered}
        \upsilon_{\rm har}^{\rm phys}>0,\qquad
        0<\eta_{\rm har}^{\rm phys}
          <\eta_{\rm har}^{\rm phys,+}
          <\eta_{\rm har}^{\rm phys,ref},\\
        0<q_{\rm har}^{\rm phys}
          <q_{\rm har}^{\rm phys,+}
          <q_{\rm har}^{\rm phys,ref}<Q_{\rm har}-1,\qquad
        0<\zeta_{\rm har}^{\rm phys}
          <\zeta_{\rm har}^{\rm phys,+}
          <\zeta_{\rm har}^{\rm phys,ref}<1,
        \end{gathered}
       \]
       such that, for every \(x\in V_a^5\), one may choose once and
       for all a pair \((u_{a,x},D_{a,x})\) for which the metric
       \(G_{A,\tau_0}\) has a
       coefficient-\(q_{\rm har}^{\rm phys,ref}\),
       domain-\(\zeta_{\rm har}^{\rm phys,ref}\) harmonic witness on
       \[
        B_{G_{A,\tau_0}}\!\left(
          x,
          (\upsilon_{\rm har}^{\rm phys}
          +2\eta_{\rm har}^{\rm phys,ref})R_a
        \right),
       \]
       and, on writing
       \[
        r_{a,\rm ref}^{\rm wit}
        :=(\upsilon_{\rm har}^{\rm phys}
           +2\eta_{\rm har}^{\rm phys,ref})R_a,\qquad
        \mathcal D_{a,x}
        :=u_{a,x}^{-1}(r_{a,\rm ref}^{\rm wit}D_{a,x}),
       \]
       the manifold preimage of its recorded Euclidean Dirichlet
       domain is compactly contained in \(W_a^{\rm har}\).
       Quantitatively, the choices have one center-independent outer
       separation
       \[
        \inf_{\substack{1\leq a\leq N_{\rm ext}\\x\in V_a^5}}
        R_a^{-1}\operatorname{dist}_{\rm ref}
        \bigl(\overline{\mathcal D_{a,x}},
              \widehat X\setminus W_a^{\rm har}\bigr)
        \geq\zeta_{\rm out}^{\rm phys}>0 .
       \]
       By restriction these charts also carry, in order, the displayed
       common \(+\)-tier and operative tier.  The reference tier is
       reserved for the static reference-to-carrier transfer, while the
       \(+\)-tier is reserved for the later carrier-to-flow transfer.
       These data satisfy
       \eqref{eq:prepared-exterior-termination} with
      \(E=U_{\rm ext}\), \(E^+=U_{\rm ext}^+\), and
      \(E^{++}=U_{\rm ext}^{++}\).  No cyclic cross-cover
      is imposed.  There are constants $c,C>0$,
      independent of $\tau_0$, such that
      \[
       \widehat X\setminus
       \iota^{-1}\bigl(\{\bar f<\Gamma/2\}\bigr)
       \subset U_{\rm ext},
      \]
      and, with
      \(R_{\rm ext}=c\sqrt A\) and
      \(R_{\min}:=\min_aR_a\geq R_{\rm ext}\),
      \begin{align}
       \min\Bigl\{
       d_{G_{A,\tau_0}}
       (U_{\rm ext},\widehat X\setminus U_{\rm ext}^+),
       d_{G_{A,\tau_0}}
       (U_{\rm ext}^+,\widehat X\setminus U_{\rm ext}^{++})
       \Bigr\}
       &\geq10R_{\rm ext},
       \label{eq:implant-buffer}\\
       \inf_{x\in U_{\rm ext}^{++}}
       r_{\rm har}(G_{A,\tau_0},x)
       &\geq cR_{\rm ext},
       \label{eq:implant-harmonic-radius}\\
       R_{\rm ext}^{2+\ell}
       |(\nabla^{G_{A,\tau_0}})^\ell
          \Rm_{G_{A,\tau_0}}|
       &\leq C_k,\qquad0\leq\ell\leq k-2,
       &&\text{on }U_{\rm ext}^{++}.
       \label{eq:implant-curvature-package}
      \end{align}
       The finite-atlas constants
       \(N_{\rm cov},C_{\rm cov},R_{\rm in},c_{\rm atl},
       \ell_{\rm Leb},\Lambda_{\rm atl}\) are fixed after \(A\) is
       chosen and are then independent of \(\tau_0\).
      The same estimates and the witnessed harmonic-coordinate
       certificates, with the same operative triples and fixed positive
       \(+\)-tier reserve margins, hold in a
      sufficiently small $C^{k,\alpha}$ neighborhood of
      $G_{A,\tau_0}$.
\end{enumerate}
The equality in \eqref{eq:implant-prepared-extension} concerns the
complete auxiliary metric on $M$; it does not assert that the closed
metric equals $S_0$ beyond the marked region.
\end{proposition}

\begin{proof}
Choose oriented coordinate balls
\[
 x_*\in B_0\Subset B_1\Subset B_{\rm imp}.
\]
In the oriented chart, replace the central ball by
$D(\mathcal O_{\mathbb P^1}(-1))$.  Use the standard local complex
blowdown on the disk bundle and choose the boundary identification
compatibly; with respect to the induced boundary orientations this
identification reverses orientation.  Thus
\[
 \widehat X
 =
 D(\mathcal O_{\mathbb P^1}(-1))
 \cup_{S^3}(X\setminus\operatorname{int}B^4)
 \cong X\#\overline{\mathbb {CP}}^{\,2}.
\]
The zero section has self-intersection $-1$.  The disk bundle together
with an open collar in $X\setminus B^4$ is diffeomorphic to $M$:
radially reparametrize a finite open collar as
$S^3\times[0,\infty)$.  This supplies $\widehat X''$ and $\iota$.
Choose the collar so that
$\iota^{-1}(\mathcal K_\Gamma)\Subset\varpi^{-1}(B_0)$.

 On $\iota^{-1}(\mathcal K_\Gamma)$ prescribe the pullback of $S_0$.
 The non-Euclidean FIK asymptotic cone \(g_{\rm C}\) and its all-order
 symbol estimate have already been fixed in
 Lemma~\ref{lem:FIK-AC-symbol}.  Combining
\eqref{eq:FIK-AC-quantitative} with the ODE for \(\varphi_\tau\) shows
 that, on every fixed normalized annulus used in the construction,
\begin{equation}\label{eq:implant-conical-scale}
 S_0=Ae^{-\tau_0}\varphi_{\tau_0}^*\bar g
\end{equation}
 has physical diameter comparable to $\sqrt A$ and, after rescaling by
 $A^{-1}$, belongs to a bounded family in every fixed $C^k$ norm as
 $\tau_0\to\infty$.  The rescaled family converges to the corresponding
 annulus in \(g_{\rm C}\), with an \(O_k(e^{-\tau_0})\) symbol error
 after the preceding fixed-annulus identification.

Choose an orientation-preserving parametrization
\[
 \kappa:B_{r_*}(0)\subset\mathbb R^4\longrightarrow B_0,
 \qquad \kappa(0)=x_*,
\]
whose differential at the origin is an oriented isometry from the
Euclidean metric to $g_X(x_*)$.
After decreasing $A_*$, the physical transition annulus, whose radius
is $O(\sqrt{A\Gamma})$, lies in this chart.  On a fixed Euclidean
annulus $\mathcal A$ define $d_A(y)=\kappa(\sqrt A\,y)$.  Then
\begin{equation}\label{eq:host-dilation-limit}
 A^{-1}d_A^*g_X\longrightarrow g_{\rm Euc}
 \quad\text{in }C^k(\mathcal A)
 \quad\text{as }A\downarrow0.
\end{equation}
Use the chosen boundary identification and radial collar map to put
the $A^{-1}$-rescaled FIK-side collar on the same annulus
$\mathcal A$.  By \eqref{eq:implant-conical-scale}, that family
converges in $C^k$ to its conical annulus as $\tau_0\to\infty$.  The
implant and collar maps may depend on $A$, but once $A$ is chosen they
are fixed independently of $\tau_0$.

 Use a fixed cutoff on $\mathcal A$ to join the two endpoint metrics.
 The limiting endpoints are deliberately not identified:
 \(g_{\rm C}\) is non-Euclidean, whereas
 \eqref{eq:host-dilation-limit} is Euclidean.  The interpolation needs
 only uniform bounded geometry and mutual equivalence on the fixed
 compact annulus, not small difference between the endpoints.  After
 restricting $A_*$ and increasing $\tau_*(A)$, compactness of the two
 limiting positive-definite families gives precisely that uniform
 equivalence.  Convexity of the
positive-definite cone preserves positivity, and the fixed cutoff
together with the uniform $C^k$ bounds gives uniform ellipticity and
derivative bounds.  Make the interpolation equal to $S_0$ on
$\mathcal K_\Gamma$ and equal to $\varpi^*g_X$ before leaving $B_1$.
Extension by $\varpi^*g_X$ proves
\eqref{eq:implant-prescribed-exterior} and
\eqref{eq:implant-exact-core}.

 Choose four nested smooth core domains
\[
 K_-\Subset K_0\Subset K_1\Subset K_2
 \Subset\iota^{-1}(\mathcal K_\Gamma)
\]
with
\[
 \iota(\overline{K_2})\Subset\{\bar f<\Gamma/2\},
 \qquad
 \overline{\Omega_\eta^{++}}\Subset\{\bar f>\Gamma/2\},
\]
whose intervening collar is mapped compactly into \(\{\eta=1\}\), and
whose boundaries lie on fixed normalized annuli.  Under
\eqref{eq:implant-conical-scale}, their successive physical
separations are comparable to \(\sqrt A\), uniformly in \(\tau_0\).
Define \(U_{\rm ext}^{++},U_{\rm ext}^+,U_{\rm ext}\) as their
complements, as in the statement.  Choose a smooth compact collar
\(\mathcal A_{\rm in}\Subset
K_1\setminus\overline{K_0}\), and choose \(\zeta\) to make its
 transition inside the interior of that collar, with \(\zeta=0\) near
 \(\partial U_{\rm ext}^{++}\) and \(\zeta=1\) near
 \(\overline{U_{\rm ext}}\).  Choose separately an interface triple
 \[
  \mathcal W_{\rm in}\Subset\mathcal W_{\rm in}^+
  \Subset\mathcal W_{\rm in}^{++}
  \Subset U_{\rm ext}^+\cap\widehat X''
 \]
 with \(\mathcal A_{\rm in}\Subset
 \operatorname{int}\mathcal W_{\rm in}\).  Since
 \(K_2\) lies strictly inside the \(\eta=1\) region, the physical graft
 transition lies in \(U_{\rm ext}\).  Choose around it the disjoint-role
 graft triple
 \[
  \mathcal W_{\rm gr}\Subset\mathcal W_{\rm gr}^+
  \Subset\mathcal W_{\rm gr}^{++}
  \Subset U_{\rm ext}\cap\widehat X'',
  \qquad
  \Omega_\eta^{++}\Subset
  \operatorname{int}\iota(\mathcal W_{\rm gr}).
 \]
 By compact containment, insert
 \(\mathcal W_{\rm gr}^{0},\ldots,\mathcal W_{\rm gr}^{5}\)
 between \(\mathcal W_{\rm gr}^+\) and
 \(\mathcal W_{\rm gr}^{++}\), with successive separations comparable
 to the local graft-collar scale.  These constants are frozen after
 \(A\) is chosen.  No collar in \(U_{\rm ext}\) is asked to contain the
 interior interface \(\mathcal A_{\rm in}\).

On the scale-\(\sqrt A\) transition choose retained triply nested
rescaled coordinate sets and, beyond every retained largest set, six
auxiliary enlargements, with each successive separation at least
\(20c\sqrt A\) and every diameter at most \(C\sqrt A\).  On the
remaining fixed-host region make the same
choice using finitely many coordinate balls with positive host-scale
separations.  By compactness, this ordinary finite family may be
chosen so that its smallest members cover
\(\overline{U_{\rm ext}}\), its retained largest members cover
\(\overline{U_{\rm ext}^+}\), and all auxiliary outer members remain
compactly contained in \(U_{\rm ext}^{++}\).  Enlarge each
\(V_a^5\) twice more to
\[
 V_a^5\Subset W_a^{\rm har}
 \Subset\widetilde W_a^{\rm har}
 \Subset U_{\rm ext}^{++}.
\]
On the scale-\(\sqrt A\) transition, the rescaled metrics on these
finitely many enlargements range in one compact bounded-geometry
family; on the fixed-host members the same is true by ordinary
compactness.  Let
\[
 b_*:=
 \min_{1\leq a\leq N_{\rm ext}}
 R_a^{-1}\operatorname{dist}_{\rm ref}
 \bigl(\overline{V_a^5},
       \widehat X\setminus W_a^{\rm har}\bigr)>0 .
\]
Let \(C_{\rm cmp}\geq1\) be a common constant for which, on these
enlargements,
\[
 d_{\rm ref}(x,z)\leq C_{\rm cmp}
 d_{G_{A,\tau_0}}(x,z),
\]
as supplied by uniform ellipticity, and let
\(\varrho_{\rm BG}>0\) be a common dimensionless outer
harmonic-coordinate radius supplied by the rescaled bounded-geometry
package.  Choose two radii
\[
 0<\varrho_{\rm in}<\varrho_{\rm out}
 <\min\{\varrho_{\rm BG},b_*/(4C_{\rm cmp})\}.
\]
For every
\((a,x)\), \(x\in V_a^5\), there is a reserved harmonic chart on
\(B_{G_{A,\tau_0}}(x,\varrho_{\rm out}R_a)\), and its restriction to
the inner ball of radius \(\varrho_{\rm in}R_a\) admits a uniformly
regular Euclidean Dirichlet domain strictly between the normalized
inner and outer coordinate images.  The comparison just fixed gives
\[
 B_{G_{A,\tau_0}}(x,\varrho_{\rm out}R_a)
 \subset
 \bigl\{z:\operatorname{dist}_{\rm ref}(z,x)<b_*R_a/4\bigr\}.
\]
This follows by applying the quantitative harmonic-coordinate and
boundary-atlas constructions after rescaling \(R_a\) to one.
 Uniform ellipticity, bounded geometry, the lower harmonic-radius
 bound, and the fixed positive gap
 \(\varrho_{\rm out}-\varrho_{\rm in}\) give a reference coefficient
 reserve \(q_{\rm har}^{\rm phys,ref}>0\), a reference domain reserve
 \(\zeta_{\rm har}^{\rm phys,ref}>0\), and all Dirichlet-domain constants
depending only on the common rescaled package, not on \(x\) or
\(\tau_0\).  Decrease the two reserves, if necessary, so that
 \(q_{\rm har}^{\rm phys,ref}<Q_{\rm har}-1\) and
 \(\zeta_{\rm har}^{\rm phys,ref}<1\).  Choose positive numbers
 \(\upsilon_{\rm har}^{\rm phys}\) and
 \(\eta_{\rm har}^{\rm phys,ref}\) such that
\[
 \upsilon_{\rm har}^{\rm phys}
  +2\eta_{\rm har}^{\rm phys,ref}
 =\varrho_{\rm in},
\]
 and then choose a common \(+\)-tier strictly componentwise below the
 reference tier and an operative tier strictly componentwise below the
 \(+\)-tier.  For every \(x\in V_a^5\), select the resulting
 reference-tier centered witness pair
\((u_{a,x},D_{a,x})\).  Its recorded manifold Dirichlet preimage is
contained in the displayed outer metric ball, and hence
\[
 R_a^{-1}\operatorname{dist}_{\rm ref}
 \bigl(\overline{\mathcal D_{a,x}},
       \widehat X\setminus W_a^{\rm har}\bigr)
 \geq\frac12b_*.
\]
Set
\[
 \zeta_{\rm out}^{\rm phys}:=\min\{b_*/2,1/2\}>0.
\]
The only minimum used to define \(b_*\) is over the finitely many cover
indices \(a\); no minimum over the center-indexed witness family is
used.  The resulting finite cover has uniform overlap, scale-comparability,
normalized Lebesgue-number, partition, and chart constants.  No
largest member is required to be covered by smaller members.  Decrease
\(c\) and then \(A_*\) so that all relevant
buffer and cutoff separations are at least \(20c\sqrt A\).  This proves
\eqref{eq:prepared-exterior-termination} and
\eqref{eq:implant-buffer}, with
\[
 R_{\min}:=\min_aR_a\geq R_{\rm ext}=c\sqrt A .
\]
The
$A^{-1}$-rescaled transition data range in a fixed bounded family, so
scaling back gives the curvature and harmonic-radius estimates there.
On the fixed-host balls the same conclusions follow from compactness
and the positive harmonic radius of $g_X$.  All constants are
independent of $\tau_0$.  Ordinary finite-jet estimates persist by
openness.  Apply the reserve-to-operative clause of
Lemma~\ref{lem:finite-physical-harmonic-openness}, with the common
\(+\)-tier as its operative tier and the reference tier as its reserve,
on
\(U_a^{++}\Subset V_a^5\Subset W_a^{\rm har}
\Subset\widetilde W_a^{\rm har}\).  Under a sufficiently small
\(C^{k,\alpha}\) perturbation it preserves, at every center
\(x\in V_a^5\), the same common \(+\)-tier and leaves a new stronger
reserve; restriction then gives the unchanged operative tier.  Thus
the reference-to-carrier and carrier-to-flow gaps remain separate.

Finally, \eqref{eq:implant-exact-core} holds wherever $\eta\ne0$,
whereas $\eta=0$ outside the marked physical graft.  Therefore
\[
 \acute G_0=\eta S_0+(1-\eta)S_0=S_0
 \quad\text{on }M.
\]
Since $\Phi_0=\varphi_{\tau_0}$,
\[
 \lambda_0^{-1}(\Phi_0^{-1})^*S_0
 =(\varphi_{-\tau_0})^*\varphi_{\tau_0}^*\bar g
 =\bar g,
\]
which proves \eqref{eq:implant-zero-perturbation}.
\end{proof}

\begin{lemma}[Canonical static exact-core pre-radius bounds]
\label{lem:static-exact-core-pre-radius-ledger}
Fix the FIK background, the radial and cutoff profiles, the fixed
pre-radius atlas based at \(\Gamma_{\rm atl}\), and
\(0<\alpha<1\).  There exist a rate-independent
\(\tau_{\rm static}\geq1\), finite canonical constants
\[
 \Lambda_{\rm ell}^0,\quad
 \Lambda_{\rm coef}^0,\quad
 \Lambda_{\rm map}^0,\quad
 \Lambda_{R,14}^0,\quad
 \Lambda_{F,6}^0,\quad
 K_{\mathcal Y,m}^0\quad(0\leq m\leq4)
\]
and positive canonical margins
\[
 \kappa_{\rm sep}^0,\qquad
 \kappa_{\rm Gram}^0,\qquad
 \kappa_{\rm map}^0,\qquad
 \kappa_{\rm har}^0,\qquad
 q_{\rm har}^0,\ q_{\rm har}^{0,+},\qquad
 \zeta_{\rm har}^0,\ \zeta_{\rm har}^{0,+}
\]
 with the following property.  For every
\(\Gamma\geq\Gamma_{\rm atl}\) and
\(\tau_0\geq\tau_{\rm static}\), the exact normalized state
\[
 h_0=0,\qquad
 R_{\tau_0}=F_0=\operatorname{Id},\qquad
 \Theta_0=\Phi_0=\varphi_{\tau_0}
\]
lies strictly inside every primitive metric, coefficient, radial,
map, inverse-map, harmonic-radius, separation, and Gram face governed
by the listed constants and margins, whenever its upper ceiling is
 strictly larger and its lower margin strictly smaller than the
 corresponding canonical value.  In particular,
\[
 \mathfrak s_{\rm sep}(\tau_0)\geq\kappa_{\rm sep}^0,
 \qquad
 \mathfrak h_{\rm har}(\tau_0)\geq\kappa_{\rm har}^0,
 \qquad
 s_{\min}M^{\rm low}(\tau_0)\geq\kappa_{\rm Gram}^0 .
\]
For every fixed
\(\kappa_{\rm har}<\kappa_{\rm har}^0\), there exist positive margins
\(\eta_{\rm har}^{0}\) and \(\eta_{\rm har}^{0,+}\), with
\[
 0<\eta_{\rm har}^{0}<\eta_{\rm har}^{0,+},
 \qquad
 \kappa_{\rm har}+2\eta_{\rm har}^{0,+}
 <\kappa_{\rm har}^0,
\]
such that the exact state carries uniform buffered harmonic witnesses
with reserve triple
\[
 (\eta_{\rm har}^{0,+},
   q_{\rm har}^{0,+},\zeta_{\rm har}^{0,+})
\]
and an operative triple
\[
 (\eta_{\rm har}^{0},q_{\rm har}^{0},\zeta_{\rm har}^{0})
\]
strictly smaller componentwise.
Its graft discrepancy is zero and
therefore has room below every fixed \(K_{\rm gr}>0\).  No scale face
is asserted here; the separate choice \(A<C_{\rm sc}\) is made in
Corollary~\ref{cor:strict-entrance-nonempty}.  Moreover,
\begin{equation}\label{eq:static-exact-core-column-ledger}
 \sup_{\Gamma\geq\Gamma_{\rm atl}}
 \sup_{\tau_0\geq\tau_{\rm static}}\sup_M
 \sum_{\ell=0}^{m}(1+\bar f)^{\ell/2}
 |\bar\nabla^\ell\mathcal Y_{j,\tau_0}|_{\bar g}
 \leq K_{\mathcal Y,m}^0,
 \qquad 0\leq m\leq4,\quad0\leq j\leq8,
\end{equation}
and the cutoff Gram matrix and support separation have the stated
canonical positive margins.  The graft face is treated separately in
Definition~\ref{def:strict-prepared-entrance} by its explicit
\(K_{\rm gr}/2\) inequality after the scale and phase choices.  These
constants are independent of the
host metric, the implantation scale \(A\), and every ordinary
physical exterior-cover constant chosen after \(\Gamma\).
\end{lemma}

\begin{proof}
At the exact state the normalized source and target metrics are both
\(\bar g\), while both relative maps are the identity.  Thus the
 ellipticity, coefficient, harmonic-radius, map, inverse-map, radial,
 and lower-singular-value constants are fixed FIK constants in the
 fixed pre-radius atlas.  Pullback and metric scaling cancel in
 \eqref{eq:dimensionless-harmonic-radius-functional}, so
 \(\mathfrak h_{\rm har}(\tau_0)\) has one positive canonical lower
 bound.  Choose the witnessing radii strictly below that lower bound.
The fixed core charts and the scale-normalized asymptotic atlas then
give a uniform coefficient reserve; restricting each chart from a
slightly larger metric ball gives a uniform domain buffer.  These are
the reserve constants
\(q_{\rm har}^{0,+},\zeta_{\rm har}^{0,+}\) above.  Choose the
operative radius, coefficient, and domain reserves strictly
componentwise below the reserve triple; these are
\(\eta_{\rm har}^{0},q_{\rm har}^{0},\zeta_{\rm har}^{0}\).  In
particular these constants do not see the
subsequent physical interpolation used to close the host manifold.

 Since \(\Phi_0=\Theta_0=\varphi_{\tau_0}\), put
\[
 c_{\Gamma,\tau_0}:=
 (1-\eta_\Gamma)\circ\varphi_{-\tau_0}.
\]
 Naturality of pullback gives, for every fixed FIK tensor \(T\),
\begin{equation}\label{eq:static-exact-core-K-formula}
 K_{\tau_0}(T)
 =(\Phi_0^{-1})^*((1-\eta_\Gamma)\Theta_0^*T)
 =c_{\Gamma,\tau_0}T .
\end{equation}
 Therefore
\begin{align}
 \mathcal Y_{0,\tau_0}
 &= (1-c_{\Gamma,\tau_0})Y_0,
 \label{eq:static-exact-core-column-zero}\\
 \mathcal Y_{j,\tau_0}
 &= (1-c_{\Gamma,\tau_0})
    \Lie_{\chi_{\tau_0}W_j}\bar g,\qquad 1\leq j\leq8.
 \label{eq:static-exact-core-column-j}
\end{align}
In particular, with the unnormalized symmetric product
\[
 u\odot v:=u\otimes v+v\otimes u,
\]
one has
\[
 \Lie_{\chi_{\tau_0}W_j}\bar g
 =\chi_{\tau_0}Y_j
  +d\chi_{\tau_0}\odot W_j^\flat,
\]
and hence the exact two-cutoff identity is
\begin{equation}\label{eq:static-exact-core-two-cutoff-tail}
 \begin{split}
 \mathcal Y_{j,\tau_0}-Y_j
 ={}&
 -\bigl[c_{\Gamma,\tau_0}
 +(1-c_{\Gamma,\tau_0})(1-\chi_{\tau_0})\bigr]Y_j\\
 &+(1-c_{\Gamma,\tau_0})
   d\chi_{\tau_0}\odot W_j^\flat .
 \end{split}
\end{equation}
This explicitly retains the geometric-column cutoff tail; it is not
absorbed into \(c_{\Gamma,\tau_0}\).

 We now estimate the two cutoffs directly, before any
package-dependent shifted-annulus lemma is available.  The radial ODE
 \(\partial_s(\bar f\circ\varphi_s)
 =(\bar f-\bar R)\circ\varphi_s\), the FIK AC symbol bounds, and the fixed-profile chain
rule give, through order five, uniform scale-normalized derivative
bounds for \(c_{\Gamma,\tau_0}\) and \(\chi_{\tau_0}\).
The former is supported in
\(\{\bar f\geq c\Gamma e^{\tau_0}\}\), while
\(1-\chi_{\tau_0}\) and \(d\chi_{\tau_0}\) are supported in
\(\{\bar f\geq(5/2)e^{\tau_0}\}\).  Hence, directly from the
FIK generator symbols and Gaussian integration,
\[
 \|\mathcal Y_{j,\tau_0}-Y_j\|_{H_\nu^m}
 \leq C_m e^{-c_m e^{\tau_0}},
 \qquad 0\leq m\leq4,\quad0\leq j\leq8,
\]
uniformly for \(\Gamma\geq\Gamma_{\rm atl}\).  The same direct symbol
calculation, without the Gaussian weight, gives
\eqref{eq:static-exact-core-column-ledger}.  No evolution, first-exit
interval, or member of \(\mathfrak P_{\rm pre}^{\rm der}\) has been
invoked.

 Choose \(\tau_{\rm static}\) sufficiently large, using only the fixed
radial and cutoff profiles, that on \(\supp\rho_{\tau_0}\) one has
\[
 c_{\Gamma,\tau_0}=0,\qquad \chi_{\tau_0}=1
 \quad
 (\Gamma\geq\Gamma_{\rm atl},\ \tau_0\geq\tau_{\rm static}).
\]
Thus the effective columns equal \(Y_j\) there exactly; the remaining
Gram error is only the Gaussian truncation tail, which is uniformly
superexponentially small by the preceding estimate.  Choose the canonical
lower margins below the background Gram least singular value, the exact support
separation, the identity-map singular value, and the exact
dimensionless harmonic-radius value.  This proves the asserted bounds.
\end{proof}

\begin{corollary}[Exact FIK cores are strict prepared entrances]
\label{cor:strict-entrance-nonempty}
Fix \(k_0\geq12\), \(0<\alpha<1\), and
$0<\sigma<\theta<\beta$.  Take the canonical static constants and
margins from
Lemma~\ref{lem:static-exact-core-pre-radius-ledger}.  Fix primitive
pre-radius data with strict room:
\[
 \begin{gathered}
 \Lambda_{\rm ell}>\Lambda_{\rm ell}^0,\qquad
 \Lambda_{\rm coef}>\Lambda_{\rm coef}^0,\qquad
 \Lambda_{R,14}^{\rm pre}>\Lambda_{R,14}^0,\qquad
 \Lambda_{F,6}^{\rm pre}>\Lambda_{F,6}^0,
 \\[2pt]
 \Lambda_{\rm map}>
 \max\{\Lambda_{\rm map}^0,
        \Lambda_{R,14}^{\rm pre},
        \Lambda_{F,6}^{\rm pre}\},
 \end{gathered}
\]
and then fix the explicit norm slacks so that
\[
 0<2\mu_R^{\rm pre}
 <\Lambda_{R,14}^{\rm pre}-\Lambda_{R,14}^0,\qquad
 0<2\mu_F^{\rm pre}
 <\Lambda_{F,6}^{\rm pre}-\Lambda_{F,6}^0.
\]
Choose each primitive lower margin strictly below its canonical
counterpart, and choose \(K_{\rm gr}>0\) and a broad upper scale
constant \(C_{\rm sc}\).  Apply
Lemma~\ref{lem:pre-radius-low-order-closure} to this primitive tuple.
It produces \(\Gamma_{\rm pre}\) and the derived forcing, column,
feedback, annulus, map, smallness, and entrance-time constants.
Before evaluating a radius functional, enlarge the derived upper
constants once, if necessary, so that
\[
 K_{\mathcal Y,m}^{\rm pre}>K_{\mathcal Y,m}^0
 \quad(0\leq m\leq4),
\]
and shrink favorable derived margins once so that the exact state has
strict room.  Replace the selected witness by this one-sided
enlarged/shrunk witness and freeze the resulting tuple; every later
reference to \(\mathfrak P_{\rm pre}^{\rm der}\) means this final
witness.  The replacement preserves every closure estimate and occurs
before either radius is evaluated.

To make the same center compatible with the later scattering
refinement, reserve the static ceilings of the common two-state
package.  This finite numerical selection uses no two-state evolution
estimate and is unnecessary for the proof of Theorem~A.  Fix
\(\Lambda_{\rm ell}^{(2)}\) and
 \(K_{\mathcal Y,0}^{(2)}\) above the corresponding common two-state
 ceilings.  Apply
Lemma~\ref{lem:static-formation-scattering-compatibility} with
\(\theta_*=\theta\), and freeze its Kato datum, reduced
\(\delta_{\rm box}^{(2)}\), positive
\(K_{\rm J}^{(2)},K_0^{(2)}\), and radius output.
Now evaluate \(\overline\Gamma_{\rm 3reg}\) and
\(\overline\Gamma_{\rm 2st}\), choose \(\Gamma\) satisfying
\eqref{eq:global-compatible-package-radius},
and construct \(\eta_\Gamma\) and its graft collars.  Now apply
Proposition~\ref{prop:exact-core-implant} at an implantation order at
least \(k_0+2\).  Choose
\[
 0<A<\min\{A_*,C_{\rm sc}\}
\]
sufficiently small, and only then choose the lower scale margin
\(c_{\rm sc}\) so that the strict entrance scale bracket is
\begin{equation}\label{eq:exact-core-scale-bracket}
 0<c_{\rm sc}<A<C_{\rm sc}<\infty .
\end{equation}
This post-\(A\) choice of \(c_{\rm sc}\) is allowed because that lower
margin does not enter either package-radius threshold.  Complete and
freeze the remaining package entries, including
\(\delta_{\rm c2}\), \(\varepsilon_{\rm ph}\),
\(\varepsilon_{\rm ent}\), and \(\tau_{\rm ad}\), in the
fixed order of
Remark~\ref{conv:authoritative-adaptive-order}.  In the choice of
\(\delta_{\rm c2}\), impose in addition
\[
 4C_{\rm har,pre}\delta_{\rm c2}
 <
 \delta_{\rm har}
 (\eta_{\rm har}^0,q_{\rm har}^0,\zeta_{\rm har}^0;
  \mathfrak P_{\rm har}^{\rm geom}),
\]
which is possible because the canonical static witnesses and their
positive reserves were fixed before this smallness choice.  The
ordinary physical
 exterior-cover, buffer, and high-regularity implantation constants
 furnished after \((\Gamma,A)\) are post-radius data: they enter
 continuation and the final entrance-time bound, but neither
 \(\mathfrak G_{\rm 3reg}\) nor \(\mathfrak G_{\rm 2st}\).

 For these fixed \((\Gamma,A)\), the scale-\(R_a\) coefficient bounds
 through order fourteen on the physical cover, its buffer constants,
 and the reference-tier harmonic witnesses furnished by
 Proposition~\ref{prop:exact-core-implant} are uniform for
 \(\tau\geq\tau_*(A)\).  Indeed, after scaling the transition members by
 \(A^{-1}\), they range in the compact bounded-geometry family used in
 the proof of that proposition, while the fixed-host members range in
 an ordinary compact family.  Apply
 Lemma~\ref{lem:reference-carrier-physical-certificate} uniformly to
 this family.  It furnishes the common numerical choices
 \[
  \varepsilon_{\rm coeff}^{\rm phys}>0,\qquad
  \Lambda_{\rm coeff}^{\rm phys}<\infty,\qquad
  \mu_{\rm coeff}^{\rm phys}>0,\qquad
  \delta_{\rm RF}>0
 \]
 independently of the final entrance time; the particular ball and
 inner locus will be instantiated only after their reference center is
 frozen.  Set
 \(\mu_{\rm RF}=1/2\).  These are post-radius package data, and they are
 now fixed before the final choice of \(\tau_0\).
 Finally,
  \(\tau_0\geq
  \max\{\tau_*(A),\tau_{\rm ad},\tau_{\rm static}\}\) may be chosen
sufficiently large so that
$G_{A,\tau_0}$, with the data
\eqref{eq:exact-core-initial-parameters}, is a strict prepared FIK
entrance of continuation order \(k_0\).  Given $\varepsilon_T>0$, the
choices may also be made so
that
\begin{equation}\label{eq:exact-core-small-time-margin}
2C_\lambda\lambda_0<\varepsilon_T,
 \qquad
 2C_\lambda\lambda_0
 <\frac12\delta_{\rm RF}R_{\min}^2.
\end{equation}
\end{corollary}

\begin{proof}
Set $\mathcal X=\widehat X$ and use the marked identification supplied
by Proposition~\ref{prop:exact-core-implant}.  We verify the six
entrance conditions in their stated order.

The scale face is strict because
\(\lambda_0=Ae^{-\tau_0}\) and
\eqref{eq:exact-core-scale-bracket} gives
\[
 c_{\rm sc}e^{-\tau_0}
 <\lambda_0
 <C_{\rm sc}e^{-\tau_0}.
\]
The same final choice of \(\tau_0\) gives
\(\tau_0\geq\tau_{\rm ad}\), as required by
\eqref{eq:strict-entrance-adaptive-time}.

First, \eqref{eq:implant-zero-perturbation} makes all nine receding
moments zero, so Proposition~\ref{prop:uniform-receding-phase} selects
the zero phase.  Moreover \(h_0=0\), \(R_{\tau_0}=F_0=\operatorname{Id}\),
and all background components are smooth with the fixed
scale-normalized bounds, so the already-sliced tuple belongs to a
common-margin ball in
\(\mathscr P_{\tau_0}^{k_0+2,\alpha}\) and uses route~(b) of
Definition~\ref{def:strict-prepared-entrance}.  Second, every
weighted, pointwise, and global $C^2$
norm of $h_0$ is zero and the normalized metric has exact ellipticity.
Since the constants in
\eqref{eq:strict-entrance-historical-constants} are already fixed,
increasing $\tau_0$ gives strict room in the explicitly typed condition
\eqref{eq:strict-entrance-historical-tail-absorption}.

Third,
\[
 R_{\tau_0}
 =\varphi_{-\tau_0}\circ\Theta_0
 =\operatorname{Id},
\]
and the extension equals $S_0$ outside the graft.  The
compatible package radius \(\Gamma\) is already frozen.  Because
\(\Phi_0=\Theta_0=\varphi_{\tau_0}\), the normalized
background map and inverse-map bounds are exact.  The support
separation in \eqref{eq:coarse-Gram-support-separation} gives
\[
 \rho_{\tau_0}K_{\tau_0}(T)=0,\qquad
    \rho_{\tau_0}\mathcal Y_{j,\tau_0}
    =\rho_{\tau_0}Y_j ,
\]
while, because \(F_0=R_{\tau_0}=\operatorname{Id}\),
\eqref{eq:static-exact-core-K-formula} computes every off-core
correction directly.  The uniform order-four column bounds and all
static map, coefficient, and lower-margin inequalities therefore follow
from Lemma~\ref{lem:static-exact-core-pre-radius-ledger}, with the strict
room built into the preceding choices.  The evolutionary closure lemma
is invoked only after its remaining hypotheses are verified below.
Finally,
Lemma~\ref{lem:cutoff-tails} and the Gaussian estimate
\eqref{eq:general-Gaussian-tail} show that increasing only the final
 entrance time \(\tau_0\) makes every remaining cutoff and Gram tail
 strictly smaller than its prescribed margin.
In particular, the direct identities
\eqref{eq:static-exact-core-column-zero}--%
\eqref{eq:static-exact-core-two-cutoff-tail} give
\eqref{eq:strict-entrance-column-faces}, while the choice
\(\tau_0\geq\tau_{\rm static}\) gives
\[
 s_{\min}M^{\rm low}(\tau_0)>\kappa_{\rm Gram},
 \qquad
 \mathfrak s_{\rm sep}(\tau_0)>\kappa_{\rm sep}.
\]
Thus \eqref{eq:strict-entrance-Gram-separation-faces} holds
quantitatively, rather than merely as a qualitative disjointness.
Fourth, on
\(\Omega_\eta^{++}\) the closed
metric equals $S_0$, so the initial graft distance and pure graft
defect vanish, including their fixed-background derivatives.

 Fifth, \(F_0=\operatorname{Id}\) gives
\[
 d_{\rm rt,sc}^{6,\alpha}(F_0,\operatorname{Id})=0
 <\varepsilon_{\rm map}^{\rm HM},\qquad
 \min\{s_{\min}(dF_0),s_{\min}(dF_0^{-1})\}=1
 >\kappa_{\rm map}.
\]
It is proper, and the complete source and
target are the same metric $S_0$.  Via $\Theta_0$, this metric is
 isometric to $\lambda_0\bar g$, so its harmonic radius, curvature
derivatives, and scale-normalized bounded-geometry quantities have
the required bounds at the pointwise intrinsic length
\(r_{{\rm sol},\tau_0}(x)\).  At the exact center its defining identity
is
\[
 r_{{\rm sol},\tau_0}(x)^2
 =Ae^{-\tau_0}
   \bigl(1+\bar f(\varphi_{\tau_0}x)\bigr)>0.
\]
No entrance-time-uniform global comparison with
\(A(1+\bar f(x))\) is asserted: the end comparison is used only on the
AC region, while on the compact core the pullback-and-scaling
covariance in the proof of
Lemma~\ref{lem:static-exact-core-pre-radius-ledger} gives the
dimensionless harmonic-radius bound directly.  Thus the canonical
static margins give
\[
 \mathfrak h_{\rm har}(\tau_0)>\kappa_{\rm har},
 \qquad
 \mathfrak C_{\rm ent}^{12,10}(\mathscr J_0)
 <\Lambda_{\rm coef}.
\]
More precisely, the global exact-core witnesses in
Lemma~\ref{lem:static-exact-core-pre-radius-ledger} supply the operative
triple
\((\eta_{\rm har}^0,q_{\rm har}^0,\zeta_{\rm har}^0)\)
and its strictly stronger \(+\)-triple required in item~(5); this is
therefore the witnessed strict face, not merely the displayed lower
bound.
At \(h_0=0\) the exact
Gram system has zero right side and hence $a(\tau_0)=b(\tau_0)=0$;
its matrix is the invertible background Gram matrix plus a
Gaussian-superexponential cutoff tail.  With $a=b=0$, $S(t)$ is the
exact self-similar Ricci flow at the initial face.  Since the closed
metric and $S_0$ agree on an open neighborhood of the graft, their
spatial Ricci jets agree there as well.  Thus the source and target
Ricci defects required at the initial face vanish.
Equivalently, \eqref{eq:S-defect} and
\eqref{eq:pure-graft-defect} give there
\[
 \partial_tS+2\Ric_S=0,\qquad
 \partial_t\acute G+2\Ric_{\acute G}
 =\mathcal G_{\rm gr}=0,
\]
and all required spatial derivatives vanish as well.

 Sixth, Proposition~\ref{prop:exact-core-implant} supplies the four
 nested core domains, the separated normalized interface, and the
 retained ordinary finite exterior atlas required in item~(6).  Its
 smallest members cover both the physical graft collar and the
 complementary outer region, its retained largest members cover
\(\overline{U_{\rm ext}^+}\), and its auxiliary six-level enlargements
 remain in \(U_{\rm ext}^{++}\).  The cutoff, atlas, Lebesgue-number,
 and buffer constants are part of the fixed prepared package, and
\[
 R_{\min}:=\min_aR_a\geq R_{\rm ext}=c\sqrt A .
\]
 The enlargements \(W_a^{\rm har}\) and the common constants
 \(\upsilon_{\rm har}^{\rm phys},\eta_{\rm har}^{\rm phys},
 q_{\rm har}^{\rm phys},\zeta_{\rm har}^{\rm phys}\), together with
 their strictly stronger \(+\)- and reference tiers, furnished by the
 same proposition give exactly the finite physical-cover witness
 certificate in item~(6).  The uniform application of
 Lemma~\ref{lem:reference-carrier-physical-certificate} made in the
 statement has already fixed the numerical choices
 \(\varepsilon_{\rm coeff}^{\rm phys}\),
 \(\Lambda_{\rm coeff}^{\rm phys}\),
 \(\mu_{\rm coeff}^{\rm phys}\), and \(\delta_{\rm RF}\), but not a
 center-dependent ball.  In the final selection of the statement,
 choose \(\tau_0\) once so that all preceding late-start requirements
 and
 \[
  2C_\lambda Ae^{-\tau_0}
  <\min\!\left\{
    \varepsilon_T,\,
    (1-\mu_{\rm RF})\delta_{\rm RF}R_{\min}^2
   \right\}
 \]
 hold.  Only then freeze
 \[
  G_{\rm ref}^{\rm phys}:=G_{A,\tau_0}.
 \]
 Instantiate the corresponding outer ball and inner locus with this
 center.  The exact carrier then belongs to the inner locus with more
 than \(4\mu_{\rm coeff}^{\rm phys}\) package-face slack.  It is
 therefore a valid actual carrier without recentering, and the uniform
 half-open quarter-window conclusion holds throughout that outer ball.
 No parameter is changed after the reference is frozen.  Keeping the
 reference fixed, the single final choice gives both inequalities in
 \eqref{eq:exact-core-small-time-margin}.  All six
conditions hold simultaneously with positive margins.  The three
parameters have separate roles: $\Gamma$ controls the normalized
graft, $A$ its physical implantation scale, and $\tau_0$ the much
smaller singular scale and remaining physical time.
\end{proof}

\begin{theorem}[Prepared high-regularity basin]
\label{thm:prepared-open-basin}
Let $(X^4,g_X)$ be a smooth closed connected oriented Riemannian
manifold, let $x_*\in X$, let $B_{\rm imp}$ be an oriented coordinate
neighborhood of $x_*$, and let $\varepsilon_T>0$.  Fix
 $0<\sigma<\theta<\beta$, an integer $k\geq14$, and
 $0<\alpha<1$.  There exist an oriented blow-up
\[
 \varpi:\widehat X\longrightarrow X,\qquad
 \widehat X\cong X\#\overline{\mathbb {CP}}^{\,2},
\]
a smooth metric $G_*$ on $\widehat X$, and a relatively
$C^{k+2,\alpha}$-open neighborhood
\[
 \mathscr U\subset\operatorname{Met}^{\infty}(\widehat X)
\]
of $G_*$ with the following properties.
\begin{enumerate}
\item
      The center is the exact-core entrance of
      Corollary~\ref{cor:strict-entrance-nonempty}; in particular its
      complete prepared extension has $h_0=0$.  It is prescribed
      outside the implantation region:
      \begin{equation}\label{eq:universal-center-exterior}
       G_*=\varpi^*g_X
       \quad\text{on }
       \widehat X\setminus\varpi^{-1}(B_{\rm imp}).
      \end{equation}
\item
      For every $G_0\in\mathscr U$, the maximal Ricci flow
      $G(t;G_0)$ has a singular time
      \begin{equation}\label{eq:universal-small-singular-time}
       0<T(G_0)<\varepsilon_T
      \end{equation}
      and develops a global Type-I singularity with the canonical FIK
      flow as its quantitative, no-subsequence marked parabolic limit.
\item
      There is $C_{\rm out}<\infty$, independent of
      $G_0\in\mathscr U$, such that
      \begin{equation}\label{eq:universal-exterior-curvature}
       \sup_{G_0\in\mathscr U}
       \sup_{0\leq t<T(G_0)}
       \sup_{\widehat X\setminus\varpi^{-1}(B_{\rm imp})}
       |\Rm_{G(t;G_0)}|
       \leq C_{\rm out}.
      \end{equation}
\item
      Every flow satisfies all quantitative conclusions of
      Theorem~\ref{thm:prepared-entrance-continuation}, including the
      receding-domain estimates and quantitative marked spacetime
      convergence.  The singular time and
      asymptotic data vary continuously by
      Lemma~\ref{lem:one-state-terminal-continuity}.
\end{enumerate}
The neighborhood $\mathscr U$ is open in the full space of metrics:
its elements need not preserve the exact core or the prescribed
exterior, and need not be $U(2)$-invariant or K\"ahler.  The
neighborhood and its constants may depend on the fixed host and
implantation data.
Here \(k\) is the phase-output order and
\[
 k_0:=k-2\geq12
\]
is the continuation output order.  The prepared phase chart uses the
two-derivative buffer
\[
 C^{k+2,\alpha}\longrightarrow
 \mathscr P_{\tau_0}^{k,\alpha}
 =\mathscr P_{\tau_0}^{k_0+2,\alpha},
\]
which is exactly the input order required for continuation and for
Lemma~\ref{lem:one-state-terminal-continuity}.
The later low-topology restart similarly asks for a prepared tuple of
order \(k_0+2\), and therefore smooths the time-\(d\) metric to order
\(k_0+4\).  Both constructions use the same two-derivative regularity
convention.
\end{theorem}

\begin{proof}
Put \(k_0=k-2\).  Thus \(k_0\geq12\), and every phase-adjusted tuple
constructed below lies in a common-margin ball of
\(\mathscr P_{\tau_0}^{k_0+2,\alpha}\).
Choose
\[
 x_*\in B_0\Subset B_1\Subset B_{\rm imp}
\]
and choose fixed open sets \(W\Subset W^+\) with
\[
 \widehat X\setminus\varpi^{-1}(B_{\rm imp})\subset W,
 \qquad
 \overline{W^+}\subset
 \widehat X\setminus\varpi^{-1}(\overline{B_1}).
\]
Having fixed \(W\Subset W^+\), make the package and implantation
choices exactly in the fixed order of
Corollary~\ref{cor:strict-entrance-nonempty}, with the interpolation
completed inside \(B_1\): fix the pre-radius ceilings, choose the
compatible \(\Gamma\), apply
Proposition~\ref{prop:exact-core-implant}, then choose \(A\) and the
lower scale margin \(c_{\rm sc}\), and finally choose \(\tau_0\).
Leave enough room in the final choice that
\begin{equation}\label{eq:universal-time-strict-margin}
 4C_\lambda\lambda_0<\varepsilon_T,
 \qquad
 4C_\lambda\lambda_0
 <\delta_{\rm RF}R_{\rm ext}^2.
\end{equation}
Corollary~\ref{cor:strict-entrance-nonempty} gives a strict exact-core
entrance.  Reset its physical time to zero, retain the normalized
label $\tau_0$, and call its closed metric $G_*$.  Equation~
\eqref{eq:implant-prescribed-exterior} gives
\eqref{eq:universal-center-exterior}.
Retain
\[
 G_{\rm ref}^{\rm phys}:=G_*
\]
as the frozen physical reference metric for this entire basin package.
Its outer coefficient ball, inner carrier locus, three physical witness
tiers, and \(\delta_{\rm RF}\) are those already fixed in
Corollary~\ref{cor:strict-entrance-nonempty}; none is recentered at a
nearby input.

For every nearby closed metric, use the fixed preparation convention
and marked chart to rebuild the auxiliary extension, adaptive target,
graft, and initial chart.  This does not replace the physical metric:
it only constructs its normalized prepared data.  The structural
equality with the target outside the graft then holds by construction.
Corollary~\ref{cor:fixed-convention-preparation-map} gives the resulting
\(C^1\), hence continuous, closed-metric-to-sliced-prepared map and
imposes the nine slice equalities.  Apart from the normalized and
finite physical-cover witnessed harmonic certificates, all other
entrance requirements are strict inequalities in finite jets or
weighted prepared norms; the last two sentences of that corollary
preserve a smaller normalized reserve and the unchanged fixed physical
\(+\)-tier.
Denote by
\(\lambda_{\rm ent}(G_0)\) the scale component of the resulting
phase-adjusted prepared tuple.  Shrink a full
$C^{k+2,\alpha}$ neighborhood
of $G_*$ so that every phase-adjusted actual carrier lies in this same
inner coefficient locus, all phase-adjusted states have common margins,
and
\[
 \lambda_{\rm ent}(G_0)\leq2\lambda_0
 \quad\text{for }G_0\in\mathscr U.
\]
  Theorem~\ref{thm:prepared-entrance-continuation}, used at continuation
  order \(k_0\), and
\eqref{eq:finite-interval-physical-width} give
\[
 T(G_0)
 \leq C_\lambda\lambda_{\rm ent}(G_0)
 \leq2C_\lambda\lambda_0<\varepsilon_T
\]
and all stated asymptotic conclusions in the prepared marked gauge.

For the exterior estimate, use the fixed sets \(W\Subset W^+\)
chosen above.
The two closed sets in this buffer have positive $G_*$-distance.
In the small-\(A\) choice above impose the requirement
\begin{equation}\label{eq:universal-exterior-buffer}
 d_{G_*}(W,\widehat X\setminus W^+)\geq20R_{\rm ext}.
\end{equation}
On $W^+$ the center metric is $\varpi^*g_X$, so it has a fixed
harmonic-radius and curvature-derivative package.  These bounds and
the weaker buffer $10R_{\rm ext}$ persist in a sufficiently small
$C^{k,\alpha}$ neighborhood.  The time margin in
\eqref{eq:universal-time-strict-margin} therefore allows
Lemma~\ref{lem:buffered-local-Ricci-control} to be applied with
\(R=R_{\rm ext}\) and \(J=0\).  It gives
\eqref{eq:universal-exterior-curvature}, uniformly in
\(G_0\in\mathscr U\).  This is only the \(J=0\) lifespan and need not
be uniform in \(J\).  Fix instead \(G_0\in\mathscr U\) and
\(K\Subset\widehat X\setminus\varpi^{-1}(B_{\rm imp})\), and choose
 strictly nested open sets
 \[
  K\Subset W_0\Subset W_1\Subset W .
 \]
 The curvature bound on the output set \(W\), together with the
 positive distance from \(W_1\) to \(\widehat X\setminus W\), and the
 metric distortion estimate following from
 \(\partial_tG=-2\Ric_G\) give uniform parabolic balls around \(W_1\)
 up to \(T(G_0)\).  The terminal local Shi estimates on these nested
 balls yield, for every \(j\geq0\) and every fixed
 \(0<d<T(G_0)\),
 \[
  \sup_{d\leq t<T(G_0)}
  \sup_{W_0}|\nabla^j\Rm_{G(t;G_0)}|
  \leq C_{j,W_0,W_1,d,G_0}<\infty .
 \]
 These constants need not be uniform in \(j\), and no higher-order
 lifespan is inferred from the \(J=0\) local lemma.  Fix
 \(d<T(G_0)\) and use \(\nabla^0=\nabla^{G(d;G_0)}\) on \(W_0\).
 Besides \(\partial_tG=-2\Ric_G\), the connection evolution is
 \[
  \partial_t(\nabla^{G(t)}-\nabla^0)
  =G^{-1}*\nabla^{G(t)}\Ric_{G(t)} .
 \]
 The Shi bounds therefore give, inductively in \(j\), uniform
 \((\nabla^0)^jG\), \((\nabla^0)^j\Ric_G\), and connection-difference
 bounds on a still smaller nested set containing \(K\).  Integrating
 \[
  \partial_t(\nabla^0)^jG=-2(\nabla^0)^j\Ric_G
 \]
 makes \(G(t;G_0)\) Cauchy in \(C^\infty(K)\) as
 \(t\uparrow T(G_0)\).  Exhausting the exterior by
 such \(K\) proves the claimed \(C^\infty_{\rm loc}\) terminal
 convergence and the localization assertion in the pointwise
 curvature-singular sense.
Apply Lemma~\ref{lem:one-state-terminal-continuity} on the sliced ball
\[
 \Sigma_{\tau_0}^{k_0+2,\alpha}
 =\Sigma_{\tau_0}^{k,\alpha}.
\]
It gives continuity of the singular time, limiting scale, phase, and
adaptive exterior data.
\end{proof}

\subsection{Positive-time smoothing and the
  \texorpdfstring{\(C^{2,\alpha}\)}{C2-alpha}-open basin}

\begin{lemma}[Gauge-fixed positive-time smoothing]
\label{lem:positive-time-smoothing}
Let $\mathcal N$ be closed, fix $d>0$, and let $\bar G(t)$ be a smooth Ricci flow on
$0\leq t\leq2d$, and fix $k\geq3$ and $0<\alpha<1$.  There is a
Banach-open neighborhood \(\widehat{\mathcal O}\) of \(\bar G(0)\) in
the positive cone of
\(h^{2,\alpha}(S^2T^*\mathcal N)\).  Put
\[
 \mathcal O:=
 \widehat{\mathcal O}\cap\operatorname{Met}^{\infty}(\mathcal N).
\]
For every \(G_0\in\widehat{\mathcal O}\), let
\(\widetilde G(t;G_0)\) be the unique Ricci--DeTurck solution on
\([0,d]\) with initial metric \(G_0\), taken relative to the
time-dependent background \(\bar G(t)\).  Then the gauge-fixed
positive-time map
\begin{equation}\label{eq:positive-time-smoothing-map}
 \widetilde{\mathfrak R}_d:
 \widehat{\mathcal O}
 \longrightarrow
 h^{k,\alpha}(S^2T^*\mathcal N),
 \qquad
 G_0\longmapsto\widetilde G(d;G_0),
\end{equation}
is \(C^1\).  Its restriction to \(\mathcal O\) is therefore locally
Lipschitz and continuous from the displayed relative
\(C^{2,\alpha}\) topology to \(C^{k,\alpha}\).  After shrinking
\(\widehat{\mathcal O}\) if necessary, for every integer \(m\geq0\)
the gauge-fixed solutions have a uniform spatial \(C^m\) bound on
\([d/2,d]\).

For every \(G_0\in\mathcal O\), there are smooth diffeomorphisms
\(\chi_t(G_0)\), with \(\chi_0(G_0)=\operatorname{Id}\), and the
corresponding smooth Ricci flow \(G(t;G_0)\) exists on \([0,d]\) and
satisfies
\begin{equation}\label{eq:DeTurck-to-Ricci-prefix}
 G(t;G_0)=\chi_t(G_0)^*\widetilde G(t;G_0).
\end{equation}
After shrinking \(\mathcal O\), the maps \(\chi_d(G_0)\) are uniformly
\(C^1\)-close to the identity, and the global curvature of these smooth
Ricci-flow prefixes is uniformly bounded on \([0,d]\).
\end{lemma}

\begin{proof}
Fix the time-dependent DeTurck gauge relative to $\bar G(t)$.  On a
small $C^{2,\alpha}$ neighborhood of $\bar G(0)$ the resulting
quasilinear system is uniformly parabolic on a common time interval.
We record the parameter spaces because the fixed positive-time gain
is stronger than bare well-posedness.  Put
\[
 \mathbb E_{\rm tr,0}=h^{0,\alpha}(S^2T^*\mathcal N),\qquad
 \mathbb E_{\rm tr,1}=h^{2,\alpha}(S^2T^*\mathcal N),
\]
and, on a short interval \(J=[s,s+\delta]\), use the initial-value
parabolic spaces
\[
 \mathbb E(J)
 =h^{1+\alpha/2,\,2+\alpha}(J\times\mathcal N),
 \qquad
 {}_0\mathbb E(J)=\{u\in\mathbb E(J):u(s)=0\},
 \qquad
 \mathbb F(J)=h^{\alpha/2,\,\alpha}(J\times\mathcal N).
\]
Here the anisotropic little-H\"older spaces are the closures of smooth
tensors in the standard parabolic norms, including the initial face;
the bounded trace map from \(\mathbb E(J)\) has target
\(\mathbb E_{\rm tr,1}\).  Fix once and for all a smooth connection
\(\nabla^\circ\) on \(S^2T^*\mathcal N\).  In these fixed
identifications the linearization has the form
\[
 L_u(t)v
 =a_u^{ij}(t,x)(\nabla^\circ)^2_{ij}v
  +b_u^i(t,x)\nabla_i^\circ v+c_u(t,x)v .
\]
Its domain is the common space
\(h^{2,\alpha}(S^2T^*\mathcal N)\), independently of \(u\).  The
principal coefficient is the scalar inverse-metric coefficient of the
Ricci--DeTurck system.  On a sufficiently small common
\(\mathbb E(J)\)-ball the coefficients are uniformly bounded in
\(h^{\alpha/2,\alpha}\), their dependence on \(u\) is \(C^1\) in
these coefficient norms, and there is one \(\lambda_{\rm par}>0\) such
that
\[
 a_u^{ij}(t,x)\xi_i\xi_j
 \geq\lambda_{\rm par}|\xi|^2\operatorname{Id}
\]
for every member of the ball.  The smooth time-dependent reference
metric \(\bar G(t)\) changes only the controlled lower-order
coefficients.  Since \(\mathcal N\) is closed, there are no lateral
boundary conditions or boundary compatibility requirements.

The nonautonomous Schauder theorem
\cite[Theorem~5.1.10]{Lunardi}, localized in one fixed finite atlas and
applied componentwise to the scalar principal part, gives
\[
 (\partial_t-L_u,\operatorname{tr}_s):
 \mathbb E(J)\longrightarrow
 \mathbb F(J)\times\mathbb E_{\rm tr,1}
\]
as an isomorphism.  The same assertion in the little-H\"older closures
follows by smooth approximation and the uniform Schauder estimate.  Its
inverse bound depends only on the fixed atlas, the displayed common
ellipticity and coefficient bounds, and \(\delta\), and is therefore
uniform in \(u\).  After extending the initial trace, the nonlinear
equation on \({}_0\mathbb E(J)\) is a \(C^1\) map with invertible
derivative.  The Banach implicit-function theorem therefore gives
\(C^1\) dependence on \(G_0\); differentiating the equation shows that
the derivative is the solution of the linearized time-dependent
Ricci--DeTurck equation.  A finite chain of such intervals covers
\([0,d]\).

For completeness, the topology gain at time \(d\) is obtained
uniformly, including for first variations.  Since the background is
the DeTurck solution with initial value \(\bar G(0)\), put
\(w=\widetilde G-\bar G\).  The mean-value formula turns the
difference equation into a linear uniformly parabolic system
\[
 \partial_tw=L_{G_0}(t)w
\]
in a fixed finite atlas.  On \([d/4,d]\), ordinary bootstrap from the
base \(\mathbb E\)-bound gives uniform coefficient bounds in every
order needed below.  Its evolution family \(U_{G_0}(t,s)\) satisfies
\begin{equation}\label{eq:positive-time-evolution-smoothing}
 \|U_{G_0}(t,s)\|_{h^{j,\alpha}\to h^{\ell,\alpha}}
 \leq C_{d,k}(t-s)^{-(\ell-j)/2},
 \qquad
 0\leq j\leq\ell\leq k,\quad d/4\leq s<t\leq d.
\end{equation}
The associated common-domain evolution-operator bounds are those of
\cite[Chapter~6, \S6.1, especially Corollary~6.1.8]{Lunardi}; their
constants are uniform for the same ellipticity and coefficient
package.  Equivalently here, the displayed estimate follows by the
interior parabolic Schauder estimate on nested time slabs, iterated
finitely many times.  Applying
\eqref{eq:positive-time-evolution-smoothing} from \(d/4\) to \(d\)
gives
\begin{equation}\label{eq:positive-time-fixed-smoothing-estimate}
 \|\widetilde G(d;G_0)-\bar G(d)\|_{C^{k,\alpha}}
 \leq C_{d,k}
 \|G_0-\bar G(0)\|_{C^{2,\alpha}}.
\end{equation}
The linearized solution satisfies the same evolution-family estimate.
We give the residual argument needed for the top H\"older topology;
boundedness in a stronger space alone is not being used to upgrade
lower-topology convergence.  Run the positive-time bootstrap on
\([d/2,d]\) one spatial order above the asserted target.  It gives,
uniformly on a smaller initial ball, high-order Lipschitz control of
the nonlinear solutions and uniform
\(h^{k+1,\alpha}\) bounds for unit-direction linearized solutions on
that slab.  Write the Ricci--DeTurck equation there as
\[
 \partial_tu=\mathscr Q(t,u),
 \qquad
 \mathscr Q(t,\cdot):h^{k+1,\alpha}\longrightarrow
 h^{k-1,\alpha},
\]
where \(\mathscr Q\) is \(C^1\) on the common positive-time ball.

Fix an initial metric \(G_0\), let \(h\to0\) in
\(h^{2,\alpha}\), put \(\eta=h/\|h\|_{h^{2,\alpha}}\), and set
\[
 q_h(t)=\frac{u(t;G_0+h)-u(t;G_0)}
                 {\|h\|_{h^{2,\alpha}}},
 \qquad
 v_\eta(t)=D_{G_0}u(t)[\eta],
 \qquad r_h=q_h-v_\eta.
\]
The base parabolic implicit-function theorem gives
\(r_h(d/2)\to0\) in \(h^{2,\alpha}\).  On \([d/2,d]\),
\[
 (\partial_t-D\mathscr Q(t,u(t;G_0)))r_h=\mathscr R_h,
\]
where the high-order Lipschitz bound and the Fr\'echet remainder for
\(\mathscr Q\) give
\[
 \|\mathscr R_h\|_
 {L^\infty([d/2,d];h^{k-1,\alpha})}\longrightarrow0.
\]
The homogeneous part is smoothed from \(h^{2,\alpha}\) to
\(h^{k,\alpha}\) across the fixed positive gap.  For the residual
source, the one-order evolution estimate gives
\[
 \|r_h(d)\|_{h^{k,\alpha}}
 \leq C_{d,k}\|r_h(d/2)\|_{h^{2,\alpha}}
 +C_{d,k}\int_{d/2}^d
   (d-s)^{-1/2}\|\mathscr R_h(s)\|_{h^{k-1,\alpha}}\,ds
 \longrightarrow0.
\]
Thus the difference quotients converge in the asserted top norm.

For continuity of the derivative, subtract the linearized equations
at two base metrics.  Besides the lower-topology initial difference,
the source is
\[
 \bigl[D\mathscr Q(t,u(t;G_0))
       -D\mathscr Q(t,u(t;G_0'))\bigr]v(t;G_0'),
\]
which tends to zero in
\(L^\infty([d/2,d];h^{k-1,\alpha})\), uniformly over unit initial
directions, by the same spare-order bounds.  The preceding estimate
therefore proves operator-norm continuity
\(h^{2,\alpha}\to h^{k,\alpha}\).  This establishes the \(C^1\)
version of \eqref{eq:positive-time-smoothing-map} without a compactness
shortcut.

The DeTurck vector field is a
first-order expression in $\widetilde G-\bar G$.  Its flow equation
therefore gives \eqref{eq:DeTurck-to-Ricci-prefix} and, from the
uniform $C^{2,\alpha}$ bound, $C^1$-closeness of $\chi_d$ to the
identity.  No gain of high coordinate regularity is asserted for
$\chi_d$: diffeomorphism modes are removed by the gauge rather than
smoothed by the ungauged Ricci flow.  Curvature is invariant under
pullback, so the gauge-fixed bounds also control the Ricci-flow
 prefix.  The \(h^{2,\alpha}\) Ricci--DeTurck well-posedness used at
 the first step is the \(k=0\) case of
 Bahuaud--Guenther--Isenberg~\cite[Theorem~6]{BahuaudGuentherIsenberg};
 after all derivatives are written relative to one fixed smooth
 connection, replacing their fixed reference metric by the smooth
 family \(\bar G(t)\) produces a nonautonomous quasilinear family with
 smoothly time-dependent controlled coefficients, while its principal
 symbol remains scalar and uniformly strongly elliptic.  The
  corresponding nonautonomous contraction therefore applies on the
  common short time interval by the coefficient and common-domain
  verification above; see \cite{GuentherIsenbergKnopf} for the
  little-H\"older Ricci--DeTurck realization and
  \cite[Theorem~5.1.10 and Chapter~6, \S6.1]{Lunardi} for the concrete
  nonautonomous Schauder and evolution-operator frameworks.  The DeTurck
pullback gives the Ricci-flow prefix, and the displayed fixed-time
estimate is the positive-time Schauder gain needed here.
\end{proof}

\begin{lemma}[Short-time persistence of strict preparation]
\label{lem:short-time-strict-preparation}
Let \(\mathbf z_{\tau_c}\) be a smooth strict prepared entrance of
continuation order \(k_0\geq12\), with all margins in
Definition~\ref{def:strict-prepared-entrance} measured in one fixed
prepared chart.  Let \(\mathbf z(\tau)\) be the local coupled feedback
evolution from Proposition~\ref{prop:coupled-local-feedback}, with the
same transported host, graft, scale, and initial-map convention.  Then
there is \(\eta>0\) such that \(\mathbf z(\tau)\) is a strict prepared
entrance, with positive (possibly smaller) common margins, for
\[
 \tau_c\leq\tau\leq\tau_c+\eta .
\]
Equivalently, after converting by \(dt/d\tau=\lambda>0\), the same
statement holds on a nontrivial physical-time interval.  The conclusion
also holds uniformly for a compact family whose initial strict margins
and prepared package constants are uniform.
\end{lemma}

\begin{proof}
Proposition~\ref{prop:coupled-local-feedback} gives continuity of the
finite-order prepared state in normalized time in the exact Banach
topology used by Definition~\ref{def:strict-prepared-entrance}.  The
weighted tensor norms, ellipticity constants, scale bracket, Gram
determinant, radial and inverse-map bounds, graft norms, and the
finitely many exterior coefficient and buffer quantities therefore
persist at \(\tau_c\).  For the normalized harmonic face, use the
cross-time graph clause of the reserve-to-operative form of
Lemma~\ref{lem:prepared-harmonic-radius-lower-stability}: continuity of
\(\bar g+h(\tau)\) in the global scale-one \(C^{2,\alpha}\) topology
compares it directly with \(\bar g+h(\tau_c)\), without identifying
the two time-typed prepared slices.  The recorded reserve triple at
\(\mathbf z_{\tau_c}\) therefore preserves the same operative
normalized triple and supplies a new strictly stronger reserve for the
restarted entrance.  On the finite physical cover, fixed-marking
continuity keeps the evolving actual carrier, for sufficiently small
\(\eta\), in the same inner locus
\eqref{eq:physical-coefficient-inner-locus} based at the unchanged
\(G_{\rm ref}^{\rm phys}\).  The reference-carrier certificate
therefore gives, at every \(x\in V_a^5\), the same common \(+\)-tier and
operative physical tier for the restarted carrier.  Since that carrier
lies in the outer ball, the precomputed half-open fixed-marking
Ricci-flow estimate \eqref{eq:witnessed-physical-quarter-modulus}
remains available with it as initial metric; no flow estimate is
reproved by continuity and no coefficient ball is recentered.
Continuity of the scale
together with the recorded \(\mu_{\rm RF}\)-slack preserves
the lifetime-width inequality after decreasing \(\eta\).  No
continuity of either bare supremal-radius functional is asserted.
Properness and
degree are unchanged under the resulting small \(C^1\) isotopy.  The
slice equalities, the graft identity, and the transported covariance
identities are preserved exactly by the coupled equations and hence
are not treated as open conditions.  Every remaining requirement is
one of the finite list of strict inequalities in the definition.
Taking the minimum of their positive initial margins and the two
witnessed harmonic thresholds proves the claim.  For a compact family,
first take a finite prepared-chart subcover carrying common positive
normalized and physical harmonic-witness reserves, then take the
minimum margin and the common local existence time over that subcover.
\end{proof}

\subsection{The invariant marked basin}

\begin{definition}[Marked FIK basin]
\label{def:marked-FIK-basin}
Let $\mathcal X$ be a closed four-manifold.  The \emph{marked FIK
basin} $\mathcal B_{\mathrm{FIK}}^{\mathrm{mark}}(\mathcal X)$ is the
set of smooth metrics $G_0$ whose maximal Ricci flow has a finite
singular time $T$, is globally Type I, meaning that for some
\(C<\infty\)
\[
 \sup_{\mathcal X}|\Rm_{G(t)}|
 \leq \frac{C}{T-t}
\]
for every sufficiently late \(t<T\), and admits exhausting FIK
markings $\Xi_t:\mathcal U_t\to\mathcal X$ with the following
property.  Here \emph{exhausting} means that each \(\mathcal U_t\) is
open, each \(\Xi_t\) is a smooth embedding, and for every
\(K\Subset M\) there is \(t_K<T\) such that
\(K\subset\mathcal U_t\) for all \(t_K<t<T\); nesting of the
\(\mathcal U_t\) is not required.  For every $t_i\uparrow T$, discard
the finitely many indices, if any, before the marking is defined on
the compact set under consideration.  With
$\delta_i=T-t_i$, the
base-time marking $\Xi_{t_i}$, frozen in the rescaled time variable,
satisfies
\begin{equation}\label{eq:marked-basin-convergence}
 \delta_i^{-1}\Xi_{t_i}^*G(T+s\delta_i)
 \longrightarrow g_{\mathrm{FIK}}(s)
 \quad\text{in }
 C^\infty_{\mathrm{loc}}(M\times(-\infty,0))
\end{equation}
without subsequence extraction.
\end{definition}

\begin{corollary}[$C^{2,\alpha}$ universal open basin]
\label{cor:low-topology-open-basin}
The construction in Theorem~\ref{thm:prepared-open-basin} may be
chosen so that its host, blow-up, and exact-core center admit a relative
$C^{2,\alpha}$-open neighborhood
\begin{equation}\label{eq:low-topology-neighborhood}
 \mathscr U_{2,\alpha}
 \subset\operatorname{Met}^{\infty}(\widehat X)
\end{equation}
of the same metric $G_*$ such that every
$G_0\in\mathscr U_{2,\alpha}$ has the singular-time, global Type-I,
exterior-curvature, and full-sequence marked-limit conclusions of
Theorem~\ref{thm:prepared-open-basin}.  After the fixed gauge restart
constructed below, its tail is generated by a strict prepared entrance
and satisfies, in the transported marking, all quantitative
conclusions of
Theorem~\ref{thm:prepared-entrance-continuation}.  In particular,
every such metric lies in
\(\mathcal B_{\mathrm{FIK}}^{\mathrm{mark}}(\widehat X)\).  The
auxiliary index \(k_0\geq12\) controls the post-smoothing continuation
argument; it does not enter the topology of the final basin.  No
optimality of this regularity threshold is claimed.
\end{corollary}

\begin{proof}
Fix the continuation output order \(k_0\geq12\) used in
Definition~\ref{def:strict-prepared-entrance}.  We first specify the
order of the geometric choices and then fix the final exact-core
center.  Let \(B_0\Subset B_{\rm imp}\) contain the
collapsing core, set
\[
 E_{\rm out}
 =\widehat X\setminus\varpi^{-1}(B_{\rm imp}),
\]
and choose fixed open sets
\[
 E_{\rm out}\subset W_-\Subset W\Subset W^+
 \Subset\widehat X\setminus\varpi^{-1}(\overline B_0).
\]
The exact-core construction has a uniform bounded-geometry package
on \(W^+\).  Choose \(R_d>0\) so that, before the final choice of
\(\tau_0\), this package includes a \(40R_d\) buffer, a harmonic-radius
lower bound at scale \(R_d\), and the required curvature jets.
In the final choice of \(\tau_0\) in
Proposition~\ref{prop:exact-core-implant}, require in addition that the
resulting center \(G_*\) satisfy
\begin{equation}\label{eq:restart-final-center-time-choice}
 T(G_*)<
 \min\left\{\frac{\varepsilon_T}{4},
             \frac18\delta_{\rm RF}R_d^2\right\}.
\end{equation}
This is possible because the center lifetime is
\(O(\lambda_0)\) and \(\lambda_0=Ae^{-\tau_0}\to0\).  The same
short-time bounded-geometry estimate makes the exterior metric
distortion as small as desired.  Freeze this final \(\tau_0\) and the
resulting center \(G_*\).

 Run that center from physical time zero.  Its entrance is strict, and
 the coupled feedback flow propagates the preparation data smoothly.
 Lemma~\ref{lem:short-time-strict-preparation} applies to this smooth
 exact-core entrance.
 Choose
\[
 0<d<\frac12T(G_*)
\]
 so small that the time-\(d\) state, with the transported preparation
 convention, is still a strict prepared entrance.  This choice is
 furnished by Lemma~\ref{lem:short-time-strict-preparation}, rather than
 being an additional persistence assumption.  Let
\(\tau_d=\tau(d;G_*)\) be its normalized-time label, and denote that
evolved prepared center by \(\mathbf z_d^*\).  By the choices above,
\begin{align}
 d_{G(d;G_*)}(W,\widehat X\setminus W^+)
 &\geq20R_d,
 \label{eq:restart-buffer-width}\\
 T(G_*)-d
 &<\frac14\delta_{\rm RF}R_d^2.
 \label{eq:restart-buffered-time-margin}
\end{align}
The time-\(d\) center also retains the harmonic-radius and
curvature-jet bounds on \(W^+\) at scale \(R_d\).  Its physical
certificate retains the original frozen reference
\(G_{\rm ref}^{\rm phys}=G_*\); the evolved carrier is an inner-locus
member, not a new reference center.

The error of \(\mathbf z_d^*\) need not vanish, but
the propagated receding slice gives
\[
 \ip{\rho_{\tau_d}h(\mathbf z_d^*)}{Z_\mu}=0,
 \qquad0\leq\mu\leq8,
\]
and its full phase-column matrix is the uniformly invertible receding
Gram matrix.  Corollary~\ref{cor:sliced-center-phase-map} therefore
provides a single centered phase convention on a neighborhood of
\(\mathbf z_d^*\), with phase \(p=0\) at the center.  The
entrance-openness argument in
Theorem~\ref{thm:prepared-open-basin}, applied with two additional
buffer derivatives, gives a \(C^{k_0+4,\alpha}\)-neighborhood
\(\mathscr V_d\) of \(G(d;G_*)\) whose members generate the full
continuation evolution and whose preparation has output order
\(k_0+2\).  Shrink it so that all of its actual restart carriers remain
in the same inner coefficient locus based at \(G_*\).  Write
\(\mathbf z_d(\widetilde G_d)\) for the resulting centered prepared
tuple.

For a prepared restart \(\mathbf z\) at \(\tau_d\), with its tail
physical clock set equal to zero at \(\tau_d\), define
\begin{equation}\label{eq:restart-tail-lifetime}
 T_{\rm tail}(\mathbf z)
 :=\int_{\tau_d}^{\infty}\lambda_{\mathbf z}(\tau)\,d\tau .
\end{equation}
The prepared continuation theorem identifies this with the maximal
tail lifetime, and
\(T_{\rm tail}(\mathbf z_d^*)=T(G_*)-d\).  The same-order continuity
of preparation and Lemma~\ref{lem:one-state-terminal-continuity} make
this quantity continuous on \(\mathscr V_d\).  Shrink
\(\mathscr V_d\), before taking any preimage, so that its metrics
retain the weaker \(10R_d\) buffer and the same bounded-geometry
package and, for every \(\widetilde G_d\in\mathscr V_d\),
\begin{equation}\label{eq:restart-tail-time-margins}
 T_{\rm tail}\bigl(\mathbf z_d(\widetilde G_d)\bigr)
 <\frac12\delta_{\rm RF}R_d^2,
 \qquad
 d+T_{\rm tail}\bigl(\mathbf z_d(\widetilde G_d)\bigr)
 <\varepsilon_T .
\end{equation}
Lemma~\ref{lem:buffered-local-Ricci-control}, restarted at time $d$, then
gives a uniform curvature bound on $W$ throughout every such tail.

Now apply Lemma~\ref{lem:positive-time-smoothing} with target order
\(k_0+4\), relative to the center flow, and, after shrinking the domain
of the gauge-fixed map, set
\begin{equation}\label{eq:gauge-fixed-low-topology-basin}
 \widehat{\mathscr U}_{2,\alpha}
 :=\widetilde{\mathfrak R}_d^{-1}(\mathscr V_d)
 \subset\operatorname{Met}^{2,\alpha}(\widehat X).
\end{equation}
Here \(\mathscr V_d\) is read in the fixed DeTurck chart at time \(d\).
For the center solution,
\[
 \widetilde{\mathfrak R}_d(G_*)=\bar G(d)=G(d;G_*),
\]
so \(\widehat{\mathscr U}_{2,\alpha}\) contains \(G_*\) and is open in
the little-H\"older \(h^{2,\alpha}\) positive cone.  Define its smooth
locus by
\begin{equation}\label{eq:smooth-low-topology-basin}
 \mathscr U_{2,\alpha}
 :=\widehat{\mathscr U}_{2,\alpha}
   \cap\operatorname{Met}^{\infty}(\widehat X).
\end{equation}
It is exactly a relative \(C^{2,\alpha}\)-open neighborhood in the
topology asserted in the corollary.

For \(G_0\in\widehat{\mathscr U}_{2,\alpha}\), put
\[
 \widetilde G_d=\widetilde{\mathfrak R}_d(G_0),
\]
reset the additive normalized-time constant to the same value
\(\tau_d\) for every nearby restart, and let
\(\mathbf z_d(\widetilde G_d)\) be the uniquely prepared tuple in the
transported convention.  Thus
\begin{equation}\label{eq:low-preparation-restart-map}
 \mathcal P_d:
 \widehat{\mathscr U}_{2,\alpha}\ni G_0\longmapsto
 \mathbf z_d\bigl(\widetilde{\mathfrak R}_d(G_0)\bigr)
 \in\Sigma_{\tau_d}^{k_0+2,\alpha}.
\end{equation}
It is \(C^1\) between the corresponding little-H\"older charts:
positive-time smoothing supplies order \(k_0+4\), and
Corollary~\ref{cor:sliced-center-phase-map} uses exactly two derivatives
to produce the order-\(k_0+2\) prepared tuple.

Now restrict to \(G_0\in\mathscr U_{2,\alpha}\), put
\(\chi_d=\chi_d(G_0)\), and let \(\widetilde H(s)\) be the Ricci flow
from \(\widetilde G_d\).  By
\eqref{eq:DeTurck-to-Ricci-prefix} and uniqueness,
\begin{equation}\label{eq:Ricci-tail-fixed-pullback}
 G(d+s;G_0)=\chi_d^*\widetilde H(s)
\end{equation}
throughout their common lifetime.  Thus diffeomorphism invariance and
saturation under a finite regular prefix transfer the Type-I and
marked-convergence conclusions already proved above to the original
Ricci flow.  The stable-profile, quadratic, and sharp
marked-spacetime refinements are proved later and transferred by the
same identity only after they have been established.  If
$\widetilde\Xi_s$ is a tail marking, the corresponding physical marking is
$\chi_d^{-1}\circ\widetilde\Xi_s$.

The terminal maps $\chi_d$ are uniformly $C^1$-close to the identity.
Shrink the low-topology neighborhood once more so that
\begin{equation}\label{eq:terminal-gauge-exterior-inclusion}
 \chi_d(E_{\rm out})\subset W_-.
\end{equation}
Equations~\eqref{eq:Ricci-tail-fixed-pullback} and
\eqref{eq:terminal-gauge-exterior-inclusion} transfer the restarted
curvature bound on $W$ to the prescribed physical exterior
$E_{\rm out}$.  Lemma~\ref{lem:positive-time-smoothing} controls the
regular prefix.  The second inequality in
\eqref{eq:restart-tail-time-margins} says exactly that the original
singular time \(d+T_{\rm tail}\) is less than \(\varepsilon_T\).

\end{proof}

\begin{theorem}[Invariant marked basin with universal open interior]
\label{thm:marked-FIK-basin}
For every closed four-manifold $\mathcal X$, the marked FIK basin is
invariant under diffeomorphism pullback, positive scaling, and the
forward Ricci semiflow.  It is also saturated under finite regular
Ricci-flow prefixes: if $G_0$ flows smoothly to a metric in the basin,
then $G_0$ is itself in the basin.  For every oriented host $(X,g_X)$
and every oriented blow-up $\widehat X$ furnished by
Theorem~\ref{thm:prepared-open-basin},
\begin{equation}\label{eq:universal-open-interior}
 \mathscr U_{2,\alpha}\subset
 \mathcal B_{\mathrm{FIK}}^{\mathrm{mark}}(\widehat X).
\end{equation}
Consequently the marked basin on $\widehat X$ has nonempty relative
$C^{2,\alpha}$ interior.
\end{theorem}

\begin{proof}
Let $\mathfrak R_t(G_0)=G(t;G_0)$.  If
$G_0\in\mathcal B_{\mathrm{FIK}}^{\mathrm{mark}}(\mathcal X)$ and
$0\leq s_0<T(G_0)$, then the flow from
$\mathfrak R_{s_0}(G_0)$ is $G(s_0+t;G_0)$ and has singular time
$T(G_0)-s_0$.  Set
\(\mathcal U'_t=\mathcal U_{s_0+t}\) and
\(\Xi'_t=\Xi_{s_0+t}\).  Base times and their distances to the singular
time are unchanged after this translation, so the Type-I bounds and
\eqref{eq:marked-basin-convergence} are unchanged.  This proves
forward invariance.

Conversely, suppose $G_0$ has a smooth Ricci flow on $[0,s_0]$ and
$G(s_0;G_0)$ lies in the marked basin.  Concatenate this regular prefix
with the basin flow.  Uniqueness identifies the concatenation with the
maximal flow from $G_0$; its singular time is shifted by $s_0$, while
all sufficiently late Type-I estimates and markings are unchanged.
Thus $G_0$ lies in the basin.  This proves saturation under regular
prefixes.

For $\vartheta\in\operatorname{Diff}(\mathcal X)$,
\[
 \mathfrak R_t(\vartheta^*G_0)=\vartheta^*\mathfrak R_t(G_0).
\]
Transport the markings by
$\Xi_t\mapsto\vartheta^{-1}\circ\Xi_t$; then every marked pullback is
identical to the original one.  For $c>0$,
\[
 \mathfrak R_t(cG_0)=c\,\mathfrak R_{t/c}(G_0).
\]
For the scaled flow set
\(\mathcal U'_t=\mathcal U_{t/c}\) and
\(\Xi'_t=\Xi_{t/c}\); its singular time is
$cT$.  If
$t_i'=ct_i$ and $\delta_i'=c\delta_i$, then
\[
 (\delta_i')^{-1}\Xi_{t_i}^*
 \bigl(cG(T+s\delta_i)\bigr)
 =
 \delta_i^{-1}\Xi_{t_i}^*G(T+s\delta_i).
\]
Thus the Type-I and marked-limit properties are invariant under
positive scaling.  The inverse diffeomorphism and inverse scaling give
equality of the corresponding basins.

  Finally, Corollary~\ref{cor:low-topology-open-basin} gives
  \eqref{eq:universal-open-interior}.  That corollary makes
$\mathscr U_{2,\alpha}$ relatively $C^{2,\alpha}$-open and supplies
the quantitative marked-limit package.
\end{proof}

\begin{proof}[Proof of Theorem~\ref{thm:intro-open-basin}]
Theorem~\ref{thm:prepared-open-basin} constructs the exact-core center
with the prescribed exterior, arbitrarily small singular time, global
Type-I control, and a uniform exterior curvature bound on a strict
high-regularity neighborhood.  Corollary~\ref{cor:low-topology-open-basin}
promotes that neighborhood to the relative \(C^{2,\alpha}\) topology
without changing the center or any of those conclusions.
Theorem~\ref{thm:marked-FIK-basin} supplies the invariant
full-sequence frozen-marking formulation.  The uniformity statements
in Theorem~\ref{thm:prepared-open-basin} and
Corollary~\ref{cor:low-topology-open-basin}, together with the
per-flow marking formulation in
Theorem~\ref{thm:marked-FIK-basin}, give precisely the quantifiers
asserted in Theorem~\ref{thm:intro-open-basin}.  Since the domains of
the marked convergence exhaust \(M\), they contain a fixed compact
tubular neighborhood of the zero section \(E\) for all sufficiently
late times.  Put
\(e_t=\Xi_t|_E\).  Restriction of
\eqref{eq:marked-pullback-identity} gives
\[
 e_t^*G(t)=\lambda(t)(\bar g+h(t))|_E .
\]
The core estimate \eqref{eq:master-core-rate}, the normalized scale
rate \eqref{eq:master-scale-rate}, and the physical scale law
\eqref{eq:master-physical-power} therefore prove
\eqref{eq:intro-marked-bolt-collapse}, with
\(\vartheta=\theta\) after decreasing the exponent if necessary.
The area asymptotic follows from smooth dependence of the volume
density on the rescaled induced metric.  Uniform bilipschitz comparison
with \(g_E\) gives the intrinsic-diameter asymptotic.  For the ambient
curvature law, use the full identity
\eqref{eq:marked-pullback-identity} on the fixed tubular neighborhood,
the \(C^2\) case of \eqref{eq:master-core-rate}, naturality of
curvature, and
\(\lvert\Rm_{c g}\rvert_{c g}=c^{-1}\lvert\Rm_g\rvert_g\).
Let
\[
 E_{\rm exc}:=\iota^{-1}(E)\subset\widehat X''.
\]
By the oriented blow-up construction in
Proposition~\ref{prop:exact-core-implant}, \(E_{\rm exc}\) is the
inserted exceptional zero section and has self-intersection \(-1\).
The prepared maps \(R_\tau\) and \(F_\tau\) lie in the degree-one
proper homotopy class fixed in the prepared-state convention, while
\(\varphi_\tau\) is a flow.  Hence
\[
 \Phi_\tau=\varphi_\tau\circ R_\tau\circ F_\tau
\]
is orientation-preserving.  On the implanted neighborhood define
\[
 \mathcal D_\tau^{\rm exc}
 :=\iota^{-1}\circ\Phi_\tau^{-1}\circ\iota .
\]
Conjugation preserves the orientation degree of a diffeomorphism.
Therefore \(\mathcal D_\tau^{\rm exc}\) is orientation-preserving;
this conclusion is
independent of any separate orientation convention for \(\iota\).
Moreover,
\[
 \Sigma_t=\Xi_t(E)=\mathcal D_{\tau(t)}^{\rm exc}(E_{\rm exc}).
\]
The differential of \(\mathcal D_{\tau(t)}^{\rm exc}\) induces an
oriented quotient-bundle isomorphism
\[
 N_{E_{\rm exc}/\widehat X}
 \cong
 \bigl(\mathcal D_{\tau(t)}^{\rm exc}|_{E_{\rm exc}}\bigr)^*
 N_{\Sigma_t/\widehat X}.
\]
Naturality of the Euler class and evaluation on the transported
fundamental class therefore give
\[
 \deg N_{\Sigma_t/\widehat X}
 =\deg N_{E_{\rm exc}/\widehat X}=-1.
\]
The ambient-diameter bound follows from the intrinsic one.
\end{proof}

\part{Nonlinear scattering and the physical first-profile foliation}

\section{Global two-state stability and asymptotic data}
\label{sec:global-two-state}

\paragraph{Roadmap for the global two-state argument.}
Theorem~C below first fixes the exact prepared domains and
the componentwise scattering target.  The proof then proceeds through
four blocks: uniform low-order hybrid propagation; coarse geometry,
map, and graft comparison; exact Kato algebra and two-state barriers
with derivative recovery; and the global continuation estimate
Theorem~\ref{thm:global-two-state-estimate}.  The terminal Lipschitz
corollary assembles the singular time, scale, phase, and stable-profile
data.  Differentiability and physical realization are deferred to the
later first-variation and transverse-disk subsections; neither is used
in the global estimate.

\begin{maintheorem}[Prepared continuation, nonlinear scattering, and
strong-stable profile foliation]
\label{thm:intro-sharp-scattering}

\emph{I. Sharp scattering on an arbitrary fixed strict-entrance ball.}
Fix \(k\geq12\), \(0<\alpha<1\), a normalized entrance time
\(\tau_0\), and first fix a rate pair
\[
0<\sigma<\theta<\beta .
\]
Next fix one numerical prepared package compatible with this pair, one
marked host, graft, scale, and initial-map convention, and finally one
entrance size
\[
 0<\varepsilon\leq\varepsilon_{\rm ent}.
\]
Let
\[
 \mathscr B\subset\Sigma_{\tau_0}^{k+2,\alpha}
\]
be a fixed common-margin sliced ball all of whose points are strict
prepared entrances of continuation order \(k\) in the exact sense of
Definition~\ref{def:strict-prepared-entrance}, with the displayed
parameters and package common on the ball.  In particular, the common
physical data include one fixed \(G_{\rm ref}^{\rm phys}\), one outer
coefficient ball, one inner carrier locus with a common
\(\mu_{\rm coeff}^{\rm phys}\)-margin, one three-tier physical witness
package, and one \(\delta_{\rm RF}\).  The actual carriers of two states
may differ, but both lie in that same inner locus and neither recenters
the certificate.  Assume also that one buffered physical cover and the
strict time-width margin
\eqref{eq:auxiliary-buffered-time-width} work throughout
\(\mathscr B\).  Let
\[
 \mathscr O\subset
 \operatorname{dom}\Pi_{\rm sl}^{\,k+4\to k+2}
 \subset\mathscr P_{\tau_0}^{k+4,\alpha}
\]
be an open buffered ambient domain such that
\[
 \overline{\mathscr O}\subset
 \operatorname{dom}\Pi_{\rm sl}^{\,k+4\to k+2},
 \qquad
 \Pi_{\rm sl}^{\,k+4\to k+2}(\mathscr O)\subset\mathscr B,
\]
and, for some uniformly interior common-margin subball
\(\mathscr B'\Subset_{\rm u}\mathscr B\),
\[
 \Pi_{\rm sl}^{\,k+4\to k+2}(\overline{\mathscr O})
 \subset\mathscr B'.
\]
Use the quantitative retraction fixed in
Proposition~\ref{prop:sliced-prepared-manifold}, and write
\(K_{\Pi,k+2}\) for the derived constant in
\eqref{eq:quantitative-phase-retraction-Lipschitz}--%
\eqref{eq:quantitative-phase-retraction-derivative}.  The sliced
constants below are the uniform constants on \(\mathscr B'\) furnished
by Theorems~\ref{thm:global-two-state-estimate} and
\ref{thm:global-two-state-scattering}.  Hence every ambient Lipschitz
or first-variation constant is bounded by the corresponding sliced
constant times \(K_{\Pi,k+2}\); the closure and interiority assumptions
retain the fixed charts and common margins.
Then every point of \(\mathscr B\) has the global continuation of
Theorem~\ref{thm:prepared-entrance-continuation}.  All estimates below
are uniform on the fixed ball and parameter package, and every
two-state assertion has a buffered ambient extension on \(\mathscr O\)
after composition with
\(\Pi_{\rm sl}^{\,k+4\to k+2}\).

In particular, for every \(\mathbf z_0\in\mathscr B\) there is a unique
\(V_\infty(\mathbf z_0)\in E_1\), and there is
\(\delta_{\rm prof}>0\), uniform on \(\mathscr B\), such that
\[
 H(\tau)
 =
 e^{-\gamma_1\tau}V_\infty
 +O_{H^1_\nu}\!\left(
   e^{-(\gamma_1+\delta_{\rm prof})\tau}\right).
\]
The marked scattering map
\[
 \mathbf z_0\longmapsto
 (T,\log\lambda_\infty,\Psi_\infty,V_\infty)
\]
is \(C^1\) on \(\mathscr B\), with buffered ambient extension on
\(\mathscr O\).  For smooth entrances the frozen base-time markings
satisfy the sharp Jacobi expansion of
Theorem~\ref{thm:sharp-marked-spacetime}; the scale and phase satisfy
the response laws of
Theorem~\ref{thm:quadratic-geometric-asymptotics}; and the combination
\(\mathfrak A_1=\lambda_\infty^{-\gamma_1}V_\infty\) has the
transported-restart-invariant physical normal form of
Corollary~\ref{cor:physical-amplitude-normal-form}.  These are
marked, gauge-fixed statements in the single transported convention;
no quotient by arbitrary re-markings is asserted.

\emph{II. Existential exact-core realization and foliation.}
Independently of the arbitrary ball and entrance time fixed in
Part~I, first fix \(k\geq12\), \(0<\alpha<1\), and a rate pair.
For the exact-core construction in Theorem~A there then exist a
rate- and construction-compatible prepared package, an entrance size
below its resulting \(\varepsilon_{\rm ent}\), and a sufficiently late
entrance time \(\widehat\tau_0\) such that
Theorem~\ref{thm:profile-realization} supplies a relative
\(C^{2,\alpha}\)-open basin \(\mathscr U_{\rm prof}\), a compactly
supported transverse physical disk, a common-margin sliced
strict-entrance ball
\[
 \mathscr B_{\rm scat}
 \subset\Sigma_{\widehat\tau_0}^{k+2,\alpha},
\]
on which the fixed exact-core buffered physical cover and the strict
time-width margin \eqref{eq:auxiliary-buffered-time-width} work
uniformly, a uniformly interior common-margin subball
\[
 \mathscr B'_{\rm scat}\Subset_{\rm u}\mathscr B_{\rm scat},
\]
and a buffered prepared neighborhood
\[
 \widetilde{\mathscr O}_{\rm scat}
 \subset\mathscr P_{\widehat\tau_0}^{k+4,\alpha}
\]
satisfying
\[
 \overline{\widetilde{\mathscr O}_{\rm scat}}
 \subset
 \operatorname{dom}\Pi_{\rm sl}^{\,k+4\to k+2},
 \qquad
 \Pi_{\rm sl}^{\,k+4\to k+2}
   (\overline{\widetilde{\mathscr O}_{\rm scat}})
 \subset\mathscr B'_{\rm scat}
 \Subset_{\rm u}\mathscr B_{\rm scat}.
\]
On this specially constructed domain, and only with this existential
quantification, the phase-retracted profile map
\[
 \mathscr S^{\rm prep}_1
 =
 V_\infty^{\rm sl}\circ
 \Pi_{\rm sl}^{\,k+4\to k+2}:
 \widetilde{\mathscr O}_{\rm scat}\longrightarrow E_1
\]
is a split \(C^1\) submersion.  Its connected fibers form the asserted
local \(C^1\) foliation, and
\[
 \mathscr W^{ss}_1
 =(\mathscr S^{\rm prep}_1)^{-1}(0)
\]
is a nonempty strong-stable submanifold.  By
Theorem~\ref{thm:profile-realization}, every sufficiently small marked
profile has a unique compactly supported transverse realization.  No
split-submersion, profile-realization, or nonempty-zero-leaf conclusion
is asserted for an arbitrary pair \((\mathscr B,\mathscr O)\) from
Part~I.
\end{maintheorem}

The proof of Theorem~\ref{thm:intro-sharp-scattering} is completed in
Section~\ref{sec:first-stable-profile}, after the one-state profile,
quadratic response, global two-state scattering, and
profile-realization theorems on which its two parts depend.

Throughout this section fix an integer \(k\geq12\) and
\(0<\alpha<1\).
The finite-horizon same-order difference and first-variation theory
used below is established in
Section~\ref{sec:finite-horizon-two-state}.

\begin{lemma}[Uniform sixth-order harmonic-map coefficient bound]
\label{lem:uniform-HMHF-C6-coefficients}
Let a global trajectory arise from a strict prepared entrance in one
fixed common-margin ball of input order \(k+2\), \(k\geq12\).  Then
\begin{equation}\label{eq:uniform-HMHF-C6-coefficients}
 \sup_{\tau\geq\tau_0}
 \|F(\tau)\|_{\operatorname{Map}_{\rm sc}^{6,\alpha}}
 \leq C ,
\end{equation}
where the typed norm includes both the map and inverse-map coordinate
representatives.  The constant is uniform on the common-margin ball
and is obtained independently of
Lemma~\ref{lem:uniform-scale-adapted-C5-bridge}.
\end{lemma}

\begin{proof}
Fix \(T>\tau_0\) and restrict the global strict trajectory to
\([\tau_0,T)\).  The physical coefficient package, target tracking,
coarse Gram and target-time package, and strict harmonic-map faces
established in the global continuation argument are precisely the
hypotheses of Lemma~\ref{lem:finite-HMHF-C6-bridge}.  Its entrance
\(C^{6,\alpha}\) norm is the fixed prepared entrance norm, and its
constant is independent of the finite terminal endpoint \(T\).
Consequently
\[
 \sup_{\tau_0\leq\tau<T}
 \|F(\tau)\|_{\operatorname{Map}_{\rm sc}^{6,\alpha}}\leq C
\]
with the same \(C\) for every \(T\).  Letting \(T\to\infty\) proves
\eqref{eq:uniform-HMHF-C6-coefficients}.
\end{proof}

\begin{lemma}[Uniform scale-adapted five-derivative bridge]
\label{lem:uniform-scale-adapted-C5-bridge}
Let a global trajectory arise from a strict prepared entrance in one
fixed common-margin ball of input order \(k+2\), \(k\geq12\).  Then,
for the weight exponent \(N\) of that ball,
\begin{equation}\label{eq:uniform-scale-adapted-C5-bridge}
 \sup_{\tau\geq\tau_0}
 \|h(\tau)\|_{\mathfrak T_{{\rm sc},N}^{5,\alpha}}
 \leq C .
\end{equation}
The constant is uniform on the common-margin ball.  The conclusion
covers, in particular, the expanding normalized inner--middle region
\[
 4\Gamma\lesssim\bar f\lesssim c\Gamma e^\tau;
\]
it is not obtained by treating that region as part of a fixed compact
core.
\end{lemma}

\begin{proof}
Put
\[
 \mathfrak q_\tau=(F_\tau^{-1})^*\acute G_\tau .
\]
Fix \(\tau\geq\tau_0\).  Apply the instantaneous construction in
Lemma~\ref{lem:C6-persistent-source-atlas} at time \(\tau\)
(equivalently, declare \(\tau\) to be a refresh time), and choose one
paired member with source buffers \(\mathcal U^j\), target buffers
\(\mathcal V^j\), and physical harmonic scale
\(r=r_{\mathcal U}\).  The endpoint-uniform source and target packages
give
\begin{equation}\label{eq:C5-algebraic-source-target-package}
 \|r^{-2}\acute G_\tau\|_
   {C^{12,\alpha}(\mathcal U^3)}
 +
 \|r^{-2}S_\tau\|_
   {C^{12,\alpha}(\mathcal V^3)}
 \leq C .
\end{equation}
The typed map-and-inverse estimate
\eqref{eq:uniform-HMHF-C6-coefficients}, the fixed forward and inverse
buffer containments
\eqref{eq:C6-paired-entry-containment}, and the coordinate
identity
\[
 r^{-2}\mathfrak q_\tau
 =(F_\tau^{-1})^*(r^{-2}\acute G_\tau)
\]
therefore give, by the buffered pullback estimate with one map
derivative to spare,
\begin{equation}\label{eq:C5-algebraic-target-metric-bound}
 \|r^{-2}(\mathfrak q_\tau-S_\tau)\|_
 {C^{5,\alpha}(\mathcal V^1,r^{-2}S_\tau)}
 \leq C .
\end{equation}
No smallness of this norm is asserted or used.

Set
\[
 \ell=\frac{r^2}{\lambda(\tau)}.
\]
The exact relative-metric identity
\eqref{eq:relative-metric-identity} gives on \(\mathcal V^1\)
\begin{equation}\label{eq:C5-algebraic-isometric-transfer}
 \Theta_\tau^*(\ell^{-1}h_\tau)
 =r^{-2}(\mathfrak q_\tau-S_\tau),\qquad
 \Theta_\tau^*(\ell^{-1}\bar g)
 =r^{-2}S_\tau .
\end{equation}
Thus \(\Theta_\tau\) is an isometry between the two scaled chart
metrics in \eqref{eq:C5-algebraic-isometric-transfer}.  In particular,
no coordinate derivative of \(\Theta_\tau\) is used: the
\(S_\tau\)-harmonic target coordinates become harmonic coordinates for
\(\ell^{-1}\bar g\) after composition with \(\Theta_\tau^{-1}\), since
constant metric scaling does not change harmonic functions.  Hence
\eqref{eq:C5-algebraic-target-metric-bound} implies
\begin{equation}\label{eq:C5-algebraic-fixed-background-bound}
 \|\ell^{-1}h_\tau\|_
 {C^{5,\alpha}(\Theta_\tau(\mathcal V^1),
                \ell^{-1}\bar g)}
 \leq C .
\end{equation}

By the output-scale clause
\eqref{eq:C6-target-output-scale},
\(\Theta_\tau(\mathcal V^0)\) is a core member or lies in a fixed
enlargement of a dyadic annulus \(A_L\), and \(L\asymp\ell\).
These inner members cover \(M\) with uniform overlap and Lebesgue
number.  Constant-rescaling equivalence in
\eqref{eq:C5-algebraic-fixed-background-bound} and the fixed FIK
dyadic atlas therefore yield the stronger estimate
\begin{equation}\label{eq:C5-strong-unweighted-bridge}
 \sup_{\tau\geq\tau_0}
 \|h(\tau)\|_{\mathfrak C_{{\rm sc},0}^{5,\alpha}}
 \leq C .
\end{equation}
Since \(N\geq0\) and every exterior label satisfies
\(L\geq\Gamma\geq1\), this controls both
\(\mathfrak C_{{\rm sc},0}^{2,\alpha}\) and
\(\mathfrak C_{{\rm sc},N}^{5,\alpha}\), proving
\eqref{eq:uniform-scale-adapted-C5-bridge}.  The proof is algebraic at
each time and depends only on the common source, target, paired-atlas,
and typed map packages.  It uses neither a transported-tube damping
estimate nor a graft-source H\"older regularization, and its constant
is uniform on the common-margin ball.
\end{proof}

\begin{remark}[Coefficient bridge]
\label{rem:acyclic-C6-C5-ledger}
The dependency just proved is
\[
 \begin{gathered}
 \text{physical/target \(C^{12,\alpha}\) packages}
 \Longrightarrow
 \text{paired-atlas \(C^{6,\alpha}\) map bridge}\\
 {}\Longrightarrow
 \text{algebraic \(C^{5,\alpha}\) tensor bridge}.
 \end{gathered}
\]
The first arrow is Lemma~\ref{lem:finite-HMHF-C6-bridge}, including
\eqref{eq:C6-remote-q2-absorption}; it uses only the raw \(C^2\) box
for \(h\), not the tensor bridge.  The second arrow is
\eqref{eq:C5-algebraic-target-metric-bound}--%
\eqref{eq:C5-strong-unweighted-bridge}.  Consequently the displayed
chain supplies the coefficient package for
Lemma~\ref{lem:uniform-weighted-Schauder-restart}, and hence for the
two-state and first-variation arguments that invoke it.
\end{remark}

\begin{lemma}[Uniform future low-order hybrid propagation]
\label{lem:uniform-weighted-Schauder-restart}
For this entrance order put
\[
 m_\#=4.
\]
Let \(\mathscr B\subset\Sigma_{\tau_0}^{k+2,\alpha}\) be a
common-margin sliced ball of strict prepared entrances for which the
fixed ordinary exterior certificate
\[
 K_-\Subset K_0\Subset K_1\Subset K_2,\qquad
 U_a\Subset U_a^+\Subset U_a^{++}\Subset V_a^0
 \Subset\cdots\Subset V_a^5\Subset E^{++}
\]
and the cutoff, atlas, and strict time-width margins in
\eqref{eq:prepared-exterior-termination} and
\eqref{eq:auxiliary-buffered-time-width} are common.  Every global
trajectory furnished by
Theorem~\ref{thm:prepared-entrance-continuation} has the following
uniform low-order coefficient package:
\begin{equation}\label{eq:uniform-future-low-coefficient-package}
 \begin{split}
 \mathfrak K_\#(\mathbf z;\tau):={}&
 \sum_{a=1}^{N_{\rm ext}}
  \|\widetilde G(t(\tau))\|_{C_{R_a}^{9,\alpha}(U_a^{++})}\\
 &+\sum_{a=1}^{N_{\rm ext}}
  \|\chi^{\pm1}(t(\tau))\|_{C_{R_a}^{8,\alpha}(U_a^{++})}\\
 &+\|\widetilde\iota(t(\tau))\|_
       {C^{8,\alpha}(\mathcal W_{\rm in}^{++})}
  +\|\widetilde\iota(t(\tau))\|_
       {C^{8,\alpha}(\mathcal W_{\rm gr}^{++})}\\
 &+\|R\|_{\operatorname{Map}_{\rm sc}^{7,\alpha}}
  +\|F\|_{\operatorname{Map}_{\rm sc}^{6,\alpha}}\\
 &+\|h\|_{\mathfrak T_{{\rm sc},N}^{5,\alpha}}
 \leq K_\# .
 \end{split}
\end{equation}
Each displayed \(\operatorname{Map}_{\rm sc}\) norm includes the map
and inverse-map representatives in the corresponding prepared charts.  The
constant \(K_\#\), as well as all ellipticity, radial-comparison, and
composition margins needed at these orders, is uniform for
\(\mathbf z_0\in\mathscr B\) and \(\tau\geq\tau_0\).

Let two such trajectories start at
\(\mathbf z_{1,0},\mathbf z_{2,0}\in\mathscr B\), and measure their
difference by
\(\mathfrak D_{m_\#}^{\rm hyb}\) from
\eqref{eq:two-state-hybrid-distance}, including the marking block
\eqref{eq:hybrid-marking-block}, the compact graft input buffer
\eqref{eq:hybrid-graft-buffer-block}, and the low-order interface block
\eqref{eq:hybrid-interface-block}.  Its \(F\)-summand is the localized
block \eqref{eq:localized-graft-F-block}, not the global finite-horizon
\(F\)-norm.  Put
\[
 d_0=\|\mathbf z_{1,0}-\mathbf z_{2,0}\|
       _{\mathscr X_{\rm prep}^{k+2,\alpha}}.
\]
For
\(I_s=[s,s+1]\), \(s\geq\tau_0\),
put \(H_i=\rho_\tau h_i\), \(i=1,2\).  Then
\[
 \begin{aligned}
 d_{\#,0}:={}&
 \mathfrak D_{m_\#}^{\rm hyb}(\tau_0)
 +d_{{\rm ext},m_\#,0}^{+}
 +d_{{\rm ext},-1,0}^{\rm corr}\\
  &+d_{{\rm Ggr},m_\#+2,0}^{0}
  +d_{{\rm Fgr},m_\#+1,0}^{++}
  +d_{{\rm F},m_\#+1,0}^{\rm glob}\\
  &+\|R_1(\tau_0)-R_2(\tau_0)\|_
        {\mathfrak X_{\rm sc}^{m_\#+2,\alpha}}\\
  &+\|H_1(\tau_0)-H_2(\tau_0)\|_{L^2_\nu}.
 \end{aligned}
\]
where \(d_{{\rm ext},m_\#,0}^{+}\) and
\(d_{{\rm ext},-1,0}^{\rm corr}\) are the buffered and separated
corridor pieces of the homogeneous exterior face, and
\(d_{{\rm Ggr},m_\#+2,0}^{0}\) is the larger-collar metric trace used
by the graft localization cutoff.  These three quantities are
original-\(\tau_0\) memories and are not inferred from the smaller
future output blocks.  The displayed order-\((m_\#+2)\) \(R\)-trace
is the auxiliary target-connection input for the Abel estimate; it is
controlled by the prepared map component of \(d_0\).
Denote by
\[
 \begin{split}
 \mathscr L_\#(s):={}&
 \sup_{\tau\in I_s}
 \left(
  \mathfrak D_{m_\#}^{\rm hyb}(\tau)
  +\mathfrak F_{m_\#+1}^{\rm glob}(\tau)
  +\|H_1(\tau)-H_2(\tau)\|_{L^2_\nu}
 \right)\\
 &+
 \left(
 \int_s^{s+1}
  \|H_1(q)-H_2(q)\|_{H^1_\nu}^2\,dq
 \right)^{1/2}.
 \end{split}
\]
The exact low-initial-data estimate is
\begin{equation}\label{eq:uniform-restart-low-initial-propagation}
 \mathscr L_\#(s)
 \leq Ce^{A(s-\tau_0)}d_{\#,0}.
\end{equation}
Since \(d_{\#,0}\leq C d_0\), in particular
\begin{equation}\label{eq:uniform-restart-low-propagation}
 \mathscr L_\#(s)
 \leq Ce^{A(s-\tau_0)}d_0.
\end{equation}
Here \(A,C\) are independent of \(s\).  The same estimate holds for a
sliced first variation with the low initial quantity replaced by
\[
\begin{aligned}
 &\mathfrak D_{m_\#}^{\rm hyb}[\dot{\mathbf z}](\tau_0)
 +d_{{\rm ext},m_\#,0}^{+}[\dot{\mathbf z}]
 +d_{{\rm ext},-1,0}^{\rm corr}[\dot{\mathbf z}]\\
 &\quad
 +d_{{\rm Ggr},m_\#+2,0}^{0}[\dot{\mathbf z}]
 +d_{{\rm Fgr},m_\#+1,0}^{++}[\dot{\mathbf z}]
 +d_{{\rm F},m_\#+1,0}^{\rm glob}[\dot{\mathbf z}]\\
 &\quad
 +\|\dot R(\tau_0)\|_{\mathfrak X_{\rm sc}^{m_\#+2,\alpha}}
 +\|\dot H(\tau_0)\|_{L^2_\nu}.
\end{aligned}
\]
Here brackets on each typed entrance quantity denote its corresponding
linearized seminorm.
The full initial tangent norm gives the stated weaker consequence.
The estimate is a global coarse propagation bound, not a cyclic
spatial restart.  No future boundedness in the full space
\(\mathscr P^{k+2,\alpha}\) is assumed or concluded.
\end{lemma}

\begin{proof}
We first derive, rather than assume,
\eqref{eq:uniform-future-low-coefficient-package}.
Corollary~\ref{cor:adaptive-auxiliary-closure} supplies the invariant
order-twelve curvature package on the auxiliary exterior chains; its
local gauges are used only to establish those one-state estimates and
are not compared across members.  Apply
Lemma~\ref{lem:anchored-exterior-interface} with \(r=9\) to the single
Dirichlet harmonic-map problem on \(E^{++}\).  Its restrictions give
the first line of
\eqref{eq:uniform-future-low-coefficient-package} in one common
Ricci--DeTurck gauge, including the effective-time coefficient modulus
used below.  The same construction on the marked graft chart gives
the transported-marking bounds.  In particular,
\[
 \widetilde\iota=\iota\circ\chi^{-1},\qquad
 (\widetilde\iota)_*\widetilde G=\iota_*G
\]
hold with the quantitative bounds in
\eqref{eq:transported-marking-definition}--%
\eqref{eq:buffered-marking-covariance}.
 The fixed \(m_{\rm ad}=13\) instance of
Proposition~\ref{prop:adaptive-target-tracking}, fixed in
Remark~\ref{conv:authoritative-adaptive-order}, gives the displayed
\(R^{\pm1}\) block, including its \(C^{7,\alpha}\) seminorm.  The fixed
order-six map and inverse-map block is
Lemma~\ref{lem:uniform-HMHF-C6-coefficients}.
The remaining \(C^{5,\alpha}\) tensor block, including the expanding
inner--middle region, is exactly
Lemma~\ref{lem:uniform-scale-adapted-C5-bridge}.  This proves the whole
package without asserting a future bound in every prepared order.

We next prove propagation.  Subtraction of the exact Gram systems and
the order-four prepared calculus give
\begin{equation}\label{eq:uniform-low-hybrid-feedback}
 |\delta c|
 \leq C\left(
  \|H_1-H_2\|_{H^1_\nu}
  +\mathfrak D_{m_\#}^{\rm hyb}\right).
\end{equation}
Every coefficient in this estimate is controlled by
\eqref{eq:uniform-future-low-coefficient-package}; hence \(C\) is
independent of the unit interval.  The exact cutoff difference
equation and stable coercivity similarly give
\begin{equation}\label{eq:uniform-low-hybrid-energy}
 \frac d{d\tau}\|H_1-H_2\|_{L^2_\nu}^2
 +c\|H_1-H_2\|_{H^1_\nu}^2
 \leq
 C\left(
  \|H_1-H_2\|_{L^2_\nu}^2
 +(\mathfrak D_{m_\#}^{\rm hyb})^2\right).
\end{equation}

For clarity separate the one-time memory blocks and the algebraically
derived interface trace from the genuinely restarted variables.  Since
\(\mathfrak D_{m_\#}^{\rm hyb}\) is the sum in
\eqref{eq:two-state-hybrid-distance}, define
\[
\begin{aligned}
 \mathfrak M_\#^{\rm mem}
 &:=
   \mathfrak G_{m_\#}+\mathfrak B_{{\rm gr},m_\#}
   +\mathfrak F_{{\rm gr},m_\#+1}^{\rm hs},\\
 \mathfrak D_{m_\#}^{\rm loc}
 &:=
   |t_1-t_2|
   +\left|\log\frac{\lambda_1}{\lambda_2}\right|
   +\mathfrak M_{{\rm gr},m_\#}\\
 &\quad
   +\|R_1-R_2\|_{\mathfrak X_{\rm sc}^{m_\#+2,\alpha}}
   +\mathfrak F_{m_\#+1}^{\rm glob}
   +\mathfrak F_{{\rm gr},m_\#+1}\\
 &\quad
   +\|h_1-h_2\|_{\mathfrak T_{{\rm sc},N}^{m_\#,\alpha}},\\
 \mathfrak W_\#
 &:=
   \|h_1-h_2\|_{\mathfrak T_{{\rm sc},N}^{4,\alpha}},\\
 \mathfrak Z_\#
 &:=
   \mathfrak D_{m_\#}^{\rm loc}-\mathfrak W_\# .
\end{aligned}
\]
\begin{equation}\label{eq:uniform-interface-algebraic-memory}
 \mathfrak I_{\rm in}
 \leq C\left(\mathfrak D_{m_\#}^{\rm loc}
             +\mathfrak M_\#^{\rm mem}\right),
 \qquad
 \mathfrak D_{m_\#}^{\rm hyb}
 \leq C\left(\mathfrak D_{m_\#}^{\rm loc}
             +\mathfrak M_\#^{\rm mem}\right).
\end{equation}
This is \eqref{eq:interface-block-algebraic-closure}; it is pointwise
and does not restart the \(P\)-equation on \(\mathcal A_{\rm in}\).
Let \(\delta_{\rm Ab}>0\) denote the short-interval constant from
Lemma~\ref{lem:effective-time-endpoint-maximal-regularity}.  Fix
\(J=[s,s+\delta]\) with
\(0<\delta\leq\min\{1,\delta_{\rm Ab}\}\).
The global right-translated \(F\)-equation, with every term linear in
\(F_1-F_2\) retained in its polarized evolution operator, gives
\begin{equation}\label{eq:uniform-global-F-Volterra}
 \begin{split}
 \sup_{s\leq q\leq\tau}
 \mathfrak F_{m_\#+1}^{\rm glob}(q)
 \leq{}&
 C\mathfrak F_{m_\#+1}^{\rm glob}(s)
 +C\int_s^\tau|\delta c(q)|\,dq\\
 &+C\delta^{1/2}
 \sup_{s\leq q\leq\tau}
 \left(\mathfrak D_{m_\#}^{\rm loc}(q)
       +\mathfrak M_\#^{\rm mem}(q)\right),
 \qquad \tau\in[s,s+\delta].
 \end{split}
\end{equation}
This is the one-derivative source-atlas Abel estimate
\eqref{eq:source-atlas-one-order-Abel-block}, applied to the
\(\mathbb F_{\rm sc}^{m_\#,\alpha}\) forcing already typed in
\eqref{eq:coupled-map-difference-forcing}.  On the collapsing core the
metric-difference forcing is retained without claiming a
\(\lambda\)-gain, and the auxiliary order-\((m_\#+2)\) \(R\)-block
types the target-connection difference.  Only the smooth
finite-dimensional phase profiles occur in the \(L^1\) term.  The
global equation is restarted from its actual global trace at \(s\);
when \(s=\tau_0\), that trace is
\(d_{{\rm F},m_\#+1,0}^{\rm glob}\).  No localized star or interface
trace is used.
Recenter on \(J\) the
same mixed contraction used in
\eqref{eq:finite-horizon-non-h-top-ledger}.  Its constants are now
independent of \(s\) by
\eqref{eq:uniform-future-low-coefficient-package}.  Propagate the
auxiliary order-\((m_\#+2)\) \(R\)-block by its ODE.  The interface
block is excluded from this
recentered system and is recovered pointwise from
\eqref{eq:uniform-interface-algebraic-memory}.
The physical exterior, graft-input, and homogeneous/separated
localized source-adapted \(F\)-blocks are not restarted; they enter as
the one-time memory
\(\mathfrak M_\#^{\rm mem}\), governed respectively by
\eqref{eq:inner-terminated-exterior-DeTurck} and
\eqref{eq:compact-graft-buffer-finite}, and
\eqref{eq:localized-graft-F-coarse-memory}.  The global \(F\)-block is
a locally restarted summand of \(\mathfrak D_{m_\#}^{\rm loc}\) and is
governed by \eqref{eq:uniform-global-F-Volterra}.
Thus, for \(s\leq\tau\leq s+\delta\),
\begin{equation}\label{eq:uniform-non-h-top-ledger}
 \begin{aligned}
 \sup_{s\leq q\leq\tau}\mathfrak Z_\#(q)
 \leq{}&
 C\mathfrak D_{m_\#}^{\rm loc}(s)
 +C\sup_{s\leq q\leq\tau}\mathfrak M_\#^{\rm mem}(q)\\
  &+C\int_s^\tau
    \left(\mathfrak D_{m_\#}^{\rm loc}(q)
          +|\delta c(q)|\right)dq\\
  &+C\bigl(\delta^{\alpha/4}+\delta^{1/2}\bigr)
    \sup_{s\leq q\leq\tau}
    \bigl(\mathfrak W_\#(q)+\mathfrak Z_\#(q)
          +\mathfrak M_\#^{\rm mem}(q)\bigr).
 \end{aligned}
\end{equation}
The inputs in this estimate are precisely the displayed exterior/graft
memory and the actual global value of the \(F\)-trace at \(s\).  The
\(\delta^{1/2}\) term is the Abel contribution from
\eqref{eq:uniform-global-F-Volterra} and from the nonseparated
localized \(F\)-source.

In the exact equation for \(h_1-h_2\), the relative top-order
transport is
\(\delta c\,\mathcal W*\bar\nabla h_2\).  The
\(C^{5,\alpha}\) tensor line of
\eqref{eq:uniform-future-low-coefficient-package} places this source at
order four, with norm \(C|\delta c|\).  The pure graft source is
in \(C_{\rm sc}^{2,\alpha}\) by
Lemma~\ref{lem:same-order-pure-graft-difference}.  Its two-derivative
endpoint recovery is \eqref{eq:h-specific-scaled-graft-Duhamel}, with
\(\lambda_*=\max\{\lambda_1,\lambda_2\}\) and \(L(q)\) the dyadic scale
of the enlarged cylinder meeting the graft at time \(q\).  Thus its
source norm is the full bracket
\[
 \mathfrak B_{{\rm gr},m_\#}
 +\left|\log\frac{\lambda_1}{\lambda_2}\right|
 +\|R_1-R_2\|_{\mathfrak X_{\rm sc}^{m_\#+1,\alpha}}
 +\mathfrak F_{{\rm gr},m_\#+1}
 \leq C\bigl(
   \mathfrak D_{m_\#}^{\rm loc}+\mathfrak M_\#^{\rm mem}\bigr),
\]
not the graft-input block alone.  On its support
\(L(q)\lambda_*(q)/\Gamma\asymp1\), so the weak
\(L(q)^{-1}\)-diffusion creates no loss.  Every remaining
nonfeedback source is placed in its time-supremum norm before weighted
Schauder is applied, whereas each feedback profile is treated by the
\(h\)-specific estimate \eqref{eq:h-feedback-L1-Duhamel} with
\(m=m_\#=4\).  In the present application the
uniform future coefficient package
\eqref{eq:uniform-future-low-coefficient-package}, together with the
fixed adaptive phase budget in \(\mathfrak P_{\rm prep}\), makes the
constant in that estimate independent of the terminal endpoint.  Hence
\begin{equation}\label{eq:uniform-h-top-ledger}
 \begin{split}
 \sup_{s\leq q\leq\tau}\mathfrak W_\#(q)
 \leq C\Bigg(&
   \mathfrak D_{m_\#}^{\rm loc}(s)
   +\sup_{s\leq q\leq\tau}\mathfrak M_\#^{\rm mem}(q)
   +\sup_{s\leq q\leq\tau}\mathfrak Z_\#(q)\\
 &+\int_s^\tau
   \bigl(\mathfrak D_{m_\#}^{\rm loc}(q)
         +|\delta c(q)|\bigr)\,dq\Bigg).
 \end{split}
\end{equation}
Choose a fixed
\(0<\delta_\#\leq\min\{1,\delta_{\rm Ab}\}\) so that both the mixed
\(\delta^{\alpha/4}\) zero-trace term and the
\(\delta^{1/2}\) global/local \(F\) Abel term in
\eqref{eq:uniform-non-h-top-ledger} are absorbed after
\eqref{eq:uniform-h-top-ledger} is substituted.  We obtain, on every
interval of length at most \(\delta_\#\),
\begin{equation}\label{eq:uniform-triangular-top-ledger}
 \begin{split}
 \sup_{s\leq q\leq\tau}\mathfrak D_{m_\#}^{\rm loc}(q)
 \leq C\Bigg(&
  \mathfrak D_{m_\#}^{\rm loc}(s)
  +\sup_{s\leq q\leq\tau}\mathfrak M_\#^{\rm mem}(q)\\
 &+\int_s^\tau
   \bigl(\mathfrak D_{m_\#}^{\rm loc}(q)
         +|\delta c(q)|\bigr)\,dq\Bigg).
 \end{split}
\end{equation}
The raw \(F\)-clock and the normalized tensor clock remain distinct
throughout this argument.  The feedback estimate is retained in its
integral form:
\[
 \int_s^\tau|\delta c(q)|\,dq
 \leq C(\tau-s)^{1/2}
 \left(\int_s^\tau
  \|H_1-H_2\|_{H^1_\nu}^2\,dq\right)^{1/2}
 +C\int_s^\tau\left(
   \mathfrak D_{m_\#}^{\rm loc}
   +\mathfrak M_\#^{\rm mem}\right)(q)\,dq .
\]
Square \eqref{eq:uniform-triangular-top-ledger}, combine it with
\eqref{eq:uniform-low-hybrid-energy}, using
\eqref{eq:uniform-interface-algebraic-memory} to eliminate the
interface block, and subdivide a unit interval
into the fixed number of pieces of length at most \(\delta_\#\).
The discrete Volterra--Gronwall inequality gives, for
\(\tau\in[s,s+1]\),
\begin{equation}\label{eq:uniform-low-hybrid-Volterra}
 \begin{split}
 &\bigl(\mathfrak D_{m_\#}^{\rm loc}(\tau)\bigr)^2
  +\|H_1(\tau)-H_2(\tau)\|_{L^2_\nu}^2
  +c\int_s^\tau\|H_1-H_2\|_{H^1_\nu}^2\,dq\\
 &\quad\leq
  C\left[
   \bigl(\mathfrak D_{m_\#}^{\rm loc}(s)\bigr)^2
   +\sup_{s\leq r\leq\tau}
      \bigl(\mathfrak M_\#^{\rm mem}(r)\bigr)^2
   +\|H_1(s)-H_2(s)\|_{L^2_\nu}^2\right]\\
 &\qquad
 +C\int_s^\tau\left(
   \bigl(\mathfrak D_{m_\#}^{\rm loc}(q)\bigr)^2
   +\|H_1-H_2\|_{L^2_\nu}^2\right)dq .
 \end{split}
\end{equation}

The physical exterior and graft-input block are now recovered once
from the initial face.
Lemma~\ref{lem:inner-terminated-exterior-DeTurck} and the algebraic
interface estimate
\eqref{eq:uniform-interface-algebraic-memory} imply the coarse
Volterra bound
\begin{equation}\label{eq:global-exterior-coarse-memory}
 \mathfrak G_{m_\#}(\tau)
 \leq Cd_{\#,0}+
 C\int_{\tau_0}^{\tau}
 \left(
   \mathfrak D_{m_\#}^{\rm loc}(q)
  +\mathfrak M_\#^{\rm mem}(q)
  +\|H_1(q)-H_2(q)\|_{L^2_\nu}\right)dq .
\end{equation}
Here the bounded Davies kernel and \(\lambda\leq Ce^{-\tau}\) have
only been discarded in the favorable direction.  More explicitly, the
kernel discarded here is bounded because, with its value at
\(\sigma=0\) understood by continuity,
\[
 \sup_{\sigma\geq0}
 \bigl(1+\sigma^{-q_{m_\#}}\bigr)
 \exp\!\left(-\frac{cR_{\rm in}^2}{\sigma}\right)<\infty .
\]
Set
\[
\begin{aligned}
 X(\tau)
 &:=
 \mathfrak D_{m_\#}^{\rm loc}(\tau)
 +\mathfrak M_\#^{\rm mem}(\tau)
 +\|H_1(\tau)-H_2(\tau)\|_{L^2_\nu},\\
 Y(\tau)
 &:=
 d_{\#,0}+\int_{\tau_0}^{\tau}
 \left(X(q)+\|H_1(q)-H_2(q)\|_{H^1_\nu}\right)dq .
\end{aligned}
\]
Thus \(Y\) is nondecreasing, \(X\) contains both the localized and
global low \(F\)-blocks, and
\eqref{eq:uniform-interface-algebraic-memory} gives
\[
 \mathfrak D_{m_\#}^{\rm hyb}
 +\mathfrak F_{m_\#+1}^{\rm glob}
 +\|H_1-H_2\|_{L^2_\nu}\leq CX.
\]
First \eqref{eq:global-exterior-coarse-memory} gives
\(\mathfrak G_{m_\#}\leq CY\).  Moreover
\(d_{{\rm gr},m_\#,0}\leq d_{\#,0}\).  By
\eqref{eq:graft-clock-source-amplitude} and
\eqref{eq:uniform-low-hybrid-feedback},
\[
 \mathfrak L_{12}(\tau)
 \leq C\left(
 d_{\#,0}+\int_{\tau_0}^{\tau}
   (X+\|H_1-H_2\|_{H^1_\nu})\,dq\right)
 \leq CY(\tau).
\]
Thus the nonseparated clock source in
\eqref{eq:compact-graft-buffer-finite} is controlled by endpoint
maximal regularity before the separated memory is estimated.
Since \(\lambda_1\leq Ce^{-q}\), in particular
\(\int_{\tau_0}^{\infty}\lambda_1\,dq\leq C\); the scalar Davies
kernel is bounded, and \(Y\) is nondecreasing.  Thus the same
bounded-kernel argument applied to that one-time estimate, with the
exterior term just obtained, gives
\[
 \mathfrak B_{{\rm gr},m_\#}(\tau)\leq CY(\tau).
\]
Finally
\[
 d_{{\rm Fgr},m_\#+1,0}^{++}\leq d_{\#,0},
\]
which is the homogeneous localized-star trace.  The only other
summand in the one-time localized map memory
\eqref{eq:localized-graft-F-hs-memory} is spatially separated.  For
that term,
\[
 \mathfrak s(q,\tau)\leq Ce^{-q},\qquad
 \lambda_1(q)
 K_{{\rm Fgr},m_\#}\bigl(\mathfrak s(q,\tau)\bigr)
 \leq Ce^{-q}(1+e^{q_{m_\#+1}q})e^{-c'e^q}
 =:\kappa_\#(q).
\]
The scalar function \(\kappa_\#\) is bounded and integrable, and the
term itself is bounded by \(\int\kappa_\#X\leq CY\).  Therefore
\[
 \mathfrak F_{{\rm gr},m_\#+1}^{\rm hs}(\tau)\leq CY(\tau).
\]
The full global and localized low \(F\)-blocks are locally restarted
components of
\(\mathfrak D_{m_\#}^{\rm loc}\subset X\), whereas only the
homogeneous/separated localized term occurs in
\(\mathfrak M_\#^{\rm mem}\).  Its corridor kernel contains the
integrable
factor \(\lambda_1K_{{\rm Fgr},m_\#}\).  Together with the exterior
and graft bounds, this gives
\begin{equation}\label{eq:uniform-one-time-memory-by-Y}
 \mathfrak M_\#^{\rm mem}(\tau)\leq CY(\tau).
\end{equation}
The resulting low-initial graft and localized-\(F\) estimates depend
only on initial contributions already contained in \(d_{\#,0}\) and
introduce no higher prepared entrance norm.

We now close the coupled Volterra system rather than estimating
\(\mathfrak M_\#^{\rm mem}\) outside it.  Put
\[
 E(\tau):=\mathfrak D_{m_\#}^{\rm loc}(\tau)
       +\|H_1(\tau)-H_2(\tau)\|_{L^2_\nu},
 \qquad X=E+\mathfrak M_\#^{\rm mem}.
\]
On \(J=[s,s+\delta]\), \(0<\delta\leq1\),
\eqref{eq:uniform-low-hybrid-Volterra} and
\eqref{eq:uniform-one-time-memory-by-Y} give
\begin{equation}\label{eq:uniform-coupled-E-estimate}
 \sup_J E^2+\int_J\|H_1-H_2\|_{H^1_\nu}^2\,dq
 \leq C\left(
  E(s)^2+\sup_JY^2+\delta\sup_JX^2\right).
\end{equation}
Here and below a supremum over \(J\) is understood to be taken only
over the portion of \(J\) under consideration.  On the other hand,
the definition of \(Y\) and Cauchy--Schwarz give
\begin{equation}\label{eq:uniform-coupled-Y-increment}
 \sup_JY^2
 \leq 2Y(s)^2
 +C\delta^2\sup_JX^2
 +C\delta\int_J\|H_1-H_2\|_{H^1_\nu}^2\,dq .
\end{equation}
Substitute \eqref{eq:uniform-coupled-E-estimate} into
\eqref{eq:uniform-coupled-Y-increment} and first choose \(\delta\) so
that the resulting \(C\delta\sup_JY^2\) term is absorbed.  Then
\[
 \sup_JY^2
 \leq C\left(
  Y(s)^2+\delta E(s)^2+\delta^2\sup_JX^2\right),
\]
and another substitution in
\eqref{eq:uniform-coupled-E-estimate} gives
\[
 \sup_JE^2+\int_J\|H_1-H_2\|_{H^1_\nu}^2\,dq
 \leq C\left(E(s)^2+Y(s)^2\right)
 +C\delta\sup_JX^2 .
\]
Finally
\eqref{eq:uniform-one-time-memory-by-Y} and
\(X=E+\mathfrak M_\#^{\rm mem}\) yield
\[
 \sup_JX^2
 \leq C\left(E(s)^2+Y(s)^2\right)+C\delta\sup_JX^2 .
\]
Choose one fixed \(\delta>0\), depending only on the common package,
so that the final term is absorbed.  Equations
\eqref{eq:uniform-coupled-E-estimate} and
\eqref{eq:uniform-coupled-Y-increment} then give
\[
 \sup_J(X^2+Y^2)
 +\int_J\|H_1-H_2\|_{H^1_\nu}^2\,dq
 \leq C\left(E(s)^2+Y(s)^2\right)
 \leq C\left(X(s)^2+Y(s)^2\right).
\]
By definition \(X(\tau_0)+Y(\tau_0)\leq Cd_{\#,0}\), including the
 buffered exterior trace, the separated exterior corridor, the
 larger-collar graft trace, the larger-star \(F\)-trace, and the global
 low initial \(F\)-trace.
Iteration over these fixed subintervals yields constants \(C,A\),
depending only on the common package, such that
\[
 X(\tau)+Y(\tau)
 \leq Ce^{A(\tau-\tau_0)}d_{\#,0},
 \qquad
 \left(\int_s^{s+1}
   \|H_1-H_2\|_{H^1_\nu}^2\,dq\right)^{1/2}
 \leq Ce^{A(s-\tau_0)}d_{\#,0}.
\]
This proves
\eqref{eq:uniform-restart-low-initial-propagation};
\eqref{eq:uniform-restart-low-propagation} follows from
\(d_{\#,0}\leq Cd_0\).  Thus the exterior and graft inputs remain
anchored at their original buffered traces.  Only the
homogeneous/separated localized \(F\)-memory remains anchored at the
original larger-star trace; its nonseparated part is recentered with
zero trace on the same star, never on a larger one.  The global
\(F\)-block is restarted only as a global equation from its actual
global trace at the left endpoint.

Sliced first variations in the finite-horizon \(C^1\) solution map
satisfy the same linear estimates, with the low initial quantity
displayed in the statement.  This conclusion is an estimate for the
exact sliced variational system.  No entrance-time-uniform,
all-block Fr\'echet-remainder estimate in the future low-order hybrid
topology is asserted here; finite-horizon Fr\'echet differentiability
is supplied separately by
Proposition~\ref{prop:two-state-prepared-evolution}, and the global
\(C^1\) scattering conclusion in
Theorem~\ref{thm:global-two-state-scattering} follows from uniform
convergence of the finite-endpoint derivatives.
\end{proof}

\begin{lemma}[Coarse two-state geometry]
\label{lem:coarse-two-state-geometry}
In the fixed exterior and marked DeTurck gauges supplied by
Proposition~\ref{prop:two-state-prepared-evolution}, consider two
global solutions from the common-margin sliced entrance ball
\(\mathscr B\) of
Lemma~\ref{lem:uniform-weighted-Schauder-restart}.  Put
\[
 d_0:=
 \|\mathbf z_{1,0}-\mathbf z_{2,0}\|
   _{\mathscr X_{\rm prep}^{k+2,\alpha}}.
\]
There are
constants \(A,C<\infty\), depending only on the common-margin ball,
such that
\begin{equation}\label{eq:coarse-two-state-growth}
 \mathfrak D_{m_\#}^{\rm hyb}(\tau)
 +\mathfrak F_{m_\#+1}^{\rm glob}(\tau)
 \leq Ce^{A(\tau-\tau_0)}d_0,
 \qquad \tau\geq\tau_0.
\end{equation}
For \(i=1,2\), write
\[
 \E_i:=(\Phi_i^{-1})^*\mathcal G_{{\rm gr},i}.
\]
Set
\[
 \Delta\mathcal C_\rho
 :=\mathcal C_\rho[h_1]-\mathcal C_\rho[h_2],
 \qquad
 \Delta\E:=\E_1-\E_2,
 \qquad
 \Delta\mathcal Y_{j,\tau}
 :=\mathcal Y_{j,\tau}(\mathbf z_1)
   -\mathcal Y_{j,\tau}(\mathbf z_2).
\]
Thus \(\Delta\mathcal C_\rho\) is the complete moving-cutoff
commutator difference, \(\Delta\E\) is the moving graft difference,
and \(\Delta\mathcal Y_{j,\tau}\) is the full effective-column
difference, including its solution-dependent cutoff and pullback.
These are exactly the outer blocks in the localized equations.
They are supported in \(\{\bar f\geq c e^\tau\}\) and obey
\begin{equation}\label{eq:coarse-two-state-Gaussian-absorption}
 \begin{split}
 &\|\Delta\mathcal C_\rho\|_{H^{-1}_\nu}
  +\|\rho_\tau\Delta\E\|_{H^{-1}_\nu}
  +\sum_{j=0}^8
    \|\rho_\tau\Delta\mathcal Y_{j,\tau}\|_{H^{-1}_\nu}\\
 &\quad+
 \sum_{\mu=0}^8
 \left(
  \left|\ip{\Delta\mathcal C_\rho}{Z_\mu}\right|
  +\left|\ip{\rho_\tau\Delta\E}{Z_\mu}\right|
  +\sum_{j=0}^8
    \left|\ip{\rho_\tau\Delta\mathcal Y_{j,\tau}}{Z_\mu}\right|
 \right)
 \leq Cd_0e^{-c'e^\tau}.
 \end{split}
\end{equation}
\end{lemma}

\begin{proof}
The supremum in
\eqref{eq:uniform-restart-low-propagation}, on the unit interval
containing \(\tau\), contains both
\(\mathfrak D_{m_\#}^{\rm hyb}\) and the separate coarse block
\(\mathfrak F_{m_\#+1}^{\rm glob}\).  It therefore gives
\eqref{eq:coarse-two-state-growth} after enlarging \(A,C\) once.  The
global block comes from the genuine global restart
\eqref{eq:uniform-global-F-Volterra} on each uniform short interval,
with the actual global trace at its left endpoint; no localized trace
and no uniform-in-time global map bound is being asserted.  The
identical conclusion holds for sliced
first variations.

The scale-one inverse, pullback, and composition estimates in
Lemma~\ref{lem:prepared-chart-calculus} make every listed outer
difference at most a fixed polynomial in \(e^\tau\)
times
\(\mathfrak D_{m_\#}^{\rm hyb}(\tau)
 +\mathfrak F_{m_\#+1}^{\rm glob}(\tau)\).
Only order-four geometric blocks and the order-five global map block
enter these second-order outer expressions.  Their support lies where
\(\bar f\geq ce^\tau\).  For some \(c_{\rm G}>0\), the Gaussian density
and the mode-growth bounds therefore contribute
\(e^{-c_{\rm G}e^\tau}\); this absorbs both the
fixed polynomial and \(e^{A\tau}\), proving
\eqref{eq:coarse-two-state-Gaussian-absorption}.  This argument uses
only the one-state bounds and the finite-horizon theory.
\end{proof}

\begin{lemma}[Buffered Davies localization]
\label{lem:buffered-Davies-localization}
Fix an integer \(m\geq4\) and \(0<\alpha<1\).
Fix an interval
\[
 I=[\tau_*,\tau^*],\qquad
 \tau_*<\tau^*\leq\infty,
\]
and let \(\lambda\in C(I;(0,\infty))\), with all assertions on an
unbounded \(I\) understood locally on compact subintervals.  Every
evolution-family estimate below is quantified for
\(\tau_*\leq\sigma\leq\tau\) in \(I\).
Let \(W\Subset W^+\) have fixed buffer width at least \(2R>0\), where
\(W^+\) is a relatively compact smooth domain.  Fix a background
metric \(g_{\rm D}\), a metric vector bundle \(E\to W^+\), and a
compatible connection \(\nabla\).  Assume that \(W^+\) has a
\(C^{m+2,\alpha}\) harmonic-radius lower bound and is covered by a
finite harmonic atlas with uniformly controlled overlap, chart
constants, and boundary charts.  Consider on \(E|_{W^+}\) a linear
system with homogeneous Dirichlet boundary condition
\[
 \partial_\tau u
 -\lambda(\tau)A^{ab}(\tau,x)\nabla_a\nabla_bu
 =
 \lambda(\tau)\bigl(B*\nabla u+C*u\bigr)+F .
\]
 The principal symbol is required to be
the symmetric scalar bundle symbol
\[
 A^{ab}(\tau,x)\xi_a\xi_b\operatorname{Id}_E,
\]
uniformly elliptic relative to \(g_{\rm D}\).  In the effective-time
variable
\[
 \mathfrak t(\tau)=\int_{\tau_*}^{\tau}\lambda(q)\,dq,
\]
assume that \(A,B,C\) have uniform spatial bounds in
\(C^{m,\alpha}\), \(C^{m-1,\alpha}\), and
\(C^{m-2,\alpha}\), respectively, in the fixed atlas and a common
\(C^{\alpha/2}\) time bound (a common parabolic time-oscillation
modulus suffices).  Put
\[
 \mathfrak s(\sigma,\tau)=\int_\sigma^\tau\lambda(q)\,dq.
\]
Here
\[
 H^{-1}(W^+;E):=\bigl(H_0^1(W^+;E)\bigr)^*
\]
with duality induced by the fixed background metric and bundle metric.
 Let \(\mathcal U(\tau,\sigma)\), \(\tau_*\leq\sigma\leq\tau\) in
 \(I\), be the evolution family for the
homogeneous equation \(F=0\).  If \(g(\sigma)\) is supported a
distance at least \(R\) from \(W\), then, for some integer \(q_m\),
\begin{equation}\label{eq:Davies-off-diagonal}
 \|\mathcal U(\tau,\sigma)g(\sigma)\|_{C^{m,\alpha}(W)}
 \leq
 C e^{C\mathfrak s(\sigma,\tau)}
 \bigl(1+\mathfrak s(\sigma,\tau)^{-q_m}\bigr)
 \exp\!\left(-\frac{cR^2}{\mathfrak s(\sigma,\tau)}\right)
 \|g(\sigma)\|_{C^{m,\alpha}(W^+)} .
\end{equation}
The right side is interpreted as zero at
\(\mathfrak s(\sigma,\tau)=0\).  More generally, if
\(g(\sigma)\in H^{-1}(W^+;E)\) has the same support separation, then
\begin{equation}\label{eq:Davies-Hminus-one}
 \|\mathcal U(\tau,\sigma)g(\sigma)\|_{L^2(W)}
 \leq
 C e^{C\mathfrak s(\sigma,\tau)}
 \mathfrak s(\sigma,\tau)^{-1/2}
 \exp\!\left(-\frac{cR^2}{\mathfrak s(\sigma,\tau)}\right)
 \|g(\sigma)\|_{H^{-1}(W^+)} .
\end{equation}
After increasing \(q_m\), interior parabolic smoothing gives
\begin{equation}\label{eq:Davies-Hminus-one-Holder}
 \|\mathcal U(\tau,\sigma)g(\sigma)\|_{C^{m,\alpha}(W)}
 \leq
 C e^{C\mathfrak s(\sigma,\tau)}
 \bigl(1+\mathfrak s(\sigma,\tau)^{-q_m}\bigr)
 \exp\!\left(-\frac{cR^2}{\mathfrak s(\sigma,\tau)}\right)
 \|g(\sigma)\|_{H^{-1}(W^+)} .
\end{equation}
Both right sides are again interpreted as zero at zero effective time.
Consequently, if \(G\in L^1([\sigma,\tau];H^{-1}(W^+;E))\)
has the same support separation for almost every \(r\), then
\begin{equation}\label{eq:Davies-Hminus-one-Duhamel}
\begin{split}
 &\left\|
   \int_\sigma^\tau\mathcal U(\tau,r)G(r)\,dr
  \right\|_{C^{m,\alpha}(W)}
 \\
 &\quad\leq
 C\int_\sigma^\tau
 e^{C\mathfrak s(r,\tau)}
 \bigl(1+\mathfrak s(r,\tau)^{-q_m}\bigr)
 \exp\!\left(-\frac{cR^2}{\mathfrak s(r,\tau)}\right)
 \|G(r)\|_{H^{-1}(W^+)}\,dr .
\end{split}
\end{equation}
More generally, suppose that \(G\) is Bochner integrable in the
displayed H\"older space and, for every \(r\),
\[
 \operatorname{dist}_{g_{\rm D}}(\supp G(r),W)\geq R.
\]
Then Duhamel's formula gives
\begin{equation}\label{eq:Davies-off-diagonal-Duhamel}
 \begin{split}
 &\left\|
   \int_\sigma^\tau\mathcal U(\tau,r)G(r)\,dr
  \right\|_{C^{m,\alpha}(W)}
 \\
 &\quad\leq
 C\int_\sigma^\tau
 e^{C\mathfrak s(r,\tau)}
 \bigl(1+\mathfrak s(r,\tau)^{-q_m}\bigr)
 \exp\!\left(-\frac{cR^2}{\mathfrak s(r,\tau)}\right)
 \|G(r)\|_{C^{m,\alpha}(W^+)}\,dr ,
 \end{split}
\end{equation}
where the kernel is set equal to zero when
\(\mathfrak s(r,\tau)=0\).  The constants depend only on the
buffer, the displayed bounded-geometry and atlas data, the bundle
rank and connection bounds, the ellipticity constants, the
coefficient bounds, and the common effective-time modulus.  The same
estimates hold on a closed double, in which case the boundary
condition is omitted.  They also apply to a cutoff of a global
solution whenever the cutoff vanishes near \(\partial W^+\).  A
divergence-form cutoff commutator is inserted in
\eqref{eq:Davies-Hminus-one-Duhamel}; only a commutator already
controlled in \(C^{m,\alpha}\) is inserted in
\eqref{eq:Davies-off-diagonal-Duhamel}.
\end{lemma}

\begin{proof}
Choose a smooth function \(\phi\) which vanishes on the support of
\(g\), equals \(R\) on a neighborhood of \(W\), and satisfies
\(|\nabla\phi|+R|\nabla^2\phi|\leq C\).  For the homogeneous evolution,
apply the energy identity to \(e^{a\phi}u\).  The homogeneous Dirichlet
condition removes the boundary term.  Uniform ellipticity
absorbs the cross term containing
\(a\lambda\nabla\phi\cdot\nabla u\), and the lower-order coefficients
give
\[
 \frac d{d\tau}\|e^{a\phi}u\|_{L^2}^2
 \leq C(1+a^2)\lambda(\tau)
       \|e^{a\phi}u\|_{L^2}^2 .
\]
Integration from \(\sigma\) to \(\tau\), followed by the choice
\(a=cR/\mathfrak s(\sigma,\tau)\), yields the \(L^2\) factor
\[
 e^{C\mathfrak s}
 \exp(-cR^2/\mathfrak s).
\]
 To pass from this base \(L^2\) estimate to the stated derivative
 estimate, choose nested relatively compact domains
 \[
   W\Subset W_1\Subset W_2\Subset W^+
 \]
 with \(W_2\) still separated from the source region.  Apply the
 preceding \(L^2\) Davies estimate with \(W_2\) in place of \(W\), and
 then apply interior parabolic H\"older estimates to the equation on
 \(W_2\), first on \(W_1\) and then on \(W\).  This avoids
differentiating the Dirichlet boundary condition.  The displayed
graded coefficient bounds, together with the
\(C^{m+2,\alpha}\) background atlas, control the interior commutators,
while
 the inverse powers of the effective time in the interior estimates
 produce a factor \(1+\mathfrak s^{-q_m}\).  Since the nested-domain
 separations are fixed fractions of \(R\), the Gaussian exponent is
 merely weakened by an absolute constant.  This proves
\eqref{eq:Davies-off-diagonal}.
For \(\psi\in L^2(W)\), duality gives
\[
 \langle\mathcal U(\tau,\sigma)g,\psi\rangle
 =\langle g,\mathcal U(\tau,\sigma)^*\psi\rangle .
\]
The same weighted energy argument for the adjoint evolution, combined
with its Caccioppoli estimate over effective time \(\mathfrak s\),
bounds the \(H^1\)-norm of
\(\mathcal U(\tau,\sigma)^*\psi\) on the source region by
\[
 C e^{C\mathfrak s}\mathfrak s^{-1/2}
   e^{-cR^2/\mathfrak s}\|\psi\|_{L^2(W)}.
\]
This proves \eqref{eq:Davies-Hminus-one}.  Applying the preceding
nested interior smoothing argument to its output proves
\eqref{eq:Davies-Hminus-one-Holder}; all additional inverse powers are
absorbed by increasing \(q_m\).
The derivatives of \(A\) produced when the nondivergence principal
term is integrated by parts are controlled by the assumed spatial
coefficient bounds and are absorbed into the same lower-order
constant.  Applying \eqref{eq:Davies-Hminus-one-Holder} and
\eqref{eq:Davies-off-diagonal} at each source time and integrating
proves \eqref{eq:Davies-Hminus-one-Duhamel} and
\eqref{eq:Davies-off-diagonal-Duhamel}, respectively.
\end{proof}

\begin{lemma}[Buffered physical two-state stability]
\label{lem:buffered-physical-two-state-stability}
Consider the two global solutions of
Lemma~\ref{lem:coarse-two-state-geometry}.  Assume that their
common-margin ball of strict prepared entrances retains one fixed
ordinary exterior certificate and the strict time-width margin
\eqref{eq:auxiliary-buffered-time-width}.  Let
\(\mathcal W_{\rm gr}\Subset\mathcal W_{\rm gr}^+\Subset E\) be the
graft collars in \eqref{eq:separated-graft-collars}, covered by the
smallest atlas members, and assume the initial exterior
curvature bounds through order eight.  If, for some \(D,\vartheta>0\),
\begin{equation}\label{eq:buffered-stability-phase-hypothesis}
 \int_\tau^\infty|c_1-c_2|(s)\,ds
 \leq Dd_0e^{-\vartheta(\tau-\tau_0)},
\end{equation}
then, in the common Ricci--DeTurck gauge,
\begin{equation}\label{eq:buffered-two-state-stability}
 \sup_{\tau\geq\tau_0}
 \|\widetilde G_1(t_1(\tau))
      -\widetilde G_2(t_2(\tau))\|
        _{C^{6,\alpha}(\mathcal W_{\rm gr}^+)}
 \leq C_Dd_0 .
\end{equation}
On the once-smaller collar
\(\mathcal W_{\rm gr}\Subset\mathcal W_{\rm gr}^+\), the
transported gauges and markings also satisfy
\begin{equation}\label{eq:buffered-two-state-marking-stability}
 \sup_{\tau\geq\tau_0}\left(
  \|\chi_1-\chi_2\|_{C^{5,\alpha}(\mathcal W_{\rm gr})}
  +\|\chi_1^{-1}-\chi_2^{-1}\|_{C^{5,\alpha}(\mathcal W_{\rm gr})}
  +\|\widetilde\iota_1-\widetilde\iota_2\|
     _{C^{5,\alpha}(\mathcal W_{\rm gr})}
 \right)
 \leq C_Dd_0 .
\end{equation}
These are exactly the larger-input/smaller-output collars in
\eqref{eq:hybrid-graft-buffer-block}.  The fixed recorded graft scale
makes the displayed ordinary and scale-normalized norms at orders six
and five uniformly equivalent.  Thus these estimates say precisely
\begin{equation}\label{eq:buffered-two-state-graft-buffer-stability}
 \sup_{\tau\geq\tau_0}
 \mathfrak B_{{\rm gr},4}(\tau)\leq C_Dd_0 .
\end{equation}
The constant is independent of the terminal normalized time.
For this lemma alone, only the consequence
\(\int_{\tau_0}^\infty|c_1-c_2|\leq Dd_0\) is needed; the exponent
\(\vartheta\) is retained because the following graft and barrier
applications use the full tail form.
\end{lemma}

\begin{proof}
The scale equations and
\eqref{eq:buffered-stability-phase-hypothesis} first give
\begin{equation}\label{eq:buffered-scale-time-difference}
 \sup_{\tau\geq\tau_0}
 \left(
  \left|\log\frac{\lambda_1}{\lambda_2}\right|
  +|t_1-t_2|
 \right)
 \leq C_Dd_0.
\end{equation}
Indeed, subtract
\((\log\lambda_i)_\tau=-(1+a_i)\), use the initial prepared distance,
 and then integrate \(t_{i,\tau}=\lambda_i\), noting that
 \(\lambda_i\leq C e^{-\tau}\).  The one-state physical scale estimate
 also shows directly that the clock limit
 \[
  T_i:=t_i(\tau_0)+\int_{\tau_0}^{\infty}\lambda_i(q)\,dq
 \]
 exists and gives, without appealing to a later extinction theorem,
 \[
  T_i-t_i(\tau)\asymp\lambda_i(\tau)\asymp e^{-\tau}.
 \]

We compare directly at the common normalized time; no synchronization
of the two physical clocks is required.  Since
\(\mathfrak I_{\rm in}\) is a summand of the hybrid distance,
\eqref{eq:coarse-two-state-growth} gives
\[
 \mathfrak I_{\rm in}(\tau)
 \leq C d_0e^{A(\tau-\tau_0)}.
\]
Combined with \eqref{eq:buffered-scale-time-difference}, this verifies
the two hypotheses in the global clause
\eqref{eq:inner-terminated-exterior-global} of the inner-terminated
exterior estimate.  Applying
 Lemma~\ref{lem:inner-terminated-exterior-DeTurck} at \(m=6\) gives
\[
 \sup_{\tau\geq\tau_0}\mathfrak G_6(\tau)\leq C_Dd_0 .
\]
Since \(\mathcal W_{\rm gr}^+\Subset E\) is covered by the recorded
smallest atlas members, this is
\eqref{eq:buffered-two-state-stability}.  Using the \(H^{-1}\) cutoff
estimate \eqref{eq:Davies-Hminus-one-Holder}, rather than
differentiating the cutoff commutator six times, avoids a
seventh-derivative requirement at this output order.

It remains only to propagate the gauge and marking on the once-smaller
collar.  In normalized time the triangular gauge and inverse-gauge
equations have forcing bounded by
\[
 C\lambda_1\|\widehat G_1-\widehat G_2\|
       _{C^{6,\alpha}(\mathcal W_{\rm gr}^+)}
 +C|\lambda_1-\lambda_2|
 \leq C_Dd_0e^{-\tau}.
\]
The inner boundary is anchored by
Lemma~\ref{lem:anchored-exterior-interface}, and the initial gauge and
marking differences are \(O(d_0)\).  Direct integration, followed by
\(\widetilde\iota_i=\iota\circ\chi_i^{-1}\) and the fixed
composition estimate, proves
\eqref{eq:buffered-two-state-marking-stability}, and the definition
\eqref{eq:hybrid-graft-buffer-block} then gives
\eqref{eq:buffered-two-state-graft-buffer-stability}.
\end{proof}

\begin{lemma}[Localized map comparison]
\label{lem:localized-graft-F-comparison}
Assume the hypotheses of
Lemma~\ref{lem:buffered-physical-two-state-stability}.  In particular,
assume the phase-tail bound
\eqref{eq:buffered-stability-phase-hypothesis}.  Then
\begin{equation}\label{eq:localized-graft-F-uniform}
 \sup_{\tau\geq\tau_0}\mathfrak F_{{\rm gr},5}(\tau)
 \leq C_Dd_0 .
\end{equation}
Only the source-adapted star in
\eqref{eq:localized-graft-F-block} is asserted here.
\end{lemma}

\begin{proof}
Apply the split cutoff-localized estimate
\eqref{eq:localized-graft-F-coarse-memory} with \(m=4\).  Its
larger-star homogeneous initial memory is bounded by \(Cd_0\).  Its
cutoff commutator uses the coarse global low \(F\)-block, but only with
the retained factor \(\lambda_1K_{{\rm Fgr},4}\).  On
\(\mathfrak S_{\rm gr}^+\) the physical metric is in the
noncollapsing exterior gauge, so the metric, target, and clock
coefficient differences occur with the physical-to-normalized factor
\(\lambda\).  Thus their contribution is bounded by
\[
 C\lambda(\tau)\left(
  \mathfrak G_6(\tau)
  +\mathfrak B_{{\rm gr},4}(\tau)
  +\left|\log\frac{\lambda_1}{\lambda_2}(\tau)\right|
  +\|R_1-R_2\|_{\mathfrak X_{\rm sc}^{6,\alpha}}\right).
\]
This estimate is localized to \(\mathfrak S_{\rm gr}^+\): on the
collapsing core, \(\lambda G^{-1}=O(1)\), so the argument does not yield
an integrable global metric forcing.

The cutoff-corridor geometric forcing and the global low lateral
\(F\)-value in \eqref{eq:localized-graft-F-coarse-memory} both carry
the scalar off-diagonal kernel and the factor \(\lambda_1\).  The
\(R\)-difference is obtained by direct
integration of its triangular ODE, and the phase profiles are
retained in their exact \(L^1\)-Duhamel form.
Now use \eqref{eq:buffered-two-state-stability},
\eqref{eq:buffered-two-state-graft-buffer-stability},
\eqref{eq:buffered-scale-time-difference}, and
\(\lambda\leq Ce^{-\tau}\).  For the nonseparated order-four source,
use the Abel factor from
\eqref{eq:effective-time-one-order-kernel}:
\[
 \int_{\tau_0}^{\tau}
 \lambda_1(q)\bigl(1+\mathfrak s(q,\tau)^{-1/2}\bigr)\,dq
 \leq C\left(\mathfrak s(\tau_0,\tau)
              +\mathfrak s(\tau_0,\tau)^{1/2}\right)\leq C .
\]
Thus the metric/clock contribution is finite at the asserted
order-five endpoint; the
factor \(\lambda_1\) times the off-diagonal kernel absorbs the coarse
growth of both the hybrid and global low \(F\)-blocks in
\eqref{eq:coarse-two-state-growth}; and
\eqref{eq:buffered-stability-phase-hypothesis} controls the phase
integral.  The order-five one-derivative Abel/Schauder estimate proves
\eqref{eq:localized-graft-F-uniform}.  Differentiating the localized
equation and using Taylor's integral remainder gives the identical
first-variation and Fr\'echet-remainder estimates.
\end{proof}

\begin{lemma}[Uniform outer-source map comparison after phase control]
\label{lem:outer-source-F-comparison}
Let \(\mathfrak S_{\rm out}\) be the uniformly locally finite union of
source-adapted atlas members meeting \(\supp(1-\eta)\), and let
\[
 \mathfrak S_{\rm out}\subset_{\rm buf}\mathfrak S_{\rm out}^+
 \subset_{\rm buf}\mathfrak S_{\rm out}^{++}
\]
denote its first and second inward/lateral buffer layers, where
\(\subset_{\rm buf}\) means positive uniform separation in the
 source-adapted atlas.  Let \(R_{\rm out}>0\) be one quarter of the
 least of these two buffer separations, and let
 \(q_{\rm out}:=p_5<\infty\) be the exponent supplied by
 Lemma~\ref{lem:fixed-bottleneck-Davies} at order five.  In the
fixed right-translated prepared charts, write
\begin{equation}\label{eq:outer-source-F-block}
 \mathfrak F_{{\rm out},5}^{\pm}(\tau)
 :=
 \|F_1-F_2\|_
   {\mathfrak X_{\rm sc}^{5,\alpha}(\mathfrak S_{\rm out}^+)}
 +\|F_1^{-1}-F_2^{-1}\|_
   {\mathfrak X_{\rm sc}^{5,\alpha}
      (F_1(\mathfrak S_{\rm out}^+)
       \cup F_2(\mathfrak S_{\rm out}^+))}.
\end{equation}
Under the hypotheses of
Lemma~\ref{lem:buffered-physical-two-state-stability}, including the
phase-tail hypothesis
\eqref{eq:buffered-stability-phase-hypothesis},
\begin{equation}\label{eq:outer-source-F-uniform}
 \sup_{\tau\geq\tau_0}
 \mathfrak F_{{\rm out},5}^{\pm}(\tau)
 \leq C_Dd_0 .
\end{equation}
This is an outer-region estimate, not a uniform bound for the global
map block \(\mathfrak F_5^{\rm glob}\).
\end{lemma}

\begin{proof}
Let \(w\) be the right-translated map difference and choose a fixed
cutoff \(\xi_{\rm out}\) equal to one on
\(\mathfrak S_{\rm out}^+\), supported in
\(\mathfrak S_{\rm out}^{++}\), with its inward derivative supported
on the deep side of one unused buffer layer.  We use the localized
unknown \(z=\xi_{\rm out}w\).  Allocate the other buffer layer to a
smooth weight \(\psi_{\rm out}\) which is zero on the collapsing core
and on \(\supp d\xi_{\rm out}\), equals \(R_{\rm out}\) on
\(\mathfrak S_{\rm out}^+\), and has gradient supported in the fixed
noncollapsing graft buffer.

For the \(F_1-F_2\) equation, take \(\mathcal L_\tau\) in
Lemma~\ref{lem:fixed-bottleneck-Davies} to be the state-\(1\)
right-translated Jacobi operator and take
\(\lambda=\lambda_1\); every state-\(2\) coefficient difference is
placed in the exact inhomogeneous polarized source.  The inverse-map,
sliced-variation, and mean-value remainder equations use their
corresponding right-translated operators, which satisfy the same
uniform package.

We verify the hypotheses of
Lemma~\ref{lem:fixed-bottleneck-Davies}.  The principal matrix of the
right-translated \(F\)-equation is
\(a_i=\lambda_i\acute G_i^{-1}\).  On
\(\supp d\psi_{\rm out}\) the anchored noncollapsing metric package
therefore gives
\[
 a_i(d\psi_{\rm out},d\psi_{\rm out})\leq C\lambda_i,
\]
and the first- and zeroth-order terms created by conjugation are
bounded by \(C\lambda_i\) there, with the required spatial jets.
Their time-oscillation modulus is the fixed-buffer modulus from the
anchored and source-adapted coefficient packages.  On the deep region
where \(\xi_{\rm out}=0\), away from \(\supp d\xi_{\rm out}\), replace
the coefficients by one fixed dissipative scalar-symbol extension.
This does not change the cutoff equation; every actual core
contribution has already become the commutator supported on
\(\supp d\xi_{\rm out}\).
Here \(a_i,b_i,c_i\) are the original coefficients.  After the
fixed-scale clock change, and after the combined spatial dilation and
clock change on a source-adapted member, the exact normalized operators
defined, respectively, in
\eqref{eq:fixed-bottleneck-normalized-operator} and
\eqref{eq:fixed-bottleneck-scale-normalized-operator} have hatted
coefficients satisfying the order-five normalized bounds and their
common effective-time modulus.  Thus the raw bottleneck inequalities
and the normalized Schauder estimates are verified for the correct,
distinct coefficient families.

It remains to estimate the genuine outer energy form.  Its covariant
rough part is dissipative.  In the pure outer region
\(\acute G_i=S_i=\lambda_i\Theta_i^*\bar g\), and radial tracking plus
the AC curvature estimates give for its Jacobi potential
\[
 \lambda_i|{\rm Rm}(S_i)|_{S_i}
 \leq
 \frac{C}{1+\bar f\circ\Theta_i}
 \leq Ce^{-\tau}.
\]
The uniform prepared \(C^1\) bounds for \(F_i\) absorb the two map
derivatives in the Jacobi endomorphism into the constant \(C\).
The positive part of the bounded transition error is
\(O(\lambda_i)\), and the right-translation error is bounded by the
integrable one-state phase tail.  Consequently
\[
 \beta_+(\tau)\leq
 C\bigl(e^{-\tau}+\lambda_1(\tau)+|c_1(\tau)|+|c_2(\tau)|\bigr)
 \in L^1([\tau_0,\infty)).
\]
Moreover, for \(i=1,2\),
\[
 \sup_{\tau_0\leq q\leq\tau<\infty}
 \int_q^\tau\lambda_i(r)\,dr
 \leq
 \sup_{q\geq\tau_0}\int_q^\infty Ce^{-r}\,dr
 \leq Ce^{-\tau_0}.
\]
This verifies the finite-total-clock hypothesis
\eqref{eq:fixed-bottleneck-total-effective-time} for each state.
The same estimate holds for the adjoint form on the boundaryless
source manifold, whose form domain is \(H^1\).
Thus the displayed form estimate verifies the energy hypothesis of
Lemma~\ref{lem:fixed-bottleneck-Davies}.  The weight
\(\psi_{\rm out}\) is constant on the collapsing core and on the pure
outer region, so the argument does not require a global estimate
\(\lambda_i\acute G_i^{-1}\simeq\lambda_i g^{-1}\).

We use the scale-adapted clause of that lemma for the unbounded output,
rather than treating its remote charts as faster than the bottleneck
clock.  If
\((\mathcal U,r_{\mathcal U})\) is a retained output member, AC radial
properness gives a smooth weight \(\psi_{\mathcal U}\) which agrees
with the fixed bottleneck weight through the graft buffer, then
increases in the pure outer region, and has value at least
\(c(R_{\rm out}+r_{\mathcal U})\) on \(\mathcal U\).  These weights
are chosen from the fixed common source atlas
and the fixed radial function \(\bar f\), independently of either
state; the identical family is therefore available for first
variations and Fr\'echet remainders.  In the radial
corridor one has \(\acute G_i=S_i\); the source-adapted ellipticity,
AC derivative bounds, and radial comparison therefore give
\eqref{eq:bottleneck-Davies-coefficient-condition}, with all order-five
jets, uniformly for \(\psi_{\mathcal U}\).  Here the relative equation
\(\partial_\tau F=\lambda\Delta_{\acute G,S}F\) has already factored
out the spatial phase transport: its spatial first-order coefficients
in this pure-outer corridor retain the factor \(\lambda\), while the
right-translation contribution to the unweighted bundle energy is the
integrable term already included in \(\beta_+\).  On an output member,
the exact state-\(i\) clock is
\[
 \vartheta_{i,\mathcal U}(q,\tau)
 :=
 \frac1{r_{\mathcal U}^{2}}
 \int_q^\tau\lambda_i(r)\,dr .
\]
This is precisely the clock in
\eqref{eq:fixed-bottleneck-scale-normalized-operator}, and it is
uniformly comparable to the frozen source-atlas clock by
\eqref{eq:reference-clock-comparison}.  Hence
\eqref{eq:fixed-bottleneck-radial-Duhamel} supplies the final
\(\mathfrak X_{\rm sc}^{5,\alpha}(\mathfrak S_{\rm out}^+)\)-norm:
the additional radial Gaussian absorbs every inverse power of the
slower remote-chart clock.  This radial \(H^{-1}\)-to-scaled-H\"older
step is used only for the separated homogeneous datum, the cutoff
commutator, and the separated geometric corridor source, all supported
where the corresponding radial weight vanishes.  The genuinely local
outer metric, target, and clock forcing has source order four for the
order-five output and is propagated by the one-derivative
source-atlas Abel estimate
\eqref{eq:source-atlas-one-order-Abel-block}; it is not inserted into
the separated radial kernel.  The smooth finite-dimensional phase
profiles remain in exact \(L^1\)-Duhamel form.
For an output member \((\mathcal U,r_{\mathcal U})\), put
\[
 \mathfrak s_{\mathcal U}(q,\tau)
 :=r_{\mathcal U}^{-2}\mathfrak s(q,\tau)
 =\vartheta_{1,\mathcal U}(q,\tau).
\]
The outer atlas is noncollapsing, so
\(r_{\mathcal U}\geq r_*>0\), and
\[
 -\partial_q\vartheta_{1,\mathcal U}(q,\tau)
 =
 \frac{\lambda_1(q)}{r_{\mathcal U}^{2}}.
\]
Consequently
\[
 \frac{\lambda_1}{r_{\mathcal U}^{2}}
 \left(1+\mathfrak s_{\mathcal U}^{-1/2}\right)
 =
 \lambda_1\left(
 r_{\mathcal U}^{-2}
 +r_{\mathcal U}^{-1}\mathfrak s^{-1/2}\right)
 \leq
 C\lambda_1\left(1+\mathfrak s^{-1/2}\right).
\]
Taking the supremum over the source atlas therefore gives precisely the
Abel kernel in the second integral below; no remote chart is treated as
having the bottleneck clock.

On the bounded transition and buffer part of
\(\mathfrak S_{\rm out}^{++}\), the anchored exterior estimate and
\eqref{eq:buffered-two-state-graft-buffer-stability} control the
source geometry.  On the pure outer part, where \(\eta=0\), one has
\(\acute G_i=S_i\).  The common Jacobi part is retained in the
covariant evolution operator just verified; after this polarization,
the remaining inhomogeneous metric, target, and clock differences
carry the integrable physical factor \(\lambda_i\).  The scale estimate
\eqref{eq:buffered-scale-time-difference}, direct integration of the
\(R\)-equation, and
Lemma~\ref{lem:localized-graft-F-comparison} give
\[
 \sup_{\tau\geq\tau_0}
 \|R_1(\tau)-R_2(\tau)\|_{\mathfrak X_{\rm sc}^{6,\alpha}}
 \leq C_Dd_0 .
\]
They therefore give the one-derivative Abel evolution-family bound
\begin{equation}\label{eq:outer-source-F-propagation}
\begin{split}
 \mathfrak F_{{\rm out},5}^{\pm}(\tau)
 \leq{}& C d_0
 +C\int_{\tau_0}^{\tau}|\delta c(q)|\,dq\\
 &+C\int_{\tau_0}^{\tau}
   \lambda_1(q)
   K_{\rm MR}^{(1)}\bigl(\mathfrak s(q,\tau)\bigr)\left(
   \mathfrak G_6+\mathfrak B_{{\rm gr},4}
   +\left|\log\frac{\lambda_1}{\lambda_2}\right|
   +\|R_1-R_2\|_{\mathfrak X_{\rm sc}^{6,\alpha}}
   \right)(q)\,dq\\
 &+C\int_{\tau_0}^{\tau}
   \lambda_1(q)
   K_{{\rm out}}\bigl(\mathfrak s(q,\tau)\bigr)
   \left(
    \mathfrak D_4^{\rm hyb}
    +\mathfrak F_5^{\rm glob}
   \right)(q)\,dq .
\end{split}
\end{equation}
Here
\[
 K_{\rm MR}^{(1)}(\sigma)=1+\sigma^{-1/2},
 \qquad \sigma>0,
\]
is used only under the locally integrable effective-time integral, and
\[
 K_{\rm out}(\sigma)
 :=(1+\sigma^{-q_{\rm out}})e^{-cR_{\rm out}^2/\sigma}.
\]
This is the scalar kernel furnished by the fixed noncollapsing
bottleneck.  The last
integral contains the divergence-form commutator
\(\lambda_1[\mathcal L,\xi_{\rm out}]w\), bounded in \(H^{-1}\) by
\(C\lambda_1\mathfrak F_5^{\rm glob}\), as well as the separated
geometric corridor source.  Thus its \(\lambda_1\)-factor is genuine.

The initial term is global and hence bounded by \(d_0\).  The first
integral is bounded by the phase-tail hypothesis.  For the second,
\[
 \int_{\tau_0}^{\tau}
 \lambda_1(q)K_{\rm MR}^{(1)}
       \bigl(\mathfrak s(q,\tau)\bigr)\,dq
 \leq C\left(
   \mathfrak s(\tau_0,\tau)
   +\mathfrak s(\tau_0,\tau)^{1/2}\right)\leq C.
\]
The phase-controlled physical estimates and the auxiliary
order-six \(R\)-bound therefore control this Abel term by \(C_Dd_0\).
Finally,
\(\mathfrak s(q,\tau)\leq Ce^{-q}\), so
\[
 K_{\rm out}\bigl(\mathfrak s(q,\tau)\bigr)
 \leq C(1+e^{q_{\rm out}q})e^{-c'e^q}.
\]
This superexponential factor, together with the displayed
\(\lambda_1\), absorbs the coarse growth of both
\(\mathfrak D_4^{\rm hyb}\) and \(\mathfrak F_5^{\rm glob}\) in
\eqref{eq:coarse-two-state-growth}, uniformly in the terminal time.
Thus \eqref{eq:outer-source-F-propagation} proves the map part of
\eqref{eq:outer-source-F-uniform}.  The inverse part follows from the
uniform lower bound for \(DF_i\), the identity
\(F_i^{-1}\circ F_i=\operatorname{Id}\), and the same buffered
composition estimate.  The argument is valid for first variations
and Fr\'echet remainders as well.
\end{proof}

\begin{lemma}[Physical graft Lipschitz estimate after phase control]
\label{lem:physical-graft-two-state-Lipschitz}
Let the two global solutions in
Lemma~\ref{lem:coarse-two-state-geometry} have initial distance \(d_0\).
Assume their common-margin ball retains the buffered physical
certificate in
Lemma~\ref{lem:buffered-physical-two-state-stability}.
Define
\[
 \mathcal E_{{\rm gr},i}
 =(\Phi_i^{-1})^*\mathcal G_{{\rm gr},i},\qquad i=1,2.
\]
Suppose, for some \(D,\vartheta>0\), that
\[
 \int_\tau^\infty|c_1-c_2|(s)\,ds
 \leq Dd_0e^{-\vartheta(\tau-\tau_0)}.
\]
Then
\begin{align}
 \supp(\mathcal E_{{\rm gr},1}-\mathcal E_{{\rm gr},2})
 &\subset
 \{c\Gamma e^\tau\leq\bar f\leq C\Gamma e^\tau\},
 \label{eq:two-graft-support}\\
 \sum_{\ell=0}^2
 |\bar\nabla^\ell
   (\mathcal E_{{\rm gr},1}-\mathcal E_{{\rm gr},2})|_{\bar g}
 &\leq C_Dd_0e^{-\tau},\label{eq:two-graft-C2}\\
 \max_{0\leq j\leq8}
 \|\mathcal Y_{j,\tau}(\mathbf z_1)
    -\mathcal Y_{j,\tau}(\mathbf z_2)\|
       _{\mathfrak C_{{\rm sc},0}^{2,\alpha}}
 &\leq C_Dd_0.\label{eq:two-column-C2}
\end{align}
\end{lemma}

\begin{proof}
The buffered physical lemma supplies precisely the input orders
required by the graft calculus:
\begin{equation}\label{eq:physical-graft-buffer-invocation}
 \sup_{\tau\geq\tau_0}
 \mathfrak B_{{\rm gr},4}(\tau)
 \leq C_Dd_0
\end{equation}
by
\eqref{eq:buffered-two-state-graft-buffer-stability}.  In particular,
the physical tensor difference is controlled in \(C^{6,\alpha}\), and
the gauge, inverse-gauge, and transported-marking differences are
controlled in \(C^{5,\alpha}\).  Thus every subsequent pullback is
taken at an order supported by these estimates.

The scale estimate
\eqref{eq:buffered-scale-time-difference} gives
\[
 \sup_{\tau\geq\tau_0}
 \left|\log\frac{\lambda_1}{\lambda_2}\right|
 \leq C_Dd_0 .
\]
Next subtract the \(R\)-equations in their fixed right-translated
scale-one charts.  Direct integration and the phase-tail hypothesis,
together with Lemma~\ref{lem:localized-graft-F-comparison}, give
\[
 \sup_{\tau\geq\tau_0}\left(
  \|R_1-R_2\|_{\mathfrak X_{\rm sc}^{5,\alpha}}
  +\mathfrak F_{{\rm gr},5}
 \right)
 \leq C_Dd_0 .
\]
For the graft pullback, only the finite source-adapted star in
Lemma~\ref{lem:localized-graft-F-comparison} is used.  For the
effective columns, whose pullback term is supported on the entire
noncompact set transported from \(\supp(1-\eta)\), invoke instead
Lemma~\ref{lem:outer-source-F-comparison}:
\[
 \sup_{\tau\geq\tau_0}
 \mathfrak F_{{\rm out},5}^{\pm}(\tau)\leq C_Dd_0.
\]
The finite graft-star estimate and the global outer-source estimate
are therefore supplied by their respective lemmas.
Together with
\eqref{eq:physical-graft-buffer-invocation}, these estimates verify
\eqref{eq:phase-controlled-graft-input-block}.

Apply the scale-sharp \(m=4\) clause
\eqref{eq:pure-graft-buffered-scale-sharp} of
Lemma~\ref{lem:same-order-pure-graft-difference}.  It gives
\eqref{eq:two-graft-C2} directly.  Its support proof uses the correct
dichotomy: after expansion, every pure-graft summand either contains a
derivative of \(\eta\), or contains a coefficient vanishing at
\(\eta=0,1\).  The two-sided annulus tracking
\eqref{eq:Phi-annulus-tracking} therefore gives
\eqref{eq:two-graft-support}.

Finally apply the \(C^{5,\alpha}\) scale-one composition estimate to
\eqref{eq:effective-column-zero}--%
\eqref{eq:effective-column-j}.  Since every
\(K_\tau(T)\) is supported in
\(\Phi_\tau(\supp(1-\eta))\), its map, inverse-map, and transported
cutoff differences are controlled by
\eqref{eq:outer-source-F-uniform}, together with the scale and
\(R\)-bounds above.  These second-order column expressions carry no
decaying tensor-rescaling factor, and hence give
\eqref{eq:two-column-C2}, with one constant common to all nine
columns.
\end{proof}

\begin{lemma}[Phase-conjugated fixed-core \(L^2\)-to-\(C^0\) estimate]
\label{lem:phase-conjugated-core-L2-C0}
Let
\[
 K_0\Subset K_1\Subset K_2\Subset M
\]
be fixed smooth relatively compact domains.  Let
\[
 u\in C([s_0,s_1];L^2(K_2))
 \cap L^2([s_0,s_1];H^1_{\rm loc}(K_2))
\]
be an energy solution on \(K_2\times[s_0,s_1]\),
\(s_1-s_0\leq1\), of the linear system
\begin{equation}\label{eq:phase-conjugated-core-linear-system}
 \partial_\tau u
 =
 a^{ij}\bar\nabla_i\bar\nabla_ju
 +A*\bar\nabla u+(B_0+qB_1)*u
 +\Lie_{V_{\rm ph}}u+F .
\end{equation}
Assume that \(a,A,B_0,B_1\) are strongly measurable in \(\tau\) with
values in the indicated spatial H\"older classes, that \(a\) is
uniformly elliptic, and that their spatial \(C^{1,\alpha}\) norms have
one common essential-supremum bound.  Assume also that
\[
 \begin{aligned}
 V_{\rm ph}(\tau)
 &=\sum_{\mu=1}^Jc_\mu(\tau)X_{\mu,\tau},&
 q(\tau)&=\sum_{\mu=1}^J|c_\mu(\tau)|,\\
 c_\mu&\in L^1([s_0,s_1])\quad(1\leq\mu\leq J).
 \end{aligned}
\]
where \(X_{\mu,\tau}\) is strongly measurable in \(\tau\) with a fixed
essential-supremum spatial \(C^{2,\alpha}\) bound on \(K_2\).  There is
a threshold \(\eta_{\rm core}>0\), depending only on the three fixed
domain buffers, the background geometry on \(K_2\), the ellipticity and
common coefficient bounds, and this fixed vector-field ceiling, such
that the following holds.  Suppose
\[
 \int_{s_0}^{s_1}q\leq\eta_{\rm core},
\]
and that \(F\) is strongly measurable with
\[
 \|F(\tau)\|_{C^{0,\alpha}(K_2)}\leq f(\tau),
 \qquad f\in L^1([s_0,s_1]).
\]
Then, whenever \(s_1-s_0\geq3/4\), one has
\begin{equation}\label{eq:phase-conjugated-core-L2-C0}
 \|u(s_1)\|_{C^0(K_0)}
 \leq C\left(
  \sup_{s_0\leq r\leq s_1}\|u(r)\|_{L^2(K_2)}
  +\int_{s_0}^{s_1}f(r)\,dr
 \right).
\end{equation}
The constant depends only on the fixed domains, background geometry,
ellipticity, and displayed spatial bounds, but not on
\(\sup_\tau q(\tau)\).
Equivalently, the estimate may first be proved for temporally smooth
coefficients and classical solutions and then extended to the stated
energy class by approximation.
\end{lemma}

\begin{proof}
Conjugate by the time-ordered flow generated by \(-V_{\rm ph}\), based at
\(s_0\).  The displacement and \(C^{2,\alpha}\) distortion of this
flow are bounded by \(C\int q\).  The choice of \(\eta_{\rm core}\) therefore
keeps the images of \(K_0\) and \(K_1\) inside the next larger domains.
Choose \(K_0\Subset K_0'\Subset K_1\) and decrease \(\eta_{\rm core}\) so that
the inverse phase image of \(K_0\) lies in \(K_0'\).  The conjugated
\(L^2(K_1)\) norm is then bounded by the original \(L^2(K_2)\) norm,
with a uniform Jacobian constant.
In the conjugated equation the Lie term is absent, while ellipticity
and the spatial coefficient bounds change by at most
\(\exp(C\int q)\).  The remaining \(qB_1\) term is retained as an
integrable zeroth-order coefficient and is handled by Gronwall.  No
pointwise bound for \(q\) enters.

We now prove the required interior bound without using a kernel gradient
estimate.  Write \(\widetilde u\) for the phase-conjugated section on
\(K_1\).  After decreasing \(\eta_{\rm core}\) as above, its equation on
\(K_1\times[s_0,s_1]\) has the form
\[
 \partial_\tau\widetilde u
 =\widetilde a^{ij}\bar\nabla_i\bar\nabla_j\widetilde u
  +\widetilde A*\bar\nabla\widetilde u
  +(\widetilde B_0+q\widetilde B_1)*\widetilde u
  +\widetilde F,
\]
where \(\widetilde a\) has the same uniform ellipticity constants, the
spatial \(C^{1,\alpha}\) norms of the displayed coefficients are
bounded by one fixed constant, and
\[
 \|\widetilde F(\tau)\|_{C^0(K_1)}\leq C f(\tau).
\]
All these assertions use only the bound for \(\int q\); no pointwise
bound for \(q\) is introduced.

We record the scalar-principal Kato calculation needed below.  Rewrite
the principal term in divergence form; the derivative of
\(\widetilde a\) is a bounded first-order coefficient.  Applying the
usual convex regularization
\(z_\varepsilon=(|\widetilde u|^2+\varepsilon^2)^{1/2}\), decomposing
\(\bar\nabla\widetilde u\) into its component parallel to
\(\widetilde u\) and its fiber-orthogonal component, and absorbing the
latter by the negative scalar-principal quadratic form gives, after
letting \(\varepsilon\downarrow0\),
\begin{equation}\label{eq:phase-core-scalar-Kato}
 \partial_\tau z
 -\bar\nabla_i(\widetilde a^{ij}\bar\nabla_j z)
 \leq C_0|\bar\nabla z|
      +(C_0+C_1q(\tau))z+|\widetilde F|
 \quad\text{weakly},
 \qquad z:=|\widetilde u|.
\end{equation}
The constants depend only on the ellipticity and the displayed spatial
bounds.  This calculation is valid for bundle-valued systems precisely
because the principal symbol is scalar on the fibers.  For an energy
solution it follows by Steklov averaging in time and the same convex
regularization, so no temporal modulus of the coefficients is used.

Put
\[
 H(\tau):=C_0(\tau-s_0)+C_1\int_{s_0}^{\tau}q(r)\,dr,
 \qquad y:=e^{-H}z,
\]
and set
\[
 g(\tau):=e^{-H(\tau)}
          \|\widetilde F(\tau)\|_{C^0(K_1)},
 \qquad
 G(\tau):=\int_{s_0}^{\tau}g(r)\,dr.
\]
Since \(H\) depends only on time,
\eqref{eq:phase-core-scalar-Kato} and the weak truncation rule imply
that
\[
 w:=(y-G)_+
\]
is a nonnegative weak subsolution of
\begin{equation}\label{eq:phase-core-homogeneous-subsolution}
 \partial_\tau w
 -\bar\nabla_i(\widetilde a^{ij}\bar\nabla_j w)
 \leq C_0|\bar\nabla w|
\end{equation}
on \(K_1\times[s_0,s_1]\).

For completeness, the interior mean-value estimate for
\eqref{eq:phase-core-homogeneous-subsolution} requires only bounded
measurable time dependence.  Indeed, test the inequality with
\(\chi^2w^{p-1}\), where \(\chi\) is a nested spatial--time cutoff.
Uniform ellipticity absorbs the term containing
\(C_0|\bar\nabla w|\) by Young's inequality.  The resulting
Caccioppoli inequality, the fixed Sobolev inequality on the buffered
charts, and the iteration \(p\mapsto p(1+2/\dim M)\) give
\begin{equation}\label{eq:phase-core-mean-value}
 \|w(s_1)\|_{C^0(K_0')}
 \leq C
 \|w\|_{L^2([s_0,s_1]\times K_1)} .
\end{equation}
Here the spatial cutoffs are chosen between
\(K_0'\Subset K_1\), and the final time cutoff is supported after
\(s_0\).  The fixed gap \(s_1-s_0\geq3/4\) makes all cutoff constants
uniform.  The estimate first holds for almost every terminal time; the
standard locally bounded representative furnished by the same
iteration, together with
\(\widetilde u\in C([s_0,s_1];L^2)\), gives it at \(s_1\).

Because \(s_1-s_0\leq1\), \(0\leq H\leq C(1+\eta_{\rm core})\), and the phase
Jacobian is uniformly bounded,
\[
 \begin{split}
 \|w\|_{L^2([s_0,s_1]\times K_1)}+G(s_1)
 &\leq C\left(
   \sup_{s_0\leq r\leq s_1}
       \|u(r)\|_{L^2(K_2)}
   +\int_{s_0}^{s_1}f(r)\,dr\right).
 \end{split}
\]
Since \(y\leq w+G\), equations
\eqref{eq:phase-core-mean-value} and the preceding display bound
\(\widetilde u(s_1)\) on \(K_0'\).  Undoing the phase map, whose inverse
image of \(K_0\) lies in \(K_0'\), proves
\eqref{eq:phase-conjugated-core-L2-C0}.  The entire argument is stable
under temporal mollification with the spatial bounds and
\(\|q\|_{L^1}\) preserved, which also proves the final approximation
clause.
\end{proof}

\begin{lemma}[Exact two-state algebra and uniform local Kato estimate]
\label{lem:uniform-two-state-Kato-ledger}
Let two prepared normalized states on a common finite interval satisfy
the exact adaptive equation \eqref{eq:adaptive-normalized}, with
coefficients \(c_i=(a_i,b_i)\), effective columns
\(\mathcal Y_{j,\tau}(\mathbf z_i)\), and normalized graft forcings
\(\E_i\), \(i=1,2\).  No two-state decay or barrier conclusion is
assumed.  Write
\[
 w=h_1-h_2,\qquad
 \delta c=c_1-c_2,\qquad q_i=|c_i|,\qquad
 \mathcal Y_j^{(i)}
 =\mathcal Y_{j,\tau}(\mathbf z_i),\qquad
 \Delta\mathcal Y_j
 =\mathcal Y_j^{(1)}-\mathcal Y_j^{(2)}.
\]
Recall the operators \(\mathscr T_j\) from
\eqref{eq:two-state-modulation-operators}, and define
\begin{align}
 \mathfrak y_{12}(\tau)
 &:=
 \max_{0\leq j\leq8}
 \|\Delta\mathcal Y_j\|_{\mathfrak C_{{\rm sc},0}^{2,\alpha}},
 \label{eq:two-state-column-size}\\
 q_{12}(\tau)
 &:=
 |\delta c(\tau)|+q_2(\tau)\mathfrak y_{12}(\tau).
 \label{eq:two-state-direct-amplitude}
\end{align}
Then the exact difference of the normalized equations is
\begin{equation}\label{eq:exact-two-state-h-equation}
 \begin{split}
 \partial_\tau w={}&
 \A w+\Q(h_1)-\Q(h_2)
 +\sum_{j=0}^8c_{1,j}\mathscr T_jw\\
 &+\sum_{j=0}^8\delta c_j
       \bigl(\mathcal Y_{j,\tau}(\mathbf z_1)
             +\mathscr T_jh_2\bigr)
 +\sum_{j=0}^8c_{2,j}\Delta\mathcal Y_j
 +(\E_1-\E_2).
 \end{split}
\end{equation}
Equation~\eqref{eq:exact-two-state-h-equation} displays every
outer-error term explicitly.

On a common scale-normalized \(C^2\) box, the prepared bounds give,
for \(0\leq\ell\leq2\),
\begin{equation}\label{eq:two-state-direct-column-bound}
 (1+\bar f)^{\ell/2}
 \left|
 \bar\nabla^\ell
 \bigl(\mathcal Y_j^{(1)}+\mathscr T_jh_2\bigr)
 \right|
 \leq C.
\end{equation}
At order zero its stipulated common ceiling is denoted
\(K_{\mathcal Y,0}^{(2)}\).  If \(g_i=\bar g+h_i\), exact
quasilinear polarization gives
\begin{equation}\label{eq:two-state-principal-polarization}
 \A w+\Q(h_1)-\Q(h_2)
 =
 g_1^{ab}\bar\nabla_a\bar\nabla_bw
 -\bar\nabla_{\bar\nabla\bar f}w
 +A_{12}*\bar\nabla w+B_{12}*w,
\end{equation}
where, on every scale-one annulus,
\begin{equation}\label{eq:two-state-polarized-coefficients}
 (1+\bar f)^{1/2}|A_{12}|
 +(1+\bar f)|B_{12}|
 \leq C,
\end{equation}
and the part of \(A_{12}\) multiplying \(\bar\nabla w\) is as small
as the common \(C^2\) box.

For every
\(\Lambda_{\rm ell}^{(2)},K_{\mathcal Y,0}^{(2)}<\infty\), there is an
admissible triple
\[
 (\delta_{\rm K}^{(2)},C_{\rm K}^{(2)},c_{\rm K}^{(2)})
 \in(0,\infty)^3,
\]
 depending only on these ceilings, the fixed FIK background and cutoff
 profiles, such that on every
 common box of size at most
\(\delta_{\rm K}^{(2)}\) and ellipticity at most
\(\Lambda_{\rm ell}^{(2)}\), whose order-zero direct tensors in
\eqref{eq:two-state-direct-column-bound} are bounded by
\(K_{\mathcal Y,0}^{(2)}\), the function
\(v=|w|_{\bar g}^2\) satisfies
\begin{equation}\label{eq:two-state-Kato}
 \begin{split}
 \mathscr P_1v
 \leq{}&
 C_{\rm K}^{(2)}q_{12}v^{1/2}
 +2v^{1/2}|\E_1-\E_2|_{\bar g}\\
 &-c_{\rm K}^{(2)}\,g_1^{ab}
   \langle\bar\nabla_aw,\bar\nabla_bw\rangle_{\bar g},
 \end{split}
\end{equation}
where
\begin{equation}\label{eq:two-state-scalar-operator}
 \mathscr P_1
 =
 \partial_\tau-g_1^{ab}\bar\nabla_a\bar\nabla_b
 +\bar\nabla_{V_1}
 -\left(
   \frac {C_{\rm K}^{(2)}}{1+\bar f}
   +C_{\rm K}^{(2)}q_1\right),
 \qquad
 V_1=(1+a_1)\bar\nabla\bar f
 -\sum_{j=1}^8b_{1,j}\chi_\tau W_j .
\end{equation}
\end{lemma}

\begin{proof}
The identity \eqref{eq:exact-two-state-h-equation} follows term by term
from
\[
 \begin{split}
 c_{1,j}\bigl(\mathcal Y_j^{(1)}+\mathscr T_jh_1\bigr)
 -c_{2,j}\bigl(\mathcal Y_j^{(2)}+\mathscr T_jh_2\bigr)
 ={}&
 c_{1,j}\mathscr T_jw\\
 &+\delta c_j
   \bigl(\mathcal Y_j^{(1)}+\mathscr T_jh_2\bigr)
 +c_{2,j}\Delta\mathcal Y_j .
 \end{split}
\]
The scale-normalized prepared bounds and
Lemma~\ref{lem:mode-growth} give
\eqref{eq:two-state-direct-column-bound}.  The normalized formulation
is essential: \(W_j=O(\bar f^{1/2})\), whereas
\(\bar\nabla h_2=O(\bar f^{-1/2})\) on scale-one annuli.

The exact formula \eqref{eq:Q-exact} and the mean-value identity
\[
 \Q(h_1)-\Q(h_2)
 =\int_0^1D\Q_{h_2+s w}[w]\,ds
\]
give \eqref{eq:two-state-principal-polarization} and
\eqref{eq:two-state-polarized-coefficients}.  For example,
\[
 (g_1^{ab}-g_2^{ab})\bar\nabla_a\bar\nabla_bh_2
 =\left(
  \int_0^1D(g^{-1})_{g_2+s(g_1-g_2)}[w]\,ds
 \right)*\bar\nabla^2h_2
\]
is zeroth order in \(w\), with
\((1+\bar f)|\bar\nabla^2h_2|\) controlled by the common box.

Take the \(\bar g\)-squared norm in the exact equation
\eqref{eq:exact-two-state-h-equation}.  Extract the transport terms as
in Lemma~\ref{lem:modulated-kato}.  Absorb the small
\(A_{12}*\bar\nabla w\) term into half the negative gradient form.
Inspection of the exact polarization shows that every nonbackground
coefficient in this squared-norm calculation is bounded by the common
\(C^2\) box and \(\Lambda_{\rm ell}^{(2)}\); the only direct
inhomogeneous tensor needed here is controlled by
\(K_{\mathcal Y,0}^{(2)}\).  No higher pre-radius coefficient ceiling
enters.  Uniform ellipticity and the fixed direct-column ceiling then
yield
\[
 (\delta_{\rm K}^{(2)},C_{\rm K}^{(2)},c_{\rm K}^{(2)})
\]
with
\eqref{eq:two-state-Kato}--\eqref{eq:two-state-scalar-operator}.
The constants involve neither a subordinate rate nor the graft
forcing.  This also proves
\(\mathscr A_{\rm K}
(\Lambda_{\rm ell}^{(2)},K_{\mathcal Y,0}^{(2)})\ne\varnothing\).
\end{proof}

\begin{lemma}[Uniform Kato estimates and package radius for subordinate
two-state rates]
\label{lem:uniform-two-state-package-radius}
Fix \(0<\theta_*<\beta\), a background ellipticity ceiling
\(\Lambda_{\rm ell}^{(2)}<\infty\), and a finite zeroth-order direct
column ceiling \(K_{\mathcal Y,0}^{(2)}<\infty\).  Fix an admissible
\[
 C_{\rm K}^{(2)}
 \in\mathscr A_{\rm K}
  (\Lambda_{\rm ell}^{(2)},K_{\mathcal Y,0}^{(2)}).
\]
Fix one witnessing pair
\((\delta_{\rm K}^{(2)},c_{\rm K}^{(2)})\)
for this fixed admissible \(C_{\rm K}^{(2)}\).  Then there exist
numerical constants
\(\delta_{\rm box}^{(2)}>0\) and
\(1\leq K_{\rm J}^{(2)},K_0^{(2)}<\infty\).
For later reference, record these new constants together with the
already fixed \(C_{\rm K}^{(2)}\) and the selected
\(c_{\rm K}^{(2)}\) as
\begin{equation}\label{eq:uniform-two-state-smallness-package}
 \delta_{\rm box}^{(2)}>0,\qquad
 C_{\rm K}^{(2)}<\infty,\qquad c_{\rm K}^{(2)}>0,\qquad
 1\leq K_{\rm J}^{(2)},K_0^{(2)}<\infty .
\end{equation}
The newly chosen constants depend only on \(\theta_*\), the two
ceilings, the fixed admissible \(C_{\rm K}^{(2)}\), its selected
witnessing pair, and the fixed background, and have the following
property.  On every common two-state prepared box of size at
most \(\delta_{\rm box}^{(2)}\), with ellipticity bounded by
\(\Lambda_{\rm ell}^{(2)}\) and the direct tensors in
\eqref{eq:two-state-direct-column-bound} bounded by
\(K_{\mathcal Y,0}^{(2)}\), the estimate
\eqref{eq:two-state-Kato}--\eqref{eq:two-state-scalar-operator}
holds with these same \(C_{\rm K}^{(2)}\) and
\(c_{\rm K}^{(2)}\).

There is also a finite number
\begin{equation}\label{eq:uniform-two-state-package-radius}
 \overline\Gamma_{\rm 2st}
 =\mathfrak G_{\rm 2st}
   (\theta_*,\Lambda_{\rm ell}^{(2)},
    K_{\mathcal Y,0}^{(2)},C_{\rm K}^{(2)})
 <\infty
\end{equation}
with the following property.  For every already fixed
\(\Gamma\geq\overline\Gamma_{\rm 2st}\) and every
\(0<\sigma'<\theta'<\theta_*\), the barrier comparison in
Lemma~\ref{lem:exact-two-state-barriers} can be carried out without
altering \(\Gamma\).  More precisely, for every finite post-radius
tuple consisting of the energy constant \(C_E\), the two-state
direct-tail constant \(C_{12}\), the graft constant
\(C_{\rm gr}^{(2)}\), finite graft-support annulus constants
\(0<c_{\rm supp}^{(2)}\leq C_{\rm supp}^{(2)}<\infty\),
and the entrance constant, there exist an outer amplitude
\(D_{\rm O}\), a delayed comparison time \(\widehat\tau_b\), and a
common amplitude multiplier for which the comparison closes.  The
 named one-state future-phase constants, evaluated at the frozen
 top rate,
 \[
  (C_P^*,c_P^*):=(C_P,c_P)\big|_{\theta=\theta_*},
 \]
 enter only the
 rate-independent base time below.  None of
\[
 C_E,\ C_{12},\ C_{\rm gr}^{(2)},\
 C_P^*,\ c_P^*,\ \text{the support constants, or the entrance constant}
\]
enters \(\overline\Gamma_{\rm 2st}\).  The named Kato and absorption
constants \(C_{\rm K}^{(2)},c_{\rm K}^{(2)},K_{\rm J}^{(2)},K_0^{(2)}\),
the radius, and the common-box, phase-budget, and base-time thresholds
are uniform over the subordinate cone.  Precisely
\[
 C_E,\ C_{12},\ C_{\rm gr}^{(2)},\ D_{\rm O},\
 \widehat\tau_b,\ \text{the common multiplier, and the final comparison
 constant}
\]
may depend on, and degenerate with, the subordinate pair.  All other
members of the post-radius input tuple are common package data.

After this \(\Gamma\) and its buffered atlas are fixed, choose smooth
domains
\[
 \{\bar f\leq\Gamma\}\Subset K_0^\Gamma
 \Subset K_1^\Gamma\Subset K_2^\Gamma
 \Subset\{\bar f<2\Gamma\}.
\]
Let \(\eta_{\rm core}^{(2)}>0\) be the threshold supplied by
Lemma~\ref{lem:phase-conjugated-core-L2-C0} for these domains and the
common coefficient ceilings, and let
\(\eta_{\rm ch}^{(2)}>0\) be a uniform phase-flow threshold preserving
the smaller/doubled scale-one chart inclusions, their inverse
inclusions, the residual-drift flow tubes used in derivative recovery,
and the required \(C^4\) distortion bounds.  Define the
post-radius, rate-independent phase threshold
\begin{equation}\label{eq:uniform-two-state-phase-threshold}
 \varepsilon_{\rm ph,*}^{(2)}
 :=
 \min\left\{
 1,\frac{\log2}{K_{\rm J}^{(2)}},
 \eta_{\rm core}^{(2)},\eta_{\rm ch}^{(2)}
 \right\}.
\end{equation}
It is fixed before any subordinate pair and is not an argument of
\(\mathfrak G_{\rm 2st}\).

For every fixed \(\Gamma\), every fixed \(c_{\rm supp}^{(2)}>0\), fixed
package phase budget
\(0<\varepsilon_{\rm ph}\leq\varepsilon_{\rm ph,*}^{(2)}\), and fixed
entrance-amplitude ceiling \(\varepsilon_{\rm ent}\), there is, in
addition, a rate-independent base-time threshold
\begin{equation}\label{eq:uniform-two-state-base-time}
 \tau_{\rm base}^{(2)}
 =\tau_{\rm base}^{(2)}
  (\theta_*,\Gamma,\Lambda_{\rm ell}^{(2)},
   K_{\mathcal Y,0}^{(2)},\varepsilon_{\rm ph},
    \varepsilon_{\rm ent},C_P^*,c_P^*,
   c_{\rm supp}^{(2)},\text{background})<\infty .
\end{equation}
The common box threshold, phase budget, and base-time threshold in
\eqref{eq:uniform-two-state-smallness-package}--%
\eqref{eq:uniform-two-state-base-time} work simultaneously for every
 \(0<\sigma'<\theta'<\theta_*\).  Only
\[
 C_E,\ C_{12},\ C_{\rm gr}^{(2)},\ D_{\rm O},\
 \widehat\tau_b,\ \text{the common multiplier, and the final comparison
 constant}
\]
may degenerate at the boundary of this open cone, including as
\(\sigma'\downarrow0\), \(\theta'-\sigma'\downarrow0\), or
\(\theta'\uparrow\theta_*\).
\end{lemma}

\begin{proof}
Select the witnessing
\((\delta_{\rm K}^{(2)},c_{\rm K}^{(2)})\) for the fixed admissible
\(C_{\rm K}^{(2)}\).  Its existence follows from
Lemma~\ref{lem:uniform-two-state-Kato-ledger} and the definition of
\(\mathscr A_{\rm K}\).  This selection precedes every barrier and
radius calculation.  Choose
\(K_{\rm J}^{(2)}\geq1\) large enough to absorb the
\(C_{\rm K}^{(2)}q_1\)-potential and every fixed phase-drift
coefficient, and then choose \(K_0^{(2)}\geq1\) large enough to absorb the
\(C_{\rm K}^{(2)}q_{12}\) direct source.

Choose once and for all
\(0<16\gamma_-<\gamma_+<1/4\) and set \(A_{\rm I}=1\).
For \(0\leq s\leq\theta_*\), let
\(F_s(r)=r^s-D_{\rm I}r^{s-1}\).  If
\(r\geq4D_{\rm I}\), then, uniformly for every metric \(g\) in the
fixed common \(\Lambda_{\rm ell}^{(2)}\)-ellipticity class,
\[
 F_s(r)\geq\frac34r^s,\qquad
 \frac{|\bar\nabla F_s(\bar f)|_g^2}{F_s(\bar f)}
 \leq C_{\rm grad}^*\bar f^{s-1},
\]
where \(C_{\rm grad}^*<\infty\) depends only on
\(\theta_*\) and the fixed ellipticity bounds, uniformly in \(s\).
Choose
\[
 D_{\rm I}>
 \theta_*(1+\theta_*)
 +8(C_{\rm K}^{(2)}+C_{\rm grad}^*+1).
\]
The background inner calculation has a strict margin after this
choice, while the coefficient of \(D_{\rm O}\) in the outer
calculation has a uniform positive margin; the latter estimate is
homogeneous in the outer amplitudes \(A_{\rm O}\) and \(D_{\rm O}\).
The scale-normalized Hessian-to-barrier ratios in both calculations are
uniformly bounded for \(s\in[0,\theta_*]\).  Compactness therefore
supplies positive thresholds
\(\delta_{\rm rad,I}^{(2)}\) and
\(\delta_{\rm rad,O}^{(2)}\), depending only on the displayed common
data, such that
\[
 (g^{-1}-\bar g^{-1})*\bar\nabla^2\mathcal B_{\rm I}^0
 \quad\hbox{and}\quad
 (g^{-1}-\bar g^{-1})*\bar\nabla^2\mathcal B_{\rm O}^0
\]
use at most one quarter of the corresponding strict radial margin
whenever the common \(C^2\) box has that size.  Freeze, before invoking
either barrier,
\begin{equation}\label{eq:two-state-Kato-radial-box-reduction}
 \delta_{\rm box}^{(2)}
 :=
 \min\{\delta_{\rm K}^{(2)},
       \delta_{\rm rad,I}^{(2)},
       \delta_{\rm rad,O}^{(2)}\}.
\end{equation}
Shrinking the Kato witness preserves the already selected
\(C_{\rm K}^{(2)}\) and \(c_{\rm K}^{(2)}\).
The calculation
\eqref{eq:radial-drift-calculation}--%
\eqref{eq:intermediate-margin} is then uniform for
\(0<s<\theta_*\).  All its remaining background coefficients are
continuous in \(s\) on the compact interval \([0,\theta_*]\), so one
finite lower bound for \(\bar f\) absorbs them.  Take
 \(\overline\Gamma_{\rm 2st}\) larger than that bound, the fixed
 critical region of \(\bar f\), and \(4D_{\rm I}\).

Fix \(\Gamma\geq\overline\Gamma_{\rm 2st}\) and its buffered atlas.
Choose the displayed domains \(K_i^\Gamma\), then select
\(\eta_{\rm core}^{(2)}\) from
Lemma~\ref{lem:phase-conjugated-core-L2-C0} and
 \(\eta_{\rm ch}^{(2)}\) from the uniform positive buffer and Lebesgue
 margins of the finite core cover and the uniformly locally finite
 scale-one dyadic atlas, including the inverse-map bounds and the
 fixed buffer fraction for the residual-drift tubes.  With
\eqref{eq:uniform-two-state-phase-threshold}, the phase flow preserves
both packages and
\(\exp(K_{\rm J}^{(2)}\varepsilon_{\rm ph})\leq2\).

With the package support constants
\(c_{\rm supp}^{(2)},C_{\rm supp}^{(2)}\) already fixed, now fix for
the subordinate pair the remaining finite inputs
\[
 (C_E,C_{12},C_{\rm gr}^{(2)},C_0).
\]
For a given \(s=\sigma'\), put
\[
 A_{\rm O}(s)=(\gamma_-\gamma_+)^{s/2}.
\]
Then
\(\gamma_-^s<A_{\rm O}(s)<\gamma_+^s\), while
\(A_{\rm O}(s)\geq
(\gamma_-\gamma_+)^{\theta_*/2}>0\).  Thus the two crossing
 inequalities are strict for every \(s>0\).  Now choose
\begin{equation}\label{eq:two-state-outer-amplitude-after-radius}
 D_{\rm O}\geq C_{\rm amp}^{(2)}
 \left(1+C_{\rm supp}^{(2)}
          \Gamma C_{\rm gr}^{(2)}\right)
\end{equation}
with \(C_{\rm amp}^{(2)}\) larger than the two annular-comparison
constants and the fixed outer radial/Kato constants.  The outer margin
\(D_{\rm O}/\bar f\) then absorbs both the background potential loss
and the graft source, while
\eqref{eq:two-state-Kato-radial-box-reduction} absorbs the
inverse-metric Hessian error.  Indeed, on the graft support
\[
 \frac{D_{\rm O}}{\bar f}
 \geq\frac{D_{\rm O}}
 {C_{\rm supp}^{(2)}\Gamma}e^{-\tau}
 \geq C\,C_{\rm gr}^{(2)}e^{-\tau}.
\]
This
choice does not enter the inner radial lower bound and hence does not
require enlarging \(\Gamma\).

With all those constants already fixed, choose one \emph{final}
\(\widehat\tau_b\geq1\) so that
\begin{equation}\label{eq:two-state-delayed-comparison-after-radius}
 \gamma_-e^{\widehat\tau_b}
 \geq4\max\{\Gamma,D_{\rm O}/A_{\rm O}(s)\},
\end{equation}
and the crossing, direct-tail, core-boundary, and phase-tail errors
are all smaller than
their corresponding strict margins.  The ratios of those errors to
the radial profiles contain either \(e^{-\widehat\tau_b}\) or
\(e^{-(\theta'-s)\widehat\tau_b}\).  Thus the choice is possible for
every \(0<s<\theta'<\theta_*\).  The support-separation, fixed-core,
and background coefficient requirements which involve the unshifted
time are rate-free.  Enlarge their finite maximum, still independently
of \(s\) and \(\theta'\), so that
\begin{equation}\label{eq:uniform-two-state-phase-tail-base-time}
 \begin{gathered}
 e^{\tau_{\rm base}^{(2)}}\geq2\Gamma,\qquad
 c_{\rm supp}^{(2)}e^{\tau_{\rm base}^{(2)}}\geq2,\\
 c_{\rm supp}^{(2)}\Gamma e^{\tau_{\rm base}^{(2)}}
 \geq2\gamma_+,\\
 C_P^*\varepsilon_{\rm ent}^{\,2}
      e^{-2\theta_*\tau_{\rm base}^{(2)}}
 +C_P^*e^{-c_P^*e^{\tau_{\rm base}^{(2)}}}
 \leq\varepsilon_{\rm ph}.
 \end{gathered}
\end{equation}
This defines the threshold in
\eqref{eq:uniform-two-state-base-time} and proves all asserted
uniformities.
\end{proof}

\begin{lemma}[Two-state equation and barriers]
\label{lem:exact-two-state-barriers}
Let the hypotheses of
Lemma~\ref{lem:physical-graft-two-state-Lipschitz} hold, and retain the
notation and exact Kato estimate of
Lemma~\ref{lem:uniform-two-state-Kato-ledger}.
For every finite \(S>\tau_0\), put
\[
 P_{12,S}(\tau)=\int_\tau^S q_{12}(s)\,ds .
\]

Fix \(0<\theta_*<\beta\).  Assume that the common box has size at most
\(\delta_{\rm box}^{(2)}\), ellipticity at most
\(\Lambda_{\rm ell}^{(2)}\), and direct-column ceiling at most
\(K_{\mathcal Y,0}^{(2)}\).  Assume also that the already frozen
package phase budget and base time satisfy
\begin{equation}\label{eq:two-state-uniform-package-inputs}
 0<\varepsilon_{\rm ph}\leq\varepsilon_{\rm ph,*}^{(2)},
 \qquad
 \tau_0\geq\tau_{\rm base}^{(2)}.
\end{equation}
For numbers
\(0<\sigma_-<\theta_-<\theta_*\), assume that the already fixed
package radius satisfies
\begin{equation}\label{eq:two-state-frozen-radius-hypothesis}
 \Gamma\geq\overline\Gamma_{\rm 2st},
\end{equation}
where the threshold is evaluated at this \(\theta_*\) and the fixed
data in
\eqref{eq:uniform-two-state-package-radius}; it is not reselected when
\((\sigma_-,\theta_-)\) is chosen.  Define
\begin{equation}\label{eq:two-state-pointwise-weight}
 \omega_{\sigma_-,\tau_0}(\tau,x)
 =
 \min\left\{
 e^{-\sigma_-(\tau-\tau_0)}
 (1+\bar f(x))^{\sigma_-},\,1
 \right\}.
\end{equation}
Fix arbitrary finite positive constants
\(C_E,C_{12},C_{\rm gr}^{(2)},C_0\) and fixed positive lower and
upper graft-support constants
\(c_{\rm supp}^{(2)},C_{\rm supp}^{(2)}\).  Assume
\begin{align}
 \|V(\tau)\|_{L^2_\nu}
 +\left(\int_\tau^S
        \|V(s)\|_{H^1_\nu}^2\,ds\right)^{1/2}
 &\leq C_Ed_0
 e^{-\theta_-(\tau-\tau_0)},
 \qquad V=\rho_\tau w,
 \label{eq:two-state-barrier-energy-tail}\\
 P_{12,S}(\tau)
 &\leq C_{12}d_0
 e^{-\theta_-(\tau-\tau_0)},
 \label{eq:two-state-direct-tail}\\
 \int_{\tau_0}^S q_1(s)\,ds
 &\leq\varepsilon_{\rm ph},
 \label{eq:two-state-background-phase-budget}\\
 \supp(\E_1-\E_2)
 &\subset
 \{c_{\rm supp}^{(2)}\Gamma e^\tau\leq\bar f
   \leq C_{\rm supp}^{(2)}\Gamma e^\tau\},
 \label{eq:two-state-graft-annulus}\\
 \sum_{\ell=0}^2|\bar\nabla^\ell(\E_1-\E_2)|
 &\leq C_{\rm gr}^{(2)}d_0e^{-\tau}.
 \label{eq:two-state-graft-size}
\end{align}
Assume also the entrance estimate
\begin{equation}\label{eq:two-state-C3-entrance}
 \sum_{\ell=0}^3|\bar\nabla^\ell w(\tau_0)|
 \leq C_0d_0
\end{equation}
in the scale-normalized sense of
\eqref{eq:scaled-tensor-holder}.  Then
\begin{equation}\label{eq:two-state-C0-barrier-result}
 |w(\tau,x)|
 \leq C d_0
 \min\left\{
 e^{-\sigma_-(\tau-\tau_0)}
 (1+\bar f(x))^{\sigma_-},\,1\right\}
\end{equation}
on \(M\times[\tau_0,S]\), with \(C\) independent of \(S\).
\end{lemma}

\begin{proof}
The exact difference identity, direct-column estimate, polarization,
and scalar Kato inequality have already been proved, with the selected
uniform constants, in
Lemma~\ref{lem:uniform-two-state-Kato-ledger}.  In particular
\(q_2\Delta\mathcal Y_j\) is measured by \(q_{12}\); it is not assigned
to \(\E_1-\E_2\), whose temporal scale is different.

If \(d_0=0\), then
\eqref{eq:two-state-barrier-energy-tail},
\eqref{eq:two-state-direct-tail},
\eqref{eq:two-state-graft-size}, and
\eqref{eq:two-state-C3-entrance} give respectively
\[
 V\equiv0,\qquad q_{12}\equiv0,\qquad
 \E_1-\E_2\equiv0,\qquad w(\tau_0)\equiv0.
\]
Uniqueness for the uniformly parabolic exact difference equation
\eqref{eq:exact-two-state-h-equation} gives \(w\equiv0\), and the
conclusion follows.  We therefore assume \(d_0>0\) in the barrier
argument below.

We give the comparison argument, taking care of the shift in
\(\omega_{\sigma_-,\tau_0}\).  Put
\(\widehat\tau=\tau-\tau_0\) and
\[
 J_1(\tau)
 =\exp\left(
   K_{\rm J}^{(2)}
   \int_{\tau_0}^{\tau}q_1(s)\,ds\right).
\]
Apply Lemma~\ref{lem:uniform-two-state-package-radius} with
\((\sigma',\theta')=(\sigma_-,\theta_-)\), after the constants
\(C_E,C_{12},C_{\rm gr}^{(2)}\) in the hypotheses have been fixed.
Choose the preliminary radial data and the \emph{final}
\(\widehat\tau_b\) in the order prescribed there, without changing the
frozen \(\Gamma\).

If \(S\leq\tau_0+\widehat\tau_b\), finite-horizon prepared stability on
\([\tau_0,S]\), with the constant for the fixed maximal length
\(\widehat\tau_b\), proves
\eqref{eq:two-state-C0-barrier-result} directly.  Hence assume
\(S>\tau_0+\widehat\tau_b\).  On
\([\tau_0,\tau_0+\widehat\tau_b]\), finite-horizon prepared stability
gives
\[
 \sum_{\ell=0}^2|\bar\nabla^\ell w|
 \leq C_{\widehat\tau_b}d_0.
\]
The shifted weight has a positive lower bound on this fixed interval.
Now multiply the preliminary
\(A_{\rm I},D_{\rm I},A_{\rm O},D_{\rm O}\) by one common factor
\(L_b\geq1\), chosen large enough that both radial branches dominate
the preceding finite-horizon bound at
\(\tau_0+\widehat\tau_b\).  This common multiplication leaves
\[
 \frac{D_{\rm I}}{A_{\rm I}},\qquad
 \frac{A_{\rm O}}{A_{\rm I}},\qquad
 \frac{D_{\rm O}}{A_{\rm O}}
\]
unchanged.  It therefore preserves the already verified radius,
positivity, and crossing inequalities, while improving the graft and
additive-tail margins.  The delayed time is not enlarged after this
step.

With these final amplitudes define
\begin{align*}
 B_{\rm I}^0
 &=d_0e^{-\sigma_-\widehat\tau}
   \left(A_{\rm I}\bar f^{\sigma_-}
         -D_{\rm I}\bar f^{\sigma_--1}\right),\\
 B_{\rm O}^0
 &=d_0\left(A_{\rm O}-D_{\rm O}\bar f^{-1}\right),\\
 B_{\rm I}
 &=J_1B_{\rm I}^0-K_0^{(2)}P_{12,S},\\
 B_{\rm O}
 &=J_1B_{\rm O}^0-K_0^{(2)}P_{12,S}.
\end{align*}
The crossover scale is \(e^{\widehat\tau}\), not \(e^\tau\).
We perform the barrier comparison on
\(\{\bar f\geq\Gamma\}\times
 [\tau_0+\widehat\tau_b,S]\).

The radial identities
\eqref{eq:radial-drift-calculation} and
\eqref{eq:outer-margin} are unchanged.  The derivative
\[
 \partial_\tau[-K_0^{(2)}P_{12,S}]=K_0^{(2)}q_{12}
\]
absorbs the first term on the right of
\eqref{eq:two-state-Kato}.  Because
\(\theta_->\sigma_-\), \eqref{eq:two-state-direct-tail} makes this
negative correction strictly smaller than both positive radial
profiles.  The radius hypothesis
\eqref{eq:two-state-frozen-radius-hypothesis} is already frozen.
Retain the crossover data and final delayed time already fixed by
Lemma~\ref{lem:uniform-two-state-package-radius}.
The term \(\E_1-\E_2\) is absorbed by the outer
\(D_{\rm O}d_0/\bar f\) margin using
\eqref{eq:two-state-graft-annulus}--%
\eqref{eq:two-state-graft-size}.  The same crossover calculation as
\eqref{eq:barrier-crossing} gives
\[
 B_{\rm O}<B_{\rm I}\quad
   \text{on }\{\bar f=\gamma_+e^{\widehat\tau}\},\qquad
 B_{\rm I}<B_{\rm O}\quad
   \text{on }\{\bar f=\gamma_-e^{\widehat\tau}\}.
\]
Indeed, on \(\bar f=\gamma e^{\widehat\tau}\),
\[
 B_{\rm I}^0
 =d_0A_{\rm I}\gamma^{\sigma_-}
  +O(d_0e^{-\widehat\tau}),\qquad
 B_{\rm O}^0
 =d_0A_{\rm O}+O(d_0e^{-\widehat\tau}),
\]
and the already fixed relation
\[
 \gamma_-^{\sigma_-}
 <\frac{A_{\rm O}}{A_{\rm I}}
 <\gamma_+^{\sigma_-}
\]
gives the two strict crossings.  This is exactly the scale at which the two branches
of \(\omega_{\sigma_-,\tau_0}\) meet.  The graft annulus lies farther
out, at \(\bar f\simeq\Gamma e^{\tau_0}e^{\widehat\tau}\), where the
outer \(d_0/\bar f\) margin has the \(d_0e^{-\tau}\) size required by
\eqref{eq:two-state-graft-size}.
Define the shifted glued barrier only on this exterior:
\begin{equation}\label{eq:two-state-glued-barrier}
 \mathcal B_{12}(\tau,x)
 =
  \begin{cases}
  B_{\rm I}(\tau,x),
   &\Gamma\leq\bar f(x)\leq\gamma_-e^{\widehat\tau},\\
  \min\{B_{\rm I}(\tau,x),B_{\rm O}(\tau,x)\},
   &\gamma_-e^{\widehat\tau}\leq\bar f(x)
       \leq\gamma_+e^{\widehat\tau},\\
  B_{\rm O}(\tau,x),
   &\bar f(x)\geq\gamma_+e^{\widehat\tau}.
 \end{cases}
\end{equation}
The inner and outer branches are positive on their stated exterior
regions.
The strict crossing inequalities give the correct one-sided ordering
at the two interfaces.  Hence \(\mathcal B_{12}\) is positive and
locally Lipschitz, and \(\mathcal B_{12}^2\) is a viscosity
supersolution of \eqref{eq:two-state-Kato}.

It remains only to supply the inner boundary value, without invoking
the conclusion to be proved.  Use the explicit hypothesis
\eqref{eq:two-state-barrier-energy-tail}.
 Use the already frozen package domains
\[
 \{\bar f\leq\Gamma\}\Subset K_0^\Gamma
 \Subset K_1^\Gamma\Subset K_2^\Gamma
 \Subset\{\bar f<2\Gamma\}.
\]
For every \(\tau\geq\tau_0\), the base-time inequalities give, on
\(K_2^\Gamma\),
\[
 \bar f<2\Gamma\leq e^{\tau_0}\leq e^\tau,
 \qquad
 \bar f<2\Gamma
 \leq c_{\rm supp}^{(2)}\Gamma e^{\tau_0}
 \leq c_{\rm supp}^{(2)}\Gamma e^\tau .
\]
Consequently \(\rho_\tau\equiv1\) and
\((\E_1-\E_2)|_{K_2^\Gamma}\equiv0\).  On \(K_2^\Gamma\) the polarized equation
\eqref{eq:two-state-principal-polarization} has the form
\eqref{eq:phase-conjugated-core-linear-system}.  The one-state phase
budget is \eqref{eq:two-state-background-phase-budget}, the direct
 source is bounded by \(Cq_{12}\) by
\eqref{eq:two-state-direct-column-bound}.  Hence
Lemma~\ref{lem:phase-conjugated-core-L2-C0}, applied on
\([\tau-1,\tau]\), gives for \(\widehat\tau\geq1\)
\begin{equation}\label{eq:two-state-core-L2-to-C0}
 \begin{split}
 \sup_{\{\bar f\leq\Gamma\}}|w(\tau)|
 \leq C\bigg(&
 \sup_{\tau-1\leq s\leq\tau}\|V(s)\|_{L^2_\nu}
 +\int_{\tau-1}^{\tau}q_{12}(s)\,ds\bigg).
 \end{split}
\end{equation}
 The first term in
\eqref{eq:two-state-core-L2-to-C0} is controlled by
\eqref{eq:two-state-barrier-energy-tail}, and the second by
\[
 \int_{\tau-1}^{\tau}q_{12}(s)\,ds
 \leq P_{12,S}(\tau-1)
 \leq Cd_0e^{-\theta_-(\tau-1-\tau_0)}.
\]
Thus
\[
 \sup_{\{\bar f\leq\Gamma\}}|w(\tau)|
 \leq Cd_0e^{-\theta_-\widehat\tau}.
\]
The defining core-boundary condition in the final choice of
\(\widehat\tau_b\) makes the preceding late-time bound no larger than
the inner barrier on \(\{\bar f=\Gamma\}\).  It therefore supplies the
inner spatial boundary.
The finite-horizon estimate at
\(\tau_0+\widehat\tau_b\) supplies the initial boundary for the
late-time comparison.  Add the exhaustion corrector
\[
 z_\delta
 =\delta J_1e^{L(\widehat\tau-\widehat\tau_b)}(1+\bar f),
\]
apply comparison on compact exhaustions, and let first the
exhaustion radius tend to infinity and then \(\delta\downarrow0\).
This gives on \(\{\bar f\geq\Gamma\}\)
\[
 |w|\leq\mathcal B_{12}
 \leq Cd_0
 \min\{e^{-\sigma_-\widehat\tau}
       (1+\bar f)^{\sigma_-},1\}.
\]
Indeed, the inner branch is bounded by the first weighted profile.  On
the outer branch the weight is \(1\) when
\(\bar f\geq e^{\widehat\tau}\), while on
\(\gamma_+e^{\widehat\tau}\leq\bar f\leq e^{\widehat\tau}\) it has the
positive lower bound \(c(\sigma_-,\gamma_+)\); hence the constant outer
profile has the same displayed bound.  On the overlap the crossing
inequalities give the smaller of the two.  Notice that neither radial
branch is asserted to be a supersolution outside the region assigned to it in
\eqref{eq:two-state-glued-barrier}.
On \(\{\bar f\leq\Gamma\}\), the separately proved core estimate gives
\[
 |w|\leq Cd_0e^{-\theta_-\widehat\tau}
 \leq Cd_0\omega_{\sigma_-,\tau_0},
\]
because \(\theta_->\sigma_-\) and this core is fixed.  Combining the
core and exterior bounds, and using the earlier finite-horizon
estimate when \(0\leq\widehat\tau\leq\widehat\tau_b\), gives the
claimed bound on all of \(M\times[\tau_0,S]\).
This is \eqref{eq:two-state-C0-barrier-result}.
\end{proof}

\begin{lemma}[Derivative recovery with an \(L^1\) difference source]
\label{lem:two-state-linear-derivative-recovery}
Under the hypotheses and conclusions of
Lemma~\ref{lem:exact-two-state-barriers}, assume that the direct
columns in \eqref{eq:exact-two-state-h-equation}, including the
normalized nonzero tensors
\(\Delta\mathcal Y_j/\mathfrak y_{12}\), have uniformly bounded
scale-normalized \(C^{2,\alpha}\) norms.  (The prepared
\(k+2\)-topology, \(k\geq12\), supplies more derivatives than this.)
When \(\mathfrak y_{12}=0\), every
\(\Delta\mathcal Y_j\) vanishes and the displayed normalized tensor is
defined to be zero.
Then
\begin{equation}\label{eq:two-state-C2-barrier-result}
 \sum_{\ell=0}^2|\bar\nabla^\ell w(\tau,x)|
 \leq C d_0
 \min\left\{
 e^{-\sigma_-(\tau-\tau_0)}
 (1+\bar f(x))^{\sigma_-},\,1\right\}
\end{equation}
on \(M\times[\tau_0,S]\), with \(C\) independent of \(S\).
\end{lemma}

\begin{proof}
The dependence on the direct difference amplitude must remain
linear.  Applying the squared Bernstein inequality
\eqref{eq:Bernstein-differential} without modification would give
\((\int q_{12})^{1/2}\), and hence only a
\(d_0^{1/2}\) estimate.  We avoid that loss by first solving off the
spatially regular \(L^1\)-in-time source.

Write the direct, nontransport source in
\eqref{eq:exact-two-state-h-equation} as
\begin{equation}\label{eq:two-state-direct-source-defined}
 \mathcal S_{12}(\tau)
 =
 \sum_{j=0}^8\delta c_j
       \bigl(\mathcal Y_{j,\tau}(\mathbf z_1)
             +\mathscr T_jh_2\bigr)
 +\sum_{j=0}^8c_{2,j}\Delta\mathcal Y_j .
\end{equation}
Thus
\(\|\mathcal S_{12}(\tau)\|_{\mathfrak C_{{\rm sc},0}^{2,\alpha}}
\leq Cq_{12}(\tau)\).
Let \(\mathcal P_{1,\tau}\), based at
\(\mathcal P_{1,\tau_0}=\operatorname{Id}\), be the time-ordered flow
generated by the negative of the first-state phase-transport field
\[
 V_{{\rm ph},1}
 =-a_1\bar\nabla\bar f
   +\sum_{j=1}^8b_{1,j}\chi_\tau W_j .
\]
Put
\[
 \widetilde w(\tau)=\mathcal P_{1,\tau}^*w(\tau),\qquad
 \widetilde{\mathcal S}_{12}(\tau)
   =\mathcal P_{1,\tau}^*\mathcal S_{12}(\tau),\qquad
 \widetilde{\Delta\E}(\tau)
   =\mathcal P_{1,\tau}^*(\E_1-\E_2)(\tau).
\]
This removes the full Lie-transport part of the first-state phase.
The undifferentiated term \(a_1w\) from
\(\mathscr T_0w=w-\Lie_{\bar\nabla\bar f}w\) remains as an
integrable-in-time zeroth-order coefficient.  The scale-normalized
\(C^4\) distortion of
\(\mathcal P_{1,\tau}\) is bounded by
\(\exp(C\int q_1)\), and no time derivative of either \(c_1\) or
\(\delta c\) is formed.

We now specify the backward rescaling, including its truncated back
face.  Given a target point \(y_*\) at time \(\tau_*\), put
\(x_0=\mathcal P_{1,\tau_*}^{-1}(y_*)\), so that \(x_0\) is its
coordinate in the phase-conjugated gauge, and set
\begin{equation}\label{eq:two-state-rescaled-time}
 \tau(s)=\tau_*-\log(1-s),\qquad
 s_-=\max\{-1,\,1-e^{\tau_*-\tau_0}\},\qquad s\in[s_-,0],
\end{equation}
and let
\[
 \phi_s=\varphi_{-\log(1-s)},\qquad
 \mathcal H(s)
 =(1-s)\phi_s^*\widetilde w(\tau(s)).
\]
In particular, the tensor being rescaled is the phase-conjugated
tensor \(\widetilde w\), not the unconjugated difference \(w\).
To fix the spatial domains, let
\[
 \widehat g_{1,\tau_*}(s)
 :=
 (1-s)\phi_s^*\mathcal P_{1,\tau(s)}^*
       \bigl(\bar g+h_1(\tau(s))\bigr),
 \qquad
 \widehat g_{1,*}:=\widehat g_{1,\tau_*}(0).
\]
Put
\begin{equation}\label{eq:two-state-local-AC-scale}
 L_*:=1+\bar f(y_*)\simeq1+\bar f(x_0),
\end{equation}
where the comparison follows from the phase threshold and radial
tracking.  Let \(\mathcal A_*^+\) be a fixed-radius buffered chart for
\(L_*^{-1}\widehat g_{1,*}\), with a fixed positive containment margin.
Inside it choose once and for all \emph{raw} radii
\(0<r_1<r_2\), independent of \(x_0,\tau_*\), and define
\begin{equation}\label{eq:two-state-nested-rescaled-balls}
 B_i:=B_{\widehat g_{1,*}}(x_0,r_i),\quad i=1,2,
 \qquad
 \overline{B_1}\Subset B_2,\qquad
 Q_2:=B_2\times[s_-,0].
\end{equation}
The balls \(B_i\) are defined using the unrescaled metric
\(\widehat g_{1,*}\).  Replacing \(\widehat g_{1,*}\) by
\(L_*^{-1}\widehat g_{1,*}\) in their definition would also require
rescaling the parabolic clock.  Instead, the entire moving tube of
\(\overline{B_2}\) is kept inside the scale-one buffer
\(\mathcal A_*^+\).  The transformed center follows the original
trajectory
\[
 y(s):=\mathcal P_{1,\tau(s)}(\phi_s(x_0)),\qquad y(0)=y_*.
\]
For every \(x\in B_2\), also put
\[
 Y_s(x):=\mathcal P_{1,\tau(s)}(\phi_s(x)).
\]
The phase distortion bound
\(\exp(C\int_{\tau_0}^{\tau_*}q_1)\), together with the soliton
scale-one chart comparison on \(-1\leq s\leq0\), ensures that
\(\mathcal P_{1,\tau(s)}(\phi_s(B_2))\) remains in the corresponding
doubled chart about \(y(s)\); conversely a fixed smaller ball about
\(y(s)\) pulls back into \(B_1\).  Thus all coefficient estimates
 below hold on the fixed nested cylinders, and
 \(\widehat\nabla\) will denote the connection of
 \(\widehat g_{1,*}\).
The same shrinker radial tracking, the fixed scale-one radius of
\(B_2\), and the chart-preserving phase threshold
\(\eta_{\rm ch}^{(2)}\) give one \(C_{\rm wt}\geq1\) such that
\begin{equation}\label{eq:two-state-transported-radial-weight}
 \begin{split}
 C_{\rm wt}^{-1}e^{-(\tau_*-\tau_0)}
 (1+\bar f(y_*))
 \leq{}&
 e^{-(\tau(s)-\tau_0)}
 (1+\bar f(Y_s(x)))\\
 \leq{}&
 C_{\rm wt}e^{-(\tau_*-\tau_0)}
 (1+\bar f(y_*))
 \end{split}
\end{equation}
for \(x\in B_2\) and \(s\in[s_-,0]\).  Since
\(0<\sigma_-<\theta_*\), this implies
\begin{equation}\label{eq:two-state-transported-weight-comparison}
 C_{\rm wt}^{-\theta_*}
 \omega_{\sigma_-,\tau_0}(\tau_*,y_*)
 \leq
 \omega_{\sigma_-,\tau_0}(\tau(s),Y_s(x))
 \leq
 C_{\rm wt}^{\theta_*}
 \omega_{\sigma_-,\tau_0}(\tau_*,y_*).
\end{equation}
Define the transformed sources, noting that no additional time factor
occurs, by
\begin{equation}\label{eq:two-state-transformed-sources}
 \widehat{\mathcal S}_{12}(s)
 =\phi_s^*\widetilde{\mathcal S}_{12}(\tau(s)),\qquad
 \widehat{\Delta\E}(s)
 =\phi_s^*\widetilde{\Delta\E}(\tau(s)).
\end{equation}
Indeed, the factor \(1-s\) in \(\mathcal H\) cancels
\(d\tau/ds=(1-s)^{-1}\) in every additive tensor source.

For completeness, phase conjugation and soliton rescaling do not
commute exactly.  With pullback of vector fields understood in the
usual sense, define the residual drift
\begin{equation}\label{eq:two-state-residual-rescaled-drift}
 \widehat Z_{1,s}
 =(1-s)^{-1}\phi_s^*
 \bigl(\bar\nabla\bar f
       -\mathcal P_{1,\tau(s)}^*\bar\nabla\bar f\bigr).
\end{equation}
It is exactly the first-order field left after the pure soliton drift
has been cancelled: differentiating \(\phi_s^*\) contributes
\(+\bar\nabla_{\bar\nabla\bar f}\), whereas the phase-conjugated
drift is
\(-\bar\nabla_{\mathcal P_{1,\tau(s)}^*\bar\nabla\bar f}\);
the difference of Lie and covariant derivatives is absorbed into the
zeroth-order coefficient below.  Differentiating the phase-flow ODE
and using the scale-normalized mode bounds gives the correctly typed
symbol estimate
\begin{equation}\label{eq:two-state-residual-drift-bound}
 \|\widehat Z_{1,s}\|_
 {C^{3,\alpha}(\mathcal A_*^+,L_*^{-1}\widehat g_{1,*})}
 \leq C\int_{\tau_0}^{\tau(s)}q_1(r)\,dr
 \leq C\varepsilon_{\rm ph},
\end{equation}
and the same calculation gives the required bounded scale-normalized
jets of all other transformed fixed lower-order coefficients.  Although
the raw zeroth-order size of \(\widehat Z_{1,s}\) may be
\(O(\varepsilon_{\rm ph}\sqrt{L_*})\), only the scale-normalized bound
\eqref{eq:two-state-residual-drift-bound} is used.

Let \(\mathcal Q_s\) be the terminally normalized residual flow
\begin{equation}\label{eq:two-state-residual-drift-flow}
 \partial_s\mathcal Q_s
 =-\widehat Z_{1,s}\circ\mathcal Q_s,\qquad
 \mathcal Q_0=\operatorname{Id}.
\end{equation}
Its relative displacement is
\(O(\varepsilon_{\rm ph}\sqrt{L_*})\) in the raw metric.  The
variational equations for \(D\mathcal Q_s\),
\(D\mathcal Q_s^{-1}\), and their required higher raw jets depend on
the scale-normalized derivatives of \(\widehat Z_{1,s}\).  In covariant
form they also contain curvature contractions with
\(\widehat Z_{1,s}\); the AC curvature-symbol decay and
\(|\widehat Z_{1,s}|=O(\varepsilon_{\rm ph}\sqrt{L_*})\) make those
contractions uniformly bounded.  Thus the jets have uniform bounds by
\eqref{eq:two-state-residual-drift-bound} and the fixed AC package.  The already fixed
phase/chart threshold was selected below this fixed buffer fraction in
the package construction and therefore keeps \(\mathcal Q_s(B_2)\) inside
\(\mathcal A_*^+\) for all \(s\in[s_-,0]\).
Define the moving sampling tube
\begin{equation}\label{eq:two-state-residual-drift-tube}
 \mathcal T_2
 :=
 \left\{(s,\mathcal Q_s(x)):
        s_-\leq s\leq0,\ x\in B_2\right\}
 \subset [s_-,0]\times\mathcal A_*^+ .
\end{equation}
After increasing \(C_{\rm wt}\) by one fixed factor,
\eqref{eq:two-state-transported-radial-weight} and
\eqref{eq:two-state-transported-weight-comparison} remain valid with
\(x\) replaced by \(\mathcal Q_s(x)\).  Indeed, the preceding
scale-normalized displacement is a fixed small fraction of the
\(L_*^{-1}\widehat g_{1,*}\)-buffer, while the scale-normalized first
derivative of \(1+\bar f\) is uniformly bounded there.  Thus both the
radial weight and the graft-annulus alternative used below are
controlled on the actual tube \(\mathcal T_2\), not merely on the
unmoved cylinder \(Q_2\).

Define
\begin{equation}\label{eq:two-state-rescaled-a1}
 \widehat a_1(s)
 :=(1-s)^{-1}a_1(\tau(s)).
\end{equation}
Then \(|\widehat a_1|\leq\widehat q_1\), where
\(\widehat q_1=(1-s)^{-1}q_1(\tau(s))\), and
\begin{equation}\label{eq:two-state-rescaled-a1-integral}
 \int_{s_-}^{0}|\widehat a_1(s)|\,ds
 \leq\int_{\tau_0}^{\tau_*}q_1(\tau)\,d\tau .
\end{equation}
Since
\(-1\leq s\leq0\), the displayed factor \((1-s)^{-1}\) is harmless.
Thus, on the cylinder \(Q_2\) in
\eqref{eq:two-state-nested-rescaled-balls}, the completely specified equation
for \(\mathcal H\) has the form
\begin{equation}\label{eq:conjugated-two-state-equation}
 \partial_s\mathcal H-\widehat a^{ab}
 \widehat\nabla_a\widehat\nabla_b\mathcal H
 =
 \Lie_{\widehat Z_{1,s}}\mathcal H
 +\widehat A*\widehat\nabla\mathcal H
 +\widehat B_0*\mathcal H
 +\widehat a_1\,\widehat B_1*\mathcal H
 +\widehat{\mathcal S}_{12}
 +\widehat{\Delta\E}.
\end{equation}
Pull back once more by the residual flow:
\begin{equation}\label{eq:two-state-drift-straightened-tensor}
 \mathcal K(s):=\mathcal Q_s^*\mathcal H(s).
\end{equation}
By \eqref{eq:two-state-residual-drift-flow},
\(\partial_s(\mathcal Q_s^*\mathcal H)
=\mathcal Q_s^*(\partial_s\mathcal H-
\Lie_{\widehat Z_{1,s}}\mathcal H)\), so the Lie term in
 \eqref{eq:conjugated-two-state-equation} cancels exactly.  Because the
sampling tube \(\mathcal T_2\) lies in the fixed scale-one buffer, the
pulled-back equation lives on the fixed raw cylinder
\(B_2\times[s_-,0]\), its principal matrix is uniformly elliptic, and
all pulled-back fixed lower-order coefficients have uniform raw
\(C^{2,\alpha}\) bounds.  The scalar
\(\widehat a_1\) is controlled only in \(L^1_s\), as in
\eqref{eq:two-state-rescaled-a1-integral}; it is not placed in a
pointwise coefficient bound.  With
\[
 \widetilde{\mathcal S}_{12}:=
 \mathcal Q_s^*\widehat{\mathcal S}_{12},\qquad
 \widetilde{\Delta\E}:=
 \mathcal Q_s^*\widehat{\Delta\E},
\]
the uniform jet bounds for \(\mathcal Q_s^{\pm1}\) give
\[
 \|\widetilde{\mathcal S}_{12}(s)\|_{C^{2,\alpha}(B_2)}
 \leq Cq_{12}(\tau(s)).
\]

On \(B_2\), let
\(\widetilde{\mathcal U}_{\tau_*}(s,r)\) be the Dirichlet evolution
family for this drift-straightened homogeneous operator.  Solve with
zero initial and lateral data for
\(\widetilde{\mathcal S}_{12}\).  If the source does not satisfy the
higher corner compatibility conditions, approximate it by smooth
sources compatible at the initial lateral corner.  The estimates below
are uniform on \(B_1\Subset B_2\), so interior compactness passes them
to the original source.  The difference between
\(\mathcal K\) and this auxiliary solution has the original lateral
data but no direct source and is estimated only on
\(B_1\Subset B_2\).  The boundary-inclusive estimate for the
auxiliary problem, followed by the interior estimate, gives
\begin{equation}\label{eq:C2-evolution-family-bound}
 \|\widetilde{\mathcal U}_{\tau_*}(s,r)F\|_{C^{2,\alpha}(B_1)}
 \leq C\|F\|_{C^{2,\alpha}(B_2)},
 \qquad s_-\leq r\leq s\leq0 .
\end{equation}
The constant depends on the first-state phase only through
\(\exp(C\int q_1)\): this follows by \(L^1\)-Gronwall for the
\(\widehat a_1\widehat B_1\) coefficient, together with the ordinary
uniform evolution estimate for the remaining coefficients.  Hence
Duhamel's formula gives the genuinely
\(L^1\)-in-time, linear estimate
\begin{equation}\label{eq:linear-direct-source-C2}
 \left\|
 \int_{s_-}^s\widetilde{\mathcal U}_{\tau_*}(s,r)
       \widetilde{\mathcal S}_{12}(r)\,dr
 \right\|_{C^{2,\alpha}(B_1)}
 \leq C\int_{s_-}^s q_{12}(\tau(r))\,dr .
\end{equation}
After subtracting this Duhamel term, the ordinary conjugated
Bernstein argument has no inhomogeneous direct column.  At the center
of the cylinder, at \(s=0\), it yields
\begin{equation}\label{eq:two-state-local-C2}
 \begin{split}
 |\widehat\nabla\mathcal K|
 +|\widehat\nabla^2\mathcal K|
 \leq C\bigg(&
 \|\mathcal K\|_{C^0(Q_2)}
 +\mathbf1_{\{s_->-1\}}
   \|\mathcal K(s_-)\|_{C^3(B_2)}\\
 &+\int_{s_-}^{0}q_{12}(\tau(r))\,dr
 +\|\widetilde{\Delta\E}\|_{C^2(Q_2)}
 \bigg).
 \end{split}
\end{equation}
The uniform jets of \(\mathcal Q_s^{\pm1}\) imply
\[
 \|\mathcal K\|_{C^0(Q_2)}
 \leq C\|\mathcal H\|_{C^0(\mathcal T_2)},\qquad
 \|\mathcal K(s_-)\|_{C^3(B_2)}
 \leq C\|\mathcal H(s_-)\|_
 {C^3(\mathcal Q_{s_-}(B_2))},
\]
and
\[
 \|\widetilde{\Delta\E}\|_{C^2(Q_2)}
 \leq C\|\widehat{\Delta\E}\|_{C^2(\mathcal T_2)}.
\]
All three right-hand norms are well-defined because
\(\mathcal T_2\subset[s_-,0]\times\mathcal A_*^+\).
Since \(\mathcal Q_0=\operatorname{Id}\),
\[
 \widehat\nabla^j\mathcal K(0,x_0)
 =\widehat\nabla^j\mathcal H(0,x_0),
 \qquad0\leq j\leq2.
\]

Return to normalized time.  On a recent cylinder,
\(dr=(1-r)\,d\tau\) and \(1\leq1-r\leq2\), so
\[
 \int_{s_-}^{0}q_{12}(\tau(r))\,dr
 \leq
 2\int_{\max\{\tau_0,\tau_*-\log2\}}^{\tau_*}
 q_{12}(\xi)\,d\xi
 \leq
 C d_0e^{-\theta_-(\tau_*-\tau_0)}.
\]
With \(\tau=\tau_*\), this is bounded by the right side of
\eqref{eq:two-state-C2-barrier-result} because
 \(\theta_->\sigma_-\).  By
\eqref{eq:two-state-C0-barrier-result} and
\eqref{eq:two-state-transported-weight-comparison},
\[
 \|\mathcal H\|_{C^0(\mathcal T_2)}
 \leq C d_0
 \omega_{\sigma_-,\tau_0}(\tau_*,y_*).
\]
The entrance term is
controlled by \eqref{eq:two-state-C3-entrance}.  If the rescaled
moving tube \(\mathcal T_2\) meets the graft annulus, the weight in
\eqref{eq:two-state-C2-barrier-result} is comparable to one and
\eqref{eq:two-state-graft-size} bounds the final term by \(Cd_0\).
 If it does not meet that annulus, the final term vanishes.  The fixed
 compact core is treated with fixed cylinders.  At the final center,
the forward and inverse phase-chart bounds give
\[
 \sum_{\ell=0}^{2}
 |\bar\nabla^\ell w|(\tau_*,y_*)
 \leq
 C\sum_{\ell=0}^{2}
 |\widehat\nabla^\ell\mathcal H|(0,x_0),
\]
with the reverse local comparison supplied by the inverse phase
chart.  Undoing the soliton and phase pullbacks proves
\eqref{eq:two-state-C2-barrier-result}.
\end{proof}

\begin{theorem}[Global two-state continuation estimate]
\label{thm:global-two-state-estimate}
Fix one normalized entrance time, one preparation convention, and a
common-margin sliced ball of strict prepared entrances
\[
 \mathscr B\subset\Sigma_{\tau_0}^{k+2,\alpha}
 \subset\mathscr P_{\tau_0}^{k+2,\alpha},
\]
chosen so that
one buffered physical cover and the strict time-width margin in
\eqref{eq:auxiliary-buffered-time-width} work throughout the ball.
Using only the fixed positive margin \(\beta-\theta\), choose this
same ball once so that every polarized stable energy form on it has
damping at least \(2\theta\).  This rate-independent shrink is
included in the package entrance ceiling \(\varepsilon_{\rm ent}\)
and is made before any subordinate pair is quantified.
 Assume that the preparation convention uses the two-state ceilings
 \(\Lambda_{\rm ell}^{(2)}\) and
\(K_{\mathcal Y,0}^{(2)}\), that the common \(C^2\) box is contained
in the ball of radius \(\delta_{\rm box}^{(2)}\), and that
\eqref{eq:two-state-uniform-package-inputs} holds.  These are one-time
 package choices made with \(\theta_*=\theta\).  In particular the
 admissible Kato bounds, \(\Gamma\), the domains \(K_i^\Gamma\),
 \(\eta_{\rm core}^{(2)},\eta_{\rm ch}^{(2)}\),
 \(\varepsilon_{\rm ph,*}^{(2)}\), the fixed-core separation
 inequalities, and \(\tau_{\rm base}^{(2)}\) are frozen before
 \(0<\sigma_-<\theta_-<\theta\) is selected; none is reselected for a
 subordinate pair.
Put the restrictions of the closed flows to \(E^{++}\) in the common
anchored exterior Ricci--DeTurck gauge relative to the fixed smooth
reference metric chosen at the center entrance, and transport the
prepared markings on the fixed graft/interface collar in that gauge.
The normalized \(h\)-block, not an exterior gauge norm, measures the
excluded collapsing core.
For two entrances \(\mathbf z_{i,0}\in\mathscr B\), let
\[
 d_0=
 \|\mathbf z_{1,0}-\mathbf z_{2,0}\|
       _{\mathscr X_{\rm prep}^{k+2,\alpha}},
 \qquad
 V=H_1-H_2,\qquad
 \E_i:=\mathcal E_{{\rm gr},i}
 =(\Phi_i^{-1})^*\mathcal G_{{\rm gr},i}.
\]
Given
\[
 0<\sigma_-<\theta_-<\theta,
\]
use the weight \(\omega_{\sigma_-,\tau_0}\) defined in
\eqref{eq:two-state-pointwise-weight}.
By Lemma~\ref{lem:uniform-two-state-package-radius}, with
\(\theta_*=\theta\), the same package radius already frozen in
\eqref{eq:global-compatible-package-radius} works for every such pair.
 The constants below may depend on the chosen
 \((\sigma_-,\theta_-)\); no uniformity at the boundary of this open
 rate cone is asserted.  More precisely, only the energy-tail constant
 \(C_E\) used below, \(C_{12}\), the two-state graft-tail constant
 \(C_{\rm gr}^{(2)}\), \(D_{\rm O}\), \(\widehat\tau_b\), the common
 amplitude multiplier, and the final comparison constant may depend on
 that subordinate pair; the Kato bounds, phase threshold, core
 separation, radius, and base time do not.
Then the global coupled solutions satisfy
\begin{align}
 \|V(\tau)\|_{L^2_\nu}
 +\left(
   \int_\tau^\infty\|V(s)\|_{H^1_\nu}^2\,ds
  \right)^{1/2}
 &\leq
 Cd_0e^{-\theta_-(\tau-\tau_0)},
 \label{eq:global-two-state-energy}\\
 \sup_{\{\bar f\leq e^\tau\}}
 \omega_{\sigma_-,\tau_0}^{-1}
 \sum_{\ell=0}^2|\bar\nabla^\ell(h_1-h_2)|
 &\leq Cd_0,
 \label{eq:global-two-state-pointwise}\\
 \int_\tau^\infty
 |(a_1,b_1)-(a_2,b_2)|(s)\,ds
 &\leq
 Cd_0e^{-\theta_-(\tau-\tau_0)}.
 \label{eq:global-two-state-phase}
\end{align}
Set
\[
 \Delta\mathcal C_\rho
 =\mathcal C_\rho[h_1]-\mathcal C_\rho[h_2],
 \qquad
 \Delta\mathcal Y_{j,\tau}
 =\mathcal Y_{j,\tau}(\mathbf z_1)
  -\mathcal Y_{j,\tau}(\mathbf z_2).
\]
The moving graft, cutoff-commutator, and effective-column differences
obey
\begin{equation}\label{eq:global-two-state-forcing}
 \begin{split}
 &\|\Delta\mathcal C_\rho\|_{H^{-1}_\nu}
  +\|\rho_\tau(\E_1-\E_2)\|_{H^{-1}_\nu}\\
 &\quad+
 \sum_{j=0}^8
 \|\rho_\tau\Delta\mathcal Y_{j,\tau}\|_{H^{-1}_\nu}\\
 &\quad+
 \left(
  \sum_{\mu=0}^8\left(
   \left|\ip{\Delta\mathcal C_\rho}{Z_\mu}\right|^2
   +
   \left|
   \ip{\rho_\tau(\E_1-\E_2)}{Z_\mu}
   \right|^2\right)
  +\sum_{\mu,j=0}^8
   \left|
    \ip{\rho_\tau\Delta\mathcal Y_{j,\tau}}{Z_\mu}
   \right|^2
 \right)^{1/2}
 \leq Cd_0e^{-ce^\tau}.
 \end{split}
\end{equation}
After \eqref{eq:global-two-state-phase} has been established, the same
forcing also satisfies the pointwise estimate
\begin{equation}\label{eq:global-two-state-forcing-C2}
 \supp(\E_1-\E_2)
 \subset
 \{c\Gamma e^\tau\leq\bar f\leq C\Gamma e^\tau\},
 \qquad
 \sum_{\ell=0}^2
 |\bar\nabla^\ell(\E_1-\E_2)|_{\bar g}
 \leq Cd_0e^{-\tau}.
\end{equation}

The terminal geometric data are locally Lipschitz:
\begin{equation}\label{eq:terminal-geometric-data-Lipschitz}
 \begin{split}
 &|T_1-T_2|
 +\left|\log\frac{\lambda_{\infty,1}}
                         {\lambda_{\infty,2}}\right|
 \\
 &\quad+
 \left\|
  \log\bigl(\Psi_{\infty,2}^{-1}
             \circ\Psi_{\infty,1}\bigr)
 \right\|_{C^m(K)}
 \leq C_{K,m}d_0
 \end{split}
\end{equation}
for every \(K\Subset M\), in an exponential chart at the identity.

All six limiting-data maps are \(C^1\) on \(\mathscr B\):
\(\mathbf z_0\mapsto T\) and
\(\mathbf z_0\mapsto\log\lambda_\infty\) are scalar-valued \(C^1\)
maps; for every \(K\Subset M\) and \(m\geq0\),
\(\mathbf z_0\mapsto\Psi_\infty\) is \(C^1\) into the fixed
\(C^m(K)\) exponential chart; and, for every integer
\(0\leq m\leq k+2\), the maps
\[
 \mathbf z_0\longmapsto
 R_\infty,\qquad R_\infty^{-1},\qquad S_\infty
\]
are \(C^1\) into the right-translated prepared
\(C^m(\Omega_\eta^+)\) charts.  Their derivatives obey the following
uniform estimates.

For every base entrance \(\mathbf z_0\in\mathscr B\) and tangent
vector
\[
 \xi\in
 T_{\mathbf z_0}\Sigma_{\tau_0}^{k+2,\alpha},
\]
write
\[
 K_\xi=DH_{\mathbf z_0}[\xi],\quad
 k_\xi=Dh_{\mathbf z_0}[\xi],\quad
 \dot c_\xi=Dc_{\mathbf z_0}[\xi].
\]
Here \(\|\xi\|\) is the tangent norm induced by
\(\mathscr E_{\rm prep}^{k+2,\alpha}\).
The first variations satisfy, uniformly on the common-margin ball,
\begin{align}
 \|K_\xi(\tau)\|_{L^2_\nu}
 +\left(\int_\tau^\infty
        \|K_\xi(s)\|_{H^1_\nu}^2\,ds\right)^{1/2}
 &\leq
 C e^{-\theta_-(\tau-\tau_0)}\|\xi\|,
 \label{eq:global-variation-energy}\\
 \int_\tau^\infty|\dot c_\xi(s)|\,ds
 &\leq
 C e^{-\theta_-(\tau-\tau_0)}\|\xi\|,
 \label{eq:global-variation-phase}\\
 \sup_{\{\bar f\leq e^\tau\}}
 \omega_{\sigma_-,\tau_0}^{-1}
 \sum_{\ell=0}^2|\bar\nabla^\ell k_\xi|
 &\leq C\|\xi\|.
 \label{eq:global-variation-pointwise}
\end{align}
If
\[
 T_S=t(\tau_0)+\int_{\tau_0}^{S}\lambda(s)\,ds,\qquad
 \ell_S=\log(\lambda(S)e^S),
\]
then, for every \(K\Subset M\) and \(m\geq0\),
\begin{equation}\label{eq:C1-geometric-tail-uniform}
 \begin{split}
 &|D(T-T_S)[\xi]|
 +|D(\log\lambda_\infty-\ell_S)[\xi]|\\
 &\qquad+
 \left\|
 D\!\left[
  \log(\Psi_S^{-1}\circ\Psi_\infty)
 \right][\xi]
 \right\|_{C^m(K)}
 \leq
 \varepsilon_{K,m}(S)\|\xi\|,
 \end{split}
\end{equation}
where \(\varepsilon_{K,m}(S)\to0\) as \(S\to\infty\), uniformly for
the base entrance in the common-margin ball.
Writing \(R_\infty\) and \(S_\infty\) for the adaptive limits in
Theorem~\ref{thm:prepared-entrance-continuation}, one also has
for every integer \(0\leq m\leq k+2\),
\begin{equation}\label{eq:C1-adaptive-tail-uniform}
 \begin{split}
 &\left\|
 D\!\left[\operatorname{Exp}_{R_S}^{-1}R_\infty\right][\xi]
 \right\|_{C^m(\Omega_\eta^+)}
 +\left\|
 D\!\left[
 \operatorname{Exp}_{R_S^{-1}}^{-1}R_\infty^{-1}
 \right][\xi]
 \right\|_{C^m(\Omega_\eta^+)}\\
 &\qquad+
 \|D(S_\infty-S(S))[\xi]\|_{C^m(\Omega_\eta^+)}
 \leq\varepsilon_m^{\rm ad}(S)\|\xi\|,
 \end{split}
\end{equation}
where \(\varepsilon_m^{\rm ad}(S)\to0\) uniformly.  The exponential
charts in this display are the right-translated prepared map charts;
any equivalent fixed bounded-geometry chart gives the same
conclusion.

On any ambient prepared neighborhood
\(\mathscr O\subset\mathscr P_{\tau_0}^{k+4,\alpha}\) satisfying
\[
 \mathscr O\subset\operatorname{dom}
    \Pi_{\rm sl}^{\,k+4\to k+2},
 \qquad
 \Pi_{\rm sl}^{\,k+4\to k+2}(\mathscr O)\subset\mathscr B,
\]
these statements apply to the phase-retracted evolution
\[
 \mathcal E_{\rm amb}
 =\mathcal E_{\rm sl}\circ
  \Pi_{\rm sl}^{\,k+4\to k+2}.
\]
For \(\mathbf z_i\in\mathscr O\), put
\(\widetilde{\mathbf z}_i
 :=\Pi_{\rm sl}^{\,k+4\to k+2}(\mathbf z_i)\).  Then
\begin{align*}
 \|\widetilde{\mathbf z}_1-\widetilde{\mathbf z}_2\|_{
    \mathscr X_{\rm prep}^{k+2,\alpha}}
 &\leq K_{\Pi,k+2}
 \|\mathbf z_1-\mathbf z_2\|_{
    \mathscr X_{\rm prep}^{k+4,\alpha}},\\
 \sup_{\mathbf z\in\mathscr O}
 \left\|
  D\Pi_{\rm sl}^{\,k+4\to k+2}(\mathbf z)
 \right\|_{
  \mathcal L(\mathscr E_{\rm prep}^{k+4,\alpha},
             \mathscr E_{\rm prep}^{k+2,\alpha})}
 &\leq K_{\Pi,k+2}.
\end{align*}
Apply every sliced two-state or first-variation estimate to the
\(\widetilde{\mathbf z}_i\).  Its ambient constant is therefore at most
the corresponding sliced constant times \(K_{\Pi,k+2}\).  The
additional conditions
\[
 \overline{\mathscr O}\subset\operatorname{dom}
 \Pi_{\rm sl}^{\,k+4\to k+2},
 \qquad
 \Pi_{\rm sl}^{\,k+4\to k+2}(\overline{\mathscr O})
 \subset\mathscr B'\Subset_{\rm u}\mathscr B
\]
retain the fixed target charts and common prepared margins; they are
not invoked as compactness hypotheses.

\end{theorem}

\begin{proof}
Proposition~\ref{prop:two-state-prepared-evolution} puts the two
solutions in one prepared chart on every finite interval and controls
the solution-dependent pullbacks.
Lemma~\ref{lem:coarse-two-state-geometry} then absorbs the possible
exponential growth of those finite-horizon comparison constants into
the Gaussian tail on the receding support.  Applied to the moving
 graft, cutoff commutator, and effective columns, it gives
\eqref{eq:global-two-state-forcing} without using any global
two-state decay estimate.

Subtract the exact cutoff equations before estimating any term.
With
\[
 w=h_1-h_2,\qquad
 \delta c=c_1-c_2,\qquad
 \mathcal Y_{j,\tau}^{(i)}
 =\mathcal Y_{j,\tau}(\mathbf z_i),\qquad
 V=H_1-H_2=\rho_\tau w,
\]
one obtains
\begin{equation}\label{eq:exact-two-state-cutoff-equation}
 \begin{split}
 \partial_\tau V={}&\A V
 +\rho_\tau\bigl(\Q(h_1)-\Q(h_2)\bigr)
 +\Delta\mathcal C_\rho+\rho_\tau(\E_1-\E_2)\\
 &+\sum_{j=0}^8c_{1,j}\rho_\tau\mathscr T_jw
 +\sum_{j=0}^8\delta c_j
      \rho_\tau\bigl(\mathcal Y_{j,\tau}^{(1)}
                      +\mathscr T_jh_2\bigr)\\
 &+\sum_{j=0}^8c_{2,j}\rho_\tau
      \Delta\mathcal Y_{j,\tau}.
 \end{split}
\end{equation}
Both sliced tensors are orthogonal to \(\mathcal Z\), so
\(V\perp\mathcal Z\).  In the term containing \(\delta c_j\), add and
subtract the fixed global column \(Y_j\).  Its \(Y_j\)-part has zero
pairing with \(V\).  The defect
\(\rho_\tau\mathcal Y_{j,\tau}^{(1)}-Y_j\) has receding support, while
\(\rho_\tau\mathscr T_jh_2\) is a core term and is estimated
separately:
\begin{align}
 \left|
 \left\langle\rho_\tau\mathscr T_jh_2,V\right\rangle
 \right|
 &\leq
 C\|H_2\|_{H^1_\nu}\|V\|_{H^1_\nu}
 +Cd_0e^{-ce^\tau},
 \label{eq:two-state-core-action-mixed}\\
 \left|
 \left\langle\rho_\tau\mathscr T_jw,V\right\rangle
 \right|
 &\leq
 C\|V\|_{H^1_\nu}^2
 +Cd_0^2e^{-ce^\tau}.
 \label{eq:two-state-core-action-difference}
\end{align}
These are the polarized forms of
\eqref{eq:B-bilinear}; the Gaussian terms account for replacing
\(\rho_\tau w\) by \(w\) on the cutoff annulus.

Polarization of \eqref{eq:localized-Q-Lip-energy}, the common
\(C^2\) box, and the coarse two-state bound on the receding annulus
give
\begin{equation}\label{eq:two-state-polarized-energy}
 \begin{split}
 \left|
 \left\langle
 \rho_\tau(\Q(h_1)-\Q(h_2)),V
 \right\rangle
 \right|
 \leq{}&
 \varepsilon_*
 \|V\|_{H^1_\nu}^2
 +Cd_0^2e^{-ce^\tau}.
 \end{split}
\end{equation}
Here one integration by parts treats the quasilinear second
derivatives, and the difference \(w-V\) is supported where
\(\bar f\geq e^\tau\).  Similarly,
\begin{equation}\label{eq:two-state-complete-outer-energy}
 \begin{split}
 &\|\Delta\mathcal C_\rho\|_{H^{-1}_\nu}
 +\|\rho_\tau(\E_1-\E_2)\|_{H^{-1}_\nu}\\
 &\quad+
 \sum_{j=0}^8
 \|\rho_\tau\Delta\mathcal Y_{j,\tau}\|_{H^{-1}_\nu}
 +\sum_{j=0}^8
 \|\rho_\tau\mathcal Y_{j,\tau}^{(1)}-Y_j\|_{H^{-1}_\nu}
 \leq Cd_0e^{-ce^\tau}+Ce^{-ce^\tau}.
 \end{split}
\end{equation}
In the last sum the term independent of \(d_0\) is always multiplied
by \(\delta c\); uniform Gram inversion absorbs that small
coefficient.  Both displays follow by splitting at the receding
support, applying the scale-one chart calculus there, and absorbing
the coarse \(e^{A\tau}\) difference growth into the Gaussian density.

Next subtract the two exact Gram systems:
\[
 M_1\delta c=-\delta d-(M_1-M_2)c_2.
\]
Use \eqref{eq:two-state-polarized-energy} in modal form,
\eqref{eq:two-state-complete-outer-energy}, and uniform inversion of
\(M_1\).  The modal analogues of
\eqref{eq:two-state-core-action-mixed}--%
\eqref{eq:two-state-core-action-difference} control the core action
terms; in particular they are not assigned to the outer package.
This gives the endpoint-independent sharp modal estimate
\begin{equation}\label{eq:two-state-sharp-modal}
 |\delta c|
 \leq
 C\bigl(\|H_1\|_{H^1_\nu}+\|H_2\|_{H^1_\nu}\bigr)
   \|V\|_{H^1_\nu}
 +Cd_0e^{-ce^\tau}.
\end{equation}
The term \((M_1-M_2)c_2\) has this size because the core difference
of \(M\) is linear in \(V\), while its outer difference is Gaussian;
the one-state estimate
\(|c_2|\leq C\|H_2\|_{H^1_\nu}^2+Ce^{-ce^\tau}\)
supplies the remaining factor.  The Gaussian summand in this bound
also produces the intermediate term
\(C\|V\|_{H^1_\nu}e^{-ce^\tau}\).  By
\eqref{eq:coarse-two-state-growth}, the tensor block in
\eqref{eq:two-state-hybrid-distance}, and Gaussian absorption of its
fixed polynomial radial loss,
\[
 \|V\|_{H^1_\nu}e^{-ce^\tau}
 \leq Cd_0e^{-c'e^\tau}.
\]
Decreasing \(c'\) and renaming it \(c\) includes this contribution in
the last term of \eqref{eq:two-state-sharp-modal}; no additional term
is omitted from that display.

Pair \eqref{eq:exact-two-state-cutoff-equation} with \(V\).
Use stable coercivity, \eqref{eq:B-bilinear},
\eqref{eq:two-state-polarized-energy}, and
\eqref{eq:two-state-sharp-modal}, together with the two core-action
bounds
\eqref{eq:two-state-core-action-mixed}--%
\eqref{eq:two-state-core-action-difference}.  The one-time common-ball
choice in the theorem statement and the fixed positive spectral margin
\(\beta-\theta\) give damping at least \(2\theta\) for the perturbed
energy form.  Now fix an
arbitrary \(\theta_-<\theta\) and weaken that coefficient to
\(2\theta_-\).  Young's inequality absorbs every occurrence of
\(\|V\|_{H^1_\nu}\) into the stable form and yields
\begin{equation}\label{eq:two-state-global-energy-differential}
 \frac d{d\tau}\|V\|_{L^2_\nu}^2
 +2\theta_-\|V\|_{L^2_\nu}^2
 +c\|V\|_{H^1_\nu}^2
 \leq Cd_0^2e^{-ce^\tau}.
\end{equation}
No finite-horizon state distance occurs on the right.  The prepared
entrance norm gives
\(\|V(\tau_0)\|_{L^2_\nu}\leq Cd_0\).  Applying the integrating factor
to \eqref{eq:two-state-global-energy-differential} from \(\tau_0\) to
\(\tau\) therefore gives
\[
 \|V(\tau)\|_{L^2_\nu}^2
 \leq
 Cd_0^2e^{-2\theta_-(\tau-\tau_0)}
 +Cd_0^2\int_{\tau_0}^{\tau}
 e^{-2\theta_-(\tau-s)}e^{-ce^s}\,ds
 \leq
 Cd_0^2e^{-2\theta_-(\tau-\tau_0)}.
\]
Next integrate
\eqref{eq:two-state-global-energy-differential} from \(\tau\) to a
finite \(S\), discard the nonnegative terminal energy, and use the
just-proved pointwise bound together with the Gaussian forcing tail.
This yields
\[
 \int_\tau^S\|V(s)\|_{H^1_\nu}^2\,ds
 \leq Cd_0^2e^{-2\theta_-(\tau-\tau_0)}
\]
with a constant independent of \(S\).  Letting \(S\to\infty\) and
taking square roots proves \eqref{eq:global-two-state-energy}.

Finally integrate \eqref{eq:two-state-sharp-modal}.  By
Cauchy--Schwarz, the one-state dissipation tail, and
\eqref{eq:global-two-state-energy},
\[
 \begin{split}
 \int_\tau^\infty|\delta c(s)|\,ds
 &\leq
 C\sum_{i=1}^2
 \left(\int_\tau^\infty\|H_i(s)\|_{H^1_\nu}^2\,ds\right)^{1/2}
 \left(\int_\tau^\infty\|V(s)\|_{H^1_\nu}^2\,ds\right)^{1/2}\\
 &\qquad+Cd_0e^{-ce^\tau}\\
 &\leq Cd_0e^{-\theta_-(\tau-\tau_0)}.
 \end{split}
\]
This proves \eqref{eq:global-two-state-phase}.
Apply Lemma~\ref{lem:physical-graft-two-state-Lipschitz} with
\(\vartheta=\theta_-\).  This proves the forcing estimate
\eqref{eq:global-two-state-forcing-C2}, and then
\eqref{eq:two-column-C2}.  This invocation follows the energy and
phase estimates, whose proof uses only the Gaussian
\(H^{-1}_\nu\) estimate.

We next obtain the pointwise estimate without suppressing the
solution-dependent column difference.  Define
\(\mathfrak y_{12}\) by \eqref{eq:two-state-column-size}.  By
Lemma~\ref{lem:physical-graft-two-state-Lipschitz},
\[
 \mathfrak y_{12}(\tau)\leq C d_0.
\]
The one-state phase estimate
\eqref{eq:master-bootstrap-phase-tail} and
\eqref{eq:global-two-state-phase} therefore give, for every finite
\(S>\tau\),
\begin{equation}\label{eq:global-two-state-augmented-tail}
 \begin{split}
 \int_\tau^S
 \bigl(|c_1-c_2|+|c_2|\mathfrak y_{12}\bigr)(s)\,ds
 &\leq
 Cd_0e^{-\theta_-(\tau-\tau_0)}
 +Cd_0\int_\tau^\infty |c_2(s)|\,ds\\
 &\leq Cd_0e^{-\theta_-(\tau-\tau_0)}.
 \end{split}
\end{equation}
Here the fixed entrance time is absorbed in the uniform common-ball
constant.  Since every entrance in \(\mathscr B\) has one-state
entrance amplitude at most \(\varepsilon_{\rm ent}\), the named
future-tail constants in \eqref{eq:named-future-tail-constants} and
\eqref{eq:master-bootstrap-phase-tail}, evaluated for the first
solution at \(\tau_0\), give
\begin{equation}\label{eq:global-two-state-background-phase-verification}
 \int_{\tau_0}^{S}q_1(s)\,ds
 \leq P_{\infty,1}(\tau_0)
 \leq C_P^*\varepsilon_{\rm ent}^{\,2}e^{-2\theta_*\tau_0}
      +C_P^*e^{-c_P^*e^{\tau_0}}
 \leq\varepsilon_{\rm ph}
 \leq\varepsilon_{\rm ph,*}^{(2)}.
\end{equation}
Here \(P_{\infty,1}(\tau)=\int_\tau^\infty q_1(s)\,ds\).
The penultimate inequality is exactly
\eqref{eq:uniform-two-state-phase-tail-base-time}, used with
\(\theta_*=\theta\), together with
\eqref{eq:two-state-uniform-package-inputs}; it is independent of the
subordinate pair.  Thus
\eqref{eq:two-state-background-phase-budget} is supplied rather than
assumed anew after \((\sigma_-,\theta_-)\) is chosen.

Fix \(C_{12}<\infty\) to dominate, uniformly on the common ball, the
constant in \eqref{eq:global-two-state-augmented-tail}.  This is a
post-radius constant and the augmented-tail estimate is precisely
\eqref{eq:two-state-direct-tail}; the
restriction of \eqref{eq:global-two-state-energy} to \([\tau,S]\) is
\eqref{eq:two-state-barrier-energy-tail}.  Thus the barrier lemma is
invoked only after its energy input has been proved.  The two-sided
support and pointwise hypotheses for the pure graft difference are
\eqref{eq:two-graft-support}--\eqref{eq:two-graft-C2}, and the
prepared entrance distance implies
\eqref{eq:two-state-C3-entrance}.  At this post-radius stage set
\(C_{\rm gr}^{(2)}\) equal to the finite constant in
\eqref{eq:two-graft-C2}; no radius threshold depends on it.  Apply
Lemma~\ref{lem:exact-two-state-barriers} and then
Lemma~\ref{lem:two-state-linear-derivative-recovery}.  Letting the
arbitrary finite endpoint \(S\) tend to infinity gives
\eqref{eq:global-two-state-pointwise}.

The equations for \(t,\log(\lambda e^\tau)\), \(R_\tau\), and
\(\Psi_\tau\), together with
\eqref{eq:global-two-state-phase}, have integrable difference tails.
Their finite-horizon values are Lipschitz by
Proposition~\ref{prop:two-state-prepared-evolution}; passing to the
limit gives \eqref{eq:terminal-geometric-data-Lipschitz}.

We now prove the asserted first-variation bounds and convergence of
derivatives.  The argument below is a separate linearized
energy--Gram--barrier estimate; in particular, no quotient of the
graph-augmented two-state distance is used.  Fix a base trajectory in
the common-margin ball, put
\[
 d_\xi:=\|\xi\|,\qquad c:=(a,b),\qquad q:=|a|+|b|,
\]
and let
\(\dot{\mathbf z}_\xi(\tau)=D\mathbf z(\tau)[\xi]\).
Bracketed hybrid seminorms below mean the corresponding componentwise
tangent seminorms.  Thus
\(\mathfrak D_{m_\#}^{\rm hyb}[\dot{\mathbf z}_\xi]\) is the
linearization of \eqref{eq:two-state-hybrid-distance}, including the
linearized gauge and transported-marking block from
\eqref{eq:hybrid-marking-block}.  The sliced first-variation clause of
Lemma~\ref{lem:uniform-weighted-Schauder-restart}, followed by unit
interval restriction, gives
\begin{equation}\label{eq:coarse-global-hybrid-variation}
 \mathfrak D_{m_\#}^{\rm hyb}[\dot{\mathbf z}_\xi](\tau)
 +\mathfrak F_{m_\#+1}^{\rm glob}[\dot{\mathbf z}_\xi](\tau)
 \leq Ce^{A(\tau-\tau_0)}d_\xi,
 \qquad \tau\geq\tau_0.
\end{equation}
This coefficient-level estimate is terminal-time independent.  It
makes the linearized metric, map, inverse-map, graft, and marking
blocks well-defined at the orders used below.  In particular the
marking variation is the derivative of
\(\widetilde\iota=\iota\circ\chi^{-1}\); the anchored exterior
interface and the triangular gauge and inverse-gauge ODEs put it one
order below the physical metric block.

We first record the exact tensors which enter the linearized equation.
Set
\begin{align}
 \mathcal C_{\rho,\xi}^{\rm lin}
 &:=D\!\left(\mathcal C_\rho[h]\right)[\xi],
 &
 \E_\xi^{\rm lin}
 &:=D\E[\xi],
 \notag\\
 \mathcal Y_{j,\xi}^{\rm lin}
 &:=D\mathcal Y_{j,\tau}[\xi],
 &
 \mathfrak y_\xi(\tau)
 &:=
 \max_{0\leq j\leq8}
 \|\mathcal Y_{j,\xi}^{\rm lin}\|_
   {\mathfrak C_{{\rm sc},0}^{2,\alpha}},
 \label{eq:linearized-column-size}\\
 s_\xi(\tau)
 &:=
 |\dot c_\xi(\tau)|+q(\tau)\mathfrak y_\xi(\tau),
 &
 \mathcal S_\xi
 &:=
 \sum_{j=0}^8\dot c_{\xi,j}
       \bigl(\mathcal Y_{j,\tau}+\mathscr T_jh\bigr)
 +\sum_{j=0}^8c_j\mathcal Y_{j,\xi}^{\rm lin}.
 \label{eq:linearized-direct-source}
\end{align}
Differentiate the exact adaptive normalized equation
\eqref{eq:adaptive-normalized} directly.  Without taking a difference
quotient, this gives
\begin{equation}\label{eq:exact-variational-h-equation}
 \partial_\tau k_\xi
 =
 \A k_\xi+D\Q_h[k_\xi]
 +\sum_{j=0}^8c_j\mathscr T_jk_\xi
 +\mathcal S_\xi+\E_\xi^{\rm lin}.
\end{equation}
Since \(K_\xi=\rho_\tau k_\xi\), its exact cutoff form is
\begin{equation}\label{eq:exact-variational-cutoff-equation}
 \begin{split}
 \partial_\tau K_\xi={}&
 \A K_\xi+\rho_\tau D\Q_h[k_\xi]
 +\mathcal C_{\rho,\xi}^{\rm lin}
 +\rho_\tau\E_\xi^{\rm lin}\\
 &+\sum_{j=0}^8c_j\rho_\tau\mathscr T_jk_\xi
 +\rho_\tau\mathcal S_\xi .
 \end{split}
\end{equation}
Because
\(\xi\in T_{\mathbf z_0}\Sigma_{\tau_0}^{k+2,\alpha}
=\ker D\mathfrak m(\mathbf z_0)\), differentiation of the propagated
slice gives
\[
 K_\xi(\tau)\perp\mathcal Z
 \quad\text{for every }\tau\geq\tau_0.
\]

We next prove the complete outer input needed by both the energy and
the linearized Gram system.  The sliced first-variation conclusion in
Lemma~\ref{lem:coarse-two-state-geometry} follows from the exact
linearized restart, not from division by \(d_0\).  Applying its
scale-one pullback and composition calculation to the three
linearized outer blocks gives, for one fixed Gaussian constant
\(c_{\rm G}>0\),
\begin{equation}\label{eq:linearized-complete-outer-package}
 \begin{split}
 &\|\mathcal C_{\rho,\xi}^{\rm lin}\|_{H^{-1}_\nu}
 +\|\rho_\tau\E_\xi^{\rm lin}\|_{H^{-1}_\nu}
 +\sum_{j=0}^8
   \|\rho_\tau\mathcal Y_{j,\xi}^{\rm lin}\|_{H^{-1}_\nu}\\
 &\quad+
 \sum_{\mu=0}^8\left(
  \left|\ip{\mathcal C_{\rho,\xi}^{\rm lin}}{Z_\mu}\right|
  +\left|\ip{\rho_\tau\E_\xi^{\rm lin}}{Z_\mu}\right|
  +\sum_{j=0}^8
   \left|\ip{\rho_\tau\mathcal Y_{j,\xi}^{\rm lin}}{Z_\mu}\right|
 \right)
 \leq Cd_\xi e^{-c_{\rm G}e^\tau},\\
 &\sum_{j=0}^8\left(
   \|\rho_\tau\mathcal Y_{j,\tau}-Y_j\|_{H^{-1}_\nu}
   +\sum_{\mu=0}^8
    \left|\ip{\rho_\tau\mathcal Y_{j,\tau}-Y_j}{Z_\mu}\right|
 \right)
 \leq Ce^{-c_{\rm G}e^\tau}.
 \end{split}
\end{equation}
Indeed, every tensor in the first two lines is supported where
\(\bar f\geq c_{\rm G}e^\tau\), is bounded by a fixed polynomial in \(e^\tau\)
times the left side of
\eqref{eq:coarse-global-hybrid-variation}, and is paired against the
Gaussian measure.  The last line is the same calculation for the
one-state column defect.  This proves
\eqref{eq:linearized-complete-outer-package} with no future high
prepared norm.

The remaining energy terms are core terms.  Polarizing
\eqref{eq:B-bilinear} in the tangent slot and using
\eqref{eq:coarse-global-hybrid-variation} on the cutoff annulus gives
\begin{align}
 \left|
  \left\langle\rho_\tau\mathscr T_jh,K_\xi\right\rangle
 \right|
 &\leq
 C\|H\|_{H^1_\nu}\|K_\xi\|_{H^1_\nu}
 +Cd_\xi e^{-c_{\rm G}e^\tau},
 \label{eq:linearized-core-action-mixed}\\
 \left|
  \left\langle\rho_\tau\mathscr T_jk_\xi,K_\xi\right\rangle
 \right|
 &\leq
 C\|K_\xi\|_{H^1_\nu}^2
 +Cd_\xi^2e^{-c_{\rm G}e^\tau}.
 \label{eq:linearized-core-action-self}
\end{align}
Likewise, differentiating the localized quasilinear form estimate
\eqref{eq:localized-Q-Lip-energy}, integrating the scalar-principal
second-order term once by parts, and using the common \(C^2\) box
gives
\begin{equation}\label{eq:linearized-polarized-energy}
 \left|
  \left\langle\rho_\tau D\Q_h[k_\xi],K_\xi\right\rangle
 \right|
 \leq
 \varepsilon_*\|K_\xi\|_{H^1_\nu}^2
 +Cd_\xi^2e^{-c_{\rm G}e^\tau}.
\end{equation}

Now differentiate the exact Gram law
\eqref{eq:receding-Gram}.  If \(Mc=-\mathbf d\) denotes that law along
the base trajectory, then
\[
 M\dot c_\xi
 =-D\mathbf d[\xi]-(DM[\xi])c.
\]
The fixed \(Y_j\)-part of the direct column has zero pairing with
\(K_\xi\) by the differentiated slice.  Uniform inversion of \(M\),
the modal versions of
\eqref{eq:linearized-core-action-mixed}--%
\eqref{eq:linearized-polarized-energy}, and
\eqref{eq:linearized-complete-outer-package} therefore give
\begin{equation}\label{eq:linearized-sharp-modal}
 |\dot c_\xi|
 \leq
 C\|H\|_{H^1_\nu}\|K_\xi\|_{H^1_\nu}
 +Cd_\xi e^{-c_{\rm G}e^\tau}.
\end{equation}
Here the core part of \(DM[\xi]\) is linear in \(K_\xi\), whereas its
outer part is Gaussian.  The factor \(c\) is controlled by the
one-state quadratic feedback estimate; thus its core contribution is
absorbed into the first term of
\eqref{eq:linearized-sharp-modal}.  The product of a Gaussian column
defect and \(K_\xi\) is also Gaussian by
\eqref{eq:coarse-global-hybrid-variation}.  Hence no unlisted
finite-horizon state distance occurs in this modal estimate.

Pair \eqref{eq:exact-variational-cutoff-equation} with \(K_\xi\).
Use stable coercivity,
\eqref{eq:linearized-core-action-mixed}--%
\eqref{eq:linearized-sharp-modal}, and the one-time common-ball
smallness.  The perturbed stable form has damping at least
\(2\theta\); after fixing \(\theta_-<\theta\), Young's inequality
gives
\begin{equation}\label{eq:linearized-global-energy-differential}
 \frac d{d\tau}\|K_\xi\|_{L^2_\nu}^2
 +2\theta_-\|K_\xi\|_{L^2_\nu}^2
 +c_{\rm en}\|K_\xi\|_{H^1_\nu}^2
 \leq Cd_\xi^2e^{-c_{\rm G}e^\tau},
 \qquad c_{\rm en}>0.
\end{equation}
The tangent entrance norm gives
\(\|K_\xi(\tau_0)\|_{L^2_\nu}\leq Cd_\xi\).
First apply the integrating factor from \(\tau_0\) to \(\tau\).  The
superexponential forcing then gives
\[
 \|K_\xi(\tau)\|_{L^2_\nu}^2
 \leq C d_\xi^2e^{-2\theta_-(\tau-\tau_0)}.
\]
Next integrate \eqref{eq:linearized-global-energy-differential}
without a weight from \(\tau\) to a finite \(S\), use this pointwise
bound at the initial face and the superexponential forcing tail, and
then let \(S\to\infty\).  This gives
\begin{equation}\label{eq:linearized-global-energy-tail}
 \|K_\xi(\tau)\|_{L^2_\nu}
 +\left(\int_\tau^\infty
        \|K_\xi(s)\|_{H^1_\nu}^2\,ds\right)^{1/2}
 \leq
 Cd_\xi e^{-\theta_-(\tau-\tau_0)}.
\end{equation}
This is \eqref{eq:global-variation-energy}.  Integrating
\eqref{eq:linearized-sharp-modal}, and using Cauchy--Schwarz together
with the one-state dissipation tail \eqref{eq:tail-diss} and
\eqref{eq:linearized-global-energy-tail}, gives
\begin{equation}\label{eq:linearized-global-phase-tail}
 \int_\tau^\infty|\dot c_\xi(s)|\,ds
 \leq
 Cd_\xi e^{-\theta_-(\tau-\tau_0)}.
\end{equation}
This is \eqref{eq:global-variation-phase}.

We next construct the pointwise inputs at the linearized level.  First
differentiate the scale and clock ODEs and use
\eqref{eq:linearized-global-phase-tail}; this bounds
\(D\log\lambda[\xi]\) and \(Dt[\xi]\) by \(Cd_\xi\) uniformly in time.
Use the sliced first-variation clauses of
Lemmas~\ref{lem:inner-terminated-exterior-DeTurck} and
\ref{lem:compact-graft-buffer-propagation}, then integrate the
linearized triangular \(R\), gauge, inverse-gauge, and marking ODEs.
Finally use the first-variation form of
Lemma~\ref{lem:localized-graft-F-coarse-memory} and the same Abel integral
which proves \eqref{eq:outer-source-F-propagation}.  The result is the
typed input estimate
\begin{equation}\label{eq:linearized-graft-input-block}
 \begin{split}
 &\mathfrak B_{{\rm gr},4}[\dot{\mathbf z}_\xi](\tau)
 +|D\log\lambda(\tau)[\xi]|
 +\|DR(\tau)[\xi]\|_{\mathfrak X_{\rm sc}^{6,\alpha}}\\
 &\qquad
 +\mathfrak F_{{\rm gr},5}[\dot{\mathbf z}_\xi](\tau)
 +\mathfrak F_{{\rm out},5}^{\pm}
       [\dot{\mathbf z}_\xi](\tau)
 \leq Cd_\xi .
 \end{split}
\end{equation}
All homogeneous memories in this calculation remain anchored at
\(\tau_0\); the nonseparated \(F\)-source is propagated by its
one-derivative Abel block, and every separated corridor term retains
its physical \(\lambda\)-factor.  Thus
\eqref{eq:linearized-graft-input-block} is a direct estimate for the
exact variational systems appearing in those formulas, not a
linearization inferred from a nonlinear Lipschitz bound.

Differentiate the exact normalized pure-graft identity
\eqref{eq:pure-graft-normalized}.  The tame calculation in the proof
of Lemma~\ref{lem:same-order-pure-graft-difference} applies directly:
each resulting summand either contains a derivative of the fixed
cutoff or a coefficient vanishing at \(\eta=0,1\).  The physical
cutoff scale and \eqref{eq:linearized-graft-input-block} therefore
give the first two assertions below.  Differentiating the exact
column formulas
\eqref{eq:effective-column-zero}--%
\eqref{eq:effective-column-j}, and using the outer-source term in
\eqref{eq:linearized-graft-input-block}, gives the third:
\begin{equation}\label{eq:linearized-graft-column-package}
 \begin{gathered}
 \supp\E_\xi^{\rm lin}
 \subset
 \{c_{\rm supp}^{(2)}\Gamma e^\tau\leq\bar f
   \leq C_{\rm supp}^{(2)}\Gamma e^\tau\},\\
 \sum_{\ell=0}^2
 |\bar\nabla^\ell\E_\xi^{\rm lin}|_{\bar g}
 \leq C_{\rm gr}^{\rm lin}d_\xi e^{-\tau},
 \qquad
 \mathfrak y_\xi(\tau)\leq C_{\rm col}^{\rm lin}d_\xi .
 \end{gathered}
\end{equation}
The constants are finite and uniform on the common-margin ball.

For a finite \(S>\tau\), define the genuinely linearized direct tail
\[
 P_{\xi,S}(\tau):=\int_\tau^S s_\xi(r)\,dr.
\]
Equations \eqref{eq:linearized-global-phase-tail},
\eqref{eq:linearized-graft-column-package}, and the one-state tail
\eqref{eq:master-bootstrap-phase-tail} give
\begin{equation}\label{eq:linearized-direct-tail}
 \begin{split}
 P_{\xi,S}(\tau)
 &\leq
 \int_\tau^\infty|\dot c_\xi(r)|\,dr
 +C_{\rm col}^{\rm lin}d_\xi
   \int_\tau^\infty q(r)\,dr\\
 &\leq
 C_{\rm dir}^{\rm lin}d_\xi
 e^{-\theta_-(\tau-\tau_0)}.
 \end{split}
\end{equation}
The already proved base-state estimate
\eqref{eq:global-two-state-background-phase-verification} supplies,
for this trajectory,
\begin{equation}\label{eq:linearized-background-phase-budget}
 \int_{\tau_0}^{S}q(r)\,dr
 \leq\varepsilon_{\rm ph}
 \leq\varepsilon_{\rm ph,*}^{(2)}.
\end{equation}
Finally, the definition of the tangent entrance norm gives
\begin{equation}\label{eq:linearized-C3-entrance}
 \sum_{\ell=0}^3
 |\bar\nabla^\ell k_\xi(\tau_0)|
 \leq C_{\rm ent}^{\rm lin}d_\xi
\end{equation}
in the scale-normalized sense of
\eqref{eq:scaled-tensor-holder}.  Thus the energy tail, direct tail,
background phase budget, graft annulus and size, and entrance trace
have all been proved with their correct linearized types.

It remains to prove the scalar inequality and perform the comparison.
Write \(g=\bar g+h\).  Direct differentiation of the exact
quasilinear expression \eqref{eq:Q-exact} gives
\begin{equation}\label{eq:linearized-principal-polarization}
 \A k_\xi+D\Q_h[k_\xi]
 =
 g^{ab}\bar\nabla_a\bar\nabla_bk_\xi
 -\bar\nabla_{\bar\nabla\bar f}k_\xi
 +A_h*\bar\nabla k_\xi+B_h*k_\xi,
\end{equation}
where
\[
 (1+\bar f)^{1/2}|A_h|+(1+\bar f)|B_h|\leq C,
\]
and the nonbackground first-order coefficient is as small as the
 common \(C^2\) box.  The uniform future coefficient package
\eqref{eq:uniform-future-low-coefficient-package}, the direct-column
bound \eqref{eq:two-state-direct-column-bound}, the
controlled-column bounds \eqref{eq:controlled-column-bounds}, the
normalized mode-growth estimate, and
\eqref{eq:linearized-column-size} give
\[
 \|\mathcal S_\xi(\tau)\|_
   {\mathfrak C_{{\rm sc},0}^{2,\alpha}}
 \leq Cs_\xi(\tau).
\]
Taking the \(\bar g\)-squared norm in
\eqref{eq:exact-variational-h-equation}, extracting the phase
transports, and absorbing the small first-order coefficient into the
negative scalar-principal form gives, weakly,
\begin{equation}\label{eq:linearized-Kato-inequality}
 \begin{split}
 \mathscr P_h|k_\xi|_{\bar g}^2
 \leq{}&
 C_{\rm K}^{(2)}s_\xi|k_\xi|_{\bar g}
 +2|k_\xi|_{\bar g}|\E_\xi^{\rm lin}|_{\bar g}\\
 &-c_{\rm K}^{(2)}g^{ab}
  \langle\bar\nabla_ak_\xi,\bar\nabla_bk_\xi\rangle_{\bar g},
 \end{split}
\end{equation}
where
\begin{equation}\label{eq:linearized-scalar-operator}
 \mathscr P_h
 =
 \partial_\tau-g^{ab}\bar\nabla_a\bar\nabla_b
 +\bar\nabla_{V_h}
 -\left(
   \frac{C_{\rm K}^{(2)}}{1+\bar f}
   +C_{\rm K}^{(2)}q\right),
 \qquad
 V_h=(1+a)\bar\nabla\bar f
 -\sum_{j=1}^8b_j\chi_\tau W_j .
\end{equation}
This is the same scalar-principal squared-norm computation which
selects the already frozen constants, now performed directly on
\eqref{eq:exact-variational-h-equation}.  In particular, no
two-state Kato or barrier conclusion has been invoked.

We spell out the barrier comparison.  If \(d_\xi=0\), uniqueness for
the linearized system gives \(k_\xi\equiv0\), so assume \(d_\xi>0\).
Put \(\widehat\tau=\tau-\tau_0\) and
\[
 J(\tau)=
 \exp\left(K_{\rm J}^{(2)}
       \int_{\tau_0}^{\tau}q(r)\,dr\right).
\]
By \eqref{eq:linearized-background-phase-budget},
\(1\leq J\leq2\).  Use the already frozen radial constants and choose
the post-radius amplitudes in the order established by
\eqref{eq:two-state-outer-amplitude-after-radius}--%
\eqref{eq:two-state-delayed-comparison-after-radius}, with
the finite energy constant \(C_E^{\rm lin}\) in
\eqref{eq:linearized-global-energy-tail},
\(C_{\rm gr}^{\rm lin}\), \(C_{\rm dir}^{\rm lin}\), and
\(C_{\rm ent}^{\rm lin}\) as the finite post-radius inputs.  Thus, for
fixed \(0<\sigma_-<\theta_-<\theta\), set
\begin{align*}
 B_{\rm I}^0
 &=d_\xi e^{-\sigma_-\widehat\tau}
   \left(A_{\rm I}\bar f^{\sigma_-}
         -D_{\rm I}\bar f^{\sigma_--1}\right),\\
 B_{\rm O}^0
 &=d_\xi\left(A_{\rm O}-D_{\rm O}\bar f^{-1}\right),\\
 B_{\rm I}
 &=JB_{\rm I}^0-K_0^{(2)}P_{\xi,S},\qquad
 B_{\rm O}
 =JB_{\rm O}^0-K_0^{(2)}P_{\xi,S}.
\end{align*}
The derivative of the final term is
\(\partial_\tau(-K_0^{(2)}P_{\xi,S})=K_0^{(2)}s_\xi\), which absorbs
the first source on the right of
\eqref{eq:linearized-Kato-inequality}.  The radial identities
\eqref{eq:radial-drift-calculation} and
\eqref{eq:outer-margin}, together with the box reduction
\eqref{eq:two-state-Kato-radial-box-reduction}, prove the strict inner
and outer supersolution inequalities.  On the graft annulus,
\eqref{eq:linearized-graft-column-package} is absorbed by the
\(D_{\rm O}d_\xi/\bar f\) margin.  Since
\(\theta_->\sigma_-\), \eqref{eq:linearized-direct-tail} keeps the
negative correction strictly below both radial profiles.  Finally the
crossing calculation \eqref{eq:barrier-crossing}, with crossover
radius \(e^{\widehat\tau}\), gives
\[
 B_{\rm O}<B_{\rm I}\quad
  \text{on }\{\bar f=\gamma_+e^{\widehat\tau}\},
 \qquad
 B_{\rm I}<B_{\rm O}\quad
  \text{on }\{\bar f=\gamma_-e^{\widehat\tau}\}.
\]
Define on \(\{\bar f\geq\Gamma\}\)
\begin{equation}\label{eq:linearized-glued-barrier}
 \mathcal B_\xi(\tau,x)
 =
 \begin{cases}
 B_{\rm I}(\tau,x),
  &\Gamma\leq\bar f(x)\leq\gamma_-e^{\widehat\tau},\\
 \min\{B_{\rm I}(\tau,x),B_{\rm O}(\tau,x)\},
  &\gamma_-e^{\widehat\tau}\leq\bar f(x)
    \leq\gamma_+e^{\widehat\tau},\\
 B_{\rm O}(\tau,x),
  &\bar f(x)\geq\gamma_+e^{\widehat\tau}.
 \end{cases}
\end{equation}
The strict crossings make \(\mathcal B_\xi\) positive and locally
Lipschitz, and
\(\mathcal B_\xi^2\) is a viscosity supersolution of
\eqref{eq:linearized-Kato-inequality} on
\(\{\bar f\geq\Gamma\}\) after the fixed delayed time.

The inner boundary is obtained independently of the desired
pointwise conclusion.  On the already frozen
\[
 \{\bar f\leq\Gamma\}\Subset K_0^\Gamma
 \Subset K_1^\Gamma\Subset K_2^\Gamma
 \Subset\{\bar f<2\Gamma\}
\]
one has \(\rho_\tau=1\) and
\(\E_\xi^{\rm lin}=0\).  Indeed, the frozen base-time inequalities
\[
 2\Gamma\leq e^{\tau_0},\qquad
 2\Gamma\leq c_{\rm supp}^{(2)}\Gamma e^{\tau_0}
\]
put \(K_2^\Gamma\) strictly inside both the unit-cutoff region and the
inner face of the graft support for every \(\tau\geq\tau_0\).
Equation
\eqref{eq:linearized-principal-polarization} is therefore of the form
\eqref{eq:phase-conjugated-core-linear-system}, with direct source
bounded by \(Cs_\xi\) and phase budget
\eqref{eq:linearized-background-phase-budget}.  Lemma
\ref{lem:phase-conjugated-core-L2-C0}, applied on
\([\tau-1,\tau]\), gives for \(\widehat\tau\geq1\)
\[
 \sup_{\{\bar f\leq\Gamma\}}|k_\xi(\tau)|
 \leq C\left(
  \sup_{\tau-1\leq r\leq\tau}\|K_\xi(r)\|_{L^2_\nu}
  +\int_{\tau-1}^{\tau}s_\xi(r)\,dr\right)
 \leq Cd_\xi e^{-\theta_-\widehat\tau}.
\]
The sliced first-variation estimate in
Proposition~\ref{prop:two-state-prepared-evolution} supplies the
initial boundary at the fixed delayed time.  Multiplying both radial
amplitudes by one common factor supplies its domination without
changing the crossings.  Add on each compact exhaustion the positive
corrector
\[
 z_\delta
 =\delta J
   e^{L(\widehat\tau-\widehat\tau_b)}(1+\bar f),
\]
let the exhaustion radius tend to infinity, and then let
\(\delta\downarrow0\).  Comparison with
\eqref{eq:linearized-glued-barrier} now proves
\begin{equation}\label{eq:linearized-C0-barrier}
 |k_\xi(\tau,x)|
 \leq
 Cd_\xi
 \min\left\{
 e^{-\sigma_-(\tau-\tau_0)}
 (1+\bar f(x))^{\sigma_-},\,1\right\}
\end{equation}
on every finite \(M\times[\tau_0,S]\), with \(C\) independent of
\(S\).  Here the fixed initial interval preceding the delayed
comparison time is covered by
\eqref{eq:coarse-global-hybrid-variation}; enlarging the common
constant once joins that estimate to the barrier bound.

We finally recover two derivatives without applying the nonlinear
recovery lemma.  The direct nontransport source is exactly
\(\mathcal S_\xi\) in \eqref{eq:linearized-direct-source}, and
\[
 \|\mathcal S_\xi(\tau)\|_
   {\mathfrak C_{{\rm sc},0}^{2,\alpha}}
 \leq Cs_\xi(\tau).
\]
Let \(\mathcal P_\tau\), based at the identity at \(\tau_0\), be the
time-ordered flow generated by the negative of the base phase field
\[
 V_{\rm ph}
 =-a\bar\nabla\bar f
  +\sum_{j=1}^8b_j\chi_\tau W_j .
\]
Its scale-normalized \(C^4\) distortion is bounded by
\(\exp(C\int q)\), hence uniformly by
\eqref{eq:linearized-background-phase-budget}.  Given a target
\((\tau_*,y_*)\), put \(x_0=\mathcal P_{\tau_*}^{-1}(y_*)\) and use
the exact backward clock and truncated face
\[
 \tau(s)=\tau_*-\log(1-s),\qquad
 s_-=\max\{-1,1-e^{\tau_*-\tau_0}\},\qquad
 \phi_s=\varphi_{-\log(1-s)}.
\]
Set
\[
 \mathcal H_\xi(s)
 =(1-s)\phi_s^*\mathcal P_{\tau(s)}^*k_\xi(\tau(s)).
\]
Choose the fixed nested raw balls
\(B_1\Subset B_2\) in the buffered scale-one chart exactly as in
\eqref{eq:two-state-nested-rescaled-balls}, put
\(Q_2=B_2\times[s_-,0]\), and let \(\widehat\nabla\) be the connection
of the terminal rescaled metric used to define those balls.
The residual drift is the base-state field
\[
 \widehat Z_s
 =(1-s)^{-1}\phi_s^*
  \bigl(\bar\nabla\bar f
       -\mathcal P_{\tau(s)}^*\bar\nabla\bar f\bigr).
\]
The proof of \eqref{eq:two-state-residual-drift-bound} uses only the
one-state phase budget and therefore gives the required
scale-normalized \(C^{3,\alpha}\) bound for \(\widehat Z_s\).
Let \(\mathcal Q_s\) solve
\(\partial_s\mathcal Q_s=-\widehat Z_s\circ\mathcal Q_s\),
\(\mathcal Q_0=\operatorname{Id}\), and put
\(\mathcal K_\xi=\mathcal Q_s^*\mathcal H_\xi\).
The phase threshold keeps its sampling tube inside the fixed
scale-one buffer, exactly as in
\eqref{eq:two-state-residual-drift-tube}.

After these two conjugations, define
\[
 \widetilde{\mathcal S}_\xi(s)
 :=\mathcal Q_s^*\phi_s^*
   \mathcal P_{\tau(s)}^*\mathcal S_\xi(\tau(s)).
\]
The factor \(1-s\) in \(\mathcal H_\xi\) cancels
\(d\tau/ds=(1-s)^{-1}\).  Consequently,
\[
 \|\widetilde{\mathcal S}_\xi(s)\|_{C^{2,\alpha}(B_2)}
 \leq Cs_\xi(\tau(s)).
\]
The drift-straightened principal matrix is uniformly elliptic; all
fixed lower-order coefficients have uniform \(C^{2,\alpha}\) bounds,
while the remaining coefficient
\((1-s)^{-1}a(\tau(s))\) is controlled in \(L^1_s\).  Let
\(\widetilde{\mathcal U}^{\rm lin}_{\tau_*}(s,r)\) be the Dirichlet
evolution family for this homogeneous drift-straightened linearized
operator on \(B_2\).  The proof of
\eqref{eq:C2-evolution-family-bound}, with the polarized operator
replaced by this base linearization, applies verbatim and gives the
linear Duhamel bound
\begin{equation}\label{eq:linearized-direct-source-C2}
 \left\|
  \int_{s_-}^s
   \widetilde{\mathcal U}^{\rm lin}_{\tau_*}(s,r)
   \widetilde{\mathcal S}_\xi(r)\,dr
 \right\|_{C^{2,\alpha}(B_1)}
 \leq C\int_{s_-}^s s_\xi(\tau(r))\,dr .
\end{equation}
Subtract this Duhamel term.  The ordinary interior Bernstein estimate
for the remaining scalar-principal system gives at \(s=0\)
\begin{equation}\label{eq:linearized-local-C2}
 \begin{split}
 |\widehat\nabla\mathcal K_\xi|
 +|\widehat\nabla^2\mathcal K_\xi|
 \leq C\bigg(&
  \|\mathcal K_\xi\|_{C^0(Q_2)}
  +\mathbf1_{\{s_->-1\}}
    \|\mathcal K_\xi(s_-)\|_{C^3(B_2)}\\
 &+\int_{s_-}^{0}s_\xi(\tau(r))\,dr
  +\|\widetilde{\E}_\xi^{\rm lin}\|_{C^2(Q_2)}
 \bigg),
 \end{split}
\end{equation}
where
\(\widetilde{\E}_\xi^{\rm lin}
=\mathcal Q_s^*\phi_s^*
\mathcal P_{\tau(s)}^*\E_\xi^{\rm lin}(\tau(s))\).

Each term in \eqref{eq:linearized-local-C2} has now been proved at the
required type.  The \(C^0\) term is controlled by
\eqref{eq:linearized-C0-barrier} and the transported-weight comparison
\eqref{eq:two-state-transported-weight-comparison}; the truncated
initial term is controlled by
\eqref{eq:linearized-C3-entrance}; and
\[
 \int_{s_-}^{0}s_\xi(\tau(r))\,dr
 \leq
 C\int_{\max\{\tau_0,\tau_*-\log2\}}^{\tau_*}
 s_\xi(r)\,dr
\]
is controlled by \eqref{eq:linearized-direct-tail}.  If the sampling
tube meets the graft annulus, the pointwise weight is comparable to
one and \eqref{eq:linearized-graft-column-package} controls the last
term after rescaling; otherwise that term vanishes.  Fixed cylinders
treat the compact core.  Undoing the residual, soliton, and phase
pullbacks therefore gives
\[
 \sum_{\ell=0}^2|\bar\nabla^\ell k_\xi(\tau,x)|
 \leq
 Cd_\xi
 \min\left\{
 e^{-\sigma_-(\tau-\tau_0)}
 (1+\bar f(x))^{\sigma_-},\,1\right\}.
\]
In particular this proves \eqref{eq:global-variation-pointwise} on
\(\{\bar f\leq e^\tau\}\).  All constants are uniform on the
common-margin ball and independent of the finite endpoint.  The
linearized energy, phase, barrier, and recovery estimates have thus
been obtained in the independent model tangent norm, without invoking
the nonlinear two-state barrier or derivative-recovery conclusions.

We prove convergence of derivatives without differentiating any
infinite-time limit before its \(C^1\) regularity is known.  Work in a
convex Banach chart whose closure lies in the common-margin ball, and
fix two finite endpoints \(S'>S\).  Proposition
\ref{prop:two-state-prepared-evolution} makes
\(T_S,\ell_S,\Psi_S\) and \(T_{S'},\ell_{S'},\Psi_{S'}\) \(C^1\) maps.
Since
\[
 \frac d{d\tau}\log(\lambda e^\tau)=-a,
\]
differentiation on the compact interval \([S,S']\) gives
\[
 D(\ell_{S'}-\ell_S)[\xi]
 =-\int_S^{S'}\dot a_\xi(s)\,ds.
\]
Equation \eqref{eq:global-variation-phase} makes this uniformly
\(o_S(1)\|\xi\|\), independently of \(S'\).  The finite-time identity
\[
 D\log\lambda(s)[\xi]
 =D\log\lambda(\tau_0)[\xi]
   -\int_{\tau_0}^{s}\dot a_\xi(\sigma)\,d\sigma
\]
therefore gives a uniform bound for \(D\log\lambda(s)[\xi]\).  Hence
\[
 |D(T_{S'}-T_S)[\xi]|
 =\left|\int_S^{S'}
       \lambda(s)D\log\lambda(s)[\xi]\,ds\right|
 \leq Ce^{-S}\|\xi\|.
\]

For the phase, use the finite time-ordered ODE
\[
 \partial_\tau\Psi_\tau
 =U_\tau\circ\Psi_\tau,\qquad
 U_\tau=\sum_{j=1}^8b_j(\tau)\chi_\tau W_j.
\]
On a compact set enlarged once to contain the relevant phase
trajectories, the variation \(Y_\xi=D\Psi_\tau[\xi]\) solves
\[
 \partial_\tau Y_\xi
 =DU_\tau(\Psi_\tau)Y_\xi
  +\sum_{j=1}^8\dot b_{\xi,j}
       (\chi_\tau W_j)(\Psi_\tau).
\]
Applying this equation to the relative flow from \(S\) to \(S'\), the
all-order bounds in Lemma~\ref{lem:mode-growth} and the radial
properness estimate \eqref{eq:phase-radial-properness} give
\begin{equation}\label{eq:C1-geometric-finite-endpoint-tail}
 \left\|D\!\left[
  \log(\Psi_S^{-1}\circ\Psi_{S'})
 \right][\xi]\right\|_{C^m(K)}
 \leq C_{K,m}\int_S^{S'}
       (|\dot b_\xi|+|b|\|\xi\|)\,ds
 =o_S(1)\|\xi\|,
\end{equation}
uniformly in \(S'\).  To convert the relative estimate into a fixed
chart estimate, choose \(S_0\) so that, on \(K\), all phase tails with
\(S\geq S_0\) lie in one exponential chart, and denote the resulting
coordinate map by \(\widehat\Psi_S\).  Smoothness of composition,
inversion, and the logarithm on that chart gives
\[
\begin{aligned}
 \|D\widehat\Psi_{S'}-D\widehat\Psi_S\|_{C^m(K)}
 &\leq C_{K,m}\bigl(
  \|D\log(\Psi_S^{-1}\circ\Psi_{S'})\|_{C^m(K)}\\
 &\qquad
  +\|\log(\Psi_S^{-1}\circ\Psi_{S'})\|_{C^{m+1}(K)}
    \|D\widehat\Psi_S\|_{C^m(K)}\bigr).
\end{aligned}
\]
The nondifferentiated phase tail tends to zero, and the variational
ODE gives a uniform bound for \(D\widehat\Psi_S\).  Thus the phase
derivatives, and hence all three finite geometric endpoint
derivatives, are uniformly Cauchy in fixed bounded-geometry charts.

We use the elementary Banach-space principle that if
\(F_S\in C^1(U,Y)\), \(F_S(z_*)\) converges at one point, and
\(DF_S\) is uniformly Cauchy on the convex set \(U\), then \(F_S\)
converges on \(U\) to a \(C^1\) map and
\(DF=\lim_SDF_S\).  This follows directly from the mean-value formula
on line segments.  The already proved, uniform nondifferentiated tails
identify the limits here as
\(T,\log\lambda_\infty,\Psi_\infty\).  Consequently these three maps
are \(C^1\).  We may therefore let \(S'\to\infty\) in the three finite
endpoint derivative estimates.  The resulting identities and bounds
are exactly \eqref{eq:C1-geometric-tail-uniform}; in particular no
derivative of a limit was used to establish its existence.

The adaptive variables require the same finite-endpoint argument.  Fix
an integer \(0\leq m\leq k+2\).  On every finite horizon, the exact
\(R\)-ODE makes \(R_S\) \(C^1\) into \(C^{k+3}\).  The identity
\[
 D(R^{-1})[\dot R]
 =-DR^{-1}\bigl(\dot R\circ R^{-1}\bigr)
\]
and the pullback formula for \(S(S)\) make \(R_S^{-1}\) and \(S(S)\)
\(C^1\) into \(C^{k+2}\).  These are finite-time ODE and composition
statements and use no regularity of an adaptive limit.  Set
\[
 Z_\tau
 =(\varphi_{-\tau})_*
   (a\bar\nabla\bar f-U_\tau),
 \qquad
 \partial_\tau R_\tau=Z_\tau\circ R_\tau
\]
as in \eqref{eq:relative-target-flow}.  On every finite interval,
\[
 \partial_\tau\dot R_\xi
 =DZ_\tau(R_\tau)\dot R_\xi
  +\dot Z_{\tau,\xi}(R_\tau),
 \qquad
 \|\dot Z_{\tau,\xi}\|_{C^m(\Omega_\eta^+)}
 \leq C_m|\dot c_\xi(\tau)|.
\]
Variation of constants for the relative map between the two finite
endpoints gives
\begin{equation}\label{eq:C1-R-tail-proof}
 \left\|
 D[\operatorname{Exp}_{R_S}^{-1}R_{S'}][\xi]
 \right\|_{C^m}
 \leq C_m\int_S^{S'}
 \bigl(|\dot c_\xi|+|c|\|\xi\|\bigr)\,d\tau
 =o_S(1)\|\xi\|,
\end{equation}
uniformly in \(S'\).  Here the mode bounds through the required order
and uniform radial properness control every composition constant.  At
order \(m+1\), the same estimate and the displayed inverse formula
give the corresponding finite-endpoint bound for
\(\operatorname{Exp}_{R_S^{-1}}^{-1}R_{S'}^{-1}\).
Choose \(S_0\) so that the adaptive tails lie in fixed
right-translated prepared charts, and write their coordinates as
\(\widehat R_S\) and \(\widehat R_S^{-1}\).  The uniform
\(C^{k+3}\) map and inverse-map bounds, smoothness of the chart
operations, and the nondifferentiated adaptive tails give
\begin{align*}
 \|D\widehat R_{S'}-D\widehat R_S\|_{C^m}
 &\leq C_m\left(
  \|D[\operatorname{Exp}_{R_S}^{-1}R_{S'}]\|_{C^m}
  +\|\operatorname{Exp}_{R_S}^{-1}R_{S'}\|_{C^{m+1}}
    \|D\widehat R_S\|_{C^m}
 \right),\\
 \|D\widehat R_{S'}^{-1}-D\widehat R_S^{-1}\|_{C^m}
 &\leq C_m\left(
  \|D[\operatorname{Exp}_{R_S^{-1}}^{-1}R_{S'}^{-1}]\|_{C^m}
  +\|\operatorname{Exp}_{R_S^{-1}}^{-1}R_{S'}^{-1}\|_{C^{m+1}}
    \|D\widehat R_S^{-1}\|_{C^m}
 \right).
\end{align*}
The finite-time variational ODE bounds the two endpoint derivatives
uniformly, so both right sides are \(o_S(1)\), uniformly in \(S'\).

Finally use the exact identity \eqref{eq:S-tau-exact}.  Differentiating
it at a finite time, using the finite-time bounds already established
for \(D\log\lambda[\xi]\), \(D R_\tau[\xi]\), and \(Dc[\xi]\), gives
\begin{equation}\label{eq:C1-S-time-derivative}
 \|D(\partial_\tau S_\tau)[\xi]\|_{C^m(\Omega_\eta^+)}
 \leq C_m\left[
  e^{-\tau}\|\xi\|
  +e^{-\tau}|\dot a_\xi|
  +|\dot b_\xi|
  +(e^{-\tau}|a|+|b|)\|\xi\|
 \right].
\end{equation}
Integrating this estimate only over \([S,S']\) shows that
\(D(S(S')-S(S))\) is uniformly \(o_S(1)\) in \(S'\), by
\eqref{eq:global-variation-phase}, the one-state phase tail, and
\(e^{-\tau}\in L^1\).

The Banach-space principle, now applied in the right-translated
prepared map charts, proves that the nondifferentiated adaptive limits
\(R_\infty,R_\infty^{-1},S_\infty\) are \(C^1\) and that their
derivatives are the uniform limits of the finite-endpoint derivatives.
Letting \(S'\to\infty\) in the last three finite-endpoint estimates
then proves \eqref{eq:C1-adaptive-tail-uniform}.  This establishes the
claimed \(C^1\) regularity of all six limiting objects in precisely the
orders stated in the theorem and
Corollary~\ref{cor:Lipschitz-asymptotic-data}.
\end{proof}

\begin{corollary}[Locally Lipschitz asymptotic data]
\label{cor:Lipschitz-asymptotic-data}
Let
\[
 \mathscr B\subset\Sigma_{\tau_0}^{k+2,\alpha}
\]
be a fixed-\(\tau_0\) sliced prepared-coordinate ball with common
strict margins.  On \(\mathscr B\), the maps
\[
 \mathbf z\longmapsto T(\mathbf z),\qquad
 \mathbf z\longmapsto\lambda_\infty(\mathbf z),\qquad
 \mathbf z\longmapsto\Psi_\infty(\mathbf z)
\]
are locally Lipschitz in the buffered prepared Banach topology.  The
first two are \(C^1\), and for every \(K\Subset M\) and \(m\geq0\) the
last is \(C^1\) into a \(C^m(K)\) exponential chart.
In every order carried by the prepared topology, the limits of
$R_\tau$, $R_\tau^{-1}$, and $S_\tau$ are locally Lipschitz and
\(C^1\) as well.

Separately, let
\[
 \mathscr O\subset\operatorname{dom}
    \Pi_{\rm sl}^{\,k+4\to k+2}
    \subset\mathscr P_{\tau_0}^{k+4,\alpha},
 \qquad
 \Pi_{\rm sl}^{\,k+4\to k+2}(\mathscr O)\subset\mathscr B,
\]
be a buffered ambient prepared neighborhood.  The ambient asymptotic
data maps are, by definition, the preceding sliced maps composed with
\(\Pi_{\rm sl}^{\,k+4\to k+2}\).  They have the same Lipschitz and
\(C^1\) conclusions, and their differentials are
\[
 D\mathcal A_{\rm amb}(\mathbf z)
 =
 D\mathcal A_{\rm sl}
  \bigl(\Pi_{\rm sl}^{\,k+4\to k+2}(\mathbf z)\bigr)
 \circ D\Pi_{\rm sl}^{\,k+4\to k+2}(\mathbf z)
\]
for each of the six asymptotic objects
\(\mathcal A\in\{T,\lambda_\infty,\Psi_\infty,
R_\infty,R_\infty^{-1},S_\infty\}\).
On each fixed \(\mathscr B'\Subset_{\rm u}\mathscr B\),
Proposition~\ref{prop:two-state-prepared-evolution},
Theorem~\ref{thm:global-two-state-estimate}, and the uniform derivative
tails \eqref{eq:C1-geometric-tail-uniform} and
\eqref{eq:C1-adaptive-tail-uniform} give, for every fixed
target component chart \((K,m)\), a finite constant
\(L_{{\rm asy},K,m}\) controlling both the sliced Lipschitz constant
and the sliced differential norm.  Consequently, if
\[
 \overline{\mathscr O}\subset\operatorname{dom}
 \Pi_{\rm sl}^{\,k+4\to k+2},
 \qquad
 \Pi_{\rm sl}^{\,k+4\to k+2}(\overline{\mathscr O})
 \subset\mathscr B'\Subset_{\rm u}\mathscr B,
\]
then the corresponding ambient constants are at most
\(L_{{\rm asy},K,m}K_{\Pi,k+2}\).  The inclusions retain the fixed
charts and margins; no compactness of \(\overline{\mathscr O}\) is
used.
\end{corollary}

\begin{proof}
On the sliced ball, at every fixed finite normalized time $S$,
Proposition~\ref{prop:two-state-prepared-evolution} and
Proposition~\ref{prop:uniform-receding-phase} give locally Lipschitz
and \(C^1\)
dependence of
\[
 t(S),\quad\lambda(S),\quad\Psi_S,\quad
 R_S,\quad R_S^{-1},\quad S(S).
\]
  Uniformly on a common-margin neighborhood,
  Theorem~\ref{thm:prepared-entrance-continuation} gives
\[
 |T-t(S)|\leq Ce^{-S},\qquad
 \left|\log\frac{\lambda(S)e^S}{\lambda_\infty}\right|
 +\|\Psi_\infty-\Psi_S\|_{C^m(K)}
 \leq C_{K,m}e^{-2\theta S},
\]
and the same tail bound for the adaptive exterior limits.  The
difference and first-variation versions of these tails are
Theorem~\ref{thm:global-two-state-estimate}, in particular
\eqref{eq:C1-geometric-tail-uniform} and
\eqref{eq:C1-adaptive-tail-uniform}.  Thus each limiting
object is a uniform \(C^1\) limit of finite-horizon maps.
The ambient assertion is then exactly the Banach-space chain rule
applied to the fixed \(C^1\) phase retraction
\(\Pi_{\rm sl}^{\,k+4\to k+2}\); no coercive estimate is being claimed
for a tangent vector transverse to the slice before this retraction.
\end{proof}

\section{The first stable profile and geometric tail}
\label{sec:first-stable-profile}

Throughout this section fix a continuation output order \(k\geq12\)
and \(0<\alpha<1\).  Thus every sliced entrance family used for
differentiation lies in
\(\Sigma_{\tau_0}^{k+2,\alpha}\), and every ambient phase retraction
is taken with the two additional input derivatives displayed in its
statement.  The relative \(C^{2,\alpha}\) topology appearing in the
physical profile theorem is obtained only after the fixed
positive-time gauge restart; it does not replace these internal
prepared regularity indices.

Theorem~\ref{thm:prepared-entrance-continuation} gives the initial
exponential rate.  The asymptotically autonomous stable equation
improves this rate to the exact spectral threshold and yields a first
stable profile together with the first nonlinear corrections to the
geometric parameters.

We proceed from one-state dynamics to geometry and then to parameter
dependence.  The stable normal form first produces \(V_\infty\) and
the quadratic response coefficient.  Quadratic modal scattering then
determines the first corrections to feedback, scale, and phase, while
the frozen-window argument converts those expansions into the sharp
marked-spacetime profile.  The subsequent two-state and variational
scattering estimates promote the asymptotic data to \(C^1\)
coordinates and yield the foliation.  The quadratic and sharp
marked-spacetime conclusions depend only on the preceding one-state
profile and quadratic coefficient.

\subsection{Stable normal form}

Let
\[
 0<\gamma_1<\gamma_2<\cdots,\qquad
 E_j=\ker(\A+\gamma_jI)\subset\mathcal Z^\perp,
\]
be the distinct stable rates and eigenspaces of $\A$, and let
$\Pi_j$ denote the $L^2_\nu$-orthogonal projection onto $E_j$.
If there is no second distinct stable rate, set $\gamma_2=\infty$.
Thus $\gamma_1$ is the exact stable gap, whereas $\beta$ in
\eqref{eq:beta} is an arbitrary strict lower bound for it.
We write
\[
 L^2_{\nu,-}:=\mathcal Z^\perp\subset L^2_\nu,\qquad
 H^1_{\nu,-}:=H^1_\nu\cap\mathcal Z^\perp,
 \qquad
 H^{-1}_{\nu,-}:=(H^1_{\nu,-})^*
\]
for the stable form space and its dual; the orthogonal complement is
understood in \(L^2_\nu\).  Equivalently, \(H^1_{\nu,-}\) is the
form-norm closure of finite sums of stable eigentensors.  This
equivalence follows from compact resolvent and the spectral
characterization of the form domain.

For the remainder of this section, fix once and for all
\begin{equation}\label{eq:section-wide-stable-kappa}
 0<\kappa<\frac{\sigma}{2},
\end{equation}
and decrease the common small prepared box once.  All inner
coefficient-decay, projected-graph, and one-state tame estimates below
use this same exponent; the same is true of the difference and
variational estimates.  The choice is otherwise arbitrary in
\((0,\sigma/2)\); every later
\(\delta_0\) and \(\zeta\) constrained by \(\kappa\) is chosen only
after \eqref{eq:section-wide-stable-kappa}.

\begin{remark}[The first stable block is not numerically identified]
\label{rem:first-stable-block-abstract}
Corollary~\ref{cor:self-contained-Hessian-families} gives a
six-real-dimensional Hessian family with
eigenvalue $1-\sqrt2$.  Hence
\[
 \gamma_1\leq\sqrt2-1<\frac12.
\]
The remaining exact block estimates prove strict negativity, but do not
give a uniform Rayleigh bound by $1-\sqrt2$ on the complement of that
family.  In particular, the stated estimates do not exclude spectrum
in $(1-\sqrt2,0)$ or additional multiplicity at $1-\sqrt2$.
Accordingly, neither $\gamma_1=\sqrt2-1$ nor a numerical
identification of $E_1$ is used below.
\end{remark}

\begin{remark}[Marked and gauge-fixed meaning of the profile]
\label{rem:stable-profile-gauge-scope}
The stable spectral spaces are those of the fixed Ricci--DeTurck
operator after the nonnegative geometric block \(\mathcal Z\) has been
removed.  They are not quotients by all stable infinitesimal
diffeomorphisms.  In particular, the stable Hessian family in
Corollary~\ref{cor:self-contained-Hessian-families} consists, up to
normalization, of tensors
\[
 \bar\nabla^2\phi=\frac12\Lie_{\bar\nabla\phi}\bar g .
\]
If its rate realizes the stable gap, these are pure-gauge directions
inside \(E_1\).  Consequently \(V_\infty\), the profile fibers, and the
zero-profile condition below are invariants of the fixed cutoffs
\((\rho,\chi,\Gamma)\), prepared marking, harmonic-map gauge, and
transported restart convention.  They
are not claimed to be invariant under an arbitrary time-dependent
re-marking, and no quotient theorem by stable Lie derivatives is used.
The nonvanishing and sharp-rate conclusions below are correspondingly
statements about the marked, gauge-fixed error.
\end{remark}

\begin{lemma}[Polynomial growth of stable eigenmodes]
\label{lem:stable-eigenmode-growth}
If $V\in E_j$, then for every $m\geq0$ and every $\delta>0$,
\begin{equation}\label{eq:stable-eigenmode-growth}
 |\bar\nabla^mV|
 \leq C_{V,m,\delta}(1+\bar f)^{\gamma_j+\delta}.
\end{equation}
In particular, every contraction of finitely many derivatives of
stable eigenmodes against a fixed polynomially growing tensor is
integrable with respect to $d\nu$.
\end{lemma}

\begin{proof}
The FIK shrinker is complete, connected, nonflat, and asymptotically
conical, so Propositions~2.34--2.35 of~\cite{Stolarski} apply to
\(\bar\Delta_{\bar f}+2\overline{\Rm}\).  Their eigenvalue
\(\lambda_j\) is \(-\gamma_j\) in our notation; hence
\(\max\{-\lambda_j,0\}=\gamma_j\), and Proposition~2.35 gives
\eqref{eq:stable-eigenmode-growth}.  The Gaussian density absorbs every
resulting polynomial.
\end{proof}

\begin{lemma}[Restricted stable semigroups and concrete projected graph recovery]
\label{lem:restricted-stable-semigroup}
For \(j\geq1\), let \(P_{\geq j}\) be the spectral projection of
\(L=-\A|_{\mathcal Z^\perp}\) onto the closed sum of the eigenspaces
with rates at least \(\gamma_j\).  If that sum is zero, all assertions
below are read trivially.  The projection \(P_{\geq j}\) is bounded on
\(H^1_{\nu,-}\) and, by duality, on \(H^{-1}_{\nu,-}\).  For \(u>0\),
\begin{align}
 \|e^{-uL}P_{\geq j}\|_{H^{-1}_{\nu,-}\to L^2_{\nu,-}}
 &\leq C(1+u^{-1/2})e^{-\gamma_ju},
 \label{eq:restricted-Hminus1-semigroup}\\
 \|e^{-uL}P_{\geq j}\|_{H^1_{\nu,-}\to H^1_{\nu,-}}
 &\leq Ce^{-\gamma_ju}.
 \label{eq:restricted-H1-semigroup}
\end{align}

The following endpoint estimate applies only to the concrete sliced
solution classes listed here; the displayed projected equation by
itself is not a hypothesis sufficient for the conclusion.
\begin{enumerate}
\item[(i)] In the one-state case,
      \(U=H=\rho_\tau h\) is the output of
      Theorem~\ref{thm:prepared-entrance-continuation} from a strict
      prepared entrance.  It satisfies the exact prepared graph and
      slice, the common global \(C^2\) box
      \eqref{eq:master-bootstrap-global}, and the feedback and
      \(H^1_\nu\) estimates
      \eqref{eq:master-bootstrap-H1}--%
      \eqref{eq:master-bootstrap-velocity}.
\item[(ii)] In the difference case,
      \(U=D=H_1-H_2\), where the two states start in one
      common-margin sliced ball, are placed in the common anchored
      exterior gauge, and are the outputs covered by
      Theorem~\ref{thm:global-two-state-estimate}.  In particular, the
      estimates
      \eqref{eq:global-two-state-energy}--%
      \eqref{eq:global-two-state-phase},
      \eqref{eq:global-two-state-forcing}--%
      \eqref{eq:global-two-state-forcing-C2}, and the order-four
      current-window propagation
      \eqref{eq:uniform-restart-low-propagation} are part of the
      admissible class.
\item[(iii)] In the variational case,
      \(U=K=DH_{\mathbf z_0}[w]\) is the first variation furnished by
      the differentiability clause of
      Theorem~\ref{thm:global-two-state-estimate}, with
      \eqref{eq:global-variation-energy}--%
      \eqref{eq:global-variation-pointwise} and the linearized
      counterparts of the preceding forcing, column, graft, and
      order-four current-window estimates.  In a late-entry
      application, \(K\) is instead an exact tangent evolution from a
      uniform late-entry \(W\)-profile disk; the estimate below is
      invoked only after
      \eqref{eq:uniform-restart-low-initial-propagation} and
      \eqref{eq:late-entry-hybrid-outer-consequence} have been
      established from the low-initial-data and separated-memory
      package.
\end{enumerate}
All graph, action, feedback, cutoff, effective-column, and graft terms
in \(F_U\) are the exact terms of the corresponding admissible sliced
equation, with the common bounds just listed; they are not replaceable
by an arbitrary \(H^{-1}_\nu\) forcing.  For any one of these three
classes put
\[
 R=P_{\geq j}U,\qquad P_{<j}=I-P_{\geq j}
 \quad\hbox{on }\mathcal Z^\perp .
\]
Let \(s_R\) denote the corresponding input size
\[
 s_R=
 \begin{cases}
  1,&\text{for a one-solution remainder},\\
  d_0,&\text{for a two-state difference, where \(d_0\) is the
             prepared entrance distance},\\
  \|w\|_{T_{\mathbf z_0}\Sigma_{\tau_0}^{k+2,\alpha}},
     &\text{for an ordinary first variation generated by \(w\)},\\
  \|\xi_{v,w}\|_{{\rm le},\tau_0},
     &\text{for the late-entry variation generated by
              \(\xi_{v,w}\)}.
 \end{cases}
\]
For the fixed \(\kappa\) in
\eqref{eq:section-wide-stable-kappa}, set
\begin{equation}\label{eq:projected-graph-lower-source}
 \begin{split}
  \alpha_R(\tau)
  &=
  \begin{cases}
   e^{-\kappa\tau}+\|H(\tau)\|_{H^1_\nu},
      &U=H,\\
   e^{-\kappa\tau}+\|H_1(\tau)\|_{H^1_\nu}
                    +\|H_2(\tau)\|_{H^1_\nu},
      &U=D,\\
   e^{-\kappa\tau}+\|H(\tau)\|_{H^1_\nu},
      &U=K,
  \end{cases}\\
  g_R(\tau)&=\alpha_R(\tau)
              \|P_{<j}U(\tau)\|_{H^1_\nu},\\
  \mathfrak e_R(\tau)&=s_R^2e^{-ce^{\tau/2}} .
 \end{split}
\end{equation}
Here \(c>0\) is fixed by the common-margin ball.  If
\[
 \partial_\tau R+LR=F_R,\qquad F_R=P_{\geq j}F_U,
\]
is the equation obtained by projecting the corresponding exact
one-state, difference, or variational sliced equation, then, for
\(\tau\geq\tau_0+1\),
\begin{equation}\label{eq:projected-stable-graph-recovery}
 \|R(\tau)\|_{H^1_\nu}^2
 \leq
 C\sup_{\tau-1\leq s\leq\tau}\|R(s)\|_{L^2_\nu}^2
 {}+C\int_{\tau-1}^{\tau}g_R(s)^2\,ds
 {}+C\mathfrak e_R(\tau).
\end{equation}
For fixed \(j\), all constants are uniform on fixed common-margin
prepared balls and for all sufficiently late entrance times.
In particular, when \(j=1\), the lower-block term \(g_R\) vanishes.
The endpoint estimate is not asserted for an arbitrary weak equation
with only \(H^{-1}_\nu\) forcing.
\end{lemma}

\begin{proof}
On \(\mathcal Z^\perp\), coercivity and boundedness of
\(\overline{\Rm}\) give
\[
 \|q\|_{H^1_\nu}^2
 \asymp
 \|(I+L)^{1/2}q\|_{L^2_\nu}^2.
\]
The dual norm is therefore equivalent to the norm defined by
\((I+L)^{-1/2}\).  The spectral theorem now reduces
\eqref{eq:restricted-Hminus1-semigroup} to
\[
 \sup_{\lambda\geq\gamma_j}
 (1+\lambda)^{1/2}e^{-u\lambda}
 \leq C(1+u^{-1/2})e^{-\gamma_ju},
\]
and gives \eqref{eq:restricted-H1-semigroup} directly.  It also shows
that every spectral projection is bounded on the form space and its
dual.

 We next isolate the additional structure which permits endpoint
 recovery.  Choose once and for all
 \(\chi_{\rm in}\in C^\infty([0,\infty);[0,1])\), equal to one on
 \([0,1]\) and zero on \([2,\infty)\), and put
 \[
  \chi_{{\rm in},\tau}(x)
  =\chi_{\rm in}\!\left(e^{-\tau/2}\bar f(x)\right).
 \]
 After increasing the entrance time, \(\rho_\tau=1\) on
 \(\supp\chi_{{\rm in},\tau}\).  Before applying the nonlocal
 projection \(P_{\geq j}\), we define the following decomposition
 globally on \(M\):
\begin{equation}\label{eq:concrete-projected-graph-splitting}
 F_U
 =
 A_U^{ij}\bar\nabla_i\bar\nabla_jR
 +\mathscr L_{U,1}R+S_R+\mathcal O_R+\mathcal Z_U .
\end{equation}
The terms are chosen so that
\begin{align}
 \|A_U^{ij}\bar\nabla_i\bar\nabla_jR
      +\mathscr L_{U,1}R\|_{L^2_\nu}
 &\leq
 \varepsilon_*\bigl(
  \|LR\|_{L^2_\nu}+\|R\|_{H^1_\nu}\bigr),
 \label{eq:concrete-graph-linear-bound}\\
 \|S_R\|_{L^2_\nu}
 &\leq Cg_R,
 \label{eq:concrete-graph-source-bound}\\
 \|\mathcal O_R\|_{L^2_\nu}
 &\leq Cs_Re^{-ce^{\tau/2}} .
 \label{eq:concrete-graph-outer-bound}
\end{align}
Here \(\mathcal Z_U(\tau)\in\mathcal Z\); consequently
\[
 P_{\geq j}\mathcal Z_U=0,
 \qquad
 \ip{LR}{\mathcal Z_U}=0.
\]
Here \(\varepsilon_*>0\) can be made arbitrarily small by decreasing
the common \(C^2\) box and taking the entrance time later.

 We verify this splitting rather than invoking an abstract regularity
 principle.  On
 \(\Omega_\tau=\{\bar f\leq e^{\tau/2}\}\) one has
 \(\chi_{{\rm in},\tau}=\rho_\tau=1\), so \(h=H\).  For one solution,
 write
\(H=R+H_{<j}\), \(H_{<j}=P_{<j}H\), in every factor of the exact
formula \eqref{eq:Q-schematic}.  In the principal term this gives
\[
 -\widehat H*\bar\nabla^2H
 =
 -\widehat H*\bar\nabla^2R
 -\widehat H*\bar\nabla^2H_{<j}.
\]
In each \(\bar\nabla H*\bar\nabla H\) term retain one actual
\(\bar\nabla H\) as a coefficient and split the other factor into
\(\bar\nabla R+\bar\nabla H_{<j}\); do the same with one factor in
\(\overline{\Rm}*H*H\).  Thus the terms containing \(R\) form a
second-order operator
\[
 A_H^{ij}\bar\nabla_i\bar\nabla_jR+
 B_H^i\bar\nabla_iR+C_HR
\]
 whose coefficients are bounded by the actual \(C^2\) size of \(h\),
 while every other term contains \(H_{<j}\).  Multiply this inner
 expression by \(\chi_{{\rm in},\tau}\): its terms containing \(R\)
 define \(A_U\) and \(\mathscr L_{U,1}\), and its terms containing
 \(H_{<j}\) define \(S_R\).  Let \(\mathcal Z_U\) be the combination
 of the global direct columns \(Y_0,\ldots,Y_8\) in the corresponding
 one-state, difference, or variational equation.  Finally define
 \[
  \mathcal O_R
  :=F_U-
   A_U^{ij}\bar\nabla_i\bar\nabla_jR
   -\mathscr L_{U,1}R-S_R-\mathcal Z_U .
 \]
 This definition includes the transition annulus of
 \(\chi_{{\rm in},\tau}\), all derivatives of \(\rho_\tau\), the graft
 term, and the effective-column defects.  It makes
 \eqref{eq:concrete-projected-graph-splitting} an exact global identity,
 not a core identity later used under a global Gaussian integral.

For a difference, the exact identity
\[
 \Q(h_1)-\Q(h_2)
 =
 \int_0^1D\Q_{h_2+s(h_1-h_2)}[D]\,ds
\]
and \(D=R+P_{<j}D\) give the identical decomposition.  The metrics on
the segment remain in the common \(C^2\) box.  The same formula with
\(D\Q_h[K]\), followed by \(K=R+P_{<j}K\), gives the variational
decomposition.  This is the polarization underlying
Remark~\ref{rem:localized-Lipschitz}.

On \(\supp\chi_{{\rm in},\tau}\), the recovered \(C^2\) bounds, the exact feedback
estimates, \eqref{eq:A-graph-estimate}, and
\eqref{eq:B-graph-bound} now give
\eqref{eq:concrete-graph-linear-bound}.  Every remaining core term
contains \(P_{<j}U\) and a coefficient bounded by \(C\alpha_R\), which
gives \eqref{eq:concrete-graph-source-bound}.  Here \(P_{<j}\) has
finite rank, so polynomial eigenmode growth and finite-dimensional norm
 equivalence control every derivative of the lower block required in
 this \(L^2_\nu\) estimate.  In the two-state case the coefficient and
 feedback differences are bounded by
\eqref{eq:global-two-state-pointwise} and
\eqref{eq:two-state-sharp-modal}.  For a first variation, use instead
the already established linearized pointwise estimate
\eqref{eq:global-variation-pointwise} and the linearized Gram-system
estimate in the direct variational argument following
\eqref{eq:exact-variational-h-equation}.  The global direct columns
\(Y_\mu\), \(0\leq\mu\leq8\), form \(\mathcal Z_U\) and vanish after
the fixed spectral projection \(P_{\geq j}\); only their
 effective-column defects remain in \(\mathcal O_R\).  By its global
definition, \(\mathcal O_R\) contains the part of the nonlinear and
action expansions outside
\(\{\bar f\leq e^{\tau/2}\}\).  All transition terms, derivatives of
\(\rho_\tau\), effective-column defects, and graft terms are supported
there as well, namely where
 \(\bar f\geq ce^{\tau/2}\); the common global \(C^2\) estimates,
 \eqref{eq:global-two-state-forcing-C2}, and
 \eqref{eq:two-column-C2}, followed by the Gaussian tail estimate, give
 \eqref{eq:concrete-graph-outer-bound}.
The outer conclusion does not require those scale-sharp pointwise
estimates as an independent hypothesis.  Indeed,
\eqref{eq:uniform-restart-low-propagation} gives the complete
order-four hybrid difference or variation on the current unit window
with at most coarse exponential growth; for late-entry variations one
uses instead the exact low-initial-data clause
\eqref{eq:uniform-restart-low-initial-propagation}.  The same-order
graft estimate
\eqref{eq:pure-graft-buffered-difference} puts the graft in its
scale-normalized \(C^{2,\alpha}\) source space, while the prepared
chart calculus applied to the same low-order \(R,F\), and marking
blocks gives that space for the effective-map defects.  Because every
such term is supported where
\(\bar f\geq ce^{\tau/2}\), the Gaussian absorbs the coarse exponential
and every chart-scale polynomial.  This gives
\eqref{eq:concrete-graph-outer-bound} directly and provides an
independent route to the outer part of the projected graph recovery.
This verification also proves
\begin{equation}\label{eq:concrete-projected-Hminus1-bound}
 \|F_R\|_{H^{-1}_\nu}
 \leq
 C\varepsilon_*\|R\|_{H^1_\nu}
 +Cg_R+Cs_Re^{-ce^{\tau/2}} .
\end{equation}

Test the weak equation against \(R\).  Coercivity, Young's inequality,
\eqref{eq:concrete-projected-Hminus1-bound}, and absorption of the
\(\varepsilon_*\)-term give
\begin{equation}\label{eq:projected-window-energy}
 \int_{\tau-1}^{\tau}\|R(s)\|_{H^1_\nu}^2\,ds
 \leq
 C\sup_{\tau-1\leq s\leq\tau}\|R(s)\|_{L^2_\nu}^2
 {}+C\int_{\tau-1}^{\tau}g_R(s)^2\,ds
 {}+C\mathfrak e_R(\tau).
\end{equation}

To recover the endpoint, put \(X_R=\ip{LR}{R}\).  Since \(LR\) belongs
to the range of \(P_{\geq j}\), self-adjointness gives
\[
 \ip{LR}{P_{\geq j}F_U}=\ip{LR}{F_U}.
\]
Thus the pairing may be estimated before the nonlocal projection.
Equations~\eqref{eq:concrete-projected-graph-splitting}--%
\eqref{eq:concrete-graph-outer-bound}, the weighted graph estimate,
and Young's inequality give the \(L^2_\nu\) pairing estimate
\begin{equation}\label{eq:projected-graph-pairing}
 2|\ip{LR}{F_R}|
 \leq
 \|LR\|_{L^2_\nu}^2+CX_R+Cg_R^2+C\mathfrak e_R .
\end{equation}
Thus this step pairs \(LR\) and \(F_R\) in \(L^2_\nu\).
Differentiating \(X_R\) along the smooth exact equation now gives
\begin{equation}\label{eq:projected-graph-differential}
 X_R'+\|LR\|_{L^2_\nu}^2
 \leq C X_R+Cg_R^2+C\mathfrak e_R .
\end{equation}

By \eqref{eq:projected-window-energy} and the form equivalence, there
is \(s_*\in[\tau-1,\tau-\tfrac12]\) such that \(X_R(s_*)\) is bounded
by twice the right-hand side of
\eqref{eq:projected-window-energy}.  Apply ordinary forward Gronwall
to \eqref{eq:projected-graph-differential} on \([s_*,\tau]\).
On this unit window the Gaussian error at time \(s\) is bounded by the
error at time \(\tau\), after decreasing \(c>0\).  This bounds
\(X_R(\tau)\) by the same right-hand side.  The form equivalence for
\(X_R\) proves \eqref{eq:projected-stable-graph-recovery}.
\end{proof}

\begin{lemma}[Stable Duhamel bootstrap with graph recovery]
\label{lem:stable-rate-bootstrap}
Let \(L=-\A|_{\mathcal Z^\perp}\), so that
\(\operatorname{spec}L\subset[\gamma_1,\infty)\).  Suppose
\[
 \partial_\tau H+LH=\mathcal N,\qquad H\perp\mathcal Z,
\]
weakly in \(H^{-1}_\nu\), and that the exact \(H\) equation gives the
time-local graph estimate
\begin{equation}\label{eq:stable-graph-recovery}
 \|H(\tau)\|_{H^1_\nu}^2
 \leq
 C\sup_{\tau-1\leq s\leq\tau}
   \|H(s)\|_{L^2_\nu}^2
 +Ce^{-ce^{\tau/2}}
\end{equation}
for \(\tau\geq\tau_0+1\).
Assume, for some \(r>0\) and \(\kappa>0\),
\[
 \|H(\tau)\|_{H^1_\nu}\leq C_re^{-r\tau}
\]
and
\[
 \|\mathcal N(\tau)\|_{H^{-1}_\nu}
 \leq
 C\left(
  e^{-\kappa\tau}\|H(\tau)\|_{H^1_\nu}
  +\|H(\tau)\|_{H^1_\nu}^2
  +e^{-ce^{\tau/2}}\right).
\]
For every \(0<\epsilon<\gamma_1\), put
\[
 \mathfrak T_\epsilon(r)
 =
 \min\{\gamma_1-\epsilon,r+\kappa,2r\}.
\]
Then
\begin{equation}\label{eq:stable-rate-improvement}
 \|H(\tau)\|_{H^1_\nu}
 \leq C_{r,\epsilon}
 e^{-\mathfrak T_\epsilon(r)\tau}.
\end{equation}
The constants are uniform on common-margin prepared families.  If
\[
 \|H(\tau_0)\|_{H^1_\nu}\leq Ce^{-\gamma_1\tau_0},
\]
the constant is uniform for all sufficiently late entrance times.
\end{lemma}

\begin{proof}
The tame estimate first gives
\[
 \|\mathcal N(\tau)\|_{H^{-1}_\nu}
 \leq C_r\left(
 e^{-(r+\kappa)\tau}
 +e^{-2r\tau}
 +e^{-ce^{\tau/2}}\right).
\]
The restricted spectral estimate
\eqref{eq:restricted-Hminus1-semigroup}, with \(j=1\), gives
\[
 \|e^{-uL}\|_{H^{-1}_{\nu,-}\to L^2_{\nu,-}}
 \leq C(1+u^{-1/2})e^{-\gamma_1u}.
\]
Duhamel's formula and the elementary convolution estimate therefore
give the \(L^2_\nu\) rate
\(\mathfrak T_\epsilon(r)\); the arbitrarily small loss
\(\epsilon\) covers equality of a forcing exponent with
\(\gamma_1\).

Since \(H=P_{\geq1}H\), the lower-block source in
\eqref{eq:projected-graph-lower-source} is zero.  Inserting the
improved \(L^2_\nu\) rate into
\eqref{eq:stable-graph-recovery} therefore proves
\eqref{eq:stable-rate-improvement}.  The graph estimate itself is the
structural estimate
\eqref{eq:projected-stable-graph-recovery} from
Lemma~\ref{lem:restricted-stable-semigroup}; it is not being inferred
from the abstract \(H^{-1}_\nu\) bound alone.  Finally,
\eqref{eq:restricted-H1-semigroup} shows that the homogeneous term in
the late-entrance case satisfies
\[
 \|e^{-(\tau-\tau_0)L}H(\tau_0)\|_{H^1_\nu}
 \leq Ce^{-\gamma_1\tau},
\]
 which gives the uniformity assertion.
\end{proof}

\begin{lemma}[Linearized stable Duhamel bootstrap]
\label{lem:linearized-stable-rate-bootstrap}
Fix a sliced prepared base state
\(\mathbf z_0\in\Sigma_{\tau_0}^{k+2,\alpha}\), a tangent datum
\[
 w\in T_{\mathbf z_0}\Sigma_{\tau_0}^{k+2,\alpha},
 \qquad
 \|w\|:=
 \|w\|_{T_{\mathbf z_0}\Sigma_{\tau_0}^{k+2,\alpha}},
\]
and numbers
\[
 \kappa>0,\qquad
 0<\zeta<\min\{\kappa,\gamma_1/2\}.
\]
Let \(H=H_{\mathbf z_0}\) be the sliced stable tensor of the base
evolution, and let \(K=DH_{\mathbf z_0}[w]\).  Since
\(w\in\ker D\mathfrak m(\mathbf z_0)\), differentiation of the exact
slice gives \(K\perp\mathcal Z\).  Suppose \(K\) solves
\[
 \partial_\tau K+LK=\mathcal N_K
\]
and satisfies the concrete \(j=1\) graph estimate
\begin{equation}\label{eq:linearized-stable-graph-recovery}
 \|K(\tau)\|_{H^1_\nu}^2
 \leq
 C\sup_{\tau-1\leq s\leq\tau}\|K(s)\|_{L^2_\nu}^2
 +C\|w\|^2e^{-ce^{\tau/2}} .
\end{equation}
Suppose
\[
 \|H(\tau)\|_{H^1_\nu}
 \leq C_\zeta e^{-(\gamma_1-\zeta)\tau}
\]
and
\begin{equation}\label{eq:linearized-stable-tame-bound}
 \|\mathcal N_K(\tau)\|_{H^{-1}_\nu}
 \leq
 C\bigl(e^{-\kappa\tau}
        +\|H(\tau)\|_{H^1_\nu}\bigr)
   \|K(\tau)\|_{H^1_\nu}
 +Ce^{-ce^{\tau/2}}\|w\|.
\end{equation}
If one has a coarse estimate
\[
 \|K(\tau)\|_{H^1_\nu}
 \leq C e^{-q(\tau-\tau_0)}\|w\|
\]
for some \(q>0\), then
\begin{equation}\label{eq:linearized-stable-rate}
 \|K(\tau)\|_{H^1_\nu}
 \leq C_\zeta
 e^{-(\gamma_1-\zeta)(\tau-\tau_0)}\|w\|.
\end{equation}
The constants are uniform on fixed common-margin sliced prepared
balls.
For the buffered ambient version, let
\[
 \mathbf z_0^{\rm amb}\in\operatorname{dom}
   \Pi_{\rm sl}^{\,k+4\to k+2},
 \qquad
 w_{\rm amb}\in
 T_{\mathbf z_0^{\rm amb}}\mathscr P_{\tau_0}^{k+4,\alpha},
\]
and put
\[
 \widetilde{\mathbf z}_0
 :=\Pi_{\rm sl}^{\,k+4\to k+2}(\mathbf z_0^{\rm amb}),
 \qquad
 \widetilde w
 :=D\Pi_{\rm sl}^{\,k+4\to k+2}(\mathbf z_0^{\rm amb})
    [w_{\rm amb}].
\]
Assume, as in the sliced statement above, that
\(\widetilde{\mathbf z}_0\) lies in the applicable fixed common-margin
sliced ball and that its base evolution satisfies the preceding
hypotheses.  Apply the lemma to the sliced pair
\((\widetilde{\mathbf z}_0,\widetilde w)\).  Then
\[
 \|\widetilde w\|_{
   \mathscr E_{\rm prep}^{k+2,\alpha}}
 \leq K_{\Pi,k+2}
 \|w_{\rm amb}\|_{\mathscr E_{\rm prep}^{k+4,\alpha}}.
\]
Thus the ambient constant is at most the sliced constant times
\(K_{\Pi,k+2}\).
\end{lemma}

\begin{proof}
The semigroup and graph-recovery argument in
Lemma~\ref{lem:stable-rate-bootstrap} applies verbatim, with the
forcing exponents generated by
\[
 q\longmapsto
 \min\{\gamma_1-\zeta,\ q+\kappa,\ q+\gamma_1-\zeta\}.
\]
Starting from any \(q>0\), finitely many iterations reach
\(\gamma_1-\zeta\).  The tail term is faster than every exponential,
and the fixed entrance time is absorbed into the uniform constant.
\end{proof}

\begin{lemma}[Stable tame normal form]
\label{lem:quadratic-stable-normal-form}
For the global evolution in
Theorem~\ref{thm:prepared-entrance-continuation}, the sliced tensor
$H=\rho_\tau h$ satisfies, weakly on $\mathcal Z^\perp$,
\begin{equation}\label{eq:stable-normal-form}
 \partial_\tau H=\A H+\mathcal N(\tau),\qquad H\perp\mathcal Z,
\end{equation}
Let \(\Pi_-\) denote the \(L^2_\nu\)-orthogonal projection onto
\(\mathcal Z^\perp\), together with its continuous dual extension to
\(\Pi_-:H^{-1}_\nu\to H^{-1}_{\nu,-}\), and define
\begin{equation}\label{eq:stable-N-definition}
 \begin{split}
 \mathcal N=\Pi_-\bigg[
  &\rho_\tau\Q(h)
   +\sum_{j=0}^8c_j\rho_\tau\mathscr B_jh
   +\mathcal C_\rho[h]+\rho_\tau\E\\
  &+\sum_{j=0}^8c_j
       \bigl(\rho_\tau\mathcal Y_{j,\tau}-Y_j\bigr)
 \bigg],
 \end{split}
\end{equation}
$c=(a,b_1,\ldots,b_8)$,
\[
 \mathscr B_0h=h-\Lie_{\bar\nabla\bar f}h,\qquad
 \mathscr B_jh=\Lie_{W_j}h\quad(1\leq j\leq8)
\]
on $\supp\rho_\tau$, and $\mathcal Y_{j,\tau}$ are the effective
columns.  For the fixed \(\kappa\) in
\eqref{eq:section-wide-stable-kappa},
\begin{equation}\label{eq:stable-tame-bound}
 \|\mathcal N(\tau)\|_{H^{-1}_\nu}
 \leq
 C e^{-\kappa\tau}\|H(\tau)\|_{H^1_\nu}
 +C\|H(\tau)\|_{H^1_\nu}^2
 +Ce^{-ce^{\tau/2}}.
\end{equation}
The section-wide choice was arbitrary in
\((0,\sigma/2)\).  Consequently, for every
$\epsilon>0$,
\begin{equation}\label{eq:near-gap-H1}
 \|H(\tau)\|_{H^1_\nu}
 \leq C_\epsilon e^{-(\gamma_1-\epsilon)\tau}.
\end{equation}
If
\begin{equation}\label{eq:profile-delta-choice}
 0<\delta_0<
 \frac14\min\{\kappa,\gamma_1,\gamma_2-\gamma_1\},
\end{equation}
with the last entry omitted when $\gamma_2=\infty$, then
\begin{equation}\label{eq:integrable-stable-remainder}
 \|\mathcal N(\tau)\|_{H^{-1}_\nu}
 \leq Ce^{-(\gamma_1+2\delta_0)\tau}.
\end{equation}
\end{lemma}

\begin{proof}
Apply $\Pi_-$ to the exact $H$ equation.  Since every global direct
column $Y_j$ lies in $\mathcal Z$, it disappears and gives
\eqref{eq:stable-normal-form}--\eqref{eq:stable-N-definition}.

To prove the tame bound, split the receding region at
$\bar f=e^{\tau/2}$.  On the inner part,
\eqref{eq:master-bootstrap-inner} and derivative recovery give
\[
 \|h\|_{C^1}
 \leq C e^{-\sigma\tau/2}.
\]
On the complement, the global $C^2$ and order-twelve bounds are
uniform, while the Gaussian density contributes
$e^{-ce^{\tau/2}}$.  The exact expansion \eqref{eq:Q-schematic},
one integration by parts, and Lemma~\ref{lem:first-moment} therefore
give
\[
 \|\rho_\tau\Q(h)\|_{H^{-1}_\nu}
 \leq
 C e^{-\kappa\tau}\|H\|_{H^1_\nu}
 +C\|H\|_{H^1_\nu}^2
 +Ce^{-ce^{\tau/2}}.
\]
The action terms are cubic because
$|c|\leq C\|H\|_{H^1_\nu}^2+Ce^{-ce^\tau}$; the cutoff, effective
column, and graft terms have the stated Gaussian tail.  This proves
\eqref{eq:stable-tame-bound}.

For later reference, the semigroup estimate used here is
\begin{equation}\label{eq:Hminus1-stable-semigroup}
 \|e^{u\A}\|_{H^{-1}_{\nu,-}\to L^2_{\nu,-}}
 \leq C(1+u^{-1/2})e^{-\gamma_1u}.
\end{equation}
Apply Lemma~\ref{lem:stable-rate-bootstrap} to
\eqref{eq:stable-tame-bound}, starting from the positive rate
\(\theta\) in \eqref{eq:master-bootstrap-H1}.  Its
finite iteration
\[
 r\longmapsto
 \min\{\gamma_1-\epsilon,\ r+\kappa,\ 2r\}
\]
reaches every exponent below \(\gamma_1\), and therefore proves
\eqref{eq:near-gap-H1}.  Choose
$\epsilon<\delta_0$.  The two terms on the right of
\eqref{eq:stable-tame-bound}, together with
\eqref{eq:profile-delta-choice}, then decay at least as
$e^{-(\gamma_1+2\delta_0)\tau}$, proving
\eqref{eq:integrable-stable-remainder}.
\end{proof}

\subsection{First profile}

\begin{theorem}[Marked, gauge-fixed first stable eigenspace profile]
\label{thm:first-stable-profile}
Fix $\delta_0$ as in \eqref{eq:profile-delta-choice}.  For every
strict prepared entrance there is a unique $V_\infty\in E_1$ such
that
\begin{equation}\label{eq:V-infinity-formula}
 V_\infty
 =\lim_{\tau\to\infty}e^{\gamma_1\tau}\Pi_1H(\tau)
 =e^{\gamma_1\tau_0}\Pi_1H(\tau_0)
  +\int_{\tau_0}^{\infty}
   e^{\gamma_1s}\Pi_1\mathcal N(s)\,ds.
\end{equation}
Here $\Pi_1\mathcal N$ is interpreted by
$H^{-1}_\nu$--$H^1_\nu$ duality on the finite-dimensional space
$E_1$.  Moreover,
\begin{align}
 \|H(\tau)-e^{-\gamma_1\tau}V_\infty\|_{H^1_\nu}
 &\leq Ce^{-(\gamma_1+\delta_0)\tau},
 \label{eq:first-profile-H1}\\
 \|h(\tau)-e^{-\gamma_1\tau}V_\infty\|_{C^m(K)}
 &\leq C_{K,m}e^{-(\gamma_1+\delta_0)\tau}
 \label{eq:first-profile-local}
\end{align}
for every $K\Subset M$, every $m\geq0$, and all sufficiently large
$\tau$.  On every fixed-$\tau_0$ common-margin prepared family,
$V_\infty$ depends continuously on the entrance data.  The profile is
defined in the fixed prepared gauge and marking, in the sense of
Remark~\ref{rem:stable-profile-gauge-scope}; in particular it is not
an unmarked geometric invariant.  A forward prepared restart at
absolute normalized time \(S\) leaves \(V_\infty\) unchanged when the
absolute \(\tau\)-coordinate is retained.  If the restarted coordinate
is instead reset to \(\widehat\tau=\tau-S\), this means a relabeling
of the \emph{same transported tail}: every time-typed object is
translated, for example
\[
 \widehat\rho_{\widehat\tau}=\rho_{\widehat\tau+S},\qquad
 \widehat\chi_{\widehat\tau}=\chi_{\widehat\tau+S},\qquad
 \widehat H(\widehat\tau)=H(\widehat\tau+S),
\]
and likewise for the grafts, effective columns, maps, and coefficient
systems.  The sliced state belongs to the translated convention
\(\widehat\Sigma_{\widehat\tau}^{[S]}
 :=\Sigma_{\widehat\tau+S}\), with all auxiliary data transported by
the same rule.  It is not thereby asserted to belong to an
independently defined standard slice at \(\widehat\tau\), nor to agree
with a fresh preparation of the restarted metric.  Put
\(\widehat\lambda(\widehat\tau)=\lambda(\widehat\tau+S)\).
Then the coefficient in
\(H=e^{-\gamma_1\widehat\tau}\widehat V_\infty+\cdots\) and the
restarted asymptotic scale satisfy
\begin{equation}\label{eq:profile-restart-rescaling}
 \widehat V_\infty=e^{-\gamma_1S}V_\infty .
\end{equation}
\begin{equation}\label{eq:scale-restart-rescaling}
 \widehat\lambda_\infty
 :=
 \lim_{\widehat\tau\to\infty}
  e^{\widehat\tau}\widehat\lambda(\widehat\tau)
 =
 e^{-S}\lambda_\infty .
\end{equation}
\end{theorem}

\begin{proof}
By \eqref{eq:integrable-stable-remainder},
$e^{\gamma_1\tau}\Pi_1\mathcal N(\tau)$ is integrable.  Projection of
\eqref{eq:stable-normal-form} onto $E_1$ proves
\eqref{eq:V-infinity-formula}.  Put
\[
 P_{\geq2}:=\Pi_--\Pi_1.
\]
Let
\[
 R(\tau)=H(\tau)-e^{-\gamma_1\tau}V_\infty.
\]
The \(E_1\) component of \(R\) is the backward tail
\[
 \Pi_1R(\tau)
 =
 -e^{-\gamma_1\tau}
 \int_\tau^\infty e^{\gamma_1s}\Pi_1\mathcal N(s)\,ds,
\]
and is \(O_{H^1_\nu}(e^{-(\gamma_1+2\delta_0)\tau})\), since \(E_1\)
is finite dimensional.  On the complementary stable subspace,
Duhamel's formula and
\eqref{eq:restricted-Hminus1-semigroup} with \(j=2\) give free rate
\(\gamma_2\) and forcing rate \(\gamma_1+2\delta_0\).
The strict inequalities in \eqref{eq:profile-delta-choice} therefore
give the \(L^2_\nu\) estimate
\[
 \|R(\tau)\|_{L^2_\nu}
 \leq Ce^{-(\gamma_1+2\delta_0)\tau}.
\]
The remainder satisfies
\(\partial_\tau R+LR=\mathcal N\).  Apply
 \eqref{eq:projected-stable-graph-recovery} to
 \(P_{\geq2}H=P_{\geq2}R\), with forcing
 \(P_{\geq2}\mathcal N\).  The \(E_1\) component already has the
 stronger \(H^1_\nu\) bound by its finite-dimensional backward-tail
 formula.  Here the lower-block source in
 \eqref{eq:projected-graph-lower-source} is
 \[
  g_R=
  \bigl(e^{-\kappa\tau}+\|H\|_{H^1_\nu}\bigr)
  \|\Pi_1H\|_{H^1_\nu}.
 \]
 Choose the loss in \eqref{eq:near-gap-H1} smaller than
 \(\delta_0\).  The inequalities in
 \eqref{eq:profile-delta-choice} then give
 \(g_R=O(e^{-(\gamma_1+2\delta_0)\tau})\).
 The preceding \(L^2_\nu\) bound and the concrete projected graph
 estimate therefore yield the
 \(H^1_\nu\) rate with \(2\delta_0\); in particular,
\eqref{eq:first-profile-H1} follows.  On a fixed compact set
\(\rho_\tau=1\) and the graft support is absent for all late times.
We include the mixed-rate pass needed to justify every local
derivative.  Put
\[
 r=h-e^{-\gamma_1\tau}V_\infty .
\]
Because \((\partial_\tau-\A)
(e^{-\gamma_1\tau}V_\infty)=0\), the exact local equation for \(r\)
is a uniformly parabolic equation whose source is the sum of:
\[
 \Q(e^{-\gamma_1\tau}V_\infty+r),\qquad
 \sum_{j=0}^8c_j\mathscr T_j
   (e^{-\gamma_1\tau}V_\infty+r),\qquad
 \sum_{j=0}^8c_j\mathcal Y_{j,\tau},
\]
together with coefficient differences multiplying \(r\).  The latter
are kept in the variable-coefficient operator.  On a fixed compact
set, the one-state compact smoothing estimate, the near-gap bound,
the feedback estimate, and the mode-growth estimates give the four
possible temporal rates
\begin{equation}\label{eq:first-profile-mixed-local-rates}
 \gamma_1+2\delta_0,\qquad
 2\gamma_1,\qquad
 \gamma_1+\kappa,\qquad
 \gamma_1+\theta .
\end{equation}
The first is the already proved \(H^1_\nu\) rate; the other three
come respectively from the quadratic term, the decaying local
coefficient defect, and the one-state compact remainder.  By
\eqref{eq:profile-delta-choice}, \(\kappa>4\delta_0\) and
\(\gamma_1>4\delta_0\), while
\(\theta>\sigma>2\kappa>8\delta_0\).  Hence every rate in
\eqref{eq:first-profile-mixed-local-rates} is strictly larger than
\(\gamma_1+\delta_0\).

Choose nested compact sets
\(K\Subset K_1\Subset\cdots\Subset K_{m+2}\).
The interior \(L^2\)-to-H\"older estimate on the last unit cylinder
first gives the \(C^0(K_{m+2})\) rate.  Differentiate the local
equation once, retain its top-order terms in the uniformly parabolic
operator, and apply the interior estimate on the next smaller
cylinder.  Repeating this finite procedure through order \(m\)
produces only products of already controlled lower derivatives, so
their rates are minima of the four numbers in
\eqref{eq:first-profile-mixed-local-rates}; there is no loss at each
step.  Compact all-order smoothing supplies the finitely many
coefficient derivatives required for the chosen \(m\).  This proves
\eqref{eq:first-profile-local} with rate
\(\gamma_1+\delta_0\), rather than merely an unspecified
compact-smoothing rate.

The profile restart law follows by writing the leading term in the
reset coordinate.  For the scale,
\[
 \lim_{\widehat\tau\to\infty}
 e^{\widehat\tau}\widehat\lambda(\widehat\tau)
 =
 e^{-S}\lim_{\tau\to\infty}e^\tau\lambda(\tau),
\]
which proves \eqref{eq:scale-restart-rescaling}.
Finally, \eqref{eq:V-infinity-formula} is the uniform limit of
finite-horizon continuous maps because its remaining integral has a
uniform exponentially decaying tail.  This proves continuity.
\end{proof}

\begin{definition}[Marked physical first amplitude]
\label{def:physical-first-amplitude}
For the fixed cutoffs, prepared marking, harmonic-map gauge, and
transported restart convention of
Remark~\ref{rem:stable-profile-gauge-scope}, define
\begin{equation}\label{eq:physical-first-amplitude}
 \mathfrak A_1
 :=
 \lambda_\infty^{-\gamma_1}V_\infty
 \in E_1 .
\end{equation}
This is a marked, gauge-fixed physical-time amplitude; it is not an
intrinsic class modulo arbitrary time-dependent re-markings.
\end{definition}

\begin{lemma}[Transported-restart covariance of the physical amplitude]
\label{lem:physical-amplitude-restart-covariance}
Under every transported forward restart in the convention of
Definition~\ref{def:physical-first-amplitude},
\begin{equation}\label{eq:physical-amplitude-restart-invariance}
 \widehat{\mathfrak A}_1
 =
 \widehat\lambda_\infty^{-\gamma_1}\widehat V_\infty
 =
 \mathfrak A_1 .
\end{equation}
Here ``transported'' includes the simultaneous relabeling of every
time-typed auxiliary object described in
Theorem~\ref{thm:first-stable-profile}.  The assertion compares two
descriptions of one prepared tail; it does not compare that tail with
an independently prepared state.
On every fixed-\(\tau_0\) common-margin prepared family,
\(\mathfrak A_1\) is continuous, and
\begin{equation}\label{eq:zero-physical-amplitude-equivalence}
 \mathfrak A_1=0
 \quad\Longleftrightarrow\quad
 V_\infty=0 .
\end{equation}
\end{lemma}

\begin{proof}
If the absolute normalized coordinate is retained, both factors in
\eqref{eq:physical-first-amplitude} are unchanged.  If it is reset,
\eqref{eq:profile-restart-rescaling} and
\eqref{eq:scale-restart-rescaling} give
\[
 (e^{-S}\lambda_\infty)^{-\gamma_1}
 e^{-\gamma_1S}V_\infty
 =
 \lambda_\infty^{-\gamma_1}V_\infty .
\]
Continuity follows from the one-state continuity of \(V_\infty\) and
\(\lambda_\infty>0\); the same positivity proves the zero-set
equivalence.
\end{proof}

\subsection{The quadratic response coefficient}

Let \(\Q_2\) be the formal symmetric quadratic Taylor coefficient of
\eqref{eq:Q-exact}.  For stable eigenmodes \(V,W\) and
\(0\leq\mu\leq8\), define its Gaussian modal moments by
\begin{equation}\label{eq:Q2-definition}
 \mathfrak q_\mu(V,W)
 :=\ip{\Q_2(V,W)}{Z_\mu}
 :=\frac12
 \left.\frac{\partial^2}{\partial s\,\partial r}\right|_{s=r=0}
 \ip{\Q(sV+rW)}{Z_\mu}.
\end{equation}
The last expression means the coefficient obtained from the formal
Taylor expansion of \eqref{eq:Q-exact}, followed by one integration by
parts.  Lemmas~\ref{lem:stable-eigenmode-growth} and
\ref{lem:mode-growth} make every resulting Gaussian integral
absolutely convergent; no global \(L^2_\nu\) assertion for
\(\Q_2(V,W)\) is needed.

\begin{lemma}[Quadratic modal form]
\label{lem:quadratic-modal-form}
For every \(0\leq\mu\leq8\), \(\mathfrak q_\mu\) extends uniquely
to a continuous symmetric bilinear form on \(H^1_{\nu,-}\), and
\begin{equation}\label{eq:quadratic-modal-bilinear-bound}
 |\mathfrak q_\mu(U,V)|
 \leq C\|U\|_{H^1_\nu}\|V\|_{H^1_\nu}.
\end{equation}
\end{lemma}

\begin{proof}
Insert the exact expression \eqref{eq:Q-exact} and integrate its
quasilinear second-derivative term once against
\(Z_\mu e^{-\bar f}\).  The resulting quadratic expression
contains only \(U,V,\bar\nabla U,\bar\nabla V\), and the fixed mode
coefficients.  The tensors \(Z_\mu\) and their first derivatives
have the bounds in Lemma~\ref{lem:mode-growth}.  When a derivative
hits the Gaussian density, the only unbounded coefficient is
\(\lvert\bar\nabla\bar f\rvert\), multiplying an undifferentiated
factor; Lemma~\ref{lem:first-moment} controls that factor in
\(L^2_\nu\).  Cauchy--Schwarz proves
\eqref{eq:quadratic-modal-bilinear-bound}.  Finite stable spectral sums
are dense in \(H^1_{\nu,-}\), so the extension is unique.
\end{proof}

Recall the fixed background matrix
\begin{equation}\label{eq:limiting-Gram-matrix}
 \mathbf G:=\mathbf G^{\rm bg},\qquad
 \mathbf G_{\mu j}=\ip{Y_j}{Z_\mu},
 \qquad0\leq\mu,j\leq8.
\end{equation}
It is invertible by \eqref{eq:background-Gram-reserve}.  For
\(V\in E_1\), put
\begin{equation}\label{eq:quadratic-feedback-coefficient}
 \mathbf c^{(2)}(V)
 =\bigl(a^{(2)}(V),b_1^{(2)}(V),\ldots,b_8^{(2)}(V)\bigr)
 :=-\mathbf G^{-1}
 \bigl(\mathfrak q_\mu(V,V)\bigr)_{\mu=0}^8.
\end{equation}
Thus \(\mathbf c^{(2)}:E_1\to\mathbb R^9\) is a continuous
finite-dimensional quadratic map, defined before it enters the
scattering package.

\subsection{Quadratic feedback, scale, and phase}

We continue the one-state branch with the quadratic modal form and
response coefficient just defined in
\eqref{eq:Q2-definition}--\eqref{eq:quadratic-feedback-coefficient};
none of the subsequent two-state, differentiability, or foliation
conclusions is an input here.

\begin{lemma}[Quadratic modal Taylor expansion]
\label{lem:quadratic-modal-Taylor}
Fix \(\mu\in\{0,\ldots,8\}\).  Let \(H=\rho_\tau h\), and put
\[
 \eta(\tau)=
 \|h(\tau)\|_{C^1(\{\bar f\leq e^{\tau/2}\})}.
\]
Whenever the smallness hypotheses of
Lemma~\ref{lem:localized-tame} hold,
\begin{equation}\label{eq:quadratic-modal-Taylor-remainder}
 \left|
  \ip{\rho_\tau\Q(h)}{Z_\mu}
  -\mathfrak q_\mu(H,H)
 \right|
 \leq
 C\eta(\tau)\|H\|_{H^1_\nu}^2
 +Ce^{-ce^{\tau/2}}.
\end{equation}
Consequently, suppose that \(V\in E_1\), \(\delta>0\),
\(\omega>0\), and
\[
 H(\tau)=e^{-\gamma_1\tau}V+R(\tau),\qquad
 \|R(\tau)\|_{H^1_\nu}
 \leq Ce^{-(\gamma_1+\delta)\tau},
\]
while \(\eta(\tau)\leq Ce^{-\omega\tau}\).  Then
\begin{equation}\label{eq:quadratic-modal-profile-expansion}
 \ip{\rho_\tau\Q(h)}{Z_\mu}
 =
 e^{-2\gamma_1\tau}\mathfrak q_\mu(V,V)
 +O\!\left(
 e^{-(2\gamma_1+\min\{\delta,\omega\})\tau}\right)
 +O(e^{-ce^{\tau/2}}).
\end{equation}
\end{lemma}

\begin{proof}
For the Taylor remainder, split the integral at
\(\bar f=e^{\tau/2}\).  On the inner region \(\rho_\tau=1\), and the
inverse-metric expansion in \eqref{eq:Q-exact} shows that every term
beyond \(\Q_2(h,h)\) gains a factor \(C\eta(\tau)\) after the same
integration by parts.  The portion where \(H\neq h\), as well as the
outer part of the integral, is supported in
\(\{\bar f\geq e^{\tau/2}\}\); the global \(C^2\) bound, polynomial
mode growth, and the Gaussian density make it
\(O(e^{-ce^{\tau/2}})\).  This proves
\eqref{eq:quadratic-modal-Taylor-remainder}.

Finally, Lemma~\ref{lem:quadratic-modal-form} and
\eqref{eq:quadratic-modal-bilinear-bound} give
\[
 \mathfrak q_\mu(H,H)
 =
 e^{-2\gamma_1\tau}\mathfrak q_\mu(V,V)
 +O(e^{-(2\gamma_1+\delta)\tau}).
\]
Combining this with
\eqref{eq:quadratic-modal-Taylor-remainder} proves
\eqref{eq:quadratic-modal-profile-expansion}.
\end{proof}

\begin{lemma}[Quadratic modal scattering]
\label{lem:quadratic-modal-scattering}
For the exact Gram system
\[
 M(\tau,h)c=-d(\tau,h),\qquad
 c=(a,b_1,\ldots,b_8),\qquad
 d=(d_0,\ldots,d_8),
\]
with \(d_\mu\) defined in \eqref{eq:receding-d}, and for the profile
$V_\infty$ in Theorem~\ref{thm:first-stable-profile},
\begin{align}
 M(\tau,h)
 &=\mathbf G+O(e^{-\gamma_1\tau}),
 \label{eq:modal-M-expansion}\\
 d_\mu(\tau,h)
 &=e^{-2\gamma_1\tau}
   \mathfrak q_\mu(V_\infty,V_\infty)
   +O(e^{-(2\gamma_1+\delta_0)\tau}).
 \label{eq:modal-d-expansion}
\end{align}
\end{lemma}

\begin{proof}
The column-tail estimate and
\eqref{eq:first-profile-H1}, inserted in
\eqref{eq:receding-M0}--\eqref{eq:receding-Mj}, prove
\eqref{eq:modal-M-expansion}.  For the right-hand side, the graft,
cutoff-commutator, and effective-column modal errors are
Gaussian-superexponential.  Apply
Lemma~\ref{lem:quadratic-modal-Taylor} with
\(V=V_\infty\), \(\delta=\delta_0\), and
\(\omega=\sigma/2\), using
\eqref{eq:first-profile-H1} and
\eqref{eq:master-bootstrap-inner}.  Since
$\delta_0<\kappa/4<\sigma/8$, these errors lie below
$e^{-(2\gamma_1+\delta_0)\tau}$.  This proves
\eqref{eq:modal-d-expansion}.
\end{proof}

\begin{theorem}[Quadratic geometric asymptotics]
\label{thm:quadratic-geometric-asymptotics}
For $V_\infty$ in Theorem~\ref{thm:first-stable-profile}, let
$\mathbf c^{(2)}=\mathbf c^{(2)}(V_\infty)$.  When an argument is
displayed, write \(a^{(2)}(V)\), \(b_j^{(2)}(V)\), and
\[
 U^{(2)}(V):=\sum_{j=1}^8b_j^{(2)}(V)W_j
\]
for the components of \(\mathbf c^{(2)}(V)\).  The quadratic
homogeneity of \(\mathbf c^{(2)}\) and
Definition~\ref{def:physical-first-amplitude} give
\begin{equation}\label{eq:physical-amplitude-quadratic-response}
 \mathbf c^{(2)}(\mathfrak A_1)
 =
 \lambda_\infty^{-2\gamma_1}
 \mathbf c^{(2)}(V_\infty).
\end{equation}
Then
\begin{equation}\label{eq:quadratic-feedback-asymptotic}
 (a,b)(\tau)
 =e^{-2\gamma_1\tau}\mathbf c^{(2)}
  +O(e^{-(2\gamma_1+\delta_0)\tau}).
\end{equation}
Consequently,
\begin{align}
 \lambda(\tau)e^\tau
 &=
 \lambda_\infty\left[
  1+\frac{a^{(2)}}{2\gamma_1}e^{-2\gamma_1\tau}
  +O(e^{-(2\gamma_1+\delta_0)\tau})
  \right],
 \label{eq:quadratic-scale-normalized}\\
 \frac{\lambda(\tau)}{T-t(\tau)}
 &=
 1+\frac{a^{(2)}}{1+2\gamma_1}e^{-2\gamma_1\tau}
  +O(e^{-(2\gamma_1+\delta_0)\tau}),
 \label{eq:quadratic-scale-physical-ratio}\\
 \lambda(t)
 &=(T-t)\left[
  1+\frac{a^{(2)}}{1+2\gamma_1}
       \lambda_\infty^{-2\gamma_1}(T-t)^{2\gamma_1}
  +o((T-t)^{2\gamma_1})
 \right].
 \label{eq:quadratic-scale-physical-time}
\end{align}
If
\(U^{(2)}=U^{(2)}(V_\infty)\), then on every compact set,
in a fixed exponential chart at the identity, write
\(\log_{\operatorname{Id}}\) for the inverse chart.  Then
\begin{equation}\label{eq:quadratic-diffeomorphism-phase}
 \log_{\operatorname{Id}}\!
 \left(\Psi_\infty\circ\Psi_\tau^{-1}\right)
 =
  \frac{e^{-2\gamma_1\tau}}{2\gamma_1}U^{(2)}
  +O_{C^m}(e^{-(2\gamma_1+\delta_0)\tau}).
\end{equation}
Equivalently, writing \(\delta=T-t\), the physical-time scale and phase
laws are
\begin{equation}\label{eq:quadratic-scale-physical-amplitude}
 \frac{\lambda(t)}{T-t}
 =
 1+
 \frac{a^{(2)}(\mathfrak A_1)}{1+2\gamma_1}
 (T-t)^{2\gamma_1}
 +o((T-t)^{2\gamma_1})
\end{equation}
and, for every \(K\Subset M\) and \(m\geq0\),
\begin{equation}\label{eq:quadratic-phase-physical-amplitude}
 \log_{\operatorname{Id}}\!
 \left(\Psi_\infty\circ\Psi_{\tau(t)}^{-1}\right)
 =
 \frac{(T-t)^{2\gamma_1}}{2\gamma_1}
   U^{(2)}(\mathfrak A_1)
 +o_{C^m(K)}((T-t)^{2\gamma_1}).
\end{equation}
The marked physical amplitude and both quadratic response maps are
continuous on common-margin prepared families.  Any or all of the
displayed order-\(2\gamma_1\) coefficients may vanish; no assertion
that one of them is the first nonzero correction is made.
\end{theorem}

\begin{proof}
Lemma~\ref{lem:quadratic-modal-scattering} and uniform inversion of
the exact Gram system prove
\eqref{eq:quadratic-feedback-asymptotic}.

Since
$\frac d{d\tau}\log(\lambda e^\tau)=-a$, integration from $\tau$ to
infinity gives, with the sign displayed explicitly,
\begin{equation}\label{eq:quadratic-scale-log-check}
 \log\frac{\lambda(\tau)e^\tau}{\lambda_\infty}
 =\int_\tau^\infty a(s)\,ds
 =\frac{a^{(2)}}{2\gamma_1}e^{-2\gamma_1\tau}
 {}+O(e^{-(2\gamma_1+\delta_0)\tau}),
\end{equation}
which proves \eqref{eq:quadratic-scale-normalized}.  For the physical
ratio use the exact identity
\[
 \frac{T-t(\tau)}{\lambda(\tau)}
 =\int_0^\infty
  \exp\left(-u-\int_\tau^{\tau+u}a(s)\,ds\right)du.
\]
The elementary integral
\[
 \int_0^\infty e^{-u}
 \frac{1-e^{-2\gamma_1u}}{2\gamma_1}\,du
 =\frac1{1+2\gamma_1}
\]
and expansion of the exponential give
\begin{equation}\label{eq:quadratic-physical-ratio-sign-check}
 \frac{T-t(\tau)}{\lambda(\tau)}
 =1-\frac{a^{(2)}}{1+2\gamma_1}e^{-2\gamma_1\tau}
  +O(e^{-(2\gamma_1+\delta_0)\tau}).
\end{equation}
Inverting this identity gives
\eqref{eq:quadratic-scale-physical-ratio}; substituting
$e^{-\tau}=\lambda_\infty^{-1}(T-t)(1+o(1))$ gives
\eqref{eq:quadratic-scale-physical-time}.  The same substitution,
the quadratic homogeneity
\eqref{eq:physical-amplitude-quadratic-response}, and division by
\(T-t\) give
\eqref{eq:quadratic-scale-physical-amplitude}.

Finally, $\chi_\tau=1$ on every fixed compact set for all sufficiently
large $\tau$, and
\[
 \int_\tau^\infty b_j(s)\,ds
 =\frac{b_j^{(2)}}{2\gamma_1}e^{-2\gamma_1\tau}
  +O(e^{-(2\gamma_1+\delta_0)\tau}).
\]
Since
\(\partial_\tau\Psi_\tau=(\sum b_j\chi_\tau W_j)\circ\Psi_\tau\),
the relative tail is \(\Psi_\infty\circ\Psi_\tau^{-1}\), so its
first-order sign is positive.  Integrating this time-ordered flow gives
\eqref{eq:quadratic-diffeomorphism-phase}; flow-composition errors are
$O(e^{-4\gamma_1\tau})$ and hence lie below the stated remainder.
Substituting
\(e^{-\tau}=\lambda_\infty^{-1}(T-t)(1+o(1))\) once more and using
\eqref{eq:physical-amplitude-quadratic-response} proves
\eqref{eq:quadratic-phase-physical-amplitude}.
Continuity of \(V_\infty\) on fixed common-margin prepared families is
already part of Theorem~\ref{thm:first-stable-profile}.  The conclusion
therefore follows from that one-state continuity, continuity of
\(\lambda_\infty>0\), and the
finite-dimensional quadratic formula
\eqref{eq:quadratic-feedback-coefficient}; it does not use the
separate two-state \(C^1\) theorem proved subsequently in this section.
\end{proof}

\subsection{Sharp marked spacetime profile}

Let
\begin{equation}\label{eq:linearized-DeTurck-vector}
 \mathfrak b_{\bar g}(V)
 :=D[B_{\bar g}]_{\bar g}[V]
\end{equation}
be the linearization of the DeTurck vector field at $\bar g$.
For $V\in E_1$, let $\mathcal W_V(u)$ be the solution of the linear
transport equation
\begin{equation}\label{eq:Jacobi-transition-transport}
 \partial_u\mathcal W_V+
 [\bar\nabla\bar f,\mathcal W_V]
 =
 -e^{-\gamma_1u}\mathfrak b_{\bar g}(V),
 \qquad
 \mathcal W_V(0)=0.
\end{equation}
Since the soliton field is complete, this equation has a global smooth
solution for every finite $u$; explicitly,
\[
 \varphi_u^*\mathcal W_V(u)
 =
 -\int_0^u
  e^{-\gamma_1r}
  \varphi_r^*\mathfrak b_{\bar g}(V)\,dr.
\]
This formulation avoids imposing any global completeness assertion on
a nonlinear perturbation of the soliton field.  For $s<0$, put
$u=-\log(-s)$ and define
\begin{equation}\label{eq:marked-Jacobi-field}
 \mathcal J_V(s)
 :=
 (-s)\varphi_u^*
 \left(
  e^{-\gamma_1u}V+
  \Lie_{\mathcal W_V(u)}\bar g
 \right).
\end{equation}
This is linear in $V$, satisfies $\mathcal J_V(-1)=V$, and is the
base-time-marked Ricci-flow linearization of the normalized
Ricci--DeTurck eigenmode $e^{-\gamma_1u}V$ along
$g_{\mathrm{FIK}}(s)$.

\begin{lemma}[Marked Jacobi equation]
\label{lem:marked-Jacobi-equation}
Let
\[
 \mathscr R(g):=-2\Ric_g .
\]
For every \(V\in E_1\), the tensor in
\eqref{eq:marked-Jacobi-field} is a smooth solution of
\begin{equation}\label{eq:marked-Jacobi-equation}
 \partial_s\mathcal J_V(s)
 =
 D\mathscr R_{g_{\mathrm{FIK}}(s)}
     [\mathcal J_V(s)],
 \qquad
 \mathcal J_V(-1)=V .
\end{equation}
\end{lemma}

\begin{proof}
Fix a compact \(u\)-interval and a compact spatial flow tube.  For a
real parameter \(\varepsilon\), set
\[
 g_\varepsilon(u)
 =\bar g+\varepsilon e^{-\gamma_1u}V,
 \qquad
 X_\varepsilon(u)
 =\bar\nabla\bar f
  -\varepsilon e^{-\gamma_1u}\mathfrak b_{\bar g}(V),
\]
and let \(D_\varepsilon(u)\) be the local flow of \(X_\varepsilon(u)\)
with \(D_\varepsilon(0)=\operatorname{Id}\).  Recall the normalized
Ricci--DeTurck operator
\[
 \mathcal R_{\bar g}(g)
 =
 -2\Ric_g+g+\Lie_{B_{\bar g}(g)}g
 -\Lie_{\bar\nabla\bar f}g .
\]
Its derivative at \(\bar g\) is \(\A\).  Since
\(\A V=-\gamma_1V\), Taylor expansion on the fixed flow tube gives
\begin{equation}\label{eq:Jacobi-normalized-residual}
 \partial_u g_\varepsilon-\mathcal R_{\bar g}(g_\varepsilon)
 =O_{C^m}(\varepsilon^2)
\end{equation}
for every \(m\).  Moreover,
\[
 X_\varepsilon
 =
 \bar\nabla\bar f-B_{\bar g}(g_\varepsilon)
 +O_{C^m}(\varepsilon^2).
\]
Consequently the normalization-and-pullback calculation, applied to
\[
 G_\varepsilon(s)
 :=
 (-s)D_\varepsilon(-\log(-s))^*
       g_\varepsilon(-\log(-s)),
\]
turns \eqref{eq:Jacobi-normalized-residual} into
\begin{equation}\label{eq:Jacobi-physical-residual}
 \partial_sG_\varepsilon-\mathscr R(G_\varepsilon)
 =O_{C^m}(\varepsilon^2)
\end{equation}
on every compact subinterval of \((-\infty,0)\).

At \(\varepsilon=0\), \(D_0(u)=\varphi_u\) and
\(G_0=g_{\mathrm{FIK}}\).  Differentiating the flow equation at
\(\varepsilon=0\) shows that the relative first-variation field is
\(\mathcal W_V\), because it satisfies
\eqref{eq:Jacobi-transition-transport}.  Hence
\[
 \left.\partial_\varepsilon\right|_{\varepsilon=0}G_\varepsilon(s)
 =
 (-s)\varphi_u^*
 \left(e^{-\gamma_1u}V
       +\Lie_{\mathcal W_V(u)}\bar g\right)
 =\mathcal J_V(s).
\]
Differentiating \eqref{eq:Jacobi-physical-residual} at
\(\varepsilon=0\) proves the evolution equation.  At \(s=-1\) one has
\(u=0\), \(D_\varepsilon(0)=\operatorname{Id}\), and
\(\mathcal W_V(0)=0\), which gives the initial value.
\end{proof}

\begin{lemma}[Frozen-window ODE first variation]
\label{lem:frozen-window-first-variation}
Fix \(m\in\mathbb N_0\), and let \(J\Subset\mathbb R\) be a compact
interval containing \(0\).
Fix compact coordinate domains
\(\mathcal K\Subset\mathcal K^+\) such that, for all sufficiently large
\(i\), every \(D_i\)- and \(\varphi\)-trajectory issuing from
\(\mathcal K\), and every comparison trajectory used below, remains in
\(\mathcal K^+\) for \(u\in J\).
All tensors \(T_i,T_0,T_1,S_i\) below are sections of one fixed
natural tensor bundle of a single type over \(\mathcal K^+\); in the
application this bundle is \(S^2T^*\mathcal K^+\).  Pullback and Lie
derivative are those of this natural bundle, and every \(C^m\) norm is
computed using one fixed background connection, equivalently in a
fixed finite coordinate trivialization on \(\mathcal K^+\).
Let \(\epsilon_i,r_i>0\) be sequences satisfying
\[
 \epsilon_i\downarrow0,\qquad
 r_i=o(\epsilon_i),\qquad
 \epsilon_i^2=O(r_i).
\]
Let
\[
 X_i(u)=X_0+\epsilon_iX_1(u)+R_i(u),
\]
where
\[
 X_0\in C_x^{m+2}(\mathcal K^+),\qquad
 X_1\in C_u^0(J;C_x^{m+2}(\mathcal K^+)),
\]
these norms are uniformly bounded, and
\[
 R_i\in C_u^0(J;C_x^{m+2}(\mathcal K^+)),\qquad
 \sup_{u\in J}\|R_i(u)\|_{C_x^{m+2}(\mathcal K^+)}
 \leq Cr_i.
\]
Let \(D_i(u)\) be the flow of \(X_i(u)\) with
\(D_i(0)=\operatorname{Id}\).
Write \(\varphi_u\) for the flow of \(X_0\).  If
\[
 T_i(u)=T_0(u)+\epsilon_iT_1(u)+S_i(u),
\]
where
\[
 T_0\in C_u^0(J;C_x^{m+2}(\mathcal K^+)),\qquad
 T_1\in C_u^0(J;C_x^{m+1}(\mathcal K^+)),
\]
\[
 \sup_{u\in J}\left(
  \|T_0(u)\|_{C_x^{m+2}(\mathcal K^+)}
  +\|T_1(u)\|_{C_x^{m+1}(\mathcal K^+)}
 \right)\leq C
\]
and
\[
 S_i\in C_u^0(J;C_x^m(\mathcal K^+)),\qquad
 \sup_{u\in J}\|S_i(u)\|_{C_x^m(\mathcal K^+)}
 \leq Cr_i,
\]
then
\begin{equation}\label{eq:frozen-window-first-variation}
 D_i(u)^*T_i(u)
 =
 \varphi_u^*\!\left[
  T_0(u)+\epsilon_i
  \bigl(T_1(u)+\Lie_{W(u)}T_0(u)\bigr)
 \right]
 +O_{C_u^0(J;C_x^m(\mathcal K))}(r_i),
\end{equation}
where \(W(0)=0\) and
\begin{equation}\label{eq:frozen-window-variation-transport}
 \partial_uW+[X_0,W]=X_1(u).
\end{equation}
The estimate is uniform for families satisfying the displayed
flow-tube and coefficient bounds.
\end{lemma}

\begin{proof}
Set \(E_i(u)=D_i(u)\circ\varphi_u^{-1}\).  Differentiating the flow
equations in local coordinates gives
\[
 \partial_uE_i
 =X_i(u)\circ E_i-DE_i\,X_0 .
\]
After substituting \(E_i=\operatorname{Id}+\epsilon_iW+F_i\), the
coefficient of \(\epsilon_i\) is precisely
\eqref{eq:frozen-window-variation-transport}.  Taylor's theorem for
the right-hand side and Gronwall on \(J\), differentiated spatially up
to order \(m+1\), give
\[
 E_i(u)=\operatorname{Id}+\epsilon_iW(u)
       +O_{C_u^0(J;C_x^{m+1}(\varphi_u(\mathcal K)))}
        (\epsilon_i^2+r_i)
 =\operatorname{Id}+\epsilon_iW(u)
       +O_{C_u^0(J;C_x^{m+1}(\varphi_u(\mathcal K)))}(r_i).
\]
The uniform \(C^{m+2}\) bound for \(T_0\), the \(C^{m+1}\) bound for
\(T_1\), and \(\epsilon_i^2=O(r_i)\) now give
\[
 E_i^*T_i
 =
 T_0+\epsilon_i\bigl(T_1+\Lie_WT_0\bigr)
 +O_{C_u^0(J;C_x^m)}(r_i).
\]
Since \(D_i=E_i\circ\varphi_u\), pulling this identity back by
\(\varphi_u\) proves \eqref{eq:frozen-window-first-variation}.
\end{proof}

\begin{theorem}[Sharp marked spacetime asymptotics]
\label{thm:sharp-marked-spacetime}
Let a strict prepared entrance generate the flow in
Theorem~\ref{thm:prepared-entrance-continuation}, and let $V_\infty$ and
$\delta_0$ be as in
Theorem~\ref{thm:first-stable-profile}.  Put
\[
 \widehat\delta=\min\{\delta_0,\gamma_1\}.
\]
For every sequence $t_i\uparrow T$, set
\[
 \delta_i=T-t_i,\qquad \tau_i=\tau(t_i),\qquad
 \Xi_i=\Xi_{t_i}.
\]
For a finite-order entrance, every norm below is understood as a norm
of the smooth pulled-back tensor furnished by
Lemma~\ref{lem:finite-order-relative-marking-cancellation}; no
\(C^\infty\)-regularity of the absolute marking \(\Xi_i\) is asserted.
Then, for every $K\Subset M$, $I\Subset(-\infty,0)$, and $m\geq0$,
\begin{equation}\label{eq:sharp-marked-spacetime-expansion}
 \left\|
  \delta_i^{-1}\Xi_i^*G(T+s\delta_i)
  -g_{\mathrm{FIK}}(s)
  -e^{-\gamma_1\tau_i}\mathcal J_{V_\infty}(s)
 \right\|_{C^m(K\times I)}
 \leq C_{K,I,m}
 e^{-(\gamma_1+\widehat\delta)\tau_i}.
\end{equation}
Equivalently, in terms of the marked physical amplitude from
Definition~\ref{def:physical-first-amplitude},
\begin{equation}\label{eq:sharp-marked-physical-amplitude-expansion}
 \left\|
  \delta_i^{-1}\Xi_i^*G(T+s\delta_i)
  -g_{\mathrm{FIK}}(s)
  -\delta_i^{\gamma_1}\mathcal J_{\mathfrak A_1}(s)
 \right\|_{C^m(K\times I)}
 \leq C_{K,I,m}
 \delta_i^{\gamma_1+\widehat\delta}.
\end{equation}
In particular,
\begin{equation}\label{eq:sharp-marked-spacetime-rate}
 \left\|
  \delta_i^{-1}\Xi_i^*G(T+s\delta_i)
  -g_{\mathrm{FIK}}(s)
 \right\|_{C^m(K\times I)}
 \leq C_{K,I,m}\delta_i^{\gamma_1},
\end{equation}
and, without subsequence extraction,
\begin{equation}\label{eq:sharp-marked-Jacobi-limit}
 \delta_i^{-\gamma_1}
 \left[
 \delta_i^{-1}\Xi_i^*G(T+s\delta_i)
  -g_{\mathrm{FIK}}(s)\right]
 \longrightarrow
 \mathcal J_{\mathfrak A_1}(s)
 =
 \lambda_\infty^{-\gamma_1}\mathcal J_{V_\infty}(s)
\end{equation}
in $C^\infty_{\mathrm{loc}}(M\times(-\infty,0))$.  If
\(\mathfrak A_1=0\), equivalently \(V_\infty=0\), the rate in
\eqref{eq:sharp-marked-spacetime-rate} improves to
$O(\delta_i^{\gamma_1+\widehat\delta})$.
If \(\mathfrak A_1\neq0\), equivalently \(V_\infty\neq0\), then on
some fixed compact spacetime window
containing \(s=-1\), the marked error is not
\(o(\delta_i^{\gamma_1})\) in \(C^0\).
These are statements for the displayed frozen base-time markings.  In
particular, when \(V_\infty\) is a stable Lie-derivative direction,
the theorem asserts a nonzero marked Jacobi field and makes no claim
that its class is nonzero modulo arbitrary re-markings.
\end{theorem}

\begin{proof}
For $s\in I$, define $\tau_i(s)$ by
$t(\tau_i(s))=T+s\delta_i$ and put
\[
 u_i(s)=\tau_i(s)-\tau_i,\qquad u(s)=-\log(-s).
\]
The quadratic physical-scale expansion
\eqref{eq:quadratic-scale-physical-ratio} gives, uniformly on $I$,
\begin{align}
 u_i(s)&=u(s)+O_I(e^{-2\gamma_1\tau_i}),
 \label{eq:sharp-time-shift}\\
 \frac{\lambda(\tau_i(s))}{\delta_i}
 &=(-s)\bigl(1+O_I(e^{-2\gamma_1\tau_i})\bigr).
 \label{eq:sharp-scale-ratio}
\end{align}
  On every compact set containing the bounded
  $\bar\nabla\bar f$-trajectories from $K$,
  Theorem~\ref{thm:first-stable-profile} and
  \eqref{eq:sharp-time-shift} give
\begin{equation}\label{eq:sharp-profile-window}
 h(\tau_i(s))
 =
 e^{-\gamma_1\tau_i}e^{-\gamma_1u(s)}V_\infty
 +
 O_{C^{m+3}}\left(
  e^{-(\gamma_1+\delta_0)\tau_i}\right).
\end{equation}
Smooth Taylor expansion of the DeTurck vector field, together with
\eqref{eq:quadratic-feedback-asymptotic}, yields
\begin{equation}\label{eq:sharp-Phi-generator}
 (\partial_\tau\Phi_\tau)\circ\Phi_\tau^{-1}
 =
 \bar\nabla\bar f
 -e^{-\gamma_1\tau_i}e^{-\gamma_1u}
   \mathfrak b_{\bar g}(V_\infty)
 +O_{C^{m+2}}\left(
   e^{-(\gamma_1+\widehat\delta)\tau_i}\right)
\end{equation}
when $\tau=\tau_i+u$ and $u$ ranges in the corresponding bounded
interval.  Here the quadratic Taylor remainder, $a\bar\nabla\bar f$,
and $U$ are $O(e^{-2\gamma_1\tau_i})$.

For \(u\) in a fixed compact interval containing all \(u_i(I)\), put
\[
 \widetilde D_i(u)
 =\Phi_{\tau_i+u}\circ\Phi_{\tau_i}^{-1}.
\]
Thus
\[
 D_{i,s}=\widetilde D_i(u_i(s))
\]
for the transition map in
\eqref{eq:master-marking-transition}.  Apply
Lemma~\ref{lem:frozen-window-first-variation} to
\(\widetilde D_i\) with
\[
 \epsilon_i=e^{-\gamma_1\tau_i},\qquad
 r_i=e^{-(\gamma_1+\widehat\delta)\tau_i},
\]
\[
 X_0=\bar\nabla\bar f,\qquad
 X_1(u)=-e^{-\gamma_1u}
          \mathfrak b_{\bar g}(V_\infty),
\]
and
\[
 T_0=\bar g,\qquad
 T_1(u)=e^{-\gamma_1u}V_\infty.
\]
The hypotheses follow from
\eqref{eq:sharp-profile-window}--\eqref{eq:sharp-Phi-generator};
the relation \(\widehat\delta\leq\gamma_1\) gives
\(\epsilon_i^2=O(r_i)\).  Equation~
\eqref{eq:frozen-window-variation-transport} is precisely
\eqref{eq:Jacobi-transition-transport}.  Evaluate the lemma at
\(u_i(s)\).  By \eqref{eq:sharp-time-shift}, replacing
\(u_i(s)\) by \(u(s)\) in the background flow, profile, and first
variation costs \(O(e^{-2\gamma_1\tau_i})=O(r_i)\).  Hence
\begin{align}
 &D_{i,s}^*
   \bigl(\bar g+h(\tau_i(s))\bigr)
 \notag\\
 &\quad=
 \varphi_{u(s)}^*
 \left[
  \bar g+
  e^{-\gamma_1\tau_i}
  \left(
   e^{-\gamma_1u(s)}V_\infty+
   \Lie_{\mathcal W_{V_\infty}(u(s))}\bar g
  \right)
 \right]
 +O_{C^0(I;C^m(K))}\left(
   e^{-(\gamma_1+\widehat\delta)\tau_i}\right).
 \label{eq:sharp-transition-expansion}
\end{align}
Insert \eqref{eq:sharp-scale-ratio} and
\eqref{eq:sharp-transition-expansion} into the exact identity
  \eqref{eq:master-frozen-pullback}.  The definition in
  \eqref{eq:marked-Jacobi-field} proves
\eqref{eq:sharp-marked-spacetime-expansion} for spatial derivatives.
To obtain the asserted mixed norms without an implicit compatibility
claim, write
\[
 G_i(s)=\delta_i^{-1}\Xi_i^*G(T+s\delta_i),\qquad
 E_i(s)=G_i(s)-g_{\mathrm{FIK}}(s)
        -e^{-\gamma_1\tau_i}\mathcal J_{V_\infty}(s).
\]
Both \(G_i\) and \(g_{\mathrm{FIK}}\) solve Ricci flow in this fixed
marking, while Lemma~\ref{lem:marked-Jacobi-equation} gives the
linearized equation for \(\mathcal J_{V_\infty}\).  Taylor expansion of
\(\mathscr R(g)=-2\Ric_g\) on a slightly larger compact spatial set
therefore gives
\[
 \partial_sE_i
 =
 D\mathscr R_{g_{\mathrm{FIK}}}[E_i]
 +O_{C^m}\!\left(
   e^{-2\gamma_1\tau_i}
   +e^{-\gamma_1\tau_i}\|E_i\|_{C^{m+2}}
   +\|E_i\|_{C^{m+2}}^2\right).
\]
The spatial estimate already proved is available at every derivative
order, and
\(e^{-2\gamma_1\tau_i}
 =O(e^{-(\gamma_1+\widehat\delta)\tau_i})\).
It follows first that
\(\partial_sE_i
=O_{C^m}(e^{-(\gamma_1+\widehat\delta)\tau_i})\).
Differentiating this equation in \(s\) and inducting, using higher
spatial estimates at each step, gives the same bound for every mixed
space--time derivative.  This proves the full
\(C^m(K\times I)\) statement.

Equations~\eqref{eq:quadratic-scale-normalized} and
\eqref{eq:quadratic-scale-physical-ratio} give the quantitative clock
conversion
\begin{equation}\label{eq:physical-amplitude-clock-conversion}
 \delta_i e^{\tau_i}
 =
 \lambda_\infty
 \bigl(1+O(e^{-2\gamma_1\tau_i})\bigr).
\end{equation}
Since \(\mathcal J_V\) is linear in \(V\),
\eqref{eq:physical-first-amplitude} and
\eqref{eq:physical-amplitude-clock-conversion} imply
\[
 e^{-\gamma_1\tau_i}\mathcal J_{V_\infty}
 =
 \delta_i^{\gamma_1}\mathcal J_{\mathfrak A_1}
 +O_{C^m(K\times I)}(\delta_i^{3\gamma_1}).
\]
Because \(\widehat\delta\leq\gamma_1\), insertion in
\eqref{eq:sharp-marked-spacetime-expansion} proves
\eqref{eq:sharp-marked-physical-amplitude-expansion}.
It also proves \eqref{eq:sharp-marked-spacetime-rate} and
\eqref{eq:sharp-marked-Jacobi-limit}.  If $V_\infty=0$, the
linear term in \eqref{eq:sharp-marked-spacetime-expansion} vanishes,
which gives the stated improved rate.
 If \(V_\infty\neq0\), then
 \(\mathcal J_{V_\infty}(-1)=V_\infty\neq0\).  Choose a compact set on
 which this tensor has nonzero \(C^0\)-norm.  The nonzero limit in
 \eqref{eq:sharp-marked-Jacobi-limit} then shows directly that the
 marked error is not \(o(\delta_i^{\gamma_1})\).  Thus both alternatives
 in the sharpness assertion are consequences of the displayed
 expansion, rather than only of its upper bound.
\end{proof}

\begin{corollary}[Transported-restart-invariant marked physical normal
form]
\label{cor:physical-amplitude-normal-form}
For every strict prepared entrance, the marked amplitude
\(\mathfrak A_1\in E_1\) is invariant under transported forward
restarts in the precise sense of
Lemma~\ref{lem:physical-amplitude-restart-covariance}.  It determines
the linear marked Jacobi term in
\eqref{eq:sharp-marked-physical-amplitude-expansion}, while its
quadratic self-interaction determines the order-\(2\gamma_1\) scale
and phase responses in
\eqref{eq:quadratic-scale-physical-amplitude} and
\eqref{eq:quadratic-phase-physical-amplitude}.  These statements are
continuous on fixed common-margin prepared families.  Moreover,
\(\mathfrak A_1=0\) holds exactly on the improved marked
strong-stable class, whereas \(\mathfrak A_1\neq0\) gives exact marked
first order \((T-t)^{\gamma_1}\) on some compact spacetime window.
All assertions retain the fixed marking and gauge of
Definition~\ref{def:physical-first-amplitude}; no quotient by arbitrary
re-markings is asserted, and an order-\(2\gamma_1\) response coefficient
is allowed to vanish.  For a finite-order entrance, the marked Jacobi
and exact-first-order statements use the relative-pullback
interpretation of
Lemma~\ref{lem:finite-order-relative-marking-cancellation}.  They
assert all-order asymptotics of the pulled-back tensor on late compact
spacetime windows, not all-order regularity of the absolute marking or
of the ODE-carried prepared maps.  On the smooth entrance subclass
this is the usual smooth frozen-marking statement.
\end{corollary}

\begin{proof}
Combine
Lemma~\ref{lem:physical-amplitude-restart-covariance},
Theorems~\ref{thm:quadratic-geometric-asymptotics} and
\ref{thm:sharp-marked-spacetime}, and
\eqref{eq:zero-physical-amplitude-equivalence}.
\end{proof}

\subsection{Transverse prepared disks and late-entry variations}

\begin{lemma}[Compactly supported transverse lift]
\label{lem:compact-profile-lift}
There is a linear map
\begin{equation}\label{eq:compact-profile-lift}
 W:E_1\longrightarrow C_c^\infty(S^2T^*M)
\end{equation}
with one common compact support and
\begin{equation}\label{eq:compact-profile-lift-properties}
 \Pi_1W(v)=v,\qquad W(v)\perp\mathcal Z.
\end{equation}
\end{lemma}

\begin{proof}
Put \(\mathcal F=E_1\oplus\mathcal Z\) and let
\(P_{\mathcal F}\) be the \(L^2_\nu\)-orthogonal projection.
Since \(C_c^\infty(S^2T^*M)\) is dense in \(L^2_\nu\), the linear
subspace
\[
 P_{\mathcal F}\bigl(C_c^\infty(S^2T^*M)\bigr)
\]
is dense in the finite-dimensional space \(\mathcal F\), and hence is
all of \(\mathcal F\).  Choose a linear right inverse of
\(P_{\mathcal F}\) and restrict it to \(E_1\).  The image of a finite
basis has a common compact support, and
\eqref{eq:compact-profile-lift-properties} follows.
\end{proof}

\begin{definition}[Prepared \(W\)-profile disk]
\label{def:prepared-W-profile-disk}
Fix the lift \(W\) in Lemma~\ref{lem:compact-profile-lift}.
A prepared \(W\)-profile disk of radius \(r\) at time \(\tau_0\) is a
\(C^1\) common-margin family
\[
 \overline B_r^{E_1}(0)\ni v\longmapsto
 \mathbf z_v\in\Sigma_{\tau_0}^{k+2,\alpha}
\]
of strict prepared entrances satisfying
\begin{equation}\label{eq:abstract-W-profile-disk}
 H_v(\tau_0)=e^{-\gamma_1\tau_0}W(v).
\end{equation}
The \(C^1\) regularity here is always understood in one fixed
same-output soliton-conjugated chart for the full prepared Banach topology
\eqref{eq:prepared-model-Banach-norm}.  Uniformity as
\(\tau_0\to\infty\), however, is measured only in the low-order hybrid
topology that is propagated by
Lemma~\ref{lem:uniform-weighted-Schauder-restart}.  More precisely,
put \(m_\#=4\), and for two entrance states at the same time define
\begin{equation}\label{eq:late-entry-profile-distance}
 \begin{aligned}
 \mathfrak d_{{\rm le},\tau_0}(\mathbf z_1,\mathbf z_2)
 :={}&
 \mathfrak D_{m_\#}^{\rm hyb}(\tau_0;\mathbf z_1,\mathbf z_2)\\
 &+d_{{\rm ext},m_\#,0}^{+}(\mathbf z_1,\mathbf z_2)\\
 &+d_{{\rm ext},-1,0}^{\rm corr}(\mathbf z_1,\mathbf z_2)\\
 &+d_{{\rm Ggr},m_\#+2,0}^{0}(\mathbf z_1,\mathbf z_2)\\
 &+d_{{\rm Fgr},m_\#+1,0}^{++}(\mathbf z_1,\mathbf z_2)\\
 &+d_{{\rm F},m_\#+1,0}^{\rm glob}(\mathbf z_1,\mathbf z_2)\\
 &+\|R_1(\tau_0)-R_2(\tau_0)\|
      _{\mathfrak X_{\rm sc}^{m_\#+2,\alpha}}\\
 &+e^{\gamma_1\tau_0}
   \|H_1(\tau_0)-H_2(\tau_0)\|
       _{\mathfrak T_{{\rm sc},N}^{m_\#,\alpha}} .
 \end{aligned}
\end{equation}
where the first term is the entrance-time version of
\eqref{eq:two-state-hybrid-distance} in the common physical gauge.
For a sliced tangent vector \(\xi\) put
\begin{equation}\label{eq:late-entry-profile-tangent}
\begin{aligned}
 \|\xi\|_{{\rm le},\tau_0}
 &:=
 \mathfrak D_{m_\#}^{\rm hyb}[\xi](\tau_0)
 +d_{{\rm ext},m_\#,0}^{+}[\xi]
 +d_{{\rm ext},-1,0}^{\rm corr}[\xi]\\
 &\quad
 +d_{{\rm Ggr},m_\#+2,0}^{0}[\xi]
 +d_{{\rm Fgr},m_\#+1,0}^{++}[\xi]
 +d_{{\rm F},m_\#+1,0}^{\rm glob}[\xi]\\
 &\quad
 +\|D R_{\mathbf z}[\xi](\tau_0)\|
      _{\mathfrak X_{\rm sc}^{m_\#+2,\alpha}}\\
 &\quad
 +e^{\gamma_1\tau_0}
   \|D H_{\mathbf z}[\xi](\tau_0)\|
       _{\mathfrak T_{{\rm sc},N}^{m_\#,\alpha}} .
\end{aligned}
\end{equation}
The two exterior terms are respectively the buffered same-order trace
and the separated \(H^{-1}\) corridor memory in
\eqref{eq:inner-terminated-exterior-DeTurck}.  The
\(d_{{\rm Ggr},m_\#+2,0}^{0}\)-term is the larger-collar homogeneous
trace required by \eqref{eq:compact-graft-buffer-finite}.  The two
\(F\)-terms are respectively the larger-source-star homogeneous trace
retained in the separated memory of
Lemma~\ref{lem:localized-graft-F-coarse-memory} and the global initial
trace used to start the global \(F\)-equation.  The former remains tied
to the larger source star, while the latter is restarted only from its
global trace.  None of these entrance quantities is inferred from a
smaller-collar or same-star trace.  Thus the collapsing core is
recorded by the normalized tensor, whereas the unscaled closed metric
is recorded only on the fixed buffered noncollapsing physical cover.
The displayed order-\((m_\#+2)\) \(R\)-trace is the auxiliary
derivative needed to type the nonseparated target-connection forcing
in the Abel estimate.  We
require
\begin{equation}\label{eq:W-profile-disk-absolute-variation}
 \mathfrak d_{{\rm le},\tau_0}(\mathbf z_v,\mathbf z_0)
 \leq C_W|v|,
 \qquad
 \|D_v\mathbf z_v[w]\|_{{\rm le},\tau_0}
 \leq C_W|w|.
\end{equation}
When a family of such disks is indexed by late entrance times
\(\tau_0\), as in Theorem~\ref{thm:profile-realization}, the phrase
\emph{uniform late-entry \(W\)-profile disk} means that \(r\), the
common strict margins, and \(C_W\) in
\eqref{eq:W-profile-disk-absolute-variation} are independent of
\(\tau_0\).  It also includes one entrance-time-independent
low-order coefficient and chart package \(K_\#\) of
\eqref{eq:uniform-future-low-coefficient-package}, with a fixed
buffered physical cover and uniform ellipticity, radial-comparison,
inverse-map, and composition margins at those orders.  Here ``common
strict margins'' means uniform fractional
slack in the scale-normalized prepared inequalities after the
displayed \(e^{-\theta\tau_0}\) and \(e^{-\sigma\tau_0}\) entrance
weights have been factored out, together with the fixed buffered
geometric inclusions.  It does not mean an entrance-time-independent
absolute gap from a face which itself tends to zero.  No
entrance-time-independent bound is asserted for the full unscaled
\(C^{k+2,\alpha}(\mathcal X)\) chart norm of the closed metric; for
each fixed \(\tau_0\), that full norm still supplies the Banach
\(C^1\) structure used by the local theory.
For such a disk we write
\[
 \mathcal V_{\tau_0}(v):=V_\infty(\mathbf z_v).
\]
\end{definition}

For a sliced state \(\mathbf z\), let
\[
 U_{\mathbf z}=H_{\mathbf z},
 \qquad
 R_{\mathbf z}=P_{\geq1}H_{\mathbf z}=H_{\mathbf z}.
\]
The last equality holds because the sliced state lies in
\(\mathcal Z^\perp\), and \(P_{<1}=0\) there.  Let
\(F_{U_{\mathbf z}}\) denote the exact right-hand side in
\[
 \partial_\tau U_{\mathbf z}+LU_{\mathbf z}
 =F_{U_{\mathbf z}};
\]
thus, in the one-state equation,
\(F_{U_{\mathbf z}}=\mathcal N_{\mathbf z}\).

Use, without reselecting any term, the \(j=1\) instance of the fixed
\(\chi_{{\rm in},\tau}\)-based global decomposition
\eqref{eq:concrete-projected-graph-splitting}.  Accordingly define the
exact outer-residual functional
\begin{equation}\label{eq:O-R-functional}
 \begin{split}
 \mathcal O_R(\tau;\mathbf z)
 :={}&F_{U_{\mathbf z}}
 -A_{U_{\mathbf z}}^{ij}
    \bar\nabla_i\bar\nabla_jR_{\mathbf z}
 -\mathscr L_{U_{\mathbf z},1}R_{\mathbf z}\\
 &-S_{R_{\mathbf z}}
 -\mathcal Z_{U_{\mathbf z}} .
 \end{split}
\end{equation}
The coefficient and source functionals in
\eqref{eq:O-R-functional} are precisely those constructed from the
inner expression multiplied by \(\chi_{{\rm in},\tau}\) in the proof of
Lemma~\ref{lem:restricted-stable-semigroup}.  Thus
\(\mathcal O_R\) is the exact global residual already used in
\eqref{eq:concrete-projected-graph-splitting}; in particular it contains
the complete transition annulus, every \(\rho_\tau\)-commutator, the
pure graft, the effective-column defects, and every terminated-interface
term.

For a sliced tangent \(\xi\), put
\[
 K=DH_{\mathbf z}[\xi].
\]
At each fixed normalized time define the complete differentiated outer
package by the Fr\'echet derivative of this exact residual:
\begin{equation}\label{eq:O-K-definition}
 \boxed{\;
 \mathcal O_K(\tau;\mathbf z,\xi)
 :=
 D_{\mathbf z}\mathcal O_R(\tau;\mathbf z)[\xi].
 \;}
\end{equation}
Equivalently, since
\(D_{\mathbf z}R_{\mathbf z}[\xi]=K\),
\begin{align*}
 \mathcal O_K
 ={}&
 D_{\mathbf z}F_{U_{\mathbf z}}[\xi]
 -\bigl(D_{\mathbf z}A_{U_{\mathbf z}}^{ij}[\xi]\bigr)
      \bar\nabla_i\bar\nabla_jR_{\mathbf z}
 -A_{U_{\mathbf z}}^{ij}\bar\nabla_i\bar\nabla_jK\\
 &-\bigl(D_{\mathbf z}\mathscr L_{U_{\mathbf z},1}[\xi]\bigr)
      R_{\mathbf z}
 -\mathscr L_{U_{\mathbf z},1}K
 -D_{\mathbf z}S_{R_{\mathbf z}}[\xi]
 -D_{\mathbf z}\mathcal Z_{U_{\mathbf z}}[\xi].
\end{align*}
All derivatives are taken in the single fixed same-output
soliton-conjugated prepared chart.  Consequently
\eqref{eq:O-K-definition} automatically differentiates every
state-dependent cutoff, graft or interface coefficient, effective
column, phase and harmonic-map pullback, and moving support occurring
in the exact equation.  The fixed cutoff
\(\chi_{{\rm in},\tau}\) itself has zero state derivative, while all of
its transition-annulus terms remain in \(\mathcal O_R\) and hence in
its differentiated residual.  No separately defined operation
\([\,\cdot\,]_{\rm out}\), support restriction, or extension convention
is used.

Define the differentiated core source by
\[
 D\mathcal N_{{\rm core},\mathbf z}[\xi]
 :=
 D_{\mathbf z}\!\left[
  A_{U_{\mathbf z}}^{ij}
    \bar\nabla_i\bar\nabla_jR_{\mathbf z}
  +\mathscr L_{U_{\mathbf z},1}R_{\mathbf z}
  +S_{R_{\mathbf z}}
  +\mathcal Z_{U_{\mathbf z}}
 \right][\xi].
\]
Then the exact differentiated stable source identity is
\[
 D\mathcal N_{\mathbf z}[\xi]
 =
 D\mathcal N_{{\rm core},\mathbf z}[\xi]
 +\mathcal O_K(\tau;\mathbf z,\xi).
\]
When \(\mathbf z,\xi\) are understood, we continue to write
\(\mathcal O_K\); for \(K=K_{v,w}\) this is the downstream notation
\(\mathcal O_{K_{v,w}}\).

\begin{lemma}[Late-entry hybrid first-variation propagation]
\label{lem:late-entry-hybrid-first-variation}
Let \(v\mapsto\mathbf z_v\) be a uniform late-entry prepared
\(W\)-profile disk and let
\(\xi_{v,w}=D_v\mathbf z_v[w]\).  Write
\[
 K_{v,w}(\tau)=D H_{\mathbf z_v(\tau)}[\xi_{v,w}(\tau)]
\]
and \(I_s=[s,s+1]\), \(s\geq\tau_0\).  There are constants \(A,C\),
independent of every sufficiently large entrance time, such that
\begin{equation}\label{eq:late-entry-hybrid-first-variation}
\begin{split}
 &\sup_{\tau\in I_s}\left(
   \mathfrak D_{4}^{\rm hyb}[\xi_{v,w}](\tau)
   +\mathfrak F_{5}^{\rm glob}[\xi_{v,w}](\tau)
   +\|K_{v,w}(\tau)\|_{L^2_\nu}\right)\\
 &\qquad+
 \left(\int_s^{s+1}
   \|K_{v,w}(q)\|_{H^1_\nu}^2\,dq\right)^{1/2}
 \leq
 C e^{A(s-\tau_0)}
   \|\xi_{v,w}\|_{{\rm le},\tau_0}.
\end{split}
\end{equation}
In particular, the left side includes the coarse global map block,
the exterior-memory block
\(\mathfrak G_4[\xi_{v,w}]\) and the compact graft input block
\(\mathfrak B_{{\rm gr},4}[\xi_{v,w}]\); neither is estimated from an
entrance-time-uniform full prepared tangent norm.  If
\(\mathcal O_{K_{v,w}}\) is the package in
\eqref{eq:O-K-definition}, then, after decreasing \(c>0\),
\begin{equation}\label{eq:late-entry-hybrid-outer-consequence}
 \|\mathcal O_{K_{v,w}}(\tau)\|_{L^2_\nu}
 +\|\mathcal O_{K_{v,w}}(\tau)\|_{H^{-1}_\nu}
 \leq
 C e^{-ce^{\tau/2}}
 \|\xi_{v,w}\|_{{\rm le},\tau_0}.
\end{equation}
\end{lemma}

\begin{proof}
Differentiate the exact low-order block systems used in
Lemma~\ref{lem:uniform-weighted-Schauder-restart}.  Their initial terms
are exactly the tangent blocks in
\eqref{eq:late-entry-profile-tangent}.  The low-initial-data clause
\eqref{eq:uniform-restart-low-initial-propagation} therefore gives the
normalized tensor, map, scale, clock, marking, and interface terms in
\eqref{eq:late-entry-hybrid-first-variation}.  The linearized
inner-terminated estimate in
Lemma~\ref{lem:inner-terminated-exterior-DeTurck} retains the entrance
exterior norm and the separated interface memory.  The compact graft
term is supplied by
Lemma~\ref{lem:compact-graft-buffer-propagation}, including its
linearized order-four clause.  The localized \(F\)-term is supplied
from the larger source-star trace by
Lemma~\ref{lem:localized-graft-F-coarse-memory}; the global \(F\)-term
starts from the global initial trace in
\eqref{eq:late-entry-profile-tangent} and is subsequently restarted
from its global value.  Thus all entrance terms are controlled by the
same low norm, and unit-interval iteration gives the
factor \(e^{A(s-\tau_0)}\).  This proves
\eqref{eq:late-entry-hybrid-first-variation}.

Every term in \(\mathcal O_{K_{v,w}}\) is a fixed low-order expression
in these blocks and is supported where
\(\bar f\geq ce^{\tau/2}\); the graft, cutoff, column-tail, and
interface pieces are supported where \(\bar f\geq ce^\tau\).
The Gaussian tail contributes \(e^{-c_{\rm G,le}e^{\tau/2}}\) for one
fixed \(c_{\rm G,le}>0\).  Since
\[
 e^{-c_{\rm G,le}e^{\tau/2}}e^{A(\tau-\tau_0)}
 \leq C e^{-ce^{\tau/2}},
 \qquad \tau\geq\tau_0,
\]
the growth in \eqref{eq:late-entry-hybrid-first-variation} is absorbed,
which proves \eqref{eq:late-entry-hybrid-outer-consequence}.
\end{proof}

\begin{lemma}[Uniform first-unit stable variation on a late-entry disk]
\label{lem:late-entry-first-unit-stable-variation}
Let
\[
 v\longmapsto\mathbf z_v,\qquad |v|\leq r,
\]
be a uniform late-entry prepared \(W\)-profile disk at a sufficiently
large entrance time \(\tau_0\).  Put
\[
 J_0=[\tau_0,\tau_0+1],\qquad
 H_v=H(\mathbf z_v),\qquad
 K_{v,w}=D_vH_v[w],\qquad
 \mathcal N'_{v,w}=D_v\mathcal N_v[w].
\]
After decreasing \(c>0\), uniformly in \(|v|\leq r\),
\begin{equation}\label{eq:late-entry-first-unit-one-state}
\begin{split}
 &\sup_{\tau\in J_0}\|H_v(\tau)\|_{L^2_\nu}
 +\left(\int_{J_0}\|H_v(s)\|_{H^1_\nu}^2\,ds\right)^{1/2}\\
 &\qquad\leq
 C_r e^{-\gamma_1\tau_0}
 +C_r e^{-ce^{\tau_0/2}},
\end{split}
\end{equation}
and
\begin{equation}\label{eq:late-entry-first-unit-variation}
\begin{split}
 &\sup_{\tau\in J_0}\|K_{v,w}(\tau)\|_{L^2_\nu}
 +\left(\int_{J_0}
       \|K_{v,w}(s)\|_{H^1_\nu}^2\,ds\right)^{1/2}\\
 &\qquad\leq
 C_r e^{-\gamma_1\tau_0}|w|
 +C_r e^{-ce^{\tau_0/2}}|w|.
\end{split}
\end{equation}
Moreover the exact one-state stable source satisfies
\begin{equation}
\label{eq:late-entry-first-unit-one-state-projected-source}
\begin{split}
 \int_{\tau_0}^{\tau_0+1}
 e^{\gamma_1s}
 \left|\Pi_1\mathcal N_v(s)\right|\,ds
 \leq C_r\bigl(
  e^{-\kappa\tau_0}
  +e^{-\gamma_1\tau_0}
  +e^{-ce^{\tau_0/2}}\bigr).
\end{split}
\end{equation}
The exact differentiated stable source satisfies
\begin{equation}\label{eq:late-entry-first-unit-projected-source}
\begin{split}
 \int_{\tau_0}^{\tau_0+1}
 e^{\gamma_1s}
 \left|\Pi_1\mathcal N'_{v,w}(s)\right|\,ds
 \leq C_r\bigl(
  e^{-\kappa\tau_0}
  +e^{-\gamma_1\tau_0}
  +e^{-ce^{\tau_0/2}}\bigr)|w|.
\end{split}
\end{equation}
The constants are independent of every sufficiently large
\(\tau_0\).
\end{lemma}

\begin{proof}
The special entrance identities and the fixed compact support of
\(W\) give
\[
 \|H_v(\tau_0)\|_{H^1_\nu}
 \leq C_r e^{-\gamma_1\tau_0},
 \qquad
 \|K_{v,w}(\tau_0)\|_{H^1_\nu}
 \leq C e^{-\gamma_1\tau_0}|w|.
\]
On \(J_0\), use the exact one-state and variational stable equations,
not the late-time estimate
\eqref{eq:variational-stable-Hminus1}.  The order-four prepared chart
calculus, integration by parts in the quasilinear principal terms,
and stable form coercivity give a constant \(c_{\rm en}>0\), chosen
after all perturbative absorptions and independent of every
sufficiently large \(\tau_0\), for which the following first-unit form
inequalities hold:
\begin{align}
 \frac d{d\tau}\|H_v\|_{L^2_\nu}^2
  +c_{\rm en}\|H_v\|_{H^1_\nu}^2
 &\leq
 C\|H_v\|_{L^2_\nu}^2
 +C e^{-ce^{\tau/2}},
 \label{eq:late-entry-first-unit-one-state-energy}\\
 \frac d{d\tau}\|K_{v,w}\|_{L^2_\nu}^2
  +c_{\rm en}\|K_{v,w}\|_{H^1_\nu}^2
 &\leq
 C\|K_{v,w}\|_{L^2_\nu}^2
 +C e^{-ce^{\tau/2}}|w|^2 .
 \label{eq:late-entry-first-unit-variation-energy}
\end{align}
We justify the two displayed inequalities term by term.  In the
collapsing core, the differentiated
second-order term is integrated once; the undifferentiated
principal-coefficient perturbation is absorbed by the strict
\(C^2\) box.  When a derivative falls on a coefficient, the common
low-order prepared package and Young's inequality contribute only a
uniform multiple of the \(L^2_\nu\) energy.  The remaining core
quadratic terms obey the same relative form bound because the global
\(C^2\) box is chosen below the fixed absorption threshold.

We record separately the scalar-feedback term which is not literally
divisible by \(K_{v,w}\).  Differentiating the exact Gram system and
using its modal core estimate together with
\eqref{eq:late-entry-hybrid-outer-consequence} gives, on \(J_0\),
with \(c_v(\tau)=c(\tau,\mathbf z_v(\tau))\),
\begin{equation}\label{eq:late-entry-first-unit-differentiated-Gram}
 |D_vc_v(\tau)[w]|
 \leq
 C\bigl(e^{-\kappa\tau}
        +\|H_v(\tau)\|_{H^1_\nu}\bigr)
   \|K_{v,w}(\tau)\|_{H^1_\nu}
 +CC_W e^{-ce^{\tau/2}}|w|.
\end{equation}
Indeed, this is the differential form of
\eqref{eq:two-state-sharp-modal}, obtained by taking the differential
quotient in the exact two-state modal identity
\(M_1(c_1-c_2)=-(d_1-d_2)-(M_1-M_2)c_2\);
the first term on the right is the differentiated core contribution and
the last is the differentiated effective-column and receding-support
contribution.  In particular,
\[
 D_v\!\left(c_{v,j}\rho_\tau\mathscr B_jh_v\right)[w]
 =
 (D_vc_{v,j}[w])\rho_\tau\mathscr B_jh_v
 +c_{v,j}\rho_\tau\mathscr B_j(D_vh_v[w]).
\]
On the core \(D_vh_v[w]=K_{v,w}\); the discrepancy created by
\(K_{v,w}=\rho_\tau D_vh_v[w]\) is supported in the cutoff annulus and
already belongs to \(\mathcal O_{K_{v,w}}\).
The bilinear
\(H^1_\nu\times H^1_\nu\to H^{-1}_\nu\) bound for
\(\mathscr B_j\), the small one-state box, and
\eqref{eq:late-entry-first-unit-differentiated-Gram} put the
\(K_{v,w}\)-dependent part into the same relative form bound as the
other core terms and leave only
\[
 CC_W e^{-ce^{\tau/2}}|w|
\]
as an additive core source.  Thus every other nonouter term contains
\(K_{v,w}\), and the exact Gram system has the required bounded
relative-form structure.

All cutoff, effective-column,
graft, moving-support, map, and interface terms are the outer package
\(\mathcal O_{K_{v,w}}\); by
\eqref{eq:late-entry-hybrid-outer-consequence} and
\eqref{eq:W-profile-disk-absolute-variation}, their \(H^{-1}_\nu\)
norm is bounded by
\[
 C C_W e^{-ce^{\tau/2}}|w|.
\]
The one-state outer terms obey the analogous bound with \(C_r\).
Pairing the outer terms and the additive scalar-feedback core source
with the corresponding stable tensor and applying Young's inequality
gives the last terms in the two energy inequalities, after decreasing
\(c\).  Thus neither inequality uses an entrance-time-uniform full
prepared tangent norm.

Integrating
\eqref{eq:late-entry-first-unit-one-state-energy} and using Gronwall
on the unit interval proves
\eqref{eq:late-entry-first-unit-one-state}.  A second Gronwall
argument applied to
\eqref{eq:late-entry-first-unit-variation-energy} therefore proves
\eqref{eq:late-entry-first-unit-variation}.

It remains to estimate the projected sources.  Let
\(\mathcal N_{v,{\rm core}}\) denote the nonouter part of the one-state
source; the complementary one-state outer package has
\(H^{-1}_\nu\)-norm at most \(C_re^{-ce^{s/2}}\) by the calculation
above.  On \(J_0\), the same termwise calculation before pairing gives
\[
 \|\mathcal N_{v,{\rm core}}(s)\|_{H^{-1}_\nu}
 \leq C\bigl(e^{-\kappa s}
              +\|H_v(s)\|_{H^1_\nu}\bigr)
        \|H_v(s)\|_{H^1_\nu}
\]
and
\[
\begin{split}
 \|\mathcal N'_{v,w,{\rm core}}(s)\|_{H^{-1}_\nu}
 &\leq C\bigl(e^{-\kappa s}
               +\|H_v(s)\|_{H^1_\nu}\bigr)\\
 &\qquad\cdot\|K_{v,w}(s)\|_{H^1_\nu}
 +CC_W e^{-ce^{s/2}}|w|.
\end{split}
\]
The finite-rank map
\(\Pi_1:H^{-1}_{\nu,-}\to E_1\) is bounded.  Hence
Cauchy--Schwarz, the two estimates just proved, and the two Gaussian
core and outer bounds yield, with
\[
 \mathsf H_0:=
 \left(\int_{J_0}\|H_v\|_{H^1_\nu}^2\right)^{1/2},
 \qquad
 \mathsf K_0:=
 \left(\int_{J_0}\|K_{v,w}\|_{H^1_\nu}^2\right)^{1/2},
\]
\[
\begin{split}
 &\int_{J_0}e^{\gamma_1s}
   |\Pi_1\mathcal N_v(s)|\,ds\\
 &\quad\leq
 C e^{\gamma_1(\tau_0+1)}
 \left[
  e^{-\kappa\tau_0}\mathsf H_0
  +\mathsf H_0^2
  +e^{-ce^{\tau_0/2}}
 \right].
\end{split}
\]
Substitution of \eqref{eq:late-entry-first-unit-one-state} proves
\eqref{eq:late-entry-first-unit-one-state-projected-source}; the
Gaussian term absorbs the harmless factor
\(e^{\gamma_1(\tau_0+1)}\).  Likewise,
\[
\begin{split}
 &\int_{J_0}e^{\gamma_1s}
   |\Pi_1\mathcal N'_{v,w}(s)|\,ds\\
 &\quad\leq
 C e^{\gamma_1(\tau_0+1)}
 \left[
  e^{-\kappa\tau_0}\mathsf K_0
  +\mathsf H_0\mathsf K_0
  +e^{-ce^{\tau_0/2}}|w|
 \right].
\end{split}
\]
Absorbing the harmless factor \(e^{\gamma_1\tau_0}\) into the
Gaussian-superexponential term proves
\eqref{eq:late-entry-first-unit-projected-source}.
\end{proof}

\subsection{Differentiable two-state dynamics}

We now turn to parameter dependence.
The one-state profile, quadratic-response, and frozen-marking branches
have already been completed.  The results below use them only as
established inputs and do not enter their proofs.

\begin{lemma}[Two-state stable normal form in \(H^{-1}_\nu\)]
\label{lem:two-state-stable-normal-form}
Let two global strict prepared evolutions belong to one fixed
common-margin ball
\(\mathscr B\subset\Sigma_{\tau_0}^{k+2,\alpha}\) satisfying the shared
buffered-cover and strict time-width hypotheses of
Theorem~\ref{thm:global-two-state-estimate}, let
\[
 D=H_1-H_2,\qquad
 d_0=\|\mathbf z_{1,0}-\mathbf z_{2,0}\|
       _{\mathscr X_{\rm prep}^{k+2,\alpha}},
\]
and let \(\mathcal N_i\) be the forcings in
\eqref{eq:stable-N-definition}.  Then, for \(\tau\geq\tau_0+1\),
\begin{equation}\label{eq:two-state-stable-Hminus1}
 \|\mathcal N_1-\mathcal N_2\|_{H^{-1}_\nu}
 \leq
 C\left(
   e^{-\kappa\tau}
   +\|H_1\|_{H^1_\nu}
   +\|H_2\|_{H^1_\nu}\right)
 \|D\|_{H^1_\nu}
 {}+Cd_0e^{-ce^{\tau/2}}.
\end{equation}
For a base entrance
\(\mathbf z_0\in\mathscr B\) and a tangent vector
\(\xi\in T_{\mathbf z_0}\Sigma_{\tau_0}^{k+2,\alpha}\), set
\(K=DH_{\mathbf z_0}[\xi]\), and give \(\xi\) the tangent norm induced
by \(\mathscr E_{\rm prep}^{k+2,\alpha}\).  Then
\begin{equation}\label{eq:variational-stable-Hminus1}
 \|D\mathcal N_{\mathbf z_0}[\xi]\|_{H^{-1}_\nu}
 \leq
 C\left(e^{-\kappa\tau}+\|H\|_{H^1_\nu}\right)
 \|K\|_{H^1_\nu}
 {}+Ce^{-ce^{\tau/2}}\|\xi\|_{\rm prep}.
\end{equation}
Both estimates include the cutoff commutator, effective columns,
solution-dependent moving support, and graft error.
\end{lemma}

\begin{proof}
Set \(w=h_1-h_2\).  On
\(\{\bar f\leq e^{\tau/2}\}\), where \(D=w\), expand the exact
quasilinear difference as
\[
 \Q(h_1)-\Q(h_2)
 =
 -\widehat h_1*\bar\nabla^2w
 -(\widehat h_1-\widehat h_2)*\bar\nabla^2h_2
 {}+\text{terms bilinear in }(w,\bar\nabla w)
\]
with coefficients depending smoothly on \(h_1,h_2,\bar\nabla h_1,
\bar\nabla h_2\).  Integrating one derivative in the first two terms
against an arbitrary unit \(H^1_\nu\) test tensor, and using
Lemma~\ref{lem:first-moment} when that derivative hits the Gaussian,
gives
\[
 \|\rho_\tau(\Q(h_1)-\Q(h_2))\|_{H^{-1}_\nu}
 \leq
 Ce^{-\kappa\tau}\|D\|_{H^1_\nu}
 {}+C(\|H_1\|_{H^1_\nu}+\|H_2\|_{H^1_\nu})
      \|D\|_{H^1_\nu}
 {}+Cd_0e^{-ce^{\tau/2}}.
\]
Here the first coefficient follows from the common inner \(C^1\)
decay.  On the complement, the common global \(C^2\) box and the
coarse two-state estimate
\eqref{eq:coarse-two-state-growth} give at most \(Cd_0e^{A\tau}\)
growth, while the region starts at
\(\bar f=e^{\tau/2}\); the Gaussian absorbs this factor.

The same support calculation, now using
Lemma~\ref{lem:prepared-chart-calculus}, gives the
explicit outer-package estimate
\begin{equation}\label{eq:two-state-complete-outer-package}
 \begin{split}
 &\|\mathcal C_\rho[h_1]-\mathcal C_\rho[h_2]\|_{H^{-1}_\nu}
 +\|\rho_\tau(\E_1-\E_2)\|_{H^{-1}_\nu}\\
 &\quad+
 \sum_{j=0}^8
 \|\rho_\tau(\mathcal Y^{(1)}_{j,\tau}
              -\mathcal Y^{(2)}_{j,\tau})\|_{H^{-1}_\nu}
 \leq Cd_0e^{-ce^{\tau/2}}.
 \end{split}
\end{equation}
Thus \eqref{eq:two-state-complete-outer-package} records the complete
outer-error contribution explicitly.

Subtracting the exact Gram systems and using the modal version of the
preceding quasilinear calculation gives
\[
 |c_1-c_2|
 \leq
 C(\|H_1\|_{H^1_\nu}+\|H_2\|_{H^1_\nu})
   \|D\|_{H^1_\nu}
 {}+Cd_0e^{-ce^{\tau/2}},
\]
whereas the one-state system gives
\[
 |c_i|\leq C\|H_i\|_{H^1_\nu}^2+Ce^{-ce^{\tau/2}}.
\]
Use the bilinear \(H^1_\nu\times H^1_\nu\to H^{-1}_\nu\) bounds for
\(\mathscr B_j\) to estimate
\[
 (c_{1,j}-c_{2,j})\mathscr B_jh_1
 +c_{2,j}\mathscr B_jw.
\]
The common small box absorbs the extra factor
\(\|H_1\|+\|H_2\|\).  The direct-column differences are controlled by
\eqref{eq:two-state-complete-outer-package}.  Combining these bounds
term by term in \eqref{eq:stable-N-definition} proves
\eqref{eq:two-state-stable-Hminus1}.

On every finite horizon the coupled solution and every term just
estimated are \(C^1\) in the buffered prepared chart.  Differentiate
the termwise proof above in an arbitrary
\(\mathscr E_{\rm prep}^{k+2,\alpha}\) model tangent direction.
The differentiated quasilinear, Gram, bilinear, cutoff-commutator, and
effective-column terms obey the displayed linearized product bounds.
The differentiated outer package remains supported on the same
receding region and is bounded by the Gaussian-tail estimate with the
model tangent norm.  This proves
\eqref{eq:variational-stable-Hminus1}.
\end{proof}

\begin{lemma}[Differentiated stable normal form]
\label{lem:differentiated-stable-normal-form}
Let \(\mathcal V\) be a \(C^1\) Banach parameter manifold, let
\(v\mapsto\mathbf z_v\) be a \(C^1\) common-margin family of strict
prepared entrances in \(\Sigma_{\tau_0}^{k+2,\alpha}\) whose image lies
in one fixed ball satisfying the shared buffered-cover and strict
time-width hypotheses of
Theorem~\ref{thm:global-two-state-estimate}, and let
\(w\in T_v\mathcal V\).  Put
\[
 \xi_{v,w}:=D_v\mathbf z_v[w]\in
 T_{\mathbf z_v}\Sigma_{\tau_0}^{k+2,\alpha}
\]
and write
\[
 K_{v,w}:=D H(\mathbf z_v)[\xi_{v,w}]
          =D_vH_v[w],\qquad
 \mathcal N'_{v,w}:=
 D\mathcal N(\mathbf z_v)[\xi_{v,w}]
          =D_v\mathcal N_v[w].
\]
We write
\(\|\xi_{v,w}\|_{\rm prep}\) for the tangent norm induced by
\(\mathscr E_{\rm prep}^{k+2,\alpha}\).  On the
finite-dimensional profile disk \(\mathcal V=B_r(E_1)\),
\(w\in E_1\) and
\(\|\xi_{v,w}\|_{\rm prep}\leq C_{r,\tau_0}|w|\) at each fixed
entrance time.  The constant in this full unscaled chart norm is not
used for late-entry uniformity.  On a uniform late-entry disk,
\eqref{eq:W-profile-disk-absolute-variation} instead gives
\(\|\xi_{v,w}\|_{{\rm le},\tau_0}\leq C_W|w|\).
For the buffered ambient variant, start instead with a \(C^1\) family
\[
 v\longmapsto\mathbf z_v^{\rm amb}\in\operatorname{dom}
 \Pi_{\rm sl}^{\,k+4\to k+2}
\]
whose phase-retracted image lies in the same fixed sliced ball, and put
\[
 \widetilde{\mathbf z}_v
 :=\Pi_{\rm sl}^{\,k+4\to k+2}(\mathbf z_v^{\rm amb}),
 \qquad
 \xi_{v,w}^{\rm amb}:=D_v\mathbf z_v^{\rm amb}[w],
 \qquad
 \widetilde\xi_{v,w}
 :=D\Pi_{\rm sl}^{\,k+4\to k+2}(\mathbf z_v^{\rm amb})
   [\xi_{v,w}^{\rm amb}].
\]
In this ambient variant, henceforth set
\[
 \mathbf z_v:=\widetilde{\mathbf z}_v,
 \qquad
 \xi_{v,w}:=\widetilde\xi_{v,w}.
\]
Thus every quantity below is evaluated at
\((\widetilde{\mathbf z}_v,\widetilde\xi_{v,w})\); in particular,
\[
 K_{v,w}:=DH(\widetilde{\mathbf z}_v)[\widetilde\xi_{v,w}],
 \qquad
 \mathcal N'_{v,w}
 :=D\mathcal N(\widetilde{\mathbf z}_v)[\widetilde\xi_{v,w}].
\]
Under this convention the prepared tangent norm satisfies
\[
 \|\xi_{v,w}\|_{\rm prep}
 =\|\widetilde\xi_{v,w}\|_{
       \mathscr E_{\rm prep}^{k+2,\alpha}}
 \leq K_{\Pi,k+2}
       \|\xi_{v,w}^{\rm amb}\|_{
       \mathscr E_{\rm prep}^{k+4,\alpha}},
\]
Consequently every general ambient variation constant below is at most
its sliced counterpart multiplied by \(K_{\Pi,k+2}\).  Any conclusion
below that is conditional on a uniform late-entry prepared
\(W\)-profile disk is asserted in this ambient variant only when the
retracted family
\(v\mapsto\widetilde{\mathbf z}_v\) is such a disk, in the sense of
Definition~\ref{def:prepared-W-profile-disk}, with the same uniform
late-entry constants.
Then
\[
 \partial_\tau K_{v,w}
 =\A K_{v,w}+\mathcal N'_{v,w},
 \qquad K_{v,w}\perp\mathcal Z,
\]
and, for \(\tau\geq\tau_0+1\),
\begin{equation}\label{eq:differentiated-tame-bound}
 \begin{split}
 \|\mathcal N'_{v,w}(\tau)\|_{H^{-1}_\nu}
 \leq\;&
 Ce^{-\kappa\tau}\|K_{v,w}(\tau)\|_{H^1_\nu}\\
  &+C\|H_v(\tau)\|_{H^1_\nu}
        \|K_{v,w}(\tau)\|_{H^1_\nu}
  +Ce^{-ce^{\tau/2}}\|\xi_{v,w}\|_{\rm prep}.
 \end{split}
\end{equation}
For a uniform late-entry \(W\)-profile disk,
Lemma~\ref{lem:late-entry-hybrid-first-variation} applies to
\(\xi_{v,w}\).  Its outer consequence gives the sharper hybrid
estimate for the outer part of the linearized source.  Precisely,
\eqref{eq:O-K-definition} gives
\(\mathcal N'_{v,w}=\mathcal N'_{v,w,{\rm core}}
+\mathcal O_{K_{v,w}}\).
Thus the complete outer package is supported where
\(\bar f\geq ce^{\tau/2}\); the latter graft, cutoff, column, and
interface pieces are in fact supported farther out, where
\(\bar f\geq ce^\tau\).  Then
\begin{equation}\label{eq:late-entry-profile-outer-source}
 \|\mathcal O_{K_{v,w}}(\tau)\|_{L^2_\nu}+{}
 \|\mathcal O_{K_{v,w}}(\tau)\|_{H^{-1}_\nu}
 \leq
 C e^{-ce^{\tau/2}}\|\xi_{v,w}\|_{{\rm le},\tau_0}
 \leq C C_W e^{-ce^{\tau/2}}|w|.
\end{equation}
Consequently, on that disk and for
\(\tau\geq\tau_0+1\),
\begin{equation}\label{eq:late-entry-differentiated-tame-bound}
 \begin{split}
 \|\mathcal N'_{v,w}(\tau)\|_{H^{-1}_\nu}
 \leq\;&
 Ce^{-\kappa\tau}\|K_{v,w}(\tau)\|_{H^1_\nu}\\
 &+C\|H_v(\tau)\|_{H^1_\nu}
       \|K_{v,w}(\tau)\|_{H^1_\nu}
 +C C_W e^{-ce^{\tau/2}}|w|.
 \end{split}
\end{equation}
In the concrete projected graph recovery
\eqref{eq:projected-stable-graph-recovery}, the same estimate replaces
the general first-variation choice of \(s_R\) by \(C_W|w|\).  Thus, on
the late-entry disk, its outer error is
\begin{equation}\label{eq:late-entry-profile-graph-error}
 \mathfrak e_R^{\rm le}(\tau)
 \leq C C_W^2|w|^2e^{-ce^{\tau/2}} .
\end{equation}
On every fixed common-margin sliced prepared ball satisfying those
shared-cover and strict time-width hypotheses, for each
\[
 0<\zeta<\min\{\kappa,\gamma_1/2\},
\]
there is a uniform variational estimate
\begin{equation}\label{eq:general-profile-variation-rate}
 \|K_{v,w}(\tau)\|_{H^1_\nu}
 \leq
 C_\zeta e^{-(\gamma_1-\zeta)(\tau-\tau_0)}
 \|\xi_{v,w}\|_{\rm prep},
 \qquad \tau\geq\tau_0+1.
\end{equation}
For a uniform late-entry prepared \(W\)-profile disk in the sense of
Definition~\ref{def:prepared-W-profile-disk},
\[
 H_v(\tau_0)=e^{-\gamma_1\tau_0}W(v)
\]
and every
\[
 0<\zeta<\min\{\kappa,\gamma_1/2\},
\]
one has, uniformly for \(|v|\leq r\),
\begin{equation}\label{eq:profile-variation-rate}
 \|K_{v,w}(\tau)\|_{H^1_\nu}
 \leq C_{r,\zeta}
 e^{-(\gamma_1-\zeta)\tau}|w|,
 \qquad \tau\geq\tau_0+1.
\end{equation}
Here and in \eqref{eq:profile-map-C1-error},
\(C_{r,\zeta}\) is independent of every sufficiently large entrance
time \(\tau_0\).
Moreover, for the \(C^1\) family
\(v\mapsto\mathbf z_v\) in the statement of the lemma, the map
\[
 v\longmapsto V_\infty(\mathbf z_v)
\]
is \(C^1\), and its derivative in the direction \(w\) is
\begin{equation}\label{eq:general-profile-map-derivative}
 \begin{split}
 D_vV_\infty(\mathbf z_v)[w]
 ={}&
 e^{\gamma_1\tau_0}\Pi_1K_{v,w}(\tau_0)
 +\int_{\tau_0}^{\tau_0+1}
  e^{\gamma_1s}\Pi_1\mathcal N'_{v,w}(s)\,ds\\
 &+\int_{\tau_0+1}^{\infty}
  e^{\gamma_1s}\Pi_1\mathcal N'_{v,w}(s)\,ds .
 \end{split}
\end{equation}
For a uniform late-entry prepared \(W\)-profile disk,
\[
 K_{v,w}(\tau_0)=e^{-\gamma_1\tau_0}W(w),
\]
so \eqref{eq:general-profile-map-derivative} reduces to
\begin{equation}\label{eq:profile-map-derivative}
 \begin{split}
 D\mathcal V_{\tau_0}(v)[w]
 ={}&
 w+\int_{\tau_0}^{\tau_0+1}
 e^{\gamma_1s}\Pi_1\mathcal N'_{v,w}(s)\,ds\\
 &+\int_{\tau_0+1}^{\infty}
 e^{\gamma_1s}\Pi_1\mathcal N'_{v,w}(s)\,ds .
 \end{split}
\end{equation}
Only for this transverse disk does the leading term equal \(w\).
Uniformly for \(|v|\leq r\),
\begin{equation}\label{eq:profile-map-C1-error}
 \begin{split}
 \sup_{|v|\leq r}
 \|D\mathcal V_{\tau_0}(v)-\operatorname{Id}_{E_1}\|
 \leq C_{r,\zeta}\bigl(
 &e^{-(\kappa-\zeta)\tau_0}\\
 &+e^{-(\gamma_1-2\zeta)\tau_0}
 +e^{-ce^{\tau_0/2}}\bigr).
 \end{split}
\end{equation}
On any \(C^1\) prepared Banach chart
\(v\mapsto\mathbf z_v\) with image in such a fixed common-margin ball
\(\mathscr B\), there is a function
\(\epsilon_{\rm prof}(S)\to0\) such that
\begin{equation}\label{eq:uniform-profile-derivative-tail}
 \sup_{\{v:\mathbf z_v\in\mathscr B\}}
 \sup_{\|\xi_{v,w}\|_{\rm prep}\leq1}
 \int_S^\infty e^{\gamma_1s}
 \left|\Pi_1\mathcal N'_{v,w}(s)\right|\,ds
 \leq\epsilon_{\rm prof}(S).
\end{equation}
\end{lemma}

\begin{proof}
The finite-horizon \(C^1\) statement in
Proposition~\ref{prop:two-state-prepared-evolution} permits
differentiation of the exact sliced equation.  The first-variation
estimate \eqref{eq:variational-stable-Hminus1}, which was obtained
term by term from that equation and includes every moving-support
term, is exactly \eqref{eq:differentiated-tame-bound} on its stated
range \(\tau\geq\tau_0+1\).

For the late-entry refinement, apply
Lemma~\ref{lem:late-entry-hybrid-first-variation}.  Its estimate
\eqref{eq:late-entry-hybrid-outer-consequence} is precisely
\eqref{eq:late-entry-profile-outer-source}; crucially, the entrance
factor is the special-disk low norm and not the full unscaled prepared
norm.  The collapsing-core variation is \(K_{v,w}\) itself and belongs
to the first two, nonouter terms of
\eqref{eq:differentiated-tame-bound}.  This also proves
\eqref{eq:late-entry-differentiated-tame-bound} on the same range.
Repeating the
\(\langle LR,\cdot\rangle_\nu\) pairing in the proof of
Lemma~\ref{lem:restricted-stable-semigroup}, now with the
\(L^2_\nu\) part of \eqref{eq:late-entry-profile-outer-source}, gives
\[
 \bigl|\langle LR,\mathcal O_{K_{v,w}}\rangle_\nu\bigr|
 \leq \varepsilon_*\|LR\|_{L^2_\nu}^2+{}
 C_{\varepsilon_*}C_W^2|w|^2e^{-ce^{\tau/2}}.
\]
All nonouter terms in that proof are unchanged, so this is precisely
\eqref{eq:late-entry-profile-graph-error} and proves the asserted
special-disk version of
\eqref{eq:projected-stable-graph-recovery}.

The sliced variational equation has no component in \(\mathcal Z\).
Fix \(0<\sigma_-<q_0<\theta\).  Apply
\eqref{eq:global-variation-energy} with \(\theta_-=q_0\) to obtain
\begin{equation}\label{eq:variational-coarse-energy}
 \|K_{v,w}(\tau)\|_{L^2_\nu}
 +\left(\int_\tau^\infty
        \|K_{v,w}(s)\|_{H^1_\nu}^2\,ds\right)^{1/2}
 \leq
 C e^{-q_0(\tau-\tau_0)}
 \|\xi_{v,w}\|_{\rm prep}.
\end{equation}
 This is an integrated \(H^1_\nu\) estimate, so we recover the
 pointwise norm before invoking the linearized bootstrap.  Apply
\eqref{eq:projected-stable-graph-recovery} with \(j=1\) to
 \(K_{v,w}\).  Its lower-block source is identically zero, and its
 outer forcing in the variational equation is linear:
\begin{equation}\label{eq:variational-outer-source-linear}
 \|\mathcal O_{K_{v,w}}(\tau)\|_{H^{-1}_\nu}
 \leq Ce^{-ce^{\tau/2}}\|\xi_{v,w}\|_{\rm prep}.
\end{equation}
When the equation is tested against \(K_{v,w}\), Young's inequality
places its contribution in the squared graph estimate as
\[
 C e^{-ce^{\tau/2}}\|\xi_{v,w}\|_{\rm prep}^2.
\]
Thus the occurrence of the squared tangent size in
 \(\mathfrak e_R=s_R^2e^{-ce^{\tau/2}}\) and in
\eqref{eq:projected-stable-graph-recovery} is required by homogeneity:
after taking the square root it gives the desired estimate linear in
\(\|\xi_{v,w}\|_{\rm prep}\).  Together with
 \eqref{eq:variational-coarse-energy}, this yields the
required pointwise coarse bound
\begin{equation}\label{eq:variational-coarse-H1}
 \|K_{v,w}(\tau)\|_{H^1_\nu}
 \leq Ce^{-q_0(\tau-\tau_0)}
 \|\xi_{v,w}\|_{\rm prep},
 \qquad \tau\geq\tau_0+1.
\end{equation}
Now all hypotheses of
Lemma~\ref{lem:linearized-stable-rate-bootstrap} are explicit, and
that lemma gives \eqref{eq:general-profile-variation-rate}.

Before specializing the variation, we establish the required
absolute-time one-state estimate without using any later
profile-realization argument.  The first-unit estimate
\eqref{eq:late-entry-first-unit-one-state}, followed by the concrete
\(j=1\) graph recovery at \(\tau_0+1\), gives
\[
 \|H_v(\tau_0+1)\|_{H^1_\nu}
 \leq C_r e^{-\gamma_1\tau_0}
      +C_r e^{-ce^{\tau_0/2}}.
\]
The entire late-entry disk consists of strict entrances for the same
fixed rate-compatible package.  Therefore the already proved estimate
\eqref{eq:master-bootstrap-H1} supplies the uniform first tail
pass
\[
 \|H_v(\tau)\|_{H^1_\nu}
 \leq C_r e^{-\theta\tau}+C_r e^{-ce^\tau},
 \qquad |v|\leq r,\quad \tau\geq\tau_0+1.
\]
Thus the absolute iteration below starts with \(q=\theta\); the
endpoint estimate above is used separately for the normalized
homogeneous contribution.
Starting at \(\tau_0+1\), apply the stable Duhamel estimate, the
unit-window energy inequality, and
\eqref{eq:projected-stable-graph-recovery}.  Given a tail pass with
uniform absolute input \(Ce^{-q\tau}\),
\eqref{eq:stable-tame-bound} gives the next absolute exponent
\[
 q\longmapsto
 \min\{\gamma_1-\zeta,\ q+\kappa,\ 2q\}.
\]
The homogeneous term is bounded by \(C_r e^{-\gamma_1\tau}\), and
every outer tail is bounded by
\(C_re^{-ce^{\tau_0/2}}\).  Consequently every constant in this
finite tail iteration is uniform for all sufficiently large
\(\tau_0\).  It terminates at
\begin{equation}\label{eq:late-disk-one-state-absolute-rate}
 \|H_v(\tau)\|_{H^1_\nu}
 \leq C_{r,\zeta}e^{-(\gamma_1-\zeta)\tau},
 \qquad |v|\leq r,\quad \tau\geq\tau_0+1 .
\end{equation}

For the uniform late-entry special disk, the first unit is supplied
by Lemma~\ref{lem:late-entry-first-unit-stable-variation}, rather than
by the late-time tame estimate.  Combining
\eqref{eq:late-entry-first-unit-variation} with the special-disk graph
recovery at the endpoint gives
\[
 \|K_{v,w}(\tau_0+1)\|_{H^1_\nu}
 \leq C_r e^{-\gamma_1\tau_0}|w|
      +C_r e^{-ce^{\tau_0/2}}|w|.
\]
On \([\tau_0+1,\infty)\), the bounds
\eqref{eq:late-entry-differentiated-tame-bound},
\eqref{eq:late-entry-profile-graph-error}, and
\eqref{eq:late-disk-one-state-absolute-rate} satisfy the hypotheses
of the finite linearized rate bootstrap with constants independent of
every sufficiently large entrance time.  Starting from the displayed
endpoint estimate proves \eqref{eq:profile-variation-rate}.

For a general tangent vector, differentiate
\eqref{eq:V-infinity-formula} first on
\([\tau_0,\tau_0+1]\) by the finite-horizon \(C^1\) theorem, and then
on \([\tau_0+1,S]\).  Passing \(S\to\infty\) after the tail estimate
below gives the explicitly split formula
\eqref{eq:general-profile-map-derivative}.  On the prepared
\(W\)-profile disk,
\(\Pi_1W(w)=w\), so the initial term reduces to \(w\), which proves
\eqref{eq:profile-map-derivative}.

To justify differentiation of the tail integral uniformly on a
prepared ball, take \(S\geq\tau_0+1\) and combine
\eqref{eq:differentiated-tame-bound},
\eqref{eq:general-profile-variation-rate}, and the uniform
one-solution near-gap estimate.  Choose
\(0<\zeta<\min\{\kappa,\gamma_1/2\}\).  After multiplication by
\(e^{\gamma_1s}\), the first, second, and outer terms are bounded by
\[
 C e^{-(\kappa-\zeta)s}\|\xi_{v,w}\|_{\rm prep},\qquad
 C e^{-(\gamma_1-2\zeta)s}\|\xi_{v,w}\|_{\rm prep},\qquad
 C e^{-ce^{s/2}}\|\xi_{v,w}\|_{\rm prep},
\]
with a constant depending only on the common-margin ball and the
prepared-chart tangent bound.  Here \(\tau_0\) is fixed, so this
constant need not be uniform in the entrance time.  Their tails tend
to zero uniformly on that fixed ball.  This proves
\eqref{eq:uniform-profile-derivative-tail}, hence the \(C^1\)
assertion by the Banach-space uniform-derivative principle used in
the proof of Theorem~\ref{thm:global-two-state-estimate}.

For the late-entry special disk, the first integral in
\eqref{eq:profile-map-derivative} is bounded by
\eqref{eq:late-entry-first-unit-projected-source}.  On the second
integral use
\eqref{eq:profile-variation-rate}, the one-state estimate
\eqref{eq:late-disk-one-state-absolute-rate}, and the hybrid outer
bound \eqref{eq:late-entry-profile-outer-source}.  After multiplication
by \(e^{\gamma_1s}\), the three corresponding integrands are bounded
by
\[
 C_{r,\zeta}e^{-(\kappa-\zeta)s}|w|,\qquad
 C_{r,\zeta}e^{-(\gamma_1-2\zeta)s}|w|,\qquad
 C C_W e^{-ce^{s/2}}|w|.
\]
Integrating the tail from \(\tau_0+1\), and adding the first-unit
bound, gives \eqref{eq:profile-map-C1-error} with constants independent
of all sufficiently large entrance times.
\end{proof}

\begin{lemma}[Initial projection and first-unit stable source]
\label{lem:first-unit-stable-source}
Let two sliced prepared solutions belong to one fixed-\(\tau_0\)
common-margin ball, put \(D=H_1-H_2\), and let \(d_0\) be their
prepared entrance distance.  Then
\begin{equation}\label{eq:first-unit-stable-source}
 e^{\gamma_1\tau_0}\|\Pi_1D(\tau_0)\|_{E_1}
 +\int_{\tau_0}^{\tau_0+1}
   e^{\gamma_1s}
   \left|\Pi_1(\mathcal N_1-\mathcal N_2)(s)\right|\,ds
 \leq C_{\tau_0}d_0 .
\end{equation}
For a sliced tangent \(\xi\), \(K=DH[\xi]\), the corresponding
linearized estimate is
\begin{equation}\label{eq:first-unit-stable-source-variation}
 e^{\gamma_1\tau_0}\|\Pi_1K(\tau_0)\|_{E_1}
 +\int_{\tau_0}^{\tau_0+1}
   e^{\gamma_1s}
   \left|\Pi_1D\mathcal N[\xi](s)\right|\,ds
 \leq C_{\tau_0}\|\xi\|_{\rm prep}.
\end{equation}
\end{lemma}

\begin{proof}
The finite-rank projection \(\Pi_1\) is bounded on \(L^2_\nu\) and,
after pairing with the smooth eigenbasis of \(E_1\), on
\(H^{-1}_\nu\).  Hence the initial term is bounded by the prepared
entrance distance.  On the compact interval
\([\tau_0,\tau_0+1]\), exact polarization of the stable equation and
the prepared chart calculus give
\[
 \|\mathcal N_1-\mathcal N_2\|_{H^{-1}_\nu}
 \leq C_{\tau_0}\left(
  \|D\|_{H^1_\nu}+\mathfrak D_k^{\rm hyb}\right).
\]
Proposition~\ref{prop:two-state-prepared-evolution} controls the
hybrid supremum and the \(L^2_\tau H^1_\nu\) norm by \(C_{\tau_0}d_0\).
Cauchy--Schwarz on this one unit interval proves
\eqref{eq:first-unit-stable-source}.  Differentiate the exact
finite-horizon system and use its all-block Fr\'echet estimate to
obtain \eqref{eq:first-unit-stable-source-variation}.
\end{proof}

\subsection{Nonlinear scattering}

\begin{theorem}[\(C^1\) marked, gauge-fixed nonlinear scattering]
\label{thm:global-two-state-scattering}
Fix \(k\geq12\), \(0<\alpha<1\), one normalized entrance time
\(\tau_0\), and first one rate pair
\[
0<\sigma<\theta<\beta,
\]
then one rate-compatible prepared package and marked preparation
convention, and finally one entrance size
\[
 0<\varepsilon\leq\varepsilon_{\rm ent}.
\]
Let
\[
 \mathscr B\subset\Sigma_{\tau_0}^{k+2,\alpha}
\]
be a common-margin sliced ball all of whose points are strict prepared
entrances of continuation order \(k\), with the parameter class and
prepared package common on the ball.  Assume, exactly as in
Theorem~\ref{thm:global-two-state-estimate}, that one buffered physical
cover and the strict time-width margin
\eqref{eq:auxiliary-buffered-time-width} work throughout
\(\mathscr B\).  Then the map
\begin{equation}\label{eq:full-scattering-map}
 \begin{aligned}
 \mathscr S_{\rm sl}:\mathbf z_0&\longmapsto
 \bigl(
  T,\log\lambda_\infty,\Psi_\infty,\\[-2pt]
 &\hspace{37mm}
  V_\infty
 \bigr).
 \end{aligned}
\end{equation}
is a \(C^1\), marked, gauge-fixed scattering package.  Precisely, its
finite-dimensional components are \(C^1\), and for every
\(K\Subset M\) and \(m\geq0\), the phase component is \(C^1\) as a map
into a \(C^m(K)\) exponential chart.  Equivalently, for each fixed
pair \((K,m)\) the displayed tuple is \(C^1\) into the corresponding
finite product Banach chart; no single Fr\'echet topology on the full
tuple is asserted.  In particular, if
\(\mathbf c^{(2)}(V_\infty)\) is the quadratic response coefficient
from \eqref{eq:quadratic-feedback-coefficient}, then
\(\mathbf c^{(2)}\circ V_\infty\) is a derived \(C^1\) response law;
it is not an additional independent scattering coordinate.  Likewise,
\begin{equation}\label{eq:physical-amplitude-C1-scattering-map}
 \mathbf z_0\longmapsto
 \mathfrak A_1(\mathbf z_0)
 =
 \exp\!\bigl(-\gamma_1\log\lambda_\infty(\mathbf z_0)\bigr)
 V_\infty(\mathbf z_0)
\end{equation}
is a derived \(C^1\) map to \(E_1\), and
\(\mathbf c^{(2)}\circ\mathfrak A_1\) is its derived physical-time
quadratic response law.  Both have the same buffered ambient extension
as the scattering package.  They retain the marked, gauge-fixed scope
of Definition~\ref{def:physical-first-amplitude}; only
\(\mathfrak A_1\), not the entire scattering tuple, has the transported
restart invariance of
Lemma~\ref{lem:physical-amplitude-restart-covariance}.  If
\[
 d_0=
 \|\mathbf z_{1,0}-\mathbf z_{2,0}\|
 _{\mathscr X_{\rm prep}^{k+2,\alpha}},
\]
then, for every \(K\Subset M\),
\begin{equation}\label{eq:scattering-data-Lipschitz}
 \begin{split}
  &|T_1-T_2|
  +\left|\log\frac{\lambda_{\infty,1}}
                         {\lambda_{\infty,2}}\right|
  +\|V_{\infty,1}-V_{\infty,2}\|_{E_1}\\
  &\quad+
  \left\|
   \log\bigl(\Psi_{\infty,2}^{-1}
             \circ\Psi_{\infty,1}\bigr)
  \right\|_{C^m(K)}
  \leq C_{K,m}d_0.
 \end{split}
\end{equation}
The stable scattering tails satisfy
\begin{equation}\label{eq:two-state-stable-scattering-tail}
 \int_{\tau_0}^{\infty}
 e^{\gamma_1s}
 \left|\Pi_1\bigl(\mathcal N_1-\mathcal N_2\bigr)(s)\right|\,ds
 \leq Cd_0.
\end{equation}
The \(V_\infty\) component and its fibers have the marked scope stated
in Remark~\ref{rem:stable-profile-gauge-scope}.

If
\[
 \mathscr O\subset\operatorname{dom}
    \Pi_{\rm sl}^{\,k+4\to k+2}
    \subset\mathscr P_{\tau_0}^{k+4,\alpha},
 \qquad
 \Pi_{\rm sl}^{\,k+4\to k+2}(\mathscr O)\subset\mathscr B,
\]
is a buffered ambient prepared neighborhood, define
\[
 \mathscr S_{\rm amb}
 :=\mathscr S_{\rm sl}\circ
   \Pi_{\rm sl}^{\,k+4\to k+2}.
\]
Then \(\mathscr S_{\rm amb}\) is \(C^1\) in the same componentwise
sense, with
\[
 D\mathscr S_{\rm amb}
 =D\mathscr S_{\rm sl}\circ
   D\Pi_{\rm sl}^{\,k+4\to k+2}.
\]
For every fixed component target chart, let \(L_{\mathscr S}\) be the
uniform sliced Lipschitz and first-variation constant on
\(\mathscr B'\Subset_{\rm u}\mathscr B\) supplied by the estimates
above.  If
\[
 \overline{\mathscr O}\subset\operatorname{dom}
 \Pi_{\rm sl}^{\,k+4\to k+2},
 \qquad
 \Pi_{\rm sl}^{\,k+4\to k+2}(\overline{\mathscr O})
 \subset\mathscr B',
\]
then
\begin{equation}\label{eq:quantitative-ambient-scattering-bound}
 \operatorname{Lip}(\mathscr S_{\rm amb}|_{\mathscr O})
 +\sup_{\mathbf z\in\mathscr O}
   \|D\mathscr S_{\rm amb}(\mathbf z)\|
 \leq2K_{\Pi,k+2}L_{\mathscr S}.
\end{equation}
Thus ambient uniformity comes from the quantitative phase retraction
and the sliced two-state/first-variation bounds; the closure and
interiority conditions retain their domains and margins.
\end{theorem}

\begin{proof}
Theorem~\ref{thm:global-two-state-estimate} gives the geometric
components of \eqref{eq:scattering-data-Lipschitz}.  We first verify
the pointwise \(H^1_\nu\) input for the stable difference rather than
deducing it from an integrated estimate.  Fix \(0<q_0<\theta\).
Equation~\eqref{eq:global-two-state-energy} gives
\[
 \|D(\tau)\|_{L^2_\nu}
 +\left(\int_\tau^\infty\|D(s)\|_{H^1_\nu}^2\,ds\right)^{1/2}
 \leq Cd_0e^{-q_0(\tau-\tau_0)}.
\]
Apply \eqref{eq:projected-stable-graph-recovery} with \(j=1\) to
 \(\partial_\tau D+LD=\mathcal N_1-\mathcal N_2\).  Its lower-block
 source is zero, while the concrete outer term is
 \(Cd_0^2e^{-ce^{\tau/2}}\).  The last unit interval of the displayed
 \(L^2_\nu\) estimate therefore gives
\[
 \|D(\tau)\|_{H^1_\nu}
 \leq Cd_0e^{-q_0(\tau-\tau_0)}.
\]
This conclusion is used only for \(\tau\geq\tau_0+1\); the first unit
interval is controlled separately by
Lemma~\ref{lem:first-unit-stable-source}.
Starting from this bound, iterate the Duhamel and graph-recovery
argument with \eqref{eq:two-state-stable-Hminus1}.  For every
\[
 0<\zeta<\min\{\kappa,\gamma_1/2\}
\]
the iteration terminates at
\[
 \|H_1(\tau)-H_2(\tau)\|_{H^1_\nu}
 \leq C_\zeta d_0
 e^{-(\gamma_1-\zeta)(\tau-\tau_0)}.
\]
Combining this estimate, the one-solution near-gap estimate, and the
explicit difference bound
\eqref{eq:two-state-stable-Hminus1} makes
\[
 e^{\gamma_1\tau}
 \|\mathcal N_1(\tau)-\mathcal N_2(\tau)\|_{H^{-1}_\nu}
\]
integrable on \([\tau_0+1,\infty)\), with integral bounded by
\(Cd_0\): the three exponents
left after multiplication by \(e^{\gamma_1\tau}\) are respectively
\(\kappa-\zeta\), \(\gamma_1-2\zeta\), and a
Gaussian-superexponential rate.  The first-unit integral is exactly
the second term in
\eqref{eq:first-unit-stable-source}; hence the two pieces prove
\eqref{eq:two-state-stable-scattering-tail}.  Finally, subtraction of
\eqref{eq:V-infinity-formula} gives
\[
 V_{\infty,1}-V_{\infty,2}
 =e^{\gamma_1\tau_0}\Pi_1D(\tau_0)
  +\int_{\tau_0}^{\infty}
    e^{\gamma_1s}\Pi_1(\mathcal N_1-\mathcal N_2)(s)\,ds .
\]
The first term in \eqref{eq:first-unit-stable-source} accounts for the
initial projection, while its second term accounts for the interval
\([\tau_0,\tau_0+1]\).  This proves the \(V_\infty\)-term in
\eqref{eq:scattering-data-Lipschitz}.

It remains to prove \(C^1\), rather than merely Lipschitz, dependence
of the limiting data.  At a finite endpoint \(S\), all components are
\(C^1\) by
Proposition~\ref{prop:two-state-prepared-evolution}.  The derivative
tails for \(T,\log\lambda_\infty,\Psi_\infty\) converge uniformly by
\eqref{eq:C1-geometric-tail-uniform}.  For the stable profile,
apply Lemma~\ref{lem:differentiated-stable-normal-form} in a Banach
chart on the split submanifold
\(\Sigma_{\tau_0}^{k+2,\alpha}\).  For the buffered ambient extension,
compose this chart with
\(\Pi_{\rm sl}^{\,k+4\to k+2}\).  The quantitative bounds
\eqref{eq:quantitative-phase-retraction-Lipschitz}--%
\eqref{eq:quantitative-phase-retraction-derivative} multiply every
sliced derivative-tail estimate by at most \(K_{\Pi,k+2}\).  Hence
uniform Cauchy convergence is preserved after composition.  Then
\eqref{eq:general-profile-map-derivative} and
\eqref{eq:uniform-profile-derivative-tail} show that the derivatives
of
\[
 e^{\gamma_1\tau_0}\Pi_1H(\tau_0)
 +\int_{\tau_0}^{S}
 e^{\gamma_1s}\Pi_1\mathcal N(s)\,ds
\]
are uniformly Cauchy on the common-margin ball.  More explicitly, the
derivative of the initial projection and the integral over
\([\tau_0,\tau_0+1]\) are controlled by
\eqref{eq:first-unit-stable-source-variation}; the uniform tail
argument is applied only on \([\tau_0+1,S]\).  Hence the limiting
map \(V_\infty\) is \(C^1\).  The coefficient
\(\mathbf c^{(2)}\) is a fixed finite-dimensional quadratic map of
\(V_\infty\).  Applying the Banach-space uniform-derivative principle
componentwise proves the stated \(C^1\) regularity of the scattering
package in the precise fixed-\((K,m)\) sense above.  The uniform
derivative-tail estimate justifies passage from finite horizons to
\(S=\infty\).
Finally, \(\lambda_\infty>0\), so the ordinary finite-dimensional
product and chain rules applied to
\eqref{eq:physical-amplitude-C1-scattering-map} prove the asserted
\(C^1\) inheritance for \(\mathfrak A_1\) and its quadratic response.
\end{proof}

\begin{remark}[Possible vanishing]
\label{rem:first-profile-resonances}
The theorem does not assert $V_\infty\ne0$.  If it vanishes, the
displayed estimate gives faster decay and all quadratic coefficients
below vanish; the first nonzero term need not automatically lie in
$E_2$ because a nonlinear resonance may intervene.

 Even for the exact-core center, $h(\tau_0)=0$ does not by itself imply
 $V_\infty=0$.  Formula~\eqref{eq:V-infinity-formula} retains the stable
 Duhamel memory of the later moving graft and cutoff forcing.  Those
 sources are Gaussian-superexponentially small at the time when they
 act, but a stable semigroup transports their projection with its
 ordinary spectral rate.  Thus no superexponential center-flow
conclusion is used here.
\end{remark}

\subsection{Profile realization and foliation}

\begin{theorem}[Marked first-profile foliation and strong-stable set]
\label{thm:profile-realization}
Equip $E_1$ with its $L^2_\nu$ norm.
For the host and implantation data of
Theorem~\ref{thm:prepared-open-basin}, the exact-core entrance time may
be chosen so late that there are $r>0$, a relative
$C^{2,\alpha}$-open basin $\mathscr U_{\rm prof}$ containing the
exact-core center $G_*$, and a \(C^1\) disk of physical metrics whose
fixed prepared lifts are strict entrances,
\[
 \overline B_r^{E_1}(0)\longrightarrow\mathscr U_{\rm prof},
 \qquad v\longmapsto G_v,
 \qquad G_0=G_*,
\]
such that every $G_v$ agrees with $G_*$ outside the implantation
region.  The basin \(\mathscr U_{\rm prof}\) satisfies every conclusion
of Theorem~A and may be taken as the existential neighborhood
\(\mathscr U\) in that theorem.  Write
\(\mathcal V_{\tau_0}(v)=V_\infty(G_v)\).  Then
\begin{equation}\label{eq:profile-realization-ball}
 B_{r/2}^{E_1}(0)
 \subset
 \mathcal V_{\tau_0}\bigl(B_r^{E_1}(0)\bigr),
 \qquad
 \{G_v:|v|\leq r\}\subset\mathscr U_{\rm prof}.
\end{equation}
More quantitatively, for every
\[
 0<\zeta<\min\{\kappa,\gamma_1/2\}
\]
one may arrange
\begin{equation}\label{eq:profile-map-degree-error}
 \sup_{|v|\leq r}|\mathcal V_{\tau_0}(v)-v|
 \leq
 C_{r,\zeta}\left(
  e^{-(\kappa-\zeta)\tau_0}
  +e^{-(\gamma_1-2\zeta)\tau_0}
  +e^{-ce^{\tau_0/2}}
 \right).
\end{equation}
Consequently every sufficiently small first stable profile occurs.
In addition,
\begin{equation}\label{eq:profile-map-degree-C1-error}
 \sup_{|v|\leq r}
 \|D\mathcal V_{\tau_0}(v)-\operatorname{Id}_{E_1}\|
 \leq
 C_{r,\zeta}\left(
  e^{-(\kappa-\zeta)\tau_0}
  +e^{-(\gamma_1-2\zeta)\tau_0}
  +e^{-ce^{\tau_0/2}}
 \right).
\end{equation}
Thus, after taking \(\tau_0\) later,
\[
 \mathcal V_{\tau_0}:
 \mathcal V_{\tau_0}^{-1}\bigl(B_{r/2}^{E_1}(0)\bigr)
 \longrightarrow B_{r/2}^{E_1}(0)
\]
 is a \(C^1\) bi-Lipschitz scattering chart: every profile in that ball
 has a unique preimage on the chosen transverse disk.

Since \(\mathscr U_{\rm prof}\) is relatively open in the smooth locus,
fix a \(C^{2,\alpha}\)-open ambient metric thickening
\[
 \widehat{\mathscr U}_{\rm prof}
 \subset\operatorname{Met}^{2,\alpha}(\widehat X),
 \qquad
 \widehat{\mathscr U}_{\rm prof}\cap
 \operatorname{Met}^{\infty}(\widehat X)
 =\mathscr U_{\rm prof}.
\]
More generally, fix for the entire disk one marked host, graft, scale,
and initial-map convention.  Let
\(\mathscr U_{\rm phys}^{k+6,\alpha}\) be a sufficiently small
\(C^{k+6,\alpha}\) physical-metric neighborhood of
\(\{G_v:|v|\leq r\}\), and let
\[
 \mathcal P_{\tau_0}:
 \mathscr U_{\rm phys}^{k+6,\alpha}\ni G
 \longmapsto\mathbf z_{\tau_0}(G)
 \in\Sigma_{\tau_0}^{k+4,\alpha}
 \hookrightarrow\Sigma_{\tau_0}^{k+2,\alpha}
\]
  be the fixed high-topology preparation map supplied, at this fixed
  entrance time and at raw order \(k+6\) and output order \(k+4\), by
  Corollary~\ref{cor:fixed-convention-preparation-map} on the compact
  already-sliced disk, and let \(\pi_G\)
  denote the closed-metric component of
a prepared tuple.  There are \(r/2<r'<r\) and a prepared open
neighborhood
\[
 \widetilde{\mathscr O}_{\rm scat}
 \subset\mathscr P_{\tau_0}^{k+4,\alpha}
\]
 of
 \(\{\mathcal P_{\tau_0}(G_v):|v|\leq r'\}\), with
 \(\pi_G(\widetilde{\mathscr O}_{\rm scat})
 \subset\widehat{\mathscr U}_{\rm prof}\).  Consequently every smooth
 closed-metric component in this prepared neighborhood belongs to
 \(\mathscr U_{\rm prof}\).  Moreover, there is one common-margin
sliced strict-entrance ball
\[
 \mathscr B_{\rm scat}\subset\Sigma_{\tau_0}^{k+2,\alpha}
\]
on which the fixed exact-core buffered physical cover and the strict
time-width margin \eqref{eq:auxiliary-buffered-time-width} work
uniformly, and a uniformly interior common-margin subball
\[
 \mathscr B'_{\rm scat}\Subset_{\rm u}\mathscr B_{\rm scat}
\]
such that
\begin{equation}\label{eq:profile-scattering-domain-inclusions}
 \overline{\widetilde{\mathscr O}_{\rm scat}}
 \subset\operatorname{dom}
   \Pi_{\rm sl}^{\,k+4\to k+2},
 \qquad
 \Pi_{\rm sl}^{\,k+4\to k+2}
   (\overline{\widetilde{\mathscr O}_{\rm scat}})
 \subset\mathscr B'_{\rm scat}
 \Subset_{\rm u}\mathscr B_{\rm scat}.
\end{equation}
If \(V_\infty^{\rm sl}\) denotes the profile component of the sliced
scattering map, then
\begin{equation}\label{eq:first-profile-scattering-map}
 \mathscr S^{\rm prep}_1:
 \widetilde{\mathscr O}_{\rm scat}\longrightarrow E_1,
 \qquad
 \mathscr S^{\rm prep}_1(\mathbf z)
 =V_\infty^{\rm sl}\!\left(
   \Pi_{\rm sl}^{\,k+4\to k+2}(\mathbf z)
  \right),
\end{equation}
is a \(C^1\) split submersion in the buffered prepared topology.  Its
level sets are \(C^1\) Banach submanifolds of codimension
\(\dim E_1\), and the connected components of these fibers form a
local \(C^1\) foliation.  In particular,
\begin{equation}\label{eq:first-strong-stable-leaf}
 \mathscr W^{ss}_1=(\mathscr S^{\rm prep}_1)^{-1}(0)
\end{equation}
is a nonempty strong-stable submanifold whose connected components are
leaves.  The condition \(V_\infty=0\) is preserved under every forward
prepared restart.  On this zero fiber the improved stable rate in
Theorem~\ref{thm:first-stable-profile} holds.  The complement of
\(\mathscr W^{ss}_1\) is open and dense in
\(\widetilde{\mathscr O}_{\rm scat}\).  All profile fibers and the
open-dense assertion in this paragraph are formed in the single gauge
and transported preparation convention fixed above; their scope is
exactly Remark~\ref{rem:stable-profile-gauge-scope}.
\end{theorem}

\begin{proof}
Use the lift \(W\) from Lemma~\ref{lem:compact-profile-lift}, and let
\(K\Subset M\) be the common compact support of its image.
Fix once and for all
\[
 0<\zeta_*<\min\{\kappa,\gamma_1/2\}
\]
and a radius \(r>0\) small enough for the fixed linear disk
\(\overline B_r^{E_1}(0)\) to lie in the chosen prepared and
positive-metric charts.  This choice of \(r\) is independent of the
entrance time and will not be changed below.

For every sufficiently large candidate entrance time \(\tau_0\), one
has \(\rho_{\tau_0}=1\) on \(K\), and the inward soliton transport gives
\begin{equation}\label{eq:profile-support-eta-one}
 \Phi_0^{-1}(K)
 \Subset\operatorname{int}\{\eta=1\}.
\end{equation}
Increase the candidate threshold, uniformly for \(|v|\leq r\), so
that this physical support is also disjoint from
\(\bigcup_a\widetilde W_a^{\rm har}\).  For each such candidate time set
\[
 u_v=e^{-\gamma_1\tau_0}W(v)
\]
and, on the marked core, define
\begin{equation}\label{eq:physical-profile-disk}
 G_v
 =
 G_*+
 \iota^*\bigl(\lambda_0\Phi_0^*u_v\bigr),
\end{equation}
extending the added compactly supported tensor by zero.  For
$|v|\leq r$ and $\tau_0$ large this is a smooth positive metric and,
by \eqref{eq:profile-support-eta-one}, on the support of the
perturbation one has the exact, uncut identity
\[
 \acute G_{0,v}
 =S_0+\lambda_0\Phi_0^*u_v.
\]
Thus the normalized prepared perturbation is exactly
\[
 h_{0,v}
 =\lambda_0^{-1}(\Phi_0^{-1})^*\acute G_{0,v}-\bar g
 =u_v,
 \qquad
 H_v(\tau_0)=u_v.
\]
All other prepared components in this construction are independent of
\(v\).  For each fixed \(\tau_0\), the pullback and extension maps make
\(v\mapsto\mathbf z_v\) \(C^1\) in the full prepared Banach topology.
The corresponding full unscaled chart constant is allowed to depend
on \(\tau_0\), since high derivatives of
\(\Phi_0=\varphi_{\tau_0}\) need not be uniformly bounded in fixed
coordinates on \(\mathcal X\).

The support of the perturbation is disjoint from the enlarged
buffered noncollapsing physical sets
\(\bigcup_{a=1}^{N_{\rm ext}}\widetilde W_a^{\rm har}\), and from the
graft
region.
Hence the physical, marking, scale, and prepared-map blocks in
\(\mathfrak D_{m_\#}^{\rm hyb}\) have zero \(v\)-variation.  In fact
the prepared map is globally independent of \(v\), so the auxiliary
order-\((m_\#+2)\) \(R\)-trace, the additional source-star trace
\(d_{{\rm Fgr},m_\#+1,0}^{++}\) and the global low initial trace
\(d_{{\rm F},m_\#+1,0}^{\rm glob}\) have zero \(v\)-variation, whereas
\[
 H_v(\tau_0)=h_{0,v}=e^{-\gamma_1\tau_0}W(v).
\]
The fixed compact support of \(W(E_1)\) therefore gives the genuine
uniform late-entry estimates
\begin{equation}\label{eq:constructed-profile-disk-absolute-variation}
 \mathfrak d_{{\rm le},\tau_0}(\mathbf z_v,\mathbf z_0)
 \leq C|v|,
 \qquad
 \|D_v\mathbf z_v[w]\|_{{\rm le},\tau_0}
 \leq C|w|,
\end{equation}
with \(C\) independent of every sufficiently large \(\tau_0\).
The buffered physical cover, graft data, scale, marking, and prepared
map blocks are the fixed ones of the exact-core construction.  On the
only varying block, finite-dimensional norm equivalence for the
compactly supported space \(W(E_1)\) bounds every normalized tensor
norm through order \(m_\#+1=5\) uniformly.  Thus the entrance
coefficient,
inverse-map, radial-comparison, and composition package has one
constant \(K_\#\) for the whole late-entry family; its future
propagation is the low-initial-data clause
\eqref{eq:uniform-restart-low-initial-propagation} of
Lemma~\ref{lem:uniform-weighted-Schauder-restart}.  The physical
exterior-memory and compact graft-buffer variations are then supplied
by Lemma~\ref{lem:late-entry-hybrid-first-variation}; no full prepared
tangent norm enters this cascade.
All nine moments vanish by
\eqref{eq:compact-profile-lift-properties}, because
\(\rho_{\tau_0}=1\) on the common support of \(W(E_1)\).
The direct effective columns are independent of \(v\), since every
prepared map and scale block is fixed.  In the exact column formulae
\eqref{eq:prepared-full-column-zero}--%
\eqref{eq:prepared-full-column-j}, the remaining terms are linear in
\(h_{0,v}\).  One Gaussian integration by parts, exactly as in the
proof of Lemma~\ref{lem:buffered-full-prepared-columns}, and
finite-dimensional norm equivalence on \(W(E_1)\) therefore give,
uniformly on the disk,
\[
 M_{\mu j}(v)
 =\ip{Y_j}{Z_\mu}
 {}+O(e^{-\gamma_1\tau_0})+O(e^{-ce^{\tau_0}}).
\]
Thus the adaptive Gram matrix is uniformly invertible for all
sufficiently large \(\tau_0\).  These tuples use the already-sliced
route~(b) of Definition~\ref{def:strict-prepared-entrance}, with
\(p=0\); no application of the full-norm implicit-function
neighborhood is needed.  Since the only varying prepared component is
the fixed compactly supported smooth tensor
\(e^{-\gamma_1\tau_0}W(v)\), finite-dimensional norm equivalence at
order \(k+2\) also puts the entire disk in one common-margin ball of
\(\mathscr P_{\tau_0}^{k+2,\alpha}\).  This verifies the
finite-regularity membership in route~(b), rather than inferring it
from smoothness alone.

The entrance norms satisfy, uniformly for $|v|\leq r$,
\begin{equation}\label{eq:profile-disk-entrance-size}
 \|H_v(\tau_0)\|_{L^2_\nu}
 +\|h_{0,v}\|_{C^3}
 \leq C_r e^{-\gamma_1\tau_0}.
\end{equation}
Since $\sigma<\theta<\beta<\gamma_1$, these bounds lie strictly below
the $e^{-\theta\tau_0}$ and $e^{-\sigma\tau_0}$ entrance faces when
$\tau_0$ is large.  On the perturbed core, finite-dimensional norm
 equivalence on the fixed smooth space \(W(E_1)\), now through order
\(\max\{12,k+2\}\), makes the scale-normalized metric-jet variation
\(O(e^{-\gamma_1\tau_0})\).  Smooth dependence of curvature and of
the Ricci defect on these jets preserves, with uniform strict slack,
the order-ten curvature and Ricci-defect margins in item~(5) of
Definition~\ref{def:strict-prepared-entrance}.  The \(C^2\) part of
the same estimate preserves ellipticity.  The exact-core bounds
supply uniform coefficient and domain reserves for the background
harmonic charts.  Compact support and finite-dimensional norm
equivalence give a \emph{global} scale-one prepared distance
\(O(e^{-\gamma_1\tau_0})\) from the exact state.
The reserve-to-operative clause of
Lemma~\ref{lem:prepared-harmonic-radius-lower-stability} may therefore
be applied simultaneously at every center with the exact-core
operative and \(+\)-triples once \(\tau_0\) is large.  It preserves the
same operative triple, its fixed modulus-compatibility inequality, and
a new common reserve triple for the whole disk; no claim is made that
a witness ball centered off the tensor support is disjoint from that
support.  The perturbation is disjoint from the graft and from every
\(\widetilde W_a^{\rm har}\), so the entire physical
certificate---including all three
witness tiers, the common coefficient ball, the quarter-modulus
future-window estimate, and the \(\mu_{\rm RF}\)-width slack---is
unchanged.  The remaining graft and exterior margins persist as well.
All these margins are scale-normalized or buffered geometric
conditions.  More explicitly, after division by the
\(e^{-\theta\tau_0}\) and \(e^{-\sigma\tau_0}\) entrance faces, the
left sides tend uniformly to zero and hence retain a fixed fractional
slack.  None of these conditions asks for a uniform high-order
unscaled chart norm on the collapsing core.  The whole disk is
therefore a common-margin family of strict prepared entrances in the
sense of Definition~\ref{def:strict-prepared-entrance}.
Together with \eqref{eq:abstract-W-profile-disk},
\eqref{eq:constructed-profile-disk-absolute-variation}, and the
uniform \(K_\#\) package above, the fixed prepared lifts of
\(v\mapsto G_v\) consequently form a uniform late-entry prepared
\(W\)-profile disk in the sense of
Definition~\ref{def:prepared-W-profile-disk}.  The notation
\(\mathcal V_{\tau_0}(v)\) used there agrees with
\(V_\infty(G_v)\) here.

The absolute-time bootstrap has already been carried out in
\eqref{eq:late-disk-one-state-absolute-rate}.  Thus
\begin{equation}\label{eq:profile-disk-near-gap}
 \|H_v(\tau)\|_{H^1_\nu}
 \leq C_{r,\zeta}e^{-(\gamma_1-\zeta)\tau}
\end{equation}
for every $|v|\leq r$, with constants independent of sufficiently
large $\tau_0$.

Formula~\eqref{eq:V-infinity-formula} and
\eqref{eq:compact-profile-lift-properties} now give
\[
 \mathcal V_{\tau_0}(v)
 =
 v+
 \int_{\tau_0}^{\tau_0+1}
 e^{\gamma_1s}\Pi_1\mathcal N_v(s)\,ds
 +\int_{\tau_0+1}^{\infty}
 e^{\gamma_1s}\Pi_1\mathcal N_v(s)\,ds.
\]
The first integral is bounded by
\eqref{eq:late-entry-first-unit-one-state-projected-source}.  On the
second integral, where \(\tau\geq\tau_0+1\), insert
\eqref{eq:profile-disk-near-gap} into the tame estimate
\eqref{eq:stable-tame-bound}.  The three tail integrals are bounded
respectively by
\[
 C_{r,\zeta} e^{-(\kappa-\zeta)\tau_0},\qquad
 C_{r,\zeta} e^{-(\gamma_1-2\zeta)\tau_0},\qquad
 C_{r,\zeta} e^{-ce^{\tau_0/2}},
\]
and the extra first-unit term
\(C_re^{-\gamma_1\tau_0}\) is absorbed by the second displayed
quantity.  This proves \eqref{eq:profile-map-degree-error}.  The map
\(\mathcal V_{\tau_0}\) is \(C^1\) by
Lemma~\ref{lem:differentiated-stable-normal-form}, whose uniform
special-disk clause uses
Lemma~\ref{lem:late-entry-hybrid-first-variation}, and
\eqref{eq:profile-map-C1-error} is exactly
\eqref{eq:profile-map-degree-C1-error}.

We now choose the entrance time once, with no subsequent change of
\(r\) or \(\tau_0\).  Let \(s_{\rm Gr}>0\) be the least singular value
of the limiting Gram matrix, let \(\delta_{\rm prep}>0\) be a common
strict-margin radius in the order-\((k+2,\alpha)\) prepared chart.
Denote by
\(\mathcal E_0(\tau)\) and \(\mathcal E_1(\tau)\) the right sides of
\eqref{eq:profile-map-degree-error} and
\eqref{eq:profile-map-degree-C1-error}, respectively, evaluated with
the fixed pair \((r,\zeta_*)\).  Choose one \(\tau_0\) so large that
\eqref{eq:profile-support-eta-one} and the buffered support separation
hold and, simultaneously,
\begin{equation}\label{eq:simultaneous-profile-disk-choice}
 \begin{gathered}
 C_r e^{-\gamma_1\tau_0}
   +Ce^{-ce^{\tau_0}}<\tfrac14s_{\rm Gr},\\
 C_r e^{-(\gamma_1-\theta)\tau_0}<\tfrac14\varepsilon,
 \qquad
 C_r e^{-(\gamma_1-\sigma)\tau_0}<\tfrac14\varepsilon,\\
 C_r e^{-\gamma_1\tau_0}<\delta_{\rm prep},\\
 \mathcal E_0(\tau_0)<\frac r2,\qquad
 \mathcal E_1(\tau_0)<\frac14 .
 \end{gathered}
\end{equation}
Every left side tends to zero for the already fixed \(r\), by the
compact support and the estimates above, so the simultaneous choice
exists.  From this point onward \(\tau_0\), the physical disk, its
marking, and its preparation convention are frozen.

The mean-value formula then gives
\begin{equation}\label{eq:profile-disk-bilipschitz}
 \frac34|v-w|
 \leq |\mathcal V_{\tau_0}(v)-\mathcal V_{\tau_0}(w)|
 \leq\frac54|v-w|,
 \qquad v,w\in\overline B_r^{E_1}(0).
\end{equation}
If
$|w|<r/2$, the homotopy
\[
 v-w+t\bigl(\mathcal V_{\tau_0}(v)-v\bigr),
 \qquad 0\leq t\leq1,
\]
does not vanish on $|v|=r$.  Hence
\[
 \deg(\mathcal V_{\tau_0}-w,B_r^{E_1},0)
 =
 \deg(\operatorname{Id}-w,B_r^{E_1},0)
 =1,
\]
and \(w\in\mathcal V_{\tau_0}(B_r^{E_1})\).  This proves
\eqref{eq:profile-realization-ball}; the injectivity in
\eqref{eq:profile-disk-bilipschitz} makes the realizing point unique
on this disk.  The inverse-function theorem makes the inverse \(C^1\)
over \(B_{r/2}^{E_1}(0)\).

For later use on the prepared neighborhood, take \(w=0\) and denote its
unique preimage by \(v_0\).  Since \(\mathcal V_{\tau_0}(v_0)=0\),
\[
 |v_0|
 \leq |\mathcal V_{\tau_0}(v_0)-v_0|
 \leq \sup_{|v|\leq r}|\mathcal V_{\tau_0}(v)-v|.
\]
The last quantity is \(<r/2\) by
\eqref{eq:simultaneous-profile-disk-choice}.  Fix
\(r'\) with \(r/2<r'<r\).  Hence \(|v_0|<r/2<r'\), and the
zero-profile point belongs to the compact
subdisk around which \(\widetilde{\mathscr O}_{\rm scat}\) is chosen.

We next place the compact physical disk in a relative
\(C^{2,\alpha}\)-open basin.  This topology-promotion step need not
use one Ricci--DeTurck chart centered at \(G_*\).  For each
\(v\in\overline B_r^{E_1}(0)\), strictness of the prepared entrance
and continuity of its transported data give a time
\(0<d_v<T(G_v)/2\) for which the time-\(d_v\) state is still a strict
prepared entrance, with positive slack in the singular-time and
buffered exterior-curvature inequalities.  Apply
Lemma~\ref{lem:positive-time-smoothing} with background the flow from
\(G_v\), and apply the finite-regularity entrance-openness argument at
time \(d_v\).  After shrinking, this gives a relative
\(C^{2,\alpha}\)-open neighborhood \(\mathscr U_v\) of \(G_v\) such
that every metric in \(\mathscr U_v\) reaches, in the corresponding
fixed gauge restart, a strict prepared entrance and hence has all
conclusions of Theorem~\ref{thm:prepared-entrance-continuation}.
Continuity of the tail lifetime and the buffered local Ricci estimate
allow \(\mathscr U_v\) to be chosen with the same prescribed
singular-time bound and an exterior-curvature bound.

The physical disk is compact in the \(C^{2,\alpha}\) topology at this
now fixed entrance time.  Choose finitely many
\(v_1,\ldots,v_{N_{\rm prof}}\) such that
\[
 \{G_v:|v|\leq r\}
 \subset\bigcup_{i=1}^{N_{\rm prof}}\mathscr U_{v_i},
\qquad
 \mathscr U_{\rm prof}
 :=\bigcup_{i=1}^{N_{\rm prof}}\mathscr U_{v_i}.
\]
Taking the maximum of the finitely many exterior constants shows that
 \(\mathscr U_{\rm prof}\) is one relative \(C^{2,\alpha}\)-open FIK
 basin containing \(G_*=G_0\) and the entire physical disk.  The
 singular-time, exterior, Type-I, localization, and full-sequence
 constants are therefore uniform after taking finite maxima of the
 upper constants and the positive minimum of the lower Type-I
 constants.  By the definition of the relative topology, choose an
 ambient \(C^{2,\alpha}\)-open metric neighborhood
 \(\widehat{\mathscr U}_{\rm prof}\) whose smooth locus is exactly
 \(\mathscr U_{\rm prof}\).  Hence
\(\mathscr U_{\rm prof}\) is an admissible choice for the existential
neighborhood \(\mathscr U\) in Theorem~A.
The
positive-time gauges in this paragraph are used only to promote the
basin topology; the disk itself retains the single marking and
prepared convention fixed at \(\tau_0\).  The support choice in
\eqref{eq:physical-profile-disk} preserves the prescribed exterior.

It remains to identify the transverse statement in the full prepared
space.  Fix the present \(\tau_0\).  The fixed marked host, graft,
scale, and map convention, together with the uniform Gram
invertibility on the compact already-sliced disk, satisfies
Corollary~\ref{cor:fixed-convention-preparation-map} at raw order
\(k+6\) and output order \(k+4\).  Hence, after shrinking a
\(C^{k+6,\alpha}\) neighborhood
\(\mathscr U_{\rm phys}^{k+6,\alpha}\) of the physical disk, there is
one fixed \(C^1\) preparation map
\[
 \mathcal P_{\tau_0}:G\longmapsto\mathbf z_{\tau_0}(G)
\]
on a high-topology neighborhood of the disk.  Exact slicing and
centered uniqueness give \(p(G_v)=0\), so it sends \(G_v\) to the
prepared tuple used throughout the preceding calculation.  The radius
and high-topology chart constants may depend on the fixed entrance
time; no \(\tau_0\)-uniform full-norm neighborhood is asserted.
 Shrink this neighborhood so that its closed-metric components lie in
 \(\widehat{\mathscr U}_{\rm prof}\).
Because the prepared disk is compact, consists of strict entrances, and
uses the one exact-core construction, the buffered physical-cover
certificates admit a finite common refinement.  Their separation and
time-width inequalities have a positive minimum on the disk.  These
strict conditions are open; hence the disk is contained in one
common-margin sliced ball
\(\mathscr B_{\rm scat}\subset
\Sigma_{\tau_0}^{k+2,\alpha}\) on which that single refined cover and
the strict time-width margin
\eqref{eq:auxiliary-buffered-time-width} work throughout.  Because the
prepared disk is compact and lies in the interior of this ball, choose
a uniformly interior common-margin subball
\(\mathscr B'_{\rm scat}\Subset_{\rm u}\mathscr B_{\rm scat}\)
containing it.  The quantitative phase retraction from
Proposition~\ref{prop:sliced-prepared-manifold} is defined on an open
neighborhood of the corresponding compact high-order disk.  By
continuity, after shrinking that neighborhood its closure is mapped
into \(\mathscr B'_{\rm scat}\).  The same proposition supplies the
uniform constant \(K_{\Pi,k+2}\); compactness is used here only for the
finite-dimensional prepared disk, not for an ambient Banach ball.

Theorem~\ref{thm:global-two-state-scattering} makes
\(\mathscr S^{\rm prep}_1
=V_\infty^{\rm sl}\circ
 \Pi_{\rm sl}^{\,k+4\to k+2}\) a
\(C^1\) map in the buffered prepared Banach chart.  The prepared
physical disk lies in \(\Sigma_{\tau_0}^{k+4,\alpha}\), and
\(D\Pi_{\rm sl}^{\,k+4\to k+2}\) is the canonical inclusion on its
tangent space.  Hence, at
\(\mathcal P_{\tau_0}(G_v)\), the restriction of
\(D\mathscr S^{\rm prep}_1\) to the disk tangent is
\(D\mathcal V_{\tau_0}(v)\).  This operator is invertible by the
derivative estimate
\eqref{eq:profile-map-degree-C1-error}.  Hence
\(D\mathscr S^{\rm prep}_1\) is onto and has the same bounded right
inverse.
Compactness of
\(\{\mathcal P_{\tau_0}(G_v):|v|\leq r'\}\) and continuity of the
differential give the prepared neighborhood
 \(\widetilde{\mathscr O}_{\rm scat}\) on which this remains true.
 Shrink it inside the phase-retraction neighborhood just chosen and
 inside the inverse image of
 \(\widehat{\mathscr U}_{\rm prof}\).  Then
\[
 \overline{\widetilde{\mathscr O}_{\rm scat}}
 \subset\operatorname{dom}
   \Pi_{\rm sl}^{\,k+4\to k+2},
 \qquad
 \Pi_{\rm sl}^{\,k+4\to k+2}
  (\overline{\widetilde{\mathscr O}_{\rm scat}})
 \subset\mathscr B'_{\rm scat}
 \Subset_{\rm u}\mathscr B_{\rm scat},
\]
which proves \eqref{eq:profile-scattering-domain-inclusions}.  These
inclusions retain strict preparation independently of the
closed-metric basin condition.  The Banach submersion theorem proves
the level-set and local-foliation assertions.

The ball realization gives a point of the zero level.  Here a
\emph{forward prepared restart} means restriction of the same coupled
trajectory at a later normalized time, with its marking, DeTurck
gauge, harmonic-map chart, and slice transported by the fixed
preparation convention; it does not mean applying an unrelated new
phase projection.  Restarting \eqref{eq:V-infinity-formula} in this
sense leaves \(V_\infty\) unchanged if the absolute
\(\tau\)-coordinate is retained.  If the coordinate is reset at
absolute time \(S\) to \(\widehat\tau=\tau-S\), then
every time-typed object and the slice are simultaneously translated
as specified in Theorem~\ref{thm:first-stable-profile}; this is only a
relabeling of the same tail.  In that translated convention
\eqref{eq:profile-restart-rescaling} gives the exact factor
\(e^{-\gamma_1S}\).  Thus the condition
\(V_\infty=0\) is preserved by every forward prepared restart.  A
submersion level has empty interior, so its complement is open and
dense.
\end{proof}

\begin{remark}[Domain of the profile foliation]
\label{rem:profile-foliation-domain}
The split submersion in
\eqref{eq:first-profile-scattering-map} is a theorem on the single
buffered prepared Banach manifold
\(\widetilde{\mathscr O}_{\rm scat}\), with its fixed marking,
Ricci--DeTurck reference metric, harmonic-map gauge, and transported
restart convention.  The physical disk \(v\mapsto G_v\) is embedded
in that manifold by the one fixed preparation map
\(\mathcal P_{\tau_0}\).  Thus the theorem gives a genuine
finite-codimensional foliation of prepared states and a transverse
physical realization of every small profile.  It does not assert
that forgetting the auxiliary prepared variables sends every leaf
injectively to the space of metrics, nor that the leaves descend to a
quotient by arbitrary diffeomorphisms.  These qualifications concern
only the domain of the foliation; they do not alter the profile
realization, split-surjectivity, or strong-stable decay conclusions.
On this ambient prepared manifold the profile map is the sliced
scattering map precomposed with the fixed buffered phase projection
\(\Pi_{\rm sl}^{\,k+4\to k+2}\).  Consequently its fibers contain the local
phase-retraction fibers; their asserted codimension is
\(\dim E_1\), not \(9+\dim E_1\).
\end{remark}

\begin{proof}[Proof of Theorem~\ref{thm:intro-sharp-scattering}]
Part~I follows from
Theorems~\ref{thm:first-stable-profile},
\ref{thm:sharp-marked-spacetime},
\ref{thm:global-two-state-scattering}, and
\ref{thm:quadratic-geometric-asymptotics}; its buffered ambient
assertion is the chain rule together with the derived bound
\(K_{\Pi,k+2}\) from
\eqref{eq:quantitative-phase-retraction-Lipschitz}--%
\eqref{eq:quantitative-phase-retraction-derivative} and the uniform
sliced constants on \(\mathscr B'\).
Part~II is exactly Theorem~\ref{thm:profile-realization}, applied to the
exact-core basin of Theorem~\ref{thm:intro-open-basin} at
\(\tau_0=\widehat\tau_0\), with the special domains in
\eqref{eq:profile-scattering-domain-inclusions}.
\end{proof}

Theorem~\ref{thm:intro-sharp-scattering} is now available before any
positive-time physical restart is introduced.  Its Part~II supplies
the exact-core transverse disk and the specially realized prepared
foliation, while its Part~I applies to every common-margin sliced ball
satisfying its hypotheses; the bridge below constructs such a ball.
It remains only to
pull this prepared scattering geometry back to an open
\(h^{2,\alpha}\)-neighborhood of actual initial metrics; no second
scattering argument is required.

\subsection{Physical-time amplitude and the fixed-restart physical foliation}

The marked physical-time amplitude \(\mathfrak A_1\) was defined in
Definition~\ref{def:physical-first-amplitude}.  The superscript
\({\rm sl}\) below records restriction to one named sliced ball, not a
new normalization.
Lemma~\ref{lem:physical-amplitude-restart-covariance} proves its exact
invariance under transported forward restart, including an additive
reset of normalized time.  It remains marked and gauge-fixed, and no
covariance under an independently chosen phase projection or an
arbitrary time-dependent re-marking is asserted.

\begin{lemma}[Uniform fixed physical restart and profile transversality]
\label{lem:uniform-physical-restart-transversality}
Fix the host, rates, prepared package, marking convention, entrance
time, and physical profile disk
\[
 v\longmapsto G_v,\qquad |v|\leq r,
\]
supplied by Part~II of
Theorem~\ref{thm:intro-sharp-scattering}.  Let \(v_0\) be
the unique point on the original prepared disk for which
\[
 V_\infty(G_{v_0})=0.
\]
Then there exist a compact parameter ball
\[
 \overline B_\varrho(v_0)\Subset B_{r'}(0),
\]
where \(r'\) is the compact-subdisk radius fixed in
Theorem~\ref{thm:profile-realization}, a sufficiently small physical
time \(d>0\), a unique parameter \(v_{\rm ss}\) near \(v_0\), and
\[
 G_{\rm ss}:=G_{v_{\rm ss}},\qquad r_{\rm ss}>0,
\]
an open neighborhood \(\mathfrak U_d^{2,\alpha}\) of
\[
 \mathfrak D_{\rm ss}
 :=\{G_v:|v-v_{\rm ss}|\leq r_{\rm ss}\}
\]
in the positive cone of the little-H\"older Banach space
\(h^{2,\alpha}(S^2T^*\widehat X)\), and common-margin sliced balls of
strict prepared entrances
\[
 \mathscr B_d\Subset_{\rm u}\mathscr B_{[0,d]}
 \subset\Sigma_{\tau_0}^{k+2,\alpha},
\]
such that the following hold.  On the larger ball set
\[
 \mathfrak A_1^{\rm sl}
 :=\left.\mathfrak A_1\right|_{\mathscr B_{[0,d]}}.
\]
\begin{enumerate}
\item
Let \(\mathbf z_v\) be the original prepared lift of \(G_v\).  There is
a single host, graft, scale, marking, initial-map, and phase convention,
independent of \(v\), and disk maps
\[
\mathcal P_s^{\rm disk}(G_v)\in\mathscr B_{[0,d]},
 \qquad
 0\leq s\leq d,\quad
 v\in\overline B_\varrho(v_0),
\]
such that
\begin{equation}\label{eq:fixed-restart-disk-C1-limit}
 \mathcal P_0^{\rm disk}(G_v)=\mathbf z_v,\qquad
 \mathcal P_s^{\rm disk}(G_\bullet)
 \longrightarrow\mathbf z_\bullet
 \quad\hbox{in }C^1_v
 \quad\text{as }s\downarrow0.
\end{equation}
For the selected \(d>0\), the corresponding fixed positive-time
preparation map
\begin{equation}\label{eq:uniform-physical-restart-map}
 \mathcal P_d:
 \mathfrak U_d^{2,\alpha}
 \longrightarrow\mathscr B_d,
 \qquad
 \mathcal P_d(G)
 =
 \mathbf z_d\!\left(
   \widetilde{\mathfrak R}_d(G)\right),
\end{equation}
is \(C^1\).  Here the source and target are the corresponding
little-H\"older Banach charts.  One buffered physical cover and the
strict time-width margin
\eqref{eq:auxiliary-buffered-time-width} work throughout
\(\mathscr B_d\).
\item
Define
\begin{equation}\label{eq:fixed-restart-amplitude-family}
 \mathcal F:
 [0,d]\times\overline B_\varrho(v_0)\longrightarrow E_1,
 \qquad
 \mathcal F(s,v):=
 \mathfrak A_1^{\rm sl}
 \bigl(\mathcal P_s^{\rm disk}(G_v)\bigr).
\end{equation}
This map is \(C^1\) down to \(s=0\), and
\begin{equation}\label{eq:fixed-restart-amplitude-base}
 \mathcal F(0,v)
 =\lambda_\infty(G_v)^{-\gamma_1}V_\infty(G_v).
\end{equation}
The selected parameter satisfies
\[
 \mathcal F(d,v_{\rm ss})=0,\qquad
 D_v\mathcal F(d,v_{\rm ss}):E_1\longrightarrow E_1
 \ \hbox{is invertible}.
\]
For a smooth \(G\in\mathfrak U_d^{2,\alpha}\), the original Ricci flow
and the tail generated by \(\mathcal P_d(G)\) are related by the
data-dependent terminal pullback produced by the one fixed DeTurck
gauge convention, as in \eqref{eq:Ricci-tail-fixed-pullback}.  The
re-prepared state \(\mathcal P_d(G_v)\) is not identified with the
transported evolved tuple from the original preparation of \(G_v\);
no cross-preparation covariance is asserted or needed.
\item
The map
\[
 \mathfrak A_{1,d}^{\rm phys}
 :=
 \mathfrak A_1^{\rm sl}\circ\mathcal P_d:
 \mathfrak U_d^{2,\alpha}\longrightarrow E_1
\]
is split \(C^1\)-submersive at \(G_{\rm ss}\).  After
\(\mathfrak U_d^{2,\alpha}\) is decreased once more, it is a split
\(C^1\) submersion throughout that neighborhood.
\end{enumerate}
The uniqueness of \(v_{\rm ss}\) is local on the selected transverse
disk for this fixed small \(d\); no disk-independent canonical metric
is claimed.
\end{lemma}

\begin{proof}
Part~II of Theorem~\ref{thm:intro-sharp-scattering}, together with
\eqref{eq:profile-disk-bilipschitz}, gives the unique \(v_0\), with
\(|v_0|<r/2\).  Choose a compact parameter ball
\(\overline B_\varrho(v_0)\Subset B_{r'}(0)\), where \(r'\) is the
compact-subdisk radius fixed in
Theorem~\ref{thm:profile-realization}.  The
original prepared lifts \(\mathbf z_v\) on this ball have one common
host, graft, scale, marking, initial-map, and phase convention; only
their compactly supported tensor component varies.

Use the smooth Ricci flow from \(G_{v_0}\) as a reference
time-dependent DeTurck background.  Choose \(d_0>0\) so that this
background exists smoothly on \([0,2d_0]\).  In fixed little-H\"older
charts, the quasilinear parabolic solution theorem gives a \(C^1\) map
\[
 [0,d_0]\times\overline B_\varrho(v_0)
 \longrightarrow h^{k+4,\alpha},
 \qquad
 (s,v)\longmapsto\widetilde{\mathfrak R}_s(G_v),
\]
including uniform \(C^1_v\) convergence to \(G_v\) as \(s\downarrow0\).
The fixed raw-lift map is \(C^1\) into the prepared
\(h^{k+4,\alpha}\) chart, and the buffered phase retraction is \(C^1\)
from that chart to \(\Sigma_{\tau_0}^{k+2,\alpha}\).  Raw-lift these
metrics using that convention, retaining the original normalized label
\(\tau_0\), and apply the centered phase projection.
Because \(\mathbf z_v\) is already sliced, centered uniqueness gives
\(\mathcal P_0^{\rm disk}(G_v)=\mathbf z_v\).  The parameterized
phase theorem and the smooth finite-time dependence give
\eqref{eq:fixed-restart-disk-C1-limit}.
The compact original disk has common strict margins, common normalized
operative and \(+\)-witness triples, and the complete finite physical
certificate.  Apply the compact-family clause of
Corollary~\ref{cor:fixed-convention-preparation-map} uniformly in
\((s,v)\), after decreasing \(d_0\).  Equivalently, the normalized
reserve-to-operative lemma preserves the normalized operative triple
and furnishes a new normalized reserve, while the actual carriers
remain in the same fixed physical inner locus.  The reference-carrier
certificate therefore supplies their unchanged common physical
\(+\)- and operative tiers, and outer-ball membership supplies the
quarter-modulus future-window estimate.  Continuity preserves the
\(\mu_{\rm RF}\)-width slack.  Hence every
\(\mathcal P_s^{\rm disk}(G_v)\) lies in one named common-margin strict
sliced ball \(\mathscr B_{[0,d_0]}\) and retains the operative
modulus-compatibility inequality.

After shrinking around this compact disk image,
Proposition~\ref{prop:intro-static-sliced-chart} supplies a buffered
ambient domain and uniformly interior subball satisfying the
hypotheses of Part~I of
Theorem~\ref{thm:intro-sharp-scattering}.  That part makes
\(\mathfrak A_1\) \(C^1\) on \(\mathscr B_{[0,d_0]}\).
Consequently \(\mathcal F\) in
\eqref{eq:fixed-restart-amplitude-family} is \(C^1\) down to \(s=0\).
At \(s=0\), Definition~\ref{def:physical-first-amplitude} on the
original prepared disk gives
\eqref{eq:fixed-restart-amplitude-base}.  Since
\(\mathcal V_{\tau_0}(v):=V_\infty(G_v)\) vanishes at \(v_0\),
\begin{equation}\label{eq:physical-amplitude-transversality}
 D_v\mathcal F(0,v_0)
 =
 \lambda_\infty(G_{v_0})^{-\gamma_1}
 D\mathcal V_{\tau_0}(v_0).
\end{equation}
This operator is invertible by
\eqref{eq:profile-map-degree-C1-error}.  The one-sided parameterized
implicit-function theorem on the manifold with boundary
\([0,d_0]\times B_\varrho(v_0)\) therefore gives
\(d_1\in(0,d_0]\) and a
unique \(C^1\) curve \(v_{\rm ss}(s)\) near \(v_0\) such that
\[
 \mathcal F(s,v_{\rm ss}(s))=0,\qquad
 D_v\mathcal F(s,v_{\rm ss}(s))\ \hbox{is invertible},
\qquad0\leq s\leq d_1.
\]
(Equivalently, extend \(\mathcal F\) to negative \(s\) by
\(\widetilde{\mathcal F}(s,v)
=2\mathcal F(0,v)-\mathcal F(-s,v)\) for \(s<0\), and apply the
ordinary parameterized theorem.)
Fix one \(d\in(0,d_1]\), put
\(v_{\rm ss}:=v_{\rm ss}(d)\) and
\(G_{\rm ss}:=G_{v_{\rm ss}}\), and then choose
\[
 0<2r_{\rm ss}<
 \operatorname{dist}(v_{\rm ss},\partial B_\varrho(v_0)).
\]

Now apply Lemma~\ref{lem:positive-time-smoothing}, with target order
\(k+4\), to the same reference background on \([0,2d]\).
After \(r_{\rm ss}\) and the \(h^{2,\alpha}\) neighborhood are
decreased, the retained disk lies in the domain of the resulting
single map
\[
 \widetilde{\mathfrak R}_d:
 h^{2,\alpha}\longrightarrow h^{k+4,\alpha}.
\]
Apply the same fixed host, marking, graft, scale, and
initial-map convention used to define
\(\mathcal P_d^{\rm disk}\), and then apply
Corollary~\ref{cor:fixed-convention-preparation-map} with raw order
\(k+4\) and output order \(k+2\).  The centered phase map uses exactly
the two available derivatives.  Its compact-family gluing clause and
centered uniqueness produce one \(C^1\) map
\eqref{eq:uniform-physical-restart-map}; by construction, its
restriction to the retained smooth disk is exactly
\(\mathcal P_d^{\rm disk}\).  Shrinking around the compact image keeps
all strict inequalities, the common buffered harmonic witnesses, the
common physical cover, and the time-width margin, and places the image
in one common-margin sliced ball
\(\mathscr B_d\Subset_{\rm u}\mathscr B_{[0,d]}\), where we set
\(\mathscr B_{[0,d]}:=\mathscr B_{[0,d_0]}\).  This proves the first
assertion.
Equation~\eqref{eq:DeTurck-to-Ricci-prefix}, uniqueness of Ricci flow,
and the definition of the transported marking give
\eqref{eq:Ricci-tail-fixed-pullback}.
The terminal diffeomorphism depends on the input metric but is
produced by this one fixed gauge convention.  More explicitly, it is
absorbed into the
physical marking as
\(\widetilde\Xi\mapsto\chi_d^{-1}\circ\widetilde\Xi\).  Pulling
\(\chi_d^*\widetilde H\) back by this marking gives exactly
\(\widetilde\Xi^*\widetilde H\), so the normalized model-side tensor
and its sliced moments are unchanged.

The transversality argument uses the \(C^1\)-nearby family
\(\mathcal F\), rather than the parameter-dependent auxiliaries
obtained by evolving the earlier preparations.  Since
\(\mathcal P_d(G_v)=\mathcal P_d^{\rm disk}(G_v)\), the derivative of
\(\mathfrak A_{1,d}^{\rm phys}\) on the finite-dimensional tangent
\(D_vG_v(E_1)|_{v=v_{\rm ss}}\) is
\(D_v\mathcal F(d,v_{\rm ss})\), hence is invertible.  Composing its
inverse with
\(D_vG_v|_{v=v_{\rm ss}}:E_1\to h^{2,\alpha}\) gives a bounded right
inverse
for
\(D\mathfrak A_{1,d}^{\rm phys}(G_{\rm ss})\).
The map \(\mathfrak A_1^{\rm sl}\) is \(C^1\) on
\(\mathscr B_d\) by the sliced conclusion in Part~I of
Theorem~\ref{thm:intro-sharp-scattering}, positivity of
\(\lambda_\infty\), and the ordinary finite-dimensional chain rule.
The same is true of its composition with \(\mathcal P_d\).
Split surjectivity is open when a fixed finite-dimensional right
inverse is retained.  A final decrease of
\(\mathfrak U_d^{2,\alpha}\) proves the last assertion.  After this
submersion neighborhood is fixed, decrease \(r_{\rm ss}\) once more
so that the retained compact transverse disk
\(\mathfrak D_{\rm ss}\) lies inside it.  This does not alter the
chosen \(d\), \(v_{\rm ss}\), or right inverse.
\end{proof}

\paragraph{Physical pullback of Theorem~C.}
Part~I of Theorem~\ref{thm:intro-sharp-scattering} supplies the
componentwise \(C^1\) scattering map and sharp asymptotics on the common
prepared ball, while Part~II supplies the exact-core transverse
realization.  Lemma~\ref{lem:uniform-physical-restart-transversality}
provides the remaining fixed positive-time
\(h^{2,\alpha}\)-to-prepared bridge.  The following proposition records
their pullback to actual initial metrics and is the final input in the
proof of Theorem~B.

\begin{proposition}[Fixed-restart pullback of prepared scattering:
first-profile moduli and quantitative marked comparison]
\label{prop:physical-amplitude-foliation}
In the setting of
Lemma~\ref{lem:uniform-physical-restart-transversality}, there is a
smaller open neighborhood
\[
 \mathfrak U_{\rm scat}^{2,\alpha}
 \subset\mathfrak U_d^{2,\alpha},
 \qquad
 \overline{\mathfrak U_{\rm scat}^{2,\alpha}}
 \subset\mathfrak U_d^{2,\alpha}
\]
of \(G_{\rm ss}\) in the positive
\(h^{2,\alpha}(S^2T^*\widehat X)\) cone with the following
properties.  Its smooth relative locus
\[
 \mathscr U_{\rm scat}
 :=
 \mathfrak U_{\rm scat}^{2,\alpha}
 \cap\operatorname{Met}^{\infty}(\widehat X)
\]
is contained in the relative \(C^{2,\alpha}\)-open basin
\(\mathscr U_{\rm prof}\) of
Theorem~\ref{thm:profile-realization}.  In particular, no assertion is
made about the whole marked basin or about initial metrics outside
this local scattering neighborhood.

Let \(T_{\rm tail}\), \(\lambda_\infty\), \(\Psi_\infty\), and
\(V_\infty\) denote the components of the Part~I scattering map in
Theorem~\ref{thm:intro-sharp-scattering}, evaluated at
\(\mathcal P_d(G)\).  Then the physical fixed-restart scattering map
\begin{equation}\label{eq:physical-amplitude-scattering-map}
 \begin{aligned}
 \mathscr S_{\rm phys,d}:
 \mathfrak U_{\rm scat}^{2,\alpha}
 &\longrightarrow
 \bigl\{\text{componentwise scattering data}\bigr\},\\
 G&\longmapsto
 \left(
  d+T_{\rm tail},\
  \log\lambda_\infty,\
  \Psi_\infty,\
  \mathfrak A_{1,d}^{\rm phys}
 \right),\\
 \mathfrak A_{1,d}^{\rm phys}(G)
 &=
 e^{-\gamma_1\log\lambda_\infty}
 V_\infty
 =
 \mathfrak A_1^{\rm sl}\!\left(\mathcal P_d(G)\right)
 \end{aligned}
\end{equation}
is \(C^1\) in the following precise sense: its scalar and
\(E_1\)-valued components are \(C^1\), and for every
\(K\Subset M\) and \(m\geq0\), its phase component is \(C^1\) into the
fixed \(C^m(K)\) exponential chart.  Thus no single Fr\'echet target
for the phase component is being asserted.

The amplitude component in
\eqref{eq:physical-amplitude-scattering-map} is a split \(C^1\)
submersion.  It vanishes at \(G_{\rm ss}\).  After decreasing
\(\mathfrak U_{\rm scat}^{2,\alpha}\), let
\[
 K_{\rm ss}^{\rm phys}
 :=\ker D\mathfrak A_{1,d}^{\rm phys}(G_{\rm ss}).
\]
There are neighborhoods \(\mathcal O_{\rm phys}\) of \(G_{\rm ss}\),
\(\mathcal O_K\) of the origin in \(K_{\rm ss}^{\rm phys}\), and
\(B_\epsilon^{E_1}(0)\), and a \(C^1\) diffeomorphism
\[
 \Phi_{\rm phys}:\mathcal O_{\rm phys}
 \longrightarrow
 \mathcal O_K\times B_\epsilon^{E_1}(0)
\]
such that
\[
 \mathfrak A_{1,d}^{\rm phys}
 \!\left(\Phi_{\rm phys}^{-1}(k,A)\right)=A.
\]
Consequently, for \(A\) in a sufficiently small ball in \(E_1\),
\begin{equation}\label{eq:physical-amplitude-fibers}
 \mathscr F_A^{\rm phys}
 :=
 \left(\mathfrak A_{1,d}^{\rm phys}\right)^{-1}(A)
 \cap\mathfrak U_{\rm scat}^{2,\alpha}
\end{equation}
is a split \(C^1\) Banach submanifold of codimension
\(\dim E_1\), and the connected components of these fibers form a
local \(C^1\) foliation of
\(\mathfrak U_{\rm scat}^{2,\alpha}\).  The zero fiber contains
\(G_{\rm ss}\) and has the improved strong-stable rate.  The numerical
value of \(\mathfrak A_{1,d}^{\rm phys}\) is unchanged when the same
prepared tail is described after a transported forward restart.  Thus,
if such a restarted description is represented in this local chart, it
lies in the same amplitude level.  No comparison with an independently
prepared state, and no assertion that the restarted metric remains
inside this chart, is made.  The complement of the zero fiber is open
and dense in
\(\mathfrak U_{\rm scat}^{2,\alpha}\).

After \(r_{\rm ss}\) is decreased, there are
\(\epsilon_{\rm A}>0\) and \(0<r_{\rm A}<r_{\rm ss}\) such that every
\(A\in B_{\epsilon_{\rm A}}^{E_1}(0)\) has a unique realization
\[
 A=\mathfrak A_{1,d}^{\rm phys}(G_v),
 \qquad |v-v_{\rm ss}|<r_{\rm A},
\]
on the chosen physical transverse disk.  These realizing metrics
agree with \(G_*\) outside the implantation region.

Finally, let \(G_1,G_2\in\mathscr U_{\rm scat}\), write
\(G_i(t)=G(t;G_i)\), let \(T_i\) be their singular times, and, for
\(\delta>0\), use the frozen base-time
markings supplied by the common transported convention to set
\[
 \mathcal G_{i,\delta}(s)
 :=
 \delta^{-1}
 \left(\Xi^{G_i}_{T_i-\delta}\right)^*
 G_i(T_i+s\delta).
\]
The two markings are thereby viewed on the same FIK model; the fixed
limiting-phase alignment is the one already built into this
convention.  Then, for every \(K\Subset M\),
\(I\Subset(-\infty,0)\), and \(m\geq0\),
\begin{equation}\label{eq:physical-quantitative-two-state}
 \delta^{-\gamma_1}
 \left(\mathcal G_{1,\delta}-\mathcal G_{2,\delta}\right)
 =
 \mathcal J_{
  \mathfrak A_{1,d}^{\rm phys}(G_1)
  -\mathfrak A_{1,d}^{\rm phys}(G_2)}
 +O_{C^m(K\times I)}(\delta^{\widehat\delta}).
\end{equation}
In particular,
\begin{equation}\label{eq:physical-first-order-equivalence}
 \mathfrak A_{1,d}^{\rm phys}(G_1)
 =
 \mathfrak A_{1,d}^{\rm phys}(G_2)
 \quad\Longleftrightarrow\quad
 \delta^{-\gamma_1}
 \left(\mathcal G_{1,\delta}
       -\mathcal G_{2,\delta}\right)
 \longrightarrow0
 \quad\text{in }C^\infty_{\rm loc}
\end{equation}
as \(\delta\downarrow0\).  This is marked first-order asymptotic
completeness in the single fixed convention.  It is not a quotient
statement under arbitrary diffeomorphisms or time-dependent
re-markings.  For the reverse implication it is enough to assume this
vanishing on every compact spatial set in one compact time window
containing \(s=-1\).  In the product chart above, the intersections of
the fibers with the smooth relative locus are therefore exactly the
local classes of marked first-order asymptotic agreement, and their
local leaf space is modeled on \(B_\epsilon^{E_1}(0)\).
\end{proposition}

\begin{proof}
Shrink the neighborhood from
Lemma~\ref{lem:uniform-physical-restart-transversality} inside the
ambient \(C^{2,\alpha}\)-open thickening
\(\widehat{\mathscr U}_{\rm prof}\) fixed in
Theorem~\ref{thm:profile-realization}, and inside the inverse image of
the common scattering ball \(\mathscr B_d\).  This gives the asserted
smooth relative locus.  After a further shrinking,
Proposition~\ref{prop:intro-static-sliced-chart} supplies the buffered
ambient phase-retraction domain required by Part~I of
Theorem~\ref{thm:intro-sharp-scattering}.  Its sliced \(C^1\)
scattering conclusion, composed with \(\mathcal P_d\), proves the
asserted regularity of
\[
 \mathscr S_{\rm sl}\circ\mathcal P_d.
\]
Adding the fixed prefix length \(d\) to the tail singular time does
not change differentiability.  Since \(\lambda_\infty>0\), the map
\((\log\lambda_\infty,V_\infty)\mapsto
e^{-\gamma_1\log\lambda_\infty}V_\infty\) is smooth.  This proves
\eqref{eq:physical-amplitude-scattering-map}.

Split submersivity is the last conclusion of
Lemma~\ref{lem:uniform-physical-restart-transversality}.  The Banach
submersion theorem gives the displayed product chart,
\eqref{eq:physical-amplitude-fibers}, and the local foliation.
The identity
\[
 \mathfrak A_{1,d}^{\rm phys}=0
 \quad\Longleftrightarrow\quad V_\infty=0
\]
follows from \(\lambda_\infty>0\).  The improved rate is the
corresponding conclusion of Part~I of
Theorem~\ref{thm:intro-sharp-scattering}
(cf. Theorems~\ref{thm:first-stable-profile} and
\ref{thm:sharp-marked-spacetime}).  Forward invariance follows
from Lemma~\ref{lem:physical-amplitude-restart-covariance}.  A submersion level
is closed in the present neighborhood and has empty interior, so its
complement is open and dense.

At \(v_{\rm ss}\), the derivative
\(D_v\mathcal F(d,v_{\rm ss})\) is invertible.  The
finite-dimensional inverse-function theorem therefore gives
\(\epsilon_{\rm A}\), \(r_{\rm A}\), and the asserted unique
realization.  Exterior agreement is inherited from
\eqref{eq:physical-profile-disk}.

It remains to verify the quantitative two-state formula and its
completeness consequence.  The sharp Jacobi normal form in Part~I of
Theorem~\ref{thm:intro-sharp-scattering}, together
with the linearity of \(V\mapsto\mathcal J_V\) and
\eqref{eq:physical-first-amplitude}, gives, for every
\(K\Subset M\), \(I\Subset(-\infty,0)\), and \(m\geq0\),
\begin{equation}\label{eq:physical-amplitude-sharp-expansion}
 \mathcal G_{i,\delta}(s)
 =
 g_{\mathrm{FIK}}(s)
 +\delta^{\gamma_1}
  \mathcal J_{\mathfrak A_{1,d}^{\rm phys}(G_i)}(s)
 +O_{C^m(K\times I)}
  \left(\delta^{\gamma_1+\widehat\delta}\right).
\end{equation}
The constants are locally uniform because
\(\mathcal P_d(\mathfrak U_{\rm scat}^{2,\alpha})\) lies in one
common-margin ball and \(\lambda_\infty\) stays in a compact subset of
\((0,\infty)\).  Subtracting
\eqref{eq:physical-amplitude-sharp-expansion} for \(i=1,2\) and dividing
by \(\delta^{\gamma_1}\) proves
\eqref{eq:physical-quantitative-two-state}.  Equality of the two
amplitudes therefore gives the vanishing limit in
\eqref{eq:physical-first-order-equivalence}.  Conversely, that
vanishing limit forces
\[
 \mathcal J_{
  \mathfrak A_{1,d}^{\rm phys}(G_1)
  -\mathfrak A_{1,d}^{\rm phys}(G_2)}=0.
\]
Since \(\mathcal J_A(-1)=A\), the two amplitudes are equal.  This also
shows that a nonzero amplitude has exact marked first order
\(\delta^{\gamma_1}\) on some compact spacetime window containing
\(s=-1\).
\end{proof}

\begin{proof}[Proof of
Theorem~\ref{thm:intro-physical-amplitude-master}]
Fix the host and implantation data of Theorem~A, and let
\(\alpha_{\rm A}\in(0,1)\) be the exponent produced by its
low-topology construction.  Instantiate every prepared result below
with \(\alpha=\alpha_{\rm A}\), choose \(k_0\geq12\), then choose
\[
 0<\sigma<\theta<\beta,
\]
and only afterward choose the compatible numerical prepared package,
entrance size, graft radius, and sufficiently late entrance time.  The
selection order is exactly the one fixed in
Remark~\ref{conv:authoritative-adaptive-order}.

We first record the regularity assignments, because the four maps used
in the argument have different domains.  On a uniformly interior
common-margin sliced ball, the finite-time solution map is
\begin{equation}\label{eq:physical-master-regularity-ledger-flow}
 \Sigma_{\tau_0}^{k_0+2,\alpha_{\rm A}}
 \longrightarrow
 \mathscr P_{\tau}^{k_0,\alpha_{\rm A}},
\end{equation}
and is \(C^1\); this is the explicit two-derivative loss in
Proposition~\ref{prop:coupled-local-feedback}.  The ambient phase
retraction is
\begin{equation}\label{eq:physical-master-regularity-ledger-phase}
 \operatorname{dom}\Pi_{\rm sl}^{\,k_0+4\to k_0+2}
 \xrightarrow{\ \Pi_{\rm sl}^{\,k_0+4\to k_0+2}\ }
 \Sigma_{\tau_0}^{k_0+2,\alpha_{\rm A}}.
\end{equation}
The raw fixed-time physical preparation map used in the transverse-disk
construction reserves two further derivatives and therefore starts at
order \(k_0+6\).  In contrast, after a fixed positive physical time
\(d>0\), parabolic smoothing and the same phase buffer give the bridge
\begin{equation}\label{eq:physical-master-regularity-ledger-restart}
 \mathcal P_d:
 h^{2,\alpha_{\rm A}}(S^2T^*\widehat X)
 \longrightarrow
 \Sigma_{\tau_0}^{k_0+2,\alpha_{\rm A}}.
\end{equation}
Here \(\tau_0\) remains the normalized label of the single fixed
positive-time re-preparation convention; no additive reset of a
parameter-dependent evolved slice is used.  It is this last map,
rather than the raw high-order preparation map, that defines the
physical neighborhood in Theorem~B.  Finally, the scattering target is
componentwise: \(T\), \(\log\lambda_\infty\), and the \(E_1\)-component
are Banach-valued \(C^1\) maps, while for each \(K\Subset M\) and
\(m\geq0\) the phase is \(C^1\) into one fixed \(C^m(K)\) exponential
chart.  No common Fr\'echet target for all phase components or
same-order continuation map is invoked.

The exact-core implant of
Proposition~\ref{prop:exact-core-implant}, the one-state continuation
theorem, and the positive-time promotion
Corollary~\ref{cor:low-topology-open-basin} produce the relative
\(C^{2,\alpha_{\rm A}}\)-open basin.  The invariant marked-basin
Theorem~\ref{thm:marked-FIK-basin} gives precisely the globally Type-I,
localized, exterior-regular, and full-sequence marked FIK conclusions
collected in Theorem~A.  Part~II of
Theorem~\ref{thm:intro-sharp-scattering} permits this same formation
basin to be chosen as
\[
 \mathscr U_0:=\mathscr U_{\rm prof},
\]
and supplies the compactly supported transverse physical profile disk
and the specially realized prepared foliation.  Together with
Theorem~A, this proves the formation assertion in Part~I of Theorem~B.
The formation conclusion alone uses no two-state comparison theorem;
the Part~II disk and prepared foliation are retained for the later
scattering synthesis.

Apply Lemma~\ref{lem:uniform-physical-restart-transversality} to that
Part~II transverse disk.  The lemma first transports its unique zero
and invertible disk derivative through a sufficiently short positive
prefix, and then fixes \(d>0\), \(G_{\rm ss}\), and a single fixed
preparation convention.  After the source is shrunk, it gives the
\(C^1\) bridge
\[
 \mathcal P_d:
 \widehat{\mathscr U}_1
 \longrightarrow\mathscr B_d
\]
from an open \(h^{2,\alpha_{\rm A}}\)-neighborhood into one
common-margin sliced ball of strict prepared entrances.  Its amplitude
derivative has a bounded right inverse.  In the notation of
\eqref{eq:physical-master-regularity-ledger-restart}, set
\[
 \widehat{\mathscr U}_1
 :=\mathfrak U_{\rm scat}^{2,\alpha_{\rm A}},
 \qquad
 \mathscr U_1
 :=\widehat{\mathscr U}_1\cap
   \operatorname{Met}^{\infty}(\widehat X)
 \subset\mathscr U_0.
\]

After one further shrinking,
Proposition~\ref{prop:intro-static-sliced-chart} supplies the buffered
ambient phase-retraction domain required in Part~I of
Theorem~\ref{thm:intro-sharp-scattering}.  That part, applied to
\(\mathscr B_d\), now gives the unique \(V_\infty\in E_1\), the
componentwise \(C^1\) scattering tuple
\[
 (T_{\rm tail},\log\lambda_\infty,\Psi_\infty,V_\infty),
\]
the sharp marked Jacobi expansion, and the quadratic scale and phase
response.  For the exponent in
\eqref{eq:intro-physical-Jacobi-normal-form} and
\eqref{eq:intro-master-quantitative-two-state}, set
\[
 \eta_1:=\widehat\delta
 =\min\{\delta_0,\gamma_1\}>0,
\]
where \(\widehat\delta\) is the remainder exponent in the sharp marked
spacetime theorem.  Under a transported forward restart, the transformation
laws for \(V_\infty\) and \(\lambda_\infty\) cancel exactly in
\[
 \mathfrak A_1=\lambda_\infty^{-\gamma_1}V_\infty
\]
by Lemma~\ref{lem:physical-amplitude-restart-covariance}.  Thus this
single prepared coordinate controls both the linear physical Jacobi
term and its first nonlinear geometric feedback.  Pulling it back
through \(\mathcal P_d\) gives the amplitude in Theorem~B; the fixed
prefix changes only \(T_{\rm tail}\) to the physical singular time
\(T=d+T_{\rm tail}\).

Proposition~\ref{prop:physical-amplitude-foliation} is precisely this
physical pullback of Part~I, with the transversality furnished by
Part~II and preserved by the fixed-restart lemma.  The Banach
submersion theorem gives the local product chart and the
codimension-\(\dim E_1\) foliation, including the nonempty zero leaf and
its open dense complement.  The retained transverse disk gives the
unique realization there of every sufficiently small amplitude.
Subtracting the two sharp marked expansions supplied by Part~I gives
the quantitative two-state Jacobi limit; marked first-order
completeness follows, with the reverse implication obtained from
\(\mathcal J_A(-1)=A\).  This proves Parts~II and~III.

Every construction in the last two paragraphs uses the same fixed
convention for positive-time re-preparation, together with the marking
and gauge transported along each resulting tail.  Consequently the
conclusions have precisely the scope stated in Part~IV: they apply on
\(\mathscr U_1\), rather than all of \(\mathscr U_0\), with this
transported marking and gauge.  This proves Part~IV and completes
Theorem~B.
\end{proof}

\appendix

\section{Self-contained certification of the FIK spectrum}
\label{app:self-contained-FIK-spectrum}

This appendix gives an independent proof of the exact spectral theorem
used in the body.  The organization follows the tensor Wigner strategy
introduced by Naff--Ozuch~\cite{NaffOzuch}, but every geometric formula,
finite certificate, and exclusion argument required for
Theorem~\ref{thm:FIK-spectrum} is proved here.


\subsection{Explicit FIK formulas and geometric Hessian modes}
\label{app:self-contained-FIK-ledger}

Following Naff--Ozuch, we use the harmonic decomposition and geometric
modes employed in their FIK stability calculation.  We reproduce all
formulas needed here so that the spectral argument is self-contained.

\paragraph{The \(SU(2)\) normalization.}
Throughout this appendix we identify \(S^3\subset\mathbb C^2\) with the
unit-quaternion group \(SU(2)\), endowed with its unit round
bi-invariant metric.  Away from a measure-zero boundary, use the Euler
parametrization
\[
 (z_1,z_2)
 =
 \left(
  \cos\frac{\theta}{2}\,e^{\frac{i}{2}(\psi+\phi)},
  \sin\frac{\theta}{2}\,e^{\frac{i}{2}(\psi-\phi)}
 \right),
 \qquad
 0<\theta<\pi,\quad
 |\phi+\psi|<2\pi,\quad |\phi-\psi|<2\pi.
\]
Fix the left-invariant coframe
\begin{align*}
 \eta_1&=\frac12\bigl(d\psi+\cos\theta\,d\phi\bigr),\\
 \eta_2&=\frac12\bigl(\cos\psi\,d\theta
                 +\sin\theta\sin\psi\,d\phi\bigr),\\
 \eta_3&=\frac12\bigl(\sin\psi\,d\theta
                 -\sin\theta\cos\psi\,d\phi\bigr).
\end{align*}
Its dual left-invariant frame is
\begin{align*}
 X_1&=2\partial_\psi,\\
 X_2&=2\cos\psi\,\partial_\theta
      +\frac{2\sin\psi}{\sin\theta}
       \bigl(\partial_\phi-\cos\theta\,\partial_\psi\bigr),\\
 X_3&=2\sin\psi\,\partial_\theta
      -\frac{2\cos\psi}{\sin\theta}
       \bigl(\partial_\phi-\cos\theta\,\partial_\psi\bigr).
\end{align*}
Thus, for every cyclic permutation \((i,j,k)\) of \((1,2,3)\),
\[
 [X_i,X_j]=-2X_k,
 \qquad
 d\eta_i=2\eta_j\wedge\eta_k.
\]
We use the orientation and Haar volume
\[
 g_{S^3}=\eta_1^2+\eta_2^2+\eta_3^2,
 \qquad
 dV_{S^3}
 =\eta_1\wedge\eta_2\wedge\eta_3
 =\frac18\sin\theta\,d\theta\wedge d\psi\wedge d\phi,
 \qquad
 \int_{S^3}dV_{S^3}=2\pi^2.
\]
The \(X_1\)-flow is the Hopf action and has least positive period
\(2\pi\); equivalently, along a Hopf fiber the Euler variable \(\psi\)
has period \(4\pi\).  Consequently, in a local trivialization of the
normal line bundle of the bolt, the normal polar angle is
\[
 \vartheta=\frac{\psi}{2}\in\mathbb R/2\pi\mathbb Z.
\]

\begin{lemma}[Exact FIK coordinate formulas]
\label{lem:self-contained-FIK-ledger}
Put \(c_0=\sqrt2-1\).  On
\((1,\infty)\times S^3\), in the normalization fixed above, one has
\begin{align}
 \bar g
 &=\frac4{F(r)}\,dr^2
   +4r^2F(r)\eta_1^2
   +4r^2(\eta_2^2+\eta_3^2),                                      \label{eq:self-FIK-metric}\\
 F(r)
 &=\frac1{\sqrt2}-\frac{c_0}{r^2}
              -\frac{c_0}{\sqrt2\,r^4}
   =\frac{(r^2-1)(r^2+c_0)}{\sqrt2\,r^4},                          \label{eq:self-FIK-F}\\
 f_{\rm NO}(r)
 &=\sqrt2(r^2-1)-\log(2c_0),                                      \label{eq:self-FIK-fNO}\\
 dV_{\bar g}
 &=16r^3\,dr\,dV_{S^3}.                                           \label{eq:self-FIK-volume}
\end{align}
The potential satisfying
\(\bar R+|\bar\nabla\bar f|^2-\bar f=0\) is
\begin{equation}
 \bar f=f_{\rm NO}+C_{\rm NO}
 =\sqrt2\,r^2+2-2\sqrt2,\qquad
 C_{\rm NO}=\log(2c_0)+\sqrt2\,c_0.                               \label{eq:self-FIK-fbar}
\end{equation}
In particular,
\begin{equation}
 d\nu=e^{-C_{\rm NO}}(4\pi)^{-2}
       e^{-f_{\rm NO}}\,dV_{\bar g},\qquad
 \int_Me^{-f_{\rm NO}}\,dV_{\bar g}=(4\pi)^2.                     \label{eq:self-FIK-measure}
\end{equation}
If
\[
 s^2=r^2F,\qquad d=\frac12rF',
\]
then, for \(r>1\),
\begin{equation}
 0<F<\frac1{\sqrt2},\qquad
 0<d<1,\qquad rF'<2,\qquad
 1<d+2F<\sqrt2.                                                   \label{eq:self-FIK-ranges}
\end{equation}
Moreover,
\begin{equation}
 \bar R=\frac{\sqrt2c_0}{r^2}>0,\qquad
 |\bar\nabla\bar f|^2=2r^2F=\bar f-\bar R.                        \label{eq:self-FIK-scalar}
\end{equation}
The curvature is bounded, and for every \(m\geq0\),
\begin{equation}
 |\bar\nabla^m\overline{\Rm}|_{\bar g}=O(r^{-2-m}).                \label{eq:self-FIK-curvature-symbol}
\end{equation}
Finally, if \(g_{\rm C}\) is obtained from
\eqref{eq:self-FIK-metric} by replacing \(F\) by \(1/\sqrt2\), then
\begin{equation}
 \left\|\delta_R^*(R^{-2}\bar g)-g_{\rm C}\right\|
 _{C^m_{g_{\rm C}}(\{1\leq r\leq2\})}
 \leq C_mR^{-2}\qquad(R\geq2).                                   \label{eq:self-FIK-AC}
\end{equation}
\end{lemma}

\begin{proof}
The orthonormal coframe on \(r>1\) is
\[
 e^0=\frac2{\sqrt F}\,dr,\qquad
 e^1=2r\sqrt F\,\eta_1,\qquad
 e^2=2r\eta_2,\qquad e^3=2r\eta_3.
\]
Its wedge product proves \eqref{eq:self-FIK-volume}.  Since
\(|S^3|=2\pi^2\), the substitution \(y=r^2-1\) gives
\[
\begin{split}
 \int_Me^{-f_{\rm NO}}dV_{\bar g}
 &=64c_0\pi^2\int_1^\infty
       r^3e^{-\sqrt2(r^2-1)}\,dr\\
 &=16c_0\pi^2(1+\sqrt2)=16\pi^2.
\end{split}
\]
This proves \eqref{eq:self-FIK-measure}.  The additive shift in
\eqref{eq:self-FIK-fbar} is obtained by subtracting the constant
\(\log(2c_0)+\sqrt2c_0\) in
\(\bar R+|\bar\nabla f_{\rm NO}|^2-f_{\rm NO}\).

Direct differentiation gives
\[
 F'=\frac{2c_0}{r^3}+\frac{2\sqrt2c_0}{r^5},\qquad
 F''=-\frac{6c_0}{r^4}-\frac{10\sqrt2c_0}{r^6}.
\]
Thus \(F>0\) by the factorization in
\eqref{eq:self-FIK-F}, while \(F<1/\sqrt2\) is immediate from its
first expression.  Also
\[
 d=\frac{c_0}{r^2}+\frac{\sqrt2c_0}{r^4}
\]
is strictly decreasing from \(1\) to \(0\), and the useful last
combination is the exact identity
\[
 d+2F=\sqrt2-\frac{c_0}{r^2}.
\]
This proves \eqref{eq:self-FIK-ranges}.

For completeness, use the convention
\[
 R_{ijkl}
 =\left\langle
  \bigl(\bar\nabla_{e_j}\bar\nabla_{e_i}
       -\bar\nabla_{e_i}\bar\nabla_{e_j}
       -\bar\nabla_{[e_j,e_i]}\bigr)e_k,e_l
 \right\rangle .
\]
The nonzero curvature entries needed here, up to the usual symmetries,
are
\begin{align*}
 R_{0101}&=-\frac18F''-\frac{3}{8r}F'
          =\frac{c_0}{\sqrt2r^6},\\
 R_{0202}=R_{0303}=R_{1212}=R_{1313}
         &=-\frac1{8r}F'
          =-\frac{c_0}{2\sqrt2r^6}-\frac{c_0}{4r^4},\\
 R_{2323}&=\frac{1-F}{r^2}
          =\frac{c_0}{\sqrt2r^6}+\frac{c_0}{r^4}
             +\frac{c_0}{\sqrt2r^2},\\
 R_{0213}&=R_{0202},\qquad
 R_{0312}=-R_{0202},\qquad R_{0123}=2R_{0202}.
\end{align*}
Their traces give
\[
 \Ric_{00}=\Ric_{11}=-\frac{c_0}{2r^4},\qquad
 \Ric_{22}=\Ric_{33}
 =\frac{c_0}{2r^4}+\frac{\sqrt2c_0}{2r^2},
\]
and hence \eqref{eq:self-FIK-scalar}.  Direct differentiation of
\eqref{eq:self-FIK-fNO} gives
\[
 \bar\nabla f_{\rm NO}=\sqrt2s\,e_0
\]
and
\[
 (\bar\nabla^2f_{\rm NO})_{00}
 =(\bar\nabla^2f_{\rm NO})_{11}
 =\frac12+\frac{c_0}{2r^4},
\qquad
 (\bar\nabla^2f_{\rm NO})_{22}
 =(\bar\nabla^2f_{\rm NO})_{33}
 =\frac12-\frac{c_0}{2r^4}
              -\frac{\sqrt2c_0}{2r^2}.
\]
Thus \(\Ric+\bar\nabla^2f_{\rm NO}=\frac12\bar g\) componentwise.
Moreover,
\[
 \bar R+|\bar\nabla f_{\rm NO}|^2-f_{\rm NO}
 =\log(2c_0)+\sqrt2c_0,
\]
which verifies the shift in \eqref{eq:self-FIK-fbar} without importing
a soliton normalization.  Since an additive shift does not change the
gradient, \(|\bar\nabla\bar f|^2=2s^2\), and substitution of
\eqref{eq:self-FIK-F} proves the second identity in
\eqref{eq:self-FIK-scalar}.

To see the bolt without importing a coordinate assertion, introduce
the normal distance
\[
 \rho=\int_1^r\frac2{\sqrt{F(t)}}\,dt.
\]
Then
\[
 r=1+\frac{\rho^2}{8}+O(\rho^4),\quad
 4r^2F=\rho^2+O(\rho^4),\quad
 4r^2=4+\rho^2+O(\rho^4).
\]
The exact function \(F\) is analytic and has a simple zero at \(r=1\).
Consequently \(r=r(\rho)\) extends as a smooth even function, while
\[
 a(\rho):=2r(\rho)\sqrt{F(r(\rho))}
 \quad\hbox{and}\quad
 b(\rho):=2r(\rho)
\]
extend respectively as a smooth odd function and a smooth even
function, with
\[
 a(0)=0,\qquad a'(0)=1,\qquad b(0)=2.
\]
Moreover, on a local Hopf trivialization,
\[
 \eta_1=d\vartheta+\frac12\cos\theta\,d\phi,
 \qquad \vartheta\in\mathbb R/2\pi\mathbb Z.
\]
These parity, slope, and period identities are precisely the smooth
polar-coefficient criterion for the normal disk bundle
\(\mathcal O_{\mathbb P^1}(-1)\).  Hence there is neither a conical nor
an orbifold singularity, and \eqref{eq:self-FIK-metric} extends
smoothly across \(\rho=0\).
The displayed curvature entries are bounded there.  At infinity,
\[
 F-\frac1{\sqrt2}=O_{\rm sym}(r^{-2}),\qquad
 |(r\partial_r)^jO_{\rm sym}(r^{-2})|\leq C_jr^{-2}.
\]
The connection coefficients are classical \(r^{-1}\)-symbols.
Differentiating the curvature formulas proves
\eqref{eq:self-FIK-curvature-symbol}, and the same symbol calculation
after dilation proves \eqref{eq:self-FIK-AC}.
\end{proof}

\begin{lemma}[Three scalar eigenpotential families]
\label{lem:self-contained-scalar-modes}
In the Euler angles fixed above, one has
\[
 X_1=2\partial_\psi,\qquad
 D^J_{M,M'}(\psi,\theta,\phi)
 =e^{-iM\psi/2}d^J_{M,M'}(\theta)e^{-iM'\phi/2}.
\]
For the only values needed in this lemma, the little Wigner matrices,
with rows and columns ordered increasingly, are
\[
 d^1=
 \begin{pmatrix}
  \cos(\theta/2)&-\sin(\theta/2)\\
  \sin(\theta/2)&\cos(\theta/2)
 \end{pmatrix}
\]
for the labels \(-1,1\), and
\[
 d^2=
 \begin{pmatrix}
 \frac{1+\cos\theta}{2}&-\frac{\sin\theta}{\sqrt2}
                       &\frac{1-\cos\theta}{2}\\
 \frac{\sin\theta}{\sqrt2}&\cos\theta
                       &-\frac{\sin\theta}{\sqrt2}\\
 \frac{1-\cos\theta}{2}&\frac{\sin\theta}{\sqrt2}
                       &\frac{1+\cos\theta}{2}
 \end{pmatrix}
\]
for the labels \(-2,0,2\).  The latter is the symmetric-square phase
convention induced by the displayed \(d^1\).  In particular, direct
differentiation with the fixed Euler frame gives the sign check
\[
 (X_3+iX_2)D^1_{-1,-1}=-2iD^1_{1,-1},
\]
and direct differentiation also gives
\[
 X_1D^J_{M,M'}=-iM D^J_{M,M'},\qquad
 (X_1^2+X_2^2+X_3^2)D^J_{M,M'}
 =-J(J+2)D^J_{M,M'}.
\]
For a radial function \(a\),
\begin{align}
 \bar\Delta_{\bar f}a
 &=\frac{e^{f_{\rm NO}}}{4r^3}
   \left(r^3Fe^{-f_{\rm NO}}a'\right)',                           \label{eq:self-radial-drift}\\
 \bar\Delta_{\bar f}(aD^J_{M,M'})
 &=\left(\bar\Delta_{\bar f}a+
 \left\{\frac{M^2-J(J+2)}{4r^2}
             -\frac{M^2}{4s^2}\right\}a\right)D^J_{M,M'}.          \label{eq:self-separated-drift}
\end{align}
Define
\begin{equation}
 \widehat u=(r^2-1)^{1/2}(r^2+c_0)^{c_0/2},\qquad
 \widehat v=r^2,\qquad
 \widehat w=\widehat u^2.                                        \label{eq:self-uvw}
\end{equation}
Then
\begin{align}
 \bar\Delta_{\bar f}(\widehat uD^1_{\pm1,M'})
 &=-\frac1{\sqrt2}\widehat uD^1_{\pm1,M'},                        \label{eq:self-u-eigen}\\
 \bar\Delta_{\bar f}(\widehat vD^2_{0,M'})
 &=-\widehat vD^2_{0,M'},                                        \label{eq:self-v-eigen}\\
 \bar\Delta_{\bar f}(\widehat wD^2_{\pm2,M'})
 &=-\sqrt2\,\widehat wD^2_{\pm2,M'}.                              \label{eq:self-w-eigen}
\end{align}
Here \(M'\) ranges over \(\{-1,1\}\) in
\eqref{eq:self-u-eigen} and over \(\{-2,0,2\}\) in the other two
displays.
\end{lemma}

\begin{proof}
Equation \eqref{eq:self-radial-drift} is the weighted divergence
formula using \eqref{eq:self-FIK-volume}; the angular part of the
inverse metric gives \eqref{eq:self-separated-drift}.  Logarithmic
differentiation and \(s^2=r^2F\) give
\[
 \widehat u'=\frac{\widehat u}{rF},\qquad
 \widehat v'=2r=\frac{2\widehat v}{r},\qquad
 \widehat w'=\frac{2\widehat w}{rF}.
\]
Substitution into \eqref{eq:self-radial-drift} yields
\begin{align*}
 \bar\Delta_{\bar f}\widehat u
 &=\left(\frac1{4s^2}+\frac1{2r^2}
                    -\frac1{\sqrt2}\right)\widehat u,\\
 \bar\Delta_{\bar f}\widehat v
 &=2-r^2=\left(\frac2{r^2}-1\right)\widehat v,\\
 \bar\Delta_{\bar f}\widehat w
 &=\left(\frac1{s^2}+\frac1{r^2}-\sqrt2\right)\widehat w.
\end{align*}
For \((J,M)=(1,\pm1),(2,0),(2,\pm2)\), respectively, the
braced angular coefficient in \eqref{eq:self-separated-drift} is
\[
 -\frac1{4s^2}-\frac1{2r^2},\qquad
 -\frac2{r^2},\qquad
 -\frac1{s^2}-\frac1{r^2}.
\]
The cancellations prove
\eqref{eq:self-u-eigen}--\eqref{eq:self-w-eigen}.
\end{proof}

\begin{lemma}[Bolt extension and real dimensions]
\label{lem:self-contained-bolt-reality}
The scalar eigenfunctions in
\eqref{eq:self-u-eigen}--\eqref{eq:self-w-eigen} extend smoothly
across the exceptional divisor and lie in every weighted Sobolev
space.  Their real spans have dimensions
\[
4\quad\text{for }\widehat uD^1_{\pm1,M'},\qquad
3\quad\text{for }\widehat vD^2_{0,M'},\qquad
6\quad\text{for }\widehat wD^2_{\pm2,M'}.
\]
\end{lemma}

\begin{proof}
With \(\rho\) as in the proof of
Lemma~\ref{lem:self-contained-FIK-ledger},
\[
 \widehat u=C_u\rho(1+O(\rho^2)),\qquad
 \widehat w=C_u^2\rho^2(1+O(\rho^2)),\qquad
 \widehat v=1+\frac{\rho^2}{4}+O(\rho^4),
\]
where every remainder is smooth and even in \(\rho\).
If \(\vartheta\) is the normal-circle angle (so that the Euler angle
used in the Wigner coefficient is \(2\vartheta\)), then its Hopf
factor is \(e^{-iM\vartheta}\).  Thus the potentially singular
products have local normal factors
\[
 \rho e^{\mp i\vartheta},\qquad
 \rho^2e^{\mp2i\vartheta},
\]
which are \(z,\bar z,z^2,\bar z^2\) in a smooth complex normal
coordinate.  The \(M=0\) coefficients are invariant under the
collapsing circle and descend smoothly to the bolt.

This fiberwise description is compatible with the global
\(U(1)\)-gluing.  Let \(s_\alpha,s_\beta\) be local Hopf sections with
\[
 s_\beta=s_\alpha\cdot e^{i\chi_{\alpha\beta}},
\]
and write
\[
 g=s_\alpha(x)\cdot e^{i\vartheta_\alpha}
  =s_\beta(x)\cdot e^{i\vartheta_\beta}.
\]
Put
\[
 d_{\alpha,M,M'}^J(x):=D_{M,M'}^J(s_\alpha(x)),
 \qquad
 d_{\beta,M,M'}^J(x):=D_{M,M'}^J(s_\beta(x)).
\]
Then
\[
 \vartheta_\beta=\vartheta_\alpha-\chi_{\alpha\beta},
 \qquad
 d_{\beta,M,M'}^J
 =e^{-iM\chi_{\alpha\beta}}d_{\alpha,M,M'}^J,
\]
the second identity following from
\(X_1D_{M,M'}^J=-iM D_{M,M'}^J\).  For the normal coordinates
\[
 z_\alpha=\rho e^{i\vartheta_\alpha}
 \quad\hbox{on}\quad \mathcal O_{\mathbb P^1}(-1),
\]
one has
\(z_\beta=e^{-i\chi_{\alpha\beta}}z_\alpha\).  Consequently,
\[
 \rho^{|M|}D_{M,M'}^J=
 \begin{cases}
  \overline z_\alpha^{\,M}d_{\alpha,M,M'}^J,&M>0,\\
  d_{\alpha,0,M'}^J,&M=0,\\
  z_\alpha^{-M}d_{\alpha,M,M'}^J,&M<0.
 \end{cases}
\]
The displayed transition laws show that these expressions agree on
overlaps.  The remaining analytic radial factors are smooth even
functions of \(\rho\), hence smooth functions of
\(|z_\alpha|^2\), and glue as well.  This proves global smooth
extension across the zero section.  At infinity all these functions
and all their derivatives have polynomial growth, whereas the measure
has Gaussian factor \(e^{-\sqrt2r^2}\); hence they belong to every
\(H^k_\nu\).

The Wigner conjugation rule is
\begin{equation}
 \overline{D^J_{M,M'}}
 =(-1)^{(M-M')/2}D^J_{-M,-M'}.                                   \label{eq:self-Wigner-reality}
\end{equation}
For \(J=1\), the four complex labels form two conjugate pairs, and
their real and imaginary parts give four real functions.  For
\(J=2,M=0\), the \(M'=0\) coefficient is real and the \(M'=\pm2\)
coefficients form one conjugate pair, giving three real functions.
For \(J=2,M=\pm2\), the six labels form three conjugate pairs, giving
six real functions.  Orthogonality of distinct matrix coefficients
proves linear independence.
\end{proof}

\begin{corollary}[The eight geometric modes and a six-real stable family]
\label{cor:self-contained-Hessian-families}
On every gradient shrinker with
\(\Ric+\nabla^2f=\frac12g\), the commutation formula
\begin{equation}
 (\Delta_f+2\Rm)\nabla^2\phi
 =\nabla^2(\Delta_f\phi+\phi)                                    \label{eq:self-Hessian-commute}
\end{equation}
holds.  Consequently the following are eigentensors of \(\A\):
\begin{align*}
 \nabla^2(\widehat uD^1_{\pm1,M'})
 &:&&1-\frac1{\sqrt2}
 &&\text{(a four-real-dimensional family)},\\
 \nabla^2(\widehat vD^2_{0,M'})
 &:&&0
 &&\text{(a three-real-dimensional family)},\\
 \nabla^2\bar f
 &:&&0
 &&\text{(one additional radial tensor)},\\
\nabla^2(\widehat wD^2_{\pm2,M'})
&:&&1-\sqrt2
&&\text{(a six-real-dimensional stable family)}.
\end{align*}
All these tensors are smooth and belong to every \(H^k_\nu\).
\end{corollary}

\begin{proof}
Commuting two covariant derivatives and inserting
\(\Ric+\nabla^2f=\frac12g\) proves
\eqref{eq:self-Hessian-commute}.  Apply that identity to
\eqref{eq:self-u-eigen}--\eqref{eq:self-w-eigen}.  In addition,
\(\Delta_{\bar f}\bar f=2-\bar f\), so taking a Hessian gives
\(\A\nabla^2\bar f=0\).

The Hessian map is injective on each of the displayed nontrivial
Wigner families.  Indeed, if
\(\Delta_{\bar f}\phi=-\mu\phi\), weighted Bochner gives
\[
 \mu^2\|\phi\|_{L^2_\nu}^2
 =\|\bar\nabla^2\phi\|_{L^2_\nu}^2
   +\frac12\|\bar\nabla\phi\|_{L^2_\nu}^2
 =\|\bar\nabla^2\phi\|_{L^2_\nu}^2
   +\frac{\mu}{2}\|\phi\|_{L^2_\nu}^2.
\]
For \(\mu\in\{1/\sqrt2,1,\sqrt2\}\), none of which equals \(1/2\),
\(\bar\nabla^2\phi=0\) therefore forces \(\phi=0\).
Hence the real dimensions from
Lemma~\ref{lem:self-contained-bolt-reality} are unchanged after
taking Hessians.  The radial tensor
\(\nabla^2\bar f\) belongs to the \(J=0\) block and is independent of
the three \(J=2\) zero modes.
\end{proof}

\begin{lemma}[All-order growth of the geometric generators]
\label{lem:self-contained-geometric-growth}
Let \(\phi\) be any one of the four real
\(\widehat uD^1_{\pm1,M'}\) potentials, any one of the three real
\(r^2D^2_{0,M'}\) potentials, or \(\bar f\).  Then, for every
\(m\geq0\),
\begin{equation}
 |\bar\nabla^{m+1}\phi|
 \leq C_m(1+r)^{1-m},\qquad
 |\bar\nabla^{m+2}\phi|
 \leq C_m(1+r)^{-m}.                                              \label{eq:self-geometric-growth}
\end{equation}
Equivalently, with
\(W=\frac12\bar\nabla\phi\) and \(Z=\bar\nabla^2\phi\),
\[
 |\bar\nabla^mW|\leq C_m(1+\bar f)^{(1-m)/2},\qquad
 |\bar\nabla^mZ|\leq C_m(1+\bar f)^{-m/2}.
\]
The scale mode obeys the stronger estimate
\begin{equation}
 |\bar\nabla^m\Ric_{\bar g}|
 \leq C_m(1+r)^{-2-m}
 \leq C_m(1+\bar f)^{-m/2}.                                      \label{eq:self-scale-mode-growth}
\end{equation}
For the six-real stable family,
\begin{equation}
 \left|\bar\nabla^m\nabla^2
   (\widehat wD^2_{\pm2,M'})\right|
 \leq C_m(1+r)^{2\sqrt2-2-m}.                                    \label{eq:self-stable-family-growth}
\end{equation}
\end{lemma}

\begin{proof}
The exact radial factors give the all-order symbols
\begin{align*}
 \widehat u
 &=r^{\sqrt2}(1-r^{-2})^{1/2}
       (1+c_0r^{-2})^{c_0/2}
   =r^{\sqrt2}\bigl(1+O_{\rm sym}(r^{-2})\bigr),\\
 \widehat w
 &=r^{2\sqrt2}(1-r^{-2})
       (1+c_0r^{-2})^{c_0}
   =r^{2\sqrt2}\bigl(1+O_{\rm sym}(r^{-2})\bigr).
\end{align*}
On the AC end, the orthonormal angular derivatives cost one factor of
\(r^{-1}\), radial unit differentiation has the same order, and the
connection is an \(r^{-1}\)-symbol.  Thus a Wigner mode with radial
order \(r^p\) has \(m\)-th covariant derivative \(O(r^{p-m})\).
Use \(p=\sqrt2\) for the positive modes and \(p=2\) for the quadratic
and radial zero modes; since \(\sqrt2<2\), this proves
\eqref{eq:self-geometric-growth}.  Use \(p=2\sqrt2\) for
\eqref{eq:self-stable-family-growth}.  Smooth bolt extension supplies
the uniform bounds on the compact core, and
\(1+\bar f\simeq1+r^2\) converts the estimates to the displayed
\(\bar f\)-weights.  Finally,
\eqref{eq:self-scale-mode-growth} is the Ricci contraction of
\eqref{eq:self-FIK-curvature-symbol}.
\end{proof}

\begin{remark}[Scope of the coordinate formulas]
\label{rem:self-contained-ledger-scope}
Corollary~\ref{cor:self-contained-Hessian-families} proves directly
the existence, eigenvalues, smoothness, real dimensions, and growth
of every explicit Hessian family used in the body of the paper.  The
assertion that these exhaust the corresponding nonnegative
eigenspaces is logically separate: it follows only after the
\(J=0\), exceptional \(J=1,2\), remaining \(J=1,\ldots,4\), and
\(J\geq5\) exclusion arguments have all been supplied.
\end{remark}

\subsection{Wigner reduction}
\label{app:self-contained-wigner}

This subsection records the geometric and representation-theoretic
calculation used in the spectral analysis.  In particular, the
matrices below are defined entry by entry from the metric and do not
rely on an external computer algebra file.  Our convention for a
complexified
tensor is
\[
 |T|^2=\langle T,\overline T\rangle ,
\]
where the metric and all differential operators are extended complex
linearly.  We use
\[
 \Delta_f=\operatorname{tr}\nabla^2-\nabla_{\nabla f},
 \qquad
 \mathrm{Rm}(h)_{ij}=R_{ikj\ell}h_{k\ell},
\]
with
\[
 R(X,Y)Z=\nabla_X\nabla_YZ-\nabla_Y\nabla_XZ-\nabla_{[X,Y]}Z,
\]
and \(R_{ijkl}=\langle R(e_j,e_i)e_k,e_l\rangle\).  In particular,
the displayed index order gives \(\mathrm{Rm}(g)=\mathrm{Ric}\), in
agreement with the curvature table below,
so that \(L_f=\Delta_f+2\mathrm{Rm}\) has quadratic form
\[
 \int_M\left(2\mathrm{Rm}(h,\overline h)-|\nabla h|^2\right)
 e^{-f}\,dV_g.
\]

\subsubsection{Frame and tensor basis}

Retain \(c_0,F,s,\bar g,f_{\rm NO}\), the orthonormal coframe from
the preceding coordinate formulas, and the left-invariant frame
\(X_1,X_2,X_3\) with dual coframe
\(\eta_1,\eta_2,\eta_3\) fixed above.  Thus, on the open orbit
\(\mathring M=(1,\infty)\times S^3\),
\[
 [X_i,X_j]=-2X_k
 \qquad\text{when }(i,j,k)\text{ is cyclic},
\]
Set
\[
 g:=\bar g,\qquad f:=f_{\rm NO},
\]
throughout this calculation.  Then
\[
 e_0=\frac{\sqrt F}{2}\partial_r=\frac{s}{2r}\partial_r,\qquad
 e_1=\frac1{2s}X_1,\qquad
 e_2=\frac1{2r}X_2,\qquad
 e_3=\frac1{2r}X_3
\]
is an orthonormal frame, with dual coframe
\[
 e^0=\frac2{\sqrt F}\,dr,\qquad
 e^1=2s\eta_1,\qquad e^2=2r\eta_2,\qquad e^3=2r\eta_3.
\]
The profile identities used in the connection calculation are
\begin{equation}
\label{eq:self-fik-profile-derivatives}
 F'=\frac{2c_0}{r^3}+\frac{2\sqrt2c_0}{r^5},
 \qquad
 F''=-\frac{6c_0}{r^4}-\frac{10\sqrt2c_0}{r^6}.
\end{equation}

For later verification, the nonzero frame brackets, up to
antisymmetry, are
\begin{align}
[e_0,e_1]&=-\frac{s'}{2r}e_1,&
[e_0,e_2]&=-\frac{s}{2r^2}e_2,&
[e_0,e_3]&=-\frac{s}{2r^2}e_3,\nonumber\\
[e_1,e_2]&=-\frac1s e_3,&
[e_1,e_3]&=\frac1s e_2,&
[e_2,e_3]&=-\frac{s}{r^2}e_1.
\label{eq:self-fik-brackets}
\end{align}
Consequently all the connection identities below can also be checked
directly from the Koszul formula.

Let
\[
 \omega_1^\pm=e^0\wedge e^1\pm e^2\wedge e^3,\quad
 \omega_2^\pm=e^0\wedge e^2\pm e^3\wedge e^1,\quad
 \omega_3^\pm=e^0\wedge e^3\pm e^1\wedge e^2,
\]
where \(e^i\wedge e^j=(e^i\otimes e^j-e^j\otimes e^i)/2\).
If \(\alpha,\beta\) are two-forms, write
\(\alpha\circ\beta\) for the symmetric tensor obtained by tracing
\(\alpha\otimes\beta\) in its middle two indices; explicitly,
\[
 (\alpha\circ\beta)_{ik}
 =\frac12\left(\alpha_{ip}\beta_{kp}
                    +\alpha_{kp}\beta_{ip}\right).
\]
Define
\begin{align*}
 b_0&=\omega_1^+\circ\omega_1^+,&
 b_1&=\omega_1^-\circ\omega_1^+,\\
 b_2&=2^{-1/2}(\omega_2^-+\omega_3^-)\circ\omega_1^+,&
 b_3&=2^{-1/2}(\omega_2^--\omega_3^-)\circ\omega_1^+,\\
 b_4&=2^{-1/2}\omega_1^-\circ(\omega_2^++\omega_3^+),&
 b_5&=2^{-1/2}\omega_1^-\circ(\omega_2^+-\omega_3^+),\\
 b_6&=2^{-1/2}(\omega_2^-\circ\omega_2^+
                     +\omega_3^-\circ\omega_3^+),&
 b_7&=2^{-1/2}(\omega_3^-\circ\omega_2^+
                     -\omega_2^-\circ\omega_3^+),\\
 b_8&=2^{-1/2}(\omega_3^-\circ\omega_2^+
                     +\omega_2^-\circ\omega_3^+),&
 b_9&=2^{-1/2}(\omega_2^-\circ\omega_2^+
                     -\omega_3^-\circ\omega_3^+).
\end{align*}
These tensors are pairwise orthogonal and
\(\langle b_p,b_q\rangle=\delta_{pq}/4\).
The first four are invariant, and the last six anti-invariant, under
the K\"ahler complex structure.  Introduce the unitary combinations
\begin{align*}
 b_\pm&=2^{-1/2}(b_2\pm i b_3),\\
 \mathbf k_1&=2^{-1/2}(b_4+i b_5),&
 \mathbf k_2&=2^{-1/2}(b_6+i b_7),&
 \mathbf k_3&=2^{-1/2}(b_8+i b_9),
\end{align*}
and put \(\mathbf k_{\bar a}=\overline{\mathbf k_a}\).  Thus a
complexified symmetric tensor has the unique form
\begin{equation}
\label{eq:self-complex-tensor-expansion}
 h=h_0b_0+h_1b_1+h_-b_++h_+b_-
   +\sum_{a=1}^3(\kappa_a\mathbf k_a+
                  \kappa_{\bar a}\mathbf k_{\bar a}).
\end{equation}
Below we usually abbreviate the coefficient \(\kappa_a\) to \(k_a\);
the boldface symbol always denotes a basis tensor.

For completeness, define
\[
 e_\pm=\frac{e_2\pm e_3}{\sqrt2},\qquad
 D_\pm=\frac{e_+\pm i e_-}{\sqrt2}
      =\frac{1\mp i}{4r}(X_3\pm iX_2),
\]
and let \(\sigma^\pm\) be the complex covectors dual to \(D_\pm\).
The three connection coefficients that enter the calculation are
\begin{equation}
\label{eq:self-gamma}
 \Gamma_1^+=\frac{rF'+4F-4}{4s},\qquad
 \Gamma_{23}^-=\frac{rF'+4}{4s},\qquad
 \Gamma_1^-=\frac Fs .
\end{equation}
The connection on the invariant basis is
\begin{align}
 \nabla b_0&=0,\nonumber\\
 \nabla b_1&=i\Gamma_1^-\sigma^+\otimes b_+
             -i\Gamma_1^-\sigma^-\otimes b_-,\nonumber\\
 \nabla b_+&=-i\Gamma_{23}^-e^1\otimes b_+
             +i\Gamma_1^-\sigma^-\otimes b_1,\nonumber\\
 \nabla b_-&= i\Gamma_{23}^-e^1\otimes b_-
             -i\Gamma_1^-\sigma^+\otimes b_1.
\label{eq:self-invariant-connection}
\end{align}
On the unbarred anti-invariant basis it is
\begin{align}
 \nabla\mathbf k_1
 &=i\Gamma_1^+e^1\otimes\mathbf k_1
   -i\Gamma_1^-\sigma^-\otimes\mathbf k_2
   +i\Gamma_1^-\sigma^+\otimes\mathbf k_3,\nonumber\\
 \nabla\mathbf k_2
 &=i(\Gamma_1^++\Gamma_{23}^-)e^1\otimes\mathbf k_2
   -i\Gamma_1^-\sigma^+\otimes\mathbf k_1,\nonumber\\
 \nabla\mathbf k_3
 &=i(\Gamma_1^+-\Gamma_{23}^-)e^1\otimes\mathbf k_3
   +i\Gamma_1^-\sigma^-\otimes\mathbf k_1,
\label{eq:self-anti-connection}
\end{align}
and the barred formulas are their complex conjugates.
The curvature endomorphism has the following action:
\begin{align}
 \mathrm{Rm}(b_0)&=\frac{\sqrt2c_0}{4r^2}b_0
                  -\frac{rF'}{4\sqrt2}b_1,\nonumber\\
 \mathrm{Rm}(b_1)&=-\frac{rF'}{4\sqrt2}b_0
 +\left(\frac{\sqrt2c_0}{4r^2}+\frac{rF'}{2r^2}\right)b_1,\nonumber\\
 \mathrm{Rm}(b_\pm)&=-\frac{rF'}{4r^2}b_\pm,\nonumber\\
 \mathrm{Rm}(\mathbf k_1)&=\frac{rF'}{2r^2}\mathbf k_1,\qquad
 \mathrm{Rm}(\mathbf k_2)
 =-\frac{rF'+2F-\sqrt2}{2r^2}\mathbf k_2,\qquad
 \mathrm{Rm}(\mathbf k_3)=-\frac{1-F}{r^2}\mathbf k_3,
\label{eq:self-curvature-action}
\end{align}
with the same diagonal entries in the barred sector.
Equations \eqref{eq:self-invariant-connection}--\eqref{eq:self-curvature-action}
follow by inserting \eqref{eq:self-fik-brackets} into the Koszul and
curvature formulas.  They are included to make every entry of the
matrices below independently checkable.

\subsubsection{Wigner functions and the orthogonal splitting}

For \(J\in\mathbb N_0\), put
\[
 \mathcal J_J=\{-J,-J+2,\ldots,J-2,J\}
\]
and, for \(M\in\mathcal J_J\), set
\begin{equation}
\label{eq:self-CJM}
 C^J_{M\pm}
 :=\sqrt{J(J+2)-M(M\pm2)}
 =\sqrt{(J+1)^2-(M\pm1)^2}.
\end{equation}
We use the convention that a Wigner function with an index outside
\(\mathcal J_J\) is zero.  There are smooth functions
\(D^J_{M,M'}\) on \(S^3\), \(M,M'\in\mathcal J_J\), satisfying
\begin{align}
 X_1D^J_{M,M'}&=-iM D^J_{M,M'},\nonumber\\
 (X_3\pm iX_2)D^J_{M,M'}&=-iC^J_{M\pm}D^J_{M\pm2,M'},\nonumber\\
 (X_1^2+X_2^2+X_3^2)D^J_{M,M'}
 &=-J(J+2)D^J_{M,M'}.
\label{eq:self-wigner-derivatives}
\end{align}
With the Haar measure \(dV_{S^3}\) fixed above, their normalization is
\begin{equation}
\label{eq:self-wigner-orthogonality}
 \int_{S^3}D^J_{M,M'}\overline{D^{K}_{N,N'}}\,dV_{S^3}
 =\frac{2\pi^2}{J+1}
   \delta_{JK}\delta_{MN}\delta_{M'N'}.
\end{equation}
They form a complete orthogonal basis of \(L^2(S^3)\), and their
phases may be chosen so that
\[
 \overline{D^J_{M,M'}}
 =(-1)^{(M-M')/2}D^J_{-M,-M'}.
\]
These statements are the Peter--Weyl theorem for \(SU(2)\), with
highest weight \(J\), in the normalization
\([X_i,X_j]=-2X_k\).

Equations \eqref{eq:self-wigner-derivatives} and the definitions of
the FIK frame give, for a radial scalar \(u\),
\begin{align}
 e_0(uD^J_{M,M'})
 &=\frac{\sqrt F}{2}u'D^J_{M,M'},&
 e_1(uD^J_{M,M'})
 &=-\frac{iM}{2s}uD^J_{M,M'},\nonumber\\
 D_+(uD^J_{M,M'})
 &=-\frac{1+i}{4r}C^J_{M+}uD^J_{M+2,M'},&
 D_-(uD^J_{M,M'})
 &= \frac{1-i}{4r}C^J_{M-}uD^J_{M-2,M'}.
\label{eq:self-fik-wigner-derivatives}
\end{align}
In particular,
\begin{equation}
\label{eq:self-berger-eigenvalue}
 \widetilde\Delta^J_M
 :=\frac{M^2-J(J+2)}{4r^2}-\frac{M^2}{4s^2}
 =-\frac1{4r^2}
  \left[J(J+2)+\left(\frac1F-1\right)M^2\right]
\end{equation}
is the angular scalar Laplacian eigenvalue on the Berger orbit.

Fix \(J\) and \(M'\).  Let \(\mathcal T^J_{M'}\) consist of tensors
whose ten coefficients in \eqref{eq:self-complex-tensor-expansion}
are sums
\[
 u(r,\theta)=\sum_{M\in\mathcal J_J}u_M(r)D^J_{M,M'}(\theta).
\]
The Peter--Weyl expansion gives the Hilbert direct sum
\begin{equation}
\label{eq:self-mode-splitting}
 L^2_f(S^2T^*M;\mathbb C)
 =\widehat{\bigoplus}_{J\geq0}
   \widehat{\bigoplus}_{M'\in\mathcal J_J}\mathcal T^J_{M'}.
\end{equation}
The same splitting holds for the weighted \(H^1\) form domain.
Indeed, \eqref{eq:self-fik-wigner-derivatives} preserves \(J,M'\),
and the tensor connection coefficients in
\eqref{eq:self-invariant-connection}--\eqref{eq:self-anti-connection}
are radial.  Hence \(\nabla\), \(\mathrm{Rm}\), and
\[
 L_f:=\Delta_f+2\mathrm{Rm}
\]
preserve every \(\mathcal T^J_{M'}\).  Orthogonality of both the
\(L^2_f\) inner product and the quadratic form follows from
\eqref{eq:self-wigner-orthogonality}.  For real tensors the
\(M'\) and \(-M'\) summands are paired by complex conjugation; it is
therefore enough to analyze one complex summand.  We first establish
these assertions for smooth tensors supported in \(\mathring M\).
Radial cutoff, angular projection, and closure in the form norm give the stated
decomposition across the collapsing orbit at \(r=1\).  The capacity
argument at that orbit is as follows.  In the smooth bolt distance
\(\rho\) from Lemma~\ref{lem:self-contained-FIK-ledger}, let
\[
 \chi_\varepsilon(\rho)=
 \begin{cases}
 0,&0\leq\rho\leq\varepsilon^2,\\[2pt]
 \displaystyle
 \frac{\log(\rho/\varepsilon^2)}{|\log\varepsilon|},
   &\varepsilon^2<\rho<\varepsilon,\\[7pt]
 1,&\rho\geq\varepsilon .
 \end{cases}
\]
After harmless smoothing at the two corners,
\[
 |d\chi_\varepsilon|
 \leq\frac{C}{\rho|\log\varepsilon|},\qquad
 \int_{\{\varepsilon^2<\rho<\varepsilon\}}
 |d\chi_\varepsilon|^2\,dV_g
 \leq\frac{C}{|\log\varepsilon|}\longrightarrow0,
\]
because \(dV_g\asymp\rho\,d\rho\,d\vartheta\,dV_{S^2}\) near the
codimension-two bolt.  A compactly supported smooth tensor is bounded,
so multiplication by \(\chi_\varepsilon\) converges to the identity in
\(H^1_f\).  Approximating a general \(H^1_f\) tensor first by compactly
supported smooth tensors and then using this logarithmic cutoff proves
that tensors supported away from the bolt are form-dense.

On that dense core the Peter--Weyl projections are mutually orthogonal
for both \(L^2_f\) and the shifted closed form
\[
 h\longmapsto
 \int_M\bigl(|\nabla h|^2-2\mathrm{Rm}(h,\overline h)
                 +C_{\rm Fr}|h|^2\bigr)e^{-f}\,dV_g,
\]
where \(C_{\rm Fr}\) is chosen so that the form is nonnegative.
Parseval and form closure therefore identify the full form as the
orthogonal Hilbert sum of its \((J,M')\) restrictions.  The
representation theorem for closed forms then shows that the
self-adjoint realization of \(L_f\) reduces every
\(\mathcal T^J_{M'}\); thus \eqref{eq:self-mode-splitting} is an
operator decomposition, not merely an \(L^2\) expansion.

\subsubsection{The exact Hermitian matrices}

It is convenient to define the zero-order terms first from \(F\):
\begin{align}
 \Lambda_{11}^{++}
 &:=-\frac{F''}{8}-\frac{7F'}{8r}-\frac F{r^2}+\frac1{r^2},
 \nonumber\\
 \Lambda_{11}^{--}
 &:=-\frac{F''}{8}+\frac{F'}{8r}-\frac{3F}{r^2}+\frac1{r^2},
 \nonumber\\
 \Lambda_{11}^{\pm}
 &:=-\frac{F''}{8}-\frac{3F'}{8r}+\frac F{r^2}-\frac1{r^2},
 \nonumber\\
 \Lambda_{1+}
 &:=\frac1F\left(
 -\frac{(F')^2}{16}-\frac{FF'}{2r}-\frac{F^2}{r^2}
 -\frac{F'}{2r}-\frac1{r^2}\right),
 \nonumber\\
 \Lambda_{1-}
 &:=\frac1F\left(
 -\frac{(F')^2}{16}+\frac{FF'}{2r}-\frac{3F^2}{r^2}
 +\frac{F'}{2r}+\frac{2F}{r^2}-\frac1{r^2}\right),
 \nonumber\\
 \Lambda_{01}
 &:=\frac1F\left(
 \frac{FF''}{4}-\frac{FF'}{4r}-\frac{(F')^2}{4}
 -\frac{2F^2}{r^2}\right),
 \qquad
 \Lambda_{23}:=\frac1F\left(\frac{2F}{r^2}-\frac4{r^2}\right).
\label{eq:self-Lambda-definitions}
\end{align}
Substitution of \eqref{eq:self-fik-profile-derivatives} gives
\begin{align}
 \Lambda_{11}^{++}
 &=\frac{\sqrt2c_0}{2r^2},\nonumber\\
 \Lambda_{11}^{--}
 &=\frac{\sqrt2}{r^2}
 \left(\frac{3c_0}{r^4}+\frac{2\sqrt2c_0}{r^2}
       -\frac{3c_0+1}{2\sqrt2}\right),\nonumber\\
 \Lambda_{11}^{\pm}
 &=-\frac{c_0}{r^4}-\frac{\sqrt2c_0}{2r^2}
 =\frac{\sqrt2}{r^2}
 \left(-\frac{c_0}{\sqrt2\,r^2}-\frac{c_0}{2}\right),
\label{eq:self-Lambda-pm-correct}\\
 \Lambda_{1+}
 &=\frac1{Fr^2}\left[
 -\frac{4\sqrt2c_0+c_0^2}{4r^4}
 -\frac{\sqrt2c_0^2}{2r^2}-\frac32\right],\nonumber\\
 \Lambda_{1-}
 &=\frac1{Fr^2}\left[
 -\frac{3c_0^2}{r^8}-\frac{5\sqrt2c_0^2}{r^6}
 +\frac{16c_0-17c_0^2}{4r^4}
 +\frac{(7\sqrt2-2)c_0}{2r^2}+c_0-\frac32\right],\nonumber\\
 \Lambda_{01}
 &=\frac1{Fr^2}\left(
 -\frac{\sqrt2c_0}{r^4}+\frac{\sqrt2c_0}{r^2}-1\right),\nonumber\\
 \Lambda_{23}
 &=\frac1{Fr^2}\left(
 -\frac{\sqrt2c_0}{r^4}-\frac{2c_0}{r^2}+\sqrt2-4\right).
\label{eq:self-Lambda-expanded}
\end{align}
The \(r^{-4}\) coefficient in
\(\Lambda_{11}^{\pm}\) is \(c_0\), not \(c_0^2\); this is forced
directly by the defining derivative identity in
\eqref{eq:self-Lambda-definitions}.  Keeping the derivative definition
next to the expanded formula is a useful convention check.

We now define the matrices in precisely the indexing used for the
finite certificates.  For \(M\in\mathcal J_J\), put
\begin{align*}
 q_0(M)&=\Lambda_{11}^{++}+\widetilde\Delta^J_M,&
 q_1(M)&=\Lambda_{11}^{--}+\widetilde\Delta^J_M,\\
 q_+(M)&=\Lambda_{1+}+\widetilde\Delta^J_M
          +\frac{M(4+rF')}{4s^2},&
 q_-(M)&=\Lambda_{1+}+\widetilde\Delta^J_M
          -\frac{M(4+rF')}{4s^2},
\end{align*}
and
\begin{align*}
 p_1(M)&=\Lambda_{1-}+\widetilde\Delta^J_M
          -\frac{M(4-4F-rF')}{4s^2},\\
 p_2(M)&=\Lambda_{01}+\widetilde\Delta^J_M
          +\frac{M(4F+2rF')}{4s^2},\\
 p_3(M)&=\Lambda_{23}+\widetilde\Delta^J_M
          -\frac{M(8-4F)}{4s^2}.
\end{align*}
On the index set
\(\{0,1,+,-\}\times\mathcal J_J\), the Hermitian matrix
\(\mathcal Q^J(r)\) has these diagonal entries and the following
nonzero upper-triangular entries:
\begin{align}
 \mathcal Q^J_{(0,M),(1,M)}&=\Lambda_{11}^{\pm},\nonumber\\
 \mathcal Q^J_{(1,M),(+,M+2)}
 &=-\frac{1+i}{2r^2}\sqrt F\,C^J_{M+},\nonumber\\
 \mathcal Q^J_{(1,M),(-,M-2)}
 &= \frac{1-i}{2r^2}\sqrt F\,C^J_{M-}.
\label{eq:self-Q-entrywise}
\end{align}
An entry with an index outside \(\mathcal J_J\) is absent, the
lower-triangular entries are the complex conjugates of those displayed,
and every other entry is zero.

On \(\{1,2,3\}\times\mathcal J_J\), the Hermitian matrix
\(\mathcal P^J(r)\) has diagonal entries
\(p_1(M),p_2(M),p_3(M)\) and nonzero upper-triangular entries
\begin{align}
 \mathcal P^J_{(1,M),(2,M+2)}
 &=-\frac{1+i}{2r^2}\sqrt F\,C^J_{M+},\nonumber\\
\mathcal P^J_{(1,M),(3,M-2)}
&= \frac{1-i}{2r^2}\sqrt F\,C^J_{M-}.
\label{eq:self-P-entrywise}
\end{align}
Let \(\widetilde{\mathcal P}^{\,J}(r)\) denote the Hermitian radial
potential matrix on the barred coefficient fiber, in the corresponding
ordered barred basis, obtained from the barred connection formulas.
Define the signed index-reversal \(U_J\) from the barred coefficient
fiber to the unbarred coefficient fiber by
\[
 (U_J\kappa)_{1,M}=-\kappa_{1,-M},\qquad
 (U_J\kappa)_{\ell,M}=\kappa_{\ell,-M},
 \quad \ell=2,3,
 \quad M\in\mathcal J_J.
\]
The barred connection formulas, which are the complex conjugates of
\eqref{eq:self-anti-connection}, give
\[
 U_J\widetilde{\mathcal P}^{\,J}(r)U_J^{-1}
 =\mathcal P^J(r),\qquad r>1.
\]
Lift \(U_J\) radially by
\[
 (\mathscr U_J\kappa)(r):=U_J\kappa(r).
\]
Because \(U_J\) is constant in \(r\), \(\mathscr U_J\) is unitary on
the full radial coefficient space with the measure induced by
\(e^{-f}dV_{\bar g}\), preserves the radial kinetic form, and maps the
smooth bolt-regular compactly supported core of the barred sector onto
the corresponding unbarred core.  It therefore maps the closed form
and Friedrichs operator domains onto one another and intertwines the
two self-adjoint realizations.  In particular, the barred and
unbarred sectors have the same spectrum.

\begin{proposition}[Exact mode identity]
\label{prop:self-exact-mode-identity}
Let \(h=h_I+h_A\in\mathcal T^J_{M'}\), and order the radial Wigner
coefficients as
\[
 \eta=(h_{0,M},h_{1,M},h_{+,M},h_{-,M})_{M\in\mathcal J_J},
 \qquad
 \kappa=(k_{1,M},k_{2,M},k_{3,M})_{M\in\mathcal J_J},
\]
with \(\bar\kappa\) denoting the barred coefficients.  Then, at every
\(r>1\),
\begin{align}
 \frac{J+1}{2\pi^2}\int_{S^3}
 \left(8\mathrm{Rm}(h_I,\overline{h_I})
 -4|\nabla h_I|^2\right)dV_{S^3}
 &=(\mathcal Q^J\eta,\eta)-\frac F4|\eta'|^2,\nonumber\\
 \frac{J+1}{2\pi^2}\int_{S^3}
 \left(8\mathrm{Rm}(h_A,\overline{h_A})
 -4|\nabla h_A|^2\right)dV_{S^3}
 &=(\mathcal P^J\kappa,\kappa)
 +(\widetilde{\mathcal P}^{\,J}\bar\kappa,\bar\kappa)
 -\frac F4\bigl(|\kappa'|^2+|\bar\kappa'|^2\bigr).
\label{eq:self-exact-mode-identity}
\end{align}
\end{proposition}

\begin{proof}
The radial term follows at once from
\(e_0=(\sqrt F/2)\partial_r\), the identity
\(\nabla_{e_0}b_p=0\), and \(|b_p|^2=1/4\).
For angular derivatives, insert
\eqref{eq:self-fik-wigner-derivatives} into
\eqref{eq:self-invariant-connection} and
\eqref{eq:self-anti-connection}, square, and use
\eqref{eq:self-wigner-orthogonality}.  For example, the cross term
between \(h_{1,M}\) and \(h_{+,M+2}\) is
\[
 -\frac{1+i}{2r^2}\sqrt F\,C^J_{M+}
 h_{1,M}\overline{h_{+,M+2}}
 \quad+\quad\text{its complex conjugate},
\]
which is the second line of \eqref{eq:self-Q-entrywise}.  The
\(h_{1,M}\)--\(h_{-,M-2}\) term and the two anti-invariant terms are
identical calculations with \(D_-\).  The uncoupled scalar angular
part is \(\widetilde\Delta^J_M\) by
\eqref{eq:self-berger-eigenvalue}.  Finally,
\eqref{eq:self-curvature-action} and the squares of the connection
coefficients \eqref{eq:self-gamma} collect into
\eqref{eq:self-Lambda-definitions}.  This gives every displayed
diagonal and off-diagonal entry and proves the identities.
\end{proof}

The sparse form makes the bounded block decomposition transparent.
For each \(M\in\mathcal J_J\), the nontrivial connected component of
\(\mathcal Q^J\) centered at \((1,M)\), in the order
\[
 (0,M),(1,M),(+,M+2),(-,M-2),
\]
is
\begin{equation}
\label{eq:self-Q-connected-block}
 \begin{pmatrix}
 q_0(M)&\Lambda_{11}^{\pm}&0&0\\
 \Lambda_{11}^{\pm}&q_1(M)&
 -\dfrac{1+i}{2r^2}\sqrt F\,C^J_{M+}&
  \dfrac{1-i}{2r^2}\sqrt F\,C^J_{M-}\\
 0&-\dfrac{1-i}{2r^2}\sqrt F\,C^J_{M+}&q_+(M+2)&0\\
 0& \dfrac{1+i}{2r^2}\sqrt F\,C^J_{M-}&0&q_-(M-2)
 \end{pmatrix},
\end{equation}
with rows and columns whose index lies outside \(\mathcal J_J\)
deleted.  The entries \((+,-J)\) and \((-,J)\) are isolated one by
one blocks.  Thus, for \(J\geq1\), \(\mathcal Q^J\) is the orthogonal
sum of two \(1\times1\), two \(3\times3\), and \(J-1\)
\(4\times4\) blocks.

Likewise, the component of \(\mathcal P^J\) centered at \((1,M)\),
ordered as \((1,M),(2,M+2),(3,M-2)\), is
\begin{equation}
\label{eq:self-P-connected-block}
 \begin{pmatrix}
 p_1(M)&-\dfrac{1+i}{2r^2}\sqrt F\,C^J_{M+}&
          \dfrac{1-i}{2r^2}\sqrt F\,C^J_{M-}\\
 -\dfrac{1-i}{2r^2}\sqrt F\,C^J_{M+}&p_2(M+2)&0\\
  \dfrac{1+i}{2r^2}\sqrt F\,C^J_{M-}&0&p_3(M-2)
 \end{pmatrix},
\end{equation}
again deleting absent indices.  The entries \((2,-J)\) and \((3,J)\)
are isolated.  Hence \(\mathcal P^J\) is the orthogonal sum of two
\(1\times1\), two \(2\times2\), and \(J-1\) \(3\times3\) blocks.
Equations \eqref{eq:self-Q-connected-block} and
\eqref{eq:self-P-connected-block}, rather than a separately transcribed
list of small matrices, are the defining block formulas.

\paragraph{Exact placement of the nonnegative Hessian modes.}
For the exceptional-block count it is useful to record the component
matching rather than infer it only from the Wigner labels.  Put
\(S^1_{\pm1,M'}=\nabla^2(\widehat uD^1_{\pm1,M'})\) and
\(S^2_{0,M'}=\nabla^2(\widehat vD^2_{0,M'})\).  Direct insertion of
\eqref{eq:self-fik-wigner-derivatives} into the connection formulas
\eqref{eq:self-invariant-connection}--\eqref{eq:self-anti-connection}
gives
\begin{align}
 S^1_{1,M'}
 & =-(1+i)\frac{s}{r}\left(\frac{\widehat u}{r}\right)'
       D^1_{-1,M'}\mathbf k_1
    +\frac{s}{r}\left(\frac{\widehat u}{s}\right)'
       D^1_{1,M'}\mathbf k_2,\nonumber\\
 S^1_{-1,M'}
 & =(1-i)\frac{s}{r}\left(\frac{\widehat u}{r}\right)'
       D^1_{1,M'}\mathbf k_{\bar1}
    +\frac{s}{r}\left(\frac{\widehat u}{s}\right)'
       D^1_{-1,M'}\mathbf k_{\bar2},\nonumber\\
 S^2_{0,M'}
 &=-c_0\frac{1+\sqrt2r^2}{r^2}D^2_{0,M'}b_0
   +\left(2+c_0\frac{r^2+\sqrt2}{r^4}\right)D^2_{0,M'}b_1\nonumber\\
 &\hspace{8mm}
   -\sqrt2(1-i)\frac{s}{r}D^2_{2,M'}b_-
   +\sqrt2(1+i)\frac{s}{r}D^2_{-2,M'}b_+ .
\label{eq:self-exceptional-mode-vectors}
\end{align}
Thus \(S^1_{1,M'}\) lies in the connected
\(\mathcal P^1[1_{-1},2_1]\) block,
\(S^1_{-1,M'}\) lies in its barred conjugate, and
\(S^2_{0,M'}\) lies in
\(\mathcal Q^2[0_0,1_0,+_2,-_{-2}]\), in precisely the coordinate
ordering of \eqref{eq:self-Q-connected-block} and
\eqref{eq:self-P-connected-block}.  Together with
Corollary~\ref{cor:self-contained-Hessian-families}, this verifies
smoothness, nonzeroness, \(L^2_f\)-integrability, and the eigenvalues
of the modes used in the exceptional-block exhaustion.

\subsubsection{Interface with the uniform high-frequency comparison}

We conclude by proving the cross-sectional inequality used for
\(J\geq5\), thereby completing the self-contained high-frequency
argument.  Write
\[
 N^J_M:=\frac{J(J+2)-M^2}{4},\qquad d:=\frac12rF',
 \qquad
 \int_{S^3}^{J}\Phi
 :=\frac{J+1}{2\pi^2}\int_{S^3}\Phi\,dV_{S^3}.
\]
For the invariant coefficients define
\begin{align*}
 \mathcal G_I^J
 :=\sum_{M\in\mathcal J_J}\bigg[
 &\left(N^J_M+\frac{M^2}{4F}-\frac72\right)|h_{0,M}|^2\\
 &+\left(N^J_M+2F+\frac{M^2}{4F}
                 -\frac9{20}-rF'\right)|h_{1,M}|^2\\
 &+\left(N^J_M+F+\frac{(M-2-d)^2}{4F}\right)|h_{+,M}|^2\\
 &+\left(N^J_M+F+\frac{(-M-2-d)^2}{4F}\right)|h_{-,M}|^2
 \bigg],
\end{align*}
and
\[
 \mathcal E_I^J
 :=\sqrt{2F}\sum_{M\in\mathcal J_J}C^J_{M+}
 \left(
 |h_{-,M}|\,|h_{1,M+2}|
 +|h_{1,M}|\,|h_{+,M+2}|
 \right).
\]
For the unbarred anti-invariant coefficients define
\begin{align*}
 \mathcal G_A^J
 :=\sum_{M\in\mathcal J_J}\bigg[
 &\left(N^J_M+2F+
 \frac{(M+2-d-2F)^2}{4F}-rF'\right)|k_{1,M}|^2\\
 &+\left(N^J_M+F+
 \frac{(M-rF'-2F)^2}{4F}\right)|k_{2,M}|^2\\
 &+\left(N^J_M+F+
 \frac{(M+4-2F)^2}{4F}\right)|k_{3,M}|^2
 \bigg],
\end{align*}
\[
 \mathcal E_A^J
 :=\sqrt{2F}\sum_{M\in\mathcal J_J}C^J_{M+}
 \left(
 |k_{3,M}|\,|k_{1,M+2}|
 +|k_{1,M}|\,|k_{2,M+2}|
 \right).
\]
Define \(\mathcal G_{\bar A}^J\) by replacing \(k_a\) by
\(k_{\bar a}\) and \(M\) by \(-M\) inside each square, and put
\[
 \mathcal E_{\bar A}^J
 :=\sqrt{2F}\sum_{M\in\mathcal J_J}C^J_{M+}
 \left(
 |k_{\bar2,M}|\,|k_{\bar1,M+2}|
 +|k_{\bar1,M}|\,|k_{\bar3,M+2}|
 \right).
\]

\begin{lemma}[Cross-sectional reduction]
\label{lem:self-cross-sectional-reduction}
For \(J\geq3\) and \(h=h_I+h_A\in\mathcal T^J_{M'}\),
\begin{align*}
 4r^2\int_{S^3}^{J}\bigl(
 2\mathrm{Rm}(h_I,\overline{h_I})
 -|\nabla_{e_1}h_I|^2
 -|\nabla_{D_+}h_I|^2
 -|\nabla_{D_-}h_I|^2\bigr)
 &\leq\mathcal E_I^J-\mathcal G_I^J,\\
 4r^2\int_{S^3}^{J}\bigl(
 2\mathrm{Rm}(h_A,\overline{h_A})
 -|\nabla_{e_1}h_A|^2
 -|\nabla_{D_+}h_A|^2
 -|\nabla_{D_-}h_A|^2\bigr)
 &\leq
 \mathcal E_A^J-\mathcal G_A^J
 +\mathcal E_{\bar A}^J-\mathcal G_{\bar A}^J.
\end{align*}
\end{lemma}

\begin{proof}
Expand the squares using
\eqref{eq:self-invariant-connection}--\eqref{eq:self-fik-wigner-derivatives}.
After multiplication by \(4r^2\), the diagonal part of the invariant
angular gradient is
\begin{align*}
 \sum_{M\in\mathcal J_J}\bigg[
 &\left(N^J_M+\frac{M^2}{4F}\right)|h_{0,M}|^2
 +\left(N^J_M+2F+\frac{M^2}{4F}\right)|h_{1,M}|^2\\
 &+\left(N^J_M+F+\frac{(M-2-d)^2}{4F}\right)|h_{+,M}|^2
 +\left(N^J_M+F+\frac{(-M-2-d)^2}{4F}\right)|h_{-,M}|^2
 \bigg].
\end{align*}
The corresponding unbarred anti-invariant diagonal is
\begin{align*}
 \sum_{M\in\mathcal J_J}\bigg[
 &\left(N^J_M+2F+
 \frac{(M+2-d-2F)^2}{4F}\right)|k_{1,M}|^2\\
 &+\left(N^J_M+F+
 \frac{(M-rF'-2F)^2}{4F}\right)|k_{2,M}|^2\\
 &+\left(N^J_M+F+
 \frac{(M+4-2F)^2}{4F}\right)|k_{3,M}|^2
 \bigg],
\end{align*}
and the barred diagonal is obtained by \(M\mapsto-M\).
These formulas follow term by term from
\eqref{eq:self-invariant-connection}--\eqref{eq:self-fik-wigner-derivatives};
for example, the \(e_1\)-square on \(h_{+,M}\) is
\((M-2-d)^2/(4F)\), while its two horizontal squares contribute
\(N^J_M+F\).  The exact off-diagonal terms are those in
\eqref{eq:self-Q-entrywise} and \eqref{eq:self-P-entrywise}; after the
factor \(4r^2\) is inserted, taking their absolute values gives the
three \(\mathcal E\)'s.

For the curvature part,
\eqref{eq:self-curvature-action},
\[
 r^3F'=2c_0\left(1+\frac{\sqrt2}{r^2}\right)\leq2,
\]
give the mixed-term bound
\[
 \left|
 -\frac{r^3F'}{2\sqrt2}
 (h_{0,M}\overline{h_{1,M}}+
  h_{1,M}\overline{h_{0,M}})
 \right|
 \leq\sqrt2\,|h_{0,M}|\,|h_{1,M}|.
\]
Put
\[
 a_0:=\frac{c_0}{\sqrt2}=1-\frac1{\sqrt2},\qquad
 A:=\frac72-a_0,\qquad B:=\frac9{20}-a_0.
\]
Both \(A,B\) are positive, and
\[
 AB=\frac{39\sqrt2-35}{40}>\frac12.
\]
Therefore Young's inequality gives
\[
 \sqrt2|xy|\leq A|x|^2+B|y|^2.
\]
Adding the
curvature diagonal \(a_0|h_{0,M}|^2+
(a_0+rF')|h_{1,M}|^2\), and discarding the favorable
\(h_\pm\) curvature terms, yields
\begin{multline*}
 4r^2\int_{S^3}^{J}2\mathrm{Rm}(h_I,\overline{h_I})\\
 \leq\sum_M
 \left(\frac72|h_{0,M}|^2+
       \left(\frac9{20}+rF'\right)|h_{1,M}|^2\right),
\end{multline*}
\begin{multline*}
 4r^2\int_{S^3}^{J}2\mathrm{Rm}(h_A,\overline{h_A})\\
 \leq\sum_M\left[
 rF'\bigl(|k_{1,M}|^2+|k_{\bar1,M}|^2\bigr)
 -\frac12\bigl(|k_{3,M}|^2+|k_{\bar3,M}|^2\bigr)
 \right].
\end{multline*}
For the second line one uses
\[
 rF'+2F-\sqrt2=\frac{\sqrt2c_0}{r^4}>0
\]
and
\(2-2F\geq2-\sqrt2>1/2\).
Dropping the remaining favorable
\(-\frac12(|k_3|^2+|k_{\bar3}|^2)\) contribution and collecting the
diagonal squares gives exactly the displayed \(\mathcal G\)'s.
\end{proof}

All profile bounds used in the subsequent uniform absorption follow
without an auxiliary lemma:
\begin{equation}
\label{eq:self-profile-ranges}
 0<F<\frac1{\sqrt2},\qquad
 0<d\leq1,\qquad
 rF'=2d\leq2,\qquad
 1<d+2F=\sqrt2-\frac{c_0}{r^2}<\sqrt2
 \quad(r>1).
\end{equation}
Thus Lemma \ref{lem:self-cross-sectional-reduction} and
\eqref{eq:self-profile-ranges} are precisely the geometric input
needed by the \(J\geq5\) star-by-star Young comparison.

\begin{lemma}[Uniform high-frequency Wigner estimate]
\label{lem:FIK-high-frequency}
In the Wigner notation of
Appendix~\ref{app:self-contained-wigner}, let \(J\geq5\),
\(M'\in\{-J,-J+2,\ldots,J\}\), and
\(h\in\mathcal T^J_{M'}\cap H^1_f\).  Then
\[
 \int_M\left(2\overline{\Rm}(h,\overline h)-|\nabla h|^2\right)
 e^{-f_{\rm NO}}\,dV_{\bar g}
 \leq-c_{\rm hf}\|h\|_{L^2_f}^2
\]
for a numerical \(c_{\rm hf}>0\), independent of \(J\), \(M'\), and
\(h\).  In particular, the inequality is strict for \(h\neq0\).
\end{lemma}

\begin{proof}[Proof of Lemma~\ref{lem:FIK-high-frequency}]
Let \(\Pi_{J,M'}\) denote the orthogonal Wigner projection onto the
fixed block.  Since it acts only in the compact angular variables,
\(\Pi_{J,M'}\) is bounded on \(H^1_f\).  It also commutes with every
smooth radial cutoff.  Hence, if \(h\in H^1_f\cap\mathcal T^J_{M'}\),
first cutting off radially and then projecting a smooth compactly
supported approximation shows that
\[
 C_c^\infty(S^2T^*M)\cap\mathcal T^J_{M'}
 \quad\text{is dense in}\quad
 H^1_f\cap\mathcal T^J_{M'}.
\]
It is therefore enough first to work with a smooth compactly supported
tensor in the fixed block.  We now give the remaining coefficient
comparison.  Put
\[
 q:=J(J+2),\qquad
 \mathcal J_J:=\{-J,-J+2,\ldots,J\},\qquad
 N_m:=\frac{q-m^2}{4},\qquad
 C_m:=\sqrt{q-m(m+2)}.
\]
Thus \(C_m=C^J_{m+}\) in
\eqref{eq:self-CJM}.  Write \(F=F(r)\) for the FIK profile and
 \(d:=\frac12rF'\).  The explicit formula for \(F\), together with
 \eqref{eq:self-profile-ranges}, gives, for \(r>1\),
\begin{equation}\label{eq:NO-F-ranges}
 0<F\leq\frac1{\sqrt2},\qquad
 0<d\leq1,\qquad
 rF'\leq2,\qquad
 1<d+2F\leq\sqrt2.
\end{equation}
For an angular integrand \(\Phi\), we use from this point onward the
normalized Wigner integral
\begin{equation}\label{eq:NO-angular-integral}
 \int_{S^3}^{J}\Phi
 :=\frac{J+1}{2\pi^2}\int_{S^3}\Phi\,dV_{S^3}.
\end{equation}

We use the explicitly defined quantities
\[
 \mathcal E_I^J,\ \mathcal G_I^J,\ \mathcal E_A^J,\
 \mathcal G_A^J,\ \mathcal E_{\bar A}^J,\ \mathcal G_{\bar A}^J
\]
from Appendix~\ref{app:self-contained-wigner}.
Lemma~\ref{lem:self-cross-sectional-reduction} reduces the desired
cross-sectional inequality to
\[
 \mathcal E_I^J\leq\mathcal G_I^J,\qquad
 \mathcal E_A^J\leq\mathcal G_A^J,\qquad
 \mathcal E_{\bar A}^J\leq\mathcal G_{\bar A}^J.
\]
We verify these inequalities uniformly in \(J\geq5\).
Set
\[
 \mathcal E^J
 :=\mathcal E_I^J+\mathcal E_A^J+\mathcal E_{\bar A}^J,
 \qquad
 \mathcal G^J
 :=\mathcal G_I^J+\mathcal G_A^J+\mathcal G_{\bar A}^J.
\]
For the fixed Wigner block \((J,M')\), let
\begin{equation}\label{eq:NO-original-block-form}
 \mathcal Q_{J,M'}[h]
 :=
 \int_M\left(
  2\overline{\Rm}(h,\overline h)-|\nabla h|^2
 \right)e^{-f_{\rm NO}}\,dV_{\bar g}
\end{equation}
denote the restriction of the original integrated stability form to
\(\mathcal T^J_{M'}\cap H^1_f\).  Thus every total form occurring in
the reduction below is defined before it is used.

For the invariant block, the coefficient in \(\mathcal G_I^J\) of
\(|h_{0,m}|^2\) is
\[
 D_m^I:=N_m+\frac{m^2}{4F}-\frac72
 \geq\frac q4-\frac72>0.
\]
For each \(m\in\mathcal J_J\), the coefficient of
\(|h_{1,m}|^2\) is
\[
 A_m^I:=N_m+2F+\frac{m^2}{4F}-\frac9{20}-rF'
 \geq\frac q4-\frac{49}{20}
 \geq\frac{63}{10}>4\sqrt2.
\]
Here we used \(F\leq1/\sqrt2<1\), so that
\((4F)^{-1}-1/4\geq0\), and \eqref{eq:NO-F-ranges}.
The two terms of \(\mathcal E_I^J\) incident to \(h_{1,m}\) are
\begin{equation}\label{eq:NO-invariant-star}
 \sqrt{2F}\,C_m|h_{1,m}|\,|h_{+,m+2}|
 +\sqrt{2F}\,C_{m-2}|h_{1,m}|\,|h_{-,m-2}|,
\end{equation}
with an absent index interpreted as an absent term.  The corresponding
leaf coefficients in \(\mathcal G_I^J\) are
\[
 \begin{split}
 B^I_{+,m}
 &=N_{m+2}+F+\frac{(m-d)^2}{4F}
 \quad (m\leq J-2),\\
 B^I_{-,m}
 &=N_{m-2}+F+\frac{(m+d)^2}{4F}
 \quad (m\geq-J+2).
 \end{split}
\]

The following elementary estimates are the uniform point in the
argument:
\begin{align}
 N_{m+2}-\frac{C_m^2}{8}
 &=\frac{(J-m-2)(J+m+4)}8\geq0
 &&(m\leq J-2),\label{eq:NO-leaf-plus}\\
 N_{m-2}-\frac{C_{m-2}^2}{8}
 &=\frac{(J+m-2)(J-m+4)}8\geq0
 &&(m\geq-J+2).\label{eq:NO-leaf-minus}
\end{align}
Consequently,
\[
 B^I_{+,m}\geq\frac{C_m^2}{8},
 \qquad
 B^I_{-,m}\geq\frac{C_{m-2}^2}{8}.
\]
Fix
\[
 \vartheta_{\rm Y}:=\frac{49}{50}.
\]
For \(B>0\), Young's inequality and either of the preceding bounds give
\begin{equation}\label{eq:NO-edge-absorption}
 \sqrt{2F}\,C\,|x|\,|y|
 \leq \vartheta_{\rm Y}B|y|^2
      +\frac{FC^2}{2\vartheta_{\rm Y}B}|x|^2
 \leq \vartheta_{\rm Y}B|y|^2
      +\frac{4F}{\vartheta_{\rm Y}}|x|^2.
\end{equation}
There are at most two edges in \eqref{eq:NO-invariant-star}; hence their
total cost at the central coefficient is at most
\[
 \frac{8F}{\vartheta_{\rm Y}}
 \leq\frac{200\sqrt2}{49}<A_m^I.
\]
Every \(h_+\)- or \(h_-\)-coefficient occurs as a leaf of at most one
such star, so summing these inequalities uses each leaf diagonal only
once and leaves the positive fraction
\((1-\vartheta_{\rm Y})B\) of it unused.  Thus every invariant star is
absorbed with the uniform central margin
\[
 \delta_I:=\frac{63}{10}-\frac{200\sqrt2}{49}>0.
\]
Together with \(D_m^I\geq21/4\), this proves the invariant estimate
with a uniform positive diagonal remainder.

For the anti-invariant block, the coefficient of
\(|k_{1,m}|^2\) is
\[
 A_m^A
 :=N_m+2F+\frac{(m+2-d-2F)^2}{4F}-rF'.
\]
Set \(a:=2-d-2F\).  By \eqref{eq:NO-F-ranges},
\(0\leq a<1\).  Completing the square gives the exact identity
\[
 N_m+\frac{(m+a)^2}{4F}
 =\frac q4+
   \frac{1-F}{4F}\left(m+\frac{a}{1-F}\right)^2
   -\frac{a^2}{4(1-F)}.
\]
Since \(F\leq1/\sqrt2\) and \(a^2\leq1\),
\[
 \frac{a^2}{4(1-F)}
 \leq\frac1{4(1-1/\sqrt2)}
 =\frac{2+\sqrt2}{4}.
\]
Using also \(2F-rF'\geq-2\), we obtain
\begin{equation}\label{eq:NO-anti-central}
 A_m^A
 \geq\frac q4-2-\frac{2+\sqrt2}{4}
 \geq\frac{25-\sqrt2}{4}>4\sqrt2.
\end{equation}
The two terms of \(\mathcal E_A^J\) incident to \(k_{1,m}\) couple it
to \(k_{2,m+2}\) with coefficient \(\sqrt{2F}\,C_m\), and to
\(k_{3,m-2}\) with coefficient \(\sqrt{2F}\,C_{m-2}\).  The
corresponding leaf coefficients in \(\mathcal G_A^J\) are
\[
 \begin{split}
 B^A_{2,m}
 &=N_{m+2}+F+
   \frac{(m+2-rF'-2F)^2}{4F},\\
 B^A_{3,m}
 &=N_{m-2}+F+
   \frac{(m+2-2F)^2}{4F},
 \end{split}
\]
whenever the indicated index belongs to \(\mathcal J_J\).
Equations~\eqref{eq:NO-leaf-plus}--\eqref{eq:NO-leaf-minus} therefore
give
\[
 B^A_{2,m}\geq\frac{C_m^2}{8},\qquad
 B^A_{3,m}\geq\frac{C_{m-2}^2}{8}.
\]
Applying \eqref{eq:NO-edge-absorption} to the at most two incident
edges costs at most \(200\sqrt2/49\), which is strictly smaller than
\((25-\sqrt2)/4\).  Hence the anti-invariant central remainder is at
least
\[
 \delta_A:=\frac{25-\sqrt2}{4}
            -\frac{200\sqrt2}{49}>0,
\]
 and every incident leaf retains the fraction
 \((1-\vartheta_{\rm Y})B\).  After \(m\mapsto-m\) and the corresponding
 interchange of the two leaf families, the barred anti-invariant block
 has the same form and the same remainders.

We also check explicitly the diagonal coefficients which are not
incident to an angular edge.  In the invariant block these are
\(h_{+,-J}\) and \(h_{-,J}\), and the displayed definition of
\(\mathcal G_I^J\) in Appendix~\ref{app:self-contained-wigner}
gives the common coefficient
\[
 D^I_{\rm unc}
 =\frac J2+F+\frac{(J+2+d)^2}{4F}
 \geq J+2+d\geq7.
\]
Here and below we use \(x+y^2/(4x)\geq y\) for \(x,y>0\).
In the unbarred anti-invariant block the uncoupled entries are
\(k_{2,-J}\) and \(k_{3,J}\), with respective coefficients
\[
 \begin{split}
 D^A_{2,-J}
 &=\frac J2+F+\frac{(J+rF'+2F)^2}{4F}
   \geq J+rF'+2F\geq J,\\
 D^A_{3,J}
 &=\frac J2+F+\frac{(J+4-2F)^2}{4F}
   \geq J+4-2F\geq J+4-\sqrt2 .
 \end{split}
\]
The barred uncoupled entries are
\(\bar k_{2,J}\) and \(\bar k_{3,-J}\); their coefficients are exactly
\(D^A_{2,-J}\) and \(D^A_{3,J}\), respectively.  Thus no endpoint
coefficient is left outside the comparison.

These remainders control the coefficient norm uniformly.  Indeed, on
 every present edge
\[
 C_m^2\geq4J\geq20
 \quad\hbox{or}\quad
 C_{m-2}^2\geq4J\geq20,
 \]
 so the retained leaf coefficient is at least
 \((1-\vartheta_{\rm Y})20/8=1/20\).  The uncoupled diagonal
 coefficients have the lower bounds just displayed.  We finally insert
 the normalization factor which is suppressed in a coefficient
 calculation.  The tensor frame of
 Appendix~\ref{app:self-contained-wigner} satisfies
 \[
  \langle b_p,b_q\rangle=\frac14\delta_{pq},
  \qquad h_p=4\langle h,b_p\rangle .
 \]
 Consequently, on the fixed \((J,M')\)-block,
 \[
 \begin{split}
  4\int_{S^3}^{J}|h|^2
  =\sum_{m\in\mathcal J_J}\biggl(
    |h_{0,m}|^2+|h_{1,m}|^2+|h_{+,m}|^2+|h_{-,m}|^2\\
    {}+\sum_{\ell=1}^3
       \bigl(|k_{\ell,m}|^2+|\bar k_{\ell,m}|^2\bigr)
  \biggr).
 \end{split}
 \]
 There is therefore a numerical \(c_{\rm ang}>0\), independent of
 \(J\geq5\), \(M'\), and \(r>1\), such that
\begin{equation}\label{eq:NO-pointwise-gap}
 \mathcal G^J-\mathcal E^J
 \geq c_{\rm ang}\int_{S^3}^{J}|h|^2 .
\end{equation}

The cross-sectional estimate in
Lemma~\ref{lem:self-cross-sectional-reduction} is therefore
nonpositive for every \(r>1\).
In the orthonormal frame of
Appendix~\ref{app:self-contained-wigner},
\(e_0=(s/2r)\partial_r\) is the radial direction, whereas \(e_1\) is
the angular Hopf direction.  The cross-sectional lemma contains the favorable
\(e_1,D_+,D_-\) derivative terms in the cross-sectional form; the
 negative \(e_0\)-radial derivative is not part of that comparison and
 remains as a separate favorable term.  For clarity, we fix the radial
 normalization entering that reduction.  Put
 \[
 d\mu_r
 :=\frac{8\pi^2}{J+1}\,
      r e^{-f_{\rm NO}(r)}\,dr .
\]
Define
\[
 \mathcal Q_{\rm cs}(r)
 :=4r^2\int_{S^3}^{J}
  \left(
   2\overline{\Rm}(h,\overline h)
   -|\nabla_{e_1}h|^2
   -|\nabla_{D_+}h|^2
   -|\nabla_{D_-}h|^2
  \right)
\]
and
\[
 \mathfrak r_{\rm rad}
 :=\int_1^\infty
    4r^2\int_{S^3}^{J}|\nabla_{e_0}h|^2\,d\mu_r .
\]
With these conventions, writing
\(\mathfrak r_{\rm ang}:=\mathcal G^J-\mathcal E^J\), the exact
reduction is
\[
 \mathcal Q_{\rm cs}\leq
 \mathcal E^J-\mathcal G^J=-\mathfrak r_{\rm ang},
 \qquad
 \mathcal Q_{J,M'}[h]=
 \int_1^\infty\mathcal Q_{\rm cs}\,d\mu_r-\mathfrak r_{\rm rad}.
\]
The retained angular margin is an \(r^{-2}\)-weighted zeroth-order
term, so one must combine it with the radial derivative rather than
identify it directly with the full \(L^2_f\) norm.  The required
one-dimensional Gaussian estimate is
\begin{equation}\label{eq:NO-radial-Hardy}
 \int_1^\infty4r^2\int_{S^3}^{J}|h|^2\,d\mu_r
 \leq C\int_1^\infty\left(
  4r^2\int_{S^3}^{J}|\nabla_{e_0}h|^2
  +\int_{S^3}^{J}|h|^2\right)d\mu_r .
\end{equation}
Indeed, choose \(R_0>1\) so that \(e_0\) is uniformly comparable with
\(\partial_r\) on \(\{r\geq R_0\}\).  On \([1,R_0]\), the factor
\(4r^2\) is bounded.  On the exterior, use
\(-\partial_r(e^{-f_{\rm NO}})=2\sqrt2\,r e^{-f_{\rm NO}}\),
integrate by parts with a cutoff, and absorb
\(r^2|u||u'|\) by
\(\varepsilon r^3|u'|^2+C_\varepsilon r|u|^2\).
This proves \eqref{eq:NO-radial-Hardy} componentwise; the radial
connection coefficients are \(O(r^{-1})\) and are absorbed by the
 last zeroth-order term.
 Multiplication by the positive radial density, integration in \(r\),
 and \eqref{eq:NO-pointwise-gap}--\eqref{eq:NO-radial-Hardy} give, for
 every smooth compactly supported tensor in the fixed block, on putting
 \[
  \mathfrak z_{\rm ang}
  :=\int_1^\infty\int_{S^3}^{J}|h|^2\,d\mu_r,
 \]
 \[
  \mathcal Q_{J,M'}[h]
  \leq-c_{\rm ang}\mathfrak z_{\rm ang}
       -\mathfrak r_{\rm rad}.
 \]
 Split both favorable terms into equal halves.  If \(C\) is the
 constant in \eqref{eq:NO-radial-Hardy}, then
 \[
  \frac{c_{\rm ang}}2\mathfrak z_{\rm ang}
  +\frac12\mathfrak r_{\rm rad}
  \geq
  \frac{\min\{c_{\rm ang},1\}}{2C}\,
  \|h\|_{L^2_f}^2.
 \]
 Hence, with
 \(c_{\rm hf}:=\min\{c_{\rm ang},1\}/(2C)>0\),
\begin{equation}\label{eq:NO-L2-gap}
 \mathcal Q_{J,M'}[h]
 \leq-c_{\rm hf}\|h\|_{L^2_f}^2
      -\frac12\mathfrak r_{\rm rad}
 \leq-c_{\rm hf}\|h\|_{L^2_f}^2 .
\end{equation}
This formulation deliberately avoids extending a separately defined
residual form across the collapsing bolt.  The original form
\(\mathcal Q_{J,M'}\) is continuous on \(H^1_f\).  For the density
passage we use only the last, pure \(L^2_f\) inequality in
\eqref{eq:NO-L2-gap}, whose right-hand side is \(L^2_f\)-continuous.
The density statement at the beginning of the proof therefore extends
that inequality directly to the full block form domain, proving the
asserted gap and strictness.  (The stronger half-radial inequality also
passes by \(H^1_f\)-continuity, but is not needed.)
\end{proof}

\subsubsection{Dictionary for exact arithmetic}

For reproducibility, set
\[
 x=r^2,\qquad z=\sqrt F,\qquad
 F=\frac1{\sqrt2}-\frac{c_0}{x}
    -\frac{c_0}{\sqrt2\,x^2},\qquad
 s^2=xF,\qquad
 rF'=2x\frac{dF}{dx}.
\]
Then \eqref{eq:self-berger-eigenvalue},
\eqref{eq:self-Lambda-expanded},
\eqref{eq:self-Q-entrywise}, and
\eqref{eq:self-P-entrywise} are exactly an entrywise specification of
the matrices over the finite algebraic extension of
\(\mathbb Q(\sqrt2)(x,z)/(z^2-F)\) obtained by adjoining \(i\) and the
numbers \(C^J_{M\pm}\).  Every connected component is a tree.  Hence
every term in any principal determinant contains each off-diagonal
edge together with its conjugate.  Since
\(\lvert1\pm i\rvert^2(C^J_{M\pm})^2
=2[J(J+2)-M(M\pm2)]\) is an integer, after \(z^2=F\) every such
determinant lies in \(\mathbb Q(\sqrt2)(x)\).  This explains why the
exact low-mode certificate may build the matrices with a formal
\(z=\sqrt F\), expand a determinant, and only then replace \(z^2\) by
\(F\).  The same reduction provides a direct entrywise verification
that the symbolic matrices coincide with the geometric matrices.

\subsection{Exact certificates for the safe low Wigner blocks}
\label{app:exact-low-mode-certificates}

We give exact certificates replacing the numerical eigenvalue plots
used to determine the signs of the low-frequency zeroth-order blocks.
The three blocks which are not pointwise negative,
\[
 \mathcal Q^0[0_0,1_0],\qquad
 \mathcal P^1[1_{-1},2_1]\quad\hbox{(and its conjugate block)},\qquad
 \mathcal Q^2[0_0,1_0,+_2,-_{-2}],
\]
are not part of the calculation below; they are treated by the
Sturm--Liouville comparison argument.  Thus ``safe'' means every other
block with \(0\leq J\leq4\).

\begin{proposition}[Exact safe-block certificate]
\label{prop:exact-safe-low-mode-blocks}
Every safe zeroth-order Wigner block with \(0\leq J\leq4\) is negative
definite at every \(r>1\).  This assertion holds for both the invariant
and anti-invariant sectors, including the barred anti-invariant sector.
\end{proposition}

\begin{proof}
Put
\[
 K=\mathbb Q(\sqrt2)\subset\mathbb R,\qquad
 x=r^2,\qquad c_0=\sqrt2-1,
\]
\[
 F(x)=\frac1{\sqrt2}-\frac{c_0}{x}
      -\frac{c_0}{\sqrt2\,x^2},
 \qquad s^2=xF(x).
\]
Here \(x>1\), and hence \(s^2>0\).  We use throughout the coefficient
obtained directly from the defining curvature formula,
\begin{equation}\label{eq:exact-lambda-pm}
\Lambda_{11}^{\pm}
   =-\frac{c_0}{x^2}-\frac{\sqrt2\,c_0}{2x}.
\end{equation}
In particular, the coefficient of \(x^{-2}\) in
\(\Lambda_{11}^{\pm}\) is \(c_0\), not \(c_0^2\).

For a safe Hermitian block \(\mathcal B\), in the ordering displayed in
the first column of the tables below, define
\[
 B_{\mathcal B}(x):=-4s^2\mathcal B(x),\qquad
 \Delta_{\mathcal B,k}(x)
 :=\det B_{\mathcal B}(x)[1\!:\!k,1\!:\!k].
\]
Although the off-diagonal entries can contain \(i\sqrt F\), the
Hermitian products eliminate both \(i\) and the odd powers of
\(\sqrt F\).  Exact expansion therefore gives
\begin{equation}\label{eq:exact-minor-normal-form}
\Delta_{\mathcal B,k}(x)
   =\frac{(x-1)^{m_{\mathcal B,k}}}
          {d_{\mathcal B,k}(x)}
      p_{\mathcal B,k}(x-1),
\end{equation}
where \(p_{\mathcal B,k}\in K[y]\) and, in the exact normalization used
throughout the coefficient calculation,
\(d_{\mathcal B,k}(x)=x^{\ell_{\mathcal B,k}}\).  Thus the denominator in
\eqref{eq:exact-minor-normal-form} is positive for \(x>1\).  The factor
\((x-1)^{m_{\mathcal B,k}}\) is the maximal bolt factor in the
numerator.

The notation in the tables is as follows.  A row of length \(n\)
records, in increasing order \(k=1,\ldots,n\), the degree
\(\deg p_{\mathcal B,k}\), the bolt multiplicity
\(m_{\mathcal B,k}\), the exact anchor \(p_{\mathcal B,k}(0)\), and
the certificate used for that minor.  The symbol \(\mathbf C\) means
that all \(\deg p_{\mathcal B,k}+1\) coefficients of
\(p_{\mathcal B,k}(y)\) are strictly positive in the ordered field
\(K\).  The three symbols \(\mathbf B_0,\mathbf B_1,\mathbf B_2\)
refer to the explicit Bernstein certificates following the tables.  Thus
every table entry certifies a particular leading principal minor; no
minor is suppressed.

\begin{equation}\label{eq:certificate-J0}
\begin{array}{c|c|c|c|c}
\multicolumn{5}{c}{J=0}\\
\mathcal B & (\deg p_k)&(m_k)&(p_k(0))&\text{certificate}\\ \hline
\mathcal Q[+_0] &(2)&(0)&(9)&(\mathbf C)\\
\mathcal P[1_0] &(4)&(0)&(1)&(\mathbf B_0)\\
\mathcal P[2_0] &(2)&(0)&(4)&(\mathbf C)\\
\mathcal P[3_0] &(2)&(0)&(16)&(\mathbf C)
\end{array}
\end{equation}

\begin{equation}\label{eq:certificate-J1}
\begin{array}{c|c|c|c|c}
\multicolumn{5}{c}{J=1}\\
\mathcal B & (\deg p_k)&(m_k)&(p_k(0))&\text{certificate}\\ \hline
\mathcal Q[0_{-1},1_{-1},+_1]
 &(2,6,8)&(0,0,0)&(1,1,4)
 &(\mathbf C,\mathbf B_1,\mathbf B_2)\\
\mathcal Q[+_{-1}] &(2)&(0)&(16)&(\mathbf C)\\
\mathcal P[1_1,3_{-1}] &(4,6)&(0,0)&(4,36)
 &(\mathbf C,\mathbf C)\\
\mathcal P[2_{-1}] &(2)&(0)&(9)&(\mathbf C)\\
\mathcal P[3_1] &(2)&(0)&(25)&(\mathbf C)
\end{array}
\end{equation}

\begin{equation}\label{eq:certificate-J2}
\begin{array}{c|c|c|c|c}
\multicolumn{5}{c}{J=2}\\
\mathcal B & (\deg p_k)&(m_k)&(p_k(0))&\text{certificate}\\ \hline
\mathcal Q[0_{-2},1_{-2},+_0]
 &(2,6,8)&(0,0,0)&(4,16,144)&(\mathbf C,\mathbf C,\mathbf C)\\
\mathcal Q[+_{-2}] &(2)&(0)&(25)&(\mathbf C)\\
\mathcal P[1_{-2},2_0] &(4,6)&(0,0)&(1,4)
 &(\mathbf C,\mathbf C)\\
\mathcal P[1_0,2_2,3_{-2}]
 &(4,5,7)&(0,1,1)&(1,12-4\sqrt2,48-16\sqrt2)
 &(\mathbf C,\mathbf C,\mathbf C)\\
\mathcal P[1_2,3_0] &(4,6)&(0,0)&(9,144)
 &(\mathbf C,\mathbf C)\\
\mathcal P[2_{-2}] &(2)&(0)&(16)&(\mathbf C)\\
\mathcal P[3_2] &(2)&(0)&(36)&(\mathbf C)
\end{array}
\end{equation}

\begin{equation}\label{eq:certificate-J3}
\begin{array}{c|c|c|c|c}
\multicolumn{5}{c}{J=3}\\
\mathcal B & (\deg p_k)&(m_k)&(p_k(0))&\text{certificate}\\ \hline
\mathcal Q[0_{-3},1_{-3},+_{-1}]
 &(2,6,8)&(0,0,0)&(9,81,1296)
 &(\mathbf C,\mathbf C,\mathbf C)\\
\mathcal Q[0_{-1},1_{-1},+_1,-_{-3}]
 &(2,6,8,9)&(0,0,0,1)&(1,1,4,40)
 &(\mathbf C,\mathbf C,\mathbf C,\mathbf C)\\
\mathcal Q[+_{-3}] &(2)&(0)&(36)&(\mathbf C)\\
\mathcal P[1_{-3},2_{-1}] &(4,6)&(0,0)&(4,36)
 &(\mathbf C,\mathbf C)\\
\mathcal P[1_{-1},2_1,3_{-3}]
 &(3,5,7)&(1,1,1)&(6,6,6)
 &(\mathbf C,\mathbf C,\mathbf C)\\
\mathcal P[1_1,2_3,3_{-1}]
 &(4,6,8)&(0,0,0)&(4,4,36)
 &(\mathbf C,\mathbf C,\mathbf C)\\
\mathcal P[1_3,3_1] &(4,6)&(0,0)&(16,400)
 &(\mathbf C,\mathbf C)\\
\mathcal P[2_{-3}] &(2)&(0)&(25)&(\mathbf C)\\
\mathcal P[3_3] &(2)&(0)&(49)&(\mathbf C)
\end{array}
\end{equation}

\[
 a_\star:=240+64\sqrt2.
\]
\begin{equation}\label{eq:certificate-J4}
\begin{array}{c|c|c|c|c}
\multicolumn{5}{c}{J=4}\\
\mathcal B & (\deg p_k)&(m_k)&(p_k(0))&\text{certificate}\\ \hline
\mathcal Q[0_{-4},1_{-4},+_{-2}]
 &(2,6,8)&(0,0,0)&(16,256,6400)
 &(\mathbf C,\mathbf C,\mathbf C)\\
\mathcal Q[0_{-2},1_{-2},+_0,-_{-4}]
 &(2,6,8,10)&(0,0,0,0)&(4,16,144,144)
 &(\mathbf C,\mathbf C,\mathbf C,\mathbf C)\\
\mathcal Q[0_0,1_0,+_2,-_{-2}]
 &(1,4,6,8)&(1,2,2,2)
 &(20+2\sqrt2,a_\star,a_\star,a_\star)
 &(\mathbf C,\mathbf C,\mathbf C,\mathbf C)\\
\mathcal Q[+_{-4}] &(2)&(0)&(49)&(\mathbf C)\\
\mathcal P[1_{-4},2_{-2}] &(4,6)&(0,0)&(9,144)
 &(\mathbf C,\mathbf C)\\
\mathcal P[1_{-2},2_0,3_{-4}]
 &(4,6,7)&(0,0,1)&(1,4,64)
 &(\mathbf C,\mathbf C,\mathbf C)\\
\mathcal P[1_0,2_2,3_{-2}]
 &(4,5,7)&(0,1,1)&(1,28-4\sqrt2,112-16\sqrt2)
 &(\mathbf C,\mathbf C,\mathbf C)\\
\mathcal P[1_2,2_4,3_0]
 &(4,6,8)&(0,0,0)&(9,36,576)
 &(\mathbf C,\mathbf C,\mathbf C)\\
\mathcal P[1_4,3_2] &(4,6)&(0,0)&(25,900)
 &(\mathbf C,\mathbf C)\\
\mathcal P[2_{-4}] &(2)&(0)&(36)&(\mathbf C)\\
\mathcal P[3_4] &(2)&(0)&(64)&(\mathbf C)
\end{array}
\end{equation}

For direct verification, write
\[
 p_{\mathcal B,k}(y)=\sum_{j=0}^{\deg p_{\mathcal B,k}}
 c_j(\mathcal B,k)y^j,\qquad
 \mathbf c(\mathcal B,k)=(c_0,\ldots,c_{\deg p_{\mathcal B,k}}).
\]
The following tables list the coefficient vector of every
\(\mathbf C\)-minor, together with the denominator exponent
\(\ell\) and bolt exponent \(m\) in
\eqref{eq:exact-minor-normal-form}.  They therefore supply a direct
coefficientwise verification of every asserted positive-coefficient
certificate.
\begingroup\small
\paragraph{The \(J=0\) coefficient vectors.}
\begin{equation*}
\begin{aligned}
\mathbf c(\mathcal Q[+_{0}],1)&=\bigl(9,\, 4 + 6 \sqrt{2},\, 6\bigr),\qquad(\ell,m)=(2,0),\quad \mathbf C.
\end{aligned}
\end{equation*}
\begin{equation*}
\begin{aligned}
\mathbf c(\mathcal P[2_{0}],1)&=\bigl(4,\, 4 \sqrt{2},\, 4\bigr),\qquad(\ell,m)=(2,0),\quad \mathbf C.
\end{aligned}
\end{equation*}
\begin{equation*}
\begin{aligned}
\mathbf c(\mathcal P[3_{0}],1)&=\bigl(16,\, 24,\, 16 - 4 \sqrt{2}\bigr),\qquad(\ell,m)=(2,0),\quad \mathbf C.
\end{aligned}
\end{equation*}
\paragraph{The \(J=1\) coefficient vectors.}
\begin{equation*}
\begin{aligned}
\mathbf c(\mathcal Q[0_{-1},1_{-1},+_{1}],1)&=\bigl(1,\, 2 \sqrt{2},\, 3 - \sqrt{2}\bigr),\qquad(\ell,m)=(2,0),\quad \mathbf C.
\end{aligned}
\end{equation*}
\begin{equation*}
\begin{aligned}
\mathbf c(\mathcal Q[+_{-1}],1)&=\bigl(16,\, 8 \sqrt{2} + 14,\, \sqrt{2} + 11\bigr),\qquad(\ell,m)=(2,0),\quad \mathbf C.
\end{aligned}
\end{equation*}
\begin{equation*}
\begin{aligned}
\mathbf c(\mathcal P[1_{1},3_{-1}],1)&=\bigl(4,\, 14 - 4 \sqrt{2},\, 59 - 25 \sqrt{2},\\[-2pt]
&\qquad 28 - 2 \sqrt{2},\, 15 - 5 \sqrt{2}\bigr),\\[-2pt]
&\qquad(\ell,m)=(4,0),\quad \mathbf C.
\end{aligned}
\end{equation*}
\begin{equation*}
\begin{aligned}
\mathbf c(\mathcal P[1_{1},3_{-1}],2)&=\bigl(36,\, 190 - 36 \sqrt{2},\, 791 - 293 \sqrt{2},\\[-2pt]
&\qquad 1298 - 468 \sqrt{2},\, 1164 - 409 \sqrt{2},\, 448 - 126 \sqrt{2},\\[-2pt]
&\qquad 145 - 68 \sqrt{2}\bigr),\\[-2pt]
&\qquad(\ell,m)=(6,0),\quad \mathbf C.
\end{aligned}
\end{equation*}
\begin{equation*}
\begin{aligned}
\mathbf c(\mathcal P[2_{-1}],1)&=\bigl(9,\, 4 + 8 \sqrt{2},\, 3 \sqrt{2} + 5\bigr),\qquad(\ell,m)=(2,0),\quad \mathbf C.
\end{aligned}
\end{equation*}
\begin{equation*}
\begin{aligned}
\mathbf c(\mathcal P[3_{1}],1)&=\bigl(25,\, 40,\, 25 - 5 \sqrt{2}\bigr),\qquad(\ell,m)=(2,0),\quad \mathbf C.
\end{aligned}
\end{equation*}
\paragraph{The \(J=2\) coefficient vectors.}
\begin{equation*}
\begin{aligned}
\mathbf c(\mathcal Q[0_{-2},1_{-2},+_{0}],1)&=\bigl(4,\, 2 \sqrt{2} + 8,\, 6\bigr),\qquad(\ell,m)=(2,0),\quad \mathbf C.
\end{aligned}
\end{equation*}
\begin{equation*}
\begin{aligned}
\mathbf c(\mathcal Q[0_{-2},1_{-2},+_{0}],2)&=\bigl(16,\, 16 \sqrt{2} + 64,\, 256,\\[-2pt]
&\qquad 64 \sqrt{2} + 416,\, 120 \sqrt{2} + 384,\, 72 \sqrt{2} + 192,\\[-2pt]
&\qquad 8 \sqrt{2} + 48\bigr),\\[-2pt]
&\qquad(\ell,m)=(6,0),\quad \mathbf C.
\end{aligned}
\end{equation*}
\begin{equation*}
\begin{aligned}
\mathbf c(\mathcal Q[0_{-2},1_{-2},+_{0}],3)&=\bigl(144,\, 240 \sqrt{2} + 768,\, 640 \sqrt{2} + 3360,\\[-2pt]
&\qquad 2464 \sqrt{2} + 7072,\, 4856 \sqrt{2} + 10240,\, 5672 \sqrt{2} + 9632,\\[-2pt]
&\qquad 3512 \sqrt{2} + 6096,\, 1760 + 1456 \sqrt{2},\, 144 \sqrt{2} + 352\bigr),\\[-2pt]
&\qquad(\ell,m)=(8,0),\quad \mathbf C.
\end{aligned}
\end{equation*}
\begin{equation*}
\begin{aligned}
\mathbf c(\mathcal Q[+_{-2}],1)&=\bigl(25,\, 10 \sqrt{2} + 28,\, 2 \sqrt{2} + 18\bigr),\qquad(\ell,m)=(2,0),\quad \mathbf C.
\end{aligned}
\end{equation*}
\begin{equation*}
\begin{aligned}
\mathbf c(\mathcal P[1_{-2},2_{0}],1)&=\bigl(1,\, -2 + 2 \sqrt{2},\, 27 - 6 \sqrt{2},\\[-2pt]
&\qquad 18 \sqrt{2},\, 2 \sqrt{2} + 6\bigr),\\[-2pt]
&\qquad(\ell,m)=(4,0),\quad \mathbf C.
\end{aligned}
\end{equation*}
\begin{equation*}
\begin{aligned}
\mathbf c(\mathcal P[1_{-2},2_{0}],2)&=\bigl(4,\, 12 \sqrt{2},\, 112 - 12 \sqrt{2},\\[-2pt]
&\qquad 112 + 132 \sqrt{2},\, 140 \sqrt{2} + 228,\, 112 + 112 \sqrt{2},\\[-2pt]
&\qquad 16 \sqrt{2} + 40\bigr),\\[-2pt]
&\qquad(\ell,m)=(6,0),\quad \mathbf C.
\end{aligned}
\end{equation*}
\begin{equation*}
\begin{aligned}
\mathbf c(\mathcal P[1_{0},2_{2},3_{-2}],1)&=\bigl(1,\, 6 - 2 \sqrt{2},\, 47 - 16 \sqrt{2},\\[-2pt]
&\qquad 10 \sqrt{2} + 16,\, 10\bigr),\\[-2pt]
&\qquad(\ell,m)=(4,0),\quad \mathbf C.
\end{aligned}
\end{equation*}
\begin{equation*}
\begin{aligned}
\mathbf c(\mathcal P[1_{0},2_{2},3_{-2}],2)&=\bigl(12 - 4 \sqrt{2},\, 96 - 50 \sqrt{2},\, 684 - 408 \sqrt{2},\\[-2pt]
&\qquad 552 - 262 \sqrt{2},\, 8 \sqrt{2} + 112,\, 80 - 36 \sqrt{2}\bigr),\\[-2pt]
&\qquad(\ell,m)=(6,1),\quad \mathbf C.
\end{aligned}
\end{equation*}
\begin{equation*}
\begin{aligned}
\mathbf c(\mathcal P[1_{0},2_{2},3_{-2}],3)&=\bigl(48 - 16 \sqrt{2},\, 528 - 248 \sqrt{2},\, 3920 - 2224 \sqrt{2},\\[-2pt]
&\qquad 9832 - 5696 \sqrt{2},\, 8432 - 4400 \sqrt{2},\, 2056 - 376 \sqrt{2},\\[-2pt]
&\qquad 800 - 176 \sqrt{2},\, 240 - 112 \sqrt{2}\bigr),\\[-2pt]
&\qquad(\ell,m)=(8,1),\quad \mathbf C.
\end{aligned}
\end{equation*}
\begin{equation*}
\begin{aligned}
\mathbf c(\mathcal P[1_{2},3_{0}],1)&=\bigl(9,\, 38 - 6 \sqrt{2},\, 99 - 30 \sqrt{2},\\[-2pt]
&\qquad 56 - 6 \sqrt{2},\, 22 - 6 \sqrt{2}\bigr),\\[-2pt]
&\qquad(\ell,m)=(4,0),\quad \mathbf C.
\end{aligned}
\end{equation*}
\begin{equation*}
\begin{aligned}
\mathbf c(\mathcal P[1_{2},3_{0}],2)&=\bigl(144,\, 896 - 96 \sqrt{2},\, 2944 - 672 \sqrt{2},\\[-2pt]
&\qquad 4608 - 1152 \sqrt{2},\, 3728 - 864 \sqrt{2},\, 1504 - 288 \sqrt{2},\\[-2pt]
&\qquad 352 - 112 \sqrt{2}\bigr),\\[-2pt]
&\qquad(\ell,m)=(6,0),\quad \mathbf C.
\end{aligned}
\end{equation*}
\begin{equation*}
\begin{aligned}
\mathbf c(\mathcal P[2_{-2}],1)&=\bigl(16,\, 12 + 12 \sqrt{2},\, 8 + 6 \sqrt{2}\bigr),\qquad(\ell,m)=(2,0),\quad \mathbf C.
\end{aligned}
\end{equation*}
\begin{equation*}
\begin{aligned}
\mathbf c(\mathcal P[3_{2}],1)&=\bigl(36,\, 60,\, 36 - 6 \sqrt{2}\bigr),\qquad(\ell,m)=(2,0),\quad \mathbf C.
\end{aligned}
\end{equation*}
\paragraph{The \(J=3\) coefficient vectors.}
\begin{equation*}
\begin{aligned}
\mathbf c(\mathcal Q[0_{-3},1_{-3},+_{-1}],1)&=\bigl(9,\, 2 \sqrt{2} + 20,\, \sqrt{2} + 11\bigr),\qquad(\ell,m)=(2,0),\quad \mathbf C.
\end{aligned}
\end{equation*}
\begin{equation*}
\begin{aligned}
\mathbf c(\mathcal Q[0_{-3},1_{-3},+_{-1}],2)&=\bigl(81,\, 36 \sqrt{2} + 450,\, 46 \sqrt{2} + 1383,\\[-2pt]
&\qquad 184 \sqrt{2} + 2148,\, 428 \sqrt{2} + 1641,\, 228 \sqrt{2} + 750,\\[-2pt]
&\qquad 34 \sqrt{2} + 155\bigr),\\[-2pt]
&\qquad(\ell,m)=(6,0),\quad \mathbf C.
\end{aligned}
\end{equation*}
\begin{equation*}
\begin{aligned}
\mathbf c(\mathcal Q[0_{-3},1_{-3},+_{-1}],3)&=\bigl(1296,\, 1224 \sqrt{2} + 9306,\, 5839 \sqrt{2} + 35295,\\[-2pt]
&\qquad 18750 \sqrt{2} + 75652,\, 37515 \sqrt{2} + 98985,\, 41988 \sqrt{2} + 84790,\\[-2pt]
&\qquad 26579 \sqrt{2} + 46503,\, 9258 \sqrt{2} + 14336,\, 1195 \sqrt{2} + 2133\bigr),\\[-2pt]
&\qquad(\ell,m)=(8,0),\quad \mathbf C.
\end{aligned}
\end{equation*}
\begin{equation*}
\begin{aligned}
\mathbf c(\mathcal Q[0_{-1},1_{-1},+_{1},-_{-3}],1)&=\bigl(1,\, 2 \sqrt{2} + 12,\, 3 + 5 \sqrt{2}\bigr),\qquad(\ell,m)=(2,0),\quad \mathbf C.
\end{aligned}
\end{equation*}
\begin{equation*}
\begin{aligned}
\mathbf c(\mathcal Q[0_{-1},1_{-1},+_{1},-_{-3}],2)&=\bigl(1,\, 4 \sqrt{2} + 18,\, 38 \sqrt{2} + 119,\\[-2pt]
&\qquad 88 \sqrt{2} + 580,\, 265 + 556 \sqrt{2},\, 164 \sqrt{2} + 414,\\[-2pt]
&\qquad 59 + 58 \sqrt{2}\bigr),\\[-2pt]
&\qquad(\ell,m)=(6,0),\quad \mathbf C.
\end{aligned}
\end{equation*}
\begin{equation*}
\begin{aligned}
\mathbf c(\mathcal Q[0_{-1},1_{-1},+_{1},-_{-3}],3)&=\bigl(4,\, 20 \sqrt{2} + 86,\, 287 \sqrt{2} + 763,\\[-2pt]
&\qquad 1498 \sqrt{2} + 4272,\, 6275 \sqrt{2} + 9237,\, 11442 + 10880 \sqrt{2},\\[-2pt]
&\qquad 6715 \sqrt{2} + 11699,\, 4124 + 3094 \sqrt{2},\, 669 + 491 \sqrt{2}\bigr),\\[-2pt]
&\qquad(\ell,m)=(8,0),\quad \mathbf C.
\end{aligned}
\end{equation*}
\begin{equation*}
\begin{aligned}
\mathbf c(\mathcal Q[0_{-1},1_{-1},+_{1},-_{-3}],4)&=\bigl(40,\, 212 \sqrt{2} + 872,\, 3188 \sqrt{2} + 7624,\\[-2pt]
&\qquad 16402 \sqrt{2} + 40779,\, 61244 \sqrt{2} + 91038,\, 127173 + 108344 \sqrt{2},\\[-2pt]
&\qquad 130244 + 93316 \sqrt{2},\, 45262 \sqrt{2} + 80291,\, 17994 + 17876 \sqrt{2},\\[-2pt]
&\qquad 1584 \sqrt{2} + 3225\bigr),\\[-2pt]
&\qquad(\ell,m)=(10,1),\quad \mathbf C.
\end{aligned}
\end{equation*}
\begin{equation*}
\begin{aligned}
\mathbf c(\mathcal Q[+_{-3}],1)&=\bigl(36,\, 12 \sqrt{2} + 46,\, 3 \sqrt{2} + 27\bigr),\qquad(\ell,m)=(2,0),\quad \mathbf C.
\end{aligned}
\end{equation*}
\begin{equation*}
\begin{aligned}
\mathbf c(\mathcal P[1_{-3},2_{-1}],1)&=\bigl(4,\, 4 \sqrt{2} + 10,\, \sqrt{2} + 43,\\[-2pt]
&\qquad 8 + 26 \sqrt{2},\, 7 + 5 \sqrt{2}\bigr),\\[-2pt]
&\qquad(\ell,m)=(4,0),\quad \mathbf C.
\end{aligned}
\end{equation*}
\begin{equation*}
\begin{aligned}
\mathbf c(\mathcal P[1_{-3},2_{-1}],2)&=\bigl(36,\, 68 \sqrt{2} + 154,\, 189 \sqrt{2} + 631,\\[-2pt]
&\qquad 802 + 704 \sqrt{2},\, 840 + 773 \sqrt{2},\, 338 \sqrt{2} + 556,\\[-2pt]
&\qquad 64 \sqrt{2} + 125\bigr),\\[-2pt]
&\qquad(\ell,m)=(6,0),\quad \mathbf C.
\end{aligned}
\end{equation*}
\begin{equation*}
\begin{aligned}
\mathbf c(\mathcal P[1_{-1},2_{1},3_{-3}],1)&=\bigl(6,\, 47 - 7 \sqrt{2},\, 12 + 22 \sqrt{2},\\[-2pt]
&\qquad 7 + 5 \sqrt{2}\bigr),\\[-2pt]
&\qquad(\ell,m)=(4,1),\quad \mathbf C.
\end{aligned}
\end{equation*}
\begin{equation*}
\begin{aligned}
\mathbf c(\mathcal P[1_{-1},2_{1},3_{-3}],2)&=\bigl(6,\, 143 - 7 \sqrt{2},\, 666 - 60 \sqrt{2},\\[-2pt]
&\qquad 364 + 365 \sqrt{2},\, 200 + 250 \sqrt{2},\, 28 \sqrt{2} + 85\bigr),\\[-2pt]
&\qquad(\ell,m)=(6,1),\quad \mathbf C.
\end{aligned}
\end{equation*}
\begin{equation*}
\begin{aligned}
\mathbf c(\mathcal P[1_{-1},2_{1},3_{-3}],3)&=\bigl(6,\, 215 - 7 \sqrt{2},\, 2292 - 114 \sqrt{2},\\[-2pt]
&\qquad 209 \sqrt{2} + 6893,\, 4010 + 5116 \sqrt{2},\, 2755 + 4469 \sqrt{2},\\[-2pt]
&\qquad 482 \sqrt{2} + 3000,\, 125 + 333 \sqrt{2}\bigr),\\[-2pt]
&\qquad(\ell,m)=(8,1),\quad \mathbf C.
\end{aligned}
\end{equation*}
\begin{equation*}
\begin{aligned}
\mathbf c(\mathcal P[1_{1},2_{3},3_{-1}],1)&=\bigl(4,\, 26 - 4 \sqrt{2},\, 83 - 19 \sqrt{2},\\[-2pt]
&\qquad 10 \sqrt{2} + 40,\, \sqrt{2} + 15\bigr),\\[-2pt]
&\qquad(\ell,m)=(4,0),\quad \mathbf C.
\end{aligned}
\end{equation*}
\begin{equation*}
\begin{aligned}
\mathbf c(\mathcal P[1_{1},2_{3},3_{-1}],2)&=\bigl(4,\, 122 - 36 \sqrt{2},\, 823 - 335 \sqrt{2},\\[-2pt]
&\qquad 2602 - 1240 \sqrt{2},\, 2008 - 719 \sqrt{2},\, 660 - 86 \sqrt{2},\\[-2pt]
&\qquad 189 - 56 \sqrt{2}\bigr),\\[-2pt]
&\qquad(\ell,m)=(6,0),\quad \mathbf C.
\end{aligned}
\end{equation*}
\begin{equation*}
\begin{aligned}
\mathbf c(\mathcal P[1_{1},2_{3},3_{-1}],3)&=\bigl(36,\, 1210 - 324 \sqrt{2},\, 10859 - 4003 \sqrt{2},\\[-2pt]
&\qquad 47072 - 20254 \sqrt{2},\, 91913 - 39259 \sqrt{2},\, 74398 - 23280 \sqrt{2},\\[-2pt]
&\qquad 27607 - 335 \sqrt{2},\, 766 \sqrt{2} + 8388,\, 25 \sqrt{2} + 1333\bigr),\\[-2pt]
&\qquad(\ell,m)=(8,0),\quad \mathbf C.
\end{aligned}
\end{equation*}
\begin{equation*}
\begin{aligned}
\mathbf c(\mathcal P[1_{3},3_{1}],1)&=\bigl(16,\, 70 - 8 \sqrt{2},\, 151 - 35 \sqrt{2},\\[-2pt]
&\qquad 92 - 10 \sqrt{2},\, 31 - 7 \sqrt{2}\bigr),\\[-2pt]
&\qquad(\ell,m)=(4,0),\quad \mathbf C.
\end{aligned}
\end{equation*}
\begin{equation*}
\begin{aligned}
\mathbf c(\mathcal P[1_{3},3_{1}],2)&=\bigl(400,\, 2582 - 200 \sqrt{2},\, 7815 - 1275 \sqrt{2},\\[-2pt]
&\qquad 11790 - 2200 \sqrt{2},\, 9264 - 1563 \sqrt{2},\, 3748 - 522 \sqrt{2},\\[-2pt]
&\qquad 761 - 168 \sqrt{2}\bigr),\\[-2pt]
&\qquad(\ell,m)=(6,0),\quad \mathbf C.
\end{aligned}
\end{equation*}
\begin{equation*}
\begin{aligned}
\mathbf c(\mathcal P[2_{-3}],1)&=\bigl(25,\, 16 \sqrt{2} + 24,\, 9 \sqrt{2} + 13\bigr),\qquad(\ell,m)=(2,0),\quad \mathbf C.
\end{aligned}
\end{equation*}
\begin{equation*}
\begin{aligned}
\mathbf c(\mathcal P[3_{3}],1)&=\bigl(49,\, 84,\, 49 - 7 \sqrt{2}\bigr),\qquad(\ell,m)=(2,0),\quad \mathbf C.
\end{aligned}
\end{equation*}
\paragraph{The \(J=4\) coefficient vectors.}
\begin{equation*}
\begin{aligned}
\mathbf c(\mathcal Q[0_{-4},1_{-4},+_{-2}],1)&=\bigl(16,\, 2 \sqrt{2} + 36,\, 2 \sqrt{2} + 18\bigr),\qquad(\ell,m)=(2,0),\quad \mathbf C.
\end{aligned}
\end{equation*}
\begin{equation*}
\begin{aligned}
\mathbf c(\mathcal Q[0_{-4},1_{-4},+_{-2}],2)&=\bigl(256,\, 64 \sqrt{2} + 1536,\, 128 \sqrt{2} + 4400,\\[-2pt]
&\qquad 448 \sqrt{2} + 6496,\, 976 \sqrt{2} + 4840,\, 536 \sqrt{2} + 2064,\\[-2pt]
&\qquad 88 \sqrt{2} + 392\bigr),\\[-2pt]
&\qquad(\ell,m)=(6,0),\quad \mathbf C.
\end{aligned}
\end{equation*}
\begin{equation*}
\begin{aligned}
\mathbf c(\mathcal Q[0_{-4},1_{-4},+_{-2}],3)&=\bigl(6400,\, 4160 \sqrt{2} + 49664,\, 23936 \sqrt{2} + 183472,\\[-2pt]
&\qquad 77344 \sqrt{2} + 385440,\, 152048 \sqrt{2} + 492936,\, 170600 \sqrt{2} + 403824,\\[-2pt]
&\qquad 108040 \sqrt{2} + 210296,\, 36544 \sqrt{2} + 63296,\, 4928 \sqrt{2} + 8688\bigr),\\[-2pt]
&\qquad(\ell,m)=(8,0),\quad \mathbf C.
\end{aligned}
\end{equation*}
\begin{equation*}
\begin{aligned}
\mathbf c(\mathcal Q[0_{-2},1_{-2},+_{0},-_{-4}],1)&=\bigl(4,\, 2 \sqrt{2} + 24,\, 6 + 8 \sqrt{2}\bigr),\qquad(\ell,m)=(2,0),\quad \mathbf C.
\end{aligned}
\end{equation*}
\begin{equation*}
\begin{aligned}
\mathbf c(\mathcal Q[0_{-2},1_{-2},+_{0},-_{-4}],2)&=\bigl(16,\, 16 \sqrt{2} + 192,\, 128 \sqrt{2} + 896,\\[-2pt]
&\qquad 448 \sqrt{2} + 2080,\, 1216 + 1432 \sqrt{2},\, 584 \sqrt{2} + 960,\\[-2pt]
&\qquad 176 + 136 \sqrt{2}\bigr),\\[-2pt]
&\qquad(\ell,m)=(6,0),\quad \mathbf C.
\end{aligned}
\end{equation*}
\begin{equation*}
\begin{aligned}
\mathbf c(\mathcal Q[0_{-2},1_{-2},+_{0},-_{-4}],3)&=\bigl(144,\, 240 \sqrt{2} + 2176,\, 2944 \sqrt{2} + 13728,\\[-2pt]
&\qquad 15392 \sqrt{2} + 46112,\, 47896 \sqrt{2} + 78400,\, 79648 + 71848 \sqrt{2},\\[-2pt]
&\qquad 44024 \sqrt{2} + 65616,\, 16432 \sqrt{2} + 24416,\, 3552 + 2640 \sqrt{2}\bigr),\\[-2pt]
&\qquad(\ell,m)=(8,0),\quad \mathbf C.
\end{aligned}
\end{equation*}
\begin{equation*}
\begin{aligned}
\mathbf c(\mathcal Q[0_{-2},1_{-2},+_{0},-_{-4}],4)&=\bigl(144,\, 5056 - 48 \sqrt{2},\, 3968 \sqrt{2} + 57152,\\[-2pt]
&\qquad 56960 \sqrt{2} + 319264,\, 323800 \sqrt{2} + 1003008,\, 1045832 \sqrt{2} + 1733760,\\[-2pt]
&\qquad 1980528 + 1693160 \sqrt{2},\, 1817792 + 1340512 \sqrt{2},\, 632992 \sqrt{2} + 1017920,\\[-2pt]
&\qquad 258816 + 204320 \sqrt{2},\, 22752 \sqrt{2} + 34752\bigr),\\[-2pt]
&\qquad(\ell,m)=(10,0),\quad \mathbf C.
\end{aligned}
\end{equation*}
\begin{equation*}
\begin{aligned}
\mathbf c(\mathcal Q[0_{0},1_{0},+_{2},-_{-2}],1)&=\bigl(2 \sqrt{2} + 20,\, 2 + 10 \sqrt{2}\bigr),\qquad(\ell,m)=(2,1),\quad \mathbf C.
\end{aligned}
\end{equation*}
\begin{equation*}
\begin{aligned}
\mathbf c(\mathcal Q[0_{0},1_{0},+_{2},-_{-2}],2)&=\bigl(64 \sqrt{2} + 240,\, 256 \sqrt{2} + 1248,\, 552 + 1328 \sqrt{2},\\[-2pt]
&\qquad 408 \sqrt{2} + 912,\, 88 \sqrt{2} + 200\bigr),\\[-2pt]
&\qquad(\ell,m)=(6,2),\quad \mathbf C.
\end{aligned}
\end{equation*}
\begin{equation*}
\begin{aligned}
\mathbf c(\mathcal Q[0_{0},1_{0},+_{2},-_{-2}],3)&=\bigl(64 \sqrt{2} + 240,\, 2016 \sqrt{2} + 6304,\, 11088 \sqrt{2} + 24456,\\[-2pt]
&\qquad 23344 + 33384 \sqrt{2},\, 17096 \sqrt{2} + 36216,\, 8256 \sqrt{2} + 13568,\\[-2pt]
&\qquad 1200 + 2080 \sqrt{2}\bigr),\\[-2pt]
&\qquad(\ell,m)=(8,2),\quad \mathbf C.
\end{aligned}
\end{equation*}
\begin{equation*}
\begin{aligned}
\mathbf c(\mathcal Q[0_{0},1_{0},+_{2},-_{-2}],4)&=\bigl(64 \sqrt{2} + 240,\, 3776 \sqrt{2} + 11360,\, 66160 \sqrt{2} + 156520,\\[-2pt]
&\qquad 348088 \sqrt{2} + 529872,\, 818984 + 801208 \sqrt{2},\, 642176 \sqrt{2} + 1128192,\\[-2pt]
&\qquad 394752 \sqrt{2} + 642176,\, 126656 + 176992 \sqrt{2},\, 6368 \sqrt{2} + 40160\bigr),\\[-2pt]
&\qquad(\ell,m)=(10,2),\quad \mathbf C.
\end{aligned}
\end{equation*}
\begin{equation*}
\begin{aligned}
\mathbf c(\mathcal Q[+_{-4}],1)&=\bigl(49,\, 14 \sqrt{2} + 68,\, 4 \sqrt{2} + 38\bigr),\qquad(\ell,m)=(2,0),\quad \mathbf C.
\end{aligned}
\end{equation*}
\begin{equation*}
\begin{aligned}
\mathbf c(\mathcal P[1_{-4},2_{-2}],1)&=\bigl(9,\, 6 \sqrt{2} + 30,\, 8 \sqrt{2} + 71,\\[-2pt]
&\qquad 24 + 34 \sqrt{2},\, 10 + 8 \sqrt{2}\bigr),\\[-2pt]
&\qquad(\ell,m)=(4,0),\quad \mathbf C.
\end{aligned}
\end{equation*}
\begin{equation*}
\begin{aligned}
\mathbf c(\mathcal P[1_{-4},2_{-2}],2)&=\bigl(144,\, 204 \sqrt{2} + 732,\, 782 \sqrt{2} + 2192,\\[-2pt]
&\qquad 2844 + 2088 \sqrt{2},\, 2440 + 2234 \sqrt{2},\, 952 \sqrt{2} + 1424,\\[-2pt]
&\qquad 172 \sqrt{2} + 304\bigr),\\[-2pt]
&\qquad(\ell,m)=(6,0),\quad \mathbf C.
\end{aligned}
\end{equation*}
\begin{equation*}
\begin{aligned}
\mathbf c(\mathcal P[1_{-2},2_{0},3_{-4}],1)&=\bigl(1,\, 2 \sqrt{2} + 14,\, 2 \sqrt{2} + 59,\\[-2pt]
&\qquad 16 + 34 \sqrt{2},\, 6 + 10 \sqrt{2}\bigr),\\[-2pt]
&\qquad(\ell,m)=(4,0),\quad \mathbf C.
\end{aligned}
\end{equation*}
\begin{equation*}
\begin{aligned}
\mathbf c(\mathcal P[1_{-2},2_{0},3_{-4}],2)&=\bigl(4,\, 12 \sqrt{2} + 80,\, 124 \sqrt{2} + 592,\\[-2pt]
&\qquad 596 \sqrt{2} + 1408,\, 964 + 1348 \sqrt{2},\, 816 + 592 \sqrt{2},\\[-2pt]
&\qquad 64 \sqrt{2} + 264\bigr),\\[-2pt]
&\qquad(\ell,m)=(6,0),\quad \mathbf C.
\end{aligned}
\end{equation*}
\begin{equation*}
\begin{aligned}
\mathbf c(\mathcal P[1_{-2},2_{0},3_{-4}],3)&=\bigl(64,\, 224 \sqrt{2} + 1280,\, 2624 \sqrt{2} + 9152,\\[-2pt]
&\qquad 12992 \sqrt{2} + 21440,\, 22144 + 26688 \sqrt{2},\, 15520 \sqrt{2} + 25408,\\[-2pt]
&\qquad 4480 \sqrt{2} + 12672,\, 256 + 1984 \sqrt{2}\bigr),\\[-2pt]
&\qquad(\ell,m)=(8,1),\quad \mathbf C.
\end{aligned}
\end{equation*}
\begin{equation*}
\begin{aligned}
\mathbf c(\mathcal P[1_{0},2_{2},3_{-2}],1)&=\bigl(1,\, 22 - 2 \sqrt{2},\, 79 - 8 \sqrt{2},\\[-2pt]
&\qquad 32 + 26 \sqrt{2},\, 10 + 8 \sqrt{2}\bigr),\\[-2pt]
&\qquad(\ell,m)=(4,0),\quad \mathbf C.
\end{aligned}
\end{equation*}
\begin{equation*}
\begin{aligned}
\mathbf c(\mathcal P[1_{0},2_{2},3_{-2}],2)&=\bigl(28 - 4 \sqrt{2},\, 640 - 138 \sqrt{2},\, 2236 - 424 \sqrt{2},\\[-2pt]
&\qquad 722 \sqrt{2} + 1224,\, 496 + 584 \sqrt{2},\, 76 \sqrt{2} + 176\bigr),\\[-2pt]
&\qquad(\ell,m)=(6,1),\quad \mathbf C.
\end{aligned}
\end{equation*}
\begin{equation*}
\begin{aligned}
\mathbf c(\mathcal P[1_{0},2_{2},3_{-2}],3)&=\bigl(112 - 16 \sqrt{2},\, 3344 - 664 \sqrt{2},\, 26896 - 5296 \sqrt{2},\\[-2pt]
&\qquad 61928 - 2368 \sqrt{2},\, 37488 + 34000 \sqrt{2},\, 22408 + 30632 \sqrt{2},\\[-2pt]
&\qquad 6352 \sqrt{2} + 16672,\, 1648 + 1680 \sqrt{2}\bigr),\\[-2pt]
&\qquad(\ell,m)=(8,1),\quad \mathbf C.
\end{aligned}
\end{equation*}
\begin{equation*}
\begin{aligned}
\mathbf c(\mathcal P[1_{2},2_{4},3_{0}],1)&=\bigl(9,\, 54 - 6 \sqrt{2},\, 131 - 22 \sqrt{2},\\[-2pt]
&\qquad 10 \sqrt{2} + 72,\, 2 \sqrt{2} + 22\bigr),\\[-2pt]
&\qquad(\ell,m)=(4,0),\quad \mathbf C.
\end{aligned}
\end{equation*}
\begin{equation*}
\begin{aligned}
\mathbf c(\mathcal P[1_{2},2_{4},3_{0}],2)&=\bigl(36,\, 576 - 132 \sqrt{2},\, 3008 - 1012 \sqrt{2},\\[-2pt]
&\qquad 7056 - 2748 \sqrt{2},\, 5524 - 1612 \sqrt{2},\, 2000 - 272 \sqrt{2},\\[-2pt]
&\qquad 424 - 80 \sqrt{2}\bigr),\\[-2pt]
&\qquad(\ell,m)=(6,0),\quad \mathbf C.
\end{aligned}
\end{equation*}
\begin{equation*}
\begin{aligned}
\mathbf c(\mathcal P[1_{2},2_{4},3_{0}],3)&=\bigl(576,\, 10944 - 2112 \sqrt{2},\, 76352 - 22240 \sqrt{2},\\[-2pt]
&\qquad 263616 - 90048 \sqrt{2},\, 451328 - 148672 \sqrt{2},\, 368000 - 80000 \sqrt{2},\\[-2pt]
&\qquad 1568 \sqrt{2} + 156160,\, 7040 \sqrt{2} + 43392,\, 1152 \sqrt{2} + 5888\bigr),\\[-2pt]
&\qquad(\ell,m)=(8,0),\quad \mathbf C.
\end{aligned}
\end{equation*}
\begin{equation*}
\begin{aligned}
\mathbf c(\mathcal P[1_{4},3_{2}],1)&=\bigl(25,\, 110 - 10 \sqrt{2},\, 215 - 40 \sqrt{2},\\[-2pt]
&\qquad 136 - 14 \sqrt{2},\, 42 - 8 \sqrt{2}\bigr),\\[-2pt]
&\qquad(\ell,m)=(4,0),\quad \mathbf C.
\end{aligned}
\end{equation*}
\begin{equation*}
\begin{aligned}
\mathbf c(\mathcal P[1_{4},3_{2}],2)&=\bigl(900,\, 5860 - 360 \sqrt{2},\, 17000 - 2150 \sqrt{2},\\[-2pt]
&\qquad 25028 - 3684 \sqrt{2},\, 19428 - 2554 \sqrt{2},\, 7840 - 840 \sqrt{2},\\[-2pt]
&\qquad 1480 - 236 \sqrt{2}\bigr),\\[-2pt]
&\qquad(\ell,m)=(6,0),\quad \mathbf C.
\end{aligned}
\end{equation*}
\begin{equation*}
\begin{aligned}
\mathbf c(\mathcal P[2_{-4}],1)&=\bigl(36,\, 20 \sqrt{2} + 40,\, 12 \sqrt{2} + 20\bigr),\qquad(\ell,m)=(2,0),\quad \mathbf C.
\end{aligned}
\end{equation*}
\begin{equation*}
\begin{aligned}
\mathbf c(\mathcal P[3_{4}],1)&=\bigl(64,\, 112,\, 64 - 8 \sqrt{2}\bigr),\qquad(\ell,m)=(2,0),\quad \mathbf C.
\end{aligned}
\end{equation*}
\endgroup

For completeness, here are the three polynomials not covered by
coefficientwise positivity:
\begin{align*}
p_{\mathbf B_0}(y)
={}&(10-4\sqrt2)y^4+(8+2\sqrt2)y^3
 +(31-20\sqrt2)y^2 \\
&\quad -(2+2\sqrt2)y+1,\\
p_{\mathbf B_1}(y)
={}&(11-2\sqrt2)y^6+(6+20\sqrt2)y^5
 +(73-20\sqrt2)y^4\\
&\quad +(-44+40\sqrt2)y^3+(23-22\sqrt2)y^2
 +(-6+4\sqrt2)y+1,\\
p_{\mathbf B_2}(y)
={}&(45-19\sqrt2)y^8+(-112+154\sqrt2)y^7
 +(491-179\sqrt2)y^6\\
&\quad +(-378+440\sqrt2)y^5+(549-331\sqrt2)y^4
 +(-348+214\sqrt2)y^3\\
&\quad +(115-103\sqrt2)y^2+(-22+20\sqrt2)y+4.
\end{align*}
They admit a short certificate with no root isolation.  For
\(0\leq y\leq\frac12\), put \(t=2y\in[0,1]\) and write
\[
 p_{\mathbf B_j}(t/2)
 =\sum_{k=0}^{n_j}b_{j,k}
   \binom{n_j}{k}t^k(1-t)^{n_j-k},
 \qquad (n_0,n_1,n_2)=(9,6,8).
\]
Exact conversion from the monomial to the Bernstein basis gives,
in increasing order of \(k\),
\begin{align}
(b_{0,k})_{k=0}^{9}
={}&\left(
1,\frac89-\frac{\sqrt2}{9},
\frac{143}{144}-\frac{13\sqrt2}{36},
\frac{445}{336}-\frac{251\sqrt2}{336},
\frac{1915}{1008}-\frac{71\sqrt2}{56},\right.\nonumber\\[-2pt]
&\left.\qquad
\frac{307}{112}-\frac{485\sqrt2}{252},
\frac{31}{8}-\frac{457\sqrt2}{168},
\frac{16}{3}-\frac{527\sqrt2}{144},
\frac{515}{72}-\frac{19\sqrt2}{4},
\frac{75}{8}-6\sqrt2
\right),\label{eq:Bernstein-B0}\\
(b_{1,k})_{k=0}^{6}
={}&\left(
1,\frac12+\frac{\sqrt2}{3},
\frac{23}{60}+\frac{3\sqrt2}{10},
\frac38+\frac{3\sqrt2}{20},
\frac{121}{240}+\frac{\sqrt2}{20},
\frac{109}{96}+\frac{3\sqrt2}{16},
\frac{203}{64}+\frac{27\sqrt2}{32}
\right),\label{eq:Bernstein-B1}\\
(b_{2,k})_{k=0}^{8}
={}&\left(
4,\frac{21}{8}+\frac{5\sqrt2}{4},
\frac{255}{112}+\frac{177\sqrt2}{112},
\frac{61}{28}+\frac{47\sqrt2}{32},
\frac{327}{160}+\frac{1229\sqrt2}{1120},\right.\nonumber\\[-2pt]
&\left.\qquad
\frac{1671}{896}+\frac{67\sqrt2}{112},
\frac{3543}{1792}+\frac{51\sqrt2}{256},
\frac{841}{256}+\frac{215\sqrt2}{512},
\frac{1977}{256}+\frac{613\sqrt2}{256}
\right).\label{eq:Bernstein-B2}
\end{align}
Every displayed Bernstein coefficient is positive in \(K\), so the
three polynomials are positive on \([0,\frac12]\).

For the remaining half-line set \(z=y-\frac12\geq0\).  Their shifted
monomial coefficient vectors, again in increasing order, are
\begin{align}
\mathbf c_0={}&\left(
\frac{75}{8}-6\sqrt2,\,
40-\frac{45\sqrt2}{2},\,
58-23\sqrt2,\,
28-6\sqrt2,\,
10-4\sqrt2\right),\label{eq:shift-B0}\\
\mathbf c_1={}&\left(
\frac{203}{64}+\frac{27\sqrt2}{32},\,
\frac{391}{16}+\frac{63\sqrt2}{8},\,
\frac{1349}{16}+\frac{249\sqrt2}{8},\,
\frac{289}{2}+45\sqrt2,\right.\nonumber\\[-2pt]
&\left.\qquad
\frac{517}{4}+\frac{45\sqrt2}{2},\,
39+14\sqrt2,\,
11-2\sqrt2\right),\label{eq:shift-B1}\\
\mathbf c_2={}&\left(
\frac{1977}{256}+\frac{613\sqrt2}{256},\,
71+\frac{1011\sqrt2}{32},\,
\frac{701}{2}+\frac{3143\sqrt2}{16},\right.\nonumber\\[-2pt]
&\left.\qquad
\frac{3465}{4}+\frac{4065\sqrt2}{8},\,
\frac{9217}{8}+\frac{5507\sqrt2}{8},\,
822+\frac{1157\sqrt2}{2},\,
\right.\nonumber\\[-2pt]
&\left.\qquad
414+227\sqrt2,\,
68+78\sqrt2,\,
45-19\sqrt2\right).\label{eq:shift-B2}
\end{align}
These are also strictly positive in \(K\).  Hence
\(p_{\mathbf B_j}(z+\frac12)>0\) for every \(z\geq0\).  Together with
the preceding Bernstein identities, this proves
\(p_{\mathbf B_j}>0\) on all of \([0,\infty)\).  Every sign just used
is exact: if \(a>0>b\), then \(a+b\sqrt2>0\) exactly when
\(a^2>2b^2\), and the analogous reversed inequality handles
\(b>0>a\).

\paragraph{Reproducibility of the finite calculation.}
The certificate can be regenerated without numerical input as follows.
For each displayed representative, assemble the block from the preceding
formulas over \(K(x,i,\sqrt F)\), form \(B_{\mathcal B}=-4s^2\mathcal
B\), and, for each \(k\), take the numerator and denominator of the
leading \(k\)-minor.  Verify that the denominator is a positive
monomial in \(x\); divide the numerator repeatedly by \(x-1\) until its
value at \(x=1\) is nonzero; then substitute \(x=y+1\).  Exact sign
tests in \(K\) give \(\mathbf C\) in \(69\) cases.  In the other three
cases, exact monomial-to-Bernstein conversion and the shift
\(y=z+\frac12\) give
\eqref{eq:Bernstein-B0}--\eqref{eq:shift-B2}.  As an independent
consistency check, the same
calculation first differentiates the displayed \(F\) and verifies all
seven curvature coefficients, including
\eqref{eq:exact-lambda-pm}, before
constructing any block.  Thus the calculation uses neither sampled
values of \(r\) nor external coefficient data.
The accompanying supplementary archive, tested with Python~3.12.13 and
SymPy~1.14.0, regenerates all \(72\)
leading minors, the \(69\) coefficient-positive certificates, and the
three Bernstein-and-shift certificates directly from the displayed
matrix formulas.  Its README records the software versions and file
manifest, together with SHA-256 hashes of the manuscript-source snapshot,
the exact-certificate subsection, and the embedded coefficient table.
The proof consists of the exact table and identities displayed
here and is independent of executing the code.

There are respectively
\[
 4,\quad8,\quad13,\quad20,\quad27
\]
certified leading minors in the five tables, for a total of \(72\).
Of these, \(69\) have certificate \(\mathbf C\), and the remaining
three have the exact Bernstein-and-shift certificates
\eqref{eq:Bernstein-B0}--\eqref{eq:shift-B2}.  It follows from
\eqref{eq:exact-minor-normal-form} that every leading principal minor
of every safe representative is positive for \(x>1\).  Sylvester's
criterion gives \(B_{\mathcal B}(x)>0\), and therefore
\(\mathcal B(x)<0\).

Finally, complex conjugation together with
\(M\mapsto-M\) and \(+\leftrightarrow-\) identifies the reflected
\(\mathcal Q\)-blocks.  In the anti-invariant sector it identifies a
barred \(\mathcal P\)-block with the correspondingly reflected
unbarred block; it does not identify the two inequivalent unbarred
\(M\)-halves.  The tables therefore display both such halves whenever
they are distinct, including the four additional \(J=3,4\) blocks
above, and omit only their antiunitarily equivalent barred partners.
Hence the conclusion holds for every safe block, not merely for the
displayed representatives.
\end{proof}

\subsection{The radial spectral block}

We record a short argument which determines the entire nonnegative
spectrum in the radial Wigner sector without importing the radial
stability theorem.  Write
\[
 c_0=\sqrt2-1,\qquad
 F(r)=\frac1{\sqrt2}-\frac{c_0}{r^2}
       -\frac{c_0}{\sqrt2\,r^4},\qquad
 f(r)=\sqrt2(r^2-1)-\log(2c_0).
\]
(Only derivatives of \(f\) enter the differential operators.)  The exact radial block
reduction has one \(2\times2\) block on
\(\operatorname{span}\{b_0,b_1\}\); every other radial
block is scalar.  The exact formulas give
\[
 \Lambda_{1+}<0,\qquad \Lambda_{1-}<0,\qquad
 \Lambda_{01}<0,\qquad \Lambda_{23}<0
 \quad\text{on }(1,\infty).
\]
(Equivalently, these four signs are the \(J=0\) rows of the exact
Sylvester--Bernstein certificate.)  Testing a scalar-block eigenvalue
equation therefore shows that every nonzero eigenvalue contribution
outside \(\operatorname{span}\{b_0,b_1\}\) is strictly
negative.

It remains to treat the latter two-component block.  Let
\(\mathcal H_{01}\subset L^2_f((S^2T^*M)_{\mathbb C})\) denote its
radial realization, let
\(\mathcal K_0\subset L^2_f((T^*M)_{\mathbb C})\) denote the radial
\(e^0\)-one-form channel, and let
\[
 D\colon\mathcal H_{01}\longrightarrow\mathcal K_0
\]
be the graph-closed realization of \(\operatorname{div}_f\) specified
in Proposition~\ref{prop:FIK-radial-domain-package} below.  The
closed-operator adjoint decomposition is
\[
 \mathcal H_{01}
 =
 \ker D
 \mathbin{\widehat\oplus}
 \overline{\operatorname{ran}D^*}.
\tag{R1}\label{eq:FIK-radial-york}
\]
The domain proposition below proves that both
summands reduce \(L_f=\Delta_f+2\mathrm{Rm}\) by promoting the core
identities
\[
 \operatorname{div}_f L_fh
 =\left(\Delta_f+\frac12\right)\operatorname{div}_fh,
 \qquad
 L_f\nabla^2u=\nabla^2(\Delta_fu+u),
\tag{R2}\label{eq:FIK-radial-commutators}
\]
and self-adjointness turns invariance into reduction.  Every smooth
one-form in the \(e^0\)-channel is exact.  If
\(\Pi_{01}\) denotes orthogonal projection onto
\(\operatorname{span}\{b_0,b_1\}\), weighted
integration by parts gives
\[
 D^*(du)=-\Pi_{01}\nabla^2u.
\]
Thus the second summand in \eqref{eq:FIK-radial-york} is the closure of
the \emph{projected} radial Hessians.  The projection is necessary
because a general radial Hessian also has a \(b_6\)-component.
Formula~\eqref{eq:FIK-radial-york} is the identity
\(\overline{\operatorname{ran}D^*}=(\ker D)^\perp\), applied to the
closed radial realization of \(D=\operatorname{div}_f\); it does not
require \(\operatorname{ran}D^*\) to be closed.

Put
\[
 \mathcal D_1a
 =\frac12\bigl(ra'+2(1-\sqrt2r^2)a\bigr),\qquad
 h_a=(a+\mathcal D_1a)b_0+(a-\mathcal D_1a)b_1.
\]
The radial divergence equation says exactly that every smooth tensor in
the kernel under consideration is \(h_a\) for a unique \(a\).  Direct
substitution in the two-component block gives
\[
 L_fh_a=h_{P(a)},\qquad
 P(a)=\frac F4\left(
 a''+\left(\frac5r-2\sqrt2r\right)a'-8\sqrt2a\right)+a.
 \tag{R3}\label{eq:FIK-radial-P}
 \]
The Ricci tensor is \(h_{-c_0r^{-4}}\).  For every bolt-regular \(a\)
and every \(R>1\), direct contraction with this tensor gives the exact
radial derivative
\[
 \left\langle h_a,\operatorname{Ric}\right\rangle_
 {L^2_f(\{1<r<R\})}
 =
 8\pi^2c_0\int_1^R
 \frac d{dr}\left[
  (1+\sqrt2r^2)a(r)e^{-f(r)}
 \right]\,dr .
\]
Since
\((1+\sqrt2)e^{-f(1)}=2\), one radial integration gives the
finite-cutoff identity
\begin{equation}\label{eq:FIK-radial-Ric-pairing-cutoff}
 \left\langle h_a,\operatorname{Ric}\right\rangle_
 {L^2_f(\{1<r<R\})}
 =
 8\pi^2c_0\left[
  (1+\sqrt2R^2)a(R)e^{-f(R)}-2a(1)
 \right].
\end{equation}

\begin{proposition}[Radial Friedrichs domains, divergence reduction,
and endpoint control]
\label{prop:FIK-radial-domain-package}
For \(j=0,1,2\), let
\(\mathscr C^{(j)}_{\mathrm{br},c}\) denote the restrictions to the open
orbit of smooth compactly supported radial sections on the completed
FIK manifold: respectively scalars, radial \(e^0\)-one-forms, and
tensors in \(\operatorname{span}\{b_0,b_1\}\).  All closures below use
the raw weighted measure \(e^{-f}\,dV_{\bar g}\).
Their ambient Hilbert spaces are the corresponding radial subspaces of
\[
 L^2_f(\mathbb C),\qquad
 L^2_f((T^*M)_{\mathbb C}),\qquad
 L^2_f((S^2T^*M)_{\mathbb C}),
\]
respectively, and their Friedrichs form domains are the closures of
these cores in the corresponding bundlewise \(H^1_f\) norms.

Define on \(\mathscr C^{(2)}_{\mathrm{br},c}\)
\[
 (D_0h)_i=\nabla^jh_{ij}-(\nabla^jf)h_{ij},
 \qquad D=\overline{D_0},
\]
and let \(D^*\) be its Hilbert-space adjoint.  Let
\[
 T_0=\Delta_f^{(0)},\qquad
 T_1=\Delta_f^{(1)}+\frac12,\qquad
 T_2=L_f\big|_{\mathcal H_{01}}
\]
be the self-adjoint realizations induced by the corresponding shifted
Friedrichs forms on the three radial channels.  Then:
\begin{enumerate}
\item
Each \(\mathscr C^{(j)}_{\mathrm{br},c}\) is a form core and an
operator core for \(T_j\);
\(\mathscr C^{(2)}_{\mathrm{br},c}\) is a graph core for \(D\), and
\(\mathscr C^{(1)}_{\mathrm{br},c}\) is a graph core for \(D^*\).
On these cores,
\[
 D^*(du)=-\Pi_{01}\nabla^2u,
 \qquad
 D_0T_2=T_1D_0.
\]
Moreover, the formal commutator has the following closed-domain form:
\begin{equation}\label{eq:FIK-radial-closed-intertwining}
 h\in D(T_2),\quad T_2h\in D(D)
 \quad\Longrightarrow\quad
 Dh\in D(T_1),\qquad T_1Dh=D(T_2h).
\end{equation}
Consequently
\[
 \mathcal H_{01}
 =\ker D\mathbin{\widehat\oplus}
   \overline{\operatorname{ran}D^*},
\]
and both summands reduce \(T_2\).

\item
Put
\[
 w=r^3e^{-f},\qquad
 p_{\rm SL}=\frac14Fr^3e^{-f}.
\]
For every bounded real-valued radial potential \(V\) that extends
smoothly across the bolt, in particular
\(V=0\) and \(V=3/(2r^2)\), the scalar expression
\[
 T_Vu=\frac1w(p_{\rm SL}u')'+Vu
\]
has the following Friedrichs endpoint realization.  At \(r=1\) the
regular solution is selected and the logarithmic solution is excluded;
at infinity the expression is limit point.  Its resolvent is compact,
its eigenvalues are simple, and, when they are ordered decreasingly,
the \(k\)-th eigenfunction has exactly \(k\) zeros in
\((1,\infty)\), for \(k=0,1,\ldots\).

Every operator-domain eigenfunction is regular at the bolt and
satisfies all endpoint limits used below.  In particular,
\[
 \lim_{r\downarrow1}p_{\rm SL}u'(r)=0,\qquad
 \lim_{r\to\infty}p_{\rm SL}(r)u'(r)\overline{u(r)}=0,
\]
and, for \(\psi_+(r)=1-\frac{3}{2r^2}\),
\[
 \lim_{r\to\infty}
 \frac{Fw}{4}\frac{\psi_+'}{\psi_+}|u|^2=0.
\]
If an eigenfunction and a comparison function have a common simple
interior zero, their quotient has a finite limit there and the
corresponding ground-state flux vanishes.

\item
Let \(P_{\rm expr}\) denote the radial differential expression \(P\)
in \eqref{eq:FIK-radial-P}.  Suppose
\(h_a\in D(T_2)\cap\ker D\) and
\(T_2h_a=\lambda h_a\), where
\[
 h_a=(a+\mathcal D_1a)b_0+(a-\mathcal D_1a)b_1,\qquad
 \mathcal D_1a
 =\frac12\bigl(ra'+2(1-\sqrt2r^2)a\bigr).
\]
Then \(a\) is smooth on the open orbit,
\(P_{\rm expr}a=\lambda a\), and
 \[
 (1+r^2)a(r)e^{-f(r)}\longrightarrow0.
 \]
Passing \(R\to\infty\) in
\eqref{eq:FIK-radial-Ric-pairing-cutoff} gives
\begin{equation*}\tag{R4}\label{eq:FIK-radial-Ric-pairing}
 \begin{split}
 \langle h_a,\operatorname{Ric}\rangle_{L^2_f}
 &=8\pi^2c_0\left[
   \lim_{R\to\infty}(1+\sqrt2R^2)a(R)e^{-f(R)}-2a(1)
  \right]\\
 &=-16\pi^2c_0a(1).
 \end{split}
\end{equation*}
In particular, \(h_a\perp\operatorname{Ric}\) if and only if
\(a(1)=0\), and in that case
\[
 a\in\mathcal Q_D(P):=
 \left\{a\in AC_{\mathrm{loc}}([1,\infty)):
 \begin{array}{l}
 a(1)=0,\\
 \displaystyle
 \int_1^\infty
 \left(
   \frac{r^5e^{-f}}{F}|a|^2+
   r^5e^{-f}|a'|^2
 \right)\,dr<\infty
 \end{array}\right\}.
\]
On \(\mathcal Q_D(P)\), define the closed upper-semibounded
sesquilinear form
\[
 \begin{split}
 \mathfrak p_D[a,b]
 :={}&-\frac14\int_1^\infty
       r^5e^{-f}a'\overline{b'}\,dr\\
 &+\int_1^\infty
       \rho(1-2\sqrt2F)a\overline b\,dr,
 \qquad
 \rho:=\frac{r^5e^{-f}}F .
 \end{split}
\]
Moreover,
\[
 C_c^\infty((1,\infty))
\]
is a form core for \(\mathcal Q_D(P)\), where these compactly
supported functions are understood as Dirichlet data at the bolt.
Let \(P_F\) be the self-adjoint Dirichlet realization determined by
this upper-semibounded closed form; equivalently, \(P_F\) is the
negative of the Friedrichs realization of \(-P_{\rm expr}\).  Its
operator domain is
\[
 \begin{split}
 D(P_F)
 =\bigl\{a\in\mathcal Q_D(P):\;&
 \text{there is }g\in L^2((1,\infty),\rho\,dr)
 \text{ such that}\\
 &\mathfrak p_D[a,b]
   =\langle g,b\rangle_{L^2(\rho\,dr)}
   \text{ for every }b\in\mathcal Q_D(P)\bigr\},
 \end{split}
\]
and \(P_Fa=g\).  Equivalently,
\[
 D(P_F)
 =
 \left\{a\in\mathcal Q_D(P):
  P_{\rm expr}a\in L^2((1,\infty),\rho\,dr)
  \text{ distributionally on }(1,\infty)\right\},
\]
where
\[
 P_{\rm expr}a
 =
 \frac1{\rho}\left(\frac14r^5e^{-f}a'\right)'
 +(1-2\sqrt2F)a .
\]
For the Ricci-orthogonal eigenfunction above,
\(P_{\rm expr}a=\lambda a\) and \(a\in\mathcal Q_D(P)\) therefore give
\[
 a\in D(P_F),\qquad P_Fa=\lambda a.
\]
The ground-state identity is the following form identity, valid for
every \(a\in\mathcal Q_D(P)\):
\[
 \mathfrak p_D[a,a]
 =
 -\int_1^\infty
 \left[
  \sqrt2F^2\left|\frac aF\right|^2+
  \frac14F^2
       \left|\left(\frac aF\right)'\right|^2
  \right]r^5e^{-f}\,dr .
\]
For \(a\in D(P_F)\), its left side is
\(\langle P_Fa,a\rangle_{L^2(\rho\,dr)}\).

\item
The radial gradient map has closed range, and its range is the entire
radial \(e^0\)-one-form channel.  More precisely,
\[
 \left\|
 u-\frac{\langle u,1\rangle_{L^2_f}}
          {\|1\|_{L^2_f}^{\,2}}\,1
 \right\|_{L^2_f}
 \leq\|du\|_{L^2_f},
 \qquad
 \|1\|_{L^2_f}^{\,2}=(4\pi)^2.
\]
Hence the nonconstant radial scalar eigenfunctions differentiate to a
complete eigenbasis in the radial one-form channel.
\end{enumerate}
\end{proposition}

\begin{proof}
For \(j=0,1,2\), let \(\widetilde T_j\) denote the ambient
self-adjoint realization of, respectively, \(\Delta_f^{(0)}\),
\(\Delta_f^{(1)}+\frac12\), and \(L_f\) on the full corresponding
bundle Hilbert space.  Conjugation by \(e^{-f/2}\) takes
\(C-\widetilde T_j\), for one fixed sufficiently large \(C\), to a
generalized Schr\"odinger operator on the complete smooth FIK manifold.
Its negative zeroth-order part is bounded, and its potential satisfies
\[
 \frac14|\nabla f|^2-\frac12\Delta f+O(1)
 =c\,r^2+O(1),\qquad c>0.
\]
Thus \cite[Corollary~2.9]{BravermanMilatovicShubin} gives
\(C_c^\infty\) as an operator core in every ambient bundle channel used
here.  The quadratic lower bound at infinity, together with the
first-moment estimate and local Rellich compactness, makes the ambient
form-domain embedding compact; equivalently, the ambient resolvent is
compact.

Let \(P_j\) denote radial averaging followed by the orthogonal
projection onto the \(j\)-th invariant block.  Compact-group averaging
and the finite-dimensional fiber projection preserve smoothness and
compact support.  By invariance of the differential expression and
the block splitting, \(P_j\) commutes on \(C_c^\infty\) with the
ambient operator, and hence with its closed form, self-adjoint
realization \(\widetilde T_j\), and resolvent.  Consequently
\(\operatorname{ran}P_j\) reduces \(\widetilde T_j\), and the
restriction of \(\widetilde T_j\) to this range is precisely \(T_j\),
by the representation theorem for the restricted closed form.  If
\(u\in D(T_j)\), choose ambient \(u_n\in C_c^\infty\) such that
\[
 u_n\longrightarrow u,
 \qquad
 \widetilde T_ju_n\longrightarrow T_ju
\]
in the ambient Hilbert space.  Then
\(P_ju_n\in\mathscr C_{\mathrm{br},c}^{(j)}\) and
\[
 P_ju_n\longrightarrow u,\qquad
 T_jP_ju_n=P_j\widetilde T_ju_n\longrightarrow T_ju.
\]
Thus \(\mathscr C_{\mathrm{br},c}^{(j)}\) is an operator core.
Compactness of the ambient resolvent passes to the reducing range of
\(P_j\).  Finally, the logarithmic bolt cutoff used in the proof of
\eqref{eq:self-mode-splitting}, followed by the same radial and block
projections, gives the asserted form core.

The graph-core assertion for \(D\) is part of the definition
\(D=\overline{D_0}\).  We prove the adjoint assertion separately.  On
\(\mathscr C_{\mathrm{br},c}^{(1)}\) put
\[
 B_0q:=-\Pi_{01}(\nabla q)_{\mathrm{sym}}.
\]
Weighted integration by parts on the completed manifold gives
\(D_0^\dagger=B_0\).  Since the Hilbert adjoint of a minimal
differential operator is its maximal distributional realization,
\begin{equation}\label{eq:FIK-radial-adjoint-maximal-domain}
 D(D^*)
 =\left\{q\in\mathcal K_0:
   B_0q\in\mathcal H_{01}\ \text{distributionally}\right\},
 \qquad D^*q=B_0q.
\end{equation}

We verify directly that the completed bolt introduces no hidden
boundary condition in this maximal domain.  Write \(q=Q(r)e^0\) and
set
\[
 A_Q=\frac{\sqrt F}{2}Q',\qquad
 B_Q=\frac{2F+rF'}{4r\sqrt F}\,Q,\qquad
 C_Q=\frac{\sqrt F}{2r}\,Q.
\]
In the radial orthonormal frame,
\[
 \nabla q=\operatorname{diag}(A_Q,B_Q,C_Q,C_Q),
\]
and the normalization of \(b_0,b_1\) gives
\[
 B_0q
 =-(A_Q+B_Q+2C_Q)b_0
  -(A_Q+B_Q-2C_Q)b_1.
\]
Let \(\varrho\) be completed normal distance from the bolt.  Choose
and fix \(0<\varrho_0\ll1\) so that the completed normal coordinate is
defined on \(0\leq\varrho<\varrho_0\), all profile and measure
comparisons below hold there with uniform constants, and the
\(O(\varrho^2Q)\) contribution in the ensuing Hardy estimate is
absorbable.  The exact profile then has, uniformly for
\(0<\varrho<\varrho_0\),
\[
 A_Q=\partial_\varrho Q,\qquad
 B_Q=(\varrho^{-1}+O(\varrho))Q,\qquad
 C_Q=O(\varrho)Q,\qquad
 w\,dr\asymp\varrho\,d\varrho.
\]
Here ``distributionally'' in
\eqref{eq:FIK-radial-adjoint-maximal-domain} means distributionally on
the completed manifold.  Put \(Z(\varrho)=\varrho Q(\varrho)\).  Since
\(A_Q+B_Q\in L^2(\varrho\,d\varrho)\) and
\(Q\in L^2(\varrho\,d\varrho)\), the identity
\[
 Z'
 =\varrho(A_Q+B_Q)+O(\varrho^2Q)
\]
shows that \(Z\) has an absolutely continuous representative with a
finite trace at \(0\).  A nonzero trace would give
\(Q(\varrho)\sim Z(0)/\varrho\), contradicting
\(Q\in L^2(\varrho\,d\varrho)\).  Thus \(Z(0)=0\), and the weighted
Hardy inequality applies.
The one-dimensional Hardy estimate applied to
\((\varrho Q)'=\varrho(\partial_\varrho Q+\varrho^{-1}Q)\) therefore
gives, for this fixed choice of \(\varrho_0\),
\begin{equation}\label{eq:FIK-radial-adjoint-bolt-Hardy}
 \int_0^{\varrho_0/2}
 \left(
   \varrho|\partial_\varrho Q|^2
   +\frac{|Q|^2}{\varrho}
 \right)d\varrho
 \leq C\left(
  \|q\|_{L^2_f(\{\varrho<\varrho_0\})}^2
  +\|B_0q\|_{L^2_f(\{\varrho<\varrho_0\})}^2
 \right).
\end{equation}
Choose a logarithmic cutoff \(\chi_\varepsilon\), equal to zero for
\(\varrho\leq\varepsilon^2\), equal to one for
\(\varrho\geq\varepsilon\), and satisfying
\(|\partial_\varrho\chi_\varepsilon|
 \leq C/(\varrho|\log\varepsilon|)\).  Then
\[
 \|[B_0,\chi_\varepsilon]q\|_{L^2_f}^2
 \leq\frac{C}{|\log\varepsilon|^2}
 \int_{\varepsilon^2}^{\varepsilon}
 \frac{|Q|^2}{\varrho}\,d\varrho
 \longrightarrow0.
\]
Together with dominated convergence and
\eqref{eq:FIK-radial-adjoint-bolt-Hardy}, this proves
\(\chi_\varepsilon q\to q\) in the graph norm of \(D^*\).
At infinity choose a radial cutoff \(\eta_R\), equal to one on
\(\{r\leq R\}\), zero on \(\{r\geq2R\}\), with
\(|\nabla\eta_R|\leq C/R\).  Then
\[
 \|[B_0,\eta_R]q\|_{L^2_f}
 \leq\frac CR\|q\|_{L^2_f(\{R<r<2R\})}
 \longrightarrow0,
\]
and the remaining graph-norm tails vanish by dominated convergence.
After these cutoffs the coefficient is supported in a compact
subinterval of \((1,\infty)\); ordinary one-dimensional mollification
gives elements of \(\mathscr C_{\mathrm{br},c}^{(1)}\) converging in
the graph norm.  Because the approximants vanish near the bolt, their
zero extensions are smooth sections on the completed manifold.  Thus
\(\mathscr C_{\mathrm{br},c}^{(1)}\) is a graph core for \(D^*\).

Weighted integration by parts on the common core gives
\[
 D^*(du)=-\Pi_{01}\nabla^2u
\]
 and the first identity in
 \eqref{eq:FIK-radial-commutators}.  The orthogonal decomposition follows
 from the closed-operator adjoint theorem.  We next justify the
 divergence-domain inclusion required for the closed commutator and
 resolvent argument, rather than inferring it from form density alone.
 Since
\(\nabla f=\nabla\bar f\), the soliton identities give
\(|\nabla f|^2\leq\bar f\).  The raw and normalized Gaussian measures
differ by one fixed positive factor, so
Lemma~\ref{lem:first-moment} implies, for
\(u\in\mathscr C_{\mathrm{br},c}^{(2)}\),
\begin{equation}\label{eq:FIK-radial-divergence-H1-bound}
 \begin{split}
 \|D_0u\|_{L^2_f}
 &\leq
 C\left(
  \|\nabla u\|_{L^2_f}
  +\||\nabla f|u\|_{L^2_f}
  +\|u\|_{L^2_f}\right)\\
 &\leq
 C\left(
  \|\nabla u\|_{L^2_f}
  +\|\sqrt{\bar f}\,u\|_{L^2_f}
  +\|u\|_{L^2_f}\right)
 \leq C\|u\|_{H^1_f}.
 \end{split}
\end{equation}
For \(u\in H^1_f\cap\mathcal H_{01}\), choose
\(u_n\in\mathscr C_{\mathrm{br},c}^{(2)}\) converging to \(u\) in
\(H^1_f\).  Estimate
\eqref{eq:FIK-radial-divergence-H1-bound} makes \(D_0u_n\) Cauchy in
\(L^2_f\).  Because \(D=\overline{D_0}\), graph closedness gives
\[
 H^1_f\cap\mathcal H_{01}\subset D(D),\qquad
 \|Du\|_{L^2_f}\leq C\|u\|_{H^1_f}.
\]
We next promote the formal commutator to the asserted closed-domain
identity.  If \(h\in D(T_2)\) and \(T_2h\in D(D)\), then
\(h\in H^1_f\cap\mathcal H_{01}\subset D(D)\), and for every
\(q\in\mathscr C_{\mathrm{br},c}^{(1)}\),
\[
\begin{aligned}
 \langle Dh,T_1q\rangle_{L^2_f}
 &=\langle h,D^*T_1q\rangle_{L^2_f}
  =\langle h,T_2D^*q\rangle_{L^2_f}\\
 &=\langle T_2h,D^*q\rangle_{L^2_f}
  =\langle D(T_2h),q\rangle_{L^2_f}.
\end{aligned}
\]
Here the middle equality is the adjoint core commutator obtained from
\(D_0T_2=T_1D_0\).  Essential self-adjointness of the displayed
\(T_1\)-core, equivalently its maximal-distributional
characterization supplied by the conjugated Schr\"odinger
realization, now gives
\[
 Dh\in D(T_1),\qquad T_1Dh=D(T_2h),
\]
which proves \eqref{eq:FIK-radial-closed-intertwining}.

To prove reduction, choose a real \(z\) larger than the upper spectral
bounds of \(T_1\) and \(T_2\).  If \(g\in\ker D\) and
\(h=(z-T_2)^{-1}g\), then \(h\in H^1_f\subset D(D)\) by the preceding
estimate, and
\(T_2h=zh-g\in D(D)\).  The closed-domain intertwining therefore gives
\[
 T_1Dh=D(zh-g)=zDh.
\]
Since \(z\in\rho(T_1)\), one has \(Dh=0\).  Thus
\((z-T_2)^{-1}\) preserves \(\ker D\).  This real resolvent is
self-adjoint and therefore also preserves
\((\ker D)^\perp=\overline{\operatorname{ran}D^*}\).  Hence both
summands reduce \(T_2\).

For the scalar endpoint statement, write \(x=r-1\).  Since
\[
 F=2x+O(x^2),\qquad w=w(1)+O(x),\qquad
 p_{\rm SL}=cx+O(x^2),
\]
the two local solutions are \(1+O(x)\) and
\(\log x+O(x\log x)\).  The latter has
\(\int p_{\rm SL}|u'|^2\,dr=\infty\), so the Friedrichs form excludes it.
At infinity the preceding unitary conjugation gives a half-line
Schr\"odinger operator with potential tending quadratically to
\(+\infty\).  Hence infinity is limit point and the spectrum is
discrete.  Simplicity follows from constancy of the Wronskian.  The
nodal ordering follows as follows.  Apply the preceding endpoint
analysis to the lower-semibounded operator \(S_V:=-T_V\).  For \(R>1\),
let \(S_{V,R}\) be its realization on \((1,R)\) with the same
Friedrichs condition at \(r=1\) and a Dirichlet condition at \(R\).
After zero extension, these form domains increase with \(R\), and their
union contains the radial scalar form core
\(\mathscr C^{(0)}_{\mathrm{br},c}\).  Monotone convergence of the
closed forms therefore gives strong-resolvent convergence
\(S_{V,R}\to S_V\).  The min--max principle and form-core density give
\[
 \mu_k(R)\downarrow\mu_k
\]
for every increasingly ordered eigenvalue of \(S_V\).

The local expansions above show that \(r=1\) is limit-circle
nonoscillatory and that the Friedrichs condition selects its principal,
regular solution; the conjugated confining potential shows that
infinity is limit point and that the spectrum is discrete.  The
singular Sturm oscillation theorem
\cite[Theorem~10.12.1]{Zettl} therefore applies directly to \(S_V\):
its \(k\)-th eigenfunction has exactly \(k\) zeros in
\((1,\infty)\).  Returning to \(T_V=-S_V\) gives the asserted
decreasing eigenvalue ordering.  Thus the nodal count follows from the
singular Sturm oscillation theorem independently of the preceding
strong-resolvent convergence.

For completeness, the endpoint estimates require no formal asymptotic
series.  On \(I_r=[r,r+r^{-1}]\), the weights
\(q_m=r^me^{-f}\) are uniformly comparable, and the
one-dimensional Sobolev inequality gives
\[
 q_m(r)|v(r)|^2
 \leq C\left(
 r\int_{I_r}q_m|v|^2+
 r^{-1}\int_{I_r}q_m|v'|^2
 \right).
\]
Let \(T_Vu=\lambda u\) and put \(y=p_{\rm SL}u'\).  The eigenvalue equation is
\[
 y'=w(\lambda-V)u.
\]
Because \(V\) is bounded and \(w\,dr\) has finite mass, the right-hand
side is in \(L^1\), so \(y\) has a finite limit at infinity.  That limit
must be zero: a nonzero limit, together with
\(p_{\rm SL}\asymp r^3e^{-f}\), would produce the non-\(L^2_f\) fast solution.
Consequently
\[
 |p_{\rm SL}(r)u'(r)|
 \leq C
 \left(\int_r^\infty w|u|^2\,ds\right)^{1/2}
 \left(\int_r^\infty w\,ds\right)^{1/2},
 \qquad
 \int_r^\infty w\,ds\leq Cr^2e^{-f(r)}.
\]
The preceding trace inequality, applied with \(m=3\), also gives
\[
 w(r)|u(r)|^2
 \leq C\left(
 r\int_r^\infty w|u|^2\,ds+
 r^{-1}\int_r^\infty p_{\rm SL}|u'|^2\,ds
 \right).
\]
Both tail integrals tend to zero; combining the last two displays
proves \(p_{\rm SL}(r)u'(r)\overline{u(r)}\to0\).  Since
\(\psi_+'/\psi_+=O(r^{-3})\), the same trace estimate gives
\[
 \frac{Fw}{4}\frac{\psi_+'}{\psi_+}|u|^2\longrightarrow0.
\]
At the bolt the fluxes vanish by regularity and \(F(1)=0\).  At a
common simple zero, ODE uniqueness gives a finite quotient and
vanishing flux.

An eigenvector \(h_a\) is smooth by elliptic regularity, and direct
substitution gives \(P_{\rm expr}a=\lambda a\).  We now justify the
polynomially weighted bounds used at the endpoints.  Put
\[
 \widetilde h_a:=e^{-f/2}h_a,\qquad
 \mathscr H_2:=e^{-f/2}(C-T_2)e^{f/2},\qquad
 \vartheta:=1+f,
\]
where \(C\) is the fixed constant chosen at the beginning of the proof.
Then
\[
 \mathscr H_2\widetilde h_a=(C-\lambda)\widetilde h_a.
\]
The zeroth-order endomorphism \(V_2\) of \(\mathscr H_2\) satisfies
\[
 V_2\geq c_1\vartheta-C_1
\]
as a quadratic form, by the confining-potential calculation above.

Choose a nonincreasing smooth cutoff \(\chi\), equal to one on
\((-\infty,1]\) and zero on \([2,\infty)\), and put
\(\chi_R=\chi(f/R)\).  Since \(|\nabla f|^2\leq f\),
\[
 |\nabla\chi_R|^2
 \leq\frac C R\,\mathbf 1_{\{R\leq f\leq2R\}},
 \qquad
 \frac{|\nabla(\vartheta^m)|^2}{\vartheta^m}
 \leq C_m\vartheta^{m-1}.
\]
Testing the eigenvalue equation against
\(\chi_R^2\vartheta^m\widetilde h_a\), integrating by parts, and using
Young's inequality gives, uniformly in \(R\),
\[
 \begin{split}
 &\int_M\chi_R^2\vartheta^m
   \left(
    |\nabla\widetilde h_a|^2
    +\vartheta|\widetilde h_a|^2
   \right)dV_{\bar g}\\
 &\qquad\leq
 C_m\left(
  \|\widetilde h_a\|_{L^2(dV_{\bar g})}^2
  +\int_M\vartheta^{m-1}|\widetilde h_a|^2\,dV_{\bar g}
 \right),\qquad m\geq1.
 \end{split}
\]
For \(m=0\), the same estimate holds with only the unweighted
\(L^2\)-norm on the right.  Indeed, the weight-derivative term is
bounded by \(C_m\vartheta^{m-1}|\widetilde h_a|^2\), while on the cutoff
annulus \(\vartheta^m|\nabla\chi_R|^2\leq C_m\vartheta^{m-1}\); the
term \(C\vartheta^m|\widetilde h_a|^2\) is absorbed by
\(c_1\vartheta^{m+1}|\widetilde h_a|^2\) outside one fixed compact set.

Letting \(R\to\infty\) and inducting on \(m\) proves every polynomial
moment of \(\widetilde h_a\) and its first derivative.  Since
\[
 e^{-f/2}\nabla h_a
 =\nabla\widetilde h_a+\frac12\,df\otimes\widetilde h_a,
 \qquad |df|^2\leq f,
\]
this is exactly the asserted family of polynomially weighted
\(H^1_f\) bounds for \(h_a\).  The coefficients \(a\) and
\(\mathcal D_1a\) are fixed radial components of \(h_a\), and the
displayed formula for \(\mathcal D_1a\) recovers \(a'\) with only
polynomial coefficients.  Bolt regularity and the preceding trace
inequality with \(m=5\) consequently give
\[
 (1+r^2)a(r)e^{-f(r)}\to0.
\]
Thus the passage to the limit in
\eqref{eq:FIK-radial-Ric-pairing-cutoff} proves
\eqref{eq:FIK-radial-Ric-pairing} and the asserted equivalence between
Ricci orthogonality and \(a(1)=0\).  Bolt regularity then gives
\(a=O(r-1)\), hence \(a^2/F=o(1)\) and
\(\int a^2/F<\infty\) near \(r=1\); at infinity the polynomially
weighted bounds give the two integrals defining \(\mathcal Q_D(P)\).

We now prove the asserted form-core and maximal-distributional
characterization independently of this particular eigenfunction.  Put
\[
 q=r^5e^{-f},\qquad A=\frac14q,\qquad
 \rho=\frac qF,\qquad V=1-2\sqrt2F.
\]
Because \(0<F<1/\sqrt2\), the shifted form norm
\[
 \|b\|_*^2
 :=2\|b\|_{L^2(\rho\,dr)}^2-\mathfrak p_D[b,b]
 =\int_1^\infty A|b'|^2\,dr
  +\int_1^\infty\rho(2-V)|b|^2\,dr
\]
is equivalent to the defining norm of \(\mathcal Q_D(P)\); indeed
\(1<2-V<3\).
We claim that \(C_c^\infty((1,\infty))\) is dense in this norm.
Write \(x=r-1\).  Since
\[
 F=2x+O(x^2),\qquad A\asymp1,\qquad\rho\asymp x^{-1}
\]
near the bolt, choose a cutoff \(\chi_\varepsilon\) which is zero for
\(x\leq\varepsilon\), one for \(x\geq2\varepsilon\), and satisfies
\(|\chi_\varepsilon'|\leq C\varepsilon^{-1}\).  The trace condition
\(b(1)=0\) and the one-dimensional Poincar\'e inequality give
\[
 \begin{split}
 \int_1^{1+2\varepsilon}
 A|\chi_\varepsilon'b|^2\,dr
 &\leq
 C\varepsilon^{-2}
 \int_0^{2\varepsilon}|b(1+x)|^2\,dx\\
 &\leq
 C\int_0^{2\varepsilon}|b'(1+x)|^2\,dx
 \longrightarrow0 .
 \end{split}
\]
The remaining derivative term and the \(L^2(\rho\,dr)\) term tend to
zero by absolute continuity of the defining integrals.  At infinity
choose \(\eta_R=1\) on \([1,R]\), \(\eta_R=0\) on
\([2R,\infty)\), and \(|\eta_R'|\leq C/R\).  Since
\(A/\rho=F/4\leq1/(4\sqrt2)\),
\[
 \int_R^{2R}A|\eta_R'b|^2\,dr
 \leq\frac C{R^2}\int_R^{2R}\rho|b|^2\,dr
 \longrightarrow0,
\]
and the other two tails vanish in the same way.  Thus
\(\chi_\varepsilon\eta_Rb\to b\) in \(\|\cdot\|_*\); ordinary
mollification on compact subintervals of \((1,\infty)\) proves the
core claim.

Now let \(a\in\mathcal Q_D(P)\) and suppose
\(P_{\rm expr}a=g\in L^2(\rho\,dr)\) distributionally on
\((1,\infty)\).  Equivalently,
\[
 (Aa')'=\rho(g-Va)
\]
in \(\mathcal D'((1,\infty))\).  Hence for every
\(\phi\in C_c^\infty((1,\infty))\),
\[
 \mathfrak p_D[a,\phi]
 =-\int_1^\infty Aa'\overline{\phi'}\,dr
   +\int_1^\infty\rho Va\overline\phi\,dr
 =\int_1^\infty\rho g\overline\phi\,dr .
\]
Both sides are continuous in \(\phi\) for the \(\|\cdot\|_*\)-norm,
so core density extends the identity to every
\(\phi\in\mathcal Q_D(P)\).  Thus \(a\in D(P_F)\) and \(P_Fa=g\).
Conversely, the variational identity defining \(D(P_F)\), tested on
\(C_c^\infty((1,\infty))\), gives the same distributional equation.
This proves the claimed equivalence and, with \(g=\lambda a\), places
the Ricci-orthogonal eigenfunction in \(D(P_F)\).

Finally, the same core density permits the ground-state calculation
to be made first on compactly supported Dirichlet data and then
passed to the full form domain.  It proves the displayed identity for
\(\mathfrak p_D\).  The form-representation theorem identifies its
left side with
\(\langle P_Fa,a\rangle_{L^2(\rho\,dr)}\) precisely when
\(a\in D(P_F)\).

Let \(\{u_k\}_{k\geq0}\) be an orthonormal radial scalar eigenbasis,
with \(T_0u_k=-\mu_ku_k\).  For \(k\geq1\), the scalar form identity
gives \(du_k\in L^2_f\), while the core commutator gives
distributionally
\[
 T_1(du_k)=d(T_0u_k+u_k)=(1-\mu_k)du_k.
\]
The maximal-distributional characterization of \(T_1\) therefore
places \(du_k\) in \(D(T_1)\).  Finite scalar spectral sums are dense
in the scalar form domain, so their gradients are dense in the
gradient range.

Finally, for \(V=0\) the constant function is the positive top
eigenfunction, while direct substitution gives
\(\Delta_f^{(0)}(r^2-2)=-(r^2-2)\).  Part~(2) and the single zero of
\(r^2-2\) therefore identify \(-1\) as the first nonzero eigenvalue.
This gives the correctly normalized Poincar\'e inequality and hence
closed range.  A radial one-form orthogonal to all gradients is weakly weighted coclosed and
therefore has the form
\[
 q(r)e^0,\qquad
 q(r)=\frac{Ce^{f(r)}}{r^3\sqrt{F(r)}}.
\]
This is not in \(L^2_f\) at either the bolt or infinity unless \(C=0\).
Thus the gradient range is the entire radial one-form channel.
\end{proof}

\begin{lemma}\label{lem:FIK-radial-divfree}
On the radial, weighted-divergence-free part of
\(\operatorname{span}\{b_0,b_1\}\), the only nonnegative eigenvalue of
\(L_f\) is \(1\), and its eigenspace is
\(\operatorname{span}\{\operatorname{Ric}\}\).
\end{lemma}

\begin{proof}
The Ricci tensor also satisfies
\(L_f\operatorname{Ric}=\operatorname{Ric}\).

On the Ricci-orthogonal subspace, Part~(3) of
Proposition~\ref{prop:FIK-radial-domain-package} identifies \(P_F\) as
the self-adjoint Dirichlet realization associated with
\(\mathfrak p_D\) in
\(L^2((1,\infty),\rho\,dr)\), where
\[
 \rho(r)=\frac{r^5e^{-f(r)}}{F(r)}.
\]
The exact supersolution identity
\[
 P_{\rm expr}(F)=-\sqrt2F^2
 \tag{R5}\label{eq:FIK-radial-PF}
\]
and the ground-state transform yield, for every
\(a\in\mathcal Q_D(P)\),
\begin{equation}\tag{R6}\label{eq:FIK-radial-ground}
\begin{split}
 \mathfrak p_D[a,a]
 ={}&-\int_1^\infty\left[
 \sqrt2F^2\left|\frac aF\right|^2
 +\frac14F^2\left|\left(\frac aF\right)'\right|^2
 \right]r^5e^{-f}\,dr .
\end{split}
\end{equation}
first for smooth Dirichlet data of compact support and then, by
Part~(3) of Proposition~\ref{prop:FIK-radial-domain-package}, on the
full Dirichlet form domain.  If \(a\in D(P_F)\), then the left side is
\(\langle P_Fa,a\rangle_{L^2(\rho\,dr)}\).  Thus, if
\(P_Fa=\lambda a\) and \(a\ne0\), then
\(\lambda\int\rho|a|^2<0\), so \(\lambda<0\).

Now let \(h_a\ne0\) be an eigenvector of \(L_f\) in this radial,
weighted-divergence-free subspace, with
\[
 L_fh_a=\lambda h_a,
 \qquad \lambda\geq0.
\]
If \(\lambda\ne1\), self-adjointness makes \(h_a\) orthogonal to
\(\operatorname{Ric}\), and the preceding argument is a contradiction.
If \(\lambda=1\), subtract the orthogonal projection of \(h_a\) onto
\(\operatorname{span}\{\operatorname{Ric}\}\).  The remainder is a
Ricci-orthogonal eigenvector with eigenvalue \(1\), and hence vanishes.
Thus \(h_a\in\operatorname{span}\{\operatorname{Ric}\}\), which proves
the lemma.
\end{proof}

\begin{lemma}\label{lem:FIK-radial-gauge}
On \(\overline{\operatorname{ran}D^*}\), the only nonnegative
eigenvalue of \(L_f\) is \(0\), and its eigenspace is
\(\operatorname{span}\{\nabla^2f\}\).
\end{lemma}

\begin{proof}
Let \(\Delta_f^{(0)}\) be the radial scalar Friedrichs realization.
Its differential expression is
\[
 \Delta_f^{(0)}u
 =\frac1{r^3e^{-f}}
   \left(\frac F4r^3e^{-f}u'\right)'.
\tag{R7}\label{eq:FIK-radial-laplacian}
\]
It has compact resolvent.  Its quadratic form is nonpositive, so its
largest eigenvalue is the simple eigenvalue \(0\), with eigenfunction
\(1\).  Moreover,
\[
 \Delta_f^{(0)}(r^2-2)=-(r^2-2).
\tag{R8}\label{eq:FIK-radial-second}
\]
Part~(2) of Proposition~\ref{prop:FIK-radial-domain-package} supplies
the regular Friedrichs bolt condition, limit-point infinity,
compactness, simplicity, and nodal ordering.  Since \(r^2-2\) has
exactly one simple zero on \((1,\infty)\), the radial scalar spectrum can be
written
\[
 0=-\mu_0>-\mu_1=-1>-\mu_2>\cdots,
 \qquad \mu_0=0,\quad\mu_1=1,\quad\mu_k>1\ (k\ge2).
\tag{R9}\label{eq:FIK-radial-scalar-spectrum}
\]

Let \(\Delta_f^{(1)}\) denote the weighted rough Laplacian on the
radial \(e^0\)-one-form channel.  The exact-gradient commutator is
\[
 \Delta_f^{(1)}(du)=d\Delta_f^{(0)}u+\frac12\,du,
 \qquad
 \left(\Delta_f^{(1)}+\frac12\right)du
 =d(\Delta_f^{(0)}u+u).
\tag{R10}\label{eq:FIK-radial-gradient-commutator}
\]
The nonzero gradients of a complete radial scalar eigenbasis are
complete in the \(e^0\)-one-form channel by Part~(4) of
Proposition~\ref{prop:FIK-radial-domain-package}.  In particular, the
scalar gap \(\mu_1=1\) gives the correctly normalized weighted
Poincar\'e estimate
\[
 \left\|
 u-\frac{\langle u,1\rangle_{L^2_f}}
          {\|1\|_{L^2_f}^{\,2}}\,1
 \right\|_{L^2_f}
 \leq\|du\|_{L^2_f},
 \qquad
 \|1\|_{L^2_f}^{\,2}=(4\pi)^2.
\]
Consequently,
\(\Delta_f^{(1)}+\frac12\) has eigenvalues
\[
 1-\mu_k,\qquad k\ge1,
\tag{R11}\label{eq:FIK-radial-oneform-spectrum}
\]
with eigenforms \(du_k\).

Now suppose \(0\ne h\in\overline{\operatorname{ran}D^*}\) and
\(L_fh=\lambda h\).  Then \(Dh\ne0\), since otherwise
\(h\in\ker D\cap(\ker D)^\perp\).  The closed-domain intertwining
\eqref{eq:FIK-radial-closed-intertwining} shows that \(Dh\) is a
\(\lambda\)-eigenform of \(\Delta_f^{(1)}+\frac12\).  If
\(\lambda\ge0\), \eqref{eq:FIK-radial-scalar-spectrum} and
\eqref{eq:FIK-radial-oneform-spectrum} force
\(\lambda=0\) and
\[
 Dh\in\operatorname{span}\{d(r^2-2)\}.
\]
On the other hand,
\[
 L_f\nabla^2(r^2-2)=0,\qquad
 D\nabla^2(r^2-2)=-\frac12d(r^2-2),
\]
and \(\nabla^2(r^2-2)=2^{-1/2}\nabla^2f\) lies in
\(\operatorname{ran}D^*\cap\mathcal H_{01}\).  Subtracting a suitable
multiple of this Hessian from \(h\) produces an element of both
\(\ker D\) and \((\ker D)^\perp\), hence zero.  This proves the
lemma.
\end{proof}

\begin{proposition}\label{prop:FIK-J0-inventory}
The nonnegative spectrum of \(L_f\) in the \(J=0\) Wigner sector is
exactly
\[
 1\quad\text{on }\operatorname{span}\{\operatorname{Ric}\},
 \qquad
 0\quad\text{on }\operatorname{span}\{\nabla^2f\}.
\]
Every other radial eigenvalue is strictly negative.
\end{proposition}

\begin{proof}
The scalar radial blocks are strictly negative by the exact signs at the
start of the subsection.  On the sole two-component block, combine the
reducing decomposition \eqref{eq:FIK-radial-york} with
Lemmas~\ref{lem:FIK-radial-divfree} and
\ref{lem:FIK-radial-gauge}.
\end{proof}

\subsection{An exact one-dimensional count for the exceptional blocks}

Put
\[
 c_0=\sqrt2-1,\qquad
 F(r)=\frac1{\sqrt2}-\frac{c_0}{r^2}
       -\frac{c_0}{\sqrt2\,r^4},\qquad
 s(r)^2=r^2F(r),\qquad
 w(r)=r^3e^{-\sqrt2 r^2}.
\]
On radial functions consider the self-adjoint Sturm--Liouville operator
\[
 \mathscr S
 :=\Delta_f+\frac{3}{2r^2},\qquad
 \Delta_f u=\frac1w\left(\frac{Fw}{4}u'\right)' ,
\]
whose closed quadratic form is
\[
 \mathfrak s[u]
 =\int_1^\infty
   \left(\frac{3}{2r^2}|u|^2-\frac F4|u'|^2\right)w\,dr.
\]
The Gaussian first-moment estimate gives a uniform tail bound on bounded
sets in the form domain, and local Rellich compactness handles compact
radial intervals.  Thus the form embedding is compact.  (Equivalently,
conjugation by \(w^{1/2}\), followed by the intrinsic radial change of
variable, produces a one-dimensional Schr\"odinger operator whose
potential tends quadratically to \(+\infty\); the regular Friedrichs
condition is imposed at the bolt.)  Since \(3/(2r^2)\) is
bounded, \(\mathscr S\) consequently has a discrete sequence of
eigenvalues, written in decreasing order as
\[
 \mu_1>\mu_2\geq\mu_3\geq\cdots\longrightarrow-\infty .
\]
The first eigenvalue is simple and has a strictly positive real
eigenfunction \(a_1\).

\begin{lemma}\label{lem:FIK-scalar-second-negative}
The second eigenvalue of \(\mathscr S\) is strictly negative:
\(\mu_2<0\).
\end{lemma}

\begin{proof}
Let \(a_2\) be a real \(\mu_2\)-eigenfunction orthogonal to \(a_1\).
Since \(a_1>0\), the function \(a_2\) changes sign and therefore has a zero
\(r_0\in(1,\infty)\).  Set \(r_*=\sqrt{3/2}\) and introduce
\[
 \psi_-(r)=\frac32-r^2,\qquad
 \psi_+(r)=1-\frac{3}{2r^2}.
\]
Direct calculation gives
\[
 \mathscr S\psi_-
 =r^2-\frac72+\frac{9}{4r^2}<0
 \quad(1<r<r_*),
\]
and
\[
 \mathscr S\psi_+
 =\frac{-3(c_0+\sqrt2)r^4+12c_0r^2+6\sqrt2c_0}{4r^8}<0
 \quad(r>r_*).
\]
For the first inequality, after multiplication by \(4r^2\), the
polynomial \(4r^4-14r^2+9\) is decreasing on
\(1\leq r^2\leq3/2\) and is already negative at \(r=1\).
For the second, the numerator is a decreasing quadratic in \(r^2\) on
\([3/2,\infty)\), and its value at \(r^2=3/2\) is
\(\frac34(1-2\sqrt2)<0\).

We use the exact ground-state identity
\[
 \int_I\left(\frac{3}{2r^2}|u|^2-\frac F4|u'|^2\right)w\,dr
 =
 \int_I\frac{\mathscr S\psi}{\psi}|u|^2w\,dr
 -\int_I\frac F4\psi^2
       \left|\left(\frac u\psi\right)'\right|^2w\,dr ,
 \tag{*}\label{eq:FIK-exceptional-ground-state}
\]
valid when \(\psi>0\) on \(I\) and the endpoint flux
\(\frac{Fw}{4}\frac{\psi'}{\psi}|u|^2\) vanishes.  It follows by writing
\(u=\psi v\) and integrating one exact derivative.  The identity extends
to the relevant form-domain restrictions by cutoff.

If \(r_0\leq r_*\), apply \eqref{eq:FIK-exceptional-ground-state} to the
nonzero restriction of \(a_2\) to \(I=(1,r_0)\), with
\(\psi=\psi_-\).  The flux vanishes at the bolt because \(F(1)=0\), and
at \(r_0\) because \(a_2(r_0)=0\).  If \(r_0\geq r_*\), apply the same
identity on \(I=(r_0,\infty)\), with \(\psi=\psi_+\); Gaussian form
cutoffs make the endpoint term at infinity vanish by Part~(2) of
Proposition~\ref{prop:FIK-radial-domain-package}.  When \(r_0=r_*\),
use either adjacent interval.
A nontrivial Sturm--Liouville eigenfunction has a simple zero there:
if \(a_2(r_*)=a_2'(r_*)=0\), ODE uniqueness would force
\(a_2\equiv0\).  The selected supersolution also has a simple zero.
Consequently \(a_2/\psi\) has a finite limit and the endpoint flux
vanishes, so the same limiting identity holds.  In either case the eigenvalue
equation, tested against the restriction, and the strict supersolution
inequality give
\[
 \mu_2\int_I|a_2|^2w\,dr
 =\mathfrak s_I[a_2]<0.
\]
Thus \(\mu_2<0\).
\end{proof}

\begin{proposition}\label{prop:FIK-exceptional-comparison}
Let \(\mathcal P,\widetilde{\mathcal P}\) be the two exceptional
\(J=1\) Hermitian \(2\times2\) potential blocks and let \(\mathcal Q\)
be the exceptional \(J=2\) Hermitian \(4\times4\) potential block
obtained from the exact Wigner reduction.  For \(r>1\),
\[
 \mathcal P,\widetilde{\mathcal P}
 \leq \frac{3}{2r^2}\operatorname{diag}(1,0),
 \qquad
 \mathcal Q
 \leq \frac{3}{2r^2}\operatorname{diag}(0,1,0,0).
\]
The inequality is strict on a nonzero vector whose distinguished
component is zero.
\end{proposition}

\begin{proof}
The two \(J=1\) blocks are unitarily conjugate, so it suffices to treat
\(\mathcal P=(\mathcal P_{ij})\).  The exact entries give
\[
 r^2\left(\mathcal P_{11}+\frac{4s^2}{r^4}\right)
 =\frac{\sqrt2c_0}{r^4}
  +\frac{c_0^2}{2(r^2+c_0)}
  +\frac1{2\sqrt2c_0}
 \leq\frac32 ,
\]
because the displayed function is decreasing in \(r^2\) and equals
\(3/2\) at \(r=1\).  Also
\[
 r^2\left(\mathcal P_{22}+\frac1{2r^2}\right)
 =
 \frac{2(c_0-1)-4(c_0-1)r^2+(2c_0-3)r^4}
      {2\sqrt2(r^2-1)(r^2+c_0)}<0.
\]
Indeed, on writing \(y=r^2-1\), the numerator becomes
\[
 (2\sqrt2-5)y^2-2y-1<0.
\]
Since the off-diagonal entry has modulus
\(\sqrt2s/r^3\), Young's inequality yields
\[
 2\frac{\sqrt2s}{r^3}|z_1||z_2|
 \leq\frac{4s^2}{r^4}|z_1|^2
     +\frac1{2r^2}|z_2|^2,
\]
which proves the first comparison and its strictness on \(z_1=0\).

For \(\mathcal Q\), interchange its first two coordinates and multiply
by \(r^2\).  The result has the block form
\[
 \begin{bmatrix}a&b\\ b^*&C\end{bmatrix},
\qquad
 C=\operatorname{diag}\left(-\frac{2+\sqrt2}{2},d,d\right),
\]
where
\[
 d=\frac{1+2\sqrt2r^2-2(1+\sqrt2)r^4}
        {2\sqrt2(r^2-1)(r^2+c_0)}<0.
\]
Indeed, the numerator, viewed as a quadratic in \(x=r^2\), equals
\(-1\) at \(x=1\) and is strictly decreasing for \(x\geq1\).
Thus \(C<0\).  With \(x=r^2\), its Schur complement is
\[
 \alpha(x):=a-bC^{-1}b^*
 =\frac{N(x)}{D(x)},
\]
where
\[
\begin{split}
 N(x)&=2\bigl(10-6\sqrt2+5(-4+3\sqrt2)x
  +(17-13\sqrt2)x^2\\
 &\hspace{35mm}{}+4(-2+\sqrt2)x^3+2x^4\bigr),\\
 D(x)&=(2+\sqrt2)x^2
       \bigl(2(1+\sqrt2)x^2-2\sqrt2x-1\bigr)>0.
\end{split}
\]
For \(y=x-1\geq0\), direct expansion gives the exact positive-coefficient
certificate
\[
\begin{split}
 \frac32D(x)-N(x)
 ={}&(8+9\sqrt2)y^4+(42+22\sqrt2)y^3\\
 &+\left(41+\frac{73}{2}\sqrt2\right)y^2
 +(28+13\sqrt2)y+1+\frac32\sqrt2>0.
\end{split}
\]
Consequently \(\alpha<3/2\), and the Schur-complement criterion proves
the \(J=2\) comparison.  Since \(C<0\), it is strict on a nonzero
vector whose distinguished coordinate vanishes.
\end{proof}

\begin{corollary}\label{cor:FIK-exceptional-count}
Each of the three exceptional vector Sturm--Liouville operators has at
most one nonnegative eigenvalue, counted with multiplicity.
\end{corollary}

\begin{proof}
Let \(\pi\) denote the distinguished scalar component in
Proposition~\ref{prop:FIK-exceptional-comparison}.  Its quadratic
form satisfies
\[
 \mathfrak q_{\mathcal B}[\xi]\leq\mathfrak s[\pi\xi],
\]
with strict inequality for nonzero \(\xi\) if \(\pi\xi=0\).
If the vector operator had a two-dimensional nonnegative spectral
subspace \(W\), choose \(0\neq\xi\in W\) so that
\(\pi\xi\perp a_1\).  If \(\pi\xi\neq0\), the min--max principle and
Lemma~\ref{lem:FIK-scalar-second-negative} give
\[
 0\leq\mathfrak q_{\mathcal B}[\xi]
 \leq\mathfrak s[\pi\xi]
 \leq\mu_2\|\pi\xi\|_{L^2(w\,dr)}^2<0,
\]
a contradiction.  If \(\pi\xi=0\), the strict part of the matrix
comparison gives the same contradiction.  Hence the nonnegative
spectral multiplicity is at most one.
\end{proof}

The component identities
\eqref{eq:self-exceptional-mode-vectors} place
\(S^1_{1,M'}\) and \(S^1_{-1,M'}\) in the two exceptional \(J=1\)
blocks and place \(S^2_{0,M'}\) in the exceptional \(J=2\) block.
Lemma~\ref{lem:self-contained-bolt-reality} and
Corollary~\ref{cor:self-contained-Hessian-families} show directly that
these tensors are smooth, nonzero, and \(L^2_f\), with eigenvalues
\(1-1/\sqrt2>0\), \(1-1/\sqrt2>0\), and \(0\), respectively.
Corollary~\ref{cor:FIK-exceptional-count} therefore shows that these
modes exhaust the nonnegative spectrum of the exceptional blocks,
completing the exceptional-block count.

\section{The diagonal scale action and physical reconstruction}
\label{app:diagonal-scale-action}

We verify the scale signs independently of the harmonic-map calculation.
Let
\[
 \mathcal R_{\bar g}(g)
 =-2\Ric_g+g+\Lie_{B_{\bar g}(g)}g
  -\Lie_{\bar\nabla\bar f}g.
\]
Then $\mathcal R_{\bar g}(\bar g)=0$, and its linearization at
$\bar g$ is $\A$.  The controlled normalized equation is
\begin{equation}\label{eq:appendix-controlled}
 \partial_\tau g
 =\mathcal R_{\bar g}(g)
  +a(g-\Lie_{\bar\nabla\bar f}g)
  +\sum_{j=1}^8b_j\Lie_{\chi_\tau W_j}g.
\end{equation}
Equivalently,
\[
 \partial_\tau g
 =-2\Ric_g+(1+a)g+\Lie_{\mathscr V}g,
\]
where
\[
 \mathscr V=B_{\bar g}(g)-(1+a)\bar\nabla\bar f
   +\sum_{j=1}^8b_j\chi_\tau W_j.
\]
Because \(\chi_\tau W_j\) is the instantaneous velocity,
differentiating the pullback produces no
\(\partial_\tau\chi_\tau\) term.
Let
\[
 \lambda_\tau=-(1+a)\lambda,\qquad
 t_\tau=\lambda,
\]
and let $\Theta_\tau$ be generated by $-\mathscr V$, with the
initial normalization $\Theta_{\tau_0}=\operatorname{Id}$.  This is the
reconstruction diffeomorphism; in the pullback convention of
Remark~\ref{rem:scale} it is
$\Theta_\tau=\Upsilon_{t(\tau)}^{-1}$ (after the same initial
normalization).  Then
\[
 G(t(\tau))=\lambda(\tau)\Theta_\tau^*g(\tau)
\]
satisfies $\partial_tG=-2\Ric_G$.  Indeed,
\[
 \partial_\tau(\Theta_\tau^*g)
 =\Theta_\tau^*(-2\Ric_g+(1+a)g),
\]
and the zeroth-order term cancels after multiplication by $\lambda$.

The finite background action corresponding to the same scale variable
is \eqref{eq:finite-scale-action}.  Its derivative is
$2\Ric_{\bar g}$, and the scale law is
\[
 \lambda(\tau)
 =\lambda(\tau_0)e^{-(\tau-\tau_0)}
  \exp\left(-\int_{\tau_0}^{\tau}a(s)\,ds\right).
\]
Thus positivity is automatic.  If $a\in L^1(d\tau)$, the remaining
physical time is asymptotic to $\lambda$ by
Lemma~\ref{lem:scale-closure}.

\section{A Gaussian tail in explicit FIK coordinates}
\label{app:FIK-Gaussian-tail}

Let $f_{\mathrm{NO}}$ denote the mass-normalized potential in
Lemma~\ref{lem:self-contained-FIK-ledger}.  In the radial coordinates
fixed there,
\[
 f_{\mathrm{NO}}(r)
 =\sqrt2(r^2-1)-\log\bigl(2(\sqrt2-1)\bigr),
\qquad
 e^{-f_{\mathrm{NO}}}dV_{\bar g}
 =16r^3e^{-f_{\mathrm{NO}}}\,dr\,dV_{S^3}.
\]
Our potential $\bar f$ differs from $f_{\mathrm{NO}}$ by an additive
constant, which only rescales the weighted measure.
We first verify the general \(\bar f\)-tail used in
\eqref{eq:general-Gaussian-tail}.  Since
\(\bar f(r)=\sqrt2\,r^2+C_0\), the substitution
\(u=\bar f(r)\) gives \(r^3\,dr\leq C(1+u)\,du\).  Hence, for every
\(R_f\geq0\),
\[
 \begin{split}
 \int_{\{\bar f\geq R_f\}}(1+\bar f)^N\,d\nu
 &\leq C_N\int_{R_f}^\infty(1+u)^{N+1}e^{-u}\,du\\
 &\leq C_N(1+R_f)^{N+1}e^{-R_f}\\
 &\leq C_N(1+R_f)^{N+2}e^{-R_f/2}.
 \end{split}
\]
The constants absorb the compact range where the radial coordinate
does not parametrize the bolt by a single product chart.  This proves
\eqref{eq:general-Gaussian-tail} directly, with room in both the
polynomial power and the exponential.

If $q$ is supported in $\{r\geq R_r\}$ and
\[
 \sum_{j=0}^m\abs{\bar\nabla^jq}\leq A(1+r)^N,
\]
then
\[
 \norm q_{H^m_\nu}^2
 \leq C A^2\int_{R_r}^\infty r^{2N+3}e^{-\sqrt2r^2}\,dr
 \leq C_NA^2(1+R_r)^{2N+2}e^{-\sqrt2R_r^2}.
\]
Taking $R_r^2\simeq e^\tau$ recovers the
$e^{-c e^\tau}$ estimates used above.


\begin{thebibliography}{99}

\bibitem{AngenentKnopf}
S.~Angenent and D.~Knopf,
\emph{An example of neckpinching for Ricci flow on $S^{n+1}$},
Math. Res. Lett. \textbf{11} (2004), no.~4, 493--518,
\doi{10.4310/MRL.2004.v11.n4.a8}.

\bibitem{Appleton}
A.~Appleton,
\emph{Eguchi--Hanson singularities in \(U(2)\)-invariant Ricci flow},
Peking Math. J. \textbf{6} (2023), no.~1, 1--141,
\doi{10.1007/s42543-022-00048-y}.

\bibitem{BahuaudGuentherIsenberg}
E.~Bahuaud, C.~Guenther, and J.~Isenberg,
\emph{Convergence stability for Ricci flow},
J. Geom. Anal. \textbf{30} (2020), no.~1, 310--336,
\doi{10.1007/s12220-018-00132-9}.

\bibitem{BamlerKleiner}
R.~H.~Bamler and B.~Kleiner,
\emph{Uniqueness and stability of Ricci flow through singularities},
Acta Math. \textbf{228} (2022), no.~1, 1--215,
\doi{10.4310/ACTA.2022.v228.n1.a1}.

\bibitem{BravermanMilatovicShubin}
M.~Braverman, O.~Milatovic, and M.~A.~Shubin,
\emph{Essential self-adjointness of Schr\"odinger-type operators on
manifolds},
Russian Math. Surveys \textbf{57} (2002), no.~4, 641--692,
\doi{10.1070/RM2002v057n04ABEH000532}.

\bibitem{CaoHamiltonIlmanen}
H.-D.~Cao, R.~S.~Hamilton, and T.~Ilmanen,
\emph{Gaussian densities and stability for some Ricci solitons},
preprint (2004), \arxiv{math/0404165}.

\bibitem{ChoiLai}
K.~Choi and Y.~Lai,
\emph{Convergence of Ricci flow and long-time existence of harmonic
map heat flow},
preprint (2025), \arxiv{2504.02804}.

\bibitem{CifarelliConlonDeruelle}
C.~Cifarelli, R.~J.~Conlon, and A.~Deruelle,
\emph{On finite time Type I singularities of the K\"ahler--Ricci flow
on compact K\"ahler surfaces},
J. Eur. Math. Soc. \textbf{28} (2026), no.~2, 463--504,
\doi{10.4171/JEMS/1485}.

\bibitem{ColdingMinicozzi}
T.~H.~Colding and W.~P.~Minicozzi,
\emph{Singularities of Ricci flow and diffeomorphisms},
Publ. Math. Inst. Hautes \'Etudes Sci. \textbf{142} (2025), 75--152,
\doi{10.1007/s10240-025-00157-1}.

\bibitem{ConlonHallgrenMa}
R.~J.~Conlon, M.~Hallgren, and Z.~Ma,
\emph{Non-collapsed finite time singularities of the Ricci flow on
compact K\"ahler surfaces are of Type I},
preprint (2025), \arxiv{2502.19804}.

\bibitem{DeTurck}
D.~M.~DeTurck,
\emph{Deforming metrics in the direction of their Ricci tensors},
J. Differential Geom. \textbf{18} (1983), no.~1, 157--162,
\doi{10.4310/jdg/1214509286}.

\bibitem{EndersMullerTopping}
J.~Enders, R.~M\"uller, and P.~M.~Topping,
\emph{On Type-I singularities in Ricci flow},
Comm. Anal. Geom. \textbf{19} (2011), no.~5, 905--922,
\doi{10.4310/CAG.2011.v19.n5.a4}.

\bibitem{FIK}
M.~Feldman, T.~Ilmanen, and D.~Knopf,
\emph{Rotationally symmetric shrinking and expanding gradient
K\"ahler--Ricci solitons},
J. Differential Geom. \textbf{65} (2003), no.~2, 169--209,
\doi{10.4310/jdg/1090511686}.

\bibitem{GarfinkleIsenbergKnopfWu}
D.~Garfinkle, J.~Isenberg, D.~Knopf, and H.~Wu,
\emph{Asymptotic behaviour of unstable perturbations of the
Fubini--Study metric in Ricci flow},
Nonlinearity \textbf{38} (2025), no.~7, 075015,
\doi{10.1088/1361-6544/addf09}.

\bibitem{GuentherIsenbergKnopf}
C.~Guenther, J.~Isenberg, and D.~Knopf,
\emph{Stability of the Ricci flow at Ricci-flat metrics},
Comm. Anal. Geom. \textbf{10} (2002), no.~4, 741--777,
\doi{10.4310/CAG.2002.v10.n4.a4}.

\bibitem{GuoSong}
B.~Guo and J.~Song,
\emph{On Feldman--Ilmanen--Knopf conjecture for the blow-up behavior
of the K\"ahler Ricci flow},
Math. Res. Lett. \textbf{23} (2016), no.~6, 1681--1719,
\doi{10.4310/MRL.2016.v23.n6.a6}.

\bibitem{Hamilton}
R.~S.~Hamilton,
\emph{Three-manifolds with positive Ricci curvature},
J. Differential Geom. \textbf{17} (1982), no.~2, 255--306,
\doi{10.4310/jdg/1214436922}.

\bibitem{Hughes}
J.~Hughes,
\emph{$L^2$-instability of the Taub--Bolt metric under the Ricci flow},
preprint (2024), \arxiv{2408.15269}.

\bibitem{IKSNeckpinch}
J.~Isenberg, D.~Knopf, and N.~\v{S}e\v{s}um,
\emph{Ricci flow neckpinches without rotational symmetry},
Comm. Partial Differential Equations \textbf{41} (2016), no.~12,
1860--1894,
\doi{10.1080/03605302.2016.1233982}.

\bibitem{IKS}
J.~Isenberg, D.~Knopf, and N.~\v{S}e\v{s}um,
\emph{Non-K\"ahler Ricci flow singularities modeled on
K\"ahler--Ricci solitons},
Pure Appl. Math. Q. \textbf{15} (2019), no.~2, 749--784,
\doi{10.4310/PAMQ.2019.v15.n2.a5}.

\bibitem{Kroencke}
K.~Kr\"oncke,
\emph{Stability of Einstein metrics under Ricci flow},
Comm. Anal. Geom. \textbf{28} (2020), no.~2, 351--394,
\doi{10.4310/CAG.2020.v28.n2.a5}.

\bibitem{LiWang}
Y.~Li and B.~Wang,
\emph{On K\"ahler Ricci shrinker surfaces},
Acta Math. \textbf{236} (2026), no.~1, 1--50,
\doi{10.4310/ACTA.2026.v236.n1.a1}.

\bibitem{Lunardi}
A.~Lunardi,
\emph{Analytic semigroups and optimal regularity in parabolic
problems},
Progress in Nonlinear Differential Equations and their Applications,
vol.~16, Birkh\"auser, Basel, 1995.

\bibitem{Maximo}
D.~M\'aximo,
\emph{On the blow-up of four-dimensional Ricci flow singularities},
J. Reine Angew. Math. \textbf{692} (2014), 153--171,
\doi{10.1515/crelle-2012-0080}.

\bibitem{NaffOzuch}
K.~Naff and T.~Ozuch,
\emph{Linear stability of the blowdown Ricci shrinker in 4D},
preprint (2025), \arxiv{2509.05470}.

\bibitem{Perelman}
G.~Perelman,
\emph{The entropy formula for the Ricci flow and its geometric
applications},
preprint (2002), \arxiv{math/0211159}.

\bibitem{PengLu}
P.~Lu,
\emph{A local curvature bound in Ricci flow},
Geom. Topol. \textbf{14} (2010), no.~2, 1095--1110,
\doi{10.2140/gt.2010.14.1095}.

\bibitem{Shi}
W.-X.~Shi,
\emph{Deforming the metric on complete Riemannian manifolds},
J. Differential Geom. \textbf{30} (1989), no.~1, 223--301,
\doi{10.4310/JDG/1214443292}.

\bibitem{SongTypeI}
J.~Song,
\emph{Some Type I solutions of Ricci flow with rotational symmetry},
Int. Math. Res. Not. IMRN (2015), no.~16, 7365--7381,
\doi{10.1093/imrn/rnu134}.

\bibitem{SongWeinkove}
J.~Song and B.~Weinkove,
\emph{The K\"ahler--Ricci flow on Hirzebruch surfaces},
J. Reine Angew. Math. \textbf{659} (2011), 141--168,
\doi{10.1515/crelle.2011.071}.

\bibitem{SongWeinkoveContraction}
J.~Song and B.~Weinkove,
\emph{Contracting exceptional divisors by the K\"ahler--Ricci flow},
Duke Math. J. \textbf{162} (2013), no.~2, 367--415,
\doi{10.1215/00127094-1962881}.

\bibitem{Stolarski}
M.~Stolarski,
\emph{Closed Ricci flows with singularities modeled on
asymptotically conical shrinkers},
Geom. Topol. \textbf{30} (2026), no.~4, 1277--1370,
\doi{10.2140/gt.2026.30.1277}.

\bibitem{Zettl}
A.~Zettl,
\emph{Sturm--Liouville theory},
Mathematical Surveys and Monographs, vol.~121,
American Mathematical Society, Providence, RI, 2005,
\doi{10.1090/surv/121}.

\enlargethispage{4\baselineskip}
\end{thebibliography}
\end{document}